\documentclass{memo-l}

\usepackage{amssymb,amsthm,enumitem,graphicx,verbatim,array,scalefnt}

\newtheorem{theorem}{Theorem}

\newtheorem{corollary}[theorem]{Corollary}
\newtheorem{thm}{Theorem}[chapter]
\newtheorem{prop}[thm]{Proposition}
\newtheorem{lem}[thm]{Lemma}
\newtheorem{cor}[thm]{Corollary}

\newtheorem{definition}[thm]{Definition}
\newcounter{fig}

\theoremstyle{remark}

\numberwithin{section}{chapter}
\numberwithin{equation}{chapter}
\numberwithin{table}{chapter}

\makeindex

\def\N{{\mathbb N}}
\def\Z{{\mathbb Z}}

\def\R{{\mathbb R}}

\def\F{{\mathbb F}}
\def\a{\mathbf{a}}
\def\b{\mathbf{b}}
\def\c{\mathbf{c}}
\def\d{\mathbf{d}}
\def\k{\mathbf{k}}

\def\codim{\mathop{\mathrm{codim}}}
\def\im{\mathop{\mathrm{im}}}
\def\Part{\mathop{\mathrm{Part}}}
\def\rank{\mathop{\mathrm{rank}}}

\def\tr{\mathop{\mathrm{tr}}}
\def\ts{\textstyle}
\def\ss{\scriptstyle}
\def\sss{\scriptscriptstyle}

\def\pmin{\phantom{-}}
\def\pa{\phantom{{}^*}}
\def\pay{\pa\hbox{yes}^*}
\def\pan{\pa\hbox{no}^*}
\def\pand{\phantom{{}^{\dagger*}}\hbox{no}^{\dagger*}}
\def\pss{{}_{\phantom{ss}}}
\def\pu{{}_{\phantom{u}}}
\def\ptw{\vbox{\hrule height0pt depth0pt width3pt}}
\def\pw{\vbox{\hrule height0pt depth0pt width1.6pt}}
\def\phmu{\pw \cdot \pw}

\def\vstrut{\vbox{\hrule height10pt depth5pt width0pt}}
\def\tbs{\vbox{\hrule height10pt depth3.75pt width0pt}}
\def\mtbs{\vbox{\hrule height9.4pt depth3.3pt width0pt}}
\def\mmtbs{\vbox{\hrule height8.8pt depth2.85pt width0pt}}
\def\tabcapsp{\vspace{-5mm}}
\def\z{\zeta}
\def\e{\epsilon}
\def\ve{\varepsilon}
\def\L{{\mathfrak L}}
\def\S{{\mathfrak S}}

\def\P{\mathcal{P}}
\def\Slambda{{\mathcal S}(\lambda)}
\def\sscon{(*)_{ss}}
\def\ssevcon{(*')_{ss}}
\def\ssdiamcon{(\diamond)_{ss}}
\def\ssdiamevcon{(\diamond')_{ss}}
\def\ssdagcon{(\dag)_{ss}}
\def\ssddagcon{(\ddag)_{ss}}
\def\ucon{(*)_u}
\def\udiamcon{(\diamond)_u}
\def\udagcon{(\dag)_u}
\def\uddagcon{(\ddag)_u}
\def\GL{\mathop{\mathrm{GL}}}

\def\GU{\mathop{\mathrm{GU}}}
\def\SL{\mathop{\mathrm{SL}}}
\def\SO{\mathop{\mathrm{SO}}}
\def\SU{\mathop{\mathrm{SU}}}
\def\PGL{\mathop{\mathrm{PGL}}}

\def\PGU{\mathop{\mathrm{PGU}}}
\def\PSL{\mathop{\mathrm{PSL}}}

\def\PSU{\mathop{\mathrm{PSU}}}
\def\Tran{\mathop{\mathrm{Tran}}}
\def\Ann{\mathop{\mathrm{Ann}}}
\def\Ad{\mathop{\mathrm{Ad}}}
\def\adj{\mathop{\mathrm{adj}}}
\def\bom{\bar\omega}
\def\diag{\mathrm{diag}}
\def\height{\mathrm{ht}}
\def\Gk{{\mathcal G}_k}
\newcommand\G[1]{{\mathcal G}_{#1}}
\def\yy{\hbox{yes} & \hbox{yes}}
\def\nn{\hbox{no} & \hbox{no}}
\def\ny{\hbox{no} & \hbox{yes}}

\newcommand\aaf[2]{\frac{a_#1}{a_#2}}
\newcommand\attopp[2]{\genfrac{}{}{0pt}{}{#1}{#2}}
\newcommand\dfourrt[4]{\textstyle{
\hbox to4pt{$\hfil\scriptstyle{#1}\hfil$}
\hbox to4pt{$\hfil\scriptstyle{#2}\hfil$}\!\vrule width0.6pt height0pt depth0pt\attopp{
{\raise-1pt\hbox to4pt{$\hfil\scriptstyle{#3}\hfil$}}}
{\raise1pt\hbox to4pt{$\hfil\scriptstyle{#4}\hfil$}}}}

\newcommand\ffourrt[4]{\textstyle{
\hbox to4pt{$\hfil\scriptstyle{#1}\hfil$}
\hbox to4pt{$\hfil\scriptstyle{#2}\hfil$}
\hbox to4pt{$\hfil\scriptstyle{#3}\hfil$}
\hbox to4pt{$\hfil\scriptstyle{#4}\hfil$}}}
\newcommand\esixrt[6]{\textstyle{\attopp{
\hbox to4pt{$\hfil\scriptstyle{#1}\hfil$}
\hbox to4pt{$\hfil\scriptstyle{#3}\hfil$}
\hbox to4pt{$\hfil\scriptstyle{#4}\hfil$}
\hbox to4pt{$\hfil\scriptstyle{#5}\hfil$}
\hbox to4pt{$\hfil\scriptstyle{#6}\hfil$}}
{\raise2pt\hbox to4pt{$\hfil\scriptstyle{#2}\hfil$}}}}
\newcommand\esevenrt[7]{\textstyle{\attopp{
\hbox to4pt{$\hfil\scriptstyle{#1}\hfil$}
\hbox to4pt{$\hfil\scriptstyle{#3}\hfil$}
\hbox to4pt{$\hfil\scriptstyle{#4}\hfil$}
\hbox to4pt{$\hfil\scriptstyle{#5}\hfil$}
\hbox to4pt{$\hfil\scriptstyle{#6}\hfil$}
\hbox to4pt{$\hfil\scriptstyle{#7}\hfil$}}
{\raise2pt\hbox to4pt{$\hfil\scriptstyle{#2}\hfil$}
\raise2pt\hbox to4pt{$\hfil\scriptstyle{\phantom{0}}\hfil$}}}}
\newcommand\eeightrt[8]{\textstyle{\attopp{
\hbox to4pt{$\hfil\scriptstyle{#1}\hfil$}
\hbox to4pt{$\hfil\scriptstyle{#3}\hfil$}
\hbox to4pt{$\hfil\scriptstyle{#4}\hfil$}
\hbox to4pt{$\hfil\scriptstyle{#5}\hfil$}
\hbox to4pt{$\hfil\scriptstyle{#6}\hfil$}
\hbox to4pt{$\hfil\scriptstyle{#7}\hfil$}
\hbox to4pt{$\hfil\scriptstyle{#8}\hfil$}}
{\raise2pt\hbox to4pt{$\hfil\scriptstyle{#2}\hfil$}
\raise2pt\hbox to4pt{$\hfil\scriptstyle{\phantom{0}}\hfil$}
\raise2pt\hbox to4pt{$\hfil\scriptstyle{\phantom{0}}\hfil$}}}}

\begin{document}

\frontmatter

\title{Generic stabilizers in actions of simple algebraic groups}


\author{R.M. Guralnick}
\address{Department of Mathematics \\
University of Southern California \\
Los Angeles \\
CA 90089-2532, USA}
\email{guralnic@usc.edu}
\thanks{The first author was partially supported by NSF grants DMS-1600056 and DMS-1901595, and a Simons Foundation Fellowship 609771.}

\author{R. Lawther}
\address{Department of Pure Mathematics and Mathematical Statistics \\
Centre for Mathematical Sciences \\
Cambridge University \\
Wilberforce Road \\
Cambridge CB3 0WB, UK}
\email{ril10@cam.ac.uk}

\date{}




\maketitle

\tableofcontents

\begin{abstract}
In this paper we treat faithful actions of simple algebraic groups on irreducible modules and on the associated Grassmannian varieties. By explicit calculation, we show that in each case, with essentially one exception, there is a dense open subset any point of which has stabilizer conjugate to a fixed subgroup, called the {\em generic stabilizer\/}. We provide tables listing generic stabilizers in the cases where they are non-trivial; in addition we decide whether or not there is a dense orbit, or a regular orbit for the action on the module.
\end{abstract}


\mainmatter
\chapter{Introduction}\label{chap: intro}

In this chapter we state our main results, establish notation which will be used throughout, recall some basic material and prove various preliminary lemmas.

\section{Statement of main results}\label{sect: statement}

Let $G$ be a simple algebraic group over an algebraically closed field $K$ of characteristic $p$; for convenience we shall take $p = \infty$ if $K$ has characteristic zero. Let $V$ be a non-trivial irreducible $G$-module of dimension $d$. Recall that for $k = 1, \dots, d$ the Grassmannian variety $\Gk(V)$ consists of the $k$-dimensional subspaces of $V$, and has dimension $k(d - k)$; as the action of $G$ on $V$ is linear, it extends naturally to $\Gk(V)$. In this paper we treat the actions of $G$ on $V$ and on the Grassmannian varieties $\Gk(V)$ for $1 \leq k \leq \frac{d}{2}$ (the reason for the upper bound is that $\Gk(V)$ is naturally isomorphic to $\G{d - k}(V^*)$, where $V^*$ is the dual of $V$).

Let $X$ be an irreducible variety on which $G$ acts, and write $G_X$ for the kernel of the action of $G$ on $X$; by definition $G/G_X$ then acts faithfully on $X$. (Observe that if in fact $X$ is a Grassmannian variety $\Gk(V)$ with $V$ non-trivial, then the kernel $G_X$ is equal to the centre of $G$, since any central element acts on $V$ as a scalar and hence fixes any line in $V$.) Our concern is with the stabilizers of points in this faithful action; note that this means that it is harmless to assume $G$ is of simply connected type, and from time to time we may do so without further comment.

To begin with, if $x \in X$ has trivial stabilizer in $G/G_X$, we say that its orbit is {\em regular\/}. We then make the following definition.

\begin{definition}
If $\hat X$ is a non-empty open set in $X$ with the property that for all $x, x' \in \hat X$ the stabilizers in $G/G_X$ of $x$ and $x'$ are isomorphic subgroups, we say that the action has a {\em semi-generic stabilizer\/}, whose isomorphism type is that shared by each such subgroup $C_{G/G_X}(x)$ for $x \in \hat X$. If moreover $\hat X$ has the property that for all $x, x' \in \hat X$ the stabilizers in $G/G_X$ of $x$ and $x'$ are conjugate subgroups, we say that the action has a {\em generic stabilizer\/}, whose conjugacy class is that containing each such subgroup $C_{G/G_X}(x)$ for $x \in \hat X$.
\end{definition}

Generic stabilizers for actions on modules have been studied extensively in the case where the field has characteristic zero. There it follows from a result of Richardson in \cite{Rich} that, when a reductive group acts on a smooth affine variety, a generic stabilizer always exists. For $G$ a complex simple Lie group acting on an irreducible module $V$, a complete list of non-trivial generic stabilizers is given in a summary table in \cite{PoVi}, which also lists references to some of the original papers. The proof of the completeness of this table makes use of Richardson's result, as well as some character theory for the Lie algebra to determine if the generic stabilizer in the group has positive dimension (which is equivalent to the generic stabilizer in the Lie algebra being non-zero). It was proved in \cite{AVE} that the generic stabilizer fails to have positive dimension if and only if, for an arbitrary non-zero element $x$ of the Lie algebra of $G$,
$$
\frac{\tr_V(x^2)}{\tr(\Ad(x))^2} > 1
$$
provided the denominator is non-zero (Dynkin's famous paper \cite{Dynk} shows that the ratio, when it exists, is independent of $x$); it follows that if $\dim V > \dim G$ the generic stabilizer is finite and hence consists of semisimple elements. In \cite{Pop1} Popov sketches the proof, which proceeds by analysing weight strings, of the fact that the generic stabilizer is in fact trivial if $\dim V$ is large enough; he also gives a table listing the cases where the generic stabilizer is finite but non-trivial.

Generic stabilizers for actions on Grassmannian varieties in characteristic zero have received rather less attention. Results on $\Gk(V)$ can be interpreted in terms of linear actions of $G \times \GL_k(K)$ on $V \otimes V'$ where $V'$ is the natural module for $\GL_k(K)$ (see Lemma~\ref{lem: reduction to projective space}); this suggests considering groups which are semisimple rather than simple. For such groups, lists corresponding to the table in \cite{PoVi} are given in \cite{Ela} and \cite[Table~1]{Pop2} (treating respectively the cases where the generic stabilizer has positive dimension, and where it is finite but non-trivial). Here \cite{AVE} gives a condition analogous to that above, but it is only sufficient and not necessary. In \cite{Pop1} Popov also treats such groups, beginning by examining the case where the semisimple group is a direct product of linear groups acting on the tensor product of the various natural modules; once these instances have been classified, he goes through various possibilities to determine the list of cases where the generic stabilizer is finite but non-trivial. However, the implications of all this work for actions on Grassmannian varieties seem not to be mentioned.

The situation in positive characteristic presents considerably greater challenges. To begin with, there is no analogue of Richardson's result which we may use; in many cases, including all those where a simple group acts on an irreducible module, we shall in fact find that generic stabilizers do exist, but this is something deduced a posteriori rather than known a priori. For actions on modules, character theory is not really of use, since the characters of irreducible modules are not available in general. There are some results for simple Lie algebras acting irreducibly on restricted modules (see for example \cite{Auld, Gue}); but these are not complete, and the connection between the result for the Lie algebra and that for the group is less close than in characteristic zero. Indeed, there are examples where the generic stabilizer for the group is trivial but that for the Lie algebra has positive dimension (see \cite{GGgen} for a complete list of such examples). In addition, finite generic stabilizers need not consist only of semisimple elements, but may also contain unipotent elements. As for actions on Grassmannian varieties, these seem to have been studied only in \cite{GLMS}, which did not consider generic stabilizers but rather classified instances of the existence of a dense orbit, or of finitely many orbits.

In the present work we settle the question of the existence of generic stabilizers, and determine them where they exist, for actions of simple algebraic groups on both modules and Grassmannian varieties, in arbitrary characteristic. The results for characteristic zero provide independent confirmation of the information given in the table in \cite{PoVi}; however, whereas some of the arguments in works mentioned above were given in sketch form only, here full details are provided.

(It should also be mentioned that generic stabilizers have applications to invariant theory, Galois cohomology and essential dimension: see \cite{GGPi, Mer, Rei} for the theory and \cite{BRV, GGspin, Karp, Loe, LMMR} for specific applications. For some of the latter, it is necessary to know the generic stabilizer as a group scheme; in \cite{GGgen} Garibaldi and the first author use the results obtained in the present work to show that, for a simple algebraic group acting on a finite-dimensional irreducible module, the generic stabilizer exists as a group scheme, and to determine it in all such cases.)

We now move on to state our results. Our most basic one is the following.

\begin{theorem}\label{thm: existence of generic stab}
Let $G$ be a simple algebraic group over an algebraically closed field of characteristic $p$, and $V$ a non-trivial irreducible $G$-module of dimension $d$.
\begin{itemize}
\item[(i)] The action of $G$ on $V$ has a generic stabilizer.
\item[(ii)] For $1 \leq k \leq \frac{d}{2}$, either the action of $G$ on $\Gk(V)$ has a generic stabilizer, or $G = B_3$ or $C_3$, $p = 2$, $V$ is the spin module for $G$, and $k = 4$, in which case the action of $G$ on $\Gk(V)$ has a semi-generic stabilizer but not a generic stabilizer.
\end{itemize}
\end{theorem}

The proof of Theorem~\ref{thm: existence of generic stab} and the determination of the generic stabilizers occupy the entirety of the present work, and involve a great deal of case analysis. In fact Lemma~\ref{lem: reduction to projective space} mentioned above, applied to the examples appearing in Theorem~\ref{thm: existence of generic stab}(ii), shows that the statement of Theorem~\ref{thm: existence of generic stab}(i) about the action on modules would not remain true if we were to allow $G$ to be a general semisimple group rather than one which is simple. It therefore seems unlikely that there is a proof of Theorem~\ref{thm: existence of generic stab} which does not involve consideration of cases. Note that if the generic stabilizer is trivial, the open set $\hat X$ in the definition above is a union of regular orbits.

In order to state our remaining results, we need a little more notation. Let $T$ be a maximal torus of $G$, and $\Phi$ be the root system of $G$ with respect to $T$; let $\Pi = \{ \alpha_1, \dots, \alpha_\ell \}$ be a simple system in $\Phi$, numbered as in \cite{Bou}, and $\omega_1, \dots, \omega_\ell$ be the corresponding fundamental dominant weights. If $\lambda$ is a dominant weight, write $L(\lambda)$ for the irreducible $G$-module with highest weight $\lambda$.

Note that throughout this paper we work modulo graph automorphisms; thus for example if $G = A_\ell$ we treat just one of the modules $L(\omega_i)$, $L(\omega_{\ell + 1 - i})$. Moreover, for $G$ of type $B_\ell$, $C_\ell$ and $D_\ell$ we normally assume $\ell \geq 2$, $\ell \geq 3$ and $\ell \geq 4$ respectively; occasionally it is convenient to relax this assumption, in which case we say so explicitly. In addition, in view of Steinberg's tensor product theorem (see Theorem~\ref{thm: Steinberg}) we may and shall always assume that the dominant weight $\lambda$ is not a multiple of $p$.

Given $G$, $\lambda$, $p$ and $k$ as above, and $V = L(\lambda)$, according as we let $X$ be $V$ or $\Gk(V)$ we say that we are considering the {\em triple\/} $(G, \lambda, p)$ or the {\em quadruple\/} $(G, \lambda, p, k)$; each quadruple $(G, \lambda, p, k)$ is said to be {\em associated\/} to the triple $(G, \lambda, p)$. The triple or quadruple is called {\em large\/} or {\em small\/} according as $\dim X > \dim G$ or $\dim X \leq \dim G$. We say that a triple or quadruple {\em has TGS\/} if the corresponding action has trivial generic stabilizer. According as $G$ is of classical or exceptional type, we say that the triple or quadruple is {\em classical\/} or {\em exceptional\/}. According as $k = 1$ or $k > 1$ we say that $(G, \lambda, p, k)$ is a {\em first quadruple\/} or a {\em higher quadruple\/}, and the variety $\Gk(V)$ is a {\em first Grassmannian variety\/} or a {\em higher Grassmannian variety\/}.

In proving that actions have trivial generic stabilizer we shall treat triples and quadruples separately (although, as we shall see, in almost all cases where a triple has TGS we can immediately conclude that all associated quadruples do as well). On the other hand, because the actions of $G$ on the module $V$ and the first Grassmannian variety $\G{1}(V)$ are so closely related, when determining generic stabilizers which are non-trivial it makes sense to treat triples and the associated first quadruples together. (We shall say more about the structure of this work later in this section.)

In addition to determining the existence of, and identifying, (semi-)generic stabilizers, we shall consider the questions of the existence of dense orbits and of regular orbits. Clearly any large triple or quadruple has no dense orbit. For each small triple or quadruple, we shall determine whether or not there is a dense orbit. On the other hand, if $X$ is a variety with $\dim X < \dim G$ then clearly there can be no regular orbit. If $X = V$, reference to \cite{Lubpaper} shows that we only have $\dim X = \dim G$ if $X$ is the Lie algebra of $G$, in which case any semisimple element is stabilized by at least a maximal torus, and the complement of the set of semisimple elements has positive codimension; thus in such a case no point can have a finite stabilizer, and it follows that any small triple has no regular orbit. If instead $X = \Gk(V)$, using \cite{Lubpaper} again we may identify the few cases where $\dim X = \dim G$; in each such case we shall find that the generic stabilizer is finite but non-trivial, and it follows that any small quadruple likewise has no regular orbit. In particular any small triple or quadruple does not have TGS. For each large triple which does not have TGS, we shall determine whether or not there is a regular orbit; we shall not however address this question for large quadruples which do not have TGS, because we have been unable to determine this in all but a very few cases.

Our main results will be given in six tables, which between them list all instances of triples and quadruples where the generic stabilizer is non-trivial. Tables~\ref{table: large triple and first quadruple non-TGS}, \ref{table: small classical triple and first quadruple generic stab} and \ref{table: small exceptional triple and first quadruple generic stab} concern triples and the associated first quadruples, while Tables~\ref{table: large higher quadruple non-TGS}, \ref{table: small classical higher quadruple generic stab} and \ref{table: small exceptional higher quadruple generic stab} concern higher quadruples; within each set of three, the first concerns triples or quadruples which are large, the second those which are small and classical, and the third those which are small and exceptional. The first few columns of each table specify the actions by listing $G$, $\lambda$, $\ell$ in the case of classical triples or quadruples, $p$, and $k$ in the case of higher quadruples. The next one or two columns give the generic stabilizers, denoted $C_X$ where $X = V$ or $\Gk(V)$ as appropriate (in Table~\ref{table: small classical higher quadruple generic stab}, the cases mentioned in the statement of Theorem~\ref{thm: existence of generic stab}(ii) as having only a semi-generic stabilizer are indicated by the presence of a symbol \lq $(*)$' beside the entry); the notation used for these groups is explained in the following section. The penultimate column of Table~\ref{table: large triple and first quadruple non-TGS} states whether or not the large triple has a regular orbit; the corresponding column in Tables~\ref{table: small classical triple and first quadruple generic stab}, \ref{table: small exceptional triple and first quadruple generic stab}, \ref{table: small classical higher quadruple generic stab} and \ref{table: small exceptional higher quadruple generic stab} states whether or not the small triple or quadruple has a dense orbit (in Tables~\ref{table: small classical triple and first quadruple generic stab} and \ref{table: small exceptional triple and first quadruple generic stab} each entry consists of two words \lq yes' or \lq no', with the first relating to the triple and the second to the first quadruple). In addition, in four rows of Table~\ref{table: large triple and first quadruple non-TGS}, indicated by asterisks in the penultimate column, the dimension of the module exceeds that of the group by one, so that the triple is large but the associated first quadruple is small; we find that in each such case the first quadruple has a dense orbit (whereas of course the triple cannot). The final column in each of the tables gives the reference to the Proposition in which the information provided is established; note that the existence or otherwise of a dense orbit follows immediately from comparing the codimension of the generic stabilizer to the dimension of the variety, and will not be mentioned in the statement of the Proposition concerned.

Throughout this work, if a parameter \lq $q$' occurs then the characteristic $p$ is finite and $q$ is a power of $p$.

\begin{table}
\caption{Large triples and associated first quadruples not having TGS}\label{table: large triple and first quadruple non-TGS}
\tabcapsp
$$
\begin{array}{|c|c|c|c|c|c|c|c|}
\hline
G      & \lambda                 & \ell      & p            & C_V                  & C_{\G{1}(V)}       & \hbox{regular?} & \hbox{reference} \mtbs \\
\hline
A_\ell & 3\omega_1               & 1         & {} \geq 5    & \Z_3                 & S_3                & \pay            & \ref{prop: A_1, 3omega_1 module} \mtbs \\
       & 3\omega_1               & 2         & {} \geq 5    & {\Z_3}^2             & {\Z_3}^2.\Z_2      & \hbox{yes}      & \ref{prop: A_2, 3omega_1, C_4, omega_4 modules} \mtbs \\
       & 4\omega_1               & 1         & {} \geq 5    & {\Z_2}^2             & {\Z_2}^2           & \hbox{yes}      & \ref{prop: B_ell or D_ell, 2omega_1 module} \mtbs \\
       & 2\omega_2               & 3         & {} \geq 3    & {\Z_2}^4             & {\Z_2}^4           & \hbox{yes}      & \ref{prop: B_ell or D_ell, 2omega_1 module} \mtbs \\
       & \omega_3                & 8         & {} \neq 3    & {\Z_3}^4.\Z_{(p, 2)} & {\Z_3}^4.\Z_2      & \hbox{yes}      & \ref{prop: A_8, omega_3, A_7, omega_4, D_8, omega_8 modules, non-special characteristic} \mtbs \\
       & \omega_3                & 8         & 3            & {\Z_3}^2             & {\Z_3}^2.\Z_2      & \hbox{yes}      & \ref{prop: A_8, omega_3, A_7, omega_4, D_8, omega_8 modules, special characteristic} \mtbs \\
       & \omega_4                & 7         & \geq 3       & {\Z_2}^6             & {\Z_2}^6           & \hbox{yes}      & \ref{prop: A_8, omega_3, A_7, omega_4, D_8, omega_8 modules, non-special characteristic} \mtbs \\
       & \omega_4                & 7         & 2            & {\Z_2}^3             & {\Z_2}^3           & \hbox{yes}      & \ref{prop: A_8, omega_3, A_7, omega_4, D_8, omega_8 modules, special characteristic} \mtbs \\
       & \omega_1 + \omega_2     & 3         & 3            & Alt_5                & S_5                & \pan            & \ref{prop: A_3, omega_1 + omega_2 module} \mtbs \\
       & \omega_1 + q\omega_1    & {} \geq 1 & {} < \infty  & \PSU_{\ell + 1}(q)   & \PGU_{\ell + 1}(q) & \pand           & \ref{prop: A_ell, omega_1 + q omega_1 and omega_1 + q omega_ell modules} \mtbs \\
       & \omega_1 + q\omega_\ell & {} \geq 2 & {} < \infty  & \PSL_{\ell + 1}(q)   & \PGL_{\ell + 1}(q) & \pan            & \ref{prop: A_ell, omega_1 + q omega_1 and omega_1 + q omega_ell modules} \mtbs \\
\hline
B_\ell & 2\omega_1               & {} \geq 2 & {} \geq 3    & {\Z_2}^{2\ell}       & {\Z_2}^{2\ell}     & \hbox{yes}      & \ref{prop: B_ell or D_ell, 2omega_1 module} \mtbs \\
       & \omega_1 + \omega_2     & 2         & 5            & \{ 1 \}              & \Z_2               & \hbox{yes}      & \ref{prop: B_2, omega_1 + omega_2 module, p = 5, C_4, omega_3 module, p = 3, k = 1} \mtbs \\
\hline
C_\ell & \omega_3                & 4         & 3            & \{ 1 \}              & \Z_2               & \hbox{yes}      & \ref{prop: B_2, omega_1 + omega_2 module, p = 5, C_4, omega_3 module, p = 3, k = 1} \mtbs \\
       & \omega_4                & 4         & {} \geq 3    & {\Z_2}^6             & {\Z_2}^6           & \hbox{yes}      & \ref{prop: A_2, 3omega_1, C_4, omega_4 modules} \mtbs \\
\hline
D_\ell & 2\omega_1               & {} \geq 4 & {} \geq 3    & {\Z_2}^{2\ell - 2}   & {\Z_2}^{2\ell - 2} & \hbox{yes}      & \ref{prop: B_ell or D_ell, 2omega_1 module} \mtbs \\
       & \omega_8                & 8         & {} \geq 3    & {\Z_2}^8             & {\Z_2}^8           & \hbox{yes}      & \ref{prop: A_8, omega_3, A_7, omega_4, D_8, omega_8 modules, non-special characteristic} \mtbs \\
       & \omega_8                & 8         & 2            & {\Z_2}^4             & {\Z_2}^4           & \hbox{yes}      & \ref{prop: A_8, omega_3, A_7, omega_4, D_8, omega_8 modules, special characteristic} \mtbs \\
\hline
\multicolumn{8}{l}{{}^\dagger \hbox{unless } \ell = 1 \hbox{ and } q \leq 3 \hbox{, in which case \lq yes'}} \\ 
\end{array}
$$
\end{table}

\begin{table}
\caption{Small classical triples and associated first quadruples}\label{table: small classical triple and first quadruple generic stab}
\tabcapsp
$$
\begin{array}{|c|c|c|c|c|c|cc|c|}
\hline
G      & \lambda                & \ell                    & p          & C_V                         & C_{\G{1}(V)}                   & \multicolumn{2}{|c|}{\hbox{dense?}}  & \hbox{reference} \mtbs \\
\hline
A_\ell & \omega_1               & {} \geq 1               & \hbox{any} & A_{\ell - 1} U_\ell         & A_{\ell - 1} T_1 U_\ell        & \yy                                  & \ref{prop: natural modules} \mtbs \\
       & 2\omega_1              & {} \geq 1, \hbox{ odd}  & {} \geq 3  & D_{\frac{1}{2}(\ell + 1)}   & D_{\frac{1}{2}(\ell + 1)}.\Z_2 & \ny                                  & \ref{prop: A_ell, 2 omega_1 and omega_2 modules} \mtbs \\
       & 2\omega_1              & {} \geq 2, \hbox{ even} & {} \geq 3  & B_{\frac{1}{2}\ell}         & B_{\frac{1}{2}\ell}            & \ny                                  & \ref{prop: A_ell, 2 omega_1 and omega_2 modules} \mtbs \\
       & \omega_2               & {} \geq 3, \hbox{ odd}  & \hbox{any} & C_{\frac{1}{2}(\ell + 1)}   & C_{\frac{1}{2}(\ell + 1)}      & \ny                                  & \ref{prop: A_ell, 2 omega_1 and omega_2 modules} \mtbs \\
       & \omega_2               & {} \geq 4, \hbox{ even} & \hbox{any} & C_{\frac{1}{2}\ell} U_\ell  & C_{\frac{1}{2}\ell} T_1 U_\ell & \yy                                  & \ref{prop: A_ell, 2 omega_1 and omega_2 modules} \mtbs \\
       & \omega_3               & 5                       & \hbox{any} & {A_2}^2.\Z_{(p, 2)}         & {A_2}^2.\Z_2                   & \ny                                  & \ref{prop: E_7, omega_7, D_6, omega_6, B_5, omega_5, A_5, omega_3, C_3, omega_3 modules} \mtbs \\
       & \omega_3               & 6                       & \hbox{any} & G_2                         & G_2                            & \ny                                  & \ref{prop: A_6, omega_3 module} \mtbs \\
       & \omega_3               & 7                       & \hbox{any} & A_2.\Z_{(p, 2)}             & A_2.\Z_2                       & \ny                                  & \ref{prop: A_7, omega_3 module} \mtbs \\
       & \omega_1 + \omega_2    & 2                       & {} \neq 3  & T_2                         & T_2                            & \nn                                  & \ref{prop: adjoint modules} \mtbs \\
       & \omega_1 + \omega_2    & 2                       & 3          & T_2.\Z_3                    & T_2.S_3                        & \ny                                  & \ref{prop: adjoint modules} \mtbs \\
       & \omega_1 + \omega_3    & 3                       & \hbox{any} & T_3.{\Z_{(p, 2)}}^2         & T_3.{\Z_{(p, 2)}}^2            & \nn                                  & \ref{prop: adjoint modules} \mtbs \\
       & \omega_1 + \omega_\ell & {} \geq 4               & \hbox{any} & T_\ell                      & T_\ell                         & \nn                                  & \ref{prop: adjoint modules} \mtbs \\
\hline
B_\ell & \omega_1               & {} \geq 2               & {} \geq 3  & D_\ell                      & D_\ell.\Z_2                    & \ny                                  & \ref{prop: natural modules} \mtbs \\
       & \omega_1               & {} \geq 2               & 2          & B_{\ell - 1} U_{2\ell - 1}  & B_{\ell - 1} T_1 U_{2\ell - 1} & \yy                                  & \ref{prop: natural module for B_ell, p = 2} \mtbs \\
       & \omega_2               & 2                       & \hbox{any} & A_1 U_3                     & A_1 T_1 U_3                    & \yy                                  & \ref{prop: natural modules} \mtbs \\
       & \omega_2               & {} \geq 3               & {} \geq 3  & T_\ell                      & T_\ell.\Z_2                    & \nn                                  & \ref{prop: adjoint modules} \mtbs \\
       & \omega_2               & 3                       & 2          & {B_1}^3                     & {B_1}^3                        & \nn                                  & \ref{prop: B_ell, omega_2 module, p = 2} \mtbs \\
       & \omega_2               & 4                       & 2          & {B_1}^4.{\Z_2}^2            & {B_1}^4.{\Z_2}^2               & \nn                                  & \ref{prop: B_ell, omega_2 module, p = 2} \mtbs \\
       & \omega_2               & {} \geq 5               & 2          & {B_1}^\ell                  & {B_1}^\ell                     & \nn                                  & \ref{prop: B_ell, omega_2 module, p = 2} \mtbs \\
       & 2\omega_2              & 2                       & {} \geq 3  & T_2                         & T_2.\Z_2                       & \nn                                  & \ref{prop: adjoint modules} \mtbs \\
       & \omega_3               & 3                       & \hbox{any} & G_2                         & G_2                            & \ny                                  & \ref{prop: B_3, omega_3 module} \mtbs \\
       & \omega_4               & 4                       & \hbox{any} & B_3                         & B_3                            & \ny                                  & \ref{prop: D_5, omega_5, B_4, omega_4 modules} \mtbs \\
       & \omega_5               & 5                       & \hbox{any} & A_4.\Z_{(p, 2)}             & A_4.\Z_2                       & \ny                                  & \ref{prop: E_7, omega_7, D_6, omega_6, B_5, omega_5, A_5, omega_3, C_3, omega_3 modules} \mtbs \\
       & \omega_6               & 6                       & \hbox{any} & {A_2}^2.{\Z_{(p, 2)}}^2     & {A_2}^2.\Z_{(p, 2)}.\Z_2       & \nn                                  & \ref{prop: D_7, omega_7, B_6, omega_6 modules} \mtbs \\
\hline
C_\ell & \omega_1               & {} \geq 3               & \hbox{any} & C_{\ell - 1} U_{2\ell - 1}  & C_{\ell - 1} T_1 U_{2\ell - 1} & \yy                                  & \ref{prop: natural modules} \mtbs \\
       & 2\omega_1              & {} \geq 3               & {} \geq 3  & T_\ell                      & T_\ell.\Z_2                    & \nn                                  & \ref{prop: adjoint modules} \mtbs \\
       & \omega_2               & 3                       & {} \neq 3  & {C_1}^3                     & {C_1}^3                        & \nn                                  & \ref{prop: C_ell, omega_2 module} \mtbs \\
       & \omega_2               & 3                       & 3          & {C_1}^3.\Z_3                & {C_1}^3.S_3                    & \ny                                  & \ref{prop: C_ell, omega_2 module} \mtbs \\
       & \omega_2               & 4                       & \hbox{any} & {C_1}^4.{\Z_{(p, 2)}}^2     & {C_1}^4.{\Z_{(p, 2)}}^2        & \nn                                  & \ref{prop: C_ell, omega_2 module} \mtbs \\
       & \omega_2               & {} \geq 5               & \hbox{any} & {C_1}^\ell                  & {C_1}^\ell                     & \nn                                  & \ref{prop: C_ell, omega_2 module} \mtbs \\
       & \omega_3               & 3                       & {} \geq 3  & \tilde A_2                  & \tilde A_2.\Z_2                & \ny                                  & \ref{prop: E_7, omega_7, D_6, omega_6, B_5, omega_5, A_5, omega_3, C_3, omega_3 modules} \mtbs \\
       & \omega_3               & 3                       & 2          & G_2                         & G_2                            & \ny                                  & \ref{prop: C_3, omega_3, C_4, omega_4, C_5, omega_5, C_6, omega_6 modules, p = 2} \mtbs \\
       & \omega_4               & 4                       & 2          & C_3                         & C_3                            & \ny                                  & \ref{prop: C_3, omega_3, C_4, omega_4, C_5, omega_5, C_6, omega_6 modules, p = 2} \mtbs \\
       & \omega_5               & 5                       & 2          & \tilde A_4.\Z_2             & \tilde A_4.\Z_2                & \ny                                  & \ref{prop: C_3, omega_3, C_4, omega_4, C_5, omega_5, C_6, omega_6 modules, p = 2} \mtbs \\
       & \omega_6               & 6                       & 2          & {{\tilde A}_2}{}^2.{\Z_2}^2 & {{\tilde A}_2}{}^2.{\Z_2}^2    & \nn                                  & \ref{prop: C_3, omega_3, C_4, omega_4, C_5, omega_5, C_6, omega_6 modules, p = 2} \mtbs \\
\hline
D_\ell & \omega_1               & {} \geq 4               & \hbox{any} & B_{\ell - 1}                & B_{\ell - 1}                   & \ny                                  & \ref{prop: natural modules} \mtbs \\
       & \omega_2               & {} \geq 4               & {} \geq 3  & T_\ell                      & T_\ell.\Z_{(2, \ell)}          & \nn                                  & \ref{prop: adjoint modules} \mtbs \\
       & \omega_2               & 4                       & 2          & T_4.{\Z_2}^3.{\Z_2}^2       & T_4.{\Z_2}^3.{\Z_2}^2          & \nn                                  & \ref{prop: adjoint modules} \mtbs \\
       & \omega_2               & {} \geq 5               & 2          & T_\ell.{\Z_2}^{\ell - 1}    & T_\ell.{\Z_2}^{\ell - 1}       & \nn                                  & \ref{prop: adjoint modules} \mtbs \\
       & \omega_5               & 5                       & \hbox{any} & B_3 U_8                     & B_3 T_1 U_8                    & \yy                                  & \ref{prop: D_5, omega_5, B_4, omega_4 modules} \mtbs \\
       & \omega_6               & 6                       & \hbox{any} & A_5.\Z_{(p, 2)}             & A_5.\Z_2                       & \ny                                  & \ref{prop: E_7, omega_7, D_6, omega_6, B_5, omega_5, A_5, omega_3, C_3, omega_3 modules} \mtbs \\
       & \omega_7               & 7                       & \hbox{any} & {G_2}^2.\Z_{(p, 2)}         & {G_2}^2.\Z_2                   & \ny                                  & \ref{prop: D_7, omega_7, B_6, omega_6 modules} \mtbs \\
\hline
\end{array}
$$
\end{table}

\begin{table}
\caption{Small exceptional triples and associated first quadruples}\label{table: small exceptional triple and first quadruple generic stab}
\tabcapsp
$$
\begin{array}{|c|c|c|c|c|cc|c|}
\hline
G      & \lambda  & p          & C_V             & C_{\G{1}(V)}    & \multicolumn{2}{|c|}{\hbox{dense?}} & \hbox{reference} \mtbs \\
\hline
E_6    & \omega_1 & \hbox{any} & F_4             & F_4             & \ny                                 & \ref{prop: E_6, omega_1, F_4, omega_4 modules} \mtbs \\
       & \omega_2 & \hbox{any} & T_6             & T_6             & \nn                                 & \ref{prop: adjoint modules} \mtbs \\
\hline
E_7    & \omega_1 & \hbox{any} & T_7.\Z_{(p, 2)} & T_7.\Z_2        & \nn                                 & \ref{prop: adjoint modules} \mtbs \\
       & \omega_7 & \hbox{any} & E_6.\Z_{(p, 2)} & E_6.\Z_2        & \ny                                 & \ref{prop: E_7, omega_7, D_6, omega_6, B_5, omega_5, A_5, omega_3, C_3, omega_3 modules} \mtbs \\
\hline
E_8    & \omega_8 & \hbox{any} & T_8.\Z_{(p, 2)} & T_8.\Z_2        & \nn                                 & \ref{prop: adjoint modules} \mtbs \\
\hline
F_4    & \omega_1 & {} \geq 3  & T_4             & T_4.\Z_2        & \nn                                 & \ref{prop: adjoint modules} \mtbs \\
       & \omega_1 & 2          & \tilde D_4      & \tilde D_4      & \nn                                 & \ref{prop: F_4, omega_1 module, p = 2} \mtbs \\
       & \omega_4 & {} \neq 3  & D_4             & D_4             & \nn                                 & \ref{prop: E_6, omega_1, F_4, omega_4 modules} \mtbs \\
       & \omega_4 & 3          & D_4.\Z_3        & D_4.S_3         & \ny                                 & \ref{prop: E_6, omega_1, F_4, omega_4 modules} \mtbs \\
\hline
G_2    & \omega_1 & {} \geq 3  & A_2             & A_2.\Z_2        & \ny                                 & \ref{prop: G_2, omega_1 module} \mtbs \\
       & \omega_1 & 2          & A_1 U_5         & A_1 T_1 U_5     & \yy                                 & \ref{prop: G_2, omega_1 module} \mtbs \\
       & \omega_2 & {} \neq 3  & T_2.\Z_{(p, 2)} & T_2.\Z_2        & \nn                                 & \ref{prop: adjoint modules} \mtbs \\
       & \omega_2 & 3          & \tilde A_2      & \tilde A_2.\Z_2 & \ny                                 & \ref{prop: G_2, omega_2 module, p = 3} \mtbs \\
\hline
\end{array}
$$
\end{table}

\begin{table}
\caption{Large higher quadruples not having TGS}\label{table: large higher quadruple non-TGS}
\tabcapsp
$$
\begin{array}{|c|c|c|c|c|c|c|}
\hline
G      & \lambda              & \ell      & p          & k & C_{\Gk(V)}                 & \hbox{reference} \mtbs \\
\hline
A_\ell & 2\omega_1            & 3         & {} \geq 3  & 2 & {\Z_2}^3.{\Z_2}^2          & \ref{prop: A_ell, 2omega_1 module, k = 2} \mtbs \\
       & 2\omega_1            & {} \geq 4 & {} \geq 3  & 2 & {\Z_2}^\ell                & \ref{prop: A_ell, 2omega_1 module, k = 2} \mtbs \\
       & 2\omega_1            & 2         & {} \geq 3  & 3 & \Z_{3/(p, 3)}.S_3          & \ref{prop: A_2, 2omega_1 module, k = 3, A_4, omega_2 module, k = 5} \mtbs \\
       & 3\omega_1            & 1         & {} \geq 5  & 2 & {\Z_2}^2                   & \ref{prop: A_1, 3omega_1 module, k = 2} \mtbs \\
       & \omega_2             & 4         & \hbox{any} & 5 & \Z_{5/(p, 5)}.Dih_{10}     & \ref{prop: A_2, 2omega_1 module, k = 3, A_4, omega_2 module, k = 5} \mtbs \\
       & \omega_2             & 5         & \hbox{any} & 3 & T_1.\Z_{3/(p, 3)}.S_3      & \ref{prop: E_6, omega_1, A_5, omega_2 modules, k = 3} \mtbs \\
       & \omega_3             & 5         & \hbox{any} & 2 & T_2.\Z_{2/(p, 2)}.\Z_2     & \ref{prop: E_7, omega_7, D_6, omega_6, A_5, omega_3, C_3, omega_3 modules, k = 2} \mtbs \\
       & \omega_1 + q\omega_1 & 1         & < \infty   & 2 & \Z_2                       & \ref{prop: A_1, omega_1 + q omega_1 module, k = 2} \mtbs \\
\hline
B_\ell & \omega_2             & 3         & 2          & 2 & T_1                        & \ref{prop: F_4, omega_1, B_3, omega_2 modules, p = 2, k = 2} \mtbs \\
       & \omega_4             & 4         & \hbox{any} & 3 & \Z_{2/(p, 2)}.\Z_2         & \ref{prop: B_4, omega_4 module, k = 3} \mtbs \\
       & \omega_5             & 5         & \hbox{any} & 2 & \Z_{2/(p, 2)}.\Z_2         & \ref{prop: B_5, omega_5 module, k = 2} \mtbs \\
\hline
C_\ell & \omega_2             & 3         & {} \neq 3  & 2 & T_1                        & \ref{prop: F_4, omega_4, C_3, omega_2 modules, k = 2} \mtbs \\
       & \omega_2             & 3         & 3          & 2 & T_1.\Z_2                   & \ref{prop: F_4, omega_4, C_3, omega_2 modules, k = 2} \mtbs \\
       & \omega_3             & 3         & {} \geq 3  & 2 & {\Z_2}^4                   & \ref{prop: E_7, omega_7, D_6, omega_6, A_5, omega_3, C_3, omega_3 modules, k = 2} \mtbs \\
       & \omega_4             & 4         & 2          & 3 & \Z_2                       & \ref{prop: C_4, omega_4 module, p = 2, k = 3} \mtbs \\
       & \omega_5             & 5         & 2          & 2 & \Z_2                       & \ref{prop: C_5, omega_5 module, p = 2, k = 2} \mtbs \\
\hline
D_\ell & \omega_5             & 5         & \hbox{any} & 4 & {\Z_{2/(p, 2)}}^2.{\Z_2}^2 & \ref{prop: D_5, omega_5 module, k = 4} \mtbs \\
\hline
\end{array}
$$
\end{table}

\begin{table}
\caption{Small classical higher quadruples}\label{table: small classical higher quadruple generic stab}
\tabcapsp
$$
\begin{array}{|c|c|c|c|c|c|c|c|}
\hline
G      & \lambda   & \ell                    & p          & k          & C_{\Gk(V)}                                                              & \hbox{dense?} & \hbox{reference} \mtbs \\
\hline
A_\ell & \omega_1  & {} \geq 1               & \hbox{any} & \hbox{any} & A_{\ell - k} A_{k - 1} T_1 U_{k(\ell + 1 - k)}                          & \hbox{yes}    & \ref{prop: natural modules, k arbitrary} \mtbs \\
       & 2\omega_1 & 2                       & {} \geq 3  & 2          & {\Z_2}^2.S_3                                                            & \hbox{yes}    & \ref{prop: A_ell, 2omega_1 module, k = 2} \mtbs \\
       & \omega_2  & 3                       & \hbox{any} & 2          & {A_1}^2 T_1.\Z_2                                                        & \hbox{yes}    & \ref{prop: natural modules, k arbitrary} \mtbs \\
       & \omega_2  & 5                       & \hbox{any} & 2          & {A_1}^3.S_3                                                             & \hbox{yes}    & \ref{prop: A_ell, ell odd, omega_2 module, k = 2} \mtbs \\
       & \omega_2  & 7                       & \hbox{any} & 2          & {A_1}^4.{\Z_2}^2                                                        & \hbox{no}     & \ref{prop: A_ell, ell odd, omega_2 module, k = 2} \mtbs \\
       & \omega_2  & {} \geq 9, \hbox{ odd}  & \hbox{any} & 2          & {A_1}^{\frac{1}{2}(\ell + 1)}                                           & \hbox{no}     & \ref{prop: A_ell, ell odd, omega_2 module, k = 2} \mtbs \\
       & \omega_2  & {} \geq 4, \hbox{ even} & \hbox{any} & 2          & A_1 T_1 U_\ell                                                          & \hbox{yes}    & \ref{prop: A_ell, ell even, omega_2 module, k = 2} \mtbs \\
       & \omega_2  & 3                       & \hbox{any} & 3          & {A_1}^2                                                                 & \hbox{yes}    & \ref{prop: natural modules, k arbitrary} \mtbs \\
       & \omega_2  & 4                       & \hbox{any} & 3          & A_1                                                                     & \hbox{yes}    & \ref{prop: A_4, omega_2 module, k = 3} \mtbs \\
       & \omega_2  & 4                       & \hbox{any} & 4          & S_5                                                                     & \hbox{yes}    & \ref{prop: A_4, omega_2 module, k = 4} \mtbs \\
\hline
B_\ell & \omega_1  & {} \geq 2               & {} \geq 3  & \hbox{odd}  & B_{\frac{1}{2}(k - 1)} D_{\ell - \frac{1}{2}(k - 1)}.\Z_2              & \hbox{yes}    & \ref{prop: natural modules, k arbitrary} \mtbs \\
       & \omega_1  & {} \geq 2               & {} \geq 3  & \hbox{even} & D_{\frac{1}{2}k} B_{\ell - \frac{1}{2}k}.\Z_2                          & \hbox{yes}    & \ref{prop: natural modules, k arbitrary} \mtbs \\
       & \omega_1  & {} \geq 2               & 2          & \hbox{odd}  & B_{\frac{1}{2}(k - 1)} B_{\ell - \frac{1}{2}(k + 1)} T_1 U_{2\ell - 1} & \hbox{yes}    & \ref{prop: natural module for B_ell, p = 2, k arbitrary} \mtbs \\
       & \omega_1  & {} \geq 2               & 2          & \hbox{even} & B_{\frac{1}{2}k} B_{\ell - \frac{1}{2}k}                               & \hbox{yes}    & \ref{prop: natural module for B_ell, p = 2, k arbitrary} \mtbs \\
       & \omega_2  & 2                       & \hbox{any} & 2           & {A_1}^2                                                                & \hbox{yes}    & \ref{prop: natural modules, k arbitrary} \mtbs \\
       & \omega_3  & 3                       & \hbox{any} & 2           & A_2 T_1.\Z_2                                                           & \hbox{yes}    & \ref{prop: B_3, omega_3 module, k = 2 or 3} \mtbs \\
       & \omega_3  & 3                       & \hbox{any} & 3           & {A_1}^2                                                                & \hbox{yes}    & \ref{prop: B_3, omega_3 module, k = 2 or 3} \mtbs \\
       & \omega_3  & 3                       & {} \geq 3  & 4           & {B_1}^2                                                                & \hbox{no}     & \ref{prop: B_3, omega_3 module, k = 4} \mtbs \\
       & \omega_3  & 3                       & 2          & 4           & {B_1}^2 (*)                                                            & \hbox{no}     & \ref{prop: B_3, omega_3 module, k = 4} \mtbs \\
       & \omega_4  & 4                       & \hbox{any} & 2           & A_2 T_1.\Z_2                                                           & \hbox{no}     & \ref{prop: D_5, omega_5, B_4, omega_4 modules, k = 2} \mtbs \\
\hline
C_\ell & \omega_1  & {} \geq 3               & \hbox{any} & \hbox{odd}  & C_{\frac{1}{2}(k - 1)} C_{\ell - \frac{1}{2}(k + 1)} T_1 U_{2\ell - 1} & \hbox{yes}    & \ref{prop: natural modules, k arbitrary} \mtbs \\
       & \omega_1  & {} \geq 3               & \hbox{any} & \hbox{even} & C_{\frac{1}{2}k} C_{\ell - \frac{1}{2}k}                               & \hbox{yes}    & \ref{prop: natural modules, k arbitrary} \mtbs \\
       & \omega_3  & 3                       & 2          & 2           & \tilde A_2 T_1.\Z_2                                                    & \hbox{yes}    & \ref{prop: C_3, omega_3 module, p = 2, k = 2 or 3} \mtbs \\
       & \omega_3  & 3                       & 2          & 3           & {\tilde A_1}{}^2                                                       & \hbox{yes}    & \ref{prop: C_3, omega_3 module, p = 2, k = 2 or 3} \mtbs \\
       & \omega_3  & 3                       & 2          & 4           & {C_1}^2 (*)                                                            & \hbox{no}     & \ref{prop: C_3, omega_3 module, p = 2, k = 4} \mtbs \\
       & \omega_4  & 4                       & 2          & 2           & \tilde A_2 T_1.\Z_2                                                    & \hbox{no}     & \ref{prop: C_4, omega_4 module, p = 2, k = 2} \mtbs \\
\hline
D_\ell & \omega_1  & {} \geq 4               & \hbox{any} & \hbox{odd}  & B_{\frac{1}{2}(k - 1)} B_{\ell - \frac{1}{2}(k + 1)}                   & \hbox{yes}    & \ref{prop: natural modules, k arbitrary} \mtbs \\
       & \omega_1  & {} \geq 4               & \hbox{any} & \hbox{even} & D_{\frac{1}{2}k} D_{\ell - \frac{1}{2}k}.\Z_2                          & \hbox{yes}    & \ref{prop: natural modules, k arbitrary} \mtbs \\
       & \omega_5  & 5                       & \hbox{any} & 2           & G_2 B_1                                                                & \hbox{yes}    & \ref{prop: D_5, omega_5, B_4, omega_4 modules, k = 2} \mtbs \\
       & \omega_5  & 5                       & \hbox{any} & 3           & {A_1}^2                                                                & \hbox{yes}    & \ref{prop: D_5, omega_5 module, k = 3} \mtbs \\
       & \omega_6  & 6                       & \hbox{any} & 2           & {A_1}^3.\Z_{2/(p, 2)}.\Z_2                                             & \hbox{no}     & \ref{prop: E_7, omega_7, D_6, omega_6, A_5, omega_3, C_3, omega_3 modules, k = 2} \mtbs \\
\hline
\end{array}
$$
\end{table}

\begin{table}
\caption{Small exceptional higher quadruples}\label{table: small exceptional higher quadruple generic stab}
\tabcapsp
$$
\begin{array}{|c|c|c|c|c|c|c|}
\hline
G   & \lambda  & p          & k & C_{\Gk(V)}             & \hbox{dense?} & \hbox{reference} \mtbs \\
\hline
E_6 & \omega_1 & \hbox{any} & 2 & D_4.S_3                & \hbox{yes}    & \ref{prop: E_6, omega_1 module, k = 2} \mtbs \\
    & \omega_1 & \hbox{any} & 3 & A_2.\Z_{3/(p, 3)}.S_3  & \hbox{no}     & \ref{prop: E_6, omega_1, A_5, omega_2 modules, k = 3} \mtbs \\
\hline
E_7 & \omega_7 & \hbox{any} & 2 & D_4.\Z_{2/(p, 2)}.\Z_2 & \hbox{no}     & \ref{prop: E_7, omega_7, D_6, omega_6, A_5, omega_3, C_3, omega_3 modules, k = 2} \mtbs \\
\hline
F_4 & \omega_1 & 2          & 2 & \tilde A_2             & \hbox{no}     & \ref{prop: F_4, omega_1, B_3, omega_2 modules, p = 2, k = 2} \mtbs \\
    & \omega_4 & {} \neq 3  & 2 & A_2                    & \hbox{no}     & \ref{prop: F_4, omega_4, C_3, omega_2 modules, k = 2} \mtbs \\
    & \omega_4 & 3          & 2 & A_2.\Z_2               & \hbox{no}     & \ref{prop: F_4, omega_4, C_3, omega_2 modules, k = 2} \mtbs \\
\hline
G_2 & \omega_1 & {} \geq 3  & 2 & A_1 T_1.\Z_2           & \hbox{yes}    & \ref{prop: G_2, omega_1 module, k = 2} \mtbs \\
    & \omega_1 & 2          & 2 & A_1 \tilde A_1         & \hbox{yes}    & \ref{prop: G_2, omega_1 module, k = 2} \mtbs \\
    & \omega_1 & {} \geq 3  & 3 & A_1                    & \hbox{no}     & \ref{prop: G_2, omega_1 module, k = 3} \mtbs \\
    & \omega_1 & 2          & 3 & A_1 U_2                & \hbox{yes}    & \ref{prop: G_2, omega_1 module, k = 3} \mtbs \\
    & \omega_2 & 3          & 2 & \tilde A_1 T_1.\Z_2    & \hbox{yes}    & \ref{prop: G_2, omega_2 module, p = 3, k = 2 or 3} \mtbs \\
    & \omega_2 & 3          & 3 & A_1                    & \hbox{no}     & \ref{prop: G_2, omega_2 module, p = 3, k = 2 or 3} \mtbs \\
\hline
\end{array}
$$
\end{table}

The theorems which we state concern large triples and associated first quadruples, small triples and associated first quadruples, large higher quadruples, and small higher quadruples respectively.

\begin{theorem}\label{thm: large triple and first quadruple generic stab}
If a large triple or associated first quadruple appears in Table~\ref{table: large triple and first quadruple non-TGS} then it has generic stabilizer as given there; in addition for a large triple the existence or otherwise of a regular orbit is indicated. If it does not appear in Table~\ref{table: large triple and first quadruple non-TGS} then it has TGS.
\end{theorem}

Note that Table~\ref{table: large triple and first quadruple non-TGS} contains two instances where the triple does have TGS, but the associated first quadruple does not.

\begin{theorem}\label{thm: small triple and first quadruple generic stab}
The generic stabilizer for a small triple or associated first quadruple is given in Table~\ref{table: small classical triple and first quadruple generic stab} or \ref{table: small exceptional triple and first quadruple generic stab} according as the triple or first quadruple is classical or exceptional; in addition the existence or otherwise of a dense orbit is indicated.
\end{theorem}

\begin{theorem}\label{thm: large higher quadruple generic stab}
If a large higher quadruple appears in Table~\ref{table: large higher quadruple non-TGS} then it has generic stabilizer as given there. If it does not appear in Table~\ref{table: large higher quadruple non-TGS} then it has TGS.
\end{theorem}

\begin{theorem}\label{thm: small higher quadruple generic stab}
The (semi-)generic stabilizer for a small higher quadruple is given in Table~\ref{table: small classical higher quadruple generic stab} or \ref{table: small exceptional higher quadruple generic stab} according as the higher quadruple is classical or exceptional; in addition the existence or otherwise of a dense orbit is indicated.
\end{theorem}

Note that in the rows of Table~\ref{table: small classical higher quadruple generic stab} corresponding to $\lambda = \omega_1$, the entry in the fifth column giving the value of $k$ is \lq any' (if $G = A_\ell$) or either \lq odd' or \lq even' (if $G = B_\ell$, $C_\ell$ or $D_\ell$); it is implicitly assumed that we restrict ourselves to values of $k$ satisfying $2 \leq k \leq \frac{1}{2} \dim L(\omega_1)$.

We find that these theorems have some interesting consequences. Inspection of Tables~\ref{table: large triple and first quadruple non-TGS} and \ref{table: large higher quadruple non-TGS} immediately reveals the following.

\begin{corollary}\label{cor: large exceptional triples and quadruples have TGS}
Any large exceptional triple or quadruple has TGS.
\end{corollary}

The next consequence follows not from the statements but from the proofs of the theorems.

\begin{corollary}\label{cor: TGS iff codim inequality}
A triple or quadruple has TGS if and only if, for any group element which is either semisimple of prime order modulo the centre, or unipotent of order $p$, the codimension of its fixed point variety is strictly greater than the dimension of its conjugacy class.
\end{corollary}

Indeed, at the start of our analysis we establish the reverse implication (see Section~\ref{sect: conditions}, where we define conditions $\ssdiamcon$ and $\udiamcon$), and thereafter use it consistently to prove that triples and quadruples have TGS; the forward implication follows from the fact that all cases which we find to have TGS are proved in this way.

The next few consequences apply to triples. Firstly we have the following.

\begin{corollary}\label{cor: large triples have finite generic stabilizers}
The triple $(G, \lambda, p)$ is large if and only if it has a finite generic stabilizer.
\end{corollary}

The forward implication follows simply from the observation that all generic stabilizers in Table~\ref{table: large triple and first quadruple non-TGS} are finite (whereas in characteristic zero it was known in advance, as stated above), while the reverse implication has already been noted. The corresponding statement does not hold for quadruples; indeed we have seen that Table~\ref{table: large triple and first quadruple non-TGS} contains instances where the triple is large but the associated first quadruple is small and has a finite generic stabilizer.

Secondly inspection of Table~\ref{table: large triple and first quadruple non-TGS} shows that in most cases (including all cases in characteristic zero) the existence of a regular orbit is linked to the finiteness of the generic stabilizer.

\begin{corollary}\label{cor: finite generic stabilizers iff regular orbit, with exceptions}
In the action of $G$ on $L(\lambda)$, if there is a regular orbit then the generic stabilizer is finite; the converse holds unless the triple $(G, \lambda, p)$ is one of the following:
\begin{itemize}
\item[(i)] $(A_3, \omega_1 + \omega_2, 3)$,
\item[(ii)] $(A_\ell, \omega_1 + q\omega_1, p)$ with $p < \infty$ and either $\ell \geq 2$, or $\ell = 1$ and $q \geq 4$,
\item[(iii)] $(A_\ell, \omega_1 + q\omega_\ell, p)$ with $p < \infty$ and $\ell \geq 2$.
\end{itemize}
\end{corollary}

The first statement of this result does not require inspection of Table~\ref{table: large triple and first quadruple non-TGS}, since a straightforward argument shows that if there is a generic stabilizer in an action on a variety then no point can have stabilizer of dimension less than that of the generic stabilizer; however we shall not need this argument in the work here.

Thirdly inspection of Tables~\ref{table: small classical triple and first quadruple generic stab} and \ref{table: small exceptional triple and first quadruple generic stab} shows that the existence of a dense orbit is linked to the structure of the generic stabilizer.

\begin{corollary}\label{cor: dense orbit for triple means non-reductive gen stab}
In the action of $G$ on $L(\lambda)$, there is a dense orbit if and only if the generic stabilizer has non-reductive connected component.
\end{corollary}

Fourthly we have the following.

\begin{corollary}\label{cor: non-TGS triple means 1-dimensional weight spaces}
If there is a non-zero weight such that the corresponding weight space in $L(\lambda)$ has dimension greater than $1$, then the triple $(G, \lambda, p)$ has TGS.
\end{corollary}

This may be seen by observing that in each case in Tables~\ref{table: large triple and first quadruple non-TGS}, \ref{table: small classical triple and first quadruple generic stab} and \ref{table: small exceptional triple and first quadruple generic stab} all weight spaces corresponding to non-zero weights are $1$-dimensional.

The next two consequences concern higher quadruples, and are immediate from inspection of Tables~\ref{table: large higher quadruple non-TGS}, \ref{table: small classical higher quadruple generic stab} and \ref{table: small exceptional higher quadruple generic stab}. The first of these is the observation that quadruples with large values of $k$ only rarely fail to have TGS.

\begin{corollary}\label{cor: k at least 4 usually implies TGS}
If $4 \leq k \leq \frac{1}{2} \dim L(\lambda)$, the quadruple $(G, \lambda, p, k)$ has TGS unless one of the following holds:
\begin{itemize}
\item[(i)] $G$ is of classical type and $\lambda = \omega_1$ (so that $V$ is the natural module for $G$);
\item[(ii)] $G = A_4$, $\lambda = \omega_2$ and $k = 4$ or $5$;
\item[(iii)] $G = B_3$ (or $C_3$ if $p = 2$), $\lambda = \omega_3$ and $k = 4$;
\item[(iv)] $G = D_5$, $\lambda = \omega_5$ and $k = 4$.
\end{itemize}
\end{corollary}

The next involves quadruples associated to a given triple.

\begin{corollary}\label{cor: k and k' - dimensions and TGS}
For a given triple $(G, \lambda, p)$, and natural numbers $k$, $k'$ satisfying $k < k' \leq \frac{1}{2} \dim L(\lambda)$, the following hold:
\begin{itemize}
\item[(i)] we have $\dim C_{\Gk(V)} \geq \dim C_{\G{k'}(V)}$;
\item[(ii)] if the associated quadruple $(G, \lambda, p, k)$ has TGS, so does the associated quadruple $(G, \lambda, p, k')$.
\end{itemize}
\end{corollary}

Finally, we note that if $G$ has finitely many orbits on an irreducible variety $X$, then one orbit must be dense. As mentioned above, cases where $G$ has finitely many orbits on $X = \Gk(V)$ were classified in \cite{GLMS}, of which Corollary~1 states that, if $k = 1$, there is a dense orbit if and only if there are finitely many orbits. Comparison of Tables~\ref{table: small classical higher quadruple generic stab} and \ref{table: small exceptional higher quadruple generic stab} here with \cite[Theorem 2]{GLMS} yields the following extension of this result to arbitrary $k$.

\begin{corollary}\label{cor: dense orbit usually means finitely many orbits}
If the action of $G$ on $\Gk(V)$ has a dense orbit, then either there are only finitely many orbits, or $G = A_\ell$ for $\ell \geq 8$ even, $\lambda = \omega_2$ and $k = 2$.
\end{corollary}

The structure of the remainder of this work is as follows. This chapter has four further sections. In Section~\ref{sect: notation} we establish notation to be used throughout this work. In Section~\ref{sect: weights and module structure} we recall some basic facts about the decomposition of modules into weight spaces, and discuss the key concept of strings of weights. In Section~\ref{sect: unipotent classes} we provide a considerable amount of detailed information which we shall require on unipotent classes and their closures. In Section~\ref{sect: prelim} we prove some preliminary results.

The next two chapters concern actions having TGS. In chapter~\ref{chap: TGS triples} we treat large triples, and show that any such which is not listed in Table~\ref{table: large triple and first quadruple non-TGS} has TGS: we start by giving a series of conditions which imply that a large triple has TGS, and then develop and apply increasingly refined methods to show that the large triples concerned satisfy them. In chapter~\ref{chap: TGS quadruples} we treat large quadruples, and show that any such which is not listed in Table~\ref{table: large higher quadruple non-TGS} has TGS: we prove a result which implies that, in all but two cases, if a large triple has TGS then all associated large quadruples also have TGS, after which we apply methods similar to those of the previous chapter to treat the remaining cases.

The final three chapters concern actions not having TGS. In chapter~\ref{chap: non-TGS methods} we explain some methods for treating such actions. In chapter~\ref{chap: non-TGS triples and first quadruples} we treat triples and first quadruples, and complete the proofs of Theorems~\ref{thm: large triple and first quadruple generic stab} and \ref{thm: small triple and first quadruple generic stab} by establishing the entries in Tables~\ref{table: large triple and first quadruple non-TGS}, \ref{table: small classical triple and first quadruple generic stab} and \ref{table: small exceptional triple and first quadruple generic stab}. Finally in chapter~\ref{chap: non-TGS higher quadruples} we treat higher quadruples, and complete the proofs of Theorems~\ref{thm: large higher quadruple generic stab} and \ref{thm: small higher quadruple generic stab} by establishing the entries in Tables~\ref{table: large higher quadruple non-TGS}, \ref{table: small classical higher quadruple generic stab} and \ref{table: small exceptional higher quadruple generic stab}.

It should be mentioned that for the work on large triples having TGS, much of the general strategy employed here is adapted from the PhD thesis \cite{Ken} of Kenneally, written under the supervision of the second author; this work tackled only part of the present problem, proving results about eigenspaces of semisimple elements but not addressing the action of unipotent elements, and considering only actions on modules but not on Grassmannian varieties. It seems rather surprising that in virtually all stages of the analysis here it proves possible to treat semisimple and unipotent elements in parallel. As a consequence there will be no need to refer to specific results obtained by Kenneally, since the calculations which he performed need to be extended to treat unipotent elements; but we acknowledge here that most of the results obtained in the present work on the action of semisimple elements on modules may be found in \cite{Ken}.

The authors are grateful to Martin Liebeck, Alexander Premet and Donna Testerman for a number of helpful conversations and discussions at various stages of this project, and to the anonymous referees for several suggestions and corrections which have led to improvements in the present work.

\section{Notation}\label{sect: notation}

In this section we establish notation to be used throughout the work.

To begin with, we let $K$ be an algebraically closed field of characteristic $p$ (writing as above $p = \infty$ if $K$ has characteristic zero), and $H$ be a simple algebraic group over $K$, of rank $\ell_H$; we write $Z(H)$ for the centre of $H$. We take a maximal torus $T_H$ of $H$, and let $N_H$ be its normalizer in $H$ and $W_H = N_H/T_H$ be the Weyl group of $H$. We let $\Phi_H$ be the irreducible root system of $H$ with respect to $T_H$; for each $\alpha \in \Phi_H$ we let $X_\alpha$ be the corresponding root subgroup of $H$, and $x_\alpha : K \to X_\alpha$ be an isomorphism of algebraic groups. As is usual, we assume that the maps $x_\alpha$ are chosen so that the Chevalley commutator relations hold, and so that for all $t \in K^*$ the element $n_\alpha(t) = x_\alpha(t) x_{-\alpha}(-t^{-1}) x_\alpha(t)$ lies in $N_H$ and $h_\alpha(t) = n_\alpha(t) n_\alpha(-1)$ lies in $T_H$; for $\alpha \in \Phi_H$ we set $n_\alpha = n_\alpha(1)$ and $w_\alpha = n_\alpha T_H \in W_H$.

We write $\Pi_H = \{ \beta_1, \dots, \beta_{\ell_H} \}$ for a simple system in $\Phi_H$, numbered as in \cite{Bou}. We let ${\Phi_H}^+$ and ${\Phi_H}^-$ be the corresponding sets of positive and negative roots in $\Phi_H$, and write $w_0$ for the long word of $W_H$, so that $w_0({\Phi_H}^+) = {\Phi_H}^-$. We let $U_H$ be the product of the root subgroups $X_\alpha$ corresponding to positive roots $\alpha$, and $B_H = U_H T_H$ be the standard Borel subgroup. We shall often represent the root $\sum m_i \beta_i$ as the $\ell_H$-tuple of coefficients $(m_1, \dots, m_{\ell_H})$ arranged as in a Dynkin diagram; thus for example if $H = E_8$ the highest root of $H$ is denoted $\eeightrt23465432$. Given $\alpha \in \Phi_H$, we write $\height(\alpha)$ for the height of $\alpha$.

We write $\L(H)$ for the Lie algebra of $H$; more generally, for a closed subgroup $H'$ of $H$ we write $\L(H')$ for the Lie algebra of $H'$, which we view as a subalgebra of $\L(H)$. For each $\alpha \in \Phi_H$ we take a root vector $e_\alpha$ in $\L(X_\alpha)$, and we write $h_\alpha$ for the vector $[e_\alpha, e_{-\alpha}]$ in $\L(T_H)$; if $\alpha \in {\Phi_H}^+$ we write $f_\alpha = e_{-\alpha}$. The structure constants of $H$ are defined by $[e_\alpha, e_\beta] = N_{\alpha \beta} e_{\alpha + \beta}$ whenever $\alpha, \beta, \alpha + \beta \in \Phi_H$; if $H = E_6$, $E_7$ or $E_8$ we shall take those given in the appendix of \cite{LSmax} unless otherwise stated.

Let $e(\Phi_H)$ be the maximum ratio of squared root lengths in $\Phi_H$, so that
$$
e(\Phi_H) = \begin{cases}
1 & \hbox{if } \Phi_H = A_\ell, D_\ell, E_6, E_7, E_8, \\
2 & \hbox{if } \Phi_H = B_\ell, C_\ell, F_4, \\
3 & \hbox{if } \Phi_H = G_2.
\end{cases}
$$
Note that if $e(\Phi_H) = 1$ we shall choose to regard all roots as short rather than long, which is not the usual convention. Given a subsystem $\Psi$ of $\Phi_H$, we write $\Psi_s$ and $\Psi_l$ respectively for the sets of short and long roots of $\Phi_H$ lying in $\Psi$.

Given $h, h' \in H$, we write $h^{h'} = {}^{{h'}^{-1}}h = {h'}^{-1}hh'$; if $h \in H$ and $A$ is a subset of $H$, we write $A^h = \{ a^h : a \in A \} = {}^{h^{-1}}A$ and $h^A = \{ h^a : a \in A \}$. Then $h^H$ is the conjugacy class of $h$ in $H$, and we write $C_H(h) = \{ h' \in H : h^{h'} = h \}$ for the centralizer of $h$ in $H$.

If $X$ is a variety on which $H$ acts, given $h \in H$ and $x \in X$ we write $h.x$ for the image of $x$ under the action of $h$, and $C_X(h) = \{ x \in X : h.x = x \}$ for the fixed point variety of $h$; given $A \leq H$ we write $A.x = \{ h.x : h \in A\}$ for the $A$-orbit containing $x$, and $C_A(x) = \{ h \in A : h.x = x \}$ for the $A$-stabilizer of $x$; given $Y \subseteq X$, we write $\overline{Y}$ for the closure of $Y$, and $\Tran_H(x, Y) = \{ h \in H : h.x \in Y \}$ for the transporter, which is closed in $H$ if $Y$ is closed in $X$.

Now let $G$ be a simple algebraic group over $K$, as in Section~\ref{sect: statement}. In the case where $H = G$ we shall mostly drop the subscript `$H$', so that the rank of $G$ is $\ell$ and we have the maximal torus $T$ with normalizer $N$, Weyl group $W$, root system $\Phi$, simple system $\Pi$, sets $\Phi^+$ and $\Phi^-$ of positive and negative roots, unipotent group $U$ and Borel subgroup $B$ (although we shall still have centralizers $C_G(h)$ and stabilizers $C_G(x)$); in addition we shall write $\Pi = \{ \alpha_1, \dots, \alpha_\ell \}$. The reason for beginning this section with $H$ rather than $G$ is that sometimes we will wish to view $G$ as a subgroup of a larger group $H$; then we may need to distinguish between maximal tori, Weyl groups, root systems and so on of the two groups.

We write $M = |\Phi| = \dim G - \ell$; the values of $M$ are as follows.
$$
\begin{array}{|c|c|c|c|c|}
\cline{1-2} \cline{4-5}
   G   &        M        & \ptw &  G  &  M  \tbs \\
\cline{1-2} \cline{4-5}
A_\ell &  \ell^2 + \ell  &      & E_6 &  72 \tbs \\
B_\ell &      2\ell^2    &      & E_7 & 126 \tbs \\
C_\ell &      2\ell^2    &      & E_8 & 240 \tbs \\
D_\ell & 2\ell^2 - 2\ell &      & F_4 &  48 \tbs \\
\cline{1-2}
\multicolumn{3}{c|}{}           & G_2 &  12 \tbs \\
\cline{4-5}
\end{array}
$$

Given $w \in W$, we write $U_w$ for the product of the root groups $X_\alpha$ for which $\alpha \in \Phi^+$ and $w(\alpha) \in \Phi^-$. The Bruhat decomposition gives each element of $G$ uniquely as $u_1 n u_2$, with $u_1 \in U$, $n \in N$ and $u_2 \in U_w$ where $w = nT \in W$. We write $G_u$ and $G_{ss}$ for the sets of unipotent and semisimple elements in $G$ respectively.

If $G$ is of classical type, we shall sometimes use the standard notation for its root system: we take an orthonormal basis $\ve_1, \dots, \ve_{\ell'}$ of $\ell'$-dimensional Euclidean space, where $\ell' = \ell + 1$ if $G = A_\ell$ and $\ell' = \ell$ if $G = B_\ell$, $C_\ell$ or $D_\ell$, and take simple roots $\alpha_i = \ve_i - \ve_{i + 1}$ for $i < \ell$ and $\alpha_\ell = \ve_\ell - \ve_{\ell + 1}$, $\ve_\ell$, $2\ve_\ell$ or $\ve_{\ell - 1} + \ve_\ell$ according as $G = A_\ell$, $B_\ell$, $C_\ell$ or $D_\ell$. Accordingly, we shall sometimes view the Weyl group $W$ as consisting of signed permutations of the set $\{ 1, \dots, \ell' \}$, where the number of minus signs is zero if $G = A_\ell$, arbitrary if $G = B_\ell$ or $C_\ell$, and even if $G = D_\ell$. In addition, we shall write $V_{nat}$ for the natural $G$-module; in what follows, we always take $1 \leq i < j \leq \ell'$, and when we describe the action of a root element $x_\alpha(t)$, any basis element whose image is not given explicitly is fixed. (See \cite[Theorem~11.3.2]{Car1}.) If $G = A_\ell$ then $V_{nat}$ has basis $v_1, \dots, v_{\ell + 1}$; root elements act by
$$
\begin{array}{rl}
x_{\ve_i - \ve_j}(t)  : & v_j \mapsto v_j + tv_i, \\
x_{-\ve_i + \ve_j}(t) : & v_i \mapsto v_i + tv_j.
\end{array}
$$
If $G = C_\ell$ then $V_{nat}$ has (hyperbolic) basis $e_1, f_1, \dots, e_\ell, f_\ell$; root elements act by
$$
\begin{array}{rl}
x_{\ve_i - \ve_j}(t)  : & e_j \mapsto e_j + te_i, \quad f_i \mapsto f_i - tf_j, \\
x_{-\ve_i + \ve_j}(t) : & e_i \mapsto e_i + te_j, \quad f_j \mapsto f_j - tf_i, \\
x_{\ve_i + \ve_j}(t)  : & f_j \mapsto f_j + te_i, \quad f_i \mapsto f_i + te_j, \\
x_{-\ve_i - \ve_j}(t) : & e_j \mapsto e_j + tf_i, \quad e_i \mapsto e_i + tf_j, \\
x_{2\ve_i}(t)        : & f_i \mapsto f_i + te_i, \\
x_{-2\ve_i}(t)       : & e_i \mapsto e_i + tf_i.
\end{array}
$$
If $G = D_\ell$ then $V_{nat}$ has basis $v_1, v_{-1}, \dots, v_\ell, v_{-\ell}$; root elements act by
$$
\begin{array}{rl}
x_{\ve_i - \ve_j}(t)  : & v_j \mapsto v_j + tv_i, \quad v_{-i} \mapsto v_{-i} - tv_{-j}, \\
x_{-\ve_i + \ve_j}(t) : & v_i \mapsto v_i + tv_j, \quad v_{-j} \mapsto v_{-j} - tv_{-i}, \\
x_{\ve_i + \ve_j}(t)  : & v_{-j} \mapsto v_{-j} + tv_i, \quad v_{-i} \mapsto v_{-i} - tv_j, \\
x_{-\ve_i - \ve_j}(t) : & v_j \mapsto v_j - tv_{-i}, \quad v_i \mapsto v_i + tv_{-j}.
\end{array}
$$
If $G = B_\ell$ then $V_{nat}$ has basis $v_0, v_1, v_{-1}, \dots, v_\ell, v_{-\ell}$; root elements act by
$$
\begin{array}{rl}
x_{\ve_i - \ve_j}(t)  : & v_j \mapsto v_j + tv_i, \quad v_{-i} \mapsto v_{-i} - tv_{-j}, \\
x_{-\ve_i + \ve_j}(t) : & v_i \mapsto v_i + tv_j, \quad v_{-j} \mapsto v_{-j} - tv_{-i}, \\
x_{\ve_i + \ve_j}(t)  : & v_{-j} \mapsto v_{-j} + tv_i, \quad v_{-i} \mapsto v_{-i} - tv_j, \\
x_{-\ve_i - \ve_j}(t) : & v_j \mapsto v_j - tv_{-i}, \quad v_i \mapsto v_i + tv_{-j}, \\
x_{\ve_i}(t)         : & v_0 \mapsto v_0 + 2tv_i, \quad v_{-i} \mapsto v_{-i} - tv_0 - t^2v_i, \\
x_{-\ve_i}(t)        : & v_0 \mapsto v_0 - 2tv_{-i}, \quad v_i \mapsto v_i + tv_0 - t^2v_{-i}.
\end{array}
$$

We write $\Lambda$ for the weight lattice of $G$ with respect to the maximal torus $T$, and let $\omega_1, \dots, \omega_\ell$ be the fundamental dominant weights of $G$ corresponding to the simple roots $\alpha_1, \dots, \alpha_\ell$ respectively. If $\lambda$ is a dominant weight of $G$, we write $L(\lambda)$ for the irreducible $G$-module with highest weight $\lambda$.

Given a $G$-module $V$, we write $\Lambda(V)$ for the set of weights in $\Lambda$ for which the weight space in $V$ is non-zero. If $\mu \in \Lambda(V)$, we write $V_\mu = \{ v \in V : \forall s \in T, \ s.v = \mu(s)v \}$ for the corresponding weight space. If $s \in G_{ss}$ and $\kappa \in K^*$, we write $V_\kappa(s) = \{ v \in V : s.v = \kappa v \}$ for the corresponding eigenspace, which is a sum of weight spaces. Given $\mu \in \Lambda(V)$, if $v \in V$ is such that the projection of $v$ on the weight space $V_\mu$ is non-zero, we say that the weight $\mu$ {\em occurs \/} in $v$.

Given $a \in \N$, we write $\Z_a$ for the cyclic group of order $a$, $Dih_{2a}$ for the dihedral group of order $2a$, and $S_a$ and $Alt_a$ for the symmetric and alternating groups of degree $a$; in addition we write $U_a$ for a connected unipotent group of dimension $a$ and $T_a$ for a torus of dimension $a$.

Given $a \in \N$ and $\kappa_1, \dots, \kappa_a \in K$, we write $\diag(\kappa_1, \dots, \kappa_a)$ for the diagonal $a \times a$ matrix whose $(i, i)$-entry is $\kappa_i$.

Given $a, b \in \N$, we write $(a, b)$ for their highest common factor, and we set
$$
\z_{a, b} =
\begin{cases}
1 & \hbox{if } a \hbox{ divides } b, \\
0 & \hbox{otherwise}.
\end{cases}
$$
If $p = \infty$ we extend this notation to cover the case where $a = p$ by setting $(p, b) = 1$ and $\z_{p, b} = 0$.

We write $\P$ for the set of primes in $\N$, and $\P'$ for $\P \setminus \{ p \}$.

Finally, given $r \in \N$, we let $\eta_r$ be a generator of the group of $r$th roots of unity in $K^*$ (so that if $(p, r) = 1$ then $\eta_r$ is a primitive $r$th root of unity); we assume this is done in such a way that whenever $r = r_1 r_2$ with $r_1, r_2 \in \N$ we have ${\eta_r}^{r_1} = \eta_{r_2}$.

\section{Weights and module structure}\label{sect: weights and module structure}

In this section we recall some basic facts about weights and modules. We start by considering weights in the abstract, and later give results linking this to the structure of modules.

Given a root system $\Phi$, its $\R$-span is a Euclidean space with an inner product $(-, -)$, on which the Weyl group $W$ acts as linear isometries via $w_\alpha(\mu) = \mu - \langle \mu, \alpha \rangle \alpha$, where we write $\langle \mu, \alpha \rangle = \frac{2(\mu, \alpha)}{(\alpha, \alpha)}$. The weights are the elements $\mu$ of this Euclidean space such that for all $\alpha \in \Phi$ we have $\langle \mu, \alpha \rangle \in \Z$; the set $\Lambda$ of weights is preserved by $W$, and thus is a union of $W$-orbits. We have a partial order $\preceq$ on $\Lambda$, where $\mu \preceq \lambda$ if and only if $\lambda - \mu$ is a sum of (zero or more) simple roots.

A weight $\lambda$ is dominant if for all $\alpha \in \Pi$ we have $\langle \lambda, \alpha \rangle \geq 0$; we write $\Lambda^+$ for the set of dominant weights, and then each $W$-orbit on $\Lambda$ contains a unique element of $\Lambda^+$. If $\lambda \in \Lambda^+$ then for all $w \in W$ we have $w(\lambda) \preceq \lambda$, as we may see by induction on the length of $w$: given $w \neq 1$ we may choose $\alpha \in \Pi$ with $w(\alpha) \in \Phi^-$, and then by \cite[Proposition~2.2.8]{Car2} we have $w = w' w_\alpha$ where $w'$ is shorter than $w$ and $w'(\alpha) \in \Phi^+$; as $\langle \lambda, \alpha \rangle \geq 0$, and by induction $w'(\lambda) \preceq \lambda$, we have
$$
w(\lambda) = w' w_\alpha(\lambda) = w'(\lambda - \langle \lambda, \alpha \rangle \alpha) = w'(\lambda) - \langle \lambda, \alpha \rangle w'(\alpha) \preceq w'(\lambda) \preceq \lambda.
$$
Note that for dominant weights the partial order is compatible with length in the Euclidean space: if $\lambda, \mu \in \Lambda^+$ with $\mu \prec \lambda$, then as $\mu \in \Lambda^+$ and $\lambda - \mu$ is a sum of simple roots, we have $(\mu, \lambda - \mu) \geq 0$, while as $\lambda - \mu \neq 0$ we have $(\lambda - \mu, \lambda - \mu) > 0$; thus $(\lambda, \lambda) - (\mu, \mu) = (\lambda + \mu, \lambda - \mu) = (\lambda - \mu, \lambda - \mu) + 2(\mu, \lambda - \mu) > 0$, so $(\mu, \mu) < (\lambda, \lambda)$.

Given $\lambda \in \Lambda^+$, set
$$
\Slambda = \{ w.\mu : w \in W, \ \mu \in \Lambda^+, \ \mu \preceq \lambda \};
$$
then $\Slambda$ is a union of $W$-orbits in $\Lambda$, and as all weights in $\Slambda$ have length at most that of $\lambda$ the set $\Slambda$ is finite. For any root $\alpha \in \Phi$ there is an equivalence relation on $\Slambda$ whereby two weights are related if and only if their difference is a multiple of $\alpha$; the equivalence classes are called {\em $\alpha$-strings\/}. Given a weight $\nu \in \Slambda$, the $\alpha$-string containing $\nu$ consists of weights of the form $\nu - t\alpha$ for $t \in \Z$. If $t_1$ and $t_2$ respectively are the maximal and minimal values of $t$ with $\nu - t\alpha \in \Slambda$ (so that $t_1 \geq 0 \geq t_2$), we may regard the $\alpha$-string as
$$
\nu - t_1\alpha \quad \dots \quad \nu \quad \dots \quad \nu - t_2\alpha
$$
where we arrange the weights in order of decreasing $t$; the reflection $w_\alpha$ acts on the $\alpha$-string by reversing the order, and we have $\langle \nu, \alpha \rangle = t_1 + t_2$. We claim that for all $t$ with $t_1 > t > t_2$ we have $\nu - t\alpha \in \Slambda$ (so that the $\alpha$-string has \lq no gaps'): by applying an appropriate element of $W$ and replacing the root $\alpha$ by its image under this element, we may assume that $\nu - t\alpha \in \Lambda^+$; by negating $\alpha$ if necessary we may assume that $\alpha \in \Phi^+$, in which case we have $\nu - t\alpha \prec \nu - t_2\alpha$; since $\nu - t_2\alpha \in \Slambda$, the dominant weight $\mu$ in its $W$-orbit satisfies $\mu \preceq \lambda$, so $\nu - t\alpha \prec \nu - t_2\alpha \preceq \mu \preceq \lambda$ as required. As a consequence we see that the set $\Slambda$ is {\em saturated\/}, meaning that for all $\nu$ in the set, all $\alpha \in \Phi$ and any $t$ between $0$ and $\langle \nu, \alpha \rangle$, the weight $\nu - t\alpha$ lies in the set.

Note that if $\langle \nu, \alpha \rangle > 1$ (so that $w_\alpha(\nu)$ is to the left of $\nu$ in the $\alpha$-string as displayed above, and there are weights lying between $w_\alpha(\nu)$ and $\nu$), then $2(\nu, \alpha) > (\alpha, \alpha)$, and so $(\nu, \nu) - (\nu - \alpha, \nu - \alpha) = 2(\nu, \alpha) - (\alpha, \alpha) > 0$, whence $(\nu - \alpha, \nu - \alpha) < (\nu, \nu)$; thus length decreases as one moves towards the centre of an $\alpha$-string.

We now consider how this relates to the structure of irreducible $G$-modules, where $\Phi$ is the root system of $G$ with respect to the maximal torus $T$. Any such module $V$ decomposes as a direct sum of weight spaces $V_\nu$ for $T$; the weights in $\Lambda(V)$ comprise a union of $W$-orbits, and for all $\nu \in \Lambda(V)$ and $w \in W$ we have $w V_\nu = V_{w \nu}$ so that $\dim V_{w \nu} = \dim V_\nu$. We shall state three theorems of central importance, in each case giving the form in which it appears in \cite{Lubpaper} (but using the notation employed here); for convenience, in each we shall assume $G$ is of simply connected type (which as mentioned at the start of Section~\ref{sect: statement} is harmless for our purposes). The first (\cite[Theorem~2.1]{Lubpaper}) is due to Chevalley, and establishes the link between irreducible $G$-modules and dominant weights.

\begin{thm}\label{thm: Chevalley}
Let $G$ be a simply connected simple algebraic group over $K$ and $V$ be a finite-dimensional irreducible $G$-module.
\begin{itemize}
\item[(i)] The set $\Lambda(V)$ contains a (unique) element $\lambda$ such that for all $\mu \in \Lambda(V)$ we have $\mu \preceq \lambda$. This $\lambda$ is called the {\em highest weight\/} of $V$, it is dominant, and we have $\dim V_\lambda = 1$.
\item[(ii)] The module $V$ is determined up to isomorphism by its highest weight.
\item[(iii)] For each $\lambda \in \Lambda^+$ there is an irreducible $G$-module $L(\lambda)$ with highest weight $\lambda$.
\end{itemize}
\end{thm}

Thus if $V = L(\lambda)$ then
$$
\Lambda(V) \subseteq \Slambda.
$$
In characteristic zero it is known that in fact we have equality (see for example \cite[Proposition~21.3]{HumLie}, where this is stated for the action of $\L(G)$, which has the same modules as $G$). There are, however, cases in positive characteristic where $\Lambda(V)$ is a proper subset of $\Slambda$.

The fundamental dominant weights $\omega_1, \dots, \omega_\ell$ are defined by $\langle \omega_i, \alpha_j \rangle = \delta_{ij}$; expressions giving the $\omega_i$ as rational linear combinations of the $\alpha_j$ appear in \cite[13.1, Table~1]{HumLie}. Any dominant weight is a sum of (zero or more) fundamental dominant weights. A dominant weight $\lambda = \sum_{i = 1}^\ell a_i \omega_i$ is called {\em $p$-restricted\/} if for all $i \leq \ell$ we have $0 \leq a_i < p$; thus if $p = \infty$ then all dominant weights are $p$-restricted. The second of our results (\cite[Theorem~2.2]{Lubpaper}) is Steinberg's tensor product theorem, which shows how an arbitrary irreducible module for $G$ is built out of ones with $p$-restricted highest weights.

\begin{thm}\label{thm: Steinberg}
Let $G$ be a simply connected simple algebraic group over $K$, and suppose $p < \infty$; write $F$ for the $p$-power Frobenius automorphism of $K$. Given a $G$-module $V$, for each $i \geq 0$ we denote by $V^{(i)}$ the $G$-module obtained by twisting the action of $G$ on $V$ by $F^i$. If $n \geq 0$ and $\lambda_0, \lambda_1, \dots, \lambda_n$ are $p$-restricted dominant weights, then
$$
L(\lambda_0 + p\lambda_1 + \cdots + p^n\lambda_n) \cong L(\lambda_0) \otimes L(\lambda_1)^{(1)} \otimes \cdots \otimes L(\lambda_n)^{(n)}.
$$
\end{thm}

(It is this result which allows us to assume that $\lambda$ is not a multiple of $p$, since otherwise the module is simply obtained by twisting.)

Recall that we define $e(\Phi)$ to be the maximum ratio of squared root lengths in $\Phi$. The third of our theorems (\cite[Theorem~4.1]{Lubpaper}) is due to Premet \cite{Prem}, and provides a condition guaranteeing equality in the containment above.

\begin{thm}\label{thm: Prem}
Let $G$ be a simply connected simple algebraic group over $K$ with root system $\Phi$, and $\lambda$ be a $p$-restricted dominant weight; write $V = L(\lambda)$. If $p > e(\Phi)$ then $\Lambda(V) = \Slambda$.
\end{thm}

However, given $V = L(\lambda)$, unless the set of dominant weights $\mu \prec \lambda$ is empty, knowledge of the set $\Lambda(V)$ alone is insufficient to determine the dimension of $V$, or the multiplicities of weights $\mu$, i.e., the dimensions of the weight spaces $V_\mu$. In characteristic zero there are formul\ae \ due to Weyl and Kostant which allow these to be computed (see for example \cite[24.2, 24.3]{HumLie}); however in positive characteristic no such formul\ae \ are known. This problem is addressed by L\"ubeck in \cite{Lubpaper}; he works with the Kostant $\Z$-form of the universal enveloping algebra of the complex Lie algebra corresponding to $G$, and we summarise his approach as follows. Let $\gamma_1, \dots, \gamma_t$ be a fixed ordering on the set $\Phi^+$. Given a $t$-tuple $\a = (a_1, \dots, a_t)$ of non-negative integers, write $f_\a = {f_{\gamma_t}}^{a_t} \dots {f_{\gamma_1}}^{a_1}$ and $e_\a = {e_{\gamma_t}}^{a_t} \dots {e_{\gamma_1}}^{a_1}$; then applying $f_\a$ or $e_\a$ to a vector in $V$ of weight $\omega$ gives a vector of weight $\omega - \sum a_i \gamma_i$ or $\omega + \sum a_i \gamma_i$ respectively. Let $v \in V$ be a vector of weight $\lambda$. Given a weight $\mu \in \Lambda(V)$, consider the set $S_\mu$ of all $\a$ such that $\lambda - \sum a_i \gamma_i = \mu$; if $\a, \b \in S_\mu$ then $e_\b f_\a v$ is again of weight $\lambda$, so there exists $n_{\a, \b} \in \Z$ such that $e_\b f_\a v = n_{\a, \b} v$. Letting $\a$ and $\b$ run through $S_\mu$ we obtain a matrix $(n_{\a, \b})$; the rank of the reduction modulo $p$ of this matrix equals $\dim V_\mu$.

Using this approach, L\"ubeck has in many cases determined all weight multiplicities in the module $V$: he treats classical root systems of bounded rank ($\ell \leq 20$ for type $A_\ell$ and $\ell \leq 11$ for other types) and exceptional root systems; for each root system he treats all modules of dimension less than some bound, in arbitrary characteristic. His results are recorded in \cite{Lubdata}; we shall make extensive use of this information. In the final result in this section we use this method to determine at least partially the structure of some particular modules for classical groups of arbitrary rank, which therefore are not given in \cite{Lubdata}.

\begin{lem}\label{lem: multiplicities in 3omega_1 and omega_1 + omega_2}
Let $G = A_\ell$ for $\ell \in [3, \infty)$, or $G = B_\ell$ or $C_\ell$ for $\ell \in [4, \infty)$; let $\lambda = 3\omega_1$ with $p \geq 5$, or $\lambda = \omega_1 + \omega_2$, and write $V = L(\lambda)$. Then
\begin{itemize}
\item[(i)] if $\lambda = 3\omega_1$, then for $\mu = \omega_1 + \omega_2$ or $\omega_3$ we have $\dim V_\mu = 1$;
\item[(ii)] if $\lambda = \omega_1 + \omega_2$, then for $\mu = \omega_3$ we have $\dim V_\mu = 2 - \z_{p, 3}$.
\end{itemize}
\end{lem}

\begin{proof}
We have $3\omega_1 - (\omega_1 + \omega_2) = \alpha_1$ and $(\omega_1 + \omega_2) - \omega_3 = \alpha_1 + \alpha_2$. We take an ordering on $\Phi^+$ such that $\gamma_1 = \alpha_2$, $\gamma_2 = \alpha_1$ and $\gamma_3 = \alpha_1 + \alpha_2$; we may assume that $[e_{\alpha_1}, e_{\alpha_2}] = e_{\alpha_1 + \alpha_2}$. Since $\lambda - \mu$ is a linear combination of $\alpha_1$ and $\alpha_2$ alone, each $t$-tuple $\a$ in $S_\mu$ has $a_i = 0$ for all $i > 3$; thus for convenience we may abbreviate $\a$ to simply $(a_1, a_2, a_3)$. Recall that for $x$ and $y$ in the Kostant $\Z$-form we have $xy = [x, y] + yx$.

First suppose $\lambda = 3\omega_1$ with $p \geq 5$, so that $h_{\alpha_1} v = 3v$ and $h_{\alpha_2} v = 0$. Take $\mu = \omega_1 + \omega_2$; then $S_\mu = \{ (0, 1, 0) \}$. We have
\begin{eqnarray*}
e_{\alpha_1} f_{\alpha_1} v & = & h_{\alpha_1} v + f_{\alpha_1} e_{\alpha_1} v = 3v + 0 = 3v.
\end{eqnarray*}
Thus the matrix $(n_{\a, \b})$ is simply $( 3 )$, whose rank is $1$. Now instead take $\mu = \omega_3$; then $S_\mu = \{ (1, 2, 0), (0, 1, 1) \}$. Since $f_{\alpha_2} v = 0$, we have $e_\b f_\a v = 0$ if $\a = (1, 2, 0)$, while
\begin{eqnarray*}
e_{\alpha_1 + \alpha_2} e_{\alpha_1} f_{\alpha_1 + \alpha_2} f_{\alpha_1} v
& = & -e_{\alpha_1 + \alpha_2} f_{\alpha_2} f_{\alpha_1} v + e_{\alpha_1 + \alpha_2} f_{\alpha_1 + \alpha_2} e_{\alpha_1} f_{\alpha_1} v \\
& = & - e_{\alpha_1} f_{\alpha_1} v - f_{\alpha_2} e_{\alpha_1 + \alpha_2} f_{\alpha_1} v + 3e_{\alpha_1 + \alpha_2} f_{\alpha_1 + \alpha_2} v \\
& = & -3v + 0 + 3(h_{\alpha_1} + h_{\alpha_2})v \\
& = & 6v.
\end{eqnarray*}
Thus the matrix $(n_{\a, \b})$ again has rank $1$. This proves (i).

Now suppose $\lambda = \omega_1 + \omega_2$, so that $h_{\alpha_1} v = h_{\alpha_2} v = v$. Take $\mu = \omega_3$; then $S_\mu = \{ (1, 1, 0), (0, 0, 1) \}$. We have
\begin{eqnarray*}
e_{\alpha_1} e_{\alpha_2} f_{\alpha_1} f_{\alpha_2} v
& = & e_{\alpha_1} f_{\alpha_1} e_{\alpha_2} f_{\alpha_2} v = e_{\alpha_1} f_{\alpha_1} h_{\alpha_2} v = e_{\alpha_1} f_{\alpha_1} v = h_{\alpha_1} v = v, \\
e_{\alpha_1 + \alpha_2} f_{\alpha_1} f_{\alpha_2} v
& = & - e_{\alpha_2} f_{\alpha_2} v + f_{\alpha_1} e_{\alpha_1 + \alpha_2} f_{\alpha_2} v = - h_{\alpha_2} v + 0 = -v, \\
e_{\alpha_1} e_{\alpha_2} f_{\alpha_1 + \alpha_2} v
& = & e_{\alpha_1} f_{\alpha_1} v + e_{\alpha_1} f_{\alpha_1 + \alpha_2} e_{\alpha_2} v = h_{\alpha_1} v + 0 = v, \\
e_{\alpha_1 + \alpha_2} f_{\alpha_1 + \alpha_2} v
& = & (h_{\alpha_1} + h_{\alpha_2}) v = v + v = 2v.
\end{eqnarray*}
Thus the matrix $(n_{\a, \b})$ is
$$
\left(
  \begin{array}{cc}
    1 & -1 \\
    1 & 2 \\
  \end{array}
\right),
$$
whose rank is $1$ if $p = 3$ and $2$ otherwise. This proves (ii).
\end{proof}

Note that if $G = A_\ell$ then $\omega_1 + \omega_2$ and $\omega_3$ are the only dominant weights lying in ${\mathcal S}(3\omega_1)$, so that the dimensions of all weight spaces in $V$ are determined; if however $G = B_\ell$ or $C_\ell$ then there are other dominant weights present, the dimensions of whose weight spaces have not been determined here.

\section{Unipotent classes}\label{sect: unipotent classes}

In this section we review some basic information about unipotent conjugacy classes of $G$, and provide some results about their dimensions and closures for later use. Our interest will be confined to classes containing elements of order $p$, as Section~\ref{sect: conditions} will make clear.

We begin with notation. In good characteristic we have the classification of Bala and Carter (extended by Pommerening), giving a bijective correspondence between unipotent classes of $G$ and conjugacy classes of pairs consisting of a Levi subgroup of $G$ and a distinguished parabolic subgroup of its semisimple part; the Bala-Carter notation labels each unipotent class by the corresponding distinguished parabolic subgroup. This notation may be extended to bad characteristic, provided additional unipotent classes are included. For $G$ of exceptional type in bad characteristic, the number of these additional classes is given in \cite[5.11]{Car2}; from \cite{LawJor} we see that such additional classes contain elements of order $p$ only if $(G, p) = (F_4, 2)$ or $(G_2, 3)$, when a single such class exists, denoted there by $\tilde A_1{}^{(p)}$. For $G$ of classical type in bad characteristic, the situation is more complicated. We shall briefly describe the classification given by Liebeck and Seitz in \cite{LSbook}; we shall also indicate the notation used by Aschbacher and Seitz in \cite{AS}, which has class representatives $a_t$ and $c_t$ for $t$ even and $b_t$ for $t$ odd.

Let $G = B_\ell$, $C_\ell$ or $D_\ell$ with $p = 2$. A class of elements of order $2$ is represented in \cite{LSbook} by an expression of the form
$$
W(1)^{a_1} + W(2)^{a_2} + V(2)^b,
$$
together with a final summand $R$ if $G = B_\ell$; here $a_1$, $a_2$ and $b$ are non-negative integers satisfying $a_1 + 2a_2 + b = \ell$, with $a_2 + b > 0$, such that $b \leq 2$, and if $G = D_\ell$ then $b$ is even. This expression gives the action on the natural module of a class representative, with $W(m)$ representing a pair of Jordan blocks of size $m$ and $V(m)$ a single Jordan block of size $m$ (and $R$ denoting the $1$-dimensional radical of the associated form if $G = B_\ell$). For some of these classes we shall use the following notation (ignoring the radical $R$ if it is present): the class corresponding to
$$
W(1)^{\ell - 2y} + W(2)^y
$$
(with representative $a_{2y}$ in \cite{AS}) will be called ${A_1}^y$ (unless $G = D_\ell$ and $y = \frac{1}{2}\ell$, in which case there are two classes $({A_1}^{\ell/2})'$ and $({A_1}^{\ell/2})''$, which are interchanged by a graph automorphism); the class corresponding to
$$
W(1)^{\ell - 2y - 1} + W(2)^y + V(2)
$$
(with representative $b_{2y + 1}$ in \cite{AS}) will be called ${A_1}^yB_1$ if $G = B_\ell$ and ${A_1}^yC_1$ if $G = C_\ell$; the class corresponding to
$$
W(1)^{\ell - 2y - 2} + W(2)^y + V(2)^2
$$
(with representative $c_{2y + 2}$ in \cite{AS}) will be called ${A_1}^yD_2$ if $G = D_\ell$ (but no notation is specified if $G = B_\ell$ or $C_\ell$). In each case the class labelled with a subsystem $\Phi'$ then contains elements regular in a subsystem subgroup of type $\Phi'$, so we recover the Bala-Carter notation for these classes.

For the remainder of this section we shall concentrate on dimensions and closures of unipotent classes. It is well known (see for example \cite[Theorem~4.2]{HumConj}) that the set $G_u$ of unipotent elements of $G$ is an irreducible closed subset of $G$, of dimension $M = |\Phi| = \dim G - \ell$; thus the closure of any unipotent class is a union of unipotent classes. Moreover by \cite[Proposition~8.3]{HumLAG} the boundary of any conjugacy class in $G$ is a union of classes of smaller dimension. Given unipotent classes ${u_1}^G$ and ${u_2}^G$, we write ${u_1}^G \leq {u_2}^G$ if ${u_1}^G \subseteq \overline{{u_2}^G}$; this gives a partial order on the set of unipotent classes of $G$. The reason for our interest in this partial order rests in the following elementary result.

\begin{lem}\label{lem: unip closure containment}
Let $V$ be a $G$-module, and take $u_1, u_2 \in G_u$ with ${u_1}^G \leq {u_2}^G$; then $\codim C_V(u_1) \leq \codim C_V(u_2)$, and if $1 \leq k \leq \frac{1}{2} \dim V$ then $\codim C_{\Gk(V)}(u_1) \leq \codim C_{\Gk(V)}(u_2)$.
\end{lem}

\begin{proof}
Let $X$ be either $V$ or $\Gk(V)$. The set
$$
\{ g \in \GL(V) : \codim C_X(g) \leq \codim C_X(u_2) \}
$$
is closed and contains ${u_2}^G$, so it contains the closure $\overline{{u_2}^G}$ and hence ${u_1}^G$; the result follows.
\end{proof}

As we shall see in Section~\ref{sect: conditions}, our method of establishing that triples and quadruples have TGS will employ conditions involving codimensions of fixed point spaces; Lemma~\ref{lem: unip closure containment} will be used frequently to limit the number of classes requiring consideration.

Here we shall first provide a brief overview of (some of) the known material concerning dimensions and closures of unipotent classes; we shall then give a number of results to be used in the work ahead.

We first consider $G$ of classical type in good characteristic. Here unipotent classes are almost entirely determined by Jordan structure, which corresponds to partitions of $n$, where $n$ is the dimension of the natural $G$-module $L(\omega_1)$ (so that $n = \ell + 1$, $2\ell + 1$, $2\ell$, or $2\ell$ according as $G = A_\ell$, $B_\ell$, $C_\ell$, or $D_\ell$); within a partition of $n$ we take the parts in decreasing order, i.e., in the partition $[a_1, a_2, \dots]$ we assume $a_1 \geq a_2 \geq \cdots$, and we shall use superscripts to indicate repeated parts. If $G = A_\ell$ there is no restriction on the partitions which occur, while if $G = C_\ell$ (respectively $G = B_\ell$ or $D_\ell$) then all odd (respectively even) parts of the partition must occur with even multiplicity. If $G = D_\ell$ and all parts of the partition are even then there are two such unipotent classes; in all other cases there is a single unipotent class corresponding to the partition. Given a unipotent class $u^G$, we shall denote the corresponding partition of $n$ by $\Part(u^G)$. The partial order on unipotent classes is given by the dominance order on partitions of $n$, whereby $[a_1, a_2, \dots] \leq [b_1, b_2, \dots]$ if and only if for all $i$ we have $a_1 + \cdots + a_i \leq b_1 + \cdots + b_i$. Given a class $u^G \in G_u$ corresponding to a partition in which the number of parts equal to $i$ is $r_i$, we have
$$
\dim u^G = \begin{cases}
(\ell + 1)^2 - \sum_i (r_i + r_{i + 1} + \cdots )^2                                                   & \hbox{if } G = A_\ell, \\
2\ell^2 + \ell - {\ts\frac{1}{2}} (\sum_i (r_i + r_{i + 1} + \cdots )^2 - \sum_{i\ \mathrm{odd}} r_i) & \hbox{if } G = B_\ell, \\
2\ell^2 + \ell - {\ts\frac{1}{2}} (\sum_i (r_i + r_{i + 1} + \cdots )^2 + \sum_{i\ \mathrm{odd}} r_i) & \hbox{if } G = C_\ell, \\
2\ell^2 - \ell - {\ts\frac{1}{2}} (\sum_i (r_i + r_{i + 1} + \cdots )^2 - \sum_{i\ \mathrm{odd}} r_i) & \hbox{if } G = D_\ell.
\end{cases}
$$
(All of this is well known; see for example \cite[13.1]{Car2} and \cite[I.2.4, I.2.5]{Sp}.)

Next we consider $G$ of classical type in bad characteristic; so $G = B_\ell$, $C_\ell$ or $D_\ell$ and $p = 2$. If we employ the Aschbacher-Seitz notation for elements, we may use \cite[Theorem~4.2]{LSbook} to see that the class dimensions are as follows:
$$
\begin{array}{cccc}
 u  &  \dim u^{B_\ell} &  \dim u^{C_\ell} &  \dim u^{D_\ell} \tbs \\
a_t &   t(2\ell - t)   &   t(2\ell - t)   & t(2\ell - 1 - t) \tbs \\
b_t & t(2\ell + 1 - t) & t(2\ell + 1 - t) &         -        \tbs \\
c_t & t(2\ell + 1 - t) & t(2\ell + 1 - t) &   t(2\ell - t)   \tbs \\
\end{array}
$$
As for the partial order, it is clear from the description above that for $x \in \{ a, b, c \}$, if $t' \leq t$ then ${x_{t'}}^G \leq {x_t}^G$, and that for $y \geq 0$ we have ${a_{2y}}^G \leq {b_{2y + 1}}^G \leq {c_{2y + 2}}^G$. We also have the following.

\begin{lem}\label{lem: a_{2y} in closure of c_{2y}}
If $G = B_\ell$, $C_\ell$ or $D_\ell$ with $p = 2$, then ${a_{2y}}^G \leq {c_{2y}}^G$ for $y \geq 1$.
\end{lem}

\begin{proof}
In the notation of \cite{LSbook}, it suffices to show that, on a $4$-dimensional space, the closure of the class containing elements acting as $V(2)^2$ contains elements acting as $W(2)$. For an element in the former class, from \cite[6.1]{LSbook} there is a basis ${v_1}^{(1)}, {v_2}^{(1)}, {v_1}^{(2)}, {v_2}^{(2)}$, such that the bilinear form satisfies $({v_i}^{(j)}, {v_{3 - i}}^{(j)}) = 1$ with the value taken at other pairs of basis vectors being zero, and if $G = B_\ell$ or $D_\ell$ the quadratic form $Q$ satisfies $Q({v_1}^{(j)}) = 1$ and $Q({v_2}^{(j)}) = 0$ for $j = 1, 2$; the element acts by fixing each ${v_1}^{(j)}$ and sending each ${v_2}^{(j)}$ to ${v_1}^{(j)} + {v_2}^{(j)}$. For $\kappa \in K^*$ write
\begin{eqnarray*}
x_{-1} & = & \kappa^{-1} ({v_1}^{(1)} + {v_2}^{(1)} + {v_2}^{(2)}), \\
   x_1 & = & \kappa^{-1} ({v_1}^{(1)} + {v_1}^{(2)}), \\
y_{-1} & = & \kappa {v_2}^{(2)}, \\
   y_1 & = & \kappa ({v_2}^{(1)} + {v_2}^{(2)});
\end{eqnarray*}
then the bilinear form satisfies $(x_i, y_{-i}) = 1$ with the value taken at other pairs of basis vectors being zero, and if $G = B_\ell$ or $D_\ell$ we find that the quadratic form $Q$ satisfies $Q(x_i) = Q(y_i) = 0$ for $i = \pm1$. Moreover the element acts as
\begin{eqnarray*}
x_{-1} & \mapsto & x_{-1} + x_1, \\
   x_1 & \mapsto & x_1, \\
y_{-1} & \mapsto & \kappa^2 x_{-1} + \kappa^2 x_1 + y_{-1} + y_1, \\
   y_1 & \mapsto & \kappa^2 x_1 + y_1;
\end{eqnarray*}
thus the closure of the class contains the element obtained from this by setting $\kappa = 0$, which acts by fixing $x_1$ and $y_1$, and sending $x_{-1}$ to $x_{-1} + x_1$ and $y_{-1}$ to $y_{-1} + y_1$. As this is exactly the description in \cite[6.1]{LSbook} of the action of an element in the latter class, the result follows.
\end{proof}

Finally we consider $G$ of exceptional type. Here the unipotent classes were originally determined by Chang in \cite{Chang} and Enomoto in \cite{Eno} for $G = G_2$, by Shinoda in \cite{Shin} and Shoji in \cite{Sho} for $G = F_4$, and by Mizuno in \cite{Miz1, Miz2} for $G = E_6$, $E_7$ and $E_8$. In \cite[II.10.4 and IV.2]{Sp} Spaltenstein provides diagrams specifying the partial ordering on unipotent classes, in all characteristics (including the additional unipotent classes occurring in bad characteristic); in the case of $E_6$, $E_7$ and $E_8$ his diagrams reproduce those appearing in \cite{Miz2}. In \cite[13.4]{Car2} Carter repeats all these diagrams in the case of characteristic zero (but using the Bala-Carter notation for classes, which neither Mizuno nor Spaltenstein employed); in \cite[13.1]{Car2} he also lists centralizer dimensions, again for characteristic zero. These dimensions are listed in all characteristics by Liebeck and Seitz in \cite[Tables~22.1.1--22.1.5]{LSbook}. Thus between them \cite{Car2}, \cite{LSbook} and \cite{Sp} give all the information we require on dimensions and closures; indeed in \cite[IV.2]{Sp} Spaltenstein also gives diagrams for some classical groups of small rank.

We now move on to the results we shall wish to use in the work here, which for unipotent classes ${u_1}^G$ and ${u_2}^G$ give conditions implying that ${u_1}^G \leq {u_2}^G$; in a few cases we allow ${u_1}^G$ to be one of two possibilities. We begin with some very general conditions, and then move on to ones which are more specific. All classes which we treat will be assumed to contain elements of order $p$; in some cases this gives a lower bound on the value of $p$, which we will not always mention.

\begin{lem}\label{lem: any class in closure of reg class}
We have ${u_1}^G \leq {u_2}^G$ if ${u_1}^G$ is any unipotent class and ${u_2}^G$ is the regular unipotent class.
\end{lem}

\begin{proof}
Since $\dim {u_2}^G = \dim G - \ell = \dim G_u$, and $G_u$ is irreducible and closed, we have $\overline{{u_2}^G} = G_u$ and hence ${u_1}^G \subseteq \overline{{u_2}^G}$.
\end{proof}

For the next result recall that $e(\Phi)$ is the maximum ratio of squared root lengths in the root system $\Phi$ of $G$.

\begin{lem}\label{lem: root elt class in closure of any non-triv class}
We have ${u_1}^G \leq {u_2}^G$ if ${u_1}^G$ contains root elements and ${u_2}^G$ is any non-trivial unipotent class, unless $e(\Phi) > 1$ and one of the following holds:
\begin{itemize}
\item[(i)] ${u_1}^G$ contains long root elements, $(G, p) = (B_\ell, 2)$, $(F_4, 2)$ or $(G_2, 3)$ and ${u_2}^G$ contains short root elements;
\item[(ii)] ${u_1}^G$ contains long root elements, $(G, p) = (C_\ell, 2)$ and ${u_2}^G = {A_1}^y$ for some $y \geq 1$;
\item[(iii)] ${u_1}^G$ contains short root elements, $G = C_\ell$, $F_4$ or $G_2$ and ${u_2}^G$ contains long root elements;
\item[(iv)] ${u_1}^G$ contains short root elements, $G = B_\ell$ and ${u_2}^G = {A_1}^y$ for some $y \geq 1$.
\end{itemize}
\end{lem}

\begin{proof}
For $G$ exceptional the result is clear from the diagrams in \cite{Sp}. For $G = A_\ell$ or $D_\ell$ the class $A_1$ of root elements has $\Part(A_1) = 2 1^{\ell - 1}$ or $2^2 1^{2\ell - 4}$ respectively, so $\Part(A_1) \leq \Part({u_2}^G)$ for any non-trivial unipotent class ${u_2}^G$ (since for $G = D_\ell$ even parts must occur with even multiplicity). Thus we may assume $G = B_\ell$ or $C_\ell$. First suppose $p \neq 2$. If $G = B_\ell$, the classes $A_1$ and $B_1$ of long and short root elements have $\Part(A_1) = 2^2 1^{2\ell - 3}$ and $\Part(B_1) = 3 1^{2\ell - 2}$, so for any non-trivial unipotent class ${u_2}^G$ we have $\Part(A_1) \leq \Part({u_2}^G)$ as for $G = D_\ell$, while $\Part(B_1) \leq \Part({u_2}^G)$ unless $\Part({u_2}^G) = 2^{2y} 1^{2\ell + 1 - 4y}$ for some $y \geq 1$, when ${u_2}^G = {A_1}^y$. If instead $G = C_\ell$, the classes $C_1$ and $A_1$ of long and short root elements have $\Part(C_1) = 2 1^{2\ell - 2}$ and $\Part(A_1) = 2^2 1^{2\ell - 4}$, so for any non-trivial unipotent class ${u_2}^G$ we have $\Part(C_1) \leq \Part({u_2}^G)$, while $\Part(A_1) \leq \Part({u_2}^G)$ unless ${u_2}^G = C_1$. Now suppose $p = 2$. For all $y \geq 0$ we have ${b_1}^G \leq {b_{2y + 1}}^G$ and ${b_1}^G \leq {c_{2y + 2}}^G$, and for all $y \geq 1$ we have ${a_2}^G \leq {a_{2y}}^G \leq {b_{2y + 1}}^G$ and ${a_2}^G \leq {c_{2y}}^G$; thus for any non-trivial unipotent class ${u_2}^G$ we have ${b_1}^G \leq {u_2}^G$ unless ${u_2}^G = {a_{2y}}^G$ for some $y \geq 1$, and ${a_2}^G \leq {u_2}^G$ unless ${u_2}^G = {b_1}^G$. Since ${a_{2y}}^G = {A_1}^y$, and ${b_1}^G = B_1$ or $C_1$ according as $G = B_\ell$ or $C_\ell$, the result follows.
\end{proof}

\begin{lem}\label{lem: A_1^2 or D_2 in D_ell}
If $G = D_\ell$ and ${u_2}^G$ is any non-trivial unipotent class apart from $A_1$, then ${u_1}^G \leq {u_2}^G$ for at least one of ${u_1}^G = {A_1}^2$ and ${u_1}^G = D_2$.
\end{lem}

\begin{proof}
First suppose $p \neq 2$; then $\Part(D_2) = 3 1^{2\ell - 3}$, so as even parts must occur with even multiplicity the only unipotent classes ${u_2}^G$ with $\Part(D_2) \not\leq \Part({u_2}^G)$ are ${A_1}^y$ with $\Part({A_1}^y) = 2^{2y} 1^{2\ell - 4y}$, and we have ${A_1}^2 \leq {A_1}^y$ if $y \geq 2$. Now suppose instead $p = 2$. We have classes ${a_{2y}}^G = {A_1}^y$ and ${c_{2y + 2}}^G = {A_1}^y D_2$; if $y \geq 2$ then ${A_1}^2 \leq {A_1}^y$, while if $y \geq 0$ then $D_2 \leq {A_1}^y D_2$. The result follows.
\end{proof}

\begin{lem}\label{lem: C_2 and A_2 in C_ell}
If $G = C_\ell$ with $p \geq 3$, we have ${u_1}^G \leq {u_2}^G$ if ${u_1}^G = C_2$ and $\Part({u_2}^G)$ has a part at least $4$, or if ${u_1}^G = A_2$ and $\Part({u_2}^G)$ has a part $3$.
\end{lem}

\begin{proof}
We have $\Part(C_2) = 4 1^{2\ell - 4}$ and $\Part(A_2) = 3^2 1^{2\ell - 6}$; the first statement is now immediate, and the second follows from the fact that odd parts must occur with even multiplicity.
\end{proof}

\begin{lem}\label{lem: unipotent classes in the two triples}
The unipotent classes containing elements of order $p$ form a totally ordered set in the following cases:
\begin{itemize}
\item[(i)] $G = C_4$ with $p = 3$, when we have $C_1 \leq A_1 \leq A_1 C_1 \leq {A_1}^2 \leq A_2 \leq A_2 C_1$, with the dimensions being $8$, $14$, $18$, $20$, $22$ and $24$ respectively;
\item[(ii)] $G = B_2$ with $p = 5$, when we have $A_1 \leq B_1 \leq B_2$, with the dimensions being $4$, $6$ and $8$ respectively.
\end{itemize}
\end{lem}

\begin{proof}
For both of these we may consult the tables in \cite[IV.2]{Sp}.
\end{proof}

In the remaining results, we take a fixed class ${u_1}^G$ and give a lower bound on $\dim {u_2}^G$ which implies that ${u_1}^G \leq {u_2}^G$. We shall proceed by considering the partially ordered set of unipotent classes $\tilde u^G$ with ${u_1}^G \not\leq \tilde u^G$; for any maximal element $\hat u^G$ of this partially ordered set, we calculate $\dim \hat u^G = \dim G - \dim C_G(\hat u)$ and observe that it does not exceed the given bound. We begin with cases where the rank $\ell$ of $G$ is unbounded.

\begin{lem}\label{lem: various classes in classical groups by dim}
We have ${u_1}^G \leq {u_2}^G$ if one of the following holds:
\begin{itemize}
\item[(i)] $G = A_\ell$ for $\ell \in [3, \infty)$, ${u_1}^G = {A_1}^2$ and $\dim {u_2}^G > 2\ell$;
\item[(ii)] $G = A_\ell$ for $\ell \in [5, \infty)$, ${u_1}^G = {A_1}^3$ and $\dim {u_2}^G > 4\ell - 2$;
\item[(iii)] $G = A_\ell$ for $\ell \in [2, \infty)$, ${u_1}^G = A_2$ and $\dim {u_2}^G > \lfloor \frac{1}{2}(\ell + 1)^2 \rfloor$;
\item[(iv)] $G = A_\ell$ for $\ell \in [5, \infty)$, ${u_1}^G = A_2 A_1$ and $\dim {u_2}^G > \lfloor \frac{1}{2}(\ell + 1)^2 \rfloor$;
\item[(v)] $G = A_\ell$ for $\ell \in [9, \infty)$, ${u_1}^G = A_2{A_1}^2$ and $\dim {u_2}^G > \lfloor \frac{1}{2}(\ell + 1)^2 \rfloor$;
\item[(vi)] $G = A_\ell$ for $\ell \in [3, \infty)$, ${u_1}^G = A_3$ and $\dim {u_2}^G > 2\lfloor \frac{1}{3}(\ell + 1)^2 \rfloor$;
\item[(vii)] $G = A_\ell$ for $\ell \in [9, \infty)$, ${u_1}^G = A_3 A_2$ and $\dim {u_2}^G > 2\lfloor \frac{1}{3}(\ell + 1)^2 \rfloor$;
\item[(viii)] $G = C_\ell$ for $\ell \in [3, \infty)$, ${u_1}^G = A_2$ and $\dim {u_2}^G > \ell(\ell + 1)$;
\item[(ix)] $G = D_\ell$ for $\ell \in [4, \infty)$, ${u_1}^G = D_2$ and $\dim {u_2}^G > \ell(\ell - 1)$.
\end{itemize}
\end{lem}

\begin{proof}
For (i) the only class $\tilde u^G$ is $A_1$, so we take $\hat u^G = A_1$; since $\Part(\hat u^G) = 2 1^{\ell - 1}$, we have $\dim \hat u^G = (\ell + 1)^2 - (\ell^2 + 1^2) = 2\ell$. For (ii) the classes $\tilde u^G$ are $A_1$, ${A_1}^2$ and $A_2$, so we take $\hat u^G = A_2$; since $\Part(\hat u^G) = 3 1^{\ell - 2}$, we have $\dim \hat u^G = (\ell + 1)^2 - ((\ell - 1)^2 + 2.1^2) = 4\ell - 2$. For (iii) the classes $\tilde u^G$ are ${A_1}^y$, so we take $\hat u^G = {A_1}^{\lfloor (\ell + 1)/2 \rfloor}$; if $\ell = 2a - 1$ is odd then $\Part(\hat u^G) = 2^a$, so $\dim \hat u^G = 4a^2 - 2a^2 = 2a^2 = \frac{1}{2}(\ell + 1)^2$, while if $\ell = 2a$ is even then $\Part(\hat u^G) = 2^a 1$, so $\dim \hat u^G = (2a + 1)^2 - ((a + 1)^2 + a^2) = 2a^2 + 2a = \lfloor \frac{1}{2}(\ell + 1)^2 \rfloor$. For (iv) the classes $\tilde u^G$ are ${A_1}^y$ and $A_2$, so we take $\hat u^G = {A_1}^{\lfloor (\ell + 1)/2 \rfloor}$ as in (iii) and $\hat u^G = A_2$; in the latter case, as in (ii) we have $\dim \hat u^G = 4\ell - 2$, which for $\ell \geq 5$ is less than or equal to $\lfloor \frac{1}{2}(\ell + 1)^2 \rfloor$. For (v) the classes $\tilde u^G$ are ${A_1}^y$, $A_2$ and $A_2 A_1$, so we take $\hat u^G = {A_1}^{\lfloor (\ell + 1)/2 \rfloor}$ as in (iii) and $\hat u^G = A_2 A_1$; in the latter case, since $\Part(\hat u^G) = 3 2 1^{\ell - 4}$, we have $\dim \hat u^G = (\ell + 1)^2 - ((\ell - 2)^2 + 2^2 + 1^2) = 6\ell - 8$, which for $\ell \geq 9$ is less than or equal to $\lfloor \frac{1}{2}(\ell + 1)^2 \rfloor$. For (vi) the classes $\tilde u^G$ are ${A_2}^z {A_1}^y$, so we take $\hat u^G = {A_2}^{\lfloor (\ell + 1)/3 \rfloor}$ or ${A_2}^{\lfloor (\ell + 1)/3 \rfloor} A_1$ according as $\ell \equiv b$ mod $3$ for $b \in \{ -1, 0 \}$ or $b = 1$; if $\ell = 3a - 1$ then $\Part(\hat u^G) = 3^a$, so $\dim \hat u^G = 9a^2 - 3a^2 = 6a^2 = \frac{2}{3}(\ell + 1)^2$, if $\ell = 3a$ then $\Part(\hat u^G) = 3^a 1$, so $\dim \hat u^G = (3a + 1)^2 - ((a + 1)^2 + 2a^2) = 6a^2 + 4a = 2\lfloor \frac{1}{3}(\ell + 1)^2 \rfloor$, while if $\ell = 3a + 1$ then $\Part(\hat u^G) = 3^a 2$, so $\dim \hat u^G = (3a + 2)^2 - (2(a + 1)^2 + a^2) = 6a^2 + 8a + 2 = 2\lfloor \frac{1}{3}(\ell + 1)^2 \rfloor$. For (vii) the classes $\tilde u^G$ are ${A_2}^z {A_1}^y$ and $A_3 {A_1}^y$, so we take $\hat u^G = {A_2}^{\lfloor (\ell + 1)/3 \rfloor}$ or ${A_2}^{\lfloor (\ell + 1)/3 \rfloor} A_1$ as in (vi) and $\hat u^G = A_3 {A_1}^{\lfloor (\ell - 3)/2 \rfloor}$; in the latter case, if $\ell = 2a - 1$ is odd then $\Part(\hat u^G) = 4 2^{a - 2}$, so $\dim \hat u^G = 4a^2 - (2(a - 1)^2 + 2.1^2) = 2a^2 + 4a - 4$, which for $\ell \geq 9$ is less than or equal to $2\lfloor \frac{1}{3}(\ell + 1)^2 \rfloor$, while if $\ell = 2a$ is even then $\Part(\hat u^G) = 4 2^{a - 2} 1$, so $\dim \hat u^G = (2a + 1)^2 - (a^2 + (a - 1)^2 + 2.1^2) = 2a^2 + 6a - 2$, which again for $\ell \geq 9$ is less than or equal to $2\lfloor \frac{1}{3}(\ell + 1)^2 \rfloor$. For (viii) (noting that $p \geq 3$ for the elements of ${u_1}^G$ to have order $p$) the classes $\tilde u^G$ are ${A_1}^y$ and ${A_1}^y C_1$, so we take $\hat u^G = {A_1}^{\ell/2}$ or ${A_1}^{(\ell - 1)/2} C_1$ according as $\ell$ is even or odd; since in either case $\Part(\hat u^G) = 2^\ell$, we have $\dim \hat u^G = (2\ell^2 + \ell) - \frac{1}{2}(2\ell^2) = \ell(\ell + 1)$. Finally for (ix) the classes $\tilde u^G$ are ${A_1}^y$, so we take $\hat u^G = {A_1}^{\lfloor \ell/2 \rfloor}$; if $p = 2$ then $\hat u^G = a_{2\lfloor \ell/2 \rfloor}$ so we have $\dim \hat u^G = 2\lfloor \frac{\ell}{2} \rfloor(2\ell - 1 - 2\lfloor \frac{\ell}{2} \rfloor) = \ell(\ell - 1)$; if instead $p \geq 3$, if $\ell = 2a$ is even then $\Part(\hat u^G) = 2^{2a}$, so $\dim \hat u^G = (8a^2 - 2a) - \frac{1}{2}(2(2a)^2) = 4a^2 - 2a = \ell(\ell - 1)$, while if $\ell = 2a + 1$ is odd then $\Part(\hat u^G) = 2^{2a} 1^2$, so $\dim \hat u^G = (8a^2 + 6a + 1) - \frac{1}{2}((2a + 2)^2 + (2a)^2 - 2) = 4a^2 + 2a = \ell(\ell - 1)$. The result follows.
\end{proof}

The remaining results in this section treat cases where the rank $\ell$ of $G$ is fixed; here the condition on ${u_2}^G$ is of the form $\dim {u_2}^G \geq m$ for some $m \in \N$. We may slightly refine the approach described above: provided the class ${u_1}^G$ has dimension at most $m$, any class lying in its boundary will have dimension strictly less than $m$, so does not require consideration; thus it suffices to consider the partially ordered set of classes $\tilde u^G$ which are not comparable to ${u_1}^G$, and we need only show that any maximal element $\hat u^G$ of this partially ordered set has dimension strictly less than $m$.

\begin{lem}\label{lem: various classes in A_ell for fixed ell}
If $G = A_\ell$, we have ${u_1}^G \leq {u_2}^G$ if one of the following holds:
\begin{itemize}
\item[(i)] $\ell = 9$, ${u_1}^G = A_4 A_1$ and $\dim {u_2}^G \geq 75$;
\item[(ii)] $\ell = 5$, ${u_1}^G = A_4$ and $\dim {u_2}^G \geq 28$.
\end{itemize}
\end{lem}

\begin{proof}
For (i) we have $\Part({u_1}^G) = 5 2 1^3$, so $\dim {u_1}^G = 100 - (5^2 + 2^2 + 3.1^2) = 68 < 75$; the only maximal element of the set of classes not comparable to ${u_1}^G$ is $\hat u^G = {A_3}^2 A_1$ with $\Part(\hat u^G) = 4^2 2$, giving $\dim \hat u^G = 100 - (2.3^2 + 2.2^2) = 74 < 75$. For (ii) we have $\Part({u_1}^G) = 5 1$, so $\dim {u_1}^G = 36 - (2^2 + 4.1^2) = 28$; the set of classes not comparable to ${u_1}^G$ is empty. The result follows.
\end{proof}

In some of the remaining results we give two possibilities for the class ${u_1}^G$; what we are claiming in these cases is that at least one of the possibilities lies in the closure of ${u_2}^G$ provided the condition on $\dim{u_2}^G$ is satisfied, and to show this we need only consider classes $\tilde u^G$ which are not comparable to either possibility.

\begin{lem}\label{lem: various classes in B_ell for fixed ell}
If $G = B_\ell$, we have ${u_1}^G \leq {u_2}^G$ if one of the following holds:
\begin{itemize}
\item[(i)] $\ell = 7$, ${u_1}^G = A_2 B_1$ and $\dim {u_2}^G \geq 63$;
\item[(ii)] $\ell = 7$, ${u_1}^G = A_3 B_1$ and $\dim {u_2}^G \geq 73$;
\item[(iii)] $\ell = 7$, ${u_1}^G = A_4 B_1$ and $\dim {u_2}^G \geq 83$;
\item[(iv)] $\ell = 4$, ${u_1}^G = A_2 B_1$ or ${u_1}^G = B_2$ and $\dim {u_2}^G \geq 24$.
\end{itemize}
\end{lem}

\begin{proof}
In each case we note that $p \geq 3$ (at least) for the elements of ${u_1}^G$ to have order $p$. For (i) we have $\Part({u_1}^G) = 3^3 1^6$, so $\dim {u_1}^G = 105 - \frac{1}{2}(9^2 + 2.3^2 - 9) = 60 < 63$; the maximal elements of the set of classes not comparable to ${u_1}^G$ are $\hat u^G = A_2{A_1}^2$ with $\Part(\hat u^G) = 3^2 2^4 1$, giving $\dim \hat u^G = 105 - \frac{1}{2}(7^2 + 6^2 + 2^2 - 3) = 62 < 63$, and $\hat u^G = B_2$ with $\Part(\hat u^G) = 5 1^{10}$, giving $\dim \hat u^G = 105 - \frac{1}{2}(11^2 + 4.1^2 - 11) = 48 < 63$. For (ii) we have $\Part({u_1}^G) = 4^2 3 1^4$, so $\dim {u_1}^G = 105 - \frac{1}{2}(7^2 + 2.3^2 + 2^2 - 5) = 72 < 73$; the maximal elements of the set of classes not comparable to ${u_1}^G$ are $\hat u^G = B_3$ with $\Part(\hat u^G) = 7 1^8$, giving $\dim \hat u^G = 105 - \frac{1}{2}(9^2 + 6.1^2 - 9) = 66 < 73$, and $\hat u^G = B_4(a_2) A_1$ with $\Part(\hat u^G) = 5 3 2^2 1^3$, giving $\dim \hat u^G = 105 - \frac{1}{2}(7^2 + 4^2 + 2^2 + 2.1^2 - 5) = 72 < 73$, and $\hat u^G = {A_2}^2 B_1$ with $\Part(\hat u^G) = 3^5$, giving $\dim \hat u^G = 105 - \frac{1}{2}(3.5^2 - 5) = 70 < 73$. For (iii) we have $\Part({u_1}^G) = 5^2 3 1^2$, so $\dim {u_1}^G = 105 - \frac{1}{2}(5^2 + 2.3^2 + 2.2^2 - 5) = 82 < 83$; the maximal elements of the set of classes not comparable to ${u_1}^G$ are $\hat u^G = B_4$ with $\Part(\hat u^G) = 9 1^6$, giving $\dim \hat u^G = 105 - \frac{1}{2}(7^2 + 8.1^2 - 7) = 80 < 83$, and $\hat u^G = B_5(a_2) A_1$ with $\Part(\hat u^G) = 7 3 2^2 1$, giving $\dim \hat u^G = 105 - \frac{1}{2}(5^2 + 4^2 + 2^2 + 4.1^2 - 3) = 82 < 83$. Finally for (iv) we may consult the tables in \cite[IV.2]{Sp}. The result follows.
\end{proof}

\begin{lem}\label{lem: various classes in D_ell for fixed ell}
If $G = D_\ell$, we have ${u_1}^G \leq {u_2}^G$ if one of the following holds:
\begin{itemize}
\item[(i)] $\ell = 9$, ${u_1}^G = D_3$ and $\dim {u_2}^G \geq 113$;
\item[(ii)] $\ell = 7$, ${u_1}^G = A_3$ or ${u_1}^G = D_3$ and $\dim {u_2}^G \geq 61$;
\item[(iii)] $\ell = 7$, ${u_1}^G = D_3$ and $\dim {u_2}^G \geq 67$;
\item[(iv)] $\ell = 6$, ${u_1}^G = A_2 A_1$ and $\dim {u_2}^G \geq 40$;
\item[(v)] $\ell = 6$, ${u_1}^G = D_3$ and $\dim {u_2}^G \geq 49$;
\item[(vi)] $\ell = 5$, ${u_1}^G = A_2 A_1$ or ${u_1}^G = D_3$ and $\dim {u_2}^G \geq 28$;
\item[(vii)] $\ell = 5$, ${u_1}^G = D_3$ and $\dim {u_2}^G \geq 33$.
\end{itemize}
\end{lem}

\begin{proof}
In each case we note that $p \geq 3$ (at least) for the elements of ${u_1}^G$ to have order $p$. For (i) we have $\Part({u_1}^G) = 5 1^{13}$, so $\dim {u_1}^G = 153 - \frac{1}{2}(14^2 + 4.1^2 - 14) = 60 < 113$; the only maximal element of the set of classes not comparable to ${u_1}^G$ is $\hat u^G = {A_3}^2$ with $\Part(\hat u^G) = 4^4 1^2$, giving $\dim \hat u^G = 153 - \frac{1}{2}(6^2 + 3.4^2 - 2) = 112 < 113$. For (ii) we have $\Part({u_1}^G) = 4^2 1^6$ or $5 1^9$, so $\dim {u_1}^G = 91 - \frac{1}{2}(8^2 + 3.2^2 - 6) = 56 < 61$ or $91 - \frac{1}{2}(10^2 + 4.1^2 - 10) = 44 < 61$; the only maximal element of the set of classes not comparable to either possibility for ${u_1}^G$ is $\hat u^G = {A_2}^2$ with $\Part(\hat u^G) = 3^4 1^2$, giving $\dim \hat u^G = 91 - \frac{1}{2}(6^2 + 2.4^2 - 6) = 60 < 61$. For (iii) we have $\Part({u_1}^G) = 5 1^9$ as in (ii); the only maximal element of the set of classes not comparable to ${u_1}^G$  is $\hat u^G = A_3 A_2$ with $\Part(\hat u^G) = 4^2 3^2$, giving $\dim \hat u^G = 91 - \frac{1}{2}(3.4^2 + 2^2 - 2) = 66 < 67$. For (iv) we have $\Part({u_1}^G) = 3^2 2^2 1^2$, so $\dim {u_1}^G = 66 - \frac{1}{2}(6^2 + 4^2 + 2^2 - 4) = 40$; the only maximal element of the set of classes not comparable to ${u_1}^G$ is $\hat u^G = D_3$ with $\Part(\hat u^G) = 5 1^7$, giving $\dim \hat u^G = 66 - \frac{1}{2}(8^2 + 4.1^2 - 8) = 36 < 40$. For (v) we have $\Part({u_1}^G) =5 1^7$, so as in (iv) $\dim {u_1}^G = 36 < 49$; the only maximal element of the set of classes not comparable to ${u_1}^G$ is $\hat u^G = A_3 D_2$ with $\Part(\hat u^G) = 4^2 3 1$, giving $\dim \hat u^G = 66 - \frac{1}{2}(4^2 + 2.3^2 + 2^2 - 2) = 48 < 49$. For (vi) we have $\Part({u_1}^G) = 3^2 2^2$ or $5 1^5$, so $\dim {u_1}^G = 45 - \frac{1}{2}(2.4^2 + 2^2 - 2) = 28$ or $45 - \frac{1}{2}(6^2 + 4.1^2 - 6) = 28$; the set of classes not comparable to either possibility for ${u_1}^G$ is empty. Finally for (vii) we have $\Part({u_1}^G) = 5 1^5$, so as in (vi) $\dim {u_1}^G = 28 < 33$; the only maximal element of the set of classes not comparable to ${u_1}^G$ is $\hat u^G = A_3$ with $\Part(\hat u^G) = 4^2 1^2$, giving $\dim \hat u^G = 45 - \frac{1}{2}(4^2 + 3.2^2 - 2) = 32 < 33$. The result follows.
\end{proof}

Each of the remaining results in this section may be proved by consulting the tables in \cite[IV.2]{Sp} in conjunction with \cite{Car2} and \cite{LSbook}.

\begin{lem}\label{lem: various classes in C_ell for fixed ell}
If $G = C_\ell$, we have ${u_1}^G \leq {u_2}^G$ if one of the following holds:
\begin{itemize}
\item[(i)] $\ell = 4$ with $p = 2$, ${u_1}^G = A_1 C_1$ or ${u_1}^G = {A_1}^2$ and $\dim {u_2}^G \geq 15$;
\item[(ii)] $\ell = 4$ with $p = 2$, ${u_1}^G = {A_1}^2$ and $\dim {u_2}^G \geq 19$;
\item[(iii)] $\ell = 4$ with $p \geq 5$, ${u_1}^G = {A_1}^2$ or ${u_1}^G = C_2$ and $\dim {u_2}^G \geq 19$;
\item[(iv)] $\ell = 4$ with $p \geq 5$, ${u_1}^G = C_2$ and $\dim {u_2}^G \geq 25$;
\item[(v)] $\ell = 3$ with $p \geq 3$, ${u_1}^G = A_1 C_1$ and $\dim {u_2}^G \geq 11$.
\end{itemize}
\end{lem}

\begin{lem}\label{lem: various classes in E_6}
If $G = E_6$, we have ${u_1}^G \leq {u_2}^G$ if one of the following holds:
\begin{itemize}
\item[(i)] ${u_1}^G = {A_1}^2$ and $\dim {u_2}^G \geq 23$;
\item[(ii)] ${u_1}^G = {A_1}^3$ and $\dim {u_2}^G \geq 33$;
\item[(iii)] ${u_1}^G = A_2$ and $\dim {u_2}^G \geq 41$;
\item[(iv)] ${u_1}^G = A_2 A_1$ and $\dim {u_2}^G \geq 43$;
\item[(v)] ${u_1}^G = A_2{A_1}^2$ and $\dim {u_2}^G \geq 49$;
\item[(vi)] ${u_1}^G = {A_2}^2$ and $\dim {u_2}^G \geq 53$;
\item[(vii)] ${u_1}^G = A_4 A_1$ and $\dim {u_2}^G \geq 61$.
\end{itemize}
\end{lem}

\begin{lem}\label{lem: various classes in E_7}
If $G = E_7$, we have ${u_1}^G \leq {u_2}^G$ if one of the following holds:
\begin{itemize}
\item[(i)] ${u_1}^G = {A_1}^2$ and $\dim {u_2}^G \geq 35$;
\item[(ii)] ${u_1}^G = ({A_1}^3)'$ and $\dim {u_2}^G \geq 55$;
\item[(iii)] ${u_1}^G = A_2 A_1$ and $\dim {u_2}^G \geq 71$;
\item[(iv)] ${u_1}^G = A_2{A_1}^2$ and $\dim {u_2}^G \geq 77$;
\item[(v)] ${u_1}^G = A_3$ and $\dim {u_2}^G \geq 91$;
\item[(vi)] ${u_1}^G = (A_3 A_1)'$ and $\dim {u_2}^G \geq 91$;
\item[(vii)] ${u_1}^G = A_4 A_1$ and $\dim {u_2}^G \geq 103$;
\item[(viii)] ${u_1}^G = A_6$ and $\dim {u_2}^G \geq 115$.
\end{itemize}
\end{lem}

\begin{lem}\label{lem: various classes in E_8}
If $G = E_8$, we have ${u_1}^G \leq {u_2}^G$ if one of the following holds:
\begin{itemize}
\item[(i)] ${u_1}^G = {A_1}^2$ and $\dim {u_2}^G \geq 59$;
\item[(ii)] ${u_1}^G = A_3$ and $\dim {u_2}^G \geq 169$.
\end{itemize}
\end{lem}

\begin{lem}\label{lem: various classes in F_4}
If $G = F_4$, we have ${u_1}^G \leq {u_2}^G$ if one of the following holds:
\begin{itemize}
\item[(i)] ${u_1}^G = A_1 \tilde A_1$ and $\dim {u_2}^G \geq 23$;
\item[(ii)] ${u_1}^G = A_2$ and $\dim {u_2}^G \geq 31$;
\item[(iii)] ${u_1}^G = A_2 \tilde A_1$ and $\dim {u_2}^G \geq 31$;
\item[(iv)] ${u_1}^G = C_3$ and $\dim {u_2}^G \geq 43$.
\end{itemize}
\end{lem}

\section{Preliminary results}\label{sect: prelim}

In this section we prove some preliminary results. The first of these will be used frequently.

\begin{lem}\label{lem: submodule and fixed points}
If $V$ is a $G$-module with submodule $V'$, then for all $g \in G$ we have $\dim C_V(g) \leq \dim C_{V'}(g) + \dim C_{V/V'}(g)$.
\end{lem}

\begin{proof}
Let $\pi : V \to V/V'$ be the quotient map; then the restriction of $\pi$ to $C_V(g)$ has kernel $C_{V'}(g)$ and image contained in $C_{V/V'}(g)$.
\end{proof}

Our next result is a technical one concerning the tensor product of two Jordan block matrices.

\begin{lem}\label{lem: Jordan block tensor product}
If $J_1$ and $J_2$ are matrices comprising single Jordan blocks with eigenvalue $1$, of sizes $r_1$ and $r_2$ respectively, then $J_1 \otimes J_2 - I$ has nullity $\min(r_1, r_2)$.
\end{lem}

\begin{proof}
We may assume $r_1 \leq r_2$. For $t = 1, 2$ let $V_t$ be a vector space of dimension $r_t$ with basis $v_1^t, \dots, v_{r_t}^t$; take the map $\theta_t : V_t \to V_t$ defined by
$$
\theta_t(v_i^t) =
\begin{cases}
v_i^t + v_{i - 1}^t & \hbox{if } i > 1, \\
v_1^t               & \hbox{if } i = 1,
\end{cases}
$$
so that $J_t$ is the matrix of $\theta_t$ with respect to the basis $v_1^t, \dots, v_{r_t}^t$. Set $V_0 = V_1 \otimes V_2$, and for $i \in [1,r_1]$ and $j \in [1, r_2]$ write $v_{ij} = v_i^1 \otimes v_j^2$. Consider the map $\phi = \theta_1 \otimes \theta_2 - 1: V_0 \to V_0$; we have
$$
\phi(v_{ij}) =
\begin{cases}
v_{i, j - 1} + v_{i - 1,j} + v_{i - 1,j - 1} & \hbox{if } i, j > 1, \\
v_{i - 1, 1} & \hbox{if } i > j = 1, \\
v_{1, j - 1} & \hbox{if } j > i = 1, \\
0 & \hbox{if } i = j = 1.
\end{cases}
$$
We claim that the vectors $\phi(v_{ij})$ with $i \in [1, r_1]$ and $j \in [2, r_2]$ form a basis of $\im \phi$. First suppose we have coefficients $\rho_{ij} \in K$ satisfying
$$
0 = \sum_{i = 1}^{r_1} \sum_{j = 2}^{r_2} \rho_{ij} \phi(v_{ij}),
$$
so that
$$
0 = \sum_{j = 2}^{r_2} \rho_{1j} v_{1, j - 1} + \sum_{i = 2}^{r_1} \sum_{j = 2}^{r_2} \rho_{ij}(v_{i, j - 1} + v_{i - 1, j} + v_{i - 1, j - 1}).
$$
For $j \in [2, r_2]$, equating coefficients of $v_{r_1, j - 1}$ shows that $\rho_{r_1, j} = 0$; now for $j \in [2, r_2]$, equating coefficients of $v_{r_1 - 1, j - 1}$ shows that $\rho_{r_1 - 1, j} = 0$; continuing in this way we see that $\rho_{ij} = 0$ for all $i \in [1, r_1]$ and $j \in [2, r_2]$. Thus the vectors specified are linearly independent; let $Z$ be their span. To show that $Z = \im \phi$, since $\phi(v_{11}) = 0$ it suffices to show that if $i \in [2, r_1]$ then $\phi(v_{i1}) = v_{i - 1, 1} \in Z$. We use induction on $i$ to show that if $i + j \leq r_2$ then $v_{ij} \in Z$: if $i = 1$ and $j \leq r_2 - 1$ we have $v_{ij} = \phi(v_{1, j + 1}) \in Z$, while if $i > 1$ and $j \leq r_2 - i$ we have $v_{ij} = \phi(v_{i, j + 1}) - v_{i - 1, j + 1} - v_{i - 1, j} \in Z$ by inductive hypothesis. Thus in particular $v_{11}, \dots, v_{r_1 - 1,1} \in Z$; so we do indeed have $Z = \im \phi$, and the result follows.
\end{proof}

In the case where $r_1 = r_2 = r$, we may view this result as saying that if $u \in A_{r - 1}$ is a regular unipotent element, then $\dim C_{L(\omega_1) \otimes L(\omega_1)}(u) = r$. Our next result treats similarly two submodules of $L(\omega_1) \otimes L(\omega_1)$.

\begin{lem}\label{lem: fixed points on L(omega_2) and L(2omega_1)}
Let $u \in A_{r - 1}$ be a regular unipotent element. Then
\begin{itemize}
\item[(i)] if $V = L(\omega_2)$, then $\dim C_V(u) = \lfloor \frac{r}{2} \rfloor$;
\item[(ii)] if $V = L(2\omega_1)$ with $p \geq 3$, then $\dim C_V(u) = \lceil \frac{r}{2} \rceil$.
\end{itemize}
\end{lem}

\begin{proof}
We prove (i); the proof of (ii) is entirely similar --- alternatively the result follows from (i) and Lemma~\ref{lem: Jordan block tensor product}, since if $p \geq 3$ we have $L(\omega_1) \otimes L(\omega_1) = L(\omega_2) \oplus L(2\omega_1)$.

Take a basis $v_1, \dots, v_r$ of the natural module for $A_{r - 1}$, such that
$$
u.v_i =
\begin{cases}
v_i + v_{i - 1} & \hbox{if } i > 1, \\
v_1             & \hbox{if } i = 1.
\end{cases}
$$
For $1 \leq i < j \leq r$ write $v_{ij} = v_i \otimes v_j - v_j \otimes v_i$, so that $V = \langle v_{ij} : 1 \leq i < j \leq r \rangle$. Let $\phi : V \to V$ be the map $v \mapsto (u - 1).v$; we have
$$
\phi(v_{ij}) =
\begin{cases}
v_{i, j - 1} + v_{i - 1, j} + v_{i - 1, j - 1} & \hbox{if } j - 1 > i > 1, \\
v_{i - 1, i + 1} + v_{i - 1, i}                & \hbox{if } j - 1 = i > 1, \\
v_{1, j - 1}                                   & \hbox{if } j - 1 > i = 1, \\
0                                              & \hbox{if } j - 1 = i = 1.
\end{cases}
$$
We claim that the vectors $\phi(v_{ij})$ with either $i < j - 1$ or $i = j - 1 > \lfloor \frac{r}{2} \rfloor$ form a basis of $\im \phi$. First suppose we have coefficients $\rho_{ij} \in K$ satisfying
$$
0 = \sum_{(i, j)} \rho_{ij} \phi(v_{ij}),
$$
where the sum runs over pairs $(i, j)$ with either $i < j - 1$ or $i = j - 1 > \lfloor \frac{r}{2} \rfloor$, so that
\begin{eqnarray*}
0 & = & \sum_{j > 2} \rho_{1j} v_{1, j - 1} + \sum_{1< i < j - 1} \rho_{ij}(v_{i, j - 1} + v_{i - 1, j} + v_{i - 1, j - 1}) \\
  &   & \qquad {} + \sum_{i > \lfloor \frac{r}{2} \rfloor} \rho_{i, i + 1}(v_{i - 1, i + 1} + v_{i - 1, i}).
\end{eqnarray*}
We show that all $\rho_{ij}$ are zero, working in order of decreasing $i + j$. If $i + j = 2r - 1$, then $(i, j) = (r - 1, r)$; equating coefficients of $v_{r - 2, r}$ shows that $\rho_{r - 1, r} = 0$. Suppose we have shown that whenever $i + j > h$ we have $\rho_{ij} = 0$. If $h > r + 1$, taking successively $i = 1, 2, \dots, \lfloor \frac{2r + 1 - h}{2} \rfloor$ and equating coefficients of $v_{h - r - 2 + i, r + 1 - i}$ shows that $\rho_{h - r - 1 + i, r + 1 - i} = 0$. If $h \leq r + 1$, taking successively $i = \lfloor \frac{h}{2} \rfloor - 1, \dots, 2, 1$ and equating coefficients of $v_{i, h - i - 1}$ shows that $\rho_{i, h - i} = 0$. Thus the vectors specified are linearly independent; let $Z$ be their span. To show that $Z = \im \phi$, since $\phi(v_{12}) = 0$ it suffices to show that if $2 \leq i \leq \lfloor \frac{r}{2} \rfloor$ then $\phi(v_{i, i + 1}) = v_{i - 1, i + 1} + v_{i - 1, i} \in Z$. We use induction on $i$ to show that if $i + j \leq r$ then $v_{ij} \in Z$: if $i = 1$ and $j \leq r - 1$ we have $v_{1j} = \phi(v_{1, j + 1}) \in Z$, while if $i > 1$ and $i < j \leq r - i$ we have $v_{ij} = \phi(v_{i, j + 1}) - v_{i - 1, j + 1} - v_{i - 1, j} \in Z$ by inductive hypothesis. Thus in particular $v_{12} + v_{13}, v_{23} + v_{24}, \dots, v_{\lfloor \frac{r}{2} \rfloor - 1, \lfloor \frac{r}{2} \rfloor} + v_{\lfloor \frac{r}{2} \rfloor - 1, \lfloor \frac{r}{2} \rfloor + 1} \in Z$; so we do indeed have $Z = \im \phi$, and the result follows.
\end{proof}

(In fact \cite[Lemma~3.4]{LSbook} proves both Lemmas~\ref{lem: Jordan block tensor product} and \ref{lem: fixed points on L(omega_2) and L(2omega_1)}; however, it requires the assumption that $p \geq 3$ for both parts of the latter, whereas we shall require Lemma~\ref{lem: fixed points on L(omega_2) and L(2omega_1)}(i) when $p = 2$.)

We may use Lemma~\ref{lem: Jordan block tensor product} to obtain the following.

\begin{lem}\label{lem: half dim bound}
If $A$ is a group of type $A_1$ defined over $K$, and $u \in A \setminus \{ 1 \}$ is unipotent, then for any non-trivial irreducible $A$-module $\tilde V$ we have $\codim C_{\tilde V}(u) \geq \frac{1}{2} \dim \tilde V$.
\end{lem}

\begin{proof}
Let $\tilde V$ be a non-trivial irreducible $A$-module with highest weight $m\omega$, where $\omega$ is the fundamental dominant weight for $A$, so that $m \in \N$. Write $m = m_0 + m_1 p + \cdots + m_t p^t$ such that for all $i$ we have $0 \leq m_i < p$, and $m_t > 0$. By Theorem~\ref{thm: Steinberg} we have $\tilde V = \tilde V_0 \otimes \tilde V_1 \otimes \cdots \otimes \tilde V_t$ where $\tilde V_i = L(m_i\omega)^{(i)}$; write $\tilde V' = \tilde V_0 \otimes \cdots \otimes \tilde V_{t - 1}$, so that $\tilde V = \tilde V' \otimes \tilde V_t$. The matrix representing the action of $u$ on $\tilde V_t$ is a single Jordan block of size $r = m_t + 1$; let the matrix representing the action of $u$ on $\tilde V'$ be a sum of Jordan blocks of sizes $r_1, \dots, r_s$. Given any such Jordan block of size $r_i$, its tensor product with the single Jordan block of size $r$ is a matrix of size $r_ir$, and by Lemma~\ref{lem: Jordan block tensor product} the fixed point space of $u$ on the underlying space has dimension $\min \{ r_i, r \} \leq \frac{1}{2}r_ir$; summing over $i$ gives the result.
\end{proof}

The next result is very straightforward.

\begin{lem}\label{lem: parabolic factorization}
Given a parabolic subgroup $P = QL$ of $G$, where $Q$ is the unipotent radical of $P$ and $L$ the Levi subgroup, let $P^- = Q^-L$ be the opposite parabolic subgroup, so that $Q \cap Q^- = 1$; then $G = P P^- P = Q L Q^- Q$.
\end{lem}

\begin{proof}
First take the case where $P = B$, so that $Q = U$ and $L = T$; write $U^-$ for the product of the root subgroups corresponding to negative roots, so that $Q^- = U^-$. Take $g \in G$, and write $g \in B^x$ for some $x \in G$; use Bruhat decomposition to write $x = bnv$ where $b \in B$, $n \in N$ and $v$ is a product of root elements corresponding to positive roots made negative by $nT$. Thus $g \in B^{nv}$, so we may write $g = (us)^{nv}$ where $u \in U$ and $s \in T$; write $u = u_1 u_2$ where $u_1$ and $u_2$ are products of root elements corresponding to positive roots such that conjugation by $n$ keeps those in $u_1$ positive and makes those in $u_2$ negative, then we have
$$
g = v^{-1}{u_1}^n.s^n.{u_2}^{sn}.v,
$$
with $v^{-1}{u_1}^n, v \in U$, ${u_2}^{sn} \in U^-$ and $s^n \in T$ as required. Now for the general case take $g \in G$ and by the above write $g = v_1 s v_2 v_3$ with $v_1, v_3 \in U$, $v_2 \in U^-$ and $s \in T$; for $i = 1, 2, 3$ write $v_i = q_i l_i$ with $q_1, q_3 \in Q$, $q_2 \in Q^-$ and $l_1, l_2, l_3 \in L$. Since $L$ normalizes both $Q$ and $Q^-$ we have
$$
g = q_1 l_1 s q_2 l_2 q_3 l_3 = q_1.l_1 s l_2 l_3.{q_2}^{l_2 l_3}.{q_3}^{l_3},
$$
with $q_1, {q_3}^{l_3} \in Q$, ${q_2}^{l_2 l_3} \in Q^-$ and $l_1 s l_2 l_3 \in L$ as required.
\end{proof}

The following result will be used repeatedly without comment.

\begin{lem}\label{lem: orbit-stab}
If the connected algebraic group $H$ acts on the variety $X$, and $x \in X$, then $\dim(\overline{H.x}) = \dim H - \dim C_H(x)$.
\end{lem}

\begin{proof}
Consider the morphism $\phi : H \to \overline{H.x}$ defined by $\phi(h) = h.x$; since $H$ is irreducible, so are $H.x$ and $\overline{H.x}$. Thus $\phi$ is a dominant morphism of irreducible varieties, so by \cite[Theorem~4.3]{HumLAG} there is a non-empty set $U \subseteq \phi(H)$ which is open in $\overline{H.x}$ such that if we take $y \in U$ then each component of $\phi^{-1}(y)$ has dimension equal to $\dim H - \dim(\overline{H.x})$; as all fibres are cosets of $C_H(x)$, the result follows.
\end{proof}

The next result is elementary.

\begin{lem}\label{lem: equal dimension fibres}
Let $r \geq 0$ be fixed, and $\phi : X \to Y$ be a dominant morphism of varieties. Suppose that for all $y \in \im\phi$ the fibre $\phi^{-1}(y)$ has dimension $r$; then $\dim X = \dim Y + r$.
\end{lem}

\begin{proof}
Let $X_1, \dots, X_s$ and $Y_1, \dots, Y_t$ be the irreducible components of $X$ and $Y$ respectively. Each set $\phi(X_i)$ is irreducible, so lies in some $Y_j$; and as $Y = \overline{\phi(X)} = \overline{\phi(X_1)} \cup \cdots \cup \overline{\phi(X_s)}$, for each $j$ there exists $i$ with $Y_j = \overline{\phi(X_i)}$. After renumbering we may assume that $\dim Y_1 \geq \dim Y_j$ for all $j > 1$, and that $Y_1 = \overline{\phi(X_1)}$; then $\dim Y = \dim Y_1$. The restriction $\phi : X_1 \to Y_1$ is then a dominant morphism of irreducible varieties, so by \cite[Theorem~4.3]{HumLAG} there is a non-empty set $U \subseteq \phi(X_1)$ which is open in $Y_1$ such that if $y \in U$ then each component of $\phi^{-1}(y)$ in $X_1$ has dimension $\dim X_1 - \dim Y_1$; as all fibres have dimension $r$, we have $\dim X_1 = \dim Y_1 + r = \dim Y + r$. Now take $i > 1$, and let $j$ be such that $\phi(X_i) \subseteq Y_j$; then $\phi : X_i \to \overline{\phi(X_i)}$ is a dominant morphism of irreducible varieties, so as before we obtain $\dim X_i = \dim \overline{\phi(X_i)} + r \leq \dim Y_j + r \leq \dim Y_1 + r = \dim Y + r$. Thus $X_1$ has maximal dimension among the irreducible components of $X$, and so $\dim X = \dim X_1 = \dim Y + r$.
\end{proof}

The next result in this section is simple, but underlies the technique which will be used to show that almost all large triples and quadruples have TGS.

\begin{lem}\label{lem: union over class of fixed points}
Let $X$ be a variety on which $G$ acts. If $g \in G$ with $C_X(g)$ non-empty, then we have
$$
\dim \overline{\bigcup_{g' \in g^G} C_X(g')} \leq \dim g^G + \dim C_X(g).
$$
Moreover if $X$ is a $G$-module $V$, then for $g \in G_{ss}$ and $\kappa \in K^*$ we have
$$
\dim \overline{\bigcup_{g' \in g^G} V_\kappa(g')} \leq \dim g^G + \dim V_\kappa(g).
$$
\end{lem}

\begin{proof}
Write $S = \{ (g', x) : g' \in g^G, \ x \in X, \ g'.x = x \}$; let $\pi_1 : S \to g^G$ and $\pi_2 : S \to X$ be the projections on the first and second components. Then $\pi_1$ is surjective, and for all $g' \in g^G$ we have ${\pi_1}^{-1}(g') = \{ (g', x) : x \in X, \ g'.x = x \} \cong C_X(g') \cong C_X(g)$, so $\dim {\pi_1}^{-1}(g') = \dim C_X(g)$; hence by Lemma~\ref{lem: equal dimension fibres} $\dim S = \dim g^G + \dim C_X(g)$. Since $\im \pi_2 = \bigcup_{g' \in g^G} C_X(g')$, the first statement follows; the proof of the second is entirely similar.
\end{proof}

The final result in this section involves subsystems of $\Phi$. Given a subsystem $\Psi$, let $m_\Psi$ be the size of the smallest possible subsystem which intersects every conjugate of $\Psi$. The values $m_\Psi$ which we will need are given in the following.

\begin{lem}\label{lem: m_Psi values}
Suppose $\Phi$ is of type $A_\ell$. If $\ell \geq 3$ then $m_{{A_1}^2} = \ell(\ell - 1) = M - 2\ell$; if $\ell \geq 5$ then $m_{{A_1}^3} = (\ell - 1)(\ell - 2) = M - (4\ell - 2)$; if $\ell \geq 2$ then $m_{A_2} = \lfloor \frac{1}{2}\ell^2 \rfloor = M - \lfloor \frac{1}{2}(\ell + 1)^2 \rfloor$; and if $\ell \geq 9$ then $m_{A_2{A_1}^2} = m_{A_2}$.
\end{lem}

\begin{proof}
We use the standard notation for the roots in $\Phi$. Let $\Phi'$ be a proper subsystem of $\Phi$; we may assume $\Phi'$ is standard. Write $c$ for the corank of $\Phi'$.

It is clear that for $\ell \geq 3$ the only subsystems $\Phi'$ which intersect every subsystem of type ${A_1}^2$ are those of type $A_{\ell - 1}$, and that for $\ell \geq 5$ the only subsystems $\Phi'$ which intersect every subsystem of type ${A_1}^3$ are those containing a subsystem of type $A_{\ell - 2}$; this gives the values claimed for $m_{{A_1}^2}$ and $m_{{A_1}^3}$. Moreover if $c = 1$ then clearly the subsystem $\Phi'$ intersects every subsystem of type $A_2$, while if $c \geq 2$ we may take $i < j$ with $\ve_i - \ve_{i + 1}, \ve_j - \ve_{j + 1} \notin \Phi'$, and then $\langle \ve_i - \ve_j, \ve_j - \ve_{j + 1} \rangle \subset \Phi \setminus \Phi'$; so for $\ell \geq 2$ the subsystems which intersect every subsystem of type $A_2$ are those of type $A_{\ell'} A_{\ell - 1 - \ell'}$, the smallest of which has $\ell' = \lfloor \frac{1}{2}(\ell - 1) \rfloor$, which gives the value claimed for $m_{A_2}$. This leaves just the value $m_{A_2 {A_1}^2}$ to determine. We claim that for $\ell \geq 6$ the only subsystems $\Phi'$ which intersect every subsystem of type $A_2{A_1}^2$ are those of corank $1$ and those containing a subsystem of type $A_{\ell - 3}$; for $\ell \geq 9$ the latter have $|\Phi'| \geq (\ell - 2)(\ell - 3) > m_{A_2}$.

Suppose $c = 2$; let $\ve_i - \ve_{i + 1}, \ve_j - \ve_{j + 1}$ be the simple roots outside $\Phi'$, with $i < j$. By applying a graph automorphism if necessary, we may assume $\ell - j \geq i - 1$. If $j = \ell$ then $i = 1$, so $\Phi'$ is of type $A_{\ell - 2}$. If $j = \ell - 1$ then $i \in \{ 1, 2 \}$; if $i = 1$ then $\Phi'$ is of type $A_{\ell - 3} A_1$, while if $i = 2$ then $\langle \ve_1 - \ve_4, \ve_4 - \ve_\ell, \ve_2 - \ve_3, \ve_5 - \ve_{\ell + 1} \rangle \subset \Phi \setminus \Phi'$. Thus we may assume $j \leq \ell - 2$. If $i \geq 3$ then $\langle \ve_i - \ve_j, \ve_j - \ve_{j + 1}, \ve_1 - \ve_{j + 2}, \ve_2 - \ve_{j + 3} \rangle \subset \Phi \setminus \Phi'$; so we may assume $i \in \{ 1, 2 \}$. If $j \geq i + 3$ then $\langle \ve_i - \ve_{i + 1}, \ve_{i + 1} - \ve_{j + 1}, \ve_{j - 1} - \ve_{j + 2}, \ve_j - \ve_{j + 3} \rangle \subset \Phi \setminus \Phi'$; so we may assume $j \in \{ i + 1, i + 2 \}$, whence $j \leq 4$. If $j = 4$ then $i = 2$, and $\langle \ve_1 - \ve_4, \ve_4 - \ve_5, \ve_2 - \ve_6, \ve_3 - \ve_7 \rangle \subset \Phi \setminus \Phi'$; if instead $j \leq 3$ then $\Phi'$ is of type $A_{\ell - 2}$ or $A_{\ell - 3} A_1$. Thus (up to graph automorphisms) the only such subsystems $\Phi'$ which intersect every subsystem of type $A_2 {A_1}^2$ are those where $(i, j) = (1, \ell)$, $(1, \ell - 1)$, $(1, 3)$, $(2, 3)$ or $(1, 2)$, which are those containing $A_{\ell - 3}$.

Now suppose $c \geq 3$; then by the above $\Phi'$ lies in a subsystem of corank $2$ which has a subsystem of type $A_2 {A_1}^2$ disjoint from it unless $c = 3$ and (up to graph automorphisms) the simple roots outside $\Phi'$ are $\ve_1 - \ve_2, \ve_2 - \ve_3, \ve_j - \ve_{j + 1}$ for some $j \in \{ 3, \ell \}$, in which case $\Phi'$ is of type $A_{\ell - 3}$. The result follows.
\end{proof}

It may be of interest to compare this result with parts of Lemma~\ref{lem: various classes in classical groups by dim}.

\chapter{Triples having TGS}\label{chap: TGS triples}

In this chapter we develop and then apply techniques to show that a triple has TGS. In Section~\ref{sect: conditions} we give a number of conditions which imply that a triple (or quadruple) has TGS. For the next six sections we concentrate on large triples $(G, \lambda, p)$ in which $\lambda$ is $p$-restricted. In Section~\ref{sect: large triple criteria} we obtain criteria which imply that a large triple satisfies the strongest conditions of Section~\ref{sect: conditions}. In Sections~\ref{sect: large triple relevance} and \ref{sect: large triple exclusion} we use these criteria, firstly in broad terms to restrict the form of the weights $\lambda$ which require consideration, and then in more detailed fashion to produce a list of large triples which must be treated. In Sections~\ref{sect: large triple weight string analysis}, \ref{sect: large triple further analysis} and \ref{sect: two triples} we employ successively more careful types of analysis of weights to show that the remaining large triples not listed in Table~\ref{table: large triple and first quadruple non-TGS} (and two of those which do) satisfy some of the weaker conditions of Section~\ref{sect: conditions}. Finally in Section~\ref{sect: large triple tensor products} we deal with large triples $(G, \lambda, p)$ in which $\lambda$ is not $p$-restricted.

\section{Conditions implying TGS}\label{sect: conditions}

In this opening section we consider both triples and quadruples; this is because the results obtained will be used both in this chapter and the next, where we prove that the large quadruples not listed in Table~\ref{table: large triple and first quadruple non-TGS} or Table~\ref{table: large higher quadruple non-TGS} have TGS. Let $(G, \lambda, p)$ or $(G, \lambda, p, k)$ be a triple or quadruple; write $V = L(\lambda)$ and set $X = V$ or $\Gk(V)$ respectively. In this section we produce a series of conditions which show that the triple or quadruple has TGS.

Recall that we write
$$
M = |\Phi| = \dim G - \rank G = \dim G_u.
$$
Let $g \mapsto \bar g = gZ(G)$ be the projection $G \to G/Z(G)$. For $r \in \N$ set
$$
G_{(r)} = \{ g \in G : o(\bar g) = r \}.
$$
In \cite{Lawdim} a lower bound $d_{\Phi, r}$ is given for $\codim G_{(r)}$, with $d_{\Phi, r} \geq \rank G$. Write
$$
M_r = \dim G - d_{\Phi, r};
$$
thus $M \geq M_r \geq \dim G_{(r)}$. For convenience we give the values $M_2$ and $M_3$ in the following table.
$$
\begin{array}{|c|c|c|c|c|c|c|}
\cline{1-3} \cline{5-7}
   G   &                   M_2                   &                     M_3                     & \ptw &  G  & M_2 & M_3 \tbs \\
\cline{1-3} \cline{5-7}
A_\ell & \lfloor \frac{1}{2}(\ell + 1)^2 \rfloor &   2\lfloor \frac{1}{3}(\ell + 1)^2 \rfloor  &      & E_6 &  40 &  54 \tbs \\
B_\ell &             \ell(\ell + 1)              & 2\lfloor \frac{1}{3}\ell(2\ell + 1) \rfloor &      & E_7 &  70 &  90 \tbs \\
C_\ell &             \ell(\ell + 1)              & 2\lfloor \frac{1}{3}\ell(2\ell + 1) \rfloor &      & E_8 & 128 & 168 \tbs \\
D_\ell &   2\lfloor \frac{1}{2}\ell^2 \rfloor    & 2\lfloor \frac{1}{3}\ell(2\ell - 1) \rfloor &      & F_4 &  28 &  36 \tbs \\
\cline{1-3}
\multicolumn{4}{c|}{}                                                                                 & G_2 &   8 &  10 \tbs \\
\cline{5-7}
\end{array}
$$
We shall also need to know that if $G = B_7$, $F_4$ or $E_7$ then $M_5 = 84$, $40$ or $106$ respectively. In addition, if $p = \infty$ we likewise define
$$
G_{(p)} = G_u \setminus \{ 1 \},
$$
and $M_p = M$; then $M_p = \dim G_{(p)}$.

The following elementary result is fundamental to our approach. Recall that we write $\P'$ for the set of primes other than $p$; thus if $r \in \P'$ then $G_{(r)} \subset G_{ss} \setminus Z(G)$.

\begin{prop}\label{prop: stabilizers meet G_{(p)} or G_{(r)}}
If $x \in X$ is such that $C_G(x) \not\leq Z(G)$, then $C_G(x)$ meets either $G_{(r)}$ for some $r \in \P'$, or $G_{(p)}$.
\end{prop}

\begin{proof}
Take $g \in C_G(x) \setminus Z(G)$; then $\langle g \rangle \leq C_G(x)$, and as the stabilizer $C_G(x)$ is closed by \cite[Proposition~8.2(b)]{HumLAG} we have $\overline{\langle g \rangle} \leq C_G(x)$. Let $g = su$ be the Jordan decomposition of $g$, with $s$ semisimple and $u$ unipotent; by \cite[Theorem~15.3(a)]{HumLAG} we have $s, u \in \overline{\langle g \rangle}$. If $u \neq 1$, then either $p$ is finite and $u$ has order $p^a$ for some $a \geq 1$, in which case $u^{p^{a - 1}} \in C_G(x) \cap G_{(p)}$, or $p = \infty$, in which case $u \in C_G(x) \cap G_{(p)}$; thus we may assume $u = 1$, whence $g = s$. Let $T'$ be a maximal torus of $G$ containing $g$; then $\langle g \rangle \leq T'$, and as $T'$ is closed we have $\overline{\langle g \rangle} \leq T'$. If $\overline{\langle g \rangle}$ is finite then $g$ has finite order, and then some power of $g$ lies in $C_G(x) \cap G_{(r)}$ for some $r \in \P'$; if instead $\overline{\langle g \rangle}$ is infinite then its connected component is a closed connected subgroup of $T'$ and so must be a torus, whence for any $r \in \P'$ it contains a non-central element $h$ such that $\bar h \in G/Z(G)$ has order $r$, so that $h \in C_G(x) \cap G_{(r)}$. The result follows.
\end{proof}

Our strategy is then to seek to show that the set of points $x$ in $X$ whose stabilizer $C_G(x)$ in $G$ contains an element as given in Proposition~\ref{prop: stabilizers meet G_{(p)} or G_{(r)}} lies in a proper subvariety of $X$; as we shall see below, this implies that the triple or quadruple has TGS. It is thus natural to subdivide the problem into consideration of semisimple elements and of unipotent elements. For each such type of element we shall obtain a hierarchy of conditions. The most basic are as follows. We say that the triple $(G, \lambda, p)$ or quadruple $(G, \lambda, p, k)$ satisfies condition $\sscon$ if
$$
\bigcup_{r \in \P'} \bigcup_{s \in G_{(r)}} C_X(s) \quad \hbox{ lies in a proper subvariety of }X,
$$
and condition $\ucon$ if
$$
\bigcup_{u \in G_{(p)}} C_X(u) \quad \hbox{ lies in a proper subvariety of }X.
$$
Our first result is then the following.

\begin{prop}\label{prop: sscon and ucon imply TGS}
If the triple $(G, \lambda, p)$ or quadruple $(G, \lambda, p, k)$ satisfies both $\sscon$ and $\ucon$, it has TGS.
\end{prop}

\begin{proof}
If the triple $(G, \lambda, p)$ or quadruple $(G, \lambda, p, k)$ satisfies both $\sscon$ and $\ucon$, then by Proposition~\ref{prop: stabilizers meet G_{(p)} or G_{(r)}} the intersection of the complements of the two proper subvarieties concerned is a non-empty open set each of whose points has stabilizer contained in $Z(G)$. For each $z \in Z(G) \setminus G_X$, the fixed point set $C_X(z)$ is a proper subvariety of $V$; since $Z(G)$ is finite, the complement of the union of these is another non-empty open set. The intersection of these two non-empty open sets is then itself a non-empty open set each of whose points has trivial stabilizer in $G/G_X$; so the triple $(G, \lambda, p)$ or quadruple $(G, \lambda, p, k)$ has TGS.
\end{proof}

We now give further conditions involving semisimple elements. The first of these concerns triples only; recall that given $s \in G_{ss}$ and $\kappa \in K^*$ we have $V_\kappa(s) = \{ v \in V : s.v = \kappa v \}$. We say that the triple $(G, \lambda, p)$ satisfies condition $\ssevcon$ if
$$
\bigcup_{r \in \P'} \bigcup_{s \in G_{(r)}} \bigcup_{\kappa \in K^*} V_\kappa(s) \quad \hbox{ lies in a proper subvariety of }V.
$$
(There is no corresponding condition for quadruples, since in actions on projective varieties the concept of eigenvalues other than $1$ does not arise.)

\begin{prop}\label{prop: ssevcon implies sscon}
If the triple $(G, \lambda, p)$ satisfies $\ssevcon$, it satisfies $\sscon$.
\end{prop}

\begin{proof}
If $s \in G_{(r)}$ for some $r \in \P'$, the eigenspace $V_1(s)$ is just the fixed point space $C_V(s)$; thus the union in condition $\sscon$ is a subset of that in condition $\ssevcon$. The result follows.
\end{proof}

In order to obtain a criterion which enables us to deduce that $\sscon$ or $\ssevcon$ holds for a given triple or quadruple, we shall employ a modified form of the approach taken by Kenneally in \cite{Ken}. We shall say that the triple $(G, \lambda, p)$ or quadruple $(G, \lambda, p, k)$ satisfies condition $\ssdiamcon$ if
$$
\hbox{for all } r \in \P' \hbox{ and all } s \in G_{(r)} \hbox{ we have } \codim C_X(s) > \dim s^G,
$$
and that the triple $(G, \lambda, p)$ satisfies condition $\ssdiamevcon$ if
$$
\hbox{for all } r \in \P', \hbox{ all } s \in G_{(r)} \hbox{ and all } \kappa \in K^* \hbox{ we have } \codim V_\kappa(s) > \dim s^G.
$$

\begin{prop}\label{prop: ssdiamevcon implies ssdiamcon}
If the triple $(G, \lambda, p)$ satisfies $\ssdiamevcon$, it satisfies $\ssdiamcon$.
\end{prop}

\begin{proof}
Again, this follows because if $s \in G_{(r)}$ for some $r \in \P'$ then we have $V_1(s) = C_V(s)$.
\end{proof}

\begin{prop}\label{prop: ssdiamcon implies sscon and ssdiamevcon implies ssevcon}
If the triple $(G, \lambda, p)$ satisfies $\ssdiamcon$ or $\ssdiamevcon$, it satisfies $\sscon$ or $\ssevcon$ respectively; likewise if the quadruple $(G, \lambda, p, k)$ satisfies $\ssdiamcon$, it satisfies $\sscon$.
\end{prop}

\begin{proof}
Set $\Gamma = \{ 1 \}$ or $K^*$. Write
$$
S = \bigcup_{r \in \P'} G_{(r)}.
$$
Recall that $T$ is a fixed maximal torus of $G$, so that each conjugacy class in $S$ meets $T$. View $V$ as a direct sum of weight spaces for $T$; since there are only finitely many ways of grouping these weight spaces into eigenspaces for an individual element of $T$, there exist $m \in \N$ and $t_1, \dots, t_m \in T \cap S$ such that if $t \in T \cap S$ then there exists $i \in [1, m]$ such that $t$ and $t_i$ have the same collection of eigenspaces with eigenvalues in $\Gamma$ (although if $\Gamma = K^*$ the eigenvalues themselves need not be the same). Note that this means that $t$ and $t_i$ have the same fixed point variety in the action on $\Gk(V)$, because any $k$-dimensional subspace of $V$ is fixed by a semisimple element if and only if it has a basis consisting of eigenvectors.

Now given $s \in S$ there exists $t \in T \cap S$ conjugate to $s$, and so $s$ has the same collection of eigenspaces with eigenvalues in $\Gamma$ as some conjugate of some $t_i$. Thus, writing $\Gamma_i$ for the finite set of eigenvalues in $\Gamma$ in the action of $t_i$ on $V$, we have
\begin{eqnarray*}
\bigcup_{r \in \P'} \bigcup_{s \in G_{(r)}} \bigcup_{\kappa \in \Gamma} V_\kappa(s)
& = & \bigcup_{i = 1}^m \bigcup_{s \in {t_i}^G} \bigcup_{\kappa \in \Gamma} V_\kappa(s) \\
& = & \bigcup_{i = 1}^m \bigcup_{\kappa \in \Gamma_i} \bigcup_{s \in {t_i}^G} V_\kappa(s);
\end{eqnarray*}
so using Lemma~\ref{lem: union over class of fixed points} we obtain
\begin{eqnarray*}
\dim \overline{\left( \bigcup_{r \in \P'} \bigcup_{s \in G_{(r)}} \bigcup_{\kappa \in \Gamma} V_\kappa(s) \right)}
& = & \max_{1 \leq i \leq m} \max_{\kappa \in \Gamma_i} \dim \overline{\left( \bigcup_{s \in {t_i}^G} V_\kappa(s) \right)} \\
& \leq & \max_{1 \leq i \leq m} \max_{\kappa \in \Gamma_i} \left( \dim {t_i}^G + \dim V_\kappa(t_i) \right) \\
& = & \dim V - \min_{1 \leq i \leq m} \left( \min_{\kappa \in \Gamma_i} \codim V_\kappa(t_i) - \dim {t_i}^G \right).
\end{eqnarray*}
Taking $\Gamma = K^*$, we see that if the triple $(G, \lambda, p)$ satisfies $\ssdiamevcon$, then
$$
\dim \overline{\left( \bigcup_{r \in \P'} \bigcup_{s \in G_{(r)}} \bigcup_{\kappa \in K^*} V_\kappa(s) \right)} < \dim V,
$$
so it satisfies $\ssevcon$; likewise taking $\Gamma = \{ 1 \}$, we see that if it satisfies $\ssdiamcon$, it satisfies $\sscon$. Finally by taking $\Gamma = \{ 1 \}$ and replacing each $V_1(s)$ by $C_{\Gk(V)}(s)$ we see that if the quadruple $(G, \lambda, p, k)$ satisfies $\ssdiamcon$, it satisfies $\sscon$.
\end{proof}

Our next result links conditions for triples and the associated first quadruples.

\begin{prop}\label{prop: ssdiamevcon for triples implies ssdiamcon for first quadruples}
If the triple $(G, \lambda, p)$ satisfies $\ssdiamevcon$, the associated first quadruple $(G, \lambda, p, 1)$ satisfies $\ssdiamcon$.
\end{prop}

\begin{proof}
Write $d = \dim V$; then $\dim \G{1}(V) = d - 1$. Take $s \in G_{(r)}$ for some $r \in \P'$, and let $d_1 = \max \{ \dim V_\kappa(s) : \kappa \in K^* \}$; 
then $\dim C_{\G{1}(V)}(s) = d_1 - 1$, so
$$
\codim C_{\G{1}(V)}(s) = (d - 1) - (d_1 - 1) = d - d_1 = \codim V_\kappa(s).
$$
If the triple $(G, \lambda, p)$ satisfies $\ssdiamevcon$, then $\codim C_{\G{1}(V)}(s) = \codim V_\kappa(s) > \dim s^G$, whence the quadruple $(G, \lambda, p, 1)$ satisfies $\ssdiamcon$.
\end{proof}

We next consider a slightly stronger condition than $\ssdiamevcon$. We say that the triple $(G, \lambda, p)$ satisfies condition $\ssdagcon$ if
$$
\hbox{for all } r \in \P', \hbox{ all } s \in G_{(r)} \hbox{ and all } \kappa \in K^* \hbox{ we have } \codim V_\kappa(s) > M_r.
$$

\begin{prop}\label{prop: ssdagcon implies ssdiamevcon}
If the triple $(G, \lambda, p)$ satisfies $\ssdagcon$, it satisfies $\ssdiamevcon$.
\end{prop}

\begin{proof}
This is immediate since if $s \in G_{(r)}$ for $r \in \P'$ then $M_r \geq \dim G_{(r)} \geq \dim s^G$.
\end{proof}

Our final condition on semisimple elements is the strongest of all. We say that the triple $(G, \lambda, p)$ satisfies condition $\ssddagcon$ if
$$
\hbox{for all } r \in \P', \hbox{ all } s \in G_{(r)} \hbox{ and all } \kappa \in K^* \hbox{ we have } \codim V_\kappa(s) > M.
$$

\begin{prop}\label{prop: ssddagcon implies ssdagcon}
If the triple $(G, \lambda, p)$ satisfies $\ssddagcon$, it satisfies $\ssdagcon$.
\end{prop}

\begin{proof}
This is immediate since for all $r \in \N$ we have $M \geq M_r$.
\end{proof}


We now turn to further conditions involving unipotent elements. We say that the triple $(G, \lambda, p)$ or quadruple $(G, \lambda, p, k)$ satisfies condition $\udiamcon$ if
$$
\hbox{for all } u \in G_{(p)} \hbox{ we have } \codim C_X(u) > \dim u^G.
$$

\begin{prop}\label{prop: udiamcon implies ucon}
If the triple $(G, \lambda, p)$ or quadruple $(G, \lambda, p, k)$ satisfies $\udiamcon$, it satisfies $\ucon$.
\end{prop}

\begin{proof}
As $G$ has only finitely many unipotent classes, there exist $m \in \N$ and $u_1, \dots, u_m \in G_{(p)}$ such that $G_{(p)} = \bigcup_{i = 1}^m {u_i}^G$; thus
$$
\bigcup_{u \in G_{(p)}} C_X(u) = \bigcup_{i = 1}^m \bigcup_{g \in {u_i}^G} C_X(g),
$$
and so using Lemma~\ref{lem: union over class of fixed points} we obtain
\begin{eqnarray*}
\dim \overline{\left( \bigcup_{u \in G_{(p)}} C_X(u) \right)}
& = & \max_{1 \leq i \leq m} \dim \overline{\left( \bigcup_{g \in {u_i}^G} C_X(g) \right)} \\
& \leq & \max_{1 \leq i \leq m} \left( \dim {u_i}^G + \dim C_X(u_i) \right) \\
& = & \dim V - \min_{1 \leq i \leq m} (\codim C_X(u_i) - \dim {u_i}^G).
\end{eqnarray*}
Thus if the triple $(G, \lambda, p)$ or quadruple $(G, \lambda, p, k)$ satisfies $\udiamcon$, it satisfies $\ucon$.
\end{proof}

Our next result links conditions for triples and the associated first quadruples, and is the analogue of Proposition~\ref{prop: ssdiamevcon for triples implies ssdiamcon for first quadruples}.

\begin{prop}\label{prop: udiamcon for triples implies udiamcon for first quadruples}
If the triple $(G, \lambda, p)$ satisfies $\udiamcon$, the associated first quadruple $(G, \lambda, p, 1)$ satisfies $\udiamcon$.
\end{prop}

\begin{proof}
Write $d = \dim V$; then $\dim \G{1}(V) = d - 1$. Take $u \in G_{(p)}$, and let $d_1 = \dim C_V(u)$. Then $\dim C_{\G{1}(V)}(u) = d_1 - 1$, so
$$
\codim C_{\G{1}(V)}(u) = (d - 1) - (d_1 - 1) = d - d_1 = \codim C_V(u).
$$
If the triple $(G, \lambda, p)$ satisfies $\udiamcon$, then $\codim C_{\G{1}(V)}(u) = \codim C_V(u) > \dim u^G$, whence the quadruple $(G, \lambda, p, 1)$ satisfies $\udiamcon$.
\end{proof}

We now produce a stronger condition than $\udiamcon$. We say that the triple $(G, \lambda, p)$ satisfies condition $\udagcon$ if
$$
\hbox{for all root elements } u \hbox{ we have } \codim C_V(u) > M_p.
$$

\begin{prop}\label{prop: udagcon implies udiamcon}
If the triple $(G, \lambda, p)$ satisfies $\udagcon$, it satisfies $\udiamcon$.
\end{prop}

\begin{proof}
By Lemma~\ref{lem: root elt class in closure of any non-triv class} we know that the closure of any non-identity unipotent class in $G$ contains root elements; since $M_p \geq \dim G_{(p)}$, the result follows from Lemma~\ref{lem: unip closure containment}.
\end{proof}

Our final condition on unipotent elements is the strongest of all. We say that the triple $(G, \lambda, p)$ satisfies condition $\uddagcon$ if
$$
\hbox{for all root elements } u \hbox{ we have } \codim C_V(u) > M.
$$

\begin{prop}\label{prop: uddagcon implies udagcon}
If the triple $(G, \lambda, p)$ satisfies $\uddagcon$, it satisfies $\udagcon$.
\end{prop}

\begin{proof}
This is immediate since $M \geq M_p$.
\end{proof}


We may summarise the relationships between the various conditions for triples in the following diagram.

\begin{center}
\begin{picture}(208,130)(000,-020)
\put(000,084){\makebox(0,0){$\pss\ssddagcon$}}%
\put(015,084){\vector(1,0){20}}%
\put(025,090){\makebox(0,0){$\ss{\ref{prop: ssddagcon implies ssdagcon}}$}}%
\put(050,084){\makebox(0,0){$\pss\ssdagcon$}}%
\put(065,084){\vector(1,0){20}}%
\put(075,090){\makebox(0,0){$\ss{\ref{prop: ssdagcon implies ssdiamevcon}}$}}%
\put(100,084){\makebox(0,0){$\pss\ssdiamevcon$}}%
\put(115,084){\vector(1,0){20}}%
\put(125,090){\makebox(0,0){$\ss{\ref{prop: ssdiamcon implies sscon and ssdiamevcon implies ssevcon}}$}}%
\put(150,084){\makebox(0,0){$\pss\ssevcon$}}%
\put(100,076){\vector(0,-1){20}}%
\put(089,066){\makebox(0,0){$\ss{\ref{prop: ssdiamevcon implies ssdiamcon}}$}}%
\put(150,076){\vector(0,-1){20}}%
\put(161,066){\makebox(0,0){$\ss{\ref{prop: ssevcon implies sscon}}$}}%
\put(100,048){\makebox(0,0){$\pss\ssdiamcon$}}%
\put(115,048){\vector(1,0){20}}%
\put(125,042){\makebox(0,0){$\ss{\ref{prop: ssdiamcon implies sscon and ssdiamevcon implies ssevcon}}$}}%
\put(150,048){\makebox(0,0){$\pss\sscon$}}%
\put(000,000){\makebox(0,0){$\pu\uddagcon$}}%
\put(015,000){\vector(1,0){20}}%
\put(025,006){\makebox(0,0){$\ss{\ref{prop: uddagcon implies udagcon}}$}}%
\put(050,000){\makebox(0,0){$\pu\udagcon$}}%
\put(065,000){\vector(1,0){20}}%
\put(075,006){\makebox(0,0){$\ss{\ref{prop: udagcon implies udiamcon}}$}}%
\put(100,000){\makebox(0,0){$\pu\udiamcon$}}%
\put(115,000){\vector(1,0){20}}%
\put(125,006){\makebox(0,0){$\ss{\ref{prop: udiamcon implies ucon}}$}}%
\put(150,000){\makebox(0,0){$\pu\ucon$}}%
\put(158,038){\line(1,-1){14}}%
\put(158,010){\line(1,1){14}}%
\put(172,024){\vector(1,0){20}}%
\put(182,030){\makebox(0,0){$\ss{\ref{prop: sscon and ucon imply TGS}}$}}%
\put(208,024){\makebox(0,0){$\hbox{TGS}$}}%
\end{picture}
\end{center}

Thus if a large triple satisfies any one of the conditions in this section concerning semisimple elements, and any one of those concerning unipotent elements, it has TGS. In the remainder of this chapter we shall show that any large triple not listed in Table~\ref{table: large triple and first quadruple non-TGS} satisfies $\ssdiamevcon$ and $\udiamcon$ (frequently by showing that it satisfies stronger conditions). It will then follow from Propositions~\ref{prop: ssdiamevcon for triples implies ssdiamcon for first quadruples} and \ref{prop: udiamcon for triples implies udiamcon for first quadruples} that the associated first quadruple satisfies $\ssdiamcon$ and $\udiamcon$, and so also has TGS; this will be of use in Chapter~\ref{chap: TGS quadruples}.

\section{Criteria involving bounds for codimensions}\label{sect: large triple criteria}

Let $(G, \lambda, p)$ be a large triple; write $V = L(\lambda)$. If $\lambda$ is a $p$-restricted dominant weight for $G$, we shall call $(G, \lambda, p)$ a {\em $p$-restricted large triple\/}. Large triples $(G, \lambda, p)$ which are not $p$-restricted will be considered in Section~\ref{sect: large triple tensor products}.

Observe that conditions $\ssddagcon$ and $\uddagcon$ of Section~\ref{sect: conditions} both require certain subspaces of $V$ to have codimension greater than $M$. In this section we shall produce a value determined by $\lambda$ which will be a lower bound for both types of codimension (if $\Phi$ has two root lengths, two values may be required); then whenever $M$ is strictly less than this value (or these values) we know that $(G, \lambda, p)$ satisfies both $\ssddagcon$ and $\uddagcon$, and thus has TGS.

We shall frequently employ the following slight abuse of terminology: given $\kappa \in K^*$ and $\mu \in \Lambda(V)$, we say that $\mu$ lies in the eigenspace $V_\kappa(s)$ if $\mu(s) = \kappa$.

Recall that we define $e(\Phi)$ as the maximum ratio of squared root lengths in $\Phi$, and that if $e(\Phi) = 1$ we choose to regard all roots as short rather than long. A subsystem of $\Phi$ which is generated by a subset of $\Pi$ will be called {\em standard\/}. Given a standard subsystem $\Psi$ of $\Phi$, we let $W(\Psi)$ be the Weyl group of $\Psi$, and define
$$
r_\Psi = \frac{|W : W(\Psi)|.|\Phi_s \setminus \Psi_s|}{2|\Phi_s|}, \qquad
{r_\Psi}' = \frac{|W : W(\Psi)|.|\Phi_l \setminus \Psi_l|}{2|\Phi_l|} \ \hbox{ if } e(\Phi) > 1.
$$

Given a dominant weight $\mu = \sum_{j = 1}^\ell a_j \omega_j$, set $\Psi = \Psi(\mu) = \langle \alpha_i : a_i = 0 \rangle$ and define
$$
r_\mu = r_\Psi, \qquad {r_\mu}' = {r_\Psi}' \ \hbox{ if } e(\Phi) > 1.
$$
Given a $p$-restricted dominant weight $\lambda$, set
$$
s_\lambda = \sum r_\mu, \qquad {s_\lambda}' = \sum {r_\mu}' \ \hbox{ if } e(\Phi) > 1,
$$
where each sum runs over the dominant weights $\mu \preceq \lambda$.

For $p > e(\Phi)$ we may apply Theorem~\ref{thm: Prem}, which enables us to prove the following.

\begin{prop}\label{prop: codim bound using Premet}
Let $(G, \lambda, p)$ be a $p$-restricted large triple; write $V = L(\lambda)$ and assume $p > e(\Phi)$. Then
\begin{itemize}
\item[(i)] for all $r \in \P'$, $s \in G_{(r)}$ and $\kappa \in K^*$ we have $\codim V_\kappa(s) \geq s_\lambda$;
\item[(ii)] for all $\alpha \in \Phi_s$ we have $\codim C_V(x_\alpha(1)) \geq s_\lambda$;
\item[(iii)] if $e(\Phi) > 1$, for all $\beta \in \Phi_l$ we have $\codim C_V(x_\beta(1)) \geq {s_\lambda}'$.
\end{itemize}
\end{prop}

\begin{proof}
Take a dominant weight $\mu \preceq \lambda$; by Theorem~\ref{thm: Prem}, the assumption on $p$ implies that $\mu \in \Lambda(V)$. Write $\Psi = \Psi(\mu)$; then $\Psi = \{ \alpha \in \Phi : \langle \mu, \alpha \rangle = 0 \}$. Thus $\Phi_s \setminus \Psi_s$ consists of the short roots in $\Phi$ not orthogonal to $\mu$. The stabilizer of $\mu$ in $W$ is $W(\Psi)$, so the orbit $W.\mu$ has size $|W : W(\Psi)|$; thus the number of pairs $(\nu, \alpha) \in W.\mu \times \Phi_s$ with $\langle \nu, \alpha \rangle \neq 0$ is $|W : W(\Psi)|.|\Phi_s \setminus \Psi_s| = 2r_\Psi |\Phi_s|$, and so for any given $\alpha \in \Phi_s$ the number of weights $\nu \in W.\mu$ not orthogonal to $\alpha$ is $2r_\Psi = 2r_\mu$. Letting $\mu$ run over the dominant weights ${} \preceq \lambda$ we see that, for a fixed $\alpha \in \Phi_s$, the number of weights in $\Lambda(V)$ not orthogonal to $\alpha$ is $2s_\lambda$. Moreover if $e(\Phi) > 1$, an exactly similar argument shows that, for a fixed $\beta \in \Phi_l$, the number of weights in $\Lambda(V)$ not orthogonal to $\beta$ is $2{s_\lambda}'$.

Now take $r \in \P'$, $s \in G_{(r)}$ and $\kappa \in K^*$; then there exists $\alpha \in \Phi_s$ with $\alpha(s) \neq 1$ (note that if $e(\Phi) > 1$ then any long root is a sum of two short roots). For this $\alpha$ we consider the $\alpha$-strings in $\Lambda(V)$; since $\alpha(s) \neq 1$, two weights which are adjacent in an $\alpha$-string cannot both lie in $V_\kappa(s)$. An $\alpha$-string of even length contains no weight orthogonal to $\alpha$, and the contribution to $\codim V_\kappa(s)$ is at least half of its length; an $\alpha$-string of odd length contains exactly one weight orthogonal to $\alpha$, and the contribution to $\codim V_\kappa(s)$ is at least half of one less than its length. Summing over the various $\alpha$-strings gives $\codim V_\kappa(s) \geq s_\lambda$, proving (i).

Now take $\alpha \in \Phi_s$ and write $A = \langle X_{\pm\alpha} \rangle \cong A_1$; again consider the $\alpha$-strings in $\Lambda(V)$. Given such an $\alpha$-string
$$
\nu - t\alpha \quad \nu - (t - 1)\alpha \quad \dots \quad \nu - \alpha \quad \nu,
$$
the sum of the corresponding weight spaces in $V$ is an $A$-module, and $V$ is the direct sum of these $A$-modules. For each such $A$-module, take a composition series, and consider one of the composition factors. If it is trivial, the weight $\nu - i\alpha$ above to which it corresponds is orthogonal to $\alpha$ (so $i = \frac{t}{2}$). If instead it is non-trivial, it is a sum of $1$-dimensional weight spaces corresponding to distinct weights $\nu - i\alpha$, and Lemma~\ref{lem: half dim bound} shows that the codimension of the fixed point space of $x_\alpha(1)$ on it is at least half of its dimension. Summing over the various composition factors in the different $\alpha$-strings, and using Lemma~\ref{lem: submodule and fixed points} repeatedly, we see that $\codim C_V(x_\alpha(1))$ is at least half of the number of weights in $\Lambda(V)$ not orthogonal to $\alpha$, i.e., $\codim C_V(x_\alpha(1)) \geq s_\lambda$, proving (ii).

Finally if $e(\Phi) > 1$, an exactly similar argument proves (iii).
\end{proof}

\begin{cor}\label{cor: TGS using Premet}
Let $(G, \lambda, p)$ be a $p$-restricted large triple; assume $p > e(\Phi)$. If $s_\lambda > M$, and also ${s_\lambda}' > M$ if $e(\Phi) > 1$, then the triple $(G, \lambda, p)$ satisfies $\ssddagcon$ and $\uddagcon$, and thus has TGS.
\end{cor}

\begin{proof}
This is immediate.
\end{proof}

However, for $p \leq e(\Phi)$ we cannot use Theorem~\ref{thm: Prem}; here a slightly different approach is required. Given a dominant weight $\mu$, define
$$
r_{\mu, p} = \frac{|W.\mu|}{|\Phi_s|}.|\{ \alpha \in \Phi_s : \langle \mu, \alpha \rangle = p^m \hbox{ for some } m \geq 0 \}|.
$$
Given a $p$-restricted dominant weight $\lambda$, set
$$
s_{\lambda, p} = \sum m_\mu r_{\mu,p}, \qquad {s_{\lambda, p}}' = \sum m_\mu {r_\mu}',
$$
where each sum runs over the dominant weights $\mu \preceq \lambda$ and $m_\mu = \dim L(\lambda)_\mu \geq 0$.

\begin{prop}\label{prop: codim bound without using Premet}
Let $(G, \lambda, p)$ be a $p$-restricted large triple; write $V = L(\lambda)$ and assume $p \leq e(\Phi)$. Then
\begin{itemize}
\item[(i)] for all $r \in \P'$, $s \in G_{(r)}$ and $\kappa \in K^*$ we have $\codim V_\kappa(s) \geq s_{\lambda, p}$;
\item[(ii)] for all $\alpha \in \Phi_s$ we have $\codim C_V(x_\alpha(1)) \geq s_{\lambda, p}$;
\item[(iii)] for all $\beta \in \Phi_l$ we have $\codim C_V(x_\beta(1)) \geq {s_{\lambda, p}}'$.
\end{itemize}
\end{prop}

\begin{proof}
An argument very similar to that in the first paragraph of the proof of Proposition~\ref{prop: codim bound using Premet} shows that, for a fixed $\alpha \in \Phi_s$, the number of weights $\nu$ in $\Lambda(V)$ (counted with multiplicity) such that $\langle \nu, \alpha \rangle = p^m$ for some $m \geq 0$ is $s_{\lambda, p}$. Likewise, for a fixed $\beta \in \Phi_l$, the number of weights $\nu$ in $\Lambda(V)$ (counted with multiplicity) such that $\langle \nu, \beta \rangle > 0$ is ${s_{\lambda, p}}'$.

Take $r \in \P'$, $s \in G_{(r)}$ and $\kappa \in K^*$; as in the proof of Proposition~\ref{prop: codim bound using Premet} there exists $\alpha \in \Phi_s$ with $\alpha(s) \neq 1$. Observe that if $\nu \in \Lambda(V)$ satisfies $\langle \nu, \alpha \rangle = p^m$ for some $m \geq 0$, then the two weights $\nu$ and $w_\alpha(\nu) = \nu - p^m \alpha$ cannot both lie in $V_\kappa(s)$, since $\alpha(s) \neq 1 \implies (p^m \alpha)(s) \neq 1$. Summing over the weights in $\Lambda(V)$ (counted with multiplicity) we see that $\codim V_\kappa(s) \geq s_{\lambda, p}$, proving (i).

Now take $\alpha \in \Phi_s$; arguing again as in the proof of Proposition~\ref{prop: codim bound using Premet}, we see that $\codim C_V(x_\alpha(1))$ is at least half of the number of weights in $\Lambda(V)$ (counted with multiplicity) which are not orthogonal to $\alpha$, which equals the number of weights $\nu$ in $\Lambda(V)$ (counted with multiplicity) such that $\langle \nu, \alpha \rangle > 0$. This number is certainly at least as great as the number of weights $\nu$ in $\Lambda(V)$ (counted with multiplicity) such that $\langle \nu, \alpha \rangle = p^m$ for some $m \geq 0$, so we have $\codim C_V(x_\alpha(1)) \geq s_{\lambda, p}$, proving (ii).

Finally take $\beta \in \Phi_l$; an exactly similar argument shows that $\codim C_V(x_\beta(1))$ is at least half of the number of weights in $\Lambda(V)$ (counted with multiplicity) which are not orthogonal to $\beta$, which equals ${s_{\lambda, p}}'$, proving (iii).
\end{proof}

\begin{cor}\label{cor: TGS without using Premet}
Let $(G, \lambda, p)$ be a $p$-restricted large triple; assume $p \leq e(\Phi)$. If $s_{\lambda, p} > M$ and ${s_{\lambda, p}}' > M$, then the triple $(G, \lambda, p)$ satisfies $\ssddagcon$ and $\uddagcon$, and thus has TGS.
\end{cor}

\begin{proof}
This is immediate.
\end{proof}

Corollaries~\ref{cor: TGS using Premet} and \ref{cor: TGS without using Premet} will form the basis of our strategy for showing that all $p$-restricted large triples which are not listed in Table~\ref{table: large triple and first quadruple non-TGS} have TGS. We call a $p$-restricted large triple $(G, \lambda, p)$ {\em excluded\/} if it satisfies the conditions of Corollary~\ref{cor: TGS using Premet} or \ref{cor: TGS without using Premet} according as $p > e(\Phi)$ or $p \leq e(\Phi)$, and {\em unexcluded\/} otherwise. Following some preliminary work on subsystems and weights in Section \ref{sect: large triple relevance}, we shall determine the unexcluded $p$-restricted large triples $(G, \lambda, p)$ in Section~\ref{sect: large triple exclusion}. These unexcluded triples will then require further investigation in the sections which follow.

\section{Relevant subsystems and dominant weights}\label{sect: large triple relevance}

Let $\Psi$ be a proper standard subsystem of the irreducible root system $\Phi$. Then $\Psi$ will be called {\em relevant\/} if $r_\Psi \leq M$, or if $e(\Phi) > 1$ and ${r_\Psi}' \leq M$; it will be called {\em irrelevant\/} if it is not relevant. Observe that if $\Psi_1$ and $\Psi_2$ are standard subsystems of $\Phi$ with $\Psi_1 \subset \Psi_2$, then $r_{\Psi_1} > r_{\Psi_2}$, and if $e(\Phi) > 1$ then ${r_{\Psi_1}}' > {r_{\Psi_2}}'$; thus if $\Psi_2$ is irrelevant, so is $\Psi_1$.

In this section, for each irreducible root system $\Phi$ we shall identify its relevant subsystems, up to automorphisms of $\Phi$. Once this is done we shall consider the associated dominant weights for a simple algebraic group having root system $\Phi$.

\begin{prop}\label{prop: rel subsystems for A_ell}
Let $\Phi$ be of type $A_\ell$; then the relevant subsystems of $\Phi$ are as follows: $A_{\ell - 1}$ for $\ell \in [1, \infty)$; $A_1A_{\ell - 2}$ for $\ell \in [3, \infty)$; $A_2A_{\ell - 3}$ for $\ell \in [5, \infty)$; $A_3A_{\ell - 4}$ for $\ell \in [7, 11]$; $A_4A_4$ for $\ell = 9$; $A_{\ell - 2}$ for $\ell \in [2, \infty)$; $A_1A_{\ell - 3}$ for $\ell \in [4, 8]$; and $\emptyset$ for $\ell = 3$.
\end{prop}

\begin{proof}
We have $M = \ell(\ell + 1)$. Let $\Psi$ be a standard subsystem of $\Phi$ of corank $c$; then we may write $\Psi = A_{j_1 - 1}A_{j_2 - 1} \dots A_{j_{c + 1} - 1}$, where $\sum j_i = \ell + 1$ and $1 \leq j_1 \leq j_2 \leq \cdots$.

First suppose $c = 1$; then $r_\Psi = \binom{\ell - 1}{j_1 - 1}$. If $j_1 \geq 6$ we have $r_\Psi \geq \binom{\ell - 1}{5} > M$; if $j_1 = 5$ we have $r_\Psi = \binom{\ell - 1}{4} > M$ for $\ell \geq 10$; if $j_1 = 4$ we have $r_\Psi = \binom{\ell - 1}{3} > M$ for $\ell \geq 12$; in all other cases we have $r_\Psi \leq M$. Thus the relevant subsystems of corank $1$ are as stated.

Next suppose $c = 2$. If $j_1 \geq 2$ then $r_\Psi \geq r_{A_1A_1A_{\ell - 4}} = (\ell - 1)(\ell - 2)^2 > M$; if $j_1 = 1$ and $j_2 \geq 3$ then $r_\Psi \geq r_{A_2A_{\ell - 4}} = \frac{1}{6}(\ell - 1)(\ell - 2)(4\ell - 9) > M$; if $j_1 = 1$ and $j_2 = 2$ then $r_\Psi = \frac{1}{2}(\ell - 1)(3\ell - 4) > M$ for $\ell \geq 9$; in all other cases we have $r_\Psi \leq M$. Thus the relevant subsystems of corank $2$ are as stated.

Finally suppose $c \geq 3$. If $j_3 \geq 2$ then $r_\Psi \geq r_{A_1A_{\ell - 4}} = \frac{1}{2}(\ell - 1)(\ell - 2)(4\ell - 7) > M$; if $j_1 = j_2 = j_3 = 1$ then $r_\Psi \geq r_{A_{\ell - 3}} = 3(\ell - 1)^2 > M$ for $\ell \geq 4$; the only other case is $\Psi = \emptyset$ for $\ell = 3$, for which $r_\Psi = 12 = M$. The result follows.
\end{proof}

\begin{prop}\label{prop: rel subsystems for D_ell}
Let $\Phi$ be of type $D_\ell$; then the relevant subsystems of $\Phi$ are as follows: $D_{\ell - 1}$, $A_1D_{\ell - 2}$ and $D_{\ell - 2}$ for $\ell \in [4, \infty)$; $A_2D_{\ell - 3}$ for $\ell \in [5, 6]$; $A_{\ell - 1}$ for $\ell \in [5, 10]$; and $A_{\ell - 2}$ for $\ell \in [4, 5]$.
\end{prop}

\begin{proof}
We have $M = 2\ell(\ell - 1)$. Let $\Psi$ be a standard subsystem of $\Phi$ of corank $c$.

First suppose $c = 1$; then $\Psi$ is either $A_{j - 1}D_{\ell - j}$ for some $j \in [1, \ell - 2]$, or $A_{\ell - 1}$. We have $r_{A_{\ell - 1}} = 2^{\ell - 3} > M$ for $\ell \geq 11$. For $j \in [1, \ell - 2]$ set $f(j) = r_{A_{j - 1}D_{\ell - j}} = 2^{j - 2}\binom{\ell - 1}{j - 1}\frac{4\ell - 3j - 1}{\ell - 1}$; then if $j < \ell - 2$ we have $f(j + 1)/f(j) = \frac{2(\ell - j)(4\ell - 3j - 4)}{j(4\ell - 3j - 1)}$, and we find that $f(j + 1)/f(j) > 1$ if and only if $j < \frac{1}{3}(2\ell - 1)$. Thus as $j$ runs from $1$ to $\ell - 2$ the values $f(j)$ increase to a maximum at $j = \lceil \frac{1}{3}(2\ell - 1) \rceil$ (provided $\ell \geq 5$) and then decrease. If $j = 1$ we have $f(j) = 2$; if $j = 2$ we have $f(j) = 4\ell - 7$; if $j = 3$ we have $f(j) = 2(\ell - 2)(2\ell - 5) > M$ for $\ell \geq 7$; if $j = 4$ we have $f(j) = \frac{2}{3}(\ell - 2)(\ell - 3)(4\ell - 13) > M$; if $j = \ell - 2$ we have $f(j) = 2^{\ell - 5}(\ell - 2)(\ell + 5) > M$ for $\ell \geq 6$. Thus the relevant subsystems of corank $1$ are as stated.

Next suppose $c = 2$. From the previous paragraph, we see that the only $\Psi$ which do not lie in an irrelevant subsystem of corank $1$ are the following: $D_{\ell - 2}$ for $\ell \in [4, \infty)$; $A_{\ell - 2}$ and $A_1A_{\ell - 3}$ for $\ell \in [4, 10]$; $A_1D_{\ell - 3}$ for $\ell \in [5, 6]$; and $A_2A_2$ for $\ell = 6$. We have $r_{D_{\ell - 2}} = 4(2\ell - 3)$; $r_{A_{\ell - 2}} = 2^{\ell - 3}(\ell + 2) > M$ for $\ell \geq 6$; $r_{A_1A_{\ell - 3}} = 2^{\ell - 4}(\ell^2 + 3\ell - 8) > M$ for $\ell \geq 5$; $r_{A_1D_{\ell - 3}} = 2(\ell - 2)(6\ell - 13) > M$; and $r_{A_2A_2} = 256 > M$ for $\ell = 6$. Thus the relevant subsystems of corank $2$ are as stated.

Finally suppose $c \geq 3$. If $\ell \geq 6$, the previous paragraph shows that $\Psi$ lies in an irrelevant subsystem of corank $2$. If $\ell = 5$, $\Psi$ must lie in $A_2$, $A_1A_1$ or $D_2$, and hence lies in $A_2A_1$ or $A_1D_2$, both of which are irrelevant. If $\ell = 4$, we have $r_{A_1} = 44 > M$. The result follows.
\end{proof}

\begin{prop}\label{prop: rel subsystems for B_ell}
Let $\Phi$ be of type $B_\ell$; then the relevant subsystems of $\Phi$ are as follows: $B_{\ell - 1}$ for $\ell \in [2, \infty)$; $A_1B_{\ell - 2}$ and $B_{\ell - 2}$ for $\ell \in [3, \infty)$; $A_2B_{\ell - 3}$ for $\ell \in [4, \infty)$; $A_{\ell - 1}$ for $\ell \in [2, 9]$; $A_{\ell - 2}B_1$ for $\ell \in [5, 6]$; and $A_{\ell - 2}$ for $\ell \in [2, 4]$.
\end{prop}

\begin{proof}
We have $M = 2\ell^2$. Let $\Psi$ be a standard subsystem of $\Phi$ of corank $c$.

First suppose $c = 1$; then $\Psi$ is $A_{j - 1}B_{\ell - j}$ for some $j \in [1, \ell]$. For $j \in [1, \ell]$ set $g(j) = r_{A_{j - 1}B_{\ell - j}} = 2^{j - 1}\binom{\ell - 1}{j - 1}$, and $f(j) = {r_{A_{j - 1}B_{\ell - j}}}' = 2^{j - 2}\binom{\ell - 1}{j - 1}\frac{4\ell - 3j - 1}{\ell - 1}$. Then $f$ is the same function as in the proof of Proposition~\ref{prop: rel subsystems for D_ell}, so as $j$ runs from $1$ to $\ell$ the values $f(j)$ increase to a maximum at $j = \lceil \frac{1}{3}(2\ell - 1) \rceil$ and then decrease. Similarly we see that $g(j + 1)/g(j) = \frac{2(\ell - j)}{j}$, so $g(j + 1)/g(j) > 1$ if and only if $j < \frac{2\ell}{3}$, and hence as $j$ runs from $1$ to $\ell$ the values $f(j)$ increase to a maximum at $j = \lceil \frac{2\ell}{3} \rceil$ and then decrease. If $j = 1$ we have $g(j) = 1$, and $f(j) = 2$; if $j = 2$ we have $g(j) = 2(\ell - 1)$, and $f(j) = 4\ell - 7$; if $j = 3$ we have $g(j) = 2(\ell - 1)(\ell - 2)$, and $f(j) = 2(\ell - 2)(2\ell - 5) > M$ for $\ell \geq 8$; if $j = 4$ we have $g(j) = \frac{4}{3}(\ell - 1)(\ell - 2)(\ell - 3) > M$ for $\ell \geq 6$, and $f(j) = \frac{2}{3}(\ell - 2)(\ell - 3)(4\ell - 13) > M$ for $\ell \geq 6$; if $j = \ell - 1$ we have $g(j) = 2^{\ell - 2}(\ell - 1) > M$ for $\ell \geq 6$, and $f(j) = 2^{\ell - 3}(\ell + 2) > M$ for $\ell \geq 7$; if $j = \ell$ we have $g(j) = 2^{\ell - 1} > M$ for $\ell \geq 9$, and $f(j) = 2^{\ell - 2} > M$ for $\ell \geq 10$. Thus the relevant subsystems of corank $1$ are as stated.

Next suppose $c = 2$. From the previous paragraph, we see that the only $\Psi$ which do not lie in an irrelevant subsystem of corank $1$ are the following: $B_{\ell - 2}$ for $\ell \in [3, \infty)$; $A_1B_{\ell - 3}$ for $\ell \in [4, \infty)$; $A_{\ell - 2}$ for $\ell \in [2, 9]$; $A_1A_{\ell - 3}$ for $\ell \in [4, 9]$; $A_2A_{\ell - 4}$ for $\ell \in [6, 9]$; and $A_{\ell - 3}B_1$ and $A_1A_{\ell - 4}B_1$ for $\ell \in [5, 6]$. We have $r_{B_{\ell - 2}} = 4(\ell - 1)$, and ${r_{B_{\ell - 2}}}' = 4(2\ell - 3)$; $r_{A_1B_{\ell - 3}} = 6(\ell - 1)(\ell - 2) > M$, and ${r_{A_1B_{\ell - 3}}}' = 2(\ell - 2)(6\ell - 13) > M$; $r_{A_{\ell - 2}} = 2^{\ell - 1}\ell > M$ for $\ell \geq 5$, and ${r_{A_{\ell - 2}}}' = 2^{\ell - 2}(\ell + 2) > M$ for $\ell \geq 5$; $r_{A_1A_{\ell - 3}} = 2^{\ell - 2}\ell(\ell - 1) > M$, and ${r_{A_1A_{\ell - 3}}}' = 2^{\ell - 3}(\ell^2 + 3\ell - 8) > M$; $r_{A_2A_{\ell - 4}} = \frac{1}{3}.2^{\ell - 2}\ell(\ell - 1)(\ell - 2) > M$, and ${r_{A_2A_{\ell - 4}}}' = \frac{1}{3}.2^{\ell - 3}(\ell - 2)(\ell^2 + 5\ell - 18) > M$; $r_{A_{\ell - 3}B_1} = 2^{\ell - 2}(\ell - 1)^2 > M$, and ${r_{A_{\ell - 3}B_1}}' = 2^{\ell - 3}(\ell^2 + 3\ell - 6) > M$; and $r_{A_1A_{\ell - 4}B_1} = 2^{\ell - 3}(\ell - 1)^2(\ell - 2) > M$, and ${r_{A_1A_{\ell - 4}B_1}}' = 2^{\ell - 4}(\ell - 2)^2(\ell + 7) > M$. Thus the relevant subsystems of corank $2$ are as stated.

Finally suppose $c \geq 3$. From the previous paragraph, we see that the only $\Psi$ which does not lie in an irrelevant subsystem of corank $2$ is $\emptyset$ for $\ell = 3$, and $r_\emptyset = {r_\emptyset}' = |W|/2 = 24 > M$. The result follows.
\end{proof}

\begin{prop}\label{prop: rel subsystems for C_ell}
Let $\Phi$ be of type $C_\ell$; then the relevant subsystems of $\Phi$ are as follows: $C_{\ell - 1}$, $A_1C_{\ell - 2}$ and $C_{\ell - 2}$ for $\ell \in [3, \infty)$; $A_2C_{\ell - 3}$ for $\ell \in [4, \infty)$; $A_{\ell - 1}$ for $\ell \in [3, 9]$; $A_{\ell - 2}C_1$ for $\ell \in [5, 6]$; and $A_{\ell - 2}$ for $\ell \in [3, 4]$.
\end{prop}

\begin{proof}
The proof may be obtained from that of Proposition~\ref{prop: rel subsystems for B_ell} by interchanging the values of $r_\Psi$ and ${r_\Psi}'$, and replacing each root system $B_r$ with $C_r$; this is because doubling the length of every short root in any root system of type $B_r$ gives a root system of type $C_r$.
\end{proof}

\begin{prop}\label{prop: rel subsystems for exceptional groups}
Let $\Phi$ be of exceptional type; then the relevant subsystems of $\Phi$ are as follows:
\begin{itemize}
\item[(i)] $D_5$ and $A_5$ if $\Phi$ is of type $E_6$;
\item[(ii)] $D_6$ and $E_6$ if $\Phi$ is of type $E_7$;
\item[(iii)] $E_7$ if $\Phi$ is of type $E_8$;
\item[(iv)] $C_3$, $\tilde A_2 A_1$, $A_2\tilde A_1$ and $B_3$ if $\Phi$ is of type $F_4$;
\item[(v)] $A_1$, $\tilde A_1$ and $\emptyset$ if $\Phi$ is of type $G_2$.
\end{itemize}
\end{prop}

\begin{proof} We have $M = 72$, $126$, $240$, $48$ or $12$ according as $\Phi$ is of type $E_6$, $E_7$, $E_8$, $F_4$ or $G_2$.

(i) If $\Phi$ is of type $E_6$, for the standard subsystems of corank $1$ we have $r_{D_5} = 6$, $r_{A_5} = 21$, $r_{A_4A_1} = 75 > M$ and $r_{A_2A_2A_1} = 290 > M$. The only standard subsystem of corank $2$ which does not lie in either $A_4A_1$ or $A_2A_2A_1$ is $D_4$, and we have $r_{D_4} = 90 > M$. The result follows.

(ii) If $\Phi$ is of type $E_7$, for the standard subsystems of corank $1$ we have $r_{D_6} = 33$, $r_{A_6} = 192 > M$, $r_{A_5A_1} = 752 > M$, $r_{A_3A_2A_1} = 4240 > M$, $r_{A_4A_2} = 1600 > M$, $r_{D_5A_1} = 252 > M$ and $r_{E_6} = 12$. Since any standard subsystem of corank $2$ then lies in an irrelevant subsystem of corank $1$, the result follows.

(iii) If $\Phi$ is of type $E_8$, for the standard subsystems of corank $1$, we have $r_{D_7} = 702 > M$, $r_{A_7} = 6624 > M$, $r_{A_6A_1} = 28224 > M$, $r_{A_4A_2A_1} = 213696 > M$, $r_{A_4A_3} = 104832 > M$, $r_{D_5A_2} = 24444 > M$, $r_{E_6A_1} = 2324 > M$ and $r_{E_7} = 57$. Since any standard subsystem of corank $2$ then lies in an irrelevant subsystem of corank $1$, the result follows.

(iv) If $\Phi$ is of type $F_4$, for the standard subsystems $\Psi$ of corank $1$, we have $(r_\Psi, {r_\Psi}') = (6, 9)$, $(36, 44)$, $(44, 36)$ and $(9, 6)$ for $\Psi = C_3$, $\tilde A_2 A_1$, $A_2\tilde A_1$ and $B_3$ respectively. For those of corank $2$, we have $(r_\Psi, {r_\Psi}') = (72, 96)$, $(96, 72)$, $(60, 60)$ and $(132, 132)$ for $\Psi = \tilde A_2$, $A_2$, $B_2$, $A_1 \tilde A_1$ respectively. The result follows.

(v) If $\Phi$ is of type $G_2$, we have $(r_\Psi, {r_\Psi}') = (3, 2)$, $(2, 3)$ and $(6, 6)$ for $\Psi = A_1$, $\tilde A_1$ and $\emptyset$ respectively. The result follows.
\end{proof}

We have thus identified the relevant subsystems for each irreducible root system $\Phi$. Now recall that $\Phi$ is the root system of the simple algebraic group $G$ over an algebraically closed field of characteristic $p$. A non-zero dominant weight $\mu$ of $G$ will be called {\em irrelevant\/} or {\em relevant\/} according as the corresponding standard subsystem $\Psi(\mu)$ is irrelevant or relevant; thus $\mu$ is relevant if $r_\mu \leq M$, or if $e(\Phi) > 1$ and ${r_\mu}' \leq M$. It is now a simple matter to identify the relevant dominant weights.

\begin{prop}\label{prop: relevant dominant weights}
Let $G$ be a simple algebraic group; then the relevant dominant weights for $G$ are as listed in Table~\ref{table: relevant dominant weights}.
\end{prop}

\begin{proof}
This is immediate from Propositions~\ref{prop: rel subsystems for A_ell}--\ref{prop: rel subsystems for exceptional groups}.
\end{proof}

\begin{table}
\caption{Relevant dominant weights}\label{table: relevant dominant weights}
\tabcapsp
$$
\begin{array}{|c|c|c|c|c|c|c|c|c|c|}
\cline{1-3} \cline{5-7} \cline{9-10}
G      & \mu                               & \ell         & \ptw & G      & \mu                               & \ell         & \ptw & G      & \mu                   \tbs \\
\cline{1-3} \cline{5-7} \cline{9-10}
A_\ell & a\omega_1                         & {} \geq 1    &      & B_\ell & a\omega_1                         & {} \geq 2    &      & E_6    & a\omega_1             \tbs \\
       & a\omega_2                         & {} \geq 3    &      &        & a\omega_2                         & {} \geq 3    &      &        & a\omega_2             \tbs \\
\cline{9-10}
       & a\omega_3                         & {} \geq 5    &      &        & a\omega_3                         & {} \geq 4    &      & E_7    & a\omega_1             \tbs \\
       & a\omega_4                         & 7, \dots, 11 &      &        & a\omega_{\ell - 1}                & 5, 6         &      &        & a\omega_7             \tbs \\
\cline{9-10}
       & a\omega_5                         & 9            &      &        & a\omega_\ell                      & 2, \dots, 9  &      & E_8    & a\omega_8             \tbs \\
\cline{9-10}
       & a\omega_1 + b\omega_\ell          & {} \geq 2    &      &        & a\omega_1 + b\omega_2             & {} \geq 3    &      & F_4    & a\omega_1             \tbs \\
       & a\omega_1 + b\omega_2             & {} \geq 3    &      &        & a\omega_1 + b\omega_\ell          & 2, 3, 4      &      &        & a\omega_2             \tbs \\
       & a\omega_2 + b\omega_\ell          & 4, \dots, 8  &      &        & a\omega_{\ell - 1} + b\omega_\ell & 3, 4         &      &        & a\omega_3             \tbs \\
\cline{5-7}
       & a\omega_1 + b\omega_3             & 5, \dots, 8  &      & C_\ell & a\omega_1                         & {} \geq 3    &      &        & a\omega_4             \tbs \\
\cline{9-10}
       & a\omega_2 + b\omega_3             & 4, \dots, 8  &      &        & a\omega_2                         & {} \geq 3    &      & G_2    & a\omega_1             \tbs \\
       & a\omega_1 + b\omega_2 + c\omega_3 & 3            &      &        & a\omega_3                         & {} \geq 4    &      &        & a\omega_2             \tbs \\
\cline{1-3}
D_\ell & a\omega_1                         & {} \geq 4    &      &        & a\omega_{\ell - 1}                & 5, 6         &      &        & a\omega_1 + b\omega_2 \tbs \\
\cline{9-10}
       & a\omega_2                         & {} \geq 4    &      &        & a\omega_\ell                      & 3, \dots, 9  & \multicolumn{3}{c}{}                  \tbs \\
       & a\omega_3                         & 5, 6         &      &        & a\omega_1 + b\omega_2             & {} \geq 3    & \multicolumn{3}{c}{}                  \tbs \\
       & a\omega_\ell                      & 5, \dots, 10 &      &        & a\omega_1 + b\omega_\ell          & 3, 4         & \multicolumn{3}{c}{}                  \tbs \\
       & a\omega_1 + b\omega_2             & {} \geq 4    &      &        & a\omega_{\ell - 1} + b\omega_\ell & 3, 4         & \multicolumn{3}{c}{}                  \tbs \\
\cline{5-7}
       & a\omega_1 + b\omega_\ell          & 4, 5         & \multicolumn{7}{c}{}                                                                                     \tbs \\
       & a\omega_4 + b\omega_5             & 5            & \multicolumn{7}{c}{}                                                                                     \tbs \\
\cline{1-3}
\end{array}
$$
\end{table}

Note that in Table~\ref{table: relevant dominant weights} the symbols $a$, $b$ and $c$ stand for arbitrary natural numbers; in particular there is no requirement that a relevant dominant weight be $p$-restricted. In Section~\ref{sect: large triple exclusion} we shall use Table~\ref{table: relevant dominant weights} to determine unexcluded $p$-restricted large triples $(G, \lambda, p)$ with $p > e(\Phi)$.

In the remainder of this section we shall assume that $\Phi$, $G$ and $p$ are as above but with $p \leq e(\Phi)$. A non-zero $p$-restricted dominant weight $\mu$ of $G$ will be called {\em $p$-relevant\/} if at least one of $r_{\mu, p} \leq M$ and ${r_\mu}' \leq M$ holds. We shall prove that the $p$-relevant dominant weights for $G$ are as listed in Table~\ref{table: p-relevant dominant weights}; this will also be used in Section~\ref{sect: large triple exclusion} to determine the corresponding unexcluded $p$-restricted large triples.

Recall that we define $r_{\mu, p} = \frac{|W.\mu|}{|\Phi_s|}.| \{ \alpha \in \Phi_s : \langle \mu, \alpha \rangle = p^m \hbox{ for some } m \geq 0 \} |$. Thus $r_{\mu, p} \leq r_\mu$, with equality precisely if $\mu$ is such that the only positive values $\langle \mu, \alpha \rangle$ for $\alpha \in \Phi_s$ are powers of $p$. In particular, any $p$-restricted dominant weight which is relevant is $p$-relevant. Thus to determine the $p$-relevant dominant weights it suffices to consider the values $r_{\mu, p}$ for weights $\mu$ such that for some $\alpha \in \Phi_s$ the value $\langle \mu, \alpha \rangle$ is positive and not a power of $p$; note that there is no need to consider $G = G_2$, since in this case all dominant weights are relevant.

\begin{table}
\caption{$p$-relevant dominant weights}\label{table: p-relevant dominant weights}
\tabcapsp
$$
\begin{array}{|c|c|c|c|c|c|c|c|}
\cline{1-4} \cline{6-8}
G      & \mu                             & \ell        & p & \ptw & G   & \mu                   & p    \tbs \\
\cline{1-4} \cline{6-8}
B_\ell & \omega_1                        & {} \geq 2   & 2 &      & F_4 & \omega_1              & 2    \tbs \\
       & \omega_2                        & {} \geq 3   & 2 &      &     & \omega_2              & 2    \tbs \\
       & \omega_3                        & {} \geq 4   & 2 &      &     & \omega_3              & 2    \tbs \\
       & \omega_{\ell - 1}               & 5, 6        & 2 &      &     & \omega_4              & 2    \tbs \\
       & \omega_\ell                     & 2, \dots, 9 & 2 &      &     & \omega_1 + \omega_2   & 2    \tbs \\
       & \omega_1 + \omega_2             & {} \geq 3   & 2 &      &     & \omega_1 + \omega_4   & 2    \tbs \\
       & \omega_1 + \omega_\ell          & 2, 3, 4     & 2 &      &     & \omega_2 + \omega_3   & 2    \tbs \\
       & \omega_2 + \omega_4             & 4           & 2 &      &     & \omega_2 + \omega_4   & 2    \tbs \\
\cline{6-8}
       & \omega_{\ell - 1} + \omega_\ell & 3, \dots, 8 & 2 &      & G_2 & \omega_1              & 2, 3 \tbs \\
       & \omega_1 + \omega_2 + \omega_3  & 3           & 2 &      &     & \omega_2              & 2, 3 \tbs \\
       & \omega_1 + \omega_3 + \omega_4  & 4           & 2 &      &     & \omega_1 + \omega_2   & 2, 3 \tbs \\
       & \omega_2 + \omega_3 + \omega_4  & 4           & 2 &      &     & 2\omega_1             & 3    \tbs \\
\cline{1-4}
C_\ell & \omega_1                        & {} \geq 3   & 2 &      &     & 2\omega_2             & 3    \tbs \\
       & \omega_2                        & {} \geq 3   & 2 &      &     & 2\omega_1 + \omega_2  & 3    \tbs \\
       & \omega_3                        & {} \geq 4   & 2 &      &     & \omega_1 + 2\omega_2  & 3    \tbs \\
       & \omega_{\ell - 1}               & 5, 6        & 2 &      &     & 2\omega_1 + 2\omega_2 & 3    \tbs \\
\cline{6-8}
       & \omega_\ell                     & 3, \dots, 9 & 2 & \multicolumn{4}{c}{}                      \tbs \\
       & \omega_1 + \omega_2             & {} \geq 3   & 2 & \multicolumn{4}{c}{}                      \tbs \\
       & \omega_1 + \omega_\ell          & 3, 4, 5     & 2 & \multicolumn{4}{c}{}                      \tbs \\
       & \omega_2 + \omega_4             & 4           & 2 & \multicolumn{4}{c}{}                      \tbs \\
       & \omega_{\ell - 1} + \omega_\ell & 3, 4, 5     & 2 & \multicolumn{4}{c}{}                      \tbs \\
       & \omega_1 + \omega_2 + \omega_3  & 3           & 2 & \multicolumn{4}{c}{}                      \tbs \\
\cline{1-4}
\end{array}
$$
\end{table}

\begin{prop}\label{prop: p-relevant dominant weights for B_ell}
Let $G = B_\ell$ and $p = 2$; then the $2$-relevant dominant weights for $G$ are as listed in Table~\ref{table: p-relevant dominant weights}.
\end{prop}

\begin{proof}
We have $M = 2\ell^2$. Let $\mu$ be a $2$-restricted dominant weight; then $\mu = \sum a_i\omega_i$ with each $a_i \in \{ 0, 1 \}$.

If $\mu = \omega_\ell$, $\omega_j$ for some $j < \ell$, or $\omega_i + \omega_j$ for some $i < j < \ell$, then all positive values $\langle \mu, \alpha \rangle$ for $\alpha \in \Phi_s$ are $1$, $2$, or either $2$ or $4$ respectively; thus for these weights we have $r_{\mu, 2} = r_\mu$ and there is no need to consider them further.

First suppose $\mu = \omega_i + \omega_\ell$ for some $i < \ell$. Then $\langle \mu, \alpha \rangle \in \{ \pm1, \pm3 \}$ for $\alpha \in \Phi_s$, and there are $\ell - i$ short roots $\alpha$ with $\langle \mu, \alpha \rangle = 1$; it follows that $r_{\mu, 2} = 2^{\ell - 1}\binom{\ell - 1}{i}$. Thus if $i < \ell - 1$ we have $r_{\mu, 2} \geq 2^{\ell - 1}(\ell - 1) > M$ for $\ell \geq 5$, while if $i = \ell - 1$ we have $r_{\mu, 2} = 2^{\ell - 1} > M$ for $\ell \geq 9$.

Next suppose $\mu = \omega_h + \cdots + \omega_i + \omega_j$ with $h < \cdots < i < j < \ell$, so that $\ell \geq 4$. Then there exist $\alpha, \alpha' \in \Phi_s$ with $\langle \mu, \alpha \rangle = 2$ and $\langle \mu, \alpha' \rangle = 4$; since $|W.\mu| \geq 2^j\binom{\ell}{j}\binom{j}{i}\binom{i}{h}$ we have $r_{\mu, 2} \geq \frac{1}{\ell}.2^j\binom{\ell}{j}\binom{j}{i}\binom{i}{h} \geq \frac{1}{\ell}.2^j\binom{\ell}{j}\binom{j}{i}i = \frac{1}{\ell}.2^j\binom{\ell}{j}\binom{j - 1}{i - 1}j = 2^j\binom{\ell - 1}{j - 1}\binom{j - 1}{i - 1} \geq 2^j\binom{\ell - 1}{j - 1}(j - 1) = 2^j(\ell - 1)\binom{\ell - 2}{j - 2} \geq 8(\ell - 1)(\ell - 2) > M$.

Finally suppose $\mu = \omega_i + \cdots + \omega_j + \omega_\ell$ with $i < \cdots < j < \ell$, so that $\ell \geq 3$. Then there exists $\alpha \in \Phi_s$ with $\langle \mu, \alpha \rangle = 1$; since $|W.\mu| \geq 2^\ell\binom{\ell}{j}\binom{j}{i}$ we have $r_{\mu, 2} \geq \frac{1}{\ell}.2^{\ell - 1}\binom{\ell}{j}\binom{j}{i} \geq \frac{1}{\ell}.2^{\ell - 1}\binom{\ell}{j}j = 2^{\ell - 1}\binom{\ell - 1}{j - 1} \geq 2^{\ell - 1}(\ell - 1) > M$ for $\ell \geq 5$. For $\ell = 4$ we have $r_{\mu, 2} = 48 > M$ for $\mu = \omega_1 + \omega_2 + \omega_4$ or $\omega_1 + \omega_2 + \omega_3 + \omega_4$, while $r_{\mu, 2} = 24$ for $\mu = \omega_1 + \omega_3 + \omega_4$ or $\omega_2 + \omega_3 + \omega_4$. For $\ell = 3$ we have $r_{\mu, 2} = 8$ for $\mu = \omega_1 + \omega_2 + \omega_3$. The result follows.
\end{proof}

\begin{prop}\label{prop: p-relevant dominant weights for C_ell}
Let $G = C_\ell$ and $p = 2$; then the $2$-relevant dominant weights for $G$ are as listed in Table~\ref{table: p-relevant dominant weights}.
\end{prop}

\begin{proof}
We have $M = 2\ell^2$. Let $\mu$ be a $2$-restricted dominant weight; then $\mu = \sum a_i\omega_i$ with each $a_i \in \{ 0, 1 \}$.

If $\mu = \omega_\ell$, or $\omega_j$ for some $j < \ell$, then all positive values $\langle \mu, \alpha \rangle$ for $\alpha \in \Phi_s$ are $2$, or either $1$ or $2$ respectively; thus for these weights we have $r_{\mu, 2} = r_\mu$ and there is no need to consider them further.

First suppose $\mu = \omega_i + \omega_j$ for some $i < j < \ell$. Then $\langle \mu, \alpha \rangle \in \{ 0, \pm1, \pm2, \pm3, \pm4 \}$ for $\alpha \in \Phi_s$, and the numbers of short roots $\alpha$ with $\langle \mu, \alpha \rangle = 1$, $2$ and $4$ are $i(j - i) + 2(j - i)(\ell - j)$, $2i(\ell - j) + \frac{1}{2}(j - i)(j - i - 1)$ and $\frac{1}{2}i(i - 1)$ respectively; it follows that $r_{\mu, 2} = 2^{j - 2}\binom{\ell - 1}{j - 1}\binom{j}{i}\frac{4\ell - 3j - 1}{\ell - 1} = \binom{j}{i}f(j)$, where $f$ is the function appearing in the proofs of Propositions~\ref{prop: rel subsystems for D_ell} and \ref{prop: rel subsystems for B_ell}. We have $3f(3) = 6(\ell  - 2)(2\ell - 5) > M$ for $\ell \geq 4$ and $3f(\ell - 1) = 3.2^{\ell - 3}(\ell + 2) > M$ for $\ell \geq 4$; thus if $j \geq 3$ we have $r_{\mu, 2} \geq jf(j) \geq 3f(j) > M$.

Next suppose $\mu = \omega_i + \omega_\ell$ for some $i < \ell$. Then $\langle \mu, \alpha \rangle \in \{ 0, \pm1, \pm2, \pm3, \pm4 \}$ for $\alpha \in \Phi_s$, and the numbers of short roots $\alpha$ with $\langle \mu, \alpha \rangle = 1$, $2$ and $4$ are $i(\ell - i)$, $\frac{1}{2}(\ell - i)(\ell - i - 1)$ and $\frac{1}{2}i(i - 1)$ respectively; it follows that $r_{\mu, 2} = 2^{\ell - 2}\binom{\ell}{i}$. Thus $r_{\mu, 2} \geq 2^{\ell - 2}\ell > M$ for $\ell \geq 6$; for $\ell = 5$ we have $r_{\mu, 2} = 80 > M$ if $i \in \{ 2, 3 \}$.

Next suppose $\mu = \omega_1 + \omega_2 + \omega_3$ and $\ell \geq 4$. Then $\langle \mu, \alpha \rangle \in \{ 0, \pm1, \pm2, \pm3, \pm4, \pm5 \}$ for $\alpha \in \Phi_s$, and the numbers of short roots $\alpha$ with $\langle \mu, \alpha \rangle = 1$, $2$ and $4$ are $2(\ell - 2)$, $2\ell - 5$ and $1$ respectively; it follows that $r_{\mu, 2} = 16(\ell - 2)^2 > M$.

Next suppose $\mu = \omega_h + \cdots + \omega_i + \omega_j$ with $h < \cdots < i < j < \ell$ and $j \geq 4$. There are at least $4$ short roots $\alpha$ with $\langle \mu, \alpha \rangle = 1$ (namely $\alpha_h$, $\alpha_i$, $\alpha_j$ and $\alpha_j + \cdots + \alpha_\ell$), and at least $2$ short roots $\alpha$ with $\langle \mu, \alpha \rangle = 2$ (namely $\alpha_i + \cdots + \alpha_j$ and $\alpha_i + \cdots + \alpha_j + \cdots + \alpha_\ell$); since $|W.\mu| \geq 2^j\binom{\ell}{j}\binom{j}{i}\binom{i}{h}$ we have $r_{\mu, 2} \geq \frac{3}{\ell(\ell - 1)}.2^j\binom{\ell}{j}\binom{j}{i}\binom{i}{h} \geq \frac{3}{\ell(\ell - 1)}.2^j\binom{\ell}{j}\binom{j}{i}i = \frac{3}{\ell(\ell - 1)}.2^j\binom{\ell}{j}\binom{j - 1}{i - 1}j = \frac{3}{\ell - 1}.2^j\binom{\ell - 1}{j - 1}\binom{j - 1}{i - 1} \geq \frac{3}{\ell - 1}.2^j\binom{\ell - 1}{j - 1}(j - 1) = 3.2^j\binom{\ell - 2}{j - 2}$. If $j = \ell - 1$ we have $r_{\mu, 2} \geq 3.2^{\ell - 1}(\ell - 2) > M$; if instead $j < \ell - 1$ we have $r_{\mu, 2} \geq 3.2^{j - 1}(\ell - 2)(\ell - 3) \geq 24(\ell - 2)(\ell - 3) > M$.

Finally suppose $\mu = \omega_i + \cdots + \omega_j + \omega_\ell$ with $i < \cdots < j < \ell$. There are at least $2$ short roots $\alpha$ with $\langle \mu, \alpha \rangle = 1$ (namely $\alpha_i$ and $\alpha_j$) and at least $2$ short roots $\alpha$ with $\langle \mu, \alpha \rangle = 2$ (namely $\alpha_\ell$ and some root $\alpha_i + \alpha_{i + 1} + \cdots$); since $|W.\mu| \geq 2^\ell\binom{\ell}{j}\binom{j}{i}$ we have $r_{\mu, 2} \geq \frac{2}{\ell(\ell - 1)}.2^\ell\binom{\ell}{j}\binom{j}{i} \geq \frac{2}{\ell(\ell - 1)}.2^\ell\binom{\ell}{j}j = 2^{\ell + 1}\binom{\ell - 2}{j - 2}(j - 1)$. If $j \geq 3$ (so that $\ell \geq 4$) we have $r_{\mu, 2} \geq 2^{\ell + 1}(\ell - 2).2 > M$; if instead $j = 2$ we have $r_{\mu, 2} \geq 2^{\ell + 1} > M$ for $\ell \geq 5$, while for $\ell = 4$ we have $r_{\mu, 2} = 64 > M$ and for $\ell = 3$ we have $r_{\mu, 2} = 16$. The result follows.
\end{proof}

\begin{prop}\label{prop: p-relevant dominant weights for F_4}
Let $G = F_4$ and $p = 2$; then the $2$-relevant dominant weights for $G$ are as listed in Table~\ref{table: p-relevant dominant weights}.
\end{prop}

\begin{proof}
We have $M = 48$. The weights $\omega_j$ for $1 \leq j \leq 4$ are all relevant; if $\mu = \omega_1 + \omega_4$ we have $r_{\mu, 2} = 36$, while if $\mu = \omega_1 + \omega_2$, $\omega_2 + \omega_3$ or $\omega_2 + \omega_4$ we have $r_{\mu, 2} = 48$. If $\mu = \omega_1 + \omega_3$ we have $r_{\mu, 2} = 84 > M$; if $\mu = \omega_3 + \omega_4$ we have $r_{\mu, 2} = 64 > M$; if $\mu = \omega_1 + \omega_2 + \omega_3$ we have $r_{\mu, 2} = 96 > M$; if $\mu = \omega_1 + \omega_2 + \omega_4$ we have $r_{\mu, 2} = 120 > M$; if $\mu = \omega_1 + \omega_3 + \omega_4$ or $\omega_2 + \omega_3 + \omega_4$ we have $r_{\mu, 2} = 144 > M$; and if $\mu = \omega_1 + \omega_2 + \omega_3 + \omega_4$ we have $r_{\mu, 2} = 192 > M$. The result follows.
\end{proof}

This completes the determination of relevant and $p$-relevant dominant weights for $G$.

\section{Exclusion of triples}\label{sect: large triple exclusion}

In this section we shall build upon the work of the previous section to determine the unexcluded $p$-restricted large triples; we list these in Table~\ref{table: unexcluded triples}. We begin with triples $(G, \lambda, p)$ with $p > e(\Phi)$, to which Corollary~\ref{cor: TGS using Premet} applies; once we have treated these we shall turn to those with $p \leq e(\Phi)$, to which Corollary~\ref{cor: TGS without using Premet} applies.

\begin{table}
\caption{Unexcluded $p$-restricted large triples}\label{table: unexcluded triples}
\tabcapsp
$$
\begin{array}{|c|c|c|c|c|c|c|c|c|c|c|c|c|}
\cline{1-4} \cline{6-9} \cline{11-13}
G      & \lambda                 & \ell         & p          & \ptw & G      & \lambda                & \ell        & p          & \ptw & G   & \lambda             & p         \tbs \\
\cline{1-4} \cline{6-9} \cline{11-13}
A_\ell & 3\omega_1               & {} \geq 1    & {} \geq 5  &      & B_\ell & 2\omega_1              & {} \geq 2   & {} \geq 3  &      & G_2 & 2\omega_1           & {} \geq 3 \tbs \\
       & 4\omega_1               & 1, 2         & {} \geq 5  &      &        & \omega_3               & {} \geq 4   & \hbox{any} &      &     & 2\omega_2           & 3         \tbs \\
       & 2\omega_2               & 3, 4, 5      & {} \geq 3  &      &        & \omega_4               & 5           & 2          &      &     & \omega_1 + \omega_2 & 3         \tbs \\
\cline{11-13}
       & \omega_3                & {} \geq 8    & \hbox{any} &      &        & \omega_\ell            & 7, 8, 9     & \hbox{any} & \multicolumn{4}{c}{}                         \tbs \\
       & \omega_4                & 7, \dots, 11 & \hbox{any} &      &        & 2\omega_\ell           & 3, 4        & {} \geq 3  & \multicolumn{4}{c}{}                         \tbs \\
       & \omega_5                & 9            & \hbox{any} &      &        & 3\omega_2              & 2           & {} \geq 5  & \multicolumn{4}{c}{}                         \tbs \\
       & 2\omega_1 + \omega_\ell & 2, 3, 4      & {} \geq 3  &      &        & \omega_1 + \omega_2    & {} \geq 3   & {} \geq 3  & \multicolumn{4}{c}{}                         \tbs \\
       & \omega_1 + \omega_2     & {} \geq 3    & \hbox{any} &      &        & \omega_1 + \omega_\ell & 2, 3        & \hbox{any} & \multicolumn{4}{c}{}                         \tbs \\
       & 2\omega_1 + \omega_2    & 3            & {} \geq 3  &      &        & \omega_1 + \omega_4    & 4           & {} \geq 3  & \multicolumn{4}{c}{}                         \tbs \\
       & \omega_2 + \omega_\ell  & 4, \dots, 8  & \hbox{any} &      &        & \omega_1 + 2\omega_2   & 2           & {} \geq 3  & \multicolumn{4}{c}{}                         \tbs \\
\cline{6-9}
       & \omega_1 + \omega_3     & 5            & \hbox{any} &      & C_\ell & 3\omega_1              & {} \geq 3   & {} \geq 5  & \multicolumn{4}{c}{}                         \tbs \\
       & \omega_2 + \omega_3     & 4            & \hbox{any} &      &        & \omega_3               & {} \geq 4   & \hbox{any} & \multicolumn{4}{c}{}                         \tbs \\
\cline{1-4}
D_\ell & 2\omega_1               & {} \geq 4    & {} \geq 3  &      &        & \omega_4               & 5           & \hbox{any} & \multicolumn{4}{c}{}                         \tbs \\
       & \omega_3                & 5, 6         & \hbox{any} &      &        & \omega_\ell            & 4, 5        & {} \geq 3  & \multicolumn{4}{c}{}                         \tbs \\
       & \omega_\ell             & 8, 9, 10     & \hbox{any} &      &        & \omega_\ell            & 7, 8, 9     & 2          & \multicolumn{4}{c}{}                         \tbs \\
       & 2\omega_5               & 5            & {} \geq 3  &      &        & \omega_1 + \omega_2    & {} \geq 3   & {} \geq 3  & \multicolumn{4}{c}{}                         \tbs \\
       & \omega_1 + \omega_\ell  & 4, 5         & \hbox{any} &      &        & \omega_1 + \omega_3    & 3           & \hbox{any} & \multicolumn{4}{c}{}                         \tbs \\
\cline{1-4} \cline{6-9}
\end{array}
$$
\end{table}

Thus assume $(G, \lambda, p)$ is a $p$-restricted large triple with $p > e(\Phi)$; for such a triple to be unexcluded, all dominant weights $\mu \preceq \lambda$ must be relevant. We shall work through the possibilities for $G$ in turn. For each $G$ we take the corresponding entries in Table~\ref{table: relevant dominant weights}; for each entry, we shall determine which if any natural numbers $a$ (or $a$ and $b$, or $a$, $b$ and $c$ as appropriate) correspond to $p$-restricted large triples $(G, \lambda, p)$ for which $s_\lambda \leq M$, or ${s_\lambda}' \leq M$ if $e(\Phi) > 1$. To show that a triple is excluded, we shall either give a single irrelevant dominant weight $\mu \preceq \lambda$, or list certain dominant weights $\mu \preceq \lambda$ and sum the corresponding values $r_\mu$ to provide a lower bound for $s_\lambda$ (and if $e(\Phi) > 1$ we shall also sum the values ${r_\mu}'$ to provide a lower bound for ${s_\lambda}'$). Note that in each case the requirement that the large triple should be $p$-restricted implies that each coefficient in $\lambda$ should be less than $p$; we state this explicitly in Table~\ref{table: unexcluded triples}, but will not mention it in the proofs in this section.

As usual we write $V = L(\lambda)$. Recall that for $(G, \lambda, p)$ to be a large triple we must have $\dim V > \dim G$; in a few cases this precludes consideration of certain small values of $a$ (and $b$ if appropriate).

\begin{prop}\label{prop: unexcluded triples for A_ell}
Let $G = A_\ell$; then the unexcluded $p$-restricted large triples $(G, \lambda, p)$ are as listed in Table~\ref{table: unexcluded triples}.
\end{prop}

\begin{proof} First suppose $\lambda = a\omega_1$ for $\ell \in [1, \infty)$. If $a \leq 2$ then $(G, \lambda, p)$ is not a large triple. If $a \geq 4$ and $\ell \geq 3$ then taking $\mu = \lambda$, $(a - 2)\omega_1 + \omega_2$, $(a - 4)\omega_1 + 2\omega_2$ and $(a - 3)\omega_1 + \omega_3$ gives $s_\lambda \geq 1 + (2\ell - 1) + (\ell - 1) + \frac{1}{2}(\ell - 1)(3\ell - 4) = \frac{1}{2}(3\ell^2 - \ell + 2) > M$. If $a \geq 5$ and $\ell = 2$ then taking $\mu = \lambda$, $(a - 2)\omega_1 + \omega_2$ and $(a - 4)\omega_1 + 2\omega_2$ gives $s_\lambda \geq 1 + 3 + 3 = 7 > M$; if $a \geq 5$ and $\ell = 1$ then taking $\mu = \lambda$, $(a - 2)\omega_1$ and $(a - 4)\omega_1$ gives $s_\lambda \geq 1 + 1 + 1 = 3 > M$. If however $a = 3$, or $a = 4$ and $\ell \in [1, 2]$, we find that $s_\lambda \leq M$.

Next suppose $\lambda = a\omega_2$ for $\ell \in [3, \infty)$. If $a = 1$ then $(G, \lambda, p)$ is not a large triple. If $a \geq 2$ and $\ell \geq 6$ then taking $\mu = \lambda$, $\omega_1 + (a - 2)\omega_2 + \omega_3$ and $(a - 2)\omega_2 + \omega_4$ gives $s_\lambda \geq (\ell - 1) + \frac{1}{2}(\ell - 1)(3\ell - 4) + \frac{1}{6}(\ell - 1)(\ell - 2)(\ell - 3) = \frac{1}{6}\ell(\ell - 1)(\ell + 4) > M$; if $a \geq 3$ and $\ell \in [3, 5]$ then taking $\mu = \lambda$ and $\omega_1 + (a - 2)\omega_2 + \omega_3$ gives $s_\lambda \geq (\ell - 1) + 3(\ell - 1)^2 = (\ell - 1)(3\ell - 2) > M$. If however $a = 2$ and $\ell \in [3, 5]$ we find that $s_\lambda \leq M$.

Next suppose $\lambda = a\omega_3$ for $\ell \in [5, \infty)$. If $a = 1$ then $(G, \lambda, p)$ is not a large triple for $\ell \in [5, 7]$, whereas for $\ell \in [8, \infty)$ we have $s_\lambda = r_\lambda = \frac{1}{2}(\ell - 1)(\ell - 2) < M$. If $a \geq 2$ then $\mu = \omega_2 + (a - 2)\omega_3 + \omega_4$ is irrelevant.

Next suppose $\lambda = a\omega_4$ for $\ell \in [7, 11]$. If $a \geq 2$ then $\mu = \omega_3 + (a - 2)\omega_4 + \omega_5$ is irrelevant. If however $a = 1$ then we have $s_\lambda = r_\lambda = \frac{1}{6}(\ell - 1)(\ell - 2)(\ell - 3) < M$.

Next suppose $\lambda = a\omega_5$ for $\ell = 9$. If $a \geq 2$ then $\mu = \omega_4 + (a - 2)\omega_5 + \omega_6$ is irrelevant. If however $a = 1$ then we have $s_\lambda = r_\lambda = 70 < M$.

Next suppose $\lambda = a\omega_1 + b\omega_\ell$ for $\ell \in [2, \infty)$; note that we may assume $a \geq b$. If $a = b = 1$ then $(G, \lambda, p)$ is not a large triple. If $a \geq 2$ and $\ell \geq 5$ then taking $\mu = \lambda$ and $(a - 2)\omega_1 + \omega_2 + b\omega_\ell$ gives $s_\lambda \geq (2\ell - 1) + \frac{1}{2}(\ell - 1)(3\ell - 4) = \frac{1}{2}(3\ell^2 - 3\ell + 2) > M$; if $a \geq 3$ and $\ell \in [3, 4]$ then taking $\mu = \lambda$ and $(a - 2)\omega_1 + \omega_2 + b\omega_\ell$ gives $s_\lambda \geq (2\ell - 1) + 3(\ell - 1)^2 = 3\ell^2 - 4\ell + 2 > M$; if $a \geq 3$ and $\ell = 2$ then taking $\mu = \lambda$, $(a - 2)\omega_1 + (b + 1)\omega_2$ and $(a - 1)\omega_1 + (b - 1)\omega_2$ gives $s_\lambda \geq 3 + 3 + 1 = 7 > M$; if $b \geq 2$ and $\ell \in [3, 4]$ then taking $\mu = \lambda$, $(a - 2)\omega_1 + \omega_2 + b\omega_\ell$ and $a\omega_1 + \omega_{\ell - 1} + (b - 2)\omega_\ell$ gives $s_\lambda \geq (2\ell - 1) + \frac{1}{2}(\ell - 1)(3\ell - 4) + \frac{1}{2}(\ell - 1)(3\ell - 4) = 3\ell^2 - 5\ell + 3 > M$; if $b \geq 2$ and $\ell = 2$ then taking $\mu = \lambda$, $(a + 1)\omega_1 + (b - 2)\omega_1$ and $(a - 1)\omega_1 + (b - 1)\omega_2$ gives $s_\lambda \geq 3 + 1 + 3 = 7 > M$. If however $a = 2$, $b = 1$ and $\ell \in [2, 4]$ we find that $s_\lambda \leq M$.

Next suppose $\lambda = a\omega_1 + b\omega_2$ for $\ell \in [3, \infty)$. If $a \geq 2$ and $\ell \geq 4$ then taking $\mu = \lambda$, $(a - 2)\omega_1 + (b + 1)\omega_2$ and $(a - 1)\omega_1 + (b - 1)\omega_2 + \omega_3$ gives $s_\lambda \geq (2\ell - 1) + (\ell - 1) + \frac{1}{2}(\ell - 1)(3\ell - 4) = \frac{1}{2}\ell(3\ell - 1) > M$; if $a \geq 3$ and $\ell = 3$ then taking $\mu = \lambda$, $(a - 2)\omega_1 + (b + 1)\omega_2$ and $(a - 1)\omega_1 + (b - 1)\omega_2 + \omega_3$ gives $s_\lambda \geq 5 + 5 + 5 = 15 > M$; if $b \geq 2$ then taking $\mu = \lambda$, $(a + 1)\omega_1 + (b - 2)\omega_2 + \omega_3$ and $(a - 1)\omega_1 + (b - 1)\omega_2 + \omega_3$ gives $s_\lambda \geq (2\ell - 1) + \frac{1}{2}(\ell - 1)(3\ell - 4) + \frac{1}{2}(\ell - 1)(3\ell - 4) = 3\ell^2 - 5\ell + 3 > M$. If however $a = b = 1$, or $a = 2$, $b = 1$ and $\ell = 3$, we find that $s_\lambda \leq M$.

Next suppose $\lambda = a\omega_2 + b\omega_\ell$ for $\ell \in [4, 8]$. If $a \geq 2$ then $\mu = \omega_1 + (a - 2)\omega_2  + \omega_3 + b\omega_\ell$ is irrelevant; if $b \geq 2$ and $\ell \in [5, 8]$ then $\mu = a\omega_2 + \omega_{\ell - 1} + (b - 2)\omega_\ell$ is irrelevant; if $b \geq 2$ and $\ell = 4$ then taking $\mu = \lambda$ and $a\omega_2 + \omega_3 + (b - 2)\omega_4$ gives $s_\lambda \geq 12 + 12 = 24 > M$. If however $a = b = 1$ we find that $s_\lambda \leq M$.

Next suppose $\lambda = a\omega_1 + b\omega_3$ for $\ell \in [5, 8]$. If $\ell \in [6, 8]$ then taking $\mu = \lambda$ and $(a - 1)\omega_1 + (b - 1)\omega_3 + \omega_4$ gives $s_\lambda \geq \frac{1}{2}(\ell - 1)(3\ell - 4) + \frac{1}{6}(\ell - 1)(\ell - 2)(\ell - 3) = \frac{1}{6}(\ell - 1)(\ell^2 + 4\ell - 6) > M$; if $a \geq 2$ and $\ell = 5$ then taking $\mu = \lambda$ and $(a - 2)\omega_1 + \omega_2 + b\omega_3$ gives $s_\lambda \geq 22 + 22 = 44 > M$; if $b \geq 2$ and $\ell = 5$ then $\mu = a\omega_1 + \omega_2 + (b - 2)\omega_3 + \omega_4$ is irrelevant. If however $a = b = 1$ and $\ell = 5$ we find that $s_\lambda \leq M$.

Next suppose $\lambda = a\omega_2 + b\omega_3$ for $\ell \in [4, 8]$. If $\ell \in [5, 8]$ then taking $\mu = \lambda$ and $\omega_1 + (a - 1)\omega_2 + (b - 1)\omega_3 + \omega_4$ gives $s_\lambda \geq \frac{1}{2}(\ell - 1)(3\ell - 4) + \frac{1}{6}(\ell - 1)(\ell - 2)(4\ell - 9) = \frac{1}{3}(\ell - 1)(2\ell^2 - 4\ell + 3) > M$; if $\ell = 4$ and either $a \geq 2$ or $b \geq 2$ then $\mu = \omega_1 + (a - 1)\omega_2 + (b - 1)\omega_3 + \omega_4$ is irrelevant. If however $a = b = 1$ and $\ell = 4$ we find that $s_\lambda \leq M$.

Finally suppose $\lambda = a\omega_1 + b\omega_2 + c\omega_3$ for $\ell = 3$. Here taking $\mu = \lambda$ and $(a - 1)\omega_1 + b\omega_2 + (c - 1)\omega_3$ gives $s_\lambda \geq 12 + 2 = 14 > M$.
\end{proof}

\begin{prop}\label{prop: unexcluded triples for D_ell}
Let $G = D_\ell$; then the unexcluded $p$-restricted large triples $(G, \lambda, p)$ are as listed in Table~\ref{table: unexcluded triples}.
\end{prop}

\begin{proof}
First suppose $\lambda = a\omega_1$ for $\ell \in [4, \infty)$. If $a = 1$ then $(G, \lambda, p)$ is not a large triple. If $a \geq 3$ then taking $\mu = \lambda$, $(a - 2)\omega_1 + \omega_2$ and either $(a - 3)\omega_1 + \omega_3$ or $(a - 3)\omega_1 + \omega_3 + \omega_4$ according as $\ell \geq 5$ or $\ell = 4$ gives $s_\lambda \geq 2 + 4(2\ell - 3) + 2(\ell - 2)(2\ell - 5) = 4\ell^2 - 10\ell + 10 > M$. If however $a = 2$ we find that $s_\lambda = 4\ell - 5 < M$.

Next suppose $\lambda = a\omega_2$ for $\ell \in [4, \infty)$. If $a = 1$ then $(G, \lambda, p)$ is not a large triple. If $a \geq 2$ then $\mu = \omega_1 + (a - 2)\omega_2 + \omega_3$ or $\omega_1 + (a - 2)\omega_2 + \omega_3 + \omega_4$ according as $\ell \geq 5$ or $\ell = 4$ is irrelevant.

Next suppose $\lambda = a\omega_3$ for $\ell \in [5, 6]$. If $a \geq 2$ then taking $\mu = \lambda$ and $2\omega_2 + (a - 2)\omega_3$ gives $s_\lambda \geq 2(\ell - 2)(2\ell - 5) + (4\ell - 7) = 4\ell^2 - 14\ell + 13 > M$. If however $a = 1$ we find that $s_\lambda = 4\ell^2 - 18\ell + 22 < M$.

Next suppose $\lambda = a\omega_\ell$ for $\ell \in [5, 10]$. If $a = 1$ then $(G, \lambda, p)$ is not a large triple for $\ell \in [5, 7]$, whereas for $\ell \in [8, 10]$ we have $s_\lambda = r_\lambda = 2^{\ell - 3} < M$. If $a \geq 2$ and $\ell \in [6, 10]$ then $\mu = \omega_{\ell - 2} + (a - 2)\omega_\ell$ is irrelevant; if $a \geq 3$ and $\ell = 5$ then $\mu = \omega_3 + (a - 2)\omega_5$ is irrelevant. If however $a = 2$ and $\ell = 5$ we find that $s_\lambda = 36 < M$.

Next suppose $\lambda = a\omega_1 + b\omega_2$ for $\ell \in [4, \infty)$. Here taking $\mu = \lambda$ and either $(a - 1)\omega_1 + (b - 1)\omega_2 + \omega_3$ or $(a - 1)\omega_1 + (b - 1)\omega_2 + \omega_3 + \omega_4$ according as $\ell \geq 5$ or $\ell = 4$ gives $s_\lambda \geq 4(2\ell - 3) + 2(\ell - 2)(2\ell - 5) = 4\ell^2 - 10\ell + 8 > M$.

Next suppose $\lambda = a\omega_1 + b\omega_\ell$ for $\ell \in [4, 5]$. If $a \geq 2$ and $\ell = 5$ then $\mu = (a - 2)\omega_1 + \omega_2 + b\omega_5$ is irrelevant; if $a \geq 2$ and $\ell = 4$ then taking $\mu = \lambda$ and $(a - 2)\omega_1 + \omega_2 + b\omega_4$ gives $s_\lambda \geq 12 + 20 = 32 > M$; if $b \geq 2$ then taking $\mu = \lambda$ and $a\omega_1 + \omega_{\ell - 2} + (b - 2)\omega_\ell$ gives $s_\lambda \geq 2^{\ell - 3}(\ell + 2) + 2^{\ell - 5}(\ell - 2)(\ell^2 + 5\ell - 16) = 2^{\ell - 5}(\ell^3 + 3\ell^2 - 22\ell + 40) > M$. If however $a = b = 1$ we find that $s_\lambda = 2^{\ell - 3}(\ell + 3) < M$.

Finally suppose $\lambda = a\omega_4 + b\omega_5$ for $\ell = 5$. Here taking $\mu = \lambda$ and $\omega_2 + (a - 1)\omega_4 + (b - 1)\omega_5$ gives $s_\lambda \geq 28 + 13 = 41 > M$.
\end{proof}

\begin{prop}\label{prop: unexcluded triples for B_ell}
Let $G = B_\ell$ and $p > 2$; then the unexcluded $p$-restricted large triples $(G, \lambda, p)$ are as listed in Table~\ref{table: unexcluded triples}.
\end{prop}

\begin{proof}
First suppose $\lambda = a\omega_1$ for $\ell \in [2, \infty)$. If $a = 1$ then $(G, \lambda, p)$ is not a large triple. If $a \geq 3$ and $\ell \geq 3$ then taking $\mu = \lambda$, $(a - 2)\omega_1 + \omega_2$, $(a - 1)\omega_1$, $(a - 3)\omega_1 + \omega_2$, $(a - 2)\omega_1$ and either $(a - 3)\omega_1 + \omega_3$ or $(a - 3)\omega_1 + 2\omega_3$ according as $\ell \geq 4$ or $\ell = 3$ gives $s_\lambda \geq 1 + 4(\ell - 1) + 1 + 2(\ell - 1) + 1 + 2(\ell - 1)(\ell - 2) = 2\ell^2 + 1 > M$ and ${s_\lambda}' \geq 2 + 4(2\ell - 3) + 2 + (4\ell - 7) + 2 + 2(\ell - 2)(2\ell - 5) = 4\ell^2 - 6\ell + 7 > M$; if $a \geq 3$ and $\ell = 2$ then taking $\mu = \lambda$, $(a - 2)\omega_1 + \omega_2$, $(a - 1)\omega_1$, $(a - 3)\omega_1 + \omega_2$ and $(a - 2)\omega_1$ gives $s_\lambda \geq 1 + 4 + 1 + 2 + 1 = 9 > M$ and ${s_\lambda}' \geq 2 + 4 + 2 + 1 + 2 = 11 > M$. If however $a = 2$ we find that $s_\lambda < M$ and ${s_\lambda}' < M$.

Next suppose $\lambda = a\omega_2$ for $\ell \in [3, \infty)$. If $a = 1$ then $(G, \lambda, p)$ is not a large triple. If $a \geq 2$ and $\ell \geq 4$ then $\mu = \omega_1 + (a - 2)\omega_2 + \omega_3$ is irrelevant; if $a \geq 2$ and $\ell = 3$ then taking $\mu = \lambda$, $\omega_1 + (a - 2)\omega_2 + 2\omega_3$ and $\omega_1 + (a - 1)\omega_2$ gives $s_\lambda \geq 4 + 12 + 8 = 24 > M$ and ${s_\lambda}' \geq 5 + 10 + 12 = 27 > M$.

Next suppose $\lambda = a\omega_3$ for $\ell \in [4, \infty)$. If $a \geq 2$ and $\ell \geq 5$ then $\mu = \omega_2 + (a - 2)\omega_3 + \omega_4$ is irrelevant; if $a \geq 2$ and $\ell = 4$ then $\mu = \omega_2 + (a - 2)\omega_3 + 2\omega_4$ is irrelevant. If however $a = 1$ we find that $s_\lambda < M$.

Next suppose $\lambda = a\omega_{\ell - 1}$ for $\ell \in [5, 6]$. If $\ell = 6$ then $\mu = \omega_4 + (a - 1)\omega_5$ is irrelevant; if $\ell = 5$ then taking $\mu = \lambda$ and $\omega_3 + (a - 1)\omega_4$ gives $s_\lambda \geq 32 + 24 = 56 > M$ and ${s_\lambda}' \geq 28 + 30 = 58 > M$.

Next suppose $\lambda = a\omega_\ell$ for $\ell \in [2, 9]$. If $a = 1$ then $(G, \lambda, p)$ is not a large triple for $\ell \in [2, 6]$, whereas for $\ell \in [7, 9]$ we have ${s_\lambda}' = {r_\lambda}' = 2^{\ell - 2} < M$. If $a = 2$ and $\ell = 2$ then $(G, \lambda, p)$ is not a large triple. If $a \geq 2$ and $\ell \in [5, 9]$ then taking $\mu = \lambda$, $\omega_{\ell - 1} + (a - 2)\omega_\ell$ and $\omega_{\ell - 2} + (a - 2)\omega_\ell$ gives $s_\lambda \geq 2^{\ell - 1} + 2^{\ell - 2}(\ell - 1) + 2^{\ell - 4}(\ell - 1)(\ell - 2) = 2^{\ell - 4}(\ell^2 + \ell + 6) > M$ and ${s_\lambda}' \geq 2^{\ell - 2} + 2^{\ell - 3}(\ell + 2) + 2^{\ell - 5}(\ell - 2)(\ell + 5) = 2^{\ell - 5}(\ell + 1)(\ell + 6) > M$; if $a \geq 3$ and $\ell = 4$ then $\mu = \omega_2 + (a - 2)\omega_4$ is irrelevant; if $a \geq 3$ and $\ell = 3$ then taking $\mu = \lambda$, $\omega_2 + (a - 2)\omega_3$ and $\omega_1 + (a - 2)\omega_3$ gives $s_\lambda \geq 4 + 12 + 12 = 28 > M$ and ${s_\lambda}' \geq 2 + 10 + 10 = 22 > M$; if $a \geq 4$ and $\ell = 2$ then taking $\mu = \lambda$, $\omega_1 + (a - 2)\omega_2$, $2\omega_1 + (a - 4)\omega_2$, $(a - 2)\omega_2$ and $\omega_1 + (a - 4)\omega_2$ gives $s_\lambda \geq 2 + 4 + 1 + 2 + 1 = 10 > M$ and ${s_\lambda}' \geq 1 + 4 + 2 + 1 + 2 = 10 > M$. If however $a = 2$ and $\ell \in [3, 4]$, or $a = 3$ and $\ell = 2$, we find that $s_\lambda \leq M$.

Next suppose $\lambda = a\omega_1 + b\omega_2$ for $\ell \in [3, \infty)$. If $a \geq 2$ and $\ell \geq 4$ then $\mu = (a - 1)\omega_1 + (b - 1)\omega_2 + \omega_3$ is irrelevant; if $a \geq 2$ and $\ell = 3$ then taking $\mu = \lambda$, $(a - 2)\omega_1 + (b + 1)\omega_2$ and $(a - 1)\omega_1 + b\omega_2$ gives $s_\lambda \geq 8 + 4 + 8 = 20 > M$ and ${s_\lambda}' \geq 12 + 5 + 12 = 29 > M$; if $b \geq 2$ and $\ell \geq 4$ then $\mu = (a - 1)\omega_1 + (b - 1)\omega_2 + \omega_3$ is irrelevant; if $b \geq 2$ and $\ell = 3$ then taking $\mu = \lambda$ and $(a + 1)\omega_2 + (b - 2)\omega_2 + 2\omega_3$ gives $s_\lambda \geq 8 + 12 = 20 > M$ and ${s_\lambda}' \geq 12 + 10 = 22 > M$. If however $a = b = 1$ we find that $s_\lambda = M$.

Next suppose $\lambda = a\omega_1 + b\omega_\ell$ for $\ell \in [2, 4]$. If $a \geq 2$ and $\ell = 4$ then $\mu = (a - 2)\omega_1 + \omega_2 + b\omega_4$ is irrelevant; if $a \geq 2$ and $\ell = 3$ then taking $\mu = \lambda$ and $(a - 2)\omega_1 + \omega_2 + b\omega_3$ gives $s_\lambda \geq 12 + 12 = 24 > M$ and ${s_\lambda}' \geq 10 + 10 = 20 > M$; if $a \geq 2$ and $\ell = 2$ then taking $\mu = \lambda$, $(a - 2)\omega_1 + (b + 2)\omega_2$ and $(a - 1)\omega_1 + b\omega_2$ gives $s_\lambda \geq 4 + 2 + 4 = 10 > M$ and ${s_\lambda}' \geq 4 + 1 + 4 = 9 > M$; if $b \geq 2$ and $\ell = 4$ then $\mu = a\omega_1 + \omega_3 + (b - 2)\omega_4$ is irrelevant; if $b \geq 2$ and $\ell = 3$ then taking $\mu = \lambda$ and $a\omega_1 + \omega_2 + (b - 2)\omega_3$ gives $s_\lambda \geq 12 + 8 = 20 > M$ and ${s_\lambda}' \geq 10 + 12 = 22 > M$; if $b \geq 3$ and $\ell = 2$ then taking $\mu = \lambda$, $(a + 1)\omega_1 + (b - 2)\omega_2$ and $(a - 1)\omega_1 + b\omega_2$ gives $s_\lambda \geq 4 + 4 + 2 = 10 > M$ and ${s_\lambda}' \geq 4 + 4 + 1 = 9 > M$. If however $a = b = 1$ we find that ${s_\lambda}' < M$, while if $a = 1$, $b = 2$ and $\ell = 2$ we find that $s_\lambda = M$.

Finally suppose $\lambda = a\omega_{\ell - 1} + b\omega_\ell$ for $\ell \in [3, 4]$. If $\ell = 4$ then $\mu = \omega_2 + (a - 1)\omega_3 + b\omega_4$ is irrelevant; if $\ell = 3$ then taking $\mu = \lambda$ and $\omega_1 + (a - 1)\omega_2 + b\omega_3$ gives $s_\lambda \geq 12 + 12 = 24 > M$ and ${s_\lambda}' \geq 10 + 10 = 20 > M$.
\end{proof}

\begin{prop}\label{prop: unexcluded triples for C_ell}
Let $G = C_\ell$ and $p > 2$; then the unexcluded $p$-restricted large triples $(G, \lambda, p)$ are as listed in Table~\ref{table: unexcluded triples}.
\end{prop}

\begin{proof}
First suppose $\lambda = a\omega_1$ for $\ell \in [3, \infty)$. If $a \leq 2$ then $(G, \lambda, p)$ is not a large triple. If $a \geq 4$ and $\ell \geq 4$ then $\mu = (a - 3)\omega_1 + \omega_3$ is irrelevant; if $a \geq 4$ and $\ell = 3$ then taking $\mu = \lambda$, $(a - 2)\omega_1 + \omega_2$ and $(a - 3)\omega_1 + \omega_3$ gives $s_\lambda \geq 2 + 12 + 10 = 24 > M$ and ${s_\lambda}' \geq 1 + 8 + 12 = 21 > M$. If however $a = 3$ we find that ${s_\lambda}' < M$.

Next suppose $\lambda = a\omega_2$ for $\ell \in [3, \infty)$. If $a = 1$ then $(G, \lambda, p)$ is not a large triple. If $a \geq 2$ and $\ell \geq 4$ then $\mu = \omega_1 + (a - 2)\omega_2 + \omega_3$ is irrelevant; if $a \geq 2$ and $\ell = 3$ then taking $\mu = \lambda$, $\omega_1 + (a - 2)\omega_2 + \omega_3$ and $(a - 1)\omega_2$ gives $s_\lambda \geq 5 + 10 + 5 = 20 > M$ and ${s_\lambda}' \geq 4 + 12 + 4 = 20 > M$.

Next suppose $\lambda = a\omega_3$ for $\ell \in [4, \infty)$. If $a \geq 2$ then $\mu = \omega_2 + (a - 2)\omega_3 + \omega_4$ is irrelevant. If however $a = 1$ we find that ${s_\lambda}' < M$.

Next suppose $\lambda = a\omega_{\ell - 1}$ for $\ell \in [5, 6]$. If $a \geq 2$ then $\mu = \omega_{\ell - 3} + (a - 1)\omega_{\ell - 1}$ is irrelevant; if $a = 1$ and $\ell = 6$ then taking $\mu = \lambda$ and $\omega_3$ gives $s_\lambda \geq 64 + 56 = 120 > M$ and ${s_\lambda}' \geq 80 + 40 = 120 > M$. If however $a = 1$ and $\ell = 5$ we find that ${s_\lambda}' < M$.

Next suppose $\lambda = a\omega_\ell$ for $\ell \in [3, 9]$. If $a = 1$ and $\ell = 3$ then $(G, \lambda, p)$ is not a large triple. If $a \geq 2$ and $\ell \in [4, 9]$ then $\mu = \omega_{\ell - 2} + (a - 1)\omega_\ell$ is irrelevant; if $a \geq 2$ and $\ell = 3$ then taking $\mu = \lambda$, $2\omega_2 + (a - 2)\omega_3$, $\omega_1 + (a - 1)\omega_3$ and $\omega_2 + (a - 2)\omega_3$ gives $s_\lambda \geq 2 + 5 + 10 + 5 = 22 > M$ and ${s_\lambda}' \geq 4 + 4 + 12 + 4 = 24 > M$; if $a = 1$ and $\ell \in [6, 9]$ then $\mu = \omega_{\ell - 2}$ is irrelevant. If however $a = 1$ and $\ell \in [4, 5]$ we find that ${s_\lambda}' < M$.

Next suppose $\lambda = a\omega_1 + b\omega_2$ for $\ell \in [3, \infty)$. If $a \geq 2$ and $\ell \geq 4$ then $\mu = (a - 1)\omega_1 + (b - 1)\omega_2 + \omega_3$ is irrelevant; if $a \geq 2$ and $\ell = 3$ then taking $\mu = \lambda$ and $(a - 1)\omega_1 + (b - 1)\omega_2 + \omega_3$ gives $s_\lambda \geq 12 + 10 = 22 > M$ and ${s_\lambda}' \geq 8 + 12 = 20 > M$; if $b \geq 2$ and $\ell \geq 4$ then $\mu = (a + 1)\omega_1 + (b - 2)\omega_2 + \omega_3$ is irrelevant; if $b \geq 2$ and $\ell = 3$ then taking $\mu = \lambda$ and $(a + 1)\omega_1 + (b - 2)\omega_2 + \omega_3$ gives $s_\lambda \geq 12 + 10 = 22 > M$ and ${s_\lambda}' \geq 8 + 12 = 20 > M$. If however $a = b = 1$ we find that ${s_\lambda}' < M$.

Next suppose $\lambda = a\omega_1 + b\omega_\ell$ for $\ell \in [3, 4]$. If $\ell = 4$ then taking $\mu = \lambda$ and $a\omega_1 + \omega_2 + (b - 1)\omega_4$ gives $s_\lambda \geq 24 + 20 = 44 > M$ and ${s_\lambda}' \geq 32 + 12 = 44 > M$; if $\ell = 3$ and $a \geq 2$ then taking $\mu = \lambda$ and $(a - 2)\omega_1 + \omega_2 + b\omega_3$ gives $s_\lambda \geq 10 + 10 = 20 > M$ and ${s_\lambda}' \geq 12 + 12 = 24 > M$; if $\ell = 3$ and $b \geq 2$ then taking $\mu = \lambda$ and $a\omega_1 + 2\omega_2 + (b - 2)\omega_3$ gives $s_\lambda \geq 10 + 12 = 22 > M$ and ${s_\lambda}' \geq 12 + 8 = 20 > M$. If however $a = b = 1$ and $\ell = 3$ we find that ${s_\lambda}' < M$.

Finally suppose $\lambda = a\omega_{\ell - 1} + b\omega_\ell$ for $\ell \in [3, 4]$. If $\ell = 4$ then $\mu = \omega_2 + a\omega_3 + (b - 1)\omega_4$ is irrelevant; if $\ell = 3$ then taking $\mu = \lambda$ and $\omega_1 + a\omega_2 + (b - 1)\omega_3$ gives $s_\lambda \geq 10 + 12 = 22 > M$ and ${s_\lambda}' \geq 12 + 8 = 20 > M$.
\end{proof}

\begin{prop}\label{prop: unexcluded triples for exceptional groups}
Let $G$ be of exceptional type and $p > e(\Phi)$; then the unexcluded $p$-restricted large triples $(G, \lambda, p)$ are as listed in Table~\ref{table: unexcluded triples}.
\end{prop}

\begin{proof}
Take $G = E_6$. First suppose $\lambda = a\omega_1$. If $a = 1$ then $(G, \lambda, p)$ is not a large triple; if $a \geq 2$ then $\mu = (a - 2)\omega_1 + \omega_3$ is irrelevant. Next suppose $\lambda = a\omega_2$. If $a = 1$ then $(G, \lambda, p)$ is not a large triple; if $a \geq 2$ then $\mu = (a - 2)\omega_2 + \omega_4$ is irrelevant.

Next take $G = E_7$. First suppose $\lambda = a\omega_1$. If $a = 1$ then $(G, \lambda, p)$ is not a large triple; if $a \geq 2$ then $\mu = (a - 2)\omega_1 + \omega_3$ is irrelevant. Next suppose $\lambda = a\omega_7$. If $a = 1$ then $(G, \lambda, p)$ is not a large triple; if $a \geq 2$ then $\mu = \omega_6 + (a - 2)\omega_7$ is irrelevant.

Next take $G = E_8$. Suppose $\lambda = a\omega_8$. If $a = 1$ then $(G, \lambda, p)$ is not a large triple; if $a \geq 2$ then $\mu = \omega_7 + (a - 2)\omega_8$ is irrelevant.

Next take $G = F_4$. First suppose $\lambda = a\omega_1$. If $a = 1$ then $(G, \lambda, p)$ is not a large triple; if $a \geq 2$ then $\mu = (a - 1)\omega_1 + \omega_4$ is irrelevant. Next suppose $\lambda = a\omega_2$. Here $\mu = \omega_1 + (a - 1)\omega_2 + \omega_4$ is irrelevant. Next suppose $\lambda = a\omega_3$. Here taking $\mu = \lambda$, $\omega_1 + (a - 1)\omega_3$ and $(a - 1)\omega_3 + \omega_4$ gives $s_\lambda \geq 44 + 6 + 9 = 59 > M$ and ${s_\lambda}' \geq 36 + 9 + 6 = 51 > M$. Finally suppose $\lambda = a\omega_4$. If $a = 1$ then $(G, \lambda, p)$ is not a large triple; if $a \geq 2$ then taking $\mu = \lambda$, $\omega_3 + (a - 2)\omega_4$ and $\omega_1 + (a - 2)\omega_4$ gives $s_\lambda \geq 9 + 44 + 6 = 59 > M$ and ${s_\lambda}' \geq 6 + 36 + 9 = 51 > M$.

Finally take $G = G_2$. First suppose $\lambda = a\omega_1$. If $a = 1$ then $(G, \lambda, p)$ is not a large triple; if $a \geq 3$ then taking $\mu = \lambda$, $(a - 2)\omega_1 + \omega_2$, $(a - 1)\omega_1$ and $(a - 3)\omega_1 + \omega_2$ gives $s_\lambda \geq 3 + 6 + 3 + 2 = 14 > M$ and ${s_\lambda}' \geq 2 + 6 + 2 + 3 = 13 > M$; if however $a = 2$ we find $s_\lambda = 8 < M$ and ${s_\lambda}' = 7 < M$. Next suppose $\lambda = a\omega_2$. If $a = 1$ then $(G, \lambda, p)$ is not a large triple; if $a \geq 2$ then taking $\mu = \lambda$, $3\omega_1 + (a - 2)\omega_2$, $\omega_1 + (a - 1)\omega_2$ and $2\omega_1 + (a - 2)\omega_2$ gives $s_\lambda \geq 2 + 3 + 6 + 3 = 14 > M$ and ${s_\lambda}' \geq 3 + 2 + 6 + 2 = 13 > M$. Finally suppose $\lambda = a\omega_1 + b\omega_2$. Here taking $\mu = \lambda$, $(a + 1)\omega_1 + (b - 1)\omega_2$, $(a - 1)\omega_1 + b\omega_2$ and $a\omega_1 + (b - 1)\omega_2$ gives $s_\lambda \geq 6 + 3 + 2 + 3 = 14 > M$ and ${s_\lambda}' \geq 6 + 2 + 3 + 2 = 13 > M$.
\end{proof}

We now assume $(G, \lambda, p)$ is a $p$-restricted large triple with $p \leq e(\Phi)$; for such a triple to be unexcluded, the weight $\lambda$ must be $p$-relevant, but the same need not be true of all dominant weights $\mu \prec \lambda$, since they need not appear in $V$ (indeed not all such dominant weights need be $p$-restricted). Here we shall make frequent use of L\"ubeck's online data \cite{Lubdata}, which for a given $G$ lists all irreducible modules of dimension less than a certain bound, and gives the weight multiplicities in each. In addition, we use Lemma~\ref{lem: multiplicities in 3omega_1 and omega_1 + omega_2} to treat the triples $(B_\ell, \omega_1 + \omega_2, 2)$ and $(C_\ell, \omega_1 + \omega_2, 2)$ for $\ell \in [4, \infty)$, and \cite{GS} to treat the triple $(G_2, 2\omega_1 + 2\omega_2, 3)$.

We shall again work through the possibilities for $G$ in turn, taking the entries in Table~\ref{table: p-relevant dominant weights}. We ignore those weights $\lambda$ for which $(G, \lambda, p)$ is not a large triple; to show that a large triple $(G, \lambda, p)$ is excluded we shall list certain dominant weights $\mu \preceq \lambda$, together with their multiplicities $m_\mu$ in $V$, and sum both the values $m_\mu r_{\mu,p}$ and the values $m_\mu {r_\mu}'$ to provide lower bounds for both $s_{\lambda,p}$ and ${s_{\lambda,p}}'$.

\begin{prop}\label{prop: unexcluded triples for B_ell when p = 2}
Let $G = B_\ell$ and $p = 2$; then the unexcluded $p$-restricted large triples $(G, \lambda, p)$ are as listed in Table~\ref{table: unexcluded triples}.
\end{prop}

\begin{proof}
If $\lambda = \omega_1$ or $\omega_2$, or $\omega_\ell$ for $\ell \in [2, 6]$, then $(G, \lambda, p)$ is not a large triple. If $\lambda = \omega_5$ for $\ell = 6$ then taking $\mu = \lambda$ and $\omega_3$ we have $m_\mu = 1$ and $2$ respectively, giving $s_{\lambda, 2} \geq 1.80 + 2.40 = 160 > M$ and ${s_{\lambda, 2}}' \geq 1.64 + 2.56 = 176 > M$. If $\lambda = \omega_1 + \omega_2$ for $\ell \in [4, \infty)$ then taking $\mu = \lambda$ and $\omega_3$ we have $m_\mu = 1$ and $2$ respectively, giving $s_{\lambda, 2} \geq 1.4(\ell - 1) + 2.2(\ell - 1)(\ell - 2) = 4(\ell - 1)^2 > M$ and ${s_{\lambda, 2}}' \geq 1.4(2\ell - 3) + 2.2(\ell - 2)(2\ell - 5) = 4(2\ell^2 - 7\ell + 7) > M$. If $\lambda = \omega_1 + \omega_2$ for $\ell = 3$ then taking $\mu = \lambda$, $2\omega_3$ and $\omega_1$ we have $m_\mu = 1$, $2$ and $4$ respectively, giving $s_{\lambda, 2} \geq 1.8 + 2.4 + 4.1 = 20 > M$ and ${s_{\lambda, 2}}' \geq 1.12 + 2.2 + 4.2 = 24 > M$. If $\lambda = \omega_1 + \omega_4$ for $\ell = 4$ then taking $\mu = \lambda$ and $\omega_4$ we have $m_\mu = 1$ and $4$ respectively, giving $s_{\lambda, 2} \geq 1.24 + 4.8 = 56 > M$ and ${s_{\lambda, 2}}' \geq 1.24 + 4.4 = 40 > M$. If $\lambda = \omega_2 + \omega_4$ for $\ell = 4$ then taking $\mu = \lambda$ and $\omega_1 + \omega_4$ we have $m_\mu = 1$ and $3$ respectively, giving $s_{\lambda, 2} \geq 1.24 + 3.24 = 96 > M$ and ${s_{\lambda, 2}}' \geq 1.40 + 3.24 = 112 > M$. If $\lambda = \omega_{\ell - 1} + \omega_\ell$ for $\ell \in [3, 8]$ then taking $\mu = \lambda$ and $\omega_{\ell - 2} + \omega_\ell$ we have $m_\mu = 1$ and $2$ respectively, giving $s_{\lambda, 2} \geq 1.2^{\ell - 1} + 2.2^{\ell - 1}(\ell - 1) = 2^{\ell - 1}(2\ell - 1) > M$ and ${s_{\lambda, 2}}' \geq 1.2^{\ell - 2}(\ell + 2) + 2.2^{\ell - 3}(\ell^2 + 3\ell - 8) = 2^{\ell - 2}(\ell^2 + 4\ell - 6) > M$. If $\lambda = \omega_1 + \omega_2 + \omega_3$ for $\ell = 3$ then taking $\mu = \lambda$ and $\omega_2 + \omega_3$ we have $m_\mu = 1$ and $4$ respectively, giving $s_{\lambda, 2} \geq 1.8 + 4.4 = 24 > M$ and ${s_{\lambda, 2}}' \geq 1.24 + 4.10 = 64 > M$. If $\lambda = \omega_1 + \omega_3 + \omega_4$ for $\ell = 4$ then taking $\mu = \lambda$ and $\omega_1 + \omega_2 + \omega_4$ we have $m_\mu = 1$ and $2$ respectively, giving $s_{\lambda, 2} \geq 1.24 + 2.48 = 120 > M$ and ${s_{\lambda, 2}}' \geq 1.88 + 2.88 = 264 > M$. If $\lambda = \omega_2 + \omega_3 + \omega_4$ for $\ell = 4$ then taking $\mu = \lambda$ and $\omega_1 + \omega_2 + \omega_4$ we have $m_\mu = 1$ and $8$ respectively, giving $s_{\lambda, 2} \geq 1.24 + 8.48 = 408 > M$ and ${s_{\lambda, 2}}' \geq 1.88 + 8.88 = 792 > M$. If however $\lambda = \omega_3$ for $\ell \in [4, \infty)$, or $\lambda = \omega_4$ for $\ell = 5$, or $\lambda = \omega_\ell$ for $\ell \in [7, 9]$, or $\lambda = \omega_1 + \omega_\ell$ for $\ell \in [2, 3]$, we find that $s_{\lambda, 2} \leq M$ or ${s_{\lambda, 2}}' \leq M$.
\end{proof}

\begin{prop}\label{prop: unexcluded triples for C_ell when p = 2}
Let $G = C_\ell$ and $p = 2$; then the unexcluded $p$-restricted large triples $(G, \lambda, p)$ are as listed in Table~\ref{table: unexcluded triples}.
\end{prop}

\begin{proof}
If $\lambda = \omega_1$ or $\omega_2$, or $\omega_\ell$ for $\ell \in [3, 6]$, then $(G, \lambda, p)$ is not a large triple. If $\lambda = \omega_5$ for $\ell = 6$ then taking $\mu = \lambda$ and $\omega_3$ we have $m_\mu = 1$ and $2$ respectively, giving $s_{\lambda, 2} \geq 1.64 + 2.56 = 176 > M$ and ${s_{\lambda, 2}}' \geq 1.80 + 2.40 = 160 > M$. If $\lambda = \omega_1 + \omega_2$ for $\ell \in [4, \infty)$ then taking $\mu = \lambda$ and $\omega_3$ we have $m_\mu = 1$ and $2$ respectively, giving $s_{\lambda, 2} \geq 1.2(4\ell - 7) + 2.2(\ell - 2)(2\ell - 5) = 2(4\ell^2 - 14\ell + 13) > M$ and ${s_{\lambda, 2}}' \geq 1.4(\ell - 1) + 2.2(\ell - 1)(\ell - 2) = 4(\ell - 1)^2 > M$. If $\lambda = \omega_1 + \omega_2$ for $\ell = 3$ then taking $\mu = \lambda$, $\omega_3$ and $\omega_1$ we have $m_\mu = 1$, $2$ and $4$ respectively, giving $s_{\lambda, 2} \geq 1.10 + 2.2 + 4.2 = 22 > M$ and ${s_{\lambda, 2}}' \geq 1.8 + 2.4 + 4.1 = 20 > M$. If $\lambda = \omega_1 + \omega_5$ for $\ell = 5$ then taking $\mu = \lambda$ and $\omega_4$ we have $m_\mu = 1$ and $2$ respectively, giving $s_{\lambda, 2} \geq 1.40 + 2.28 = 96 > M$ and ${s_{\lambda, 2}}' \geq 1.80 + 2.32 = 144 > M$. If $\lambda = \omega_1 + \omega_4$ for $\ell = 4$ then taking $\mu = \lambda$ and $\omega_3$ we have $m_\mu = 1$ and $2$ respectively, giving $s_{\lambda, 2} \geq 1.16 + 2.12 = 40 > M$ and ${s_{\lambda, 2}}' \geq 1.32 + 2.12 = 56 > M$. If $\lambda = \omega_2 + \omega_4$ for $\ell = 4$ then taking $\mu = \lambda$ and $\omega_1 + \omega_3$ we have $m_\mu = 1$ and $2$ respectively, giving $s_{\lambda, 2} \geq 1.24 + 2.36 = 96 > M$ and ${s_{\lambda, 2}}' \geq 1.48 + 2.36 = 120 > M$. If $\lambda = \omega_4 + \omega_5$ for $\ell = 5$ then taking $\mu = \lambda$ and $\omega_1 + \omega_2$ we have $m_\mu = 1$ and $8$ respectively, giving $s_{\lambda, 2} \geq 1.40 + 8.26 = 248 > M$ and ${s_{\lambda, 2}}' \geq 1.80 + 8.16 = 208 > M$. If $\lambda = \omega_3 + \omega_4$ for $\ell = 4$ then taking $\mu = \lambda$ and $\omega_2 + \omega_3$ we have $m_\mu = 1$ and $2$ respectively, giving $s_{\lambda, 2} \geq 1.16 + 2.36 = 88 > M$ and ${s_{\lambda, 2}}' \geq 1.32 + 2.36 = 104 > M$. If $\lambda = \omega_2 + \omega_3$ for $\ell = 3$ then taking $\mu = \lambda$ and $\omega_1 + \omega_2$ we have $m_\mu = 1$ and $2$ respectively, giving $s_{\lambda, 2} \geq 1.6 + 2.10 = 26 > M$ and ${s_{\lambda, 2}}' \geq 1.12 + 2.8 = 28 > M$. If $\lambda = \omega_1 + \omega_2 + \omega_3$ for $\ell = 3$ then taking $\mu = \lambda$ and $\omega_2$ we have $m_\mu = 1$ and $12$ respectively, giving $s_{\lambda, 2} \geq 1.16 + 12.5 = 76 > M$ and ${s_{\lambda, 2}}' \geq 1.24 + 12.4 = 72 > M$. If however $\lambda = \omega_3$ for $\ell \in [4, \infty)$, or $\lambda = \omega_4$ for $\ell = 5$, or $\lambda = \omega_\ell$ for $\ell \in [7, 9]$, or $\lambda = \omega_1 + \omega_3$ for $\ell = 3$, we find that $s_{\lambda, 2} \leq M$ or ${s_{\lambda, 2}}' \leq M$.
\end{proof}

\begin{prop}\label{prop: unexcluded triples for F_4 when p = 2}
Let $G = F_4$ and $p = 2$; then there are no unexcluded $p$-restricted large triples $(G, \lambda, p)$.
\end{prop}

\begin{proof}
If $\lambda = \omega_1$ or $\omega_4$ then $(G, \lambda, p)$ is not a large triple. If $\lambda = \omega_2$ then taking $\mu = \lambda$ and $\omega_1$ we have $m_\mu = 1$ and $4$ respectively, giving $s_{\lambda, 2} \geq 1.36 + 4.6 = 60 > M$ and ${s_{\lambda, 2}}' \geq 1.44 + 4.9 = 80 > M$. If $\lambda = \omega_3$ then taking $\mu = \lambda$ and $\omega_4$ we have $m_\mu = 1$ and $4$ respectively, giving $s_{\lambda, 2} \geq 1.36 + 4.9 = 72 > M$ and ${s_{\lambda, 2}}' \geq 1.36 + 4.6 = 60 > M$. If $\lambda = \omega_1 + \omega_2$ then taking $\mu = \lambda$ and $\omega_2$ we have $m_\mu = 1$ and $14$ respectively, giving $s_{\lambda, 2} \geq 1.48 + 14.36 = 552 > M$ and ${s_{\lambda, 2}}' \geq 1.96 + 14.44 = 712 > M$. If $\lambda = \omega_1 + \omega_4$ then taking $\mu = \lambda$ and $\omega_3$ we have $m_\mu = 1$ and $3$ respectively, giving $s_{\lambda, 2} \geq 1.36 + 3.36 = 144 > M$ and ${s_{\lambda, 2}}' \geq 1.60 + 3.36 = 168 > M$. If $\lambda = \omega_2 + \omega_3$ then taking $\mu = \lambda$ and $\omega_1 + \omega_3 + \omega_4$ we have $m_\mu = 1$ and $2$ respectively, giving $s_{\lambda, 2} \geq 1.48 + 2.144 = 336 > M$ and ${s_{\lambda, 2}}' \geq 1.132 + 2.264 = 660 > M$. If $\lambda = \omega_2 + \omega_4$ then taking $\mu = \lambda$ and $\omega_1 + \omega_3$ we have $m_\mu = 1$ and $2$ respectively, giving $s_{\lambda, 2} \geq 1.48 + 2.84 = 216 > M$ and ${s_{\lambda, 2}}' \geq 1.132 + 2.132 = 396 > M$.
\end{proof}

\begin{prop}\label{prop: unexcluded triples for G_2 when p = 2 or 3}
Let $G = G_2$ and $p = 2$ or $3$; then the unexcluded $p$-restricted large triples $(G, \lambda, p)$ are as listed in Table~\ref{table: unexcluded triples}.
\end{prop}

\begin{proof}
If $\lambda = \omega_1$ or $\omega_2$ then $(G, \lambda, p)$ is not a large triple. Suppose $p = 2$. If $\lambda = \omega_1 + \omega_2$ then taking $\mu = \lambda$ and $\omega_1$ we have $m_\mu = 1$ and $4$ respectively, giving $s_{\lambda, 2} \geq 1.4 + 4.3 = 16 > M$ and ${s_{\lambda, 2}}' \geq 1.6 + 4.2 = 14 > M$. Suppose $p = 3$. If $\lambda = 2\omega_1 + \omega_2$ then taking $\mu = \lambda$ and $\omega_1$ we have $m_\mu = 1$ and $8$ respectively, giving $s_{\lambda, 3} \geq 1.0 + 8.2 = 16 > M$ and ${s_{\lambda, 3}}' \geq 1.6 + 8.2 = 22 > M$. If $\lambda = \omega_1 + 2\omega_2$ then taking $\mu = \lambda$ and $\omega_1$ we have $m_\mu = 1$ and $7$ respectively, giving $s_{\lambda, 3} \geq 1.2 + 7.2 = 16 > M$ and ${s_{\lambda, 3}}' \geq 1.6 + 7.2 = 20 > M$. If $\lambda = 2\omega_1 + 2\omega_2$ then taking $\mu = \lambda$ and $\omega_1$ we have $m_\mu = 1$ and $19$ respectively, giving $s_{\lambda, 3} \geq 1.0 + 19.2 = 38 > M$ and ${s_{\lambda, 3}}' \geq 1.6 + 19.6 = 120 > M$. If however $\lambda = \omega_1 + \omega_2$, or $\lambda = 2\omega_1$, or $\lambda = 2\omega_2$, we find that $s_{\lambda, 3} \leq M$.
\end{proof}

This concludes the application of Corollaries~\ref{cor: TGS using Premet} and \ref{cor: TGS without using Premet} to the task of proving that all $p$-restricted large triples which are not listed in Table~\ref{table: large triple and first quadruple non-TGS} have TGS. In the next two sections we shall complete this task by dealing with the unexcluded $p$-restricted large triples which are listed in Table~\ref{table: unexcluded triples} but not in Table~\ref{table: large triple and first quadruple non-TGS}.

\section{Weight string analysis}\label{sect: large triple weight string analysis}

In this section we shall treat some of the unexcluded $p$-restricted large triples listed in Table~\ref{table: unexcluded triples}. Our approach is to consider weight strings in more detail than we have done up to this point. In some cases we shall see that it is still possible to show that the triple satisfies both $\ssddagcon$ and $\uddagcon$; in others we shall instead use one or two of the weaker conditions given in Section~\ref{sect: conditions}, but in each case we shall show that it at least satisfies both $\ssdiamevcon$ and $\udiamcon$, and thus has TGS.

We shall use the following notation throughout. Given a triple $(G, \lambda, p)$, we write $V = L(\lambda)$. We let $s$ be an element of $G_{(r)}$ for some $r \in \P'$, and take $\kappa \in K^*$; we assume $s$ lies in $T$. We write $\Phi(s) = \{ \alpha \in \Phi : \alpha(s) = 1 \}$, so that $C_G(s)^\circ = \langle T, X_\alpha : \alpha \in \Phi(s) \rangle$. We take $\alpha \in \Phi_s$ and write $u_\alpha = x_\alpha(1)$; if $e(\Phi) > 1$, we take $\beta \in \Phi_l$ and write $u_\beta = x_\beta(1)$.

We start with triples $(G, \lambda, p)$ where $G = A_\ell$ or $D_\ell$, in which $e(\Phi) = 1$. For each triple we shall give two tables. The first is the weight table, which lists the dominant weights $\mu \in \Lambda(V)$, and gives the sizes of their $W$-orbits and their multiplicities $m_\mu$; this information is taken from \cite{Lubdata}. The first column of this table numbers the $W$-orbits, in an order compatible with length in the Euclidean space containing $\Lambda$, and thus with the usual partial ordering on dominant weights, as seen in Section~\ref{sect: weights and module structure}; thereafter we let $\mu_i$ stand for any weight in the $i$th $W$-orbit. The second table is the $\alpha$-string table, whose rows correspond to the different types of $\alpha$-string which appear among the weights in $\Lambda(V)$. In each row the entries are as follows: the first column gives the type of $\alpha$-string; the second column gives the number $m$ of such $\alpha$-strings; the remaining columns give lower bounds $c(s)$ and $c(u_\alpha)$ for the contributions to $\codim V_\kappa(s)$ and $\codim C_V(u_\alpha)$ respectively, where for the former we assume (as we saw in the proof of Proposition~\ref{prop: codim bound using Premet} that we may) that $\alpha \notin \Phi(s)$.

Note that if $\mu_j$ is a weight lying between two weights $\mu_i$ in an $\alpha$-string, then as seen in Section~\ref{sect: weights and module structure} the length of $\mu_j$ is less than that of $\mu_i$, so our ordering of $W$-orbits ensures that $j < i$. Thus in any given $\alpha$-string the outermost weights lie in one $W$-orbit, with any internal weights lying in \lq lower' $W$-orbits; moreover exactly one of the outermost weights $\mu$ has $\langle \mu, \alpha \rangle \geq 0$, and for this choice of $\mu$ the number of weights in the $\alpha$-string is $\langle \mu, \alpha \rangle + 1$.

For some types of $\alpha$-string, the lower bound $c(s)$ which we are able to obtain will depend on the order $r$ of $\bar s = sZ(G)$, since two weights in the same $\alpha$-string can only lie in the same eigenspace $V_\kappa(s)$ if they differ by a multiple of $r\alpha$. Similarly, the lower bound $c(u_\alpha)$ may depend on $p$, since this may affect the way in which the sum of the weight spaces corresponding to a given $\alpha$-string decomposes into composition factors for $\langle X_{\pm\alpha} \rangle$. For this reason, the $c(s)$ and $c(u_\alpha)$ columns may often be subdivided according to the values of $r$ and $p$ respectively.

We give an example to show how the entries in the $\alpha$-string table may be calculated. Let $G = A_3$ and $\lambda = 2\omega_1 + \omega_2$ with $p \geq 3$. From \cite{Lubdata} we see that the weight table is as follows.
$$
\begin{array}{|*{4}{>{\ss}c|}}
\hline
i & \mu & |W.\mu| & m_\mu \\
\hline
3 & 2\omega_1 + \omega_2 & 12 & 1 \\
2 &       2\omega_2      &  6 & 1 \\
1 &  \omega_1 + \omega_3 & 12 & 2 \\
0 &           0          &  1 & 3 \\
\hline
\end{array}
$$
Recall from Section~\ref{sect: notation} that the root system lies in a $4$-dimensional Euclidean space with orthonormal basis $\ve_1, \ve_2, \ve_3, \ve_4$; the simple roots are $\ve_1 - \ve_2, \ve_2 - \ve_3, \ve_3 - \ve_4$, and the Weyl group acts by permuting the vectors $\ve_i$. From \cite[13.1, Table~1]{HumLie} we see that in this notation we have $2\omega_1 + \omega_2 = 2\ve_1 - \ve_3 - \ve_4$, $2\omega_2 = \ve_1 + \ve_2 - \ve_3 - \ve_4$ and $\omega_1 + \omega_3 = \ve_1 - \ve_4$. We shall represent $a_1\ve_1 + a_2\ve_2 + a_3\ve_3 + a_4\ve_4$ as $a_1a_2a_3a_4$; for convenience we write $\bar 1$ for $-1$. Thus the weights in $\Lambda(V)$ are obtained from $20{\bar 1}{\bar 1}$, $11{\bar 1}{\bar 1}$, $100{\bar 1}$, $0000$ by permuting symbols.

Let $\alpha = \alpha_1 = \ve_1 - \ve_2$; write $\mu = a_1a_2a_3a_4$, then $\langle \mu, \alpha \rangle = a_1 - a_2$. Hence any $\alpha$-string has a unique outermost weight $\mu$ with $a_1 \geq a_2$, and it contains $a_1 - a_2 + 1$ weights. We start by determining the $\alpha$-strings with outermost weights lying in the $W$-orbit containing $\lambda$ itself; we then move to the next $W$-orbit, treating only the remaining weights, and continue until all weights have been dealt with.

Initially then we consider the weights $\mu$ of the form $\mu_3$; for these we have
$$
\langle \mu, \alpha \rangle = \begin{cases}
3 & \hbox{if } \mu = 2{\bar 1}0{\bar 1} \hbox{ or } 2{\bar 1}{\bar 1}0, \\
2 & \hbox{if } \mu = 20{\bar 1}{\bar 1}, \\
1 & \hbox{if } \mu = 0{\bar 1}2{\bar 1} \hbox{ or } 0{\bar 1}{\bar 1}2, \\
0 & \hbox{if } \mu = {\bar 1}{\bar 1}20 \hbox{ or } {\bar 1}{\bar 1}02.
\end{cases}
$$
In the first possibility the two internal weights are $100{\bar 1}, 010{\bar 1}$ or $10{\bar 1}0, 01{\bar 1}0$, giving two $\alpha$-strings $\mu_3 \ \mu_1 \ \mu_1 \ \mu_3$; in the second the internal weight is $11{\bar 1}{\bar 1}$, giving one $\alpha$-string $\mu_3 \ \mu_2 \ \mu_3$; the third and fourth give two $\alpha$-strings $\mu_3 \ \mu_3$ and two $\alpha$-strings $\mu_3$. Next we consider the remaining weights $\mu$ of the form $\mu_2$; for these we have
$$
\langle \mu, \alpha \rangle = \begin{cases}
2 & \hbox{if } \mu = 1{\bar 1}1{\bar 1} \hbox{ or } 1{\bar 1}{\bar 1}1, \\
0 & \hbox{if } \mu = {\bar 1}{\bar 1}11.
\end{cases}
$$
In the first possibility the internal weight is $001{\bar 1}$ or $00{\bar 1}1$, giving two $\alpha$-strings $\mu_2 \ \mu_1 \ \mu_2$; the second gives one $\alpha$-string $\mu_2$. Now we consider the remaining weights $\mu$ of the form $\mu_1$; for these we have
$$
\langle \mu, \alpha \rangle = \begin{cases}
2 & \hbox{if } \mu = 1{\bar 1}00, \\
1 & \hbox{if } \mu = 0{\bar 1}10 \hbox{ or } 0{\bar 1}01.
\end{cases}
$$
In the first possibility the internal weight is $0000$, giving one $\alpha$-string $\mu_1 \ \mu_0 \ \mu_1$; the second gives two $\alpha$-strings $\mu_1 \ \mu_1$. As there are now no remaining weights, this completes the determination of the $\alpha$-strings.

We now turn to the lower bounds $c(s)$ and $c(u_\alpha)$; take an $\alpha$-string of type
$$
\mu_3 \ \mu_1 \ \mu_1 \ \mu_3
$$
and note that the multiplicities of the weights $\mu_1$ and $\mu_3$ are $2$ and $1$ respectively. First consider $c(s)$. If $r = 2$, the eigenspace $V_\kappa(s)$ may contain at most the first and third, or the second and fourth weights; thus we may take $c(s) = 3$. If $r = 3$, $V_\kappa(s)$ may contain at most one of the inner weights, or both of the outer weights; in either case we may take $c(s) = 4$. If however $r \geq 5$, $V_\kappa(s)$ may contain at most one weight; again we may take $c(s) = 4$. Now consider $c(u_\alpha)$; write $A$ for the $A_1$ subgroup $\langle X_{\pm\alpha} \rangle$, and regard the sum of the weight spaces as a $6$-dimensional $A$-module with weights $3\bom, \bom, \bom, -\bom, -\bom, -3\bom$, where $\bom$ is the fundamental dominant weight for $A$. If $p = 3$, there are three composition factors, with high weights $3\bom$, $\bom$ and $\bom$, on each of which $u_\alpha$ has a $1$-dimensional fixed point space; thus by Lemma~\ref{lem: submodule and fixed points} we may take $c(u_\alpha) = 3$. If however $p \geq 5$, there are two composition factors, with high weights $3\bom$ and $\bom$, on each of which $u_\alpha$ has a $1$-dimensional fixed point space; thus this time we may take $c(u_\alpha) = 4$.

The bottom row of the $\alpha$-string table sums the values $c(s)$ and $c(u_\alpha)$ to give lower bounds for $\codim V_\kappa(s)$ and $\codim C_V(u_\alpha)$. Provided the lower bounds on $\codim V_\kappa(s)$ all exceed $M$ the triple $(G, \lambda, p)$ satisfies $\ssddagcon$, while if those for various $r$ all exceed $M_r$ it satisfies $\ssdagcon$. Likewise provided the lower bound on $\codim C_V(u_\alpha)$ exceeds $M$ the triple $(G, \lambda, p)$ satisfies $\uddagcon$, while if the bound exceeds the appropriate value $M_p$ it satisfies $\udagcon$.

\begin{prop}\label{prop: A_2, 4omega_1, strings}
Let $G = A_2$ and $\lambda = 4\omega_1$ with $p \geq 5$; then the triple $(G, \lambda, p)$ satisfies $\ssdagcon$ and $\uddagcon$.
\end{prop}

\begin{proof}
The tables described above are as follows.
$$
\begin{array}{|*{4}{>{\ss}c|}}
\hline
i & \mu & |W.\mu| & m_\mu \\
\hline
4 &       4\omega_1      & 3 & 1 \\
3 & 2\omega_1 + \omega_2 & 6 & 1 \\
2 &       2\omega_2      & 3 & 1 \\
1 &        \omega_1      & 3 & 1 \\
\hline
\end{array}
\quad
\begin{array}{|*{6}{>{\ss}c|}}
\hline
 & & \multicolumn{3}{|>{\ss}c|}{c(s)} & \multicolumn{1}{|>{\ss}c|}{c(u_\alpha)} \\
\cline{3-6}
     \ss{\alpha\mathrm{-strings}}     & m & r = 2 & r = 3 & r \geq 5 & p \geq 5 \\
\hline
\mu_4 \ \mu_3 \ \mu_2 \ \mu_3 \ \mu_4 & 1 &   2   &   3   &     4    &     4    \\
                \mu_4                 & 1 &       &       &          &          \\
    \mu_3 \ \mu_1 \ \mu_1 \ \mu_3     & 1 &   2   &   2   &     3    &     3    \\
            \mu_3 \ \mu_3             & 1 &   1   &   1   &     1    &     1    \\
        \mu_2 \ \mu_1 \ \mu_2         & 1 &   1   &   2   &     2    &     2    \\
\hline
\multicolumn{2}{c|}{}                     &   6   &   8   &    10    &    10    \\
\cline{3-6}
\end{array}
$$
We have $M = 6$ and $M_2 = 4$. Thus $\codim C_V(u_\alpha) > M$, and $\codim V_\kappa(s) > M$ unless $r = 2$, in which case $\codim V_\kappa(s) > M_r$; so the triple $(G, \lambda, p)$ satisfies $\ssdagcon$ and $\uddagcon$.
\end{proof}

\begin{prop}\label{prop: A_ell, 2omega_2, strings}
Let $G = A_\ell$ for $\ell \in [4, 5]$ and $\lambda = 2\omega_2$ with $p \geq 3$; then the triple $(G, \lambda, p)$ satisfies $\ssdagcon$ and $\uddagcon$.
\end{prop}

\begin{proof}
Write $\z = \z_{p, 3}$. First suppose $\ell = 5$. In this case the tables are as follows.
$$
\begin{array}{|*{4}{>{\ss}c|}}
\hline
i & \mu & |W.\mu| & m_\mu \\
\hline
3 &       2\omega_2      & 15 &    1   \\
2 &  \omega_1 + \omega_3 & 60 &    1   \\
1 &        \omega_4      & 15 & 2 - \z \\
\hline
\end{array}
\quad
\begin{array}{|*{6}{>{\ss}c|}}
\hline
 & & \multicolumn{2}{|>{\ss}c|}{c(s)} & \multicolumn{2}{|>{\ss}c|}{c(u_\alpha)} \\
\cline{3-6}
 \ss{\alpha\mathrm{-strings}} &  m &   r = 2   & r \geq 3 & p = 3 & p \geq 5 \\
\hline
    \mu_3 \ \mu_2 \ \mu_3     &  4 &     4     &     8    &   8   &     8    \\
            \mu_3             &  7 &           &          &       &          \\
    \mu_2 \ \mu_1 \ \mu_2     &  6 & 12 - 6\z  &    12    &  12   &    12    \\
        \mu_2 \ \mu_2         & 16 &    16     &    16    &  16   &    16    \\
            \mu_2             & 12 &           &          &       &          \\
        \mu_1 \ \mu_1         &  4 &  8 - 4\z  &  8 - 4\z &   4   &     8    \\
            \mu_1             &  1 &           &          &       &          \\
\hline
\multicolumn{2}{c|}{}              & 40 - 10\z & 44 - 4\z &  40   &    44    \\
\cline{3-6}
\end{array}
$$
We have $M = 30$ and $M_2 = 18$. Thus $\codim C_V(u_\alpha) > M$, and $\codim V_\kappa(s) > M$ unless $p = 3$ and $r = 2$, in which case $\codim V_\kappa(s) > M_r$; so the triple $(G, \lambda, p)$ satisfies $\ssdagcon$ and $\uddagcon$.

Now suppose $\ell = 4$. In this case the tables are as follows.
$$
\begin{array}{|*{4}{>{\ss}c|}}
\hline
i & \mu & |W.\mu| & m_\mu \\
\hline
3 &      2\omega_2      & 10 &    1   \\
2 & \omega_1 + \omega_3 & 30 &    1   \\
1 &       \omega_4      &  5 & 2 - \z \\
\hline
\end{array}
\quad
\begin{array}{|*{6}{>{\ss}c|}}
\hline
 & & \multicolumn{2}{|>{\ss}c|}{c(s)} & \multicolumn{2}{|>{\ss}c|}{c(u_\alpha)} \\
\cline{3-6}
 \ss{\alpha\mathrm{-strings}} & m &   r = 2  & r \geq 3 & p = 3 & p \geq 5 \\
\hline
    \mu_3 \ \mu_2 \ \mu_3     & 3 &     3    &     6    &   6   &     6    \\
            \mu_3             & 4 &          &          &       &          \\
    \mu_2 \ \mu_1 \ \mu_2     & 3 &  6 - 3\z &     6    &   6   &     6    \\
        \mu_2 \ \mu_2         & 9 &     9    &     9    &   9   &     9    \\
            \mu_2             & 3 &          &          &       &          \\
        \mu_1 \ \mu_1         & 1 &  2 -  \z &  2 - \z  &   1   &     2    \\
\hline
\multicolumn{2}{c|}{}             & 20 - 4\z & 23 - \z  &  22   &    23    \\
\cline{3-6}
\end{array}
$$
We have $M = 20$ and $M_2 = 12$. Thus $\codim C_V(u_\alpha) > M$, and $\codim V_\kappa(s) > M$ unless $r = 2$, in which case $\codim V_\kappa(s) > M_r$; so the triple $(G, \lambda, p)$ satisfies $\ssdagcon$ and $\uddagcon$.
\end{proof}

\begin{prop}\label{prop: A_ell, 2omega_1 + omega_ell, strings}
Let $G = A_\ell$ for $\ell \in [2, 4]$ and $\lambda = 2\omega_1 + \omega_\ell$ with $p \geq 3$; then the triple $(G, \lambda, p)$ satisfies $\ssddagcon$ and $\uddagcon$.
\end{prop}

\begin{proof}
First suppose $\ell = 4$; write $\z = \z_{p, 3}$. In this case the tables are as follows.
$$
\begin{array}{|*{4}{>{\ss}c|}}
\hline
i & \mu & |W.\mu| & m_\mu \\
\hline
3 & 2\omega_1 + \omega_4 & 20 &    1   \\
2 &  \omega_2 + \omega_4 & 30 &    1   \\
1 &        \omega_1      &  5 & 4 - \z \\
\hline
\end{array}
\quad
\begin{array}{|*{6}{>{\ss}c|}}
\hline
 & & \multicolumn{2}{|>{\ss}c|}{c(s)} & \multicolumn{2}{|>{\ss}c|}{c(u_\alpha)} \\
\cline{3-6}
 \ss{\alpha\mathrm{-strings}} & m &  r = 2  & r \geq 3 & p = 3 & p \geq 5 \\
\hline
\mu_3 \ \mu_1 \ \mu_1 \ \mu_3 & 1 &  5 - \z &  6 - \z  &   4   &     6    \\
    \mu_3 \ \mu_2 \ \mu_3     & 3 &    3    &     6    &   6   &     6    \\
        \mu_3 \ \mu_3         & 3 &    3    &     3    &   3   &     3    \\
            \mu_3             & 6 &         &          &       &          \\
    \mu_2 \ \mu_1 \ \mu_2     & 3 &    6    &     6    &   6   &     6    \\
        \mu_2 \ \mu_2         & 9 &    9    &     9    &   9   &     9    \\
            \mu_2             & 3 &         &          &       &          \\
\hline
\multicolumn{2}{c|}{}             & 26 - \z & 30 - \z  &  28   &    30    \\
\cline{3-6}
\end{array}
$$
We have $M = 20$. Thus $\codim V_\kappa(s), \ \codim C_V(u_\alpha) > M$; so the triple $(G, \lambda, p)$ satisfies $\ssddagcon$ and $\uddagcon$.

Now suppose $\ell = 3$; write $\z = \z_{p, 5}$. In this case the tables are as follows.
$$
\begin{array}{|*{4}{>{\ss}c|}}
\hline
i & \mu & |W.\mu| & m_\mu \\
\hline
3 & 2\omega_1 + \omega_3 & 12 &    1   \\
2 &  \omega_2 + \omega_3 & 12 &    1   \\
1 &        \omega_1      &  4 & 3 - \z \\
\hline
\end{array}
\quad
\begin{array}{|*{7}{>{\ss}c|}}
\hline
 & & \multicolumn{2}{|>{\ss}c|}{c(s)} & \multicolumn{3}{|>{\ss}c|}{c(u_\alpha)} \\
\cline{3-7}
 \ss{\alpha\mathrm{-strings}} & m &  r = 2  & r \geq 3 & p = 3 & p = 5 & p \geq 7 \\
\hline
\mu_3 \ \mu_1 \ \mu_1 \ \mu_3 & 1 &  4 - \z &  5 - \z  &   4   &   4   &     5    \\
    \mu_3 \ \mu_2 \ \mu_3     & 2 &    2    &     4    &   4   &   4   &     4    \\
        \mu_3 \ \mu_3         & 2 &    2    &     2    &   2   &   2   &     2    \\
            \mu_3             & 2 &         &          &       &       &          \\
    \mu_2 \ \mu_1 \ \mu_2     & 2 &    4    &     4    &   4   &   4   &     4    \\
        \mu_2 \ \mu_2         & 3 &    3    &     3    &   3   &   3   &     3    \\
\hline
\multicolumn{2}{c|}{}             & 15 - \z & 18 - \z  &  17   &  17   &    18    \\
\cline{3-7}
\end{array}
$$
We have $M = 12$. Thus $\codim V_\kappa(s), \ \codim C_V(u_\alpha) > M$; so the triple $(G, \lambda, p)$ satisfies $\ssddagcon$ and $\uddagcon$.

Finally suppose $\ell = 2$. In this case the tables are as follows.
$$
\begin{array}{|*{4}{>{\ss}c|}}
\hline
i & \mu & |W.\mu| & m_\mu \\
\hline
3 & 2\omega_1 + \omega_2 & 6 & 1 \\
2 &       2\omega_2      & 3 & 1 \\
1 &        \omega_1      & 3 & 2 \\
\hline
\end{array}
\quad
\begin{array}{|*{6}{>{\ss}c|}}
\hline
 & & \multicolumn{2}{|>{\ss}c|}{c(s)} & \multicolumn{2}{|>{\ss}c|}{c(u_\alpha)} \\
\cline{3-6}
 \ss{\alpha\mathrm{-strings}} & m & r = 2 & r \geq 3 & p = 3 & p \geq 5 \\
\hline
\mu_3 \ \mu_1 \ \mu_1 \ \mu_3 & 1 &   3   &     4    &   3   &     4    \\
    \mu_3 \ \mu_2 \ \mu_3     & 1 &   1   &     2    &   2   &     2    \\
        \mu_3 \ \mu_3         & 1 &   1   &     1    &   1   &     1    \\
    \mu_2 \ \mu_1 \ \mu_2     & 1 &   2   &     2    &   2   &     2    \\
\hline
\multicolumn{2}{c|}{}             &   7   &     9    &   8   &     9    \\
\cline{3-6}
\end{array}
$$
We have $M = 6$. Thus $\codim V_\kappa(s), \ \codim C_V(u_\alpha) > M$; so the triple $(G, \lambda, p)$ satisfies $\ssddagcon$ and $\uddagcon$.
\end{proof}

\begin{prop}\label{prop: A_3, 2omega_1 + omega_2, strings}
Let $G = A_3$ and $\lambda = 2\omega_1 + \omega_2$ with $p \geq 3$; then the triple $(G, \lambda, p)$ satisfies $\ssddagcon$ and $\uddagcon$.
\end{prop}

\begin{proof}
The tables are as follows.
$$
\begin{array}{|*{4}{>{\ss}c|}}
\hline
i & \mu & |W.\mu| & m_\mu \\
\hline
3 & 2\omega_1 + \omega_2 & 12 & 1 \\
2 &       2\omega_2      &  6 & 1 \\
1 &  \omega_1 + \omega_3 & 12 & 2 \\
0 &           0          &  1 & 3 \\
\hline
\end{array}
\quad
\begin{array}{|*6{>{\ss}c|}}
\hline
 & & \multicolumn{2}{|>{\ss}c|}{c(s)} & \multicolumn{2}{|>{\ss}c|}{c(u_\alpha)} \\
\cline{3-6}
 \ss{\alpha\mathrm{-strings}} & m & r = 2 & r \geq 3 & p = 3 & p \geq 5 \\
\hline
\mu_3 \ \mu_1 \ \mu_1 \ \mu_3 & 2 &   6   &     8    &   6   &     8    \\
    \mu_3 \ \mu_2 \ \mu_3     & 1 &   1   &     2    &   2   &     2    \\
        \mu_3 \ \mu_3         & 2 &   2   &     2    &   2   &     2    \\
            \mu_3             & 2 &       &          &       &          \\
    \mu_2 \ \mu_1 \ \mu_2     & 2 &   4   &     4    &   4   &     4    \\
            \mu_2             & 1 &       &          &       &          \\
    \mu_1 \ \mu_0 \ \mu_1     & 1 &   3   &     4    &   4   &     4    \\
        \mu_1 \ \mu_1         & 2 &   4   &     4    &   4   &     4    \\
\hline
\multicolumn{2}{c|}{}             &  20   &    24    &  22   &    24    \\
\cline{3-6}
\end{array}
$$
We have $M = 12$. Thus $\codim V_\kappa(s), \ \codim C_V(u_\alpha) > M$; so the triple $(G, \lambda, p)$ satisfies $\ssddagcon$ and $\uddagcon$.
\end{proof}

\begin{prop}\label{prop: A_ell, omega_2 + omega_ell, strings}
Let $G = A_\ell$ for $\ell \in [6, 8]$ and $\lambda = \omega_2 + \omega_\ell$; then the triple $(G, \lambda, p)$ satisfies $\ssddagcon$ and $\udagcon$ for $\ell = 6$, and $\ssddagcon$ and $\uddagcon$ for $\ell \in [7, 8]$.
\end{prop}

\begin{proof}
First suppose $\ell = 8$; write $\z = \z_{p, 2}$. The tables are as follows.
$$
\begin{array}{|*4{>{\ss}c|}}
\hline
i & \mu & |W.\mu| & m_\mu \\
\hline
2 & \omega_2 + \omega_8 & 252 &    1   \\
1 &       \omega_1      &  9  & 7 - \z \\
\hline
\end{array}
\quad
\begin{array}{|*5{>{\ss}c|}}
\hline
 & & \multicolumn{1}{|>{\ss}c|}{c(s)} & \multicolumn{2}{|>{\ss}c|}{c(u_\alpha)} \\
\cline{3-5}
 \ss{\alpha\mathrm{-strings}} &  m  & r \geq 2 & p = 2 & p \geq 3 \\
\hline
    \mu_2 \ \mu_1 \ \mu_2     &   7 &    14    &   7   &    14    \\
        \mu_2 \ \mu_2         &  63 &    63    &  63   &    63    \\
            \mu_2             & 112 &          &       &          \\
        \mu_1 \ \mu_1         &   1 &   7 - \z &   6   &     7    \\
\hline
\multicolumn{2}{c|}{}               &  84 - \z &  76   &    84    \\
\cline{3-5}
\end{array}
$$
We have $M = 72$. Thus $\codim V_\kappa(s), \ \codim C_V(u_\alpha) > M$; so the triple $(G, \lambda, p)$ satisfies $\ssddagcon$ and $\uddagcon$.

Next suppose $\ell = 7$; write $\z = \z_{p, 7}$. The tables are as follows.
$$
\begin{array}{|*4{>{\ss}c|}}
\hline
i & \mu & |W.\mu| & m_\mu \\
\hline
2 & \omega_2 + \omega_7 & 168 &    1   \\
1 &       \omega_1      &  8  & 6 - \z \\
\hline
\end{array}
\quad
\begin{array}{|*5{>{\ss}c|}}
\hline
 & & \multicolumn{1}{|>{\ss}c|}{c(s)} & \multicolumn{2}{|>{\ss}c|}{c(u_\alpha)} \\
\cline{3-5}
 \ss{\alpha\mathrm{-strings}} &  m & r \geq 2 & p = 2 & p \geq 3 \\
\hline
    \mu_2 \ \mu_1 \ \mu_2     &  6 &    12    &   6   &    12    \\
        \mu_2 \ \mu_2         & 45 &    45    &  45   &    45    \\
            \mu_2             & 66 &          &       &          \\
        \mu_1 \ \mu_1         &  1 &   6 - \z &   6   &   6 - \z \\
\hline
\multicolumn{2}{c|}{}              &  63 - \z &  57   &  63 - \z \\
\cline{3-5}
\end{array}
$$
We have $M = 56$. Thus $\codim V_\kappa(s), \ \codim C_V(u_\alpha) > M$; so the triple $(G, \lambda, p)$ satisfies $\ssddagcon$ and $\uddagcon$.

Finally suppose $\ell = 6$; write $\z = \z_{p, 6}$. The tables are as follows.
$$
\begin{array}{|*4{>{\ss}c|}}
\hline
i & \mu & |W.\mu| & m_\mu \\
\hline
2 & \omega_2 + \omega_6 & 105 &    1   \\
1 &       \omega_1      &  7  & 5 - \z \\
\hline
\end{array}
\quad
\begin{array}{|*5{>{\ss}c|}}
\hline
 & & \multicolumn{1}{|>{\ss}c|}{c(s)} & \multicolumn{2}{|>{\ss}c|}{c(u_\alpha)} \\
\cline{3-5}
 \ss{\alpha\mathrm{-strings}} &  m & r \geq 2 & p = 2 & p \geq 3 \\
\hline
    \mu_2 \ \mu_1 \ \mu_2     &  5 &    10    &   5   &    10    \\
        \mu_2 \ \mu_2         & 30 &    30    &  30   &    30    \\
            \mu_2             & 35 &          &       &          \\
        \mu_1 \ \mu_1         &  1 &   5 - \z &   4   &   5 - \z \\
\hline
\multicolumn{2}{c|}{}              &  45 - \z &  39   &  45 - \z \\
\cline{3-5}
\end{array}
$$
We have $M = 42$ and $M_2 = 24$. Thus $\codim V_\kappa(s) > M$, and $\codim C_V(u_\alpha) > M$ unless $p = 2$, in which case $\codim C_V(u_\alpha) > M_p$; so the triple $(G, \lambda, p)$ satisfies $\ssddagcon$ and $\udagcon$.
\end{proof}

\begin{prop}\label{prop: A_5, omega_1 + omega_3, strings}
Let $G = A_5$ and $\lambda = \omega_1 + \omega_3$; then the triple $(G, \lambda, p)$ satisfies $\ssddagcon$ and $\udagcon$.
\end{prop}

\begin{proof}
Write $\z = \z_{p, 2}$. The tables are as follows.
$$
\begin{array}{|*4{>{\ss}c|}}
\hline
i & \mu & |W.\mu| & m_\mu \\
\hline
2 & \omega_1 + \omega_3 & 60 &    1   \\
1 &       \omega_4      & 15 & 3 - \z \\
\hline
\end{array}
\quad
\begin{array}{|*6{>{\ss}c|}}
\hline
 & & \multicolumn{2}{|>{\ss}c|}{c(s)} & \multicolumn{2}{|>{\ss}c|}{c(u_\alpha)} \\
\cline{3-6}
 \ss{\alpha\mathrm{-strings}} &  m & r = 2 & r \geq 3 & p = 2 & p \geq 3 \\
\hline
    \mu_2 \ \mu_1 \ \mu_2     &  6 &  12   &    12    &   6   &    12    \\
        \mu_2 \ \mu_2         & 16 &  16   &    16    &  16   &    16    \\
            \mu_2             & 16 &       &          &       &          \\
        \mu_1 \ \mu_1         &  4 &  12   & 12 - 4\z &   8   &    12    \\
            \mu_1             &  1 &       &          &       &          \\
\hline
\multicolumn{2}{c|}{}              &  40   & 40 - 4\z &  30   &    40    \\
\cline{3-6}
\end{array}
$$
We have $M = 30$ and $M_2 = 18$. Thus $\codim V_\kappa(s) > M$, and $\codim C_V(u_\alpha) > M$ unless $p = 2$, in which case $\codim C_V(u_\alpha) > M_p$; so the triple $(G, \lambda, p)$ satisfies $\ssddagcon$ and $\udagcon$.
\end{proof}

\begin{prop}\label{prop: A_4, omega_2 + omega_3, strings}
Let $G = A_4$ and $\lambda = \omega_2 + \omega_3$; then the triple $(G, \lambda, p)$ satisfies $\ssdagcon$ and $\uddagcon$.
\end{prop}

\begin{proof}
Write $\z = \z_{p, 3}$ and $\z' = \z_{p, 2}$. The tables are as follows.
$$
\begin{array}{|*4{>{\ss}c|}}
\hline
i & \mu & |W.\mu| & m_\mu \\
\hline
2 & \omega_2 + \omega_3 & 30 &       1       \\
1 & \omega_1 + \omega_4 & 20 &     2 - \z    \\
0 &          0          &  1 & 5 - 4\z - \z' \\
\hline
\end{array}
\quad
\begin{array}{|*7{>{\ss}c|}}
\hline
 & & \multicolumn{2}{|>{\ss}c|}{c(s)} & \multicolumn{3}{|>{\ss}c|}{c(u_\alpha)} \\
\cline{3-7}
 \ss{\alpha\mathrm{-strings}} & m &   r = 2   & r \geq 3 & p = 2 & p = 3 & p \geq 5 \\
\hline
    \mu_2 \ \mu_1 \ \mu_2     & 6 & 12 - 6\z  &    12    &   6   &  12   &    12    \\
        \mu_2 \ \mu_2         & 6 &     6     &     6    &   6   &   6   &     6    \\
            \mu_2             & 6 &           &          &       &       &          \\
    \mu_1 \ \mu_0 \ \mu_1     & 1 &  4 - 3\z  &  4 - 2\z &   2   &   2   &     4    \\
        \mu_1 \ \mu_1         & 6 & 12 - 6\z  & 12 - 6\z &  12   &   6   &    12    \\
\hline
\multicolumn{2}{c|}{}             & 34 - 15\z & 34 - 8\z &  26   &  26   &    34    \\
\cline{3-7}
\end{array}
$$
We have $M = 20$ and $M_2 = 12$. Thus $\codim C_V(u_\alpha) > M$, and $\codim V_\kappa(s) > M$ unless $p = 3$ and $r = 2$, in which case $\codim V_\kappa(s) > M_r$; so the triple $(G, \lambda, p)$ satisfies $\ssdagcon$ and $\uddagcon$.
\end{proof}

\begin{prop}\label{prop: D_6, omega_3, strings}
Let $G = D_6$ and $\lambda = \omega_3$; then the triple $(G, \lambda, p)$ satisfies $\ssddagcon$ and $\uddagcon$.
\end{prop}

\begin{proof}
Write $\z = \z_{p, 2}$. The tables are as follows.
$$
\begin{array}{|*4{>{\ss}c|}}
\hline
i & \mu & |W.\mu| & m_\mu \\
\hline
2 & \omega_3 & 160 &    1   \\
1 & \omega_1 &  12 & 5 - \z \\
\hline
\end{array}
\quad
\begin{array}{|*5{>{\ss}c|}}
\hline
 & & \multicolumn{1}{|>{\ss}c|}{c(s)} & \multicolumn{2}{|>{\ss}c|}{c(u_\alpha)} \\
\cline{3-5}
 \ss{\alpha\mathrm{-strings}} &  m & r \geq 2 & p = 2 & p \geq 3 \\
\hline
    \mu_2 \ \mu_1 \ \mu_2     &  8 &    16    &   8   &    16    \\
        \mu_2 \ \mu_2         & 48 &    48    &  48   &    48    \\
            \mu_2             & 48 &          &       &          \\
        \mu_1 \ \mu_1         &  2 & 10 - 2\z &   8   &    10    \\
\hline
\multicolumn{2}{c|}{}              & 74 - 2\z &  64   &    74    \\
\cline{3-5}
\end{array}
$$
We have $M = 60$. Thus $\codim V_\kappa(s), \ \codim C_V(u_\alpha) > M$; so the triple $(G, \lambda, p)$ satisfies $\ssddagcon$ and $\uddagcon$.
\end{proof}

\begin{prop}\label{prop: D_5, omega_3, strings}
Let $G = D_5$ and $\lambda = \omega_3$ with $p \geq 3$; then the triple $(G, \lambda, p)$ satisfies $\ssddagcon$ and $\uddagcon$.
\end{prop}

\begin{proof}
The tables are as follows.
$$
\begin{array}{|*4{>{\ss}c|}}
\hline
i & \mu & |W.\mu| & m_\mu \\
\hline
2 & \omega_3 & 80 & 1 \\
1 & \omega_1 & 10 & 4 \\
\hline
\end{array}
\quad
\begin{array}{|*4{>{\ss}c|}}
\hline
 & & \multicolumn{1}{|>{\ss}c|}{c(s)} & \multicolumn{1}{|>{\ss}c|}{c(u_\alpha)} \\
\cline{3-4}
 \ss{\alpha\mathrm{-strings}} &  m & r \geq 2 & p \geq 3 \\
\hline
    \mu_2 \ \mu_1 \ \mu_2     &  6 &    12    &    12    \\
        \mu_2 \ \mu_2         & 24 &    24    &    24    \\
            \mu_2             & 20 &          &          \\
        \mu_1 \ \mu_1         &  2 &     8    &     8    \\
\hline
\multicolumn{2}{c|}{}              &    44    &    44    \\
\cline{3-4}
\end{array}
$$
We have $M = 40$. Thus $\codim V_\kappa(s), \ \codim C_V(u_\alpha) > M$; so the triple $(G, \lambda, p)$ satisfies $\ssddagcon$ and $\uddagcon$.
\end{proof}

\begin{prop}\label{prop: D_5, 2omega_5, strings}
Let $G = D_5$ and $\lambda = 2\omega_5$ with $p \geq 3$; then the triple $(G, \lambda, p)$ satisfies $\ssddagcon$ and $\uddagcon$.
\end{prop}

\begin{proof}
The tables are as follows.
$$
\begin{array}{|*4{>{\ss}c|}}
\hline
i & \mu & |W.\mu| & m_\mu \\
\hline
3 & 2\omega_5 & 16 & 1 \\
2 &  \omega_3 & 80 & 1 \\
1 &  \omega_1 & 10 & 3 \\
\hline
\end{array}
\quad
\begin{array}{|*5{>{\ss}c|}}
\hline
 & & \multicolumn{2}{|>{\ss}c|}{c(s)} & \multicolumn{1}{|>{\ss}c|}{c(u_\alpha)} \\
\cline{3-5}
 \ss{\alpha\mathrm{-strings}} &  m & r = 2 & r \geq 3 & p \geq 3 \\
\hline
    \mu_3 \ \mu_2 \ \mu_3     &  4 &   4   &     8    &     8    \\
            \mu_3             &  8 &       &          &          \\
    \mu_2 \ \mu_1 \ \mu_2     &  6 &  12   &    12    &    12    \\
        \mu_2 \ \mu_2         & 24 &  24   &    24    &    24    \\
            \mu_2             & 16 &       &          &          \\
        \mu_1 \ \mu_1         &  2 &   6   &     6    &     6    \\
\hline
\multicolumn{2}{c|}{}              &  46   &    50    &    50    \\
\cline{3-5}
\end{array}
$$
We have $M = 40$. Thus $\codim V_\kappa(s), \ \codim C_V(u_\alpha) > M$; so the triple $(G, \lambda, p)$ satisfies $\ssddagcon$ and $\uddagcon$.
\end{proof}

\begin{prop}\label{prop: D_5, omega_1 + omega_5, strings}
Let $G = D_5$ and $\lambda = \omega_1 + \omega_5$; then the triple $(G, \lambda, p)$ satisfies $\ssddagcon$ and $\uddagcon$.
\end{prop}

\begin{proof}
Write $\z = \z_{p, 5}$. The tables are as follows.
$$
\begin{array}{|*4{>{\ss}c|}}
\hline
i & \mu & |W.\mu| & m_\mu \\
\hline
2 & \omega_1 + \omega_5 & 80 &    1   \\
1 &       \omega_4      & 16 & 4 - \z \\
\hline
\end{array}
\quad
\begin{array}{|*5{>{\ss}c|}}
\hline
 & & \multicolumn{1}{|>{\ss}c|}{c(s)} & \multicolumn{2}{|>{\ss}c|}{c(u_\alpha)} \\
\cline{3-5}
 \ss{\alpha\mathrm{-strings}} &  m & r \geq 2 & p = 2 & p \geq 3 \\
\hline
    \mu_2 \ \mu_1 \ \mu_2     &  8 &    16    &   8   &    16    \\
        \mu_2 \ \mu_2         & 20 &    20    &  20   &    20    \\
            \mu_2             & 24 &          &       &          \\
        \mu_1 \ \mu_1         &  4 & 16 - 4\z &  16   & 16 - 4\z \\
\hline
\multicolumn{2}{c|}{}              & 52 - 4\z &  44   & 52 - 4\z \\
\cline{3-5}
\end{array}
$$
We have $M = 40$. Thus $\codim V_\kappa(s), \ \codim C_V(u_\alpha) > M$; so the triple $(G, \lambda, p)$ satisfies $\ssddagcon$ and $\uddagcon$.
\end{proof}

We now turn to triples $(G, \lambda, p)$ where $G = B_\ell$, $C_\ell$ or $G_2$, in which $e(\Phi) > 1$; here we must consider both short and long root elements. We therefore give three tables for each case: the weight table, the $\alpha$-string table and the $\beta$-string table, of which the second and third between them provide lower bounds $c(s)$, $c(u_\alpha)$ and $c(u_\beta)$ for the contributions to the codimensions of $V_\kappa(s)$, $C_V(u_\alpha)$ and $C_V(u_\beta)$ respectively. We proceed much as in the previous cases. Note however that any short root is $\frac{1}{e(\Phi)}$ times the sum of two long roots. Thus if $p = e(\Phi)$, then for any $s \in G_{(r)}$ we may assume not only that $\alpha \notin \Phi(s)$ but also that $\beta \notin \Phi(s)$; as a result we sometimes give the $c(s)$ calculations in the $\beta$-string table rather than the $\alpha$-string table, since these may lead to better lower bounds on $\codim V_\kappa(s)$. Note also that if $p \leq e(\Phi)$ then the set $\Lambda(V)$ need not be saturated, so that some $\alpha$-strings or $\beta$-strings may have missing entries.

Again, provided the lower bounds on $\codim V_\kappa(s)$ all exceed $M$ the triple $(G, \lambda, p)$ satisfies $\ssddagcon$, while if those for various $r$ all exceed $M_r$ it satisfies $\ssdagcon$. Likewise provided the lower bounds on both $\codim C_V(u_\alpha)$ and $\codim C_V(u_\beta)$ exceed $M$ the triple $(G, \lambda, p)$ satisfies $\uddagcon$, while if they exceed the appropriate values $M_p$ it satisfies $\udagcon$. Here however there may be cases in which one of $\codim C_V(u_\alpha)$ and $\codim C_V(u_\beta)$ exceeds the appropriate bound but the other does not. If this is so, it may be possible to argue using the partial ordering on unipotent classes and Lemma~\ref{lem: unip closure containment} that the triple $(G, \lambda, p)$ satisfies $\udiamcon$.

We shall treat a few families of cases where the rank $\ell$ is unbounded. Here we will begin with values of $\ell$ up to $11$, for which the tables in \cite{Lubdata} give weight multiplicities; for these values we may proceed as before. For larger values of $\ell$ our knowledge is less complete, so we cannot provide precise tables. Instead we shall focus on one type each of $\alpha$-strings and $\beta$-strings where we do have information on the multiplicities (sometimes using Theorem~\ref{thm: Prem} or Lemma~\ref{lem: multiplicities in 3omega_1 and omega_1 + omega_2}). The corresponding entries in the tables will give lower bounds on the codimensions of $V_\kappa(s)$, $C_V(u_\alpha)$ and $C_V(u_\beta)$; as we are ignoring other types of weight string, in some cases the conditions from Section~\ref{sect: conditions} which we prove for $\ell \geq 12$ are weaker than those for $\ell \leq 11$, but they suffice for our purposes.

\begin{prop}\label{prop: B_ell, omega_3, strings}
Let $G = B_\ell$ for $\ell \in [4, \infty)$ and $\lambda = \omega_3$; then if $p \geq 3$ and $\ell \in [5, \infty)$ the triple $(G, \lambda, p)$ satisfies $\ssdagcon$ and $\uddagcon$, while if $p \geq 3$ and $\ell = 4$ it satisfies $\ssdagcon$ and $\udiamcon$; if instead $p = 2$ and $\ell \in [7, \infty)$ it satisfies $\ssddagcon$ and $\udagcon$.
\end{prop}

\begin{proof}
First suppose $p \geq 3$. If $\ell \in [4, 11]$ the tables are as follows.
$$
\begin{array}{|*4{>{\ss}c|}}
\hline
i & \mu & |W.\mu| & m_\mu \\
\hline
3 & \omega_3 & \frac{4}{3}\ell(\ell - 1)(\ell - 2) &     1    \\
2 & \omega_2 &           2\ell(\ell - 1)           &     1    \\
1 & \omega_1 &                2\ell                & \ell - 1 \\
0 &    0     &                  1                  &   \ell   \\
\hline
\end{array}
\quad
\begin{array}{|*3{>{\ss}c|}}
\hline
 & & \multicolumn{1}{|>{\ss}c|}{c(u_\beta)} \\
\cline{3-3}
 \ss{\beta\mathrm{-strings}}  &                      m                     &        p \geq 3       \\
\hline
    \mu_3 \ \mu_1 \ \mu_3     &                 2(\ell - 2)                &      4(\ell - 2)      \\
        \mu_3 \ \mu_3         &            4(\ell - 2)(\ell - 3)           & 4(\ell - 2)(\ell - 3) \\
            \mu_3             & \frac{4}{3}(\ell - 2)(\ell^2 - 7\ell + 15) &                       \\
    \mu_2 \ \mu_0 \ \mu_2     &                      1                     &           2           \\
        \mu_2 \ \mu_2         &                 4(\ell - 2)                &      4(\ell - 2)      \\
            \mu_2             &            2(\ell^2 - 5\ell + 7)           &                       \\
        \mu_1 \ \mu_1         &                      2                     &      2(\ell - 1)      \\
\hline
\multicolumn{2}{c|}{}                                                      &  4\ell^2 - 10\ell + 8 \\
\cline{3-3}
\end{array}
$$
$$
\begin{array}{|*5{>{\ss}c|}}
\hline
 & & \multicolumn{2}{|>{\ss}c|}{c(s)} & \multicolumn{1}{|>{\ss}c|}{c(u_\alpha)} \\
\cline{3-5}
 \ss{\alpha\mathrm{-strings}} &                     m                     &         r = 2         &        r \geq 3       &        p \geq 3       \\
\hline
    \mu_3 \ \mu_2 \ \mu_3     &           2(\ell - 1)(\ell - 2)           & 2(\ell - 1)(\ell - 2) & 4(\ell - 1)(\ell - 2) & 4(\ell - 1)(\ell - 2) \\
            \mu_3             & \frac{4}{3}(\ell - 1)(\ell - 2)(\ell - 3) &                       &                       &                       \\
    \mu_2 \ \mu_1 \ \mu_2     &                2(\ell - 1)                &      4(\ell - 1)      &      4(\ell - 1)      &      4(\ell - 1)      \\
    \mu_1 \ \mu_0 \ \mu_1     &                     1                     &          \ell         &      2(\ell - 1)      &      2(\ell - 1)      \\
\hline
\multicolumn{2}{c|}{}                                                     &    2\ell^2 - \ell     &  4\ell^2 - 6\ell + 2  &  4\ell^2 - 6\ell + 2  \\
\cline{3-5}
\end{array}
$$
We have $M = 2\ell^2$ and $M_2 = \ell^2 + \ell$. Thus $\codim V_\kappa(s) > M$ unless $r = 2$, in which case $\codim V_\kappa(s) > M_r$; so the triple $(G, \lambda, p)$ satisfies $\ssdagcon$. Moreover if $\ell \geq 5$ then $\codim C_V(u_\alpha), \ \codim C_V(u_\beta) > M$; so the triple $(G, \lambda, p)$ satisfies $\uddagcon$. If however $\ell = 4$ we have $\codim C_V(u_\alpha) > M$ and $\codim C_V(u_\beta) \geq M$ --- by Lemma~\ref{lem: unip closure containment}, for any unipotent class $u^G$ we have $\codim C_V(u) \geq M$, and the only unipotent class $u^G$ with $\dim u^G \geq M$ is the regular unipotent class, whose closure contains $u_\alpha$ by Lemma~\ref{lem: any class in closure of reg class}; so the triple $(G, \lambda, p)$ satisfies $\udiamcon$.

If instead $\ell \in [12, \infty)$, we consider $\alpha$-strings and $\beta$-strings of types
$$
\mu_3 \ \mu_2 \ \mu_3 \qquad \hbox{and} \qquad \mu_3 \ \mu_3
$$
respectively; note that weights $\mu_3$ and $\mu_2$ have multiplicity at least $1$ by Theorem~\ref{thm: Prem}. We have $\codim V_\kappa(s) \geq 4(\ell - 1)(\ell - 2) > M$ unless $r = 2$, in which case $\codim V_\kappa(s) \geq 2(\ell - 1)(\ell - 2) > M_r$; so the triple $(G, \lambda, p)$ satisfies $\ssdagcon$. Moreover, $\codim C_V(u_\alpha) \geq 4(\ell - 1)(\ell - 2) > M$ and $\codim C_V(u_\beta) \geq 4(\ell - 2)(\ell - 3) > M$; so the triple $(G, \lambda, p)$ satisfies $\uddagcon$.

Now suppose $p = 2$; write $\z = \z_{2, \ell - 1}$. If $\ell \in [7, 11]$  the tables are as follows.
$$
\begin{array}{|*4{>{\ss}c|}}
\hline
i & \mu & |W.\mu| & m_\mu \\
\hline
2 & \omega_3 & \frac{4}{3}\ell(\ell - 1)(\ell - 2) &       1       \\
1 & \omega_1 &                2\ell                & \ell - 2 - \z \\
\hline
\end{array}
\quad
\begin{array}{|*3{>{\ss}c|}}
\hline
 & & \multicolumn{1}{|>{\ss}c|}{c(u_\alpha)} \\
\cline{3-3}
 \ss{\alpha\mathrm{-strings}} &                     m                     &           p = 2          \\
\hline
    \mu_2 \ \phmu \ \mu_2     &           2(\ell - 1)(\ell - 2)           &  2(\ell - 1)(\ell - 2)   \\
            \mu_2             & \frac{4}{3}(\ell - 1)(\ell - 2)(\ell - 3) &                          \\
    \mu_1 \ \phmu \ \mu_1     &                     1                     &       \ell - 2 - \z      \\
            \mu_1             &                2(\ell - 1)                &                          \\
\hline
\multicolumn{2}{c|}{}                                                     & 2\ell^2 - 5\ell + 2 - \z \\
\cline{3-3}
\end{array}
$$
$$
\begin{array}{|*4{>{\ss}c|}}
\hline
 & & \multicolumn{1}{|>{\ss}c|}{c(s)} & \multicolumn{1}{|>{\ss}c|}{c(u_\beta)} \\
\cline{3-4}
 \ss{\beta\mathrm{-strings}}  &                      m                     &           r \geq 3          &            p = 2            \\
\hline
    \mu_2 \ \mu_1 \ \mu_2     &                 2(\ell - 2)                &         4(\ell - 2)         &         2(\ell - 2)         \\
        \mu_2 \ \mu_2         &            4(\ell - 2)(\ell - 3)           &    4(\ell - 2)(\ell - 3)    &    4(\ell - 2)(\ell - 3)    \\
            \mu_2             & \frac{4}{3}(\ell - 2)(\ell^2 - 7\ell + 15) &                             &                             \\
        \mu_1 \ \mu_1         &                      2                     &       2\ell - 4 - 2\z       &       2\ell - 4 - 2\z       \\
\hline
\multicolumn{2}{c|}{}                                                      & 4\ell^2 - 14\ell + 12 - 2\z & 4\ell^2 - 16\ell + 16 - 2\z \\
\cline{3-4}
\end{array}
$$
We have $M = 2\ell^2$ and $M_2 = \ell^2 +\ell$. Thus $\codim V_\kappa(s), \ \codim C_V(u_\beta) > M$, while $\codim C_V(u_\alpha) > M_p$; so the triple $(G, \lambda, p)$ satisfies $\ssddagcon$ and $\udagcon$.

If instead $\ell \in [12, \infty)$,  we consider $\alpha$-strings and $\beta$-strings of types
$$
\mu_2 \ \phmu \ \mu_2 \qquad \hbox{and} \qquad \mu_2 \ \mu_2
$$
respectively; note that weights $\mu_2$ have multiplicity $1$. We have $\codim V_\kappa(s)$, $\codim C_V(u_\beta) \geq 4(\ell - 2)(\ell - 3) > M$, while $\codim C_V(u_\alpha) \geq 2(\ell - 1)(\ell - 2) > M_p$; so the triple $(G, \lambda, p)$ satisfies $\ssddagcon$ and $\udagcon$.
\end{proof}

\begin{prop}\label{prop: B_5, omega_4, strings}
Let $G = B_5$ and $\lambda = \omega_4$ with $p = 2$; then the triple $(G, \lambda, p)$ satisfies $\ssddagcon$ and $\udagcon$.
\end{prop}

\begin{proof}
The tables are as follows.
$$
\begin{array}{|*4{>{\ss}c|}}
\hline
i & \mu & |W.\mu| & m_\mu \\
\hline
2 & \omega_4 & 80 & 1 \\
1 & \omega_2 & 40 & 2 \\
0 &    0     &  1 & 4 \\
\hline
\end{array}
\quad
\begin{array}{|*3{>{\ss}c|}}
\hline
 & & \multicolumn{1}{|>{\ss}c|}{c(u_\alpha)} \\
\cline{3-3}
 \ss{\alpha\mathrm{-strings}} &  m & p = 2 \\
\hline
    \mu_2 \ \phmu \ \mu_2     & 32 &  32   \\
            \mu_2             & 16 &       \\
    \mu_1 \ \phmu \ \mu_1     &  8 &  16   \\
            \mu_1             & 24 &       \\
            \mu_0             &  1 &       \\
\hline
\multicolumn{2}{c|}{}              &  48   \\
\cline{3-3}
\end{array}
\quad
\begin{array}{|*4{>{\ss}c|}}
\hline
 & & \multicolumn{1}{|>{\ss}c|}{c(s)} & \multicolumn{1}{|>{\ss}c|}{c(u_\beta)} \\
\cline{3-4}
 \ss{\beta\mathrm{-strings}}  &  m & r \geq 3 & p = 2 \\
\hline
    \mu_2 \ \mu_1 \ \mu_2     & 12 &    24    &  12   \\
        \mu_2 \ \mu_2         & 16 &    16    &  16   \\
            \mu_2             & 24 &          &       \\
    \mu_1 \ \mu_0 \ \mu_1     &  1 &     4    &   2   \\
        \mu_1 \ \mu_1         & 12 &    24    &  24   \\
            \mu_1             &  2 &          &       \\
\hline
\multicolumn{2}{c|}{}              &    68    &  54   \\
\cline{3-4}
\end{array}
$$
We have $M = 50$ and $M_2 = 30$. Thus $\codim V_\kappa(s), \ \codim C_V(u_\beta) > M$, while $\codim C_V(u_\alpha) > M_p$; so the triple $(G, \lambda, p)$ satisfies $\ssddagcon$ and $\udagcon$.
\end{proof}

\begin{prop}\label{prop: B_ell, 2omega_ell, strings}
Let $G = B_\ell$ for $\ell \in [3, 4]$ and $\lambda = 2\omega_\ell$ with $p \geq 3$; then if $\ell = 4$ the triple $(G, \lambda, p)$ satisfies $\ssddagcon$ and $\uddagcon$, while if $\ell = 3$ it satisfies $\ssdagcon$ and $\udiamcon$.
\end{prop}

\begin{proof}
First suppose $\ell = 4$. In this case the tables are as follows.
$$
\begin{array}{|*4{>{\ss}c|}}
\hline
i & \mu & |W.\mu| & m_\mu \\
\hline
4 & 2\omega_4 & 16 & 1 \\
3 &  \omega_3 & 32 & 1 \\
2 &  \omega_2 & 24 & 2 \\
1 &  \omega_1 &  8 & 3 \\
0 &     0     &  1 & 6 \\
\hline
\end{array}
\quad
\begin{array}{|*5{>{\ss}c|}}
\hline
 & & \multicolumn{2}{|>{\ss}c|}{c(s)} & \multicolumn{1}{|>{\ss}c|}{c(u_\alpha)} \\
\cline{3-5}
 \ss{\alpha\mathrm{-strings}} &  m & r = 2 & r \geq 3 & p \geq 3 \\
\hline
    \mu_4 \ \mu_3 \ \mu_4     &  8 &   8   &    16    &    16    \\
    \mu_3 \ \mu_2 \ \mu_3     & 12 &  24   &    24    &    24    \\
    \mu_2 \ \mu_1 \ \mu_2     &  6 &  18   &    24    &    24    \\
    \mu_1 \ \mu_0 \ \mu_1     &  1 &   6   &     6    &     6    \\
\hline
\multicolumn{2}{c|}{}              &  56   &    70    &    70    \\
\cline{3-5}
\end{array}
\quad
\begin{array}{|*3{>{\ss}c|}}
\hline
 & & \multicolumn{1}{|>{\ss}c|}{c(u_\beta)} \\
\cline{3-3}
 \ss{\beta\mathrm{-strings}}  & m & p \geq 3 \\
\hline
    \mu_4 \ \mu_2 \ \mu_4     & 4 &     8    \\
            \mu_4             & 8 &          \\
    \mu_3 \ \mu_1 \ \mu_3     & 4 &     8    \\
        \mu_3 \ \mu_3         & 8 &     8    \\
            \mu_3             & 8 &          \\
    \mu_2 \ \mu_0 \ \mu_2     & 1 &     4    \\
        \mu_2 \ \mu_2         & 8 &    16    \\
            \mu_2             & 2 &          \\
        \mu_1 \ \mu_1         & 2 &     6    \\
\hline
\multicolumn{2}{c|}{}             &    50    \\
\cline{3-3}
\end{array}
$$
We have $M = 32$. Thus $\codim V_\kappa(s), \ \codim C_V(u_\alpha), \ \codim C_V(u_\beta) > M$; so the triple $(G, \lambda, p)$ satisfies $\ssddagcon$ and $\uddagcon$.

Now suppose $\ell = 3$. In this case the tables are as follows.
$$
\begin{array}{|*4{>{\ss}c|}}
\hline
i & \mu & |W.\mu| & m_\mu \\
\hline
3 & 2\omega_3 &  8 & 1 \\
2 &  \omega_2 & 12 & 1 \\
1 &  \omega_1 &  6 & 2 \\
0 &     0     &  1 & 3 \\
\hline
\end{array}
\quad
\begin{array}{|*5{>{\ss}c|}}
\hline
 & & \multicolumn{2}{|>{\ss}c|}{c(s)} & \multicolumn{1}{|>{\ss}c|}{c(u_\alpha)} \\
\cline{3-5}
 \ss{\alpha\mathrm{-strings}} & m & r = 2 & r \geq 3 & p \geq 3 \\
\hline
    \mu_3 \ \mu_2 \ \mu_3     & 4 &   4   &     8    &     8    \\
    \mu_2 \ \mu_1 \ \mu_2     & 4 &   8   &     8    &     8    \\
    \mu_1 \ \mu_0 \ \mu_1     & 1 &   3   &     4    &     4    \\
\hline
\multicolumn{2}{c|}{}             &  15   &    20    &    20    \\
\cline{3-5}
\end{array}
\quad
\begin{array}{|*3{>{\ss}c|}}
\hline
 & & \multicolumn{1}{|>{\ss}c|}{c(u_\beta)} \\
\cline{3-3}
 \ss{\beta\mathrm{-strings}}  & m & p \geq 3 \\
\hline
    \mu_3 \ \mu_1 \ \mu_3     & 2 &     4    \\
            \mu_3             & 4 &          \\
    \mu_2 \ \mu_0 \ \mu_2     & 1 &     2    \\
        \mu_2 \ \mu_2         & 4 &     4    \\
            \mu_2             & 2 &          \\
        \mu_1 \ \mu_1         & 2 &     4    \\
\hline
\multicolumn{2}{c|}{}             &    14    \\
\cline{3-3}
\end{array}
$$
We have $M = 18$ and $M_2 = 12$. Thus $\codim V_\kappa(s) > M$ unless $r = 2$, in which case $\codim V_\kappa(s) > M_r$; so the triple $(G, \lambda, p)$ satisfies $\ssdagcon$. Moreover $\codim C_V(u_\alpha) > M$, and $\codim C_V(u_\beta) > 8 = \dim {u_\beta}^G$ --- Lemma~\ref{lem: root elt class in closure of any non-triv class} shows that all non-identity unipotent classes other than ${u_\beta}^G$ contain $u_\alpha$ in their closure; so the triple $(G, \lambda, p)$ satisfies $\udiamcon$.
\end{proof}

\begin{prop}\label{prop: B_2, 3omega_2, strings}
Let $G = B_2$ and $\lambda = 3\omega_2$ with $p \geq 5$; then the triple $(G, \lambda, p)$ satisfies $\ssddagcon$ and $\uddagcon$.
\end{prop}

\begin{proof}
The tables are as follows.
$$
\begin{array}{|*4{>{\ss}c|}}
\hline
i & \mu & |W.\mu| & m_\mu \\
\hline
3 &      3\omega_2      & 4 & 1 \\
2 & \omega_1 + \omega_2 & 8 & 1 \\
1 &       \omega_2      & 4 & 2 \\
\hline
\end{array}
\quad
\begin{array}{|*6{>{\ss}c|}}
\hline
 & & \multicolumn{3}{|>{\ss}c|}{c(s)} & \multicolumn{1}{|>{\ss}c|}{c(u_\alpha)} \\
\cline{3-6}
 \ss{\alpha\mathrm{-strings}} & m & r = 2 & r = 3 & r \geq 5 & p \geq 5 \\
\hline
\mu_3 \ \mu_2 \ \mu_2 \ \mu_3 & 2 &   4   &   4   &     6    &     6    \\
\mu_2 \ \mu_1 \ \mu_1 \ \mu_2 & 2 &   6   &   8   &     8    &     8    \\
\hline
\multicolumn{2}{c|}{}             &  10   &  12   &    14    &    14    \\
\cline{3-6}
\end{array}
\quad
\begin{array}{|*3{>{\ss}c|}}
\hline
 & & \multicolumn{1}{|>{\ss}c|}{c(u_\beta)} \\
\cline{3-3}
 \ss{\beta\mathrm{-strings}}  & m & p \geq 5 \\
\hline
\mu_3 \ \mu_1 \ \mu_1 \ \mu_3 & 1 &     4    \\
            \mu_3             & 2 &          \\
    \mu_2 \ \mu_1 \ \mu_2     & 2 &     4    \\
        \mu_2 \ \mu_2         & 2 &     2    \\
\hline
\multicolumn{2}{c|}{}             &    10    \\
\cline{3-3}
\end{array}
$$
We have $M = 8$. Thus $\codim V_\kappa(s), \ \codim C_V(u_\alpha), \ \codim C_V(u_\beta) > M$; so the triple $(G, \lambda, p)$ satisfies $\ssddagcon$ and $\uddagcon$.
\end{proof}

\begin{prop}\label{prop: B_ell, omega_1 + omega_2, strings}
Let $G = B_\ell$ for $\ell \in [3, \infty)$ and $\lambda = \omega_1 + \omega_2$ with $p \geq 3$; then if $\ell \in [12, \infty)$ and $p = 3$ the triple $(G, \lambda, p)$ satisfies $\ssdagcon$ and $\uddagcon$, while otherwise it satisfies $\ssddagcon$ and $\uddagcon$.
\end{prop}

\begin{proof}
First suppose $p \geq 5$; write $\z = \z_{p, \ell}$. If $\ell \in [3, 11]$ the tables are as follows.
$$
\begin{array}{|*{4}{>{\ss}c|}}
\hline
i & \mu & |W.\mu| & m_\mu \\
\hline
5 & \omega_1 + \omega_2 &           4\ell(\ell - 1)           &        1       \\
4 &      2\omega_1      &                2\ell                &        1       \\
3 &   \pa \omega_3 {}^* & \frac{4}{3}\ell(\ell - 1)(\ell - 2) &        2       \\
2 &       \omega_2      &           2\ell(\ell - 1)           &        2       \\
1 &       \omega_1      &                2\ell                & 2\ell - 1 - \z \\
0 &          0          &                  1                  & 2\ell - 1 - \z \\
\hline
\multicolumn{4}{>{\ss}l}{\quad {}^* 2\omega_3 \ \mathrm{if} \ \ell = 3} \\
\end{array}
\quad
\begin{array}{|*{3}{>{\ss}c|}}
\hline
 & & \multicolumn{1}{|>{\ss}c|}{c(u_\beta)} \\
\cline{3-3}
 \ss{\beta\mathrm{-strings}}  &                     m                     &          p \geq 5         \\
\hline
\mu_5 \ \mu_1 \ \mu_1 \ \mu_5 &                     2                     &      4\ell + 2 - 2\z      \\
    \mu_5 \ \mu_3 \ \mu_5     &                4(\ell - 2)                &        8(\ell - 2)        \\
        \mu_5 \ \mu_5         &                2(2\ell - 3)               &        2(2\ell - 3)       \\
            \mu_5             &           4(\ell - 2)(\ell - 3)           &                           \\
    \mu_4 \ \mu_2 \ \mu_4     &                     2                     &             4             \\
            \mu_4             &                2(\ell - 2)                &                           \\
    \mu_3 \ \mu_1 \ \mu_3     &                2(\ell - 2)                &        8(\ell - 2)        \\
        \mu_3 \ \mu_3         &           4(\ell - 2)(\ell - 3)           &   8(\ell - 2)(\ell - 3)   \\
            \mu_3             & \frac{4}{3}(\ell - 2)(\ell - 3)(\ell - 4) &                           \\
    \mu_2 \ \mu_0 \ \mu_2     &                     1                     &             4             \\
        \mu_2 \ \mu_2         &                4(\ell - 2)                &        8(\ell - 2)        \\
            \mu_2             &           2(\ell - 2)(\ell - 3)           &                           \\
\hline
\multicolumn{2}{c|}{}                                                     & 8\ell^2 - 8\ell + 4 - 2\z \\
\cline{3-3}
\end{array}
$$
\vspace{-0.5mm}
$$
\begin{array}{|*{6}{>{\ss}c|}}
\hline
 & & \multicolumn{3}{|>{\ss}c|}{c(s)} & \multicolumn{1}{|>{\ss}c|}{c(u_\alpha)} \\
\cline{3-6}
     \ss{\alpha\mathrm{-strings}}     &                     m                     &         r = 2         &           r = 3           &        r \geq 5       &        p \geq 5       \\
\hline
\mu_5 \ \mu_2 \ \mu_1 \ \mu_2 \ \mu_5 &                2(\ell - 1)                &      8(\ell - 1)      &        12(\ell - 1)       &      12(\ell - 1)     &      12(\ell - 1)     \\
        \mu_5 \ \mu_4 \ \mu_5         &                2(\ell - 1)                &      2(\ell - 1)      &        4(\ell - 1)        &      4(\ell - 1)      &      4(\ell - 1)      \\
                \mu_5                 &           4(\ell - 1)(\ell - 2)           &                       &                           &                       &                       \\
\mu_4 \ \mu_1 \ \mu_0 \ \mu_1 \ \mu_4 &                     1                     &     2\ell + 1 - \z    &      4\ell - 1 - 2\z      &      4\ell - 2\z      &      4\ell - 2\z      \\
        \mu_3 \ \mu_2 \ \mu_3         &           2(\ell - 1)(\ell - 2)           & 4(\ell - 1)(\ell - 2) &   8(\ell - 1)(\ell - 2)   & 8(\ell - 1)(\ell - 2) & 8(\ell - 1)(\ell - 2) \\
                \mu_3                 & \frac{4}{3}(\ell - 1)(\ell - 2)(\ell - 3) &                       &                           &                       &                       \\
\hline
\multicolumn{2}{c|}{}                                                             &    4\ell^2 - 1 - \z   & 8\ell^2 - 4\ell - 1 - 2\z & 8\ell^2 - 4\ell - 2\z & 8\ell^2 - 4\ell - 2\z \\
\cline{3-6}
\end{array}
$$
\vspace{-0.5mm}
We have $M = 2\ell^2$. Thus $\codim V_\kappa(s), \ \codim C_V(u_\alpha), \ \codim C_V(u_\beta) > M$; so the triple $(G, \lambda, p)$ satisfies $\ssddagcon$ and $\uddagcon$.

\vspace{-0.5mm}

If instead $\ell \in [12, \infty)$, we consider $\alpha$-strings and $\beta$-strings of types
$$
\mu_3 \ \mu_2 \ \mu_3 \qquad \hbox{and} \qquad \mu_3 \ \mu_3
$$
respectively; note that weights $\mu_3$ have multiplicity $2$ by Lemma~\ref{lem: multiplicities in 3omega_1 and omega_1 + omega_2}, and as the weight spaces corresponding to these $\alpha$-strings must decompose into composition \ factors \ for \ $\langle X_{\pm\alpha} \rangle$, \ weights \ $\mu_2$ \ must \ have \ multiplicity \ at \ least \ $2$. \ We \ have $\codim V_\kappa(s) \geq 4(\ell - 1)(\ell - 2) > M$, while $\codim C_V(u_\alpha) \geq 8(\ell - 1)(\ell - 2) > M$ and $\codim C_V(u_\beta) \geq 8(\ell - 2)(\ell - 3) > M$; so the triple $(G, \lambda, p)$ satisfies $\ssddagcon$ and $\uddagcon$.

\vspace{-0.5mm}

Now suppose $p = 3$; write $\z = \z_{3, \ell}$. If $\ell \in [3, 11]$ the tables are as follows.
\vspace{-0.5mm}
$$
\begin{array}{|*4{>{\ss}c|}}
\hline
i & \mu & |W.\mu| & m_\mu \\
\hline
5 & \omega_1 + \omega_2 &           4\ell(\ell - 1)           &       1       \\
4 &      2\omega_1      &                2\ell                &       1       \\
3 &   \pa \omega_3 {}^* & \frac{4}{3}\ell(\ell - 1)(\ell - 2) &       1       \\
2 &       \omega_2      &           2\ell(\ell - 1)           &       1       \\
1 &       \omega_1      &                2\ell                &   \ell - \z   \\
0 &          0          &                  1                  & \ell - 1 - \z \\
\hline
\multicolumn{4}{>{\ss}l}{\quad {}^* 2\omega_3 \ \mathrm{if} \ \ell = 3} \\
\end{array}
\quad
\begin{array}{|*3{>{\ss}c|}}
\hline
 & & \multicolumn{1}{|>{\ss}c|}{c(u_\beta)} \\
\cline{3-3}
 \ss{\beta\mathrm{-strings}}  &                     m                     &           p = 3           \\
\hline
\mu_5 \ \mu_1 \ \mu_1 \ \mu_5 &                     2                     &      2\ell + 2 - 2\z      \\
    \mu_5 \ \mu_3 \ \mu_5     &                4(\ell - 2)                &        8(\ell - 2)        \\
        \mu_5 \ \mu_5         &                2(2\ell - 3)               &        2(2\ell - 3)       \\
            \mu_5             &           4(\ell - 2)(\ell - 3)           &                           \\
    \mu_4 \ \mu_2 \ \mu_4     &                     2                     &             4             \\
            \mu_4             &                2(\ell - 2)                &                           \\
    \mu_3 \ \mu_1 \ \mu_3     &                2(\ell - 2)                &        4(\ell - 2)        \\
        \mu_3 \ \mu_3         &           4(\ell - 2)(\ell - 3)           &   4(\ell - 2)(\ell - 3)   \\
            \mu_3             & \frac{4}{3}(\ell - 2)(\ell - 3)(\ell - 4) &                           \\
    \mu_2 \ \mu_0 \ \mu_2     &                     1                     &             2             \\
        \mu_2 \ \mu_2         &                4(\ell - 2)                &        4(\ell - 2)        \\
            \mu_2             &           2(\ell - 2)(\ell - 3)           &                           \\
\hline
\multicolumn{2}{c|}{}                                                     & 4\ell^2 + 2\ell - 6 - 2\z \\
\cline{3-3}
\end{array}
$$
$$
\begin{array}{|*5{>{\ss}c|}}
\hline
 & & \multicolumn{2}{|>{\ss}c|}{c(s)} & \multicolumn{1}{|>{\ss}c|}{c(u_\alpha)} \\
\cline{3-5}
     \ss{\alpha\mathrm{-strings}}     &                     m                     &          r = 2          &          r \geq 5         &          p = 3        \\
\hline
\mu_5 \ \mu_2 \ \mu_1 \ \mu_2 \ \mu_5 &                2(\ell - 1)                &       4(\ell - 1)       &        8(\ell - 1)        &      4(\ell - 1)      \\
        \mu_5 \ \mu_4 \ \mu_5         &                2(\ell - 1)                &       2(\ell - 1)       &        4(\ell - 1)        &      4(\ell - 1)      \\
                \mu_5                 &           4(\ell - 1)(\ell - 2)           &                         &                           &                       \\
\mu_4 \ \mu_1 \ \mu_0 \ \mu_1 \ \mu_4 &                     1                     &      \ell + 1 - \z      &      2\ell + 1 - 2\z      &      2\ell - 2\z      \\
        \mu_3 \ \mu_2 \ \mu_3         &           2(\ell - 1)(\ell - 2)           &  2(\ell - 1)(\ell - 2)  &   4(\ell - 1)(\ell - 2)   & 4(\ell - 1)(\ell - 2) \\
                \mu_3                 & \frac{4}{3}(\ell - 1)(\ell - 2)(\ell - 3) &                         &                           &                       \\
\hline
\multicolumn{2}{c|}{}                                                             & 2\ell^2 + \ell - 1 - \z & 4\ell^2 + 2\ell - 3 - 2\z & 4\ell^2 - 2\ell - 2\z \\
\cline{3-5}
\end{array}
$$
We have $M = 2\ell^2$. Thus $\codim V_\kappa(s), \ \codim C_V(u_\alpha), \ \codim C_V(u_\beta) > M$; so the triple $(G, \lambda, p)$ satisfies $\ssddagcon$ and $\uddagcon$.

If instead $\ell \in [12, \infty)$, we consider $\alpha$-strings and $\beta$-strings of types
$$
\mu_3 \ \mu_2 \ \mu_3 \qquad \hbox{and} \qquad \mu_3 \ \mu_3
$$
respectively; note that weights $\mu_3$ have multiplicity $1$ by Lemma~\ref{lem: multiplicities in 3omega_1 and omega_1 + omega_2}, and weights $\mu_2$ have multiplicity at least $1$ by Theorem~\ref{thm: Prem}. We have $\codim V_\kappa(s) \geq 4(\ell - 1)(\ell - 2) > M$ unless $r = 2$, in which case $\codim V_\kappa(s) \geq 2(\ell - 1)(\ell - 2) > M_r$; so the triple $(G, \lambda, p)$ satisfies $\ssdagcon$. Moreover, $\codim C_V(u_\alpha) \geq 4(\ell - 1)(\ell - 2) > M$ and $\codim C_V(u_\beta) \geq 4(\ell - 2)(\ell - 3) > M$; so the triple $(G, \lambda, p)$ satisfies $\uddagcon$.
\end{proof}

\begin{prop}\label{prop: B_2, omega_1 + omega_2, strings}
Let $G = B_2$ and $\lambda = \omega_1 + \omega_2$ with $p \neq 5$; then the triple $(G, \lambda, p)$ satisfies $\ssdagcon$ and $\udiamcon$.
\end{prop}

\begin{proof}
The tables are as follows.
$$
\begin{array}{|*4{>{\ss}c|}}
\hline
i & \mu & |W.\mu| & m_\mu \\
\hline
2 & \omega_1 + \omega_2 & 8 & 1 \\
1 &       \omega_2      & 4 & 2 \\
\hline
\end{array}
$$
$$
\begin{array}{|*7{>{\ss}c|}}
\hline
 & & \multicolumn{2}{|>{\ss}c|}{c(s)} & \multicolumn{3}{|>{\ss}c|}{c(u_\alpha)} \\
\cline{3-7}
 \ss{\alpha\mathrm{-strings}} & m & r = 2 & r \geq 3 & p = 2 & p = 3 & p \geq 7 \\
\hline
\mu_2 \ \mu_1 \ \mu_1 \ \mu_2 & 2 &   6   &     8    &   6   &   6   &     8    \\
        \mu_2 \ \mu_2         & 2 &   2   &     2    &   2   &   2   &     2    \\
\hline
\multicolumn{2}{c|}{}             &   8   &    10    &   8   &   8   &    10    \\
\cline{3-7}
\end{array}
\quad
\begin{array}{|*4{>{\ss}c|}}
\hline
 & & \multicolumn{2}{|>{\ss}c|}{c(u_\beta)} \\
\cline{3-4}
 \ss{\beta\mathrm{-strings}}  & m & p = 2 & p \geq 3 \\
\hline
    \mu_2 \ \mu_1 \ \mu_2     & 2 &   2   &     4    \\
        \mu_2 \ \mu_2         & 2 &   2   &     2    \\
        \mu_1 \ \mu_1         & 1 &   2   &     2    \\
\hline
\multicolumn{2}{c|}{}             &   6   &     8    \\
\cline{3-4}
\end{array}
$$
We have $M = 8$ and $M_2 = M_3 = 6$. Thus $\codim V_\kappa(s) > M$ unless $r = 2$, in which case $\codim V_\kappa(s) > M_r$; so the triple $(G, \lambda, p)$ satisfies $\ssdagcon$. Moreover $\codim C_V(u_\alpha) > M$ unless $p \in \{ 2, 3 \}$, in which case $\codim C_V(u_\alpha) > M_p$, while $\codim C_V(u_\beta) > 4 = \dim {u_\beta}^G$ --- Lemma~\ref{lem: root elt class in closure of any non-triv class} shows that all non-identity unipotent classes other than ${u_\beta}^G$ contain $u_\alpha$ in their closure; so the triple $(G, \lambda, p)$ satisfies $\udiamcon$.
\end{proof}

\begin{prop}\label{prop: B_3, omega_1 + omega_3, strings}
Let $G = B_3$ and $\lambda = \omega_1 + \omega_3$; then if $p \neq 7$ the triple $(G, \lambda, p)$ satisfies $\ssddagcon$ and $\udagcon$, while if $p = 7$ it satisfies $\ssddagcon$ and $\udiamcon$.
\end{prop}

\begin{proof}
Write $\z = \z_{p, 7}$. The tables are as follows.
$$
\begin{array}{|*4{>{\ss}c|}}
\hline
i & \mu & |W.\mu| & m_\mu \\
\hline
2 & \omega_1 + \omega_3 & 24 & 1 \\
1 &       \omega_3      &  8 & 3 - \z \\
\hline
\end{array}
\quad
\begin{array}{|*4{>{\ss}c|}}
\hline
 & & \multicolumn{2}{|>{\ss}c|}{c(u_\beta)} \\
\cline{3-4}
 \ss{\beta\mathrm{-strings}}  & m & p = 2 & p \geq 3 \\
\hline
    \mu_2 \ \mu_1 \ \mu_2     & 4 &   4   &     8    \\
        \mu_2 \ \mu_2         & 6 &   6   &     6    \\
            \mu_2             & 4 &       &          \\
        \mu_1 \ \mu_1         & 2 &   6   &  6 - 2\z \\
\hline
\multicolumn{2}{c|}{}             &  16   & 20 - 2\z \\
\cline{3-4}
\end{array}
$$
$$
\begin{array}{|*7{>{\ss}c|}}
\hline
 & & \multicolumn{2}{|>{\ss}c|}{c(s)} & \multicolumn{3}{|>{\ss}c|}{c(u_\alpha)} \\
\cline{3-7}
 \ss{\alpha\mathrm{-strings}} & m &   r = 2  & r \geq 3 & p = 2 & p = 3 & p \geq 5 \\
\hline
\mu_2 \ \mu_1 \ \mu_1 \ \mu_2 & 4 & 16 - 4\z & 20 - 4\z &  16   &  16   & 20 - 4\z \\
        \mu_2 \ \mu_2         & 8 &     8    &     8    &   8   &   8   &     8    \\
\hline
\multicolumn{2}{c|}{}             & 24 - 4\z & 28 - 4\z &  24   &  24   & 28 - 4\z \\
\cline{3-7}
\end{array}
$$
We have $M = 18$ and $M_2 = 12$. Thus $\codim V_\kappa(s) > M$; so the triple $(G, \lambda, p)$ satisfies $\ssddagcon$. Moreover $\codim C_V(u_\alpha) > M$, while $\codim C_V(u_\beta) > M$ unless either $p = 2$, in which case $\codim C_V(u_\beta) > M_p$, or $p = 7$, in which case $\codim C_V(u_\beta) \geq M$ --- if $p = 7$, by Lemma~\ref{lem: unip closure containment}, for any unipotent class $u^G$ we have $\codim C_V(u) \geq M$, and the only unipotent class $u^G$ with $\dim u^G \geq M$ is the regular unipotent class, whose closure contains $u_\alpha$ by Lemma~\ref{lem: any class in closure of reg class}; so if $p \neq 7$ then the triple $(G, \lambda, p)$ satisfies $\udagcon$, while if $p = 7$ it satisfies $\udiamcon$.
\end{proof}

\begin{prop}\label{prop: B_4, omega_1 + omega_4, strings}
Let $G = B_4$ and $\lambda = \omega_1 + \omega_4$ with $p \geq 3$; then the triple $(G, \lambda, p)$ satisfies $\ssddagcon$ and $\uddagcon$.
\end{prop}

\begin{proof}
Write $\z = \z_{p, 3}$. The tables are as follows.
$$
\begin{array}{|*4{>{\ss}c|}}
\hline
i & \mu & |W.\mu| & m_\mu \\
\hline
2 & \omega_1 + \omega_4 & 64 & 1 \\
1 &       \omega_4      & 16 & 4 - \z \\
\hline
\end{array}
$$
$$
\begin{array}{|*6{>{\ss}c|}}
\hline
 & & \multicolumn{2}{|>{\ss}c|}{c(s)} & \multicolumn{2}{|>{\ss}c|}{c(u_\alpha)} \\
\cline{3-6}
 \ss{\alpha\mathrm{-strings}} &  m &   r = 2  &  r \geq 3 & p = 3 & p \geq 5 \\
\hline
\mu_2 \ \mu_1 \ \mu_1 \ \mu_2 &  8 & 40 - 8\z &  48 - 8\z &  32   &    48    \\
        \mu_2 \ \mu_2         & 24 &    24    &     24    &  24   &    24    \\
\hline
\multicolumn{2}{c|}{}              & 64 - 8\z &  72 - 8\z &  56   &    72    \\
\cline{3-6}
\end{array}
\quad
\begin{array}{|*4{>{\ss}c|}}
\hline
 & & \multicolumn{2}{|>{\ss}c|}{c(u_\beta)} \\
\cline{3-4}
 \ss{\beta\mathrm{-strings}}  &  m & p = 3 & p \geq 5 \\
\hline
    \mu_2 \ \mu_1 \ \mu_2     &  8 &  16   &    16    \\
        \mu_2 \ \mu_2         & 16 &  16   &    16    \\
            \mu_2             & 16 &       &          \\
        \mu_1 \ \mu_1         &  4 &  12   &    16    \\
\hline
\multicolumn{2}{c|}{}              &  44   &    48    \\
\cline{3-4}
\end{array}
$$
We have $M = 32$. Thus $\codim V_\kappa(s), \ \codim C_V(u_\alpha), \ \codim C_V(u_\beta) > M$; so the triple $(G, \lambda, p)$ satisfies $\ssddagcon$ and $\uddagcon$.
\end{proof}

\begin{prop}\label{prop: B_2, omega_1 + 2omega_2, strings}
Let $G = B_2$ and $\lambda = \omega_1 + 2\omega_2$ with $p \geq 3$; then the triple $(G, \lambda, p)$ satisfies $\ssddagcon$ and $\uddagcon$.
\end{prop}

\begin{proof}
First suppose $p \geq 5$. In this case the tables are as follows.
$$
\begin{array}{|*4{>{\ss}c|}}
\hline
i & \mu & |W.\mu| & m_\mu \\
\hline
4 & \omega_1 + 2\omega_2 & 8 & 1 \\
3 &      2\omega_1       & 4 & 1 \\
2 &      2\omega_2       & 4 & 2 \\
1 &       \omega_1       & 4 & 3 \\
0 &          0           & 1 & 3 \\
\hline
\end{array}
$$
$$
\begin{array}{|*6{>{\ss}c|}}
\hline
 & & \multicolumn{3}{|>{\ss}c|}{c(s)} & \multicolumn{1}{|>{\ss}c|}{c(u_\alpha)} \\
\cline{3-6}
     \ss{\alpha\mathrm{-strings}}     & m & r = 2 & r = 3 & r \geq 5 & p \geq 5 \\
\hline
\mu_4 \ \mu_2 \ \mu_1 \ \mu_2 \ \mu_4 & 2 &   8   &  12   &    12    &    12    \\
        \mu_4 \ \mu_3 \ \mu_4         & 2 &   2   &   4   &     4    &     4    \\
\mu_3 \ \mu_1 \ \mu_0 \ \mu_1 \ \mu_3 & 1 &   5   &   7   &     8    &     8    \\
\hline
\multicolumn{2}{c|}{}                     &  15   &  23   &    24    &    24    \\
\cline{3-6}
\end{array}
\quad
\begin{array}{|*3{>{\ss}c|}}
\hline
 & & \multicolumn{1}{|>{\ss}c|}{c(u_\beta)} \\
\cline{3-3}
 \ss{\beta\mathrm{-strings}}  & m & p \geq 5 \\
\hline
\mu_4 \ \mu_1 \ \mu_1 \ \mu_4 & 2 &    10    \\
        \mu_4 \ \mu_4         & 2 &     2    \\
    \mu_3 \ \mu_2 \ \mu_3     & 2 &     4    \\
    \mu_2 \ \mu_0 \ \mu_2     & 1 &     4    \\
\hline
\multicolumn{2}{c|}{}             &    20    \\
\cline{3-3}
\end{array}
$$
We have $M = 8$. Thus $\codim V_\kappa(s), \ \codim C_V(u_\alpha), \ \codim C_V(u_\beta) > M$; so the triple $(G, \lambda, p)$ satisfies $\ssddagcon$ and $\uddagcon$.

Now suppose $p = 3$. In this case the tables are as follows.
$$
\begin{array}{|*4{>{\ss}c|}}
\hline
i & \mu & |W.\mu| & m_\mu \\
\hline
4 & \omega_1 + 2\omega_2 & 8 & 1 \\
3 &      2\omega_1       & 4 & 1 \\
2 &      2\omega_2       & 4 & 1 \\
1 &       \omega_1       & 4 & 2 \\
0 &          0           & 1 & 1 \\
\hline
\end{array}
\quad
\begin{array}{|*5{>{\ss}c|}}
\hline
 & & \multicolumn{2}{|>{\ss}c|}{c(s)} & \multicolumn{1}{|>{\ss}c|}{c(u_\alpha)} \\
\cline{3-5}
     \ss{\alpha\mathrm{-strings}}     & m & r = 2 & r \geq 5 & p = 3 \\
\hline
\mu_4 \ \mu_2 \ \mu_1 \ \mu_2 \ \mu_4 & 2 &   4   &     8    &   4   \\
        \mu_4 \ \mu_3 \ \mu_4         & 2 &   2   &     4    &   4   \\
\mu_3 \ \mu_1 \ \mu_0 \ \mu_1 \ \mu_3 & 1 &   3   &     5    &   4   \\
\hline
\multicolumn{2}{c|}{}                     &   9   &    17    &  12   \\
\cline{3-5}
\end{array}
\quad
\begin{array}{|*3{>{\ss}c|}}
\hline
 & & \multicolumn{1}{|>{\ss}c|}{c(u_\beta)} \\
\cline{3-3}
 \ss{\beta\mathrm{-strings}}  & m & p = 3 \\
\hline
\mu_4 \ \mu_1 \ \mu_1 \ \mu_4 & 2 &   6   \\
        \mu_4 \ \mu_4         & 2 &   2   \\
    \mu_3 \ \mu_2 \ \mu_3     & 2 &   4   \\
    \mu_2 \ \mu_0 \ \mu_2     & 1 &   2   \\
\hline
\multicolumn{2}{c|}{}             &  14   \\
\cline{3-3}
\end{array}
$$
We have $M = 8$. Thus $\codim V_\kappa(s), \ \codim C_V(u_\alpha), \ \codim C_V(u_\beta) > M$; so the triple $(G, \lambda, p)$ satisfies $\ssddagcon$ and $\uddagcon$.
\end{proof}

\begin{prop}\label{prop: C_ell, 3omega_1, strings}
Let $G = C_\ell$ for $\ell \in [3, \infty)$ and $\lambda = 3\omega_1$ with $p \geq 5$; then if $\ell \in [3, 11]$ the triple $(G, \lambda, p)$ satisfies $\ssddagcon$ and $\uddagcon$, while if $\ell \in [12, \infty)$ it satisfies $\ssddagcon$ and $\udiamcon$.
\end{prop}

\begin{proof}
If $\ell \in [3, 11]$ the tables are as follows.
$$
\begin{array}{|*4{>{\ss}c|}}
\hline
i & \mu & |W.\mu| & m_\mu \\
\hline
4 &      3\omega_1      &                2\ell                &   1  \\
3 & \omega_1 + \omega_2 &           4\ell(\ell - 1)           &   1  \\
2 &       \omega_3      & \frac{4}{3}\ell(\ell - 1)(\ell - 2) &   1  \\
1 &       \omega_1      &                2\ell                & \ell \\
\hline
\end{array}
\quad
\begin{array}{|*3{>{\ss}c|}}
\hline
 & & \multicolumn{1}{|>{\ss}c|}{c(u_\beta)} \\
\cline{3-3}
 \ss{\beta\mathrm{-strings}}  &                     m                     &        p \geq 5       \\
\hline
\mu_4 \ \mu_1 \ \mu_1 \ \mu_4 &                     1                     &        \ell + 2       \\
            \mu_4             &                2(\ell - 1)                &                       \\
    \mu_3 \ \mu_1 \ \mu_3     &                2(\ell - 1)                &      4(\ell - 1)      \\
        \mu_3 \ \mu_3         &                2(\ell - 1)                &      2(\ell - 1)      \\
            \mu_3             &           4(\ell - 1)(\ell - 2)           &                       \\
        \mu_2 \ \mu_2         &           2(\ell - 1)(\ell - 2)           & 2(\ell - 1)(\ell - 2) \\
            \mu_2             & \frac{4}{3}(\ell - 1)(\ell - 2)(\ell - 3) &                       \\
\hline
\multicolumn{2}{c|}{}                                                     &     2\ell^2 + \ell    \\
\cline{3-3}
\end{array}
$$
$$
\begin{array}{|*6{>{\ss}c|}}
\hline
 & & \multicolumn{3}{|>{\ss}c|}{c(s)} & \multicolumn{1}{|>{\ss}c|}{c(u_\alpha)} \\
\cline{3-6}
 \ss{\alpha\mathrm{-strings}} &                     m                     &         r = 2         &         r = 3         &        r \geq 5       &        p \geq 5       \\
\hline
\mu_4 \ \mu_3 \ \mu_3 \ \mu_4 &                     2                     &           4           &           4           &           6           &           6           \\
            \mu_4             &                2(\ell - 2)                &                       &                       &                       &                       \\
\mu_3 \ \mu_1 \ \mu_1 \ \mu_3 &                     2                     &      2(\ell + 1)      &      2(\ell + 2)      &      2(\ell + 2)      &      2(\ell + 2)      \\
    \mu_3 \ \mu_2 \ \mu_3     &                4(\ell - 2)                &      4(\ell - 2)      &      8(\ell - 2)      &      8(\ell - 2)      &      8(\ell - 2)      \\
        \mu_3 \ \mu_3         &                4(\ell - 2)                &      4(\ell - 2)      &      4(\ell - 2)      &      4(\ell - 2)      &      4(\ell - 2)      \\
            \mu_3             &           4(\ell - 2)(\ell - 3)           &                       &                       &                       &                       \\
    \mu_2 \ \mu_1 \ \mu_2     &                2(\ell - 2)                &      4(\ell - 2)      &      4(\ell - 2)      &      4(\ell - 2)      &      4(\ell - 2)      \\
        \mu_2 \ \mu_2         &           4(\ell - 2)(\ell - 3)           & 4(\ell - 2)(\ell - 3) & 4(\ell - 2)(\ell - 3) & 4(\ell - 2)(\ell - 3) & 4(\ell - 2)(\ell - 3) \\
            \mu_2             & \frac{4}{3}(\ell - 2)(\ell - 3)(\ell - 4) &                       &                       &                       &                       \\
\hline
\multicolumn{2}{c|}{}                                                     &  4\ell^2 - 6\ell + 6  &    4\ell^2 - 2\ell    &  4\ell^2 - 2\ell + 2  &  4\ell^2 - 2\ell + 2  \\
\cline{3-6}
\end{array}
$$
We have $M = 2\ell^2$. Thus $\codim V_\kappa(s), \ \codim C_V(u_\alpha), \ \codim C_V(u_\beta) > M$; so the triple $(G, \lambda, p)$ satisfies $\ssddagcon$ and $\uddagcon$.

If instead $\ell \in [12, \infty)$, we consider both $\alpha$-strings and $\beta$-strings of type
$$
\mu_2 \ \mu_2;
$$
note that weights $\mu_2$ have multiplicity $1$ by Lemma~\ref{lem: multiplicities in 3omega_1 and omega_1 + omega_2}. We have $\codim V_\kappa(s) \geq 4(\ell - 2)(\ell - 3) > M$; so the triple $(G, \lambda, p)$ satisfies $\ssddagcon$. Moreover, $\codim C_V(u_\alpha) \geq 4(\ell - 2)(\ell - 3) > M$ and $\codim C_V(u_\beta) \geq 2(\ell - 1)(\ell - 2) > 2\ell = \dim {u_\beta}^G$ --- Lemma~\ref{lem: root elt class in closure of any non-triv class} shows that all non-identity unipotent classes other than ${u_\beta}^G$ contain $u_\alpha$ in their closure; so the triple $(G, \lambda, p)$ satisfies $\udiamcon$.
\end{proof}

\begin{prop}\label{prop: C_ell, omega_3, strings}
Let $G = C_\ell$ for $\ell \in [7, \infty)$ and $\lambda = \omega_3$; then the triple $(G, \lambda, p)$ satisfies $\ssddagcon$ and $\udiamcon$.
\end{prop}

\begin{proof}
Write $\z = \z_{p, \ell - 1}$. If $\ell \in [7, 11]$ the tables are as follows.
$$
\begin{array}{|*4{>{\ss}c|}}
\hline
i & \mu & |W.\mu| & m_\mu \\
\hline
2 & \omega_3 & \frac{4}{3}\ell(\ell - 1)(\ell - 2) &       1       \\
1 & \omega_1 &                2\ell                & \ell - 2 - \z \\
\hline
\end{array}
\quad
\begin{array}{|*3{>{\ss}c|}}
\hline
 & & \multicolumn{1}{|>{\ss}c|}{c(u_\beta)} \\
\cline{3-3}
 \ss{\beta\mathrm{-strings}}  &                     m                     &         p \geq 2         \\
\hline
        \mu_2 \ \mu_2         &           2(\ell - 1)(\ell - 2)           &   2(\ell - 1)(\ell - 2)  \\
            \mu_2             & \frac{4}{3}(\ell - 1)(\ell - 2)(\ell - 3) &                          \\
        \mu_1 \ \mu_1         &                     1                     &       \ell - 2 - \z      \\
            \mu_1             &                2(\ell - 1)                &                          \\
\hline
\multicolumn{2}{c|}{}                                                     & 2\ell^2 - 5\ell + 2 - \z \\
\cline{3-3}
\end{array}
$$
$$
\begin{array}{|*5{>{\ss}c|}}
\hline
 & & \multicolumn{1}{|>{\ss}c|}{c(s)} & \multicolumn{2}{|>{\ss}c|}{c(u_\alpha)} \\
\cline{3-5}
 \ss{\alpha\mathrm{-strings}} &                      m                     &           r \geq 2          &            p = 2            &           p \geq 3          \\
\hline
    \mu_2 \ \mu_1 \ \mu_2     &                 2(\ell - 2)                &         4(\ell - 2)         &         2(\ell - 2)         &         4(\ell - 2)         \\
        \mu_2 \ \mu_2         &            4(\ell - 2)(\ell - 3)           &    4(\ell - 2)(\ell - 3)    &    4(\ell - 2)(\ell - 3)    &    4(\ell - 2)(\ell - 3)    \\
            \mu_2             & \frac{4}{3}(\ell - 2)(\ell^2 - 7\ell + 15) &                             &                             &                             \\
        \mu_1 \ \mu_1         &                      2                     &       2\ell - 4 - 2\z       &       2\ell - 4 - 2\z       &       2\ell - 4 - 2\z       \\
\hline
\multicolumn{2}{c|}{}                                                      & 4\ell^2 - 14\ell + 12 - 2\z & 4\ell^2 - 16\ell + 16 - 2\z & 4\ell^2 - 14\ell + 12 - 2\z \\
\cline{3-5}
\end{array}
$$
We have $M = 2\ell^2$ and $M_2 = \ell^2 + \ell$. Thus $\codim V_\kappa(s) > M$; so the triple $(G, \lambda, p)$ satisfies $\ssddagcon$. Moreover $\codim C_V(u_\alpha) > M$ unless $\ell = 7$ and $p = 2$, in which case $\codim C_V(u_\alpha) > M_p$, and $\codim C_V(u_\beta) > 2\ell = \dim {u_\beta}^G$ --- Lemma~\ref{lem: root elt class in closure of any non-triv class} shows that all non-identity unipotent classes other than ${u_\beta}^G$ contain $u_\alpha$ in their closure; so the triple $(G, \lambda, p)$ satisfies $\udiamcon$.

If instead $\ell \in [12, \infty)$, we consider both $\alpha$-strings and $\beta$-strings of type
$$
\mu_2 \ \mu_2;
$$
note that weights $\mu_2$ have multiplicity $1$. We have $\codim V_\kappa(s) \geq 4(\ell - 2)(\ell - 3) > M$; so the triple $(G, \lambda, p)$ satisfies $\ssddagcon$. Moreover, $\codim C_V(u_\alpha) \geq 4(\ell - 2)(\ell - 3) > M$ and $\codim C_V(u_\beta) \geq 2(\ell - 1)(\ell - 2) > 2\ell = \dim {u_\beta}^G$ --- Lemma~\ref{lem: root elt class in closure of any non-triv class} shows that all non-identity unipotent classes other than ${u_\beta}^G$ contain $u_\alpha$ in their closure; so the triple $(G, \lambda, p)$ satisfies $\udiamcon$.
\end{proof}

\begin{prop}\label{prop: C_5, omega_4, strings}
Let $G = C_5$ and $\lambda = \omega_4$; then if $p \neq 3$ the triple $(G, \lambda, p)$ satisfies $\ssddagcon$ and $\udiamcon$, while if $p = 3$ it satisfies $\ssdagcon$ and $\udiamcon$.
\end{prop}

\begin{proof}
First suppose $p \neq 3$; write $\z = \z_{p, 2}$. In this case the tables are as follows.
$$
\begin{array}{|*4{>{\ss}c|}}
\hline
i & \mu & |W.\mu| & m_\mu \\
\hline
2 & \omega_4 & 80 & 1 \\
1 & \omega_2 & 40 & 2 \\
0 &    0     &  1 & 5 - \z \\
\hline
\end{array}
\quad
\begin{array}{|*5{>{\ss}c|}}
\hline
 & & \multicolumn{1}{|>{\ss}c|}{c(s)} & \multicolumn{2}{|>{\ss}c|}{c(u_\alpha)} \\
\cline{3-5}
 \ss{\alpha\mathrm{-strings}} &  m & r \geq 2 & p = 2 & p \geq 5 \\
\hline
    \mu_2 \ \mu_1 \ \mu_2     & 12 &    24    &   12  &    24    \\
        \mu_2 \ \mu_2         & 16 &    16    &   16  &    16    \\
            \mu_2             & 24 &          &       &          \\
    \mu_1 \ \mu_0 \ \mu_1     &  1 &     4    &    2  &     4    \\
        \mu_1 \ \mu_1         & 12 &    24    &   24  &    24    \\
            \mu_1             &  2 &          &       &          \\
\hline
\multicolumn{2}{c|}{}              &    68    &   54  &    68    \\
\cline{3-5}
\end{array}
\quad
\begin{array}{|*3{>{\ss}c|}}
\hline
 & & \multicolumn{1}{|>{\ss}c|}{c(u_\beta)} \\
\cline{3-3}
 \ss{\beta\mathrm{-strings}}  &  m & p \neq 3 \\
\hline
        \mu_2 \ \mu_2         & 32 &    32    \\
            \mu_2             & 16 &          \\
       \mu_1 \ \mu_1         &  8 &    16    \\
            \mu_1             & 24 &          \\
             \mu_0             &  1 &          \\
\hline
\multicolumn{2}{c|}{}              &    48    \\
\cline{3-3}
\end{array}
$$
We have $M = 50$. Thus $\codim V_\kappa(s) > M$; so the triple $(G, \lambda, p)$ satisfies $\ssddagcon$. Moreover $\codim C_V(u_\alpha) > M$ and $\codim C_V(u_\beta) > 10 = \dim {u_\beta}^G$ --- Lemma~\ref{lem: root elt class in closure of any non-triv class} shows that all non-identity unipotent classes other than ${u_\beta}^G$ contain $u_\alpha$ in their closure; so the triple $(G, \lambda, p)$ satisfies $\udiamcon$.

Now suppose $p = 3$. In this case the tables are as follows.
$$
\begin{array}{|*4{>{\ss}c|}}
\hline
i & \mu & |W.\mu| & m_\mu \\
\hline
2 & \omega_4 & 80 & 1 \\
1 & \omega_2 & 40 & 1 \\
0 &    0     &  1 & 1 \\
\hline
\end{array}
\quad
\begin{array}{|*5{>{\ss}c|}}
\hline
 & & \multicolumn{2}{|>{\ss}c|}{c(s)} & \multicolumn{1}{|>{\ss}c|}{c(u_\alpha)} \\
\cline{3-5}
 \ss{\alpha\mathrm{-strings}} &  m & r = 2 & r \geq 5 & p = 3 \\
\hline
    \mu_2 \ \mu_1 \ \mu_2     & 12 &  12   &    24    &   24  \\
        \mu_2 \ \mu_2         & 16 &  16   &    16    &   16  \\
            \mu_2             & 24 &       &          &       \\
    \mu_1 \ \mu_0 \ \mu_1     &  1 &   1   &     2    &    2  \\
        \mu_1 \ \mu_1         & 12 &  12   &    12    &   12  \\
            \mu_1             &  2 &       &          &       \\
\hline
\multicolumn{2}{c|}{}              &  41   &    54    &   54  \\
\cline{3-5}
\end{array}
\quad
\begin{array}{|*3{>{\ss}c|}}
\hline
 & & \multicolumn{1}{|>{\ss}c|}{c(u_\beta)} \\
\cline{3-3}
 \ss{\beta\mathrm{-strings}}  &  m & p = 3 \\
\hline
        \mu_2 \ \mu_2         & 32 &   32  \\
            \mu_2             & 16 &       \\
        \mu_1 \ \mu_1         &  8 &    8  \\
            \mu_1             & 24 &       \\
            \mu_0             &  1 &       \\
\hline
\multicolumn{2}{c|}{}              &   40  \\
\cline{3-3}
\end{array}
$$
We have $M = 50$ and $M_2 = 30$. Thus $\codim V_\kappa(s) > M$ unless $r = 2$, in which case $\codim V_\kappa(s) > M_r$; so the triple $(G, \lambda, p)$ satisfies $\ssdagcon$. Moreover $\codim C_V(u_\alpha) > M$ and $\codim C_V(u_\beta) > 10 = \dim {u_\beta}^G$ --- Lemma~\ref{lem: root elt class in closure of any non-triv class} shows that all non-identity unipotent classes other than ${u_\beta}^G$ contain $u_\alpha$ in their closure; so the triple $(G, \lambda, p)$ satisfies $\udiamcon$.
\end{proof}

\begin{prop}\label{prop: C_5, omega_5, strings}
Let $G = C_5$ and $\lambda = \omega_5$ with $p \geq 3$; then the triple $(G, \lambda, p)$ satisfies $\ssdagcon$ and $\udiamcon$.
\end{prop}

\begin{proof}
Write $\z = \z_{p, 3}$. The tables are as follows.
$$
\begin{array}{|*4{>{\ss}c|}}
\hline
i & \mu & |W.\mu| & m_\mu \\
\hline
3 & \omega_5 & 32 & 1 \\
2 & \omega_3 & 80 & 1 \\
1 & \omega_1 & 10 & 2 - \z \\
\hline
\end{array}
\quad
\begin{array}{|*6{>{\ss}c|}}
\hline
 & & \multicolumn{2}{|>{\ss}c|}{c(s)} & \multicolumn{2}{|>{\ss}c|}{c(u_\alpha)} \\
\cline{3-6}
 \ss{\alpha\mathrm{-strings}} &  m &   r = 2  & r \geq 3 & p = 3 & p \geq 5 \\
\hline
    \mu_3 \ \mu_2 \ \mu_3     &  8 &     8    &    16    &  16   &    16    \\
            \mu_3             & 16 &          &          &       &          \\
    \mu_2 \ \mu_1 \ \mu_2     &  6 & 12 - 6\z &    12    &  12   &    12    \\
        \mu_2 \ \mu_2         & 24 &    24    &    24    &  24   &    24    \\
            \mu_2             & 12 &          &          &       &          \\
        \mu_1 \ \mu_1         &  2 &  4 - 2\z &  4 - 2\z &   2   &     4    \\
\hline
\multicolumn{2}{c|}{}              & 48 - 8\z & 56 - 2\z &  54   &    56    \\
\cline{3-6}
\end{array}
\quad
\begin{array}{|*4{>{\ss}c|}}
\hline
 & & \multicolumn{2}{|>{\ss}c|}{c(u_\beta)} \\
\cline{3-4}
 \ss{\beta\mathrm{-strings}}  &  m & p = 3 & p \geq 5 \\
\hline
        \mu_3 \ \mu_3         & 16 &  16  &    16    \\
        \mu_2 \ \mu_2         & 24 &  24  &    24    \\
            \mu_2             & 32 &      &          \\
        \mu_1 \ \mu_1         &  1 &   1  &     2    \\
            \mu_1             &  8 &      &          \\
\hline
\multicolumn{2}{c|}{}              &  41  &    42    \\
\cline{3-4}
\end{array}
$$
We have $M = 50$ and $M_2 = 30$. Thus $\codim V_\kappa(s) > M$ unless $r = 2$, in which case $\codim V_\kappa(s) > M_r$; so the triple $(G, \lambda, p)$ satisfies $\ssdagcon$. Moreover $\codim C_V(u_\alpha) > M$ and $\codim C_V(u_\beta) > 10 = \dim {u_\beta}^G$ --- Lemma~\ref{lem: root elt class in closure of any non-triv class} shows that all non-identity unipotent classes other than ${u_\beta}^G$ contain $u_\alpha$ in their closure; so the triple $(G, \lambda, p)$ satisfies $\udiamcon$.
\end{proof}

\begin{prop}\label{prop: C_ell, omega_1 + omega_2, strings}
Let $G = C_\ell$ for $\ell \in [3, \infty)$ and $\lambda = \omega_1 + \omega_2$ with $p \geq 3$; then if $\ell \in [12, \infty)$ and $p = 3$ the triple $(G, \lambda, p)$ satisfies $\ssddagcon$ and $\udiamcon$, while otherwise it satisfies $\ssddagcon$ and $\uddagcon$.
\end{prop}

\begin{proof}
First suppose $p \geq 5$; write $\z = \z_{p, 2\ell + 1}$. If $\ell \in [3, 11]$ the tables are as follows.
$$
\begin{array}{|*4{>{\ss}c|}}
\hline
i & \mu & |W.\mu| & m_\mu \\
\hline
3 & \omega_1 + \omega_2 &           4\ell(\ell - 1)           &        1       \\
2 &       \omega_3      & \frac{4}{3}\ell(\ell - 1)(\ell - 2) &        2       \\
1 &       \omega_1      &                2\ell                & 2\ell - 2 - \z \\
\hline
\end{array}
\quad
\begin{array}{|*3{>{\ss}c|}}
\hline
 & & \multicolumn{1}{|>{\ss}c|}{c(u_\beta)} \\
\cline{3-3}
 \ss{\beta\mathrm{-strings}}  &                     m                     &        p \geq 5       \\
\hline
    \mu_3 \ \mu_1 \ \mu_3     &                2(\ell - 1)                &      4(\ell - 1)      \\
        \mu_3 \ \mu_3         &                2(\ell - 1)                &      2(\ell - 1)      \\
            \mu_3             &           4(\ell - 1)(\ell - 2)           &                       \\
        \mu_2 \ \mu_2         &           2(\ell - 1)(\ell - 2)           & 4(\ell - 1)(\ell - 2) \\
            \mu_2             & \frac{4}{3}(\ell - 1)(\ell - 2)(\ell - 3) &                       \\
        \mu_1 \ \mu_1         &                     1                     &     2\ell - 2 - \z    \\
\hline
\multicolumn{2}{c|}{}                                                     & 4\ell^2 - 4\ell - \z  \\
\cline{3-3}
\end{array}
$$
$$
\begin{array}{|*5{>{\ss}c|}}
\hline
 & & \multicolumn{2}{|>{\ss}c|}{c(s)} & \multicolumn{1}{|>{\ss}c|}{c(u_\alpha)} \\
\cline{3-5}
 \ss{\alpha\mathrm{-strings}} &                     m                     &            r = 2           &           r \geq 3          &           p \geq 5          \\
\hline
\mu_3 \ \mu_1 \ \mu_1 \ \mu_3 &                     2                     &       4\ell - 2 - 2\z      &         4\ell - 2\z         &         4\ell - 2\z         \\
    \mu_3 \ \mu_2 \ \mu_3     &                4(\ell - 2)                &         8(\ell - 2)        &         8(\ell - 2)         &         8(\ell - 2)         \\
        \mu_3 \ \mu_3         &                2(2\ell - 3)               &        2(2\ell - 3)        &         2(2\ell - 3)        &         2(2\ell - 3)        \\
            \mu_3             &           4(\ell - 2)(\ell - 3)           &                            &                             &                             \\
    \mu_2 \ \mu_1 \ \mu_2     &                2(\ell - 2)                &      8(\ell - 2) - 2\z     &         8(\ell - 2)         &         8(\ell - 2)         \\
        \mu_2 \ \mu_2         &           4(\ell - 2)(\ell - 3)           &    8(\ell - 2)(\ell - 3)   &    8(\ell - 2)(\ell - 3)    &    8(\ell - 2)(\ell - 3)    \\
            \mu_2             & \frac{4}{3}(\ell - 2)(\ell - 3)(\ell - 4) &                            &                             &                             \\
\hline
\multicolumn{2}{c|}{}                                                     & 8\ell^2 - 16\ell + 8 - 4\z & 8\ell^2 - 16\ell + 10 - 2\z & 8\ell^2 - 16\ell + 10 - 2\z \\
\cline{3-5}
\end{array}
$$
We have $M = 2\ell^2$. Thus $\codim V_\kappa(s), \ \codim C_V(u_\alpha), \ \codim C_V(u_\beta) > M$; so the triple $(G, \lambda, p)$ satisfies $\ssddagcon$ and $\uddagcon$.

If instead $\ell \in [12, \infty)$, we consider both $\alpha$-strings and $\beta$-strings of type
$$
\mu_2 \ \mu_2;
$$
note that weights $\mu_2$ have multiplicity $2$ by Lemma~\ref{lem: multiplicities in 3omega_1 and omega_1 + omega_2}. We have $\codim V_\kappa(s)$, $\codim C_V(u_\alpha) \geq 8(\ell - 2)(\ell - 3) > M$, and $\codim C_V(u_\beta) \geq 4(\ell - 1)(\ell - 2) > M$; so the triple $(G, \lambda, p)$ satisfies $\ssddagcon$ and $\uddagcon$.

Now suppose $p = 3$. If $\ell \in [3, 11]$ the tables are as follows.
$$
\begin{array}{|*4{>{\ss}c|}}
\hline
i & \mu & |W.\mu| & m_\mu \\
\hline
3 & \omega_1 + \omega_2 &           4\ell(\ell - 1)           &   1  \\
2 &       \omega_3      & \frac{4}{3}\ell(\ell - 1)(\ell - 2) &   1  \\
1 &       \omega_1      &                2\ell                & \ell \\
\hline
\end{array}
\quad
\begin{array}{|*3{>{\ss}c|}}
\hline
 & & \multicolumn{1}{|>{\ss}c|}{c(u_\beta)} \\
\cline{3-3}
 \ss{\beta\mathrm{-strings}}  &                     m                     &         p = 3         \\
\hline
    \mu_3 \ \mu_1 \ \mu_3     &                2(\ell - 1)                &      4(\ell - 1)      \\
        \mu_3 \ \mu_3         &                2(\ell - 1)                &      2(\ell - 1)      \\
            \mu_3             &           4(\ell - 1)(\ell - 2)           &                       \\
        \mu_2 \ \mu_2         &           2(\ell - 1)(\ell - 2)           & 2(\ell - 1)(\ell - 2) \\
            \mu_2             & \frac{4}{3}(\ell - 1)(\ell - 2)(\ell - 3) &                       \\
        \mu_1 \ \mu_1         &                     1                     &          \ell         \\
\hline
\multicolumn{2}{c|}{}                                                     &   2\ell^2 + \ell - 2  \\
\cline{3-3}
\end{array}
$$
$$
\begin{array}{|*5{>{\ss}c|}}
\hline
 & & \multicolumn{2}{|>{\ss}c|}{c(s)} & \multicolumn{1}{|>{\ss}c|}{c(u_\alpha)} \\
\cline{3-5}
 \ss{\alpha\mathrm{-strings}} &                     m                     &         r = 2         &        r \geq 5       &         p = 3         \\
\hline
\mu_3 \ \mu_1 \ \mu_1 \ \mu_3 &                     2                     &      2(\ell + 1)      &      2(\ell + 2)      &      2(\ell + 1)      \\
    \mu_3 \ \mu_2 \ \mu_3     &                4(\ell - 2)                &      4(\ell - 2)      &      8(\ell - 2)      &      8(\ell - 2)      \\
        \mu_3 \ \mu_3         &                2(2\ell - 3)               &      2(2\ell - 3)     &      2(2\ell - 3)     &      2(2\ell - 3)     \\
            \mu_3             &           4(\ell - 2)(\ell - 3)           &                       &                       &                       \\
    \mu_2 \ \mu_1 \ \mu_2     &                2(\ell - 2)                &      4(\ell - 2)      &      4(\ell - 2)      &      4(\ell - 2)      \\
        \mu_2 \ \mu_2         &           4(\ell - 2)(\ell - 3)           & 4(\ell - 2)(\ell - 3) & 4(\ell - 2)(\ell - 3) & 4(\ell - 2)(\ell - 3) \\
            \mu_2             & \frac{4}{3}(\ell - 2)(\ell - 3)(\ell - 4) &                       &                       &                       \\
\hline
\multicolumn{2}{c|}{}                                                     &  4\ell^2 - 6\ell + 4  &  4\ell^2 - 2\ell - 2  &  4\ell^2 - 2\ell - 4  \\
\cline{3-5}
\end{array}
$$
We have $M = 2\ell^2$. Thus $\codim V_\kappa(s), \ \codim C_V(u_\alpha), \ \codim C_V(u_\beta) > M$; so the triple $(G, \lambda, p)$ satisfies $\ssddagcon$ and $\uddagcon$.

If instead $\ell \in [12, \infty)$, we consider both $\alpha$-strings and $\beta$-strings of type
$$
\mu_2 \ \mu_2;
$$
note that weights $\mu_2$ have multiplicity $1$ by Lemma~\ref{lem: multiplicities in 3omega_1 and omega_1 + omega_2}. We have $\codim V_\kappa(s) \geq 4(\ell - 2)(\ell - 3) > M$; so the triple $(G, \lambda, p)$ satisfies $\ssddagcon$. Moreover, $\codim C_V(u_\alpha) \geq 4(\ell - 2)(\ell - 3) > M$ and $\codim C_V(u_\beta) \geq 2(\ell - 1)(\ell - 2) > 2\ell = \dim {u_\beta}^G$ --- Lemma~\ref{lem: root elt class in closure of any non-triv class} shows that all non-identity unipotent classes other than ${u_\beta}^G$ contain $u_\alpha$ in their closure; so the triple $(G, \lambda, p)$ satisfies $\udiamcon$.
\end{proof}

\begin{prop}\label{prop: C_3, omega_1 + omega_3, strings}
Let $G = C_3$ and $\lambda = \omega_1 + \omega_3$; then the triple $(G, \lambda, p)$ satisfies $\ssddagcon$ and $\uddagcon$.
\end{prop}

\begin{proof}
First suppose $p \geq 3$; write $\z = \z_{p, 3}$. In this case the tables are as follows.
$$
\begin{array}{|*4{>{\ss}c|}}
\hline
i & \mu & |W.\mu| & m_\mu \\
\hline
3 & \omega_1 + \omega_3 & 24 & 1 \\
2 &      2\omega_1      &  6 & 1 \\
1 &       \omega_2      & 12 & 3 - \z \\
0 &          0          &  1 & 4 - \z \\
\hline
\end{array}
\quad
\begin{array}{|*6{>{\ss}c|}}
\hline
 & & \multicolumn{2}{|>{\ss}c|}{c(s)} & \multicolumn{2}{|>{\ss}c|}{c(u_\alpha)} \\
\cline{3-6}
 \ss{\alpha\mathrm{-strings}} & m &   r = 2  & r \geq 3 & p = 3 & p \geq 5 \\
\hline
\mu_3 \ \mu_1 \ \mu_1 \ \mu_3 & 4 & 16 - 4\z & 20 - 4\z &  12   &    20    \\
    \mu_3 \ \mu_2 \ \mu_3     & 2 &     2    &     4    &   4   &     4    \\
        \mu_3 \ \mu_3         & 4 &     4    &     4    &   4   &     4    \\
            \mu_3             & 4 &          &          &       &          \\
    \mu_2 \ \mu_1 \ \mu_2     & 2 &     4    &     4    &   4   &     4    \\
    \mu_1 \ \mu_0 \ \mu_1     & 1 &  4 - \z  &  6 - 2\z &   4   &     6    \\
\hline
\multicolumn{2}{c|}{}             & 30 - 5\z & 38 - 6\z &  28   &    38    \\
\cline{3-6}
\end{array}
\quad
\begin{array}{|*3{>{\ss}c|}}
\hline
 & & \multicolumn{1}{|>{\ss}c|}{c(u_\beta)} \\
\cline{3-3}
 \ss{\beta\mathrm{-strings}}  & m & p \geq 3 \\
\hline
    \mu_3 \ \mu_1 \ \mu_3     & 4 &     8    \\
        \mu_3 \ \mu_3         & 8 &     8    \\
    \mu_2 \ \mu_0 \ \mu_2     & 1 &     2    \\
            \mu_2             & 4 &          \\
        \mu_1 \ \mu_1         & 4 & 12 - 4\z \\
\hline
\multicolumn{2}{c|}{}             & 30 - 4\z \\
\cline{3-3}
\end{array}
$$
We have $M = 18$. Thus $\codim V_\kappa(s), \ \codim C_V(u_\alpha), \ \codim C_V(u_\beta) > M$; so the triple $(G, \lambda, p)$ satisfies $\ssddagcon$ and $\uddagcon$.

Now suppose $p = 2$. In this case the tables are as follows.
$$
\begin{array}{|*4{>{\ss}c|}}
\hline
i & \mu & |W.\mu| & m_\mu \\
\hline
2 & \omega_1 + \omega_3 & 24 & 1 \\
1 &       \omega_2      & 12 & 2 \\
\hline
\end{array}
\quad
\begin{array}{|*4{>{\ss}c|}}
\hline
 & & \multicolumn{1}{|>{\ss}c|}{c(s)} & \multicolumn{1}{|>{\ss}c|}{c(u_\alpha)} \\
\cline{3-4}
 \ss{\alpha\mathrm{-strings}} & m & r \geq 3 & p = 2 \\
\hline
\mu_2 \ \mu_1 \ \mu_1 \ \mu_2 & 4 &    16    &  12   \\
    \mu_2 \ \phmu \ \mu_2     & 2 &     2    &   2   \\
        \mu_2 \ \mu_2         & 4 &     4    &   4   \\
            \mu_2             & 4 &          &       \\
    \mu_1 \ \phmu \ \mu_1     & 1 &     2    &   2   \\
            \mu_1             & 2 &          &       \\
\hline
\multicolumn{2}{c|}{}             &    24    &  20   \\
\cline{3-4}
\end{array}
\quad
\begin{array}{|*3{>{\ss}c|}}
\hline
 & & \multicolumn{1}{|>{\ss}c|}{c(u_\beta)} \\
\cline{3-3}
 \ss{\beta\mathrm{-strings}}  & m & p = 2 \\
\hline
    \mu_2 \ \mu_1 \ \mu_2     & 4 &   4   \\
        \mu_2 \ \mu_2         & 8 &   8   \\
        \mu_1 \ \mu_1         & 4 &   8   \\
\hline
\multicolumn{2}{c|}{}             &  20   \\
\cline{3-3}
\end{array}
$$
We have $M = 18$. Thus $\codim V_\kappa(s), \ \codim C_V(u_\alpha), \ \codim C_V(u_\beta) > M$; so the triple $(G, \lambda, p)$ satisfies $\ssddagcon$ and $\uddagcon$.
\end{proof}

\begin{prop}\label{prop: G_2, 2omega_1, strings}
Let $G = G_2$ and $\lambda = 2\omega_1$ with $p \geq 3$; then the triple $(G, \lambda, p)$ satisfies $\ssdagcon$ and $\udiamcon$.
\end{prop}

\begin{proof}
Write $\z = \z_{p, 7}$. The tables are as follows.
$$
\begin{array}{|*4{>{\ss}c|}}
\hline
i & \mu & |W.\mu| & m_\mu \\
\hline
3 & 2\omega_1 &  6 & 1 \\
2 &  \omega_2 &  6 & 1 \\
1 &  \omega_1 &  6 & 2 \\
0 &     0     &  1 & 3 - \z \\
\hline
\end{array}
\quad
\begin{array}{|*7{>{\ss}c|}}
\hline
 & & \multicolumn{3}{|>{\ss}c|}{c(s)} & \multicolumn{2}{|>{\ss}c|}{c(u_\alpha)} \\
\cline{3-7}
     \ss{\alpha\mathrm{-strings}}     & m & r = 2 &  r = 3  & r \geq 5 & p = 3 & p \geq 5 \\
\hline
\mu_3 \ \mu_1 \ \mu_0 \ \mu_1 \ \mu_3 & 1 &   4   &  6 - \z &     6    &   4   &     6    \\
        \mu_3 \ \mu_2 \ \mu_3         & 2 &   2   &    4    &     4    &   4   &     4    \\
    \mu_2 \ \mu_1 \ \mu_1 \ \mu_2     & 2 &   6   &    8    &     8    &   6   &     8    \\
\hline
\multicolumn{2}{c|}{}                     &  12   & 18 - \z &    18    &  14   &    18    \\
\cline{3-7}
\end{array}
\quad
\begin{array}{|*3{>{\ss}c|}}
\hline
 & & \multicolumn{1}{|>{\ss}c|}{c(u_\beta)} \\
\cline{3-3}
 \ss{\beta\mathrm{-strings}}  & m & p \geq 3 \\
\hline
    \mu_3 \ \mu_1 \ \mu_3     & 2 &     4    \\
            \mu_3             & 2 &          \\
    \mu_2 \ \mu_0 \ \mu_2     & 1 &     2    \\
        \mu_2 \ \mu_2         & 2 &     2    \\
        \mu_1 \ \mu_1         & 2 &     4    \\
\hline
\multicolumn{2}{c|}{}             &    12    \\
\cline{3-3}
\end{array}
$$
We have $M = 12$ and $M_2 = 8$. Thus $\codim V_\kappa(s) > M$ unless $r = 2$, in which case $\codim V_\kappa(s) > M_r$; so the triple $(G, \lambda, p)$ satisfies $\ssdagcon$. Moreover $\codim C_V(u_\alpha) > M$, and $\codim C_V(u_\beta) \geq M$ --- by Lemma~\ref{lem: unip closure containment}, for any unipotent class $u^G$ we have $\codim C_V(u) \geq M$, and the only unipotent class $u^G$ with $\dim u^G \geq M$ is the regular unipotent class, whose closure contains $u_\alpha$ by Lemma~\ref{lem: any class in closure of reg class}; so the triple $(G, \lambda, p)$ satisfies $\udiamcon$.
\end{proof}

\begin{prop}\label{prop: G_2, 2omega_2, strings}
Let $G = G_2$ and $\lambda = 2\omega_2$ with $p = 3$; then the triple $(G, \lambda, p)$ satisfies $\ssdagcon$ and $\udiamcon$.
\end{prop}

\begin{proof}
The tables are as follows.
$$
\begin{array}{|*4{>{\ss}c|}}
\hline
i & \mu & |W.\mu| & m_\mu \\
\hline
3 & 2\omega_2 &  6 & 1 \\
2 & 3\omega_1 &  6 & 1 \\
1 &  \omega_2 &  6 & 2 \\
0 &     0     &  1 & 3 \\
\hline
\end{array}
\quad
\begin{array}{|*3{>{\ss}c|}}
\hline
 & & \multicolumn{1}{|>{\ss}c|}{c(u_\alpha)} \\
\cline{3-3}
             \ss{\alpha\mathrm{-strings}}             & m & p = 3 \\
\hline
\mu_3 \ \phmu \ \phmu \ \mu_1 \ \phmu \ \phmu \ \mu_3 & 2 &   4   \\
                        \mu_3                         & 2 &       \\
\mu_2 \ \phmu \ \phmu \ \mu_0 \ \phmu \ \phmu \ \mu_2 & 1 &   2   \\
            \mu_2 \ \phmu \ \phmu \ \mu_2             & 2 &   2   \\
            \mu_1 \ \phmu \ \phmu \ \mu_1             & 2 &   4   \\
\hline
\multicolumn{2}{c|}{}                                     &  12   \\
\cline{3-3}
\end{array}
\quad
\begin{array}{|*5{>{\ss}c|}}
\hline
 & & \multicolumn{2}{|>{\ss}c|}{c(s)} & \multicolumn{1}{|>{\ss}c|}{c(u_\beta)} \\
\cline{3-5}
     \ss{\beta\mathrm{-strings}}      & m & r = 2 & r \geq 5 & p = 3 \\
\hline
\mu_3 \ \mu_1 \ \mu_0 \ \mu_1 \ \mu_3 & 1 &   4   &     6    &   4   \\
        \mu_3 \ \mu_2 \ \mu_3         & 2 &   2   &     4    &   4   \\
    \mu_2 \ \mu_1 \ \mu_1 \ \mu_2     & 2 &   6   &     8    &   6   \\
\hline
\multicolumn{2}{c|}{}                     &  12   &    18    &  14   \\
\cline{3-5}
\end{array}
$$
We have $M = 12$ and $M_2 = 8$. Thus $\codim V_\kappa(s) > M$ unless $r = 2$, in which case $\codim V_\kappa(s) > M_r$; so the triple $(G, \lambda, p)$ satisfies $\ssdagcon$. Moreover $\codim C_V(u_\beta) > M$, and $\codim C_V(u_\alpha) \geq M$ --- by Lemma~\ref{lem: unip closure containment}, for any unipotent class $u^G$ we have $\codim C_V(u) \geq M$, and the only unipotent class $u^G$ with $\dim u^G \geq M$ is the regular unipotent class, whose closure contains $u_\beta$ by Lemma~\ref{lem: any class in closure of reg class}; so the triple $(G, \lambda, p)$ satisfies $\udiamcon$.
\end{proof}

\begin{prop}\label{prop: G_2, omega_1 + omega_2, strings}
Let $G = G_2$ and $\lambda = \omega_1 + \omega_2$ with $p = 3$; then the triple $(G, \lambda, p)$ satisfies $\ssddagcon$ and $\uddagcon$.
\end{prop}

\begin{proof}
The tables are as follows.
$$
\begin{array}{|*4{>{\ss}c|}}
\hline
i & \mu & |W.\mu| & m_\mu \\
\hline
4 & \omega_1 + \omega_2 & 12 & 1 \\
3 &      2\omega_1      &  6 & 2 \\
2 &       \omega_2      &  6 & 1 \\
1 &       \omega_1      &  6 & 3 \\
0 &          0          &  1 & 1 \\
\hline
\end{array}
\quad
\begin{array}{|*5{>{\ss}c|}}
\hline
 & & \multicolumn{2}{|>{\ss}c|}{c(s)} & \multicolumn{1}{|>{\ss}c|}{c(u_\alpha)} \\
\cline{3-5}
         \ss{\alpha\mathrm{-strings}}         & m & r = 2 & r \geq 5 & p = 3 \\
\hline
\mu_4 \ \mu_2 \ \mu_1 \ \mu_1 \ \mu_2 \ \mu_4 & 2 &  10   &    14    &  12   \\
    \mu_4 \ \mu_3 \ \mu_2 \ \mu_3 \ \mu_4     & 2 &   6   &    10    &   8   \\
                \mu_4 \ \mu_4                 & 2 &   2   &     2    &   2   \\
    \mu_3 \ \mu_1 \ \mu_0 \ \mu_1 \ \mu_3     & 1 &   5   &     8    &   6   \\
\hline
\multicolumn{2}{c|}{}                             &  23   &    34    &  28   \\
\cline{3-5}
\end{array}
\quad
\begin{array}{|*3{>{\ss}c|}}
\hline
 & & \multicolumn{1}{|>{\ss}c|}{c(u_\beta)} \\
\cline{3-3}
 \ss{\beta\mathrm{-strings}}  & m & p = 3 \\
\hline
\mu_4 \ \mu_1 \ \mu_1 \ \mu_4 & 2 &   8   \\
    \mu_4 \ \mu_3 \ \mu_4     & 2 &   4   \\
        \mu_4 \ \mu_4         & 2 &   2   \\
    \mu_3 \ \mu_1 \ \mu_3     & 2 &   8   \\
    \mu_2 \ \mu_0 \ \mu_2     & 1 &   2   \\
        \mu_2 \ \mu_2         & 2 &   2   \\
\hline
\multicolumn{2}{c|}{}             &  26   \\
\cline{3-3}
\end{array}
$$
We have $M = 12$. Thus $\codim V_\kappa(s), \ \codim C_V(u_\alpha), \ \codim C_V(u_\beta) > M$; so the triple $(G, \lambda, p)$ satisfies $\ssddagcon$ and $\uddagcon$.
\end{proof}

Combining the results in this section with Table~\ref{table: unexcluded triples} we have the following.

\begin{prop}\label{prop: remaining p-restricted triples}
Any $p$-restricted large triple which does not appear in Table~\ref{table: remaining triples} satisfies both $\ssdiamevcon$ and $\udiamcon$, and thus has TGS.
\end{prop}

\begin{table}[ht]
\caption{Remaining $p$-restricted large triples}\label{table: remaining triples}
\tabcapsp
$$
\begin{array}{|*4{c|}}
\hline
G      & \lambda                & \ell         & p          \tbs \\
\hline
A_\ell & 3\omega_1              & {} \geq 1    & {} \geq 5  \tbs \\
       & 4\omega_1              & 1            & {} \geq 5  \tbs \\
       & 2\omega_2              & 3            & {} \geq 3  \tbs \\
       & \omega_3               & {} \geq 8    & \hbox{any} \tbs \\
       & \omega_4               & 7, \dots, 11 & \hbox{any} \tbs \\
       & \omega_5               & 9            & \hbox{any} \tbs \\
       & \omega_1 + \omega_2    & {} \geq 3    & \hbox{any} \tbs \\
       & \omega_2 + \omega_\ell & 4, 5         & \hbox{any} \tbs \\
\hline
B_\ell & 2\omega_1              & {} \geq 2    & {} \geq 3  \tbs \\
       & \omega_3               & 4, 5, 6      & 2          \tbs \\
       & \omega_\ell            & 7, 8, 9      & \hbox{any} \tbs \\
       & \omega_1 + \omega_2    & 2            & 5          \tbs \\
\hline
C_\ell & \omega_3               & 4, 5, 6      & \hbox{any} \tbs \\
       & \omega_4               & 4            & {} \geq 3  \tbs \\
       & \omega_\ell            & 7, 8, 9      & 2          \tbs \\
\hline
D_\ell & 2\omega_1              & {} \geq 4    & {} \geq 3  \tbs \\
       & \omega_3               & 5            & 2          \tbs \\
       & \omega_\ell            & 8, 9, 10     & \hbox{any} \tbs \\
       & \omega_1 + \omega_4    & 4            & \hbox{any} \tbs \\
\hline
\end{array}
$$
\end{table}

In the following section we shall treat each of the triples which are listed in Table~\ref{table: remaining triples} but not in Table~\ref{table: large triple and first quadruple non-TGS}.

\section{Further analysis}\label{sect: large triple further analysis}

In this section we shall show that each of the $p$-restricted large triples listed in Table~\ref{table: remaining triples} which does not appear in Table~\ref{table: large triple and first quadruple non-TGS} satisfies $\ssdiamevcon$ and $\udiamcon$, and thus has TGS. Our approach is to generalize the strategy employed in Section~\ref{sect: large triple weight string analysis}, since consideration of weight strings alone will be insufficient for our purposes.

We continue with much of the notation of Section~\ref{sect: large triple weight string analysis}. Given a triple $(G, \lambda, p)$ listed in Table~\ref{table: remaining triples}, we let $s$ be an element of $G_{(r)}$ for some $r \in \P'$, and $\kappa$ be an element of $K^*$; we may assume $s$ lies in $T$. Write $\Phi(s) = \{ \alpha \in \Phi : \alpha(s) = 1 \}$, so that $C_G(s)^\circ = \langle T, X_\alpha : \alpha \in \Phi(s) \rangle$; observe that $\dim s^G = |\Phi| - |\Phi(s)| = M - |\Phi(s)|$. In addition, given a subsystem $\Psi$ of $\Phi$, we write $G_\Psi = \langle X_\alpha : \alpha \in \Psi \rangle$ for the corresponding subsystem subgroup of $G$.

For our generalization, let $\Psi$ be a standard subsystem of $\Phi$. We define an equivalence relation on $\Lambda(V)$ by saying that two weights are related if and only if their difference is a sum of roots in $\Psi$; we call the equivalence classes {\em $\Psi$-nets\/}. Moreover, we write $\Psi$ as a disjoint union of irreducible subsystems $\Psi_i$, any two of which are orthogonal to each other; then each $\Psi_i$ is standard, and $G_\Psi$ is the product of the subsystem subgroups $G_{\Psi_i}$. For each $i$ let $u_{\Psi_i}$ be a regular unipotent element of $G_{\Psi_i}$; let $u_\Psi$ be the product of the $u_{\Psi_i}$, so that $u_\Psi$ is regular unipotent in $G_\Psi$. Observe that if $\Psi = \langle \alpha \rangle$, then $\Psi$-nets are simply $\alpha$-strings, and we may take $u_\Psi = u_\alpha$.

Now on the one hand, if we assume $\Psi$ is disjoint from $\Phi(s)$, then in a given $\Psi$-net any two weights whose difference is a root must lie in different eigenspaces for $s$; we may use this observation to obtain a lower bound $c(s)$ for the contribution to $\codim V_\kappa(s)$. On the other hand, for the same $\Psi$-net the sum of the weight spaces corresponding to the set of weights therein is a $G_\Psi$-module; if we assume $u_\Psi \in G_{(p)}$, we may determine a lower bound $c(u_\Psi)$ for the contribution to $\codim C_V(u_\Psi)$.

As with the $\alpha$-string tables in Section~\ref{sect: large triple weight string analysis}, we provide a $\Psi$-net table, whose rows correspond to the different types of $\Psi$-net which appear among the weights in $\Lambda(V)$. In each row of this table the entries are as follows: the first column gives the type of $\Psi$-net, using notation explained below; the next few columns give the numbers $n_i$ of weights in the $\Psi$-net which lie in the $W$-orbit numbered $i$ in the weight table; the next column gives the number $m$ of such $\Psi$-nets; and the remaining columns give the lower bounds $c(s)$ and $c(u_\Psi)$ (which may depend on $r$ or $p$ respectively). The bottom row of the table sums the values $c(s)$ and $c(u_\Psi)$ to give lower bounds $c(\Psi)_{ss}$ and $c(\Psi)_u$ for $\codim V_\kappa(s)$ and $\codim C_V(u_\Psi)$ respectively.

Our strategy is then as follows. We first give the weight table, as in Section~\ref{sect: large triple weight string analysis}. We then begin by giving the $\Psi$-net table for $\Psi$ of rank $1$ (if $e(\Phi) > 1$ we need to consider $\Psi = \langle \alpha \rangle$ both for $\alpha$ short and for $\alpha$ long). If $s^G$ is any semisimple class with $\dim s^G < c(\Psi)_{ss}$ such that $\Psi$ is disjoint from $\Phi(s)$ (which we may assume if $\Psi = \langle \alpha \rangle$, unless $p \neq e(\Phi) > 1$ and $\alpha$ is long), then for all $\kappa \in K^*$ we have
$$
\codim V_\kappa(s) \geq c(\Psi)_{ss} > \dim s^G
$$
as required for $\ssdiamevcon$; thus we may henceforth ignore all such semisimple classes. Likewise if $u^G$ is any unipotent class in $G_{(p)}$ with $\dim u^G < c(\Psi)_u$ such that the closure of $u^G$ contains ${u_\Psi}^G$ (which by Lemma~\ref{lem: root elt class in closure of any non-triv class} we may assume if $\Psi$ has rank $1$, provided we treat both long and short roots if $e(\Phi) > 1$), then by Lemma~\ref{lem: unip closure containment} we have
$$
\codim C_V(u) \geq \codim C_V(u_\Psi) \geq c(\Psi)_u > \dim u^G
$$
as required for $\udiamcon$; thus we may henceforth ignore all such unipotent classes. The remaining semisimple classes $s^G$ satisfy $\dim s^G \geq c(\Psi)_{ss}$, i.e., $|\Phi(s)| \leq M - c(\Psi)_{ss}$; we identify a larger subsystem $\Psi$ such that each of these remaining $\Phi(s)$ has a conjugate of $\Psi$ disjoint from it, and such that all remaining unipotent classes in $G_{(p)}$ have ${u_\Psi}^G$ in their closure. We take this $\Psi$ and repeat the procedure to obtain improved lower bounds $c(\Psi)_{ss}$ and $c(\Psi)_u$. Continuing thus, in most cases we eventually show that the triple $(G, \lambda, p)$ satisfies $\ssdiamevcon$ and $\udiamcon$.

In order to explain our notation for $\Psi$-nets, we begin with Weyl $G_\Psi$-modules. We have $\Psi = \langle \alpha_j : j \in S \rangle$ for some subset $S$ of $\{ 1, \dots, \ell \}$. For each $j \in S$, there exists $i$ with $\alpha_j \in \Psi_i$; we let $\bom_j$ be the fundamental dominant weight of $G_{\Psi_i}$ corresponding to $\alpha_j$. The highest weight of any Weyl $G_\Psi$-module may then be written as $\bar\nu = \sum_{j \in S} a_j \bom_j$ for some non-negative integers $a_j$; we write $W_{G_\Psi}(\bar\nu)$ for the Weyl $G_\Psi$-module with highest weight $\bar\nu$.

Note that, given $\alpha \in \Psi$, any $\Psi$-net is a union of $\alpha$-strings; as we saw in Section~\ref{sect: large triple weight string analysis}, in any $\alpha$-string the outermost weights lie in one $W$-orbit and any internal weights lie in \lq lower' $W$-orbits. Thus for a given $\Psi$-net, if $i$ is maximal such that it contains weights of the form $\mu_i$, then for each $\alpha \in \Psi$ each weight of the form $\mu_i$ lying therein must be outermost in the $\alpha$-string containing it; we may regard these weights as forming the \lq outer layer' of the $\Psi$-net concerned.

We find that in the cases treated here, the set of weights (ignoring multiplicities) appearing in a given $\Psi$-net is that of a single Weyl $G_\Psi$-module $W_{G_\Psi}(\bar\nu)$. The corresponding entry in the first column of the $\Psi$-net table is the weight $\bar\nu$.

We also find that in the cases treated here, each $\Psi$-net contains weights from either one or two $W$-orbits. Given $\alpha \in \Psi$, any $\alpha$-string of length $3$ or $4$ must be of the form $\mu_i \ \mu_j \ \mu_i$ or $\mu_i \ \mu_j \ \mu_j \ \mu_i$ with $j < i$. Thus for any row of the $\Psi$-net table, the entry $\bar\nu$ in the first column and the values $n_i$ in the next few suffice to determine the exact form of the $\Psi$-net concerned.

We give an example to show both how the entries in the $\Psi$-net table may be calculated, and also how they determine the forms of the $\Psi$-nets. Let $G = A_\ell$ for $\ell \in [4, \infty)$ and $\lambda = \omega_1 + \omega_2$. By Lemma~\ref{lem: multiplicities in 3omega_1 and omega_1 + omega_2} the weight table is as follows, where we write $\z = \z_{p, 3}$.
$$
\begin{array}{|*4{>{\ss}c|}}
\hline
i & \mu & |W.\mu| & m_\mu \\
\hline
2 & \omega_1 + \omega_2 &        \ell(\ell + 1)       &    1   \\
1 &       \omega_3      & \frac{1}{6}\ell(\ell^2 - 1) & 2 - \z \\
\hline
\end{array}
$$
Recall from Section~\ref{sect: notation} that the root system lies in an $(\ell + 1)$-dimensional Euclidean space with orthonormal basis $\ve_1, \dots, \ve_{\ell + 1}$; the simple roots are $\ve_1 - \ve_2, \dots, \ve_\ell - \ve_{\ell + 1}$, and the Weyl group acts by permuting the vectors $\ve_i$. From \cite[13.1, Table~1]{HumLie} we see that in this notation we have $\omega_1 + \omega_2 = 2\ve_1 + \ve_2 - \xi$ and $\omega_3 = \ve_1 + \ve_2 + \ve_3 - \xi$, where $\xi = \frac{3}{\ell + 1}\sum_{i = 1}^{\ell + 1}\ve_i$. We shall represent $a_1\ve_1 + \cdots + a_{\ell + 1}\ve_{\ell + 1} - \xi$ as $a_1 \dots a_{\ell + 1}$. Thus the weights in $\Lambda(V)$ are obtained from $210 \dots 0$, $1110 \dots 0$ by permuting symbols.

Let $\Psi = \langle \alpha_1, \alpha_2 \rangle$ of type $A_2$; write $\mu = a_1 \dots a_{\ell + 1}$, then $\langle \mu, \alpha_1 \rangle = a_1 - a_2$ and $\langle \mu, \alpha_2 \rangle = a_2 - a_3$. Any $\Psi$-net is a union of $\langle w_{\alpha_1}, w_{\alpha_2} \rangle$-orbits, each of which contains a single weight $\mu$ with $\langle \mu, \alpha_1 \rangle \geq 0$ and $\langle \mu, \alpha_2 \rangle \geq 0$, i.e., with $a_1 \geq a_2 \geq a_3$. Since all weights $\mu$ in a given $\Psi$-net have the same final segment $a_4 \dots a_{\ell + 1}$, we shall concentrate attention on the initial segment $a_1a_2a_3$, bearing in mind that each way of completing $a_1a_2a_3$ to $a_1a_2a_3a_4 \dots a_{\ell + 1}$ gives a distinct $\Psi$-net, and that all permutations of the initial segment give weights in the same $\Psi$-net. We start by determining the $\Psi$-nets with outer layer comprising weights in the $W$-orbit containing $\lambda$ itself; we then move to the next $W$-orbit, treating only the remaining weights, and continue until all weights have been dealt with.

Initially then we consider the weights $\mu$ of the form $\mu_2$; for these we have
$$
(\langle \mu, \alpha_1 \rangle, \langle \mu, \alpha_2 \rangle) = \begin{cases}
(1, 1) & \hbox{if } \mu = 210\dots, \\
(2, 0) & \hbox{if } \mu = 200\dots, \\
(1, 0) & \hbox{if } \mu = 100\dots, \\
(0, 0) & \hbox{if } \mu = 000\dots.
\end{cases}
$$
In the first possibility there are six weights in the outer layer, and a single internal weight $111\dots$; the final segment is simply $0 \dots 0$, giving one $\Psi$-net with $\bar\nu = \bom_1 + \bom_2$ and $(n_1, n_2) = (1, 6)$. In the second there are three weights in the outer layer, and three internal weights $110\dots$; the final segment is a permutation of $10 \dots 0$, giving $\ell - 2$ $\Psi$-nets with $\bar\nu = 2\bom_1$ and $(n_1, n_2) = (3, 3)$. In the third there are just three weights; the final segment is a permutation of $20 \dots 0$, giving $\ell - 2$ $\Psi$-nets with $\bar\nu = \bom_1$ and $(n_1, n_2) = (0, 3)$. In the fourth there is a single weight; the final segment is a permutation of $210 \dots 0$, giving $(\ell - 2)(\ell - 3)$ $\Psi$-nets with $\bar\nu = 0$ and $(n_1, n_2) = (0, 1)$. Next we consider the remaining weights $\mu$ of the form $\mu_1$; for these we have
$$
(\langle \mu, \alpha_1 \rangle, \langle \mu, \alpha_2 \rangle) = \begin{cases}
(1, 0) & \hbox{if } \mu = 100\dots, \\
(0, 0) & \hbox{if } \mu = 000\dots.
\end{cases}
$$
In the first possibility there are three weights; the final segment is a permutation of $110 \dots 0$, giving $\frac{1}{2}(\ell - 2)(\ell - 3)$ $\Psi$-nets with $\bar\nu = \bom_1$ and $(n_1, n_2) = (3, 0)$. In the second there is a single weight; the final segment is a permutation of $1110 \dots 0$, giving $\frac{1}{6}(\ell - 2)(\ell - 3)(\ell - 4)$ $\Psi$-nets with $\bar\nu = 0$ and $(n_1, n_2) = (1, 0)$. This completes the determination of the $\Psi$-nets.

Consider now the first two types of $\Psi$-net described; as we have seen, the entries in the first three columns of the $\Psi$-net table are
$$
\phantom{2\bom_1} \bom_1 + \bom_2 \quad 1 \quad 6 \qquad\qquad \hbox{or} \qquad\qquad 2\bom_1 \quad 3 \quad 3 \phantom{\bom_1 + \bom_2}
$$
and it follows that they correspond to $\Psi$-nets of the form
$$
\begin{array}{c}
       \mu_2 \ \ \mu_2        \\
  \mu_2 \ \ \mu_1 \ \ \mu_2   \\
       \mu_2 \ \ \mu_2
\end{array}
\qquad\qquad \hbox{or} \qquad\qquad
\begin{array}{c}
            \mu_2             \\
       \mu_1 \ \ \mu_1        \\
  \mu_2 \ \ \mu_1 \ \ \mu_2
\end{array}
$$
respectively, in which the lines sloping upwards and to the right are $\alpha_1$-strings and those sloping downwards and to the right are $\alpha_2$-strings. Thus the weights are
$$
\mu, \mu - \alpha_1, \mu - \alpha_2, \mu - \alpha_1 - \alpha_2, \mu - 2\alpha_1 - \alpha_2, \mu - \alpha_1 - 2\alpha_2, \mu - 2\alpha_1 - 2\alpha_2
$$
or
$$
\mu, \mu - \alpha_1, \mu - 2\alpha_1, \mu - \alpha_1 - \alpha_2, \mu - 2\alpha_1 - \alpha_2, \mu - 2\alpha_1 - 2\alpha_2
$$
respectively.

It will next be helpful to say something about how the values $c(s)$ and $c(u_\Psi)$ are obtained for a given $\Psi$-net; we begin with the former. We may take any given weight $\mu$ in the $\Psi$-net and suppose $\mu$ lies in $V_\kappa(s)$; since $r$ is the order of $\bar s = sZ(G)$, this implies that any other weight in the $\Psi$-net of the form $\mu - t \alpha$, where $\alpha \in \Psi$ and $t$ is not a multiple of $r$, does not lie in $V_\kappa(s)$. If there are weights not of this form, we may take any of them and repeat the process until we obtain a maximal set of weights which could all lie in $V_\kappa(s)$. After working through all possibilities we may let $c(s)$ be the smallest possible contribution to $\codim V_\kappa(s)$ obtained. Of course, we may use the action of the Weyl group to reduce the number of possibilities to be considered; for example, when choosing the initial weight $\mu$ we need only treat one from each $W$-orbit represented in the $\Psi$-net. Note that if $\Psi_i$ is one of the irreducible subsystems comprising $\Psi$, then we may decompose the $\Psi$-net into $\Psi_i$-nets and add together the lower bounds obtained from them to obtain a lower bound for the $\Psi$-net.

For example, consider the two $\Psi$-nets above for $G_\Psi$ of type $A_2$. In the first, where $\bar \nu = \bom_1 + \bom_2$, up to the action of $W$ the maximal sets of weights which could all lie in $V_\kappa(s)$ are as follows:
$$
\{ \mu, \mu - 2\alpha_1 - \alpha_2, \mu - \alpha_1 - 2\alpha_2 \}; \{ \mu - \alpha_1 - \alpha_2 \}; \hbox{ and } \{ \mu, \mu - 2\alpha_1 - 2\alpha_2 \} \hbox{ if } r = 2.
$$
In the second, where $\bar \nu = 2\bom_1$, they are as follows:
$$
\{ \mu, \mu - 2\alpha_1 - \alpha_2 \}; \hbox{ and } \{ \mu, \mu - 2\alpha_1, \mu - 2\alpha_1 - 2\alpha_2 \} \hbox{ if } r = 2.
$$
Using the multiplicities given in the weight table, for the first $\Psi$-net we have $c(s) = \min \{ 3 + (2 - \z), 6 \} = 5 - \z$ (even if $r = 2$), while for the second we have $c(s) = 2 + 2(2 - \z) = 6 - 2\z$, unless $r = 2$ in which case $c(s) = \min \{ 2 + 2(2 - \z), 3(2 - \z) \} = 6 - 3\z$.

We now turn to the value $c(u_\Psi)$ for a given $\Psi$-net, where we assume $u_\Psi \in G_{(p)}$. For each $i$, let $J_{\Psi_i}$ be an $A_1$ subgroup of $G_{\Psi_i}$ containing the regular unipotent element $u_{\Psi_i}$. Let $J_\Psi$ be the product of the subgroups $J_{\Psi_i}$; then $u_\Psi$ is regular in $J_\Psi$. As observed above, the sum of the weight spaces corresponding to the weights in the $\Psi$-net is a $G_\Psi$-module; we may decompose this into $J_\Psi$-composition factors. Any such $J_\Psi$-composition factor $X$ is then a tensor product of $J_{\Psi_i}$-composition factors $X_i$. For each $i$ we may compute $\codim C_{X_i}(u_{\Psi_i})$; multiplying by $\prod_{j \neq i} \dim X_j$ gives $\codim C_X(u_{\Psi_i})$. Since the closure of ${u_\Psi}^G$ contains each $u_{\Psi_i}$, we may then take the maximum of the values $\codim C_X(u_{\Psi_i})$ as a lower bound for the contribution to $\codim C_V(u_\Psi)$ from the $J_\Psi$-composition factor $X$; summing over the $J_\Psi$-composition factors in the $G_\Psi$-module gives the value $c(u_\Psi)$ for the $\Psi$-net.

We next discuss further the process of successively increasing the size of the subsystem $\Psi$ used in the calculations. At a given stage, using $\Psi$ we obtain a lower bound $c(\Psi)_{ss}$ for $\codim V_\kappa(s)$, and hence an upper bound $M - c(\Psi)_{ss}$ for the size of the subsystems $\Phi(s)$ for the semisimple classes $s^G$ which still require consideration. We then wish to take a certain larger subsystem $\Psi$ at the next stage, which requires us to know that each of these remaining $\Phi(s)$ has a conjugate of $\Psi$ disjoint from it. If the upper bound for $|\Phi(s)|$ is relatively small, this may be seen by inspection. For the other instances which arise, we recall that for a subsystem $\Psi$ we write $m_\Psi$ for the size of the smallest possible subsystem which intersects every conjugate of $\Psi$, and use the values $m_\Psi$ given in Lemma~\ref{lem: m_Psi values}; provided the upper bound for $|\Phi(s)|$ is smaller than $m_\Psi$, we may conclude that for each of the semisimple classes $s^G$ still under consideration there is indeed a conjugate of $\Psi$ disjoint from $\Phi(s)$.

We now work through the triples. As in Section~\ref{sect: large triple weight string analysis}, we shall begin with those in which all roots in $\Phi$ have the same length.

\begin{prop}\label{prop: A_ell, 3omega_1, nets}
Let $G = A_\ell$ for $\ell \in [3, \infty)$ and $\lambda = 3\omega_1$ with $p \geq 5$; then the triple $(G, \lambda, p)$ satisfies $\ssdiamevcon$ and $\udiamcon$.
\end{prop}

\begin{proof}
By Lemma~\ref{lem: multiplicities in 3omega_1 and omega_1 + omega_2} the weight table is as follows.
$$
\begin{array}{|*4{>{\ss}c|}}
\hline
i & \mu & |W.\mu| & m_\mu \\
\hline
3 &      3\omega_1      &           \ell + 1          & 1 \\
2 & \omega_1 + \omega_2 &        \ell(\ell + 1)       & 1 \\
1 &       \omega_3      & \frac{1}{6}\ell(\ell^2 - 1) & 1 \\
\hline
\end{array}
$$
We have $M = \ell(\ell + 1)$, $M_3 = 2\lfloor \frac{1}{3}(\ell + 1)^2 \rfloor$ and $M_2 = \lfloor \frac{1}{2}(\ell + 1)^2 \rfloor$.

Take $\Psi = \langle \alpha_1 \rangle$ of type $A_1$. The $\Psi$-net table is as follows.
$$
\begin{array}{|*9{>{\ss}c|}}
\hline
\multicolumn{4}{|>{\ss}c|}{\Psi\mathrm{-nets}} & & \multicolumn{3}{|>{\ss}c|}{c(s)} & \multicolumn{1}{|>{\ss}c|}{c(u_\Psi)} \\
\cline{1-4} \cline{6-9}
    \bar\nu    & n_1 & n_2 & n_3 &                     m                     &              r = 2              &              r = 3              &             r \geq 5            &             p \geq 5            \\
\hline
    3\bom_1    &  0  &  2  &  2  &                     1                     &                2                &                2                &                3                &                3                \\
    2\bom_1    &  1  &  2  &  0  &                  \ell - 1                 &             \ell - 1            &           2(\ell - 1)           &           2(\ell - 1)           &           2(\ell - 1)           \\
     \bom_1    &  0  &  2  &  0  &                  \ell - 1                 &             \ell - 1            &             \ell - 1            &             \ell - 1            &             \ell - 1            \\
     \bom_1    &  2  &  0  &  0  &      \frac{1}{2}(\ell - 1)(\ell - 2)      & \frac{1}{2}(\ell - 1)(\ell - 2) & \frac{1}{2}(\ell - 1)(\ell - 2) & \frac{1}{2}(\ell - 1)(\ell - 2) & \frac{1}{2}(\ell - 1)(\ell - 2) \\
       0       &  0  &  0  &  1  &                  \ell - 1                 &                                 &                                 &                                 &                                 \\
       0       &  0  &  1  &  0  &            (\ell - 1)(\ell - 2)           &                                 &                                 &                                 &                                 \\
       0       &  1  &  0  &  0  & \frac{1}{6}(\ell - 1)(\ell - 2)(\ell - 3) &                                 &                                 &                                 &                                 \\
\hline
\multicolumn{5}{c|}{}                                                        &  \frac{1}{2}(\ell^2 + \ell + 2) &     \frac{1}{2}\ell(\ell + 3)   & \frac{1}{2}(\ell^2 + 3\ell + 2) & \frac{1}{2}(\ell^2 + 3\ell + 2) \\
\cline{6-9}
\end{array}
$$
Thus $\codim V_\kappa(s), \codim C_V(u_\Psi) > 2\ell = \dim {u_\Psi}^G$. We therefore need only consider semisimple classes $s^G$ with $|\Phi(s)| < M - 2\ell = m_{{A_1}^2}$, each of which has a subsystem of type ${A_1}^2$ disjoint from $\Phi(s)$, and unipotent classes of dimension greater than $2\ell$, each of which has the class ${A_1}^2$ in its closure by Lemma~\ref{lem: various classes in classical groups by dim}(i).

Now take $\Psi = \langle \alpha_1, \alpha_3 \rangle$ of type ${A_1}^2$. The $\Psi$-net table is as follows.
$$
\begin{array}{|*9{>{\ss}c|}}
\hline
\multicolumn{4}{|>{\ss}c|}{\Psi\mathrm{-nets}} & & \multicolumn{3}{|>{\ss}c|}{c(s)} & \multicolumn{1}{|>{\ss}c|}{c(u_\Psi)} \\
\cline{1-4} \cline{6-9}
      \bar\nu      & n_1 & n_2 & n_3 &                     m                     &              r = 2              &              r = 3              &             r \geq 5            &             p \geq 5            \\
\hline
      3\bom_1      &  0  &  2  &  2  &                     1                     &                2                &                2                &                3                &                3                \\
      3\bom_3      &  0  &  2  &  2  &                     1                     &                2                &                2                &                3                &                3                \\
 2\bom_1 + \bom_3  &  2  &  4  &  0  &                     1                     &                3                &                4                &                4                &                4                \\
  \bom_1 + 2\bom_3 &  2  &  4  &  0  &                     1                     &                3                &                4                &                4                &                4                \\
      2\bom_1      &  1  &  2  &  0  &                  \ell - 3                 &             \ell - 3            &           2(\ell - 3)           &           2(\ell - 3)           &           2(\ell - 3)           \\
      2\bom_3      &  1  &  2  &  0  &                  \ell - 3                 &             \ell - 3            &           2(\ell - 3)           &           2(\ell - 3)           &           2(\ell - 3)           \\
  \bom_1 + \bom_3  &  4  &  0  &  0  &                  \ell - 3                 &           2(\ell - 3)           &           2(\ell - 3)           &           2(\ell - 3)           &           2(\ell - 3)           \\
       \bom_1      &  0  &  2  &  0  &                  \ell - 3                 &             \ell - 3            &             \ell - 3            &             \ell - 3            &             \ell - 3            \\
       \bom_3      &  0  &  2  &  0  &                  \ell - 3                 &             \ell - 3            &             \ell - 3            &             \ell - 3            &             \ell - 3            \\
       \bom_1      &  2  &  0  &  0  &      \frac{1}{2}(\ell - 3)(\ell - 4)      & \frac{1}{2}(\ell - 3)(\ell - 4) & \frac{1}{2}(\ell - 3)(\ell - 4) & \frac{1}{2}(\ell - 3)(\ell - 4) & \frac{1}{2}(\ell - 3)(\ell - 4) \\
       \bom_3      &  2  &  0  &  0  &      \frac{1}{2}(\ell - 3)(\ell - 4)      & \frac{1}{2}(\ell - 3)(\ell - 4) & \frac{1}{2}(\ell - 3)(\ell - 4) & \frac{1}{2}(\ell - 3)(\ell - 4) & \frac{1}{2}(\ell - 3)(\ell - 4) \\
         0         &  0  &  0  &  1  &                  \ell - 3                 &                                 &                                 &                                 &                                 \\
         0         &  0  &  1  &  0  &            (\ell - 3)(\ell - 4)           &                                 &                                 &                                 &                                 \\
         0         &  1  &  0  &  0  & \frac{1}{6}(\ell - 3)(\ell - 4)(\ell - 5) &                                 &                                 &                                 &                                 \\
\hline
\multicolumn{5}{c|}{}                                                            &        \ell^2 - \ell + 4        &          \ell(\ell + 1)         &        \ell^2 + \ell + 2        &        \ell^2 + \ell + 2        \\
\cline{6-9}
\end{array}
$$
Thus $\codim C_V(u_\Psi) > M$, and $\codim V_\kappa(s) > M$ unless $r = 2$ or $r = 3$, in which case $\codim V_\kappa(s) > M_r$; so the triple $(G, \lambda, p)$ satisfies $\ssdiamevcon$ and $\udiamcon$.
\end{proof}

\begin{prop}\label{prop: A_ell, omega_3, nets}
Let $G = A_\ell$ for $\ell \in [9, \infty)$ and $\lambda = \omega_3$; then the triple $(G, \lambda, p)$ satisfies $\ssdiamevcon$ and $\udiamcon$.
\end{prop}

\begin{proof}
The weight table is as follows.
$$
\begin{array}{|*4{>{\ss}c|}}
\hline
i & \mu & |W.\mu| & m_\mu \\
\hline
1 & \omega_3 & \frac{1}{6}\ell(\ell^2 - 1) & 1 \\
\hline
\end{array}
$$
We have $M = \ell(\ell + 1)$, $M_3 = 2\lfloor \frac{1}{3}(\ell + 1)^2 \rfloor$ and $M_2 = \lfloor \frac{1}{2}(\ell + 1)^2 \rfloor$.

Take $\Psi = \langle \alpha_1 \rangle$ of type $A_1$. The $\Psi$-net table is as follows.
$$
\begin{array}{|*5{>{\ss}c|}}
\hline
\multicolumn{2}{|>{\ss}c|}{\Psi\mathrm{-nets}} & & \multicolumn{1}{|>{\ss}c|}{c(s)} & \multicolumn{1}{|>{\ss}c|}{c(u_\Psi)} \\
\cline{1-2} \cline{4-5}
    \bar\nu    & n_1 &                      m                     &             r \geq 2            &             p \geq 2            \\
\hline
     \bom_1    &  2  &       \frac{1}{2}(\ell - 1)(\ell - 2)      & \frac{1}{2}(\ell - 1)(\ell - 2) & \frac{1}{2}(\ell - 1)(\ell - 2) \\
       0       &  1  & \frac{1}{6}(\ell - 1)(\ell^2 - 5\ell + 12) &                                 &                                 \\
\hline
\multicolumn{3}{c|}{}                                             & \frac{1}{2}(\ell - 1)(\ell - 2) & \frac{1}{2}(\ell - 1)(\ell - 2) \\
\cline{4-5}
\end{array}
$$
Thus $\codim V_\kappa(s), \codim C_V(u_\Psi) > 2\ell = \dim {u_\Psi}^G$. We therefore need only consider semisimple classes $s^G$ with $|\Phi(s)| < M - 2\ell = m_{{A_1}^2}$, each of which has a subsystem of type ${A_1}^2$ disjoint from $\Phi(s)$, and unipotent classes of dimension greater than $2\ell$, each of which has the class ${A_1}^2$ in its closure by Lemma~\ref{lem: various classes in classical groups by dim}(i).

Now take $\Psi = \langle \alpha_1, \alpha_3 \rangle$ of type ${A_1}^2$. The $\Psi$-net table is as follows.
$$
\begin{array}{|*5{>{\ss}c|}}
\hline
\multicolumn{2}{|>{\ss}c|}{\Psi\mathrm{-nets}} & & \multicolumn{1}{|>{\ss}c|}{c(s)} & \multicolumn{1}{|>{\ss}c|}{c(u_\Psi)} \\
\cline{1-2} \cline{4-5}
      \bar\nu      & n_1 &                      m                     &             r \geq 2             &             p \geq 2             \\
\hline
  \bom_1 + \bom_3  &  4  &                  \ell - 3                  &            2(\ell - 3)           &            2(\ell - 3)           \\
       \bom_1      &  2  &      \frac{1}{2}(\ell^2 - 7\ell + 14)      & \frac{1}{2}(\ell^2 - 7\ell + 14) & \frac{1}{2}(\ell^2 - 7\ell + 14) \\
       \bom_3      &  2  &      \frac{1}{2}(\ell^2 - 7\ell + 14)      & \frac{1}{2}(\ell^2 - 7\ell + 14) & \frac{1}{2}(\ell^2 - 7\ell + 14) \\
         0         &  1  & \frac{1}{6}(\ell - 3)(\ell^2 - 9\ell + 32) &                                  &                                  \\
\hline
\multicolumn{3}{c|}{}                                                 &        \ell^2 - 5\ell + 8        &        \ell^2 - 5\ell + 8        \\
\cline{4-5}
\end{array}
$$
Thus $\codim V_\kappa(s), \codim C_V(u_\Psi) > 4\ell - 2 > 4\ell - 4 = \dim {u_\Psi}^G$. We therefore need only consider semisimple classes $s^G$ with $|\Phi(s)| < M - (4\ell - 2) = m_{{A_1}^3}$, each of which has a subsystem of type ${A_1}^3$ disjoint from $\Phi(s)$, and unipotent classes of dimension greater than $4\ell - 2$, each of which has the class ${A_1}^3$ in its closure by Lemma~\ref{lem: various classes in classical groups by dim}(ii).

Now take $\Psi = \langle \alpha_1, \alpha_3, \alpha_5 \rangle$ of type ${A_1}^3$. The $\Psi$-net table is as follows.
$$
\begin{array}{|*5{>{\ss}c|}}
\hline
\multicolumn{2}{|>{\ss}c|}{\Psi\mathrm{-nets}} & & \multicolumn{1}{|>{\ss}c|}{c(s)} & \multicolumn{1}{|>{\ss}c|}{c(u_\Psi)} \\
\cline{1-2} \cline{4-5}
          \bar\nu          & n_1 &                      m                      &              r \geq 2              &              p \geq 2              \\
\hline
  \bom_1 + \bom_3 + \bom_5 &  8  &                      1                      &                  4                 &                  4                 \\
      \bom_1 + \bom_3      &  4  &                   \ell - 5                  &             2(\ell - 5)            &             2(\ell - 5)            \\
      \bom_1 + \bom_5      &  4  &                   \ell - 5                  &             2(\ell - 5)            &             2(\ell - 5)            \\
      \bom_3 + \bom_5      &  4  &                   \ell - 5                  &             2(\ell - 5)            &             2(\ell - 5)            \\
           \bom_1          &  2  &      \frac{1}{2}(\ell^2 - 11\ell + 34)      &  \frac{1}{2}(\ell^2 - 11\ell + 34) &  \frac{1}{2}(\ell^2 - 11\ell + 34) \\
           \bom_3          &  2  &      \frac{1}{2}(\ell^2 - 11\ell + 34)      &  \frac{1}{2}(\ell^2 - 11\ell + 34) &  \frac{1}{2}(\ell^2 - 11\ell + 34) \\
           \bom_5          &  2  &      \frac{1}{2}(\ell^2 - 11\ell + 34)      &  \frac{1}{2}(\ell^2 - 11\ell + 34) &  \frac{1}{2}(\ell^2 - 11\ell + 34) \\
             0             &  1  & \frac{1}{6}(\ell - 5)(\ell^2 - 13\ell + 60) &                                    &                                    \\
\hline
\multicolumn{3}{c|}{}                                                          & \frac{1}{2}(3\ell^2 - 21\ell + 50) & \frac{1}{2}(3\ell^2 - 21\ell + 50) \\
\cline{4-5}
\end{array}
$$
Thus $\codim V_\kappa(s), \codim C_V(u_\Psi) > M_2 > 6\ell - 12 = \dim {u_\Psi}^G$; we may therefore assume from now on that $r \geq 3$, and that $p \geq 3$ when we treat unipotent classes. Moreover if $\ell \in [21, \infty)$ then $\codim V_\kappa(s), \codim C_V(u_\Psi) > M$; so the triple $(G, \lambda, p)$ satisfies $\ssdiamevcon$ and $\udiamcon$. We may therefore assume from now on that $\ell \in [9, 20]$. We need only consider semisimple classes $s^G$ with $|\Phi(s)| < M - M_2 = m_{A_2{A_1}^2}$, each of which has a subsystem of type $A_2{A_1}^2$ disjoint from $\Phi(s)$, and unipotent classes of dimension greater than $M_2$, each of which has the class $A_2{A_1}^2$ in its closure by Lemma~\ref{lem: various classes in classical groups by dim}(v).

Now take $\Psi = \langle \alpha_1, \alpha_2, \alpha_4, \alpha_6 \rangle$ of type $A_2{A_1}^2$. The $\Psi$-net table is as follows.
$$
\begin{array}{|*5{>{\ss}c|}}
\hline
\multicolumn{2}{|>{\ss}c|}{\Psi\mathrm{-nets}} & & \multicolumn{1}{|>{\ss}c|}{c(s)} & \multicolumn{1}{|>{\ss}c|}{c(u_\Psi)} \\
\cline{1-2} \cline{4-5}
          \bar\nu          & n_1 &                        m                       &              r \geq 3             &              p \geq 3             \\
\hline
  \bom_1 + \bom_4 + \bom_6 & 12  &                        1                       &                 8                 &                 8                 \\
      \bom_1 + \bom_4      &  6  &                    \ell - 6                    &            4(\ell - 6)            &            4(\ell - 6)            \\
      \bom_1 + \bom_6      &  6  &                    \ell - 6                    &            4(\ell - 6)            &            4(\ell - 6)            \\
      \bom_2 + \bom_4      &  6  &                        1                       &                 4                 &                 4                 \\
      \bom_2 + \bom_6      &  6  &                        1                       &                 4                 &                 4                 \\
      \bom_4 + \bom_6      &  4  &                    \ell - 6                    &            2(\ell - 6)            &            2(\ell - 6)            \\
           \bom_1          &  3  &        \frac{1}{2}(\ell^2 - 13\ell + 46)       &        \ell^2 - 13\ell + 46       &        \ell^2 - 13\ell + 46       \\
           \bom_2          &  3  &                    \ell - 6                    &            2(\ell - 6)            &            2(\ell - 6)            \\
           \bom_4          &  2  &        \frac{1}{2}(\ell^2 - 13\ell + 44)       & \frac{1}{2}(\ell^2 - 13\ell + 44) & \frac{1}{2}(\ell^2 - 13\ell + 44) \\
           \bom_6          &  2  &        \frac{1}{2}(\ell^2 - 13\ell + 44)       & \frac{1}{2}(\ell^2 - 13\ell + 44) & \frac{1}{2}(\ell^2 - 13\ell + 44) \\
             0             &  1  & \frac{1}{6}(\ell^3 - 21\ell^2 + 158\ell - 402) &                                   &                                   \\
\hline
\multicolumn{3}{c|}{}                                                             &       2\ell^2 - 14\ell + 34       &       2\ell^2 - 14\ell + 34       \\
\cline{4-5}
\end{array}
$$
Thus $\codim V_\kappa(s), \codim C_V(u_\Psi) > M_3 > 8\ell - 18 = \dim {u_\Psi}^G$; we may therefore assume from now on that $r \geq 5$, and that $p \geq 5$ when we treat unipotent classes. Moreover if $\ell \in [13, 20]$ then $\codim V_\kappa(s), \codim C_V(u_\Psi) > M$; so the triple $(G, \lambda, p)$ satisfies $\ssdiamevcon$ and $\udiamcon$. We may therefore assume from now on that $\ell \in [9, 12]$. We need only consider semisimple classes $s^G$ with $|\Phi(s)| \leq M - (2\ell^2 - 14\ell + 34) \leq 20$, each of which by inspection has a subsystem of type $A_3A_2$ disjoint from $\Phi(s)$, and unipotent classes of dimension greater than $M_3$, each of which has the class $A_3A_2$ in its closure by Lemma~\ref{lem: various classes in classical groups by dim}(vii).

Now take $\Psi = \langle \alpha_1, \alpha_2, \alpha_3, \alpha_5, \alpha_6 \rangle$ of type $A_3A_2$. The $\Psi$-net table is as follows.
$$
\begin{array}{|*5{>{\ss}c|}}
\hline
\multicolumn{2}{|>{\ss}c|}{\Psi\mathrm{-nets}} & & \multicolumn{1}{|>{\ss}c|}{c(s)} & \multicolumn{1}{|>{\ss}c|}{c(u_\Psi)} \\
\cline{1-2} \cline{4-5}
      \bar\nu      & n_1 &                        m                       &              r \geq 5              &              p \geq 5              \\
\hline
  \bom_2 + \bom_5  & 18  &                        1                       &                 12                 &                 14                 \\
  \bom_1 + \bom_5  & 12  &                    \ell - 6                    &             9(\ell - 6)            &             9(\ell - 6)            \\
  \bom_1 + \bom_6  & 12  &                        1                       &                  9                 &                  9                 \\
       \bom_1      &  4  &         \frac{1}{2}(\ell - 6)(\ell - 7)        &  \frac{3}{2}(\ell - 6)(\ell - 7)   &  \frac{3}{2}(\ell - 6)(\ell - 7)   \\
       \bom_2      &  6  &                    \ell - 6                    &             4(\ell - 6)            &             4(\ell - 6)            \\
       \bom_3      &  4  &                        1                       &                  3                 &                  3                 \\
       \bom_5      &  3  &         \frac{1}{2}(\ell - 6)(\ell - 7)        &        (\ell - 6)(\ell - 7)        &        (\ell - 6)(\ell - 7)        \\
       \bom_6      &  3  &                    \ell - 6                    &             2(\ell - 6)            &             2(\ell - 6)            \\
         0         &  1  & \frac{1}{6}(\ell^3 - 21\ell^2 + 146\ell - 330) &                                    &                                    \\
\hline
\multicolumn{3}{c|}{}                                                     & \frac{1}{2}(5\ell^2 - 35\ell + 78) & \frac{1}{2}(5\ell^2 - 35\ell + 82) \\
\cline{4-5}
\end{array}
$$
Thus if $\ell \in [10, 12]$ we have $\codim V_\kappa(s), \codim C_V(u_\Psi) > M$; so the triple $(G, \lambda, p)$ satisfies $\ssdiamevcon$ and $\udiamcon$. We may therefore assume from now on that $\ell = 9$. We have $\codim V_\kappa(s) \geq 84$ and $\codim C_V(u_\Psi) \geq 86$, whereas $M = 90$. We need only consider semisimple classes $s^G$ with $|\Phi(s)| \leq 6$, each of which by inspection has a subsystem of type $A_4A_1$ disjoint from $\Phi(s)$, and unipotent classes of dimension at least $86$, each of which has the class $A_4A_1$ in its closure by Lemma~\ref{lem: various classes in A_ell for fixed ell}(i).

Now take $\Psi = \langle \alpha_1, \alpha_2, \alpha_3, \alpha_4, \alpha_6 \rangle$ of type $A_4A_1$. The $\Psi$-net table is as follows.
$$
\begin{array}{|*5{>{\ss}c|}}
\hline
\multicolumn{2}{|>{\ss}c|}{\Psi\mathrm{-nets}} & & \multicolumn{1}{|>{\ss}c|}{c(s)} & \multicolumn{1}{|>{\ss}c|}{c(u_\Psi)} \\
\cline{1-2} \cline{4-5}
      \bar\nu      & n_1 & m & r \geq 5 & p \geq 5 \\
\hline
  \bom_2 + \bom_6  & 20  & 1 &    16    &    16    \\
  \bom_1 + \bom_6  & 10  & 3 &    24    &    24    \\
       \bom_1      &  5  & 4 &    16    &    16    \\
       \bom_2      & 10  & 3 &    24    &    24    \\
       \bom_3      & 10  & 1 &     8    &     8    \\
       \bom_6      &  2  & 3 &     3    &     3    \\
         0         &  1  & 4 &          &          \\
\hline
\multicolumn{3}{c|}{}        &    91    &    91    \\
\cline{4-5}
\end{array}
$$
Thus $\codim V_\kappa(s), \codim C_V(u_\Psi) > M$; so the triple $(G, \lambda, p)$ satisfies $\ssdiamevcon$ and $\udiamcon$.
\end{proof}

\begin{prop}\label{prop: A_ell, omega_4, nets}
Let $G = A_\ell$ for $\ell \in [8, 11]$ and $\lambda = \omega_4$; then the triple $(G, \lambda, p)$ satisfies $\ssdiamevcon$ and $\udiamcon$.
\end{prop}

\begin{proof}
The weight table is as follows.
$$
\begin{array}{|*4{>{\ss}c|}}
\hline
i & \mu & |W.\mu| & m_\mu \\
\hline
1 & \omega_4 & \frac{1}{24}\ell(\ell^2 - 1)(\ell - 2) & 1 \\
\hline
\end{array}
$$
We have $M = \ell(\ell + 1)$, $M_3 = 2\lfloor \frac{1}{3}(\ell + 1)^2 \rfloor$ and $M_2 = \lfloor \frac{1}{2}(\ell + 1)^2 \rfloor$.

Take $\Psi = \langle \alpha_1 \rangle$ of type $A_1$. The $\Psi$-net table is as follows.
$$
\begin{array}{|*5{>{\ss}c|}}
\hline
\multicolumn{2}{|>{\ss}c|}{\Psi\mathrm{-nets}} & & \multicolumn{1}{|>{\ss}c|}{c(s)} & \multicolumn{1}{|>{\ss}c|}{c(u_\Psi)} \\
\cline{1-2} \cline{4-5}
    \bar\nu    & n_1 &                           m                           &                  r \geq 2                 &                  p \geq 2                 \\
\hline
     \bom_1    &  2  &       \frac{1}{6}(\ell - 1)(\ell - 2)(\ell - 3)       & \frac{1}{6}(\ell - 1)(\ell - 2)(\ell - 3) & \frac{1}{6}(\ell - 1)(\ell - 2)(\ell - 3) \\
       0       &  1  & \frac{1}{24}(\ell - 1)(\ell - 2)(\ell^2 - 7\ell + 24) &                                           &                                           \\
\hline
\multicolumn{3}{c|}{}                                                        & \frac{1}{6}(\ell - 1)(\ell - 2)(\ell - 3) & \frac{1}{6}(\ell - 1)(\ell - 2)(\ell - 3) \\
\cline{4-5}
\end{array}
$$
Thus $\codim V_\kappa(s), \codim C_V(u_\Psi) > 2\ell = \dim {u_\Psi}^G$. We therefore need only consider semisimple classes $s^G$ with $|\Phi(s)| < M - 2\ell = m_{{A_1}^2}$, each of which has a subsystem of type ${A_1}^2$ disjoint from $\Phi(s)$, and unipotent classes of dimension greater than $2\ell$, each of which has the class ${A_1}^2$ in its closure by Lemma~\ref{lem: various classes in classical groups by dim}(i).

Now take $\Psi = \langle \alpha_1, \alpha_3 \rangle$ of type ${A_1}^2$. The $\Psi$-net table is as follows.
$$
\begin{array}{|*5{>{\ss}c|}}
\hline
\multicolumn{2}{|>{\ss}c|}{\Psi\mathrm{-nets}} & & \multicolumn{1}{|>{\ss}c|}{c(s)} & \multicolumn{1}{|>{\ss}c|}{c(u_\Psi)} \\
\cline{1-2} \cline{4-5}
      \bar\nu      & n_1 &                              m                              &                  r \geq 2                  &                  p \geq 2                  \\
\hline
  \bom_1 + \bom_3  &  4  &               \frac{1}{2}(\ell - 3)(\ell - 4)               &            (\ell - 3)(\ell - 4)            &            (\ell - 3)(\ell - 4)            \\
       \bom_1      &  2  &          \frac{1}{6}(\ell - 3)(\ell^2 - 9\ell + 26)         & \frac{1}{6}(\ell - 3)(\ell^2 - 9\ell + 26) & \frac{1}{6}(\ell - 3)(\ell^2 - 9\ell + 26) \\
       \bom_3      &  2  &          \frac{1}{6}(\ell - 3)(\ell^2 - 9\ell + 26)         & \frac{1}{6}(\ell - 3)(\ell^2 - 9\ell + 26) & \frac{1}{6}(\ell - 3)(\ell^2 - 9\ell + 26) \\
         0         &  1  & \frac{1}{24}(\ell^4 - 18\ell^3 + 143\ell^2 - 510\ell + 672) &                                            &                                   \\
\hline
\multicolumn{3}{c|}{}                                                                  & \frac{1}{3}(\ell - 3)(\ell^2 - 6\ell + 14) & \frac{1}{3}(\ell - 3)(\ell^2 - 6\ell + 14) \\
\cline{4-5}
\end{array}
$$
Thus $\codim V_\kappa(s), \codim C_V(u_\Psi) > M_2 > 4\ell - 4 = \dim {u_\Psi}^G$; we may therefore assume from now on that $r \geq 3$, and that $p \geq 3$ when we treat unipotent classes. Moreover if $\ell \in [10, 11]$ then $\codim V_\kappa(s), \codim C_V(u_\Psi) > M$; so the triple $(G, \lambda, p)$ satisfies $\ssdiamevcon$ and $\udiamcon$. We may therefore assume from now on that $\ell \in [8, 9]$. We need only consider semisimple classes $s^G$ with $|\Phi(s)| < M - M_2 = m_{A_2}$, each of which has a subsystem of type $A_2$ disjoint from $\Phi(s)$, and unipotent classes of dimension greater than $M_2$, each of which has the class $A_2$ in its closure by Lemma~\ref{lem: various classes in classical groups by dim}(iii).

Now take $\Psi = \langle \alpha_1, \alpha_2 \rangle$ of type $A_2$. The $\Psi$-net table is as follows.
$$
\begin{array}{|*5{>{\ss}c|}}
\hline
\multicolumn{2}{|>{\ss}c|}{\Psi\mathrm{-nets}} & & \multicolumn{1}{|>{\ss}c|}{c(s)} & \multicolumn{1}{|>{\ss}c|}{c(u_\Psi)} \\
\cline{1-2} \cline{4-5}
      \bar\nu      & n_1 &                            m                            &                  r \geq 3                 &                  p \geq 3                 \\
\hline
       \bom_1      &  3  &        \frac{1}{6}(\ell - 2)(\ell - 3)(\ell - 4)        & \frac{1}{3}(\ell - 2)(\ell - 3)(\ell - 4) & \frac{1}{3}(\ell - 2)(\ell - 3)(\ell - 4) \\
       \bom_2      &  3  &             \frac{1}{2}(\ell - 2)(\ell - 3)             &            (\ell - 2)(\ell - 3)           &            (\ell - 2)(\ell - 3)           \\
         0         &  1  & \frac{1}{24}(\ell - 2)(\ell^3 - 12\ell^2 + 47\ell - 36) &                                           &                                           \\
\hline
\multicolumn{3}{c|}{}                                                              & \frac{1}{3}(\ell - 1)(\ell - 2)(\ell - 3) & \frac{1}{3}(\ell - 1)(\ell - 2)(\ell - 3) \\
\cline{4-5}
\end{array}
$$
Thus $\codim V_\kappa(s), \codim C_V(u_\Psi) > M_3 > 4\ell - 2 = \dim {u_\Psi}^G$; we may therefore assume from now on that $r \geq 5$, and that $p \geq 5$ when we treat unipotent classes. Moreover if $\ell = 9$ then $\codim V_\kappa(s), \codim C_V(u_\Psi) > M$; so the triple $(G, \lambda, p)$ satisfies $\ssdiamevcon$ and $\udiamcon$. We may therefore assume from now on that $\ell = 8$. We have $\codim V_\kappa(s), \codim C_V(u_\Psi) \geq 70$ while $M = 72$. We therefore need only consider semisimple classes $s^G$ with $|\Phi(s)| \leq 2$, each of which has a subsystem of type $A_3$ disjoint from $\Phi(s)$, and unipotent classes of dimension at least $70$, each of which has the class $A_3$ in its closure by Lemma~\ref{lem: various classes in classical groups by dim}(vi).

Now take $\Psi = \langle \alpha_1, \alpha_2, \alpha_3 \rangle$ of type $A_3$. The $\Psi$-net table is as follows.
$$
\begin{array}{|*5{>{\ss}c|}}
\hline
\multicolumn{2}{|>{\ss}c|}{\Psi\mathrm{-nets}} & & \multicolumn{1}{|>{\ss}c|}{c(s)} & \multicolumn{1}{|>{\ss}c|}{c(u_\Psi)} \\
\cline{1-2} \cline{4-5}
      \bar\nu      & n_1 &  m & r \geq 5 & p \geq 5 \\
\hline
       \bom_1      &  4  & 10 &    30    &    30    \\
       \bom_2      &  6  & 10 &    40    &    40    \\
       \bom_3      &  4  &  5 &    15    &    15    \\
         0         &  1  &  6 &          &          \\
\hline
\multicolumn{3}{c|}{}         &    85    &    85    \\
\cline{4-5}
\end{array}
$$
Thus $\codim V_\kappa(s), \codim C_V(u_\Psi) > M$; so the triple $(G, \lambda, p)$ satisfies $\ssdiamevcon$ and $\udiamcon$.
\end{proof}

\begin{prop}\label{prop: A_9, omega_5, nets}
Let $G = A_9$ and $\lambda = \omega_5$; then the triple $(G, \lambda, p)$ satisfies $\ssdiamevcon$ and $\udiamcon$.
\end{prop}

\begin{proof}
The weight table is as follows.
$$
\begin{array}{|*4{>{\ss}c|}}
\hline
i & \mu & |W.\mu| & m_\mu \\
\hline
1 & \omega_5 & 252 & 1 \\
\hline
\end{array}
$$
We have $M = 90$.

Take $\Psi = \langle \alpha_1 \rangle$ of type $A_1$. The $\Psi$-net table is as follows.
$$
\begin{array}{|*5{>{\ss}c|}}
\hline
\multicolumn{2}{|>{\ss}c|}{\Psi\mathrm{-nets}} & & \multicolumn{1}{|>{\ss}c|}{c(s)} & \multicolumn{1}{|>{\ss}c|}{c(u_\Psi)} \\
\cline{1-2} \cline{4-5}
    \bar\nu    & n_1 &  m  & r \geq 2 & p \geq 2 \\
\hline
     \bom_1    &  2  &  70 &    70    &    70    \\
       0       &  1  & 112 &          &          \\
\hline
\multicolumn{3}{c|}{}      &    70    &    70    \\
\cline{4-5}
\end{array}
$$
Thus $\codim V_\kappa(s), \codim C_V(u_\Psi) \geq 70 > 18 = \dim {u_\Psi}^G$. We therefore need only consider semisimple classes $s^G$ with $|\Phi(s)| \leq 20 < 72 = m_{{A_1}^2}$, each of which has a subsystem of type ${A_1}^2$ disjoint from $\Phi(s)$, and unipotent classes of dimension at least $70$, each of which has the class ${A_1}^2$ in its closure by Lemma~\ref{lem: various classes in classical groups by dim}(i).

Now take $\Psi = \langle \alpha_1, \alpha_3 \rangle$ of type ${A_1}^2$. The $\Psi$-net table is as follows.
$$
\begin{array}{|*5{>{\ss}c|}}
\hline
\multicolumn{2}{|>{\ss}c|}{\Psi\mathrm{-nets}} & & \multicolumn{1}{|>{\ss}c|}{c(s)} & \multicolumn{1}{|>{\ss}c|}{c(u_\Psi)} \\
\cline{1-2} \cline{4-5}
     \bar\nu     & n_1 &  m & r \geq 2 & p \geq 2 \\
\hline
 \bom_1 + \bom_3 &  4  & 20 &    40    &    40    \\
      \bom_1     &  2  & 30 &    30    &    30    \\
      \bom_3     &  2  & 30 &    30    &    30    \\
         0       &  1  & 52 &          &          \\
\hline
\multicolumn{3}{c|}{}       &   100    &   100    \\
\cline{4-5}
\end{array}
$$
Thus $\codim V_\kappa(s), \codim C_V(u_\Psi) > M$; so the triple $(G, \lambda, p)$ satisfies $\ssdiamevcon$ and $\udiamcon$.
\end{proof}

\begin{prop}\label{prop: A_ell, omega_1 + omega_2, nets}
Let $G = A_\ell$ for $\ell \in [4, \infty)$ and $\lambda = \omega_1 + \omega_2$; then the triple $(G, \lambda, p)$ satisfies $\ssdiamevcon$ and $\udiamcon$.
\end{prop}

\begin{proof}
Write $\z = \z_{p, 3}$. By Lemma~\ref{lem: multiplicities in 3omega_1 and omega_1 + omega_2} the weight table is as follows.
$$
\begin{array}{|*4{>{\ss}c|}}
\hline
i & \mu & |W.\mu| & m_\mu \\
\hline
2 & \omega_1 + \omega_2 &        \ell(\ell + 1)       &   1    \\
1 &       \omega_3      & \frac{1}{6}\ell(\ell^2 - 1) & 2 - \z \\
\hline
\end{array}
$$
We have $M = \ell(\ell + 1)$, $M_3 = 2\lfloor \frac{1}{3}(\ell + 1)^2 \rfloor$ and $M_2 = \lfloor \frac{1}{2}(\ell + 1)^2 \rfloor$.

First suppose $p \neq 3$. Take $\Psi = \langle \alpha_1 \rangle$ of type $A_1$. The $\Psi$-net table is as follows.
$$
\begin{array}{|*7{>{\ss}c|}}
\hline
\multicolumn{3}{|>{\ss}c|}{\Psi\mathrm{-nets}} & & \multicolumn{1}{|>{\ss}c|}{c(s)} & \multicolumn{2}{|>{\ss}c|}{c(u_\Psi)} \\
\cline{1-3} \cline{5-7}
    \bar\nu    & n_1 & n_2 &                     m                     &       r \geq 2       &         p = 2        &       p \geq 5       \\
\hline
    2\bom_1    &  1  &  2  &                  \ell - 1                 &      2(\ell - 1)     &       \ell - 1       &      2(\ell - 1)     \\
     \bom_1    &  0  &  2  &                    \ell                   &         \ell         &         \ell         &         \ell         \\
     \bom_1    &  2  &  0  &      \frac{1}{2}(\ell - 1)(\ell - 2)      & (\ell - 1)(\ell - 2) & (\ell - 1)(\ell - 2) & (\ell - 1)(\ell - 2) \\
       0       &  0  &  1  &            (\ell - 1)(\ell - 2)           &                      &                      &                      \\
       0       &  1  &  0  & \frac{1}{6}(\ell - 1)(\ell - 2)(\ell - 3) &                      &                      &                      \\
\hline
\multicolumn{4}{c|}{}                                                  &        \ell^2        &   \ell^2 - \ell + 1  &        \ell^2        \\
\cline{5-7}
\end{array}
$$
Thus $\codim V_\kappa(s), \codim C_V(u_\Psi) > M_2 > 2\ell = \dim {u_\Psi}^G$; we may therefore assume from now on that $r \geq 3$, and that $p \geq 5$ when we treat unipotent classes. We need only consider semisimple classes $s^G$ with $|\Phi(s)| < M - M_2 = m_{A_2}$, each of which has a subsystem of type $A_2$ disjoint from $\Phi(s)$, and unipotent classes of dimension greater than $M_2$, each of which has the class $A_2$ in its closure by Lemma~\ref{lem: various classes in classical groups by dim}(iii).

Now take $\Psi = \langle \alpha_1, \alpha_2 \rangle$ of type $A_2$. The $\Psi$-net table is as follows.
$$
\begin{array}{|*6{>{\ss}c|}}
\hline
\multicolumn{3}{|>{\ss}c|}{\Psi\mathrm{-nets}} & & \multicolumn{1}{|>{\ss}c|}{c(s)} & \multicolumn{1}{|>{\ss}c|}{c(u_\Psi)} \\
\cline{1-3} \cline{5-6}
      \bar\nu      & n_1 & n_2 &                     m                     &        r \geq 3       &        p \geq 5       \\
\hline
  \bom_1 + \bom_2  &  1  &  6  &                     1                     &           5           &           6           \\
      2\bom_1      &  3  &  3  &                  \ell - 2                 &      6(\ell - 2)      &      6(\ell - 2)      \\
       \bom_1      &  0  &  3  &                  \ell - 2                 &      2(\ell - 2)      &      2(\ell - 2)      \\
       \bom_1      &  3  &  0  &      \frac{1}{2}(\ell - 2)(\ell - 3)      & 2(\ell - 2)(\ell - 3) & 2(\ell - 2)(\ell - 3) \\
         0         &  0  &  1  &            (\ell - 2)(\ell - 3)           &                       &                       \\
         0         &  1  &  0  & \frac{1}{6}(\ell - 2)(\ell - 3)(\ell - 4) &                       &                       \\
\hline
\multicolumn{4}{c|}{}                                                      &  2\ell^2 - 2\ell + 1  &  2\ell^2 - 2\ell + 2  \\
\cline{5-6}
\end{array}
$$
Thus $\codim V_\kappa(s), \codim C_V(u_\Psi) > M$; so the triple $(G, \lambda, p)$ satisfies $\ssdiamevcon$ and $\udiamcon$.

Now suppose $p = 3$. Take $\Psi = \langle \alpha_1 \rangle$ of type $A_1$. The $\Psi$-net table is as follows.
$$
\begin{array}{|*7{>{\ss}c|}}
\hline
\multicolumn{3}{|>{\ss}c|}{\Psi\mathrm{-nets}} & & \multicolumn{2}{|>{\ss}c|}{c(s)} & \multicolumn{1}{|>{\ss}c|}{c(u_\Psi)} \\
\cline{1-3} \cline{5-7}
    \bar\nu    & n_1 & n_2 &                     m                     &              r = 2              &             r \geq 5            &              p = 3              \\
\hline
    2\bom_1    &  1  &  2  &                  \ell - 1                 &             \ell - 1            &           2(\ell - 1)           &           2(\ell - 1)           \\
     \bom_1    &  0  &  2  &                    \ell                   &               \ell              &               \ell              &               \ell              \\
     \bom_1    &  2  &  0  &      \frac{1}{2}(\ell - 1)(\ell - 2)      & \frac{1}{2}(\ell - 1)(\ell - 2) & \frac{1}{2}(\ell - 1)(\ell - 2) & \frac{1}{2}(\ell - 1)(\ell - 2) \\
       0       &  0  &  1  &            (\ell - 1)(\ell - 2)           &                                 &                                 &                                 \\
       0       &  1  &  0  & \frac{1}{6}(\ell - 1)(\ell - 2)(\ell - 3) &                                 &                                 &                                 \\
\hline
\multicolumn{4}{c|}{}                                                  &    \frac{1}{2}\ell(\ell + 1)    & \frac{1}{2}(\ell^2 + 3\ell - 2) & \frac{1}{2}(\ell^2 + 3\ell - 2) \\
\cline{5-7}
\end{array}
$$
Thus $\codim V_\kappa(s), \codim C_V(u_\Psi) > 2\ell = \dim {u_\Psi}^G$. We therefore need only consider semisimple classes $s^G$ with $|\Phi(s)| < M - 2\ell = m_{{A_1}^2}$, each of which has a subsystem of type ${A_1}^2$ disjoint from $\Phi(s)$, and unipotent classes of dimension greater than $2\ell$, each of which has the class ${A_1}^2$ in its closure by Lemma~\ref{lem: various classes in classical groups by dim}(i).

Now take $\Psi = \langle \alpha_1, \alpha_3 \rangle$ of type ${A_1}^2$. The $\Psi$-net table is as follows.
$$
\begin{array}{|*7{>{\ss}c|}}
\hline
\multicolumn{3}{|>{\ss}c|}{\Psi\mathrm{-nets}} & & \multicolumn{2}{|>{\ss}c|}{c(s)} & \multicolumn{1}{|>{\ss}c|}{c(u_\Psi)} \\
\cline{1-3} \cline{5-7}
      \bar\nu      & n_1 & n_2 &                     m                     &              r = 2              &             r \geq 5            &              p = 3              \\
\hline
  2\bom_1 + \bom_3 &  2  &  4  &                     1                     &                3                &                4                &                4                \\
  \bom_1 + 2\bom_3 &  2  &  4  &                     1                     &                3                &                4                &                4                \\
  \bom_1 + \bom_3  &  4  &  0  &                  \ell - 3                 &           2(\ell - 3)           &           2(\ell - 3)           &           2(\ell - 3)           \\
      2\bom_1      &  1  &  2  &                  \ell - 3                 &             \ell - 3            &           2(\ell - 3)           &           2(\ell - 3)           \\
      2\bom_3      &  1  &  2  &                  \ell - 3                 &             \ell - 3            &           2(\ell - 3)           &           2(\ell - 3)           \\
       \bom_1      &  0  &  2  &                  \ell - 2                 &             \ell - 2            &             \ell - 2            &             \ell - 2            \\
       \bom_3      &  0  &  2  &                  \ell - 2                 &             \ell - 2            &             \ell - 2            &             \ell - 2            \\
       \bom_1      &  2  &  0  &      \frac{1}{2}(\ell - 3)(\ell - 4)      & \frac{1}{2}(\ell - 3)(\ell - 4) & \frac{1}{2}(\ell - 3)(\ell - 4) & \frac{1}{2}(\ell - 3)(\ell - 4) \\
       \bom_3      &  2  &  0  &      \frac{1}{2}(\ell - 3)(\ell - 4)      & \frac{1}{2}(\ell - 3)(\ell - 4) & \frac{1}{2}(\ell - 3)(\ell - 4) & \frac{1}{2}(\ell - 3)(\ell - 4) \\
         0         &  0  &  1  &            (\ell - 3)(\ell - 4)           &                                 &                                 &                                 \\
         0         &  1  &  0  & \frac{1}{6}(\ell - 3)(\ell - 4)(\ell - 5) &                                 &                                 &                                 \\
\hline
\multicolumn{4}{c|}{}                                                      &        \ell^2 - \ell + 2        &        \ell^2 + \ell - 2        &        \ell^2 + \ell - 2        \\
\cline{5-7}
\end{array}
$$
Thus $\codim V_\kappa(s) \geq M - 2$ unless $r = 2$, in which case $\codim V_\kappa(s) > M_r$, and $\codim C_V(u_\Psi) > M_p$; so the triple $(G, \lambda, p)$ satisfies $\udiamcon$, and we may assume from now on that $r \geq 5$. We need only consider semisimple classes $s^G$ with $|\Phi(s)| \leq 2$, each of which has a subsystem of type $A_3$ disjoint from $\Phi(s)$.

Now take $\Psi = \langle \alpha_1, \alpha_2, \alpha_3 \rangle$ of type $A_3$. The $\Psi$-net table is as follows.
$$
\begin{array}{|*5{>{\ss}c|}}
\hline
\multicolumn{3}{|>{\ss}c|}{\Psi\mathrm{-nets}} & & \multicolumn{1}{|>{\ss}c|}{c(s)} \\
\cline{1-3} \cline{5-5}
      \bar\nu      & n_1 & n_2 &                     m                     &             r \geq 5            \\
\hline
  \bom_1 + \bom_2  &  4  & 12  &                     1                     &                12               \\
      2\bom_1      &  6  &  4  &                  \ell - 3                 &           8(\ell - 3)           \\
       \bom_1      &  0  &  4  &                  \ell - 3                 &           3(\ell - 3)           \\
       \bom_1      &  4  &  0  &      \frac{1}{2}(\ell - 3)(\ell - 4)      & \frac{3}{2}(\ell - 3)(\ell - 4) \\
         0         &  0  &  1  &            (\ell - 3)(\ell - 4)           &                                 \\
         0         &  1  &  0  & \frac{1}{6}(\ell - 3)(\ell - 4)(\ell - 5) &                                 \\
\hline
\multicolumn{4}{c|}{}                                                      & \frac{1}{2}(3\ell^2 + \ell - 6) \\
\cline{5-5}
\end{array}
$$
Thus $\codim V_\kappa(s) > M$; so the triple $(G, \lambda, p)$ satisfies $\ssdiamevcon$.
\end{proof}

\begin{prop}\label{prop: A_3, omega_1 + omega_2, nets}
Let $G = A_3$ and $\lambda = \omega_1 + \omega_2$ with $p \neq 3$; then the triple $(G, \lambda, p)$ satisfies $\ssdiamevcon$ and $\udiamcon$.
\end{prop}

\begin{proof}
The weight table is as follows.
$$
\begin{array}{|*4{>{\ss}c|}}
\hline
i & \mu & |W.\mu| & m_\mu \\
\hline
2 & \omega_1 + \omega_2 & 12 & 1 \\
1 &       \omega_3      &  4 & 2 \\
\hline
\end{array}
$$
We have $M = 12$ and $M_2 = 8$.

Take $\Psi = \langle \alpha_1 \rangle$ of type $A_1$. The $\Psi$-net table is as follows.
$$
\begin{array}{|*8{>{\ss}c|}}
\hline
\multicolumn{3}{|>{\ss}c|}{\Psi\mathrm{-nets}} & & \multicolumn{2}{|>{\ss}c|}{c(s)} & \multicolumn{2}{|>{\ss}c|}{c(u_\Psi)} \\
\cline{1-3} \cline{5-8}
    \bar\nu    & n_1 & n_2 & m & r = 2 & r \geq 3 & p = 2 & p \geq 5 \\
\hline
    2\bom_1    &  1  &  2  & 2 &   4   &     4    &   2   &     4    \\
     \bom_1    &  0  &  2  & 3 &   3   &     3    &   3   &     3    \\
     \bom_1    &  2  &  0  & 1 &   2   &     2    &   2   &     2    \\
       0       &  0  &  1  & 2 &       &          &       &          \\
\hline
\multicolumn{4}{c|}{}          &   9   &     9    &   7   &     9    \\
\cline{5-8}
\end{array}
$$
Thus $\codim V_\kappa(s) \geq 9 > M_2$, and $\codim C_V(u_\Psi) > 6 = \dim {u_\Psi}^G$; we may therefore assume from now on that $r \geq 3$. We need only consider semisimple classes $s^G$ with $|\Phi(s)| \leq 3 < 6 = m_{{A_1}^2}$, each of which has a subsystem of type ${A_1}^2$ disjoint from $\Phi(s)$, and unipotent classes of dimension greater than $6$, each of which has the class ${A_1}^2$ in its closure by Lemma~\ref{lem: various classes in classical groups by dim}(i).

Now take $\Psi = \langle \alpha_1, \alpha_3 \rangle$ of type ${A_1}^2$. The $\Psi$-net table is as follows.
$$
\begin{array}{|*7{>{\ss}c|}}
\hline
\multicolumn{3}{|>{\ss}c|}{\Psi\mathrm{-nets}} & & \multicolumn{1}{|>{\ss}c|}{c(s)} & \multicolumn{2}{|>{\ss}c|}{c(u_\Psi)} \\
\cline{1-3} \cline{5-7}
      \bar\nu      & n_1 & n_2 & m & r \geq 3 & p = 2 & p \geq 5 \\
\hline
 2\bom_1 + \bom_3  &  2  &  4  & 1 &     5    &   4   &     5    \\
 \bom_1 + 2\bom_3  &  2  &  4  & 1 &     5    &   4   &     5    \\
       \bom_1      &  0  &  2  & 1 &     1    &   1   &     1    \\
       \bom_3      &  0  &  2  & 1 &     1    &   1   &     1    \\
\hline
\multicolumn{4}{c|}{}              &    12    &  10   &    12    \\
\cline{5-7}
\end{array}
$$
Thus $\codim V_\kappa(s) \geq M$, and $\codim C_V(u_\Psi) \geq M$ unless $p = 2$, in which case $\codim C_V(u_\Psi) > M_p$; we may therefore assume that $p \geq 5$ when we treat unipotent classes. We need only consider semisimple classes $s^G$ with $\Phi(s) = \emptyset$, each of which has a subsystem of type $A_2$ disjoint from $\Phi(s)$, and unipotent classes of dimension $12$, of which the only one is the regular class $A_3$, which has the class $A_2$ in its closure by Lemma~\ref{lem: any class in closure of reg class}.

Now take $\Psi = \langle \alpha_1, \alpha_2 \rangle$ of type $A_2$. The $\Psi$-net table is as follows.
$$
\begin{array}{|*6{>{\ss}c|}}
\hline
\multicolumn{3}{|>{\ss}c|}{\Psi\mathrm{-nets}} & & \multicolumn{1}{|>{\ss}c|}{c(s)} & \multicolumn{1}{|>{\ss}c|}{c(u_\Psi)} \\
\cline{1-3} \cline{5-6}
      \bar\nu      & n_1 & n_2 & m & r \geq 3 & p \geq 5 \\
\hline
  \bom_1 + \bom_2  &  1  &  6  & 1 &     5    &     6    \\
      2\bom_1      &  3  &  3  & 1 &     6    &     6    \\
       \bom_1      &  0  &  3  & 1 &     2    &     2    \\
\hline
\multicolumn{4}{c|}{}              &    13    &    14    \\
\cline{5-6}
\end{array}
$$
Thus $\codim V_\kappa(s), \codim C_V(u_\Psi) > M$; so the triple $(G, \lambda, p)$ satisfies $\ssdiamevcon$ and $\udiamcon$.
\end{proof}

\begin{prop}\label{prop: A_ell, omega_2 + omega_ell, nets}
Let $G = A_\ell$ for $\ell \in [4, 5]$ and $\lambda = \omega_2 + \omega_\ell$; then the triple $(G, \lambda, p)$ satisfies $\ssdiamevcon$ and $\udiamcon$.
\end{prop}

\begin{proof}
Write $\z = \z_{p, \ell}$ and $\z' = \z\z_{p, 2}$. The weight table is as follows.
$$
\begin{array}{|*4{>{\ss}c|}}
\hline
i & \mu & |W.\mu| & m_\mu \\
\hline
2 & \omega_2 + \omega_\ell & \frac{1}{2}\ell(\ell^2 - 1) &       1       \\
1 &        \omega_1        &           \ell + 1          & \ell - 1 - \z \\
\hline
\end{array}
$$
We have $M = \ell(\ell + 1)$ and $M_2 = \lfloor \frac{1}{2}(\ell + 1)^2 \rfloor$.

Take $\Psi = \langle \alpha_1 \rangle$ of type $A_1$. The $\Psi$-net table is as follows.
$$
\begin{array}{|*7{>{\ss}c|}}
\hline
\multicolumn{3}{|>{\ss}c|}{\Psi\mathrm{-nets}} & & \multicolumn{1}{|>{\ss}c|}{c(s)} & \multicolumn{2}{|>{\ss}c|}{c(u_\Psi)} \\
\cline{1-3} \cline{5-7}
    \bar\nu    & n_1 & n_2 &                     m                     &             r \geq 2            &                 p = 2                 &             p \geq 3            \\
\hline
    2\bom_1    &  1  &  2  &                  \ell - 1                 &           2(\ell - 1)           &                \ell - 1               &           2(\ell - 1)           \\
     \bom_1    &  0  &  2  &      \frac{3}{2}(\ell - 1)(\ell - 2)      & \frac{3}{2}(\ell - 1)(\ell - 2) &    \frac{3}{2}(\ell - 1)(\ell - 2)    & \frac{3}{2}(\ell - 1)(\ell - 2) \\
     \bom_1    &  2  &  0  &                     1                     &          \ell - 1 - \z          &             \ell - 1 - \z             &          \ell - 1 - \z          \\
       0       &  0  &  1  & \frac{1}{2}(\ell - 1)(\ell^2 - 5\ell + 8) &                                 &                                       &                                 \\
\hline
\multicolumn{4}{c|}{}                                                  &  \frac{3}{2}\ell(\ell - 1) - \z & \frac{1}{2}(\ell - 1)(3\ell - 2) - \z & \frac{3}{2}\ell(\ell - 1) - \z  \\
\cline{5-7}
\end{array}
$$
Thus $\codim V_\kappa(s), \codim C_V(u_\Psi) > M_2 > 2\ell = \dim {u_\Psi}^G$; we may therefore assume from now on that $r \geq 3$, and that $p \geq 3$ when we treat unipotent classes. We need only consider semisimple classes $s^G$ with $|\Phi(s)| < M - M_2 = m_{A_2}$, each of which has a subsystem of type $A_2$ disjoint from $\Phi(s)$, and unipotent classes of dimension greater than $M_2$, each of which has the class $A_2$ in its closure by Lemma~\ref{lem: various classes in classical groups by dim}(iii).

Now take $\Psi = \langle \alpha_1, \alpha_2 \rangle$ of type $A_2$. The $\Psi$-net table is as follows.
$$
\begin{array}{|*7{>{\ss}c|}}
\hline
\multicolumn{3}{|>{\ss}c|}{\Psi\mathrm{-nets}} & & \multicolumn{1}{|>{\ss}c|}{c(s)} & \multicolumn{2}{|>{\ss}c|}{c(u_\Psi)} \\
\cline{1-3} \cline{5-7}
      \bar\nu      & n_1 & n_2 &                     m                     &             r \geq 3             &         p = 3         &          p \geq 5         \\
\hline
  \bom_1 + \bom_2  &  1  &  6  &                  \ell - 2                 &    6(\ell - 2) - (\ell - 2)\z'   &      4(\ell - 2)      &        6(\ell - 2)        \\
      2\bom_2      &  3  &  3  &                     1                     &            2\ell - 2\z           &       2\ell - 2       &        2\ell - 2\z        \\
       \bom_1      &  0  &  3  &            (\ell - 2)(\ell - 3)           &      2(\ell - 2)(\ell - 3)       & 2(\ell - 2)(\ell - 3) &   2(\ell - 2)(\ell - 3)   \\
       \bom_2      &  0  &  3  &      \frac{1}{2}(\ell - 1)(\ell - 2)      &       (\ell - 1)(\ell - 2)       &  (\ell - 1)(\ell - 2) &    (\ell - 1)(\ell - 2)   \\
         0         &  0  &  1  & \frac{1}{2}(\ell - 2)(\ell - 3)(\ell - 4) &                                  &                       &                           \\
\hline
\multicolumn{4}{c|}{}                                                      & 3\ell^2 - 5\ell + 2 - 2\z - 2\z' &  3\ell^2 - 7\ell + 4  & 3\ell^2 - 5\ell + 2 - 2\z \\
\cline{5-7}
\end{array}
$$
Thus $\codim V_\kappa(s), \codim C_V(u_\Psi) > M$; so the triple $(G, \lambda, p)$ satisfies $\ssdiamevcon$ and $\udiamcon$.
\end{proof}

\begin{prop}\label{prop: D_5, omega_3, nets}
Let $G = D_5$ and $\lambda = \omega_3$ with $p = 2$; then the triple $(G, \lambda, p)$ satisfies $\ssdiamevcon$ and $\udagcon$.
\end{prop}

\begin{proof}
The weight table is as follows.
$$
\begin{array}{|*4{>{\ss}c|}}
\hline
i & \mu & |W.\mu| & m_\mu \\
\hline
2 & \omega_3 & 80 & 1 \\
1 & \omega_1 & 10 & 2 \\
\hline
\end{array}
$$
We have $M = 40$ and $M_2 = 24$.

Take $\Psi = \langle \alpha_1 \rangle$ of type $A_1$. The $\Psi$-net table is as follows.
$$
\begin{array}{|*6{>{\ss}c|}}
\hline
\multicolumn{3}{|>{\ss}c|}{\Psi\mathrm{-nets}} & & \multicolumn{1}{|>{\ss}c|}{c(s)} & \multicolumn{1}{|>{\ss}c|}{c(u_\Psi)} \\
\cline{1-3} \cline{5-6}
    \bar\nu    & n_1 & n_2 &  m & r \geq 3 & p = 2 \\
\hline
    2\bom_1    &  1  &  2  &  6 &    12    &   6   \\
     \bom_1    &  0  &  2  & 24 &    24    &  24   \\
     \bom_1    &  2  &  0  &  2 &     4    &   4   \\
       0       &  0  &  1  & 20 &          &       \\
\hline
\multicolumn{4}{c|}{}           &    40    &  34   \\
\cline{5-6}
\end{array}
$$
Thus $\codim V_\kappa(s) \geq M$, and $\codim C_V(u_\Psi) > M_2$; so the triple $(G, \lambda, p)$ satisfies $\udagcon$. We need only consider semisimple classes $s^G$ with $\Phi(s) = \emptyset$, each of which has a subsystem of type $A_2$ disjoint from $\Phi(s)$.

Now take $\Psi = \langle \alpha_1, \alpha_2 \rangle$ of type $A_2$. The $\Psi$-net table is as follows.
$$
\begin{array}{|*5{>{\ss}c|}}
\hline
\multicolumn{3}{|>{\ss}c|}{\Psi\mathrm{-nets}} & & \multicolumn{1}{|>{\ss}c|}{c(s)} \\
\cline{1-3} \cline{5-5}
       \bar\nu       & n_1 & n_2 & m & r \geq 3 \\
\hline
   \bom_1 + \bom_2   &  1  &  6  & 4 &    20    \\
       2\bom_1       &  3  &  3  & 1 &     6    \\
       2\bom_2       &  3  &  3  & 1 &     6    \\
        \bom_1       &  0  &  3  & 8 &    16    \\
        \bom_2       &  0  &  3  & 8 &    16    \\
           0         &  0  &  1  & 2 &          \\
\hline
\multicolumn{4}{c|}{}                &    64    \\
\cline{5-5}
\end{array}
$$
Thus $\codim V_\kappa(s) > M$; so the triple $(G, \lambda, p)$ satisfies $\ssdiamevcon$.
\end{proof}

\begin{prop}\label{prop: D_4, omega_1 + omega_4, nets}
Let $G = D_4$ and $\lambda = \omega_1 + \omega_4$; then the triple $(G, \lambda, p)$ satisfies $\ssdiamevcon$ and $\udiamcon$.
\end{prop}

\begin{proof}
Write $\z = \z_{p, 2}$. The weight table is as follows.
$$
\begin{array}{|*4{>{\ss}c|}}
\hline
i & \mu & |W.\mu| & m_\mu \\
\hline
2 & \omega_1 + \omega_4 & 32 &   1    \\
1 &       \omega_3      &  8 & 3 - \z \\
\hline
\end{array}
$$
We have $M = 24$ and $M_2 = 16$.

Take $\Psi = \langle \alpha_1 \rangle$ of type $A_1$. The $\Psi$-net table is as follows.
$$
\begin{array}{|*8{>{\ss}c|}}
\hline
\multicolumn{3}{|>{\ss}c|}{\Psi\mathrm{-nets}} & & \multicolumn{2}{|>{\ss}c|}{c(s)} & \multicolumn{2}{|>{\ss}c|}{c(u_\Psi)} \\
\cline{1-3} \cline{5-8}
    \bar\nu    & n_1 & n_2 & m & r = 2 & r \geq 3 & p = 2 & p \geq 3 \\
\hline
    2\bom_1    &  1  &  2  & 4 &   8   &     8    &   4   &     8    \\
     \bom_1    &  0  &  2  & 8 &   8   &     8    &   8   &     8    \\
     \bom_1    &  2  &  0  & 2 &   6   &  6 - 2\z &   4   &     6    \\
       0       &  0  &  1  & 8 &       &          &       &          \\
\hline
\multicolumn{4}{c|}{}          &  22   & 22 - 2\z &  16   &    22    \\
\cline{5-8}
\end{array}
$$
Thus $\codim V_\kappa(s) \geq 22 - 2\z$, and $\codim C_V(u_\Psi) \geq 22 > 10 = \dim {u_\Psi}^G$ unless $p = 2$, in which case $\codim C_V(u_\Psi) \geq M_p$. We therefore need only consider semisimple classes $s^G$ with $|\Phi(s)| \leq 2 + 2\z$, each of which has a subsystem of type $D_2$ disjoint from $\Phi(s)$, and unipotent classes of dimension at least $16$ or $22$ according as $p = 2$ or $p \geq 3$, each of which has the class $D_2$ in its closure by Lemma~\ref{lem: various classes in classical groups by dim}(ix).

Now take $\Psi = \langle \alpha_3, \alpha_4 \rangle$ of type $D_2$. The $\Psi$-net table is as follows.
$$
\begin{array}{|*8{>{\ss}c|}}
\hline
\multicolumn{3}{|>{\ss}c|}{\Psi\mathrm{-nets}} & & \multicolumn{2}{|>{\ss}c|}{c(s)} & \multicolumn{2}{|>{\ss}c|}{c(u_\Psi)} \\
\cline{1-3} \cline{5-8}
      \bar\nu      & n_1 & n_2 & m & r = 2 & r \geq 3 & p = 2 & p \geq 3 \\
\hline
 2\bom_3 + \bom_4  &  2  &  4  & 2 &  10   & 12 - 2\z &   8   &    12    \\
 \bom_3 + 2\bom_4  &  2  &  4  & 2 &  10   & 12 - 2\z &   8   &    12    \\
       \bom_3      &  0  &  2  & 4 &   4   &     4    &   4   &     4    \\
       \bom_4      &  0  &  2  & 4 &   4   &     4    &   4   &     4    \\
\hline
\multicolumn{4}{c|}{}              &  28   & 32 - 4\z &  24   &    32    \\
\cline{5-8}
\end{array}
$$
Thus $\codim V_\kappa(s) > M$, and $\codim C_V(u_\Psi) > M$ unless $p = 2$, in which case $\codim C_V(u_\Psi) > M_p$; so the triple $(G, \lambda, p)$ satisfies $\ssdiamevcon$ and $\udiamcon$.
\end{proof}

\begin{prop}\label{prop: D_ell, omega_ell, nets}
Let $G = D_\ell$ for $\ell \in [9, 10]$ and $\lambda = \omega_\ell$; then the triple $(G, \lambda, p)$ satisfies $\ssdiamevcon$ and $\udiamcon$.
\end{prop}

\begin{proof}
The weight table is as follows.
$$
\begin{array}{|*4{>{\ss}c|}}
\hline
i & \mu & |W.\mu| & m_\mu \\
\hline
1 & \omega_\ell & 2^{\ell - 1} & 1 \\
\hline
\end{array}
$$
We have $M = 2\ell(\ell - 1)$, $M_3 = 2\lfloor \frac{1}{3}\ell(2\ell - 1) \rfloor$ and $M_2 = 2\lfloor \frac{1}{2}\ell^2 \rfloor$.

Take $\Psi = \langle \alpha_1 \rangle$ of type $A_1$. The $\Psi$-net table is as follows.
$$
\begin{array}{|*5{>{\ss}c|}}
\hline
\multicolumn{2}{|>{\ss}c|}{\Psi\mathrm{-nets}} & & \multicolumn{1}{|>{\ss}c|}{c(s)} & \multicolumn{1}{|>{\ss}c|}{c(u_\Psi)} \\
\cline{1-2} \cline{4-5}
    \bar\nu    & n_1 &       m      &   r \geq 2   &   p \geq 2   \\
\hline
     \bom_1    &  2  & 2^{\ell - 3} & 2^{\ell - 3} & 2^{\ell - 3} \\
       0       &  1  & 2^{\ell - 2} &              &              \\
\hline
\multicolumn{3}{c|}{}               & 2^{\ell - 3} & 2^{\ell - 3} \\
\cline{4-5}
\end{array}
$$
Thus $\codim V_\kappa(s), \codim C_V(u_\Psi) \geq 2^{\ell - 3} > 4\ell - 6 = \dim {u_\Psi}^G$. We therefore need only consider semisimple classes $s^G$ with $|\Phi(s)| \leq M - 2^{\ell - 3}$, each of which has a subsystem of type ${A_1}^2$ or $D_2$ disjoint from $\Phi(s)$, and unipotent classes of dimension at least $2^{\ell - 3}$, each of which has the class ${A_1}^2$ or $D_2$ in its closure by Lemma~\ref{lem: A_1^2 or D_2 in D_ell}.

Now take $\Psi = \langle \alpha_1, \alpha_3 \rangle$ of type ${A_1}^2$, and $\langle \alpha_{\ell - 1}, \alpha_\ell \rangle$ of type $D_2$. The $\Psi$-net tables are as follows.
$$
\begin{array}{|*5{>{\ss}c|}}
\hline
\multicolumn{2}{|>{\ss}c|}{\Psi\mathrm{-nets}} & & \multicolumn{1}{|>{\ss}c|}{c(s)} & \multicolumn{1}{|>{\ss}c|}{c(u_\Psi)} \\
\cline{1-2} \cline{4-5}
      \bar\nu      & n_1 &       m      &    r \geq 2    &    p \geq 2    \\
\hline
  \bom_1 + \bom_3  &  4  & 2^{\ell - 5} &  2^{\ell - 4}  &  2^{\ell - 4}  \\
       \bom_1      &  2  & 2^{\ell - 4} &  2^{\ell - 4}  &  2^{\ell - 4}  \\
       \bom_3      &  2  & 2^{\ell - 4} &  2^{\ell - 4}  &  2^{\ell - 4}  \\
         0         &  1  & 2^{\ell - 3} &                &                \\
\hline
\multicolumn{3}{c|}{}                   & 3.2^{\ell - 4} & 3.2^{\ell - 4} \\
\cline{4-5}
\end{array}
\quad
\begin{array}{|*5{>{\ss}c|}}
\hline
\multicolumn{2}{|>{\ss}c|}{\Psi\mathrm{-nets}} & & \multicolumn{1}{|>{\ss}c|}{c(s)} & \multicolumn{1}{|>{\ss}c|}{c(u_\Psi)} \\
\cline{1-2} \cline{4-5}
         \bar\nu        & n_1 &       m      &   r \geq 2   &   p \geq 2   \\
\hline
     \bom_{\ell - 1}    &  2  & 2^{\ell - 3} & 2^{\ell - 3} & 2^{\ell - 3} \\
        \bom_\ell       &  2  & 2^{\ell - 3} & 2^{\ell - 3} & 2^{\ell - 3} \\
\hline
\multicolumn{3}{c|}{}                        & 2^{\ell - 2} & 2^{\ell - 2} \\
\cline{4-5}
\end{array}
$$
Thus $\codim V_\kappa(s), \codim C_V(u_\Psi) > M_2 > 8\ell - 20 = \dim {u_\Psi}^G$ if $\Psi = {A_1}^2$, and $\codim V_\kappa(s), \codim C_V(u_\Psi) > M_3 > 4\ell - 4 = \dim {u_\Psi}^G$ if $\Psi = D_2$; taking the smaller of the two lower bounds, we see that we may assume from now on that $r \geq 3$, and that $p \geq 3$ when we treat unipotent classes. Moreover if $\ell = 10$ then $\codim V_\kappa(s), \codim C_V(u_\Psi) > M$ for either choice of $\Psi$; so the triple $(G, \lambda, p)$ satisfies $\ssdiamevcon$ and $\udiamcon$. We may therefore assume from now on that $\ell = 9$. We have $\codim V_\kappa(s), \codim C_V(u_\Psi) \geq 96$ if $\Psi = {A_1}^2$, and $\codim V_\kappa(s), \codim C_V(u_\Psi) \geq 128$ if $\Psi = D_2$, while $M = 144$. Again taking the smaller of the two lower bounds, we see that we need only consider semisimple classes $s^G$ with $|\Phi(s)| \leq 48$, and unipotent classes of dimension at least $96$; since each of the former has a subsystem of type $D_2$ disjoint from $\Phi(s)$, and each of the latter has the class $D_2$ in its closure by Lemma~\ref{lem: various classes in classical groups by dim}(ix), we may actually take the larger of the two lower bounds. We may therefore assume from now on that $r \geq 5$, and that $p \geq 5$ when we treat unipotent classes; we need only consider semisimple classes $s^G$ with $|\Phi(s)| \leq 16$, each of which by inspection has a subsystem $D_3$ disjoint from $\Phi(s)$, and unipotent classes of dimension at least $128$, each of which has the class $D_3$ in its closure by Lemma~\ref{lem: various classes in D_ell for fixed ell}(i).

Now take $\Psi = \langle \alpha_7, \alpha_8, \alpha_9 \rangle$ of type $D_3$. The $\Psi$-net table is as follows.
$$
\begin{array}{|*5{>{\ss}c|}}
\hline
\multicolumn{2}{|>{\ss}c|}{\Psi\mathrm{-nets}} & & \multicolumn{1}{|>{\ss}c|}{c(s)} & \multicolumn{1}{|>{\ss}c|}{c(u_\Psi)} \\
\cline{1-2} \cline{4-5}
    \bar\nu    & n_1 &  m & r \geq 5 & p \geq 5 \\
\hline
     \bom_8    &  4  & 32 &    96    &    96    \\
     \bom_9    &  4  & 32 &    96    &    96    \\
\hline
\multicolumn{3}{c|}{}     &   192    &   192    \\
\cline{4-5}
\end{array}
$$
Thus $\codim V_\kappa(s), \codim C_V(u_\Psi) > M$; so the triple $(G, \lambda, p)$ satisfies $\ssdiamevcon$ and $\udiamcon$.
\end{proof}

\begin{prop}\label{prop: B_ell, omega_ell, nets}
Let $G = B_\ell$ for $\ell \in [7, 9]$ and $\lambda = \omega_\ell$; then the triple $(G, \lambda, p)$ satisfies $\ssdiamevcon$ and $\udiamcon$.
\end{prop}

\begin{proof}
For $\ell \in [8, 9]$ this follows from Proposition~\ref{prop: D_ell, omega_ell, nets}, since $B_\ell$ is a subgroup of $D_{\ell + 1}$ and the spin module for $B_\ell$ is the restriction of the half-spin module for $D_{\ell + 1}$. We may therefore assume that $\ell = 7$.

The weight table is as follows.
$$
\begin{array}{|*4{>{\ss}c|}}
\hline
i & \mu & |W.\mu| & m_\mu \\
\hline
1 & \omega_7 & 128 & 1 \\
\hline
\end{array}
$$
We have $M = 98$, $M_5 = 84$, $M_3 = 70$ and $M_2 = 56$.

Take $\Psi = \langle \alpha_1 \rangle$ of type $A_1$, $\Psi = \langle \alpha_1, \alpha_3 \rangle$ of type ${A_1}^2$, and $\Psi = \langle \alpha_1, \alpha_3, \alpha_5 \rangle$ of type ${A_1}^3$. The $\Psi$-net tables are as follows.
$$
\begin{array}{|*4{>{\ss}c|}}
\hline
\multicolumn{2}{|>{\ss}c|}{\Psi\mathrm{-nets}} & & \multicolumn{1}{|>{\ss}c|}{c(u_\Psi)} \\
\cline{1-3} \cline{4-4}
      \bar\nu      & n_1 &  m & p \geq 2 \\
\hline
       \bom_1      &  2  & 32 &    32    \\
         0         &  1  & 64 &          \\
\hline
\multicolumn{3}{c|}{}         &    32    \\
\cline{4-4}
\end{array}
\quad
\begin{array}{|*4{>{\ss}c|}}
\hline
\multicolumn{2}{|>{\ss}c|}{\Psi\mathrm{-nets}} & & \multicolumn{1}{|>{\ss}c|}{c(u_\Psi)} \\
\cline{1-3} \cline{4-4}
        \bar\nu        & n_1 &  m & p \geq 2 \\
\hline
    \bom_1 + \bom_3    &  4  &  8 &    16    \\
         \bom_1        &  2  & 16 &    16    \\
         \bom_3        &  2  & 16 &    16    \\
           0           &  1  & 32 &          \\
\hline
\multicolumn{3}{c|}{}             &    48    \\
\cline{4-4}
\end{array}
\quad
\begin{array}{|*4{>{\ss}c|}}
\hline
\multicolumn{2}{|>{\ss}c|}{\Psi\mathrm{-nets}} & & \multicolumn{1}{|>{\ss}c|}{c(u_\Psi)} \\
\cline{1-3} \cline{4-4}
          \bar\nu          & n_1 &  m & p \geq 2 \\
\hline
  \bom_1 + \bom_3 + \bom_5 &  8  &  2 &     8    \\
      \bom_1 + \bom_3      &  4  &  4 &     8    \\
      \bom_1 + \bom_5      &  4  &  4 &     8    \\
      \bom_3 + \bom_5      &  4  &  4 &     8    \\
           \bom_1          &  2  &  8 &     8    \\
           \bom_3          &  2  &  8 &     8    \\
           \bom_5          &  2  &  8 &     8    \\
             0             &  1  & 16 &          \\
\hline
\multicolumn{3}{c|}{}                 &    56    \\
\cline{4-4}
\end{array}
$$
Thus $\codim C_V(u_\Psi) > 24 = \dim {u_\Psi}^G$ if $\Psi$ is of type $A_1$, $\codim C_V(u_\Psi) > 40 = \dim {u_\Psi}^G$ if $\Psi$ is of type ${A_1}^2$, and $\codim C_V(u_\Psi) > 48 = \dim {u_\Psi}^G$ if $\Psi$ is of type ${A_1}^3$. Each of the remaining non-trivial unipotent classes has the class $B_1$ in its closure by Lemma~\ref{lem: root elt class in closure of any non-triv class}.

Now take $\Psi = \langle \alpha_7 \rangle$ of type $B_1$. The $\Psi$-net table is as follows.
$$
\begin{array}{|*5{>{\ss}c|}}
\hline
\multicolumn{2}{|>{\ss}c|}{\Psi\mathrm{-nets}} & & \multicolumn{1}{|>{\ss}c|}{c(s)} & \multicolumn{1}{|>{\ss}c|}{c(u_\Psi)} \\
\cline{1-2} \cline{4-5}
    \bar\nu    & n_1 &  m & r \geq 2 & p \geq 2 \\
\hline
     \bom_7    &  2  & 64 &    64    &    64    \\
\hline
\multicolumn{3}{c|}{}     &    64    &    64    \\
\cline{4-5}
\end{array}
$$
Thus $\codim V_\kappa(s), \codim C_V(u_\Psi) \geq 64 > M_2 > 26 - 12\delta_{p, 2} = \dim {u_\Psi}^G$; we may therefore assume from now on that $r \geq 3$, and that $p \geq 3$ when we treat unipotent classes. We need only consider semisimple classes $s^G$ with $|\Phi(s)| \leq 34$, each of which by inspection has a subsystem of type $A_2B_1$ disjoint from $\Phi(s)$, and unipotent classes of dimension at least $64$, each of which has the class $A_2B_1$ in its closure by Lemma~\ref{lem: various classes in B_ell for fixed ell}(i).

Now take $\Psi = \langle \alpha_1, \alpha_2, \alpha_7 \rangle$ of type $A_2B_1$. The $\Psi$-net table is as follows.
$$
\begin{array}{|*5{>{\ss}c|}}
\hline
\multicolumn{2}{|>{\ss}c|}{\Psi\mathrm{-nets}} & & \multicolumn{1}{|>{\ss}c|}{c(s)} & \multicolumn{1}{|>{\ss}c|}{c(u_\Psi)} \\
\cline{1-2} \cline{4-5}
      \bar\nu      & n_1 &  m & r \geq 3 & p \geq 3 \\
\hline
  \bom_1 + \bom_7  &  6  &  8 &    32    &    32    \\
  \bom_2 + \bom_7  &  6  &  8 &    32    &    32    \\
       \bom_7      &  2  & 16 &    16    &    16    \\
\hline
\multicolumn{3}{c|}{}         &    80    &    80    \\
\cline{4-5}
\end{array}
$$
Thus $\codim V_\kappa(s), \codim C_V(u_\Psi) \geq 80 > M_3 > 60 = \dim {u_\Psi}^G$; we may therefore assume from now on that $r \geq 5$, and that $p \geq 5$ when we treat unipotent classes. We need only consider semisimple classes $s^G$ with $|\Phi(s)| \leq 18$, each of which by inspection has a subsystem of type $A_3B_1$ disjoint from $\Phi(s)$, and unipotent classes of dimension at least $80$, each of which has the class $A_3B_1$ in its closure by Lemma~\ref{lem: various classes in B_ell for fixed ell}(ii).

Now take $\Psi = \langle \alpha_1, \alpha_2, \alpha_3, \alpha_7 \rangle$ of type $A_3B_1$. The $\Psi$-net table is as follows.
$$
\begin{array}{|*5{>{\ss}c|}}
\hline
\multicolumn{2}{|>{\ss}c|}{\Psi\mathrm{-nets}} & & \multicolumn{1}{|>{\ss}c|}{c(s)} & \multicolumn{1}{|>{\ss}c|}{c(u_\Psi)} \\
\cline{1-2} \cline{4-5}
      \bar\nu      & n_1 & m & r \geq 5 & p \geq 5 \\
\hline
  \bom_1 + \bom_7  &  8  & 4 &    24    &    24    \\
  \bom_2 + \bom_7  & 12  & 4 &    32    &    36    \\
  \bom_3 + \bom_7  &  8  & 4 &    24    &    24    \\
       \bom_7      &  2  & 8 &     8    &     8    \\
\hline
\multicolumn{3}{c|}{}        &    88    &    92    \\
\cline{4-5}
\end{array}
$$
Thus $\codim V_\kappa(s) \geq 88 > M_5$, and $\codim C_V(u_\Psi) \geq 92 > M_5 > 72 = \dim {u_\Psi}^G$; we may therefore assume from now on that $r \geq 7$, and that $p \geq 7$ when we treat unipotent classes. We therefore need only consider semisimple classes $s^G$ with $|\Phi(s)| \leq 10$, each of which by inspection has a subsystem of type $A_4B_1$ disjoint from $\Phi(s)$, and unipotent classes of dimension at least $92$, each of which has the class $A_4B_1$ in its closure by Lemma~\ref{lem: various classes in B_ell for fixed ell}(iii).

Now take $\Psi = \langle \alpha_1, \alpha_2, \alpha_3, \alpha_4, \alpha_7 \rangle$ of type $A_4B_1$. The $\Psi$-net table is as follows.
$$
\begin{array}{|*5{>{\ss}c|}}
\hline
\multicolumn{2}{|>{\ss}c|}{\Psi\mathrm{-nets}} & & \multicolumn{1}{|>{\ss}c|}{c(s)} & \multicolumn{1}{|>{\ss}c|}{c(u_\Psi)} \\
\cline{1-2} \cline{4-5}
      \bar\nu      & n_1 & m & r \geq 7 & p \geq 7 \\
\hline
  \bom_1 + \bom_7  & 10  & 2 &    16    &    16    \\
  \bom_2 + \bom_7  & 20  & 2 &    32    &    32    \\
  \bom_3 + \bom_7  & 20  & 2 &    32    &    32    \\
  \bom_4 + \bom_7  & 10  & 2 &    16    &    16    \\
       \bom_7      &  2  & 4 &     4    &     4    \\
\hline
\multicolumn{3}{c|}{}        &   100    &   100    \\
\cline{4-5}
\end{array}
$$
Thus $\codim V_\kappa(s), \codim C_V(u_\Psi) > M$; so the triple $(G, \lambda, p)$ satisfies $\ssdiamevcon$ and $\udiamcon$.
\end{proof}

\begin{prop}\label{prop: C_ell, omega_ell, nets}
Let $G = C_\ell$ for $\ell \in [7, 9]$ and $\lambda = \omega_\ell$ with $p = 2$; then the triple $(G, \lambda, p)$ satisfies $\ssdiamevcon$ and $\udiamcon$.
\end{prop}

\begin{proof}
This is an immediate consequence of Proposition~\ref{prop: B_ell, omega_ell, nets}, using the exceptional isogeny $B_\ell \to C_\ell$ which exists in characteristic $2$.
\end{proof}

\begin{prop}\label{prop: C_ell, omega_3, nets}
Let $G = C_\ell$ for $\ell \in [5, 6]$ and $\lambda = \omega_3$; then the triple $(G, \lambda, p)$ satisfies $\ssdiamevcon$ and $\udiamcon$.
\end{prop}

\begin{proof}
Write $\z = \z_{p, \ell - 1}$ and $\z' = \z\z_{p, 2}$. The weight table is as follows.
$$
\begin{array}{|*4{>{\ss}c|}}
\hline
i & \mu & |W.\mu| & m_\mu \\
\hline
2 & \omega_3 & \frac{4}{3}\ell(\ell - 1)(\ell - 2) &       1       \\
1 & \omega_1 &                2\ell                & \ell - 2 - \z \\
\hline
\end{array}
$$
We have $M = 2\ell^2$ and $M_2 = \ell(\ell + 1)$.

Take $\Psi = \langle \alpha_\ell \rangle$ of type $C_1$. The $\Psi$-net table is as follows.
$$
\begin{array}{|*5{>{\ss}c|}}
\hline
\multicolumn{3}{|>{\ss}c|}{\Psi\mathrm{-nets}} & & \multicolumn{1}{|>{\ss}c|}{c(u_\Psi)} \\
\cline{1-3} \cline{5-5}
      \bar\nu      & n_1 & n_2 &                     m                     &          p \geq 2          \\
\hline
     \bom_\ell     &  0  &  2  &           2(\ell - 1)(\ell - 2)           &    2(\ell - 1)(\ell - 2)   \\
     \bom_\ell     &  2  &  0  &                     1                     &        \ell - 2 - \z       \\
         0         &  0  &  1  & \frac{4}{3}(\ell - 1)(\ell - 2)(\ell - 3) &                            \\
         0         &  1  &  0  &                2(\ell - 1)                &                            \\
\hline
\multicolumn{4}{c|}{}                                                      & (\ell - 2)(2\ell - 1) - \z \\
\cline{5-5}
\end{array}
$$
Thus $\codim C_V(u_\Psi) > 2\ell = \dim {u_\Psi}^G$. Each of the remaining non-trivial unipotent classes has the class $A_1$ in its closure by Lemma~\ref{lem: root elt class in closure of any non-triv class}.

Now take $\Psi = \langle \alpha_1 \rangle$ of type $A_1$. The $\Psi$-net table is as follows.
$$
\begin{array}{|*7{>{\ss}c|}}
\hline
\multicolumn{3}{|>{\ss}c|}{\Psi\mathrm{-nets}} & & \multicolumn{1}{|>{\ss}c|}{c(s)} & \multicolumn{2}{|>{\ss}c|}{c(u_\Psi)} \\
\cline{1-3} \cline{5-7}
    \bar\nu    & n_1 & n_2 &                      m                     &           r \geq 2           &         p = 2         &           p \geq 3           \\
\hline
    2\bom_1    &  1  &  2  &                 2(\ell - 2)                &          4(\ell - 2)         &      2(\ell - 2)      &          4(\ell - 2)         \\
     \bom_1    &  0  &  2  &            4(\ell - 2)(\ell - 3)           &     4(\ell - 2)(\ell - 3)    & 4(\ell - 2)(\ell - 3) &     4(\ell - 2)(\ell - 3)    \\
     \bom_1    &  2  &  0  &                      2                     &       2(\ell - 2) - 2\z      &   2(\ell - 2) - 2\z   &       2(\ell - 2) - 2\z      \\
       0       &  0  &  1  & \frac{4}{3}(\ell - 2)(\ell^2 - 7\ell + 15) &                              &                       &                              \\
\hline
\multicolumn{4}{c|}{}                                                   & 2(\ell - 2)(2\ell - 3) - 2\z &  4(\ell - 2)^2 - 2\z  & 2(\ell - 2)(2\ell - 3) - 2\z \\
\cline{5-7}
\end{array}
$$
Thus $\codim V_\kappa(s), \codim C_V(u_\Psi) > M_2 > 4\ell - 2 - 2\delta_{p, 2} = \dim {u_\Psi}^G$; we may therefore assume from now on that $r \geq 3$, and that $p \geq 3$ when we treat unipotent classes. We need only consider semisimple classes $s^G$ with $|\Phi(s)| < M - M_2$, each of which by inspection has a subsystem of type $A_2$ disjoint from $\Phi(s)$, and unipotent classes of dimension greater than $M_2$, each of which has the class $A_2$ in its closure by Lemma~\ref{lem: various classes in classical groups by dim}(viii).

Now take $\Psi = \langle \alpha_1, \alpha_2 \rangle$ of type $A_2$. The $\Psi$-net table is as follows.
$$
\begin{array}{|*7{>{\ss}c|}}
\hline
\multicolumn{3}{|>{\ss}c|}{\Psi\mathrm{-nets}} & & \multicolumn{1}{|>{\ss}c|}{c(s)} & \multicolumn{2}{|>{\ss}c|}{c(u_\Psi)} \\
\cline{1-3} \cline{5-7}
      \bar\nu      & n_1 & n_2 &                        m                       &            r \geq 3           &         p = 3         &       p \geq 5      \\
\hline
  \bom_1 + \bom_2  &  1  &  6  &                   2(\ell - 3)                  & 12(\ell - 3) - 2(\ell - 3)\z' &      8(\ell - 3)      &     12(\ell - 3)    \\
      2\bom_1      &  3  &  3  &                        1                       &       2(\ell - 1) - 2\z       &      2(\ell - 2)      &  2(\ell - 1) - 2\z  \\
      2\bom_2      &  3  &  3  &                        1                       &       2(\ell - 1) - 2\z       &      2(\ell - 2)      &  2(\ell - 1) - 2\z  \\
       \bom_1      &  0  &  3  &                  2(\ell - 3)^2                 &         4(\ell - 3)^2         &     4(\ell - 3)^2     &    4(\ell - 3)^2    \\
       \bom_2      &  0  &  3  &                  2(\ell - 3)^2                 &         4(\ell - 3)^2         &     4(\ell - 3)^2     &    4(\ell - 3)^2    \\
         0         &  0  &  1  & \frac{2}{3}(2\ell^3 - 24\ell^2 + 94\ell - 117) &                               &                       &                     \\
\hline
\multicolumn{4}{c|}{}                                                           &   8(\ell - 2)^2 - 4\z - 4\z'  & 8\ell^2 - 36\ell + 40 & 8(\ell - 2)^2 - 4\z \\
\cline{5-7}
\end{array}
$$
Thus $\codim V_\kappa(s), \codim C_V(u_\Psi) > M$; so the triple $(G, \lambda, p)$ satisfies $\ssdiamevcon$ and $\udiamcon$.
\end{proof}

\begin{prop}\label{prop: C_4, omega_3, nets}
Let $G = C_4$ and $\lambda = \omega_3$ with $p \neq 3$; then the triple $(G, \lambda, p)$ satisfies $\ssdiamevcon$ and $\udiamcon$.
\end{prop}

\begin{proof}
The weight table is as follows.
$$
\begin{array}{|*4{>{\ss}c|}}
\hline
i & \mu & |W.\mu| & m_\mu \\
\hline
2 & \omega_3 & 32 & 1 \\
1 & \omega_1 &  8 & 2 \\
\hline
\end{array}
$$
We have $M = 32$ and $M_2 = 20$.

Take $\Psi = \langle \alpha_4 \rangle$ of type $C_1$. The $\Psi$-net table is as follows.
$$
\begin{array}{|*6{>{\ss}c|}}
\hline
\multicolumn{3}{|>{\ss}c|}{\Psi\mathrm{-nets}} & & \multicolumn{2}{|>{\ss}c|}{c(u_\Psi)} \\
\cline{1-3} \cline{5-6}
      \bar\nu      & n_1 & n_2 &  m & p = 2 & p \geq 5 \\
\hline
       \bom_4      &  0  &  2  & 12 &  12   &    12    \\
       \bom_4      &  2  &  0  &  1 &   2   &     2    \\
         0         &  0  &  1  &  8 &       &          \\
         0         &  1  &  0  &  6 &       &          \\
\hline
\multicolumn{4}{c|}{}               &  14   &    14    \\
\cline{5-6}
\end{array}
$$
Thus $\codim C_V(u_\Psi) \geq 14 > 8 = \dim {u_\Psi}^G$. Each of the remaining non-trivial unipotent classes has the class $A_1$ in its closure by Lemma~\ref{lem: root elt class in closure of any non-triv class}.

Now take $\Psi = \langle \alpha_1 \rangle$ of type $A_1$. The $\Psi$-net table is as follows.
$$
\begin{array}{|*7{>{\ss}c|}}
\hline
\multicolumn{3}{|>{\ss}c|}{\Psi\mathrm{-nets}} & & \multicolumn{1}{|>{\ss}c|}{c(s)} & \multicolumn{2}{|>{\ss}c|}{c(u_\Psi)} \\
\cline{1-3} \cline{5-7}
    \bar\nu    & n_1 & n_2 & m & r \geq 2 & p = 2 & p \geq 5 \\
\hline
    2\bom_1    &  1  &  2  & 4 &     8    &   4   &     8    \\
     \bom_1    &  0  &  2  & 8 &     8    &   8   &     8    \\
     \bom_1    &  2  &  0  & 2 &     4    &   4   &     4    \\
       0       &  0  &  1  & 8 &          &       &          \\
\hline
\multicolumn{4}{c|}{}          &    20    &  16   &    20    \\
\cline{5-7}
\end{array}
$$
Thus $\codim V_\kappa(s) \geq 20$, and $\codim C_V(u_\Psi) \geq 20 > 14 = \dim {u_\Psi}^G$ unless $p = 2$, in which case $\codim C_V(u_\Psi) \geq 16 > 12 = \dim {u_\Psi}^G$. We therefore need only consider semisimple classes $s^G$ with $|\Phi(s)| \leq 12$, each of which by inspection has a subsystem of type ${A_1}^2$ disjoint from $\Phi(s)$, and unipotent classes of dimension at least $20$ if $p \geq 5$, each of which has the class ${A_1}^2$ or $C_2$ in its closure by Lemma~\ref{lem: various classes in C_ell for fixed ell}(iii), or at least $16$ if $p = 2$, each of which has the class $A_1C_1$ or ${A_1}^2$ in its closure by Lemma~\ref{lem: various classes in C_ell for fixed ell}(i).

Now take $\Psi = \langle \alpha_1, \alpha_4 \rangle$ of type $A_1C_1$ with $p = 2$. The $\Psi$-net table is as follows.
$$
\begin{array}{|*5{>{\ss}c|}}
\hline
\multicolumn{3}{|>{\ss}c|}{\Psi\mathrm{-nets}} & & \multicolumn{1}{|>{\ss}c|}{c(u_\Psi)} \\
\cline{1-3} \cline{5-5}
      \bar\nu      & n_1 & n_2 & m & p = 2 \\
\hline
  2\bom_1 + \bom_4 &  2  &  4  & 1 &   4   \\
  \bom_1 + \bom_4  &  0  &  4  & 4 &   8   \\
      2\bom_1      &  1  &  2  & 2 &   2   \\
       \bom_1      &  2  &  0  & 2 &   4   \\
       \bom_4      &  0  &  2  & 2 &   2   \\
         0         &  0  &  1  & 4 &       \\
\hline
\multicolumn{4}{c|}{}              &  20   \\
\cline{5-5}
\end{array}
$$
Thus $\codim C_V(u_\Psi) \geq 20 > 18 = \dim {u_\Psi}^G$. Each of the remaining unipotent classes requiring consideration with $p = 2$ has the class ${A_1}^2$ in its closure by Lemma~\ref{lem: various classes in C_ell for fixed ell}(ii).

Now take $\Psi = \langle \alpha_1, \alpha_3 \rangle$ of type ${A_1}^2$. The $\Psi$-net table is as follows.
$$
\begin{array}{|*8{>{\ss}c|}}
\hline
\multicolumn{3}{|>{\ss}c|}{\Psi\mathrm{-nets}} & & \multicolumn{2}{|>{\ss}c|}{c(s)} & \multicolumn{2}{|>{\ss}c|}{c(u_\Psi)} \\
\cline{1-3} \cline{5-8}
      \bar\nu      & n_1 & n_2 & m & r = 2 & r \geq 3 & p = 2 & p \geq 5 \\
\hline
 2\bom_1 + \bom_3  &  2  &  4  & 2 &   8   &    10    &   8   &    10    \\
 \bom_1 + 2\bom_3  &  2  &  4  & 2 &   8   &    10    &   8   &    10    \\
       \bom_1      &  0  &  2  & 4 &   4   &     4    &   4   &     4    \\
       \bom_3      &  0  &  2  & 4 &   4   &     4    &   4   &     4    \\
\hline
\multicolumn{4}{c|}{}              &  24   &    28    &  24   &    28    \\
\cline{5-8}
\end{array}
$$
Thus $\codim V_\kappa(s), \codim C_V(u_\Psi) > M_2 \geq 20 - 4\delta_{p, 2} = \dim {u_\Psi}^G$; we may therefore assume from now on that $r \geq 3$, and that $p \geq 5$ when treating unipotent classes. We need only consider semisimple classes $s^G$ with $|\Phi(s)| \leq 4$, each of which by inspection has a subsystem of type $C_2$ disjoint from $\Phi(s)$, and unipotent classes of dimension at least $28$, each of which has the class $C_2$ in its closure by Lemma~\ref{lem: various classes in C_ell for fixed ell}(iv).

Now take $\Psi = \langle \alpha_3, \alpha_4 \rangle$ of type $C_2$. The $\Psi$-net table is as follows.
$$
\begin{array}{|*6{>{\ss}c|}}
\hline
\multicolumn{3}{|>{\ss}c|}{\Psi\mathrm{-nets}} & & \multicolumn{1}{|>{\ss}c|}{c(s)} & \multicolumn{1}{|>{\ss}c|}{c(u_\Psi)} \\
\cline{1-3} \cline{5-6}
      \bar\nu      & n_1 & n_2 & m & r \geq 3 & p \geq 5 \\
\hline
       \bom_4      &  1  &  4  & 4 &    16    &    16    \\
       \bom_3      &  0  &  4  & 4 &    12    &    12    \\
       \bom_3      &  4  &  0  & 1 &     6    &     6    \\
\hline
\multicolumn{4}{c|}{}              &    34    &    34    \\
\cline{5-6}
\end{array}
$$
Thus $\codim V_\kappa(s), \codim C_V(u_\Psi) > M$; so the triple $(G, \lambda, p)$ satisfies $\ssdiamevcon$ and $\udiamcon$.
\end{proof}

\begin{prop}\label{prop: B_ell, omega_3, nets}
Let $G = B_\ell$ for $\ell \in [4, 6]$ and $\lambda = \omega_3$ with $p = 2$; then the triple $(G, \lambda, p)$ satisfies $\ssdiamevcon$ and $\udiamcon$.
\end{prop}

\begin{proof}
This is an immediate consequence of Propositions~\ref{prop: C_ell, omega_3, nets} and \ref{prop: C_4, omega_3, nets}, using the exceptional isogeny $B_\ell \to C_\ell$ which exists in characteristic $2$.
\end{proof}

This completes the treatment of the $p$-restricted large triples listed in Table~\ref{table: remaining triples} which do not appear in Table~\ref{table: large triple and first quadruple non-TGS}. As a consequence of this section and the preceding four we have proved the following.

\begin{prop}\label{prop: p-restricted large triples with TGS}
Any $p$-restricted large triple which is not listed in Table~\ref{table: large triple and first quadruple non-TGS} satisfies $\ssdiamevcon$ and $\udiamcon$, and so has TGS.
\end{prop}

\section{The triples $(C_4, \omega_3, 3)$ and $(B_2, \omega_1 + \omega_2, 5)$}\label{sect: two triples}

In this section we treat the two triples $(G, \lambda, p)$ which appear in Table~\ref{table: large triple and first quadruple non-TGS} but have TGS, namely $(C_4, \omega_3, 3)$ and $(B_2, \omega_1 + \omega_2, 5)$. The reason for handling them separately is that the $\Psi$-net analysis used in Section~\ref{sect: large triple further analysis} is insufficient for our purposes. Indeed in each case we shall be unable to show that $\ssdiamevcon$ is satisfied, and shall instead show that $(G, \lambda, p)$ satisfies $\ssdiamcon$ and $\udiamcon$; this is sufficient to show that the triple has TGS, but not the associated first quadruple. We shall proceed as follows. As before we begin with the weight table. We then list the possible subsystems $\Phi(s)$ corresponding to semisimple classes $s^G$, and the unipotent classes $u^G$ lying in $G_{(p)}$. Next we use $\Psi$-nets to dispose of all of the latter, and all but a few of the former; for each choice of $\Psi$ we shall of course only be able to dispose of subsystems $\Phi(s)$ which are disjoint from a conjugate of $\Psi$, while for unipotent classes we shall sometimes be forced to build appropriate representations and calculate using Jordan blocks to obtain strong enough lower bounds $c(u_\Psi)$. Finally we treat the remaining possibilities for $\Phi(s)$, arguing more closely using weights.

For this more detailed analysis of weights, we follow the approach of \cite{Ken}. Assume we have the semisimple element $s$ with corresponding subsystem $\Phi(s)$. We consider the equivalence relation on the set of weights defined by setting two weights to be equivalent if their difference is a sum of roots in $\Phi(s)$; the equivalence classes are called {\em clusters\/}. Thus if two weights lie in the same cluster then they must lie in the same eigenspace for $s$. One cluster {\em excludes\/} another if there exist two weights, one in the first cluster and one in the second, whose difference is a root (necessarily not in $\Phi(s)$). A {\em clique\/} is a set of clusters each of which excludes all of the others. All clusters in a clique must then lie in different eigenspaces.

\begin{prop}\label{prop: C_4, omega_3, p = 3, nets}
Let $G = C_4$ and $\lambda = \omega_3$ with $p = 3$; then the triple $(G, \lambda, p)$ satisfies $\ssdiamcon$ and $\udiamcon$.
\end{prop}

\begin{proof}
The weight table is as follows.
$$
\begin{array}{|*4{>{\ss}c|}}
\hline
i & \mu & |W.\mu| & m_\mu \\
\hline
2 & \omega_3 & 32 & 1 \\
1 & \omega_1 &  8 & 1 \\
\hline
\end{array}
$$
We have $M = 32$ and $M_2 = 20$. Since the order $r$ of $\bar s$ is prime, the possibilities for the subsystem $\Phi(s)$ are $C_3C_1$, ${C_2}^2$, $A_3$, $A_2C_1$, $A_1C_2$, $C_3$, $A_2$, ${A_1}^2$, $A_1C_1$, $C_2$, $A_1$, $C_1$ and $\emptyset$, of which only the first three apply if $r = 2$. By Lemma~\ref{lem: unipotent classes in the two triples}(i) the unipotent classes lying in $G_{(p)}$ are $C_1 \leq A_1 \leq A_1 C_1 \leq {A_1}^2 \leq A_2 \leq A_2 C_1$, with the dimensions being $8$, $14$, $18$, $20$, $22$ and $24$ respectively.

Take $\Psi = \langle \alpha_4 \rangle$ of type $C_1$. The $\Psi$-net table is as follows.
$$
\begin{array}{|*5{>{\ss}c|}}
\hline
\multicolumn{3}{|>{\ss}c|}{\Psi\mathrm{-nets}} & & \multicolumn{1}{|>{\ss}c|}{c(u_\Psi)} \\
\cline{1-3} \cline{5-5}
      \bar\nu      & n_1 & n_2 &  m & p = 3 \\
\hline
       \bom_4      &  0  &  2  & 12 &  12   \\
       \bom_4      &  2  &  0  &  1 &   1   \\
         0         &  0  &  1  &  8 &       \\
         0         &  1  &  0  &  6 &       \\
\hline
\multicolumn{4}{c|}{}               &  13   \\
\cline{5-5}
\end{array}
$$
Thus $\codim C_V(u_\Psi) \geq 13 > 8 = \dim {u_\Psi}^G$, which disposes of the unipotent class $C_1$.

Now take $\Psi = \langle \alpha_1 \rangle$ of type $A_1$. The $\Psi$-net table is as follows.
$$
\begin{array}{|*7{>{\ss}c|}}
\hline
\multicolumn{3}{|>{\ss}c|}{\Psi\mathrm{-nets}} & & \multicolumn{2}{|>{\ss}c|}{c(s)} & \multicolumn{1}{|>{\ss}c|}{c(u_\Psi)} \\
\cline{1-3} \cline{5-7}
    \bar\nu    & n_1 & n_2 & m & r = 2 & r \geq 5 & p = 3 \\
\hline
    2\bom_1    &  1  &  2  & 4 &   4   &     8    &   8   \\
     \bom_1    &  0  &  2  & 8 &   8   &     8    &   8   \\
     \bom_1    &  2  &  0  & 2 &   2   &     2    &   2   \\
       0       &  0  &  1  & 8 &       &          &       \\
\hline
\multicolumn{4}{c|}{}          &  14   &    18    &  18   \\
\cline{5-7}
\end{array}
$$
Thus $\codim V_\kappa(s) \geq 18$ unless $r = 2$, in which case $\codim V_\kappa(s) \geq 14$; this disposes of the possibilities $\Phi(s) = C_3C_1$ and $C_3$. Moreover $\codim C_V(u_\Psi) \geq 18 > 14 = \dim {u_\Psi}^G$, which disposes of the unipotent class $A_1$.

Now take $\Psi = \langle \alpha_1, \alpha_4 \rangle$ of type $A_1C_1$. The $\Psi$-net table is as follows.
$$
\begin{array}{|*5{>{\ss}c|}}
\hline
\multicolumn{3}{|>{\ss}c|}{\Psi\mathrm{-nets}} & & \multicolumn{1}{|>{\ss}c|}{c(u_\Psi)} \\
\cline{1-3} \cline{5-5}
      \bar\nu      & n_1 & n_2 & m & p = 3 \\
\hline
  2\bom_1 + \bom_4 &  2  &  4  & 1 &   4   \\
  \bom_1 + \bom_4  &  0  &  4  & 4 &   8   \\
      2\bom_1      &  1  &  2  & 2 &   4   \\
       \bom_1      &  2  &  0  & 2 &   2   \\
       \bom_4      &  0  &  2  & 2 &   2   \\
         0         &  0  &  1  & 4 &       \\
\hline
\multicolumn{4}{c|}{}              &  20   \\
\cline{5-5}
\end{array}
$$
Thus $\codim C_V(u_\Psi) \geq 20 > 18 = \dim {u_\Psi}^G$, which disposes of the unipotent class $A_1C_1$.

Now take $\Psi = \langle \alpha_1, \alpha_3 \rangle$ of type ${A_1}^2$. The $\Psi$-net table is as follows.
$$
\begin{array}{|*7{>{\ss}c|}}
\hline
\multicolumn{3}{|>{\ss}c|}{\Psi\mathrm{-nets}} & & \multicolumn{2}{|>{\ss}c|}{c(s)} & \multicolumn{1}{|>{\ss}c|}{c(u_\Psi)} \\
\cline{1-3} \cline{5-7}
      \bar\nu      & n_1 & n_2 & m & r = 2 & r \geq 5 & p = 3 \\
\hline
 2\bom_1 + \bom_3  &  2  &  4  & 2 &   6   &     8    &   8   \\
 \bom_1 + 2\bom_3  &  2  &  4  & 2 &   6   &     8    &   8   \\
       \bom_1      &  0  &  2  & 4 &   4   &     4    &   4   \\
       \bom_3      &  0  &  2  & 4 &   4   &     4    &   4   \\
\hline
\multicolumn{4}{c|}{}              &  20   &    24    &  24   \\
\cline{5-7}
\end{array}
$$
Thus $\codim V_\kappa(s) \geq 24$ unless $r = 2$, in which case $\codim V_\kappa(s) \geq 20$; this disposes of the possibilities $\Phi(s) = {C_2}^2$ and $A_1C_2$. Moreover $\codim C_V(u_\Psi) \geq 24 > 20 = \dim {u_\Psi}^G$, which disposes of the unipotent classes ${A_1}^2$ and $A_2$.

Now take $\Psi = \langle \alpha_3, \alpha_4 \rangle$ of type $C_2$. The $\Psi$-net table is as follows.
$$
\begin{array}{|*6{>{\ss}c|}}
\hline
\multicolumn{3}{|>{\ss}c|}{\Psi\mathrm{-nets}} & & \multicolumn{2}{|>{\ss}c|}{c(s)} \\
\cline{1-3} \cline{5-6}
      \bar\nu      & n_1 & n_2 & m & r = 2 & r \geq 5 \\
\hline
       \bom_4      &  1  &  4  & 4 &   12  &    16    \\
       \bom_3      &  0  &  4  & 4 &   12  &    12    \\
       \bom_3      &  4  &  0  & 1 &    3  &     3    \\
\hline
\multicolumn{4}{c|}{}              &   27  &    31    \\
\cline{5-6}
\end{array}
$$
Thus $\codim V_\kappa(s) \geq 31$ unless $r = 2$, in which case $\codim V_\kappa(s) > M_r$; this disposes of the possibilities $\Phi(s) = A_2$, ${A_1}^2$, $A_1C_1$, $C_2$, $A_1$ and $C_1$.

Now take $\Psi = \langle \alpha_1, \alpha_2, \alpha_4 \rangle$ of type $A_2C_1$. The $\Psi$-net table is as follows.
$$
\begin{array}{|*5{>{\ss}c|}}
\hline
\multicolumn{3}{|>{\ss}c|}{\Psi\mathrm{-nets}} & & \multicolumn{1}{|>{\ss}c|}{c(u_\Psi)} \\
\cline{1-3} \cline{5-5}
          \bar\nu          & n_1 & n_2 & m & p = 3 \\
\hline
  \bom_1 + \bom_2 + \bom_4 &  2  & 12  & 1 &   9   \\
       \bom_1 + \bom_4     &  0  &  6  & 1 &   4   \\
       \bom_2 + \bom_4     &  0  &  6  & 1 &   4   \\
           2\bom_1         &  3  &  3  & 1 &   4   \\
           2\bom_2         &  3  &  3  & 1 &   4   \\
              0            &  0  &  1  & 2 &       \\
\hline
\multicolumn{4}{c|}{}                      &  25   \\
\cline{5-5}
\end{array}
$$
The values $c(u_\Psi)$ in the first, fourth and fifth rows are obtained as follows. For the fourth, we have the Weyl $G_\Psi$-module with highest weight $2\bom_1$. We may take root elements corresponding to roots $\alpha_1$ and $\alpha_2$ to act on it as
$$
\left(
  \begin{array}{cccccc}
    1 & 2 & 1 &   &   &   \\
      & 1 & 1 &   &   &   \\
      &   & 1 &   &   &   \\
      &   &   & 1 & 1 &   \\
      &   &   &   & 1 &   \\
      &   &   &   &   & 1 \\
  \end{array}
\right)
\qquad \hbox{and} \qquad
\left(
  \begin{array}{cccccc}
    1 &   &   &   &   &   \\
      & 1 &   & 1 &   &   \\
      &   & 1 &   & 2 & 1 \\
      &   &   & 1 &   &   \\
      &   &   &   & 1 & 1 \\
      &   &   &   &   & 1 \\
  \end{array}
\right)
$$
respectively; we may then take the product of these to represent $u_\Psi$, since the $C_1$ factor of $G_\Psi$ acts trivially, and subtracting $I$ leaves a matrix of rank $4$, so we may take $c(u_\Psi) = 4$. Likewise we have $c(u_\Psi) = 4$ in the fifth row. Finally, for the first we have the Weyl $G_\Psi$-module with highest weight $\bom_1 + \bom_2 + \bom_4$. We may proceed similarly to obtain a matrix
$$
\left(
  \begin{array}{ccccccc}
    1 & 2 & 1 & 2 & 1 &   &   \\
      & 1 &   & 1 & 2 &   &   \\
      &   & 1 & 1 & 1 & 2 & 1 \\
      &   &   & 1 & 1 & 1 & 2 \\
      &   &   &   & 1 &   & 1 \\
      &   &   &   &   & 1 & 2 \\
      &   &   &   &   &   & 1 \\
  \end{array}
\right)
\otimes
\left(
  \begin{array}{cc}
    1 & 1 \\
      & 1 \\
  \end{array}
\right)
$$
which we may take to represent the action of $u_\Psi$; subtracting $I$ leaves a matrix of rank $9$, so we may take $c(u_\Psi) = 9$. Thus $\codim C_V(u_\Psi) \geq 25 > 24 = \dim {u_\Psi}^G$, which disposes of the unipotent class $A_2C_1$. We have now disposed of all the unipotent classes lying in $G_{(p)}$; so the triple $(G, \lambda, p)$ satisfies $\udiamcon$.

At this point we are left with just the possibilities $\Phi(s) = A_3$, $A_2C_1$ and $\emptyset$ to consider. We shall treat each of these in turn. We first explain the notation we shall use for the weights here. We shall use the standard notation for roots in $\Phi$ as in Section~\ref{sect: notation}. The weights are integer linear combinations of the $\ve_i$. We shall represent $a_1\ve_1 + a_2\ve_2 + a_3\ve_3 + a_4\ve_4$ as $a_1a_2a_3a_4$; then $\omega_3 = 1110$ and $\omega_1 = 1000$, and the weights of the form $\mu_2$ are strings of three $\pm 1$s and one $0$, while those of the form $\mu_1$ are strings of one $\pm 1$ and three $0$s. For convenience we write $\bar 1$ for $-1$. In addition, if we enclose part of a string in brackets it means that all possible permutations of the symbols inside are to be taken; thus for example $(1 0 0) 0$ stands for the three weights $1 0 0 0$, $0 1 0 0$ and $0 0 1 0$.

We start with $\Phi(s) = \langle \alpha_1, \alpha_2, \alpha_3 \rangle$ of type $A_3$; then $\dim s^G = 20$, and we have $(\ve_1 - \ve_2)(s) = (\ve_2 - \ve_3)(s) = (\ve_3 - \ve_4)(s) = 1$, so $\ve_1(s) = \ve_2(s) = \ve_3(s) = \ve_4(s)$. The clusters are as follows:
$$
\{ (1 1 1 0) \}, \{ (1 1 0 \bar 1), (1 0 0 0) \}, \{ (\bar 1 \bar 1 0 1), (\bar 1 0 0 0) \}, \{ (\bar 1 \bar 1 \bar 1 0) \}.
$$
The second and third clusters are of size $16$ and the first and fourth are of size $4$. If neither cluster of size $16$ is in $V_\kappa(s)$ then $\codim V_\kappa(s) \geq 32 > \dim s^G$, so we may assume (without loss of generality) the second cluster is in $V_\kappa(s)$; this excludes the first and third clusters, so $\codim V_\kappa(s) \geq 4 + 16 = 20 = \dim s^G$. For equality we must have the second and fourth clusters in $V_\kappa(s)$; then $\ve_1(s) = (-\ve_1 - \ve_2 - \ve_3)(s) = \kappa$, so $(4\ve_1)(s) = 1$, and as $2\ve_1 \notin \Phi(s)$ we must have $(2\ve_1)(s) \neq 1$, so $(2\ve_1)(s) = -1$ and hence $\kappa$ is a square root of $-1$. Thus we do have $\codim C_V(s) > \dim s^G$ here; but if we set $s = h_{\alpha_1}(\eta_4) h_{\alpha_2}(-1) h_{\alpha_3}(-\eta_4)$, then $s$ has eigenvalues $\eta_4$ and $-\eta_4$ on $V$, and for $\kappa \in \{ \pm\eta_4 \}$ we have $\codim V_\kappa(s) = 20 = \dim s^G$.

Next we take $\Phi(s) = \langle \alpha_1, \alpha_2, \alpha_4 \rangle$ of type $A_2C_1$; then $\dim s^G = 24$, and we have $(\ve_1 - \ve_2)(s) = (\ve_2 - \ve_3)(s) = (2\ve_4)(s) = 1$, so $\ve_1(s) = \ve_2(s) = \ve_3(s)$ and $\ve_4(s) = \pm 1$. The clusters are as follows:
\begin{eqnarray*}
& \{ (1 0 \bar 1) 1, (1 0 \bar 1) \bar 1, 0 0 0 1, 0 0 0 \bar 1 \}; & \\
& \{ (1 1 \bar 1) 0, (1 0 0) 0 \}, \{ (1 1 0) 1, (1 1 0) \bar 1 \}, \{ 1 1 1 0 \}; & \\
& \{ (\bar 1 \bar 1 1) 0, (\bar 1 0 0) 0 \}, \{ (\bar 1 \bar 1 0) 1, (\bar 1 \bar 1 0) \bar 1 \}, \{ \bar 1 \bar 1 \bar 1 0 \}. &
\end{eqnarray*}
The clusters within each row form a clique; the first clique is a single cluster of size $14$, while the second and third cliques each comprise three clusters, of sizes $6$, $6$ and $1$. The contribution to $c(s)$ from each of the second and third cliques is therefore at least $6 + 1 = 7$. Thus if the first cluster is not in $V_\kappa(s)$ then $\codim V_\kappa(s) \geq 14 + 7 + 7 = 28 > \dim s^G$, so we may assume the first cluster is in $V_\kappa(s)$; this excludes the four clusters of size $6$, so $\codim V_\kappa(s) \geq 4.6 = 24 = \dim s^G$. For equality we must have the first cluster and the two of size $1$ in $V_\kappa(s)$; then $(\ve_1 + \ve_2 + \ve_3)(s) = \ve_4(s) = \kappa$, so $\kappa = \pm 1$, and $3\ve_1(s) = \kappa$, whence $\ve_1(s) = \kappa$, but then $(\ve_1 - \ve_4)(s) = 1$ contrary to $\ve_1 - \ve_4 \notin \Phi(s)$. Thus equality is impossible, and we have $\codim V_\kappa(s) > \dim s^G$ here.

Finally we take $\Phi(s) = \emptyset$; then $\dim s^G = 32$. All clusters are single weights. First suppose some weight of the form $\mu_1$ is in $V_\kappa(s)$; using the Weyl group $W$ we may assume $1 0 0 0$ is in $V_\kappa(s)$, which excludes all weights except those of the form $\mu_2$ with first coefficient either $0$ or $\bar 1$. If some weight with first coefficient $0$ is in $V_\kappa(s)$, using the stabilizer in $W$ of $1 0 0 0$ we may assume $0 1 1 1$ is in $V_\kappa(s)$, which excludes all but the following four cliques:
$$
\{ \bar 1 1 \bar 1 0, \bar 1 1 0 \bar 1, 0 1 \bar 1 \bar 1 \};
\{ \bar 1 \bar 1 1 0, \bar 1 0 1 \bar 1, 0 \bar 1 1 \bar 1 \};
\{ \bar 1 \bar 1 0 1, \bar 1 0 \bar 1 1, 0 \bar 1 \bar 1 1 \};
\{ \bar 1 \bar 1 \bar 1 0, \bar 1 \bar 1 0 \bar 1, \bar 1 0 \bar 1 \bar 1, 0 \bar 1 \bar 1 \bar 1 \}.
$$
Thus at most $6$ weights can lie in $V_\kappa(s)$, so $\codim V_\kappa(s) \geq 34 > \dim s^G$. We may therefore assume no weight with first coefficient $0$ is in $V_\kappa(s)$; but then the remaining $12$ weights form the following six cliques:
$$
\{ \bar 1 1 1 0, \bar 1 1 \bar 1 0 \}; \{ \bar 1 \bar 1 1 0, \bar 1 \bar 1 \bar 1 0 \}; \{ \bar 1 1 0 1, \bar 1 1 0 \bar 1 \}; \{ \bar 1 \bar 1 0 1, \bar 1 \bar 1 0 \bar 1 \}; \{ \bar 1 0 1 1, \bar 1 0 1 \bar 1 \}; \{ \bar 1 0 \bar 1 1, \bar 1 0 \bar 1 \bar 1 \}.
$$
Thus $\codim V_\kappa(s) \geq 33 > \dim s^G$. We may therefore assume that no weight of the form $\mu_1$ is in $V_\kappa(s)$.

Using $W$ we may then assume $1 1 1 0$ is in $V_\kappa(s)$, which excludes all but the following seven cliques:
\begin{eqnarray*}
& \{ \bar 1 \bar 1 \bar 1 0, \bar 1 \bar 1 0 \bar 1, \bar 1 0 \bar 1 \bar 1, 0 \bar 1 \bar 1 \bar 1 \}; & \\
& \{ \bar 1 \bar 1 1 0, \bar 1 \bar 1 0 1, \bar 1 0 1 1, 0 \bar 1 1 1 \};
\{ \bar 1 1 \bar 1 0, \bar 1 1 0 1, \bar 1 0 \bar 1 1, 0 1 \bar 1 1 \};
\{ 1 \bar 1 \bar 1 0, 1 \bar 1 0 1, 1 0 \bar 1 1, 0 \bar 1 \bar 1 1 \}; & \\
& \{ 1 0 \bar 1 \bar 1, 0 1 \bar 1 \bar 1 \};
\{ 1 \bar 1 0 \bar 1, 0 \bar 1 1 \bar 1 \};
\{ \bar 1 1 0 \bar 1, \bar 1 0 1 \bar 1 \}. &
\end{eqnarray*}
Thus $\codim V_\kappa(s) \geq 32 = \dim s^G$. For equality we must have exactly one weight from each clique in $V_\kappa(s)$. Since the weight $0 \bar 1 \bar 1 \bar 1$ from the first clique excludes both weights in the fifth clique, it then cannot lie in $V_\kappa(s)$; using the stabilizer in $W$ of $1 1 1 0$, we see that neither $\bar 1 0 \bar 1 \bar 1$ nor $\bar 1 \bar 1 0 \bar 1$ can lie in $V_\kappa(s)$. Thus from the first clique we must have $\bar 1 \bar 1 \bar 1 0$ in $V_\kappa(s)$, which excludes $\bar 1 \bar 1 1 0$ and $\bar 1 \bar 1 0 1$ from the second clique, $\bar 1 1 \bar 1 0$ and $\bar 1 0 \bar 1 1$ from the third, and $1 \bar 1 \bar 1 0$ and $0 \bar 1 \bar 1 1$ from the fourth; so we are left with six cliques of size $2$. Again using the stabilizer in $W$ of $1 1 1 0$ we may assume from the fifth clique we have $1 0 \bar 1 \bar 1$ in $V_\kappa(s)$; this excludes $1 0 \bar 1 1$ from the fourth and $1 \bar 1 0 \bar 1$ from the sixth, giving $1 \bar 1 0 1$ and $0 \bar 1 1 \bar 1$ in $V_\kappa(s)$, which exclude $0 \bar 1 1 1$ from the second and $\bar 1 0 1 \bar 1$ from the seventh, giving $\bar 1 0 1 1$ and $\bar 1 1 0 \bar 1$ in $V_\kappa(s)$, which exclude $\bar 1 1 0 1$ from the third, giving $0 1 \bar 1 1$ in $V_\kappa(s)$. Thus the weights in $V_\kappa(s)$ are $1 1 1 0$, $1 \bar 1 0 1$, $0 1 \bar 1 1$, $\bar 1 0 1 1$ and their negatives; then $(\ve_1 + \ve_2 + \ve_3)(s) = (\ve_1 - \ve_2 + \ve_4)(s) = (\ve_2 - \ve_3 + \ve_4)(s) = (-\ve_1 + \ve_3 + \ve_4)(s) = \kappa$, and so $(2\ve_2 + \ve_3 - \ve_4)(s) = (\ve_1 + 2\ve_3 - \ve_4)(s) = (2\ve_1 + \ve_2 - \ve_4)(s) = 1$, whence $(-\ve_1 + 2\ve_2 - \ve_3)(s) = (2\ve_1 - \ve_2 - \ve_3)(s) = 1$, and so $(3\ve_1 - 3\ve_2)(s) = 1$, which forces $(\ve_1 - \ve_2)(s) = 1$, contrary to $\ve_1 - \ve_2 \notin \Phi(s)$. Thus equality is impossible, and we have $\codim V_\kappa(s) > \dim s^G$ here.

Therefore the triple $(G, \lambda, p)$ satisfies $\ssdiamcon$, but not $\ssdiamevcon$.
\end{proof}

\begin{prop}\label{prop: B_2, omega_1 + omega_2, p = 5, nets}
Let $G = B_2$ and $\lambda = \omega_1 + \omega_2$ with $p = 5$; then the triple $(G, \lambda, p)$ satisfies $\ssdiamcon$ and $\udiamcon$.
\end{prop}

\begin{proof}
The weight table is as follows.
$$
\begin{array}{|*4{>{\ss}c|}}
\hline
i & \mu & |W.\mu| & m_\mu \\
\hline
2 & \omega_1 + \omega_2 & 8 & 1 \\
1 &      \omega_2       & 4 & 1 \\
\hline
\end{array}
$$

\pagebreak

\noindent We have $M = 8$ and $M_2 = 6$. The possibilities for the subsystem $\Phi(s)$ are ${A_1}^2$, $B_1$, $A_1$ and $\emptyset$, of which only the first two apply if $r = 2$. By Lemma~\ref{lem: unipotent classes in the two triples}(ii) the unipotent classes lying in $G_{(p)}$ are $A_1 \leq B_1 \leq B_2$, with the dimensions being $4$, $6$ and $8$ respectively.

Take $\Psi = \langle \alpha_1 \rangle$ of type $A_1$. The $\Psi$-net table is as follows.
$$
\begin{array}{|*5{>{\ss}c|}}
\hline
\multicolumn{3}{|>{\ss}c|}{\Psi\mathrm{-nets}} & & \multicolumn{1}{|>{\ss}c|}{c(u_\Psi)} \\
\cline{1-3} \cline{5-5}
      \bar\nu      & n_1 & n_2 & m & p = 5 \\
\hline
      2\bom_1      &  1  &  2  & 2 &   4   \\
       \bom_1      &  0  &  2  & 2 &   2   \\
       \bom_1      &  2  &  0  & 1 &   1   \\
\hline
\multicolumn{4}{c|}{}              &   7   \\
\cline{5-5}
\end{array}
$$
Thus $\codim C_V(u_\Psi) \geq 7 > 4 = \dim {u_\Psi}^G$, which disposes of the unipotent class $A_1$.

Now take $\Psi = \langle \alpha_2 \rangle$ of type $B_1$. The $\Psi$-net table is as follows.
$$
\begin{array}{|*8{>{\ss}c|}}
\hline
\multicolumn{3}{|>{\ss}c|}{\Psi\mathrm{-nets}} & & \multicolumn{3}{|>{\ss}c|}{c(s)} & \multicolumn{1}{|>{\ss}c|}{c(u_\Psi)} \\
\cline{1-3} \cline{5-8}
    \bar\nu    & n_1 & n_2 & m & r = 2 & r = 3 & r \geq 7 & p = 5 \\
\hline
    3\bom_2    &  2  &  2  & 2 &   4   &   4   &     6    &   6   \\
     \bom_2    &  0  &  2  & 2 &   2   &   2   &     2    &   2   \\
\hline
\multicolumn{4}{c|}{}          &   6   &   6   &     8    &   8   \\
\cline{5-8}
\end{array}
$$
Thus $\codim V_\kappa(s) \geq 8$ unless $r = 2$ or $r = 3$, in which case $\codim V_\kappa(s) \geq 6$; this disposes of the possibility $\Phi(s) = {A_1}^2$. Moreover $\codim C_V(u_\Psi) \geq 8 > 6 = \dim {u_\Psi}^G$, which disposes of the unipotent class $B_1$.

Now take $\Psi = \langle \alpha_1, \alpha_2 \rangle$ of type $B_2$. The $\Psi$-net table is as follows.
$$
\begin{array}{|*5{>{\ss}c|}}
\hline
\multicolumn{3}{|>{\ss}c|}{\Psi\mathrm{-nets}} & & \multicolumn{1}{|>{\ss}c|}{c(u_\Psi)} \\
\cline{1-3} \cline{5-5}
      \bar\nu      & n_1 & n_2 & m & p = 5 \\
\hline
  \bom_1 + \bom_2  &  4  &  8  & 1 &   9   \\
\hline
\multicolumn{4}{c|}{}              &   9   \\
\cline{5-5}
\end{array}
$$
Here we take $x_{\alpha_1}(1)$ and $x_{\alpha_2}(1)$ to act on $V$ as
$$
\left(
  \begin{array}{cccccccccccc}
    1 & 1 &   &   &   &   &   &   &   &   &   &   \\
      & 1 &   &   &   &   &   &   &   &   &   &   \\
      &   & 1 & 4 & 1 &   &   &   &   &   &   &   \\
      &   &   & 1 & 3 &   &   &   &   &   &   &   \\
      &   &   &   & 1 &   &   &   &   &   &   &   \\
      &   &   &   &   & 1 & 1 &   &   &   &   &   \\
      &   &   &   &   &   & 1 &   &   &   &   &   \\
      &   &   &   &   &   &   & 1 & 1 & 4 &   &   \\
      &   &   &   &   &   &   &   & 1 & 3 &   &   \\
      &   &   &   &   &   &   &   &   & 1 &   &   \\
      &   &   &   &   &   &   &   &   &   & 1 & 4 \\
      &   &   &   &   &   &   &   &   &   &   & 1 \\
  \end{array}
\right)
\quad\hbox{and}\quad
\left(
  \begin{array}{cccccccccccc}
    1 &   & 1 &   &   &   &   &   &   &   &   &   \\
      & 1 &   & 2 &   & 2 &   & 4 &   &   &   &   \\
      &   & 1 &   &   &   &   &   &   &   &   &   \\
      &   &   & 1 &   & 2 &   & 1 &   &   &   &   \\
      &   &   &   & 1 &   & 3 &   & 2 &   & 1 &   \\
      &   &   &   &   & 1 &   & 1 &   &   &   &   \\
      &   &   &   &   &   & 1 &   & 3 &   & 1 &   \\
      &   &   &   &   &   &   & 1 &   &   &   &   \\
      &   &   &   &   &   &   &   & 1 &   & 4 &   \\
      &   &   &   &   &   &   &   &   & 1 &   & 4 \\
      &   &   &   &   &   &   &   &   &   & 1 &   \\
      &   &   &   &   &   &   &   &   &   &   & 1 \\
  \end{array}
\right)
$$
respectively; taking the product of the two matrices and subtracting $I$ leaves a matrix of rank $9$, so we may take $c(u_\Psi) = 9$. Thus $\codim C_V(u_\Psi) \geq 9 > 8 = \dim {u_\Psi}^G$, which disposes of the unipotent class $B_2$. We have now disposed of all the unipotent classes lying in $G_{(p)}$; so the triple $(G, \lambda, p)$ satisfies $\udiamcon$.

At this point we are left with just the possibilities $\Phi(s) = B_1$, $A_1$ and $\emptyset$ to consider. We shall treat each of these in turn. We first explain the notation we shall use for the weights here. We shall use the standard notation for roots in $\Phi$ as in Section~\ref{sect: notation}. The weights are half-integer linear combinations of the $\ve_i$. We shall represent $\frac{1}{2}(a_1\ve_1 + a_2\ve_2)$ as $a_1a_2$; then $\omega_1 + \omega_2 = 31$ and $\omega_2 = 11$, and the weights of the form $\mu_2$ are strings of one $\pm 3$ and one $\pm 1$, while those of the form $\mu_1$ are strings of two $\pm 1$s. For convenience we write $\bar 1$ for $-1$ and $\bar 3$ for $-3$.

We start with $\Phi(s) = \langle \alpha_2 \rangle$ of type $B_1$; then $\dim s^G = 6$, and we have $\ve_2(s) = 1$. The clusters are as follows:
$$
\{ 3 1, 3 \bar 1 \}, \{ 1 3, 1 1, 1 \bar 1, 1 \bar 3 \}, \{ \bar 1 3, \bar 1 1, \bar 1 \bar 1, \bar 1 \bar 3 \}, \{ \bar 3 1, \bar 3 \bar 1 \}.
$$
If neither cluster of size $4$ is in $V_\kappa(s)$ we have $\codim V_\kappa(s) \geq 8 > \dim s^G$, so we may assume (without loss of generality) the second cluster is in $V_\kappa(s)$; this excludes the first and third clusters, so $\codim V_\kappa(s) \geq 2 + 4 = 6 = \dim s^G$. For equality we must have the second and fourth clusters in $V_\kappa(s)$; then $(\frac{1}{2}(\ve_1 + \ve_2))(s) = (\frac{1}{2}(-3\ve_1 + \ve_2))(s) = \kappa$, so $(2\ve_1)(s) = 1$, and as $\ve_1 \notin \Phi(s)$ we must have $\ve_1(s) \neq 1$, so $\ve_1(s) = -1$ and hence $\kappa$ is a square root of $-1$. Thus we do have $\codim C_V(s) > \dim s^G$ here; but if we set $s = h_{\alpha_1}(-1) h_{\alpha_2}(\eta_4)$, then $s$ has eigenvalues $\eta_4$ and $-\eta_4$ on $V$, and for $\kappa \in \{ \pm\eta_4 \}$ we have $\codim V_\kappa(s) = 6 = \dim s^G$.

Next we take $\Phi(s) = \langle \alpha_1 \rangle$ of type $A_1$; then $\dim s^G = 6$, and we have $(\ve_1 - \ve_2)(s) = 1$, so $\ve_1(s) = \ve_2(s)$. The clusters are as follows:
$$
\{ 3 1, 1 3 \}, \{ 3 \bar 1, 1 1, \bar 1 3 \}, \{ 1 \bar 1, \bar 1 1 \}, \{ 1 \bar 3, \bar 1 \bar 1, \bar 3 1 \}, \{ \bar 1 \bar 3, \bar 3 \bar 1 \}.
$$
Each excludes its neighbours, and as $r \neq 2$ we cannot have either the two clusters of size $3$ or the three of size $2$ in $V_\kappa(s)$; thus we have $\codim V_\kappa(s) \! \geq \! 7 \! > \! \dim s^G$ here.

Finally take $\Phi(s) = \emptyset$; then $\dim s^G = 8$. All clusters are single weights; again $r \neq 2$. If $V_\kappa(s)$ contains a weight of the form $\mu_1$, we may assume it contains $1 1$; this excludes all the other weights except $1 \bar 3$, $\bar 1 \bar 3$, $\bar 3 1$ and $\bar 3 \bar 1$, of which the first two and the second two form two cliques, so $\codim V_\kappa(s) \geq 9 > \dim s^G$. Thus we may assume $V_\kappa(s)$ contains no weight of the form $\mu_1$; we may then assume it contains $3 1$, which excludes all other weights of the form $\mu_2$ except $\bar 1 3$, $1 \bar 3$, $\bar 1 \bar 3$, $\bar 3 1$ and $\bar 3 \bar 1$, of which the second and third form a clique, as do the fourth and fifth, so $\codim V_\kappa(s) \geq 8 = \dim s^G$. For equality we must have $\bar 1 3$ in $V_\kappa(s)$, which excludes $\bar 3 1$, so we must have $\bar 3 \bar 1$ in $V_\kappa(s)$, which excludes $\bar 1 \bar 3$, so we must have $1 \bar 3$ in $V_\kappa(s)$. Thus the weights in $V_\kappa(s)$ are $3 1$, $\bar 1 3$ and their negatives; then $(\pm\frac{1}{2}(3\ve_1 + \ve_2))(s) = (\pm\frac{1}{2}(\ve_1 - 3\ve_2))(s) = \kappa$, so $(3\ve_1 + \ve_2)(s) = (2\ve_1 - \ve_2)(s) = 1$, whence $(5\ve_1)(s) = 1$, which forces $\ve_1(s) = 1$, contrary to $\ve_1 \notin \Phi(s)$. Thus equality is impossible, and we have $\codim V_\kappa(s) > \dim s^G$ here.

Therefore the triple $(G, \lambda, p)$ satisfies $\ssdiamcon$, but not $\ssdiamevcon$.
\end{proof}

We have therefore shown that the two triples $(C_4, \omega_3, 3)$ and $(B_2, \omega_1 + \omega_2, 5)$ have TGS; as a result the only $p$-restricted large triples which have not yet been treated are those listed in Table~\ref{table: large triple and first quadruple non-TGS} as not having TGS. In the final section of this chapter we turn to large triples which are not $p$-restricted.

\section{Tensor products}\label{sect: large triple tensor products}

Let $(G, \lambda, p)$ be a large triple; as usual write $V = L(\lambda)$. In this section we assume that $\lambda$ is not $p$-restricted (so that in particular $p$ is finite); thus by Theorem~\ref{thm: Steinberg} we have $V = V_1 \otimes V_2$ with $\dim V_1, \ \dim V_2 > 1$. We shall show that if $(G, \lambda, p)$ is not listed in Table~\ref{table: large triple and first quadruple non-TGS} then it satisfies $\ssdiamevcon$ and $\udiamcon$, and thus has TGS.

We begin with a couple of elementary results which between them imply that if either $\dim V_1 > M$ or $\dim V_2 > M$ then the triple $(G, \lambda, p)$ satisfies $\ssddagcon$ and $\uddagcon$. The first, concerning semisimple elements, is essentially \cite[Proposition~4.3]{Ken}.

\begin{lem}\label{lem: ss codim bound for tensor products}
Let $G$ be an algebraic group, $s \in G$ be semisimple, and $V = V_1 \otimes V_2$ be a $G$-module; then if each eigenspace of $s$ on $V_1$ has codimension at least $c$, then for all $\kappa \in K^*$ we have $\codim V_\kappa(s) \geq c.\dim V_2$.
\end{lem}

\begin{proof}
Write
$$
V_1 = \bigoplus_{\kappa_1 \in K^*} (V_1)_{\kappa_1}(s), \qquad V_2 = \bigoplus_{\kappa_2 \in K^*} (V_2)_{\kappa_2}(s);
$$
then $V_\kappa(s)$ is the sum of the spaces $(V_1)_{\kappa_1}(s) \otimes (V_2)_{\kappa_2}(s)$ with $\kappa_1\kappa_2 = \kappa$. For each $\kappa_2 \in K^*$ we have $\dim (V_1)_{\kappa{\kappa_2}^{-1}}(s) \leq \dim V_1 - c$, so
$$
\dim ((V_1)_{\kappa{\kappa_2}^{-1}}(s) \otimes (V_2)_{\kappa_2}(s)) \leq (\dim V_1 - c).\dim (V_2)_{\kappa_2}(s);
$$
summing over $\kappa_2$ gives $\dim V_\kappa(s) \leq (\dim V_1 - c).\dim V_2$ as required.
\end{proof}

The second is an analogous result concerning unipotent elements.

\begin{lem}\label{lem: unip codim bound for tensor products}
Let $G$ be an algebraic group, $u \in G$ be unipotent, and $V = V_1 \otimes V_2$ be a $G$-module; then $\codim C_{V_1 \otimes V_2}(u) \geq \codim C_{V_1}(u).\dim V_2$.
\end{lem}

\begin{proof}
Let $u$ have Jordan block sizes $r_1^1, \dots, r_1^{d_1}$ on $V_1$ and $r_2^1, \dots, r_2^{d_2}$ on $V_2$, where $d_1 = \dim C_{V_1}(u)$ and $d_2 = \dim C_{V_2}(u)$. Then summing over the various Jordan blocks in $V_1 \otimes V_2$, by Lemma~\ref{lem: Jordan block tensor product} we have
$$
\dim C_{V_1 \otimes V_2}(u) = \sum_{i = 1}^{d_1} \sum_{j = 1}^{d_2} \min(r_1^i, r_2^j) \leq \sum_{i = 1}^{d_1} \sum_{j = 1}^{d_2} r_2^j = d_1 \dim V_2;
$$
the result follows.
\end{proof}

Thus if either $\dim V_1 > M$ or $\dim V_2 > M$ then the triple $(G, \lambda, p)$ satisfies both $\ssddagcon$ and $\uddagcon$. We therefore need only consider large triples $(G, \lambda, p)$ with both $\dim V_1 \leq M$ and $\dim V_2 \leq M$.

\begin{table}[ht]
\caption{Modules $L(\lambda')$ of dimension at most $M$}\label{table: small modules}
\tabcapsp
$$
\begin{array}{|c|c|c|c|c|c|c|c|c|c|}
\cline{1-5} \cline{7-10}
G      & \lambda'     & \ell        & p          & \dim L(\lambda')                 & \ptw & G   & \lambda' & p          & \dim L(\lambda') \tbs \\
\cline{1-5} \cline{7-10}
A_\ell & \omega_1     & {} \geq 1   & \hbox{any} & \ell + 1                         &      & E_6 & \omega_1 & \hbox{any} & 27               \tbs \\
\cline{7-10}
       & 2\omega_1    & {} \geq 2   & {} \geq 3  & \frac{1}{2}(\ell + 1)(\ell + 2)  &      & E_7 & \omega_7 & \hbox{any} & 56               \tbs \\
\cline{7-10}
       & \omega_2     & {} \geq 3   & \hbox{any} & \frac{1}{2}\ell(\ell + 1)        &      & F_4 & \omega_4 & \hbox{any} & 26 - \z_{p, 3}   \tbs \\
       & \omega_3     & 5, 6, 7     & \hbox{any} & \frac{1}{6}\ell(\ell^2 - 1)      &      &     & \omega_1 & 2          & 26               \tbs \\
\cline{1-5} \cline{7-10}
B_\ell & \omega_1     & {} \geq 2   & \hbox{any} & 2\ell + 1 - \z_{p, 2}            &      & G_2 & \omega_1 & \hbox{any} & 7 - \z_{p, 2}    \tbs \\
       & \omega_2     & {} \geq 3   & 2          & 2\ell^2 - \ell - 1 -\z_{2, \ell} &      &     & \omega_2 & 3          & 7                \tbs \\
\cline{7-10}
       & \omega_\ell  & 2, \dots, 6 & \hbox{any} & 2^\ell                           & \multicolumn{5}{c}{}                                  \tbs \\
\cline{1-5}
C_\ell & \omega_1     & {} \geq 3   & \hbox{any} & 2\ell                            & \multicolumn{5}{c}{}                                  \tbs \\
       & \omega_2     & {} \geq 3   & \hbox{any} & 2\ell^2 - \ell - 1 -\z_{p, \ell} & \multicolumn{5}{c}{}                                  \tbs \\
       & \omega_3     & 3           & {} \geq 3  & 14                               & \multicolumn{5}{c}{}                                  \tbs \\
       & \omega_\ell  & 3, \dots, 6 & 2          & 2^\ell                           & \multicolumn{5}{c}{}                                  \tbs \\
\cline{1-5}
D_\ell & \omega_1     & {} \geq 4   & \hbox{any} & 2\ell                            & \multicolumn{5}{c}{}                                  \tbs \\
       & \omega_\ell  & 5, 6, 7     & \hbox{any} & 2^{\ell - 1}                     & \multicolumn{5}{c}{}                                  \tbs \\
\cline{1-5}
\end{array}
$$
\end{table}

Reference to \cite{Lubpaper} shows that, up to graph automorphisms, the irreducible modules $L(\lambda')$ with $p$-restricted $\lambda'$ which are of dimension at most $M$ are those listed in Table~\ref{table: small modules}. In particular we see that no non-trivial tensor product has dimension at most $M$. Thus we need only consider triples $(G, \lambda, p)$ with
$$
\lambda = \lambda_1 + q\lambda_2,
$$
where both $\lambda_1$ and $\lambda_2$ are $p$-restricted and $q = p^i$ for some $i \geq 1$; then we have
$$
V = V_1 \otimes V_2
$$
with $V_1 = L(\lambda_1)$ and $V_2 = L(\lambda_2)^{(i)}$.

We shall proceed as follows. As in Section~\ref{sect: large triple weight string analysis}, we let $s \in T$ be an element of $G_{(r)}$ for some $r \in \P'$, and $\kappa$ be an element of $K^*$; we take $\alpha \in \Phi_s$ with $\alpha(s) \neq 1$, and write $u_\alpha = x_\alpha(1)$; if $e(\Phi) > 1$, we take $\beta \in \Phi_l$, and write $u_\beta = x_\beta(1)$. For each group $G$, we take the various modules $V' = L(\lambda')$ listed in Table~\ref{table: small modules}; we again provide tables enabling us to compute lower bounds $c(s)$ for $\codim (V')_\kappa(s)$ and $c(u_\alpha)$ for $\codim C_{V'}(u_\alpha)$ (and $c(u_\beta)$ for $\codim C_{V'}(u_\beta)$ if $e(\Phi) > 1$). We then consider the possible pairs of weights $(\lambda_1, \lambda_2)$ (as usual, working up to graph automorphisms); note that the order of the weights within a pair is immaterial. In most cases Lemmas~\ref{lem: ss codim bound for tensor products} and \ref{lem: unip codim bound for tensor products} immediately show that the triple $(G, \lambda, p)$ satisfies both $\uddagcon$ and $\ssddagcon$; in some cases we can obtain an improved lower bound for $\codim C_{V'}(u_\alpha)$ or $\codim C_{V'}(u_\beta)$ by applying Lemma~\ref{lem: Jordan block tensor product} as in the proof of Lemma~\ref{lem: unip codim bound for tensor products}. We then investigate further the few remaining cases.

For the classical groups $B_\ell$, $C_\ell$ and $D_\ell$, we postpone until the end of this section the consideration of the case where $(\lambda_1, \lambda_2) = (\omega_1, \omega_1)$; note that for the group $A_\ell$ both this case and that where $(\lambda_1, \lambda_2) = (\omega_1, \omega_\ell)$ are not being considered here, since they are listed in Table~\ref{table: large triple and first quadruple non-TGS}. Since we work modulo graph automorphisms, when we come to treat the group $D_4$ the assumption that the pair $(\lambda_1, \lambda_2)$ is not $(\omega_1, \omega_1)$ also rules out the pairs $(\omega_3, \omega_3)$ and $(\omega_4, \omega_4)$.

\begin{prop}\label{prop: A_ell, tensor products}
Let $G = A_\ell$; suppose $(\lambda_1, \lambda_2) \neq (\omega_1, \omega_1)$ or $(\omega_1, \omega_\ell)$. Then if $\ell = 3$ and $(\lambda_1, \lambda_2) = (\omega_2, \omega_2)$, or if $\ell \in [3, \infty)$ and $(\lambda_1, \lambda_2) = (\omega_2, \omega_1)$ or $(\omega_2, \omega_\ell)$, the triple $(G, \lambda, p)$ satisfies $\ssdiamevcon$ and $\udiamcon$; if $\ell \in [2, \infty)$ and $(\lambda_1, \lambda_2) = (2\omega_1, \omega_1)$ or $(2\omega_1, \omega_\ell)$ with $p \geq 3$, the triple $(G, \lambda, p)$ satisfies $\ssdagcon$ and $\uddagcon$; in all other cases the triple $(G, \lambda, p)$ satisfies $\ssddagcon$ and $\uddagcon$.
\end{prop}

\begin{proof}
First suppose $\ell \in [1, \infty)$ and $\lambda' = \omega_1$. In this case the tables are as follows.
$$
\begin{array}{|*4{>{\ss}c|}}
\hline
i & \mu & |W.\mu| & m_\mu \\
\hline
1 & \omega_1 & \ell + 1 & 1 \\
\hline
\end{array}
\quad
\begin{array}{|*4{>{\ss}c|}}
\hline
 & & \multicolumn{1}{|>{\ss}c|}{c(s)} & \multicolumn{1}{|>{\ss}c|}{c(u_\alpha)} \\
\cline{3-4}
 \ss{\alpha\mathrm{-strings}} &     m    & r \geq 2 & p \geq 2 \\
\hline
            \mu_1             & \ell - 1 &          &          \\
        \mu_1 \ \mu_1         &     1    &     1    &     1    \\
\hline
\multicolumn{2}{c|}{}                    &     1    &     1    \\
\cline{3-4}
\end{array}
$$
Next suppose $\ell \in [2, \infty)$ and $\lambda' = 2\omega_1$ with $p \geq 3$. In this case the tables are as follows.
$$
\begin{array}{|*4{>{\ss}c|}}
\hline
i & \mu & |W.\mu| & m_\mu \\
\hline
2 & 2\omega_1 &          \ell + 1         & 1 \\
1 &  \omega_2 & \frac{1}{2}\ell(\ell + 1) & 1 \\
\hline
\end{array}
\quad
\begin{array}{|*5{>{\ss}c|}}
\hline
 & & \multicolumn{2}{|>{\ss}c|}{c(s)} & \multicolumn{1}{|>{\ss}c|}{c(u_\alpha)} \\
\cline{3-5}
 \ss{\alpha\mathrm{-strings}} &                m                &   r = 2  & r \geq 3 & p \geq 3 \\
\hline
            \mu_2             &             \ell - 1            &          &          &          \\
    \mu_2 \ \mu_1 \ \mu_2     &                1                &     1    &     2    &     2    \\
            \mu_1             & \frac{1}{2}(\ell - 1)(\ell - 2) &          &          &          \\
        \mu_1 \ \mu_1         &             \ell - 1            & \ell - 1 & \ell - 1 & \ell - 1 \\
\hline
\multicolumn{2}{c|}{}                                           &   \ell   & \ell + 1 & \ell + 1 \\
\cline{3-5}
\end{array}
$$
Next suppose $\ell \in [3, \infty)$ and $\lambda' = \omega_2$. In this case the tables are as follows.
$$
\begin{array}{|*4{>{\ss}c|}}
\hline
i & \mu & |W.\mu| & m_\mu \\
\hline
1 & \omega_2 & \frac{1}{2}\ell(\ell + 1) & 1 \\
\hline
\end{array}
\quad
\begin{array}{|*4{>{\ss}c|}}
\hline
 & & \multicolumn{1}{|>{\ss}c|}{c(s)} & \multicolumn{1}{|>{\ss}c|}{c(u_\alpha)} \\
\cline{3-4}
 \ss{\alpha\mathrm{-strings}} &                m                & r \geq 2 & p \geq 2 \\
\hline
            \mu_1             & \frac{1}{2}(\ell^2 - 3\ell + 4) &          &          \\
        \mu_1 \ \mu_1         &             \ell - 1            & \ell - 1 & \ell - 1 \\
\hline
\multicolumn{2}{c|}{}                                           & \ell - 1 & \ell - 1 \\
\cline{3-4}
\end{array}
$$
Finally suppose $\ell \in [5, 7]$ and $\lambda' = \omega_3$. In this case the tables are as follows.
$$
\begin{array}{|*4{>{\ss}c|}}
\hline
i & \mu & |W.\mu| & m_\mu \\
\hline
1 & \omega_3 & \frac{1}{6}\ell(\ell^2 - 1) & 1 \\
\hline
\end{array}
\quad
\begin{array}{|*4{>{\ss}c|}}
\hline
 & & \multicolumn{1}{|>{\ss}c|}{c(s)} & \multicolumn{1}{|>{\ss}c|}{c(u_\alpha)} \\
\cline{3-4}
 \ss{\alpha\mathrm{-strings}} &                      m                     &             r \geq 2            &             p \geq 2            \\
\hline
            \mu_1             & \frac{1}{6}(\ell - 1)(\ell^2 - 5\ell + 12) &                                 &                                 \\
        \mu_1 \ \mu_1         &       \frac{1}{2}(\ell - 1)(\ell - 2)      & \frac{1}{2}(\ell - 1)(\ell - 2) & \frac{1}{2}(\ell - 1)(\ell - 2) \\
\hline
\multicolumn{2}{c|}{}                                                      & \frac{1}{2}(\ell - 1)(\ell - 2) & \frac{1}{2}(\ell - 1)(\ell - 2) \\
\cline{3-4}
\end{array}
$$
We have $M = \ell(\ell + 1)$ and $M_2 = \lfloor \frac{1}{2}(\ell + 1)^2\rfloor$.

Now if $\ell \in [2, \infty)$ and $(\lambda_1, \lambda_2) = (2\omega_1, 2\omega_1)$ or $(2\omega_1, 2\omega_\ell)$ with $p \geq 3$, then $\codim V_\kappa(s) \geq \frac{1}{2}\ell(\ell + 1)(\ell + 2) > M$ and $\codim C_V(u_\alpha) \geq \frac{1}{2}(\ell + 1)^2(\ell + 2) > M$. If $\ell \in [3, \infty)$ and $(\lambda_1, \lambda_2) = (\omega_2, 2\omega_1)$ or $(\omega_2, 2\omega_\ell)$ with $p \geq 3$, then $\codim V_\kappa(s)$, $\codim C_V(u_\alpha) \geq \frac{1}{2}(\ell^2 - 1)(\ell + 2) > M$. If $\ell \in [4, \infty)$ and $(\lambda_1, \lambda_2) = (\omega_2, \omega_2)$ or $(\omega_2, \omega_{\ell - 1})$, then $\codim V_\kappa(s), \ \codim C_V(u_\alpha) \geq \frac{1}{2}\ell(\ell^2 - 1) > M$. If $\ell \in [5, 7]$ and $(\lambda_1, \lambda_2) = (\omega_3, \omega_1)$ or $(\omega_3, \omega_\ell)$, then $\codim V_\kappa(s), \ \codim C_V(u_\alpha) \geq \frac{1}{2}(\ell^2 - 1)(\ell - 2) > M$. If $\ell \in [5, 7]$ and $(\lambda_1, \lambda_2) = (\omega_3, 2\omega_1)$ or $(\omega_3, 2\omega_\ell)$ with $p \geq 3$, then $\codim V_\kappa(s), \ \codim C_V(u_\alpha) \geq \frac{1}{4}(\ell^2 - 1)(\ell^2 - 4) > M$. If $\ell \in [5, 7]$ and $(\lambda_1, \lambda_2) = (\omega_3, \omega_2)$ or $(\omega_3, \omega_{\ell - 1})$, then $\codim V_\kappa(s), \ \codim C_V(u_\alpha) \geq \frac{1}{4}\ell(\ell^2 - 1)(\ell - 2) > M$. If $\ell \in [5, 7]$ and $(\lambda_1, \lambda_2) = (\omega_3, \omega_3)$ or $(\omega_3, \omega_{\ell - 2})$, then $\codim V_\kappa(s), \ \codim C_V(u_\alpha) \geq \frac{1}{12}\ell(\ell^2 - 1)(\ell - 1)(\ell - 2) > M$. Thus in these cases the triple $(G, \lambda, p)$ satisfies $\ssddagcon$ and $\uddagcon$.

Next if $\ell \in [2, \infty)$ and $(\lambda_1, \lambda_2) = (2\omega_1, \omega_1)$ or $(2\omega_1, \omega_\ell)$ with $p \geq 3$, then $\codim C_V(u_\alpha) \geq (\ell + 1)^2 > M$, and $\codim V_\kappa(s) \geq (\ell + 1)^2 > M$ unless $r = 2$, in which case $\codim V_\kappa(s) \geq \ell(\ell + 1) > M_2$. Thus in these cases the triple $(G, \lambda, p)$ satisfies $\ssdagcon$ and $\uddagcon$.

Next if $\ell = 3$ and $(\lambda_1, \lambda_2) = (\omega_2, \omega_2)$, then $\codim V_\kappa(s), \ \codim C_V(u_\alpha) \geq 12 = M$; thus we need only consider regular classes. If $s \in G_{(r)}$ is regular, then for each weight $\mu$ in $W.\omega_2$ there are $4$ other weights $\mu'$ with $\mu - \mu' \in \Phi$ (e.g. if $\mu = \omega_2$ then we may take $\mu' = \mu - \alpha$ for $\alpha \in \{ \alpha_2, \alpha_1 + \alpha_2, \alpha_2 + \alpha_3, \alpha_1 + \alpha_2 + \alpha_3 \}$); thus $\codim (V_1)_\kappa(s) \geq 4$ and hence $\codim V_\kappa(s) \geq 4 \dim V_2 = 24 > M$. If $u \in G_{(p)}$ is regular, then in the natural representation of $G$ with basis $v_1, v_2, v_3, v_4$ we may assume $u$ maps $v_1 \mapsto v_1$ and $v_i \mapsto v_{i - 1} + v_i$ for $i = 2, 3, 4$; then writing $v_{ij}$ for $v_i \wedge v_j$, on $L(\omega_2)$ we see that $u$ maps
$$
\begin{array}{l}
v_{12} \mapsto v_{12}, \\
v_{13} \mapsto v_{12} + v_{13}, \\
v_{14} \mapsto v_{13} + v_{14}, \\
v_{23} \mapsto v_{12} + v_{13} + v_{23}, \\
v_{24} \mapsto v_{13} + v_{14} + v_{23} + v_{24}, \\
v_{34} \mapsto v_{23} + v_{24} + v_{34}, \\
\end{array}
\qquad\qquad \hbox{and so }
u =
\left(
\begin{array}{cccccc}
 1 & 1 &   & 1 &   &   \\
   & 1 & 1 & 1 & 1 &   \\
   &   & 1 &   & 1 &   \\
   &   &   & 1 & 1 & 1 \\
   &   &   &   & 1 & 1 \\
   &   &   &   &   & 1 \\
\end{array}
\right),
$$
so $\codim C_{V_1}(u) = 4$ and hence $\codim C_V(u) \geq 4 \dim V_2 = 24 > M$. Thus in this case the triple $(G, \lambda, p)$ satisfies $\ssdiamevcon$ and $\udiamcon$.

Finally if $\ell \in [3, \infty)$ and $(\lambda_1, \lambda_2) = (\omega_2, \omega_1)$ or $(\omega_2, \omega_\ell)$, then $\codim V_\kappa(s)$, $\codim C_V(u_\alpha) \geq \ell^2 - 1$; thus we need only consider classes of dimension at least $\ell^2 - 1$. Take $s \in G_{(r)}$. If its centralizer is $A_{\ell - 1}$ then $\dim s^G = 2\ell < \ell^2 - 1$; if not then $\codim (V_2)_\kappa(s) \geq 2$, so $\codim V_\kappa(s) \geq 2 \dim V_1 = \ell(\ell + 1) = M$, so we need only consider regular semisimple classes; if $s$ is regular then any eigenspace in $V_2$ has codimension at least $\ell$, so $\codim V_\kappa(s) \geq \ell \dim V_1 = \frac{1}{2}\ell^2(\ell + 1) > M$. Now take $u \in G_{(p)}$. If $u$ is a root element then $\dim u^G = 2\ell < \ell^2 - 1$; if not then $\codim C_{V_2}(u) \geq 2$, so $\codim C_V(u) \geq 2 \dim V_1 = \ell(\ell + 1) = M$, so we need only consider regular unipotent elements; if $u$ is regular then $\codim C_{V_2}(u) = \ell$, so $\codim C_V(u) \geq \ell \dim V_1 = \frac{1}{2}\ell^2(\ell + 1) > M$. Thus in these cases the triple $(G, \lambda, p)$ satisfies $\ssdiamevcon$ and $\udiamcon$.
\end{proof}

\begin{prop}\label{prop: B_ell, tensor products}
Let $G = B_\ell$; suppose $(\lambda_1, \lambda_2) \neq (\omega_1, \omega_1)$. Then if $\ell = 2$ and $(\lambda_1, \lambda_2) = (\omega_2, \omega_1)$ or $(\omega_2, \omega_2)$, the triple $(G, \lambda, p)$ satisfies $\ssdiamevcon$ and $\udiamcon$; in all other cases the triple $(G, \lambda, p)$ satisfies $\ssddagcon$ and $\uddagcon$.
\end{prop}

\begin{proof}
First suppose $\ell \in [2, \infty)$ and $\lambda' = \omega_1$; write $\z = \z_{p, 2}$. In this case the tables are as follows.
$$
\begin{array}{|*4{>{\ss}c|}}
\hline
i & \mu & |W.\mu| & m_\mu \\
\hline
1 & \omega_1 & 2\ell &    1   \\
0 &    0     &   1   & 1 - \z \\
\hline
\end{array}
\quad
\begin{array}{|*6{>{\ss}c|}}
\hline
 & & \multicolumn{2}{|>{\ss}c|}{c(s)} & \multicolumn{2}{|>{\ss}c|}{c(u_\alpha)} \\
\cline{3-6}
 \ss{\alpha\mathrm{-strings}} &      m      & r = 2 & r \geq 3 & p = 2 & p \geq 3 \\
\hline
            \mu_1             & 2(\ell - 1) &       &          &       &          \\
    \mu_1 \ \mu_0 \ \mu_1     &      1      &   1   &  2 - \z  &   1   &     2    \\
\hline
\multicolumn{2}{c|}{}                       &   1   &  2 - \z  &   1   &     2    \\
\cline{3-6}
\end{array}
\quad
\begin{array}{|*3{>{\ss}c|}}
\hline
 & & \multicolumn{1}{|>{\ss}c|}{c(u_\beta)} \\
\cline{3-3}
 \ss{\beta\mathrm{-strings}} &      m      & p \geq 2 \\
\hline
            \mu_1             & 2(\ell - 2) &          \\
        \mu_1 \ \mu_1         &      2      &     2    \\
            \mu_0             &      1      &          \\
\hline
\multicolumn{2}{c|}{}                       &     2    \\
\cline{3-3}
\end{array}
$$
Next suppose $\ell \in [3, \infty)$ and $\lambda' = \omega_2$ with $p = 2$; write $\z = \z_{2, \ell}$. In this case the tables are as follows.
$$
\begin{array}{|*4{>{\ss}c|}}
\hline
i & \mu & |W.\mu| & m_\mu \\
\hline
1 & \omega_2 & 2\ell(\ell - 1) &       1       \\
0 &    0     &         1       & \ell - 1 - \z \\
\hline
\end{array}
$$
$$
\begin{array}{|*4{>{\ss}c|}}
\hline
 & & \multicolumn{1}{|>{\ss}c|}{c(s)} & \multicolumn{1}{|>{\ss}c|}{c(u_\alpha)} \\
\cline{3-4}
 \ss{\alpha\mathrm{-strings}} &           m           &   r \geq 3  &    p = 2    \\
\hline
            \mu_1             & 2(\ell - 1)(\ell - 2) &             &             \\
    \mu_1 \ \phmu \ \mu_1     &      2(\ell - 1)      & 2(\ell - 1) & 2(\ell - 1) \\
            \mu_0             &           1           &             &             \\
\hline
\multicolumn{2}{c|}{}                                 & 2(\ell - 1) & 2(\ell - 1) \\
\cline{3-4}
\end{array}
\quad
\begin{array}{|*3{>{\ss}c|}}
\hline
 & & \multicolumn{1}{|>{\ss}c|}{c(u_\beta)} \\
\cline{3-3}
 \ss{\beta\mathrm{-strings}} &           m           &    p = 2    \\
\hline
            \mu_1             & 2\ell^2 - 10\ell + 14 &             \\
        \mu_1 \ \mu_1         &      4(\ell - 2)      & 4(\ell - 2) \\
    \mu_1 \ \mu_0 \ \mu_1     &           1           &      1      \\
\hline
\multicolumn{2}{c|}{}                                 &  4\ell - 7  \\
\cline{3-3}
\end{array}
$$
Finally suppose $\ell \in [2, 6]$ and $\lambda' = \omega_\ell$. In this case the tables are as follows.
$$
\begin{array}{|*4{>{\ss}c|}}
\hline
i & \mu & |W.\mu| & m_\mu \\
\hline
1 & \omega_\ell & 2^\ell & 1 \\
\hline
\end{array}
\quad
\begin{array}{|*4{>{\ss}c|}}
\hline
 & & \multicolumn{1}{|>{\ss}c|}{c(s)} & \multicolumn{1}{|>{\ss}c|}{c(u_\alpha)} \\
\cline{3-4}
 \ss{\alpha\mathrm{-strings}} &       m      &   r \geq 2   &   p \geq 2   \\
\hline
        \mu_1 \ \mu_1         & 2^{\ell - 1} & 2^{\ell - 1} & 2^{\ell - 1} \\
\hline
\multicolumn{2}{c|}{}                        & 2^{\ell - 1} & 2^{\ell - 1} \\
\cline{3-4}
\end{array}
\quad
\begin{array}{|*3{>{\ss}c|}}
\hline
 & & \multicolumn{1}{|>{\ss}c|}{c(u_\beta)} \\
\cline{3-3}
 \ss{\beta\mathrm{-strings}} &       m      &    p \geq 2   \\
\hline
            \mu_1             & 2^{\ell - 1} &               \\
        \mu_1 \ \mu_1         & 2^{\ell - 2} &  2^{\ell - 2} \\
\hline
\multicolumn{2}{c|}{}                        &  2^{\ell - 2} \\
\cline{3-3}
\end{array}
$$
We have $M = 2\ell^2$.

Now if $\ell \in [3, \infty)$ and $(\lambda_1, \lambda_2) = (\omega_2, \omega_1)$ with $p = 2$, then $\codim V_\kappa(s)$, $\codim C_V(u_\alpha) \geq 4\ell(\ell - 1) > M$, and $\codim C_V(u_\beta) \geq 2\ell(4\ell - 7) > M$. If $\ell \in [5, 6]$ and $(\lambda_1, \lambda_2) = (\omega_\ell, \omega_1)$, then $\codim V_\kappa(s), \ \codim C_V(u_\alpha) \geq 2^\ell \ell > M$, and $\codim C_V(u_\beta) \geq 2^{\ell - 1}\ell > M$. If $\ell \in [3, \infty)$ and $(\lambda_1, \lambda_2) = (\omega_2, \omega_2)$ with $p = 2$, then $\codim V_\kappa(s), \ \codim C_V(u_\alpha) \geq 2(\ell - 1)(2\ell^2 - \ell - 2) > M$, and $\codim C_V(u_\beta) \geq (4\ell - 7)(2\ell^2 - \ell - 2) > M$. If $\ell \in [3, 6]$ and $(\lambda_1, \lambda_2) = (\omega_\ell, \omega_2)$ with $p = 2$, then $\codim V_\kappa(s), \ \codim C_V(u_\alpha) \geq 2^{\ell - 1}(2\ell^2 - \ell - 2) > M$, and $\codim C_V(u_\beta) \geq 2^{\ell - 2}(2\ell^2 - \ell - 2) > M$. If $\ell \in [4, 6]$ and $(\lambda_1, \lambda_2) = (\omega_\ell, \omega_\ell)$, then $\codim V_\kappa(s), \ \codim C_V(u_\alpha) \geq 2^{2\ell - 1} > M$, and $\codim C_V(u_\beta) \geq 2^{2\ell - 2} > M$. Thus in these cases the triple $(G, \lambda, p)$ satisfies $\ssddagcon$ and $\uddagcon$.

Next if $\ell = 4$ and $(\lambda_1, \lambda_2) = (\omega_\ell, \omega_1)$, then $\codim V_\kappa(s), \ \codim C_V(u_\alpha) \geq 64 > M$, and $\codim C_V(u_\beta) \geq 32 = M$. In fact $u_\beta$ has Jordan block sizes $2^4, 1^8$ on $V_1$ and $2^2, 1^4$ or $2^2, 1^5$ on $V_2$ according as $p = 2$ or $p \geq 3$, so by Lemma~\ref{lem: Jordan block tensor product} we have $\dim C_V(u_\beta) = 8.2 + 64.1 = 80$ or $8.2 + 76.1 = 92$, whence $\codim C_V(u_\beta) = 48 > M$ or $52 > M$. Thus in this case the triple $(G, \lambda, p)$ satisfies $\ssddagcon$ and $\uddagcon$.

Next if $\ell = 3$ and $(\lambda_1, \lambda_2) = (\omega_\ell, \omega_1)$, then $\codim V_\kappa(s), \ \codim C_V(u_\alpha) \geq 24 > M$, and $\codim C_V(u_\beta) \geq 12$. In fact $u_\beta$ has Jordan block sizes $2^2, 1^4$ on $V_1$ and $2^2, 1^2$ or $2^2, 1^3$ on $V_2$ according as $p = 2$ or $p \geq 3$, so by Lemma~\ref{lem: Jordan block tensor product} we have $\dim C_V(u_\beta) = 4.2 + 20.1 = 28$ or $4.2 + 26.1 = 34$, whence $\codim C_V(u_\beta) = 20 > M$ or $22 > M$. Thus in this case the triple $(G, \lambda, p)$ satisfies $\ssddagcon$ and $\uddagcon$.

Next if $\ell = 3$ and $(\lambda_1, \lambda_2) = (\omega_\ell, \omega_\ell)$, then $\codim V_\kappa(s), \ \codim C_V(u_\alpha) \geq 32 > M$, and $\codim C_V(u_\beta) \geq 16$. In fact $u_\beta$ has Jordan block sizes $2^2, 1^4$ on both $V_1$ and $V_2$, so by Lemma~\ref{lem: Jordan block tensor product} we have $\dim C_V(u_\beta) = 4.2 + 32.1 = 40$, whence $\codim C_V(u_\beta) = 24 > M$. Thus in this case the triple $(G, \lambda, p)$ satisfies $\ssddagcon$ and $\uddagcon$.

Next if $\ell = 2$ and $(\lambda_1, \lambda_2) = (\omega_\ell, \omega_1)$, then $\codim V_\kappa(s), \ \codim C_V(u_\alpha) \geq 10 - 2\z_{p, 2} \geq M > 6 - 2\z_{p, 2} = \dim{u_\alpha}^G$, and $\codim C_V(u_\beta) \geq 5 - \z_{p, 2}$; thus we need only consider regular semisimple classes with $p = 2$ and unipotent classes other than ${u_\alpha}^G$. If $s \in G_{(r)}$ is regular with $p = 2$, for each weight $\mu$ in $W.\omega_\ell$ the other $3$ weights $\mu'$ all satisfy $\mu - \mu' \in \Phi$; thus $\codim (V_1)_\kappa(s) \geq 3$ and hence $\codim V_\kappa(s) \geq 3 \dim V_2 = 12 > M$. Now take $u \in G_{(p)}$. If $u = u_\beta$ then $\codim C_{V_2}(u) = 2$ and so $\codim C_V(u) \geq 2 \dim V_1 = 10 - 2\z_{p, 2} \geq M > 4 = \dim u^G$; thus we need only consider $u$ regular, in which case $\codim C_{V_1}(u) = 3$ and hence $\codim C_V(u) \geq 3 \dim V_2 = 15 - 3\z_{p, 2} > M$. Thus in this case the triple $(G, \lambda, p)$ satisfies $\ssdiamevcon$ and $\udiamcon$.

Finally if $\ell = 2$ and $(\lambda_1, \lambda_2) = (\omega_\ell, \omega_\ell)$, then $\codim V_\kappa(s), \ \codim C_V(u_\alpha) \geq 8 = M > 6 - 2\z_{p, 2} = \dim{u_\alpha}^G$, and $\codim C_V(u_\beta) \geq 4$; thus we need only consider regular semisimple classes and unipotent classes other than ${u_\alpha}^G$. If $s \in G_{(r)}$ is regular, for each weight $\mu$ in $W.\omega_\ell$ the other $3$ weights $\mu'$ all satisfy $\mu - \mu' \in \Phi$; thus $\codim (V_1)_\kappa(s) \geq 3$ and hence $\codim V_\kappa(s) \geq 3 \dim V_2 = 12 > M$. Now take $u \in G_{(p)}$. If $u = u_\beta$ then $u$ has Jordan block sizes $2, 1^2$ on both $V_1$ and $V_2$, so by Lemma~\ref{lem: Jordan block tensor product} we have $\dim C_V(u) = 2 + 8.1 = 10$, whence $\codim C_V(u) = 6 > 4 = \dim u^G$; if $p = 2$ there is an additional class, but it has $u_\alpha$ in its closure by Lemma~\ref{lem: root elt class in closure of any non-triv class} and dimension $6$; finally if $u$ is regular then $\codim C_{V_1}(u) = 3$ and hence $\codim C_V(u) \geq 3 \dim V_2 = 12 > M$. Thus in this case the triple $(G, \lambda, p)$ satisfies $\ssdiamevcon$ and $\udiamcon$.
\end{proof}

\begin{prop}\label{prop: C_ell, tensor products}
Let $G = C_\ell$; suppose $(\lambda_1, \lambda_2) \neq (\omega_1, \omega_1)$. Then the triple $(G, \lambda, p)$ satisfies $\ssddagcon$ and $\uddagcon$.
\end{prop}

\begin{proof}
First suppose $\ell \in [3, \infty)$ and $\lambda' = \omega_1$. In this case the tables are as follows.
$$
\begin{array}{|*4{>{\ss}c|}}
\hline
i & \mu & |W.\mu| & m_\mu \\
\hline
1 & \omega_1 & 2\ell & 1 \\
\hline
\end{array}
\quad
\begin{array}{|*4{>{\ss}c|}}
\hline
 & & \multicolumn{1}{|>{\ss}c|}{c(s)} & \multicolumn{1}{|>{\ss}c|}{c(u_\alpha)} \\
\cline{3-4}
 \ss{\alpha\mathrm{-strings}} &      m      & r \geq 2 & p \geq 2 \\
\hline
            \mu_1             & 2(\ell - 2) &          &          \\
        \mu_1 \ \mu_1         &      2      &     2    &     2    \\
\hline
\multicolumn{2}{c|}{}                       &     2    &     2    \\
\cline{3-4}
\end{array}
\quad
\begin{array}{|*3{>{\ss}c|}}
\hline
 & & \multicolumn{1}{|>{\ss}c|}{c(u_\beta)} \\
\cline{3-3}
 \ss{\beta\mathrm{-strings}} &      m      & p \geq 2 \\
\hline
            \mu_1             & 2(\ell - 1) &          \\
        \mu_1 \ \mu_1         &      1      &     1    \\
\hline
\multicolumn{2}{c|}{}                       &     1    \\
\cline{3-3}
\end{array}
$$
Next suppose $\ell \in [3, \infty)$ and $\lambda' = \omega_2$; write $\z = \z_{p, \ell}$ and $\z' = \z\z_{\ell, 3}$. In this case the tables are as follows.
$$
\begin{array}{|*4{>{\ss}c|}}
\hline
i & \mu & |W.\mu| & m_\mu \\
\hline
1 & \omega_2 & 2\ell(\ell - 1) &       1       \\
0 &    0     &         1       & \ell - 1 - \z \\
\hline
\end{array}
\quad
\begin{array}{|*3{>{\ss}c|}}
\hline
 & & \multicolumn{1}{|>{\ss}c|}{c(u_\beta)} \\
\cline{3-3}
 \ss{\beta\mathrm{-strings}} &           m           &   p \geq 2  \\
\hline
            \mu_1             & 2(\ell - 1)(\ell - 2) &             \\
        \mu_1 \ \mu_1         &      2(\ell - 1)      & 2(\ell - 1) \\
            \mu_0             &           1           &             \\
\hline
\multicolumn{2}{c|}{}                                 & 2(\ell - 1) \\
\cline{3-3}
\end{array}
$$
$$
\begin{array}{|*6{>{\ss}c|}}
\hline
 & & \multicolumn{2}{|>{\ss}c|}{c(s)} & \multicolumn{2}{|>{\ss}c|}{c(u_\alpha)} \\
\cline{3-6}
 \ss{\alpha\mathrm{-strings}} &           m           &      r = 2      &   r \geq 3  &    p = 2    &   p \geq 3  \\
\hline
            \mu_1             & 2\ell^2 - 10\ell + 14 &                 &             &             &             \\
        \mu_1 \ \mu_1         &      4(\ell - 2)      &   4(\ell - 2)   & 4(\ell - 2) & 4(\ell - 2) & 4(\ell - 2) \\
    \mu_1 \ \mu_0 \ \mu_1     &           1           &     2 - \z'     &      2      &      1      &      2      \\
\hline
\multicolumn{2}{c|}{}                                 & 4\ell - 6 - \z' &  4\ell - 6  &  4\ell - 7  &  4\ell - 6  \\
\cline{3-6}
\end{array}
$$
Next suppose $\ell = 3$ and $\lambda' = \omega_3$ with $p \geq 3$. In this case the tables are as follows.
$$
\begin{array}{|*4{>{\ss}c|}}
\hline
i & \mu & |W.\mu| & m_\mu \\
\hline
2 & \omega_3 & 8 & 1 \\
1 & \omega_1 & 6 & 1 \\
\hline
\end{array}
\quad
\begin{array}{|*5{>{\ss}c|}}
\hline
 & & \multicolumn{2}{|>{\ss}c|}{c(s)} & \multicolumn{1}{|>{\ss}c|}{c(u_\alpha)} \\
\cline{3-5}
 \ss{\alpha\mathrm{-strings}} & m & r = 2 & r \geq 3 & p \geq 3 \\
\hline
            \mu_2             & 4 &       &          &          \\
    \mu_2 \ \mu_1 \ \mu_2     & 2 &   2   &     4    &     4    \\
        \mu_1 \ \mu_1         & 2 &   2   &     2    &     2    \\
\hline
\multicolumn{2}{c|}{}             &   4   &     6    &     6    \\
\cline{3-5}
\end{array}
\quad
\begin{array}{|*3{>{\ss}c|}}
\hline
 & & \multicolumn{1}{|>{\ss}c|}{c(u_\beta)} \\
\cline{3-3}
 \ss{\beta\mathrm{-strings}} & m & p \geq 3 \\
\hline
        \mu_2 \ \mu_2         & 4 &     4    \\
            \mu_1             & 4 &          \\
        \mu_1 \ \mu_1         & 1 &     1    \\
\hline
\multicolumn{2}{c|}{}             &     5    \\
\cline{3-3}
\end{array}
$$
Finally suppose $\ell \in [3, 6]$ and $\lambda' = \omega_\ell$ with $p = 2$. In this case the tables are as follows.
$$
\begin{array}{|*4{>{\ss}c|}}
\hline
i & \mu & |W.\mu| & m_\mu \\
\hline
1 & \omega_\ell & 2^\ell & 1 \\
\hline
\end{array}
\quad
\begin{array}{|*3{>{\ss}c|}}
\hline
 & & \multicolumn{1}{|>{\ss}c|}{c(u_\alpha)} \\
\cline{3-3}
 \ss{\alpha\mathrm{-strings}} &       m      &     p = 2    \\
\hline
            \mu_1             & 2^{\ell - 1} &              \\
        \mu_1 \ \mu_1         & 2^{\ell - 2} & 2^{\ell - 2} \\
\hline
\multicolumn{2}{c|}{}                        & 2^{\ell - 2} \\
\cline{3-3}
\end{array}
\quad
\begin{array}{|*4{>{\ss}c|}}
\hline
 & & \multicolumn{1}{|>{\ss}c|}{c(s)} & \multicolumn{1}{|>{\ss}c|}{c(u_\beta)} \\
\cline{3-4}
 \ss{\beta\mathrm{-strings}} &       m      &    r \geq 3   &     p = 2    \\
\hline
        \mu_1 \ \mu_1         & 2^{\ell - 1} &  2^{\ell - 1} & 2^{\ell - 1} \\
\hline
\multicolumn{2}{c|}{}                        &  2^{\ell - 1} & 2^{\ell - 1} \\
\cline{3-4}
\end{array}
$$
We have $M = 2\ell^2$.

Now if $\ell \in [3, \infty)$ and $(\lambda_1, \lambda_2) = (\omega_2, \omega_1)$, then $\codim V_\kappa(s), \ \codim C_V(u_\alpha) \geq 2\ell(4\ell - 7) > M$, and $\codim C_V(u_\beta) \geq 4\ell(\ell - 1) > M$. If $\ell = 3$ and $(\lambda_1, \lambda_2) = (\omega_3, \omega_1)$ with $p \geq 3$, then $\codim V_\kappa(s) \geq 24 > M$, $\codim C_V(u_\alpha) \geq 36 > M$, and $\codim C_V(u_\beta) \geq 30 > M$. If $\ell \in [5, 6]$ and $(\lambda_1, \lambda_2) = (\omega_\ell, \omega_1)$ with $p = 2$, then $\codim C_V(u_\alpha) \geq 2^{\ell - 1} \ell > M$, and $\codim V_\kappa(s), \ \codim C_V(u_\beta) \geq 2^\ell \ell > M$. If $\ell \in [3, \infty)$ and $(\lambda_1, \lambda_2) = (\omega_2, \omega_2)$, then $\codim V_\kappa(s), \ \codim C_V(u_\alpha) \geq (4\ell - 7)(2\ell^2 - \ell - 2) > M$, and $\codim C_V(u_\beta) \geq 2(\ell - 1)(2\ell^2 - \ell - 2) > M$. If $\ell = 3$ and $(\lambda_1, \lambda_2) = (\omega_3, \omega_2)$ with $p \geq 3$, then $\codim V_\kappa(s) \geq 52 > M$, $\codim C_V(u_\alpha) \geq 78 > M$, and $\codim C_V(u_\beta) \geq 65 > M$. If $\ell \in [3, 6]$ and $(\lambda_1, \lambda_2) = (\omega_\ell, \omega_2)$ with $p = 2$, then $\codim C_V(u_\alpha) \geq 2^{\ell - 2}(2\ell^2 - \ell - 2) > M$, and $\codim V_\kappa(s), \ \codim C_V(u_\beta) \geq 2^{\ell - 1}(2\ell^2 - \ell - 2) > M$. If $\ell = 3$ and $(\lambda_1, \lambda_2) = (\omega_3, \omega_3)$ with $p \geq 3$, then $\codim V_\kappa(s) \geq 56 > M$, $\codim C_V(u_\alpha) \geq 84 > M$, and $\codim C_V(u_\beta) \geq 70 > M$. If $\ell \in [4, 6]$ and $(\lambda_1, \lambda_2) = (\omega_\ell, \omega_\ell)$ with $p = 2$, then $\codim C_V(u_\alpha) \geq 2^{2\ell - 2} > M$, and $\codim V_\kappa(s), \ \codim C_V(u_\beta) \geq 2^{2\ell - 1} > M$. Thus in these cases the triple $(G, \lambda, p)$ satisfies $\ssddagcon$ and $\uddagcon$.

If $\ell = 4$ and $(\lambda_1, \lambda_2) = (\omega_\ell, \omega_1)$ with $p = 2$, then $\codim V_\kappa(s), \ \codim C_V(u_\beta) \geq 64 > M$, and $\codim C_V(u_\alpha) \geq 32 = M$. In fact $u_\alpha$ has Jordan block sizes $2^4, 1^8$ on $V_1$ and $2^2, 1^4$ on $V_2$, so by Lemma~\ref{lem: Jordan block tensor product} we have $\dim C_V(u_\alpha) = 8.2 + 64.1 = 80$, whence $\codim C_V(u_\alpha) = 48 > M$. Thus in this case the triple $(G, \lambda, p)$ satisfies $\ssddagcon$ and $\uddagcon$.

If $\ell = 3$ and $(\lambda_1, \lambda_2) = (\omega_\ell, \omega_1)$ with $p = 2$, then $\codim V_\kappa(s), \ \codim C_V(u_\beta) \geq 24 > M$, and $\codim C_V(u_\alpha) \geq 12$. In fact $u_\alpha$ has Jordan block sizes $2^2, 1^4$ on $V_1$ and $2^2, 1^2$ on $V_2$, so by Lemma~\ref{lem: Jordan block tensor product} we have $\dim C_V(u_\alpha) = 4.2 + 20.1 = 28$, whence $\codim C_V(u_\alpha) = 20 > M$. Thus in this case the triple $(G, \lambda, p)$ satisfies $\ssddagcon$ and $\uddagcon$.

If $\ell = 3$ and $(\lambda_1, \lambda_2) = (\omega_\ell, \omega_\ell)$ with $p = 2$, then $\codim V_\kappa(s), \ \codim C_V(u_\beta) \geq 32 > M$, and $\codim C_V(u_\alpha) \geq 16$. In fact $u_\alpha$ has Jordan block sizes $2^2, 1^4$ on both $V_1$ and $V_2$, so by Lemma~\ref{lem: Jordan block tensor product} we have $\dim C_V(u_\alpha) = 4.2 + 32.1 = 40$, whence $\codim C_V(u_\alpha) = 24 > M$. Thus in this case the triple $(G, \lambda, p)$ satisfies $\ssddagcon$ and $\uddagcon$.
\end{proof}

\begin{prop}\label{prop: D_ell, tensor products}
Let $G = D_\ell$; suppose $(\lambda_1, \lambda_2) \neq (\omega_1, \omega_1)$. Then if $\ell \in [4, 5]$ and $(\lambda_1, \lambda_2) = (\omega_\ell, \omega_1)$, the triple $(G, \lambda, p)$ satisfies $\ssdiamevcon$ and $\udiamcon$; in all other cases the triple $(G, \lambda, p)$ satisfies $\ssddagcon$ and $\uddagcon$.
\end{prop}

\begin{proof}
First suppose $\ell \in [4, \infty)$ and $\lambda' = \omega_1$. In this case the tables are as follows.
$$
\begin{array}{|*4{>{\ss}c|}}
\hline
i & \mu & |W.\mu| & m_\mu \\
\hline
1 & \omega_1 & 2\ell & 1 \\
\hline
\end{array}
\quad
\begin{array}{|*4{>{\ss}c|}}
\hline
 & & \multicolumn{1}{|>{\ss}c|}{c(s)} & \multicolumn{1}{|>{\ss}c|}{c(u_\alpha)} \\
\cline{3-4}
 \ss{\alpha\mathrm{-strings}} &      m      & r \geq 2 & p \geq 2 \\
\hline
            \mu_1             & 2(\ell - 2) &          &          \\
        \mu_1 \ \mu_1         &      2      &     2    &     2    \\
\hline
\multicolumn{2}{c|}{}                       &     2    &     2    \\
\cline{3-4}
\end{array}
$$
Next suppose $\ell \in [5, 7]$ and $\lambda' = \omega_\ell$. In this case the tables are as follows.
$$
\begin{array}{|*4{>{\ss}c|}}
\hline
i & \mu & |W.\mu| & m_\mu \\
\hline
1 & \omega_\ell & 2^{\ell - 1} & 1 \\
\hline
\end{array}
\quad
\begin{array}{|*4{>{\ss}c|}}
\hline
 & & \multicolumn{1}{|>{\ss}c|}{c(s)} & \multicolumn{1}{|>{\ss}c|}{c(u_\alpha)} \\
\cline{3-4}
 \ss{\alpha\mathrm{-strings}} &       m      &    r \geq 2   &    p \geq 2   \\
\hline
            \mu_1             & 2^{\ell - 2} &               &               \\
        \mu_1 \ \mu_1         & 2^{\ell - 3} &  2^{\ell - 3} &  2^{\ell - 3} \\
\hline
\multicolumn{2}{c|}{}                        &  2^{\ell - 3} &  2^{\ell - 3} \\
\cline{3-4}
\end{array}
$$
We have $M = 2\ell(\ell - 1)$.

Now if $\ell \in [6, 7]$ and $(\lambda_1, \lambda_2) = (\omega_\ell, \omega_1)$, then $\codim V_\kappa(s), \ \codim C_V(u_\alpha) \geq 2^{\ell - 2}\ell > M$. If $\ell \in [5, 7]$ and $(\lambda_1, \lambda_2) = (\omega_\ell, \omega_\ell)$ or $(\omega_\ell, \omega_{\ell - 1})$, then $\codim V_\kappa(s)$, $\codim C_V(u_\alpha) \geq 2^{2\ell - 4} > M$. Thus in these cases the triple $(G, \lambda, p)$ satisfies $\ssddagcon$ and $\uddagcon$.

Next if $\ell = 5$ and $(\lambda_1, \lambda_2) = (\omega_\ell, \omega_1)$, then $\codim V_\kappa(s), \ \codim C_V(u_\alpha) \geq 40 = M$; thus we need only consider regular classes. If $s \in G_{(r)}$ is regular, for each weight $\mu$ in $W.\omega_\ell$ there are $10$ other weights $\mu'$ with $\mu - \mu' \in \Phi$ (e.g. if $\mu = \omega_5$ then we may take $\mu' = \mu - \alpha$ for $\alpha \in \Phi^+ \setminus \langle \alpha_1, \alpha_2, \alpha_3, \alpha_4 \rangle$); thus $\codim (V_1)_\kappa(s) \geq 10$ and hence $\codim V_\kappa(s) \geq 10 \dim V_2 = 100 > M$. If $u \in G_{(p)}$ is regular, then $\codim C_{V_2}(u) = 8$ and hence $\codim C_V(u) \geq 8 \dim V_1 = 128 > M$. Thus in this case the triple $(G, \lambda, p)$ satisfies $\ssdiamevcon$ and $\udiamcon$.

Finally if $\ell = 4$ and $(\lambda_1, \lambda_2) = (\omega_\ell, \omega_1)$, then $\codim V_\kappa(s), \ \codim C_V(u_\alpha) \geq 16$; thus we need only consider classes of dimension at least $16$. If the centralizer of $s \in G_{(r)}$ is $A_3$ or $D_3$ then $\dim s^G = 12 < 16$; if not, there are $3$ mutually orthogonal roots $\alpha$ with $\alpha(s) \neq 1$, which we may assume are $\alpha_1$, $\alpha_3$ and $\alpha_4$. The $8$ weights in $W.\omega_\ell$ may then be divided into $4$ pairs
\begin{eqnarray*}
& \{ \omega_4, \omega_4 - \alpha_4 \}, & \\
& \{ \omega_4 - \alpha_2 - \alpha_4, \omega_4 - \alpha_1 - \alpha_2 - \alpha_4 \}, & \\
& \{ \omega_4 - \alpha_2 - \alpha_3 - \alpha_4, \omega_4 - \alpha_1 - \alpha_2 - \alpha_3 - \alpha_4 \}, & \\
& \{ \omega_4 - \alpha_1 - 2\alpha_2 - \alpha_3 - \alpha_4, \omega_4 - \alpha_1 - 2\alpha_2 - \alpha_3 - 2\alpha_4 \}, &
\end{eqnarray*}
with the two weights in each pair differing by $\alpha_1$, $\alpha_3$ or $\alpha_4$; thus any eigenspace for $s$ in $V_1$ has codimension at least $4$, so $\codim V_\kappa(s) \geq 4 \dim V_2 = 32 > M$. Now take $u \in G_{(p)}$. If $u = u_\alpha$ then $u$ has Jordan block sizes $2^2, 1^4$ on both $V_1$ and $V_2$, so by Lemma~\ref{lem: Jordan block tensor product} we have $\dim C_V(u) = 4.2 + 32.1 = 40$, whence $\codim C_V(u) = 24 = M$, so we need only consider regular unipotent elements; if $u$ is regular then $\codim C_{V_1}(u) = 6$ and hence $\codim C_V(u) \geq 6 \dim V_2 = 48 > M$. Thus in this case the triple $(G, \lambda, p)$ satisfies $\ssdiamevcon$ and $\udiamcon$.
\end{proof}

\begin{prop}\label{prop: exceptional groups, tensor products}
Let $G$ be of exceptional type. Then the triple $(G, \lambda, p)$ satisfies $\ssddagcon$ and $\uddagcon$.
\end{prop}

\begin{proof}
Take $G = E_6$. Suppose $\lambda' = \omega_1$. In this case the tables are as follows.
$$
\begin{array}{|*4{>{\ss}c|}}
\hline
i & \mu & |W.\mu| & m_\mu \\
\hline
1 & \omega_1 & 27 & 1 \\
\hline
\end{array}
\quad
\begin{array}{|*4{>{\ss}c|}}
\hline
 & & \multicolumn{1}{|>{\ss}c|}{c(s)} & \multicolumn{1}{|>{\ss}c|}{c(u_\alpha)} \\
\cline{3-4}
 \ss{\alpha\mathrm{-strings}} &  m & r \geq 2 & p \geq 2 \\
\hline
            \mu_1             & 15 &          &          \\
        \mu_1 \ \mu_1         &  6 &     6    &     6    \\
\hline
\multicolumn{2}{c|}{}              &     6    &     6    \\
\cline{3-4}
\end{array}
$$
We have $M = 72$. Now if $(\lambda_1, \lambda_2) = (\omega_1, \omega_1)$ or $(\omega_1, \omega_6)$, then $\codim V_\kappa(s)$, $\codim C_V(u_\alpha) \geq 162 > M$. Thus in these cases the triple $(G, \lambda, p)$ satisfies $\ssddagcon$ and $\uddagcon$.

Next take $G = E_7$. Suppose $\lambda' = \omega_7$. In this case the tables are as follows.
$$
\begin{array}{|*4{>{\ss}c|}}
\hline
i & \mu & |W.\mu| & m_\mu \\
\hline
1 & \omega_7 & 56 & 1 \\
\hline
\end{array}
\quad
\begin{array}{|*4{>{\ss}c|}}
\hline
 & & \multicolumn{1}{|>{\ss}c|}{c(s)} & \multicolumn{1}{|>{\ss}c|}{c(u_\alpha)} \\
\cline{3-4}
 \ss{\alpha\mathrm{-strings}} &  m & r \geq 2 & p \geq 2 \\
\hline
            \mu_1             & 32 &          &          \\
        \mu_1 \ \mu_1         & 12 &    12    &    12    \\
\hline
\multicolumn{2}{c|}{}              &    12    &    12    \\
\cline{3-4}
\end{array}
$$
We have $M = 126$. Now if $(\lambda_1, \lambda_2) = (\omega_7, \omega_7)$, then $\codim V_\kappa(s), \ \codim C_V(u_\alpha) \geq 672 > M$. Thus in this case the triple $(G, \lambda, p)$ satisfies $\ssddagcon$ and $\uddagcon$.

Next take $G = F_4$. Suppose $\lambda' = \omega_4$; write $\z = \z_{p, 3}$. In this case the tables are as follows.
$$
\begin{array}{|*4{>{\ss}c|}}
\hline
i & \mu & |W.\mu| & m_\mu \\
\hline
1 & \omega_4 & 24 & 1 \\
0 &    0     &  1 & 2 - \z \\
\hline
\end{array}
\quad
\begin{array}{|*6{>{\ss}c|}}
\hline
 & & \multicolumn{2}{|>{\ss}c|}{c(s)} & \multicolumn{2}{|>{\ss}c|}{c(u_\alpha)} \\
\cline{3-6}
 \ss{\alpha\mathrm{-strings}} & m &  r = 2  & r \geq 3 & p = 2 & p \geq 3 \\
\hline
            \mu_1             & 6 &         &          &       &          \\
        \mu_1 \ \mu_1         & 8 &    8    &     8    &   8   &     8    \\
    \mu_1 \ \mu_0 \ \mu_1     & 1 &  2 - \z &     2    &   1   &     2    \\
\hline
\multicolumn{2}{c|}{}             & 10 - \z &    10    &   9   &    10    \\
\cline{3-6}
\end{array}
\quad
\begin{array}{|*3{>{\ss}c|}}
\hline
 & & \multicolumn{1}{|>{\ss}c|}{c(u_\beta)} \\
\cline{3-3}
 \ss{\beta\mathrm{-strings}} &  m & p \geq 2 \\
\hline
            \mu_1             & 12 &          \\
        \mu_1 \ \mu_1         &  6 &     6    \\
            \mu_0             &  1 &          \\
\hline
\multicolumn{2}{c|}{}              &     6    \\
\cline{3-3}
\end{array}
$$
Now suppose $\lambda' = \omega_1$ with $p = 2$. In this case the tables are as follows.
$$
\begin{array}{|*4{>{\ss}c|}}
\hline
i & \mu & |W.\mu| & m_\mu \\
\hline
1 & \omega_1 & 24 & 1 \\
0 &    0     &  1 & 2 \\
\hline
\end{array}
\quad
\begin{array}{|*4{>{\ss}c|}}
\hline
 & & \multicolumn{1}{|>{\ss}c|}{c(s)} & \multicolumn{1}{|>{\ss}c|}{c(u_\alpha)} \\
\cline{3-4}
 \ss{\alpha\mathrm{-strings}} &  m & r \geq 3 & p = 2 \\
\hline
            \mu_1             & 12 &          &       \\
    \mu_1 \ \phmu \ \mu_1     &  6 &     6    &   6   \\
            \mu_0             &  1 &          &       \\
\hline
\multicolumn{2}{c|}{}              &     6    &   6   \\
\cline{3-4}
\end{array}
\quad
\begin{array}{|*3{>{\ss}c|}}
\hline
 & & \multicolumn{1}{|>{\ss}c|}{c(u_\beta)} \\
\cline{3-3}
 \ss{\beta\mathrm{-strings}} & m & p = 2 \\
\hline
            \mu_1             & 6 &       \\
        \mu_1 \ \mu_1         & 8 &   8   \\
    \mu_1 \ \mu_0 \ \mu_1     & 1 &   1   \\
\hline
\multicolumn{2}{c|}{}             &   9   \\
\cline{3-3}
\end{array}
$$
We have $M = 48$. Now if $(\lambda_1, \lambda_2) = (\omega_4, \omega_4)$, then $\codim V_\kappa(s), \ \codim C_V(u_\alpha) \geq 225 > M$, and $\codim C_V(u_\beta) \geq 150 > M$. If $(\lambda_1, \lambda_2) = (\omega_1, \omega_4)$ or $(\omega_1, \omega_1)$ with $p = 2$, then $\codim V_\kappa(s), \ \codim C_V(u_\alpha) \geq 156 > M$, and $\codim C_V(u_\beta) \geq 234 > M$. Thus in these cases the triple $(G, \lambda, p)$ satisfies $\ssddagcon$ and $\uddagcon$.

Finally take $G = G_2$. Suppose $\lambda' = \omega_1$; write $\z = \z_{p, 2}$. In this case the tables are as follows.
$$
\begin{array}{|*4{>{\ss}c|}}
\hline
i & \mu & |W.\mu| & m_\mu \\
\hline
1 & \omega_1 & 6 & 1 \\
0 &    0     & 1 & 1 - \z \\
\hline
\end{array}
\quad
\begin{array}{|*6{>{\ss}c|}}
\hline
 & & \multicolumn{2}{|>{\ss}c|}{c(s)} & \multicolumn{2}{|>{\ss}c|}{c(u_\alpha)} \\
\cline{3-6}
 \ss{\alpha\mathrm{-strings}} & m & r = 2 & r \geq 3 & p = 2 & p \geq 3 \\
\hline
        \mu_1 \ \mu_1         & 2 &   2   &     2    &   2   &     2    \\
    \mu_1 \ \mu_0 \ \mu_1     & 1 &   1   &  2 - \z  &   1   &     2    \\
\hline
\multicolumn{2}{c|}{}             &   3   &  4 - \z  &   3   &     4    \\
\cline{3-6}
\end{array}
\quad
\begin{array}{|*3{>{\ss}c|}}
\hline
 & & \multicolumn{1}{|>{\ss}c|}{c(u_\beta)} \\
\cline{3-3}
 \ss{\beta\mathrm{-strings}} & m & p \geq 2 \\
\hline
            \mu_1             & 2 &          \\
        \mu_1 \ \mu_1         & 2 &     2    \\
            \mu_0             & 1 &          \\
\hline
\multicolumn{2}{c|}{}             &     2    \\
\cline{3-3}
\end{array}
$$
Now suppose $\lambda' = \omega_2$ with $p = 3$. In this case the tables are as follows.
$$
\begin{array}{|*4{>{\ss}c|}}
\hline
i & \mu & |W.\mu| & m_\mu \\
\hline
1 & \omega_2 & 6 & 1 \\
0 &    0     & 1 & 1 \\
\hline
\end{array}
\quad
\begin{array}{|*4{>{\ss}c|}}
\hline
 & & \multicolumn{1}{|>{\ss}c|}{c(s)} & \multicolumn{1}{|>{\ss}c|}{c(u_\alpha)} \\
\cline{3-4}
 \ss{\alpha\mathrm{-strings}} & m & r \neq 3 & p = 3 \\
\hline
            \mu_1             & 2 &          &       \\
\mu_1 \ \phmu \ \phmu \ \mu_1 & 2 &     2    &   2   \\
            \mu_0             & 1 &          &       \\
\hline
\multicolumn{2}{c|}{}             &     2    &   2   \\
\cline{3-4}
\end{array}
\quad
\begin{array}{|*3{>{\ss}c|}}
\hline
 & & \multicolumn{1}{|>{\ss}c|}{c(u_\beta)} \\
\cline{3-3}
 \ss{\beta\mathrm{-strings}} & m & p = 3 \\
\hline
        \mu_1 \ \mu_1         & 2 &   2   \\
    \mu_1 \ \mu_0 \ \mu_1     & 1 &   2   \\
\hline
\multicolumn{2}{c|}{}             &   4   \\
\cline{3-3}
\end{array}
$$
We have $M = 12$ and $M_2 = 8$. Now if $(\lambda_1, \lambda_2) = (\omega_1, \omega_1)$ with $p \geq 3$, then $\codim V_\kappa(s) \geq 21 > M$, $\codim C_V(u_\alpha) \geq 28 > M$, and $\codim C_V(u_\beta) \geq 14 > M$. If $(\lambda_1, \lambda_2) = (\omega_2, \omega_1)$ or $(\omega_2, \omega_2)$ with $p = 3$, then $\codim V_\kappa(s), \ \codim C_V(u_\alpha) \geq 14 > M$, and $\codim C_V(u_\beta) \geq 28 > M$. Thus in these cases the triple $(G, \lambda, p)$ satisfies $\ssddagcon$ and $\uddagcon$.

Finally if $(\lambda_1, \lambda_2) = (\omega_1, \omega_1)$ with $p = 2$, then $\codim V_\kappa(s), \ \codim C_V(u_\alpha) \geq 18 > M$, and $\codim C_V(u_\beta) \geq 12$. In fact $u_\beta$ has Jordan block sizes $2^2, 1^2$ on both $V_1$ and $V_2$, so by Lemma~\ref{lem: Jordan block tensor product} we have $\dim C_V(u_\beta) = 4.2 + 12.1 = 20$, whence $\codim C_V(u_\beta) = 16 > M$. Thus in this case the triple $(G, \lambda, p)$ satisfies $\ssddagcon$ and $\uddagcon$.
\end{proof}

We now turn to the postponed cases. For convenience we give separate results treating semisimple and unipotent classes.

\begin{prop}\label{prop: B_ell, C_ell, D_ell, omega_1 + q omega_1, ss}
Let $G = B_\ell$, $C_\ell$ or $D_\ell$ and $(\lambda_1, \lambda_2) = (\omega_1, \omega_1)$. Then the triple $(G, \lambda, p)$ satisfies $\ssdiamevcon$.
\end{prop}

\begin{proof}
Take $s \in G_{(r)}$. First suppose $G = C_\ell$. Take a basis 
$$
v_{-\ell}, \dots, v_{-1}, v_1, \dots, v_\ell
$$
of $V_1 = L(\omega_1)$ such that $sv_i = \delta_i v_i$ for $i = -\ell, \dots, -1, 1, \dots, \ell$, where $\delta_{-i} = {\delta_i}^{-1}$. Applying a suitable Weyl group element we may assume the eigenvalues $\delta_1, \dots, \delta_\ell$ are
$$
1^a, (-1)^b, {\kappa_1}^{a_1}, (-\kappa_1)^{b_1}, \dots, {\kappa_t}^{a_t}, (-\kappa_t)^{b_t},
$$
where the multiplicities $a, b, a_1, b_1, \dots, a_t, b_t$ are all non-negative integers such that $a, b < \ell$, for each $h$ we have $a_h + b_h > 0$ and ${\kappa_h}^2 \neq 1$, and for each $i, j$ with $i \neq j$ we have $\kappa_i{\kappa_j}^{\pm 1} \neq \pm 1$. Then $C_G(s)$ is of type
$$
C_a C_b A_{a_1 - 1} A_{b_1 - 1} \dots A_{a_t - 1} A_{b_t - 1}
$$
(where we ignore terms $A_{-1}$), and so
$$
\dim s^G = 2\ell^2 - \left[2a^2 + 2b^2 + \sum (a_h(a_h - 1) + b_h(b_h - 1))\right].
$$
Now the vectors $v_i \otimes v_j$ for $i, j \in \{ -\ell, \dots, -1, 1, \dots, \ell \}$ form a basis of $V = V_1 \otimes V_2$, and we have $s(v_i \otimes v_j) = \kappa_{i, j} v_i \otimes v_j$, where $\kappa_{i, j} = \delta_i{\delta_j}^q$. We consider the multiplicity in $V$ of a given eigenvalue $\kappa$; we shall show that for all choices of $\kappa$ we have $\codim V_\kappa(s) > \dim s^G$.

First assume $\kappa \neq \pm 1, {\pm \kappa_h}^{\pm1}$ for each $h$. Given $i$ and $j$, at most one of $\kappa_{i, j}$ and $\kappa_{i, -j}$ can equal $\kappa$; thus for each $i$ there are at least $\ell$ values $j$ with $\kappa_{i, j} \neq \kappa$, so $\codim V_\kappa(s) \geq 2\ell.\ell = 2\ell^2 \geq \dim s^G$. If $\dim s^G = 2\ell^2$ then $a = b = 0$ and all $a_h$ and $b_h$ are at most $1$, so all $\delta_j$ are distinct, as are all $\kappa_{i, j}$ for any fixed $i$; thus $\codim V_\kappa(s) \geq 2\ell(2\ell - 1) > \dim s^G$.

Next assume $\kappa = {\pm \kappa_h}^{\pm1}$ for some $h$; without loss of generality we may assume $\kappa = \kappa_h$. Given $i$, if $\delta_i \neq \pm \kappa_h$ then for each $j$ at most one of $\kappa_{i, j}$ and $\kappa_{i, -j}$ can equal $\kappa$, which gives at least $(2\ell - (a_h + b_h)).\ell$ pairs $(i, j)$ with $\kappa_{i, j} \neq \kappa$; if instead $\delta_i = \kappa_h$ (which occurs for $a_h$ values $i$) then there are $2a$ values $j$ for which $\kappa_{i, j} = \kappa$, while if $\delta_i = -\kappa_h$ (which occurs for $b_h$ values $i$) then there are $2b$ values $j$ for which $\kappa_{i, j} = \kappa$. Thus $\codim V_\kappa(s) \geq (2\ell - (a_h + b_h)).\ell + a_h(2\ell - 2a) + b_h(2\ell - 2b)$; so
\begin{eqnarray*}
\codim V_\kappa(s) - \dim s^G
& \geq & (2\ell - (a_h + b_h)).\ell + a_h(2\ell - 2a) + b_h(2\ell - 2b) \\
&      & \qquad {} - (2\ell^2 - 2a^2 - 2b^2) \\
&   =  & a_h\ell - 2aa_h + 2a^2 + b_h\ell -2bb_h + 2b^2 \\
&   >  & 0
\end{eqnarray*}
(observe that if $a < \frac{\ell}{2}$ then certainly $a_h\ell - 2aa_h + 2a^2 \geq 0$ with equality only if $a_h = a = 0$, while if $a \geq \frac{\ell}{2}$ then $a_h\ell - 2aa_h + 2a^2 = \frac{\ell^2}{2} + 2(a - \frac{\ell}{2})(\ell - a_h) + 2(a - \frac{\ell}{2})^2 > 0$; likewise $b_h\ell - 2bb_h + 2b^2 \geq 0$ with equality only if $b_h = b = 0$ --- but we cannot have both $a_h$ and $b_h$ equal to $0$).

Next assume $\kappa = 1$. Given $i$, if $\delta_i = 1$ (which occurs for $2a$ values $i$) then there are $2a$ values $j$ for which $\kappa_{i, j} = \kappa$, while if $\delta_i = -1$ (which occurs for $2b$ values $i$) then there are $2b$ values $j$ for which $\kappa_{i, j} = \kappa$; if instead $\delta_i \neq \pm 1$, then for each $j$ if $\delta_j = \pm 1$ then $\kappa_{i, j} \neq \kappa$ while if $\delta_j \neq \pm 1$ then at most one of $\kappa_{i, j}$ and $\kappa_{i, -j}$ can equal $\kappa$, so there are at least $\ell + a + b$ values $j$ with $\kappa_{i, j} \neq \kappa$. Thus $\codim V_\kappa(s) \geq 2a(2\ell - 2a) + 2b(2\ell - 2b) + (2\ell - 2a - 2b)(\ell + a + b)$; so
\begin{eqnarray*}
\codim V_\kappa(s) - \dim s^G
& \geq & 2a(2\ell - 2a) + 2b(2\ell - 2b) + (2\ell - 2a - 2b)(\ell + a + b) \\
&      & \qquad {} - (2\ell^2 - 2a^2 - 2b^2) \\
&   =  & 4a\ell - 4a^2 + 4b\ell - 4b^2 - 4ab \\
&   =  & 2a(\ell - a) + 2b(\ell - b) + 2(a + b)(\ell - (a + b)) \\
& \geq & 0.
\end{eqnarray*}
For equality we must have $a = b = 0$; but then if some $a_h$ or $b_h$ is greater than $1$ then $\codim V_\kappa(s) \geq 2\ell^2 > \dim s^G$, while if all $a_h$ and $b_h$ are at most $1$ then for each $i$ there can be at most one $j$ with $\kappa_{i, j} = \kappa$, so $\codim V_\kappa(s) \geq 2\ell(2\ell - 1) > \dim s^G$.

Finally assume $\kappa = -1$. The calculation here is identical to that for $\kappa = 1$, except that in the expression for the lower bound on $\codim V_\kappa(s)$ we must replace $2a(2\ell - 2a) + 2b(2\ell - 2b)$ by $2a(2\ell - 2b) + 2b(2\ell - 2a)$, which cannot decrease the value since the difference is $4a^2 - 8ab + 4b^2 = (2a - 2b)^2 \geq 0$.

Thus for all $\kappa$ we have $\codim V_\kappa(s) > \dim s^G$; so the triple $(G, \lambda, p)$ satisfies $\ssdiamevcon$.

Next suppose $G = D_\ell$. Again take a basis $v_{-\ell}, \dots, v_{-1}, v_1, \dots, v_\ell$ of $V_1 = L(\omega_1)$; we may need to interchange $v_\ell$ and $v_{-\ell}$ to ensure $\delta_1, \dots, \delta_\ell$ are as given above. The calculations for $\codim V_\kappa(s)$ are identical to those for $G = C_\ell$. On the other hand, in the expression for $\dim s^G$ we must replace $2\ell^2 - 2a^2 - 2b^2$ by $2\ell(\ell - 1) - 2a(a - 1) - 2b(b - 1)$ because the simple factors $C_a$ and $C_b$ in $C_G(s)$ are replaced by $D_a$ and $D_b$ respectively; the difference is $2\ell - 2a - 2b \geq 0$, so the value of $\dim s^G$ for $G = D_\ell$ is no larger than it is for $G = C_\ell$. Thus for all $\kappa$ we again have $\codim V_\kappa(s) > \dim s^G$; so the triple $(G, \lambda, p)$ satisfies $\ssdiamevcon$.

Finally suppose $G = B_\ell$. If $p = 2$ the details are exactly as for $G = C_\ell$, so assume $p \geq 3$. Here we take a basis $v_{-\ell}, \dots, v_{-1}, v_0, v_1, \dots, v_\ell$ of $V_1 = L(\omega_1)$ such that $sv_i = \delta_i v_i$ for $i = -\ell, \dots, -1, 0, 1, \dots, \ell$, where $\delta_{-i} = {\delta_i}^{-1}$ and $\delta_0 = 1$. We again have $\delta_1, \dots, \delta_\ell$ as above; this time the simple factors $C_a$ and $C_b$ in $C_G(s)$ are replaced by $B_a$ and $D_b$ respectively, so we obtain
$$
\dim s^G = 2\ell^2 - \left[2a^2 + 2b(b - 1) + \sum (a_h(a_h - 1) + b_h(b_h - 1))\right].
$$

First assume $\kappa \neq \pm 1, {\pm \kappa_h}^{\pm1}$ for each $h$. Given $i$, there are at least $\ell$ non-zero values $j$ with $\kappa_{i, j} \neq \kappa$, together with the value $0$, so $\codim V_\kappa(s) \geq (2\ell + 1)(\ell + 1) > 2\ell^2 \geq \dim s^G$.

Next assume $\kappa = {\pm \kappa_h}^{\pm1}$ for some $h$; without loss of generality we may assume $\kappa = \kappa_h$. Arguing as in the $G = C_\ell$ case and allowing for the extra terms $\kappa_{i, j}$ with $ij = 0$, we have
\begin{eqnarray*}
\codim V_\kappa(s) - \dim s^G
& \geq & (2\ell + 1 - (a_h + b_h))(\ell + 1) + a_h(2\ell + 1 - 2a) \\
&      & \qquad {} + b_h(2\ell + 1 - 2b) - (2\ell^2 - 2a^2 - 2b(b - 1)) \\
&   =  & a_h\ell - 2aa_h + 2a^2 + b_h\ell -2bb_h + 2b^2 + 3\ell + 1 - 2b\\
&   >  & 0.
\end{eqnarray*}

Next assume $\kappa = 1$. Arguing again as in the $G = C_\ell$ case and allowing for the extra terms $\kappa_{i, j}$, we have
\begin{eqnarray*}
\codim V_\kappa(s) - \dim s^G
& \geq & (2a + 1)(2\ell + 1 - 2a - 1) + 2b(2\ell + 1 - 2b) \\
&      & \qquad {} + (2\ell + 1 - 2a - 1 - 2b)(\ell + a + b) \\
&      & \qquad {} - (2\ell^2 - 2a^2 - 2b(b - 1)) \\
&   =  & 4a\ell - 4a^2 + 2\ell - 2a + 4b\ell - 4b^2 - 4ab \\
&   =  & 2(a + 1)(\ell - a) + 2b(\ell - b) + 2(a + b)(\ell - (a + b)) \\
&   >  & 0.
\end{eqnarray*}

Finally assume $\kappa = -1$. The calculation here is identical to that for $\kappa = 1$, except that in the expression for the lower bound on $\codim V_\kappa(s)$ we must replace $(2a + 1)(2\ell + 1 - 2a - 1) + 2b(2\ell + 1 - 2b)$ by $(2a + 1)(2\ell + 1 - 2b) + 2b(2\ell + 1 - 2a - 1)$, which cannot decrease the value since the difference is $(2a + 1)^2 - 2(2a + 1)2b + 4b^2 = (2a + 1 - 2b)^2 \geq 0$.

Thus for all $\kappa$ we have $\codim V_\kappa(s) > \dim s^G$; so the triple $(G, \lambda, p)$ satisfies $\ssdiamevcon$.
\end{proof}

\begin{prop}\label{prop: B_ell, C_ell, D_ell, omega_1 + q omega_1, unip}
Let $G = B_\ell$, $C_\ell$ or $D_\ell$ and $(\lambda_1, \lambda_2) = (\omega_1, \omega_1)$. Then the triple $(G, \lambda, p)$ satisfies $\udiamcon$.
\end{prop}

\begin{proof}
Write
$$
d = \dim V_1 =
\left\{
\begin{array}{ll}
2\ell + 1 - \z_{p, 2}, & \hbox{if $G = B_\ell$;} \\
2\ell,                 & \hbox{if $G = C_\ell$ or $D_\ell$.}
\end{array}
\right.
$$
Take $u \in G_{(p)}$ and suppose $u$ has Jordan block sizes $m_1, \dots, m_r$ on $V_1$, where $m_1 \geq m_2 \geq \cdots \geq m_r$ and $\sum m_i = d$. By Lemma~\ref{lem: Jordan block tensor product} we have
\begin{eqnarray*}
\codim C_V(u) & = & \sum_{i, j = 1}^r (m_im_j - \min(m_i, m_j)) \\
              & = & d^2 - (m_1 + 3m_2 + \cdots + (2r - 1)m_r) \\
              & = & d^2 + d - 2\sum_{i = 1}^r im_i.
\end{eqnarray*}
On the other hand from \cite{LSbook} we see that $\dim C_G(u) = \sum_{i = 1}^r (im_i - \chi(m_i))$, where the precise definition of the function $\chi$ depends on the type of $G$; thus
$$
\codim C_V(u) - \dim u^G = d^2 + d - \dim G - \sum_{i = 1}^r (im_i + \chi(m_i)).
$$
We now consider separately the possibilities for $G$.

If $G = C_\ell$ then $d = 2\ell$, $\dim G = 2\ell^2 + \ell$, and $\chi(m) \leq \frac{1}{2}m$; since $u \neq 1$ we have $\sum_{i = 1}^r im_i \leq 1.2 + 2.1 + 3.1 + \cdots + (2\ell - 1).1 = 2\ell^2 - \ell + 1$, and $\sum_{i = 1}^r \chi(m_i) \leq \sum_{i = 1}^r \frac{1}{2}m_i = \ell$. Thus
\begin{eqnarray*}
\codim C_V(u) - \dim u^G & \geq & 4\ell^2 + 2\ell - (2\ell^2 + \ell) - (2\ell^2 - \ell + 1) - \ell \\
                         &   =  & \ell - 1 \\
                         &   >  & 0.
\end{eqnarray*}

If $G = D_\ell$ then $d = 2\ell$, $\dim G = 2\ell^2 - \ell$, and $\chi(m) \leq \frac{1}{2}(m + 2)$; since $u \neq 1$ we have $\sum_{i = 1}^r im_i \leq 1.2 + 2.2 + 3.1 + \cdots + (2\ell - 2).1 = 2\ell^2 - 3\ell + 4$, and $\sum_{i = 1}^r \chi(m_i) \leq \sum_{i = 1}^r \frac{1}{2}(m_i + 2) \leq 3\ell - 2$. Thus
\begin{eqnarray*}
\codim C_V(u) - \dim u^G & \geq & 4\ell^2 + 2\ell - (2\ell^2 - \ell) - (2\ell^2 - 3\ell + 4) - (3\ell - 2) \\
                         &   =  & 3\ell - 2 \\
                         &   >  & 0.
\end{eqnarray*}

If $G = B_\ell$ then for $p = 2$ the calculation is identical to that for the case $G = C_\ell$, so we may assume $p \geq 3$; then $d = 2\ell + 1$, $\dim G = 2\ell^2 + \ell$, and $\chi(m) = \lceil \frac{1}{2}m \rceil \leq \frac{1}{2}(m + 1)$; since $u \neq 1$ we have $\sum_{i = 1}^r im_i \leq 1.2 + 2.1 + 3.1 + \cdots + 2\ell.1 = 2\ell^2 + \ell + 1$, and $\sum_{i = 1}^r \chi(m_i) \leq \sum_{i = 1}^r \frac{1}{2}(m_i + 1) \leq 2\ell + 1$. Thus
\begin{eqnarray*}
\codim C_V(u) - \dim u^G & \geq & 4\ell^2 + 6\ell + 2 - (2\ell^2 + \ell) - (2\ell^2 + \ell + 1) - (2\ell + 1) \\
                         &   =  & 2\ell \\
                         &   >  & 0.
\end{eqnarray*}

Thus in all cases we have $\codim C_V(u) > \dim u^G$; so the triple $(G, \lambda, p)$ satisfies $\udiamcon$.
\end{proof}

The results proved in this section have established the following.

\begin{prop}\label{prop: large triples which are not p-restricted with TGS}
Any large triple which is not $p$-restricted and is not listed in Table~\ref{table: large triple and first quadruple non-TGS} satisfies $\ssdiamevcon$ and $\udiamcon$, and so has TGS.
\end{prop}

\chapter{Quadruples having TGS}\label{chap: TGS quadruples}

In this chapter we develop techniques for showing that a large quadruple has TGS, and then apply them to prove that any large quadruple not listed in Table~\ref{table: large triple and first quadruple non-TGS} or \ref{table: large higher quadruple non-TGS} has TGS. The structure of this chapter is as follows. In Section~\ref{sect: prelim tuples} we do some initial work on integer tuples. In Section~\ref{sect: reduction of quadruples} we substantially reduce the number of large quadruples requiring consideration. In the remaining two sections we treat the remaining large quadruples which have TGS: Sections~\ref{sect: large quadruple individual cases} and \ref{sect: large quadruple infinite families} concern individual quadruples and infinite families of quadruples respectively.

\section{Preliminary results on tuples}\label{sect: prelim tuples}

In this section we prove some preliminary results on integer tuples, the significance of which will become clear in the following section.

Firstly, given a tuple $\a = (a_1, \dots, a_t)$ of integers, we shall write
$$
|\a| = a_1 + \cdots + a_t;
$$
we say that $\a$ is {\em decreasing\/} if $a_1 \geq \cdots \geq a_t$.

Now let $\d = (d_1, \dots, d_t)$ be a tuple of natural numbers. Given a tuple $\k = (k_1, \dots, k_t)$ of integers which has the same length as $\d$, if for each $i \in [1, t]$ we have $0 \leq k_i \leq d_i$ then we say that $\k$ is {\em $\d$-feasible\/}, and define
$$
B_{\d, \k} = |\k|(|\d| - |\k|) - \sum k_i(d_i - k_i).
$$
Given a natural number $k$ with $k \leq |\d|$, we set
$$
B_{\d, k} = \min \{ B_{\d, \k} : \k \hbox{ is $\d$-feasible}, \ |\k| = k \}.
$$
Our first result here shows that the value $B_{\d, k}$ increases with $k$ up to $\frac{1}{2}|\d|$.

\begin{prop}\label{prop: B value increases with k}
Given a tuple $\d$ of natural numbers and a natural number $k$ with $2 \leq k \leq \frac{1}{2}|\d|$, we have $B_{\d, k - 1} \leq B_{\d, k}$.
\end{prop}

\begin{proof}
Write $\d = (d_1, \dots, d_t)$ and set $d = |\d|$; take a $\d$-feasible tuple $\k = (k_1, \dots, k_t)$ with $|\k| = k$ such that $B_{\d, \k} = B_{\d, k}$. Observe that $\sum_{i = 1}^t (d_i - 2k_i) = d - 2k \geq 0$. If for all $i \in [1, t]$ we have $d_i \geq 2k_i$, choose $j$ such that $k_j > 0$; if not, choose $j$ such that $d_j < 2k_j$ (so certainly $k_j > 0$). In either case the choice of $j$ then means that $\sum_{i \neq j} (d_i - 2k_i) \geq 0$. Define a new tuple $\k' = ({k_1}', \dots, {k_t}')$ by setting
$$
{k_i}' =
\begin{cases}
k_j - 1 & \hbox{if } i = j, \\
k_i & \hbox{if } i \neq j;
\end{cases}
$$
then $|\k'| = k - 1$, and $\k'$ is $\d$-feasible. We have
\begin{eqnarray*}
B_{\d, \k} - B_{\d, \k'}
& = & k(d - k) - k_j(d_j - k_j) - (k - 1)(d - (k - 1)) \\
&   & \quad {} + (k_j - 1)(d_j - (k_j - 1)) \\
& = & dk - k^2 - d_jk_j + {k_j}^2 - dk + d + k^2 - 2k + 1 \\
&   & \quad {} + d_jk_j - d_j - {k_j}^2 + 2k_j - 1 \\
& = & (d - d_j) - 2(k - k_j) \\
& = & \sum_{i \neq j} (d_i - 2k_i) \\
& \geq & 0;
\end{eqnarray*}
thus $B_{\d, k - 1} \leq B_{\d, \k'} \leq B_{\d, \k} = B_{\d, k}$ as required.
\end{proof}

Note that the value $B_{\d, k}$ is unaffected by permutations of the parts of $\d$ (since corresponding permutations may be applied to the parts of the $\d$-feasible tuples $\k$); from now on we shall assume that $\d$ is decreasing. Our next result shows that the value $B_{\d, k}$ is then attained by a decreasing tuple $\k$.

\begin{prop}\label{prop: decreasing tuples are best}
Given a decreasing tuple $\d$ of natural numbers and a non-negative integer $k$ with $k \leq |\d|$, there exists a decreasing $\d$-feasible tuple $\k$ with $|\k| = k$ such that $B_{\d, \k} = B_{\d, k}$.
\end{prop}

\begin{proof}
Write $\d = (d_1, \dots, d_t)$, and take a $\d$-feasible tuple $\k = (k_1, \dots, k_t)$ with $|\k| = k$ such that $B_{\d, \k} = B_{\d, k}$. Suppose $\k$ is not decreasing; then there exist $j_1, j_2 \leq t$ with $j_1 < j_2$ and $k_{j_1} < k_{j_2}$. Define a new tuple $\k' = ({k_1}', \dots, {k_t}')$ by setting
$$
{k_i}' =
\begin{cases}
k_{j_2} & \hbox{if } i = j_1, \\
k_{j_1} & \hbox{if } i = j_2, \\
k_i & \hbox{if } i \neq j_1, j_2;
\end{cases}
$$
then $|\k'| = k$, and as ${k_{j_1}}' = k_{j_2} \leq d_{j_2} \leq d_{j_1}$ and ${k_{j_2}}' = k_{j_1} < k_{j_2} \leq d_{j_2}$ we see that $\k'$ is $\d$-feasible. As
\begin{eqnarray*}
B_{\d, \k} - B_{\d, \k'}
& = & -k_{j_1}(d_{j_1} - k_{j_1}) - k_{j_2}(d_{j_2} - k_{j_2}) + {k_{j_1}}'(d_{j_1} - {k_{j_1}}') \\
&   & \quad {} + {k_{j_2}}'(d_{j_2} - {k_{j_2}}') \\
& = & -d_{j_1}k_{j_1} + {k_{j_i}}^2 - d_{j_2}k_{j_2} + {k_{j_2}}^2 + d_{j_1}k_{j_2} - {k_{j_2}}^2 + d_{j_2}k_{j_1} - {k_{j_1}}^2 \\
& = & (d_{j_1} - d_{j_2})(k_{j_2} - k_{j_1}) \\
& \geq & 0,
\end{eqnarray*}
we have $B_{\d, \k'} \leq B_{\d, \k} = B_{\d, k}$; thus by definition we must have $B_{\d, \k'} = B_{\d, k}$. Iterating this procedure gives the result.
\end{proof}

For small values of $k$ this has the following consequence.

\begin{cor}\label{cor: B values for small k}
Given a decreasing tuple $\d = (d_1, \dots, d_t)$ of natural numbers, write $d = |\d|$; then we have
\begin{itemize}
\item[(i)] $B_{\d, 1} = d - d_1$;
\item[(ii)] $B_{\d, 2} =
\begin{cases}
2d - 2d_1 & \hbox{if } d_1 \geq d_2 + 2, \\
2d - d_1 - d_2 - 2 & \hbox{if } d_1 < d_2 + 2;
\end{cases}$
\item[(iii)] $B_{\d, 3} =
\begin{cases}
3d - 3d_1 & \hbox{if } d_1 \geq d_2 + 4, \\
3d - 2d_1 - d_2 - 4 & \hbox{if } d_2 + 4 > d_1 \geq d_3 + 2, \\
3d - d_1 - d_2 - d_3 - 6 & \hbox{if } d_1 < d_3 + 2.
\end{cases}$
\end{itemize}
\end{cor}

\begin{proof}
By Proposition~\ref{prop: decreasing tuples are best} we need only consider decreasing tuples $\k$. If $|\k| = 1$ we must have $\k = (1, 0, \dots, 0)$; if $|\k| = 2$ we must have $\k = (2, 0, \dots, 0)$ or $(1, 1, 0, \dots, 0)$; if $|\k| = 3$ we must have $\k = (3, 0, \dots, 0)$, $(2, 1, 0, \dots, 0)$ or $(1, 1, 1, 0, \dots, 0)$. The values $B_{\d, \k}$ are then as shown; the inequalities stated in parts (ii) and (iii) are the conditions which must hold for the relevant value $B_{\d, \k}$ to be minimal.
\end{proof}

Another special case is when $\d$ is a $2$-tuple.

\begin{prop}\label{prop: B value when t = 2}
Given a decreasing tuple $\d = (d_1, d_2)$ of natural numbers, write $d = |\d| = d_1 + d_2$; then if $1 \leq k \leq \frac{d}{2}$ we have
$$
B_{\d, k} =
\begin{cases}
d_2k & \hbox{if } d_2 + k \leq \ts{\frac{d}{2}}, \\
\lceil \ts{\frac{1}{2}} d_1d_2 - \ts{\frac{1}{8}} (d - 2k)^2 \rceil & \hbox{if } d_2 + k > \ts{\frac{d}{2}}.
\end{cases}
$$
\end{prop}

\begin{proof}
Given $\k = (k_1, k_2)$ with $k = |\k| = k_1 + k_2$, we have
\begin{eqnarray*}
B_{\d, \k}
& = & k(d - k) - k_1(d_1 - k_1) - k_2(d_2 - k_2) \\
& = & k(d - k) - (k - k_2)((d - d_2) - (k - k_2)) - k_2(d_2 - k_2) \\
& = & dk - k^2 - dk + dk_2 + d_2k - d_2k_2 + k^2 - 2kk_2 + {k_2}^2 - d_2k_2 + {k_2}^2 \\
& = & 2{k_2}^2 - (2k - d + 2d_2)k_2 + d_2k.
\end{eqnarray*}
The minimum value of the quadratic $2x^2 - (2k - d + 2d_2)x + d_2k$ occurs when $x = \frac{1}{4}(2k - d + 2d_2)$. Thus if $d_2 + k \leq \frac{d}{2}$ we take $k_2 = 0$ to give $B_{\d, k} = d_2k$. If instead $d_2 + k > \frac{d}{2}$, we take $k_2$ to be the nearest integer to $\frac{1}{4}(2k - d + 2d_2)$, say $k_2 = \frac{1}{4}(2k - d + 2d_2) + \e$ where $|\e| \leq \frac{1}{2}$, to give $B_{\d, k} = 2\e^2 - \frac{1}{8}(2k - d + 2d_2)^2 + d_2k = 2\e^2 - \frac{1}{8}(2k - d)^2 - \frac{1}{2}d_2(2k - d) - \frac{1}{2}{d_2}^2 + d_2k = \frac{1}{2}d_1d_2 - \frac{1}{8}(d - 2k)^2 + 2\e^2$; since the value must be an integer and $|2\e^2| \leq \frac{1}{2}$, we must have $B_{\d, k} = \lceil \frac{1}{2}d_1d_2 - \frac{1}{8}(d - 2k)^2 \rceil$ as required.
\end{proof}

So far we have been considering a fixed tuple $\d$; we now consider ranging over tuples with fixed value of $|\d|$. Given natural numbers $d$ and $b$, we define
$$
{\mathcal T}_d^b = \{ \d: |\d| = d, \ \d = (d_1, \dots, d_t), \ b \geq d_1 \geq \cdots \geq d_t \};
$$
given additionally a natural number $k$ with $k \leq d$ we define
$$
B_{d,k}^b = \min \{ B_{\d, k} : \d \in {\mathcal T}_d^b \}.
$$
Our result here is the following.

\begin{prop}\label{prop: min B value for varying bounded tuple d}
Given natural numbers $d$, $k$ and $b$ with $k \leq \frac{d}{2}$, we have $B_{d, k}^b = B_{\d_0, k}$, where $\d_0 = (d_1, \dots, d_t)$ with $t = \lceil \frac{d}{b} \rceil$ and $d_1 = \cdots = d_{t-1} = b$, $d_t = d - (t - 1)b$.
\end{prop}

\begin{proof}
Take $\d = (d_1, \dots, d_t) \in {\mathcal T}_d^b$ (for an arbitrary $t$) and suppose $d_{t-1} < b$; let $j$ be minimal with $d_j < b$, so that $j < t$. We shall show that there exists $\d' \in {\mathcal T}_d^b$ such that $\d'$ strictly precedes $\d$ in the standard partial ordering and $B_{\d', k} \leq B_{\d, k}$, from which the result follows.

First suppose $d_j + d_t \leq b$. In this case we define $\d' = ({d_1}', \dots, {d_{t - 1}}')$ by
$$
{d_i}' =
\begin{cases}
d_j + d_t & \hbox{if } i = j, \\
d_i & \hbox{if } i \neq j;
\end{cases}
$$
then $\d' \in {\mathcal T}_d^b$ strictly precedes $\d$. Given a $\d$-feasible tuple $\k = (k_1, \dots, k_t)$ of integers with $|\k| = k$ and $B_{\d, \k} = B_{\d, k}$, we obtain a tuple $\k' = ({k_1}', \dots, {k_{t - 1}}')$ by setting
$$
{k_i}' =
\begin{cases}
k_j + k_t & \hbox{if } i = j, \\
k_i & \hbox{if } i \neq j;
\end{cases}
$$
then $\k'$ is $\d'$-feasible with $|\k'| = k$, and we have
\begin{eqnarray*}
B_{\d, \k} - B_{\d', \k'}
& = & -k_j(d_j - k_j) - k_t(d_t - k_t) + {k_j}'({d_j}' - {k_j}') \\
& = & -k_j(d_j - k_j) - k_t(d_t - k_t) + (k_j + k_t)((d_j + d_t) - (k_j + k_t)) \\
& = & k_j(d_t - k_t) + k_t(d_j - k_j) \\
& \geq & 0.
\end{eqnarray*}
Thus $B_{\d', k} \leq B_{\d', \k'} \leq B_{\d, \k} = B_{\d, k}$ as required.

Now suppose $d_j + d_t > b$; note that as $d_j < b$ we must have $d_t \geq 2$. In this case we define $\d' = ({d_1}', \dots, {d_t}')$ by
$$
{d_i}' =
\begin{cases}
d_j + 1 & \hbox{if } i = j, \\
d_t - 1 & \hbox{if } i = t, \\
d_i & \hbox{if } i \neq j, t;
\end{cases}
$$
then $\d' \in {\mathcal T}_d^b$ strictly precedes $\d$. We must show that $B_{\d', k} \leq B_{\d, k}$.

Take a $\d$-feasible tuple $\k = (k_1, \dots, k_t)$ of integers with $|\k| = k$ and $B_{\d, \k} = B_{\d, k}$; by Proposition~\ref{prop: decreasing tuples are best} we may assume $\k$ is decreasing. Suppose if possible that $k_t = d_t$; then as $k_t > \frac{d_t}{2}$ and $k \leq \frac{d}{2}$, there must exist $a$ with $k_a < \frac{d_a}{2}$. Now if we define $\k' = ({k_1}', \dots, {k_t}')$ by setting
$$
{k_i}' =
\begin{cases}
k_a + 1 & \hbox{if } i = a, \\
k_t - 1 & \hbox{if } i = t, \\
k_i & \hbox{if } i \neq a, t;
\end{cases}
$$
then $\k'$ is $\d$-feasible with $|\k'| = k$, and we have
\begin{eqnarray*}
B_{\d, \k} - B_{\d, \k'}
& = & -k_a(d_a - k_a) - k_t(d_t - k_t) + {k_a}'(d_a - {k_a}') + {k_t}'(d_t - {k_t}') \\
& = & -k_a(d_a - k_a) + (k_a + 1)(d_a - (k_a + 1)) + (d_t - 1) \\
& = & (d_a - 2k_a - 1) + (d_t - 1) \\
& \geq & 1,
\end{eqnarray*}
contrary to the assumption that $B_{\d, \k} = B_{\d, k}$. Thus we must have $k_t < d_t$, and hence $k_t \leq {d_t}'$; so $\k$ is $\d'$-feasible. Now
\begin{eqnarray*}
B_{\d, \k} - B_{\d', \k}
& = & -k_j(d_j - k_j) - k_t(d_t - k_t) + k_j({d_j}' - k_j) + k_t({d_t}' - k_t) \\
& = & k_j({d_j}' - d_j) + k_t({d_t}' - d_t) \\
& = & k_j - k_t \\
& \geq & 0,
\end{eqnarray*}
so $B_{\d', k} \leq B_{\d', \k} \leq B_{\d, \k} = B_{\d, k}$ as required.
\end{proof}

\section{Reduction of quadruples requiring consideration}\label{sect: reduction of quadruples}

We now return to the situation where we have a simple algebraic group $G$ over an algebraically closed field $K$ of characteristic $p$, and a dominant weight $\lambda$; as before we write $V = L(\lambda)$ and $d = \dim V$. We take an integer $k$ satisfying $1 \leq k \leq \frac{d}{2}$, and suppose that $(G, \lambda, p, k)$ is a large quadruple; we write $X = \Gk(V)$. In this chapter we wish to show that a large quadruple not listed in Table~\ref{table: large triple and first quadruple non-TGS} or \ref{table: large higher quadruple non-TGS} has TGS; in this section we shall substantially reduce the number of large quadruples requiring consideration.

In order to show that a large quadruple has TGS, we shall prove that it satisfies both $\ssdiamcon$ and $\udiamcon$ of Section~\ref{sect: conditions}. We therefore require information on $\codim C_X(g)$ for $g$ either semisimple or unipotent. Our first result in this section links this codimension to the work on tuples in the previous section.

\begin{prop}\label{prop: codim formula for quadruples}
Let $(G, \lambda, p, k)$ be a quadruple, and take $g \in G$ either semi-\linebreak simple or unipotent; define the tuple $\d = (d_1, \dots, d_t)$ as follows.
\begin{itemize}
\item[(i)] If $g = s$ is semisimple, let $\kappa_1, \dots, \kappa_t$ be the distinct eigenvalues of $s$ in its action on $V$; for $i \in [1, t]$ let $d_i = \dim V_{\kappa_i}(s)$.
\item[(ii)] If $g = u$ is unipotent, let $1^{b_1}, 2^{b_2}, \dots, t^{b_t}$ be the sizes of the Jordan blocks of $u$ in its action on $V$, where $b_t > 0$; for $i \in [1, t]$ let $d_i = b_i + \cdots + b_t$.
\end{itemize}
Then $\codim C_{\Gk(V)}(g) = B_{\d, k}$.
\end{prop}

\begin{proof}
Suppose $\bar V$ is a $k$-dimensional subspace of $V$ which is fixed by $g$. If $g = s$ we have $V = \bigoplus_{i = 1}^t V_{\kappa_i}(s)$, and as $\bar V$ must have a basis of eigenvectors for $s$ we see that $\bar V = \bigoplus_{i = 1}^t \bar V_{\kappa_i}(s)$; for $i \in [1, t]$ let $k_i = \dim \bar V_{\kappa_i}(s)$. If instead $g = u$ the sizes of the Jordan blocks of $u$ in its action on $\bar V$ must be $1^{a_1}, 2^{a_2}, \dots, t^{a_t}$, where for $i \in [1, t]$ we have $a_i + \cdots + a_t \leq b_i + \cdots + b_t$; for $i \in [1, t]$ let $k_i = a_i + \cdots + a_t$. In either case write $\k = (k_1, \dots, k_t)$; then $\k$ is $\d$-feasible and $|\k| = k$.

For a fixed $\d$-feasible tuple $\k$, let ${\mathcal V}_{\k}$ be the variety of such subspaces $\bar V$. If $g = s$ then clearly $\dim {\mathcal V}_{\k} = \sum_{i = 1}^t k_i(d_i - k_i)$; if instead $g = u$ then $\dim {\mathcal V}_{\k} = \sum_{j=1}^t a_j((d_1 + \cdots + d_j) - (k_1 + \cdots + k_j)) = \sum_{j = 1}^t a_j \sum_{i = 1}^j (d_i - k_i) = \sum_{i = 1}^t (d_i - k_i) \sum_{j = i}^t a_j = \sum_{i = 1}^t k_i(d_i - k_i)$. Thus in either case $\codim {\mathcal V}_{\k} = B_{\d, \k}$. Taking the union of the varieties ${\mathcal V}_{\k}$ as $\k$ ranges through the finite set of $\d$-feasible tuples with $|\k| = k$, we see that $\codim C_{\Gk(V)}(g) = B_{\d, k}$.
\end{proof}

Note that if $g = u$ then the tuple $\d$ is automatically decreasing, while if $g = s$ it is harmless to assume this. Our next result in this section is the following.

\begin{prop}\label{prop: large quadruple having TGS implies same for increased k}
Let $(G, \lambda, p, k)$ be a large quadruple, and take $k' \in [k, \frac{d}{2}]$; if $(G, \lambda, p, k)$ satisfies $\ssdiamcon$ or $\udiamcon$, so does $(G, \lambda, p, k')$.
\end{prop}

\begin{proof}
If $g$ is either semisimple or unipotent, by Propositions~\ref{prop: B value increases with k} and \ref{prop: codim formula for quadruples} we have $\codim C_{\G{k'}(V)}(g) = B_{\d, k'} \geq B_{\d, k} = \codim C_{\Gk(V)}(g)$. Thus if $(G, \lambda, p, k)$ satisfies $\ssdiamcon$ or $\udiamcon$, then whenever $g$ is $s \in G_{(r)}$ for some $r \in \P'$ or $u \in G_{(p)}$ respectively we have $\codim C_{\Gk(V)}(g) > \dim g^G$, and so $\codim C_{\G{k'}(V)}(g) > \dim g^G$, so that $(G, \lambda, p, k')$ also satisfies $\ssdiamcon$ or $\udiamcon$ respectively.
\end{proof}

As a consequence we are able to prove the following.

\begin{prop}\label{prop: large quadruples with TGS}
If $(G, \lambda, p)$ is a triple which does not appear in Tables~\ref{table: large triple and first quadruple non-TGS}, \ref{table: small classical triple and first quadruple generic stab} or \ref{table: small exceptional triple and first quadruple generic stab}, then any associated quadruple $(G, \lambda, p, k)$ has TGS.
\end{prop}

\begin{proof}
Let $(G, \lambda, p)$ be a triple as in the statement; then it must be large, and by Propositions~\ref{prop: p-restricted large triples with TGS} and \ref{prop: large triples which are not p-restricted with TGS} it satisfies $\ssdiamevcon$ and $\udiamcon$. By Propositions~\ref{prop: ssdiamevcon for triples implies ssdiamcon for first quadruples} and \ref{prop: udiamcon for triples implies udiamcon for first quadruples} the associated first quadruple $(G, \lambda, p, 1)$ satisfies $\ssdiamcon$ and $\udiamcon$. Proposition~\ref{prop: large quadruple having TGS implies same for increased k} therefore shows that any associated quadruple $(G, \lambda, p, k)$ also satisfies $\ssdiamcon$ and $\udiamcon$, and so has TGS.
\end{proof}

We are therefore left only to consider quadruples associated to triples appearing in Tables~\ref{table: large triple and first quadruple non-TGS}, \ref{table: small classical triple and first quadruple generic stab} and \ref{table: small exceptional triple and first quadruple generic stab}. We list in Table~\ref{table: remaining large quadruples} the remaining large quadruples $(G, \lambda, p, k)$ which we shall show have TGS. Note that the final column of Table~\ref{table: remaining large quadruples} is headed `$k_0$'; the entry here is the least value of $k$ for which we claim that the quadruple $(G, \lambda, p, k)$ has TGS. By Proposition~\ref{prop: large quadruple having TGS implies same for increased k} it suffices to prove that $(G, \lambda, p, k_0)$ satisfies $\ssdiamcon$ and $\udiamcon$; to do this, by Proposition~\ref{prop: codim formula for quadruples} it suffices to establish certain inequalities involving $B_{\d, k_0}$ for appropriate tuples $\d$.

\begin{table}
\caption{Remaining large quadruples}\label{table: remaining large quadruples}
\tabcapsp
$$
\begin{small}
\begin{array}{|c|c|c|c|c|c|c|c|c|c|c|c|c|c|c|c|}
\cline{1-5} \cline{7-11} \cline{13-16}
G      & \lambda                 & \ell      & p          & k_0 & \ptw & G      & \lambda             & \ell      & p          & k_0 & \ptw & G   & \lambda  & p          & k_0 \mmtbs \\
\cline{1-5} \cline{7-11} \cline{13-16}
A_\ell & 2\omega_1               & {} \geq 3 & {} \geq 3  & 3   &      & B_\ell & 2\omega_1           & {} \geq 2 & {} \geq 3  & 2   &      & E_6 & \omega_1 & \hbox{any} & 4   \mmtbs \\
       & \omega_2                & 5         & \hbox{any} & 4   &      &        & \omega_2            & {} \geq 3 & {} \geq 3  & 2   &      &     & \omega_2 & \hbox{any} & 2   \mmtbs \\
\cline{13-16}
       & \omega_2                & {} \geq 6 & \hbox{any} & 3   &      &        & \omega_2            & 3         & 2          & 3   &      & E_7 & \omega_1 & \hbox{any} & 2   \mmtbs \\
       & \omega_3                & 5         & \hbox{any} & 3   &      &        & \omega_2            & {} \geq 4 & 2          & 2   &      &     & \omega_7 & \hbox{any} & 3   \mmtbs \\
\cline{13-16}
       & \omega_3                & 6, 7, 8   & \hbox{any} & 2   &      &        & \omega_1 + \omega_2 & 2         & 5          & 2   &      & E_8 & \omega_8 & \hbox{any} & 2   \mmtbs \\
\cline{13-16}
       & 3\omega_1               & 2         & {} \geq 5  & 2   &      &        & 2\omega_2           & 2         & {} \geq 3  & 2   &      & F_4 & \omega_1 & {} \geq 3  & 2   \mmtbs \\
       & 4\omega_1               & 1         & {} \geq 5  & 2   &      &        & \omega_4            & 4         & \hbox{any} & 4   &      &     & \omega_1 & 2          & 3   \mmtbs \\
       & 2\omega_2               & 3         & {} \geq 3  & 2   &      &        & \omega_5            & 5         & \hbox{any} & 3   &      &     & \omega_4 & \hbox{any} & 3   \mmtbs \\
\cline{13-16}
       & \omega_4                & 7         & \hbox{any} & 2   &      &        & \omega_6            & 6         & \hbox{any} & 2   &      & G_2 & \omega_2 & {} \neq 3  & 2   \mmtbs \\
\cline{7-11} \cline{13-16}
       & \omega_1 + \omega_2     & 3         & 3          & 2   &      & C_\ell & 2\omega_1           & {} \geq 3 & {} \geq 3  & 2   & \multicolumn{5}{c}{}                     \mmtbs \\
       & \omega_1 + \omega_\ell  & {} \geq 2 & \hbox{any} & 2   &      &        & \omega_2            & 3         & \hbox{any} & 3   & \multicolumn{5}{c}{}                     \mmtbs \\
       & \omega_1 + q\omega_1    & {} \geq 2 & < \infty   & 2   &      &        & \omega_2            & {} \geq 4 & \hbox{any} & 2   & \multicolumn{5}{c}{}                     \mmtbs \\
       & \omega_1 + q\omega_\ell & {} \geq 2 & < \infty   & 2   &      &        & \omega_3            & 4         & 3          & 2   & \multicolumn{5}{c}{}                     \mmtbs \\
\cline{1-5}
D_\ell & 2\omega_1               & {} \geq 4 & {} \geq 3  & 2   &      &        & \omega_3            & 3         & {} \geq 3  & 3   & \multicolumn{5}{c}{}                     \mmtbs \\
       & \omega_2                & {} \geq 4 & \hbox{any} & 2   &      &        & \omega_4            & 4         & {} \geq 3  & 2   & \multicolumn{5}{c}{}                     \mmtbs \\
       & \omega_5                & 5         & \hbox{any} & 5   &      &        & \omega_4            & 4         & 2          & 4   & \multicolumn{5}{c}{}                     \mmtbs \\
       & \omega_6                & 6         & \hbox{any} & 3   &      &        & \omega_5            & 5         & 2          & 3   & \multicolumn{5}{c}{}                     \mmtbs \\
       & \omega_\ell             & 7, 8      & \hbox{any} & 2   &      &        & \omega_6            & 6         & 2          & 2   & \multicolumn{5}{c}{}                     \mmtbs \\
\cline{1-5} \cline{7-11}
\end{array}
\end{small}
$$
\end{table}

It will be seen that in some rows in Table~\ref{table: remaining large quadruples} the rank $\ell$ takes a (small) finite number of values, while in others it takes infinitely many values. We shall find that the two require different strategies. In Section~\ref{sect: large quadruple individual cases} we shall treat the former cases, taking each individually; in Section~\ref{sect: large quadruple infinite families} we shall then consider the infinite families.

\section{Analysis of individual cases}\label{sect: large quadruple individual cases}

We continue with the notation of the previous section. For a given quadruple $(G, \lambda, p, k_0)$ listed in Table~\ref{table: remaining large quadruples} which does not occur in an infinite family, we wish to show that it satisfies $\ssdiamcon$ and $\udiamcon$; to do this we must take elements $s \in G_{(r)}$ for some $r \in \P'$ and $u \in G_{(p)}$ and consider the codimensions of their fixed point varieties, which by Proposition~\ref{prop: codim formula for quadruples} we know to be equal to $B_{\d, k_0}$ for the appropriate tuples $\d$. We assume $s$ lies in $T$, and write $\Phi(s) = \{ \alpha \in \Phi : \alpha(s) = 1 \}$, so that $C_G(s)^\circ = \langle T, X_\alpha : \alpha \in \Phi(s) \rangle$; we take $\kappa \in K^*$. Our approach in this section is based on that employed in Section~\ref{sect: large triple further analysis}, as follows.

We start by giving the weight table. We then take a standard subsystem $\Psi$ of $\Phi$ (initially we take $\Psi$ of rank $1$) and give the $\Psi$-net table. This provides lower bounds $c(\Psi)_{ss}$ on $\codim V_\kappa(s)$ for any $s \in G_{(r)}$ with $\Phi(s)$ disjoint from $\Psi$, and $c(\Psi)_u$ on $\codim C_V(u_\Psi)$; these lower bounds may depend on $r$ or $p$ respectively. Write $c$ for either $c(\Psi)_{ss}$ or $c(\Psi)_u$. At this point we may not know precisely the tuple $\d$ associated to the element $s$ or $u_\Psi$. However, we have $d - d_1 \geq c$, and hence $d_1 \leq d - c$; thus if we write $b = d - c$, then $\d \in {\mathcal T}_d^b$. Proposition~\ref{prop: min B value for varying bounded tuple d} now shows that $B_{\d, k_0} \geq B_{\d_0, k_0}$, where $\d_0$ has all parts equal to $b$ except possibly the last.

We therefore compute $B_{\d_0, k_0}$. For this, often we have $c \leq \frac{d}{2}$, so that $b \geq \frac{d}{2}$ and hence $\d_0 = (b, c)$, in which case we can employ Proposition~\ref{prop: B value when t = 2}; if $c > \frac{d}{2}$ and $k_0 \leq 3$ we can employ Corollary~\ref{cor: B values for small k}; if $c > \frac{d}{2}$ and $k_0 \geq 4$ we may need to consider the various possibilities for $\k$, but Proposition~\ref{prop: decreasing tuples are best} means that we need only treat those which are decreasing.

As in Section~\ref{sect: large triple further analysis}, if $c = c(\Psi)_{ss}$ we may henceforth ignore all semisimple classes $s^G$ with $\dim s^G < B_{\d_0, k_0}$ such that $\Psi$ is disjoint from $\Phi(s)$, while if $c = c(\Psi)_u$ we may henceforth ignore all unipotent classes $u^G$ in $G_{(p)}$ with $\dim u^G < B_{\d_0, k_0}$ such that the closure of $u^G$ contains ${u_\Psi}^G$. The remaining semisimple classes $s^G$ satisfy $\dim s^G \geq B_{\d_0, k_0}$, i.e., $|\Phi(s)| \leq M - B_{\d_0, k_0}$; we identify a larger subsystem $\Psi$ such that each of these remaining $\Phi(s)$ has a conjugate of $\Psi$ disjoint from it, and such that all remaining unipotent classes in $G_{(p)}$ have ${u_\Psi}^G$ in their closure. We take this $\Psi$ and repeat the procedure to obtain improved lower bounds $c(\Psi)_{ss}$ and $c(\Psi)_u$, and hence larger values $B_{\d_0, k_0}$. Continuing thus, we eventually show that the quadruple $(G, \lambda, p, k_0)$ satisfies $\ssdiamcon$ and $\udiamcon$.

In two cases it is not true that the set of weights (ignoring multiplicities) appearing in a given $\Psi$-net is that of a single Weyl $G_\Psi$-module $W_{G_\Psi}(\bar\nu)$. In the cases concerned we have $G = E_7$, $\lambda = \omega_1$ and $G = E_8$, $\lambda = \omega_8$; the subsystem $\Psi$ concerned is $\langle \alpha_1, \alpha_4 \rangle$ of type ${A_1}^2$, and the $\Psi$-net consists of the weights in the Lie algebra of $G_\Psi$. The corresponding entry in the first column of the $\Psi$-net table is `$2\bom_1/2\bom_4$', indicating that the weights concerned are those lying in the union of the sets of weights of the Weyl $G_\Psi$-modules $W_{G_\Psi}(2\bom_1)$ and $W_{G_\Psi}(2\bom_4)$.

We now work through the quadruples. As in sections~\ref{sect: large triple weight string analysis} and \ref{sect: large triple further analysis}, we shall begin with those in which all roots in $\Phi$ have the same length.

\begin{prop}\label{prop: A_5, omega_2, k geq 4, nets}
Let $G = A_5$ and $\lambda = \omega_2$; then for $k \in [4, \frac{d}{2}]$ the quadruple $(G, \lambda, p, k)$ satisfies $\ssdiamcon$ and $\udiamcon$.
\end{prop}

\begin{proof}
The weight table is as follows.
$$
\begin{array}{|*4{>{\ss}c|}}
\hline
i & \mu & |W.\mu| & m_\mu \\
\hline
1 & \omega_2 & 15 & 1 \\
\hline
\end{array}
$$
We have $M = 30$, $M_3 = 24$ and $M_2 = 18$; we take $k_0 = 4$.

Take $\Psi = \langle \alpha_1 \rangle$ of type $A_1$. The $\Psi$-net table is as follows.
$$
\begin{array}{|*5{>{\ss}c|}}
\hline
\multicolumn{2}{|>{\ss}c|}{\Psi\mathrm{-nets}} & & \multicolumn{1}{|>{\ss}c|}{c(s)} & \multicolumn{1}{|>{\ss}c|}{c(u_\Psi)} \\
\cline{1-2} \cline{4-5}
    \bar\nu    & n_1 & m & r \geq 2 & p \geq 2 \\
\hline
     \bom_1    &  2  & 4 &     4    &     4    \\
       0       &  1  & 7 &          &          \\
\hline
\multicolumn{3}{c|}{}    &     4    &     4    \\
\cline{4-5}
\end{array}
$$
Thus $c(\Psi)_{ss} = c(\Psi)_u = 4$, so we take $\d_0 = (11, 4)$; using Proposition~\ref{prop: B value when t = 2} we then have $B_{\d_0, 4} = 16 > 10 = \dim {u_\Psi}^G$. We need only consider semisimple classes $s^G$ with $|\Phi(s)| \leq M - 16 = 14 < 20 = m_{{A_1}^2}$, each of which has a subsystem of type ${A_1}^2$ disjoint from $\Phi(s)$, and unipotent classes of dimension at least $16$, each of which has the class ${A_1}^2$ in its closure by Lemma~\ref{lem: various classes in classical groups by dim}(i).

Now take $\Psi = \langle \alpha_1, \alpha_3 \rangle$ of type ${A_1}^2$. The $\Psi$-net table is as follows.
$$
\begin{array}{|*5{>{\ss}c|}}
\hline
\multicolumn{2}{|>{\ss}c|}{\Psi\mathrm{-nets}} & & \multicolumn{1}{|>{\ss}c|}{c(s)} & \multicolumn{1}{|>{\ss}c|}{c(u_\Psi)} \\
\cline{1-2} \cline{4-5}
     \bar\nu     & n_1 & m & r \geq 2 & p \geq 2 \\
\hline
 \bom_1 + \bom_3 &  4  & 1 &     2    &     2    \\
      \bom_1     &  2  & 2 &     2    &     2    \\
      \bom_3     &  2  & 2 &     2    &     2    \\
        0        &  1  & 3 &          &          \\
\hline
\multicolumn{3}{c|}{}      &     6    &     6    \\
\cline{4-5}
\end{array}
$$
Thus $c(\Psi)_{ss} = c(\Psi)_u = 6$, so we take $\d_0 = (9, 6)$; using Proposition~\ref{prop: B value when t = 2} we then have $B_{\d_0, 4} = 21 > M_2 > 16 = \dim {u_\Psi}^G$. We may therefore assume from now on that $r \geq 3$, and that $p \geq 3$ when we treat unipotent classes. We need only consider semisimple classes $s^G$ with $|\Phi(s)| \leq M - 21 = 9$, each of which has a subsystem of type $A_2A_1$ disjoint from $\Phi(s)$, and unipotent classes of dimension at least $21$, each of which has the class $A_2A_1$ in its closure by Lemma~\ref{lem: various classes in classical groups by dim}(iv).

Now take $\Psi = \langle \alpha_1, \alpha_2, \alpha_4 \rangle$ of type $A_2A_1$. The $\Psi$-net table is as follows.
$$
\begin{array}{|*5{>{\ss}c|}}
\hline
\multicolumn{2}{|>{\ss}c|}{\Psi\mathrm{-nets}} & & \multicolumn{1}{|>{\ss}c|}{c(s)} & \multicolumn{1}{|>{\ss}c|}{c(u_\Psi)} \\
\cline{1-2} \cline{4-5}
     \bar\nu     & n_1 & m & r \geq 3 & p \geq 3 \\
\hline
 \bom_1 + \bom_4 &  6  & 1 &     4    &     4    \\
      \bom_1     &  3  & 1 &     2    &     2    \\
      \bom_2     &  3  & 1 &     2    &     2    \\
      \bom_4     &  2  & 1 &     1    &     1    \\
        0        &  1  & 1 &          &          \\
\hline
\multicolumn{3}{c|}{}       &     9    &     9    \\
\cline{4-5}
\end{array}
$$
Thus $c(\Psi)_{ss} = c(\Psi)_u = 9$, so we take $\d_0 = (6, 6, 3)$; according as $\k_0 = (4, 0, 0)$, $(3, 1, 0)$, $(2, 2, 0)$ or $(2, 1, 1)$ we have $B_{\d_0, \k_0} = 36$, $30$, $28$ or $29$, so $B_{\d_0, 4} = 28 > M_3 > 22 = \dim {u_\Psi}^G$. We may therefore assume from now on that $r \geq 5$, and that $p \geq 5$ when we treat unipotent classes. We need only consider semisimple classes $s^G$ with $|\Phi(s)| \leq M - 28 = 2$, each of which has a subsystem of type $A_4$ disjoint from $\Phi(s)$, and unipotent classes of dimension at least $28$, each of which has the class $A_4$ in its closure by Lemma~\ref{lem: various classes in A_ell for fixed ell}(ii).

Now take $\Psi = \langle \alpha_1, \alpha_2, \alpha_3, \alpha_4 \rangle$ of type $A_4$. The $\Psi$-net table is as follows.
$$
\begin{array}{|*5{>{\ss}c|}}
\hline
\multicolumn{2}{|>{\ss}c|}{\Psi\mathrm{-nets}} & & \multicolumn{1}{|>{\ss}c|}{c(s)} & \multicolumn{1}{|>{\ss}c|}{c(u_\Psi)} \\
\cline{1-2} \cline{4-5}
     \bar\nu     & n_1 & m & r \geq 5 & p \geq 5 \\
\hline
      \bom_1     &  5  & 1 &     4    &     4    \\
      \bom_2     & 10  & 1 &     8    &     8    \\
\hline
\multicolumn{3}{c|}{}      &    12    &    12    \\
\cline{4-5}
\end{array}
$$
Thus $c(\Psi)_{ss} = c(\Psi)_u = 12$, so we take $\d_0 = (3, 3, 3, 3, 3)$; according as $\k_0 = (3, 1, 0, 0, 0)$, $(2, 2, 0, 0, 0)$, $(2, 1, 1, 0, 0)$ or $(1, 1, 1, 1, 0)$ we have $B_{\d_0, \k_0} = 42$, $40$, $38$ or $36$, so $B_{\d_0, 4} = 36 > M$. Therefore if $k \in [4, \frac{d}{2}]$ the quadruple $(G, \lambda, p, k)$ satisfies $\ssdiamcon$ and $\udiamcon$.
\end{proof}

\begin{prop}\label{prop: A_5, omega_3, k geq 3, nets}
Let $G = A_5$ and $\lambda = \omega_3$; then for $k \in [3, \frac{d}{2}]$ the quadruple $(G, \lambda, p, k)$ satisfies $\ssdiamcon$ and $\udiamcon$.
\end{prop}

\begin{proof}
The weight table is as follows.
$$
\begin{array}{|*4{>{\ss}c|}}
\hline
i & \mu & |W.\mu| & m_\mu \\
\hline
1 & \omega_3 & 20 & 1 \\
\hline
\end{array}
$$
We have $M = 30$, $M_3 = 24$ and $M_2 = 18$; we take $k_0 = 3$.

Take $\Psi = \langle \alpha_1 \rangle$ of type $A_1$. The $\Psi$-net table is as follows.
$$
\begin{array}{|*5{>{\ss}c|}}
\hline
\multicolumn{2}{|>{\ss}c|}{\Psi\mathrm{-nets}} & & \multicolumn{1}{|>{\ss}c|}{c(s)} & \multicolumn{1}{|>{\ss}c|}{c(u_\Psi)} \\
\cline{1-2} \cline{4-5}
    \bar\nu    & n_1 & m & r \geq 2 & p \geq 2 \\
\hline
     \bom_1    &  2  & 6 &     6    &     6    \\
       0       &  1  & 8 &          &          \\
\hline
\multicolumn{3}{c|}{}    &     6    &     6    \\
\cline{4-5}
\end{array}
$$
Thus $c(\Psi)_{ss} = c(\Psi)_u = 6$, so we take $\d_0 = (14, 6)$; using Proposition~\ref{prop: B value when t = 2} we then have $B_{\d_0, 3} = 18 > 10 = \dim {u_\Psi}^G$. We need only consider semisimple classes $s^G$ with $|\Phi(s)| \leq M - 18 = 12 < 20 = m_{{A_1}^2}$, each of which has a subsystem of type ${A_1}^2$ disjoint from $\Phi(s)$, and unipotent classes of dimension at least $18$, each of which has the class ${A_1}^2$ in its closure by Lemma~\ref{lem: various classes in classical groups by dim}(i).

Now take $\Psi = \langle \alpha_1, \alpha_3 \rangle$ of type ${A_1}^2$. The $\Psi$-net table is as follows.
$$
\begin{array}{|*5{>{\ss}c|}}
\hline
\multicolumn{2}{|>{\ss}c|}{\Psi\mathrm{-nets}} & & \multicolumn{1}{|>{\ss}c|}{c(s)} & \multicolumn{1}{|>{\ss}c|}{c(u_\Psi)} \\
\cline{1-2} \cline{4-5}
     \bar\nu     & n_1 & m & r \geq 2 & p \geq 2 \\
\hline
 \bom_1 + \bom_3 &  4  & 2 &     4    &     4    \\
      \bom_1     &  2  & 2 &     2    &     2    \\
      \bom_3     &  2  & 2 &     2    &     2    \\
        0        &  1  & 4 &          &          \\
\hline
\multicolumn{3}{c|}{}      &     8    &     8    \\
\cline{4-5}
\end{array}
$$
Thus $c(\Psi)_{ss} = c(\Psi)_u = 8$, so we take $\d_0 = (12, 8)$; using Proposition~\ref{prop: B value when t = 2} we then have $B_{\d_0, 3} = 24 > M_2 > 16 = \dim {u_\Psi}^G$. We may therefore assume from now on that $r \geq 3$, and that $p \geq 3$ when we treat unipotent classes. We need only consider semisimple classes $s^G$ with $|\Phi(s)| \leq M - 24 = 6$, each of which has a subsystem of type $A_2A_1$ disjoint from $\Phi(s)$, and unipotent classes of dimension at least $24$, each of which has the class $A_2A_1$ in its closure by Lemma~\ref{lem: various classes in classical groups by dim}(iv).

Now take $\Psi = \langle \alpha_1, \alpha_2, \alpha_4 \rangle$ of type $A_2A_1$. The $\Psi$-net table is as follows.
$$
\begin{array}{|*5{>{\ss}c|}}
\hline
\multicolumn{2}{|>{\ss}c|}{\Psi\mathrm{-nets}} & & \multicolumn{1}{|>{\ss}c|}{c(s)} & \multicolumn{1}{|>{\ss}c|}{c(u_\Psi)} \\
\cline{1-2} \cline{4-5}
     \bar\nu     & n_1 & m & r \geq 3 & p \geq 3 \\
\hline
 \bom_1 + \bom_4 &  6  & 1 &     4    &     4    \\
 \bom_2 + \bom_4 &  6  & 1 &     4    &     4    \\
      \bom_1     &  3  & 1 &     2    &     2    \\
      \bom_2     &  3  & 1 &     2    &     2    \\
        0        &  1  & 2 &          &          \\
\hline
\multicolumn{3}{c|}{}      &    12    &    12    \\
\cline{4-5}
\end{array}
$$
Thus $c(\Psi)_{ss} = c(\Psi)_u = 12$, so we take $\d_0 = (8, 8, 4)$; using Corollary~\ref{cor: B values for small k} we then have $B_{\d_0, 3} = 32 > M$. Therefore if $k \in [3, \frac{d}{2}]$ the quadruple $(G, \lambda, p, k)$ satisfies $\ssdiamcon$ and $\udiamcon$.
\end{proof}

\begin{prop}\label{prop: A_6, omega_3, k geq 2, nets}
Let $G = A_6$ and $\lambda = \omega_3$; then for $k \in [2, \frac{d}{2}]$ the quadruple $(G, \lambda, p, k)$ satisfies $\ssdiamcon$ and $\udiamcon$.
\end{prop}

\begin{proof}
The weight table is as follows.
$$
\begin{array}{|*4{>{\ss}c|}}
\hline
i & \mu & |W.\mu| & m_\mu \\
\hline
1 & \omega_3 & 35 & 1 \\
\hline
\end{array}
$$
We have $M = 42$, $M_3 = 32$ and $M_2 = 24$; we take $k_0 = 2$.

Take $\Psi = \langle \alpha_1 \rangle$ of type $A_1$. The $\Psi$-net table is as follows.
$$
\begin{array}{|*5{>{\ss}c|}}
\hline
\multicolumn{2}{|>{\ss}c|}{\Psi\mathrm{-nets}} & & \multicolumn{1}{|>{\ss}c|}{c(s)} & \multicolumn{1}{|>{\ss}c|}{c(u_\Psi)} \\
\cline{1-2} \cline{4-5}
    \bar\nu    & n_1 &  m & r \geq 2 & p \geq 2 \\
\hline
     \bom_1    &  2  & 10 &    10    &    10    \\
       0       &  1  & 15 &          &          \\
\hline
\multicolumn{3}{c|}{}     &    10    &    10    \\
\cline{4-5}
\end{array}
$$
Thus $c(\Psi)_{ss} = c(\Psi)_u = 10$, so we take $\d_0 = (25, 10)$; using Proposition~\ref{prop: B value when t = 2} we then have $B_{\d_0, 2} = 20 > 12 = \dim {u_\Psi}^G$. We need only consider semisimple classes $s^G$ with $|\Phi(s)| \leq M - 20 = 22 < 30 = m_{{A_1}^2}$, each of which has a subsystem of type ${A_1}^2$ disjoint from $\Phi(s)$, and unipotent classes of dimension at least $20$, each of which has the class ${A_1}^2$ in its closure by Lemma~\ref{lem: various classes in classical groups by dim}(i).

Now take $\Psi = \langle \alpha_1, \alpha_3 \rangle$ of type ${A_1}^2$. The $\Psi$-net table is as follows.
$$
\begin{array}{|*5{>{\ss}c|}}
\hline
\multicolumn{2}{|>{\ss}c|}{\Psi\mathrm{-nets}} & & \multicolumn{1}{|>{\ss}c|}{c(s)} & \multicolumn{1}{|>{\ss}c|}{c(u_\Psi)} \\
\cline{1-2} \cline{4-5}
     \bar\nu     & n_1 & m & r \geq 2 & p \geq 2 \\
\hline
 \bom_1 + \bom_3 &  4  & 3 &     6    &     6    \\
      \bom_1     &  2  & 4 &     4    &     4    \\
      \bom_3     &  2  & 4 &     4    &     4    \\
        0        &  1  & 7 &          &          \\
\hline
\multicolumn{3}{c|}{}      &    14    &    14    \\
\cline{4-5}
\end{array}
$$
Thus $c(\Psi)_{ss} = c(\Psi)_u = 14$, so we take $\d_0 = (21, 14)$; using Proposition~\ref{prop: B value when t = 2} we then have $B_{\d_0, 2} = 28 > M_2 > 20 = \dim {u_\Psi}^G$. We may therefore assume from now on that $r \geq 3$, and that $p \geq 3$ when we treat unipotent classes. We need only consider semisimple classes $s^G$ with $|\Phi(s)| \leq M - 28 = 14 < 18 = m_{A_2}$, each of which has a subsystem of type $A_2$ disjoint from $\Phi(s)$, and unipotent classes of dimension at least $28$, each of which has the class $A_2$ in its closure by Lemma~\ref{lem: various classes in classical groups by dim}(iii).

Now take $\Psi = \langle \alpha_1, \alpha_2 \rangle$ of type $A_2$. The $\Psi$-net table is as follows.
$$
\begin{array}{|*5{>{\ss}c|}}
\hline
\multicolumn{2}{|>{\ss}c|}{\Psi\mathrm{-nets}} & & \multicolumn{1}{|>{\ss}c|}{c(s)} & \multicolumn{1}{|>{\ss}c|}{c(u_\Psi)} \\
\cline{1-2} \cline{4-5}
     \bar\nu     & n_1 & m & r \geq 3 & p \geq 3 \\
\hline
      \bom_1     &  3  & 6 &    12    &    12    \\
      \bom_2     &  3  & 4 &     8    &     8    \\
        0        &  1  & 5 &          &          \\
\hline
\multicolumn{3}{c|}{}      &    20    &    20    \\
\cline{4-5}
\end{array}
$$
Thus $c(\Psi)_{ss} = c(\Psi)_u = 20$, so we take $\d_0 = (15, 15, 5)$; using Corollary~\ref{cor: B values for small k} we then have $B_{\d_0, 2} = 38 > M_3 > 22 = \dim {u_\Psi}^G$. We may therefore assume from now on that $r \geq 5$, and that $p \geq 5$ when we treat unipotent classes. We need only consider semisimple classes $s^G$ with $|\Phi(s)| \leq M - 38 = 4$, each of which has a subsystem of type $A_3$ disjoint from $\Phi(s)$, and unipotent classes of dimension at least $38$, each of which has the class $A_3$ in its closure by Lemma~\ref{lem: various classes in classical groups by dim}(vi).

Now take $\Psi = \langle \alpha_1, \alpha_2, \alpha_3 \rangle$ of type $A_3$. The $\Psi$-net table is as follows.
$$
\begin{array}{|*5{>{\ss}c|}}
\hline
\multicolumn{2}{|>{\ss}c|}{\Psi\mathrm{-nets}} & & \multicolumn{1}{|>{\ss}c|}{c(s)} & \multicolumn{1}{|>{\ss}c|}{c(u_\Psi)} \\
\cline{1-2} \cline{4-5}
     \bar\nu     & n_1 & m & r \geq 5 & p \geq 5 \\
\hline
      \bom_1     &  4  & 3 &     9    &     9    \\
      \bom_2     &  6  & 3 &    12    &    12    \\
      \bom_3     &  4  & 1 &     3    &     3    \\
        0        &  1  & 1 &          &          \\
\hline
\multicolumn{3}{c|}{}      &    24    &    24    \\
\cline{4-5}
\end{array}
$$
Thus $c(\Psi)_{ss} = c(\Psi)_u = 24$, so we take $\d_0 = (11, 11, 11, 2)$; using Corollary~\ref{cor: B values for small k} we then have $B_{\d_0, 2} = 46 > M$. Therefore if $k \in [2, \frac{d}{2}]$ the quadruple $(G, \lambda, p, k)$ satisfies $\ssdiamcon$ and $\udiamcon$.
\end{proof}

\begin{prop}\label{prop: A_7, omega_3, k geq 2, nets}
Let $G = A_7$ and $\lambda = \omega_3$; then for $k \in [2, \frac{d}{2}]$ the quadruple $(G, \lambda, p, k)$ satisfies $\ssdiamcon$ and $\udiamcon$.
\end{prop}

\begin{proof}
The weight table is as follows.
$$
\begin{array}{|*4{>{\ss}c|}}
\hline
i & \mu & |W.\mu| & m_\mu \\
\hline
1 & \omega_3 & 56 & 1 \\
\hline
\end{array}
$$
We have $M = 56$, $M_3 = 42$ and $M_2 = 32$; we take $k_0 = 2$.

Take $\Psi = \langle \alpha_1 \rangle$ of type $A_1$. The $\Psi$-net table is as follows.
$$
\begin{array}{|*5{>{\ss}c|}}
\hline
\multicolumn{2}{|>{\ss}c|}{\Psi\mathrm{-nets}} & & \multicolumn{1}{|>{\ss}c|}{c(s)} & \multicolumn{1}{|>{\ss}c|}{c(u_\Psi)} \\
\cline{1-2} \cline{4-5}
    \bar\nu    & n_1 &  m & r \geq 2 & p \geq 2 \\
\hline
     \bom_1    &  2  & 15 &    15    &    15    \\
       0       &  1  & 26 &          &          \\
\hline
\multicolumn{3}{c|}{}     &    15    &    15    \\
\cline{4-5}
\end{array}
$$
Thus $c(\Psi)_{ss} = c(\Psi)_u = 15$, so we take $\d_0 = (41, 15)$; using Proposition~\ref{prop: B value when t = 2} we then have $B_{\d_0, 2} = 30 > 14 = \dim {u_\Psi}^G$. We need only consider semisimple classes $s^G$ with $|\Phi(s)| \leq M - 30 = 26 < 30 = m_{{A_1}^3}$, each of which has a subsystem of type ${A_1}^3$ disjoint from $\Phi(s)$, and unipotent classes of dimension at least $30$, each of which has the class ${A_1}^3$ in its closure by Lemma~\ref{lem: various classes in classical groups by dim}(ii).

Now take $\Psi = \langle \alpha_1, \alpha_3, \alpha_5 \rangle$ of type ${A_1}^3$. The $\Psi$-net table is as follows.
$$
\begin{array}{|*5{>{\ss}c|}}
\hline
\multicolumn{2}{|>{\ss}c|}{\Psi\mathrm{-nets}} & & \multicolumn{1}{|>{\ss}c|}{c(s)} & \multicolumn{1}{|>{\ss}c|}{c(u_\Psi)} \\
\cline{1-2} \cline{4-5}
          \bar\nu          & n_1 & m & r \geq 2 & p \geq 2 \\
\hline
  \bom_1 + \bom_3 + \bom_5 &  8  & 1 &     4    &     4    \\
      \bom_1 + \bom_3      &  4  & 2 &     4    &     4    \\
      \bom_1 + \bom_3      &  4  & 2 &     4    &     4    \\
      \bom_1 + \bom_3      &  4  & 2 &     4    &     4    \\
           \bom_1          &  2  & 3 &     3    &     3    \\
           \bom_1          &  2  & 3 &     3    &     3    \\
           \bom_1          &  2  & 3 &     3    &     3    \\
             0             &  1  & 6 &          &          \\
\hline
\multicolumn{3}{c|}{}                &    25    &    25    \\
\cline{4-5}
\end{array}
$$
Thus $c(\Psi)_{ss} = c(\Psi)_u = 25$, so we take $\d_0 = (31, 25)$; using Proposition~\ref{prop: B value when t = 2} we then have $B_{\d_0, 2} = 50 > M_3 > 30 = \dim {u_\Psi}^G$. We may therefore assume from now on that $r \geq 5$, and that $p \geq 5$ when we treat unipotent classes. We need only consider semisimple classes $s^G$ with $|\Phi(s)| \leq M - 50 = 6$, each of which has a subsystem of type $A_3$ disjoint from $\Phi(s)$, and unipotent classes of dimension at least $50$, each of which has the class $A_3$ in its closure by Lemma~\ref{lem: various classes in classical groups by dim}(vi).

Now take $\Psi = \langle \alpha_1, \alpha_2, \alpha_3 \rangle$ of type $A_3$. The $\Psi$-net table is as follows.
$$
\begin{array}{|*5{>{\ss}c|}}
\hline
\multicolumn{2}{|>{\ss}c|}{\Psi\mathrm{-nets}} & & \multicolumn{1}{|>{\ss}c|}{c(s)} & \multicolumn{1}{|>{\ss}c|}{c(u_\Psi)} \\
\cline{1-2} \cline{4-5}
     \bar\nu     & n_1 & m & r \geq 5 & p \geq 5 \\
\hline
      \bom_1     &  4  & 6 &    18    &    18    \\
      \bom_2     &  6  & 4 &    16    &    16    \\
      \bom_3     &  4  & 1 &     3    &     3    \\
        0        &  1  & 4 &          &          \\
\hline
\multicolumn{3}{c|}{}      &    37    &    37    \\
\cline{4-5}
\end{array}
$$
Thus $c(\Psi)_{ss} = c(\Psi)_u = 37$, so we take $\d_0 = (19, 19, 18)$; using Corollary~\ref{cor: B values for small k} we then have $B_{\d_0, 2} = 72 > M$. Therefore if $k \in [2, \frac{d}{2}]$ the quadruple $(G, \lambda, p, k)$ satisfies $\ssdiamcon$ and $\udiamcon$.
\end{proof}

\begin{prop}\label{prop: A_8, omega_3, k geq 2, nets}
Let $G = A_8$ and $\lambda = \omega_3$; then for $k \in [2, \frac{d}{2}]$ the quadruple $(G, \lambda, p, k)$ satisfies $\ssdiamcon$ and $\udiamcon$.
\end{prop}

\begin{proof}
The weight table is as follows.
$$
\begin{array}{|*4{>{\ss}c|}}
\hline
i & \mu & |W.\mu| & m_\mu \\
\hline
1 & \omega_3 & 84 & 1 \\
\hline
\end{array}
$$
We have $M = 72$ and $M_2 = 40$; we take $k_0 = 2$.

Take $\Psi = \langle \alpha_1 \rangle$ of type $A_1$. The $\Psi$-net table is as follows.
$$
\begin{array}{|*5{>{\ss}c|}}
\hline
\multicolumn{2}{|>{\ss}c|}{\Psi\mathrm{-nets}} & & \multicolumn{1}{|>{\ss}c|}{c(s)} & \multicolumn{1}{|>{\ss}c|}{c(u_\Psi)} \\
\cline{1-2} \cline{4-5}
    \bar\nu    & n_1 &  m & r \geq 2 & p \geq 2 \\
\hline
     \bom_1    &  2  & 21 &    21    &    21    \\
       0       &  1  & 42 &          &          \\
\hline
\multicolumn{3}{c|}{}     &    21    &    21    \\
\cline{4-5}
\end{array}
$$
Thus $c(\Psi)_{ss} = c(\Psi)_u = 21$, so we take $\d_0 = (63, 21)$; using Proposition~\ref{prop: B value when t = 2} we then have $B_{\d_0, 2} = 42 > M_2 > 16 = \dim {u_\Psi}^G$. We may therefore assume from now on that $r \geq 3$, and that $p \geq 3$ when we treat unipotent classes. We need only consider semisimple classes $s^G$ with $|\Phi(s)| \leq M - 42 = 30$, each of which has a subsystem of type $A_2A_1$ disjoint from $\Phi(s)$, and unipotent classes of dimension at least $42$, each of which has the class $A_2A_1$ in its closure by Lemma~\ref{lem: various classes in classical groups by dim}(iv).

Now take $\Psi = \langle \alpha_1, \alpha_2, \alpha_4 \rangle$ of type $A_2A_1$. The $\Psi$-net table is as follows.
$$
\begin{array}{|*5{>{\ss}c|}}
\hline
\multicolumn{2}{|>{\ss}c|}{\Psi\mathrm{-nets}} & & \multicolumn{1}{|>{\ss}c|}{c(s)} & \multicolumn{1}{|>{\ss}c|}{c(u_\Psi)} \\
\cline{1-2} \cline{4-5}
     \bar\nu     & n_1 & m & r \geq 3 & p \geq 3 \\
\hline
 \bom_1 + \bom_4 &  6  & 4 &    16    &    16    \\
 \bom_2 + \bom_4 &  6  & 1 &     4    &     4    \\
      \bom_1     &  3  & 7 &    14    &    14    \\
      \bom_2     &  3  & 4 &     8    &     8    \\
      \bom_4     &  2  & 6 &     6    &     6    \\
        0        &  1  & 9 &          &          \\
\hline
\multicolumn{3}{c|}{}      &    48    &    48    \\
\cline{4-5}
\end{array}
$$
Thus $c(\Psi)_{ss} = c(\Psi)_u = 48$, so we take $\d_0 = (36, 36, 12)$; using Corollary~\ref{cor: B values for small k} we then have $B_{\d_0, 2} = 94 > M$. Therefore if $k \in [2, \frac{d}{2}]$ the quadruple $(G, \lambda, p, k)$ satisfies $\ssdiamcon$ and $\udiamcon$.
\end{proof}

\begin{prop}\label{prop: A_2, 3omega_1, k geq 2, nets}
Let $G = A_2$ and $\lambda = 3\omega_1$ with $p \geq 5$; then for $k \in [2, \frac{d}{2}]$ the quadruple $(G, \lambda, p, k)$ satisfies $\ssdiamcon$ and $\udiamcon$.
\end{prop}

\begin{proof}
The weight table is as follows.
$$
\begin{array}{|*4{>{\ss}c|}}
\hline
i & \mu & |W.\mu| & m_\mu \\
\hline
2 &      3\omega_1      & 3 & 1 \\
1 & \omega_1 + \omega_2 & 6 & 1 \\
0 &          0          & 1 & 1 \\
\hline
\end{array}
$$
We have $M = 6$; we take $k_0 = 2$.

Take $\Psi = \langle \alpha_1 \rangle$ of type $A_1$. The $\Psi$-net table is as follows.
$$
\begin{array}{|*9{>{\ss}c|}}
\hline
\multicolumn{4}{|>{\ss}c|}{\Psi\mathrm{-nets}} & & \multicolumn{3}{|>{\ss}c|}{c(s)} & \multicolumn{1}{|>{\ss}c|}{c(u_\Psi)} \\
\cline{1-4} \cline{6-9}
    \bar\nu    & n_0 & n_1 & n_2 & m & r = 2 & r = 3 & r \geq 5 & p \geq 5 \\
\hline
    3\bom_1    &  0  &  2  &  2  & 1 &   2   &   2   &     3    &     3    \\
    2\bom_1    &  1  &  2  &  0  & 1 &   1   &   2   &     2    &     2    \\
     \bom_1    &  0  &  2  &  0  & 1 &   1   &   1   &     1    &     1    \\
       0       &  0  &  0  &  1  & 1 &       &       &          &          \\
\hline
\multicolumn{5}{c|}{}                &   4   &   5   &     6    &     6    \\
\cline{6-9}
\end{array}
$$
Thus $c(\Psi)_{ss}, c(\Psi)_u \geq 4$, so we may take $\d_0 = (6, 4)$; using Proposition~\ref{prop: B value when t = 2} we then have $B_{\d_0, 2} = 8 > M$. Therefore if $k \in [2, \frac{d}{2}]$ the quadruple $(G, \lambda, p, k)$ satisfies $\ssdiamcon$ and $\udiamcon$.
\end{proof}

\begin{prop}\label{prop: A_1, 4omega_1, k geq 2, nets}
Let $G = A_1$ and $\lambda = 4\omega_1$ with $p \geq 5$; then for $k \in [2, \frac{d}{2}]$ the quadruple $(G, \lambda, p, k)$ satisfies $\ssdiamcon$ and $\udiamcon$.
\end{prop}

\begin{proof}
The weight table is as follows.
$$
\begin{array}{|*4{>{\ss}c|}}
\hline
i & \mu & |W.\mu| & m_\mu \\
\hline
2 & 4\omega_1 & 2 & 1 \\
1 & 2\omega_1 & 2 & 1 \\
0 &     0     & 1 & 1 \\
\hline
\end{array}
$$
We have $M = 2$; we take $k_0 = 2$.

Take $\Psi = \langle \alpha_1 \rangle$ of type $A_1$. The $\Psi$-net table is as follows.
$$
\begin{array}{|*9{>{\ss}c|}}
\hline
\multicolumn{4}{|>{\ss}c|}{\Psi\mathrm{-nets}} & & \multicolumn{3}{|>{\ss}c|}{c(s)} & \multicolumn{1}{|>{\ss}c|}{c(u_\Psi)} \\
\cline{1-4} \cline{6-9}
    \bar\nu    & n_0 & n_1 & n_2 & m & r = 2 & r = 3 & r \geq 5 & p \geq 5 \\
\hline
    4\bom_1    &  1  &  2  &  2  & 1 &   2   &   3   &     4    &     4    \\
\hline
\multicolumn{5}{c|}{}                &   2   &   3   &     4    &     4    \\
\cline{6-9}
\end{array}
$$
Thus $c(\Psi)_{ss}, c(\Psi)_u \geq 2$, so we may take $\d_0 = (3, 2)$; using Proposition~\ref{prop: B value when t = 2} we then have $B_{\d_0, 2} = 3 > M$. Therefore if $k \in [2, \frac{d}{2}]$ the quadruple $(G, \lambda, p, k)$ satisfies $\ssdiamcon$ and $\udiamcon$.
\end{proof}

\begin{prop}\label{prop: A_3, 2omega_2, k geq 2, nets}
Let $G = A_3$ and $\lambda = 2\omega_2$ with $p \geq 3$; then for $k \in [2, \frac{d}{2}]$ the quadruple $(G, \lambda, p, k)$ satisfies $\ssdiamcon$ and $\udiamcon$.
\end{prop}

\begin{proof}
Write $\z = \z_{p, 3}$. The weight table is as follows.
$$
\begin{array}{|*4{>{\ss}c|}}
\hline
i & \mu & |W.\mu| & m_\mu \\
\hline
2 &      2\omega_2      &  6 &    1   \\
1 & \omega_1 + \omega_3 & 12 &    1   \\
0 &          0          &  1 & 2 - \z \\
\hline
\end{array}
$$
We have $M = 12$; we take $k_0 = 2$.

Take $\Psi = \langle \alpha_1 \rangle$ of type $A_1$. The $\Psi$-net table is as follows.
$$
\begin{array}{|*8{>{\ss}c|}}
\hline
\multicolumn{4}{|>{\ss}c|}{\Psi\mathrm{-nets}} & & \multicolumn{2}{|>{\ss}c|}{c(s)} & \multicolumn{1}{|>{\ss}c|}{c(u_\Psi)} \\
\cline{1-4} \cline{6-8}
    \bar\nu    & n_0 & n_1 & n_2 & m &  r = 2 & r \geq 3 & p \geq 3 \\
\hline
    2\bom_1    &  0  &  1  &  2  & 2 &    2   &     4    &     4    \\
    2\bom_1    &  1  &  2  &  0  & 1 & 2 - \z &     2    &     2    \\
     \bom_1    &  0  &  2  &  0  & 4 &    4   &     4    &     4    \\
       0       &  0  &  0  &  1  & 2 &        &          &          \\
\hline
\multicolumn{5}{c|}{}                & 8 - \z &    10    &    10    \\
\cline{6-8}
\end{array}
$$
Thus $c(\Psi)_{ss}, c(\Psi)_u \geq 8 - \z$, so we may take $\d_0 = (12, 8 - \z)$; using Proposition~\ref{prop: B value when t = 2} we then have $B_{\d_0, 2} = 16 - 2\z > M$. Therefore if $k \in [2, \frac{d}{2}]$ the quadruple $(G, \lambda, p, k)$ satisfies $\ssdiamcon$ and $\udiamcon$.
\end{proof}

\begin{prop}\label{prop: A_7, omega_4, k geq 2, nets}
Let $G = A_7$ and $\lambda = \omega_4$; then for $k \in [2, \frac{d}{2}]$ the quadruple $(G, \lambda, p, k)$ satisfies $\ssdiamcon$ and $\udiamcon$.
\end{prop}

\begin{proof}
The weight table is as follows.
$$
\begin{array}{|*4{>{\ss}c|}}
\hline
i & \mu & |W.\mu| & m_\mu \\
\hline
1 & \omega_4 & 70 & 1 \\
\hline
\end{array}
$$
We have $M = 56$ and $M_2 = 32$; we take $k_0 = 2$.

Take $\Psi = \langle \alpha_1 \rangle$ of type $A_1$. The $\Psi$-net table is as follows.
$$
\begin{array}{|*5{>{\ss}c|}}
\hline
\multicolumn{2}{|>{\ss}c|}{\Psi\mathrm{-nets}} & & \multicolumn{1}{|>{\ss}c|}{c(s)} & \multicolumn{1}{|>{\ss}c|}{c(u_\Psi)} \\
\cline{1-2} \cline{4-5}
    \bar\nu    & n_1 &  m & r \geq 2 & p \geq 2 \\
\hline
     \bom_1    &  2  & 20 &    20    &    20    \\
       0       &  1  & 30 &          &          \\
\hline
\multicolumn{3}{c|}{}     &    20    &    20    \\
\cline{4-5}
\end{array}
$$
Thus $c(\Psi)_{ss} = c(\Psi)_u = 20$, so we take $\d_0 = (50, 20)$; using Proposition~\ref{prop: B value when t = 2} we then have $B_{\d_0, 2} = 40 > M_2 > 14 = \dim {u_\Psi}^G$. We may therefore assume from now on that $r \geq 3$, and that $p \geq 3$ when we treat unipotent classes. We need only consider semisimple classes $s^G$ with $|\Phi(s)| \leq M - 40 = 16 < 24 = m_{A_2}$, each of which has a subsystem of type $A_2$ disjoint from $\Phi(s)$, and unipotent classes of dimension at least $40$, each of which has the class $A_2$ in its closure by Lemma~\ref{lem: various classes in classical groups by dim}(iii).

Now take $\Psi = \langle \alpha_1, \alpha_2 \rangle$ of type $A_2$. The $\Psi$-net table is as follows.
$$
\begin{array}{|*5{>{\ss}c|}}
\hline
\multicolumn{2}{|>{\ss}c|}{\Psi\mathrm{-nets}} & & \multicolumn{1}{|>{\ss}c|}{c(s)} & \multicolumn{1}{|>{\ss}c|}{c(u_\Psi)} \\
\cline{1-2} \cline{4-5}
     \bar\nu     & n_1 &  m & r \geq 3 & p \geq 3 \\
\hline
      \bom_1     &  3  & 10 &    20    &    20    \\
      \bom_2     &  3  & 10 &    20    &    20    \\
        0        &  1  & 10 &          &          \\
\hline
\multicolumn{3}{c|}{}       &    40    &    40    \\
\cline{4-5}
\end{array}
$$
Thus $c(\Psi)_{ss} = c(\Psi)_u = 40$, so we take $\d_0 = (30, 30, 10)$; using Corollary~\ref{cor: B values for small k} we then have $B_{\d_0, 2} = 78 > M$. Therefore if $k \in [2, \frac{d}{2}]$ the quadruple $(G, \lambda, p, k)$ satisfies $\ssdiamcon$ and $\udiamcon$.
\end{proof}

\begin{prop}\label{prop: A_3, omega_1 + omega_2, k geq 2, nets}
Let $G = A_3$ and $\lambda = \omega_1 + \omega_2$ with $p = 3$; then for $k \in [2, \frac{d}{2}]$ the quadruple $(G, \lambda, p, k)$ satisfies $\ssdiamcon$ and $\udiamcon$.
\end{prop}

\begin{proof}
The weight table is as follows.
$$
\begin{array}{|*4{>{\ss}c|}}
\hline
i & \mu & |W.\mu| & m_\mu \\
\hline
1 & \omega_1 + \omega_2 & 12 &    1   \\
2 &       \omega_3      &  4 &    1   \\
\hline
\end{array}
$$
We have $M = 12$ and $M_2 = 8$; we take $k_0 = 2$.

Take $\Psi = \langle \alpha_1 \rangle$ of type $A_1$. The $\Psi$-net table is as follows.
$$
\begin{array}{|*7{>{\ss}c|}}
\hline
\multicolumn{3}{|>{\ss}c|}{\Psi\mathrm{-nets}} & & \multicolumn{2}{|>{\ss}c|}{c(s)} & \multicolumn{1}{|>{\ss}c|}{c(u_\Psi)} \\
\cline{1-3} \cline{5-7}
    \bar\nu    & n_1 & n_2 & m & r = 2 & r \geq 5 & p = 3 \\
\hline
    2\bom_1    &  1  &  2  & 2 &   2   &     4    &   4   \\
     \bom_1    &  0  &  2  & 3 &   3   &     3    &   3   \\
     \bom_1    &  2  &  0  & 1 &   1   &     1    &   1   \\
       0       &  0  &  1  & 2 &       &          &       \\
\hline
\multicolumn{4}{c|}{}          &   6   &     8    &   8   \\
\cline{5-7}
\end{array}
$$
Thus if $r = 2$ then $c(\Psi)_{ss} = 6$, so we take $\d_0 = (10, 6)$; using Proposition~\ref{prop: B value when t = 2} we then have $B_{\d_0, 2} = 12 > M_2$. If instead $r \geq 5$ then $c(\Psi)_{ss} = c(\Psi)_u = 8$, so we take $\d_0 = (8, 8)$; using Proposition~\ref{prop: B value when t = 2} we then have $B_{\d_0, 2} = 14 > M$. Therefore if $k \in [2, \frac{d}{2}]$ the quadruple $(G, \lambda, p, k)$ satisfies $\ssdiamcon$ and $\udiamcon$.
\end{proof}

\begin{prop}\label{prop: D_5, omega_5, k geq 5, nets}
Let $G = D_5$ and $\lambda = \omega_5$; then for $k \in [5, \frac{d}{2}]$ the quadruple $(G, \lambda, p, k)$ satisfies $\ssdiamcon$ and $\udiamcon$.
\end{prop}

\begin{proof}
The weight table is as follows.
$$
\begin{array}{|*4{>{\ss}c|}}
\hline
i & \mu & |W.\mu| & m_\mu \\
\hline
1 & \omega_5 & 16 & 1 \\
\hline
\end{array}
$$
We have $M = 40$, $M_3 = 30$ and $M_2 = 24$; we take $k_0 = 5$.

Take $\Psi = \langle \alpha_1 \rangle$ of type $A_1$. The $\Psi$-net table is as follows.
$$
\begin{array}{|*5{>{\ss}c|}}
\hline
\multicolumn{2}{|>{\ss}c|}{\Psi\mathrm{-nets}} & & \multicolumn{1}{|>{\ss}c|}{c(s)} & \multicolumn{1}{|>{\ss}c|}{c(u_\Psi)} \\
\cline{1-2} \cline{4-5}
    \bar\nu    & n_1 & m & r \geq 2 & p \geq 2 \\
\hline
     \bom_1    &  2  & 4 &     4    &     4    \\
       0       &  1  & 8 &          &          \\
\hline
\multicolumn{3}{c|}{}    &     4    &     4    \\
\cline{4-5}
\end{array}
$$
Thus $c(\Psi)_{ss} = c(\Psi)_u = 4$, so we take $\d_0 = (12, 4)$; using Proposition~\ref{prop: B value when t = 2} we then have $B_{\d_0, 5} = 20 > 14 = \dim {u_\Psi}^G$. We need only consider semisimple classes $s^G$ with $|\Phi(s)| \leq M - 20 = 20$, each of which has a subsystem of type ${A_1}^2$ or a subsystem of type $D_2$ disjoint from $\Phi(s)$, and unipotent classes of dimension at least $20$, each of which has the class ${A_1}^2$ or the class $D_2$ in its closure by Lemma~\ref{lem: A_1^2 or D_2 in D_ell}.

Now take $\Psi = \langle \alpha_1, \alpha_3 \rangle$ of type ${A_1}^2$, and $\Psi = \langle \alpha_4, \alpha_5 \rangle$ of type $D_2$. The $\Psi$-net tables are as follows.
$$
\begin{array}{|*5{>{\ss}c|}}
\hline
\multicolumn{2}{|>{\ss}c|}{\Psi\mathrm{-nets}} & & \multicolumn{1}{|>{\ss}c|}{c(s)} & \multicolumn{1}{|>{\ss}c|}{c(u_\Psi)} \\
\cline{1-2} \cline{4-5}
      \bar\nu      & n_1 & m & r \geq 2 & p \geq 2 \\
\hline
  \bom_1 + \bom_3  &  4  & 1 &     2    &     2    \\
       \bom_1      &  2  & 2 &     2    &     2    \\
       \bom_3      &  2  & 2 &     2    &     2    \\
         0         &  1  & 4 &          &          \\
\hline
\multicolumn{3}{c|}{}        &     6    &     6    \\
\cline{4-5}
\end{array}
\qquad
\begin{array}{|*5{>{\ss}c|}}
\hline
\multicolumn{2}{|>{\ss}c|}{\Psi\mathrm{-nets}} & & \multicolumn{1}{|>{\ss}c|}{c(s)} & \multicolumn{1}{|>{\ss}c|}{c(u_\Psi)} \\
\cline{1-2} \cline{4-5}
     \bar\nu     & n_1 & m & r \geq 2 & p \geq 2 \\
\hline
       \bom_4    &  2  & 4 &     4    &     4    \\
       \bom_5    &  2  & 4 &     4    &     4    \\
\hline
\multicolumn{3}{c|}{}      &     8    &     8    \\
\cline{4-5}
\end{array}
$$
Thus according as $\Psi = {A_1}^2$ or $D_2$ we have $c(\Psi)_{ss} = c(\Psi)_u = 6$ or $8$, so we take $\d_0 = (10, 6)$ or $(8, 8)$; using Proposition~\ref{prop: B value when t = 2} we then have $B_{\d_0, 5} = 26 > M_2 > 20 = \dim {u_\Psi}^G$ or $B_{\d_0, 5} = 28 > M_2 > 16 = \dim {u_\Psi}^G$. Taking the smaller of the two lower bounds, we see that we may therefore assume from now on that $r \geq 3$, and that $p \geq 3$ when we treat unipotent classes. Moreover we need only consider semisimple classes $s^G$ with $|\Phi(s)| \leq M - 26 = 14$, and unipotent classes of dimension at least $26$; since each of the former has a subsystem of type $D_2$ disjoint from $\Phi(s)$, and each of the latter has the class $D_2$ in its closure by Lemma~\ref{lem: various classes in classical groups by dim}(ix), we may actually take the larger of the two lower bounds. We need therefore only consider semisimple classes $s^G$ with $|\Phi(s)| \leq M - 28 = 12$, each of which has a subsystem of type $A_2A_1$ or a subsystem of type $D_3$ disjoint from $\Phi(s)$, and unipotent classes of dimension at least $28$, each of which has the class $A_2A_1$ or the class $D_3$ in its closure by Lemma~\ref{lem: various classes in D_ell for fixed ell}(vi) (and for the unipotent class $D_3$ to lie in $G_{(p)}$ we need $p \geq 5$).

Now take $\Psi = \langle \alpha_1, \alpha_2, \alpha_4 \rangle$ of type $A_2A_1$, and $\Psi = \langle \alpha_3, \alpha_4, \alpha_5 \rangle$ of type $D_3$. The $\Psi$-net tables are as follows.
$$
\begin{array}{|*5{>{\ss}c|}}
\hline
\multicolumn{2}{|>{\ss}c|}{\Psi\mathrm{-nets}} & & \multicolumn{1}{|>{\ss}c|}{c(s)} & \multicolumn{1}{|>{\ss}c|}{c(u_\Psi)} \\
\cline{1-2} \cline{4-5}
      \bar\nu      & n_1 & m & r \geq 3 & p \geq 3 \\
\hline
  \bom_2 + \bom_4  &  6  & 1 &     4    &     4    \\
       \bom_1      &  3  & 2 &     4    &     4    \\
       \bom_4      &  2  & 1 &     1    &     1    \\
         0         &  1  & 2 &          &          \\
\hline
\multicolumn{3}{c|}{}        &     9    &     9    \\
\cline{4-5}
\end{array}
\qquad
\begin{array}{|*5{>{\ss}c|}}
\hline
\multicolumn{2}{|>{\ss}c|}{\Psi\mathrm{-nets}} & & \multicolumn{1}{|>{\ss}c|}{c(s)} & \multicolumn{1}{|>{\ss}c|}{c(u_\Psi)} \\
\cline{1-2} \cline{4-5}
     \bar\nu     & n_1 & m & r \geq 3 & p \geq 5 \\
\hline
       \bom_4    &  4  & 2 &     6    &     6    \\
       \bom_5    &  4  & 2 &     6    &     6    \\
\hline
\multicolumn{3}{c|}{}      &    12    &    12    \\
\cline{4-5}
\end{array}
$$
Thus according as $\Psi = A_2A_1$ or $D_3$ we have $c(\Psi)_{ss} = c(\Psi)_u = 9$ or $12$, so we take $\d_0 = (7, 7, 2)$ or $(4, 4, 4, 4)$. In the former case, according as $\k_0 = (5, 0, 0)$, $(4, 1, 0)$, $(3, 2, 0)$, $(3, 1, 1)$ or $(2, 2, 1)$ we have $B_{\d_0, \k_0} = 45$, $37$, $33$, $36$ or $34$, so $B_{\d_0, 5} = 33 > M_3 > 28 = \dim {u_\Psi}^G$; in the latter case, according as $\k_0 = (4, 1, 0, 0)$, $(3, 2, 0, 0)$, $(3, 1, 1, 0)$, $(2, 2, 1, 0)$ or $(2, 1, 1, 1)$ we have $B_{\d_0, \k_0} = 52$, $48$, $46$, $44$ or $42$, so $B_{\d_0, 5} = 42 > M$. Taking the smaller of the two lower bounds, we see that we need only consider semisimple classes $s^G$ with $|\Phi(s)| \leq M - 33 = 7$, and unipotent classes of dimension at least $33$; since each of the former has a subsystem of type $D_3$ disjoint from $\Phi(s)$, and each of the latter has the class $D_3$ in its closure by Lemma~\ref{lem: various classes in D_ell for fixed ell}(vii), we may actually take the larger of the two lower bounds. Therefore if $k \in [5, \frac{d}{2}]$ the quadruple $(G, \lambda, p, k)$ satisfies $\ssdiamcon$ and $\udiamcon$.
\end{proof}

\begin{prop}\label{prop: D_6, omega_6, k geq 3, nets}
Let $G = D_6$ and $\lambda = \omega_6$; then for $k \in [3, \frac{d}{2}]$ the quadruple $(G, \lambda, p, k)$ satisfies $\ssdiamcon$ and $\udiamcon$.
\end{prop}

\begin{proof}
The weight table is as follows.
$$
\begin{array}{|*4{>{\ss}c|}}
\hline
i & \mu & |W.\mu| & m_\mu \\
\hline
1 & \omega_6 & 32 & 1 \\
\hline
\end{array}
$$
We have $M = 60$, $M_3 = 44$ and $M_2 = 36$; we take $k_0 = 3$.

Take $\Psi = \langle \alpha_1 \rangle$ of type $A_1$. The $\Psi$-net table is as follows.
$$
\begin{array}{|*5{>{\ss}c|}}
\hline
\multicolumn{2}{|>{\ss}c|}{\Psi\mathrm{-nets}} & & \multicolumn{1}{|>{\ss}c|}{c(s)} & \multicolumn{1}{|>{\ss}c|}{c(u_\Psi)} \\
\cline{1-2} \cline{4-5}
    \bar\nu    & n_1 &  m & r \geq 2 & p \geq 2 \\
\hline
     \bom_1    &  2  &  8 &     8    &     8    \\
       0       &  1  & 16 &          &          \\
\hline
\multicolumn{3}{c|}{}     &     8    &     8    \\
\cline{4-5}
\end{array}
$$
Thus $c(\Psi)_{ss} = c(\Psi)_u = 8$, so we take $\d_0 = (24, 8)$; using Proposition~\ref{prop: B value when t = 2} we then have $B_{\d_0, 3} = 24 > 18 = \dim {u_\Psi}^G$. We need only consider semisimple classes $s^G$ with $|\Phi(s)| \leq M - 24 = 36$, each of which has a subsystem of type ${A_1}^2$ or a subsystem of type $D_2$ disjoint from $\Phi(s)$, and unipotent classes of dimension at least $24$, each of which has the class ${A_1}^2$ or the class $D_2$ in its closure by Lemma~\ref{lem: A_1^2 or D_2 in D_ell}.

Now take $\Psi = \langle \alpha_1, \alpha_3 \rangle$ of type ${A_1}^2$, and $\Psi = \langle \alpha_5, \alpha_6 \rangle$ of type $D_2$. The $\Psi$-net tables are as follows.
$$
\begin{array}{|*5{>{\ss}c|}}
\hline
\multicolumn{2}{|>{\ss}c|}{\Psi\mathrm{-nets}} & & \multicolumn{1}{|>{\ss}c|}{c(s)} & \multicolumn{1}{|>{\ss}c|}{c(u_\Psi)} \\
\cline{1-2} \cline{4-5}
      \bar\nu      & n_1 & m & r \geq 2 & p \geq 2 \\
\hline
  \bom_1 + \bom_3  &  4  & 2 &     4    &     4    \\
       \bom_1      &  2  & 4 &     4    &     4    \\
       \bom_3      &  2  & 4 &     4    &     4    \\
         0         &  1  & 8 &          &          \\
\hline
\multicolumn{3}{c|}{}        &    12    &    12    \\
\cline{4-5}
\end{array}
\qquad
\begin{array}{|*5{>{\ss}c|}}
\hline
\multicolumn{2}{|>{\ss}c|}{\Psi\mathrm{-nets}} & & \multicolumn{1}{|>{\ss}c|}{c(s)} & \multicolumn{1}{|>{\ss}c|}{c(u_\Psi)} \\
\cline{1-2} \cline{4-5}
     \bar\nu     & n_1 & m & r \geq 2 & p \geq 2 \\
\hline
       \bom_5    &  2  & 8 &     8    &     8    \\
       \bom_6    &  2  & 8 &     8    &     8    \\
\hline
\multicolumn{3}{c|}{}      &    16    &    16    \\
\cline{4-5}
\end{array}
$$
Thus according as $\Psi = {A_1}^2$ or $D_2$ we have $c(\Psi)_{ss} = c(\Psi)_u = 12$ or $16$, so we take $\d_0 = (20, 12)$ or $(16, 16)$; using Proposition~\ref{prop: B value when t = 2} we then have $B_{\d_0, 3} = 36 > 28 = \dim {u_\Psi}^G$ or $B_{\d_0, 3} = 44 > M_2 > 20 = \dim {u_\Psi}^G$. Taking the smaller of the two lower bounds, we need only consider semisimple classes $s^G$ with $|\Phi(s)| \leq M - 36 = 24$, and unipotent classes of dimension at least $36$; since each of the former has a subsystem of type $D_2$ disjoint from $\Phi(s)$, and each of the latter has the class $D_2$ in its closure by Lemma~\ref{lem: various classes in classical groups by dim}(ix), we may actually take the larger of the two lower bounds. We may therefore assume from now on that $r \geq 3$, and that $p \geq 3$ when we treat unipotent classes; and we need only consider semisimple classes $s^G$ with $|\Phi(s)| \leq M - 44 = 16$, each of which has a subsystem of type $A_2A_1$ disjoint from $\Phi(s)$, and unipotent classes of dimension at least $44$, each of which has the class $A_2A_1$ in its closure by Lemma~\ref{lem: various classes in D_ell for fixed ell}(iv).

Now take $\Psi = \langle \alpha_1, \alpha_2, \alpha_4 \rangle$ of type $A_2A_1$. The $\Psi$-net table is as follows.
$$
\begin{array}{|*5{>{\ss}c|}}
\hline
\multicolumn{2}{|>{\ss}c|}{\Psi\mathrm{-nets}} & & \multicolumn{1}{|>{\ss}c|}{c(s)} & \multicolumn{1}{|>{\ss}c|}{c(u_\Psi)} \\
\cline{1-2} \cline{4-5}
      \bar\nu      & n_1 & m & r \geq 2 & p \geq 2 \\
\hline
  \bom_1 + \bom_4  &  6  & 1 &     4    &     4    \\
  \bom_2 + \bom_4  &  6  & 1 &     4    &     4    \\
       \bom_1      &  3  & 2 &     4    &     4    \\
       \bom_2      &  3  & 2 &     4    &     4    \\
       \bom_4      &  2  & 2 &     2    &     2    \\
         0         &  1  & 4 &          &          \\
\hline
\multicolumn{3}{c|}{}        &    18    &    18    \\
\cline{4-5}
\end{array}
$$
Thus $c(\Psi)_{ss} = c(\Psi)_u = 18$, so we take $\d_0 = (14, 14, 4)$; using Corollary~\ref{cor: B values for small k} we then have $B_{\d_0, 3} = 50 > M_3 > 40 = \dim {u_\Psi}^G$. We need only consider semisimple classes $s^G$ with $|\Phi(s)| \leq M - 50 = 10$, each of which has a subsystem of type $D_3$ disjoint from $\Phi(s)$, and unipotent classes of dimension at least $50$, each of which has the class $D_3$ in its closure by Lemma~\ref{lem: various classes in D_ell for fixed ell}(v).

Now take $\Psi = \langle \alpha_3, \alpha_4, \alpha_5 \rangle$ of type $D_3$. The $\Psi$-net table is as follows.
$$
\begin{array}{|*5{>{\ss}c|}}
\hline
\multicolumn{2}{|>{\ss}c|}{\Psi\mathrm{-nets}} & & \multicolumn{1}{|>{\ss}c|}{c(s)} & \multicolumn{1}{|>{\ss}c|}{c(u_\Psi)} \\
\cline{1-2} \cline{4-5}
     \bar\nu     & n_1 & m & r \geq 2 & p \geq 2 \\
\hline
       \bom_5    &  4  & 4 &    12    &    12    \\
       \bom_6    &  4  & 4 &    12    &    12    \\
\hline
\multicolumn{3}{c|}{}      &    24    &    24    \\
\cline{4-5}
\end{array}
$$
Thus $c(\Psi)_{ss} = c(\Psi)_u = 24$, so we take $\d_0 = (8, 8, 8, 8)$; using Corollary~\ref{cor: B values for small k} we then have $B_{\d_0, 3} = 66 > M$. Therefore if $k \in [3, \frac{d}{2}]$ the quadruple $(G, \lambda, p, k)$ satisfies $\ssdiamcon$ and $\udiamcon$.
\end{proof}

\begin{prop}\label{prop: D_7, omega_7, k geq 2, nets}
Let $G = D_7$ and $\lambda = \omega_7$; then for $k \in [2, \frac{d}{2}]$ the quadruple $(G, \lambda, p, k)$ satisfies $\ssdiamcon$ and $\udiamcon$.
\end{prop}

\begin{proof}
The weight table is as follows.
$$
\begin{array}{|*4{>{\ss}c|}}
\hline
i & \mu & |W.\mu| & m_\mu \\
\hline
1 & \omega_7 & 64 & 1 \\
\hline
\end{array}
$$
We have $M = 84$, $M_3 = 60$ and $M_2 = 48$; we take $k_0 = 2$.

Take $\Psi = \langle \alpha_1 \rangle$ of type $A_1$. The $\Psi$-net table is as follows.
$$
\begin{array}{|*5{>{\ss}c|}}
\hline
\multicolumn{2}{|>{\ss}c|}{\Psi\mathrm{-nets}} & & \multicolumn{1}{|>{\ss}c|}{c(s)} & \multicolumn{1}{|>{\ss}c|}{c(u_\Psi)} \\
\cline{1-2} \cline{4-5}
    \bar\nu    & n_1 &  m & r \geq 2 & p \geq 2 \\
\hline
     \bom_1    &  2  & 16 &    16    &    16    \\
       0       &  1  & 32 &          &          \\
\hline
\multicolumn{3}{c|}{}     &    16    &    16    \\
\cline{4-5}
\end{array}
$$
Thus $c(\Psi)_{ss} = c(\Psi)_u = 16$, so we take $\d_0 = (48, 16)$; using Proposition~\ref{prop: B value when t = 2} we then have $B_{\d_0, 2} = 32 > 22 = \dim {u_\Psi}^G$. We need only consider semisimple classes $s^G$ with $|\Phi(s)| \leq M - 32 = 52$, each of which has a subsystem of type ${A_1}^2$ or a subsystem of type $D_2$ disjoint from $\Phi(s)$, and unipotent classes of dimension at least $32$, each of which has the class ${A_1}^2$ or the class $D_2$ in its closure by Lemma~\ref{lem: A_1^2 or D_2 in D_ell}.

Now take $\Psi = \langle \alpha_1, \alpha_3 \rangle$ of type ${A_1}^2$, and $\Psi = \langle \alpha_6, \alpha_7 \rangle$ of type $D_2$. The $\Psi$-net tables are as follows.
$$
\begin{array}{|*5{>{\ss}c|}}
\hline
\multicolumn{2}{|>{\ss}c|}{\Psi\mathrm{-nets}} & & \multicolumn{1}{|>{\ss}c|}{c(s)} & \multicolumn{1}{|>{\ss}c|}{c(u_\Psi)} \\
\cline{1-2} \cline{4-5}
      \bar\nu      & n_1 &  m & r \geq 2 & p \geq 2 \\
\hline
  \bom_1 + \bom_3  &  4  &  4 &     8    &     8    \\
       \bom_1      &  2  &  8 &     8    &     8    \\
       \bom_3      &  2  &  8 &     8    &     8    \\
         0         &  1  & 16 &          &          \\
\hline
\multicolumn{3}{c|}{}         &    24    &    24    \\
\cline{4-5}
\end{array}
\qquad
\begin{array}{|*5{>{\ss}c|}}
\hline
\multicolumn{2}{|>{\ss}c|}{\Psi\mathrm{-nets}} & & \multicolumn{1}{|>{\ss}c|}{c(s)} & \multicolumn{1}{|>{\ss}c|}{c(u_\Psi)} \\
\cline{1-2} \cline{4-5}
     \bar\nu     & n_1 &  m & r \geq 2 & p \geq 2 \\
\hline
       \bom_6    &  2  & 16 &    16    &    16    \\
       \bom_7    &  2  & 16 &    16    &    16    \\
\hline
\multicolumn{3}{c|}{}       &    32    &    32    \\
\cline{4-5}
\end{array}
$$
Thus according as $\Psi = {A_1}^2$ or $D_2$ we have $c(\Psi)_{ss} = c(\Psi)_u = 24$ or $32$, so we take $\d_0 = (40, 24)$ or $(32, 32)$; using Proposition~\ref{prop: B value when t = 2} we then have $B_{\d_0, 2} = 48 > 36 = \dim {u_\Psi}^G$ or $B_{\d_0, 2} = 62 > M_3 > 24 = \dim {u_\Psi}^G$. Taking the smaller of the two lower bounds, we see that we need only consider semisimple classes $s^G$ with $|\Phi(s)| \leq M - 48 = 36$, and unipotent classes of dimension at least $48$; since each of the former has a subsystem of type $D_2$ disjoint from $\Phi(s)$, and each of the latter has the class $D_2$ in its closure by Lemma~\ref{lem: various classes in classical groups by dim}(ix), we may actually take the larger of the two lower bounds. We may therefore assume from now on that $r \geq 5$, and that $p \geq 5$ when we treat unipotent classes. Moreover we need therefore only consider semisimple classes $s^G$ with $|\Phi(s)| \leq M - 62 = 22$, each of which has a subsystem of type $A_3$ or a subsystem of type $D_3$ disjoint from $\Phi(s)$, and unipotent classes of dimension at least $62$, each of which has the class $A_3$ or the class $D_3$ in its closure by Lemma~\ref{lem: various classes in D_ell for fixed ell}(ii).

Now take $\Psi = \langle \alpha_1, \alpha_2, \alpha_3 \rangle$ of type $A_3$, and $\Psi = \langle \alpha_3, \alpha_4, \alpha_5 \rangle$ of type $D_3$. The $\Psi$-net tables are as follows.
$$
\begin{array}{|*5{>{\ss}c|}}
\hline
\multicolumn{2}{|>{\ss}c|}{\Psi\mathrm{-nets}} & & \multicolumn{1}{|>{\ss}c|}{c(s)} & \multicolumn{1}{|>{\ss}c|}{c(u_\Psi)} \\
\cline{1-2} \cline{4-5}
      \bar\nu      & n_1 & m & r \geq 5 & p \geq 5 \\
\hline
       \bom_1      &  4  & 4 &    12    &    12    \\
       \bom_2      &  6  & 4 &    16    &    16    \\
       \bom_3      &  4  & 4 &    12    &    12    \\
         0         &  1  & 8 &          &          \\
\hline
\multicolumn{3}{c|}{}        &    40    &    40    \\
\cline{4-5}
\end{array}
\qquad
\begin{array}{|*5{>{\ss}c|}}
\hline
\multicolumn{2}{|>{\ss}c|}{\Psi\mathrm{-nets}} & & \multicolumn{1}{|>{\ss}c|}{c(s)} & \multicolumn{1}{|>{\ss}c|}{c(u_\Psi)} \\
\cline{1-2} \cline{4-5}
     \bar\nu     & n_1 & m & r \geq 5 & p \geq 5 \\
\hline
       \bom_6    &  4  & 8 &    24    &    24    \\
       \bom_7    &  4  & 8 &    24    &    24    \\
\hline
\multicolumn{3}{c|}{}      &    48    &    48    \\
\cline{4-5}
\end{array}
$$
Thus according as $\Psi = A_3$ or $D_3$ we have $c(\Psi)_{ss} = c(\Psi)_u = 40$ or $48$, so we take $\d_0 = (24, 24, 16)$ or $(16, 16, 16, 16)$; using Corollary~\ref{cor: B values for small k} we then have $B_{\d_0, 2} = 78 > 56 = \dim {u_\Psi}^G$ or $B_{\d_0, 2} = 94 > M$.  Taking the smaller of the two lower bounds, we see that we need only consider semisimple classes $s^G$ with $|\Phi(s)| \leq M - 78 = 6$, and unipotent classes of dimension at least $78$; since each of the former has a subsystem of type $D_3$ disjoint from $\Phi(s)$, and each of the latter has the class $D_3$ in its closure by Lemma~\ref{lem: various classes in D_ell for fixed ell}(iii), we may actually take the larger of the two lower bounds. Therefore if $k \in [2, \frac{d}{2}]$ the quadruple $(G, \lambda, p, k)$ satisfies $\ssdiamcon$ and $\udiamcon$.
\end{proof}

\begin{prop}\label{prop: D_8, omega_8, k geq 2, nets}
Let $G = D_8$ and $\lambda = \omega_8$; then for $k \in [2, \frac{d}{2}]$ the quadruple $(G, \lambda, p, k)$ satisfies $\ssdiamcon$ and $\udiamcon$.
\end{prop}

\begin{proof}
The weight table is as follows.
$$
\begin{array}{|*4{>{\ss}c|}}
\hline
i & \mu & |W.\mu| & m_\mu \\
\hline
1 & \omega_8 & 128 & 1 \\
\hline
\end{array}
$$
We have $M = 112$; we take $k_0 = 2$.

Take $\Psi = \langle \alpha_1 \rangle$ of type $A_1$. The $\Psi$-net table is as follows.
$$
\begin{array}{|*5{>{\ss}c|}}
\hline
\multicolumn{2}{|>{\ss}c|}{\Psi\mathrm{-nets}} & & \multicolumn{1}{|>{\ss}c|}{c(s)} & \multicolumn{1}{|>{\ss}c|}{c(u_\Psi)} \\
\cline{1-2} \cline{4-5}
    \bar\nu    & n_1 &  m & r \geq 2 & p \geq 2 \\
\hline
     \bom_1    &  2  & 32 &    32    &    32    \\
       0       &  1  & 64 &          &          \\
\hline
\multicolumn{3}{c|}{}     &    32    &    32    \\
\cline{4-5}
\end{array}
$$
Thus $c(\Psi)_{ss} = c(\Psi)_u = 32$, so we take $\d_0 = (96, 32)$; using Proposition~\ref{prop: B value when t = 2} we then have $B_{\d_0, 2} = 64 > 26 = \dim {u_\Psi}^G$. We need only consider semisimple classes $s^G$ with $|\Phi(s)| \leq M - 64 = 48$, each of which has a subsystem of type ${A_1}^2$ or a subsystem of type $D_2$ disjoint from $\Phi(s)$, and unipotent classes of dimension at least $64$, each of which has the class ${A_1}^2$ or the class $D_2$ in its closure by Lemma~\ref{lem: A_1^2 or D_2 in D_ell}.

Now take $\Psi = \langle \alpha_1, \alpha_3 \rangle$ of type ${A_1}^2$, and $\Psi = \langle \alpha_7, \alpha_8 \rangle$ of type $D_2$. The $\Psi$-net tables are as follows.
$$
\begin{array}{|*5{>{\ss}c|}}
\hline
\multicolumn{2}{|>{\ss}c|}{\Psi\mathrm{-nets}} & & \multicolumn{1}{|>{\ss}c|}{c(s)} & \multicolumn{1}{|>{\ss}c|}{c(u_\Psi)} \\
\cline{1-2} \cline{4-5}
      \bar\nu      & n_1 &  m & r \geq 2 & p \geq 2 \\
\hline
  \bom_1 + \bom_3  &  4  &  8 &    16    &    16    \\
       \bom_1      &  2  & 16 &    16    &    16    \\
       \bom_3      &  2  & 16 &    16    &    16    \\
         0         &  1  & 32 &          &          \\
\hline
\multicolumn{3}{c|}{}         &    48    &    48    \\
\cline{4-5}
\end{array}
\qquad
\begin{array}{|*5{>{\ss}c|}}
\hline
\multicolumn{2}{|>{\ss}c|}{\Psi\mathrm{-nets}} & & \multicolumn{1}{|>{\ss}c|}{c(s)} & \multicolumn{1}{|>{\ss}c|}{c(u_\Psi)} \\
\cline{1-2} \cline{4-5}
     \bar\nu     & n_1 &  m & r \geq 2 & p \geq 2 \\
\hline
       \bom_7    &  2  & 32 &    32    &    32    \\
       \bom_8    &  2  & 32 &    32    &    32    \\
\hline
\multicolumn{3}{c|}{}       &    64    &    64    \\
\cline{4-5}
\end{array}
$$
Thus according as $\Psi = {A_1}^2$ or $D_2$ we have $c(\Psi)_{ss} = c(\Psi)_u = 48$ or $64$, so we take $\d_0 = (80, 48)$ or $(64, 64)$; using Proposition~\ref{prop: B value when t = 2} we then have $B_{\d_0, 2} = 96 > 44 = \dim {u_\Psi}^G$ or $B_{\d_0, 2} = 126 > M$. Taking the smaller of the two lower bounds, we see that we need only consider semisimple classes $s^G$ with $|\Phi(s)| \leq M - 96 = 16$, and unipotent classes of dimension at least $96$; since each of the former has a subsystem of type $D_2$ disjoint from $\Phi(s)$, and each of the latter has the class $D_2$ in its closure by Lemma~\ref{lem: various classes in classical groups by dim}(ix), we may actually take the larger of the two lower bounds. Therefore if $k \in [2, \frac{d}{2}]$ the quadruple $(G, \lambda, p, k)$ satisfies $\ssdiamcon$ and $\udiamcon$.
\end{proof}

\begin{prop}\label{prop: E_6, omega_1, k geq 4, nets}
Let $G = E_6$ and $\lambda = \omega_1$; then for $k \in [4, \frac{d}{2}]$ the quadruple $(G, \lambda, p, k)$ satisfies $\ssdiamcon$ and $\udiamcon$.
\end{prop}

\begin{proof}
The weight table is as follows.
$$
\begin{array}{|*4{>{\ss}c|}}
\hline
i & \mu & |W.\mu| & m_\mu \\
\hline
1 & \omega_1 & 27 & 1 \\
\hline
\end{array}
$$
We have $M = 72$, $M_3 = 54$ and $M_2 = 40$; we take $k_0 = 4$.

Take $\Psi = \langle \alpha_1 \rangle$ of type $A_1$. The $\Psi$-net table is as follows.
$$
\begin{array}{|*5{>{\ss}c|}}
\hline
\multicolumn{2}{|>{\ss}c|}{\Psi\mathrm{-nets}} & & \multicolumn{1}{|>{\ss}c|}{c(s)} & \multicolumn{1}{|>{\ss}c|}{c(u_\Psi)} \\
\cline{1-2} \cline{4-5}
    \bar\nu    & n_1 &  m & r \geq 2 & p \geq 2 \\
\hline
     \bom_1    &  2  &  6 &     6    &     6    \\
       0       &  1  & 15 &          &          \\
\hline
\multicolumn{3}{c|}{}     &     6    &     6    \\
\cline{4-5}
\end{array}
$$
Thus $c(\Psi)_{ss} = c(\Psi)_u = 6$, so we take $\d_0 = (21, 6)$; using Proposition~\ref{prop: B value when t = 2} we then have $B_{\d_0, 4} = 24 > 22 = \dim {u_\Psi}^G$. We need only consider semisimple classes $s^G$ with $|\Phi(s)| \leq M - 24 = 48$, each of which has a subsystem of type ${A_1}^2$ disjoint from $\Phi(s)$, and unipotent classes of dimension at least $24$, each of which has the class ${A_1}^2$ in its closure by Lemma~\ref{lem: various classes in E_6}(i).

Now take $\Psi = \langle \alpha_1, \alpha_4 \rangle$ of type ${A_1}^2$. The $\Psi$-net table is as follows.
$$
\begin{array}{|*5{>{\ss}c|}}
\hline
\multicolumn{2}{|>{\ss}c|}{\Psi\mathrm{-nets}} & & \multicolumn{1}{|>{\ss}c|}{c(s)} & \multicolumn{1}{|>{\ss}c|}{c(u_\Psi)} \\
\cline{1-2} \cline{4-5}
     \bar\nu     & n_1 & m & r \geq 2 & p \geq 2 \\
\hline
 \bom_1 + \bom_4 &  4  & 1 &     2    &     2    \\
      \bom_1     &  2  & 4 &     4    &     4    \\
      \bom_4     &  2  & 4 &     4    &     4    \\
        0        &  1  & 7 &          &          \\
\hline
\multicolumn{3}{c|}{}      &    10    &    10    \\
\cline{4-5}
\end{array}
$$
Thus $c(\Psi)_{ss} = c(\Psi)_u = 10$, so we take $\d_0 = (17, 10)$; using Proposition~\ref{prop: B value when t = 2} we then have $B_{\d_0, 4} = 40 > 32 = \dim {u_\Psi}^G$. We need only consider semisimple classes $s^G$ with $|\Phi(s)| \leq M - 40 = 32$, each of which has a subsystem of type ${A_1}^3$ disjoint from $\Phi(s)$, and unipotent classes of dimension at least $40$, each of which has the class ${A_1}^3$ in its closure by Lemma~\ref{lem: various classes in E_6}(ii).

Now take $\Psi = \langle \alpha_1, \alpha_4, \alpha_6 \rangle$ of type ${A_1}^3$. The $\Psi$-net table is as follows.
$$
\begin{array}{|*5{>{\ss}c|}}
\hline
\multicolumn{2}{|>{\ss}c|}{\Psi\mathrm{-nets}} & & \multicolumn{1}{|>{\ss}c|}{c(s)} & \multicolumn{1}{|>{\ss}c|}{c(u_\Psi)} \\
\cline{1-2} \cline{4-5}
     \bar\nu     & n_1 & m & r \geq 2 & p \geq 2 \\
\hline
 \bom_1 + \bom_4 &  4  & 1 &     2    &     2    \\
 \bom_1 + \bom_6 &  4  & 1 &     2    &     2    \\
 \bom_4 + \bom_6 &  4  & 1 &     2    &     2    \\
      \bom_1     &  2  & 2 &     2    &     2    \\
      \bom_4     &  2  & 2 &     2    &     2    \\
      \bom_6     &  2  & 2 &     2    &     2    \\
        0        &  1  & 3 &          &          \\
\hline
\multicolumn{3}{c|}{}      &    12    &    12    \\
\cline{4-5}
\end{array}
$$
Thus $c(\Psi)_{ss} = c(\Psi)_u = 12$, so we take $\d_0 = (15, 12)$; using Proposition~\ref{prop: B value when t = 2} we then have $B_{\d_0, 4} = 45 > M_2 = 40 = \dim {u_\Psi}^G$. We may therefore assume from now on that $r \geq 3$, and that $p \geq 3$ when we treat unipotent classes. We need only consider semisimple classes $s^G$ with $|\Phi(s)| \leq M - 45 = 27$, each of which has a subsystem of type $A_2A_1$ disjoint from $\Phi(s)$, and unipotent classes of dimension at least $45$, each of which has the class $A_2A_1$ in its closure by Lemma~\ref{lem: various classes in E_6}(iv).

Now take $\Psi = \langle \alpha_1, \alpha_3, \alpha_6 \rangle$ of type $A_2A_1$. The $\Psi$-net table is as follows.
$$
\begin{array}{|*5{>{\ss}c|}}
\hline
\multicolumn{2}{|>{\ss}c|}{\Psi\mathrm{-nets}} & & \multicolumn{1}{|>{\ss}c|}{c(s)} & \multicolumn{1}{|>{\ss}c|}{c(u_\Psi)} \\
\cline{1-2} \cline{4-5}
     \bar\nu     & n_1 & m & r \geq 3 & p \geq 3 \\
\hline
 \bom_3 + \bom_6 &  6  & 1 &     4    &     4    \\
      \bom_1     &  3  & 3 &     6    &     6    \\
      \bom_3     &  3  & 1 &     2    &     2    \\
      \bom_6     &  2  & 3 &     3    &     3    \\
        0        &  1  & 3 &          &          \\
\hline
\multicolumn{3}{c|}{}      &    15    &    15    \\
\cline{4-5}
\end{array}
$$
Thus $c(\Psi)_{ss} = c(\Psi)_u = 15$, so we take $\d_0 = (12, 12, 3)$; according as $\k_0 = (4, 0, 0)$, $(3, 1, 0)$, $(2, 2, 0)$ or $(2, 1, 1)$ we have $B_{\d_0, \k_0} = 60$, $54$, $52$ or $59$, so $B_{\d_0, 4} = 52 > 46 = \dim {u_\Psi}^G$. We need only consider semisimple classes $s^G$ with $|\Phi(s)| \leq M - 52 = 20$, each of which has a subsystem of type $A_2{A_1}^2$ disjoint from $\Phi(s)$, and unipotent classes of dimension at least $52$, each of which has the class $A_2{A_1}^2$ in its closure by Lemma~\ref{lem: various classes in E_6}(v).

Now take $\Psi = \langle \alpha_1, \alpha_2, \alpha_4, \alpha_6 \rangle$ of type $A_2{A_1}^2$. The $\Psi$-net table is as follows.
$$
\begin{array}{|*5{>{\ss}c|}}
\hline
\multicolumn{2}{|>{\ss}c|}{\Psi\mathrm{-nets}} & & \multicolumn{1}{|>{\ss}c|}{c(s)} & \multicolumn{1}{|>{\ss}c|}{c(u_\Psi)} \\
\cline{1-2} \cline{4-5}
     \bar\nu     & n_1 & m & r \geq 3 & p \geq 3 \\
\hline
 \bom_1 + \bom_4 &  6  & 1 &     4    &     4    \\
 \bom_1 + \bom_6 &  4  & 1 &     2    &     2    \\
 \bom_2 + \bom_6 &  6  & 1 &     4    &     4    \\
      \bom_1     &  2  & 1 &     1    &     1    \\
      \bom_2     &  3  & 1 &     2    &     2    \\
      \bom_4     &  3  & 1 &     2    &     2    \\
      \bom_6     &  2  & 1 &     1    &     1    \\
        0        &  1  & 1 &          &          \\
\hline
\multicolumn{3}{c|}{}      &    16    &    16    \\
\cline{4-5}
\end{array}
$$
Thus $c(\Psi)_{ss} = c(\Psi)_u = 16$, so we take $\d_0 = (11, 11, 5)$; according as $\k_0 = (4, 0, 0)$, $(3, 1, 0)$, $(2, 2, 0)$ or $(2, 1, 1)$ we have $B_{\d_0, \k_0} = 64$, $58$, $56$ or $60$, so $B_{\d_0, 4} = 56 > M_3 > 50 = \dim {u_\Psi}^G$. We may therefore assume from now on that $r \geq 5$, and that $p \geq 5$ when we treat unipotent classes. We need only consider semisimple classes $s^G$ with $|\Phi(s)| \leq M - 56 = 16$, each of which has a subsystem of type ${A_2}^2$ disjoint from $\Phi(s)$, and unipotent classes of dimension at least $56$, each of which has the class ${A_2}^2$ in its closure by Lemma~\ref{lem: various classes in E_6}(vi).

Now take $\Psi = \langle \alpha_1, \alpha_3, \alpha_5, \alpha_6 \rangle$ of type ${A_2}^2$. The $\Psi$-net table is as follows.
$$
\begin{array}{|*5{>{\ss}c|}}
\hline
\multicolumn{2}{|>{\ss}c|}{\Psi\mathrm{-nets}} & & \multicolumn{1}{|>{\ss}c|}{c(s)} & \multicolumn{1}{|>{\ss}c|}{c(u_\Psi)} \\
\cline{1-2} \cline{4-5}
     \bar\nu     & n_1 & m & r \geq 5 & p \geq 5 \\
\hline
 \bom_3 + \bom_6 &  9  & 1 &     6    &     6    \\
      \bom_1     &  3  & 3 &     6    &     6    \\
      \bom_5     &  3  & 3 &     6    &     6    \\
\hline
\multicolumn{3}{c|}{}      &    18    &    18    \\
\cline{4-5}
\end{array}
$$
Thus $c(\Psi)_{ss} = c(\Psi)_u = 18$, so we take $\d_0 = (9, 9, 9)$; according as $\k_0 = (4, 0, 0)$, $(3, 1, 0)$, $(2, 2, 0)$ or $(2, 1, 1)$ we have $B_{\d_0, \k_0} = 72$, $66$, $64$ or $62$, so $B_{\d_0, 4} = 62 > 48 = \dim {u_\Psi}^G$. We need only consider semisimple classes $s^G$ with $|\Phi(s)| \leq M - 62 = 10$, each of which has a subsystem of type $A_4A_1$ disjoint from $\Phi(s)$, and unipotent classes of dimension at least $62$, each of which has the class $A_4A_1$ in its closure by Lemma~\ref{lem: various classes in E_6}(vii).

Now take $\Psi = \langle \alpha_1, \alpha_2, \alpha_4, \alpha_5, \alpha_6 \rangle$ of type $A_4A_1$. The $\Psi$-net table is as follows.
$$
\begin{array}{|*5{>{\ss}c|}}
\hline
\multicolumn{2}{|>{\ss}c|}{\Psi\mathrm{-nets}} & & \multicolumn{1}{|>{\ss}c|}{c(s)} & \multicolumn{1}{|>{\ss}c|}{c(u_\Psi)} \\
\cline{1-2} \cline{4-5}
     \bar\nu     & n_1 & m & r \geq 5 & p \geq 5 \\
\hline
 \bom_1 + \bom_6 & 10  & 1 &     8    &     8    \\
      \bom_1     &  2  & 1 &     1    &     1    \\
      \bom_2     &  5  & 1 &     4    &     4    \\
      \bom_4     & 10  & 1 &     8    &     8    \\
\hline
\multicolumn{3}{c|}{}      &    21    &    21    \\
\cline{4-5}
\end{array}
$$
Thus $c(\Psi)_{ss} = c(\Psi)_u = 21$, so we take $\d_0 = (6, 6, 6, 6, 3)$; according as $\k_0 = (4, 0, 0, 0, 0)$, $(3, 1, 0, 0, 0)$, $(2, 2, 0, 0, 0)$, $(2, 1, 1, 0, 0)$ or $(1, 1, 1, 1, 0)$ we have $B_{\d_0, \k_0} = 84$, $78$, $76$, $74$ or $72$, so $B_{\d_0, 4} = 72  = M > 62 = \dim {u_\Psi}^G$. We need only consider regular semisimple classes $s^G$, each of which has a subsystem of type $A_5$ disjoint from $\Phi(s)$, and the regular unipotent class, which has the class $A_5$ in its closure by Lemma~\ref{lem: any class in closure of reg class} (and for the unipotent class $A_5$ to lie in $G_{(p)}$ we need $p \geq 7$).

Now take $\Psi = \langle \alpha_1, \alpha_3, \alpha_4, \alpha_5, \alpha_6 \rangle$ of type $A_5$. The $\Psi$-net table is as follows.
$$
\begin{array}{|*5{>{\ss}c|}}
\hline
\multicolumn{2}{|>{\ss}c|}{\Psi\mathrm{-nets}} & & \multicolumn{1}{|>{\ss}c|}{c(s)} & \multicolumn{1}{|>{\ss}c|}{c(u_\Psi)} \\
\cline{1-2} \cline{4-5}
     \bar\nu     & n_1 & m & r \geq 5 & p \geq 7 \\
\hline
      \bom_1     &  6  & 2 &    10    &    10    \\
      \bom_5     & 15  & 1 &    12    &    12    \\
\hline
\multicolumn{3}{c|}{}      &    22    &    22    \\
\cline{4-5}
\end{array}
$$
Thus $c(\Psi)_{ss} = c(\Psi)_u = 22$, so we take $\d_0 = (5, 5, 5, 5, 5, 2)$; according as $\k_0 = (4, 0, 0, 0, 0, 0)$, $(3, 1, 0, 0, 0, 0)$, $(2, 2, 0, 0, 0, 0)$, $(2, 1, 1, 0, 0, 0)$ or $(1, 1, 1, 1, 0, 0)$ we have $B_{\d_0, \k_0} = 88$, $82$, $80$, $78$ or $76$, so $B_{\d_0, 4} = 76 > M$. Therefore if $k \in [4, \frac{d}{2}]$ the quadruple $(G, \lambda, p, k)$ satisfies $\ssdiamcon$ and $\udiamcon$.
\end{proof}

\begin{prop}\label{prop: E_6, omega_2, k geq 2, nets}
Let $G = E_6$ and $\lambda = \omega_2$; then for $k \in [2, \frac{d}{2}]$ the quadruple $(G, \lambda, p, k)$ satisfies $\ssdiamcon$ and $\udiamcon$.
\end{prop}

\begin{proof}
Write $\z = \z_{p, 3}$. The weight table is as follows.
$$
\begin{array}{|*4{>{\ss}c|}}
\hline
i & \mu & |W.\mu| & m_\mu \\
\hline
1 & \omega_2 & 72 &   1    \\
0 &    0     &  1 & 6 - \z \\
\hline
\end{array}
$$
We have $M = 72$ and $M_2 = 40$; we take $k_0 = 2$.

Take $\Psi = \langle \alpha_1 \rangle$ of type $A_1$. The $\Psi$-net table is as follows.
$$
\begin{array}{|*7{>{\ss}c|}}
\hline
\multicolumn{3}{|>{\ss}c|}{\Psi\mathrm{-nets}} & & \multicolumn{1}{|>{\ss}c|}{c(s)} & \multicolumn{2}{|>{\ss}c|}{c(u_\Psi)} \\
\cline{1-3} \cline{5-7}
    \bar\nu    & n_0 & n_1 &  m & r \geq 2 & p = 2 & p \geq 3 \\
\hline
    2\bom_1    &  1  &  2  &  1 &     2    &   1   &     2    \\
     \bom_1    &  0  &  2  & 20 &    20    &  20   &    20    \\
       0       &  0  &  1  & 30 &          &       &          \\
\hline
\multicolumn{4}{c|}{}           &    22    &  21   &    22    \\
\cline{5-7}
\end{array}
$$
Thus $c(\Psi)_{ss}, c(\Psi)_u \geq 21$, so we may take $\d_0 = (57 - \z, 21)$; using Proposition~\ref{prop: B value when t = 2} we then have $B_{\d_0, 2} = 42 > M_2 > 22 = \dim {u_\Psi}^G$. We may therefore assume from now on that $r \geq 3$, and that $p \geq 3$ when we treat unipotent classes. We need only consider semisimple classes $s^G$ with $|\Phi(s)| \leq M - 42 = 30$, each of which has a subsystem of type $A_2$ disjoint from $\Phi(s)$, and unipotent classes of dimension at least $42$, each of which has the class $A_2$ in its closure by Lemma~\ref{lem: various classes in E_6}(iii).

Now take $\Psi = \langle \alpha_1, \alpha_3 \rangle$ of type $A_2$. The $\Psi$-net table is as follows.
$$
\begin{array}{|*7{>{\ss}c|}}
\hline
\multicolumn{3}{|>{\ss}c|}{\Psi\mathrm{-nets}} & & \multicolumn{1}{|>{\ss}c|}{c(s)} & \multicolumn{2}{|>{\ss}c|}{c(u_\Psi)} \\
\cline{1-3} \cline{5-7}
     \bar\nu     & n_0 & n_1 &  m & r \geq 3 & p = 3 & p \geq 5 \\
\hline
 \bom_1 + \bom_3 &  1  &  6  &  1 &     6    &   4   &     6    \\
      \bom_1     &  0  &  3  &  9 &    18    &  18   &    18    \\
      \bom_3     &  0  &  3  &  9 &    18    &  18   &    18    \\
        0        &  0  &  1  & 12 &          &       &          \\
\hline
\multicolumn{4}{c|}{}             &    42    &  40   &    42    \\
\cline{5-7}
\end{array}
$$
Thus $c(\Psi)_{ss}, c(\Psi)_u \geq 40$, so we may take $\d_0 = (38 - \z, 38 - \z, 2 + \z)$; using Corollary~\ref{cor: B values for small k} we then have $B_{\d_0, 2} = 78 > M$. Therefore if $k \in [2, \frac{d}{2}]$ the quadruple $(G, \lambda, p, k)$ satisfies $\ssdiamcon$ and $\udiamcon$.
\end{proof}

\begin{prop}\label{prop: E_7, omega_1, k geq 2, nets}
Let $G = E_7$ and $\lambda = \omega_1$; then for $k \in [2, \frac{d}{2}]$ the quadruple $(G, \lambda, p, k)$ satisfies $\ssdiamcon$ and $\udiamcon$.
\end{prop}

\begin{proof}
Write $\z = \z_{p, 2}$. The weight table is as follows.
$$
\begin{array}{|*4{>{\ss}c|}}
\hline
i & \mu & |W.\mu| & m_\mu \\
\hline
1 & \omega_1 & 126 &   1    \\
0 &    0     &   1 & 7 - \z \\
\hline
\end{array}
$$
We have $M = 126$, $M_3 = 90$ and $M_2 = 70$; we take $k_0 = 2$.

Take $\Psi = \langle \alpha_1 \rangle$ of type $A_1$. The $\Psi$-net table is as follows.
$$
\begin{array}{|*7{>{\ss}c|}}
\hline
\multicolumn{3}{|>{\ss}c|}{\Psi\mathrm{-nets}} & & \multicolumn{1}{|>{\ss}c|}{c(s)} & \multicolumn{2}{|>{\ss}c|}{c(u_\Psi)} \\
\cline{1-3} \cline{5-7}
    \bar\nu    & n_0 & n_1 &  m & r \geq 2 & p = 2 & p \geq 3 \\
\hline
    2\bom_1    &  1  &  2  &  1 &     2    &   1   &     2    \\
     \bom_1    &  0  &  2  & 32 &    32    &  32   &    32    \\
       0       &  0  &  1  & 60 &          &       &          \\
\hline
\multicolumn{4}{c|}{}           &    34    &  33   &    34    \\
\cline{5-7}
\end{array}
$$
Thus $c(\Psi)_{ss}, c(\Psi)_u \geq 33$, so we may take $\d_0 = (100 - \z, 33)$; using Proposition~\ref{prop: B value when t = 2} we then have $B_{\d_0, 2} = 66 > 34 = \dim {u_\Psi}^G$. We need only consider semisimple classes $s^G$ with $|\Phi(s)| \leq M - 66 = 60$, each of which has a subsystem of type ${A_1}^2$ disjoint from $\Phi(s)$, and unipotent classes of dimension at least $66$, each of which has the class ${A_1}^2$ in its closure by Lemma~\ref{lem: various classes in E_7}(i).

Now take $\Psi = \langle \alpha_1, \alpha_4 \rangle$ of type ${A_1}^2$. The $\Psi$-net table is as follows.
$$
\begin{array}{|*7{>{\ss}c|}}
\hline
\multicolumn{3}{|>{\ss}c|}{\Psi\mathrm{-nets}} & & \multicolumn{1}{|>{\ss}c|}{c(s)} & \multicolumn{2}{|>{\ss}c|}{c(u_\Psi)} \\
\cline{1-3} \cline{5-7}
     \bar\nu     & n_0 & n_1 &  m & r \geq 2 & p = 2 & p \geq 3 \\
\hline
 2\bom_1/2\bom_4 &  1  &  4  &  1 &     4    &   2   &     4    \\
 \bom_1 + \bom_4 &  0  &  4  &  8 &    16    &  16   &    16    \\
      \bom_1     &  0  &  2  & 16 &    16    &  16   &    16    \\
      \bom_4     &  0  &  2  & 16 &    16    &  16   &    16    \\
        0        &  0  &  1  & 26 &          &       &          \\
\hline
\multicolumn{4}{c|}{}             &    52    &  50   &    52    \\
\cline{5-7}
\end{array}
$$
Thus $c(\Psi)_{ss}, c(\Psi)_u \geq 50$, so we may take $\d_0 = (83 - \z, 50)$; using Proposition~\ref{prop: B value when t = 2} we then have $B_{\d_0, 2} = 100 > M_3 > 52 = \dim {u_\Psi}^G$. We may therefore assume from now on that $r \geq 5$, and that $p \geq 5$ when we treat unipotent classes. We need only consider semisimple classes $s^G$ with $|\Phi(s)| \leq M - 100 = 26$, each of which has a subsystem of type $A_3$ disjoint from $\Phi(s)$, and unipotent classes of dimension at least $100$, each of which has the class $A_3$ in its closure by Lemma~\ref{lem: various classes in E_7}(v).

Now take $\Psi = \langle \alpha_1, \alpha_3, \alpha_4 \rangle$ of type $A_3$. The $\Psi$-net table is as follows.
$$
\begin{array}{|*6{>{\ss}c|}}
\hline
\multicolumn{3}{|>{\ss}c|}{\Psi\mathrm{-nets}} & & \multicolumn{1}{|>{\ss}c|}{c(s)} & \multicolumn{1}{|>{\ss}c|}{c(u_\Psi)} \\
\cline{1-3} \cline{5-6}
     \bar\nu     & n_0 & n_1 &  m & r \geq 5 & p \geq 5 \\
\hline
 \bom_1 + \bom_4 &  1  & 12  &  1 &    12    &    12    \\
      \bom_1     &  0  &  4  &  8 &    24    &    24    \\
      \bom_3     &  0  &  6  &  6 &    24    &    24    \\
      \bom_4     &  0  &  4  &  8 &    24    &    24    \\
        0        &  0  &  1  & 14 &          &          \\
\hline
\multicolumn{4}{c|}{}             &    84    &    84    \\
\cline{5-6}
\end{array}
$$
Thus $c(\Psi)_{ss} = c(\Psi)_u = 84$, so we take $\d_0 = (49 - \z, 49 - \z, 35 + \z)$; using Corollary~\ref{cor: B values for small k} we then have $B_{\d_0, 2} = 166 > M$. Therefore if $k \in [2, \frac{d}{2}]$ the quadruple $(G, \lambda, p, k)$ satisfies $\ssdiamcon$ and $\udiamcon$.
\end{proof}

\begin{prop}\label{prop: E_7, omega_7, k geq 3, nets}
Let $G = E_7$ and $\lambda = \omega_7$; then for $k \in [3, \frac{d}{2}]$ the quadruple $(G, \lambda, p, k)$ satisfies $\ssdiamcon$ and $\udiamcon$.
\end{prop}

\begin{proof}
The weight table is as follows.
$$
\begin{array}{|*4{>{\ss}c|}}
\hline
i & \mu & |W.\mu| & m_\mu \\
\hline
1 & \omega_7 & 56 & 1 \\
\hline
\end{array}
$$
We have $M = 126$, $M_5 = 106$, $M_3 = 90$ and $M_2 = 70$; we take $k_0 = 3$.

Take $\Psi = \langle \alpha_1 \rangle$ of type $A_1$. The $\Psi$-net table is as follows.
$$
\begin{array}{|*5{>{\ss}c|}}
\hline
\multicolumn{2}{|>{\ss}c|}{\Psi\mathrm{-nets}} & & \multicolumn{1}{|>{\ss}c|}{c(s)} & \multicolumn{1}{|>{\ss}c|}{c(u_\Psi)} \\
\cline{1-2} \cline{4-5}
    \bar\nu    & n_1 &  m & r \geq 2 & p \geq 2 \\
\hline
     \bom_1    &  2  & 12 &    12    &    12    \\
       0       &  1  & 32 &          &          \\
\hline
\multicolumn{3}{c|}{}     &    12    &    12    \\
\cline{4-5}
\end{array}
$$
Thus $c(\Psi)_{ss} = c(\Psi)_u = 12$, so we take $\d_0 = (44, 12)$; using Proposition~\ref{prop: B value when t = 2} we then have $B_{\d_0, 3} = 36 > 34 = \dim {u_\Psi}^G$. We need only consider semisimple classes $s^G$ with $|\Phi(s)| \leq M - 36 = 90$, each of which has a subsystem of type ${A_1}^2$ disjoint from $\Phi(s)$, and unipotent classes of dimension at least $36$, each of which has the class ${A_1}^2$ in its closure by Lemma~\ref{lem: various classes in E_7}(i).

Now take $\Psi = \langle \alpha_1, \alpha_4 \rangle$ of type ${A_1}^2$. The $\Psi$-net table is as follows.
$$
\begin{array}{|*5{>{\ss}c|}}
\hline
\multicolumn{2}{|>{\ss}c|}{\Psi\mathrm{-nets}} & & \multicolumn{1}{|>{\ss}c|}{c(s)} & \multicolumn{1}{|>{\ss}c|}{c(u_\Psi)} \\
\cline{1-2} \cline{4-5}
     \bar\nu     & n_1 &  m & r \geq 2 & p \geq 2 \\
\hline
 \bom_1 + \bom_4 &  4  &  2 &     4    &     4    \\
      \bom_1     &  2  &  8 &     8    &     8    \\
      \bom_4     &  2  &  8 &     8    &     8    \\
        0        &  1  & 16 &          &          \\
\hline
\multicolumn{3}{c|}{}       &    20    &    20    \\
\cline{4-5}
\end{array}
$$
Thus $c(\Psi)_{ss} = c(\Psi)_u = 20$, so we take $\d_0 = (36, 20)$; using Proposition~\ref{prop: B value when t = 2} we then have $B_{\d_0, 3} = 60 > 52 = \dim {u_\Psi}^G$. We need only consider semisimple classes $s^G$ with $|\Phi(s)| \leq M - 60 = 66$, each of which has a subsystem of type $({A_1}^3)'$ disjoint from $\Phi(s)$, and unipotent classes of dimension at least $60$, each of which has the class $({A_1}^3)'$ in its closure by Lemma~\ref{lem: various classes in E_7}(ii).

Now take $\Psi = \langle \alpha_1, \alpha_4, \alpha_6 \rangle$ of type $({A_1}^3)'$. The $\Psi$-net table is as follows.
$$
\begin{array}{|*5{>{\ss}c|}}
\hline
\multicolumn{2}{|>{\ss}c|}{\Psi\mathrm{-nets}} & & \multicolumn{1}{|>{\ss}c|}{c(s)} & \multicolumn{1}{|>{\ss}c|}{c(u_\Psi)} \\
\cline{1-2} \cline{4-5}
     \bar\nu     & n_1 & m & r \geq 2 & p \geq 2 \\
\hline
 \bom_1 + \bom_4 &  4  & 2 &     4    &     4    \\
 \bom_1 + \bom_6 &  4  & 2 &     4    &     4    \\
 \bom_4 + \bom_6 &  4  & 2 &     4    &     4    \\
      \bom_1     &  2  & 4 &     4    &     4    \\
      \bom_4     &  2  & 4 &     4    &     4    \\
      \bom_6     &  2  & 4 &     4    &     4    \\
        0        &  1  & 8 &          &          \\
\hline
\multicolumn{3}{c|}{}      &    24    &    24    \\
\cline{4-5}
\end{array}
$$
Thus $c(\Psi)_{ss} = c(\Psi)_u = 24$, so we take $\d_0 = (32, 24)$; using Proposition~\ref{prop: B value when t = 2} we then have $B_{\d_0, 3} = 72 > M_2 > 64 = \dim {u_\Psi}^G$. We may therefore assume from now on that $r \geq 3$, and that $p \geq 3$ when we treat unipotent classes. We need only consider semisimple classes $s^G$ with $|\Phi(s)| \leq M - 72 = 54$, each of which has a subsystem of type $A_2A_1$ disjoint from $\Phi(s)$, and unipotent classes of dimension at least $72$, each of which has the class $A_2A_1$ in its closure by Lemma~\ref{lem: various classes in E_7}(iii).

Now take $\Psi = \langle \alpha_1, \alpha_3, \alpha_5 \rangle$ of type $A_2A_1$. The $\Psi$-net table is as follows.
$$
\begin{array}{|*5{>{\ss}c|}}
\hline
\multicolumn{2}{|>{\ss}c|}{\Psi\mathrm{-nets}} & & \multicolumn{1}{|>{\ss}c|}{c(s)} & \multicolumn{1}{|>{\ss}c|}{c(u_\Psi)} \\
\cline{1-2} \cline{4-5}
     \bar\nu     & n_1 & m & r \geq 3 & p \geq 3 \\
\hline
 \bom_1 + \bom_5 &  6  & 1 &     4    &     4    \\
 \bom_3 + \bom_5 &  6  & 1 &     4    &     4    \\
      \bom_1     &  3  & 4 &     8    &     8    \\
      \bom_3     &  3  & 4 &     8    &     8    \\
      \bom_5     &  2  & 6 &     6    &     6    \\
        0        &  1  & 8 &          &          \\
\hline
\multicolumn{3}{c|}{}      &    30    &    30    \\
\cline{4-5}
\end{array}
$$
Thus $c(\Psi)_{ss} = c(\Psi)_u = 30$, so we take $\d_0 = (26, 26, 4)$; using Corollary~\ref{cor: B values for small k} we then have $B_{\d_0, 3} = 86 > 76 = \dim {u_\Psi}^G$. We need only consider semisimple classes $s^G$ with $|\Phi(s)| \leq M - 86 = 40$, each of which has a subsystem of type $A_2{A_1}^2$ disjoint from $\Phi(s)$, and unipotent classes of dimension at least $86$, each of which has the class $A_2{A_1}^2$ in its closure by Lemma~\ref{lem: various classes in E_7}(iv).

Now take $\Psi = \langle \alpha_1, \alpha_3, \alpha_5, \alpha_7 \rangle$ of type $A_2{A_1}^2$. The $\Psi$-net table is as follows.
$$
\begin{array}{|*5{>{\ss}c|}}
\hline
\multicolumn{2}{|>{\ss}c|}{\Psi\mathrm{-nets}} & & \multicolumn{1}{|>{\ss}c|}{c(s)} & \multicolumn{1}{|>{\ss}c|}{c(u_\Psi)} \\
\cline{1-2} \cline{4-5}
     \bar\nu     & n_1 & m & r \geq 3 & p \geq 3 \\
\hline
 \bom_1 + \bom_5 &  6  & 1 &     4    &     4    \\
 \bom_1 + \bom_7 &  6  & 1 &     4    &     4    \\
 \bom_3 + \bom_5 &  6  & 1 &     4    &     4    \\
 \bom_3 + \bom_7 &  6  & 1 &     4    &     4    \\
 \bom_5 + \bom_7 &  4  & 2 &     4    &     4    \\
      \bom_1     &  3  & 2 &     4    &     4    \\
      \bom_3     &  3  & 2 &     4    &     4    \\
      \bom_5     &  2  & 2 &     2    &     2    \\
      \bom_7     &  2  & 2 &     2    &     2    \\
        0        &  1  & 4 &          &          \\
\hline
\multicolumn{3}{c|}{}      &    32    &    32    \\
\cline{4-5}
\end{array}
$$
Thus $c(\Psi)_{ss} = c(\Psi)_u = 32$, so we take $\d_0 = (24, 24, 8)$; using Corollary~\ref{cor: B values for small k} we then have $B_{\d_0, 3} = 92 > M_3 > 82 = \dim {u_\Psi}^G$. We may therefore assume from now on that $r \geq 5$, and that $p \geq 5$ when we treat unipotent classes. We need only consider semisimple classes $s^G$ with $|\Phi(s)| \leq M - 92 = 34$, each of which has a subsystem of type $(A_3A_1)'$ disjoint from $\Phi(s)$, and unipotent classes of dimension at least $92$, each of which has the class $(A_3A_1)'$ in its closure by Lemma~\ref{lem: various classes in E_7}(vi).

Now take $\Psi = \langle \alpha_1, \alpha_5, \alpha_6, \alpha_7 \rangle$ of type $(A_3A_1)'$. The $\Psi$-net table is as follows.
$$
\begin{array}{|*5{>{\ss}c|}}
\hline
\multicolumn{2}{|>{\ss}c|}{\Psi\mathrm{-nets}} & & \multicolumn{1}{|>{\ss}c|}{c(s)} & \multicolumn{1}{|>{\ss}c|}{c(u_\Psi)} \\
\cline{1-2} \cline{4-5}
     \bar\nu     & n_1 & m & r \geq 5 & p \geq 5 \\
\hline
 \bom_1 + \bom_5 &  8  & 1 &     6    &     6    \\
 \bom_1 + \bom_7 &  8  & 1 &     6    &     6    \\
      \bom_1     &  2  & 4 &     4    &     4    \\
      \bom_5     &  4  & 2 &     6    &     6    \\
      \bom_6     &  6  & 2 &     8    &     8    \\
      \bom_7     &  4  & 2 &     6    &     6    \\
        0        &  1  & 4 &          &          \\
\hline
\multicolumn{3}{c|}{}      &    36    &    36    \\
\cline{4-5}
\end{array}
$$
Thus $c(\Psi)_{ss} = c(\Psi)_u = 36$, so we take $\d_0 = (20, 20, 16)$; using Corollary~\ref{cor: B values for small k} we then have $B_{\d_0, 3} = 104 > 92 = \dim {u_\Psi}^G$. We need only consider semisimple classes $s^G$ with $|\Phi(s)| \leq M - 104 = 22$, each of which has a subsystem of type $A_4A_1$ disjoint from $\Phi(s)$, and unipotent classes of dimension at least $104$, each of which has the class $A_4A_1$ in its closure by Lemma~\ref{lem: various classes in E_7}(vii).

Now take $\Psi = \langle \alpha_1, \alpha_4, \alpha_5, \alpha_6, \alpha_7 \rangle$ of type $A_4A_1$. The $\Psi$-net table is as follows.
$$
\begin{array}{|*5{>{\ss}c|}}
\hline
\multicolumn{2}{|>{\ss}c|}{\Psi\mathrm{-nets}} & & \multicolumn{1}{|>{\ss}c|}{c(s)} & \multicolumn{1}{|>{\ss}c|}{c(u_\Psi)} \\
\cline{1-2} \cline{4-5}
     \bar\nu     & n_1 & m & r \geq 5 & p \geq 5 \\
\hline
 \bom_1 + \bom_4 & 10  & 1 &     8    &     8    \\
 \bom_1 + \bom_7 & 10  & 1 &     8    &     8    \\
      \bom_1     &  2  & 2 &     2    &     2    \\
      \bom_4     &  5  & 1 &     4    &     4    \\
      \bom_5     & 10  & 1 &     8    &     8    \\
      \bom_6     & 10  & 1 &     8    &     8    \\
      \bom_7     &  5  & 1 &     4    &     4    \\
        0        &  1  & 2 &          &          \\
\hline
\multicolumn{3}{c|}{}      &    42    &    42    \\
\cline{4-5}
\end{array}
$$
Thus $c(\Psi)_{ss} = c(\Psi)_u = 42$, so we take $\d_0 = (14, 14, 14, 14)$; using Corollary~\ref{cor: B values for small k} we then have $B_{\d_0, 3} = 120 > M_5 > 104 = \dim {u_\Psi}^G$. We may therefore assume from now on that $r \geq 7$, and that $p \geq 7$ when we treat unipotent classes. We need only consider semisimple classes $s^G$ with $|\Phi(s)| \leq M - 120 = 6$, each of which has a subsystem of type $A_6$ disjoint from $\Phi(s)$, and unipotent classes of dimension at least $120$, each of which has the class $A_6$ in its closure by Lemma~\ref{lem: various classes in E_7}(viii).

Now take $\Psi = \langle \alpha_1, \alpha_3, \alpha_4, \alpha_5, \alpha_6, \alpha_7 \rangle$ of type $A_6$. The $\Psi$-net table is as follows.
$$
\begin{array}{|*5{>{\ss}c|}}
\hline
\multicolumn{2}{|>{\ss}c|}{\Psi\mathrm{-nets}} & & \multicolumn{1}{|>{\ss}c|}{c(s)} & \multicolumn{1}{|>{\ss}c|}{c(u_\Psi)} \\
\cline{1-2} \cline{4-5}
     \bar\nu     & n_1 & m & r \geq 7 & p \geq 7 \\
\hline
      \bom_1     &  7  & 1 &     6    &     6    \\
      \bom_3     & 21  & 1 &    18    &    18    \\
      \bom_6     & 21  & 1 &    18    &    18    \\
      \bom_7     &  7  & 1 &     6    &     6    \\
\hline
\multicolumn{3}{c|}{}      &    48    &    48    \\
\cline{4-5}
\end{array}
$$
Thus $c(\Psi)_{ss} = c(\Psi)_u = 48$, so we take $\d_0 = (8, 8, 8, 8, 8, 8, 8)$; using Corollary~\ref{cor: B values for small k} we then have $B_{\d_0, 3} = 138 > M$. Therefore if $k \in [3, \frac{d}{2}]$ the quadruple $(G, \lambda, p, k)$ satisfies $\ssdiamcon$ and $\udiamcon$.
\end{proof}

\begin{prop}\label{prop: E_8, omega_8, k geq 2, nets}
Let $G = E_8$ and $\lambda = \omega_8$; then for $k \in [2, \frac{d}{2}]$ the quadruple $(G, \lambda, p, k)$ satisfies $\ssdiamcon$ and $\udiamcon$.
\end{prop}

\begin{proof}
The weight table is as follows.
$$
\begin{array}{|*4{>{\ss}c|}}
\hline
i & \mu & |W.\mu| & m_\mu \\
\hline
1 & \omega_8 & 240 & 1 \\
0 &    0     &   1 & 8 \\
\hline
\end{array}
$$
We have $M = 240$, $M_3 = 168$ and $M_2 = 128$; we take $k_0 = 2$.

Take $\Psi = \langle \alpha_1 \rangle$ of type $A_1$. The $\Psi$-net table is as follows.
$$
\begin{array}{|*7{>{\ss}c|}}
\hline
\multicolumn{3}{|>{\ss}c|}{\Psi\mathrm{-nets}} & & \multicolumn{1}{|>{\ss}c|}{c(s)} & \multicolumn{2}{|>{\ss}c|}{c(u_\Psi)} \\
\cline{1-3} \cline{5-7}
    \bar\nu    & n_0 & n_1 &  m  & r \geq 2 & p = 2 & p \geq 3 \\
\hline
    2\bom_1    &  1  &  2  &   1 &     2    &   1   &     2    \\
     \bom_1    &  0  &  2  &  56 &    56    &  56   &    56    \\
       0       &  0  &  1  & 126 &          &       &          \\
\hline
\multicolumn{4}{c|}{}            &    58    &  57   &    58    \\
\cline{5-7}
\end{array}
$$
Thus $c(\Psi)_{ss}, c(\Psi)_u \geq 57$, so we may take $\d_0 = (191, 57)$; using Proposition~\ref{prop: B value when t = 2} we then have $B_{\d_0, 2} = 114 > 58 = \dim {u_\Psi}^G$. We need only consider semisimple classes $s^G$ with $|\Phi(s)| \leq M - 114 = 126$, each of which has a subsystem of type ${A_1}^2$ disjoint from $\Phi(s)$, and unipotent classes of dimension at least $114$, each of which has the class ${A_1}^2$ in its closure by Lemma~\ref{lem: various classes in E_8}(i).

Now take $\Psi = \langle \alpha_1, \alpha_4 \rangle$ of type ${A_1}^2$. The $\Psi$-net table is as follows.
$$
\begin{array}{|*7{>{\ss}c|}}
\hline
\multicolumn{3}{|>{\ss}c|}{\Psi\mathrm{-nets}} & & \multicolumn{1}{|>{\ss}c|}{c(s)} & \multicolumn{2}{|>{\ss}c|}{c(u_\Psi)} \\
\cline{1-3} \cline{5-7}
     \bar\nu     & n_0 & n_1 &  m & r \geq 2 & p = 2 & p \geq 3 \\
\hline
 2\bom_1/2\bom_4 &  1  &  4  &  1 &     4    &   2   &     4    \\
 \bom_1 + \bom_4 &  0  &  4  & 12 &    24    &  24   &    24    \\
      \bom_1     &  0  &  2  & 32 &    32    &  32   &    32    \\
      \bom_4     &  0  &  2  & 32 &    32    &  32   &    32    \\
        0        &  0  &  1  & 60 &          &       &          \\
\hline
\multicolumn{4}{c|}{}             &    92    &  90   &    92    \\
\cline{5-7}
\end{array}
$$
Thus $c(\Psi)_{ss}, c(\Psi)_u \geq 90$, so we may take $\d_0 = (158, 90)$; using Proposition~\ref{prop: B value when t = 2} we then have $B_{\d_0, 2} = 180 > M_3 > 92 = \dim {u_\Psi}^G$. We may therefore assume from now on that $r \geq 5$, and that $p \geq 5$ when we treat unipotent classes. We need only consider semisimple classes $s^G$ with $|\Phi(s)| \leq M - 180 = 60$, each of which has a subsystem of type $A_3$ disjoint from $\Phi(s)$, and unipotent classes of dimension at least $180$, each of which has the class $A_3$ in its closure by Lemma~\ref{lem: various classes in E_8}(ii).

Now take $\Psi = \langle \alpha_1, \alpha_3, \alpha_4 \rangle$ of type $A_3$. The $\Psi$-net table is as follows.
$$
\begin{array}{|*6{>{\ss}c|}}
\hline
\multicolumn{3}{|>{\ss}c|}{\Psi\mathrm{-nets}} & & \multicolumn{1}{|>{\ss}c|}{c(s)} & \multicolumn{1}{|>{\ss}c|}{c(u_\Psi)} \\
\cline{1-3} \cline{5-6}
     \bar\nu     & n_0 & n_1 &  m & r \geq 5 & p \geq 5 \\
\hline
 \bom_1 + \bom_4 &  1  & 12  &  1 &    12    &    12    \\
      \bom_1     &  0  &  4  & 16 &    48    &    48    \\
      \bom_3     &  0  &  6  & 10 &    40    &    40    \\
      \bom_4     &  0  &  4  & 16 &    48    &    48    \\
        0        &  0  &  1  & 40 &          &          \\
\hline
\multicolumn{4}{c|}{}             &   148    &   148    \\
\cline{5-6}
\end{array}
$$
Thus $c(\Psi)_{ss} = c(\Psi)_u = 148$, so we take $\d_0 = (100, 100, 48)$; using Corollary~\ref{cor: B values for small k} we then have $B_{\d_0, 2} = 294 > M$. Therefore if $k \in [2, \frac{d}{2}]$ the quadruple $(G, \lambda, p, k)$ satisfies $\ssdiamcon$ and $\udiamcon$.
\end{proof}

\begin{prop}\label{prop: B_2, omega_1 + omega_2, k geq 2, nets}
Let $G = B_2$ and $\lambda = \omega_1 + \omega_2$ with $p = 5$; then for $k \in [2, \frac{d}{2}]$ the quadruple $(G, \lambda, p, k)$ satisfies $\ssdiamcon$ and $\udiamcon$.
\end{prop}

\begin{proof}
The weight table is as follows.
$$
\begin{array}{|*4{>{\ss}c|}}
\hline
i & \mu & |W.\mu| & m_\mu \\
\hline
2 & \omega_1 + \omega_2 & 8 & 1 \\
1 &       \omega_2      & 4 & 1 \\
\hline
\end{array}
$$
We have $M = 8$; we take $k_0 = 2$.

Take $\Psi = \langle \alpha_1 \rangle$ of type $A_1$. The $\Psi$-net table is as follows.
$$
\begin{array}{|*5{>{\ss}c|}}
\hline
\multicolumn{3}{|>{\ss}c|}{\Psi\mathrm{-nets}} & & \multicolumn{1}{|>{\ss}c|}{c(u_\Psi)} \\
\cline{1-3} \cline{5-5}
    \bar\nu    & n_1 & n_2 & m & p = 5 \\
\hline
    2\bom_1    &  1  &  2  & 2 &   4   \\
     \bom_1    &  0  &  2  & 2 &   2   \\
     \bom_1    &  2  &  0  & 1 &   1   \\
\hline
\multicolumn{4}{c|}{}          &   7   \\
\cline{5-5}
\end{array}
$$
Thus $c(\Psi)_u = 7$, so we take $\d_0 = (5, 5, 2)$; using Corollary~\ref{cor: B values for small k} we then have $B_{\d_0, 2} = 12 > M$. Each of the remaining non-trivial unipotent classes has $B_1$ in its closure by Lemma~\ref{lem: root elt class in closure of any non-triv class}.

Now take $\Psi = \langle \alpha_2 \rangle$ of type $B_1$. The $\Psi$-net table is as follows.
$$
\begin{array}{|*8{>{\ss}c|}}
\hline
\multicolumn{3}{|>{\ss}c|}{\Psi\mathrm{-nets}} & & \multicolumn{3}{|>{\ss}c|}{c(s)} & \multicolumn{1}{|>{\ss}c|}{c(u_\Psi)} \\
\cline{1-3} \cline{5-8}
    \bar\nu    & n_1 & n_2 & m & r = 2 & r = 3 & r \geq 7 & p = 5 \\
\hline
    3\bom_2    &  2  &  2  & 2 &   4   &   4   &     6    &   6   \\
     \bom_2    &  0  &  2  & 2 &   2   &   2   &     2    &   2   \\
\hline
\multicolumn{4}{c|}{}          &   6   &   6   &     8    &   8   \\
\cline{5-8}
\end{array}
$$
Thus $c(\Psi)_{ss}, c(\Psi)_u \geq 6$, so we may take $\d_0 = (6, 6)$; using Proposition~\ref{prop: B value when t = 2} we then have $B_{\d_0, 2} = 10 > M$. Therefore if $k \in [2, \frac{d}{2}]$ the quadruple $(G, \lambda, p, k)$ satisfies $\ssdiamcon$ and $\udiamcon$.
\end{proof}

\begin{prop}\label{prop: B_2, 2omega_2, k geq 2, nets}
Let $G = B_2$ and $\lambda = 2\omega_2$ with $p \geq 3$; then for $k \in [2, \frac{d}{2}]$ the quadruple $(G, \lambda, p, k)$ satisfies $\ssdiamcon$ and $\udiamcon$.
\end{prop}

\begin{proof}
The weight table is as follows.
$$
\begin{array}{|*4{>{\ss}c|}}
\hline
i & \mu & |W.\mu| & m_\mu \\
\hline
2 & 2\omega_2 & 4 & 1 \\
1 &  \omega_1 & 4 & 1 \\
0 &     0     & 1 & 2 \\
\hline
\end{array}
$$
We have $M = 8$ and $M_2 = 6$; we take $k_0 = 2$.

Take $\Psi = \langle \alpha_1 \rangle$ of type $A_1$. The $\Psi$-net table is as follows.
$$
\begin{array}{|*6{>{\ss}c|}}
\hline
\multicolumn{4}{|>{\ss}c|}{\Psi\mathrm{-nets}} & & \multicolumn{1}{|>{\ss}c|}{c(u_\Psi)} \\
\cline{1-4} \cline{6-6}
    \bar\nu    & n_0 & n_1 & n_2 & m & p \geq 3 \\
\hline
    2\bom_1    &  1  &  0  &  2  & 1 &     2    \\
     \bom_1    &  0  &  2  &  0  & 2 &     2    \\
       0       &  0  &  0  &  1  & 2 &          \\
\hline
\multicolumn{5}{c|}{}                &     4    \\
\cline{6-6}
\end{array}
$$
Thus $c(\Psi)_u = 4$, so we take $\d_0 = (6, 4)$; using Proposition~\ref{prop: B value when t = 2} we then have $B_{\d_0, 2} = 8 > 4 = \dim {u_\Psi}^G$. Each of the remaining non-trivial unipotent classes has $B_1$ in its closure by Lemma~\ref{lem: root elt class in closure of any non-triv class}.

Now take $\Psi = \langle \alpha_2 \rangle$ of type $B_1$. The $\Psi$-net table is as follows.
$$
\begin{array}{|*8{>{\ss}c|}}
\hline
\multicolumn{4}{|>{\ss}c|}{\Psi\mathrm{-nets}} & & \multicolumn{2}{|>{\ss}c|}{c(s)} & \multicolumn{1}{|>{\ss}c|}{c(u_\Psi)} \\
\cline{1-4} \cline{6-8}
    \bar\nu    & n_0 & n_1 & n_2 & m & r = 2 & r \geq 3 & p \geq 3 \\
\hline
    2\bom_2    &  0  &  1  &  2  & 2 &   2   &     4    &     4    \\
    2\bom_2    &  1  &  2  &  0  & 1 &   2   &     2    &     2    \\
\hline
\multicolumn{5}{c|}{}                &   4   &     6    &     6    \\
\cline{6-8}
\end{array}
$$
Thus if $r = 2$ then $c(\Psi)_{ss} = 4$, so we take $\d_0 = (6, 4)$; using Proposition~\ref{prop: B value when t = 2} we then have $B_{\d_0, 2} = 8 > M_2$. If instead $r \geq 3$ then $c(\Psi)_{ss} = c(\Psi)_u = 6$, so we take $\d_0 = (4, 4, 2)$; using Corollary~\ref{cor: B values for small k} we then have $B_{\d_0, 2} = 10 > M$. Therefore if $k \in [2, \frac{d}{2}]$ the quadruple $(G, \lambda, p, k)$ satisfies $\ssdiamcon$ and $\udiamcon$.
\end{proof}

\begin{prop}\label{prop: B_4, omega_4, k geq 4, nets}
Let $G = B_4$ and $\lambda = \omega_4$; then for $k \in [4, \frac{d}{2}]$ the quadruple $(G, \lambda, p, k)$ satisfies $\ssdiamcon$ and $\udiamcon$.
\end{prop}

\begin{proof}
The weight table is as follows.
$$
\begin{array}{|*4{>{\ss}c|}}
\hline
i & \mu & |W.\mu| & m_\mu \\
\hline
1 &  \omega_4 & 16 & 1 \\
\hline
\end{array}
$$
We have $M = 32$ and $M_2 = 20$; we take $k_0 = 4$.

Take $\Psi = \langle \alpha_1 \rangle$ of type $A_1$, and $\Psi = \langle \alpha_1, \alpha_3 \rangle$ of type ${A_1}^2$. The $\Psi$-net tables are as follows.
$$
\begin{array}{|*4{>{\ss}c|}}
\hline
\multicolumn{2}{|>{\ss}c|}{\Psi\mathrm{-nets}} & & \multicolumn{1}{|>{\ss}c|}{c(u_\Psi)} \\
\cline{1-2} \cline{4-4}
    \bar\nu    & n_1 & m & p \geq 2 \\
\hline
     \bom_1    &  2  & 4 &     4    \\
       0       &  1  & 8 &          \\
\hline
\multicolumn{3}{c|}{}    &     4    \\
\cline{4-4}
\end{array}
\qquad
\begin{array}{|*4{>{\ss}c|}}
\hline
\multicolumn{2}{|>{\ss}c|}{\Psi\mathrm{-nets}} & & \multicolumn{1}{|>{\ss}c|}{c(u_\Psi)} \\
\cline{1-3} \cline{4-4}
        \bar\nu        & n_1 & m & p \geq 2 \\
\hline
    \bom_1 + \bom_3    &  4  & 1 &     2    \\
         \bom_1        &  2  & 2 &     2    \\
         \bom_3        &  2  & 2 &     2    \\
           0           &  1  & 4 &          \\
\hline
\multicolumn{3}{c|}{}            &     6    \\
\cline{4-4}
\end{array}
$$
Thus according as $\Psi = A_1$ or ${A_1}^2$ we have $c(\Psi)_u = 4$ or $6$, so we take $\d_0 = (12, 4)$ or $(10, 6)$; using Proposition~\ref{prop: B value when t = 2} we then have $B_{\d_0, 4} = 16 > 12 = \dim {u_\Psi}^G$ or $B_{\d_0, 4} = 22 > 16 = \dim {u_\Psi}^G$. Each of the remaining non-trivial unipotent classes has $B_1$ in its closure by Lemma~\ref{lem: root elt class in closure of any non-triv class}.

Now take $\Psi = \langle \alpha_4 \rangle$ of type $B_1$. The $\Psi$-net table is as follows.
$$
\begin{array}{|*5{>{\ss}c|}}
\hline
\multicolumn{2}{|>{\ss}c|}{\Psi\mathrm{-nets}} & & \multicolumn{1}{|>{\ss}c|}{c(s)} & \multicolumn{1}{|>{\ss}c|}{c(u_\Psi)} \\
\cline{1-2} \cline{4-5}
    \bar\nu    & n_1 & m & r \geq 2 & p \geq 2 \\
\hline
     \bom_4    &  2  & 8 &     8    &     8    \\
\hline
\multicolumn{3}{c|}{}    &     8    &     8    \\
\cline{4-5}
\end{array}
$$
Thus $c(\Psi)_{ss} = c(\Psi)_u = 8$, so we take $\d_0 = (8, 8)$; using Proposition~\ref{prop: B value when t = 2} we then have $B_{\d_0, 4} = 24 > M_2 > 14 = \dim {u_\Psi}^G$. We may therefore assume from now on that $r \geq 3$, and that $p \geq 3$ when we treat unipotent classes. We need only consider semisimple classes $s^G$ with $|\Phi(s)| \leq M - 24 = 8$, each of which has a subsystem of type $A_2B_1$ or a subsystem of type $B_2$ disjoint from $\Phi(s)$, and unipotent classes of dimension at least $24$, each of which has the class $A_2B_1$ or the class $B_2$ in its closure by Lemma~\ref{lem: various classes in B_ell for fixed ell}(iv) (and for the unipotent class $B_2$ to lie in $G_{(p)}$ we need $p \geq 5$).

Now take $\Psi = \langle \alpha_1, \alpha_2, \alpha_4 \rangle$ of type $A_2B_1$, and $\Psi = \langle \alpha_3, \alpha_4 \rangle$ of type $B_2$. The $\Psi$-net tables are as follows.
$$
\begin{array}{|*5{>{\ss}c|}}
\hline
\multicolumn{2}{|>{\ss}c|}{\Psi\mathrm{-nets}} & & \multicolumn{1}{|>{\ss}c|}{c(s)} & \multicolumn{1}{|>{\ss}c|}{c(u_\Psi)} \\
\cline{1-2} \cline{4-5}
     \bar\nu     & n_1 & m & r \geq 3 & p \geq 3 \\
\hline
 \bom_1 + \bom_4 &  6  & 1 &     4    &     4    \\
 \bom_2 + \bom_4 &  6  & 1 &     4    &     4    \\
      \bom_4     &  2  & 2 &     2    &     2    \\
\hline
\multicolumn{3}{c|}{}      &    10    &    10    \\
\cline{4-5}
\end{array}
\qquad
\begin{array}{|*5{>{\ss}c|}}
\hline
\multicolumn{2}{|>{\ss}c|}{\Psi\mathrm{-nets}} & & \multicolumn{1}{|>{\ss}c|}{c(s)} & \multicolumn{1}{|>{\ss}c|}{c(u_\Psi)} \\
\cline{1-2} \cline{4-5}
    \bar\nu    & n_1 & m & r \geq 3 & p \geq 5 \\
\hline
     \bom_4    &  4  & 4 &    12    &    12    \\
\hline
\multicolumn{3}{c|}{}    &    12    &    12    \\
\cline{4-5}
\end{array}
$$
Thus according as $\Psi = A_2B_1$ or $B_2$ we have $c(\Psi)_{ss} = c(\Psi)_u = 10$ or $12$, so we take $\d_0 = (6, 6, 4)$ or $(4, 4, 4, 4)$. In the former case, according as $\k_0 = (4, 0, 0)$, $(3, 1, 0)$, $(2, 2, 0)$ or $(2, 1, 1)$ we have $B_{\d_0, \k_0} = 40$, $34$, $32$ or $32$, so $B_{\d_0, 4} = 32 = M > 24 = \dim {u_\Psi}^G$; in the latter case, according as $\k_0 = (4, 0, 0, 0)$, $(3, 1, 0, 0)$, $(2, 2, 0, 0)$, $(2, 1, 1, 0)$ or $(1, 1, 1, 1)$ we have $B_{\d_0, \k_0} = 48$, $42$, $40$, $38$ or $36$, so $B_{\d_0, 4} = 36 > M$. Taking the smaller of the two lower bounds, we see that we need only consider regular semisimple classes $s^G$, and the regular unipotent class; since each of the former has a subsystem of type $B_2$ disjoint from $\Phi(s)$, and the latter has the class $B_2$ in its closure by Lemma~\ref{lem: any class in closure of reg class}, we may actually take the larger of the two lower bounds. Therefore if $k \in [4, \frac{d}{2}]$ the quadruple $(G, \lambda, p, k)$ satisfies $\ssdiamcon$ and $\udiamcon$.
\end{proof}

\begin{prop}\label{prop: B_5 and B_6, omega_5 and omega_6, k geq 3 and 2, nets}
Let $G = B_5$ and $\lambda = \omega_5$, or $G = B_6$ and $\lambda = \omega_6$; then for $k \in [3, \frac{d}{2}]$ or $k \in [2, \frac{d}{2}]$ respectively the quadruple $(G, \lambda, p, k)$ satisfies $\ssdiamcon$ and $\udiamcon$.
\end{prop}

\begin{proof}
These follow from Propositions~\ref{prop: D_6, omega_6, k geq 3, nets} and \ref{prop: D_7, omega_7, k geq 2, nets}, since $B_\ell$ is a subgroup of $D_{\ell + 1}$ and the spin module for $B_\ell$ is the restriction of the half-spin module for $D_{\ell + 1}$.
\end{proof}

\begin{prop}\label{prop: C_4, C_5 and C_6, omega_4, omega_5 and omega_6, k geq 4, 3 and 2, nets}
Let $G = C_4$ and $\lambda = \omega_4$, or $G = C_5$ and $\lambda = \omega_5$, or $G = C_6$ and $\lambda = \omega_6$, all with $p = 2$; then for $k \in [4, \frac{d}{2}]$, or $k \in [3, \frac{d}{2}]$ or $k \in [2, \frac{d}{2}]$ respectively the quadruple $(G, \lambda, p, k)$ satisfies $\ssdiamcon$ and $\udiamcon$.
\end{prop}

\begin{proof}
These are immediate consequences of Propositions~\ref{prop: B_4, omega_4, k geq 4, nets} and \ref{prop: B_5 and B_6, omega_5 and omega_6, k geq 3 and 2, nets}, using the exceptional isogeny $B_\ell \to C_\ell$ which exists in characteristic $2$.
\end{proof}

\begin{prop}\label{prop: C_3, omega_2, k geq 3, nets}
Let $G = C_3$ and $\lambda = \omega_2$; then for $k \in [3, \frac{d}{2}]$ the quadruple $(G, \lambda, p, k)$ satisfies $\ssdiamcon$ and $\udiamcon$.
\end{prop}

\begin{proof}
Write $\z = \z_{p, 3}$. The weight table is as follows.
$$
\begin{array}{|*4{>{\ss}c|}}
\hline
i & \mu & |W.\mu| & m_\mu \\
\hline
1 & \omega_2 & 12 &   1    \\
0 &    0     &  1 & 2 - \z \\
\hline
\end{array}
$$
We have $M = 18$, $M_3 = 14$ and $M_2 = 12$; we take $k_0 = 3$.

Take $\Psi = \langle \alpha_3 \rangle$ of type $C_1$. The $\Psi$-net table is as follows.
$$
\begin{array}{|*5{>{\ss}c|}}
\hline
\multicolumn{3}{|>{\ss}c|}{\Psi\mathrm{-nets}} & & \multicolumn{1}{|>{\ss}c|}{c(u_\Psi)} \\
\cline{1-3} \cline{5-5}
    \bar\nu    & n_0 & n_1 & m & p \geq 2 \\
\hline
     \bom_3    &  0  &  2  & 4 &     4    \\
       0       &  0  &  1  & 4 &          \\
       0       &  1  &  0  & 1 &          \\
\hline
\multicolumn{4}{c|}{}          &     4    \\
\cline{5-5}
\end{array}
$$
Thus $c(\Psi)_u = 4$, so we take $\d_0 = (10 - \z, 4)$; using Proposition~\ref{prop: B value when t = 2} we then have $B_{\d_0, 3} = 12 > 6 = \dim {u_\Psi}^G$. Each of the remaining non-trivial unipotent classes has $A_1$ in its closure by Lemma~\ref{lem: root elt class in closure of any non-triv class}.

Now take $\Psi = \langle \alpha_1 \rangle$ of type $A_1$. The $\Psi$-net table is as follows.
$$
\begin{array}{|*8{>{\ss}c|}}
\hline
\multicolumn{3}{|>{\ss}c|}{\Psi\mathrm{-nets}} & & \multicolumn{2}{|>{\ss}c|}{c(s)} & \multicolumn{2}{|>{\ss}c|}{c(u_\Psi)} \\
\cline{1-3} \cline{5-8}
    \bar\nu    & n_0 & n_1 & m &  r = 2  & r \geq 3 & p = 2 & p \geq 3 \\
\hline
    2\bom_1    &  1  &  2  & 1 &  2 - \z &     2    &   1  &   2   \\
     \bom_1    &  0  &  2  & 4 &    4    &     4    &   4  &   4   \\
       0       &  0  &  1  & 2 &         &          &      &       \\
\hline
\multicolumn{4}{c|}{}          &  6 - \z &     6    &   5  &   6   \\
\cline{5-8}
\end{array}
$$
Thus if $r = 2$ then $c(\Psi)_{ss} = 6 - \z$, so we take $\d_0 = (8, 6 - \z)$; using Proposition~\ref{prop: B value when t = 2} we then have $B_{\d_0, 3} = 16 - 2\z > M_2$. If $p = 2$ then $c(\Psi)_u = 5$, so we take $\d_0 = (9, 5)$; using Proposition~\ref{prop: B value when t = 2} we then have $B_{\d_0, 3} = 15 > M_2$. If instead $r \geq 3$ and $p \geq 3$ then $c(\Psi)_{ss} = c(\Psi)_u = 6$, so we take $\d_0 = (8 - \z, 6)$; using Proposition~\ref{prop: B value when t = 2} we then have $B_{\d_0, 3} = 16 - \z > M_3 > 10 = \dim {u_\Psi}^G$. We may therefore assume from now on that $r \geq 5$, and that $p \geq 5$ when we treat unipotent classes. We need only consider semisimple classes $s^G$ with $|\Phi(s)| \leq M - (16 - \z) = 2 + \z$, each of which has a subsystem of type $A_2$ disjoint from $\Phi(s)$, and unipotent classes of dimension at least $16$, each of which has the class $A_2$ in its closure by Lemma~\ref{lem: various classes in classical groups by dim}(viii).

Now take $\Psi = \langle \alpha_1, \alpha_2 \rangle$ of type $A_2$. The $\Psi$-net table is as follows.
$$
\begin{array}{|*6{>{\ss}c|}}
\hline
\multicolumn{3}{|>{\ss}c|}{\Psi\mathrm{-nets}} & & \multicolumn{1}{|>{\ss}c|}{c(s)} & \multicolumn{1}{|>{\ss}c|}{c(u_\Psi)} \\
\cline{1-3} \cline{5-6}
     \bar\nu     & n_0 & n_1 & m & r \geq 5 & p \geq 5 \\
\hline
 \bom_1 + \bom_2 &  1  &  6  & 1 &  5 - \z  &     6    \\
      \bom_1     &  0  &  3  & 1 &     2    &     2    \\
      \bom_2     &  0  &  3  & 1 &     2    &     2    \\
\hline
\multicolumn{4}{c|}{}            &  9 - \z  &    10    \\
\cline{5-6}
\end{array}
$$
Thus $c(\Psi)_{ss}, c(\Psi)_u \geq 9 - \z$, so we may take $\d_0 = (5, 5, 4 - \z)$; using Corollary~\ref{cor: B values for small k} we then have $B_{\d_0, 3} = 22 - 2\z > M$. Therefore if $k \in [3, \frac{d}{2}]$ the quadruple $(G, \lambda, p, k)$ satisfies $\ssdiamcon$ and $\udiamcon$.
\end{proof}

\begin{prop}\label{prop: B_3, omega_2, k geq 3, nets}
Let $G = B_3$ and $\lambda = \omega_2$ with $p = 2$; then for $k \in [3, \frac{d}{2}]$ the quadruple $(G, \lambda, p, k)$ satisfies $\ssdiamcon$ and $\udiamcon$.
\end{prop}

\begin{proof}
This is an immediate consequence of Proposition~\ref{prop: C_3, omega_2, k geq 3, nets}, using the exceptional isogeny $B_\ell \to C_\ell$ which exists in characteristic $2$.
\end{proof}

\begin{prop}\label{prop: C_4, omega_3, k geq 2, nets}
Let $G = C_4$ and $\lambda = \omega_3$ with $p = 3$; then for $k \in [2, \frac{d}{2}]$ the quadruple $(G, \lambda, p, k)$ satisfies $\ssdiamcon$ and $\udiamcon$.
\end{prop}

\begin{proof}
The weight table is as follows.
$$
\begin{array}{|*4{>{\ss}c|}}
\hline
i & \mu & |W.\mu| & m_\mu \\
\hline
2 & \omega_3 & 32 & 1 \\
1 & \omega_1 &  8 & 1 \\
\hline
\end{array}
$$
We have $M = 32$ and $M_2 = 20$; we take $k_0 = 2$.

Take $\Psi = \langle \alpha_4 \rangle$ of type $C_1$. The $\Psi$-net table is as follows.
$$
\begin{array}{|*5{>{\ss}c|}}
\hline
\multicolumn{3}{|>{\ss}c|}{\Psi\mathrm{-nets}} & & \multicolumn{1}{|>{\ss}c|}{c(u_\Psi)} \\
\cline{1-3} \cline{5-5}
    \bar\nu    & n_1 & n_2 &  m & p = 3 \\
\hline
     \bom_4    &  0  &  2  & 12 &  12   \\
     \bom_4    &  2  &  0  &  1 &   1   \\
       0       &  0  &  1  &  8 &       \\
       0       &  1  &  0  &  6 &       \\
\hline
\multicolumn{4}{c|}{}           &  13   \\
\cline{5-5}
\end{array}
$$
Thus $c(\Psi)_u = 13$, so we take $\d_0 = (27, 13)$; using Proposition~\ref{prop: B value when t = 2} we then have $B_{\d_0, 2} = 26 > 8 = \dim {u_\Psi}^G$. Each of the remaining non-trivial unipotent classes has $A_1$ in its closure by Lemma~\ref{lem: root elt class in closure of any non-triv class}.

Now take $\Psi = \langle \alpha_1 \rangle$ of type $A_1$. The $\Psi$-net table is as follows.
$$
\begin{array}{|*7{>{\ss}c|}}
\hline
\multicolumn{3}{|>{\ss}c|}{\Psi\mathrm{-nets}} & & \multicolumn{2}{|>{\ss}c|}{c(s)} & \multicolumn{1}{|>{\ss}c|}{c(u_\Psi)} \\
\cline{1-3} \cline{5-7}
    \bar\nu    & n_1 & n_2 & m & r = 2 & r \geq 5 & p = 3 \\
\hline
    2\bom_1    &  1  &  2  & 4 &   4   &     8    &   8   \\
     \bom_1    &  0  &  2  & 8 &   8   &     8    &   8   \\
     \bom_1    &  2  &  0  & 2 &   2   &     2    &   2   \\
       0       &  0  &  1  & 8 &       &          &       \\
\hline
\multicolumn{4}{c|}{}          &  14   &    18    &  18   \\
\cline{5-7}
\end{array}
$$
Thus if $r = 2$ then $c(\Psi)_{ss} = 14$, so we take $\d_0 = (26, 14)$; using Proposition~\ref{prop: B value when t = 2} we then have $B_{\d_0, 2} = 28 > M_2$. If instead $r \geq 5$ then $c(\Psi)_{ss} = c(\Psi)_u = 18$, so we take $\d_0 = (22, 18)$; using Proposition~\ref{prop: B value when t = 2} we then have $B_{\d_0, 2} = 36 > M$. Therefore if $k \in [2, \frac{d}{2}]$ the quadruple $(G, \lambda, p, k)$ satisfies $\ssdiamcon$ and $\udiamcon$.
\end{proof}

\begin{prop}\label{prop: C_3, omega_3, k geq 3, nets}
Let $G = C_3$ and $\lambda = \omega_3$ with $p \geq 3$; then for $k \in [3, \frac{d}{2}]$ the quadruple $(G, \lambda, p, k)$ satisfies $\ssdiamcon$ and $\udiamcon$.
\end{prop}

\begin{proof}
The weight table is as follows.
$$
\begin{array}{|*4{>{\ss}c|}}
\hline
i & \mu & |W.\mu| & m_\mu \\
\hline
2 & \omega_3 & 8 & 1 \\
1 & \omega_1 & 6 & 1 \\
\hline
\end{array}
$$
We have $M = 18$ and $M_2 = 12$; we take $k_0 = 3$.

Take $\Psi = \langle \alpha_3 \rangle$ of type $C_1$. The $\Psi$-net table is as follows.
$$
\begin{array}{|*5{>{\ss}c|}}
\hline
\multicolumn{3}{|>{\ss}c|}{\Psi\mathrm{-nets}} & & \multicolumn{1}{|>{\ss}c|}{c(u_\Psi)} \\
\cline{1-3} \cline{5-5}
    \bar\nu    & n_1 & n_2 & m & p \geq 3 \\
\hline
     \bom_3    &  0  &  2  & 4 &     4    \\
     \bom_3    &  2  &  0  & 1 &     1    \\
       0       &  1  &  0  & 4 &          \\
\hline
\multicolumn{4}{c|}{}          &     5    \\
\cline{5-5}
\end{array}
$$
Thus $c(\Psi)_u = 5$, so we take $\d_0 = (9, 5)$; using Proposition~\ref{prop: B value when t = 2} we then have $B_{\d_0, 3} = 15 > 6 = \dim {u_\Psi}^G$. Each of the remaining non-trivial unipotent classes has $A_1$ in its closure by Lemma~\ref{lem: root elt class in closure of any non-triv class}.

Now take $\Psi = \langle \alpha_1 \rangle$ of type $A_1$. The $\Psi$-net table is as follows.
$$
\begin{array}{|*7{>{\ss}c|}}
\hline
\multicolumn{3}{|>{\ss}c|}{\Psi\mathrm{-nets}} & & \multicolumn{2}{|>{\ss}c|}{c(s)} & \multicolumn{1}{|>{\ss}c|}{c(u_\Psi)} \\
\cline{1-3} \cline{5-7}
    \bar\nu    & n_1 & n_2 & m & r = 2 & r \geq 3 & p \geq 3 \\
\hline
    2\bom_1    &  1  &  2  & 2 &   2   &     4    &     4    \\
     \bom_1    &  2  &  0  & 2 &   2   &     2    &     2    \\
       0       &  0  &  1  & 4 &       &          &          \\
\hline
\multicolumn{4}{c|}{}          &   4   &     6    &     6    \\
\cline{5-7}
\end{array}
$$
Thus $c(\Psi)_{ss}, c(\Psi)_u \geq 4$, so we may take $\d_0 = (10, 4)$; using Proposition~\ref{prop: B value when t = 2} we then have $B_{\d_0, 3} = 12 > 10 = \dim {u_\Psi}^G$. We need only consider semisimple classes $s^G$ with $|\Phi(s)| \leq M - 12 = 6$, each of which has a subsystem of type $A_1C_1$ disjoint from $\Phi(s)$, and unipotent classes of dimension at least $12$, each of which has the class $A_1C_1$ in its closure by Lemma~\ref{lem: various classes in C_ell for fixed ell}(v).

Now take $\Psi = \langle \alpha_1, \alpha_3 \rangle$ of type $A_1C_1$. The $\Psi$-net table is as follows.
$$
\begin{array}{|*7{>{\ss}c|}}
\hline
\multicolumn{3}{|>{\ss}c|}{\Psi\mathrm{-nets}} & & \multicolumn{2}{|>{\ss}c|}{c(s)} & \multicolumn{1}{|>{\ss}c|}{c(u_\Psi)} \\
\cline{1-3} \cline{5-7}
     \bar\nu     & n_1 & n_2 & m & r = 2 & r \geq 3 & p \geq 3 \\
\hline
 2\bom_1 + \bom_3 &  2  &  4  & 1 &   3   &     4    &     4    \\
       \bom_1     &  2  &  0  & 2 &   2   &     2    &     2    \\
       \bom_3     &  0  &  2  & 2 &   2   &     2    &     2    \\
\hline
\multicolumn{4}{c|}{}             &   7   &     8    &     8    \\
\cline{5-7}
\end{array}
$$
Thus if $r = 2$ then $c(\Psi)_{ss} = 7$, so we take $\d_0 = (7, 7)$; using Proposition~\ref{prop: B value when t = 2} we then have $B_{\d_0, 3} = 17 > M_2$. If instead $r \geq 3$ then $c(\Psi)_{ss} = c(\Psi)_u = 8$, so we take $\d_0 = (6, 6, 2)$; using Corollary~\ref{cor: B values for small k} we then have $B_{\d_0, 3} = 20 > M$. Therefore if $k \in [3, \frac{d}{2}]$ the quadruple $(G, \lambda, p, k)$ satisfies $\ssdiamcon$ and $\udiamcon$.
\end{proof}

\begin{prop}\label{prop: C_4, omega_4, k geq 2, nets}
Let $G = C_4$ and $\lambda = \omega_4$ with $p \geq 3$; then for $k \in [2, \frac{d}{2}]$ the quadruple $(G, \lambda, p, k)$ satisfies $\ssdiamcon$ and $\udiamcon$.
\end{prop}

\begin{proof}
Write $\z = \z_{p, 3}$. The weight table is as follows.
$$
\begin{array}{|*4{>{\ss}c|}}
\hline
i & \mu & |W.\mu| & m_\mu \\
\hline
2 & \omega_4 & 16 &   1    \\
1 & \omega_2 & 24 &   1    \\
0 &    0     &  1 & 2 - \z \\
\hline
\end{array}
$$
We have $M = 32$ and $M_2 = 20$; we take $k_0 = 2$.

Take $\Psi = \langle \alpha_4 \rangle$ of type $C_1$. The $\Psi$-net table is as follows.
$$
\begin{array}{|*6{>{\ss}c|}}
\hline
\multicolumn{4}{|>{\ss}c|}{\Psi\mathrm{-nets}} & & \multicolumn{1}{|>{\ss}c|}{c(u_\Psi)} \\
\cline{1-4} \cline{6-6}
    \bar\nu    & n_0 & n_1 & n_2 &  m & p \geq 3 \\
\hline
     \bom_1    &  0  &  0  &  2  &  8 &     8    \\
     \bom_1    &  0  &  2  &  0  &  6 &     6    \\
       0       &  0  &  1  &  0  & 12 &          \\
       0       &  1  &  0  &  0  &  1 &          \\
\hline
\multicolumn{5}{c|}{}                 &    14    \\
\cline{6-6}
\end{array}
$$
Thus $c(\Psi)_u = 14$, so we take $\d_0 = (28 - \z, 14)$; using Proposition~\ref{prop: B value when t = 2} we then have $B_{\d_0, 2} = 28 > 8 = \dim {u_\Psi}^G$. Each of the remaining non-trivial unipotent classes has $A_1$ in its closure by Lemma~\ref{lem: root elt class in closure of any non-triv class}.

Now take $\Psi = \langle \alpha_1 \rangle$ of type $A_1$. The $\Psi$-net table is as follows.
$$
\begin{array}{|*8{>{\ss}c|}}
\hline
\multicolumn{4}{|>{\ss}c|}{\Psi\mathrm{-nets}} & & \multicolumn{2}{|>{\ss}c|}{c(s)} & \multicolumn{1}{|>{\ss}c|}{c(u_\Psi)} \\
\cline{1-4} \cline{6-8}
    \bar\nu    & n_0 & n_1 & n_2 & m &  r = 2  & r \geq 3 & p \geq 3 \\
\hline
    2\bom_1    &  0  &  1  &  2  & 4 &    4    &     8    &     8    \\
    2\bom_1    &  1  &  2  &  0  & 1 &  2 - \z &     2    &     2    \\
     \bom_1    &  0  &  2  &  0  & 8 &    8    &     8    &     8    \\
     \bom_1    &  0  &  0  &  1  & 8 &         &          &          \\
       0       &  0  &  1  &  0  & 2 &         &          &          \\
\hline
\multicolumn{5}{c|}{}                & 14 - \z &    18    &    18    \\
\cline{6-8}
\end{array}
$$
Thus if $r = 2$ then $c(\Psi)_{ss} = 14 - \z$, so we take $\d_0 = (28, 14 - \z)$; using Proposition~\ref{prop: B value when t = 2} we then have $B_{\d_0, 2} = 28 - 2\z > M_2$. If instead $r \geq 5$ then $c(\Psi)_{ss} = c(\Psi)_u = 18$, so we take $\d_0 = (24 - \z, 18)$; using Corollary~\ref{cor: B values for small k} we then have $B_{\d_0, 2} = 36 > M$. Therefore if $k \in [2, \frac{d}{2}]$ the quadruple $(G, \lambda, p, k)$ satisfies $\ssdiamcon$ and $\udiamcon$.
\end{proof}

\begin{prop}\label{prop: F_4, omega_1, p = 2, k geq 3, nets}
Let $G = F_4$ and $\lambda = \omega_1$ with $p = 2$; then for $k \in [3, \frac{d}{2}]$ the quadruple $(G, \lambda, p, k)$ satisfies $\ssdiamcon$ and $\udiamcon$.
\end{prop}

\begin{proof}
The weight table is as follows.
$$
\begin{array}{|*4{>{\ss}c|}}
\hline
i & \mu & |W.\mu| & m_\mu \\
\hline
1 & \omega_1 & 24 & 1 \\
0 &    0     &  1 & 2 \\
\hline
\end{array}
$$
We have $M = 48$, $M_3 = 36$ and $M_2 = 28$; we take $k_0 = 3$.

Take $\Psi = \langle \alpha_4 \rangle$ of type $\tilde A_1$. The $\Psi$-net table is as follows.
$$
\begin{array}{|*5{>{\ss}c|}}
\hline
\multicolumn{3}{|>{\ss}c|}{\Psi\mathrm{-nets}} & & \multicolumn{1}{|>{\ss}c|}{c(u_\Psi)} \\
\cline{1-3} \cline{5-5}
    \bar\nu    & n_0 & n_1 &  m & p = 2 \\
\hline
    2\bom_4    &  0  &  2  &  6 &   6   \\
       0       &  0  &  1  & 12 &       \\
       0       &  1  &  0  &  1 &       \\
\hline
\multicolumn{4}{c|}{}           &   6   \\
\cline{5-5}
\end{array}
$$
Thus $c(\Psi)_u = 6$, so we take $\d_0 = (20, 6)$; using Proposition~\ref{prop: B value when t = 2} we then have $B_{\d_0, 3} = 18 > 16 = \dim {u_\Psi}^G$. Each of the remaining non-trivial unipotent classes has $A_1$ in its closure by Lemma~\ref{lem: root elt class in closure of any non-triv class}.

Now take $\Psi = \langle \alpha_1 \rangle$ of type $A_1$. The $\Psi$-net table is as follows.
$$
\begin{array}{|*6{>{\ss}c|}}
\hline
\multicolumn{3}{|>{\ss}c|}{\Psi\mathrm{-nets}} & & \multicolumn{1}{|>{\ss}c|}{c(s)} & \multicolumn{1}{|>{\ss}c|}{c(u_\Psi)} \\
\cline{1-3} \cline{5-6}
    \bar\nu    & n_0 & n_1 & m & r \geq 3 & p = 2 \\
\hline
    2\bom_1    &  1  &  2  & 1 &     2    &   1   \\
     \bom_1    &  0  &  2  & 8 &     8    &   8   \\
       0       &  0  &  1  & 6 &          &       \\
\hline
\multicolumn{4}{c|}{}          &    10    &   9   \\
\cline{5-6}
\end{array}
$$
Thus $c(\Psi)_{ss}, c(\Psi)_u \geq 9$, so we may take $\d_0 = (17, 9)$; using Proposition~\ref{prop: B value when t = 2} we then have $B_{\d_0, 3} = 27 > 16 = \dim {u_\Psi}^G$. We need only consider semisimple classes $s^G$ with $|\Phi(s)| \leq M - 27 = 21$, each of which has a subsystem of type $A_1\tilde A_1$ disjoint from $\Phi(s)$, and unipotent classes of dimension at least $27$, each of which has the class $A_1\tilde A_1$ in its closure by Lemma~\ref{lem: various classes in F_4}(i).

Now take $\Psi = \langle \alpha_1, \alpha_4 \rangle$ of type $A_1\tilde A_1$. The $\Psi$-net table is as follows.
$$
\begin{array}{|*6{>{\ss}c|}}
\hline
\multicolumn{3}{|>{\ss}c|}{\Psi\mathrm{-nets}} & & \multicolumn{1}{|>{\ss}c|}{c(s)} & \multicolumn{1}{|>{\ss}c|}{c(u_\Psi)} \\
\cline{1-3} \cline{5-6}
      \bar\nu      & n_0 & n_1 & m & r \geq 3 & p = 2 \\
\hline
  \bom_1 + 2\bom_4 &  0  &  4  & 2 &     4    &   4   \\
      2\bom_1      &  1  &  2  & 1 &     2    &   1   \\
       \bom_1      &  0  &  2  & 4 &     4    &   4   \\
      2\bom_4      &  0  &  2  & 2 &     2    &   2   \\
         0         &  0  &  1  & 2 &          &       \\
\hline
\multicolumn{4}{c|}{}              &    12    &  11   \\
\cline{5-6}
\end{array}
$$
Thus $c(\Psi)_{ss}, c(\Psi)_u \geq 11$, so we may take $\d_0 = (15, 11)$; using Proposition~\ref{prop: B value when t = 2} we then have $B_{\d_0, 3} = 33 > M_2 = 28 = \dim {u_\Psi}^G$. Therefore if $k \in [3, \frac{d}{2}]$ the quadruple $(G, \lambda, p, k)$ satisfies $\udiamcon$. We need only consider semisimple classes $s^G$ with $|\Phi(s)| \leq M - 33 = 15$, each of which has a subsystem of type $A_2$ disjoint from $\Phi(s)$.

Now take $\Psi = \langle \alpha_1, \alpha_2 \rangle$ of type $A_2$. The $\Psi$-net table is as follows.
$$
\begin{array}{|*5{>{\ss}c|}}
\hline
\multicolumn{3}{|>{\ss}c|}{\Psi\mathrm{-nets}} & & \multicolumn{1}{|>{\ss}c|}{c(s)} \\
\cline{1-3} \cline{5-5}
      \bar\nu      & n_0 & n_1 & m & r \geq 3 \\
\hline
  \bom_1 + \bom_2  &  1  &  6  & 1 &     5    \\
       \bom_1      &  0  &  3  & 3 &     6    \\
       \bom_2      &  0  &  3  & 3 &     6    \\
\hline
\multicolumn{4}{c|}{}              &    17    \\
\cline{5-5}
\end{array}
$$
Thus $c(\Psi)_{ss} = 17$, so we take $\d_0 = (9, 9, 8)$; using Corollary~\ref{cor: B values for small k} we then have $B_{\d_0, 3} = 46 > M_3$. We may therefore assume from now on that $r \geq 5$. We need only consider semisimple classes $s^G$ with $|\Phi(s)| \leq M - 46 = 2$, each of which has a subsystem of type $B_3$ disjoint from $\Phi(s)$.

Now take $\Psi = \langle \alpha_1, \alpha_2, \alpha_3 \rangle$ of type $B_3$. The $\Psi$-net table is as follows.
$$
\begin{array}{|*5{>{\ss}c|}}
\hline
\multicolumn{3}{|>{\ss}c|}{\Psi\mathrm{-nets}} & & \multicolumn{1}{|>{\ss}c|}{c(s)} \\
\cline{1-3} \cline{5-5}
      \bar\nu      & n_0 & n_1 & m & r \geq 5 \\
\hline
       \bom_1      &  0  &  6  & 2 &    10    \\
       \bom_2      &  1  & 12  & 1 &    11    \\
\hline
\multicolumn{4}{c|}{}              &    21    \\
\cline{5-5}
\end{array}
$$
Thus $c(\Psi)_{ss} = 21$, so we take $\d_0 = (5, 5, 5, 5, 5, 1)$; using Corollary~\ref{cor: B values for small k} we then have $B_{\d_0, 3} = 57 > M$. Therefore if $k \in [3, \frac{d}{2}]$ the quadruple $(G, \lambda, p, k)$ satisfies $\ssdiamcon$.
\end{proof}

\begin{prop}\label{prop: F_4, omega_1, p geq 3, k geq 2, nets}
Let $G = F_4$ and $\lambda = \omega_1$ with $p \geq 3$; then for $k \in [2, \frac{d}{2}]$ the quadruple $(G, \lambda, p, k)$ satisfies $\ssdiamcon$ and $\udiamcon$.
\end{prop}

\begin{proof}
The weight table is as follows.
$$
\begin{array}{|*4{>{\ss}c|}}
\hline
i & \mu & |W.\mu| & m_\mu \\
\hline
2 & \omega_1 & 24 & 1 \\
1 & \omega_4 & 24 & 1 \\
0 &    0     &  1 & 4 \\
\hline
\end{array}
$$
We have $M = 48$, $M_3 = 36$ and $M_2 = 28$; we take $k_0 = 2$.

Take $\Psi = \langle \alpha_1 \rangle$ of type $A_1$. The $\Psi$-net table is as follows.
$$
\begin{array}{|*6{>{\ss}c|}}
\hline
\multicolumn{4}{|>{\ss}c|}{\Psi\mathrm{-nets}} & & \multicolumn{1}{|>{\ss}c|}{c(u_\Psi)} \\
\cline{1-4} \cline{6-6}
    \bar\nu    & n_0 & n_1 & n_2 &  m & p \geq 3 \\
\hline
    2\bom_1    &  1  &  0  &  2  &  1 &     2    \\
     \bom_1    &  0  &  0  &  2  &  8 &     8    \\
     \bom_1    &  0  &  2  &  0  &  6 &     6    \\
       0       &  0  &  0  &  1  &  6 &          \\
       0       &  0  &  1  &  0  & 12 &          \\
\hline
\multicolumn{5}{c|}{}                 &    16    \\
\cline{6-6}
\end{array}
$$
Thus $c(\Psi)_u = 16$, so we take $\d_0 = (36, 16)$; using Proposition~\ref{prop: B value when t = 2} we then have $B_{\d_0, 2} = 32 > 16 = \dim {u_\Psi}^G$. Each of the remaining non-trivial unipotent classes has $\tilde A_1$ in its closure by Lemma~\ref{lem: root elt class in closure of any non-triv class}.

Now take $\Psi = \langle \alpha_4 \rangle$ of type $\tilde A_1$. The $\Psi$-net table is as follows.
$$
\begin{array}{|*8{>{\ss}c|}}
\hline
\multicolumn{4}{|>{\ss}c|}{\Psi\mathrm{-nets}} & & \multicolumn{2}{|>{\ss}c|}{c(s)} & \multicolumn{1}{|>{\ss}c|}{c(u_\Psi)} \\
\cline{1-4} \cline{6-8}
    \bar\nu    & n_0 & n_1 & n_2 &  m & r = 2 & r \geq 3 & p \geq 3 \\
\hline
    2\bom_4    &  0  &  1  &  2  &  6 &   6   &    12    &    12    \\
    2\bom_4    &  1  &  2  &  0  &  1 &   2   &     2    &     2    \\
     \bom_4    &  0  &  2  &  0  &  8 &   8   &     8    &     8    \\
       0       &  0  &  0  &  1  & 12 &       &          &          \\
\hline
\multicolumn{5}{c|}{}                 &  16   &    22    &    22    \\
\cline{6-8}
\end{array}
$$
Thus if $r = 2$ then $c(\Psi)_{ss} = 16$, so we take $\d_0 = (36, 16)$; using Proposition~\ref{prop: B value when t = 2} we then have $B_{\d_0, 2} = 32 > M_2$. If instead $r \geq 3$ then $c(\Psi)_{ss} = c(\Psi)_u = 22$, so we take $\d_0 = (30, 22)$; using Proposition~\ref{prop: B value when t = 2} we then have $B_{\d_0, 2} = 44 > M_3 > 22 = \dim {u_\Psi}^G$. We may therefore assume from now on that $r \geq 5$, and that $p \geq 5$ when we treat unipotent classes. We need only consider semisimple classes $s^G$ with $|\Phi(s)| \leq M - 44 = 4$, each of which has a subsystem of type $A_2$ disjoint from $\Phi(s)$, and unipotent classes of dimension at least $44$, each of which has the class $A_2$ in its closure by Lemma~\ref{lem: various classes in F_4}(ii).

Now take $\Psi = \langle \alpha_1, \alpha_2 \rangle$ of type $A_2$. The $\Psi$-net table is as follows.
$$
\begin{array}{|*7{>{\ss}c|}}
\hline
\multicolumn{4}{|>{\ss}c|}{\Psi\mathrm{-nets}} & & \multicolumn{1}{|>{\ss}c|}{c(s)} & \multicolumn{1}{|>{\ss}c|}{c(u_\Psi)} \\
\cline{1-4} \cline{6-7}
     \bar\nu     & n_0 & n_1 & n_2 & m & r \geq 5 & p \geq 5 \\
\hline
 \bom_1 + \bom_2 &  1  &  0  &  6  & 1 &     6    &     6    \\
      \bom_1     &  0  &  0  &  3  & 3 &     6    &     6    \\
      \bom_1     &  0  &  3  &  0  & 3 &     6    &     6    \\
      \bom_2     &  0  &  0  &  3  & 3 &     6    &     6    \\
      \bom_2     &  0  &  3  &  0  & 3 &     6    &     6    \\
        0        &  0  &  1  &  0  & 6 &          &          \\
\hline
\multicolumn{5}{c|}{}                  &    30    &    30    \\
\cline{6-7}
\end{array}
$$
Thus $c(\Psi)_{ss} = c(\Psi)_u = 30$, so we take $\d_0 = (22, 22, 8)$; using Corollary~\ref{cor: B values for small k} we then have $B_{\d_0, 2} = 58 > M$. Therefore if $k \in [2, \frac{d}{2}]$ the quadruple $(G, \lambda, p, k)$ satisfies $\ssdiamcon$ and $\udiamcon$.
\end{proof}

\begin{prop}\label{prop: F_4, omega_4, k geq 3, nets}
Let $G = F_4$ and $\lambda = \omega_4$; then for $k \in [3, \frac{d}{2}]$ the quadruple $(G, \lambda, p, k)$ satisfies $\ssdiamcon$ and $\udiamcon$.
\end{prop}

\begin{proof}
Write $\z = \z_{p, 3}$. The weight table is as follows.
$$
\begin{array}{|*4{>{\ss}c|}}
\hline
i & \mu & |W.\mu| & m_\mu \\
\hline
1 & \omega_4 & 24 &   1    \\
0 &    0     &  1 & 2 - \z \\
\hline
\end{array}
$$
We have $M = 48$, $M_5 = 40$, $M_3 = 36$ and $M_2 = 28$; we take $k_0 = 3$.

Take $\Psi = \langle \alpha_1 \rangle$ of type $A_1$. The $\Psi$-net table is as follows.
$$
\begin{array}{|*5{>{\ss}c|}}
\hline
\multicolumn{3}{|>{\ss}c|}{\Psi\mathrm{-nets}} & & \multicolumn{1}{|>{\ss}c|}{c(u_\Psi)} \\
\cline{1-3} \cline{5-5}
    \bar\nu    & n_0 & n_1 &  m & p = 2 \\
\hline
     \bom_1    &  0  &  2  &  6 &   6   \\
       0       &  0  &  1  & 12 &       \\
       0       &  1  &  0  &  1 &       \\
\hline
\multicolumn{4}{c|}{}           &   6   \\
\cline{5-5}
\end{array}
$$
Thus $c(\Psi)_u = 6$, so we take $\d_0 = (20 - \z, 6)$; using Proposition~\ref{prop: B value when t = 2} we then have $B_{\d_0, 3} = 18 > 16 = \dim {u_\Psi}^G$. Each of the remaining non-trivial unipotent classes has $\tilde A_1$ in its closure by Lemma~\ref{lem: root elt class in closure of any non-triv class}.

Now take $\Psi = \langle \alpha_4 \rangle$ of type $\tilde A_1$. The $\Psi$-net table is as follows.
$$
\begin{array}{|*8{>{\ss}c|}}
\hline
\multicolumn{3}{|>{\ss}c|}{\Psi\mathrm{-nets}} & & \multicolumn{2}{|>{\ss}c|}{c(s)} & \multicolumn{2}{|>{\ss}c|}{c(u_\Psi)} \\
\cline{1-3} \cline{5-8}
    \bar\nu    & n_0 & n_1 & m &  r = 2  & r \geq 3 & p = 2 & p \geq 3 \\
\hline
    2\bom_4    &  1  &  2  & 1 &  2 - \z &     2    &   1   &     2    \\
     \bom_4    &  0  &  2  & 8 &    8    &     8    &   8   &     8    \\
       0       &  0  &  1  & 6 &         &          &       &          \\
\hline
\multicolumn{4}{c|}{}          & 10 - \z &    10    &   9   &    10    \\
\cline{5-8}
\end{array}
$$
Thus $c(\Psi)_{ss}, c(\Psi)_u \geq 9$, so we may take $\d_0 = (17 - \z, 9)$; using Proposition~\ref{prop: B value when t = 2} we then have $B_{\d_0, 3} = 27 > 22 - 6\delta_{p, 2} = \dim {u_\Psi}^G$. We need only consider semisimple classes $s^G$ with $|\Phi(s)| \leq M - 27 = 21$, each of which has a subsystem of type $A_1\tilde A_1$ disjoint from $\Phi(s)$, and unipotent classes of dimension at least $27$, each of which has the class $A_1\tilde A_1$ in its closure by Lemma~\ref{lem: various classes in F_4}(i).

Now take $\Psi = \langle \alpha_1, \alpha_4 \rangle$ of type $A_1\tilde A_1$. The $\Psi$-net table is as follows.
$$
\begin{array}{|*8{>{\ss}c|}}
\hline
\multicolumn{3}{|>{\ss}c|}{\Psi\mathrm{-nets}} & & \multicolumn{2}{|>{\ss}c|}{c(s)} & \multicolumn{2}{|>{\ss}c|}{c(u_\Psi)} \\
\cline{1-3} \cline{5-8}
      \bar\nu      & n_0 & n_1 & m &  r = 2  & r \geq 3 & p = 2 & p \geq 3 \\
\hline
  \bom_1 + \bom_4  &  0  &  4  & 2 &    4    &     4    &   4   &     4    \\
       \bom_1      &  0  &  2  & 2 &    2    &     2    &   2   &     2    \\
      2\bom_4      &  1  &  2  & 1 &  2 - \z &     2    &   1   &     2    \\
       \bom_4      &  0  &  2  & 4 &    4    &     4    &   4   &     4    \\
         0         &  0  &  1  & 2 &         &          &       &          \\
\hline
\multicolumn{4}{c|}{}              & 12 - \z &    12    &  11   &    12    \\
\cline{5-8}
\end{array}
$$
Thus $c(\Psi)_{ss}, c(\Psi)_u \geq 11$, so we may take $\d_0 = (15 - \z, 11)$; using Proposition~\ref{prop: B value when t = 2} we then have $B_{\d_0, 3} = 33 - \z > M_2 = 28 = \dim {u_\Psi}^G$. We may therefore assume from now on that $r \geq 3$, and that $p \geq 3$ when we treat unipotent classes. We need only consider semisimple classes $s^G$ with $|\Phi(s)| \leq M - (33 - \z) = 15 + \z$, each of which has a subsystem of type $A_2\tilde A_1$ disjoint from $\Phi(s)$, and unipotent classes of dimension at least $33 - \z$, each of which has the class $A_2\tilde A_1$ in its closure by Lemma~\ref{lem: various classes in F_4}(iii).

Now take $\Psi = \langle \alpha_1, \alpha_2, \alpha_4 \rangle$ of type $A_2\tilde A_1$. The $\Psi$-net table is as follows.
$$
\begin{array}{|*6{>{\ss}c|}}
\hline
\multicolumn{3}{|>{\ss}c|}{\Psi\mathrm{-nets}} & & \multicolumn{1}{|>{\ss}c|}{c(s)} & \multicolumn{1}{|>{\ss}c|}{c(u_\Psi)} \\
\cline{1-3} \cline{5-6}
      \bar\nu      & n_0 & n_1 & m & r \geq 3 & p \geq 3 \\
\hline
  \bom_1 + \bom_4  &  0  &  6  & 1 &     4    &     4    \\
  \bom_2 + \bom_4  &  0  &  6  & 1 &     4    &     4    \\
       \bom_1      &  0  &  3  & 1 &     2    &     2    \\
       \bom_2      &  0  &  3  & 1 &     2    &     2    \\
      2\bom_4      &  1  &  2  & 1 &     2    &     2    \\
       \bom_4      &  0  &  2  & 2 &     2    &     2    \\
\hline
\multicolumn{4}{c|}{}              &    16    &    16    \\
\cline{5-6}
\end{array}
$$
Thus $c(\Psi)_{ss} = c(\Psi)_u = 16$, so we take $\d_0 = (10 - \z, 10 - \z, 6 + \z)$; using Corollary~\ref{cor: B values for small k} we then have $B_{\d_0, 3} = 44 > M_5$. We may therefore assume from now on that $r \geq 7$, and that $p \geq 7$ when we treat unipotent classes. We need only consider semisimple classes $s^G$ with $|\Phi(s)| \leq M - 44 = 4$, each of which has a subsystem of type $C_3$ disjoint from $\Phi(s)$, and unipotent classes of dimension at least $44$, each of which has the class $C_3$ in its closure by Lemma~\ref{lem: various classes in F_4}(iv).

Now take $\Psi = \langle \alpha_2, \alpha_3, \alpha_4 \rangle$ of type $C_3$. The $\Psi$-net table is as follows.
$$
\begin{array}{|*6{>{\ss}c|}}
\hline
\multicolumn{3}{|>{\ss}c|}{\Psi\mathrm{-nets}} & & \multicolumn{1}{|>{\ss}c|}{c(s)} & \multicolumn{1}{|>{\ss}c|}{c(u_\Psi)} \\
\cline{1-3} \cline{5-6}
      \bar\nu      & n_0 & n_1 & m & r \geq 7 & p \geq 7 \\
\hline
       \bom_3      &  1  & 12  & 1 & 11 - \z  &     12   \\
       \bom_4      &  0  &  6  & 2 &    10    &     10   \\
\hline
\multicolumn{4}{c|}{}              & 21 - \z  &     22   \\
\cline{5-6}
\end{array}
$$
Thus $c(\Psi)_u = 22$, so we may take $\d_0 = (4, 4, 4, 4, 4, 4, 2)$; using Corollary~\ref{cor: B values for small k} we then have $B_{\d_0, 3} = 60 > M$. Also $c(\Psi)_{ss} = 21 - \z$, so according as $\z = 0$ or $1$ we take $\d_0 = (5, 5, 5, 5, 5, 1)$ or $(5, 5, 5, 5, 5)$; using Corollary~\ref{cor: B values for small k} we then have $B_{\d_0, 3} = 57 > M$ or $B_{\d_0, 3} = 54 > M$. Therefore if $k \in [3, \frac{d}{2}]$ the quadruple $(G, \lambda, p, k)$ satisfies $\ssdiamcon$ and $\udiamcon$.
\end{proof}

\begin{prop}\label{prop: G_2, omega_2, k geq 2, nets}
Let $G = G_2$ and $\lambda = \omega_2$ with $p \neq 3$; then for $k \in [2, \frac{d}{2}]$ the quadruple $(G, \lambda, p, k)$ satisfies $\ssdiamcon$ and $\udiamcon$.
\end{prop}

\begin{proof}
The weight table is as follows.
$$
\begin{array}{|*4{>{\ss}c|}}
\hline
i & \mu & |W.\mu| & m_\mu \\
\hline
2 & \omega_2 & 6 & 1 \\
1 & \omega_1 & 6 & 1 \\
0 &    0     & 1 & 2 \\
\hline
\end{array}
$$
We have $M = 12$, $M_3 = 10$ and $M_2 = 8$; we take $k_0 = 2$.

Take $\Psi = \langle \alpha_2 \rangle$ of type $A_1$. The $\Psi$-net table is as follows.
$$
\begin{array}{|*7{>{\ss}c|}}
\hline
\multicolumn{4}{|>{\ss}c|}{\Psi\mathrm{-nets}} & & \multicolumn{2}{|>{\ss}c|}{c(u_\Psi)} \\
\cline{1-4} \cline{6-7}
    \bar\nu    & n_0 & n_1 & n_2 & m & p = 2 & p \geq 5 \\
\hline
    2\bom_2    &  1  &  0  &  2  & 1 &   1   &     2    \\
     \bom_2    &  0  &  0  &  2  & 2 &   2   &     2    \\
     \bom_2    &  0  &  2  &  0  & 2 &   2   &     2    \\
       0       &  0  &  1  &  0  & 2 &       &          \\
\hline
\multicolumn{5}{c|}{}                &   5   &     6    \\
\cline{6-7}
\end{array}
$$
Thus $c(\Psi)_u \geq 5$, so we may take $\d_0 = (9, 5)$; using Proposition~\ref{prop: B value when t = 2} we then have $B_{\d_0, 2} = 10 > 6 = \dim {u_\Psi}^G$. Each of the remaining non-trivial unipotent classes has $\tilde A_1$ in its closure by Lemma~\ref{lem: root elt class in closure of any non-triv class}.

Now take $\Psi = \langle \alpha_1 \rangle$ of type $\tilde A_1$. The $\Psi$-net table is as follows.
$$
\begin{array}{|*{10}{>{\ss}c|}}
\hline
\multicolumn{4}{|>{\ss}c|}{\Psi\mathrm{-nets}} & & \multicolumn{3}{|>{\ss}c|}{c(s)} & \multicolumn{2}{|>{\ss}c|}{c(u_\Psi)} \\
\cline{1-4} \cline{6-10}
    \bar\nu    & n_0 & n_1 & n_2 & m & r = 2 & r = 3 & r \geq 5 & p = 2 & p \geq 5 \\
\hline
    3\bom_1    &  0  &  2  &  2  & 2 &   4   &   4   &     6    &   4   &     6    \\
     \bom_1    &  1  &  2  &  0  & 1 &   2   &   2   &     2    &   1   &     2    \\
       0       &  0  &  0  &  1  & 2 &       &       &          &       &          \\
\hline
\multicolumn{5}{c|}{}                &   6   &   6   &     8    &   5   &     8    \\
\cline{6-10}
\end{array}
$$
Thus if $r = 2$ or $3$ then $c(\Psi)_{ss} = 6$, so we take $\d_0 = (8, 6)$; using Proposition~\ref{prop: B value when t = 2} we then have $B_{\d_0, 2} = 12 > M_r$. If $p = 2$ then $c(\Psi)_u = 5$, so we take $\d_0 = (9, 5)$; using Proposition~\ref{prop: B value when t = 2} we then have $B_{\d_0, 2} = 10 > M_2 = 8 = \dim {u_\Psi}^G$. If instead $r \geq 5$ and $p \geq 5$ then $c(\Psi)_{ss} = c(\Psi)_u = 8$, so we take $\d_0 = (6, 6, 2)$; using Corollary~\ref{cor: B values for small k} we then have $B_{\d_0, 2} = 14 > M$. Therefore if $k \in [2, \frac{d}{2}]$ the quadruple $(G, \lambda, p, k)$ satisfies $\ssdiamcon$ and $\udiamcon$.
\end{proof}

This completes the treatment of the individual cases listed in Table~\ref{table: remaining large quadruples}.

\section{Analysis of infinite families}\label{sect: large quadruple infinite families}

Again we continue with the notation of Section~\ref{sect: conditions}. Our strategy for handling the infinite families of cases listed in Table~\ref{table: remaining large quadruples} will usually be more direct than that of Section~\ref{sect: large quadruple individual cases}: we shall simply take $g$ to be either an arbitrary $u \in G_{(p)}$ or an arbitrary $s \in G_{(r)}$ for some $r \in \P'$, let $\d = (d_1, d_2, \dots)$ be the tuple associated to $g$ in Proposition~\ref{prop: codim formula for quadruples}, and seek to show that $B_{\d, k_0} > \dim g^G$. We will however on occasion use weight tables and $\Psi$-nets as in Section~\ref{sect: large quadruple individual cases}. Once more we assume $s$ lies in $T$, and write $\Phi(s) = \{ \alpha \in \Phi : \alpha(s) = 1 \}$, so that $C_G(s)^\circ = \langle T, X_\alpha : \alpha \in \Phi(s) \rangle$.

We shall use the standard notation given in Section~\ref{sect: notation} for the roots of a root system of classical type, and extend it to the weights in $\Lambda(V)$.

We begin with a lemma on eigenspace dimensions for the action of $s \in G_{(r)}$ on $\L(G)$. Recall that we let $\eta_r$ be a generator of the group of $r$th roots of unity in $K^*$.

\begin{lem}\label{lem: eigenspace dimension bound on Lie algebra}
Let $G$ be a classical group, and $s \in G_{(r)}$ for some prime $r \in \P'$ which is good for $G$. Then if $G = A_\ell$ and $r | \ell + 1$ we have $\dim \L(G)_{\eta_r}(s) \leq \dim \L(G)_1(s) + 1$, while in all other cases we have $\dim \L(G)_{\eta_r}(s) \leq \dim \L(G)_1(s)$.
\end{lem}

\begin{proof}
We have $\dim \L(G)_1(s) = \ell + |\{ \alpha \in \Phi : \alpha(s) = 1 \}|$ and $\dim \L(G)_{\eta_r}(s) = |\{ \alpha \in \Phi : \alpha(s) = \eta_r \}|$. Write $\delta = \dim \L(G)_1(s) - \dim \L(G)_{\eta_r}(s)$; so we must show that $\delta \geq -1$ if $G = A_\ell$ and $r | \ell + 1$, and $\delta \geq 0$ otherwise. For all $\alpha \in \Phi$ the value $\alpha(s)$ is a power of $\eta_r$; we consider what this implies for the values $\ve_i(s)$.

First take $G = A_\ell$. Write $\ve_1(s) = \xi$; then for all $i$ there exists $j$ with $\ve_i(s) = \xi{\eta_r}^j$. For $j = 0, 1, \dots, r - 1$ set $m_j = | \{ i : \ve_i(s) = \xi{\eta_r}^j \}|$; then $\sum_{j = 0}^{r - 1}m_j = \ell + 1$. We then have
\begin{eqnarray*}
\dim \L(G)_1(s)        & = & m_0(m_0 - 1) + m_1(m_1 - 1) + \cdots + m_{r - 1}(m_{r - 1} - 1) + \ell, \\
\dim \L(G)_{\eta_r}(s) & = & m_0m_1 + m_1m_2 + \cdots + m_{r - 2}m_{r - 1} + m_{r - 1}m_0.
\end{eqnarray*}
Thus
\begin{eqnarray*}
\delta &   =  & {\ts{\frac{1}{2}}}\left[(m_0 - m_1)^2 + (m_1 - m_2)^2 + \cdots + (m_{r - 1} - m_0)^2\right] - 1 \\
       & \geq & -1;
\end{eqnarray*}
we have equality if and only if $m_0 = m_1 = \cdots = m_{r - 1}$, which forces $r$ to divide $\ell + 1$.

Next take $G = D_\ell$. Since $(2\ve_1)(s) = (\ve_1 - \ve_2)(s).(\ve_1 + \ve_2)(s)$ is a power of $\eta_r$, and $r$ is odd, there exists $\xi \in \{ \pm 1 \}$ such that for some $j$ we have $\ve_1(s) = \xi{\eta_r}^j$; then for all $i$ there exists $j$ with $\ve_i(s) = \xi{\eta_r}^j$. For $j = 0, 1, \dots, r - 1$ set $m_j = | \{ i : \ve_i(s) = \xi{\eta_r}^j \}|$; then $\sum_{j = 0}^{r - 1}m_j = \ell$. We then have
\begin{eqnarray*}
\dim \L(G)_1(s)        & = & m_0(m_0 - 1) + m_1(m_1 - 1) + \cdots + m_{r - 1}(m_{r - 1} - 1) \\
                       &   & {} + m_0(m_0 - 1) + m_1m_{r - 1} + m_2m_{r - 2} + \cdots + m_{r - 1}m_1 + \ell, \\
\dim \L(G)_{\eta_r}(s) & = & m_0m_1 + m_1m_2 + \cdots + m_{r - 2}m_{r - 1} + m_{r - 1}m_0 \\
                       &   & {} + m_0m_1 + m_{r - 1}m_2 + \cdots + m_{\frac{r + 3}{2}}m_{\frac{r - 1}{2}} + {\ts{\frac{1}{2}}}m_{\frac{r + 1}{2}}(m_{\frac{r + 1}{2}} - 1) \\
                       &   & {} + m_0m_{r - 1} + m_1m_{r - 2} + \cdots + m_{\frac{r - 3}{2}}m_{\frac{r + 1}{2}} + {\ts{\frac{1}{2}}}m_{\frac{r - 1}{2}}(m_{\frac{r - 1}{2}} - 1).
\end{eqnarray*}
Thus
\begin{eqnarray*}
\delta & = & {\ts{\frac{1}{2}}}\left[\left((2m_0 - m_1 - m_{r - 1} - {\ts{\frac{1}{2}}})^2 - {\ts{\frac{1}{4}}}\right) + \left((m_1 - m_2 - m_{r - 2} + m_{r - 1} - {\ts{\frac{1}{2}}})^2 - {\ts{\frac{1}{4}}}\right) \right. \\
       &   & \quad {} + \cdots + \left.\left((m_{\frac{r - 3}{2}} - m_{\frac{r - 1}{2}} - m_{\frac{r + 1}{2}} + m_{\frac{r + 3}{2}} - {\ts{\frac{1}{2}}})^2 - {\ts{\frac{1}{4}}}\right)\right] \\
       & \geq & 0.
\end{eqnarray*}

Next take $G = B_\ell$. For all $i$ there exists $j$ with $\ve_i(s) = {\eta_r}^j$. For $j = 0, 1, \dots, r - 1$ set $m_j = | \{ i : \ve_i(s) = {\eta_r}^j \}|$; then $\sum_{j = 0}^{r - 1}m_j = \ell$. We then have
\begin{eqnarray*}
\dim \L(G)_1(s)        & = & m_0(m_0 - 1) + m_1(m_1 - 1) + \cdots + m_{r - 1}(m_{r - 1} - 1) \\
                       &   & {} + m_0(m_0 - 1) + m_1m_{r - 1} + m_2m_{r - 2} + \cdots + m_{r - 1}m_1 \\
                       &   & {} + 2m_0 + \ell, \\
\dim \L(G)_{\eta_r}(s) & = & m_0m_1 + m_1m_2 + \cdots + m_{r - 2}m_{r - 1} + m_{r - 1}m_0 \\
                       &   & {} + m_0m_1 + m_{r - 1}m_2 + \cdots + m_{\frac{r + 3}{2}}m_{\frac{r - 1}{2}} + {\ts{\frac{1}{2}}}m_{\frac{r + 1}{2}}(m_{\frac{r + 1}{2}} - 1) \\
                       &   & {} + m_0m_{r - 1} + m_1m_{r - 2} + \cdots + m_{\frac{r - 3}{2}}m_{\frac{r + 1}{2}} + {\ts{\frac{1}{2}}}m_{\frac{r - 1}{2}}(m_{\frac{r - 1}{2}} - 1) \\
                       &   & {} + m_1 + m_{r - 1}.
\end{eqnarray*}
Thus
\begin{eqnarray*}
\delta & = & {\ts{\frac{1}{2}}}\left[\left((2m_0 - m_1 - m_{r - 1} + {\ts{\frac{1}{2}}})^2 - {\ts{\frac{1}{4}}}\right) + \left((m_1 - m_2 - m_{r - 2} + m_{r - 1} - {\ts{\frac{1}{2}}})^2 - {\ts{\frac{1}{4}}}\right) \right. \\
       &   & \quad {} + \cdots + \left.\left((m_{\frac{r - 3}{2}} - m_{\frac{r - 1}{2}} - m_{\frac{r + 1}{2}} + m_{\frac{r + 3}{2}} - {\ts{\frac{1}{2}}})^2 - {\ts{\frac{1}{4}}}\right)\right] \\
       & \geq & 0.
\end{eqnarray*}

Finally take $G = C_\ell$. As in the $D_\ell$ case, there exists $\xi \in \{ \pm 1 \}$ such that for all $i$ there exists $j$ with $\ve_i(s) = \xi{\eta_r}^j$. For $j = 0, 1, \dots, r - 1$ set $m_j = | \{ i : \ve_i(s) = \xi{\eta_r}^j \}|$; then $\sum_{j = 0}^{r - 1}m_j = \ell$. We then have
\begin{eqnarray*}
\dim \L(G)_1(s)        & = & m_0(m_0 - 1) + m_1(m_1 - 1) + \cdots + m_{r - 1}(m_{r - 1} - 1) \\
                       &   & {} + m_0(m_0 - 1) + m_1m_{r - 1} + m_2m_{r - 2} + \cdots + m_{r - 1}m_1 \\
                       &   & {} + 2m_0 + \ell, \\
\dim \L(G)_{\eta_r}(s) & = & m_0m_1 + m_1m_2 + \cdots + m_{r - 2}m_{r - 1} + m_{r - 1}m_0 \\
                       &   & {} + m_0m_1 + m_{r - 1}m_2 + \cdots + m_{\frac{r + 3}{2}}m_{\frac{r - 1}{2}} + {\ts{\frac{1}{2}}}m_{\frac{r + 1}{2}}(m_{\frac{r + 1}{2}} - 1) \\
                       &   & {} + m_0m_{r - 1} + m_1m_{r - 2} + \cdots + m_{\frac{r - 3}{2}}m_{\frac{r + 1}{2}} + {\ts{\frac{1}{2}}}m_{\frac{r - 1}{2}}(m_{\frac{r - 1}{2}} - 1) \\
                       &   & {} + m_{\frac{r + 1}{2}} + m_{\frac{r - 1}{2}}.
\end{eqnarray*}
Thus
\begin{eqnarray*}
\delta & = & {\ts{\frac{1}{2}}}\left[\left((2m_0 - m_1 - m_{r - 1} + {\ts{\frac{1}{2}}})^2 - {\ts{\frac{1}{4}}}\right) + \left((m_1 - m_2 - m_{r - 2} + m_{r - 1} + {\ts{\frac{1}{2}}})^2 - {\ts{\frac{1}{4}}}\right) \right. \\
       &   & \quad {} + \cdots + \left.\left((m_{\frac{r - 3}{2}} - m_{\frac{r - 1}{2}} - m_{\frac{r + 1}{2}} + m_{\frac{r + 3}{2}} + {\ts{\frac{1}{2}}})^2 - {\ts{\frac{1}{4}}}\right)\right] \\
       & \geq & 0.
\end{eqnarray*}
The result follows.
\end{proof}

We now consider our infinite families. As before, given a quadruple $(G, \lambda, p, k)$ we write $V = L(\lambda)$. In most of the cases we have $k_0 = 2$; note that Corollary~\ref{cor: B values for small k} allows us to write $B_{\d, 2} = 2d - 2d_1 - x$ where
$$
x =
\begin{cases}
0             & \hbox{if } d_1 \geq d_2 + 2, \\
d_2 + 2 - d_1 & \hbox{if } d_1 < d_2 + 2,
\end{cases}
$$
so that $x \in [0, 2]$.

We begin with those families where $V$ is the quotient of the Lie algebra of $G$ by its centre.

\begin{prop}\label{prop: classical Lie algebras, k = 2}
Let $G = A_\ell$ for $\ell \in [2, \infty)$ and $\lambda = \omega_1 + \omega_\ell$, or $G = B_\ell$ for $\ell \in [3, \infty)$ and $\lambda = \omega_2$ with $p \geq 3$, or $G = C_\ell$ for $\ell \in [3, \infty)$ and $\lambda = 2\omega_1$ with $p \geq 3$, or $G = D_\ell$ for $\ell \in [4, \infty)$ and $\lambda = \omega_2$; then for $k \in [2, \frac{d}{2}]$ the quadruple $(G, \lambda, p, k)$ satisfies $\ssdiamcon$ and $\udiamcon$.
\end{prop}

\begin{proof}
In all these cases we have $V \!=\! \L(G)/Z(\L(G))$. Write $z \!=\! \dim Z(\L(G))$; then
$$
z =
\begin{cases}
\z_{p, \ell + 1}            & \hbox{if } G = A_\ell, \\
\z_{p, 2}(1 + \z_{2, \ell}) & \hbox{if } G = D_\ell, \\
0                           & \hbox{otherwise}.
\end{cases}
$$
In particular we have $z \leq 2$.

First take $u \in G_{(p)}$; then $\dim C_{\L(G)}(u) = \dim C_G(u) + z'$ where
$$
\begin{cases}
0 \leq z' \leq z    & \hbox{if } G = A_\ell, \\
0 \leq z' \leq \ell & \hbox{if } G = D_\ell \hbox{ and } p = 2, \\
z' = 0              & \hbox{otherwise}.
\end{cases}
$$
Write $c = \codim C_V(u)$; then we have
\begin{eqnarray*}
c & = & \dim V - \dim C_V(u) \\
  & = & (\dim \L(G) - z) - (\dim C_{\L(G)}(u) - z) \\
  & = & \dim \L(G) - \dim C_{\L(G)}(u) \\
  & = & \dim G - (\dim C_G(u) + z') \\
  & = & \dim u^G - z'.
\end{eqnarray*}
We have $c = d - d_1$. Thus
$$
2d - 2d_1 - x = 2c - x = \dim u^G + (\dim u^G - 2z' - x).
$$
We clearly have $\dim u^G - 2z' - x > 0$ unless either $G = A_2$ and $\dim u^G = 4$ with $z' = 1$, or $G = D_4$ and $\dim u^G = 10$ with $z' = 4$; in the former case we have $d = 7$, so $d_1 = d - \dim u^G + z' = 4$, and then $d_1 > d_2$, whence $x < 2$, while in the latter case we have $d = 26$, so $d_1 = d - \dim u^G + z' = 20$, and then $d_1 > d_2 + 2$, whence $x = 0$. In all cases here we therefore have $2d - 2d_1 - x > \dim u^G$.

Thus $B_{\d, 2} > \dim u^G$. Therefore if $k \in [2, \frac{d}{2}]$ the quadruple $(G, \lambda, p, k)$ satisfies $\udiamcon$.

Now take $s \in G_{(r)}$ for $r \in \P'$; then $C_{\L(G)}(s) = \L(C_G(s))$, so $\dim C_{\L(G)}(s) = \dim C_G(s)$. Write $c = \codim C_V(s)$; then we have
\begin{eqnarray*}
c & = & \dim V - \dim C_V(s) \\
  & = & (\dim \L(G) - z) - (\dim C_{\L(G)}(s) - z) \\
  & = & \dim G - \dim C_G(s) \\
  & = & \dim s^G.
\end{eqnarray*}
We have $c = d - d_i$ for some $i \geq 1$; then $d_i = (\dim \L(G) - z) - (\dim \L(G) - \dim C_{\L(G)}(s)) = \dim C_{\L(G)}(s) - z = \dim \L(G)_1(s) - z$. Write
$$
a =
\begin{cases}
1 & \hbox{if } G = A_\ell \hbox{ and } r | \ell + 1, \\
0 & \hbox{otherwise};
\end{cases}
$$
observe that $z + a \leq 2$.

First suppose either $r \geq 3$, or $r = 2$ and $G = A_\ell$. If $j \neq i$ there is a primitive $r$th root of unity $\eta$ such that $d_j = \dim V_\eta(s) = \dim \L(G)_\eta(s)$; using Lemma~\ref{lem: eigenspace dimension bound on Lie algebra} we see that $d_j \leq \dim \L(G)_1(s) + a = d_i + z + a$. Thus
\begin{eqnarray*}
2d - 2d_1 - x & \geq & 2(d - d_i - z - a) - x \\
              &   =  & (d - d_i) + (d - d_i - 2z - 2a - x) \\
              &   =  & \dim s^G + (|\Phi| - |\Phi(s)| - 2(z + a) - x);
\end{eqnarray*}
since $|\Phi| - |\Phi(s)| \geq 2\ell$ we have $(2d - 2d_1 - x) - \dim s^G \geq 2\ell - 2(z + a) - x$, which is clearly positive unless $G = A_2$ or $A_3$, in which case we cannot have both $p | \ell + 1$ and $r | \ell + 1$, so $z + a \leq 1$; thus we need only consider $G = A_2$ and $\Phi(s)$ of type $A_1$, in which case the eigenspaces of $s$ on $\L(G)$ have dimensions $4$, $2$ and $2$, so $d_1 = 4 - z > 2 = d_2$ and hence $x \leq 1$. In all cases here we therefore have $2d - 2d_1 - x > \dim s^G$.

Now suppose $r = 2$ and $G = B_\ell$, $C_\ell$ or $D_\ell$; note that then $a = z = 0$, and $\d = (d_1, d_2)$. We have $c = \dim s^G \leq M_2$, so $d_i = d - c \geq d - M_2 = d_{\Phi, 2} \geq \frac{1}{2}(d - \ell)$ by \cite[Lemma~1.2]{Lawdim}. Thus if $i = 1$ then we have
$$
(2d - 2d_1 - x) - \dim s^G = 2d_2 - x - d_2 = d_2 - x > 0
$$
(because $d_2 \geq 2$, and if $d_2 = 2$ then $d_1 > d_2 + 2$ so $x = 0$); if instead $i = 2$ then we have
\begin{eqnarray*}
(2d - 2d_1 - x) - \dim s^G &   =  & (2d_2 - x) - (d - d_2) \\
                           &   =  & 3d_2 - d - x \\
                           & \geq & {\ts{\frac{3}{2}}}(d - \ell) - d - x \\
                           &   =  & {\ts{\frac{1}{2}}}(d - 3\ell - 2x) \\
                           &   >  & 0.
\end{eqnarray*}
In all cases here we therefore have $2d - 2d_1 - x > \dim s^G$.

Thus $B_{\d, 2} > \dim s^G$. Therefore if $k \in [2, \frac{d}{2}]$ the quadruple $(G, \lambda, p, k)$ satisfies $\ssdiamcon$.
\end{proof}

Next we treat the two families of cases which are not $p$-restricted.

\begin{prop}\label{prop: A_ell, omega_1 + q omega_1 and omega_1 + q omega_ell, k geq 2}
Let $G = A_\ell$ for $\ell \in [2, \infty)$ and $\lambda = \omega_1 + q\omega_1$ or $\omega_1 + q\omega_\ell$; then for $k \in [2, \frac{d}{2}]$ the quadruple $(G, \lambda, p, k)$ satisfies $\ssdiamcon$ and $\udiamcon$.
\end{prop}

\begin{proof}
We take $G = \SL_{\ell + 1}(K)$. Recall that $V_{nat} = L(\omega_1) = \langle v_1, \dots, v_{\ell + 1} \rangle$. Take $A \in G$, so that for each $i$ we have $A.v_i = \sum_{i' = 1}^{\ell + 1} a_{i'i} v_{i'}$. In the case $\lambda = \omega_1 + q\omega_1$, we have $V = L(\omega_1) \otimes L(\omega_1)^{(q)}$; we see that $A$ maps $v_i \otimes v_j \mapsto \sum_{i', j' = 1}^{\ell + 1} a_{i'i}{a_{j'j}}^q v_{i'} \otimes v_{j'}$. We may then identify $V$ with the space of $(\ell + 1) \times (\ell + 1)$ matrices $D$ over $K$, and the matrix unit $E_{ij}$ with $v_i \otimes v_j$; then $A$ maps $E_{ij} \mapsto \sum_{i', j' = 1}^{\ell + 1} a_{i'i}{a_{j'j}}^q E_{i'j'}$, and so $\sum_{i, j = 1}^{\ell + 1} d_{ij}E_{ij} \mapsto \sum_{i', j' = 1}^{\ell + 1} \left( \sum_{i, j = 1}^{\ell + 1} a_{i'i}d_{ij}{a_{j'j}}^q \right) E_{i'j'}$, or $D \mapsto AD(A^{(q)})^T$. Similarly in the case $\lambda = \omega_1 + q\omega_\ell$, we may again identify $V$ with the space of $(\ell + 1) \times (\ell + 1)$ matrices $D$ over $K$, but such that $A$ maps $D \mapsto AD(A^{(q)})^{-1}$.

First take $u \in G_{(p)}$; let $A$ be the matrix representing $u$. We may assume
$$
A = \left(
      \begin{array}{cccc}
        J_1 &     &        &     \\
            & J_2 &        &     \\
            &     & \ddots &     \\
            &     &        & J_t \\
      \end{array}
    \right),
$$
where $J_i$ is a single Jordan block of size $m_i$, with $m_1 \geq m_2 \geq \cdots \geq m_t$; then $\dim C_G(u) = 1.m_1 + 3.m_2 + 5.m_3 + \cdots + (2t - 1)m_t - 1$. Note that $A^{(q)} = A$. Given $D \in V$, write
$$
D = \left(
      \begin{array}{cccc}
        D_{11} & D_{12} & \cdots & D_{1t} \\
        D_{21} & D_{22} & \cdots & D_{2t} \\
        \vdots & \vdots & \ddots & \vdots \\
        D_{t1} & D_{t2} & \cdots & D_{tt} \\
      \end{array}
    \right),
$$
where $D_{ij}$ is an $m_i \times m_j$ matrix. If $\lambda = \omega_1 + q\omega_1$, then
$$
D \in C_V(u) \iff ADA^T = D \iff \hbox{for all $i$ and $j$ we have } J_i D_{ij} {J_j}^T = D_{ij};
$$
if instead $\lambda = \omega_1 + q\omega_\ell$, then
$$
D \in C_V(u) \iff ADA^{-1} = D \iff \hbox{for all $i$ and $j$ we have } J_i D_{ij} = D_{ij} J_j.
$$
In both cases it is easy to check that for a fixed pair $(i, j)$ the set of such matrices $D_{ij}$ has dimension $\min(m_i, m_j)$. Thus $\dim C_V(u) = \sum_{i, j} \min(m_i, m_j) = 1.m_1 + 3.m_2 + 5.m_3 + \cdots + (2t - 1)m_t = \dim C_G(u) + 1$. Write $c = \codim C_V(u)$; then we have
$$
c = \dim V - \dim C_V(u) = (\dim G + 1) - (\dim C_G(u) + 1) = \dim u^G.
$$
We have $c = d - d_1$, so as $\dim u^G \geq 2\ell \geq 4$ we have
$$
2d - 2d_1 - x = 2c - x = 2\dim u^G - x > \dim u^G.
$$
Thus $B_{\d, 2} > \dim u^G$. Therefore if $k \in [2, \frac{d}{2}]$ the quadruple $(G, \lambda, p, k)$ satisfies $\udiamcon$.

Now take $s \in G_{(r)}$ for $r \in \P'$; let $A$ be the matrix representing $s$. We may assume $A = (a_{ij})$ is diagonal. We have $A^r = \rho I$ for some $\rho \in K$ with $\rho^{\ell + 1} = 1$; choose $\kappa \in K$ with $\kappa^r = \rho$, then each diagonal entry of $A$ has the form $\kappa{\eta_r}^j$ for some $j$. For $j = 0, 1, \dots, r - 1$ set $m_j = |\{ i : a_{ii} = \kappa {\eta_r}^j \}|$; then $\sum_{j = 0}^{r - 1} m_j = \ell + 1$. We have $\dim C_G(s) = \sum_{j = 0}^{r - 1} {m_j}^2 - 1$, so
$$
\dim s^G = \dim G - \dim C_G(s) = (\ell + 1)^2 - \sum_{i = 0}^{r - 1} {m_i}^2.
$$

Write $\xi = 1$ or $-1$ according as $\lambda = \omega_1 + q\omega_1$ or $\omega_1 + q\omega_\ell$. Each matrix unit $E_{ij}$ is an eigenvector for $s$, with eigenvalue $a_{ii}{a_{jj}}^{\xi q}$; so the eigenvalues are $\kappa^{1 + \xi q} {\eta_r}^h$ for various values of $h$. For a fixed $h$, we have $\dim V_{\kappa^{1 + \xi q} {\eta_r}^h}(s) = \sum_{(i, j)} m_im_j$, where the sum runs over all pairs $(i, j)$ such that $i + \xi q j \equiv h$ (mod $r$); note that $j \mapsto h - \xi q j$ is a permutation $\pi_h$, say, of $\Z/r\Z$, and then we have $\dim V_{\kappa^{1 + \xi q} {\eta_r}^h}(s) = \sum_{j = 0}^{r - 1} m_{\pi_h(j)} m_j$. Thus the various dimensions $d_i$ are the various values $\sum_{j = 0}^{r - 1} m_{\pi_h(j)} m_j$ as $h$ runs from $0$ to $r - 1$. Therefore $d - d_1 = (\ell + 1)^2 - \sum_{j = 0}^{r - 1} m_{\pi_h(j)} m_j$ for some $h$, and so
$$
(d - d_1) - \dim s^G = \sum_{j = 0}^{r - 1} {m_j}^2 - \sum_{j = 0}^{r - 1} m_{\pi_h(j)} m_j = {\ts{\frac{1}{2}}}\sum_{j = 0}^{r - 1} (m_j - m_{\pi_h(j)})^2 \geq 0,
$$
i.e., $d - d_1 \geq \dim s^G$; so as $\dim s^G \geq 2\ell \geq 4$ we have
$$
2d - 2d_1 - x \geq 2\dim s^G - x > \dim s^G.
$$
Thus $B_{\d, 2} > \dim s^G$. Therefore if $k \in [2, \frac{d}{2}]$ the quadruple $(G, \lambda, p, k)$ satisfies $\ssdiamcon$.
\end{proof}

The next few results treat the remaining infinite families where $k_0 = 2$.

\begin{prop}\label{prop: B_ell and D_ell, 2omega_1, k geq 2}
Let $G = B_\ell$ for $\ell \in [2, \infty)$ or $D_\ell$ for $\ell \in [4, \infty)$, and $\lambda = 2\omega_1$ with $p \geq 3$; then for $k \in [2, \frac{d}{2}]$ the quadruple $(G, \lambda, p, k)$ satisfies $\ssdiamcon$ and $\udiamcon$.
\end{prop}

\begin{proof}
First take $u \in G_{(p)}$. Write $\ell' = 2\ell$ or $2\ell - 1$ according as $G = B_\ell$ or $D_\ell$, and $\z = \z_{p, \ell' + 1}$; then $\dim G = \frac{1}{2}\ell'(\ell' + 1)$. Let $H$ be the simply connected group of type $A_{\ell'}$ over $K$; then $\dim Z(\L(H)) = \z$. As we shall see in the proof of Proposition~\ref{prop: B_ell or D_ell, 2omega_1 module}, we have
$$
\L(H) = \L(G) \oplus \tilde V,
$$
where $Z(\L(H)) \leq \tilde V$ and $V = \tilde V/Z(\L(H))$. We have $d = \frac{1}{2}\ell'(\ell' + 3) - \z$.

Let $1^{r_1}, 2^{r_2}, \dots$ be the sizes of the Jordan blocks of $u$ on the natural module for $H$, so that $\sum_i ir_i = \ell' + 1$. Write $n_0 = \sum_i (r_i + r_{i + 1} + \cdots)^2$ and $n_1 = \sum_{i \ \rm{odd}} r_i$, and observe that $n_1 \leq \ell' - 1$. Then $\dim C_H(u) = n_0 - 1$, and $\dim C_{\L(H)}(u) = \dim C_H(u) + z$ where $0 \leq z \leq \z$; moreover $\dim C_G(u) = \frac{1}{2}n_0 - \frac{1}{2}n_1$, and $C_{\L(G)}(u) = \L(C_G(u))$. Thus
\begin{eqnarray*}
d_1 & = & \dim C_V(u) \\
    & = & \dim C_{\tilde V}(u) - \z \\
    & = & \dim C_{\L(H)}(u) - \dim C_{\L(G)}(u) - \z \\
    & = & \dim C_H(u) + z - \dim C_G(u) - \z \\
    & = & n_0 - 1 + z - {\ts{\frac{1}{2}}}n_0 + {\ts{\frac{1}{2}}}n_1 - \z \\
    & = & {\ts{\frac{1}{2}}}n_0 + {\ts{\frac{1}{2}}}n_1 + z - \z - 1,
\end{eqnarray*}
while
$$
\dim u^G = \dim G - \dim C_G(u) = {\ts{\frac{1}{2}}}\ell'(\ell' + 1) - {\ts{\frac{1}{2}}}n_0 + {\ts{\frac{1}{2}}}n_1.
$$
Hence
\begin{eqnarray*}
(2d - 2d_1 - x) - 2\dim u^G &   =  & \ell'(\ell' + 3) - 2\z - n_0 - n_1 - 2z + 2\z + 2 - x \\
                            &      & \quad {} - \ell'(\ell' + 1) + n_0 - n_1 \\
                            &   =  & 2\ell' - 2n_1 - 2z + 2 - x \\
                            & \geq & 2 - 2z + 2 - x \\
                            & \geq & 0,
\end{eqnarray*}
and so $2d - 2d_1 - x \geq 2\dim u^G > \dim u^G$. Thus $B_{\d, 2} > \dim u^G$. Therefore if $k \in [2, \frac{d}{2}]$ the quadruple $(G, \lambda, p, k)$ satisfies $\udiamcon$.

Now take $s \in G_{(r)}$ for $r \in \P'$. We analyse the weights in $V$ as in the proof of Lemma~\ref{lem: eigenspace dimension bound on Lie algebra}; we set $\delta = \dim V_1(s) - \dim V_{\eta_r}(s)$. We shall treat the cases $G = B_\ell$ and $G = D_\ell$ separately.

We begin with $G = B_\ell$. The weights are $\pm 2\ve_i$ for $1 \leq i \leq \ell$, $\pm \ve_i \pm \ve_j$ for $1 \leq i < j \leq \ell$, $\pm \ve_i$ for $1 \leq i \leq \ell$, and $0$; all have multiplicity $1$ except the last, which has multiplicity $\ell - \z$. For all $i$ there exists $j$ with $\ve_i(s) = {\eta_r}^j$. For $j = 0, 1, \dots, r - 1$ set $m_j = | \{ i : \ve_i(s) = {\eta_r}^j \}|$; then $\sum_{j = 0}^{r - 1}m_j = \ell$.

First suppose $r \geq 3$. We then have
\begin{eqnarray*}
\dim V_1(s)        & = & 2m_0 + m_0(m_0 - 1) + m_1(m_1 - 1) + \cdots + m_{r - 1}(m_{r - 1} - 1) \\
                   &   & {} + m_0(m_0 - 1) + m_1m_{r - 1} + m_2m_{r - 2} + \cdots + m_{r - 1}m_1 \\
                   &   & {} + 2m_0 + \ell - \z, \\
\dim V_{\eta_r}(s) & = & m_{\frac{r + 1}{2}} + m_{\frac{r - 1}{2}} + m_0m_1 + m_1m_2 + \cdots + m_{r - 2}m_{r - 1} + m_{r - 1}m_0 \\
                   &   & {} + m_0m_1 + m_{r - 1}m_2 + \cdots + m_{\frac{r + 3}{2}}m_{\frac{r - 1}{2}} + {\ts{\frac{1}{2}}}m_{\frac{r + 1}{2}}(m_{\frac{r + 1}{2}} - 1) \\
                   &   & {} + m_0m_{r - 1} + m_1m_{r - 2} + \cdots + m_{\frac{r - 3}{2}}m_{\frac{r + 1}{2}} + {\ts{\frac{1}{2}}}m_{\frac{r - 1}{2}}(m_{\frac{r - 1}{2}} - 1) \\
                   &   & {} + m_1 + m_{r - 1}.
\end{eqnarray*}
Thus
\begin{eqnarray*}
\delta & = & {\ts{\frac{1}{2}}}\left[\left((2m_0 - m_1 - m_{r - 1} + {\ts{\frac{1}{2}}})^2 - {\ts{\frac{1}{4}}}\right) + \left((m_1 - m_2 - m_{r - 2} + m_{r - 1} + {\ts{\frac{1}{2}}})^2 - {\ts{\frac{1}{4}}}\right) \right. \\
       &   & \quad {} + \cdots + \left.\left((m_{\frac{r - 3}{2}} - m_{\frac{r - 1}{2}} - m_{\frac{r + 1}{2}} + m_{\frac{r + 3}{2}} + {\ts{\frac{1}{2}}})^2 - {\ts{\frac{1}{4}}}\right)\right] \\
       &   & \quad {} + 2m_0 - m_1 - m_{r - 1} - \z \\
       & \geq & 2m_0 - m_1 - m_{r - 1} - \z.
\end{eqnarray*}
Hence $\dim V_{\eta_r}(s) \leq \dim V_1(s) - (2m_0 - m_1 - m_{r - 1} - \z) \leq \dim V_1(s) - 2m_0 + \ell + \z$. Also $|\Phi(s)| = \dim V_1(s) - (2m_0 + \ell - \z)$, so $\dim s^G = |\Phi| - |\Phi(s)| = 2\ell^2 - \dim V_1(s) + (2m_0 + \ell - \z)$. Thus if $d_1 > \dim V_1(s)$, then $d_1 \leq \dim V_1(s) - 2m_0 + \ell + \z$, and so $d_1 + \dim s^G \leq 2\ell^2 + 2\ell$, whence
$$
(2d - 2d_1 - x) - 2\dim s^G \geq 2(2\ell^2 + 3\ell - \z - (2\ell^2 + 2\ell)) - x = 2\ell - 2\z - x \geq 0.
$$
If instead $d_1 = \dim V_1(s)$, then $d_1 + \dim s^G = 2\ell^2 + \ell + 2m_0 - \z \leq 2\ell^2 + \ell + 2(\ell - 1) - \z = 2\ell^2 + 3\ell - 2 - \z = d - 2$, and so
$$
(2d - 2d_1 - x) - 2\dim s^G \geq 4 - x > 0.
$$
In all cases here we therefore have $2d - 2d_1 - x \geq 2\dim s^G > \dim s^G$.

Now suppose $r = 2$. We then have
\begin{eqnarray*}
\dim V_1(s)    & = & 2\ell + 2m_0(m_0 - 1) + 2m_1(m_1 - 1) + 2m_0 + \ell - \z, \\
\dim V_{-1}(s) & = & 4m_0m_1 + 2m_1.
\end{eqnarray*}
Thus
$$
\delta = 2[(m_0 - m_1 + {\ts{\frac{1}{2}}})^2 - {\ts{\frac{1}{4}}}] + \ell - \z > 0.
$$
Hence $d_1 = \dim V_1(s)$ and $d_2 = \dim V_{-1}(s)$. Also $|\Phi(s)| = \dim V_1(s) - (2\ell + \ell - \z)$, so $\dim s^G = |\Phi| - |\Phi(s)| = 2\ell^2 - \dim V_1(s) + (2\ell + \ell - \z) = d - \dim V_1(s) = d - d_1$. Therefore
$$
(2d - 2d_1 - x) - 2\dim s^G = - x,
$$
and so $2d - 2d_1 - x = 2\dim s^G - x = \dim s^G + (\dim s^G - x) > \dim s^G$.

Thus $B_{\d, 2} > \dim s^G$. Therefore if $G = B_\ell$ and $k \in [2, \frac{d}{2}]$ the quadruple $(G, \lambda, p, k)$ satisfies $\ssdiamcon$.

We now take $G = D_\ell$. The weights are $\pm 2\ve_i$ for $1 \leq i \leq \ell$, $\pm \ve_i \pm \ve_j$ for $1 \leq i < j \leq \ell$, and $0$; all have multiplicity $1$ except the last, which has multiplicity $\ell - 1 - \z$. Thus the non-zero weights form a root system of type $C_\ell$.

First suppose $r \geq 3$. As in the proof of Lemma~\ref{lem: eigenspace dimension bound on Lie algebra}, there exists $\xi \in \{ \pm 1 \}$ such that for all $i$ there exists $j$ with $\ve_i(s) = \xi {\eta_r}^j$. For $j = 0, 1, \dots, r - 1$ set $m_j = | \{ i : \ve_i(s) = \xi {\eta_r}^j \}|$; then $\sum_{j = 0}^{r - 1}m_j = \ell$. We then have
\begin{eqnarray*}
\dim V_1(s)        & = & 2m_0 + m_0(m_0 - 1) + m_1(m_1 - 1) + \cdots + m_{r - 1}(m_{r - 1} - 1) \\
                   &   & {} + m_0(m_0 - 1) + m_1m_{r - 1} + m_2m_{r - 2} + \cdots + m_{r - 1}m_1 + \ell - 1 - \z, \\
\dim V_{\eta_r}(s) & = & m_{\frac{r + 1}{2}} + m_{\frac{r - 1}{2}} + m_0m_1 + m_1m_2 + \cdots + m_{r - 2}m_{r - 1} + m_{r - 1}m_0 \\
                   &   & {} + m_0m_1 + m_{r - 1}m_2 + \cdots + m_{\frac{r + 3}{2}}m_{\frac{r - 1}{2}} + {\ts{\frac{1}{2}}}m_{\frac{r + 1}{2}}(m_{\frac{r + 1}{2}} - 1) \\
                   &   & {} + m_0m_{r - 1} + m_1m_{r - 2} + \cdots + m_{\frac{r - 3}{2}}m_{\frac{r + 1}{2}} + {\ts{\frac{1}{2}}}m_{\frac{r - 1}{2}}(m_{\frac{r - 1}{2}} - 1).
\end{eqnarray*}
Thus
\begin{eqnarray*}
\delta & = & {\ts{\frac{1}{2}}}\left[\left((2m_0 - m_1 - m_{r - 1} + {\ts{\frac{1}{2}}})^2 - {\ts{\frac{1}{4}}}\right) + \left((m_1 - m_2 - m_{r - 2} + m_{r - 1} + {\ts{\frac{1}{2}}})^2 - {\ts{\frac{1}{4}}}\right) \right. \\
       &   & \quad {} + \cdots + \left.\left((m_{\frac{r - 3}{2}} - m_{\frac{r - 1}{2}} - m_{\frac{r + 1}{2}} + m_{\frac{r + 3}{2}} + {\ts{\frac{1}{2}}})^2 - {\ts{\frac{1}{4}}}\right)\right] - 1 - \z \\
       & \geq & -1 - \z.
\end{eqnarray*}
Hence $\dim V_{\eta_r}(s) \leq \dim V_1(s) + 1 + \z$. Also $|\Phi(s)| = \dim V_1(s) - (2m_0 + \ell - 1 - \z)$, so $\dim s^G = |\Phi| - |\Phi(s)| = 2\ell(\ell - 1) - \dim V_1(s) + (2m_0 + \ell - 1 - \z)$. Thus $d_1 + \dim s^G \leq 2\ell^2 - 2\ell + 2m_0 + \ell = 2\ell^2 - \ell + 2m_0$, and so
\begin{eqnarray*}
(2d - 2d_1 - x) - 2\dim s^G & \geq & 2(2\ell^2 + \ell - 1 - \z - (2\ell^2 - \ell + 2m_0)) - x \\
                            &   =  & 2(2\ell - 2m_0 - 1 - \z) - x \\
                            & \geq & -x.
\end{eqnarray*}
Therefore $(2d - 2d_1 - x) - \dim s^G \geq \dim s^G - x > 0$, so $2d - 2d_1 - x > \dim s^G$.

Now suppose $r = 2$. Here either all $\ve_i(s)$ are $\pm 1$, or all $\ve_i(s)$ are $\pm \eta_4$.

If all $\ve_i(s)$ are $\pm 1$, set $m_0 = |\{ i : \ve_i(s) = 1 \}|$ and $m_1 = |\{ i : \ve_i(s) = -1 \}|$; then $m_0 + m_1 = \ell$. We then have
\begin{eqnarray*}
\dim V_1(s)    & = & 2\ell + 2m_0(m_0 - 1) + 2m_1(m_1 - 1) + \ell - 1 - \z, \\
\dim V_{-1}(s) & = & 4m_0m_1.
\end{eqnarray*}
Thus
$$
\delta = 2(m_0 - m_1)^2 + \ell - 1 - \z > 0.
$$
Hence $d_1 = \dim V_1(s)$ and $d_2 = \dim V_{-1}(s)$. Also $|\Phi(s)| = \dim V_1(s) - (2\ell + \ell - 1 - \z)$, so $\dim s^G = |\Phi| - |\Phi(s)| = 2\ell(\ell - 1) - \dim V_1(s) + (2\ell + \ell - 1 - \z) = d - \dim V_1(s) = d - d_1$. Therefore
$$
(2d - 2d_1 - x) - 2\dim s^G = -x,
$$
and so $2d - 2d_1 - x = 2\dim s^G - x = \dim s^G + (\dim s^G - x) > \dim s^G$.

If instead all $\ve_i(s)$ are $\pm \eta_4$, set $m_0 = |\{ i : \ve_i(s) = \eta_4 \}|$ and $m_1 = |\{ i : \ve_i(s) = -\eta_4 \}|$; then $m_0 + m_1 = \ell$. We then have
\begin{eqnarray*}
\dim V_1(s)    & = & m_0(m_0 - 1) + m_1(m_1 - 1) + 2m_0m_1 + \ell - 1 - \z = \ell^2 - 1 - \z, \\
\dim V_{-1}(s) & = & 2\ell + 2m_0m_1 + m_0(m_0 - 1) + m_1(m_1 - 1) = \ell^2 + \ell.
\end{eqnarray*}
Thus
$$
\delta = -\ell - 1 - \z < 0.
$$
Hence $d_1 = \dim V_{-1}(s)$ and $d_2 = \dim V_1(s)$. Also $|\Phi(s)| = \frac{1}{2}|\Phi|$, so $\dim s^G = |\Phi| - |\Phi(s)| = \frac{1}{2}|\Phi| = \ell(\ell - 1)$. Therefore
\begin{eqnarray*}
(2d - 2d_1 - x) - \dim s^G & = & 2d_2 - x - \dim s^G \\
                           & = & 2\ell^2 - 2 - 2\z - x - \ell^2 + \ell \\
                           & = & \ell^2 + \ell - 2 - 2\z - x \\
                           & > & 0,
\end{eqnarray*}
and so $2d - 2d_1 - x > \dim s^G$.

Thus $B_{\d, 2} > \dim s^G$. Therefore if $G = D_\ell$ and $k \in [2, \frac{d}{2}]$ the quadruple $(G, \lambda, p, k)$ satisfies $\ssdiamcon$.
\end{proof}

\begin{prop}\label{prop: C_ell, omega_2, k geq 2}
Let $G = C_\ell$ for $\ell \in [4, \infty)$ and $\lambda = \omega_2$; then for $k \in [2, \frac{d}{2}]$ the quadruple $(G, \lambda, p, k)$ satisfies $\ssdiamcon$ and $\udiamcon$.
\end{prop}

\begin{proof}
Write $\z = \z_{p, \ell}$. The weight table is as follows.
$$
\begin{array}{|*4{>{\ss}c|}}
\hline
i & \mu & |W.\mu| & m_\mu \\
\hline
1 & \omega_2 & 2\ell(\ell - 1) &       1       \\
0 &    0     &        1        & \ell - 1 - \z \\
\hline
\end{array}
$$
Thus $d = 2\ell^2 - \ell - 1 - \z$.

First take $u \in G_{(p)}$. We treat the cases $p \geq 3$ and $p = 2$ separately.

Begin by assuming $p \geq 3$. Much as in the proof of Proposition~\ref{prop: B_ell and D_ell, 2omega_1, k geq 2}, we let $H$ be the simply connected group of type $A_{2\ell - 1}$ over $K$; then $\dim Z(\L(H)) = \z$, and we have
$$
\L(H) = \L(G) \oplus \tilde V,
$$
where $Z(\L(H)) \leq \tilde V$ and $V = \tilde V/Z(\L(H))$.

Let $1^{r_1}, 2^{r_2}, \dots$ be the sizes of the Jordan blocks of $u$ on the natural module for $H$, so that $\sum_i ir_i = 2\ell$. Write $n_0 = \sum_i (r_i + r_{i + 1} + \cdots)^2$ and $n_1 = \sum_{i \ \rm{odd}} r_i$. Then $\dim C_H(u) = n_0 - 1$, and $\dim C_{\L(H)}(u) = \dim C_H(u) + z$ where $0 \leq z \leq \z$; moreover $\dim C_G(u) = \frac{1}{2}n_0 + \frac{1}{2}n_1$, and $C_{\L(G)}(u) = \L(C_G(u))$. Thus
\begin{eqnarray*}
d_1 & = & \dim C_V(u) \\
    & = & \dim C_{\tilde V}(u) - \z \\
    & = & \dim C_{\L(H)}(u) - \dim C_{\L(G)}(u) - \z \\
    & = & \dim C_H(u) + z - \dim C_G(u) - \z \\
    & = & n_0 - 1 + z - {\ts{\frac{1}{2}}}n_0 - {\ts{\frac{1}{2}}}n_1 - \z \\
    & = & {\ts{\frac{1}{2}}}n_0 - {\ts{\frac{1}{2}}}n_1 + z - \z - 1,
\end{eqnarray*}
while
$$
\dim u^G = \dim G - \dim C_G(u) = 2\ell^2 + \ell - {\ts{\frac{1}{2}}}n_0 - {\ts{\frac{1}{2}}}n_1.
$$
Hence
\begin{eqnarray*}
(2d - 2d_1 - x) - 2\dim u^G & = & 4\ell^2 - 2\ell - 2 - 2\z - n_0 + n_1 - 2z + 2\z + 2 - x \\
                            &   & \quad {} - 4\ell^2 - 2\ell + n_0 + n_1 \\
                            & = & -4\ell + 2n_1 - 2z - x,
\end{eqnarray*}
and so
\begin{eqnarray*}
(2d - 2d_1 - x) - \dim u^G &   =  & \dim u^G - 4\ell + 2n_1 - 2z - x \\
                           & \geq & \dim u^G - (4\ell + 4).
\end{eqnarray*}
Thus $2d - 2d_1 - x > \dim u^G$ provided $\dim u^G > 4\ell + 4$; so we must consider the unipotent classes of dimension at most $4\ell + 4$. We analyse the possibilities using Jordan normal form and the partial order on unipotent classes given by containment of closures.

If $r_i > 0$ for some $i \geq 4$, then by Lemma~\ref{lem: C_2 and A_2 in C_ell}(i) the class has $C_2$ in its closure, whose dimension is $6\ell - 4$; if $\ell \in [5, \infty)$ we have $6\ell - 4 > 4\ell + 4$, while if $\ell = 4$ the only such class needing consideration is $C_2$ itself, for which $n_1 = 4$ and hence
\begin{eqnarray*}
\dim u^G - 4\ell + 2n_1 - 2z - x & = & 20 - 16 + 8 - 2z - x \\
                                 & = & 12 - 2z - x \\
                                 & > & 0.
\end{eqnarray*}
If $r_3 > 0$, then by Lemma~\ref{lem: C_2 and A_2 in C_ell}(ii) the class has $A_2$ in its closure, whose dimension is $8\ell - 10$, which is greater than $4\ell + 4$. Thus we may assume $r_i = 0$ for $i \geq 3$. Write $y = r_2$, then $r_1 = 2\ell - 2y$; so $n_0 = 4\ell^2 - 4\ell y + 2y^2$ and $n_1 = 2\ell - 2y$, whence $\dim u^G = 2\ell y - y^2 + y$. Therefore
\begin{eqnarray*}
\dim u^G - 4\ell + 2n_1 - 2z - x & = & 2\ell y - y^2 + y - 4\ell + 4\ell - 4y - 2z - x \\
                                 & = & 2\ell y - y^2 - 3y - 2z - x;
\end{eqnarray*}
let the expression on the right hand side be $f(y)$, say. Then $f(y)$ is increasing for $y < \ell - \frac{3}{2}$ and decreasing for $y > \ell - \frac{3}{2}$. We have $f(2) = 4\ell - 10 - 2z - x > 0$ and $f(\ell - 1) = \ell^2 - 3\ell + 2 - 2z - x > 0$, while $f(1) = 2\ell - 4 - 2z - x$ and $f(\ell) = \ell^2 - 3\ell - 2z - x$, each of which is positive unless $\ell = 4$, $z = 1$ and $x = 2$. Thus we may assume $\ell = 4$ and $z = 1$, and need only consider $y = 1$ and $4$, corresponding to the classes $C_1$ and ${A_1}^2$; we must have $\z = 1$ and hence $d = 26$ while $d_1 = \frac{1}{2}n_0 - \frac{1}{2}n_1 - 1$. If $y = 1$ then $n_0 = 50$ and $n_1 = 6$, so $d_1 = 21$ and $d_2 = 5$; if instead $y = 4$ then $n_0 = 32$ and $n_1 = 0$, so $d_1 = 15$ and $d_2 = 11$. In both cases $d_1 \geq d_2 + 2$, so by Corollary~\ref{cor: B values for small k} we in fact have $x = 0$; so in all cases $f(y) > 0$.

Thus $B_{\d, 2} > \dim u^G$. Therefore if $p \geq 3$ and $k \in [2, \frac{d}{2}]$ the quadruple $(G, \lambda, p, k)$ satisfies $\udiamcon$.

We now assume instead $p = 2$. As in Section~\ref{sect: unipotent classes}, the unipotent classes in $G_{(2)}$ are
\begin{eqnarray*}
{a_{2y}}^G     & \hbox{ for } & y \in [1, {\ts\frac{\ell}{2}}], \\
{b_{2y + 1}}^G & \hbox{ for } & y \in [0, {\ts\frac{\ell - 1}{2}}], \\
{c_{2y + 2}}^G & \hbox{ for } & y \in [0, {\ts\frac{\ell - 2}{2}}];
\end{eqnarray*}
where ${a_{2y}}^G = {A_1}^y$ and ${b_{2y + 1}}^G = {A_1}^yC_1$, and we have
\begin{eqnarray*}
\dim {a_{2y}}^G     & = & 2y(2\ell - 2y), \\
\dim {b_{2y + 1}}^G & = & (2y + 1)(2\ell - 2y), \\
\dim {c_{2y + 2}}^G & = & (2y + 2)(2\ell - 2y - 1).
\end{eqnarray*}
We take each type of class in turn; we shall analyse the first two using an appropriate $\Psi$-net. Much as in the cases $G = E_7$, $\lambda = \omega_1$ and $G = E_8$, $\lambda = \omega_8$ in Section~\ref{sect: large quadruple individual cases}, we shall write `$2\bom_1/2\bom_3/\cdots/2\bom_{2y - 1}$' to denote a $\Psi$-net whose weights are those lying in the union of the Weyl $G_\Psi$-modules $W_{G_\Psi}(2\bom_1)$, $W_{G_\Psi}(2\bom_3)$, \dots, $W_{G_\Psi}(2\bom_{2y - 1})$.

First take $u \in {a_{2y}}^G$. Here we take $\Psi = \langle \alpha_1, \alpha_3, \dots, \alpha_{2y - 1} \rangle$ of type ${A_1}^y$; then we may assume $u = u_\Psi$. The $\Psi$-net table is as follows.
$$
\begin{array}{|*5{>{\ss}c|}}
\hline
\multicolumn{3}{|>{\ss}c|}{\Psi\mathrm{-nets}} & & \multicolumn{1}{|>{\ss}c|}{c(u_\Psi)} \\
\cline{1-3} \cline{5-5}
                \bar\nu                & n_0 & n_1 &                   m                   &       p = 2       \\
\hline
 2\bom_1/2\bom_3/\cdots/2\bom_{2y - 1} &  1  & 2y  &                   1                   &         y         \\
     \bom_{2i - 1} + \bom_{2j - 1}     &  0  &  4  &               2y(y - 1)               &     4y(y - 1)     \\
             \bom_{2i - 1}             &  0  &  2  &             4y(\ell - 2y)             &   4y(\ell - 2y)   \\
                   0                   &  0  &  1  & 2\ell^2 - 8\ell y + 8y^2 - 2\ell + 6y &                   \\
\hline
\multicolumn{4}{c|}{}                                                                      & y(4\ell - 4y - 3) \\
\cline{5-5}
\end{array}
$$
Thus $d - d_1 \geq y(4\ell - 4y - 3)$, so
\begin{eqnarray*}
(2d - 2d_1 - x) - \dim u^G & \geq & 2y(4\ell - 4y - 3) - x - 2y(2\ell - 2y) \\
                           &   =  & 2y(2\ell - 2y - 3) - x;
\end{eqnarray*}
let the expression on the right hand side be $f(y)$, say. Then $f(y)$ is increasing for $y < \frac{2\ell - 3}{4}$ and decreasing for $y > \frac{2\ell - 3}{4}$. We have $f(1) = 2(2\ell - 5) - x > 0$; if $y = \lfloor \frac{\ell}{2} \rfloor$, then according as $\ell$ is odd or even we have $\ell = 2y + 1$ or $2y$, so $f(y) = 2y(2y - 1) - x$ or $2y(2y - 3) - x$, each of which is positive. Therefore for all values of $y$ we have $2d - 2d_1 - x > \dim u^G$.

Now take $u \in {b_{2y + 1}}^G$. Here we take $\Psi = \langle \alpha_1, \alpha_3, \dots, \alpha_{2y - 1}, \alpha_\ell \rangle$ of type ${A_1}^y C_1$; then we may assume $u = u_\Psi$. The $\Psi$-net table is as follows.
$$
\begin{array}{|*5{>{\ss}c|}}
\hline
\multicolumn{3}{|>{\ss}c|}{\Psi\mathrm{-nets}} & & \multicolumn{1}{|>{\ss}c|}{c(u_\Psi)} \\
\cline{1-3} \cline{5-5}
                \bar\nu                & n_0 & n_1 &                      m                     &              p = 2              \\
\hline
 2\bom_1/2\bom_3/\cdots/2\bom_{2y - 1} &  1  & 2y  &                      1                     &                y                \\
     \bom_{2i - 1} + \bom_{2j - 1}     &  0  &  4  &                  2y(y - 1)                 &            4y(y - 1)            \\
             \bom_{2i - 1}             &  0  &  2  &              4y(\ell - 2y - 1)             &        4y(\ell - 2y - 1)        \\
       \bom_{2i - 1} + \bom_\ell       &  0  &  4  &                     2y                     &                4y               \\
               \bom_\ell               &  0  &  2  &              2(\ell - 2y - 1)              &         2(\ell - 2y - 1)        \\
                   0                   &  0  &  1  & 2\ell^2 - 8\ell y + 8y^2 - 6\ell + 14y + 4 &                                 \\
\hline
\multicolumn{4}{c|}{}                                                                           & 4\ell y - 4y^2 + 2\ell - 7y - 2 \\
\cline{5-5}
\end{array}
$$
Thus $d - d_1 \geq 4\ell y - 4y^2 + 2\ell - 7y - 2$, so
\begin{eqnarray*}
(2d - 2d_1 - x) - \dim u^G & \geq & 8\ell y - 8y^2 + 4\ell - 14y - 4 - x - (2y + 1)(2\ell - 2y) \\
                           &   =  & 4\ell y - 4y^2 + 2\ell - 12y - 4 - x;
\end{eqnarray*}
let the expression on the right hand side be $f(y)$, say. Then $f(y)$ is increasing for $y < \frac{\ell - 3}{2}$ and decreasing for $y > \frac{\ell - 3}{2}$. We have $f(0) = 2\ell - 4 - x > 0$; if $y = \lfloor \frac{\ell - 1}{2} \rfloor$, then according as $\ell$ is odd or even we have $\ell = 2y + 1$ or $2y + 2$, so $f(y) = 4y^2 - 4y - 2 - x$ or $4y^2 - x$, each of which is positive. Therefore for all values of $y$ we have $2d - 2d_1 - x > \dim u^G$.

Finally take $u \in {c_{2y + 2}}^G$. By Lemma~\ref{lem: a_{2y} in closure of c_{2y}} we have ${a_{2y + 2}}^G \leq {c_{2y + 2}}^G$, so we may use the bound obtained above to see that
\begin{eqnarray*}
(2d - 2d_1 - x) - \dim u^G & \geq & 2(y + 1)(4\ell - 4(y + 1) - 3) - x \\
                           &      & \quad {} - (2y + 2)(2\ell - 2y - 1) \\
                           &   =  & 2(y + 1)(2\ell - 2y - 6) - x;
\end{eqnarray*}
let the expression on the right hand side be $f(y)$, say. Then $f(y)$ is increasing for $y < \frac{\ell - 4}{2}$ and decreasing for $y > \frac{\ell - 4}{2}$. We have $f(0) = 4\ell - 12 - x > 0$; if $y = \lfloor \frac{\ell - 2}{2} \rfloor$, then according as $\ell$ is odd or even we have $\ell = 2y + 3$ or $2y + 2$, so $f(y) = 2(y + 1)(2y) - x$ or $2(y + 1)(2y - 2) - x$, each of which is positive with the exception of the latter when $y = 1$. Therefore for all values of $y$ we have $2d - 2d_1 - x > \dim u^G$, unless $(\ell, y) = (4, 1)$.

We are thus left to consider the unipotent class ${c_4}^G$ in $G = C_4$; then $\dim u^G = 20$ and $d = 26$. We may take $u = x_{\alpha_1}(1) x_{\alpha_4}(1) x_{2\alpha_3 + \alpha_4}(1)$; if as before we write the non-zero weights in $V$ as $\pm \ve_i \pm \ve_j$ for $1 \leq i < j \leq 4$, then $u = x_{\ve_1 - \ve_2}(1) x_{\ve_3}(1) x_{\ve_4}(1)$, and we have weight nets
\begin{eqnarray*}
& \{ \pm \ve_3 \pm \ve_4 \}, & \\
& \{ \ve_1 \pm \ve_3 \}, \{ -\ve_1 \pm \ve_3 \}, \{ \ve_2 \pm \ve_3 \}, \{ -\ve_2 \pm \ve_3 \}, & \\
& \{ \ve_1 \pm \ve_4 \}, \{ -\ve_1 \pm \ve_4 \}, \{ \ve_2 \pm \ve_4 \}, \{ -\ve_2 \pm \ve_4 \}, & \\
& \{ \ve_1 - \ve_2, 0, -\ve_1 + \ve_2 \}, & \\
& \{ \ve_1 + \ve_2 \}, \{ -\ve_1 - \ve_2 \}, &
\end{eqnarray*}
giving $d - d_1 \geq 2 + 8 + 1 = 11$. If $d - d_1 = 11$ then $d_1 = 15 > 13 = d_2 + 2$ so that $x = 0$; if not then $d - d_1 \geq 12$. In either case we have $2d - 2d_1 - x \geq 22 > \dim u^G$.

Thus $B_{\d, 2} > \dim u^G$. Therefore if $p = 2$ and $k \in [2, \frac{d}{2}]$ the quadruple $(G, \lambda, p, k)$ satisfies $\udiamcon$.

Now take $s \in G_{(r)}$ for $r \in \P'$. The weights are $\pm \ve_i \pm \ve_j$ for $1 \leq i < j \leq \ell$, and $0$. Thus the non-zero weights form a root system of type $D_\ell$.

First suppose $r \geq 3$. Again we set $\delta = \dim V_1(s) - \dim V_{\eta_r}(s)$. As in the proof of Lemma~\ref{lem: eigenspace dimension bound on Lie algebra}, there exists $\xi \in \{ \pm 1 \}$ such that for all $i$ there exists $j$ with $\ve_i(s) = \xi {\eta_r}^j$. For $j = 0, 1, \dots, r - 1$ set $m_j = | \{ i : \ve_i(s) = \xi {\eta_r}^j \}|$; then $\sum_{j = 0}^{r - 1}m_j = \ell$. We then have
\begin{eqnarray*}
\dim V_1(s)        & = & m_0(m_0 - 1) + m_1(m_1 - 1) + \cdots + m_{r - 1}(m_{r - 1} - 1) \\
                   &   & {} + m_0(m_0 - 1) + m_1m_{r - 1} + m_2m_{r - 2} + \cdots + m_{r - 1}m_1 + \ell - 1 - \z, \\
\dim V_{\eta_r}(s) & = & m_0m_1 + m_1m_2 + \cdots + m_{r - 2}m_{r - 1} + m_{r - 1}m_0 \\
                   &   & {} + m_0m_1 + m_{r - 1}m_2 + \cdots + m_{\frac{r + 3}{2}}m_{\frac{r - 1}{2}} + {\ts{\frac{1}{2}}}m_{\frac{r + 1}{2}}(m_{\frac{r + 1}{2}} - 1) \\
                   &   & {} + m_0m_{r - 1} + m_1m_{r - 2} + \cdots + m_{\frac{r - 3}{2}}m_{\frac{r + 1}{2}} + {\ts{\frac{1}{2}}}m_{\frac{r - 1}{2}}(m_{\frac{r - 1}{2}} - 1).
\end{eqnarray*}
Thus
\begin{eqnarray*}
\delta & = & {\ts{\frac{1}{2}}}\left[\left((2m_0 - m_1 - m_{r - 1} - {\ts{\frac{1}{2}}})^2 - {\ts{\frac{1}{4}}}\right) + \left((m_1 - m_2 - m_{r - 2} + m_{r - 1} - {\ts{\frac{1}{2}}})^2 - {\ts{\frac{1}{4}}}\right) \right. \\
       &   & \quad {} + \cdots + \left.\left((m_{\frac{r - 3}{2}} - m_{\frac{r - 1}{2}} - m_{\frac{r + 1}{2}} + m_{\frac{r + 3}{2}} - {\ts{\frac{1}{2}}})^2 - {\ts{\frac{1}{4}}}\right)\right] - 1 - \z \\
       & \geq & -1 - \z.
\end{eqnarray*}
Hence $\dim V_{\eta_r}(s) \leq \dim V_1(s) + 1 + \z$. Also $|\Phi(s)| = \dim V_1(s) + 2m_0 - (\ell - 1 - \z)$, so $\dim s^G = |\Phi| - |\Phi(s)| = 2\ell^2 - \dim V_1(s) - 2m_0 + (\ell - 1 - \z)$. Thus $d_1 + \dim s^G \leq 2\ell^2 + \ell - 2m_0$, and so
\begin{eqnarray*}
(2d - 2d_1 - x) - 2\dim s^G & \geq & 2(2\ell^2 - \ell - 1 - \z - (2\ell^2 + \ell - 2m_0)) - x \\
                            &   =  & 2(2m_0 - 2\ell - 1 - \z) - x.
\end{eqnarray*}
Thus $(2d - 2d_1 - x) - \dim s^G \geq \dim s^G + 4m_0 - 4\ell - 2 - 2\z - x \geq \dim s^G - 4\ell - 6$. Thus $2d - 2d_1 - x > \dim s^G$ provided $\dim s^G > 4\ell + 6$; so we must consider the semisimple classes of dimension at most $4\ell + 6$, i.e., the classes $s^G$ with $|\Phi(s)| \geq |\Phi| - 4\ell - 6 = 2\ell^2 - 4\ell - 6$.

Since $r \geq 3$, the subsystem $\Phi(s)$ can have at most one simple factor of type $C$. Suppose $\Phi(s)$ has a factor $C_{\ell - y}$ for $1 \leq y \leq \ell$, then $\Phi(s) \subseteq A_{y - 1}C_{\ell - y}$, and so
\begin{eqnarray*}
(2\ell^2 - 4\ell - 6) - |\Phi(s)| & \geq & (2\ell^2 - 4\ell - 6) - (y(y - 1) + 2(\ell - y)^2) \\
                                  &   =  & 4\ell y - 4\ell - 3y^2 + y - 6;
\end{eqnarray*}
let the expression on the right hand side be $f(y)$, say. Then $f(y)$ is increasing for $y < \frac{4\ell + 1}{6}$ and decreasing for $y > \frac{4\ell + 1}{6}$. We have $f(3) = 8\ell - 30 > 0$ and $f(\ell - 1) = \ell^2 - \ell - 10 > 0$, while $f(2) = 4\ell - 16$ and $f(\ell) = \ell^2 - 3\ell - 6$, each of which is positive unless $\ell = 4$, and $f(1) = -8 < 0$. Thus we may assume either $\Phi(s) = C_{\ell - 1}$, or $\ell = 4$ and $\Phi(s) \subseteq A_1C_2$ or $A_3$.

If $\Phi(s) = C_{\ell - 1}$ then $m_0 = \ell - 1$ and $\dim s^G = 2\ell^2 - 2(\ell - 1)^2 = 4\ell - 2$; so
\begin{eqnarray*}
(2d - 2d_1 - x) - \dim s^G & \geq & \dim s^G + 4m_0 - 4\ell - 2 - 2\z - x  \\
                           &   =  & 4\ell - 2 + 4\ell - 4 - 4\ell - 2 - 2\z - x \\
                           &   =  & 4\ell - 8 - 2\z - x \\
                           &   >  & 0.
\end{eqnarray*}
Thus we may assume $\ell = 4$; so we need only consider semisimple classes $s^G$ with $|\Phi(s)| \geq 32 - 16 - 6 = 10$. Thus if $\Phi(s) \subseteq A_1C_2$ we need only consider $\Phi(s) = A_1C_2$; then $\dim s^G = 32 - 10 = 22$. We may assume $\ve_1(s) = \ve_2(s) = \xi$, $\ve_3(s) = \ve_4(s) = \xi\eta_r$, so $m_0 = 2$, and
\begin{eqnarray*}
(2d - 2d_1 - x) - \dim s^G & \geq & \dim s^G + 4m_0 - 4\ell - 2 - 2\z - x \\
                           &   =  & 22 + 8 - 16 - 2 - 2\z - x \\
                           &   =  & 12 - 2\z - x \\
                           &   >  & 0.
\end{eqnarray*}
If instead $\Phi(s) \subseteq A_3$ we need only consider $\Phi(s) = A_3$; then $\dim s^G = 32 - 12 = 20$. We may assume $\ve_1(s) = \ve_2(s) = \ve_3(s) = \ve_4(s) = \xi\eta_r$; here we have $\dim V_1(s) = 15 - \z$ while $\dim V_{{\eta_r}^2}(s) = \dim V_{{\eta_r}^{-2}}(s) = 6$, so
$$
(2d - 2d_1 - x) - \dim s^G = 2(27 - \z) - 2(15 - \z) - x - 20 = 4 - x > 0.
$$
Thus if $r \geq 3$ we have $2d - 2d_1 - x > \dim s^G$.

Now suppose $r = 2$. Here either all $\ve_i(s)$ are $\pm 1$, or all $\ve_i(s)$ are $\pm \eta_4$.

If all $\ve_i(s)$ are $\pm 1$, set $m_0 = |\{ i : \ve_i(s) = 1 \}|$ and $m_1 = |\{ i : \ve_i(s) = -1 \}|$; then $m_0 + m_1 = \ell$. We then have
\begin{eqnarray*}
\dim V_1(s)    & = & 2m_0(m_0 - 1) + 2m_1(m_1 - 1) + \ell - 1 - \z, \\
\dim V_{-1}(s) & = & 4m_0m_1.
\end{eqnarray*}
Thus
$$
\delta = 2(m_0 - m_1)^2 - \ell - 1 - \z;
$$
so $d_1$ may be $\dim V_1(s)$ or $\dim V_{-1}(s)$. Also $|\Phi(s)| = 2{m_0}^2 + 2{m_1}^2$, so $\dim s^G = |\Phi| - |\Phi(s)| = 4m_0m_1$. If $d_1 = \dim V_1(s)$ and $d_2 = \dim V_{-1}(s)$ then
\begin{eqnarray*}
(2d - 2d_1 - x) - \dim s^G & = & 2d_2 - x - \dim s^G \\
                           & = & 8m_0m_1 - x - 4m_0m_1 \\
                           & = & 4m_0m_1 - x \\
                           & > & 0.
\end{eqnarray*}
If on the other hand $d_1 = \dim V_{-1}(s)$ and $d_2 = \dim V_1(s)$ then
\begin{eqnarray*}
(2d - 2d_1 - x) - \dim s^G & = & 2d_2 - x - \dim s^G \\
                           & = & 4{m_0}^2 + 4{m_1}^2 - 2\ell - 2 - 2\z - x - 4m_0m_1 \\
                           & = & 3(m_0 - m_1)^2 + \ell^2 - 2\ell - 2 - 2\z - x \\
                           & > & 0.
\end{eqnarray*}
Hence $2d - 2d_1 - x > \dim s^G$.

If instead all $\ve_i(s)$ are $\pm \eta_4$, set $m_0 = |\{ i : \ve_i(s) = \eta_4 \}|$ and $m_1 = |\{ i : \ve_i(s) = -\eta_4 \}|$; then $m_0 + m_1 = \ell$. We then have
\begin{eqnarray*}
\dim V_1(s)    & = & m_0(m_0 - 1) + m_1(m_1 - 1) + 2m_0m_1 + \ell - 1 - \z = \ell^2 - 1 - \z, \\
\dim V_{-1}(s) & = & 2m_0m_1 + m_0(m_0 - 1) + m_1(m_1 - 1) = \ell^2 - \ell.
\end{eqnarray*}
Thus
$$
\delta = \ell - 1 - \z > 0.
$$
Hence $d_1 = \dim V_1(s)$ and $d_2 = \dim V_{-1}(s)$. Also $|\Phi(s)| = \ell(\ell - 1)$, so $\dim s^G = |\Phi| - |\Phi(s)| = \ell(\ell + 1)$. Therefore
\begin{eqnarray*}
(2d - 2d_1 - x) - \dim s^G & = & 2d_2 - x - \dim s^G \\
                           & = & 2\ell^2 - 2\ell - x - (\ell^2 + \ell) \\
                           & = & \ell^2 - 3\ell - x \\
                           & > & 0,
\end{eqnarray*}
and so $2d - 2d_1 - x > \dim s^G$.

Thus $B_{\d, 2} > \dim s^G$. Therefore if $k \in [2, \frac{d}{2}]$ the quadruple $(G, \lambda, p, k)$ satisfies $\ssdiamcon$.
\end{proof}

\begin{prop}\label{prop: B_ell, omega_2, p = 2, k geq 2}
Let $G = B_\ell$ for $\ell \in [4, \infty)$ and $\lambda = \omega_2$ with $p = 2$; then for $k \in [2, \frac{d}{2}]$ the quadruple $(G, \lambda, p, k)$ satisfies $\ssdiamcon$ and $\udiamcon$.
\end{prop}

\begin{proof}
This is an immediate consequence of Proposition~\ref{prop: C_ell, omega_2, k geq 2}, using the exceptional isogeny $B_\ell \to C_\ell$ which exists in characteristic $2$.
\end{proof}

Finally we take the two infinite families where $k_0 = 3$; note that Corollary~\ref{cor: B values for small k} allows us to write $B_{\d, 3} = 3d - 3d_1 - x$ where
$$
x =
\begin{cases}
0                    & \hbox{if } d_1 \geq d_2 + 4, \\
d_2 + 4 - d_1        & \hbox{if } d_2 + 4 > d_1  \geq d_3 + 2, \\
d_2 + d_3 + 6 - 2d_1 & \hbox{if } d_1 < d_3 + 2,
\end{cases}
$$
so that $x \in [0, 6]$.

\begin{prop}\label{prop: A_ell, omega_2 and 2omega_1, k geq 3}
Let $G = A_\ell$ for $\ell \in [6, \infty)$ and $\lambda = \omega_2$, or $G = A_\ell$ for $\ell \in [3, \infty)$ and $\lambda = 2\omega_1$ with $p \geq 3$; then for $k \in [3, \frac{d}{2}]$ the quadruple $(G, \lambda, p, k)$ satisfies $\ssdiamcon$ and $\udiamcon$.
\end{prop}

\begin{proof}
Number the cases (i) and (ii) according as $\lambda = \omega_2$ or $2\omega_1$; in what follows, whenever we give two choices followed by the word `respectively', we are taking the cases in the order (i), (ii).

First take $u \in G_{(p)}$. We may assume $u = u_\Psi$ for $\Psi$ of type $A_{m_1 - 1} A_{m_2 - 1} \dots$, where $\sum_i m_i = \ell + 1$ and $p \geq m_1 \geq m_2 \geq \cdots$; we have
$$
\dim u^G = (\ell + 1)(\ell + 2) - 2 \sum_i im_i.
$$
For each $i$ write $l_i = m_1 + \cdots + m_{i - 1}$; then we may assume the simple roots of the $i$th factor $A_{m_i - 1}$ are $\alpha_{l_i + 1}, \alpha_{l_i + 2}, \dots, \alpha_{l_i + m_i - 1}$. Each $\Psi$-net then corresponds to a weight $\bar\nu$ which is of one of two forms: $\bom_{l_i + 1} + \bom_{l_j + 1}$ with $i < j$, and $\bom_{l_i + 2}$ (where $i$ is such that $m_i \geq 2$) or $2\bom_{l_i + 1}$ respectively. If $\bar\nu$ is of the first form, the weights in the $\Psi$-net are $\bar\mu_i + \bar\mu_j$ where $\bar\mu_i \in W(A_{m_i - 1}).\bom_{l_i + 1}$ and $\bar\mu_j \in W(A_{m_j - 1}).\bom_{l_j + 1}$; for any fixed $\bar\mu_j$, the sum of the weight spaces corresponding to the weights $\bar\mu_i + \bar\mu_j$ is a natural module for $A_{m_i - 1}$ on which $u_\Psi$ acts regularly, so we may take the contribution $c(u_\Psi)$ from the $\Psi$-net to be $(m_i - 1)m_j$. If $\bar\nu$ is of the second form, the sum of the weight spaces corresponding to the weights in the $\Psi$-net is a module $L(\omega_2)$ or $L(2\omega_1)$ respectively for $A_{m_i - 1}$ on which $u_\Psi$ acts regularly, so by Lemma~\ref{lem: fixed points on L(omega_2) and L(2omega_1)} we may take the contribution $c(u_\Psi)$ from the $\Psi$-net to be $\frac{1}{2}m_i(m_i - 1) - \lfloor \frac{m_i}{2} \rfloor$ or $\frac{1}{2}m_i(m_i + 1) - \lceil \frac{m_i}{2} \rceil$ respectively. Thus
\begin{eqnarray*}
d - d_1 & \geq & \sum_{i < j} (m_i - 1)m_j +
\begin{cases}
\sum_i ({\ts{\frac{1}{2}}}m_i(m_i - 1) - \lfloor {\ts{\frac{m_i}{2}}} \rfloor) & \hbox{in case~(i),}  \\
\sum_i ({\ts{\frac{1}{2}}}m_i(m_i + 1) - \lceil {\ts{\frac{m_i}{2}}} \rceil)   & \hbox{in case~(ii),}
\end{cases} \\
        & = & \sum_{i < j} m_i m_j - \sum_{i < j} m_j + {\ts{\frac{1}{2}}} \sum_i {m_i}^2 +
\begin{cases}
\sum_i (-{\ts{\frac{m_i}{2}}} - \lfloor {\ts{\frac{m_i}{2}}} \rfloor) & \hbox{in case~(i),}  \\
\sum_i ({\ts{\frac{m_i}{2}}} - \lceil {\ts{\frac{m_i}{2}}} \rceil)    & \hbox{in case~(ii).}
\end{cases}
\end{eqnarray*}
Observe that $2 \sum_{i < j} m_im_j = \sum_{i \neq j}m_im_j = (\sum_i m_i)^2 - \sum_i {m_i}^2 = (\ell + 1)^2 - \sum_i {m_i}^2$, so $\sum_{i < j} m_im_j = \frac{1}{2}(\ell + 1)^2 - \frac{1}{2} \sum_i {m_i}^2$; also $\sum_{i < j} m_j = \sum_j (j - 1)m_j$. Thus if we write $v$ for the number of odd $m_i$, we have
\begin{eqnarray*}
d - d_1 & \geq & {\ts{\frac{1}{2}}}(\ell + 1)^2 - \sum_i (i - 1)m_i +
\begin{cases}
\sum_i (-{\ts{\frac{m_i}{2}}} - \lfloor {\ts{\frac{m_i}{2}}} \rfloor) & \hbox{in case~(i),}  \\
\sum_i ({\ts{\frac{m_i}{2}}} - \lceil {\ts{\frac{m_i}{2}}} \rceil)    & \hbox{in case~(ii),}
\end{cases} \\
        & = &
\begin{cases}
{\ts{\frac{1}{2}}}(\ell + 1)^2         - \sum_i im_i + {\ts{\frac{1}{2}}} v & \hbox{in case~(i),}  \\
{\ts{\frac{1}{2}}}(\ell + 1)(\ell + 3) - \sum_i im_i - {\ts{\frac{1}{2}}} v & \hbox{in case~(ii).}
\end{cases}
\end{eqnarray*}
Hence
$$
(3d - 3d_1 - x) - \dim u^G \geq
\begin{cases}
{\ts{\frac{1}{2}}}(\ell + 1)(\ell - 1) - \sum_i im_i + {\ts{\frac{3}{2}}} v - x & \hbox{in case~(i),} \\
{\ts{\frac{1}{2}}}(\ell + 1)(\ell + 5) - \sum_i im_i - {\ts{\frac{3}{2}}} v - x & \hbox{in case~(ii).}
\end{cases}
$$

First assume we are in case~(ii). If $\Psi$ is of type $A_1$ then $m_1 = 2$, $m_2 = \cdots = m_\ell = 1$, so $\sum_i im_i = \frac{1}{2}\ell^2 + \frac{1}{2}\ell + 1$ while $v = \ell - 1$; thus we have $(3d - 3d_1 - x) - \dim u^G \geq \ell + 3 - x$. Since $x \leq 6$ this is positive for $\ell \in [4, \infty)$. For $\ell = 3$ we have $d = 10$ and $d - d_1 \geq 4$, so we may take $\d_0 = (6, 4)$; using Proposition~\ref{prop: B value when t = 2} we then have $B_{\d_0, 3} = 10 > 6 = \dim u^G$. For any other $\Psi$, the class ${u_\Psi}^G$ contains $A_1$ in its closure by Lemma~\ref{lem: root elt class in closure of any non-triv class}, so $\sum_i im_i < \frac{1}{2} \ell^2 + \frac{1}{2} \ell + 1$, and $v \leq \ell - 1$ as $v \equiv \ell + 1$ (mod $2$); thus we have $(3d - 3d_1 - x) - \dim u^G > \ell + 3 - x \geq 0$ for $\ell \in [3, \infty)$. Thus in all cases here $B_{\d, 3} > \dim u^G$. Therefore in case~(ii) if $k \in [3, \frac{d}{2}]$ the quadruple $(G, \lambda, p, k)$ satisfies $\udiamcon$.

Now assume we are in case~(i). We begin with some small subsystems $\Psi$. If $\Psi$ is of type $A_1$ then as above $\sum_i im_i = \frac{1}{2}\ell^2 + \frac{1}{2}\ell + 1$ and $v = \ell - 1$, so $d - d_1 \geq \ell - 1$; as $\ell - 1 \leq \frac{1}{2}d$ we may take $\d_0 = (d - (\ell - 1), \ell - 1)$, and then Proposition~\ref{prop: B value when t = 2} shows that we have $B_{\d_0, 3} = 3(\ell - 1) > 2\ell = \dim u^G$. If $\Psi$ is of type ${A_1}^2$ then $m_1 = m_2 = 2$, $m_3 = \cdots = m_{\ell - 1} = 1$, whence $\sum_i im_i = \frac{1}{2}\ell^2 - \frac{1}{2}\ell + 3$ and $v = \ell - 3$, so $d - d_1 \geq 2\ell - 4$; as $2\ell - 4 \leq \frac{1}{2}d$ we may take $\d_0 = (d - (2\ell - 4), 2\ell - 4)$, and then Proposition~\ref{prop: B value when t = 2} shows that we have $B_{\d_0, 3} = 3(2\ell - 4) > 4\ell - 4 = \dim u^G$. If $\Psi$ is of type $A_2$ then $u^G$ has dimension $4\ell - 2$ and contains the class ${A_1}^2$ in its closure; we have $3(2\ell - 4) > 4\ell - 2$, so the bound for the class ${A_1}^2$ suffices. If $\Psi$ is of type ${A_1}^3$ then $m_1 = m_2 = m_3 = 2$, $m_4 = \cdots = m_{\ell - 2} = 1$, whence $\sum_i im_i = \frac{1}{2}\ell^2 - \frac{3}{2}\ell + 7$ and $v = \ell - 5$, so $d - d_1 \geq 3\ell - 9$; as $3\ell - 9 \leq \frac{1}{2}d$ we may take $\d_0 = (d - (3\ell - 9), 3\ell - 9)$, and then Proposition~\ref{prop: B value when t = 2} shows that for $\ell \in [7, \infty)$ we have $B_{\d_0, 3} = 3(3\ell - 9) > 6\ell - 12 = \dim u^G$, while for $\ell = 6$ we have $B_{\d_0, 3} = 26 > 24 = \dim u^G$. If $\Psi$ is of type $A_2A_1$ then $u^G$ has dimension $6\ell - 8$ and contains the class ${A_1}^3$ in its closure. For $\ell \in [7, \infty)$ we have $3(3\ell - 9) > 6\ell - 8$, so the bound for the class ${A_1}^3$ suffices. For $\ell = 6$ we have $m_1 = 3$, $m_2 = 2$, $m_3 = m_4 = 1$, whence $\sum_i im_i = 14$ and $v = 3$, so $d - d_1 \geq 12$; we may take $\d_0 = (9, 9, 3)$, and then Corollary~\ref{cor: B values for small k} shows that we have $B_{\d_0, 3} = 32 > 28 = \dim u^G$. If $\Psi$ is of type $A_3$ then $u^G$ has dimension $6\ell - 6$ and again contains the class ${A_1}^3$ in its closure. For $\ell \in [8, \infty)$ we have $3(3\ell - 9) > 6\ell - 6$, so the bound for the class ${A_1}^3$ suffices. For $\ell \in [6, 7]$ we have $m_1 = 4$, $m_2 = \cdots = m_{\ell - 2} = 1$, whence $\sum_i im_i = \frac{1}{2}\ell^2 - \frac{3}{2}\ell + 4$ and $v = \ell - 3$, so $d - d_1 \geq 3\ell - 5$; according as $\ell = 6$ or $7$ we may take $\d_0 = (8, 8, 5)$ or $(12, 12, 4)$, and then Corollary~\ref{cor: B values for small k} shows that we have $B_{\d_0, 3} = 35 > 30 = \dim u^G$ or $B_{\d_0, 3} = 44 > 36 = \dim u^G$.

We now treat the remaining subsystems $\Psi$. First suppose $\ell = 6$. If $\Psi$ is of type $A_2{A_1}^2$ then $m_1 = 3$, $m_2 = m_3 = 2$, whence $\sum_i im_i = 13$ and $v = 1$, so $d - d_1 \geq 12$; we may take $\d_0 = (9, 9, 3)$, and then Corollary~\ref{cor: B values for small k} shows that we have $B_{\d_0, 3} = 32 > 30 = \dim u^G$. All remaining classes in $G_{(p)}$ (assuming $p \geq 3$) have ${A_2}^2$ in their closure, for which $m_1 = m_2 = 3$, $m_3 = 1$, whence $\sum_i im_i = 12$; noting that $v$ is odd we have $(3d - 3d_1 - x) - \dim u^G \geq \frac{35}{2} - 12 + \frac{3}{2}v - x \geq 7 - x > 0$. Now instead suppose $\ell \in [7, \infty)$. Here all remaining classes in $G_{(p)}$ have ${A_1}^4$ in their closure, for which $m_1 = \cdots = m_4 = 2$, $m_5 = \cdots = m_{\ell - 3} = 1$, whence $\sum_i im_i = \frac{1}{2}\ell^2 - \frac{5}{2}\ell + 13$; thus we have $(3d - 3d_1 - x) - \dim u^G \geq \frac{5}{2}\ell - \frac{27}{2} + \frac{3}{2}v - x$, which is positive unless $\ell = 7$ and $v = 0$. Thus we may suppose $\ell = 7$, in which case $d = 28$. If $\Psi$ is of type ${A_1}^4$ then as we have seen $\sum_i im_i = 20$, so $d - d_1 \geq 12$; we may take $\d_0 = (16, 12)$, and then Proposition~\ref{prop: B value when t = 2} shows that we have $B_{\d_0, 3} = 36 > 32 = \dim u^G$. The other classes in $G_{(p)}$ with $v = 0$ (assuming $p$ is sufficiently large for each) are $A_3{A_1}^2$, ${A_3}^2$, $A_5A_1$ and $A_7$, each of which has $A_3{A_1}^2$ in its closure, for which $m_1 = 4$, $m_2 = m_3 = 2$, whence $\sum_i im_i = 14$; thus we have $(3d - 3d_1 - x) - \dim u^G \geq 10 - x > 0$. Thus in all cases here $B_{\d, 3} > \dim u^G$. Therefore in case~(i) if $k \in [3, \frac{d}{2}]$ the quadruple $(G, \lambda, p, k)$ satisfies $\udiamcon$.

Now take $s \in G_{(r)}$ for $r \in \P'$. Write $\Phi(s) = A_{a_1 - 1} A_{a_2 - 1} \dots$, where $\sum_i a_i = \ell + 1$ and $a_1 \geq a_2 \geq \cdots$; then $(a_1, a_2, \dots)$ is a partition of $\ell + 1$. Let $(m_1, m_2, \dots)$ be the partition of $\ell + 1$ dual to $(a_1, a_2, \dots)$, where $m_1 \geq m_2 \geq \cdots$; take $\Psi$ of type $A_{m_1 - 1} A_{m_2 - 1} \dots$. For each $i$ write $l_i = m_1 + \cdots + m_{i - 1}$, and let $r_i$ be the number of $m_j$ equal to $i$, so that $\sum_i ir_i = \ell + 1$. Then for each $i$ we have $a_i = r_i + r_{i + 1} + \cdots$; moreover $u_\Psi$ has Jordan block sizes $1^{r_1}, 2^{r_2}, \dots$, and as $\sum_i (2i - 1)m_i = \sum_i (r_i + r_{i + 1} + \cdots)^2 = \sum_i {a_i}^2$ we have $\dim {u_\Psi}^G = (\ell + 1)^2 - \sum_i {a_i}^2$. Thus
$$
\dim s^G = |\Phi| - |\Phi(s)| = \ell(\ell + 1) - \sum_i a_i(a_i - 1) = \ell(\ell + 1) - \sum_i {a_i}^2 + \sum_i a_i = \dim {u_\Psi}^G.
$$
Moreover, if we take the Young tableau whose $i$th row has length $a_i$, and number its boxes from $1$ to $\ell + 1$ in order, working down the first column, then down the second and so on, then the roots $\ve_i - \ve_j$ where $i$ and $j$ lie in the same row form a subsystem of type $\Phi(s)$, while those where $i$ and $j$ lie in the same column form a subsystem of type $\Psi$, and the two are clearly disjoint; so we may assume $\Psi$ is disjoint from $\Phi(s)$. Moreover the simple roots of the $i$th factor $A_{m_i - 1}$ of $\Psi$ are $\alpha_{l_i + 1}, \alpha_{l_i + 2}, \dots, \alpha_{l_i + m_i - 1}$; thus $\Psi$ is as given at the beginning of this proof.

Now consider the $\Psi$-nets; recall that these correspond to weights $\bar\nu$ which are of one of two forms. If $\bar\nu = \bom_{l_i + 1} + \bom_{l_j + 1}$ with $i < j$, as stated above the weights in the $\Psi$-net are $\bar\mu_i + \bar\mu_j$ where $\bar\mu_i \in W(A_{m_i - 1}).\bom_{l_i + 1}$ and $\bar\mu_j \in W(A_{m_j - 1}).\bom_{l_j + 1}$; for any fixed $\bar\mu_j$, the weights $\bar\mu_i + \bar\mu_j$ are such that any two differ by a root in $\Phi(s)$, so that at most one can lie in any given eigenspace for $s$. We may therefore take the contribution $c(s)$ from the $\Psi$-net to be $(m_i - 1)m_j$, which is what we took $c(u_\Psi)$ to be. The other possibility is that $\bar\nu = \bom_{l_i + 2}$ (with $m_i \geq 2$) or $2\bom_{l_i + 1}$ respectively; the root system of the relevant factor consists of roots $\pm(\ve_{j_1} - \ve_{j_2})$ for $l_i + 1 \leq j_1 < j_2 \leq l_i + m_i$. If we are in case~(i), the weights in the $\Psi$-net are $\ve_{j_1} + \ve_{j_2}$ where $l_i + 1 \leq j_1 < j_2 \leq l_i + m_i$. In any given eigenspace we then cannot have two weights $\ve_{j_1} + \ve_{j_2}$ and $\ve_{{j_1}'} + \ve_{{j_2}'}$ for which the sets $\{ j_1, j_2 \}$ and $\{ {j_1}', {j_2}' \}$ have non-empty intersection; so at most $\lfloor \frac{m_i}{2} \rfloor$ weights can lie in an eigenspace, and hence we may take $c(s)$ to be $\frac{1}{2}m_i(m_i - 1) - \lfloor \frac{m_i}{2} \rfloor$, which is again what we took $c(u_\Psi)$ to be. If instead we are in case~(ii), the weights in the $\Psi$-net are $2\ve_j$ where $l_i + 1 \leq j \leq l_i + m_i$ and $\ve_{j_1} + \ve_{j_2}$ where $l_i + 1 \leq j_1 < j_2 \leq l_i + m_i$. First suppose $r \geq 3$. In any given eigenspace, if no weights $2\ve_j$ are present then as in case~(i) at most $\lfloor \frac{m_i}{2} \rfloor$ weights can be present; if instead some weight $2\ve_j$ is present, then no other weight $2\ve_{j'}$ or $\ve_{j_1} + \ve_{j_2}$ with $j \in \{ j_1, j_2 \}$ can be, so that at most $1 + \lfloor \frac{m_i - 1}{2} \rfloor = \lceil \frac{m_i}{2} \rceil$ weights can lie in the eigenspace. Hence we may take $c(s)$ to be $\frac{1}{2}m_i(m_i + 1) - \lceil \frac{m_i}{2} \rceil$, which is once more what we took $c(u_\Psi)$ to be. If however $r = 2$, then $\Phi(s)$ has at most two simple factors, so each simple factor of $\Psi$ has rank at most one, and hence $\Psi$ is of type ${A_1}^y$ for some $y$; here if $m_i = 1$ the $\Psi$-net contains just the one weight $2\ve_{l_i + 1}$ so that $c(s)$ may be $0$, while if $m_i = 2$ it contains the three weights $2\ve_{l_i + 1}$, $\ve_{l_i + 1} + \ve_{l_i + 2}$ and $2\ve_{l_i + 2}$, of which the first and third lie in the same eigenspace so that we may take $c(s) = 1$.

Thus if we are in case~(i), or case~(ii) with $r \geq 3$, the sum of the contributions $c(s)$ is the same as the sum of the contributions $c(u_\Psi)$. If instead we are in case~(ii) with $r = 2$, and $\Psi$ is of type ${A_1}^y$, we obtain
$$
d - d_1 \geq \sum_{i < j} (m_i - 1)m_j + y = (\ell - 1) + (\ell - 3) + \cdots + (\ell - (2y - 1)) + y = y(\ell + 1 - y),
$$
while $\dim s^G = 2y(\ell + 1 - y)$; if we write $c = y(\ell + 1 - y)$, then as $c \leq \frac{d}{2}$ we may take $\d_0 = (d - c, c)$, and then Corollary~\ref{cor: B values for small k} shows that according as $d \geq 2c + 4$ or $d < 2c + 4$ we have $B_{\d_0, 3} = 3c > 2c = \dim s^G$ or $B_{\d_0, 3} = 3d - 2(d - c) - c - 4 = 2c + (d - c - 4) \geq 2c + (\frac{d}{2} - 4) > 2c = \dim s^G$.

Thus in all cases $B_{\d, 3} > \dim s^G$. Therefore if $k \in [3, \frac{d}{2}]$ the quadruple $(G, \lambda, p, k)$ satisfies $\ssdiamcon$.
\end{proof}

This completes the treatment of the infinite families of cases listed in Table~\ref{table: remaining large quadruples}.

\chapter{Methods for treating cases not having TGS}\label{chap: non-TGS methods}

In this chapter we develop methods for treating cases which do not have trivial generic stabilizers. The structure of this chapter is as follows. In Section~\ref{sect: localization} we explain the key idea of localization to a subvariety. In Section~\ref{sect: semisimple auts} we consider a particular type of situation where the group and the module respectively occur within a larger simple algebraic group and its Lie algebra. In Section~\ref{sect: Lie algebra annihilators} we detail an approach involving annihilators in the Lie algebra. In Section~\ref{sect: invariants} we treat a very specific set-up which allows us to make use of an invariant. In Section~\ref{sect: gen ht fns} we describe a technique which greatly simplifies the determination of stabilizers (and more generally transporters) provided certain conditions hold. Finally in Section~\ref{sect: reduction} we give a result which links actions on higher Grassmannian varieties to those on projective spaces.

\section{Localization to a subvariety}\label{sect: localization}

This section concerns a basic approach which will be used in almost all cases where a triple $(G, \lambda, p)$ or quadruple $(G, \lambda, p, k)$ does not have TGS.

Let $X$ be a variety on which $G$ acts. Given a subvariety $Y$ of $X$, the morphism $\phi: G \times Y \to X$ defined by $\phi((g, y)) = g.y$ is known as the {\em orbit map\/}; clearly its image is the union of those $G$-orbits in $X$ which meet $Y$. Our first result here concerns dimensions of transporters.

\begin{lem}\label{lem: transporter dimension}
With the notation established, given $y \in Y$ we have
\begin{itemize}
\item[(i)] $\dim \Tran_G(y, Y) = \dim \phi^{-1}(y)$;
\item[(ii)] $\codim \Tran_G(y, Y) = \dim(\overline{G.y}) - \dim(\overline{G.y \cap Y})$.
\end{itemize}
\end{lem}

\begin{proof}
We have
\begin{eqnarray*}
\phi^{-1}(y) &   =   & \{ (g, y') : g \in G, \ y' \in Y, \ g.y' = y \} \\
             &   =   & \{ (g, g^{-1}.y) : g \in G, \ g^{-1}.y \in Y \} \\
             & \cong & \{ g^{-1} : g \in G, \ g^{-1}.y \in Y \} \\
             &   =   & {\textstyle\Tran_G(y, Y)};
\end{eqnarray*}
this proves (i). The fibre $\phi^{-1}(y)$ is closed in $G \times Y$; let $\pi_2 : \phi^{-1}(y) \to Y$ be the projection on the second component. Then
$$
\im \pi_2 = \{ y' \in Y : \exists g \in G \hbox{ with } g.y' = y \} = G.y \cap Y,
$$
so that the morphism $\pi_2 : \phi^{-1}(y) \to \overline{G.y \cap Y}$ is dominant; for each $y' \in \im \pi_2$ we have ${\pi_2}^{-1}(y') = \{ (g, y') : g \in G,\ g.y' = y \} \cong \{ g \in G : g.y' = y \}$ which is a coset of $C_G(y)$, so all fibres of $\pi_2$ have dimension equal to $\dim C_G(y)$. Thus by Lemma~\ref{lem: equal dimension fibres} we have $\dim \phi^{-1}(y) = \dim(\overline{G.y \cap Y}) + \dim C_G(y)$. Since $\dim(\overline{G.y}) = \dim G - \dim C_G(y)$, using (i) we have
\begin{eqnarray*}
\codim {\textstyle\Tran_G(y, Y)} & = & \dim G - \dim {\textstyle\Tran_G(y, Y)} \\
                                 & = & (\dim(\overline{G.y}) + \dim C_G(y)) - \dim \phi^{-1}(y) \\
                                 & = & \dim(\overline{G.y}) - \dim(\overline{G.y \cap Y});
\end{eqnarray*}
this proves (ii).
\end{proof}

We shall be interested in subvarieties $Y$ which are \lq sufficiently representative', in the sense that almost all orbits in $X$ meet them, with the intersections having the appropriate dimensions. In order to give a condition for this, we make the following definition.

\begin{definition}
Given a subvariety $Y$ of $X$, a point $y \in Y$ is called {\em $Y$-exact\/} if
$$
\codim {\textstyle\Tran_G(y, Y)} = \codim Y.
$$
\end{definition}

Note that by Lemma~\ref{lem: transporter dimension}(ii) a point $y \in Y$ is $Y$-exact if and only if
$$
\dim X - \dim(\overline{G.y}) = \dim Y - \dim(\overline{G.y \cap Y});
$$
in other words, the codimension in $X$ of the closure of the orbit containing $y$ is equal to that in $Y$ of the closure of the orbit's intersection with $Y$.

Our result is then as follows.

\begin{lem}\label{lem: exactness condition}
Let $Y$ be a subvariety of $X$, and $\hat Y$ be a dense open subset of $Y$; suppose that all points in $\hat Y$ are $Y$-exact. Then $\phi(G \times \hat Y)$ contains a dense open subset of $X$.
\end{lem}

\begin{proof}
Take $y \in \hat Y$; by assumption and Lemma~\ref{lem: transporter dimension}(i) we have
\begin{eqnarray*}
\dim \phi^{-1}(y) & = & \dim {\textstyle\Tran_G(y, Y)} \\
                  & = & \dim G - (\dim X - \dim Y) \\
                  & = & \dim(G \times Y) - \dim X.
\end{eqnarray*}
Let $X' = \overline{\phi(G \times Y)}$, and regard $\phi$ as a morphism $G \times Y \to X'$; both $G \times Y$ and $X'$ are irreducible, and here $\phi$ is dominant. By \cite[Theorem 4.1]{HumLAG}, each component of $\phi^{-1}(y)$ thus has dimension at least $\dim(G \times Y) - \dim X'$; so by the above $\dim X' \geq \dim X$, and as $X$ is irreducible and contains the closed set $X'$ we must have $X' = X$. Thus the morphism $\phi : G \times Y \to X$ is dominant.

Since $\hat Y$ is a dense open subset of $Y$, we see that $G \times \hat Y$ is a dense open subset of $G \times Y$, and hence constructible; as morphisms send constructible sets to constructible sets by \cite[Theorem 4.4]{HumLAG}, $\phi(G \times \hat Y)$ is constructible. Moreover the closure of $\phi(G \times \hat Y)$ contains $\phi(\overline{G \times \hat Y}) = \phi(G \times Y)$, so it contains $\overline{\phi(G \times Y)} = X$. As any constructible set contains a dense open subset of its closure, we see that $\phi(G \times \hat Y)$ contains a dense open subset of $X$.
\end{proof}

As a consequence we have the following.

\begin{lem}\label{lem: generic stabilizer from exact subset}
Let $Y$ be a subvariety of $X$, and $\hat Y$ be a dense open subset of $Y$; let $C$ be a subgroup of $G$ containing $G_X$. Suppose that for each $y \in \hat Y$ the following are true:
\begin{itemize}
\item[(i)] $y$ is $Y$-exact;
\item[(ii)] the stabilizer $C_G(y)$ is a conjugate of $C$.
\end{itemize}
Then $C/G_X$ is the generic stabilizer in the action of $G$ on $X$.
\end{lem}

\begin{proof}
From (i) we know by Lemma~\ref{lem: exactness condition} that $\phi(G \times \hat Y)$ contains a dense open subset $\hat X$ of $X$. As elements of $X$ lying in the same orbit have conjugate stabilizers, from (ii) we know that each element of $\hat X$ has stabilizer equal to a conjugate of $C$; taking the quotient by the kernel $G_X$ proves the result.
\end{proof}

This result may be seen as localizing the problem: we seek a subvariety $Y$ of $X$, and a dense open subset $\hat Y$ of $Y$ all of whose points are $Y$-exact and have conjugate stabilizers. In practice we want $Y$ to be a relatively small subvariety, since we need to determine both transporters and stabilizers of all points in the dense open subset $\hat Y$; of course, since the transporter $\Tran_G(y, Y)$ contains the stabilizer $C_G(y)$, identifying the former takes us some way towards finding the latter. In fact, often we are able to arrange things such that the stabilizer of each point in $\hat Y$ is conjugate to the subgroup $C$ by an element of $T$, as opposed to a general element of $G$.

\section{Semisimple automorphisms}\label{sect: semisimple auts}

In many of the cases where a triple $(G, \lambda, p)$ or quadruple $(G, \lambda, p, k)$ fails to have TGS, it turns out that we may locate $G$ inside a larger simple algebraic group $H$ and the module $V = L(\lambda)$ inside $\L(H)$. Usually we do so by taking a maximal parabolic subgroup $P$ of $H$, such that $G$ is the derived group of the Levi subgroup of $P$ while $V$ lies in the Lie algebra of the unipotent radical of $P$. In this section, however, we discuss a slightly different set-up.

As in Section~\ref{sect: notation}, take a simple algebraic group $H$ over the algebraically closed field $K$ of characteristic $p$, with maximal torus $T_H$, Lie algebra $\L(H)$ and so on; we shall assume that $H$ is of simply connected type. Let $\theta$ be a semisimple automorphism of $H$ of order $r$ coprime to $p$; we may assume $\theta$ preserves the torus $T_H$. Then $\L(H)$ decomposes as the direct sum of $r$ eigenspaces for $\theta$; for $0 \leq i < r$ denote the eigenspace corresponding to the eigenvalue ${\eta_r}^i$ by $\L(H)_{(i)}$. We have $\L(H)_{(0)} = \L(C_H(\theta))$; we shall focus on the eigenspace $\L(H)_{(1)}$, which is clearly a $C_H(\theta)$-module. Set $\L(T_H)_{(1)} = \L(T_H) \cap \L(H)_{(1)}$ and $Z(\L(H))_{(1)} = Z(\L(H)) \cap \L(H)_{(1)}$, and let
$$
Y = \L(T_H)_{(1)}/Z(\L(H))_{(1)}.
$$
Write
\begin{eqnarray*}
({W_H}^\ddagger)_{(1)} & = & \{ w \in W_H : \exists \xi \in K^*, \ \forall y \in Y, \ w.y = \xi y \}, \\
({W_H}^\dagger)_{(1)}  & = & \{ w \in W_H : \forall y \in Y, \ w.y = y \};
\end{eqnarray*}
let $({N_H}^\ddagger)_{(1)}$ and $({N_H}^\dagger)_{(1)}$ be the pre-images of $({W_H}^\ddagger)_{(1)}$ and $({W_H}^\dagger)_{(1)}$ respectively under the quotient map $N_H \to W_H$. We then have the following result.

\begin{lem}\label{lem: semisimple auts}
With the notation established above, write $G = C_H(\theta)$ and $V = \L(H)_{(1)}/Z(\L(H))_{(1)}$.
\begin{itemize}
\item[(i)] Suppose $G_V = G \cap Z(H)$, and $v \in \L(U_H) \cap \L(H)_{(1)}$ is a regular nilpotent element such that $G \cap C_{U_H}(v) = \{ 1 \}$; then in the action of $G$ on $V$ the orbit containing $v + Z(\L(H))_{(1)}$ is regular.
\item[(ii)] Suppose $\dim \L(H)_{(1)} - \dim \L(T_H)_{(1)} = \dim G - \dim (G \cap T_H)$, and $\L(T_H)_{(1)}$ contains regular semisimple elements; then the generic stabilizer for the action of $G$ on $V$ is $C_{({N_H}^\dagger)_{(1)}}(\theta)/G_V$, while that for the action of $G$ on $\G{1}(V)$ is $C_{({N_H}^\ddagger)_{(1)}}(\theta)/Z(G)$.
\end{itemize}
\end{lem}

\begin{proof}
(i) As $v \in \L(U_H)$ is regular nilpotent, we have $C_H(v) = C_{U_H}(v) Z(H)$; thus as $G \cap C_{U_H}(v) = 1$ we have $C_G(v) = G \cap Z(H) = G_V$. Moreover the only nilpotent element in the coset $v + Z(\L(H))_{(1)}$ is $v$ itself; so the $G$-orbit containing $v$ must meet $v + Z(\L(H))_{(1)}$ simply in $v$, and thus $C_G(v + Z(\L(H)_{(1)}) = C_G(v) = G_V$. Therefore the stabilizer in $G/G_V$ of $v + Z(\L(H))_{(1)}$ is trivial as required.

(ii) Suppose $v \in \L(T_H)_{(1)}$ is regular semisimple. Since any $H$-orbit in $\L(H)$ has finite intersection with $\L(T_H)$, there are only finitely many elements $z \in Z(\L(H))_{(1)}$ such that $v + z$ lies in $H.v$; thus $C_H(v + Z(\L(H))_{(1)})$ is a finite union of cosets of $C_H(v)$, and so $\dim C_H(v + Z(\L(H))_{(1)}) = \dim C_H(v)$. Hence $v + Z(\L(H))_{(1)} \in Y$ is also regular semisimple. Let $\hat Y_1$ be the set of regular semisimple elements in $Y$. Now given $w \in W_H \setminus ({W_H}^\dagger)_{(1)}$, take $n \in N_H$ with $n T_H = w$; by assumption the set of elements of $Y$ fixed by $n$ is a proper closed subvariety of $Y$. Let $\hat Y_2$ be the complement of the union of these subvarieties as $w$ runs over $W_H \setminus ({W_H}^\dagger)_{(1)}$. Set $\hat Y = \hat Y_1 \cap \hat Y_2$; as both $\hat Y_1$ and $\hat Y_2$ are dense open subsets of $Y$, the same is true of $\hat Y$.

Take $y \in \hat Y$. We have $C_H(y) = ({N_H}^\dagger)_{(1)}$, and hence
\begin{eqnarray*}
C_G(y) & = & G \cap ({N_H}^\dagger)_{(1)} \\
       & = & C_H(\theta) \cap ({N_H}^\dagger)_{(1)} \\
       & = & C_{({N_H}^\dagger)_{(1)}}(\theta),
\end{eqnarray*}
which is the union of a finite number of cosets of $C_{T_H}(\theta) = G \cap T_H$. Thus we have $\dim(\overline{G.y}) = \dim G - \dim(G \cap T_H)$, while $\dim(\overline{G.y \cap Y}) = 0$ because of the observation above about orbits having finite intersection with $\L(T_H)$; therefore
\begin{eqnarray*}
\dim V - \dim(\overline{G.y}) & = & (\dim \L(H)_{(1)} - \dim Z(\L(H))_{(1)}) - (\dim G - \dim(G \cap T_H)) \\
                              & = & \dim \L(T_H)_{(1)} - \dim Z(\L(H))_{(1)} \\
                              & = & \dim Y - \dim(\overline{G.y \cap Y}).
\end{eqnarray*}
Hence $y$ is $Y$-exact. Thus the conditions of Lemma~\ref{lem: generic stabilizer from exact subset} hold; so the generic stabilizer for the action of $G$ on $V$ is $C_{({N_H}^\dagger)_{(1)}}(\theta)/G_V$ as required. Replacing $Y$ by $\G{1}(Y)$, and $({N_H}^\dagger)_{(1)}$ and $({W_H}^\dagger)_{(1)}$ by $({N_H}^\ddagger)_{(1)}$ and $({W_H}^\ddagger)_{(1)}$ respectively, in an exactly similar fashion we see that the generic stabilizer for the action of $G$ on $\G{1}(V)$ is $C_{({N_H}^\ddagger)_{(1)}}(\theta)/Z(G)$.
\end{proof}

In the remainder of this section we address two issues involved in applying this result, one relating to each part.

For the first part we require information about $C_{U_H}(v)$ for a regular nilpotent element $v$ lying in $\L(U_H)$. We begin by observing that we have a filtration of $U_H$ given by the heights of roots: for $i \in \N$ we write
$$
{U_H}^{(i)} = \prod_{\height(\alpha) \geq i} X_\alpha,
$$
and then we have $U_H = {U_H}^{(1)} > {U_H}^{(2)} > \cdots > {U_H}^{(m)} > {U_H}^{(m + 1)} = \{ 1 \}$ where $m$ is the height of the highest root. We wish to know the relationship between $C_{U_H}(v)$ and the subgroups ${U_H}^{(i)}$. We shall be interested in only a few possibilities for $H$.

\begin{lem}\label{lem: m_i for H}
Let $H = A_{\ell_H}$, or $D_4$ with $p \geq 3$, or $E_6$ with $p \geq 3$, or $E_7$, or $E_8$. Write
$$
m_1, \dots, m_{\ell_H} =
\begin{cases}
1, 2, \dots, \ell_H           & \hbox{if } H = A_{\ell_H}, \\
1, 3, 3, 5                    & \hbox{if } H = D_4 \hbox{ with } p \geq 3, \\
1, 4, 5, 7, 8, 11             & \hbox{if } H = E_6 \hbox{ with } p \geq 5, \\
3, 4, 5, 7, 8, 11             & \hbox{if } H = E_6 \hbox{ with } p = 3, \\
1, 5, 7, 9, 11, 13, 17        & \hbox{if } H = E_7 \hbox{ with } p \geq 5, \\
3, 5, 7, 9, 11, 13, 17        & \hbox{if } H = E_7 \hbox{ with } p = 3, \\
5, 7, 8, 9, 11, 13, 17        & \hbox{if } H = E_7 \hbox{ with } p = 2, \\
1, 7, 11, 13, 17, 19, 23, 29  & \hbox{if } H = E_8 \hbox{ with } p \geq 7, \\
5, 7, 11, 13, 17, 19, 23, 29  & \hbox{if } H = E_8 \hbox{ with } p = 5, \\
7, 9, 11, 13, 17, 19, 23, 29  & \hbox{if } H = E_8 \hbox{ with } p = 3, \\
8, 11, 13, 14, 17, 19, 23, 29 & \hbox{if } H = E_8 \hbox{ with } p = 2. \\
\end{cases}
$$
With the notation established, if $v \in \L(U_H)$ is a regular nilpotent element, then $C_{U_H}(v) = \{ y_1(c_1) \dots y_{\ell_H}(c_{\ell_H}) : c_i \in K \}$, where each $y_i$ is an injection from $K$ into ${U_H}^{(m_i)}$ whose image does not lie in ${U_H}^{(m_i + 1)}$.
\end{lem}

\begin{proof}
Most of this is proved in \cite[Chapter~13]{LSbook}; see Tables~13.4 and 13.6, and the proof (not merely the statement) of Proposition~13.5. The exception is the case of $E_7$ with $p = 2$, for which the argument is the same as that for $E_8$ with $p = 2$ given in the first half of the penultimate paragraph of Proposition~13.5.
\end{proof}

Indeed, from the proof of \cite[Proposition~13.5]{LSbook} we see that each $y_i(c)$ is of the form $\left(\prod_{\height(\alpha) = m_i} x_\alpha(n_\alpha c)\right)x$, where $x \in {U_H}^{(m_i + 1)}$, and the $n_\alpha \in K$ satisfy $\sum_{\height(\alpha) = m_i} n_\alpha e_\alpha \in C_{\L(G)}(v)$. It is thus a straightforward calculation to determine the coset $y_i(c) {U_H}^{(m_i + 1)}$ in the quotient group ${U_H}^{(m_i)}/{U_H}^{(m_i + 1)}$. In some cases we shall require this additional information; we shall deal with these as they arise.

Note that if $p$ is not a bad prime for $H$, then the values $m_1, \dots, m_{\ell_H}$ are those listed in \cite[Proposition~10.2.5]{Car1} as the integers obtained by subtracting $1$ from the degrees of the basic polynomial invariants of the Weyl group $W_H$; the value $j$ appears in the list once (respectively twice) if the difference between the numbers of roots of heights $j$ and $j + 1$ is one (respectively two).


To apply the second part of Lemma~\ref{lem: semisimple auts} we need to determine the subgroup $G \cap T_H = C_{T_H}(\theta)$ of $T_H$, the subspace $\L(T_H)_{(1)}$ of $\L(T_H)$, and the subgroups $({W_H}^\ddagger)_{(1)}$ and $({W_H}^\dagger)_{(1)}$ of $W_H$. The first two involve straightforward calculations; and to see that $\L(T_H)_{(1)}$ contains regular semisimple elements it is enough to check that for each $\alpha \in \Phi_H$ there exists $v \in \L(T_H)_{(1)}$ with $[v e_\alpha] \neq 0$, which is routine (and needed only if $\L(T_H)_{(1)} \neq \L(T_H)$). However, identifying the subgroups of $W_H$ can be more involved; the remainder of the present section is devoted to this issue.

Observe that the set $\{ h_\beta : \beta \in \Phi_H \}$ is a root system dual to $\Phi_H$, with simple system $\{ h_{\beta_1}, \dots, h_{\beta_{\ell_H}} \}$; given $\beta = \sum a_i \beta_i$ with the coefficients $a_i \in \Z$, we have $h_\beta = \sum \frac{a_i \langle \beta_i, \beta_i \rangle}{\langle \beta, \beta \rangle} h_{\beta_i}$ (note that \cite[Lemma 1.2]{Lawmaxab} shows that if $e(\Phi_H) > 1$ then $\beta$ is long if and only if $e(\Phi_H)$ divides each $a_i$ for which $\beta_i$ is short, so the fraction always gives an integer, which may then be regarded as an element of $K$). Moreover given $w \in W_H$ we have $w.h_\beta = h_{w(\beta)}$.

If $\L(T_H)_{(1)} = \L(T_H)$ (which occurs if either $r = 1$, or $r = 2$ and the automorphism $\theta$ acts on $\L(T_H)$ by negation), we shall abbreviate $({W_H}^\ddagger)_{(1)}$ and $({W_H}^\dagger)_{(1)}$ to ${W_H}^\ddagger$ and ${W_H}^\dagger$ respectively. Here we have the following result.

\begin{lem}\label{lem: W_H on L(T_H)/Z(L(H))}
With the notation established, we have the following.
\begin{itemize}
\item[(i)] If $H = A_1$ with $p = 2$, then ${W_H}^\ddagger = {W_H}^\dagger = W_H \cong \Z_2$.
\item[(ii)] If $H = A_2$ with $p = 3$, then ${W_H}^\ddagger = W_H \cong S_3$, and ${W_H}^\dagger = \langle w_{\beta_1} w_{\beta_2} \rangle \cong \Z_3$.
\item[(iii)] If $H = A_3$ with $p = 2$, then ${W_H}^\ddagger = {W_H}^\dagger = \langle w_{\beta_1} w_{\beta_3}, w_{\beta_1 + \beta_2} w_{\beta_2 + \beta_3} \rangle \cong {\Z_2}^2$.
\item[(iv)] If $H = B_2$ with $p = 2$, then ${W_H}^\ddagger = {W_H}^\dagger = W_H \cong Dih_8$.
\item[(v)] If $H = B_\ell$ or $C_\ell$ for $\ell \in [3, \infty)$ with $p = 2$, then ${W_H}^\ddagger = {W_H}^\dagger \cong {\Z_2}^\ell$.
\item[(vi)] If $H = D_4$ with $p = 2$, then ${W_H}^\ddagger = {W_H}^\dagger \cong {\Z_2}^3.{\Z_2}^2$.
\item[(vii)] If $H = D_\ell$ for $\ell \in [5, \infty)$ with $p = 2$, then ${W_H}^\ddagger = {W_H}^\dagger \cong {\Z_2}^{\ell - 1}$.
\item[(viii)] If $H = A_1$, $B_\ell$ for $\ell \in [2, \infty)$, $C_\ell$ for $\ell \in [3, \infty)$, or $D_\ell$ for even $\ell \in [4, \infty)$, with $p \neq 2$, then ${W_H}^\ddagger = \langle w_0 \rangle \cong \Z_2$, and ${W_H}^\dagger = \{ 1 \}$.
\item[(ix)] If $H = E_7$, $E_8$, $F_4$ or $G_2$, then ${W_H}^\ddagger = \langle w_0 \rangle \cong \Z_2$, and ${W_H}^\dagger = \{ 1 \}$ or $\langle w_0 \rangle$ according as $p \geq 3$ or $p = 2$.
\item[(x)] If $H = A_2$ with $p \neq 3$, or $A_3$ with $p \neq 2$, or $A_\ell$ for $\ell \in [4, \infty)$, or $D_\ell$ for odd $\ell \in [5, \infty)$ with $p \neq 2$, or $E_6$, then ${W_H}^\ddagger = {W_H}^\dagger = \{ 1 \}$.
\end{itemize}
\end{lem}

\begin{proof}
Note that as $\L(T_H)_{(1)} = \L(T_H)$ we have $Z(\L(H))_{(1)} = Z(\L(H))$; and if $Z(\L(H)) \neq \{ 0 \}$ then $H = A_\ell$ with $p$ a factor of $\ell + 1$, or $B_\ell$, $C_\ell$, $D_\ell$ or $E_7$ with $p = 2$, or $E_6$ with $p = 3$. We shall start with the cases where $Z(\L(H)) = \{ 0 \}$.

First suppose $e(\Phi_H) = 1$, so that $H = A_\ell$, $D_\ell$, $E_6$, $E_7$ or $E_8$, and we assume $p$ is not a factor of $\ell + 1$ in the first of these cases, $p \neq 2$ in the second and fourth, and $p \neq 3$ in the third. Here the root system $\{ h_\beta : \beta \in \Phi_H \}$ is isomorphic to $\Phi_H$. The result is clear if $H = A_1$, so assume this is not the case. We claim that, for each $j$, if $\beta = \sum a_i \beta_i$ is any root other than $\pm \beta_j$ then the highest common factor in $\Z$ of the coefficients $a_i$ for $i \neq j$ is $1$. It suffices to consider $\beta$ positive. If $H = A_2$ and $i \neq j$ then the only positive roots other than $\beta_j$ have $n_i = 1$. If $H = A_\ell$ for $\ell \geq 3$, or $D_\ell$ for $\ell \geq 4$, or $E_6$, then given any $j$ there exists $i \neq j$ such that the coefficient of $\beta_i$ in the highest root is $1$; then $a_i$ is either $0$, in which case $\beta$ lies in a proper subsystem and the claim follows by induction, or $1$, in which case the highest common factor is certainly $1$. If $H = E_7$, then for $j \neq 7$ we may take $i = 7$ and the same argument applies; if instead $j = 7$ we take $i = 1$ and observe that $a_i \in \{ 0, 1, 2 \}$ --- the cases $a_i = 0$ and $a_i = 1$ are as before, while $a_i = 2$ only occurs if $\beta = \esevenrt2234321$, when the highest common factor is $1$ by inspection. Finally if $H = E_8$ and $j \neq 8$ we may take $i = 8$ and observe that $a_i \in \{ 0, 1, 2 \}$ --- again the cases $a_i = 0$ and $a_i = 1$ are as before, while $a_i = 2$ only occurs if $\beta = \eeightrt23465432$, when the highest common factor is $1$ by inspection; if instead $j = 8$ we take $i = 1$ and again observe that $a_i \in \{ 0, 1, 2 \}$ --- once more the cases $a_i = 0$ and $a_i = 1$ are as before, while this time $a_i = 2$ implies $a_3 \in \{ 3, 4 \}$, and $a_3 = 4$ implies $a_4 \in \{ 5, 6 \}$, and $a_4 = 6$ implies $a_2 = 3$, so that again the highest common factor is $1$. Now that the claim has been proved, it follows that any $w \in {W_H}^\ddagger$ must send each $\beta_j$ to $\pm \beta_j$, since $p$ must divide $w.h_{\beta_j} - \xi h_{\beta_j}$ for some $\xi$; the connectedness of the Dynkin diagram forces all the signs to be the same, since if $\beta_j$ and $\beta_{j'}$ correspond to adjacent nodes then $\beta_j + \beta_{j'}$ is a root while $\beta_j - \beta_{j'}$ is not; thus $w$ must be either $1$ or $-1$, with the latter occurring only if $H = A_1$, $D_\ell$ for $\ell$ even, $E_7$ or $E_8$. Hence in these cases we have ${W_H}^\ddagger = \langle w_0 \rangle$, and so ${W_H}^\dagger = \{ 1 \}$ or $\langle w_0 \rangle$ according as $p \geq 3$ or $p = 2$; in the cases where $H = A_\ell$ for $\ell \geq 2$, or $D_\ell$ for $\ell$ odd, or $E_6$, we have ${W_H}^\ddagger = {W_H}^\dagger = \{ 1 \}$.

Next suppose $e(\Phi_H) = 3$, so that $H = G_2$. Given $w \in W_H$ we have $w.h_{\beta_2} \in \{ \pm h_{\beta_2}, \pm (h_{\beta_1} + h_{\beta_2}), \pm (h_{\beta_1} + 2h_{\beta_2}) \}$; thus if $w \in {W_H}^\ddagger$ we must have $w(\beta_2) = \pm \beta_2$, so that $w \in \{ 1, w_{\beta_2}, w_0, w_0 w_{\beta_2} \}$. Certainly $w_0 = -1 \in {W_H}^\ddagger$. However if $w = w_{\beta_2}$ then $w.h_{\beta_1} = h_{\beta_1} + 3h_{\beta_2}$ while $w.h_{\beta_2} = -h_{\beta_2}$; thus the condition $w.h_{\beta_1} = \xi h_{\beta_1}$ forces $p = 3$ and $\xi = 1$, so we do not have $w.h_{\beta_2} = \xi h_{\beta_2}$. Therefore ${W_H}^\ddagger = \langle w_0 \rangle$, and so ${W_H}^\dagger = \{ 1 \}$ or $\langle w_0 \rangle$ according as $p \geq 3$ or $p = 2$.

Now suppose $e(\Phi_H) = 2$, so that $H = B_\ell$, $C_\ell$ or $F_4$, and we assume $p \neq 2$ in the first two cases; take $w \in {W_H}^\ddagger$. If $H = B_\ell$ for $\ell \geq 3$, then for $j < \ell$ the coefficient of $h_{\beta_\ell}$ in $w.h_{\beta_j}$ lies in $\{ 0, \pm 1 \}$; as before it cannot be $\pm 1$, so it must be $0$ and then arguing in the subsystem of type $A_{\ell - 1}$ we see that $w(\beta_j) = \pm \beta_j$. Again the signs must all be the same, so as $w_0 = -1$ we see that either $w$ or $w_0 w$ must fix each $\beta_j$ for $j < \ell$; but the only element of $W_H$ which does this is $1$, so ${W_H}^\ddagger = \langle w_0 \rangle$ and ${W_H}^\dagger = \{ 1 \}$. If $H = C_\ell$ for $\ell \geq 2$, the elements $h_\beta$ for $\beta$ long are of the form
$$
h_{\beta_i} + h_{\beta_{i + 1}} + \cdots + h_{\beta_\ell}
$$
for $i \leq \ell$; thus $w(\beta_\ell)$ must be $\pm \beta_\ell$. For $\beta$ short the coefficient of $h_{\beta_\ell}$ in $h_\beta$ lies in $\{ 0, \pm 2 \}$. Thus for each $j < \ell$ we must have $w(\beta_j) \in \langle \beta_1, \dots, \beta_{\ell - 1} \rangle$; arguing in the subsystem of type $A_{\ell - 1}$ we see that we must have $w(\beta_j) = \pm \beta_j$, from which it follows as before that $w \in \{ 1, w_0 \}$, so that ${W_H}^\ddagger = \langle w_0 \rangle$ and ${W_H}^\dagger = \{ 1 \}$. Finally if $H = F_4$, for $j \in \{ 1, 2 \}$ the coefficient of $h_{\beta_4}$ in $w.h_{\beta_j}$ lies in $\{ 0, \pm 1, \pm 2 \}$, and if it is $2\e$ for $\e = \pm 1$ then the coefficient of $h_{\beta_3}$ is $3\e$ --- so again $w(\beta_j) = \pm \beta_j$, and the signs for $j = 1$ and $j = 2$ must be the same. As $w_0 = -1$, either $w$ or $w_0 w$ must fix both $\beta_1$ and $\beta_2$, so must lie in $\langle w_{\beta_4}, w_{\beta_1 + 2\beta_2 + 3\beta_3 + \beta_4} \rangle \cong S_3$; of these six elements, four map $h_{\beta_3}$ to either $h_{\beta_3} + h_{\beta_4}$ or $-(2h_{\beta_1} + 4h_{\beta_2} + 2h_{\beta_3} + h_{\beta_4})$, and a fifth fixes $h_{\beta_3}$ but maps $h_{\beta_4}$ to $-(2h_{\beta_1} + 4h_{\beta_2} + 3h_{\beta_3} + h_{\beta_4})$, so the only one lying in ${W_H}^\ddagger$ is $1$. Thus ${W_H}^\ddagger = \langle w_0 \rangle$, and ${W_H}^\dagger = \{ 1 \}$ or $\langle w_0 \rangle$ according as $p \geq 3$ or $p = 2$.

We now turn to the cases where $Z(\L(H)) \neq \{ 0 \}$. Note that ${W_H}^\ddagger \lhd W_H$.

First suppose $H$ is of exceptional type. If $H = E_6$ with $p = 3$, then $W_H \cong \mathrm{S}_4(3).\Z_2$; as ${W_H}^\ddagger$ clearly does not contain $\mathrm{S}_4(3)$, we must have ${W_H}^\ddagger = {W_H}^\dagger = \{ 1 \}$. If $H = E_7$ with $p = 2$, then $W_H \cong \mathrm{S}_6(2) \times \Z_2$, where the $\Z_2$ is $\langle w_0 \rangle$; as $w_0 = -1$, and ${W_H}^\ddagger$ clearly does not contain $\mathrm{S}_6(2)$, we must have ${W_H}^\ddagger = {W_H}^\dagger = \langle w_0 \rangle$.

Next suppose $H = A_\ell$ with $p$ a factor of $\ell + 1$; then $W_H \cong S_{\ell + 1}$. If $\ell = 1$ and $p = 2$ we have $\L(T_H) = \langle h_{\beta_1} \rangle = Z(\L(H))$, so ${W_H}^\ddagger = {W_H}^\dagger = W_H$. If $\ell = 2$ and $p = 3$ we have $\L(T_H) = \langle h_{\beta_1}, h_{\beta_2} \rangle$ and $Z(\L(H)) = \langle z_1 \rangle$, where $z_1 = h_{\beta_1} - h_{\beta_2}$. Since $\dim \L(T_H)/Z(\L(H)) = 1$, all elements of $W_H$ act on $\L(T_H)/Z(\L(H))$ as scalars, so ${W_H}^\ddagger = W_H$; the transposition $w_{\beta_1}$ negates $h_{\beta_1}$, and sends $h_{\beta_2}$ to $h_{\beta_1} + h_{\beta_2} = -h_{\beta_2} + z_1$, so it acts on $\L(T_H)/Z(\L(H))$ as $-1$; similarly each of the other two transpositions acts as $-1$, and so the $3$-cycles act as $1$, whence ${W_H}^\dagger = \langle w_{\beta_1} w_{\beta_2} \rangle \cong \Z_3$. If $\ell = 3$ and $p = 2$ we have $\L(T_H) = \langle h_{\beta_1}, h_{\beta_2}, h_{\beta_3} \rangle$ and $Z(\L(H)) = \langle z_1 \rangle$, where $z_1 = h_{\beta_1} + h_{\beta_3}$. The double transposition $w_{\beta_1} w_{\beta_3}$ negates both $h_{\beta_1}$ and $h_{\beta_3}$, and sends $h_{\beta_2}$ to $h_{\beta_2} + z_1$, so it acts on $\L(T_H)/Z(\L(H))$ as $1$; similarly each of the other two double transpositions acts as $1$, and as the $3$-cycle $w_{\beta_1} w_{\beta_2}$ sends $h_{\beta_1}$ to $h_{\beta_2}$ and so does not act as a scalar, we see that ${W_H}^\ddagger = {W_H}^\dagger = {\Z_2}^2$. If $\ell \geq 4$ then the only proper non-trivial normal subgroup of $S_{\ell + 1}$ is $Alt_{\ell + 1}$; as ${W_H}^\ddagger$ clearly does not contain $Alt_{\ell + 1}$, we must have ${W_H}^\ddagger = {W_H}^\dagger = \{ 1 \}$.

Next suppose $H = B_\ell$ for $\ell \geq 2$ with $p = 2$; then $W_H \cong {\Z_2}^\ell.S_\ell$. We have $\L(T_H) = \langle h_{\beta_1}, \dots, h_{\beta_\ell} \rangle$ and $Z(\L(H)) = \langle z_1 \rangle$, where $z_1 = h_{\beta_\ell}$. If $\ell = 2$ then $\dim \L(T_H)/Z(\L(H)) = 1$, so all elements of $W_H$ act on $\L(T_H)/Z(\L(H))$ as scalars, and hence ${W_H}^\ddagger = W_H \cong Dih_8$; indeed any $w \in W_H$ fixes $h_{\beta_2}$ and maps $h_{\beta_1}$ to either $\pm h_{\beta_1} = h_{\beta_1}$ or $\pm h_{\beta_1} + h_{\beta_2} = h_{\beta_1} + z_1$, so ${W_H}^\dagger = W_H$. Now assume $\ell \geq 3$ and take $w \in {W_H}^\ddagger$. Given $j < \ell$, if $j > 1$ the coefficient of $h_{\beta_1}$ in $w.h_{\beta_j}$ lies in $\{ 0, \pm 1 \}$, so must be $0$; arguing similarly we see that for each $i < j$ the coefficient of $h_{\beta_i}$ in $w.h_{\beta_j}$ must be $0$, so that $w(\beta_j) \in \langle \beta_j, \beta_{j + 1}, \dots, \beta_\ell \rangle$; this gives $w(\beta_{\ell - 1}) \in \{ \pm \beta_{\ell - 1}, \pm (\beta_{\ell - 1} + 2\beta_\ell) \}$, while if $j < \ell - 1$ the coefficient of $h_{\beta_{j + 1}}$ in $w.h_{\beta_j}$ cannot be $\pm 1$, so must be $0$ or $\pm 2$, whence $w(\beta_j) \in \{ \pm \beta_j, \pm (\beta_j + 2\beta_{j + 1} + \cdots + 2\beta_{\ell - 1} + 2\beta_\ell) \}$. Now if we identify the dual root system $\{ h_\beta : \beta \in \Phi_H \}$ with the standard root system of type $C_\ell$, then for each $j < \ell$ we must map $\ve_j - \ve_{j + 1}$ to $\pm \ve_j \pm \ve_{j + 1}$; the elements concerned are those which map each $\ve_i$ to $\pm \ve_i$ with independent choice of signs, so we have ${W_H}^\ddagger = {W_H}^\dagger \cong {\Z_2}^\ell$.

Now suppose $H = C_\ell$ for $\ell \geq 3$ with $p = 2$; then $W_H \cong {\Z_2}^\ell.S_\ell$. We have $\L(T_H) = \langle h_{\beta_1}, \dots, h_{\beta_\ell} \rangle$ and $Z(\L(H)) = \langle z_1 \rangle$, where $z_1 = \sum_{i = 1}^{\lceil \ell/2 \rceil} h_{\beta_{2i - 1}}$. Take $w \in {W_H}^\ddagger$. Here the elements $h_\beta$ for $\beta$ long are of the form
$$
h_{\beta_i} + h_{\beta_{i + 1}} + \cdots + h_{\beta_\ell}
$$
for $i \leq \ell$, while those for $\beta$ short are of the form
$$
h_{\beta_i} + h_{\beta_{i + 1}} + \cdots + h_{\beta_j}
$$
for $i < j < \ell$; thus $w(\beta_\ell)$ must be $\pm \beta_\ell$. Provided $\ell \neq 4$, we see that for each $j < \ell$ and $\kappa \in K^*$ the element $\kappa h_{\beta_j} + z_1$ is not of the form $\kappa' h_{\beta}$ for any root $\beta$ and any $\kappa' \in K^*$, so we must have $w.h_{\beta_j} = h_{\beta_j}$; arguing just as in the previous paragraph shows that $w(\beta_{\ell - 1}) \in \{ \pm \beta_{\ell - 1}, \pm (\beta_{\ell - 1} + \beta_\ell) \}$, while if $j < \ell - 1$ then $w(\beta_j) \in \{ \pm \beta_j, \pm (\beta_j + 2\beta_{j + 1} + \cdots + 2\beta_{\ell - 1} + \beta_\ell) \}$. If however $\ell = 4$ then as $z_1 = h_{\beta_1} + h_{\beta_3}$ we have $h_{\beta_3} = h_{\beta_1} + z_1$ and $h_{\beta_1 + \beta_2 + \beta_3} = h_{\beta_2} + z_1$; as a result we can initially conclude only that $w(\beta_1), w(\beta_3) \in \{ \pm \beta_3, \pm(\beta_3 + \beta_4), \pm \beta_1, \pm(\beta_1 + 2\beta_2 + 2\beta_3 + \beta_4) \}$ and $w(\beta_2) \in \{ \pm \beta_2, \pm(\beta_2 + 2\beta_3 + \beta_4), \pm(\beta_1 + \beta_2 + \beta_3), \pm(\beta_1 + \beta_2 + \beta_3 + \beta_4) \}$. However, the fact that $w(\beta_1)$ and $w(\beta_2)$ must be orthogonal to $w(\beta_4)$, while $w(\beta_3) + w(\beta_4)$ is a root, reduces to the possibilities given before. Now if we identify the dual root system $\{ h_\beta : \beta \in \Phi_H \}$ with the standard root system of type $B_\ell$, then we must map $\ve_\ell$ to $\pm \ve_\ell$, and for each $j < \ell$ we must map $\ve_j - \ve_{j + 1}$ to $\pm \ve_j \pm \ve_{j + 1}$; the elements concerned are those which map each $\ve_i$ to $\pm \ve_i$ with independent choice of signs, so we have ${W_H}^\ddagger = {W_H}^\dagger \cong {\Z_2}^\ell$.

Finally suppose $H = D_\ell$ for $\ell \geq 4$ with $p = 2$; then $W_H \cong {\Z_2}^{\ell - 1}.S_\ell$. We have $\L(T_H) = \langle h_{\beta_1}, \dots, h_{\beta_\ell} \rangle$ and $Z(\L(H)) = \langle z_1 \rangle$ or $\langle z_1, z_2 \rangle$ according as $\ell$ is odd or even, where $z_1 = h_{\beta_{\ell - 1}} + h_{\beta_\ell}$ and if $\ell$ is even $z_2 = \sum_{i = 1}^{\ell/2} h_{\beta_{2i - 1}}$. Clearly $w_{\beta_{\ell - 1}} w_{\beta_\ell}$ fixes $h_{\beta_i}$ for $i < \ell - 2$, negates both $h_{\beta_{\ell - 1}}$ and $h_{\beta_\ell}$, and sends $h_{\beta_{\ell - 2}}$ to $h_{\beta_{\ell - 2}} + h_{\beta_{\ell - 1}} + h_{\beta_\ell} = h_{\beta_{\ell - 2}} + z_1$, so it acts on $\L(T_H)/Z(\L(H))$ as $1$; similarly each of the other elements of the normal subgroup ${\Z_2}^{\ell - 1}$ acts as $1$. Since the $3$-cycle $w_{\beta_1} w_{\beta_2}$ sends $h_{\beta_1}$ to $h_{\beta_2}$ and so does not act as a scalar, we see that if $\ell \geq 5$ then we must have ${W_H}^\ddagger = {W_H}^\dagger = {\Z_2}^{\ell - 1}$. If however $\ell = 4$, then $w_{\beta_1} w_{\beta_3}$ negates both $h_{\beta_1}$ and $h_{\beta_3}$, fixes $h_{\beta_4}$ and sends $h_{\beta_2}$ to $h_{\beta_1} + h_{\beta_2} + h_{\beta_3} = h_{\beta_2} + z_2$, so it acts on $\L(T_H)/Z(\L(H))$ as $1$; similarly each of the other double transpositions acts as $1$, and so ${W_H}^\ddagger = {W_H}^\dagger = {\Z_2}^3.{\Z_2}^2$.
\end{proof}

In each of the remaining cases it will turn out that $Z(\L(H)) = \{ 0 \}$, so we shall be concerned simply with the action of $W_H$ on $\L(T_H)_{(1)}$. Rather than presenting a series of seemingly unmotivated results here, we shall include the determination of the subgroups $({W_H}^\ddagger)_{(1)}$ and $({W_H}^\dagger)_{(1)}$ of $W_H$ within the proofs of the results giving the generic stabilizers concerned. However, we make some general comments here on an approach which may often be applied.

Take $w \in W_H$ with the property that there exists $\xi \in K^*$ such that for all $y \in \L(T_H)_{(1)}$ we have $w.y = \xi y$. Suppose $\L(T_H)_{(1)}$ contains a vector of the form $y = \kappa h_{\beta_j} + \kappa' h_{\beta_{j'}}$ for $\kappa, \kappa' \in K^*$ and two simple roots $\beta_j$ and $\beta_{j'}$ of the same length; then $\kappa h_{w(\beta_j)} + \kappa' h_{w(\beta_{j'})} = \xi\kappa h_{\beta_j} + \xi\kappa' h_{\beta_{j'}}$. As above we may write $w(\beta_j) = \sum a_i \beta_i$ and $w(\beta_{j'}) = \sum {a_i}' \beta_i$ where all $a_i, {a_i}' \in \Z$; then $h_{w(\beta_j)} = \sum \frac{a_i \langle \beta_i, \beta_i \rangle}{\langle \beta_j, \beta_j \rangle} h_{\beta_i}$ and $h_{w(\beta_{j'})} = \sum \frac{{a_i}' \langle \beta_i, \beta_i \rangle}{\langle \beta_{j'}, \beta_{j'} \rangle} h_{\beta_i}$ (where we must now regard the coefficients as lying in $K$). Thus for $i \neq j, j'$ we must have $\kappa \frac{a_i \langle \beta_i, \beta_i \rangle}{\langle \beta_j, \beta_j \rangle} + \kappa' \frac{{a_i}' \langle \beta_i, \beta_i \rangle}{\langle \beta_{j'}, \beta_{j'} \rangle} = 0$, so as $\langle \beta_j, \beta_j \rangle = \langle \beta_{j'}, \beta_{j'} \rangle$ we have ${a_i}' = -\frac{\kappa}{\kappa'} a_i$; we say that the roots $w(\beta_j)$ and $w(\beta_{j'})$ are {\em proportional outside $\{ \beta_j, \beta_{j'} \}$\/}. Often inspection of the root system (regarding coefficients as lying in $K$) reveals that the only possibility is that there exists $\e \in \{ \pm1 \}$ such that for all $i \neq j, j'$ we have ${a_i}' = \e a_i$; and if $\frac{\kappa}{\kappa'} \neq \pm1$, then for all $i \neq j, j'$ we must have $a_i = {a_i}' = 0$, so that $w$ preserves $\langle \beta_j, \beta_{j'} \rangle$.

\section{Use of Lie algebra annihilators}\label{sect: Lie algebra annihilators}

In this section we describe an approach which in certain circumstances may be applied in combination with Lemma~\ref{lem: generic stabilizer from exact subset} to determine the generic stabilizer for a triple $(G, \lambda, p)$ or quadruple $(G, \lambda, p, k)$; the argument is in essence the work of Alexander Premet. Write $V = L(\lambda)$ as usual, and let $X$ be the variety $V$ or $\Gk(V)$ as appropriate.

We begin by defining the annihilator in $\L(G)$ of an element $y$ of $X$. If $X = V$, so that $y$ is a vector in $V$, we write
$$
{\ts\Ann_{\L(G)}(y)} = \{ v \in \L(G) : v.y = 0 \};
$$
if instead $X = \Gk(V)$, so that $y$ is a $k$-dimensional subspace of $V$, we write
$$
{\ts\Ann_{\L(G)}(y)} = \{ v \in \L(G) : v.y \leq y \}.
$$

We now suppose that we have an irreducible subvariety $Y$ of $X$, a dense open subset $\hat Y$ of $Y$, and a subalgebra $\S$ of $\L(T)$, which between them satisfy a number of conditions.

\begin{lem}\label{lem: exactness via Premet}
Suppose $C_{\L(G)}(\S) = \L(T)$, and $\S$ lies in the annihilator of each $y \in Y$.
\begin{itemize}
\item[(i)] Given $y \in Y$, if $\Ann_{\L(G)}(y) = \S$ then $\Tran_G(y, Y) \subseteq N$ (so in particular $C_G(y) \leq N$).
\item[(ii)] Suppose $\codim Y = \dim G$, and for all $y \in \hat Y$ we have
\begin{itemize}
\item[(a)] $\Ann_{\L(G)}(y) = \S$,
\item[(b)] $C_T(y) = \{ 1 \}$,
\item[(c)] $|T.y \cap Y| < \infty$, and
\item[(d)] $N.y \cap Y \subset \hat Y$;
\end{itemize}
then each $y \in \hat Y$ is $Y$-exact.
\end{itemize}
\end{lem}

\begin{proof}
(i) If $g \in \Tran_G(y, Y)$, then as $g.y \in Y$ we have $\S \leq \Ann_{\L(G)}(g.y) = \Ad(g).\Ann_{\L(G)}(y) = \Ad(g).\S$, whence $\Ad(g).\S = \S$. Therefore $\Ad(g).C_{\L(G)}(\S) = C_{\L(G)}(\S)$, so $\Ad(g).\L(T) = \L(T)$, giving $g \in N_G(T) = N$ as required.

(ii) Take $y \in \hat Y$. As $T.y \cap Y$ is finite and contained in $\hat Y$, and $C_T(y) = \{ 1 \}$, it follows that $\Tran_T(y, Y)$ is finite. For each $w \in W$ choose $n_w \in N$ with $n_wT = w$; then the set $\Tran_T(n_w.y, Y) = \{ t \in T : tn_w.y \in Y \}$ is finite (because if it is non-empty and $tn_w.y$ lies in $Y$ then $tn_w.y$ lies in $\hat Y$), as is thus $\Tran_G(y, Y) = \bigcup_{w \in W} \Tran_T(n_w.y, Y)n_w$. Therefore we have $\codim \Tran_G(y, Y) = \dim G = \codim Y$, whence $y$ is $Y$-exact as required.
\end{proof}

Thus if the conditions of Lemma~\ref{lem: exactness via Premet} hold, provided one can show that for each $y \in \hat Y$ we have $C_N(y) = C$ where $C$ is a fixed subgroup of $G$, Lemma~\ref{lem: generic stabilizer from exact subset} may be applied to identify the generic stabilizer as $C/G_X$. We will do this several times in Sections~\ref{sect: non-TGS large triples and first quadruples} and \ref{sect: non-TGS large higher quadruples}; on each occasion, once appropriate $Y$, $\hat Y$ and $\S$ have been defined, the bulk of the proof will therefore consist of calculations showing that all the conditions hold, and the identification of the subgroup $C$. In fact we shall find that there is one case in Section~\ref{sect: non-TGS large higher quadruples} where only some of the conditions hold; although Lemma~\ref{lem: exactness via Premet} cannot therefore be applied as it stands, it will turn out that the basic strategy of the proof does go through, although the argument is considerably more complicated than in the other cases.

\section{Invariants}\label{sect: invariants}

In this section we consider a rather special situation. Let $(G, \lambda, p)$ be a triple and as usual write $V = L(\lambda)$. Suppose that $G$ is of type $A_\ell$ for some $\ell$, and $\dim V = \dim G + 1$ (so that the triple is large and the associated first quadruple $(G, \lambda, p, 1)$ is small). We may take $G = \SL_{\ell + 1}(K)$; set $G^+ = \GL_{\ell + 1}(K)$. Suppose also that the action of $G$ on $V$ extends to the group $G^+$, and that there is a non-trivial invariant in $K[V]$ for the action of $G$; let $f$ be an invariant of minimal positive degree, which we may take to have zero constant term. For $a \in K$ write $V(a) = \{ v \in V : f(v) = a \}$, so that $V$ is the disjoint union of the varieties $V(a)$, each of which is preserved by $G$ and of dimension $\dim V - 1 = \dim G$.

\begin{lem}\label{lem: invariant is homogeneous and irreducible}
With the notation established, the invariant $f$ is homogeneous; if $a \in K$ then $f - a$ is irreducible, so that $V(a)$ is an irreducible variety.
\end{lem}

\begin{proof}
Since homogeneous components of an invariant are also invariants, the first statement follows from minimality of degree. Given $a \in K$, if $f - a$ were a product of irreducible factors, then each element of $G$ would have to permute and scale them; as $G$ is connected the permutation must be trivial, and as $G$ is perfect so must the scaling, so again the minimality of degree implies irreducibility of $f - a$, and therefore of its zero set $V(a)$.
\end{proof}

Now suppose additionally that there exists $c \in \Z \setminus \{ 0 \}$ such that for all $\kappa \in K^*$ and $v \in V$ we have $(\kappa I).v = \kappa^{c}v$; then given $g \in G^+$ there exist $\kappa \in K^*$ and $g' \in G$ such that $g = g'(\kappa I)$, so by the homogeneity of $f$, for all $v \in V$ we have $f(g.v) = f(g'.(\kappa I).v) = f((\kappa I).v) = f(\kappa^{c}v) = \kappa^{c\deg f}f(v)$. It follows that $f$ is a relative invariant for the action of $G^+$ with associated (linear) character $\chi$, where $\chi(g) = \det(g)^{c \deg f/(\ell + 1)}$. In particular $G^+$ preserves $V(0)$ and hence $V \setminus V(0)$. In this situation, although the triple is large, identifying a single stabilizer may suffice to determine the generic stabilizer, and identifying a second may settle the question of the existence or otherwise of a regular orbit.

\begin{lem}\label{lem: gen stab in invariant situation}
With the notation established, suppose there exists $y_0 \in V$ such that $C_{G^+}(y_0)$ is finite. Then
\begin{itemize}
\item[(i)] the generic stabilizers for the actions of $G$ on $V$ and $\G{1}(V)$ are $C_G(y_0)/G_V$ and $C_G(\langle y_0 \rangle)/Z(G) \cong C_{G^+}(\langle y_0 \rangle)/Z(G^+)$ respectively;
\item[(ii)] if moreover $C_G(y_0) \neq G_V$, and there exists $y_1 \in V(0)$ such that $C_G(y_1)$ is finite, then in the action of $G$ on $V$ there is a regular orbit if and only if $C_G(y_1) = G_V$.
\end{itemize}
\end{lem}

\begin{proof}
By assumption the orbit $G^+.y_0$ is dense in $V$, and hence must lie in $V \setminus V(0)$. Given $a \in K^*$ there exists $\kappa \in K^*$ such that $\kappa y_0 \in V(a)$; as the stabilizer $C_G(\kappa y_0) = C_G(y_0)$ is finite, and the variety $V(a)$ is irreducible by Lemma~\ref{lem: invariant is homogeneous and irreducible}, it follows that the orbit $G.\kappa y_0$ is dense in $V(a)$. Since the union of the orbits $G.\kappa y_0$ as $\kappa$ runs through $K^*$ is the dense set $G^+.y_0$, (i) follows. Now suppose $C_G(y_0)/G_V$ is non-trivial. If there is a regular orbit in the action of $G$ on $V$, it must lie in some $V(a)$. If $a \neq 0$ there is a dense orbit $G.\kappa y_0$ in $V(a)$ which is not regular, and $V(a) \setminus G.\kappa y_0$ has dimension less than $\dim G$; thus any regular orbit must lie in $V(0)$. Since $y_1 \in V(0)$ has finite stabilizer $C_G(y_1)$, and the variety $V(0)$ is irreducible by Lemma~\ref{lem: invariant is homogeneous and irreducible}, the orbit $G.y_1$ is dense in $V(0)$, and $V(0) \setminus G.y_1$ has dimension less than $\dim G$; thus the only possible regular orbit is $G.y_1$, and (ii) follows.
\end{proof}

We shall use Lemma~\ref{lem: gen stab in invariant situation} to prove several results in Section~\ref{sect: non-TGS large triples and first quadruples}.

\section{Generalized height functions}\label{sect: gen ht fns}

We begin this section by recalling that \cite[Lemma~2.1]{GLMS} shows that if $G$ acts on a module $V$, and $v, v'$ are two vectors lying in the zero weight space of $V$, then $v$ and $v'$ lie in the same $G$-orbit if and only if they lie in the same $N$-orbit. The proof is straightforward: suppose $g \in G$ is such that $g.v = v'$; write $g$ in Bruhat decomposition as $g = unu'$ where $n \in N$, $u \in U$ and $u' \in U_w$ where $w = nT \in W$, then $nu'n^{-1}.(n.v) = u^{-1}.v'$; observe that the weights occurring on the left hand side are all zero or sums of negative roots, while those on the right are all zero or sums of positive roots; thus in each case the only weight present must be zero, and we have $n.v = v'$.

The basic idea of taking two points lying in the \lq middle' of a variety, and using Bruhat decomposition and comparison of weights to gain information about group elements which send one to the other, turns out to be very relevant to the work undertaken here. We shall provide a general setting and prove results which generalize \cite[Lemma~2.1]{GLMS}; these will be of use in determining transporters for suitably chosen subvarieties of either a $G$-module or a Grassmannian variety. Recall that $\Lambda$ is the weight lattice of $G$.

\begin{definition}
A {\em generalized height function\/} is a linear function $\Lambda \to \Z$ whose value at each simple root is non-negative; we refer to the value at any weight as the {\em generalized height\/} of the weight. A generalized height function is {\em strictly positive\/} if the generalized height of each simple root is in fact positive.
\end{definition}

There are of course many generalized height functions; in a given context we shall define the particular one being considered. Recall from Section~\ref{sect: weights and module structure} that the set of weights of $G$ has the partial order $\prec$ defined by $\mu \prec \nu$ if and only if $\nu - \mu$ is a non-empty sum of positive roots; thus if $\mu \prec \nu$, then for any strictly positive generalized height function the generalized height of $\mu$ is strictly less than that of $\nu$.

Let $V$ be a $G$-module. Given a generalized height function, for $i \in \Z$ we write $\Lambda(V)_{[i]}$ for the set of weights in $\Lambda(V)$ whose generalized height is $i$, and set
$$
V_{[i]} = \bigoplus_{\nu \in \Lambda(V)_{[i]}} V_\nu;
$$
we write $\Lambda(V)_{[-]} = \bigcup_{i < 0} \Lambda(V)_{[i]}$ and $\Lambda(V)_{[+]} = \bigcup_{i > 0} \Lambda(V)_{[i]}$, and set
$$
V_{[-]} = \bigoplus_{i < 0} V_{[i]} = \bigoplus_{\nu \in \Lambda(V)_{[-]}} V_\nu, \qquad V_{[+]} = \bigoplus_{i > 0} V_{[i]} = \bigoplus_{\nu \in \Lambda(V)_{[+]}} V_\nu,
$$
so that
$$
V = V_{[-]} \oplus V_{[0]} \oplus V_{[+]}.
$$
In addition we let $\Phi_{[0]}$ be the set of roots of generalized height $0$, and write
$$
G_{[0]} = \langle T, X_\alpha : \alpha \in \Phi_{[0]} \rangle,
$$
and
$$
U_{[0]} = \prod_{\alpha \in \Phi^+ \cap \Phi_{[0]}} X_\alpha, \qquad U_{[+]} = \prod_{\alpha \in \Phi^+ \setminus \Phi_{[0]}} X_\alpha,
$$
so that $U = U_{[0]} U_{[+]} = U_{[+]} U_{[0]}$ and $U_{[+]} \cap U_{[0]} = \{ 1 \}$ (and if the generalized height function is strictly positive we have $G_{[0]} = T$, $U_{[0]} = \{ 1 \}$ and $U_{[+]} = U$). Finally we write $W_{\Lambda(V)_{[0]}}$ for the stabilizer in $W$ of $\Lambda(V)_{[0]}$, and $N_{\Lambda(V)_{[0]}}$ for the preimage in $N$ of $W_{\Lambda(V)_{[0]}}$.

\begin{definition}
With the notation established, a subset $\Delta$ of $\Lambda(V)_{[0]}$ has {\em ZLC\/} (denoting \lq zero linear combination\rq) if there is a linear combination $\sum_{\nu \in \Delta} c_\nu \nu = 0$ in which for all $\nu \in \Delta$ we have $c_\nu \in \N$. A subset $\Delta$ of $\Lambda(V)_{[0]}$ has {\em ZLCE\/} (denoting \lq zero linear combination extended\rq) if all subsets $\Delta'$ of $\Lambda(V)_{[0]}$ with $\Delta \subseteq \Delta'$ have ZLC.
\end{definition}

Note that for a subset $\Delta$ of $\Lambda(V)_{[0]}$ to have ZLCE it is sufficient merely that all subsets $\Delta'$ of $\Lambda(V)_{[0]}$ with $\Delta \subseteq \Delta'$ and $|\Delta' \setminus \Delta| \leq 1$ have ZLC: if this weaker condition holds, then given any subset $\Delta'$ of $\Lambda(V)_{[0]}$ with $\Delta \subseteq \Delta'$, for each weight $\nu$ in $\Delta' \setminus \Delta$ we may take the corresponding linear combination of the weights in $\Delta \cup \{ \nu \}$; summing them all together then gives a linear combination of the weights in $\Delta'$ as required.

Our first result here gives partial information on certain transporters.

\begin{lem}\label{lem: gen height zero not strictly positive}
Let $X$ be either $V$ or $\Gk(V)$ for some $k \in \N$. Given a generalized height function on the weight lattice of $G$, let $Y$ be a subvariety of either $V_{[0]}$ or $\Gk(V_{[0]})$. Assume $W_{\Lambda(V)_{[0]}}$ stabilizes $\Phi_{[0]}$. Suppose $y \in Y$ has the property that for all $u \in U_{[0]}$ the set of weights occurring in $u.y$ has ZLCE; take $g \in \Tran_G(y, Y)$ and set $y' = g.y \in Y$. Then we may write $g = u_1 g' u_2$ with $u_1 \in C_{U_{[+]}}(y')$, $u_2 \in C_{U_{[+]}}(y)$, and $g' \in G_{[0]} N_{\Lambda(V)_{[0]}}$ with $g'.y = y'$. In particular $G.y \cap Y = G_{[0]} N_{\Lambda(V)_{[0]}}.y \cap Y$, and $C_G(y) = C_{U_{[+]}}(y) C_{G_{[0]} N_{\Lambda(V)_{[0]}}}(y) C_{U_{[+]}}(y)$.
\end{lem}

\begin{proof}
Suppose $y$, $g$ and $y'$ are as given; use the Bruhat decomposition and the factorization $U = U_{[0]} U_{[+]} = U_{[+]} U_{[0]}$ above to write $g = u_1u_{1, [0]}nu_{2, [0]}u_2$ with $u_1, u_2 \in U_{[+]}$, $u_{1, [0]}, u_{2, [0]} \in U_{[0]}$ and $n \in N$, such that if we write $w = nT \in W$ then $u_2, u_{2, [0]} \in U_w$. We have $n.(u_{2, [0]} u_2.y) = {u_{1, [0]}}^{-1}{u_1}^{-1}.y'$; write ${u_1}' = {u_1}^{u_{1, [0]}}$ and ${u_2}' = {u_2}^{{u_{2, [0]}}^{-1}}$, so that ${u_1}', {u_2}' \in U_{[+]}$, and set $y_1 = u_{2, [0]}.y$ and ${y_1}' = {u_{1, [0]}}^{-1}.y'$, then we have $n.({u_2}'.y_1) = {{u_1}'}^{-1}.{y_1}'$. Let $\Delta$ be the set of weights occurring in $y_1$; by assumption $\Delta$ has ZLCE. Write $\Delta = \{ \nu_1, \dots, \nu_r \}$, and let $c_1, \dots, c_r \in \N$ be such that $c_1\nu_1 + \cdots + c_r\nu_r = 0$.

First suppose $X = V$. We have $y_1, {y_1}' \in V_{[0]}$; since adding a non-empty sum of positive roots to any $\nu_i$ gives a weight of positive generalized height, we see that both ${u_2}'.y_1 - y_1$ and ${{u_1}'}^{-1}.{y_1}' - {y_1}'$ must lie in $V_{[+]}$. Since all weights $\nu_i$ occur in $y_1$, they therefore occur in ${u_2}'.y_1$; thus all weights $w(\nu_i)$ occur in $n.({u_2}'.y_1)$, and as ${u_2}' \in U_w$ we see that each term in $n.({u_2}'.y_1) - n.y_1$ corresponds to a weight $\nu'$ such that for some $i$ we have $\nu' \prec w(\nu_i)$. As $n.({u_2}'.y_1) = {{u_1}'}^{-1}.{y_1}'$, each weight $w(\nu_i)$ occurs in ${{u_1}'}^{-1}.{y_1}'$, so lies in $\Lambda(V)_{[0]} \cup \Lambda(V)_{[+]}$; as $c_1w(\nu_1) + \cdots + c_rw(\nu_r) = w(c_1\nu_1 + \cdots + c_r\nu_r) = 0$, for each $i$ we must have $w(\nu_i) \in \Lambda(V)_{[0]}$. Since then $n.({u_2}'.y_1) - n.y_1 \in V_{[-]}$, we must have $n.({u_2}'.y_1) - n.y_1 = 0 = {{u_1}'}^{-1}.{y_1}' - {y_1}'$, so that ${u_1}' \in C_{U_{[+]}}({y_1}')$ and ${u_2}' \in C_{U_{[+]}}(y_1)$, while $n.y_1 = {y_1}'$. Since ${u_1}' = {u_1}^{u_{1, [0]}}$ and ${y_1}' = {u_{1, [0]}}^{-1}.y'$ we have $u_1 \in C_{U_{[+]}}(y')$; likewise as ${u_2}' = {u_2}^{{u_{2, [0]}}^{-1}}$ and $y_1 = u_{2, [0]}.y$ we have $u_2 \in C_{U_{[+]}}(y)$; and $u_{1, [0]} n u_{2, [0]}.y = y'$. Moreover, as $\Delta$ has ZLCE, for each $\nu \in \Lambda(V)_{[0]} \setminus \Delta$ there exist ${c_1}^*, \dots, {c_r}^*, c^* \in \N$ such that ${c_1}^*\nu_1 + \cdots + {c_r}^*\nu_r + c^*\nu = 0$, and so ${c_1}^*w(\nu_1) + \cdots + {c_r}^*w(\nu_r) + c^*w(\nu) = 0$; as each $w(\nu_i) \in \Lambda(V)_{[0]}$ we must also have $w(\nu) \in \Lambda(V)_{[0]}$. Thus $w \in W_{\Lambda(V)_{[0]}}$; since $N_{\Lambda(V)_{[0]}}$ normalises $G_{[0]}$ because $W_{\Lambda(V)_{[0]}}$ stabilizes $\Phi_{[0]}$, if we set $g' = u_{1, [0]} n u_{2, [0]}$ then $g' \in U_{[0]} N_{\Lambda(V)_{[0]}} U_{[0]} = G_{[0]} N_{\Lambda(V)_{[0]}}$ and $g'.y = y'$. The result follows.

Now suppose $X = \Gk(V)$ for some $k \in \N$. Write $y_1 = \langle x_1, \dots, x_k \rangle$; since each weight $\nu_i$ occurs in some basis vector $x_j$, by changing basis if necessary we may ensure that each $\nu_i$ occurs in each $x_j$. For each $j$ the argument of the previous paragraph now applies to the vectors $x_j$ and ${u_1}'n{u_2}'.x_j$; the result follows.
\end{proof}

Note that the assumption that $W_{\Lambda(V)_{[0]}}$ stabilizes $\Phi_{[0]}$ does not always hold. For example, let $G = A_6$ and $V = L(\omega_1)$; then $\Lambda(V) = \{ \omega_1, \omega_1 - \alpha_1, \omega_1 - \alpha_1 - \alpha_2, \dots, \omega_1 - \alpha_1 - \cdots - \alpha_6 \}$. If we take the generalized height function on the weight lattice of $G$ whose value at $\alpha_1$ and $\alpha_6$ is $0$ and at each other simple root $\alpha_i$ is $1$, then the generalized height of $\omega_1 = \frac{1}{7}(6 \alpha_1 + 5\alpha_2 + 4\alpha_3 + 3\alpha_4 + 2\alpha_5 + \alpha_6)$ is $2$, and as $\omega_1$ and $\Phi$ generate the weight lattice it follows that the generalized height of any weight is an integer; the generalized heights of the weights in $\Lambda(V)$ are $2$, $2$, $1$, $0$, $-1$, $-2$, $-2$, so $\Lambda(V)_{[0]} = \{ \omega_1 - \alpha_1 - \alpha_2 - \alpha_3 \}$, and hence $W_{\Lambda(V)_{[0]}} = \langle w_{\alpha_1}, w_{\alpha_2}, w_{\alpha_5}, w_{\alpha_6} \rangle$, which does not stabilize $\Phi_{[0]} = \langle \alpha_1, \alpha_6 \rangle$. However, we shall see that in the cases where we wish to apply Lemma~\ref{lem: gen height zero not strictly positive} the assumption does hold.

In the case of a strictly positive generalized height function we can go further; here of course $\Phi_{[0]}$ is empty so the assumption automatically holds.

\begin{lem}\label{lem: gen height zero}
Let $X$ be either $V$ or $\Gk(V)$ for some $k \in \N$. Given a strictly positive generalized height function on the weight lattice of $G$ such that $\Lambda(V)_{[0]}$ has ZLC, let $Y$ be a subvariety of either $V_{[0]}$ or $\Gk(V_{[0]})$. Suppose $y \in Y$ has the property that each weight in $\Lambda(V)_{[0]}$ occurs in $y$; take $g \in \Tran_G(y, Y)$ and set $y' = g.y \in Y$. Then we may write $g = u_1 n u_2$ with $u_1 \in C_U(y')$, $u_2 \in C_U(y)$, and $n \in N_{\Lambda(V)_{[0]}}$ with $n.y = y'$. In particular $G.y \cap Y = N_{\Lambda(V)_{[0]}}.y \cap Y$, and $C_G(y) = C_U(y) C_{N_{\Lambda(V)_{[0]}}}(y) C_U(y)$.
\end{lem}

\begin{proof}
Clearly if $\Lambda(V)_{[0]}$ has ZLC it has ZLCE. Since the generalized height function is strictly positive, we have $U_{[0]} = \{ 1 \}$ and $U_{[+]} = U$; the result thus follows from Lemma~\ref{lem: gen height zero not strictly positive}.
\end{proof}

Note that in the case where $X = V$, $\Lambda(V)_{[0]} = \{ 0 \}$ and $Y = V_0$ (the zero weight space), we conclude that two elements of $Y$ lie in the same $G$-orbit if and only if they lie in the same $N$-orbit, which as we said at the beginning of this section is the statement of \cite[Lemma~2.1]{GLMS}. It should be observed that Lemma~\ref{lem: gen height zero} reduces the often challenging problem of determining stabilizers in $G$ to the considerably simpler problems of identifying stabilizers in $U$ and in $N_{\Lambda(V)_{[0]}}$.

Many of the results of Section~\ref{sect: small triples and first quadruples}, where we determine generic stabilizers for small triples and associated first quadruples, will use this approach. In most cases the generalized height function chosen will be strictly positive, so we can apply Lemma~\ref{lem: gen height zero}. The details will of course vary from case to case, but typically we proceed as follows.

We identify the set $\Lambda(V)_{[0]}$, and show that it has ZLC and find its stabilizer $W_{\Lambda(V)_{[0]}}$; we set $Y = V_{[0]}$, choose a dense open subset $\hat Y$ of $Y$, and pick $y_0 \in \hat Y$. We define a subgroup $C$ of $G$ such that $C \leq C_G(y_0)$; our aim is to show that we have equality. Using $W_{\Lambda(V)_{[0]}}$ we determine $N_{\Lambda(V)_{[0]}}.y_0$ as a small number of cosets of $T.y_0$ (often it is just $T.y_0$ itself); we then show that $C_{N_{\Lambda(V)_{[0]}}}(y_0) = C \cap N$. Next we choose a subset $\Xi$ of $\Phi^+$ and set $U' = \prod_{\alpha \in \Xi} X_\alpha$ such that $U'$ is a complement to $C \cap U$ in $U$; by considering sums $\nu_i + \alpha$ for $\nu_i \in \Lambda(V)_{[0]}$ and $\alpha \in \Xi$ we argue that $C_{U'}(y_0) = \{ 1 \}$, whence $C_U(y_0) = C \cap U$. Lemma~\ref{lem: gen height zero} now shows that we do indeed have $C_G(y_0) = C$, and identifies $G.y_0 \cap Y$. Finally given an arbitrary $y \in \hat Y$ we observe that there exists $h \in T$ with $C_G(y) = {}^h C$, and find $G.y \cap Y$; by comparing dimensions we see that $y$ is $Y$-exact, and then Lemma~\ref{lem: generic stabilizer from exact subset} gives the result for the triple $(G, \lambda, p)$. Indeed, as we explain at the start of Chapter~\ref{chap: non-TGS triples and first quadruples}, we may simultaneously obtain the result for the associated first quadruple $(G, \lambda, p, 1)$, by taking similarly a subgroup $C'$ of $G$ such that $C' \leq C_G(\langle y_0 \rangle)$ and arguing in parallel to show that with $y$ and $h$ as above we have $C_G(\langle y \rangle) = {}^h C'$.

In fact, in Section~\ref{sect: small triples and first quadruples} there is just one case where we have to work with a generalized height function which is not strictly positive, and thus cannot apply Lemma~\ref{lem: gen height zero} but must instead employ the more complicated result Lemma~\ref{lem: gen height zero not strictly positive}. However, when we deal with higher quadruples in Section~\ref{sect: small higher quadruples}, there are several instances where we make use of Lemma~\ref{lem: gen height zero not strictly positive}.

\section{Reduction from higher Grassmannian varieties}\label{sect: reduction}

In this section we give a general result which links the generic stabilizer for the action of $G$ on a Grassmannian variety $\Gk(V)$ with $k > 1$ to that for the action of a larger group on an appropriate projective space. The basic idea may be found in \cite[Proposition~3.2]{GLMS}, but that result is stated only in terms of the existence of finitely many orbits.

\begin{lem}\label{lem: reduction to projective space}
Let $V$ be a $G$-module; take $k > 1$, and let $V_{nat}$ be the natural module for the group $A_{k - 1}$, so that $G \times A_{k - 1}$ acts on the module $V \otimes V_{nat}$. If there is a generic stabilizer in the action of $G \times A_{k - 1}$ on $\G{1}(V \otimes V_{nat})$, then there is also one in the action of $G$ on $\G{k}(V)$, and the two are isomorphic.
\end{lem}

\begin{proof}
Write $G^+ = G \times A_{k - 1}$, and let $\pi_1 : G^+ \to G$ be the projection on the first component; write $V^+ = V \otimes V_{nat}$. Set $X = \Gk(V)$ and $X^+ = \G{1}(V^+)$; recall that $v_1, \dots, v_k$ is the natural basis of $V_{nat}$. Each element of $X^+$ is of the form $\langle v^+ \rangle$ where $v^+ = \sum_{i = 1}^k v^{(i)} \otimes v_i \in V^+$ for some $v^{(1)}, \dots, v^{(k)} \in V$; let $\tilde X^+$ be the dense open subset of $X^+$ consisting of such points where the vectors $v^{(1)}, \dots, v^{(k)}$ are linearly independent. We may define a surjection $\psi : \tilde X^+ \to X$ sending such a point $\langle \sum_{i = 1}^k v^{(i)} \otimes v_i \rangle$ to $\langle v^{(1)}, \dots, v^{(k)} \rangle$.

Take $x^+ = \langle v^+ \rangle \in \tilde X^+$ where $v^+ = \sum_{i = 1}^k v^{(i)} \otimes v_i$; let $x = \langle v^{(1)}, \dots, v^{(k)} \rangle = \psi(x^+)$. Given $g^+ = (g, a) \in C_{G^+}(x^+)$, for $i = 1, \dots, k$ write $a.v_i = \sum_{j = 1}^k c_{ij}v_j$; let $(d_{ij})$ be the inverse of the matrix $(c_{ij})$. We have
$$
g^+.v^+ = \sum_{i = 1}^k (g.v^{(i)}) \otimes (a.v_i) = \sum_{i, j = 1}^k (g.v^{(i)}) \otimes c_{ij} v_j = \sum_{j = 1}^k g.\left( \sum_{i = 1}^k c_{ij} v^{(i)} \right) \otimes v_j.
$$
Since $g^+.x^+ = x^+$ there exists $\kappa \in K^*$ with $g^+.v^+ = \kappa v^+$; thus for $j = 1, \dots, k$ we have $g.\left(\sum_{i = 1}^k c_{ij} v^{(i)}\right) = \kappa v^{(j)}$, and so $g.v^{(i)} = \kappa \sum_{j = 1}^k d_{ji} v^{(j)}$. Since $g$ maps each $v^{(i)}$ into $\langle v^{(1)}, \dots, v^{(k)} \rangle$, we have $g \in C_G(x)$. Conversely given $g \in C_G(x)$, for $i = 1, \dots, k$ write $g.v^{(i)} = \sum_{j = 1}^k d_{ji} v^{(j)}$. Take $\kappa \in K^*$ satisfying $\kappa^k = \det (d_{ij})$, and define $a \in A_{k - 1}$ by $a^{-1}.v_j = \kappa^{-1} \sum_{i = 1}^k d_{ji} v_i$; then with $g^+ = (g, a) \in G^+$ we have
$$
g^+.v^+ = \sum_{i, j = 1}^k d_{ji} v^{(j)} \otimes (a.v_i) = \sum_{j = 1}^k v^{(j)} \otimes a.\left( \sum_{i = 1}^k d_{ji}v_i \right) = \sum_{j = 1}^k v^{(j)} \otimes \kappa v_j = \kappa v^+,
$$
and so $g^+.x^+ = x^+$, i.e., $g^+ \in C_{G^+}(x^+)$. Thus $C_G(x) = \pi_1(C_{G^+}(x^+))$.

Now suppose the action of $G^+$ on $X^+$ has generic stabilizer $C/Z(G^+)$. Let $\hat X^+$ be a dense open subset of $X^+$ all of whose points have $G^+$-stabilizer equal to a $G^+$-conjugate of $C$; by replacing $\hat X^+$ by its intersection with $\tilde X^+$ we may assume each point of $\hat X^+$ is of the form $\langle \sum_{i = 1}^k v^{(i)} \otimes v_i \rangle$ with $\langle v^{(1)}, \dots, v^{(k)} \rangle \in X$. As $\hat X^+$ is a dense open subset of $\tilde X^+$, we see that $\psi(\hat X^+)$ is a dense open subset $\hat X$ of $X$. By the above, for all $x \in \hat X$ there exists $h^+ \in G^+$ such that $C_G(x) = \pi_1(C^{h^+}) = \pi_1(C)^{\pi_1(h^+)}$. Thus all points in $\hat X$ have $G$-stabilizer a $G$-conjugate of $\pi_1(C)$; so the action of $G$ on $X$ has generic stabilizer $\pi_1(C)/Z(G)$.

It remains to determine $\pi_1(C)$. Take $v^+ = \sum_{i = 1}^k v^{(i)} \otimes v_i \in V^+$ such that $\langle v^+ \rangle$ has $G^+$-stabilizer $C$. If $g^+ \in C \cap \ker \pi_1$, then $g^+ = (1, a)$ for some $a \in A_{k - 1}$, and there exists $\kappa \in K^*$ with $\sum_{i = 1}^k v^{(i)} \otimes a.v_i = \kappa \sum_{i = 1}^k v^{(i)} \otimes v_i$; as the $v^{(i)}$ are linearly independent, for all $i$ we must have $a.v_i = \kappa v_i$, so that $a \in Z(A_{k - 1})$. Therefore $C \cap \ker \pi_1 \subseteq \{ (1, a) : a \in Z(A_{k - 1}) \}$; as the reverse inclusion is obvious, we have $\pi_1(C)/Z(G) \cong C/Z(G \times A_{k - 1}) = C/Z(G^+)$. The result follows.
\end{proof}

Note that although this result is stated for the action of the direct product $G \times A_{k - 1}$, it is harmless to replace $G \times A_{k - 1}$ by any quotient by a subgroup of the centre, since in any action on a Grassmannian variety the kernel is the full centre of the group concerned.

Although the result holds generally, we shall apply it only in certain situations, where there is a simple algebraic group $H$ having a maximal rank subgroup $GA_{k - 1}$. In the cases concerned, we shall see that provided $(p, k) = 1$ there is a semisimple element of $H$ whose centralizer is $GA_{k - 1}$; we may then apply Lemma~\ref{lem: semisimple auts}.

We conclude this section by observing that Lemma~\ref{lem: reduction to projective space} links the existence of generic stabilizers for a higher quadruple and a related first quadruple, but proves an implication in one direction only. In fact the converse is false in general, as the following counterexample shows. Consider the higher quadruple $(A_7, \omega_2, 2, 2)$. In the proof of Proposition~\ref{prop: A_ell, ell odd, omega_2 module, k = 2} we shall obtain a family of $2$-dimensional subspaces $\langle {v^{(1)}}', {v^{(2)}}' \rangle$ of the $A_7$-module $V = L(\omega_2)$, each having $G$-stabilizer equal to $A \langle n^*, n^{**} \rangle$ where $A$ is a fixed ${A_1}^4$ subgroup and $n^*, n^{**}$ are fixed commuting involutions in $N$; using Lemma~\ref{lem: generic stabilizer from exact subset} we shall conclude that the generic stabilizer is ${A_1}^4.{\Z_2}^2$. Each such $2$-dimensional subspace is determined by a single parameter $a$ which may take any value in $K \setminus \{ 0, 1 \}$; on the subspace parametrized by $a$, the subgroup $A$ acts trivially, while with respect to the basis ${v^{(1)}}', {v^{(2)}}'$ the elements $n^*$ and $n^{**}$ act as the matrices
$$
J =
\left(
  \begin{array}{cc}
    0 & 1 \\
    1 & 0 \\
  \end{array}
\right)
\quad \hbox{and} \quad
M_a =
\left(
  \begin{array}{cc}
    a + 1 &   a   \\
      a   & a + 1 \\
  \end{array}
\right)
$$
respectively. We have $G^+ = G \times A_1$ and $V^+ = V \otimes V_{nat}$. The preimage under the map $\psi$ of the subspace $\langle {v^{(1)}}', {v^{(2)}}' \rangle$ contains the line $\langle {v^{(1)}}' \otimes v_1 + {v^{(2)}}' \otimes v_2 \rangle$, and the union of the $G^+$-orbits containing these lines contains a dense open subset of $\G{1}(V^+)$. If we take the subspace parametrized by $a$ and consider the corresponding line, its $G^+$-stabilizer is isomorphic to ${A_1}^4.{\Z_2}^2$; the connected component is $A \times \{ 1 \}$, while corresponding to $n^*$ and $n^{**}$ we have elements $(n^*, J)$ and $(n^{**}, M_a)$. Thus the projection of the $G^+$-stabilizer on the second factor is $\langle J, M_a \rangle \cong {\Z_2}^2$. If infinitely many of these ${\Z_2}^2$ subgroups of $A_1$ were conjugate, then certainly there would exist $a, b \in K \setminus \{ 0, 1 \}$ distinct and a conjugating element of $A_1$ which fixed $J$ and sent $M_a$ to $M_b$; but the $A_1$-centralizer of $J$ is equal to
$$
\left\{
\left(
  \begin{array}{cc}
    \kappa + 1 &   \kappa   \\
      \kappa   & \kappa + 1 \\
  \end{array}
\right)
: \kappa \in K \right\},
$$
which is an abelian group containing $M_a$, so no such conjugating element can exist. It follows that the first quadruple $(A_7 \times A_1, \omega_2 \otimes \omega_1, 2, 1)$ has no generic stabilizer (but there is a semi-generic stabilizer ${A_1}^4.{\Z_2}^2$). 
\chapter{Triples and first quadruples not having TGS}\label{chap: non-TGS triples and first quadruples}

In this chapter we consider triples and first quadruples which do not have TGS, and establish the entries in Tables~\ref{table: large triple and first quadruple non-TGS}, \ref{table: small classical triple and first quadruple generic stab} and \ref{table: small exceptional triple and first quadruple generic stab}. In Sections~\ref{sect: non-TGS large triples and first quadruples} and \ref{sect: small triples and first quadruples} we treat triples which are large and small respectively, together with the associated first quadruples. Throughout, given a triple $(G, \lambda, p)$ or quadruple $(G, \lambda, p, 1)$ we write $V = L(\lambda)$.

In many cases our approach will be to apply Lemma~\ref{lem: generic stabilizer from exact subset} to determine the required generic stabilizers. We consider the action of $G$ on the module $V = L(\lambda)$; we choose a subspace $Y$ of $V$, and take a dense open subset $\hat Y$ of $Y$ which is closed under taking non-zero scalar multiples. For all $y \in \hat Y$, we show that the stabilizers $C_G(y)$ and $C_G(\langle y \rangle)$ in the actions on $V$ and $\G{1}(V)$ are conjugates of fixed subgroups $C$ and $C'$ respectively, and that $y$ is $Y$-exact; since clearly $\Tran_G(\langle y \rangle, \G{1}(Y)) = \Tran_G(y, Y)$ and $\codim \G{1}(Y) = \dim \G{1}(V) - \dim \G{1}(Y) = \dim V - \dim Y = \codim Y$, it follows that $\langle y \rangle$ is $\G{1}(Y)$-exact. By Lemma~\ref{lem: generic stabilizer from exact subset} applied to $Y$ and $\G{1}(Y)$ we may now conclude that the triple $(G, \lambda, p)$ and the associated first quadruple $(G, \lambda, p, 1)$ have generic stabilizer $C/G_V$ and $C'/Z(G)$ respectively.

\section{Large triples and associated first quadruples}\label{sect: non-TGS large triples and first quadruples}

In this section we shall treat the large triples not dealt with in Chapter~\ref{chap: TGS triples} and the associated first quadruples, and establish the entries in Table~\ref{table: large triple and first quadruple non-TGS}, thus proving Theorem~\ref{thm: large triple and first quadruple generic stab}. We shall begin with the $p$-restricted triples and associated first quadruples, including the two cases treated in Propositions~\ref{prop: C_4, omega_3, p = 3, nets} and \ref{prop: B_2, omega_1 + omega_2, p = 5, nets}, where we showed that the triples have TGS but were unable to deduce the same of the associated first quadruples; we shall conclude by considering the two cases where the triple is not $p$-restricted.

The first three results in this section will be proved using the approach of Section~\ref{sect: semisimple auts}. Recall that we take a simple algebraic group $H$ of simply connected type over $K$, and let $\theta$ be a semisimple automorphism of $H$ of order $r$ coprime to $p$; then $\theta$ acts on $\L(H)$, and for $i = 0, 1, \dots, r - 1$ we write $\L(H)_{(i)}$ for the eigenspace corresponding to the eigenvalue ${\eta_r}^i$. Writing $Z(\L(H))_{(1)} = Z(\L(H)) \cap \L(H)_{(1)}$, we consider the action of the group $C_H(\theta)$ on the module $\L(H)_{(1)}/Z(\L(H))_{(1)}$.

We shall group together triples $(G, \lambda, p)$ and associated first quadruples for which the arguments are similar, although the details may vary. We shall begin by choosing $H$, and shall in fact give two semisimple automorphisms $\theta_1$ and $\theta_2$, which we shall show are conjugate. Taking $\theta = \theta_1$, we shall see that we may take $G = C_H(\theta)$ and $V = \L(H)_{(1)}/Z(\L(H))_{(1)}$, and using Lemma~\ref{lem: semisimple auts}(i) we shall prove the existence of regular orbits. Taking $\theta = \theta_2$ and using Lemma~\ref{lem: semisimple auts}(ii) we shall then determine the generic stabilizers.

We start with those cases in which $\theta$ is an inner automorphism, when we shall identify it with an element of $H$. In these cases the assumption that $r$ is coprime to $p$ excludes one choice of $p$ requiring treatment; we shall postpone dealing with these triples and first quadruples until later in this section.

\begin{prop}\label{prop: A_8, omega_3, A_7, omega_4, D_8, omega_8 modules, non-special characteristic}
Let $G = A_7$ and $\lambda = \omega_4$ with $p \geq 3$, or $G = D_8$ and $\lambda = \omega_8$ with $p \geq 3$, or $G = A_8$ and $\lambda = \omega_3$ with $p \neq 3$. Then the triple $(G, \lambda, p)$ has generic stabilizer ${\Z_2}^6$, or ${\Z_2}^8$, or ${\Z_3}^4.\Z_{(p, 2)}$, respectively, and there is a regular orbit; the associated first quadruple $(G, \lambda, p, 1)$ has generic stabilizer ${\Z_2}^6$, or ${\Z_2}^8$, or ${\Z_3}^4.\Z_2$, respectively.
\end{prop}

\begin{proof}
Number the cases (i), (ii) and (iii) according as $G = A_7$, $D_8$ or $A_8$; whenever we give three choices followed by the word `respectively' we are taking the cases in the order (i), (ii), (iii).

Let $H$ be the simply connected group defined over $K$ of type $E_7$, $E_8$ or $E_8$ respectively (so that in each case $\ell_H = \ell$), with simple roots $\beta_1, \dots, \beta_\ell$, and let $r$ be $2$, $2$ or $3$ respectively; assume $p \neq r$. We have $Z(\L(H)) = \{ 0 \}$.

Define $\theta_1 \in T_H$ to be
$$
\begin{array}{ll}
h_{\beta_1}(-1) h_{\beta_2}(\eta_4) h_{\beta_3}(-1) h_{\beta_5}(-\eta_4) h_{\beta_7}(-\eta_4)                   & \vstrut \hbox{in case~(i),} \\
h_{\beta_3}(-1) h_{\beta_4}(-1) h_{\beta_7}(-1) h_{\beta_8}(-1)                                                 & \vstrut \hbox{in case~(ii),} \\
h_{\beta_1}(\eta_3) h_{\beta_2}({\eta_3}^2) h_{\beta_3}(\eta_3) h_{\beta_5}({\eta_3}^2) h_{\beta_8}({\eta_3}^2) & \vstrut \hbox{in case~(iii).}
\end{array}
$$
Then $\langle {\theta_1}^r \rangle = Z(H)$, and $\theta_1$ sends $x_\alpha(t)$ to $x_\alpha({\eta_r}^{\height(\alpha)}t)$; so $X_\alpha < C_H(\theta_1)$ if and only if $\height(\alpha) \equiv 0$ (mod $r$). It follows that $C_H(\theta_1)$ is a connected group of type $A_7$, $D_8$ or $A_8$ respectively, with simple root elements $x_{\alpha_i}(t)$, where $\alpha_1, \dots, \alpha_\ell$ are
$$
\begin{array}{ll}
\esevenrt0111100, \esevenrt0000011, \esevenrt0001100, \esevenrt1010000, \esevenrt0101000, \esevenrt0000110, \esevenrt0011000                           & \vstrut \hbox{in case~(i),} \\
\eeightrt01111000, \eeightrt00000110, \eeightrt00011000, \eeightrt10100000, \eeightrt01010000, \eeightrt00001100, \eeightrt00110000, \eeightrt00000011 & \vstrut \hbox{in case~(ii),} \\
\eeightrt00011100, \eeightrt01110000, \eeightrt00001110, \eeightrt10110000, \eeightrt01011000, \eeightrt00000111, \eeightrt00111000, \eeightrt11111100 & \vstrut \hbox{in case~(iii);}
\end{array}
$$
in each case we see that $Z(C_H(\theta_1)) = \langle \theta_1 \rangle$.

Now let $\delta_1, \dots, \delta_\ell$ be
$$
\begin{array}{ll}
\esevenrt0010000, \esevenrt0100000, \esevenrt0000100, \esevenrt0112100, \esevenrt0000001, \esevenrt0112221, \esevenrt2234321                           & \vstrut \hbox{in case~(i),} \\
\eeightrt00100000, \eeightrt01000000, \eeightrt00001000, \eeightrt01121000, \eeightrt00000010, \eeightrt01122210, \eeightrt22343210, \eeightrt23465432 & \vstrut \hbox{in case~(ii),} \\
\eeightrt10000000, \eeightrt00100000, \eeightrt00000100, \eeightrt00001000, \eeightrt01000000, \eeightrt11232100, \eeightrt00000001, \eeightrt23465431 & \vstrut \hbox{in case~(iii);}
\end{array}
$$
then $\langle \delta_1, \dots, \delta_\ell \rangle$ is a subsystem of type ${A_{r - 1}}^{\ell/(r - 1)}$, i.e., ${A_1}^7$, ${A_1}^8$ or ${A_2}^4$ respectively. Set $\theta_2 = n_{\delta_1} \dots n_{\delta_\ell}$; then $\langle {\theta_2}^r \rangle = Z(H)$, and indeed in cases (i) and (ii) the element of $W_H$ corresponding to $\theta_2$ is the long word. We find that $\theta_2$ acts fixed-point-freely on both $\Phi_H$ and $\L(T_H)$ (these are now obvious in cases~(i) and (ii), while in case~(iii) they are both easy calculations). Thus $\dim C_{\L(H)}(\theta_2) = |\Phi_H|/r = 63$, $120$ or $80$ respectively; the classification of semisimple elements of $H$ (see e.g. \cite[Tables~4.3.1 and 4.7.1]{GLS}) now shows that $\theta_2$ must be a conjugate of $\theta_1$.

First set $\theta = \theta_1$; then we may take $G = C_H(\theta)$. We see that $e_\alpha \in \L(H)_{(i)}$ if and only if $\height(\alpha) \equiv i$ (mod $r$). Thus in $\L(H)_{(1)}$ we have a highest weight vector $e_\beta$ for $\beta = \esevenrt2234321$, $\eeightrt23465432$ or $\eeightrt23465431$ respectively; the expressions above for the simple root elements of $G$ show that $\L(H)_{(1)}$ is the Weyl $G$-module with high weight $\omega_4$, $\omega_8$ or $\omega_3$ respectively. We have $Z(\L(H))_{(1)} = Z(\L(H)) = \{ 0 \}$; we may take $V = \L(H)_{(1)}/Z(\L(H))_{(1)}$, and then $G_V = Z(H)$.

Take $v = e_{\beta_1} + \cdots + e_{\beta_\ell} \in \L(U_H) \cap \L(H)_{(1)}$; then $v$ is a regular nilpotent element. From Lemma~\ref{lem: m_i for H} we see that there are natural numbers $m_1, \dots, m_\ell$, which are listed there, such that we may write
$$
C_{U_H}(v) = \{ y_1(c_1) \dots y_\ell(c_\ell) : c_i \in K \},
$$
with each $y_i(c)$ of the form $\left(\prod_{\height(\alpha) = m_i} x_\alpha(n_\alpha c)\right)x$, where $x$ is a product of root elements corresponding to roots of height greater than $m_i$, and the $n_\alpha \in K$ are not all zero and satisfy $\sum_{\height(\alpha) = m_i} n_\alpha e_\alpha \in C_{\L(G)}(v)$. Since no $m_i$ is divisible by $r$, we have $G \cap C_{U_H}(v) = \{ 1 \}$; since $G_V = Z(H) = G \cap Z(H)$, Lemma~\ref{lem: semisimple auts}(i) shows that the orbit containing $v + Z(\L(H))_{(1)}$ is regular.

Now set $\theta = \theta_2$, and again take $G = C_H(\theta)$ and $V = \L(H)_{(1)}/Z(\L(H))_{(1)}$. We have $G \cap T_H = C_{T_H}(\theta) \cong {\Z_r}^{\ell/(r - 1)}$, i.e., ${\Z_2}^7$, ${\Z_2}^8$ or ${\Z_3}^4$ respectively (this is obvious in cases~(i) and (ii) as then $\theta$ acts on $T_H$ as inversion; in case~(iii) we have $C_{T_H}(\theta) = \langle h_{\beta_1}(\eta_3) h_{\beta_3}({\eta_3}^2), h_{\beta_5}({\eta_3}^2) h_{\beta_6}(\eta_3), h_{\beta_1}(\eta_3) h_{\beta_2}(\eta_3) h_{\beta_6}(\eta_3), h_{\beta_1}({\eta_3}^2) h_{\beta_6}(\eta_3) h_{\beta_8}(\eta_3) \rangle$). Moreover, in cases~(i) and (ii) we have $\L(T_H)_{(1)} = \L(T_H)$, while in case~(iii) we find that $\L(T_H)_{(1)} = \langle h_{\delta_{2i - 1}} - \eta_3 h_{\delta_{2i}} : i = 1, 2, 3, 4 \rangle$. Thus $\dim \L(H)_{(1)} - \dim \L(T_H)_{(1)} = \dim G - \dim (G \cap T_H)$ in each case; and in case (iii) a routine check shows that $\L(T_H)_{(1)}$ contains regular semisimple elements.

By Lemma~\ref{lem: W_H on L(T_H)/Z(L(H))}, in cases~(i) and (ii) we have $({W_H}^\ddagger)_{(1)} = {W_H}^\ddagger = \langle w_0 \rangle = \langle \theta T_H \rangle$; we claim that in case~(iii) we have $({W_H}^\ddagger)_{(1)} = \langle \theta T_H, w_0 \rangle$. Thus suppose $w \in W_H$ and there exists $\xi \in K^*$ such that for all $y \in \L(T_H)_{(1)}$ we have $w.y = \xi y$. For $i = 1, 2, 3, 4$ write $\Psi_i = \langle \delta_{2i - 1}, \delta_{2i} \rangle$ and $y_i = h_{\delta_{2i - 1}} - \eta_3 h_{\delta_{2i}} \in \L(T_H)_{(1)}$. Taking $y = y_1$ and arguing as at the end of Section~\ref{sect: semisimple auts} shows that $w(\beta_1)$ and $w(\beta_3)$ must be proportional outside $\{ \beta_1, \beta_3 \}$, and as $\eta_3 \neq \pm1$ that $w$ must preserve $\Psi_1$. Now take $i \in \{ 2, 3, 4 \}$. There exists $w' \in W_H$ with $w'(\delta_1) = \delta_{2i - 1}$ and $w'(\delta_2) = \delta_{2i}$, and so $w'.y_1 = y_i$, whence $w.y_i = \xi y_i$ gives $w^{w'}.y_1 = y_1$; by the above $w^{w'}$ preserves $\Psi_1$, so $w$ preserves $\Psi_i$. Thus $w = w_1w_2w_3w_4.{w_0}^j$ where each $w_i$ lies in $W(\Psi_i)$ and $j \in \{ 0, 1 \}$. For each $i$, the three elements in $W(\Psi_i)$ of odd length send $y_i$ to a scalar multiple of $\eta_3 h_{\delta_{2i - 1}} - h_{\delta_{2i}}$, so we must have $w_i \in \langle w_{\delta_{2i - 1}} w_{\delta_{2i}} \rangle$; since $w$ must multiply each of the four vectors $y_i$ by the same scalar, we must have $w_1w_2w_3w_4 \in \langle \theta T_H \rangle$, so that $w \in \langle \theta T_H, w_0 \rangle$ as required. Note that in this case if we write $n_0 = n_{\delta_1} \dots n_{\delta_8}$ where $\delta_1, \dots, \delta_8$ are as in case~(ii), then $n_0$ is an involution in $N_H$ corresponding to $w_0$ which commutes with $\theta$.

Now $w_0$ acts on $\L(T_H)_{(1)}$ as negation, and in case~(iii) $\theta$ acts on $\L(T_H)_{(1)}$ as multiplication by $\eta_3$. Thus if we are in case~(iii) with $p = 2$ then $({W_H}^\dagger)_{(1)} = \langle w_0 \rangle$, and so $C_{({N_H}^\dagger)_{(1)}}(\theta) = C_{T_H}(\theta) \langle n_0 \rangle$; if instead we are in case~(iii) with $p \geq 3$, or in case~(i) or (ii), then $({W_H}^\dagger)_{(1)} = \{ 1 \}$, and so $C_{({N_H}^\dagger)_{(1)}}(\theta) = C_{T_H}(\theta)$. Also in cases~(i) and (ii) we have $({N_H}^\ddagger)_{(1)} = T_H \langle \theta \rangle$, so $C_{({N_H}^\ddagger)_{(1)}}(\theta) = C_{T_H}(\theta) \langle \theta \rangle$, while in case~(iii) we have $({N_H}^\ddagger)_{(1)} = T_H \langle \theta, n_0 \rangle$, so $C_{({N_H}^\ddagger)_{(1)}}(\theta) = C_{T_H}(\theta) \langle \theta, n_0 \rangle$. Since $G_V = Z(H) \cong \Z_2$, $\{ 1 \}$ or $\{ 1 \}$ respectively, and $Z(G) = \langle \theta \rangle$ in each case, Lemma~\ref{lem: semisimple auts}(ii) shows that the triple $(G, \lambda, p)$ has generic stabilizer $C_{({N_H}^\dagger)_{(1)}}(\theta)/G_V \cong {\Z_2}^6$, or ${\Z_2}^8$, or ${\Z_3}^4.\Z_{(p, 2)}$, respectively, while the quadruple $(G, \lambda, p, 1)$ has generic stabilizer $C_{({N_H}^\ddagger)_{(1)}}(\theta)/Z(G) \cong {\Z_2}^6$, or ${\Z_2}^8$, or ${\Z_3}^4.\Z_2$, respectively.
\end{proof}

We now turn to the cases where $\theta$ is an outer automorphism. We begin with two individual cases.

\begin{prop}\label{prop: A_2, 3omega_1, C_4, omega_4 modules}
Let $G = A_2$ and $\lambda = 3\omega_1$ with $p \geq 5$, or $G = C_4$ and $\lambda = \omega_4$ with $p \geq 3$. Then the triple $(G, \lambda, p)$ has generic stabilizer ${\Z_3}^2$ or ${\Z_2}^6$ respectively, and there is a regular orbit; the associated first quadruple $(G, \lambda, p, 1)$ has generic stabilizer ${\Z_3}^2.\Z_2$ or ${\Z_2}^6$ respectively.
\end{prop}

\begin{proof}
Number the cases (i) and (ii) according as $G = A_2$ or $C_4$; whenever we give two choices followed by the word `respectively' we are taking the cases in the order (i), (ii).

Let $H$ be the simply connected group of type $D_4$ or $E_6$ respectively over $K$, with simple roots $\beta_1, \dots, \beta_{\ell_H}$, and set $r = 3$ or $2$ respectively; assume $p > r$. We have $Z(\L(H)) = \{ 0 \}$ unless we are in case~(ii) with $p = 3$, in which case $Z(\L(H)) = \langle h_{\beta_1} - h_{\beta_3} + h_{\beta_5} - h_{\beta_6} \rangle$.

Let $\tau$ be the automorphism of $\Phi_H$ which preserves $\Pi_H$ and permutes simple roots as follows:
$$
\begin{array}{ll}
(\dfourrt1000 \ \dfourrt0010 \ \dfourrt0001) (\dfourrt0100)                                    & \vstrut \hbox{in case~(i),} \\
(\esixrt100000 \ \esixrt000001)(\esixrt001000 \ \esixrt000010) (\esixrt010000) (\esixrt000100) & \vstrut \hbox{in case~(ii).}
\end{array}
$$
We claim that we may assume that the isomorphisms $x_\alpha : K \to X_\alpha$ are chosen such that the structure constants are preserved by $\tau$, i.e., for all $\alpha, \beta \in \Phi_H$ we have $N_{\alpha, \beta} = N_{\tau(\alpha), \tau(\beta)}$. The map $x_\alpha(t) \mapsto x_{\tau(\alpha)}(t)$ then gives rise to a graph automorphism of $H$, which by slight abuse of notation we also call $\tau$.

To obtain these structure constants we use the method explained in \cite[4.2]{Car1}. We begin by defining a total ordering on the set of positive roots as follows: we take $j_1, \dots, j_{\ell_H} = 2, 4, 3, 1$ or $2, 4, 5, 3, 6, 1$ respectively, and then given two positive roots $\sum a_i \beta_i$ and $\sum b_i \beta_i$ we say that $\sum a_i \beta_i$ precedes $\sum b_i \beta_i$ if there exists $i'$ such that for $i < i'$ we have $a_{j_i} = b_{j_i}$, while $a_{j_{i'}} < b_{j_{i'}}$. The choice of the $j_i$ means that this total ordering respects $\tau$-orbits, in the sense that no two roots in the same $\tau$-orbit are separated by a root in a different $\tau$-orbit. The total ordering then determines a set of extraspecial pairs $(\alpha, \beta)$, one for each non-simple positive root. We define $N_{\alpha, \beta} = 1$ for each such extraspecial pair; this then suffices to determine the full collection of structure constants, and a direct check (or a proof using induction on the height of a root) shows that it has the property of preservation by $\tau$ stated above. In case~(i) we find that the pairs $(\alpha, \beta)$ for which $N_{\alpha, \beta} = 1$ are
$$
(\dfourrt1000, \dfourrt0100), \ (\dfourrt1000, \dfourrt0110), \ (\dfourrt1000, \dfourrt0101), \ (\dfourrt1000, \dfourrt0111), \ (\dfourrt1110, \dfourrt0101), \ (\dfourrt0100, \dfourrt1111)
$$
together with their images under $\tau$ and $\tau^2$. In case~(ii) there are too many pairs to list conveniently, but the structure constants may be obtained from those given in the appendix of \cite{LSmax} by negating the root vectors $e_\alpha$ for the following roots $\alpha$:
\begin{eqnarray*}
& \esixrt000011, \esixrt000110, \esixrt010100, \esixrt011100, \esixrt001110, \esixrt111100, \esixrt010111, \esixrt101110, \esixrt011111, & \\
& \esixrt011211, \esixrt111111, \esixrt011221, \esixrt111211, \esixrt112211, \esixrt111221, \esixrt112221, \esixrt112321, \esixrt122321. &
\end{eqnarray*}

Define $\theta_1$ to be
$$
\begin{array}{ll}
\tau h_{\beta_2}({\eta_3}^2)                                         & \vstrut \hbox{in case~(i),} \\
\tau h_{\beta_2}(-1) h_{\beta_3}(-1) h_{\beta_4}(-1) h_{\beta_5}(-1) & \vstrut \hbox{in case~(ii).}
\end{array}
$$
Then ${\theta_1}^r = 1$, and $\theta_1$ sends $x_\alpha(t)$ to $x_{\tau(\alpha)}({\eta_r}^{\height(\alpha)}t)$; so if $\tau(\alpha) = \alpha$ then $X_\alpha < C_H(\theta_1)$ if and only if $\height(\alpha) \equiv 0$ (mod $r$), while if $\tau(\alpha) \neq \alpha$ then the intersection of $C_H(\theta_1)$ with $X_\alpha X_{\tau(\alpha)} X_{\tau^2(\alpha)}$ or $X_\alpha X_{\tau(\alpha)}$ is the $1$-dimensional group
$$
\left\{ x_\alpha(t) x_{\tau(\alpha)}({\eta_3}^{\height(\alpha)} t) x_{\tau^2(\alpha)}({\eta_3}^{2\height(\alpha)} t) : t \in K \right\}
$$
or
$$
\left\{ x_\alpha(t) x_{\tau(\alpha)}((-1)^{\height(\alpha)} t) : t \in K \right\}
$$
respectively. It follows that $C_H(\theta_1)$ is a group of type $A_2$ or $C_4$ respectively, with simple root elements
$$
x_{\alpha_1}(t) = x_{\dfourrt1000}(t) x_{\dfourrt0010}(\eta_3 t) x_{\dfourrt0001}({\eta_3}^2 t), \quad
x_{\alpha_2}(t) = x_{\dfourrt1100}(t) x_{\dfourrt0110}({\eta_3}^2 t) x_{\dfourrt0101}(\eta_3 t)
$$
or
\begin{eqnarray*}
x_{\alpha_1}(t) = x_{\esixrt001100}(t) x_{\esixrt000110}(t), \phantom{(-t)} & &
x_{\alpha_2}(t) = x_{\esixrt100000}(t) x_{\esixrt000001}(-t), \\
x_{\alpha_3}(t) = x_{\esixrt001000}(t) x_{\esixrt000010}(-t), \phantom{(t)} & &
x_{\alpha_4}(t) = x_{\esixrt010100}(t)
\end{eqnarray*}
respectively; in each case we see that $Z(C_H(\theta_1)) = \{ 1 \}$.

Now define $\theta_2$ to be
$$
\begin{array}{ll}
\tau n_{\dfourrt1000} n_{\dfourrt0010} n_{\dfourrt0001} n_{\dfourrt1100} n_{\dfourrt0110} n_{\dfourrt0101} & \vstrut \hbox{in case~(i),} \\
\tau n_{\esixrt000100} n_{\esixrt001110} n_{\esixrt101111} n_{\esixrt122321}                               & \vstrut \hbox{in case~(ii);} \\
\end{array}
$$
then ${\theta_2}^r = 1$, and indeed in case~(ii) $\theta_2$ sends each root subgroup $X_\alpha$ to $X_{-\alpha}$. We find that $\theta_2$ acts fixed-point-freely on both $\Phi_H$ and $\L(T_H)$ (these are now obvious in case~(ii), while in case~(i) they are both easy calculations). Thus $\dim C_{\L(H)}(\theta_2) = |\Phi_H|/r = 8$ or $36$ respectively; the classification of outer automorphisms of $H$ (see e.g. \cite[Tables~4.3.1 and 4.7.1]{GLS}) now shows that $\theta_2$ must be a conjugate of $\theta_1$.

First set $\theta = \theta_1$; then we may take $G = C_H(\theta)$. We see that if $\tau(\alpha) = \alpha$ then $e_\alpha \in \L(H)_{(i)}$ if and only if $\height(\alpha) \equiv i$ (mod $r$), while if $\tau(\alpha) \neq \alpha$ then $\L(H)_{(i)}$ contains the vector $e_\alpha + {\eta_3}^{\height(\alpha) - i} e_{\tau(\alpha)} + {\eta_3}^{2\height(\alpha) - 2i} e_{\tau^2(\alpha)}$ or $e_\alpha + (-1)^{\height(\alpha)-i} e_{\tau(\alpha)}$ respectively. Thus in $\L(H)_{(1)}$ we have a highest weight vector $e_\beta$ for $\beta = \dfourrt1111$ or $\esixrt122321$ respectively; the expressions above for the simple root elements of $G$ show that $\L(H)_{(1)}$ is the Weyl $G$-module with high weight $3\omega_1$ or $\omega_4$ respectively. We have $Z(\L(H))_{(1)} = Z(\L(H))$, since if $Z(\L(H))$ is non-zero its generator given above is negated by $\theta$; we may take $V = \L(H)_{(1)}/Z(\L(H))_{(1)}$, and then $G_V = \{ 1 \}$.

Take $v = e_{\beta_1} + \cdots + e_{\beta_{\ell_H}} \in \L(U_H)$; then $v$ is a regular nilpotent element, and by the previous paragraph we have $v \in \L(H)_{(1)}$. From Lemma~\ref{lem: m_i for H} we see that there are natural numbers $m_1, \dots, m_{\ell_H}$, which are listed there, such that we may write
$$
C_{U_H}(v) = \{ y_1(c_1) \dots y_{\ell_H}(c_{\ell_H}) : c_i \in K \},
$$
with each $y_i(c)$ of the form $\left(\prod_{\height(\alpha) = m_i} x_\alpha(n_\alpha c)\right)x$, where $x$ is a product of root elements corresponding to roots of height greater than $m_i$, and the $n_\alpha \in K$ are not all zero and satisfy $\sum_{\height(\alpha) = m_i} n_\alpha e_\alpha \in C_{\L(G)}(v)$. Suppose $g = y_1(c_1) \dots y_{\ell_H}(c_{\ell_H}) \in G \cap C_{U_H}(v)$; we shall prove that $g = 1$. Since $G_V = \{ 1 \} = G \cap Z(H)$, Lemma~\ref{lem: semisimple auts}(i) will then show that the orbit containing $v + Z(\L(H))_{(1)}$ is regular.

First suppose we are in case (i); here we have $m_1 = 1$, $m_2 = m_3 = 3$ and $m_4 = 5$. Write
$$
v_1 = v, \quad v_2 = e_{\dfourrt0111} - e_{\dfourrt1101}, \quad v_3 = e_{\dfourrt0111} - e_{\dfourrt1110}, \quad v_4 = e_{\dfourrt1211};
$$
then each $v_i$ is the vector lying in $C_{\L(G)}(v)$ corresponding to the element $y_i(c_i)$. Since by the above $G \cap U_H$ contains no element with non-trivial projection on the root subgroup corresponding to the root $\dfourrt0100$ or $\dfourrt1211$, we must have $c_1 = c_4 = 0$. Moreover, for $\alpha = \dfourrt1110$ the projection of $G \cap U_H$ on the product of the root groups corresponding to roots $\alpha$, $\tau(\alpha)$ and $\tau^2(\alpha)$ consists of elements $x_\alpha(t) x_{\tau(\alpha)}(t) x_{\tau^2(\alpha)}(t)$; since the vector
$$
e_{\dfourrt1110} + e_{\dfourrt0111} + e_{\dfourrt1101}
$$
is not a linear combination of $v_2$ and $v_3$, we must also have $c_2 = c_3 = 0$. Therefore $g = 1$ as required.

Now suppose we are in case (ii); here we have $m_1 = 1$ or $3$ according as $p \geq 5$ or $p = 3$, $m_2 = 4$, $m_3 = 5$, $m_4 = 7$, $m_5 = 8$ and $m_6 = 11$. Write
\begin{eqnarray*}
v_1    & = & \begin{cases}
v,                                                                                                 & \hbox{if } p \geq 5, \\
e_{\esixrt101100} + e_{\esixrt000111} + e_{\esixrt011100} + e_{\esixrt010110} - e_{\esixrt001110}, & \hbox{if } p = 3, \\
\end{cases} \\
v_2    & = & e_{\esixrt111100} - e_{\esixrt010111} + e_{\esixrt101110} - e_{\esixrt001111}, \\
v_3    & = & e_{\esixrt111110} + e_{\esixrt011111} + 2e_{\esixrt101111} - e_{\esixrt011210}, \\
v_4    & = & e_{\esixrt112210} + e_{\esixrt011221} - e_{\esixrt111211}, \\
v_5    & = & e_{\esixrt112211} - e_{\esixrt111221}, \\
v_6    & = & e_{\esixrt122321};
\end{eqnarray*}
then each $v_i$ is the vector lying in $C_{\L(G)}(v)$ corresponding to the element $y_i(c_i)$. Since by the above $G \cap U_H$ contains no element with non-trivial projection on the root subgroup corresponding to the root $\esixrt010000$, $\esixrt001110$, $\esixrt101111$, $\esixrt111211$ or $\esixrt122321$, we must have $c_1 = c_3 = c_4 = c_6 = 0$. Moreover, for $\alpha = \esixrt111100$, $\esixrt101110$ or $\esixrt112211$ the projection of $G \cap U_H$ on the product of the root groups corresponding to roots $\alpha$ and $\tau(\alpha)$ consists of elements $x_\alpha(t) x_{\tau(\alpha)}(t)$ rather than $x_\alpha(t) x_{\tau(\alpha)}(-t)$; thus we must also have $c_2 = c_5 = 0$. Therefore $g = 1$ as required.

Now set $\theta = \theta_2$, and again take $G = C_H(\theta)$ and $V = \L(H)_{(1)}/Z(\L(H))_{(1)}$. We have $G \cap T_H = C_{T_H}(\theta) \cong {\Z_r}^{\ell_H/(r - 1)}$, i.e., ${\Z_3}^2$ or ${\Z_2}^6$ respectively (this is obvious in case~(ii) as then $\theta$ acts on $T_H$ as inversion; in case~(i) we have $C_{T_H}(\theta) = \langle h_{\beta_1}(\eta_3) h_{\beta_3}(\eta_3) h_{\beta_4}(\eta_3), h_{\beta_2}(\eta_3) h_{\beta_3}({\eta_3}^2) h_{\beta_4}(\eta_3) \rangle$). Moreover, in case~(ii) we have $\L(T_H)_{(1)} = \L(T_H)$, while in case~(i) we find that $\L(T_H)_{(1)} = \langle h_{\beta_1} + {\eta_3}^2 h_{\beta_3} + \eta_3 h_{\beta_4}, h_{\beta_2} - {\eta_3}^2 h_{\beta_3} + h_{\beta_4} \rangle$. Thus $\dim \L(H)_{(1)} - \dim \L(T_H)_{(1)} = \dim G - \dim (G \cap T_H)$ in each case; and in case (i) a routine check shows that $\L(T_H)_{(1)}$ contains regular semisimple elements.

By Lemma~\ref{lem: W_H on L(T_H)/Z(L(H))}, in case~(ii) we have $({W_H}^\ddagger)_{(1)} = {W_H}^\ddagger = \{ 1 \}$; we claim that in case~(i) we have $({W_H}^\ddagger)_{(1)} = \langle w_0 \rangle$. Thus suppose $w \in W_H$ and there exists $\xi \in K^*$ such that for all $y \in \L(T_H)_{(1)}$ we have $w.y = \xi y$. We use the standard notation for the roots of $\Phi_H$; then the elements of $W_H$ act as signed permutations of $\{ 1, 2, 3, 4 \}$. Taking $y = h_{\beta_2} - {\eta_3}^2 h_{\beta_3} + h_{\beta_4}$ we see that the permutation involved in $w$ must be some power of the $3$-cycle $(2 \ 3 \ 4)$; taking $y = h_{\beta_1} + {\eta_3}^2 h_{\beta_3} + \eta_3 h_{\beta_4}$ then forces the permutation to be the identity, and all signs to be equal, so $w \in \langle w_0 \rangle$ as required. Note that in this case if we write $n_0 = n_{\dfourrt1000} n_{\dfourrt0010} n_{\dfourrt0001} n_{\dfourrt1211}$, then $n_0$ is an involution in $N_H$ corresponding to $w_0$ which commutes with $\theta$.

Now $w_0$ acts on $\L(T_H)_{(1)}$ as negation. Thus in both cases $({W_H}^\dagger)_{(1)} = \{ 1 \}$, and so $C_{({N_H}^\dagger)_{(1)}}(\theta) = C_{T_H}(\theta)$. Also we have $({N_H}^\ddagger)_{(1)} = T_H \langle n_0 \rangle$ or $T_H$ respectively, so $C_{({N_H}^\ddagger)_{(1)}}(\theta) = C_{T_H}(\theta) \langle n_0 \rangle$ or $C_{T_H}(\theta)$ respectively. Since in each case $G_V = Z(G) = \{ 1 \}$, Lemma~\ref{lem: semisimple auts}(ii) shows that the triple $(G, \lambda, p)$ has generic stabilizer $C_{({N_H}^\dagger)_{(1)}}(\theta)/G_V \cong {\Z_3}^2$ or ${\Z_2}^6$ respectively, while the quadruple $(G, \lambda, p, 1)$ has generic stabilizer $C_{({N_H}^\ddagger)_{(1)}}(\theta)/Z(G) \cong {\Z_3}^2.\Z_2$ or ${\Z_2}^6$ respectively.
\end{proof}

Next we consider two infinite families of cases, in which $G$ is an orthogonal group. In the statement of the following result, for convenience we refer to the cases where $G = B_1$, $\lambda = 2\omega_1$ and $G = D_3$, $\lambda = 2\omega_1$; these appear in Table~\ref{table: large triple and first quadruple non-TGS} as $G = A_1$, $\lambda = 4\omega_1$ and $G = A_3$, $\lambda = 2\omega_2$ respectively.

\begin{prop}\label{prop: B_ell or D_ell, 2omega_1 module}
Let $G = B_\ell$ for $\ell \in [1, \infty)$ or $D_\ell$ for $\ell \in [3, \infty)$, and $\lambda = 2\omega_1$ with $p \geq 3$ (and if $G = B_1$ then $p \neq 3$). Then the triple $(G, \lambda, p)$ has generic stabilizer ${\Z_2}^{2\ell}$ or ${\Z_2}^{2\ell - 2}$ respectively, and there is a regular orbit; the associated first quadruple $(G, \lambda, p, 1)$ has generic stabilizer ${\Z_2}^{2\ell}$ or ${\Z_2}^{2\ell - 2}$ respectively.
\end{prop}

\begin{proof}
Number the cases (i) and (ii) according as $G = B_\ell$ or $D_\ell$; whenever we give two choices followed by the word `respectively' we are taking the cases in the order (i), (ii). As the proof here is so similar to those of Propositions~\ref{prop: A_8, omega_3, A_7, omega_4, D_8, omega_8 modules, non-special characteristic} and \ref{prop: A_2, 3omega_1, C_4, omega_4 modules}, we shall be brief in places.

Let $H$ be the simply connected group of type $A_{\ell_H}$ over $K$, with simple roots $\beta_1, \dots, \beta_{\ell_H}$, where $\ell_H = 2\ell$ or $2\ell - 1$ respectively, so that $H = \SL_{\ell_H + 1}(K)$; assume $p \geq 3$ (and if $\ell_H = 2$ then $p \geq 5$). We have $Z(\L(H)) = \{ 0 \}$ unless $p$ divides $\ell_H + 1$, in which case $Z(\L(H)) = \langle h_{\beta_1} + 2h_{\beta_2} + \cdots + \ell_H h_{\beta_{\ell_H}} \rangle$.

Let $\tau$ be the automorphism of $\Phi_H$ which preserves $\Pi_H$ and permutes simple roots by sending $\beta_i$ to $\beta_{\ell_H + 1 - i}$. This time we shall not assume that $\tau$ preserves the structure constants (indeed it cannot in case~(i), since then it interchanges $\beta_\ell$ and $\beta_{\ell + 1}$); rather we shall assume that for all $h < i < j$ we have $N_{\beta_h + \cdots + \beta_{i - 1}, \beta_i + \cdots + \beta_{j - 1}} = 1$. Again we obtain a graph automorphism of $H$, which by slight abuse of notation we also call $\tau$.

Multiplying $\tau$ by an appropriate element of $T_H$ we obtain $\theta_1$, such that ${\theta_1}^2 = 1$ and $\theta_1$ sends $x_{\beta_i}(t)$ to $x_{\beta_{\ell_H + 1 - i}}(-t)$, unless we are in case~(i) and $i \in \{ \ell, \ell + 1 \}$, when it sends $x_{\beta_\ell}(t)$ to $x_{\beta_{\ell + 1}}(-\frac{1}{2}t)$ and $x_{\beta_{\ell + 1}}(t)$ to $x_{\beta_\ell}(-2t)$. We find that if $\tau(\alpha) = \alpha$ then $X_\alpha \not < C_H(\theta_1)$; if $\tau(\alpha)$ is orthogonal to $\alpha$ then the intersection of $C_H(\theta_1)$ with $X_\alpha X_{\tau(\alpha)}$ is the $1$-dimensional group $\{ x_\alpha(t) x_{\tau(\alpha)}(-t) : t \in K \}$; if we are in case~(i) and $\alpha = \beta_i + \cdots + \beta_\ell$ then the intersection of $C_H(\theta_1)$ with $X_\alpha X_{\tau(\alpha)} X_{\alpha + \tau(\alpha)}$ is the $1$-dimensional group $\{ x_\alpha(2t) x_{\tau(\alpha)}(-t) x_{\alpha + \tau(\alpha)}(t^2) : t \in K \}$. It follows that $C_H(\theta_1) = \SO_{\ell_H + 1}(K)$ is a group of type $B_\ell$ or $D_\ell$ respectively, with simple root elements
\begin{eqnarray*}
x_{\alpha_i}(t)    & = & x_{\beta_i}(t) x_{\beta_{\ell_H + 1 - i}}(-t) \quad \hbox{for } i < \ell, \\
x_{\alpha_\ell}(t) & = &
\begin{cases}
x_{\beta_\ell}(2t) x_{\beta_{\ell + 1}}(-t) x_{\beta_\ell + \beta_{\ell + 1}}(t^2) & \hbox{in case~(i),} \\
x_{\beta_{\ell - 1} + \beta_\ell}(t) x_{\beta_\ell + \beta_{\ell + 1}}(-t)         & \hbox{in case~(ii).}
\end{cases}
\end{eqnarray*}
Regarding the elements of $H$ as matrices, we may take $x_{\beta_i + \cdots + \beta_{j - 1}}(t)$ as $I_{\ell_H} + tE_{ij}$ where $E_{ij}$ is the matrix unit with $(i, j)$-entry $1$ and all other entries $0$; we then recover the action of $C_H(\theta_1)$ on its natural module described in Section~\ref{sect: notation}.

Now for $i = 1, \dots, \ell$ let
$$
\delta_i = \beta_i + \beta_{i + 1} + \cdots + \beta_{\ell_H + 1 - i};
$$
then $\langle \delta_1, \dots, \delta_\ell \rangle$ is a subsystem of type ${A_1}^\ell$. Set $\theta_2 = \tau n_{\delta_1} \dots n_{\delta_\ell}$; then ${\theta_2}^2 = 1$, and indeed $\theta_2$ sends each root subgroup $X_\alpha$ to $X_{-\alpha}$. We find that $\theta_2$ acts fixed-point-freely on both $\Phi_H$ and $\L(T_H)$. Thus $\dim C_{\L(H)}(\theta_2) = |\Phi_H|/2 = \frac{1}{2}\ell_H(\ell_H + 1)$; the classification of outer automorphisms of $H$ (see e.g. \cite[Table~4.3.1]{GLS}) now shows that $\theta_2$ must be a conjugate of $\theta_1$.

First set $\theta = \theta_1$; then we may take $G = C_H(\theta)$. We see that $\L(H)_{(1)}$ is the Weyl $G$-module with high weight $2\omega_1$. We have $Z(\L(H))_{(1)} = Z(\L(H))$, since if $Z(\L(H))$ is non-zero its generator given above is negated by $\theta$; we may take $V = \L(H)_{(1)}/Z(\L(H))_{(1)}$, and then $G_V = Z(G)$.

Take $v = e_{\beta_1} + \cdots + e_{\beta_{\ell - 1}} + 2e_{\beta_\ell} + e_{\beta_{\ell + 1}} + \cdots + e_{\beta_{\ell_H}} \in \L(U_H)$; then $v$ is a regular nilpotent element lying in $\L(H)_{(1)}$. As before we may refer to Lemma~\ref{lem: m_i for H} to see the structure of $C_{U_H}(v)$, but here we can be more explicit: by taking the known group $C_{U_H}(v')$ where $v'$ is obtained from $v$ by changing the coefficient of $e_{\beta_\ell}$ from $2$ to $1$, and conjugating by a suitable element of $T_H$, we see that $C_{U_H}(v)$ comprises upper unitriangular matrices $g$ with the property that there exist $c_1, \dots, c_{\ell_H} \in K$ such that the $(i, j)$-entry is $c_{j - i}$ if either $j \leq \ell$ or $i > \ell$, and $2c_{j - i}$ if $i \leq \ell < j$. By comparing with the description above of the root groups in $G$, we see that if $g \in G \cap C_{U_H}(v)$ we must have all $c_i = 0$, so that $g = 1$. Since $G_V = Z(G) = G \cap Z(H)$, Lemma~\ref{lem: semisimple auts}(i) shows that the orbit containing $v + Z(\L(H))$ is regular.

Now set $\theta = \theta_2$, and again take $G = C_H(\theta)$ and $V = \L(H)_{(1)}/Z(\L(H))_{(1)}$. We have $G \cap T_H = C_{T_H}(\theta) \cong {\Z_2}^{\ell_H}$, i.e., ${\Z_2}^{2\ell}$ or ${\Z_2}^{2\ell - 1}$ respectively (as $\theta$ acts on $T_H$ as inversion). Moreover, we have $\L(T_H)_{(1)} = \L(T_H)$. Thus $\dim \L(H)_{(1)} - \dim \L(T_H)_{(1)} = \dim G - \dim (G \cap T_H)$ in each case. By Lemma~\ref{lem: W_H on L(T_H)/Z(L(H))}, we have $({W_H}^\ddagger)_{(1)} = {W_H}^\ddagger = \{ 1 \} = ({W_H}^\dagger)_{(1)} = {W_H}^\dagger$, and so $C_{({N_H}^\ddagger)_{(1)}}(\theta) = C_{({N_H}^\dagger)_{(1)}}(\theta) = C_{T_H}(\theta)$. Since $G_V = Z(G) = \{ 1 \}$ or $\Z_2$ respectively, Lemma~\ref{lem: semisimple auts}(ii) shows that the triple $(G, \lambda, p)$ has generic stabilizer $C_{({N_H}^\dagger)_{(1)}}(\theta)/G_V \cong {\Z_2}^{2\ell}$ or ${\Z_2}^{2\ell - 2}$ respectively, while the quadruple $(G, \lambda, p, 1)$ has generic stabilizer $C_{({N_H}^\ddagger)_{(1)}}(\theta)/Z(G) \cong {\Z_2}^{2\ell}$ or ${\Z_2}^{2\ell - 2}$ respectively.
\end{proof}

We now treat the three postponed cases, where the approach using Section~\ref{sect: semisimple auts} does not apply; we shall instead employ that of Section~\ref{sect: Lie algebra annihilators}.

\begin{prop}\label{prop: A_8, omega_3, A_7, omega_4, D_8, omega_8 modules, special characteristic}
Let $G = A_7$ and $\lambda = \omega_4$ with $p = 2$, or $G = D_8$ and $\lambda = \omega_8$ with $p = 2$, or $G = A_8$ and $\lambda = \omega_3$ with $p = 3$. Then the triple $(G, \lambda, p)$ has generic stabilizer ${\Z_2}^3$, or ${\Z_2}^4$, or ${\Z_3}^2$, respectively, and there is a regular orbit; the associated first quadruple $(G, \lambda, p, 1)$ has generic stabilizer ${\Z_2}^3$, or ${\Z_2}^4$, or ${\Z_3}^2.\Z_2$, respectively.
\end{prop}

\begin{proof}
Although the approach using Lemma~\ref{lem: semisimple auts} does not apply in these cases, there are points of contact with the proof of Proposition~\ref{prop: A_8, omega_3, A_7, omega_4, D_8, omega_8 modules, non-special characteristic}. As there, number the cases (i), (ii) and (iii) according as $G = A_7$, $D_8$ or $A_8$; whenever we give three choices followed by the word `respectively' we are taking the cases in the order (i), (ii), (iii). Note that in each case $Z(G) = \{ 1 \}$, so that $G_V = \{ 1 \}$.

First let $H$ be the simply connected group defined over $K$ of type $E_7$, $E_8$ or $E_8$ respectively. Although we cannot take $\theta_1 \in T_H$ of order $p$, we may still let $G = \langle T_H, X_\alpha : \alpha \in \Phi_H, \ \height(\alpha) \equiv 0 \hbox{ (mod $p$)} \rangle$, and then $G$ is a connected group of type $A_7$, $D_8$ or $A_8$ respectively, with simple roots $\alpha_1, \dots, \alpha_\ell$ as listed in the proof of Proposition~\ref{prop: A_8, omega_3, A_7, omega_4, D_8, omega_8 modules, non-special characteristic}; likewise we may let $V = \langle e_\alpha : \alpha \in \Phi_H, \ \height(\alpha) \equiv 1 \hbox{ (mod $p$)} \rangle$ (note that in case (i) $Z(\L(H)) \neq \{ 0 \}$, but we choose to define $Z(\L(H))_{(1)} = \{ 0 \}$). We still have the regular nilpotent element $v = e_{\beta_1} + \cdots + e_{\beta_\ell} \in \L(U_H) \cap V$, but here Lemma~\ref{lem: m_i for H} does not show that it lies in a regular orbit since at least one of the values $m_i$ is divisible by $p$ (indeed we find that $C_G(v)$ is non-trivial). We shall show the existence of regular orbits in a different way.

It will in fact prove convenient to use a different notation for elements of $V$. In case (i) we may view $V$ as the exterior power $\bigwedge^4(V_{nat})$; for $i_1, i_2, i_3, i_4 \leq 8$ we write $v_{i_1i_2i_3i_4} = v_{i_1} \wedge v_{i_2} \wedge v_{i_3} \wedge v_{i_4}$, and then $V = \langle v_{i_1i_2i_3i_4} : 1 \leq i_1 < i_2 < i_3 < i_4 \leq 8 \rangle$. Likewise in case (iii) we may view $V$ as the exterior cube $\bigwedge^3(V_{nat})$; for $i_1, i_2, i_3 \leq 9$ we write $v_{i_1i_2i_3} = v_{i_1} \wedge v_{i_2} \wedge v_{i_3}$, and then $V = \langle v_{i_1i_2i_3} : 1 \leq i_1 < i_2 < i_3 \leq 9 \rangle$. In case (ii), we use the standard notation for the roots in $\Phi$, and then each weight $\nu \in \Lambda(V)$ is of the form $\frac{1}{2}\sum_{i = 1}^8 \pm\ve_i$, where the number of minus signs is even; we shall represent such a weight as a string of $8$ plus or minus signs, and write $v_\nu$ for the corresponding weight vector, so that $V = \langle v_\nu : \nu \in \Lambda(V) \rangle$ and each element $n_\alpha$ for $\alpha \in \Phi$ permutes the vectors $v_\nu$. In addition, in cases (i) and (iii) we shall identify $W$ with the symmetric group $S_8$ or $S_9$.

Before proceeding it is worth noting that cases (i) and (ii) are linked: we have the obvious $A_7$ subgroup of $D_8$ with simple roots $\alpha_1, \dots, \alpha_7$, and the $A_7$-module may be identified with the span in the $D_8$-module of the vectors $v_\nu$ where the weight $\nu$ has $4$ plus and $4$ minus signs; indeed the vector $v_{i_1i_2i_3i_4}$ is then equal to $v_\nu$ where the weight $\nu$ has plus signs in positions $i_1$, $i_2$, $i_3$, $i_4$ and minus signs elsewhere. We may use either notation for vectors in $V$ in case (i).

Writing $h_i$ for $h_{\alpha_i}$, let $\S \leq \L(T)$ be
$$
\begin{array}{ll}
\langle h_1 + h_3, h_3 + h_5, h_5 + h_7, h_2 + h_6 \rangle               & \vstrut \hbox{in cases~(i) and (ii),} \\
\langle h_1 - h_2 - h_7 + h_8, h_1 + h_2 + h_4 + h_5 + h_7 + h_8 \rangle & \vstrut \hbox{in case~(iii).} \\
\end{array}
$$
It is easy to see that if $\alpha \in \Phi$ there exists $h \in \S$ with $[h e_\alpha] \neq 0$, so $C_{\L(G)}(\S) = \L(T)$.

For $1 \leq i \leq \frac{\ell}{p - 1}$ and $1 \leq j \leq p$ we define vectors $x_{ij}$ in $V$ as follows: in cases (i) and (ii) we set
\begin{eqnarray*}
x_{11} = v_{\sss{+--+-++-}}, & & x_{12} = v_{\sss{-++-+--+}}, \\
x_{21} = v_{\sss{+--++--+}}, & & x_{22} = v_{\sss{-++--++-}}, \\
x_{31} = v_{\sss{+-+--+-+}}, & & x_{32} = v_{\sss{-+-++-+-}}, \\
x_{41} = v_{\sss{+-+-+-+-}}, & & x_{42} = v_{\sss{-+-+-+-+}}, \\
x_{51} = v_{\sss{++----++}}, & & x_{52} = v_{\sss{--++++--}}, \\
x_{61} = v_{\sss{++--++--}}, & & x_{62} = v_{\sss{--++--++}}, \\
x_{71} = v_{\sss{++++----}}, & & x_{72} = v_{\sss{----++++}}, \\
x_{81} = v_{\sss{++++++++}}, & & x_{82} = v_{\sss{--------}}
\end{eqnarray*}
(where in case (i) we ignore the vectors $x_{81}$ and $x_{82}$); in case (iii) we set
\begin{eqnarray*}
& x_{11} = v_{348}, \quad x_{12} = v_{267}, \quad x_{13} = v_{159}, & \\
& x_{21} = v_{168}, \quad x_{22} = v_{357}, \quad x_{23} = v_{249}, & \\
& x_{31} = v_{258}, \quad x_{32} = v_{147}, \quad x_{33} = v_{369}, & \\
& x_{41} = v_{456}, \quad x_{42} = v_{123}, \quad x_{43} = v_{789}. &
\end{eqnarray*}
(If we regard $V$ as a submodule of $\L(H)$ as above and use the notation of the proof of Proposition~\ref{prop: A_8, omega_3, A_7, omega_4, D_8, omega_8 modules, non-special characteristic}, then in cases (i) and (ii) we have $x_{i1} = e_{\delta_i}$ and $x_{i2} = e_{-\delta_i}$, while in case (iii) we have $x_{i1} = e_{\delta_{2i - 1}}$, $x_{i2} = e_{\delta_{2i}}$ and $x_{i3} = e_{-(\delta_{2i - 1} + \delta_{2i})}$.) For each pair $(i, j)$ let $\nu_{ij}$ be the weight such that $V_{\nu_{ij}} = \langle x_{ij} \rangle$; thus for each $i$ we have $\sum_{j = 1}^p \nu_{ij} = 0$. Let
$$
Y' = \langle x_{ij} : 1 \leq i \leq {\ts\frac{\ell}{p - 1}}, \ 1 \leq j \leq p \rangle;
$$
a straightforward calculation shows that $Y'$ is the subspace of $V$ annihilated by the subalgebra $\S$.

Define $\Upsilon$ to be the following set of subsets of $\{ 1, \dots, \ell \}$, where in the interests of brevity we write simply \lq $i_1i_2\dots$' for \lq $\{ i_1, i_2, \dots \}$':
$$
\begin{array}{ll}
\phantom{ \left.\right. } \left\{ 1234, 1256, 1357, 1467, 2367, 2457, 3456 \right\} & \vstrut \hbox{in case~(i),} \\
\phantom{ \left.\right. } \left\{ 1234, 1256, 1357, 1467, 2367, 2457, 3456, \right. & \vstrut \\
\phantom{ \left\{\right. } \left. 5678, 3478, 2468, 2358, 1458, 1368, 1278 \right\} & \vstrut \hbox{in case~(ii),} \\
\phantom{ \left.\right. } \left\{ 123, 124, 134, 234 \right\}                       & \vstrut \hbox{in case~(iii).} \\
\end{array}
$$
Note that in cases (i) and (ii) the set $\Upsilon$ may be characterised as follows: given a subset $S$ of $\{ 1, \dots, \ell \}$, we have $S \in \Upsilon$ if and only if there exists $\alpha \in \Phi$ such that the weights $\nu_{ij}$ not orthogonal to $\alpha$ are precisely those for which $i \in S$ (for example, the weights $\nu_{ij}$ not orthogonal to $\alpha = \ve_1 - \ve_2$ are those with $i \in \{ 1, 2, 3, 4 \}$). Set
$$
\hat{Y'} = \left\{ {\ts\sum}_{i, j} a_{ij} x_{ij} : \forall i \ (a_{ij}, a_{ij'}) \neq (0, 0) \hbox{ for } j \neq j', \ \forall S \in \Upsilon, \ {\ts\sum}_{i \in S} (\pm{\ts\prod}_j a_{ij}) \neq 0 \right\};
$$
then $\hat{Y'}$ is a dense open subset of $Y'$. Note that as for each $i$ we have $\sum_j \nu_{ij} = 0$, applying an element of $T$ to an element of $Y'$ has no effect on the values $\prod_j a_{ij}$, so $T$ preserves $\hat{Y'}$. Moreover if $y \in \hat{Y'}$ and $s \in T$ with $s.y \in \langle y \rangle$, writing $s.y = \kappa y$ and considering the coefficients of $x_{ij}$ for fixed $i$ such that $\prod_j a_{ij} \neq 0$ we see that $\kappa^p = 1$, whence $\kappa = 1$ so that in fact $s.y = y$; a straightforward calculation now shows that $s = 1$, so $C_T(y) = C_T(\langle y \rangle) = \{ 1 \}$. Take
$$
y = {\ts\sum}_{i, j} a_{ij} x_{ij} \in \hat{Y'}.
$$

First suppose $x \in \Ann_{\L(G)}(y)$; write $x = h + e$ where $h \in \L(T)$ and $e \in \langle e_\alpha : \alpha \in \Phi \rangle$. Clearly $h.y \in Y'$; since the difference of two weights $\nu_{ij}$ is never a root, for each pair $(i, j)$ we see that $e.y$ contains no term $x_{ij}$. Thus we must have $h.y = e.y = 0$. A quick calculation shows that we must have $h \in \S$. Now write $e = \sum_{\alpha \in \Phi} t_\alpha e_\alpha$; then the equation $e.y = 0$ may be expressed in matrix form as $A{\bf t} = {\bf 0}$, where $A$ is an $M \times M$ matrix and ${\bf t}$ is a column vector whose entries are the various coefficients $t_\alpha$. We find that if the rows and columns of $A$ are suitably ordered then it becomes block diagonal, having $7$, $14$ or $8$ blocks respectively, with each block being an $8 \times 8$, $8 \times 8$ or $9 \times 9$ matrix respectively. In cases (i) and (ii) each block may be written in the form
$$
\left(
  \begin{array}{cc|cc|cc|cc}
             & a_{i_11} & a_{i_21} &          & a_{i_31} &          & a_{i_42} &          \\
    a_{i_12} &          &          & a_{i_22} &          & a_{i_32} &          & a_{i_41} \\
  \hline
    a_{i_21} &          &          & a_{i_11} & a_{i_41} &          & a_{i_32} &          \\
             & a_{i_22} & a_{i_12} &          &          & a_{i_42} &          & a_{i_31} \\
  \hline
    a_{i_31} &          & a_{i_41} &          &          & a_{i_11} & a_{i_22} &          \\
             & a_{i_32} &          & a_{i_42} & a_{i_12} &          &          & a_{i_21} \\
  \hline
    a_{i_42} &          & a_{i_32} &          & a_{i_22} &          &          & a_{i_11} \\
             & a_{i_41} &          & a_{i_31} &          & a_{i_21} & a_{i_12} &          \\
  \end{array}
\right)
$$
where $S = \{ i_1, i_2, i_3, i_4 \} \in \Upsilon$; calculation shows that the determinant of this $8 \times 8$ matrix is $\sum_{i \in S}(a_{i1}a_{i2})^4 = (\sum_{i \in S} (\prod_j a_{ij}))^4$. In case (iii), each block or its transpose may be written in the form
$$
\left(
  \begin{array}{ccc|ccc|ccc}
                     &                   &   a_{i_1j_{11}}   &  \e a_{i_2j_{21}} &                   &                   & \e' a_{i_3j_{31}} &                   &                   \\
     a_{i_1j_{12}}   &                   &                   &                   &  \e a_{i_2j_{22}} &                   &                   & \e' a_{i_3j_{32}} &                   \\
                     &   a_{i_1j_{13}}   &                   &                   &                   &  \e a_{i_2j_{23}} &                   &                   & \e' a_{i_3j_{33}} \\
  \hline
   \e' a_{i_3j_{33}} &                   &                   &                   &                   &   a_{i_1j_{11}}   &  \e a_{i_2j_{22}} &                   &                   \\
                     & \e' a_{i_3j_{31}} &                   &   a_{i_1j_{12}}   &                   &                   &                   &  \e a_{i_2j_{23}} &                   \\
                     &                   & \e' a_{i_3j_{32}} &                   &   a_{i_1j_{13}}   &                   &                   &                   &  \e a_{i_2j_{21}} \\
  \hline
    \e a_{i_2j_{23}} &                   &                   & \e' a_{i_3j_{32}} &                   &                   &                   &                   &   a_{i_1j_{11}}   \\
                     &  \e a_{i_2j_{21}} &                   &                   & \e' a_{i_3j_{33}} &                   &   a_{i_1j_{12}}   &                   &                   \\
                     &                   &  \e a_{i_2j_{22}} &                   &                   & \e' a_{i_3j_{31}} &                   &   a_{i_1j_{13}}   &                   \\
  \end{array}
\right)
$$
where $S = \{ i_1, i_2, i_3 \} \in \Upsilon$, $\e, \e' \in \{ \pm1 \}$, and for each $i$ the values $j_{i1}$, $j_{i2}$, $j_{i3}$ are $1$, $2$, $3$ in some order; calculation shows that the determinant of this $9 \times 9$ matrix is $\sum_{i \in S} (\pm a_{i1}a_{i2}a_{i3})^3 = (\sum_{i \in S} (\pm \prod_j a_{ij}))^3$. Thus in each case the final condition in the definition of the set $\hat{Y'}$ implies that each block of $A$ is non-singular, as therefore is $A$ itself; so ${\bf t}$ must be the zero vector and hence $e = 0$. Thus $x = h + e \in \S$; so $\Ann_{\L(G)}(y) = \S$. By Lemma~\ref{lem: exactness via Premet}(i) we have $\Tran_G(y, Y') \subseteq N$, so $C_G(y) \leq C_G(\langle y \rangle) \leq N$.

Now take $y \in \hat{Y'}$ as follows: in cases (i) and (ii) let
$$
y = x_{11} + x_{21} + x_{31} + a_4(x_{41} + x_{42}) + a_5(x_{51} + x_{52}) + a_6(x_{61} + x_{62}) + a_7(x_{71} + x_{72}) + x_{81}
$$
for $a_4, a_5, a_6, a_7 \in K^*$ distinct (where in case (i) we ignore the vector $x_{81}$); in case (iii) let
$$
y = x_{11} + x_{12} + x_{22} + x_{23} + a_3(x_{31} + x_{32} + x_{33}) + a_4(x_{41} + x_{42} + x_{43})
$$
for $a_3, a_4 \in K^*$ distinct up to sign. Write $I = \{ 4, 5, 6, 7 \}$ in cases (i) and (ii) and $I = \{ 3, 4 \}$ in case (iii). Take $n \in C_G(\langle y \rangle)$; then $n$ must permute the weights occurring in $y$, and as the minimal sets of such weights summing to zero are the $\{ \nu_{i1}, \dots, \nu_{ip} \}$ for $i \in I$ it must permute these sets. Indeed, we may write $n = n^* s$, where $n^*$ is a product of various elements $n_\alpha$ and $s \in T$. If $i \in I$ and $n$ sends the set $\{ \nu_{i1}, \dots, \nu_{ip} \}$ to the set $\{ \nu_{i'1}, \dots, \nu_{i'p} \}$, we see that $n^*$ must send each $x_{ij}$ to some $\pm x_{i'j'}$; since applying $s$ has no effect on the products of the coefficients of the $x_{i'j}$ for $1 \leq j \leq p$, the distinctness up to sign of the $a_i$ implies that $n$ must in fact fix each set $\{ \nu_{i1}, \dots, \nu_{ip} \}$ for $i \in I$, and permute the other weights occurring in $y$.

In case (i) $n$ must permute the weights $\nu_{11}, \nu_{21}, \nu_{31}$, so $nT \in W$ must permute the sets $\{ 1, 4, 6, 7 \}$, $\{ 1, 4, 5, 8 \}$ and $\{ 1, 3, 6, 8 \}$. As $1$ and $2$ are the only numbers appearing in all or none of these sets respectively, $nT$ must fix $1$ and $2$; as $n$ must also fix the sets $\{ \nu_{i1}, \nu_{i2} \}$ for $i = 4, 5, 6, 7$ we see that $nT$ must fix the sets $\{ 3, 5, 7 \}$, $\{ 7, 8 \}$, $\{ 5, 6 \}$ and $\{ 3, 4 \}$, so it fixes each of $3, 4, 5, 6, 7, 8$ and hence $nT = 1$. Thus the orbit containing $y$ is regular in case (i).

In case (ii) $n$ must permute the weights $\nu_{11}, \nu_{21}, \nu_{31}, \nu_{81}$. If $n$ fixes $\nu_{81}$ it must lie in $W(A_7)$, so by the previous paragraph $nT = 1$. If $n$ sends $\nu_{31}$ to $\nu_{81}$ then we must have $n = n' n''$ where $n'' = n_{2 - 4} n_{2 + 4} n_{5 - 7} n_{5 + 7}$ (where we use the standard notation for the roots in $\Phi$, and write $n_{i \pm j}$ for $n_{\ve_i \pm \ve_j}$) and $n'T \in W(A_7)$; then $n'T$ must send the sets $\{ 1, 2, 5, 6 \}$, $\{ 1, 2, 7, 8 \}$ and $\{ 1, 3, 6, 8 \}$ to the sets $\{ 1, 4, 6, 7 \}$, $\{ 1, 4, 5, 8 \}$ and $\{ 1, 3, 6, 8 \}$ in some order, so as before it must fix $1$ and send $4$ to $2$. However, $n''$ sends $\nu_{71}$ to $\nu_{41}$, and $n'$ cannot now send $\{ \nu_{41}, \nu_{42} \}$ to $\{ \nu_{71}, \nu_{72} \}$, because $1$ and $4$ are in different parts of the partition given by the former pair whereas $1$ and $2$ are in the same part of the partition given by the latter pair. Entirely similar arguments show that $n$ cannot send $\nu_{21}$ or $\nu_{11}$ to $\nu_{81}$. Thus the orbit containing $y$ is regular in case (ii).

In case (iii) $n$ must permute the weights $\nu_{11}, \nu_{12}, \nu_{22}, \nu_{23}$, so $nT \in W$ must permute the sets $\{ 3, 4, 8 \}$, $\{ 2, 6, 7 \}$, $\{ 1, 6, 8 \}$ and $\{ 2, 4, 9 \}$. As $5$ is the only number appearing in none of these sets, $nT$ must fix $5$; as $n$ must also fix $\{ \nu_{31}, \nu_{32}, \nu_{33} \}$ it must fix the weight $\nu_{31}$ and hence $nT$ must fix the set $\{ 2, 8 \}$. Write
$$
n_0 = \left(
        \begin{array}{ccc}
           &   & J \\
           & J &   \\
         J &   &   \\
        \end{array}
      \right)
\quad \hbox{and} \quad
n' = -\left(
        \begin{array}{ccc}
         J &   &   \\
           & J &   \\
           &   & J \\
        \end{array}
      \right)
\quad \hbox{for} \quad
J =
\left(
  \begin{array}{ccc}
      &   & 1 \\
      & 1 &   \\
    1 &   &   \\
  \end{array}
\right),
$$
so that $n_0T = (1\ 9)(2\ 8)(3\ 7)(4\ 6)$ is the long word of $W$ and $n'T = (1\ 3)(4\ 6)(7\ 9)$. Since $n_0T$ also fixes the above sets, ${n_0}^j nT$ for some $j \in  \{ 0, 1 \}$ must also fix $2$ and hence $8$; as ${n_0}^j n$ fixes $\{ \nu_{41}, \nu_{42}, \nu_{43} \}$ we see that ${n_0}^j nT$ must fix the sets $\{ 1, 3 \}$, $\{ 4, 6 \}$ and $\{ 7, 9 \}$, and then according as ${n_0}^j n$ fixes or interchanges $\nu_{32}$ and $\nu_{33}$ we must have ${n_0}^j nT = 1$ or $n'T$, so that $nT = {n_0}^j{n'}^{j'}T$ for some $j, j' \in \{ 0, 1 \}$. However, both $n_0$ and $n'$ send the vectors $x_{31}, x_{32}, x_{33}$ to $-x_{31}, -x_{33}, -x_{32}$ respectively, while $n_0$ sends the vectors $x_{41}, x_{42}, x_{43}$ to $-x_{41}, -x_{43}, -x_{42}$ respectively and $n'$ fixes all three of these vectors; so for $n.y = y$ we must have $j = j' = 0$, giving $C_G(y) = \{ 1 \}$. Thus the orbit containing $y$ is regular in case (iii).

We now consider generic stabilizers. For $1 \! \leq \! i \! \leq \! \frac{\ell}{p - 1}$ write $y_i = \sum_j x_{ij}$, and set
$$
Y = \langle y_i : 1 \leq i \leq {\ts\frac{\ell}{p - 1}} \rangle,
$$
so that $Y$ is a subvariety of $Y'$; note that in each case $\codim Y = \dim G$. Write
$$
\hat Y = \left\{ {\ts\sum}_i a_i y_i \in Y \cap \hat{Y'} : {\ts\prod}_i a_i \neq 0, \ a_i \neq \pm a_{i'} \hbox{ for } i \neq i'\right\};
$$
then $\hat Y$ is a dense open subset of $Y$. Take
$$
y = {\ts\sum}_i a_i y_i \in \hat Y.
$$

A straightforward calculation shows that $T.y \cap Y = \{ y \}$. Moreover, we have $\Tran_G(y, Y) \subseteq N$; if $n \in \Tran_G(y, Y)$ then $n$ must permute the sets $\{ \nu_{i1}, \dots, \nu_{ip} \}$ for $1 \leq i \leq \frac{\ell}{p - 1}$, and arguing as above we see that $n.y \in \hat Y$ (note that in cases (i) and (ii) the characterisation above of the set $\Upsilon$ shows that $n$ must preserve it). Thus by Lemma~\ref{lem: exactness via Premet}(ii) $y$ is $Y$-exact. We have $C_G(y) \leq C_G(\langle y \rangle) \leq N$. Note that as $C_T(\langle y \rangle) = \{ 1 \}$, each coset of $T$ in $N$ can contain at most one element of $C_G(\langle y \rangle)$.

In case (i) set
\begin{eqnarray*}
n_1 & = & n_{\alpha_1} n_{\alpha_3} n_{\alpha_5} n_{\alpha_7}, \\
n_2 & = & n_{\alpha_1 + \alpha_2} n_{\alpha_2 + \alpha_3} n_{\alpha_5 + \alpha_6} n_{\alpha_6 + \alpha_7}, \\
n_3 & = & n_{\alpha_1 + \alpha_2 + \alpha_3 + \alpha_4} n_{\alpha_2 + \alpha_3 + \alpha_4 + \alpha_5} n_{\alpha_3 + \alpha_4 + \alpha_5 + \alpha_6} n_{\alpha_4 + \alpha_5 + \alpha_6 + \alpha_7},
\end{eqnarray*}
so that we have $n_1T = (1\ 2)(3\ 4)(5\ 6)(7\ 8)$, $n_2T = (1\ 3)(2\ 4)(5\ 7)(6\ 8)$, $n_3T = (1\ 5)(2\ 6)(3\ 7)(4\ 8)$; in case (ii) set $n_1$, $n_2$, $n_3$ as given and also set
\begin{eqnarray*}
n_0 & = & n_{\alpha_1} n_{\alpha_1 + 2\alpha_2 + 2\alpha_3 + 2\alpha_4 + 2\alpha_5 + 2\alpha_6 + \alpha_7 + \alpha_8} \\
    &   & {} \times n_{\alpha_3} n_{\alpha_3 + 2\alpha_4 + 2\alpha_5 + 2\alpha_6 + \alpha_7 + \alpha_8} n_{\alpha_5} n_{\alpha_5 + 2\alpha_6 + \alpha_7 + \alpha_8} n_{\alpha_7} n_{\alpha_8},
\end{eqnarray*}
so that $n_0 T$ is the long word of the Weyl group; in case (iii) set
\begin{eqnarray*}
n_1 & = & n_{\alpha_1} n_{\alpha_2} n_{\alpha_4} n_{\alpha_5} n_{\alpha_7} n_{\alpha_8} s_1, \\
n_2 & = & n_{\alpha_1 + \alpha_2 + \alpha_3} n_{\alpha_4 + \alpha_5 + \alpha_6} n_{\alpha_2 + \alpha_3 + \alpha_4} n_{\alpha_5 + \alpha_6 + \alpha_7} n_{\alpha_3 + \alpha_4 + \alpha_5} n_{\alpha_6 + \alpha_7 + \alpha_8} s_2
\end{eqnarray*}
(where $s_1, s_2 \in T$ are chosen so that $n_1$ and $n_2$ are permutation matrices), so that we have $n_1T = (1\ 2\ 3)(4\ 5\ 6)(7\ 8\ 9)$, $n_2T = (1\ 4\ 7)(2\ 5\ 8)(3\ 6\ 9)$, and in addition let $n_0$ be as defined above with $n_0T = (1\ 9)(2\ 8)(3\ 7)(4\ 6)$. Let $C$ be the subgroup
$$
\begin{array}{ll}
\langle n_1, n_2, n_3 \rangle      & \vstrut \hbox{in case~(i),} \\
\langle n_1, n_2, n_3, n_0 \rangle & \vstrut \hbox{in case~(ii),} \\
\langle n_1, n_2 \rangle           & \vstrut \hbox{in case~(iii),} \\
\end{array}
$$
and let $C' = C$, $C$ or $C \langle n_0 \rangle$ respectively. Clearly we then have $C \leq C_G(y)$ and $C' \leq C_G(\langle y \rangle)$ (in case (iii), for each $i$ we see that $n_0$ negates one vector $x_{ij}$ and interchanges and negates the other two such vectors, so that $n_0.y = -y$); we shall show that in fact $C_G(y) = C$ and $C_G(\langle y \rangle) = C'$. Take $n \in C_G(\langle y \rangle)$; as before we see that $n$ must fix each set $\{ \nu_{i1}, \dots, \nu_{ip} \}$.

In case (i) the projection on $W$ of the group $C$ acts transitively on $\{ 1, \dots, 8 \}$, so there exists $c \in C$ such that $cnT \in W$ fixes $1$; as $cn$ fixes each set $\{ \nu_{i1}, \nu_{i2} \}$ we see that $cnT$ must preserve each of the $7$ corresponding partitions of $\{ 1, \dots, 8 \}$ into two sets of size $4$, and it immediately follows that $cnT = 1$, whence $n = c^{-1} \in C$. In case (ii), ${n_0}^j n$ for some $j \in  \{ 0, 1 \}$ fixes $\nu_{81}$, and therefore ${n_0}^j nT$ lies in $W(A_7)$, so by the previous sentence we again have $n \in C$. Finally in case (iii) the projection on $W$ of the group $C$ acts transitively on $\{ 1, \dots, 9 \}$, so there exists $c \in C$ such that $cnT \in W$ fixes $5$; as $cn$ fixes each set $\{ \nu_{i1}, \nu_{i2}, \nu_{i3} \}$ we see that $cnT$ must preserve the sets $\{ 1, 9 \}$, $\{ 3, 7 \}$, $\{ 2, 8 \}$, $\{ 4, 6 \}$. Since $n_0T$ also fixes these four sets, ${n_0}^j cnT$ for some $j \in  \{ 0, 1 \}$ must also fix $1$, and hence $9$, and then as ${n_0}^j cn$ fixes $\{ \nu_{41}, \nu_{42}, \nu_{43} \}$ we see that ${n_0}^j cnT$ must fix $\{ 2, 3 \}$, so it must fix $2$ and $3$ and hence $8$ and $7$, and finally as ${n_0}^j cn$ fixes $\{ \nu_{31}, \nu_{32}, \nu_{33} \}$ we see that ${n_0}^j cnT$ must fix both $4$ and $6$, so ${n_0}^j cnT = 1$, whence $n = c^{-1} {n_0}^{-j} \in C'$; moreover, if in fact $n \in C_G(y)$ we must have $n \in C$. Therefore in each case we do indeed have $C_G(y) = C$ and $C_G(\langle y \rangle) = C'$.

Thus the conditions of Lemma~\ref{lem: generic stabilizer from exact subset} hold; so the triple $(G, \lambda, p)$ has generic stabilizer $C/G_V \cong {\Z_2}^3$, or ${\Z_2}^4$, or ${\Z_3}^2$, respectively, while the quadruple $(G, \lambda, p, 1)$ has generic stabilizer $C'/Z(G) \cong {\Z_2}^3$, or ${\Z_2}^4$, or ${\Z_3}^2.\Z_2$, respectively.
\end{proof}

Now we employ a similar approach to deal with the first quadruples associated to the two triples treated in Section~\ref{sect: two triples}.

\begin{prop}\label{prop: B_2, omega_1 + omega_2 module, p = 5, C_4, omega_3 module, p = 3, k = 1}
Let $G = B_2$ and $\lambda = \omega_1 + \omega_2$ with $p = 5$, or $G = C_4$ and $\lambda = \omega_3$ with $p = 3$. Then the first quadruple $(G, \lambda, p, 1)$ has generic stabilizer $\Z_2$.
\end{prop}

\begin{proof}
The argument is very similar to that of Proposition~\ref{prop: A_8, omega_3, A_7, omega_4, D_8, omega_8 modules, special characteristic}; as there we shall work in $V$ and deduce the result for $\G{1}(V)$. Number the cases (i) and (ii) according as $G = B_2$ or $C_4$; whenever we give two choices followed by the word `respectively' we are taking the cases in the order (i), (ii). Note that $Z(G) = \langle z \rangle$ where $z = h_{\alpha_2}(-1)$ or $h_{\alpha_1}(-1) h_{\alpha_3}(-1)$ respectively, and $z$ acts on $V$ as negation.

First suppose $G = B_2$. Here we begin with the tensor product $L(\omega_1) \otimes L(\omega_2)$; the first factor is the natural module $V_{nat}$ for $B_2$, with basis $v_0, v_1, v_{-1}, v_2, v_{-2}$, while we may regard the second factor as the natural module for $C_2$, with basis $e_1, f_1, e_2, f_2$, where the numbering of the simple roots of $B_2$ means that it is $x_{\alpha_1}(t)$ which sends $f_2 \mapsto f_2 + te_2$ and fixes $e_1$, $e_2$ and $f_1$, and $x_{\alpha_2}(t)$ which sends $e_2 \mapsto e_2 + te_1$ and $f_1 \mapsto f_1 - tf_2$ and fixes $e_1$ and $f_2$. This tensor product contains two submodules $X_1$ and $X_2$, where $X_1$ comprises those vectors such that the coefficients of the vectors in each of the sets
\begin{eqnarray*}
& & \{ v_0 \otimes e_1, v_1 \otimes f_2, -v_2 \otimes e_2 \}, \\
& & \{ v_0 \otimes e_2, v_1 \otimes f_1, v_{-2} \otimes e_1 \}, \\
& & \{ v_0 \otimes f_2, -v_{-1} \otimes e_1, v_2 \otimes f_1 \}, \\
& & \{ v_0 \otimes f_1, -v_{-1} \otimes e_2, -v_{-2} \otimes f_2 \}
\end{eqnarray*}
sum to $0$, and
\begin{eqnarray*}
X_2 & = & \langle v_0 \otimes e_1 + 2v_1 \otimes f_2 - 2v_2 \otimes e_2, \\
    &   & \phantom{\langle} v_0 \otimes e_2 + 2v_1 \otimes f_1 + 2v_{-2} \otimes e_1, \\
    &   & \phantom{\langle} v_0 \otimes f_2 - 2v_{-1} \otimes e_1 + 2v_2 \otimes f_1, \\
    &   & \phantom{\langle} v_0 \otimes f_1 - 2v_{-1} \otimes e_2 - 2v_{-2} \otimes f_2 \rangle.
\end{eqnarray*}
The module $X_1$ is the Weyl $G$-module of high weight $\omega_1 + \omega_2$; since $p = 5$ we have $X_2 \subset X_1$, and $V = X_1/X_2$. Indeed, the matrices given in the proof of Proposition~\ref{prop: B_2, omega_1 + omega_2, p = 5, nets} may be obtained by taking the following ordered basis of $V$:
\begin{eqnarray*}
v_1 \otimes e_1 + X_2, & v_2 \otimes e_1 + X_2, & v_1 \otimes e_2 + X_2, \\
v_0 \otimes e_1 - v_1 \otimes f_2 + X_2, & v_2 \otimes f_2 + X_2, & v_0 \otimes e_2 - v_1 \otimes f_1 + X_2, \\
v_0 \otimes f_2 - v_2 \otimes f_1 + X_2, & v_{-2} \otimes e_2 + X_2, & v_0 \otimes f_1 + v_{-2} \otimes f_2 + X_2, \\
v_{-1} \otimes f_2 + X_2, & v_{-2} \otimes f_1 + X_2, & v_{-1} \otimes f_1 + X_2.
\end{eqnarray*}

Now suppose instead $G = C_4$. Here we begin with the exterior power $\bigwedge^3(V_{nat})$; this contains two submodules $X_1$ and $X_2$, where $X_1$ comprises those vectors such that for $i = 1, \dots, 4$ the sum as $j$ varies of the coefficients of both $e_i \wedge e_j \wedge f_j$ and $f_i \wedge f_j \wedge e_j$ is $0$, and
$$
X_2 = {\ts \langle e_i \wedge \sum_{j \neq i} (e_j \wedge f_j), f_i \wedge \sum_{j \neq i} (f_j \wedge e_j) : i = 1, \dots, 4 \rangle}.
$$
The module $X_1$ is the Weyl $G$-module of high weight $\omega_3$; since $p = 3$ we have $X_2 \subset X_1$, and $V = X_1/X_2$. Thus the basis vectors of $V$ are of the form $v + X_2$ where $v$ is either $v_1 \wedge v_2 \wedge v_3$ such that each $v_i$ is either $e_{j_i}$ or $f_{j_i}$ for distinct $j_1, j_2, j_3$, or $e_i \wedge (e_{i + 1} \wedge f_{i + 1} - e_{i + 2} \wedge f_{i + 2})$ or $f_i \wedge (f_{i + 1} \wedge e_{i + 1} - f_{i + 2} \wedge e_{i + 2})$ where $i \in \{ 1, \dots, 4 \}$ and subscripts are taken mod $4$.

For $1 \leq i \leq \ell$ and $j = 1, 2$ we define vectors $x_{ij}$ as follows: in case (i) we set
\begin{eqnarray*}
x_{11} = v_1 \otimes e_2 + X_2, & & x_{12} = v_{-1} \otimes f_2 + X_2, \\
x_{21} = v_2 \otimes e_1 + X_2, & & x_{22} = v_{-2} \otimes f_1 + X_2;
\end{eqnarray*}
in case (ii) we set
\begin{eqnarray*}
& x_{11} = e_1 \wedge e_2 \wedge e_3 + X_2, \quad x_{12} = f_1 \wedge f_2 \wedge f_3 + X_2, & \\
& x_{21} = f_1 \wedge e_2 \wedge f_4 + X_2, \quad x_{22} = e_1 \wedge f_2 \wedge e_4 + X_2, & \\
& x_{31} = e_1 \wedge f_3 \wedge f_4 + X_2, \quad x_{32} = f_1 \wedge e_3 \wedge e_4 + X_2, & \\
& x_{41} = f_2 \wedge e_3 \wedge f_4 + X_2, \quad x_{42} = e_2 \wedge f_3 \wedge e_4 + X_2. &
\end{eqnarray*}
For each pair $(i, j)$ let $\nu_{ij}$ be the weight such that $V_{\nu_{ij}} = \langle x_{ij} \rangle$; thus for each $i$ we have $\nu_{i1} + \nu_{i2} = 0$. For $1 \leq i \leq \ell$ set $y_i = x_{i1} + \eta_4 x_{i2}$. Define
$$
Y = \langle y_i : 1 \leq i \leq \ell \rangle;
$$
note that in each case $\codim Y = \dim G$.

Writing $h_i$ for $h_{\alpha_i}$, let $\S \leq \L(T)$ be
$$
\begin{array}{ll}
\langle h_1 + 2h_2 \rangle           & \vstrut \hbox{in case~(i),} \\
\langle h_1 + h_4, h_2 + h_4 \rangle & \vstrut \hbox{in case~(ii).} \\
\end{array}
$$
It is easy to see that if $\alpha \in \Phi$ then there exists $h \in \S$ with $[h e_\alpha] \neq 0$, so $C_{\L(G)}(\S) = \L(T)$. A straightforward calculation shows that the span of the vectors $x_{ij}$ is the subspace of $V$ annihilated by the subalgebra $\S$.

Set
\begin{eqnarray*}
\hat Y & = & \left\{ {\ts\sum}_i a_i y_i : {\ts\prod}_i a_i \neq 0, \ {a_i}^2 \neq \pm {a_{i'}}^2, \pm 2{a_{i'}}^2 \hbox{ for } i \neq i' \right\}
\end{eqnarray*}
or
\begin{eqnarray*}
\hat Y & = & \left\{ {\ts\sum}_i a_i y_i : {\ts\prod}_i a_i \neq 0, \ {a_i}^2 \neq \pm {a_{i'}}^2 \hbox{ for } i \neq i', \right. \\
       &   & \left.\phantom{\left\{ {\ts\sum}_i a_i y_i : \right.} {\ts\sum}_{j \in \F_3} ({a_{i_j}}^2 \pm {a_{i_{j + 1}}}^2)^2 \neq 0 \hbox{ for } i_1, i_2, i_3 \hbox{ distinct} \right\}
\end{eqnarray*}
respectively (with $\F_3$ the field of size $3$); then $\hat Y$ is a dense open subset of $Y$. Take
$$
y = {\ts\sum}_{i, j} a_i y_i \in \hat Y.
$$

First suppose $x \in \Ann_{\L(G)}(y)$; write $x = h + e$ where $h \in \L(T)$ and $e \in \langle e_\alpha : \alpha \in \Phi \rangle$. Clearly $h.y \in Y$; since the difference of two weights $\nu_{ij}$ is never a root, for each pair $(i, j)$ we see that $e.y$ contains no term $x_{ij}$. Thus we must have $h.y = e.y = 0$. A quick calculation shows that we must have $h \in \S$. Now write $e = \sum_{\alpha \in \Phi} t_\alpha e_\alpha$; then the equation $e.y = 0$ may be expressed in matrix form as $A{\bf t} = {\bf 0}$, where $A$ is an $M \times M$ matrix and ${\bf t}$ is a column vector whose entries are the various coefficients $t_\alpha$. We find that if the rows and columns of $A$ are suitably ordered then it becomes block diagonal, having $4$ or $8$ blocks respectively, with each block being a $2 \times 2$ or $4 \times 4$ matrix respectively. In case (i) the blocks are
$$
\left(
  \begin{array}{cc}
    a_1 & 2\eta_4 a_2 \\
    a_2 & -\eta_4 a_1 \\
  \end{array}
\right),
\quad
\left(
  \begin{array}{cc}
      a_1 & \eta_4 a_2 \\
    -2a_2 & \eta_4 a_1 \\
  \end{array}
\right),
\quad
\left(
  \begin{array}{cc}
        a_1    &     a_2     \\
    \eta_4 a_2 & 2\eta_4 a_1 \\
  \end{array}
\right),
\quad
\left(
  \begin{array}{cc}
       2a_1    &     a_2    \\
    \eta_4 a_2 & \eta_4 a_1 \\
  \end{array}
\right),
$$
each of which has determinant a scalar multiple of $2{a_1}^2 - {a_2}^2$. In case (ii), after some negation of columns each block may be written in the form
$$
\left(
  \begin{array}{cccc}
                       &       a_{i_1}       &    \e_1 a_{i_2}    & \e_2\eta_4 a_{i_3} \\
         a_{i_1}       &                     &    \e_3 a_{i_3}    & \e_4\eta_4 a_{i_2} \\
       \e_1 a_{i_2}    &    -\e_3 a_{i_3}    &                    &   \eta_4 a_{i_1}   \\
   -\e_2\eta_4 a_{i_3} &  \e_4\eta_4 a_{i_2} &   \eta_4 a_{i_1}   &                    \\
  \end{array}
\right)
\quad \hbox{or} \quad
\left(
  \begin{array}{cccc}
                       &    \eta_4 a_{i_1}   & \e_1\eta_4 a_{i_2} & \e_2 a_{i_3} \\
      \eta_4 a_{i_1}   &                     & \e_3\eta_4 a_{i_3} & \e_4 a_{i_2} \\
    \e_1\eta_4 a_{i_2} & -\e_3\eta_4 a_{i_3} &                    &    a_{i_1}   \\
      -\e_2 a_{i_3}    &     \e_4 a_{i_2}    &       a_{i_1}      &              \\
  \end{array}
\right)
$$
where $i_1, i_2, i_3$ are distinct and $\e_1, \dots, \e_4 \in \{ \pm1 \}$; each of these matrices has determinant $\sum_{j \in \F_3} ({a_{i_j}}^2 \pm {a_{i_{j + 1}}}^2)^2$. Thus in each case the final condition in the definition of the set $\hat Y$ implies that each block of $A$ is non-singular, as therefore is $A$ itself; so ${\bf t}$ must be the zero vector and hence $e = 0$. Thus $x = h + e \in \S$; so $\Ann_{\L(G)}(y) = \S$.

A straightforward calculation shows that $C_T(y) = \{ 1 \}$, and $T.y \cap Y = \{ \pm y \}$. Take $w \in W$, and suppose there exists $n \in \Tran_G(y, Y)$ with $nT = w$; as the minimal sets of weights $\nu_{ij}$ summing to zero are the sets $\{ \nu_{i1}, \nu_{i2} \}$ for $i = 1, \dots, \ell$, we see that $n$ must permute these sets. We may write $n = s n^*$, where $n^*$ is a product of elements $n_\alpha$ for various roots $\alpha$, and $s \in T$. Take $i \leq \ell$, then there exists $i' \leq \ell$ such that $n.y_i = cy_{i'}$ for some $c \in K^*$; as for each $j$ the element $n^*$ must send $x_{ij}$ to $\pm x_{i'j'}$ for some $j'$, we must have $c(x_{i'1} + \eta_4 x_{i'2}) = n.(x_{i1} + \eta_4 x_{i2})$ which is either $s.(\pm x_{i'1} \pm \eta_4 x_{i'2}) = \pm \nu_{i'1}(s) x_{i'1} \pm \eta_4 \nu_{i'2}(s) x_{i'2}$ or $s.(\pm x_{i'2} \pm \eta_4 x_{i'1}) = \pm \nu_{i'2}(s) x_{i'2} \pm \eta_4 \nu_{i'1}(s) x_{i'1}$, and taking the product of the coefficients on both sides gives $\eta_4 c^2 = \pm \eta_4 \nu_{i'1}(s) \nu_{i'2}(s) = \pm \eta_4 (\nu_{i'1} + \nu_{i'2})(s) = \pm \eta_4$, whence $c^2 = \pm 1$. Thus $n$ permutes and possibly scales by a power of $\eta_4$ the vectors $y_i$, so sends $y = \sum a_iy_i$ to $\sum {\eta_4}^{b_i} a_{\pi(i)} y_i$ for some permutation $\pi$ of $\{ 1, \dots, \ell \}$ and some integers $b_1, \dots, b_\ell$. It now follows that $n.y \in \hat Y$. Thus $N.y \cap Y \subset \hat Y$. Since all the conditions of Lemma~\ref{lem: exactness via Premet} hold, it follows that $\Tran_G(y, Y) \subseteq N$, and $y$ is $Y$-exact.

We thus have $C_G(\langle y \rangle) \leq N$ (since $\Tran_G(\langle y \rangle, \G{1}(Y)) = \Tran_G(y, Y)$); as $C_T(\langle y \rangle) = \langle z \rangle$, each coset of $T$ in $N$ can contain at most two elements of $C_G(\langle y \rangle)$. Given $n \in C_G(\langle y \rangle)$, since by the above $n$ permutes and possibly scales by a power of $\eta_4$ the vectors $y_i$, the final or penultimate condition respectively in the definition of the set $\hat Y$ shows that $n$ must in fact fix each line $\langle y_i \rangle$, and thus each set $\{ \nu_{i1}, \nu_{i2} \}$.

Set $n_0$ to be
$$
\begin{array}{ll}
n_{\alpha_1} n_{\alpha_1 + 2\alpha_2}                                                                                       & \vstrut \hbox{in case~(i),} \\
n_{\alpha_4} n_{2\alpha_3 + \alpha_4} n_{2\alpha_2 + 2\alpha_3 + \alpha_4} n_{2\alpha_1 + 2\alpha_2 + 2\alpha_3 + \alpha_4} & \vstrut \hbox{in case~(ii),} \\
\end{array}
$$
so that $n_0 T$ is the long word of the Weyl group, and ${n_0}^2 = z$. Let $C = \langle n_0 \rangle$. Since for each $i$ we have $n_0.x_{i1} = -x_{i2}$ and $n_0.x_{i2} = x_{i1}$, we see that $n_0.y = \eta_4 y$, and so $C \leq C_G(\langle y \rangle)$; we shall show that in fact $C_G(\langle y \rangle) = C$.

Take $n \in C_G(\langle y \rangle)$. In case (i) the $W$-stabilizer of the weight $\nu_{11}$ is trivial; in case (ii) it is $\langle w_{\alpha_1}, w_{\alpha_2}, w_{\alpha_4} \rangle$, of which elements only $1$ and $w_{\alpha_1} w_{\alpha_4}$ stabilize the set $\{ \nu_{21}, \nu_{22} \}$, and the latter element interchanges the sets $\{ \nu_{31}, \nu_{32} \}$ and $\{ \nu_{41}, \nu_{42} \}$. In both cases we thus have $nT = T$ or $n_0T$, and so $n \in \{ 1, z, n_0, n_0z \} = \langle n_0 \rangle = C$. Therefore in each case we do indeed have $C_G(\langle y \rangle) = C$.

Thus the conditions of Lemma~\ref{lem: generic stabilizer from exact subset} hold; so the quadruple $(G, \lambda, p, 1)$ has generic stabilizer $C/Z(G) \cong \Z_2$.
\end{proof}

The remaining cases in this section may be treated using the material of Section~\ref{sect: invariants}. The first is very straightforward.

\begin{prop}\label{prop: A_1, 3omega_1 module}
Let $G = A_1$ and $\lambda = 3\omega_1$ with $p \geq 5$. Then the triple $(G, \lambda, p)$ has generic stabilizer $\Z_3$, and there is a regular orbit; the associated first quadruple $(G, \lambda, p, 1)$ has generic stabilizer $S_3$.
\end{prop}

\begin{proof}
We take $G = \SL_2(K)$. Recall that $V_{nat} = \langle v_1, v_2 \rangle$. We may identify $V$ with $S^3(V_{nat})$, the space of homogeneous polynomials in $v_1$ and $v_2$ of degree $3$, so that $V = \langle {v_1}^3, {v_1}^2v_2, v_1{v_2}^2, {v_2}^3 \rangle$. We write $G^+ = \GL_2(K)$ and extend the action of $G$ on $V$ to $G^+$ in the obvious way; clearly for all $\kappa \in K^*$ and $v \in V$ we have $(\kappa I).v = \kappa^3v$, so that $G_V = \{ 1 \}$.

Given $v = a_3{v_1}^3 + a_2{v_1}^2v_2 + a_1v_1{v_2}^2 + a_0{v_2}^3 \in V$, define $A_v$ to be the symmetric $2 \times 2$ matrix
$$
\left(
  \begin{array}{cc}
    6a_3a_1 - 2{a_2}^2 &  9a_3a_0 - a_2a_1  \\
     9a_3a_0 - a_2a_1  & 6a_2a_0 - 2{a_1}^2 \\
  \end{array}
\right).
$$
We then find that for all $g \in G$ we have
$$
A_{g.v} = [(\det g)g] A_v [(\det g)g]^T = (\det g)^2 gA_vg^T;
$$
indeed it suffices to check this for root elements and scalar multiples of $I$. Thus $\det A_{g.v} = (\det g)^6 \det A_v$; so the map $f : v \mapsto \det(A_v)$ is the relative invariant, and the associated character $\chi$ is given by $\chi(g) = (\det g)^6$.

Take $y_0 = v_1v_2(v_1 + v_2)$; then any $g \in C_{G^+}(\langle y_0 \rangle)$ must permute and scale the three linear factors of $y_0$. Write
$$
g_1 = \left(
   \begin{array}{cc}
     0 & -1 \\
     1 & -1 \\
   \end{array}
\right),
\quad
g_2 = \left(
   \begin{array}{cc}
     -\eta_4 & \eta_4 \\
        0    & \eta_4 \\
   \end{array}
\right),
$$
so that $g_1, g_2 \in G$. Then up to scaling $g_1$ cycles the three linear factors, while $g_2$ fixes the first linear factor while interchanging the second and third; moreover $g_1.y_0 = y_0$ while $g_2.y_0 = \eta_4 y_0$. Since $G/Z(G)$ acts sharply $3$-transitively on lines in $V_{nat}$, we have $C_{G^+}(\langle y_0 \rangle) = Z(G^+)\langle g_1, g_2 \rangle$; so $C_{G^+}(y_0) = \langle g_1, \eta_4 g_2 \rangle$, which is finite, and $C_G(\langle y_0 \rangle) = \langle g_1, g_2 \rangle$ while $C_G(y_0) = \langle g_1 \rangle$. Thus by Lemma~\ref{lem: gen stab in invariant situation}(i) the triple $(G, \lambda, p)$ has generic stabilizer $C_G(y_0)/G_V \cong\Z_3$, while the quadruple $(G, \lambda, p, 1)$ has generic stabilizer $C_{G^+}(\langle y_0 \rangle)/Z(G^+) \cong S_3$.

Now take $y_1 = {v_1}^2v_2 \in V(0)$. Any element of $C_G(y_1)$ must fix each of the lines $\langle v_1 \rangle$ and $\langle v_2 \rangle$, so must be diagonal; since $\diag(\kappa, \kappa^{-1}).y_1 = \kappa y_1$, we see that $C_G(y_1) = \{ 1 \}$. Thus there is a regular orbit in the action of $G$ on $V$.
\end{proof}

We next treat a case where the result may be easily obtained from the literature; the triple $(G, \lambda, p)$ concerned is $(A_3, \omega_1 + \omega_2, 3)$. This action was first investigated by Chen in \cite{Z}: he showed that there is a dense $G^+$-orbit and obtained representatives of it and several others, in each case giving the stabilizers in both $G^+$ and $\L(G^+)$; moreover he proved that in this action there must be a relative invariant of degree $8$, with associated character $\chi$ given by $\chi(g) = (\det g)^6$. Subsequently Cohen and Wales in \cite{CW} built upon Chen's work and obtained a complete set of orbit representatives, together with their stabilizers in $\L(G^+)$ and in some cases in $G^+$ as well. However, although the Lie algebra calculations are straightforward, Chen gave no proof that the group stabilizers were as stated, while Cohen and Wales employed computer calculations to obtain their results. Moreover, Chen did not find the invariant, saying \lq The determination of this relative invariant must be very interesting'; Cohen and Wales referred to it in passing, mentioning a computer calculation showing the existence of a relative invariant of degree $8$ on $S^3(V_{nat})$ over a field of arbitrary characteristic other than $2$, but did not give it explicitly. For the sake of both completeness and independence of computer calculations, we provide full details here.

\begin{prop}\label{prop: A_3, omega_1 + omega_2 module}
Let $G = A_3$ and $\lambda = \omega_1 + \omega_2$ with $p = 3$. Then the triple $(G, \lambda, p)$ has generic stabilizer $Alt_5$, and there is no regular orbit; the associated first quadruple $(G, \lambda, p, 1)$ has generic stabilizer $S_5$.
\end{prop}

\begin{proof}
We take $G = \SL_4(K)$. Recall that $V_{nat} = \langle v_1, v_2, v_3, v_4 \rangle$. As in \cite{CW}, we may identify $V$ with the quotient of $S^3(V_{nat})$, the space of homogeneous polynomials in $v_1, v_2, v_3, v_4$ of degree $3$, by the subspace $\langle {v_1}^3, {v_2}^3, {v_3}^3, {v_4}^3 \rangle$. Given $1 \leq i_1, i_2, i_3 \leq 4$, write $x_{i_1i_2i_3}$ for the image of $v_{i_1}v_{i_2}v_{i_3}$ in $V$; thus the subscripts in a vector $x_{i_1i_2i_3}$ may be freely permuted. The $16$-dimensional module $V$ then has a basis consisting of $4$ vectors $x_{i_1i_2i_3}$ with $i_1, i_2, i_3$ all distinct, and $12$ vectors $x_{i_1i_1i_2}$ with $i_1, i_2$ distinct. We write $G^+ = \GL_4(K)$ and extend the action of $G$ on $V$ to $G^+$ in the obvious way; clearly for all $\kappa \in K^*$ and $v \in V$ we have $(\kappa I).v = \kappa^3v$, so that $G_V = \{ 1 \}$.

Given $v = \sum a_{i_1i_2i_3} x_{i_1i_2i_3} + \sum a_{i_1i_1i_2} x_{i_1i_1i_2} \in V$, define $A_v$ to be the $4 \times 4$ symmetric matrix with $(i_1, i_1)$-entry
$$
a_{i_2i_2i_4}a_{i_3i_3i_4} + a_{i_2i_2i_3}a_{i_4i_4i_3} + a_{i_3i_3i_2}a_{i_4i_4i_2} - {a_{i_2i_3i_4}}^2
$$
and $(i_1, i_2)$-entry
$$
a_{i_1i_2i_3}a_{i_4i_4i_3} + a_{i_1i_2i_4}a_{i_3i_3i_4} + a_{i_3i_3i_1}a_{i_4i_4i_2} + a_{i_4i_4i_1}a_{i_3i_3i_2} + a_{i_3i_4i_1}a_{i_3i_4i_2},
$$
where we write $\{ i_1, i_2, i_3, i_4 \} = \{ 1, 2, 3, 4 \}$. Recall that if $g \in G^+$ the adjugate matrix $\adj g$ satisfies $g(\adj g) = (\det g)I$, so that $\adj g = (\det g)g^{-1}$. We then find that for all $g \in G^+$ we have
$$
A_{g.v} = (\adj g)^T A_v (\adj g) = (\det g)^2 (g^{-1})^T A_v g^{-1};
$$
indeed it suffices to check this for root elements and scalar multiples of $I$. Thus $\det A_{g.v} = (\det g)^6 \det A_v$; so the map $f : v \mapsto \det A_v$ is the relative invariant, and the associated character $\chi$ is given by $\chi(g) = (\det g)^6$.

Take
$$
y_0 = x_{123} + x_{124} + x_{134} + x_{234}
$$
and suppose $g \in C_{G^+}(\langle y_0 \rangle)$; write $g = (a_{ij})$. Define the matrix $x$ whose $(i, j)$-entry is $a_{ji'}a_{ji''} + a_{ji''}a_{ji'''} + a_{ji'''}a_{ji'}$ where $\{ i, i', i'', i''' \} = \{ 1, 2, 3, 4 \}$. Then in the product $gx$ the $(i, i)$-entry is
$$
\sum_j a_{ij}\left(\sum_{\genfrac{}{}{0pt}{}{j',j'' \neq j,}{j' \neq j''}} a_{ij'}a_{ij''}\right) = 3\sum_{j, j', j'' \ \ss{\mathrm{distinct}}} a_{ij}a_{ij'}a_{ij''} = 0,
$$
while for $i \neq j$ the $(i, j)$-entry is the coefficient of ${v_j}^2v_i$ in $g.y_0$, which is $0$ as $g.y_0 \in \langle y_0 \rangle$. As $g$ is non-singular, all entries in the matrix $x$ must be zero; so for fixed $j$ we have
\begin{eqnarray*}
a_{j2}a_{j3} + a_{j3}a_{j4} + a_{j4}a_{j2} & = & 0, \\
a_{j1}a_{j3} + a_{j3}a_{j4} + a_{j4}a_{j1} & = & 0, \\
a_{j1}a_{j2} + a_{j2}a_{j4} + a_{j4}a_{j1} & = & 0, \\
a_{j1}a_{j2} + a_{j2}a_{j3} + a_{j3}a_{j1} & = & 0.
\end{eqnarray*}
These equations are certainly satisfied if $a_{j1} = a_{j2} = a_{j3} = a_{j4}$, so suppose (say) $a_{j1} \neq a_{j2}$. Subtracting the first from the second gives $(a_{j1} - a_{j2})(a_{j3} + a_{j4}) = 0$, so we must have $a_{j3} + a_{j4} = 0$; the first now reduces to $a_{j3}a_{j4} = 0$, so that $a_{j3} = a_{j4} = 0$, while the sum of the third and fourth gives $2a_{j1}a_{j2} = 0$, so that (say) $a_{j2} = 0$. Thus in each row of $g$ either all four entries are equal, or three of the four entries are zero. Requiring the four basis vectors occurring in $y_0$ to have equal coefficient in $g.y_0$ now shows that there exists $\kappa \in K^*$ such that the rows of $\kappa^{-1} g$ are four of $(-1 \ -1 \ -1 \ -1)$, $(1 \ 0 \ 0 \ 0)$, $(0 \ 1 \ 0 \ 0)$, $(0 \ 0 \ 1 \ 0)$, $(0 \ 0 \ 0 \ 1)$. Thus $C_{G^+}(\langle y_0 \rangle) = Z(G^+)S_5$ where the $S_5$ is generated by the permutation matrices together with
$$
\left(
  \begin{array}{cccc}
     -1   &   -1   &   -1   & -1 \\
   \pmin1 &        &        &    \\
          & \pmin1 &        &    \\
          &        & \pmin1 &    \\
  \end{array}
\right);
$$
so $C_{G^+}(y_0) = S_5$, which is finite, and $C_G(\langle y_0 \rangle) = S_5$ while $C_G(y_0) = Alt_5$. Thus by Lemma~\ref{lem: gen stab in invariant situation}(i) the triple $(G, \lambda, p)$ has generic stabilizer $C_G(y_0)/G_V \cong Alt_5$, while the quadruple $(G, \lambda, p, 1)$ has generic stabilizer $C_{G^+}(\langle y_0 \rangle)/Z(G^+) \cong S_5$.

Now take
$$
y_1 = x_{133} + x_{224} + x_{114} \in V(0),
$$
and suppose $g \in C_G(y_1)$; write $g = u_1 n u_2$, with $u_1 \in U$, $n \in N$ and $u_2 \in U_w$ where $w = nT \in W$, then we have ${u_1}^{-1}.y_1 = nu_2.y_1$. Let $\nu_1, \nu_2, \nu_3$ be the weights with $x_{133} \in V_{\nu_1}$, $x_{224} \in V_{\nu_2}$, $x_{114} \in V_{\nu_3}$, and write $\Lambda' = \{ \mu \in \Lambda(V) : \exists i \hbox{ with } \nu_i \preceq \mu \}$; then each weight occurring in ${u_1}^{-1}.y_1$ lies in $\Lambda'$. Now $u_2.y_1$ contains $x_{133}$ and $x_{224}$, and if it does not contain $x_{114}$ then $u_2$ must involve a root element in $X_{1 - 2}$ (where we write $X_{i - j}$ for the root group corresponding to the root $\ve_i - \ve_j$), so that it must contain $x_{124}$, the weight corresponding to which is $\frac{1}{2}(\nu_2 + \nu_3)$. Thus $w(\nu_1)$, $w(\nu_2)$, and at least one of $w(\nu_3)$ and $w(\frac{1}{2}(\nu_2 + \nu_3))$ must all lie in $\Lambda'$; a straightforward check shows that this forces $w = 1$, so that $u_2 = 1$ and $g = u_1 t$ with $t \in T$. Equating coefficients of weight vectors in the order $x_{124}$, $x_{223}$, $x_{123}$, $x_{113}$, $x_{122}$, $x_{112}$ shows that the projection of $u_1$ must be trivial on each of the root groups $X_{1 - 2}$, $X_{3 - 4}$, $X_{2 - 3}$, $X_{1 - 3}$, $X_{1 - 4}$, $X_{2 - 4}$ in turn; thus $u_1 = 1$. Finally write $t = \diag(\kappa_1, \kappa_2, \kappa_3, \kappa_4)$ with $\kappa_1\kappa_2\kappa_3\kappa_4 = 1$; then $\kappa_1{\kappa_3}^2 = {\kappa_2}^2\kappa_4 = {\kappa_1}^2\kappa_4 = 1$, so $\kappa_1 = \kappa_4 = 1$, $\kappa_2 = \kappa_3 = \pm 1$. Hence $C_G(y_1) = \langle \diag(-1, 1, -1, 1) \rangle$, which is finite and non-trivial. Thus by Lemma~\ref{lem: gen stab in invariant situation}(ii) there is no regular orbit in the action of $G$ on $V$.
\end{proof}

Finally we turn to the two cases where the triple is not $p$-restricted.

\begin{prop}\label{prop: A_ell, omega_1 + q omega_1 and omega_1 + q omega_ell modules}
Let $G = A_\ell$ and $\lambda = \omega_1 + q\omega_1$ or $\omega_1 + q\omega_\ell$. Then the triple $(G, \lambda, p)$ has generic stabilizer $\PSU_{\ell + 1}(q)$ or $\PSL_{\ell + 1}(q)$ respectively, and there is no regular orbit unless $\ell = 1$ and $q \leq 3$; the associated first quadruple $(G, \lambda, p, 1)$ has generic stabilizer $\PGU_{\ell + 1}(q)$ or $\PGL_{\ell + 1}(q)$ respectively.
\end{prop}

\begin{proof}
As in the proof of Proposition~\ref{prop: A_ell, omega_1 + q omega_1 and omega_1 + q omega_ell, k geq 2}, we take $G = \SL_{\ell + 1}(K)$ and identify $V$ with the space of $(\ell + 1) \times (\ell + 1)$ matrices $D$ over $K$, where $A \in G$ maps $D \mapsto AD(A^{(q)})^T$ or $D \mapsto AD(A^{(q)})^{-1}$ according as $\lambda = \omega_1 + q\omega_1$ or $\omega_1 + q\omega_\ell$. We write $G^+ = \GL_{\ell + 1}(K)$, and extend the action of $G$ on $V$ to $G^+$ in the obvious way. We have the invariant $\det \in K[V]$, and $\kappa I \in G^+$ maps $D \mapsto \kappa^{1 + q}D$ or $D \mapsto \kappa^{1 - q}D$ respectively. In each case we shall take $y_0 = I$.

First assume $\lambda = \omega_1 + q\omega_1$; then $G_V = \{ \kappa I \in G : \kappa^{q + 1} = 1 \} = Z(\SU_{\ell + 1}(q))$. Also $C_{G^+}(y_0) = \{ A \in G^+ : AI(A^{(q)})^T = I \} = \GU_{\ell + 1}(q)$, so $C_G(y_0) = C_{G^+}(y_0) \cap G = \SU_{\ell + 1}(q)$; and $C_{G^+}(\langle y_0 \rangle) = \{ A \in G^+ : AI(A^{(q)})^T \in \langle I \rangle \} = \{ A \in G^+ : A^{(q)} = \kappa (A^{-1})^T \hbox{ for some } \kappa \in K^* \} = Z(G^+)\GU_{\ell + 1}(q)$. Thus by Lemma~\ref{lem: gen stab in invariant situation}(i) the triple $(G, \lambda, p)$ has generic stabilizer $\SU_{\ell + 1}(q)/Z(\SU_{\ell + 1}(q)) = \PSU_{\ell + 1}(q)$, and the quadruple $(G, \lambda, p, 1)$ has generic stabilizer $Z(G^+)\GU_{\ell + 1}(q)/Z(G^+) \cong \PGU_{\ell + 1}(q)$. Set
$$
y_1 = \left(
\begin{array}{ccc}
 I_{\ell - 1} & 0 & 0 \\
 0            & 0 & 0 \\
 0            & 1 & 0 \\
\end{array}
\right) \in V(0).
$$
We regard $V_{nat}$ as the space of column vectors, with standard basis $v_1, \dots, v_{\ell + 1}$. We have $\{ v \in V_{nat} : y_1v = 0 \} = \langle v_{\ell + 1} \rangle$, and
$\{ v \in V_{nat} : v^Ty_1 = 0 \} = \langle v_\ell \rangle$. Suppose $A \in G^+$ stabilizes $y_1$. Then $A y_1(A^{(q)})^T v_{\ell + 1} = y_1v_{\ell + 1} = 0$, so $y_1(A^{(q)})^T v_{\ell + 1} = 0$, and hence $(A^{(q)})^T v_{\ell + 1} \in \langle v_{\ell + 1} \rangle$; similarly ${v_\ell}^T A y_1(A^{(q)})^T = {v_\ell}^T y_1 = 0$, so ${v_\ell}^T A y_1 = 0$, and hence $A^T v_\ell \in \langle v_\ell \rangle$. Thus
$$
A = \left(
\begin{array}{ccc}
 R & x & x' \\
 0 & \kappa & 0  \\
 0 & 0 & \kappa' \\
\end{array}
\right),
$$
where $\kappa, \kappa' \in K^*$, the column vectors $x, x'$ have length $\ell - 1$, and $R \in \GL_{\ell - 1}(K)$. Equating $A y_1(A^{(q)})^T$ and $y_1$ now shows that $x = x' = 0$, while $\kappa'\kappa^q = 1$ and $R(R^{(q)})^T = I_{\ell - 1}$, so that
$$
A = \left(
\begin{array}{ccc}
 R & 0 & 0      \\
 0 & \kappa & 0      \\
 0 & 0 & \kappa^{-q} \\
\end{array}
\right)
$$
where $\kappa \in K^*$ and $R \in GU_{\ell - 1}(q)$; the $G$-stabilizer of $y_1$ therefore consists of such matrices where $\det R = \kappa^{q - 1}$, and thus is finite and not $G_V$, unless $\ell = 1$ and $q \leq 3$. Thus by Lemma~\ref{lem: gen stab in invariant situation}(ii) there is no regular orbit in the action of $G$ on $V$ in this case, unless $\ell = 1$ and $q \leq 3$.

Now assume $\lambda = \omega_1 + q\omega_\ell$; then $G_V = \{ \kappa I \in G : \kappa^{q - 1} = 1 \} = Z(\SL_{\ell + 1}(q))$. Also $C_{G^+}(y_0) = \{ A \in G^+ : AI(A^{(q)})^{-1} = I \} = \GL_{\ell + 1}(q)$, so $C_G(y_0) = C_{G^+}(y_0) \cap G = \SL_{\ell + 1}(q)$; and $C_{G^+}(\langle y_0 \rangle) = \{ A \in G^+ : AI(A^{(q)})^{-1} \in \langle I \rangle \} = \{ A \in G^+ : A^{(q)} = \kappa A \hbox{ for some } \kappa \in K^* \} = Z(G^+)\GL_{\ell + 1}(q)$. Thus by Lemma~\ref{lem: gen stab in invariant situation}(i) the triple $(G, \lambda, p)$ has generic stabilizer $\SL_{\ell + 1}(q)/Z(\SL_{\ell + 1}(q)) = \PSL_{\ell + 1}(q)$, and the quadruple $(G, \lambda, p, 1)$ has generic stabilizer $Z(G^+)\GL_{\ell + 1}(q)/Z(G^+) \cong \PGL_{\ell + 1}(q)$. Set
$$
y_1 = \left(
\begin{array}{cc}
 I_\ell & 0 \\
 0      & 0 \\
\end{array}
\right) \in V(0).
$$
It is easy to see that the $G$-stabilizer of $y_1$ consists of matrices
$$
A = \left(
\begin{array}{cc}
 R & 0 \\
 0 & \kappa \\
\end{array}
\right)
$$
where $R \in \GL_\ell(q)$ and $\kappa = (\det R)^{-1}$, and thus is finite and not $G_V$. Thus by Lemma~\ref{lem: gen stab in invariant situation}(ii) there is no regular orbit in the action of $G$ on $V$ in this case.
\end{proof}

This completes the justification of the entries in Table~\ref{table: large triple and first quadruple non-TGS}, and hence the proof of Theorem~\ref{thm: large triple and first quadruple generic stab}.

\section{Small triples and associated first quadruples}\label{sect: small triples and first quadruples}

In this section we shall treat small triples and associated first quadruples, and establish the entries in Tables~\ref{table: small classical triple and first quadruple generic stab} and \ref{table: small exceptional triple and first quadruple generic stab}, thus proving Theorem~\ref{thm: small triple and first quadruple generic stab}.

We begin with those triples and first quadruples where the module $V$ is the unique non-trivial composition factor of the Lie algebra $\L(G)$.

\begin{prop}\label{prop: adjoint modules}
Let $G = A_1$ and $\lambda = 2\omega_1$ with $p \geq 3$, or $G = A_\ell$ for $\ell \in [2, \infty)$ and $\lambda = \omega_1 + \omega_\ell$, or $G = B_2$ and $\lambda = 2\omega_2$ with $p \geq 3$, or $G = B_\ell$ for $\ell \in [3, \infty)$ and $\lambda = \omega_2$ with $p \geq 3$, or $G = C_\ell$ for $\ell \in [3, \infty)$ and $\lambda = 2\omega_1$ with $p \geq 3$, or $G = D_\ell$ for $\ell \in [4, \infty)$ and $\lambda = \omega_2$, or $G = E_6$ and $\lambda = \omega_2$, or $G = E_7$ and $\lambda = \omega_1$, or $G = E_8$ and $\lambda = \omega_8$, or $G = F_4$ and $\lambda = \omega_1$ with $p \geq 3$, or $G = G_2$ and $\lambda = \omega_2$ with $p \neq 3$. Then the triple $(G, \lambda, p)$ and the associated first quadruple $(G, \lambda, p, 1)$ have generic stabilizers $C_V$ and $C_{\G{1}(V)}$ respectively, where
\begin{itemize}
\item[(i)] $C_V = T_2.\Z_3$ and $C_{\G{1}(V)} = T_2.S_3$ if $G = A_2$ with $p = 3$;
\item[(ii)] $C_V = C_{\G{1}(V)} = T_3.{\Z_2}^2$ if $G = A_3$ with $p = 2$;
\item[(iii)] $C_V = C_{\G{1}(V)} = T_4.{\Z_2}^3.{\Z_2}^2$ if $G = D_4$ with $p = 2$;
\item[(iv)] $C_V = C_{\G{1}(V)} = T_\ell.{\Z_2}^{\ell - 1}$ if $G = D_\ell$ for $\ell \in [5, \infty)$ with $p = 2$;
\item[(v)] $C_V = C_{\G{1}(V)} = T_\ell$ if $G = A_2$ with $p \neq 3$, or $A_3$ with $p \geq 3$, or $A_\ell$ for $\ell \in [4, \infty)$, or $D_\ell$ for odd $\ell \in [5, \infty)$ with $p \geq 3$, or $E_6$;
\item[(vi)] $C_V = T_\ell$ and $C_{\G{1}(V)} = T_\ell.\Z_2$ if $G = A_1$, or $B_2$, or $B_\ell$, or $C_\ell$, or $D_\ell$ for even $\ell \in [4, \infty)$ with $p \geq 3$, or $E_7$ with $p \geq 3$, or $E_8$ with $p \geq 3$, or $F_4$, or $G_2$ with $p \geq 5$;
\item[(vii)] $C_V = C_{\G{1}(V)} = T_\ell.\Z_2$ if $G = E_7$ with $p = 2$, or $E_8$ with $p = 2$, or $G_2$ with $p = 2$.
\end{itemize}
\end{prop}

\begin{proof}
Take $G$ to be of simply connected type. In all these cases we have $V = \L(G)/Z(\L(G))$, so that $G_V = Z(G)$. We apply Lemma~\ref{lem: semisimple auts}, taking $H = G$ and $\theta = 1$; then $\L(T_H)_{(1)} = \L(T_H)$, so that $({W_H}^\ddagger)_{(1)} = {W_H}^\ddagger$ and $({W_H}^\dagger)_{(1)} = {W_H}^\dagger$, as given by Lemma~\ref{lem: W_H on L(T_H)/Z(L(H))}. As the quotient of a torus by a finite group is still a torus of the same rank, the triple $(G, \lambda, p)$ has generic stabilizer $T_\ell.W^\dagger$, while the quadruple $(G, \lambda, p, 1)$  has generic stabilizer $T_\ell.W^\ddagger$.
\end{proof}

Next we take the cases where $G$ is a classical group and $V$ is the natural module. In the statement of the following result, for convenience we refer to the case where $G = C_2$, $\lambda = \omega_1$; this appears in Table~\ref{table: small classical triple and first quadruple generic stab} as $G = B_2$, $\lambda = \omega_2$.

\begin{prop}\label{prop: natural modules}
Let $G = A_\ell$ for $\ell \in [1, \infty)$, or $G = B_\ell$ for $\ell \in [2, \infty)$ with $p \geq 3$, or $G = C_\ell$ for $\ell \in [2, \infty)$, or $G = D_\ell$ for $\ell \in [4, \infty)$, and $\lambda = \omega_1$. Then the triple $(G, \lambda, p)$ has generic stabilizer $A_{\ell - 1} U_\ell$, or $D_\ell$, or $C_{\ell - 1} U_{2\ell - 1}$, or $B_{\ell - 1}$, respectively; the associated first quadruple $(G, \lambda, p, 1)$ has generic stabilizer $A_{\ell - 1} T_1 U_\ell$, or $D_\ell.\Z_2$, or $C_{\ell - 1} T_1 U_{2\ell - 1}$, or $B_{\ell - 1}$, respectively.
\end{prop}

\begin{proof}
In all these cases $V = V_{nat}$. If $G = A_\ell$ or $C_\ell$, then $G$ acts transitively on $V \setminus \{ 0 \}$; if we write $QL$ for the maximal parabolic subgroup corresponding to the first simple root, where $Q$ is the unipotent radical and $L$ the Levi subgroup, then the stabilizer of the first vector of the standard basis for $V$ is $QL'$ (where $L'$ is the derived group of $L$), which is of form $A_{\ell - 1} U_{\ell - 1}$ or $C_{\ell - 1} U_{2\ell - 1}$ respectively, while that of the corresponding line is $QL$, which is of form $A_{\ell - 1} T_1 U_{\ell - 1}$ or $C_{\ell - 1} T_1 U_{2\ell - 1}$ respectively. If $G = B_\ell$ with $p \geq 3$ or $D_\ell$, then $G$ acts transitively on the sets of singular and non-singular vectors in $V \setminus \{ 0 \}$, and the latter is dense in $V$; the stabilizer of a non-singular vector is an orthogonal group on a space of dimension $\dim V - 1$, which is of form $D_\ell$ or $B_{\ell - 1}$ respectively, while that of the corresponding line is of form $D_\ell.\Z_2$ or $B_{\ell - 1}$ respectively.
\end{proof}

\begin{prop}\label{prop: natural module for B_ell, p = 2}
Let $G = B_\ell$ for $\ell \in [2, \infty)$ and $\lambda = \omega_1$ with $p = 2$. Then the triple $(G, \lambda, p)$ has generic stabilizer $B_{\ell - 1} U_{2\ell - 1}$; the associated first quadruple $(G, \lambda, p, 1)$ has generic stabilizer $B_{\ell - 1} T_1 U_{2\ell - 1}$.
\end{prop}

\begin{proof}
This is an immediate consequence of Proposition~\ref{prop: natural modules}, using the exceptional isogeny $B_\ell \to C_\ell$ which exists in characteristic $2$.
\end{proof}

For most of the remaining results of this section we shall employ the approach of Section~\ref{sect: gen ht fns}. We first consider the remaining cases which occur in infinite families.

\begin{prop}\label{prop: A_ell, 2 omega_1 and omega_2 modules}
Let $G = A_\ell$ for $\ell \in [1, \infty)$ and $\lambda = 2\omega_1$ with $p \geq 3$, or $G = A_\ell$ for $\ell \in [3, \infty)$ and $\lambda = \omega_2$. Then the triple $(G, \lambda, p)$ has generic stabilizer $D_{\frac{1}{2}(\ell + 1)}$ or $C_{\frac{1}{2}(\ell + 1)}$ respectively if $\ell$ is odd, and $B_{\frac{1}{2}\ell}$ or $C_{\frac{1}{2}\ell} U_\ell$ respectively if $\ell$ is even; the associated first quadruple $(G, \lambda, p, 1)$ has generic stabilizer $D_{\frac{1}{2}(\ell + 1)}.\Z_2$ or $C_{\frac{1}{2}(\ell + 1)}$ respectively if $\ell$ is odd, and $B_{\frac{1}{2}\ell}$ or $C_{\frac{1}{2}\ell} T_1 U_\ell$ respectively if $\ell$ is even.
\end{prop}

\begin{proof}
In each case we may view $V$ as a submodule or quotient of the tensor square $V_{nat} \otimes V_{nat}$ of the natural module, and identify $W$ with the symmetric group $S_{\ell + 1}$. Write $\ell_1 = \lceil \frac{1}{2}\ell \rceil$, so that $\ell = 2\ell_1 - 1$ or $2\ell_1$. Note that $Z(G) = \langle z \rangle$ where $z = \prod_{i = 1}^\ell h_{\alpha_i}({\eta_{\ell + 1}}^i)$; as $z$ acts on $V_{nat}$ as multiplication by $\eta_{\ell + 1}$, it acts on $V_{nat} \otimes V_{nat}$ and hence on $V$ as multiplication by ${\eta_{\ell + 1}}^2$, so $G_V = \langle z^{\ell_1} \rangle$ or $\{ 1 \}$ according as $\ell = 2\ell_1 - 1$ or $2\ell_1$.

We take the strictly positive generalized height function on the weight lattice of $G$ whose value at each simple root $\alpha_i$ is $2$; then the generalized height of $\omega_1 = \frac{1}{\ell + 1}(\ell\alpha_1 + (\ell - 1)\alpha_2 + (\ell - 2)\alpha_3 + \cdots + \alpha_\ell)$ is $\ell$, and as $\omega_1$ and $\Phi$ generate the weight lattice we see that the generalized height of any weight is an integer.

First suppose $\lambda = 2\omega_1$ with $p \geq 3$; then we may view $V$ as the symmetric square $S^2(V_{nat})$. For convenience, for $1 \leq i, j \leq \ell + 1$ write $v_{i, j} = v_i \otimes v_j + v_j \otimes v_i$; thus $V = \langle v_{i, j} : 1 \leq i \leq j \leq \ell + 1 \rangle$. The details here differ slightly depending on whether $\ell = 2\ell_1 - 1$ or $2\ell_1$; whenever we give two choices followed by \lq respectively' we are taking the two possibilities in this order.

Since $V_\lambda = \langle v_{1, 1} \rangle$, and by the above the generalized height of $\lambda$ is $2\ell$, we see that if $\mu \in \Lambda(V)$ and $v_{i, j} \in V_\mu$ then the generalized height of $\mu$ is $2(\ell + 2 - i - j)$. Thus $\Lambda(V)_{[0]} = \{ \nu_1, \dots, \nu_{\ell + 1 - \ell_1} \}$, where we write
$$
x_1 = v_{1, \ell + 1}, \quad x_2 = v_{2, \ell}, \quad \cdots, \quad x_{\ell + 1 - \ell_1} = v_{\ell + 1 - \ell_1, \ell_1 + 1},
$$
and for each $i$ we let $\nu_i$ be the weight such that $x_i \in V_{\nu_i}$. Observe that if we take $s = \prod_{i = 1}^\ell h_{\alpha_i}(\kappa_i) \in T$ then $\nu_1(s) = \frac{\kappa_1}{\kappa_\ell}$, and for $i = 2, \dots, \ell + 1 - \ell_1$ we have $\nu_i(s) = \frac{\kappa_i \kappa_{\ell + 2 - i}}{\kappa_{i - 1} \kappa_{\ell + 1 - i}}$; thus $\nu_1 + \cdots + \nu_{\ell_1} = 0$ or $2\nu_1 + \cdots + 2\nu_{\ell_1} + \nu_{\ell_1 + 1} = 0$ respectively, and so $\Lambda(V)_{[0]}$ has ZLC. Set $Y = V_{[0]} = \langle x_1, \dots, x_{\ell + 1 - \ell_1} \rangle$, and
$$
\hat Y = \{ a_1 x_1 + \cdots + a_{\ell + 1 - \ell_1} x_{\ell + 1 - \ell_1} : a_1 \dots a_{\ell + 1 - \ell_1} \neq 0 \},
$$
so that $\hat Y$ is a dense open subset of $Y$. Write
$$
y_0 = x_1 + \cdots + x_{\ell + 1 - \ell_1} \in \hat Y.
$$

Here $W$ acts on $\Lambda(V)$ such that if $w \in W$ and $\mu \in \Lambda(V)$ with $v_{i, j} \in V_\mu$ then $v_{w(i), w(j)} \in V_{w(\mu)}$. The pointwise stabilizer in $W$ of $\Lambda(V)_{[0]}$ is the subgroup $\langle (1 \ \ell + 1), (2 \ \ell), \dots, (\ell_1 \ \ell + 2 - \ell_1) \rangle = \langle w_{\alpha_1 + \cdots + \alpha_\ell}, w_{\alpha_2 + \cdots + \alpha_{\ell - 1}}, \dots, w_{\alpha'} \rangle$, where we set $\alpha' = \alpha_{\ell_1}$ or $\alpha_{\ell_1} + \alpha_{\ell_1 + 1}$ respectively. Note that there are two $W$-orbits on weights here: in the notation of earlier sections, weights of the form $\mu_2$ lie in $W.\lambda$ while those of the form $\mu_1$ lie in $W.\omega_2$. If $\ell = 2\ell_1 - 1$ then all the weights $\nu_i$ are of the form $\mu_1$; if however $\ell = 2\ell_1$ then the last weight $\nu_{\ell + 1 - \ell_1}$ is of the form $\mu_2$ while the remaining $\nu_i$ are of the form $\mu_1$. Since for $1 \leq i < \ell_1$ the element $w_{\alpha_i} w_{\alpha_{\ell + 1 - i}}$ interchanges $\nu_i$ and $\nu_{i + 1}$ while fixing the remaining $\nu_j$, we see that the setwise stabilizer in $W$ of $\Lambda(V)_{[0]}$ is
\begin{eqnarray*}
&   & \langle w_{\alpha_1 + \cdots + \alpha_\ell}, w_{\alpha_2 + \cdots + \alpha_{\ell - 1}}, \dots, w_{\alpha'}, w_{\alpha_1} w_{\alpha_\ell}, \dots, w_{\alpha_{\ell_1 - 1}} w_{\alpha_{\ell + 2 - \ell_1}} \rangle \\
& = & \langle w_{\alpha'}, w_{\alpha_{\ell_1 - 1}} w_{\alpha_{\ell + 2 - \ell_1}}, \dots, w_{\alpha_1} w_{\alpha_\ell} \rangle.
\end{eqnarray*}

Let $A$ be the $D_{\ell_1}$ or $B_{\ell_1}$ subgroup respectively whose first $\ell_1 - 1$ simple root groups are $\{ x_{\alpha_i}(t) x_{\alpha_{\ell + 1 - i}}(-t) : t \in K \}$ for $i = 1, \dots, \ell_1 - 1$, and whose last is $\{ x_{\alpha_{\ell_1 - 1} + \alpha_{\ell_1}}(t) x_{\alpha_{\ell_1} + \alpha_{\ell_1 + 1}}(-t) : t \in K \}$ or $\{ x_{\alpha_{\ell_1}}(t) x_{\alpha_{\ell_1 + 1}}(-2t) x_{\alpha_{\ell_1} + \alpha_{\ell_1 + 1}}(t^2) : t \in K \}$ respectively; then $Z(A) = \langle z^{\ell_1} \rangle$ or $\{ 1 \}$ respectively. If $\ell = 2\ell_1 - 1$ write $n^* = n_{\alpha_{\ell_1}} \prod_{i = 1}^{\ell_1 - 1} h_{\alpha_i}({\eta_{2\ell_1}}^i) \in N$, so that $(n^*)^2 = zs$ where
\begin{eqnarray*}
s &  =  & \prod_{i = 1}^{\ell_1 - 2} h_{\alpha_i}({\eta_{2\ell_1}}^i) h_{\alpha_{2\ell_1 - i}}({\eta_{2\ell_1}}^i) \times \\
  &     & h_{\alpha_{\ell_1 - 1}}({\eta_{4\ell_1}}^{\ell_1 - 1}) h_{\alpha_{\ell_1 + 1}}({\eta_{4\ell_1}}^{\ell_1 - 1}).h_{\alpha_{\ell_1 - 1} + \alpha_{\ell_1}}({\eta_{4\ell_1}}^{\ell_1 - 1}) h_{\alpha_{\ell_1} + \alpha_{\ell_1 + 1}}({\eta_{4\ell_1}}^{\ell_1 - 1}) \\
  & \in & A \cap T,
\end{eqnarray*}
and conjugation by $n^*$ induces a graph automorphism of $A$; then for $i = 1, \dots, \ell_1$ we have $n^*.x_i = \eta_{2\ell_1} x_i$. Set $C = A$, and $C' = Z(G) A \langle n^* \rangle$ or $Z(G) A$ respectively. Clearly we then have $C \leq C_G(y_0)$ and $C' \leq C_G(\langle y_0 \rangle)$; we shall show that in fact $C_G(y_0) = C$ and $C_G(\langle y_0 \rangle) = C'$.

By Lemma~\ref{lem: gen height zero}, if we take $g \in \Tran_G(y_0, Y)$ and set $y' = g.y_0 \in Y$, then we have $g = u_1 n u_2$ with $u_1 \in C_U(y')$, $u_2 \in C_U(y_0)$, and $n \in N_{\Lambda(V)_{[0]}}$ with $n.y_0 = y'$. In particular $G.y_0 \cap Y = N_{\Lambda(V)_{[0]}}.y_0 \cap Y$, and $C_G(y_0) = C_U(y_0) C_{N_{\Lambda(V)_{[0]}}}(y_0) C_U(y_0)$ while $C_G(\langle y_0 \rangle) = C_U(y_0) C_{N_{\Lambda(V)_{[0]}}}(\langle y_0 \rangle) C_U(y_0)$.

First, from the above the elements of $W$ which preserve $\Lambda(V)_{[0]}$ are those corresponding to elements of $C \langle n^* \rangle \cap N$ or $C \cap N$ respectively; so we have $N_{\Lambda(V)_{[0]}}.y_0 = T.y_0 \cup n^*T.y_0$ or $T.y_0$ respectively. Since any element of $T$ may be written as $\prod_{i = 1}^{\ell - \ell_1}h_{\alpha_i}(\kappa_i) t$ where $\kappa_1, \dots, \kappa_{\ell - \ell_1} \in K^*$ and $t \in C \cap T$, by the above if $\ell = 2\ell_1 - 1$ we have
\begin{eqnarray*}
T.y_0    &\! = \!& \left\{ \kappa_1 x_1 + {\ts\frac{\kappa_2}{\kappa_1}} x_2 + \cdots + {\ts\frac{\kappa_{\ell_1 - 1}}{\kappa_{\ell_1 - 2}}} x_{\ell_1 - 1} + {\ts\frac{1}{\kappa_{\ell_1 - 1}}} x_{\ell_1} : \kappa_1, \dots, \kappa_{\ell_1 - 1} \in K^* \right\}, \\
n^*T.y_0 &\! = \!& \left\{ \eta_{2\ell_1}(\kappa_1 x_1 + {\ts\frac{\kappa_2}{\kappa_1}} x_2 + \cdots + {\ts\frac{\kappa_{\ell_1 - 1}}{\kappa_{\ell_1 - 2}}} x_{\ell_1 - 1} + {\ts\frac{1}{\kappa_{\ell_1 - 1}}} x_{\ell_1}) : \kappa_1, \dots, \kappa_{\ell_1 - 1} \in K^* \right\}
\end{eqnarray*}
while if $\ell = 2\ell_1$ we have
$$
T.y_0 = \left\{ \kappa_1 x_1 + {\ts\frac{\kappa_2}{\kappa_1}} x_2 + \cdots + {\ts\frac{\kappa_{\ell_1}}{\kappa_{\ell_1 - 1}}} x_{\ell_1} + {\ts\frac{1}{{\kappa_{\ell_1}}^2}} x_{\ell_1 + 1} : \kappa_1, \dots, \kappa_{\ell_1} \in K^* \right\}.
$$
Hence $C_{N_{\Lambda(V)_{[0]}}}(y_0) = C \cap N$; also $N_{\Lambda(V)_{[0]}}.y_0 \subseteq \hat Y$, and $N_{\Lambda(V)_{[0]}}.y_0 \cap \langle y_0 \rangle = \{ {\eta_{\ell + 1}}^i y_0 : i = 0, \dots, \ell \} = Z(G) \langle n^* \rangle.y_0$ or $Z(G).y_0$ respectively, so $C_{N_{\Lambda(V)_{[0]}}}(\langle y_0 \rangle) = C' \cap N$.

Next, let $\Xi = \{ \alpha_i + \alpha_{i + 1} + \cdots + \alpha_j : i \leq j, \ i + j \leq \ell + 1 \}$, and set $U' = \prod_{\alpha \in \Xi} X_\alpha$; then $U = U'.(C \cap U)$ and $U' \cap (C \cap U) = \{ 1 \}$. We now observe that if $\alpha \in \Xi$ then $\nu_i + \alpha$ is a weight in $V$ for exactly one value of $i$; moreover each weight in $V$ of positive generalized height is of the form $\nu_i + \alpha$ for exactly one such root $\alpha$. Thus if we take $u = \prod x_\alpha(t_\alpha) \in U'$ satisfying $u.y_0 = y_0$, and equate coefficients of weight vectors, taking them in an order compatible with increasing generalized height, we see that for all $\alpha$ we must have $t_\alpha = 0$, so that $u = 1$; so $C_U(y_0) = C \cap U$.

Thus $C_U(y_0), C_{N_{\Lambda(V)_{[0]}}}(y_0) \leq C$ and $C_{N_{\Lambda(V)_{[0]}}}(\langle y_0 \rangle) \leq C'$, so we do indeed have $C_G(y_0) = C$ and $C_G(\langle y_0 \rangle) = C'$. Moreover $G.y_0 \cap Y = \{ b_1 x_1 + \cdots + b_{\ell + 1 - \ell_1} x_{\ell + 1 - \ell_1} : (b_1 \dots b_{\ell - \ell_1})^2 {b_{\ell + 1 - \ell_1}}^{2/(2, \ell)} = 1 \}$.

Take $y = a_1 x_1 + \cdots + a_{\ell + 1 - \ell_1} x_{\ell + 1 - \ell_1} \in \hat Y$. By the above, if we choose $\kappa \in K^*$ satisfying $\kappa^{\ell_1} = a_1 \dots a_{\ell_1}$ or $\kappa^{2\ell_1 + 1} = (a_1 \dots a_{\ell_1})^2 a_{\ell_1 + 1}$ respectively, then $\kappa^{-1}y \in T.y_0$, so there exists $h \in T$ with $h.y_0 = \kappa^{-1} y$; so $C_G(y) = C_G(\kappa^{-1} y) = C_G(h.y_0) = {}^h C$ and likewise $C_G(\langle y \rangle) = {}^h C'$. Moreover, we see that $G.y \cap Y = G.h.\kappa y_0 \cap Y = \kappa(G.y_0 \cap Y) = \{ b_1 x_1 + \cdots + b_{\ell + 1 - \ell_1} x_{\ell + 1 - \ell_1} : (b_1 \dots b_{\ell - \ell_1})^2 {b_{\ell + 1 - \ell_1}}^{2/(2, \ell)} = (a_1 \dots a_{\ell - \ell_1})^2 {a_{\ell + 1 - \ell_1}}^{2/(2, \ell)} \}$. Since $\dim C = 2{\ell_1}^2 - \ell_1$ or $2{\ell_1}^2 + \ell_1$ respectively, for all $y \in \hat Y$ we have $\dim(\overline{G.y}) = \dim G - \dim C = (4{\ell_1}^2 - 1) - (2{\ell_1}^2 - \ell_1) = 2{\ell_1}^2 + \ell_1 - 1$ or $(4{\ell_1}^2 + 4{\ell_1}) - (2{\ell_1}^2 + \ell_1) = 2{\ell_1}^2 + 3{\ell_1}$ respectively, while $\dim(\overline{G.y \cap Y}) = \ell_1 - 1$ or $\ell_1$ respectively; therefore if $\ell = 2\ell_1 - 1$ then
$$
\dim V - \dim(\overline{G.y}) = (2{\ell_1}^2 + \ell_1) - (2{\ell_1}^2 + \ell_1 - 1) = 1
$$
and
$$
\dim Y - \dim(\overline{G.y \cap Y}) = \ell_1 - (\ell_1 - 1) = 1,
$$
while if $\ell = 2\ell_1$ then
$$
\dim V - \dim(\overline{G.y}) = (2{\ell_1}^2 + 3\ell_1 + 1) - (2{\ell_1}^2 + 3\ell_1) = 1
$$
and
$$
\dim Y - \dim(\overline{G.y \cap Y}) = (\ell_1 + 1) - \ell_1 = 1.
$$
Hence $y$ is $Y$-exact. Thus the conditions of Lemma~\ref{lem: generic stabilizer from exact subset} hold; so the triple $(G, \lambda, p)$ has generic stabilizer $C/G_V \ \cong \ D_{\ell_1}$ or $B_{\ell_1}$ respectively, while the quadruple $(G, \lambda, p, 1)$ has generic stabilizer $C'/Z(G) \cong D_{\ell_1}.\Z_2$ or $B_{\ell_1}$ respectively, where the $D_{\ell_1}$ or $B_{\ell_1}$ is of adjoint type.

Now suppose instead $\lambda = \omega_2$; then we may view $V$ as the exterior square $\bigwedge^2(V_{nat})$. For convenience, for $1 \leq i, j \leq \ell + 1$ with $i \neq j$ write $\bar v_{i, j} = v_i \wedge v_j$; thus $V = \langle \bar v_{i, j} : 1 \leq i < j \leq \ell + 1 \rangle$.

Write
$$
x_1 = \bar v_{1, 2\ell_1}, \quad x_2 = \bar v_{2, 2\ell_1 - 1}, \quad \cdots, \quad x_{\ell_1} = \bar v_{\ell_1, \ell_1 + 1},
$$
and for each $i$ let $\nu_i$ be the weight such that $x_i \in V_{\nu_i}$. Observe that if we take $s = \prod_{i = 1}^\ell h_{\alpha_i}(\kappa_i) \in T$ then $\nu_1(s) = \frac{\kappa_1}{\kappa_{2\ell_1 - 1}}$ or $\frac{\kappa_1 \kappa_{2\ell_1}}{\kappa_{2\ell_1 - 1}}$ according as $\ell = 2\ell_1 - 1$ or $2\ell_1$, and for $i = 2, \dots, \ell_1$ we have $\nu_i(s) = \frac{\kappa_i \kappa_{2\ell_1 + 1 - i}}{\kappa_{i - 1} \kappa_{2\ell_1 - i}}$. Set $Y = \langle x_1, \dots, x_{\ell_1} \rangle$, and
$$
\hat Y = \{ a_1 x_1 + \cdots + a_{\ell_1} x_{\ell_1} : a_1 \dots a_{\ell_1} \neq 0 \},
$$
so that $\hat Y$ is a dense open subset of $Y$. Write
$$
y_0 = x_1 + \cdots + x_{\ell_1} \in \hat Y.
$$

Here $W$ acts on $\Lambda(V)$ such that if $w \in W$ and $\mu \in \Lambda(V)$ with $\bar v_{i, j} \in V_\mu$ then $\bar v_{w(i), w(j)} \in V_{w(\mu)}$. The pointwise stabilizer in $W$ of $\{ \nu_1, \dots, \nu_{\ell_1} \}$ is the subgroup $\langle (1 \ 2\ell_1), (2 \ 2\ell_1 - 1), \dots, (\ell_1 \ \ell_1 + 1) \rangle = \langle w_{\alpha_1 + \cdots + \alpha_{2\ell_1 - 1}}, w_{\alpha_2 + \cdots + \alpha_{2\ell_1 - 2}}, \dots, w_{\alpha_{\ell_1}} \rangle$. Since for $1 \leq i < \ell_1$ the element $w_{\alpha_i} w_{\alpha_{2\ell_1 - i}}$ interchanges $\nu_i$ and $\nu_{i + 1}$ while fixing the remaining $\nu_j$, we see that the setwise stabilizer in $W$ of $\{ \nu_1, \dots, \nu_{\ell_1} \}$ is
\begin{eqnarray*}
   &   & \langle w_{\alpha_1 + \cdots + \alpha_{2\ell_1 - 1}}, w_{\alpha_2 + \cdots + \alpha_{2\ell_1 - 2}}, \dots, w_{\alpha_{\ell_1}}, w_{\alpha_1} w_{\alpha_{2\ell_1 - 1}}, \dots, w_{\alpha_{\ell_1 - 1}} w_{\alpha_{\ell_1 + 1}} \rangle \\
   & = & \langle w_{\alpha_{\ell_1}}, w_{\alpha_{\ell_1 - 1}} w_{\alpha_{\ell_1 + 1}}, \dots, w_{\alpha_1} w_{\alpha_{2\ell_1 - 1}} \rangle.
\end{eqnarray*}

Let $A$ be the $C_{\ell_1}$ subgroup with short simple root groups $\{ x_{\alpha_i}(t) x_{\alpha_{2\ell_1 - i}}(-t) : t \in K \}$ for $i = 1, \dots, \ell_1 - 1$, and long simple root group $X_{\alpha_{\ell_1}}$; then $Z(A) = \langle z^{\ell_1} \rangle$.

First assume $\ell = 2\ell_1 - 1$, and set $C = A$ and $C' = Z(G) A$. Clearly we then have $C \leq C_G(y_0)$ and $C' \leq C_G(\langle y_0 \rangle)$; we shall show that in fact $C_G(y_0) = C$ and $C_G(\langle y_0 \rangle) = C'$.

Since $V_\lambda = \langle \bar v_{1, 2} \rangle$, and by the above the generalized height of $\lambda$ is $2\ell - 2$, we see that if $\mu \in \Lambda(V)$ and $\bar v_{i, j} \in V_\mu$, the generalized height of $\mu$ is $2(\ell + 2 - i - j)$. Thus $\Lambda(V)_{[0]} = \{ \nu_1, \dots, \nu_{\ell_1} \}$, and so $Y = V_{[0]}$; from the above we have $\nu_1 + \cdots + \nu_{\ell_1} = 0$, so that $\Lambda(V)_{[0]}$ has ZLC. By Lemma~\ref{lem: gen height zero}, if we take $g \in \Tran_G(y_0, Y)$ and set $y' = g.y_0 \in Y$, then we have $g = u_1 n u_2$ with $u_1 \in C_U(y')$, $u_2 \in C_U(y_0)$, and $n \in N_{\Lambda(V)_{[0]}}$ with $n.y_0 = y'$. In particular $G.y_0 \cap Y = N_{\Lambda(V)_{[0]}}.y_0 \cap Y$, and $C_G(y_0) = C_U(y_0) C_{N_{\Lambda(V)_{[0]}}}(y_0) C_U(y_0)$ while $C_G(\langle y_0 \rangle) = C_U(y_0) C_{N_{\Lambda(V)_{[0]}}}(\langle y_0 \rangle) C_U(y_0)$.

First, from the above the elements of $W$ which preserve $\Lambda(V)_{[0]}$ are those corresponding to elements of $C \cap N$; so we have $N_{\Lambda(V)_{[0]}}.y_0 = T.y_0$. Since any element of $T$ may be written as $\prod_{i = 1}^{\ell_1 - 1}h_{\beta_i}(\kappa_i) t$ where $\kappa_1, \dots, \kappa_{\ell_1 - 1} \in K^*$ and $t \in C \cap T$, by the above we have
$$
T.y_0 = \left\{ \kappa_1 x_1 + {\ts\frac{\kappa_2}{\kappa_1}} x_2 + \cdots + {\ts\frac{\kappa_{\ell_1 - 1}}{\kappa_{\ell_1 - 2}}} x_{\ell_1 - 1} + {\ts\frac{1}{\kappa_{\ell_1 - 1}}} x_{\ell_1} : \kappa_1, \dots, \kappa_{\ell_1 - 1} \in K^* \right\}.
$$
Hence $C_{N_{\Lambda(V)_{[0]}}}(y_0) = C \cap N$; also $N_{\Lambda(V)_{[0]}}.y_0 \subseteq \hat Y$, and $N_{\Lambda(V)_{[0]}}.y_0 \cap \langle y_0 \rangle = \{ {\eta_{\ell_1}}^i y_0 : i = 0, \dots, \ell_1 - 1 \} = Z(G).y_0$, so $C_{N_{\Lambda(V)_{[0]}}}(\langle y_0 \rangle) = C' \cap N$.

Next, let $\Xi = \{ \alpha_i + \alpha_{i + 1} + \cdots + \alpha_j : i \leq j, \ i + j \leq \ell \}$, and set $U' = \prod_{\alpha \in \Xi} X_\alpha$; then $U = U'.(C \cap U)$ and $U' \cap (C \cap U) = \{ 1 \}$. We now observe that if $\alpha \in \Xi$ then $\nu_i + \alpha$ is a weight in $V$ for exactly one value of $i$; moreover each weight in $V$ of positive generalized height is of the form $\nu_i + \alpha$ for exactly one such root $\alpha$. Thus if we take $u = \prod x_\alpha(t_\alpha) \in U'$ satisfying $u.y_0 = y_0$, and equate coefficients of weight vectors, taking them in an order compatible with increasing generalized height, we see that for all $\alpha$ we must have $t_\alpha = 0$, so that $u = 1$; so $C_U(y_0) = C \cap U$.

Thus $C_U(y_0), C_{N_{\Lambda(V)_{[0]}}}(y_0) \leq C$ and $C_{N_{\Lambda(V)_{[0]}}}(\langle y_0 \rangle) \leq C'$, so we do indeed have $C_G(y_0) = C$ and $C_G(\langle y_0 \rangle) = C'$. Moreover $G.y_0 \cap Y = \{ b_1 x_1 + \cdots + b_{\ell_1} x_{\ell_1} : b_1 \dots b_{\ell_1} = 1 \}$.

Take $y = a_1 x_1 + \cdots + a_{\ell_1} x_{\ell_1} \in \hat Y$. By the above, if we choose $\kappa \in K^*$ satisfying $\kappa^{\ell_1} = a_1 \dots a_{\ell_1}$, then $\kappa^{-1}y \in T.y_0$, so there exists $h \in T$ with $h.y_0 = \kappa^{-1} y$; so $C_G(y) = C_G(\kappa^{-1} y) = C_G(h.y_0) = {}^h C$ and likewise $C_G(\langle y \rangle) = {}^h C'$. Moreover, we see that $G.y \cap Y = G.h.\kappa y_0 \cap Y = \kappa(G.y_0 \cap Y) = \{ b_1 x_1 + \cdots + b_{\ell_1} x_{\ell_1} : b_1 \dots b_{\ell_1} = a_1 \dots a_{\ell_1} \}$. Since $\dim C = 2{\ell_1}^2 + \ell_1$, for all $y \in \hat Y$ we have $\dim(\overline{G.y}) = \dim G - \dim C = (4{\ell_1}^2 - 1) - (2{\ell_1}^2 + \ell_1) = 2{\ell_1}^2 - \ell_1 - 1$, while $\dim(\overline{G.y \cap Y}) = \ell_1 - 1$; therefore
$$
\dim V - \dim(\overline{G.y}) = (2{\ell_1}^2 - \ell_1) - (2{\ell_1}^2 - \ell_1 - 1) = 1
$$
and
$$
\dim Y - \dim(\overline{G.y \cap Y}) = \ell_1 - (\ell_1 - 1) = 1.
$$
Hence $y$ is $Y$-exact. Thus the conditions of Lemma~\ref{lem: generic stabilizer from exact subset} hold; so the triple $(G, \lambda, p)$ has generic stabilizer $C/G_V \cong C_{\ell_1}$, while the quadruple $(G, \lambda, p, 1)$ has generic stabilizer $C'/Z(G) \cong C_{\ell_1}$, where the $C_{\ell_1}$ is of adjoint type.

Now assume $\ell = 2\ell_1$. Let $P = QL$ be the standard $A_{\ell - 1}$ parabolic subgroup of $G$ corresponding to the last simple root, with Levi subgroup $L = \langle T, X_\alpha : \alpha = \sum m_i \alpha_i, \ m_\ell = 0 \rangle$ and $\ell$-dimensional unipotent radical $Q = \langle X_\alpha : \alpha = \sum m_i \alpha_i, \ m_\ell = 1 \rangle$; then each element of $Q$ fixes each element of $Y$. Write $P^- = Q^- L$ for the opposite parabolic subgroup, so that $Q^- = \langle X_\alpha : \alpha = \sum m_i \alpha_i, \ m_\ell = -1 \rangle$. Note that $A \leq L$; set $C = QA$ and $C' = Z(L)QA$. Clearly we then have $C \leq C_G(y_0)$ and $C' \leq C_G(\langle y_0 \rangle)$; we shall show that in fact $C_G(y_0) = C$ and $C_G(\langle y_0 \rangle) = C'$.

Suppose first that $g \in G$ satisfies $g.y_0 \in Y$; write $g.y_0 = y$. Using Lemma~\ref{lem: parabolic factorization} we may write $g = q_1xq_2q_3$, where $q_1, q_3 \in Q$, $q_2 \in Q^-$ and $x \in L$; then we have $xq_2q_3.y_0 = {q_1}^{-1}.y$, whence $xq_2.y_0 = y$. Now the root element $x_{-(\alpha_i + \cdots + \alpha_\ell)}(t)$ sends $v_i$ to $v_i + tv_{\ell + 1}$ and fixes all other basis vectors of $V_{nat}$. Thus if $q_2 \neq 1$ then $q_2.y_0$ has at least one term $\bar v_{i, \ell + 1}$, as therefore does $xq_2.y_0$, contrary to $y \in Y$; so we must have $q_2 = 1$, and hence $g = q_1xq_3$ and $x.y_0 = y$. Now write $x = tx'$ where $x' \in L'$ and $t = \prod_{i = 1}^\ell h_{\alpha_i}(\kappa^i) \in Z(L)$ for some $\kappa \in K^*$; then $t.y_0 = \kappa^2 y_0$ and so we have $x'.y_0 = \kappa^{-2}y$.

By the previous case we see that $x' = u_1nu_2$ where $u_1 \in C_{L' \cap U}(\kappa^{-2} y) = C_{L' \cap U}(y)$, $u_2 \in C_{L' \cap U}(y_0)$ and $n \in L' \cap N$ with $n.y_0 = \kappa^{-2} y$. From the above the elements of $W$ which preserve $\{ \nu_1, \dots, \nu_{\ell_1} \}$ are those corresponding to elements of $A \cap N$, so $n \in  \{ \prod_{i = 1}^{\ell_1 - 1} h_{\alpha_i}(\kappa_i) : \kappa_1, \dots, \kappa_{\ell_1 - 1} \in K^* \} (A \cap N)$; then we may write $n = sn'$ where $s = \prod_{i = 1}^{\ell_1 - 1} h_{\alpha_i}(\kappa_i)$ for some $\kappa_1, \dots, \kappa_{\ell_1 - 1} \in K^*$ and $n' \in A \cap N$, and so $\kappa^{-2}y = sn'.y_0 = s.y_0$. From the previous case again we have $C_{L' \cap U}(y_0) = A \cap U$, and so $C_{L' \cap U}(y) = C_{L' \cap U}(\kappa^{-2} y) = C_{L' \cap U}(s.y_0) = {}^s (A \cap U)$; therefore $x' \in {}^s (A \cap U).s(A \cap N).(A \cap U) = s(A \cap U)(A \cap N)(A \cap U)$, and so we have $x' = sa$ for some $a \in A$.

Hence $x = tx' = tsa$; so $g = q_1tsaq_3 = ts.({q_1}^{ts})({}^a{q_3}).a \in TQA = TC$. In particular, if $g.y_0 \in \langle y_0 \rangle$ we must have $s.y_0 \in \langle y_0 \rangle$, so as $s.y_0 = \kappa_1 x_1 + \frac{\kappa_2}{\kappa_1} x_2 + \cdots + \frac{\kappa_{\ell_1 - 1}}{\kappa_{\ell_1 - 2}} x_{\ell_1 - 1} + \frac{1}{\kappa_{\ell_1 - 1}} x_{\ell_1}$ we must have $\kappa_1 = \frac{\kappa_2}{\kappa_1} = \cdots = \frac{\kappa_{\ell_1 - 1}}{\kappa_{\ell_1 - 2}} = \frac{1}{\kappa_{\ell_1 - 1}}$; this implies ${\kappa_1}^{\ell_1} = 1$ and $\kappa_i = {\kappa_1}^i$ for $1 < i < \ell_1$, so that if we take $\kappa' \in K^*$ satisfying ${\kappa'}^2 = \kappa_1$ then $s = s_1s_2$ where $s_1 = \prod_{i = 1}^\ell h_{\alpha_i}({\kappa'}^i) \in Z(L)$ and $s_2 = \left( \prod_{i = 1}^{\ell_1 - 1} h_{\alpha_i}({\kappa'}^i) h_{\alpha_{2\ell_1 - i}}({\kappa'}^i) \right) h_{\alpha_{\ell_1}}({\kappa'}^{\ell_1}) \in A \cap T$. Therefore $C_G(\langle y_0 \rangle) = C'$; and as $C_{Z(L)}(y_0) = \langle \prod_{i = 1}^\ell h_{\alpha_i}((-1)^i) \rangle < A$ we also have $C_G(y_0) = C$. Moreover we see that $G.y_0 \cap Y = TC.y_0 = T.y_0 = \hat Y$, since given $y = a_1 x_1 + \cdots + a_{\ell_1} x_{\ell_1} \in \hat Y$ we have $y = h.y_0$ for $h = h_{\alpha_{\ell_1 + 1}}(a_{\ell_1}) h_{\alpha_{\ell_1 + 2}}(a_{\ell_1 - 1} a_{\ell_1}) \dots h_{\alpha_{2\ell_1}}(a_1 \dots a_{\ell_1})$.

Take $y \in \hat Y$. By the above, there exists $h \in T$ with $h.y_0 = y$; so $C_G(y) = C_G(h.y_0) = {}^h C$ and likewise $C_G(\langle y \rangle) = {}^h C'$. Since $\dim C = 2{\ell_1}^2 + 3\ell_1$, for all $y \in \hat Y$ we have $\dim(\overline{G.y}) = \dim G - \dim C = (4{\ell_1}^2 + 4\ell_1) - (2{\ell_1}^2 + 3\ell_1) = 2{\ell_1}^2 + \ell_1$, while $\dim(\overline{G.y \cap Y}) = \ell_1$; therefore
$$
\dim V - \dim(\overline{G.y}) = (2{\ell_1}^2 + \ell_1) - (2{\ell_1}^2 + \ell_1) = 0
$$
and
$$
\dim Y - \dim(\overline{G.y \cap Y}) = \ell_1 - \ell_1 = 0.
$$
Hence $y$ is $Y$-exact. Thus the conditions of Lemma~\ref{lem: generic stabilizer from exact subset} hold; so the triple $(G, \lambda, p)$ has generic stabilizer $C/G_V \cong C_{\ell_1} U_\ell$, while the quadruple $(G, \lambda, p, 1)$ has generic stabilizer $C'/Z(G) \cong C_{\ell_1} T_1 U_\ell$, where the $C_{\ell_1}$ is of simply connected type.
\end{proof}

\begin{prop}\label{prop: C_ell, omega_2 module}
Let $G = C_\ell$ for $\ell \in [3, \infty)$ and $\lambda = \omega_2$. Then if $\ell = p = 3$ the triple $(G, \lambda, p)$ has generic stabilizer ${C_1}^3.\Z_3$ and the associated first quadruple $(G, \lambda, p, 1)$ has generic stabilizer ${C_1}^3.S_3$; if $\ell = 4$ and $p = 2$ both the triple $(G, \lambda, p)$ and the associated first quadruple $(G, \lambda, p, 1)$ have generic stabilizer ${C_1}^4.{\Z_2}^2$; in all other cases both the triple $(G, \lambda, p)$ and the associated first quadruple $(G, \lambda, p, 1)$ have generic stabilizer ${C_1}^\ell$.
\end{prop}

\begin{proof}
Inside the exterior square $\bigwedge^2(V_{nat})$ of the natural module we have the submodules $X_1 = \{ \sum_{i < j} \rho_{ij} e_i \wedge e_j + \sum_{i < j} \sigma_{ij} f_i \wedge f_j + \sum_{i, j} \tau_{ij} e_i \wedge f_j : \sum_i \tau_{ii} = 0 \}$ and $X_2 = \langle \sum_{i = 1}^\ell e_i \wedge f_i \rangle$, with the latter being trivial. If $p$ is coprime to $\ell$ then $\bigwedge^2(V_{nat}) = X_1 \oplus X_2$, and $V = X_1$; if however $p$ divides $\ell$ then $X_2 < X_1$, and $V = X_1/X_2$. Thus in all cases $V = X_1/(X_1 \cap X_2)$, where $X_1 \cap X_2$ is either zero or the trivial $G$-module. Moreover $Z(G) = \langle z \rangle$ where $z = \prod_{i = 1}^{\lceil \frac{\ell}{2} \rceil} h_{\alpha_{2i - 1}}(-1)$ which acts as negation on $V_{nat}$ and therefore trivially on $V$, so $G_V = Z(G)$.

We take the strictly positive generalized height function on the weight lattice of $G$ whose value at $\alpha_\ell$ is $2$ and at each other simple root $\alpha_i$ is $1$; then the generalized height of $\lambda = \alpha_1 + 2\alpha_2 + \cdots + 2\alpha_{\ell - 1} + \alpha_\ell$ is $2\ell - 1$, and as $\frac{1}{2}\alpha_\ell$ and $\Phi$ generate the weight lattice it follows that the generalized height of any weight is an integer. Since $V_\lambda = \langle e_1 \wedge e_2 \rangle$, we see that if $\mu \in \Lambda(V)$ and $v \in V_\mu$ where $v = e_i \wedge e_j$, $e_i \wedge f_j$ or $f_i \wedge f_j$, then the generalized height of $\mu$ is $2\ell + 2 - (i + j)$, $j - i$ or $i + j - (2\ell + 2)$ respectively. Thus $\Lambda(V)_{[0]} = \{ 0 \}$; so trivially $\Lambda(V)_{[0]}$ has ZLC, and $N_{\Lambda(V)_{[0]}} = N$. For $i = 1, \dots, \ell$ write $x_i = e_i \wedge f_i$. Set
$$
Y = V_{[0]} = \left\{ {\ts\sum} a_i x_i + (X_1 \cap X_2) : {\ts\sum} a_i = 0 \right\}.
$$
If $\ell = 3$ set
$$
\hat Y = \left\{ {\ts\sum} a_i x_i + (X_1 \cap X_2) : {\ts\sum} a_i = 0, \ a_1 \dots a_\ell \neq 0, \ (\ts{\frac{a_i}{a_j}})^3 \neq 1 \hbox{ if } i \neq j \right\};
$$
if $\ell = 4$ and $p = 2$ set
\begin{eqnarray*}
\hat Y & = & \left\{ {\ts\sum} a_i x_i + (X_1 \cap X_2) : {\ts\sum} a_i = 0, \ a_1 \dots a_\ell \neq 0, \ a_i \neq a_j \hbox{ if } i \neq j, \right. \\
       &   & \left. \phantom{\left\{ {\ts\sum} a_i x_i + (X_1 \cap X_2) : \right.} \ (\ts{\frac{a_i - a_{i'}}{a_{i'} - a_{i''}}})^{12} \neq 1 \hbox{ if } i, i', i'' \hbox{ distinct} \right\};
\end{eqnarray*}
otherwise set
\begin{eqnarray*}
\hat Y & = & \left\{ {\ts\sum} a_i x_i + (X_1 \cap X_2) : {\ts\sum} a_i = 0, \ a_1 \dots a_\ell \neq 0, \ a_i \neq a_j \hbox{ if } i \neq j, \right. \\
       &   & \left. \phantom{\left\{ {\ts\sum} a_i x_i + (X_1 \cap X_2) : \right.} \ (\ts{\frac{a_i - a_{i'}}{a_{i'} - a_{i''}}})^{\ell!} \neq 1 \hbox{ if } i, i', i'' \hbox{ distinct}, \right. \\
       &   & \left. \phantom{\left\{ {\ts\sum} a_i x_i + (X_1 \cap X_2) : \right.} \ a_i + a_j \neq a_{i'} + a_{j'} \hbox{ if } i, j, i', j' \hbox{ distinct} \right\}.
\end{eqnarray*}
Thus $\hat Y$ is a dense open subset of $Y$. Take
$$
y = {\ts\sum} a_i x_i + (X_1 \cap X_2) \in \hat Y.
$$

Let $A$ be the ${C_1}^\ell$ subgroup with simple roots $2\ve_1$, \dots, $2\ve_\ell$; as $Z(A) = \langle h_\alpha(-1) : \alpha \in \Phi_l \cap \Phi^+ \rangle$, we see that $A$ is of simply connected type. If $\ell = p = 3$, write $n^* = n_{\alpha_1}$ and $n^{**} = n_{\alpha_2}$, and set $C = A \langle n^* n^{**} \rangle \cong {C_1}^3.\Z_3$ and $C' = A \langle n^*, n^{**} \rangle \cong {C_1}^3.S_3$; if $\ell = 4$ and $p = 2$, write $n^\dagger = n_{\alpha_1} n_{\alpha_3}$ and $n^\ddagger = n_{\alpha_1 + \alpha_2} n_{\alpha_2 + \alpha_3}$, and set $C = C' = A \langle n^\dagger, n^\ddagger \rangle \cong {C_1}^4.{\Z_2}^2$; otherwise set $C = C' = A$. Note that if $\ell = p = 3$ then $\dim Y = 1$, and both $n^*$ and $n^{**}$ act on $Y$ as negation since for example $n^*.(\sum a_i x_i) = a_1 x_2 + a_2 x_1 + a_3 x_3 = -(\sum a_i x_i) - a_3 \sum x_i$ as $\sum a_i = 0$; if instead $\ell = 4$ and $p = 2$ then both $n^\dagger$ and $n^\ddagger$ act on $Y$ as the identity since for example $n^\dagger.(\sum a_i x_i) = a_1 x_2 + a_2 x_1 + a_3 x_4 + a_4 x_3 = (\sum a_i x_i) + (a_1 + a_2)\sum x_i$ as $\sum a_i = 0$. Thus we have $C \leq C_G(y)$ and $C' \leq C_G(\langle y \rangle)$; we shall show that in fact $C_G(y) = C$ and $C_G(\langle y \rangle) = C'$.

By Lemma~\ref{lem: gen height zero}, if we take $g \in \Tran_G(y, Y)$ and set $y' = g.y \in Y$, then we have $g = u_1 n u_2$ with $u_1 \in C_U(y')$, $u_2 \in C_U(y)$, and $n \in N$ with $n.y = y'$. In particular $G.y \cap Y = N.y$, and $C_G(y) = C_U(y) C_N(y) C_U(y)$ while $C_G(\langle y \rangle) = C_U(y) C_N(\langle y \rangle) C_U(y)$.

First, we note that any element of $N$ may be written as $n's$, where $n'$ is a product of elements $n_\alpha$ for various roots $\alpha$, and $s \in T$; since $s$ stabilizes $y$, and each element $n_\alpha$ can only permute the individual vectors $x_i$, we see that
$$
N.y = \left\{ {\ts\sum} a_{\pi(i)} x_i + (X_1 \cap X_2) : \pi \in S_\ell \right\}.
$$
Moreover in the case where $\ell = 3$ and $p \neq 3$, suppose $\pi \in S_3$ satisfies $\sum a_{\pi(i)} x_i = \kappa \sum a_i x_i$ for some $\kappa \in K^*$. If $\pi$ is a transposition, say $(1 \ 2)$, then equating coefficients of $x_3$ and $x_1$ gives $\kappa = 1$ and then $a_1 = a_2$; if instead $\pi$ is a $3$-cycle, say $(1 \ 2 \ 3)$, then equating coefficients gives $\kappa = \frac{a_2}{a_1} = \frac{a_3}{a_2} = \frac{a_1}{a_3}$, so $\kappa^3 = 1$ and $(\frac{a_2}{a_1})^3 = 1$. The definition of $\hat Y$ rules out both possibilities, so we must have $\pi = 1$. Now assume instead $\ell \geq 4$, and suppose $\pi \in S_\ell$ satisfies $\sum a_{\pi(i)} x_i = \kappa \sum a_i x_i + \nu \sum x_i$ for some $\kappa \in K^*$ and $\nu \in K$. If $\pi$ contains an $r$-cycle for some $r \geq 3$, say $(1 \ 2 \ \dots \ r)$, then equating coefficients of $x_1, x_2, \dots, x_r$ gives $\nu = a_2 - \kappa a_1 = a_3 - \kappa a_2 = \cdots = a_r - \kappa a_{r - 1} = a_1 - \kappa a_r$, so $\kappa = \frac{a_3 - a_2}{a_2 - a_1} = \frac{a_4 - a_3}{a_3 - a_2} = \cdots = \frac{a_1 - a_r}{a_r - a_{r - 1}} = \frac{a_2 - a_1}{a_1 - a_r}$ and hence $(\frac{a_3 - a_2}{a_2 - a_1})^r = \kappa^r = 1$; if instead $\pi$ has order $2$ and contains at least two transpositions, say $(1 \ 2)$ and $(3 \ 4)$, then if $p \neq 2$ equating coefficients of $x_1, \dots, x_4$ gives $\nu = a_2 - \kappa a_1 = a_1 - \kappa a_2 = a_4 - \kappa a_3 = a_3 - \kappa a_4$, so $(1 + \kappa)(a_1 - a_2) = 0 = (1 + \kappa)(a_3 - a_4)$, whence either $a_1 = a_2$ and $a_3 = a_4$, or $\kappa = -1$ and $a_1 + a_2 = a_3 + a_4$; finally if $\pi$ is a transposition, say $(1 \ 2)$, then equating coefficients of $x_1, \dots, x_4$ gives $\nu = a_2 - \kappa a_1 = a_1 - \kappa a_2 = (1 - \kappa)a_3 = (1 - \kappa)a_4$, whence either $a_3 = a_4$, or $\kappa = 1$ and $a_1 = a_2$. Again the definition of $\hat Y$ rules out each of these possibilities, so we must have $\pi = 1$. Thus the only elements of $N$ which stabilize $y$ or $\langle y \rangle$ are those in $C$ or $C'$ respectively. Hence $C_N(y) = C \cap N$ and $C_N(\langle y \rangle) = C' \cap N$.

Next, let $\Xi = \Phi_s \cap \Phi^+$, and set $U' = \prod_{\alpha \in \Xi} X_\alpha$; then $U = U'.(C \cap U)$ and $U' \cap (C \cap U) = \{ 1 \}$. Observe that for $i < j$ the short root elements $x_{\ve_i - \ve_j}(t)$ and $x_{\ve_i + \ve_j}(t)$ send $y$ to $y + t(a_j - a_i)e_i \wedge f_j$ and $y + t(a_i - a_j)e_i \wedge e_j$ respectively. Thus if we take $u = \prod x_\alpha(t_\alpha) \in U'$ satisfying $u.y = y$, and equate coefficients of weight vectors, taking them in an order compatible with increasing generalized height, we see that for all $\alpha$ we must have $t_\alpha = 0$, so that $u = 1$; so $C_U(y) = C \cap U$.

Thus $C_U(y), C_N(y) \leq C$ and $C_N(\langle y \rangle) \leq C'$, so we do indeed have $C_G(y) = C$ and $C_G(\langle y \rangle) = C'$.

Since $\dim C = 3\ell$, we have $\dim(\overline{G.y}) = \dim G - \dim C = 2\ell^2 - 2\ell$, while $\dim(\overline{G.y \cap Y}) = 0$ because any $N$-orbit on $Y = V_0$ is finite; therefore
$$
\dim V - \dim(\overline{G.y}) = (2\ell^2 - \ell - 1 - \z_{p, \ell}) - (2\ell^2 - 2\ell) = \ell - 1 - \z_{p, \ell}
$$
and
$$
\dim Y - \dim(\overline{G.y \cap Y}) = (\ell - 1 - \z_{p, \ell}) - 0 = \ell - 1 - \z_{p, \ell}.
$$
Hence $y$ is $Y$-exact. Thus the conditions of Lemma~\ref{lem: generic stabilizer from exact subset} hold; so if $\ell = p = 3$ the triple $(G, \lambda, p)$ has generic stabilizer $C/G_V \cong {C_1}^3.\Z_3$ and the quadruple $(G, \lambda, p, 1)$ has generic stabilizer $C'/Z(G) \cong {C_1}^3.S_3$, if $\ell = 4$ and $p = 2$ the triple $(G, \lambda, p)$ has generic stabilizer $C/G_V \cong {C_1}^4.{\Z_2}^2$ and the associated first quadruple $(G, \lambda, p, 1)$ has generic stabilizer $C'/Z(G) \cong {C_1}^4.{\Z_2}^2$, while in all other cases the triple $(G, \lambda, p)$ has generic stabilizer $C/G_V \cong {C_1}^\ell$ and the quadruple $(G, \lambda, p, 1)$ has generic stabilizer $C'/Z(G) \cong {C_1}^\ell$, where each ${C_1}^\ell$ is a central product.
\end{proof}

\begin{prop}\label{prop: B_ell, omega_2 module, p = 2}
Let $G = B_\ell$ for $\ell \in [3, \infty)$ and $\lambda = \omega_2$ with $p = 2$. Then if $\ell = 4$ the triple $(G, \lambda, p)$ and the associated first quadruple $(G, \lambda, p, 1)$ both have generic stabilizer ${B_1}^4.{\Z_2}^2$, while otherwise the triple $(G, \lambda, p)$ and the associated first quadruple $(G, \lambda, p, 1)$ both have generic stabilizer ${B_1}^\ell$.
\end{prop}

\begin{proof}
This is an immediate consequence of Proposition~\ref{prop: C_ell, omega_2 module}, using the exceptional isogeny $B_\ell \to C_\ell$ which exists in characteristic $2$.
\end{proof}

This concludes the treatment of the cases occurring in infinite families. Although the remaining cases must be treated individually, it will be seen that there are connections between some of them which significantly reduce the amount of work involved.

\begin{prop}\label{prop: E_7, omega_7, D_6, omega_6, B_5, omega_5, A_5, omega_3, C_3, omega_3 modules}
Let $G = E_7$ and $\lambda = \omega_7$, or $G = D_6$ and $\lambda = \omega_6$, or $G = B_5$ and $\lambda = \omega_5$, or $G = A_5$ and $\lambda = \omega_3$, or $G = C_3$ and $\lambda = \omega_3$ with $p \geq 3$. Then the triple $(G, \lambda, p)$ has generic stabilizer $E_6.\Z_{(p, 2)}$, or $A_5.\Z_{(p, 2)}$, or $A_4.\Z_{(p, 2)}$, or ${A_2}^2.\Z_{(p, 2)}$, or $\tilde A_2$, respectively; the associated first quadruple $(G, \lambda, p, 1)$ has generic stabilizer $E_6.\Z_2$, or $A_5.\Z_2$, or $A_4.\Z_2$, or ${A_2}^2.\Z_2$, or $\tilde A_2.\Z_2$, respectively.
\end{prop}

\begin{proof}
Throughout this proof we take $H$ to be the (simply connected) group defined over $K$ of type $E_8$, with simple roots $\beta_1, \dots, \beta_8$.

We begin with the case where $G = E_7$ and $\lambda = \omega_7$. Let $G$ have simple roots $\alpha_i = \beta_i$ for $i \leq 7$, so that $G = \langle X_\alpha : \alpha = \sum m_i \beta_i, \ m_8 = 0 \rangle < H$; then we may take $V = \langle e_\alpha : \alpha = \sum m_i \beta_i, \ m_8 = 1 \rangle < \L(H)$. Note that $Z(G) = \langle z \rangle$ where $z = h_{\beta_2}(-1) h_{\beta_5}(-1) h_{\beta_7}(-1)$; since $z$ acts on $V$ as multiplication by $-1$, we have $G_V = \{ 1 \}$.

We take the strictly positive generalized height function on the weight lattice of $G$ whose value at $\alpha_5$ is $2$, and at each other simple root $\alpha_i$ is $1$; then the generalized height of $\lambda = \frac{1}{2}(2\alpha_1 + 3\alpha_2 + 4\alpha_3 + 6\alpha_4 + 5\alpha_5 + 4\alpha_6 + 3\alpha_7)$ is $16$, and as $\lambda$ and $\Phi$ generate the weight lattice it follows that the generalized height of any weight is an integer. Since $V_\lambda = \langle e_\delta \rangle$ where $\delta = \eeightrt23465431$, we see that if $\mu \in \Lambda(V)$ and $e_\alpha \in V_\mu$ where $\alpha = \sum m_i \beta_i$ with $m_8 = 1$, then the generalized height of $\mu$ is $\sum_{i = 1}^7 m_i + m_5 - 16$. Thus $\Lambda(V)_{[0]} = \{ \nu_1, \nu_2 \}$, where we write
$$
\gamma_1 = \eeightrt12232221, \quad \gamma_2 = \eeightrt11233211,
$$
and for each $i$ we let $\nu_i$ be the weight such that $V_{\nu_i} = \langle e_{\gamma_i} \rangle$. Observe that if we take $s = \prod_{i = 1}^7 h_{\beta_i}(\kappa_i) \in T$ then $\nu_1(s) = \frac{\kappa_2 \kappa_7}{\kappa_5}$ and $\nu_2(s) = \frac{\kappa_5}{\kappa_2 \kappa_7}$; thus $\nu_1 + \nu_2 = 0$, and so $\Lambda(V)_{[0]}$ has ZLC. Set $Y = V_{[0]} = \langle e_{\gamma_1}, e_{\gamma_2} \rangle$, and
$$
\hat Y = \{ a_1 e_{\gamma_1} + a_2 e_{\gamma_2} : a_1a_2 \neq 0 \},
$$
so that $\hat Y$ is a dense open subset of $Y$. Write
$$
y_0 = e_{\gamma_1} + e_{\gamma_2} \in \hat Y.
$$

Note that $W$ acts transitively on the set $\Sigma$ of roots $\alpha$ of $H$ corresponding to the root vectors $e_\alpha$ spanning $V$. Thus if we write $W_1$ for the stabilizer in $W$ of $\gamma_1$, then $|W_1| = \frac{|W|}{|\Sigma|} = \frac{|W|}{56} = |W(E_6)|$; we then see that $W_1 = \langle w_{\beta_6}, w_{\beta_1}, w_{\beta_2 + \beta_4 + \beta_5}, w_{\beta_3},$ $w_{\beta_4}, w_{\beta_5 + \beta_6 + \beta_7} \rangle$. As each generator of $W_1$ stabilizes $\gamma_2$, the pointwise stabilizer in $W$ of $\{ \gamma_1, \gamma_2 \}$ is $W_1$. Now write $w^* = w_{\beta_2} w_{\beta_5} w_{\beta_7}$; as $w^*$ interchanges $\gamma_1$ and $\gamma_2$, the setwise stabilizer in $W$ of $\{ \gamma_1, \gamma_2 \}$, and hence of $\Lambda(V)_{[0]}$, is $W_1 \langle w^* \rangle$.

Let $A$ be the $E_6$ subgroup having simple roots $\beta_6$, $\beta_1$, $\beta_2 + \beta_4 + \beta_5$, $\beta_3$, $\beta_4$ and $\beta_5 + \beta_6 + \beta_7$; since $Z(A) = \langle z' \rangle$ where $z' = h_{\beta_2}({\eta_3}^2) h_{\beta_5}(\eta_3) h_{\beta_7}({\eta_3}^2)$, we see that $A$ is of simply connected type. Write $n^* = n_{\beta_2} n_{\beta_5} n_{\beta_7} h_{\beta_5}(-\eta_4) \in N$, so that $(n^*)^2 = z$ and conjugation by $n^*$ induces a graph automorphism of $A$; then $n^*.e_{\gamma_1} = \eta_4 e_{\gamma_2}$ and $n^*.e_{\gamma_2} = \eta_4 e_{\gamma_1}$. Set $C = A$ or $A \langle n^* \rangle$ according as $p \geq 3$ or $p = 2$, and $C' = A \langle n^* \rangle$. Clearly we then have $C \leq C_G(y_0)$ and $C' \leq C_G(\langle y_0 \rangle)$; we shall show that in fact $C_G(y_0) = C$ and $C_G(\langle y_0 \rangle) = C'$.

By Lemma~\ref{lem: gen height zero}, if we take $g \in \Tran_G(y_0, Y)$ and set $y' = g.y_0 \in Y$, then we have $g = u_1 n u_2$ with $u_1 \in C_U(y')$, $u_2 \in C_U(y_0)$, and $n \in N_{\Lambda(V)_{[0]}}$ with $n.y_0 = y'$. In particular $G.y_0 \cap Y = N_{\Lambda(V)_{[0]}}.y_0 \cap Y$, and $C_G(y_0) = C_U(y_0) C_{N_{\Lambda(V)_{[0]}}}(y_0) C_U(y_0)$ while $C_G(\langle y_0 \rangle) = C_U(y_0) C_{N_{\Lambda(V)_{[0]}}}(\langle y_0 \rangle) C_U(y_0)$.

First, from the above the elements of $W$ which preserve $\Lambda(V)_{[0]}$ are those corresponding to elements of $A \langle n^* \rangle \cap N$; so we have $N_{\Lambda(V)_{[0]}}.y_0 = T.y_0 \cup n^*T.y_0$. Since any element of $T$ may be written as $h_{\beta_5}(\kappa_5) t$ where $\kappa_5 \in K^*$ and $t \in A \cap T$, by the above we have
\begin{eqnarray*}
T.y_0    & = & \left\{ {\ts\frac{1}{\kappa_5}} e_{\gamma_1} + \kappa_5 e_{\gamma_2} : \kappa_5 \in K^* \right\}, \\
n^*T.y_0 & = & \left\{ \eta_4({\ts\frac{1}{\kappa_5}} e_{\gamma_2} + \kappa_5 e_{\gamma_1}) : \kappa_5 \in K^* \right\}.
\end{eqnarray*}
Hence $C_{N_{\Lambda(V)_{[0]}}}(y_0) = C \cap N$; also $N_{\Lambda(V)_{[0]}}.y_0 \subseteq \hat Y$, and $N_{\Lambda(V)_{[0]}}.y_0 \cap \langle y_0 \rangle = \{ {\eta_4}^i y_0 : i = 0, 1, 2, 3 \} = \langle n^* \rangle.y_0$, so $C_{N_{\Lambda(V)_{[0]}}}(\langle y_0 \rangle) = C' \cap N$.

Next, let $\Xi = \Phi^+ \setminus \Phi_A$, and set $U' = \prod_{\alpha \in \Xi} X_\alpha$; then $U = U'.(C \cap U)$ and $U' \cap (C \cap U) = \{ 1 \}$. We now observe that if $\alpha \in \Xi$ then $\nu_i + \alpha$ is a weight in $V$ for exactly one value of $i$; moreover each weight in $V$ of positive generalized height is of the form $\nu_i + \alpha$ for exactly one such root $\alpha$. Thus if we take $u = \prod x_\alpha(t_\alpha) \in U'$ satisfying $u.y_0 = y_0$, and equate coefficients of weight vectors, taking them in an order compatible with increasing generalized height, we see that for all $\alpha$ we must have $t_\alpha = 0$, so that $u = 1$; so $C_U(y_0) = C \cap U$.

Thus $C_U(y_0), C_{N_{\Lambda(V)_{[0]}}}(y_0) \leq C$ and $C_{N_{\Lambda(V)_{[0]}}}(\langle y_0 \rangle) \leq C'$, so we do indeed have $C_G(y_0) = C$ and $C_G(\langle y_0 \rangle) = C'$. Moreover $G.y_0 \cap Y = \{ b_1 e_{\gamma_1} + b_2 e_{\gamma_2} : (b_1b_2)^2 = 1 \}$.

Take $y = a_1 e_{\gamma_1} + a_2 e_{\gamma_2} \in \hat Y$. By the above, if we choose $\kappa \in K^*$ satisfying $\kappa^2 = a_1a_2$, then $\kappa^{-1}y \in T.y_0$, so there exists $h \in T$ with $h.y_0 = \kappa^{-1} y$; so $C_G(y) = C_G(\kappa^{-1} y) = C_G(h.y_0) = {}^h C$ and likewise $C_G(\langle y \rangle) = {}^h C'$. Moreover, we see that $G.y \cap Y = G.h.\kappa y_0 \cap Y = \kappa(G.y_0 \cap Y) = \{ b_1 e_{\gamma_1} + b_2 e_{\gamma_2} : (b_1b_2)^2 = (a_1a_2)^2 \}$. Since $\dim C = 78$, we have $\dim(\overline{G.y}) = \dim G - \dim C = 133 - 78 = 55$, while $\dim(\overline{G.y \cap Y}) = 1$; therefore
$$
\dim V - \dim(\overline{G.y}) = 56 - 55 = 1 \quad \hbox{and} \quad \dim Y - \dim(\overline{G.y \cap Y}) = 2 - 1 = 1.
$$
Hence $y$ is $Y$-exact. Thus the conditions of Lemma~\ref{lem: generic stabilizer from exact subset} hold; so the triple $(G, \lambda, p)$ has generic stabilizer $C/G_V \cong E_6.\Z_{(p, 2)}$, while the quadruple $(G, \lambda, p, 1)$ has generic stabilizer $C'/Z(G) \cong E_6.\Z_2$, where the $E_6$ is of simply connected type.

Now if we take the $D_6$ subgroup $\langle X_\alpha : \alpha = \sum m_i \beta_i, \ m_1 = m_8 = 0 \rangle$ of $E_7$, then $V|_{D_6} = V^{(0)} \oplus V^{(1)} \oplus V^{(2)}$, where $V^{(j)} = \langle e_\alpha : \alpha = \sum m_i \beta_i, \ m_8 = 1, \ m_1 = j \rangle$ for $j = 0, 1, 2$; of these three summands, $V^{(0)}$ and $V^{(2)}$ are natural $D_6$-modules, while $V^{(1)}$ is a half-spin $D_6$-module and contains $Y$. To treat the case where $G = D_6$ and $\lambda = \omega_6$, we may therefore replace $G$ by $D_6$ and $V$ by $V^{(1)}$. We then have $Z(G) = \langle z_1, z_2 \rangle$ where $z_1 = h_{\beta_3}(-1) h_{\beta_5}(-1) h_{\beta_7}(-1)$ and $z_2 = h_{\beta_2}(-1) h_{\beta_3}(-1)$; since $z_1$ and $z_2$ act on $V$ as multiplication by $1$ and $-1$ respectively, we have $G_V = \langle z_1 \rangle$. We replace $A$ by the intersection of that above with $G$, which is the $A_5$ subgroup having simple roots $\beta_6$, $\beta_2 + \beta_4 + \beta_5$, $\beta_3$, $\beta_4$ and $\beta_5 + \beta_6 + \beta_7$; since $Z(A) = \langle z' \rangle$ where $z' = h_{\beta_2}({\eta_6}^2) h_{\beta_3}(-1) h_{\beta_5}(\eta_6) h_{\beta_7}({\eta_6}^5)$, we see that $A$ is of simply connected type. We again set $n^* = n_{\beta_2} n_{\beta_5} n_{\beta_7} h_{\beta_5}(-\eta_4) \in N$, and then $(n^*)^2 = z_1z_2$ and conjugation by $n^*$ still induces a graph automorphism of $A$; we again set $C = A$ or $A \langle n^* \rangle$ according as $p \geq 3$ or $p = 2$, and $C' = A \langle n^* \rangle$. Take $y \in \hat Y$ and $h$ as above; again we have $C_G(y) = {}^h C$ and $C_G(\langle y \rangle) = {}^h C'$, and as we still have $\{ h_{\beta_5}(\kappa_5) : \kappa_5 \in K^* \} \subset T$ we see that $G.y \cap Y$ is as before. Since $\dim C = 35$, we have $\dim(\overline{G.y}) = \dim G - \dim C = 66 - 35 = 31$, while $\dim(\overline{G.y \cap Y}) = 1$; therefore
$$
\dim V - \dim(\overline{G.y}) = 32 - 31 = 1 \quad \hbox{and} \quad \dim Y - \dim(\overline{G.y \cap Y}) = 2 - 1 = 1.
$$
Hence $y$ is $Y$-exact. Thus the conditions of Lemma~\ref{lem: generic stabilizer from exact subset} hold; so the triple $(G, \lambda, p)$ has generic stabilizer $C/G_V \cong A_5.\Z_{(p, 2)}$, while the quadruple $(G, \lambda, p, 1)$ has generic stabilizer $C'/Z(G) \cong A_5.\Z_2$, where the $A_5$ has centre of order $3/(p, 3)$.

To treat the case where $G = B_5$ and $\lambda = \omega_5$ we leave $V$, $Y$ and $\hat Y$ unchanged, but replace $G$ by the $B_5$ subgroup of $D_6$ which has simple root groups $X_{\beta_7}$, $X_{\beta_6}$, $X_{\beta_5}$, $X_{\beta_4}$ and $\{ x_{\beta_2}(t) x_{\beta_3}(t) : t \in K \}$. Here we have $Z(G) = \langle z_2 \rangle$ where $z_2$ is as above, so $G_V = \{ 1 \}$. We replace $A$ by the intersection of that above with $G$, which is the $A_4$ subgroup having simple roots $\beta_6$, $\beta_2 + \beta_3 + \beta_4 + \beta_5$, $\beta_4$ and $\beta_5 + \beta_6 + \beta_7$; since $Z(A) = \langle z' \rangle$ where $z' = h_{\beta_2}({\eta_5}^2) h_{\beta_3}({\eta_5}^2) h_{\beta_5}(\eta_5) h_{\beta_7}({\eta_5}^4)$, we see that $A$ is of simply connected type. We also replace $n^*$ by $n_{\beta_2} n_{\beta_3} n_{\beta_5} n_{\beta_7} h_{\beta_5}(-\eta_4) \in N$, and then $(n^*)^2 = z_2s$ where $s = h_{\beta_6}(-1) h_{\beta_5 + \beta_6 + \beta_7}(-1) \in A \cap T$, and conjugation by $n^*$ still induces a graph automorphism of $A$; again let $C = A$ or $A \langle n^* \rangle$ according as $p \geq 3$ or $p = 2$, and $C' = A \langle n^* \rangle$. Take $y \in \hat Y$ and $h$ as above; again we have $C_G(y) = {}^h C$ and $C_G(\langle y \rangle) = {}^h C'$, and as we still have $\{ h_{\beta_5}(\kappa_5) : \kappa_5 \in K^* \} \subset T$ we see that $G.y \cap Y$ is as before. Since $\dim C = 24$, we have $\dim(\overline{G.y}) = \dim G - \dim C = 55 - 24 = 31$, while $\dim(\overline{G.y \cap Y}) = 1$; therefore
$$
\dim V - \dim(\overline{G.y}) = 32 - 31 = 1 \quad \hbox{and} \quad \dim Y - \dim(\overline{G.y \cap Y}) = 2 - 1 = 1.
$$
Hence $y$ is $Y$-exact. Thus the conditions of Lemma~\ref{lem: generic stabilizer from exact subset} hold; so the triple $(G, \lambda, p)$ has generic stabilizer $C/G_V \cong A_4.\Z_{(p, 2)}$, while the quadruple $(G, \lambda, p, 1)$ has generic stabilizer $C'/Z(G) \cong A_4.\Z_2$, where the $A_4$ is of simply connected type.

Now if we take the $A_5$ subgroup $\langle X_\alpha : \alpha = \sum m_i \beta_i, \ m_1 = m_3 = m_8 = 0 \rangle$ of $D_6$, then $V|_{A_5} = V^{(1)} \oplus V^{(2)} \oplus V^{(3)}$, where $V^{(j)} = \langle e_\alpha : \alpha = \sum m_i \beta_i, \ m_1 = m_8 = 1 , \ m_3 = j \rangle$ for $j = 1, 2, 3$; of these three summands, $V^{(1)}$ and $V^{(3)}$ are natural $A_5$-modules or their duals, while $V^{(2)}$ is the exterior cube of the natural $A_5$-module and contains $Y$. To treat the case where $G = A_5$ and $\lambda = \omega_3$, we may therefore replace $G$ by $A_5$ and $V$ by $V^{(2)}$. We then have $Z(G) = \langle z_3 \rangle$ where $z_3 = h_{\beta_2}(\eta_6) h_{\beta_4}({\eta_6}^2) h_{\beta_5}(-1) h_{\beta_6}({\eta_6}^4) h_{\beta_7}({\eta_6}^5)$; since $z_3$ acts on $V$ as multiplication by $-1$, we have $G_V = \langle {z_3}^2 \rangle$ or $\langle z_3 \rangle$ according as $p \geq 3$ or $p = 2$. We replace $A$ by the intersection of that from the $D_6$ case with $G$, which is the ${A_2}^2$ subgroup with one factor having simple roots $\beta_6$ and $\beta_2 + \beta_4 + \beta_5$, and the other factor having simple roots $\beta_4$ and $\beta_5 + \beta_6 + \beta_7$; since $Z(A) = \langle {z_1}', {z_2}' \rangle$ where ${z_1}' = h_{\beta_2}({\eta_3}^2) h_{\beta_4}({\eta_3}^2) h_{\beta_5}({\eta_3}^2) h_{\beta_6}(\eta_3)$ and ${z_2}' = h_{\beta_4}(\eta_3) h_{\beta_5}({\eta_3}^2) h_{\beta_6}({\eta_3}^2) h_{\beta_7}({\eta_3}^2)$, we see that $A$ is of simply connected type. We again set $n^* = n_{\beta_2} n_{\beta_5} n_{\beta_7} h_{\beta_5}(-\eta_4) \in N$, and then $(n^*)^2 = {z_3}^3$ and conjugation by $n^*$ still induces a graph automorphism of $A$, which here interchanges the simple factors; we again set $C = A$ or $A \langle n^* \rangle$ according as $p \geq 3$ or $p = 2$, and $C' = A \langle n^* \rangle$. Take $y \in \hat Y$ and $h$ as above; again we have $C_G(y) = {}^h C$ and $C_G(\langle y \rangle) = {}^h C'$, and as we still have $\{ h_{\beta_5}(\kappa_5) : \kappa_5 \in K^* \} \subset T$ we see that $G.y \cap Y$ is as before. Since $\dim C = 16$, we have $\dim(\overline{G.y}) = \dim G - \dim C = 35 - 16 = 19$, while $\dim(\overline{G.y \cap Y}) = 1$; therefore
$$
\dim V - \dim(\overline{G.y}) = 20 - 19 = 1 \quad \hbox{and} \quad \dim Y - \dim(\overline{G.y \cap Y}) = 2 - 1 = 1.
$$
Hence $y$ is $Y$-exact. Thus the conditions of Lemma~\ref{lem: generic stabilizer from exact subset} hold; so the triple $(G, \lambda, p)$ has generic stabilizer $C/G_V \cong {A_2}^2.\Z_{(p, 2)}$, while the quadruple $(G, \lambda, p, 1)$ has generic stabilizer $C'/Z(G) \cong {A_2}^2.\Z_2$, where the ${A_2}^2$ has centre of order $3/(p, 3)$.

Finally if we take the $C_3$ subgroup of $A_5$ with simple root groups $\{ x_{\beta_2}(t) x_{\beta_7}(t) : t \in K \}$, $\{ x_{\beta_4}(t) x_{\beta_6}(t) : t \in K \}$ and $X_{\beta_5}$, then $V|_{C_3} = V' \oplus V''$, where $V'$ is a natural $C_3$-module, while $V''$ has highest weight $\omega_3$ and contains $Y$. To treat the case where $G = C_3$ and $\lambda = \omega_3$ with $p \geq 3$, we may therefore replace $G$ by $C_3$ and $V$ by $V''$. Here we have $Z(G) = \langle {z_3}^3 \rangle$ where $z_3$ is as above, so $G_V = \{ 1 \}$. We replace $A$ by the intersection of that above with $G$, which is the $\tilde A_2$ subgroup having simple root groups $\{ x_{\beta_4}(t) x_{\beta_6}(t) : t \in K \}$ and $\{ x_{\beta_2 + \beta_4 + \beta_5}(t) x_{\beta_5 + \beta_6 + \beta_7}(t) : t \in K \}$; since $Z(A) = \langle z' \rangle$ where $z' = h_{\beta_2}({\eta_3}^2) h_{\beta_5}(\eta_3) h_{\beta_7}({\eta_3}^2)$, we see that $A$ is of simply connected type. We again set $n^* = n_{\beta_2} n_{\beta_5} n_{\beta_7} h_{\beta_5}(-\eta_4) \in N$, and then $(n^*)^2 = {z_3}^3$ and conjugation by $n^*$ still induces a graph automorphism of $A$; we again set $C = A$ and $C' = A \langle n^* \rangle$. Take $y \in \hat Y$ and $h$ as above; again we have $C_G(y) = {}^h C$ and $C_G(\langle y \rangle) = {}^h C'$, and as we still have $\{ h_{\beta_5}(\kappa_5) : \kappa_5 \in K^* \} \subset T$ we see that $G.y \cap Y$ is as before. Since $\dim C = 8$, we have $\dim(\overline{G.y}) = \dim G - \dim C = 21 - 8 = 13$, while $\dim(\overline{G.y \cap Y}) = 1$; therefore
$$
\dim V - \dim(\overline{G.y}) = 14 - 13 = 1 \quad \hbox{and} \quad \dim Y - \dim(\overline{G.y \cap Y}) = 2 - 1 = 1.
$$
Hence $y$ is $Y$-exact. Thus the conditions of Lemma~\ref{lem: generic stabilizer from exact subset} hold; so the triple $(G, \lambda, p)$ has generic stabilizer $C/G_V \cong \tilde A_2$, while the quadruple $(G, \lambda, p, 1)$ has generic stabilizer $C'/Z(G) \cong \tilde A_2.\Z_2$, where the $\tilde A_2$ is of simply connected type.
\end{proof}

\begin{prop}\label{prop: D_5, omega_5, B_4, omega_4 modules}
Let $G = D_5$ and $\lambda = \omega_5$, or $G = B_4$ and $\lambda = \omega_4$. Then the triple $(G, \lambda, p)$ has generic stabilizer $B_3 U_8$ or $B_3$ respectively; the associated first quadruple $(G, \lambda, p, 1)$ has generic stabilizer $B_3 T_1 U_8$ or $B_3$ respectively.
\end{prop}

\begin{proof}
Throughout this proof we take $H$ to be the simply connected group defined over $K$ of type $E_6$, with simple roots $\beta_1, \dots, \beta_6$.

We begin with the case where $G = D_5$ and $\lambda = \omega_5$. Let $G$ have simple roots $\alpha_1 = \beta_1$, $\alpha_2 = \beta_3$, $\alpha_3 = \beta_4$, $\alpha_4 = \beta_5$, $\alpha_5 = \beta_2$, so that $G = \langle X_\alpha : \alpha = \sum m_i \beta_i, \ m_6 = 0 \rangle < H$; then we may take $V = \langle e_\alpha : \alpha = \sum m_i \beta_i, \ m_6 = 1 \rangle < \L(H)$. Note that $Z(G) = \langle z \rangle$ where $z = h_{\beta_1}(-1) h_{\beta_2}(\eta_4) h_{\beta_4}(-1) h_{\beta_5}(-\eta_4)$; since $z$ acts on $V$ as multiplication by $\eta_4$, we have $G_V = \{ 1 \}$.

Write
$$
\gamma_1 = \esixrt111221, \quad \gamma_2 = \esixrt112211.
$$
Set $Y = \langle e_{\gamma_1}, e_{\gamma_2} \rangle$, and let
$$
\hat Y = \{ a_1 e_{\gamma_1} + a_2 e_{\gamma_2} : a_1a_2 \neq 0 \},
$$
so that $\hat Y$ is a dense open subset of $Y$. Write
$$
y_0 = e_{\gamma_1} + e_{\gamma_2} \in \hat Y.
$$

Note that $W$ acts transitively on the set $\Sigma$ of roots $\alpha$ of $H$ corresponding to the root vectors $e_\alpha$ spanning $V$. Thus if we write $W_1$ for the stabilizer in $W$ of $\gamma_1$, then $|W_1| = \frac{|W|}{|\Sigma|} = \frac{|W|}{16} = |W(A_4)|$; we then see that $W_1 = \langle w_{\beta_3 + \beta_4 + \beta_5}, w_{\beta_2}, w_{\beta_4},$ $w_{\beta_1 + \beta_3} \rangle$. Now the stabilizer in $W$ of any $\alpha \in \Sigma$ acts transitively on the set $\Sigma'$ of roots $\alpha' \in \Sigma$ orthogonal to $\alpha$ (this is evident if we take $\alpha = \esixrt000001$, as then its stabilizer in $W$ is $\langle w_{\beta_1}, w_{\beta_3}, w_{\beta_4}, w_{\beta_2} \rangle$, which acts transitively on the set of roots $\alpha' = \sum m_i \beta_i$ with $m_5 = 2$ and $m_6 = 1$). Thus if we write $W_2$ for the stabilizer in $W_1$ of $\gamma_2$, then $|W_2| = \frac{|W_1|}{|\Sigma'|} = \frac{|W_1|}{5} = |W(A_3)|$; we then see that $W_2 = \langle w_{\beta_3 + \beta_4 + \beta_5}, w_{\beta_2}, w_{\beta_4} \rangle$. Thus the pointwise stabilizer in $W$ of $\{ \gamma_1, \gamma_2 \}$ is $W_2$. Now write $w^* = w_{\beta_3} w_{\beta_5}$; as $w^*$ interchanges $\gamma_1$ and $\gamma_2$, the setwise stabilizer in $W$ of $\{ \gamma_1, \gamma_2 \}$ is $W_2\langle w^* \rangle = \langle w_{\beta_2}, w_{\beta_4}, w_{\beta_3} w_{\beta_5} \rangle$.

Let $P = QL$ be the standard $D_4$ parabolic subgroup of $G$, with Levi subgroup $L = \langle T, X_\alpha : \alpha = \sum m_i \beta_i, \ m_1 = m_6 = 0 \rangle$ and $8$-dimensional unipotent radical $Q = \langle X_\alpha : \alpha = \sum m_i \beta_i, \ m_1 = 1, \ m_6 = 0 \rangle$; then each element of $Q$ fixes each element of $Y$. Write $P^- = Q^-L$ for the opposite parabolic subgroup, so that $Q^- = \langle X_\alpha : \alpha = \sum m_i \beta_i, \ m_1 = -1, \ m_6 = 0 \rangle$. Let $A$ be the $B_3$ subgroup of $L$ with simple root groups $X_{\beta_2}$, $X_{\beta_4}$ and $\{ x_{\beta_3}(t) x_{\beta_5}(-t) : t \in K \}$; as $Z(A) = \langle z' \rangle$ where $z' = h_{\beta_3}(-1) h_{\beta_5}(-1)$, we see that $A$ is of simply connected type. Set $C = QA$ and $C' = Z(L)QA$. Clearly we then have $C \leq C_G(y_0)$ and $C' \leq C_G(\langle y_0 \rangle)$; we shall show that in fact $C_G(y_0) = C$ and $C_G(\langle y_0 \rangle) = C'$.

Suppose first that $g \in G$ satisfies $g.y_0 \in Y$; write $g.y_0 = y$. Using Lemma~\ref{lem: parabolic factorization} we may write $g = q_1xq_2q_3$, where $q_1, q_3 \in Q$, $q_2 \in Q^-$ and $x \in L$; then we have $xq_2q_3.y_0 = {q_1}^{-1}.y$, whence $xq_2.y_0 = y$. Now if $q_2 \neq 1$ then $q_2.y_0$ has at least one term $e_\alpha$ for a root $\alpha$ of the form $\sum m_i \beta_i$ with $m_1 = 0$, as therefore does $xq_2.y_0$, contrary to $y \in Y$; so we must have $q_2 = 1$, and hence $g = q_1xq_3$ and $x.y_0 = y$. Now write $x = tx'$ where $x' \in L'$ and $t = h_{\beta_1}(\kappa^2) h_{\beta_2}(\kappa) h_{\beta_3}(\kappa^2) h_{\beta_4}(\kappa^2) h_{\beta_5}(\kappa) \in Z(L)$ for some $\kappa \in K^*$; then $t.y_0 = \kappa y_0$ and so we have $x'.y_0 = \kappa^{-1} y$.

Now consider the action of the $D_4$ subgroup $L'$ on $V$: the subspace $V' = \langle e_\alpha : \alpha = \sum m_i \beta_i, \ m_1 = m_6 = 1 \rangle$ is an irreducible $L'$-module of highest weight $\lambda' = \frac{1}{2}(\alpha_2 + 2\alpha_3 + \alpha_4 + 2\alpha_5)$, and contains $Y$. We take the strictly positive generalized height function on the weight lattice of $L'$ whose value at each $\alpha_i$ (with $i > 1$) is $1$; then the generalized height of $\lambda'$ is $3$, and as $\lambda'$, $\lambda' + \frac{1}{2}\alpha_4 - \frac{1}{2}\alpha_5$ and $\Phi(L')$ generate the weight lattice it follows that the generalized height of any weight is an integer. Since ${V'}_{\lambda'} = \langle e_\delta \rangle$ where $\delta = \esixrt122321$, we see that if $\mu \in \Lambda(V')$ and $e_\alpha \in {V'}_\mu$, where $\alpha = \sum m_i \beta_i$ with $m_1 = m_6 = 1$, then the generalized height of $\mu$ is $\sum_{i = 2}^5 m_i - 6$. Thus if for each $i$ we let $\nu_i$ be the weight such that ${V'}_{\nu_i} = \langle e_{\gamma_i} \rangle$, then $\Lambda(V')_{[0]} = \{ \nu_1, \nu_2 \}$. Observe that if we take $s = \prod_{i = 2}^5 h_{\beta_i}(\kappa_i) \in L' \cap T$, then $\nu_1(s) = \frac{\kappa_5}{\kappa_3}$ and $\nu_2(s) = \frac{\kappa_3}{\kappa_5}$; thus $\nu_1 + \nu_2 = 0$, and hence $\Lambda(V')_{[0]}$ has ZLC.

By Lemma~\ref{lem: gen height zero} we see that $x' = u_1nu_2$ where $u_1 \in C_{L' \cap U}(\kappa^{-1} y) = C_{L' \cap U}(y)$, $u_2 \in C_{L' \cap U}(y_0)$ and $n \in L' \cap N$ with $n.y_0 = \kappa^{-1} y$. From the above the elements of $W$ which preserve $\{ \gamma_1, \gamma_2 \}$ are those corresponding to elements of $A \cap N$, so $n \in \{ h_{\beta_3}(\kappa_3) : \kappa_3 \in K^* \} (A \cap N)$; then we may write $n = sn'$ where $s = h_{\beta_3}(\kappa_3)$ for some $\kappa_3 \in K^*$ and $n' \in A \cap N$, and so $\kappa^{-1}y = sn'.y_0 = s.y_0$. Clearly $C_{L' \cap U}(y_0) = A \cap U$, and so $C_{L' \cap U}(y) = C_{L' \cap U}(\kappa^{-1} y) = C_{L' \cap U}(s.y_0) = {}^s (A \cap U)$; therefore $x' \in {}^s (A \cap U).s(A \cap N).(A \cap U) = s(A \cap U)(A \cap N)(A \cap U)$, and so we have $x' = sa$ for some $a \in A$.

Hence $x = tx' = tsa$; so $g = q_1tsaq_3 = ts.({q_1}^{ts})({}^a{q_3}).a \in TQA = TC$. In particular, if $g.y_0 \in \langle y_0 \rangle$ we must have $s.y_0 \in \langle y_0 \rangle$, so as $s.y_0 = {\kappa_3}^{-1} e_{\gamma_1} + \kappa_3 e_{\gamma_2}$ we must have $\kappa_3 = \pm1$; then $s = s_1s_2$ where $s_1 = h_{\beta_2}(\kappa_3) h_{\beta_5}(\kappa_3) \in Z(L)$ and $s_2 = h_{\beta_2}(\kappa_3) h_{\beta_3}(\kappa_3) h_{\beta_5}(\kappa_3) \in A \cap T$. Therefore $C_G(\langle y_0 \rangle) = C'$; and as $C_{Z(L)}(y_0) = \{ 1 \}$ we also have $C_G(y_0) = C$. Moreover we see that $G.y_0 \cap Y = TC.y_0 = T.y_0 = \hat Y$, since given $y = a_1 e_{\gamma_1} + a_2 e_{\gamma_2} \in \hat Y$ we have $y = h.y_0$ for $h = h_{\alpha_1}(a_1a_2) h_{\alpha_3}(a_2)$.

Take $y \in \hat Y$. By the above, there exists $h \in T$ with $h.y_0 = y$; so $C_G(y) = C_G(h.y_0) = {}^h C$ and likewise $C_G(\langle y \rangle) = {}^h C'$. Since $\dim C = 29$, we have $\dim(\overline{G.y}) = \dim G - \dim C = 45 - 29 = 16$, while $\dim(\overline{G.y \cap Y}) = 2$; therefore
$$
\dim V - \dim(\overline{G.y}) = 16 - 16 = 0 \quad \hbox{and} \quad \dim Y - \dim(\overline{G.y \cap Y}) = 2 - 2 = 0.
$$
Hence $y$ is $Y$-exact. Thus the conditions of Lemma~\ref{lem: generic stabilizer from exact subset} hold; so the triple $(G, \lambda, p)$ has generic stabilizer $C/G_V \cong B_3U_8$, while the quadruple $(G, \lambda, p, 1)$ has generic stabilizer $C'/Z(G) \cong B_3 T_1 U_8$, where the $B_3$ is of simply connected type.

To treat the case where $G = B_4$ and $\lambda = \omega_4$, we leave $H$, $V$, $Y$ and $\hat Y$ unchanged, but replace $G$ by the $B_4$ subgroup of $D_5$ having simple root groups $X_{\beta_5}$, $X_{\beta_4}$, $X_{\beta_3}$ and $\{ x_{\beta_1}(t) x_{-(\beta_1 + \beta_2 + 2\beta_3 + 2\beta_4 + \beta_5)}(t) : t \in K \}$. Here we have $Z(G) = \langle z^2 \rangle$ where $z$ is as above, so again $G_V = \{ 1 \}$. We have $Q \cap G = \{ 1 \}$ while $A < G$, and $Z(L) \cap G = Z(G)$; so we replace $C$ by $A$ and $C'$ by $Z(G) A$.

Take $y = a_1 e_{\gamma_1} + a_2 e_{\gamma_2} \in \hat Y$. If we choose $\kappa \in K^*$ satisfying $\kappa^2 = a_1a_2$ and set $h = h_{\beta_5}(\kappa^{-1} a_1)$, then $h.y_0 = \kappa^{-1} y$; so $C_G(y) = C_G(\kappa^{-1} y) = C_G(h.y_0) = {}^h C$ and likewise $C_G(\langle y \rangle) = {}^h C'$. Now if $g \in G$ satisfies $g.y \in Y$, by the above we must certainly have $g.y \in \hat Y$, whence $g.h.\kappa y_0 \in \hat Y$, so $gh.y_0 \in \hat Y$; then we must have $gh \in TC = Th^{-1}C_G(y)h = TC_G(y)h$, so $g \in TC_G(y)$, and so $g.y \in TC_G(y).y = T.y$. Thus $G.y \cap Y \subset T.y$; the reverse inclusion is obvious. Moreover, since any element of $T$ may be written as $h_{\beta_3}(\kappa_3) t$ where $\kappa_3 \in K^*$ and $t \in C \cap T$, by the above we have
$$
T.y = \left\{ {\ts\frac{1}{\kappa_3}} a_1 e_{\gamma_1} + \kappa_3 a_2 e_{\gamma_2} : \kappa_3 \in K^* \right\}.
$$
Hence $G.y \cap Y = \{ b_1 e_{\gamma_1} + b_2 e_{\gamma_2} : b_1b_2 = a_1a_2 \}$. Since $\dim C = 21$, we have $\dim(\overline{G.y}) = \dim G - \dim C = 36 - 21 = 15$, while $\dim(\overline{G.y \cap Y}) = 1$; therefore
$$
\dim V - \dim(\overline{G.y}) = 16 - 15 = 1 \quad \hbox{and} \quad \dim Y - \dim(\overline{G.y \cap Y}) = 2 - 1 = 1.
$$
Hence $y$ is $Y$-exact. Thus the conditions of Lemma~\ref{lem: generic stabilizer from exact subset} hold; so the triple $(G, \lambda, p)$ has generic stabilizer $C/G_V \cong B_3$, while the quadruple $(G, \lambda, p, 1)$ has generic stabilizer $C'/Z(G) \cong B_3$, where the $B_3$ is of simply connected type.
\end{proof}

\begin{prop}\label{prop: D_7, omega_7, B_6, omega_6 modules}
Let $G = D_7$ and $\lambda = \omega_7$, or $G = B_6$ and $\lambda = \omega_6$. Then the triple $(G, \lambda, p)$ has generic stabilizer ${G_2}^2.\Z_{(p, 2)}$ or ${A_2}^2.{\Z_{(p, 2)}}^2$ respectively; the associated first quadruple $(G, \lambda, p, 1)$ has generic stabilizer ${G_2}^2.\Z_2$ or ${A_2}^2.\Z_{(p, 2)}.\Z_2$ respectively.
\end{prop}

\begin{proof}
Throughout this proof we take $H$ to be the (simply connected) group defined over $K$ of type $E_8$, with simple roots $\beta_1, \dots, \beta_8$.

We begin with the case where $G = D_7$ and $\lambda = \omega_7$. Let $G$ have simple roots $\alpha_i = \beta_{9 - i}$ for $i \leq 7$, so that $G = \langle X_\alpha : \alpha = \sum m_i \beta_i, \ m_1 = 0 \rangle < H$; then we may take $V = \langle e_\alpha : \alpha = \sum m_i \beta_i, \ m_1 = 1 \rangle < \L(H)$. Note that $Z(G) = \langle z \rangle$ where $z = h_{\beta_2}(\eta_4) h_{\beta_3}(-\eta_4) h_{\beta_4}(-1) h_{\beta_6}(-1) h_{\beta_8}(-1)$; since $z$ acts on $V$ as multiplication by $\eta_4$, we have $G_V = \{ 1 \}$.

We take the strictly positive generalized height function on the weight lattice of $G$ whose value at $\alpha_1$ is $4$, and at each other simple root $\alpha_i$ is $1$; then the generalized height of $\lambda = \frac{1}{2}(\alpha_1 + 2\alpha_2 + 3\alpha_3 + 4\alpha_4 + 5\alpha_5 + \frac{5}{2}\alpha_6 + \frac{7}{2}\alpha_7)$ is $12$, and as $\lambda$, $\omega_6 = \lambda + \frac{1}{2}\alpha_6 - \frac{1}{2}\alpha_7$ and $\Phi$ generate the weight lattice it follows that the generalized height of any weight is an integer. Since $V_\lambda = \langle e_\delta \rangle$ where $\delta = \eeightrt13354321$, we see that if $\mu \in \Lambda(V)$ and $e_\alpha \in V_\mu$ where $\alpha = \sum m_i \beta_i$ with $m_1 = 1$, then the generalized height of $\mu$ is $\sum_{i = 2}^7 m_i + 4m_8 - 12$. Thus $\Lambda(V)_{[0]} = \{ \nu_1, \nu_2, \nu_3, \nu_4 \}$, where we write
$$
\gamma_1 = \eeightrt12232210, \quad \gamma_2 = \eeightrt11122111, \quad \gamma_3 = \eeightrt11233210, \quad \gamma_4 = \eeightrt11221111,
$$
and for each $i$ we let $\nu_i$ be the weight such that $V_{\nu_i} = \langle e_{\gamma_i} \rangle$. Observe that if we take $s = \prod_{i = 2}^8 h_{\beta_i}(\kappa_i) \in T$, then $\nu_1(s) = \frac{\kappa_2 \kappa_6}{\kappa_5 \kappa_8}$, $\nu_2(s) = \frac{\kappa_5 \kappa_8}{\kappa_3 \kappa_6}$, $\nu_3(s) = \frac{\kappa_5}{\kappa_2 \kappa_8}$ and $\nu_4(s) = \frac{\kappa_3 \kappa_8}{\kappa_5}$; thus $\nu_1 + \nu_2 + \nu_3 + \nu_4 = 0$, and so $\Lambda(V)_{[0]}$ has ZLC. Set $Y = V_{[0]} = \langle e_{\gamma_1}, e_{\gamma_2}, e_{\gamma_3}, e_{\gamma_4} \rangle$, and
$$
\hat Y = \{ a_1 e_{\gamma_1} + a_2 e_{\gamma_2} + a_3 e_{\gamma_3} + a_4 e_{\gamma_4} : a_1a_2a_3a_4 \neq 0 \},
$$
so that $\hat Y$ is a dense open subset of $Y$. Write
$$
y_0 = e_{\gamma_1} + e_{\gamma_2} + e_{\gamma_3} + e_{\gamma_4} \in \hat Y.
$$

Note that $W$ acts transitively on the set $\Sigma$ of roots $\alpha$ of $H$ corresponding to the root vectors $e_\alpha$ spanning $V$. Thus if we write $W_1$ for the stabilizer in $W$ of $\gamma_1$, then $|W_1| = \frac{|W|}{|\Sigma|} = \frac{|W|}{64} = |W(A_6)|$; we then see that $W_1 = \langle w_{\beta_7}, w_{\beta_5 + \beta_6}, w_{\beta_4}, w_{\beta_3},$ $w_{\beta_2 + \beta_4 + \beta_5}, w_{\beta_6 + \beta_7 + \beta_8} \rangle$. Now the stabilizer in $W$ of any $\alpha \in \Sigma$ acts transitively on the set $\Sigma'$ of roots $\alpha' \in \Sigma$ such that $\alpha + \alpha' \in \Phi_H$ (this is evident if we take $\alpha = \eeightrt10000000$, as then its stabilizer in $W$ is $\langle w_{\beta_2}, w_{\beta_4}, w_{\beta_5}, w_{\beta_6}, w_{\beta_7}, w_{\beta_8} \rangle$, which acts transitively on the set of roots $\alpha' = \sum m_i \beta_i$ with $m_1 = 1$ and $m_3 = 3$). Thus if we write $W_2$ for the stabilizer in $W_1$ of $\gamma_2$, then $|W_2| = \frac{|W_1|}{|\Sigma'|} = \frac{|W_1|}{7} = |W(A_5)|$; we then see that $W_2 = \langle w_{\beta_7}, w_{\beta_5 + \beta_6}, w_{\beta_4}, w_{\beta_2 + \beta_3 + \beta_4 + \beta_5}, w_{\beta_6 + \beta_7 + \beta_8} \rangle$. Next the joint stabilizer in $W$ of any pair of roots $\alpha, \alpha' \in \Sigma$ such that $\alpha + \alpha' \in \Phi_H$ acts transitively on the set $\Sigma''$ of roots $\alpha'' \in \Sigma$ orthogonal to both $\alpha$ and $\alpha'$ (this is evident if we take $\alpha = \eeightrt10000000$ and $\alpha' = \eeightrt12343210$, as then the joint stabilizer in $W$ is $\langle w_{\beta_2}, w_{\beta_4}, w_{\beta_5}, w_{\beta_6}, w_{\beta_7} \rangle$, which acts transitively on the set of roots $\alpha'' = \sum m_i \beta_i$ with $m_1 = 1$, $m_3 = 2$ and $m_8 = 1$). Thus if we write $W_3$ for the stabilizer in $W_2$ of $\gamma_3$, then $|W_3| = \frac{|W_2|}{|\Sigma''|} = \frac{|W_2|}{20} = |W({A_2}^2)|$; we then see that $W_3 = \langle w_{\beta_4}, w_{\beta_2 + \beta_3 + \beta_4 + \beta_5}, w_{\beta_7}, w_{\beta_2 + \beta_3 + 2\beta_4 + 2\beta_5 + 2\beta_6 + \beta_7 + \beta_8} \rangle$. As each generator of $W_3$ stabilizes $\gamma_4$, the pointwise stabilizer in $W$ of $\{ \gamma_1, \gamma_2, \gamma_3, \gamma_4 \}$ is $W_3$. Now write
\begin{eqnarray*}
w^*     & = & w_{\beta_4 + \beta_5 + \beta_6} w_{\beta_5 + \beta_6 + \beta_7} w_{\beta_6 + \beta_7 + \beta_8}, \\
w^{**}  & = & w_{\beta_2 + \beta_3 + \beta_4 + \beta_5 + \beta_6} w_{\beta_6 + \beta_7} w_{\beta_5 + \beta_6 + \beta_7 + \beta_8}, \\
w^{***} & = & w_{\beta_2} w_{\beta_3} w_{\beta_5};
\end{eqnarray*}
then $w^*$ interchanges $\gamma_3$ and $\gamma_4$ while fixing both $\gamma_1$ and $\gamma_2$, and similarly $w^{**}$ interchanges $\gamma_1$ and $\gamma_2$ while fixing both $\gamma_3$ and $\gamma_4$, while $w^{***}$ interchanges $\gamma_1$ and $\gamma_3$, and also $\gamma_2$ and $\gamma_4$. Thus as $\gamma_1$ is orthogonal to $\gamma_3$ and $\gamma_4$ but not $\gamma_2$, the setwise stabilizer in $W$ of $\{ \gamma_1, \gamma_2, \gamma_3, \gamma_4 \}$, and hence of $\Lambda(V)_{[0]}$, is
\begin{eqnarray*}
W_3 \langle w^*, w^{**}, w^{***} \rangle & = & \langle w_{\beta_4}, w_{\beta_2} w_{\beta_3} w_{\beta_5}, w_{\beta_7}, w_{\beta_2 + \beta_4 + \beta_5 + \beta_6} w_{\beta_3 + \beta_4 + \beta_5 + \beta_6} w_{\beta_8}, \\
                                         &   & \phantom{\langle} w_{\beta_4 + \beta_5 + \beta_6} w_{\beta_5 + \beta_6 + \beta_7} w_{\beta_6 + \beta_7 + \beta_8} \rangle.
\end{eqnarray*}

Set ${\beta_2}' = \beta_2 + \beta_4 + \beta_5 + \beta_6$ and ${\beta_3}' = \beta_3 + \beta_4 + \beta_5 + \beta_6$; let $A$ be the ${G_2}^2$ subgroup with one factor having simple root groups $\{ x_{\beta_2}(-t) x_{\beta_3}(-t) x_{\beta_5}(t) : t \in K \}$ and $X_{\beta_4}$, and the other factor having simple root groups $\{ x_{{\beta_2}'}(-t) x_{{\beta_3}'}(t) x_{\beta_8}(t) : t \in K \}$ and $X_{\beta_7}$. Write $n^* = n_{\beta_4 + \beta_5 + \beta_6} n_{\beta_5 + \beta_6 + \beta_7} n_{\beta_6 + \beta_7 + \beta_8} h_{\beta_2}({\eta_8}^3) h_{\beta_3}(\eta_8) h_{\beta_6}({\eta_8}^6) \in N$, and then $(n^*)^2 = zs$ where
$$
s = h_{\beta_2}(\eta_4) h_{\beta_3}(\eta_4) h_{\beta_5}(\eta_4).h_{\beta_4}(\eta_4).h_{{\beta_2}'}(\eta_4) h_{{\beta_3}'}(\eta_4) h_{\beta_8}(\eta_4).h_{\beta_7}(\eta_4) \in A \cap T,
$$
and conjugation by $n^*$ interchanges the two factors of $A$; then $n^*.e_{\gamma_1} = \eta_8 e_{\gamma_1}$, $n^*.e_{\gamma_2} = \eta_8 e_{\gamma_2}$, $n^*.e_{\gamma_3} = \eta_8 e_{\gamma_4}$ and $n^*.e_{\gamma_4} = \eta_8 e_{\gamma_3}$. Set $C = A$ or $A \langle n^* \rangle$ according as $p \geq 3$ or $p = 2$, and $C' = A \langle n^* \rangle$. Clearly we then have $C \leq C_G(y_0)$ and $C' \leq C_G(\langle y_0 \rangle)$; we shall show that in fact $C_G(y_0) = C$ and $C_G(\langle y_0 \rangle) = C'$.

By Lemma~\ref{lem: gen height zero}, if we take $g \in \Tran_G(y_0, Y)$ and set $y' = g.y_0 \in Y$, then we have $g = u_1 n u_2$ with $u_1 \in C_U(y')$, $u_2 \in C_U(y_0)$, and $n \in N_{\Lambda(V)_{[0]}}$ with $n.y_0 = y'$. In particular $G.y_0 \cap Y = N_{\Lambda(V)_{[0]}}.y_0 \cap Y$, and $C_G(y_0) = C_U(y_0) C_{N_{\Lambda(V)_{[0]}}}(y_0) C_U(y_0)$ while $C_G(\langle y_0 \rangle) = C_U(y_0) C_{N_{\Lambda(V)_{[0]}}}(\langle y_0 \rangle) C_U(y_0)$.

First, from the above the elements of $W$ which preserve $\Lambda(V)_{[0]}$ are those corresponding to elements of $A \langle n^* \rangle \cap N$; so we have $N_{\Lambda(V)_{[0]}}.y_0 = T.y_0 \cup n^*T.y_0$. Since any element of $T$ may be written as $h_{\beta_2}(\kappa_2) h_{\beta_3}(\kappa_3) h_{\beta_6}(\kappa_6) t$ where $\kappa_2, \kappa_3, \kappa_6 \in K^*$ and $t \in A \cap T$, by the above we have
\begin{eqnarray*}
T.y_0    & = & \left\{ \kappa_2 \kappa_6 e_{\gamma_1} + {\ts\frac{1}{\kappa_3 \kappa_6}} e_{\gamma_2} + {\ts\frac{1}{\kappa_2}} e_{\gamma_3} + \kappa_3 e_{\gamma_4} : \kappa_2, \kappa_3, \kappa_6 \in K^* \right\}, \\
n^*T.y_0 & = & \left\{ \eta_8(\kappa_2 \kappa_6 e_{\gamma_1} + {\ts\frac{1}{\kappa_3 \kappa_6}} e_{\gamma_2} + {\ts\frac{1}{\kappa_2}} e_{\gamma_4} + \kappa_3 e_{\gamma_3}) : \kappa_2, \kappa_3, \kappa_6 \in K^* \right\}.
\end{eqnarray*}
Hence $C_{N_{\Lambda(V)_{[0]}}}(y_0) = C \cap N$; also $N_{\Lambda(V)_{[0]}}.y_0 \subseteq \hat Y$, and $N_{\Lambda(V)_{[0]}}.y_0 \cap \langle y_0 \rangle = \{ {\eta_8}^i y_0 : i = 0, 1, \dots, 7 \} = \langle n^* \rangle.y_0$, so $C_{N_{\Lambda(V)_{[0]}}}(\langle y_0 \rangle) = C' \cap N$.

Next, take the ${A_3}^2$ subsystem $\Psi$ of $\Phi$ with one factor having simple roots $\beta_4$, $\beta_5$ and $\beta_2 + \beta_3 + \beta_4$, and the other factor having simple roots $\beta_7$, $\beta_8$ and $\beta_2 + \beta_3 + 2\beta_4 + 2\beta_5 + 2\beta_6 + \beta_7$; then each of the long root subgroups in $C$ is $X_\alpha$ for some $\alpha \in \Psi$, and each of the short root subgroups in $C$ is diagonally embedded in $X_\alpha X_{\alpha'} X_{\alpha''}$ for some $\alpha \in \Psi$ and $\alpha', \alpha'' \notin \Psi$. Therefore let $\Xi = \Phi^+ \setminus \Psi$, and set $U' = \prod_{\alpha \in \Xi} X_\alpha$; then $U = U'.(C \cap U)$ and $U' \cap (C \cap U) = \{ 1 \}$. We now observe that if $\alpha \in \Xi$ then $\nu_i + \alpha$ is a weight in $V$ for exactly one value of $i$; moreover each weight in $V$ of positive generalized height is of the form $\nu_i + \alpha$ for exactly one such root $\alpha$. Thus if we take $u = \prod x_\alpha(t_\alpha) \in U'$ satisfying $u.y_0 = y_0$, and equate coefficients of weight vectors, taking them in an order compatible with increasing generalized height, we see that for all $\alpha$ we must have $t_\alpha = 0$, so that $u = 1$; so $C_U(y_0) = C \cap U$.

Thus $C_U(y_0), C_{N_{\Lambda(V)_{[0]}}}(y_0) \leq C$ and $C_{N_{\Lambda(V)_{[0]}}}(\langle y_0 \rangle) \leq C'$, so we do indeed have $C_G(y_0) = C$ and $C_G(\langle y_0 \rangle) = C'$.  Moreover $G.y_0 \cap Y = \{ b_1 e_{\gamma_1} + b_2 e_{\gamma_2} + b_3 e_{\gamma_3} + b_4 e_{\gamma_4} : (b_1b_2b_3b_4)^2 = 1 \}$.

Take $y = a_1 e_{\gamma_1} + a_2 e_{\gamma_2} + a_3 e_{\gamma_3} + a_4 e_{\gamma_4} \in \hat Y$. By the above, if we choose $\kappa \in K^*$ satisfying $\kappa^4 = a_1a_2a_3a_4$, then $\kappa^{-1}y \in T.y_0$, so there exists $h \in T$ with $h.y_0 = \kappa^{-1}y$; so $C_G(y) = C_G(\kappa^{-1}y) = C_G(h.y_0) = {}^h C$ and likewise $C_G(\langle y \rangle) = {}^h C'$. Moreover, we see that $G.y \cap Y = G.h.\kappa y_0 \cap Y = \kappa(G.y_0 \cap Y) = \{ b_1 e_{\gamma_1} + b_2 e_{\gamma_2} + b_3 e_{\gamma_3} + b_4 e_{\gamma_4} : (b_1b_2b_3b_4)^2 = (a_1a_2a_3a_4)^2 \}$. Since $\dim C = 28$, we have $\dim(\overline{G.y}) = \dim G - \dim C = 91 - 28 = 63$, while $\dim(\overline{G.y \cap Y}) = 3$; therefore
$$
\dim V - \dim(\overline{G.y}) = 64 - 63 = 1 \quad \hbox{and} \quad \dim Y - \dim(\overline{G.y \cap Y}) = 4 - 3 = 1.
$$
Hence $y$ is $Y$-exact. Thus the conditions of Lemma~\ref{lem: generic stabilizer from exact subset} hold; so the triple $(G, \lambda, p)$ has generic stabilizer $C/G_V \cong {G_2}^2.\Z_{(p, 2)}$, while the quadruple $(G, \lambda, p, 1)$ has generic stabilizer $C'/Z(G) \cong {G_2}^2.\Z_2$.

Before continuing, we note that with $y \in \hat Y$ and $h$ as above the short simple root groups in $C_G(\langle y \rangle) = {}^h C'$ are $\{ x_{\beta_2}(-\frac{a_1}{a_3}t) x_{\beta_3}(-\frac{a_4}{a_2}t) x_{\beta_5}(t) : t \in K \}$ and $\{ x_{{\beta_2}'}(-\frac{a_1}{a_4}t) x_{{\beta_3}'}(\frac{a_3}{a_2}t) x_{\beta_8}(t) : t \in K \}$. Moreover $C_G(\langle y \rangle) \cap T = T_A$ is a $4$-dimensional torus, and the intersection of $N$ with $({}^h C')^\circ = {}^h A$ comprises cosets $n_1n_2T_A$, with $n_1$ and $n_2$ corresponding to Weyl group elements arising from the first and second $G_2$ factors respectively; for $i = 1, 2$ we may write $n_i \in \{ {n_i}'', {n_i}' {n_i}'' \}$, with ${n_1}' = h_{\beta_2}(-\frac{a_1}{a_3}) h_{\beta_3}(-\frac{a_4}{a_2}) n_{\beta_2} n_{\beta_3} n_{\beta_5}$ and ${n_2}' = h_{{\beta_2}'}(-\frac{a_1}{a_4}) h_{{\beta_3}'}(\frac{a_3}{a_2}) n_{{\beta_2}'} n_{{\beta_3}'} n_{\beta_8}$, and ${n_1}''$ and ${n_2}''$ lying in $\langle n_{\beta_4}, n_{\beta_2 + \beta_3 + \beta_4 + \beta_5} \rangle$ and $\langle n_{\beta_7}, n_{{\beta_2}' + {\beta_3}' + \beta_7 + \beta_8} \rangle$ respectively. Note that we have ${n_1}' {n_2}' = h_{\beta_2}(-\frac{a_1}{a_2}) h_{\beta_3}(\frac{a_1}{a_2}) h_{\beta_4 + \beta_5 + \beta_6}(-\frac{a_1a_3}{a_2a_4}) n_{\beta_2} n_{\beta_3} n_{\beta_5} n_{{\beta_2}'} n_{{\beta_3}'} n_{\beta_8}$. Write $n^{**} = n_{\beta_2 + \beta_3 + \beta_4 + \beta_5 + \beta_6} n_{\beta_6 + \beta_7} n_{\beta_5 + \beta_6 + \beta_7 + \beta_8} h_{\beta_2}({\eta_8}^7) h_{\beta_3}(\eta_8) h_{\beta_6}({\eta_8}^2) \in N$, and then $(n^{**})^2 = zs'$ where
$$
s' = h_{\beta_2}({\eta_4}^3) h_{\beta_3}({\eta_4}^3) h_{\beta_5}({\eta_4}^3).h_{\beta_4}({\eta_4}^3).h_{\beta_7}({\eta_4}^3) \in A \cap T,
$$
and conjugation by $n^{**}$ interchanges the two factors of $A$; then $n^{**}.e_{\gamma_1} = \eta_8 e_{\gamma_2}$, $n^{**}.e_{\gamma_2} = \eta_8 e_{\gamma_1}$, $n^{**}.e_{\gamma_3} = \eta_8 e_{\gamma_3}$ and $n^{**}.e_{\gamma_4} = \eta_8 e_{\gamma_4}$.

To treat the case where $G = B_6$ and $\lambda = \omega_6$ we leave $V$ and $Y$ unchanged, but replace $G$ by the $B_6$ subgroup of $D_7$ having long simple roots $\beta_8$, $\beta_7$, $\beta_6$, $\beta_5$ and $\beta_4$ and short simple root group $\{ x_{\beta_2}(t) x_{\beta_3}(t) : t \in K \}$. Here we have $Z(G) = \langle z^2 \rangle$ where $z$ is as above, so $G_V = \{ 1 \}$. We replace $\hat Y$ by
$$
\{ a_1 e_{\gamma_1} + a_2 e_{\gamma_2} + a_3 e_{\gamma_3} + a_4 e_{\gamma_4} : a_1a_2a_3a_4 \neq 0, \ (a_1a_2)^2 \neq (a_3a_4)^2 \}.
$$
For $y \in \hat Y$ we then see that the short root subgroups of the previous paragraph now meet $G$ trivially; moreover $G$ contains neither ${n_1}'$ nor ${n_2}'$, and contains ${n_1}' {n_2}'$ only if $p = 2$. We therefore replace $A$ by the ${A_2}^2$ subgroup with one factor having simple roots $\beta_4$ and $\beta_2 + \beta_3 + \beta_4 + \beta_5$, and the other factor having simple roots $\beta_7$ and $\beta_2 + \beta_3 + 2\beta_4 + 2\beta_5 + 2\beta_6 + \beta_7 + \beta_8$; since $Z(A) = \langle {z_1}', {z_2}' \rangle$ where ${z_1}' = h_{\beta_2}(\eta_3) h_{\beta_3}(\eta_3) h_{\beta_5}(\eta_3)$ and ${z_2}' = h_{\beta_2}(\eta_3) h_{\beta_3}(\eta_3) h_{\beta_4}({\eta_3}^2) h_{\beta_5}({\eta_3}^2) h_{\beta_6}({\eta_3}^2) h_{\beta_8}(\eta_3)$, we see that $A$ is of simply connected type. We note that $n^*n^{**} \in G$, and conjugation by $n^*n^{**}$ acts as a graph automorphism on each factor (and $n^*n^{**}$ and ${n_1}'{n_2}'$ correspond to the same Weyl group element); but $n^*, n^{**} \in G$ only if $p = 2$. Thus according as $p \geq 3$ or $p = 2$ set $C = A$ or $A \langle n^*, n^{**} \rangle$, and $C' = A \langle n^*n^{**} \rangle$ or $A \langle n^*, n^{**} \rangle$. Take $y \in \hat Y$ as above. Here we cannot take $h$ as above as it must now lie in the torus of $B_6$ rather than that of $D_7$, but if instead we take $\kappa, \kappa' \in K^*$ satisfying $\kappa^2 = \frac{a_3}{a_4}$ and ${\kappa'}^2 = \frac{a_2}{a_1}$, and set $h = h_{\beta_6}(\frac{1}{\kappa\kappa'}) h_{\beta_8}(\frac{1}{\kappa})$, then we see that ${}^h n^*.y = {}^h n^{**}.y = \eta_8 y$; since ${}^h A = A$, we have $C_G(y) = {}^h C$ and $C_G(\langle y \rangle) = {}^h C'$. Here we have $G.y \cap Y = N.y \cap Y = T.y \cup n^*T.y \cup n^{**}T.y \cup n^*n^{**}T.y$. This time any element of $T$ may be written as $h_{\beta_2}(\kappa_2) h_{\beta_3}(\kappa_2) h_{\beta_6}(\kappa_6) t$ where $\kappa_2, \kappa_6 \in K^*$ and $t \in A \cap T$, so we have
\begin{eqnarray*}
T.y          & = & \left\{ \kappa_2 \kappa_6 a_1 e_{\gamma_1} + {\ts\frac{1}{\kappa_2 \kappa_6}} a_2 e_{\gamma_2} + {\ts\frac{1}{\kappa_2}} a_3 e_{\gamma_3} + \kappa_2 a_4 e_{\gamma_4} : \kappa_2, \kappa_6 \in K^* \right\}, \\
n^*T.y       & = & \left\{ \eta_8(\kappa_2 \kappa_6 a_1 e_{\gamma_1} + {\ts\frac{1}{\kappa_2 \kappa_6}} a_2 e_{\gamma_2} + {\ts\frac{1}{\kappa_2}} a_3 e_{\gamma_4} + \kappa_2 a_4 e_{\gamma_3}) : \kappa_2, \kappa_6 \in K^* \right\}, \\
n^{**}T.y    & = & \left\{ \eta_8(\kappa_2 \kappa_6 a_1 e_{\gamma_2} + {\ts\frac{1}{\kappa_2 \kappa_6}} a_2 e_{\gamma_1} + {\ts\frac{1}{\kappa_2}} a_3 e_{\gamma_3} + \kappa_2 a_4 e_{\gamma_4}) : \kappa_2, \kappa_6 \in K^* \right\}, \\
n^*n^{**}T.y & = & \left\{ \eta_4(\kappa_2 \kappa_6 a_1 e_{\gamma_2} + {\ts\frac{1}{\kappa_2 \kappa_6}} a_2 e_{\gamma_1} + {\ts\frac{1}{\kappa_2}} a_3 e_{\gamma_4} + \kappa_2 a_4 e_{\gamma_3}) : \kappa_2, \kappa_6 \in K^* \right\}.
\end{eqnarray*}
Hence $G.y \cap Y = \{ b_1 e_{\gamma_1} + b_2 e_{\gamma_2} + b_3 e_{\gamma_3} + b_4 e_{\gamma_4} : (b_1b_2)^2 = (a_1a_2)^2, \ (b_3b_4)^2 = (a_3a_4)^2 \}$. Since $\dim C = 16$, we have $\dim(\overline{G.y}) = \dim G - \dim C = 78 - 16 = 62$, while $\dim(\overline{G.y \cap Y}) = 2$; therefore
$$
\dim V - \dim(\overline{G.y}) = 64 - 62 = 2 \quad \hbox{and} \quad \dim Y - \dim(\overline{G.y \cap Y}) = 4 - 2 = 2.
$$
Hence $y$ is $Y$-exact. Thus the conditions of Lemma~\ref{lem: generic stabilizer from exact subset} hold; so the triple $(G, \lambda, p)$ has generic stabilizer $C/G_V \cong {A_2}^2.{\Z_{(p, 2)}}^2$, while the quadruple $(G, \lambda, p, 1)$ has generic stabilizer $C'/Z(G) \cong {A_2}^2.\Z_{(p, 2)}.\Z_2$, where the ${A_2}^2$ is of simply connected type.
\end{proof}

\begin{prop}\label{prop: B_3, omega_3 module}
Let $G = B_3$ and $\lambda = \omega_3$. Then the triple $(G, \lambda, p)$ and the associated first quadruple $(G, \lambda, p, 1)$ both have generic stabilizer $G_2$.
\end{prop}

\begin{proof}
Take $H$ to be the (simply connected) group defined over $K$ of type $F_4$, with simple roots $\beta_1, \beta_2, \beta_3, \beta_4$. Let $G$ have simple roots $\alpha_i = \beta_i$ for $i \leq 3$, so that $G = \langle X_\alpha : \alpha = \sum m_i \beta_i, \ m_4 = 0 \rangle < H$; then we may take $V = \langle e_\alpha : \alpha = \sum m_i \beta_i, \ m_4 = 1 \rangle < \L(H)$. Note that $Z(G) = \langle z \rangle$ where $z = h_{\beta_3}(-1)$; since $z$ acts on $V$ as multiplication by $-1$, we have $G_V = \{ 1 \}$.

We take the strictly positive generalized height function on the weight lattice of $G$ whose value at each simple root $\alpha_i$ is $1$; then the generalized height of $\lambda = \frac{1}{2}(\alpha_1 + 2\alpha_2 + 3\alpha_3)$ is $3$, and as $\lambda$ and $\Phi$ generate the weight lattice it follows that the generalized height of any weight is an integer. Since $V_\lambda = \langle e_\delta \rangle$ where $\delta = \ffourrt1231$, we see that if $\mu \in \Lambda(V)$ and $e_\alpha \in V_\mu$ where $\alpha = \sum m_i \beta_i$ with $m_4 = 1$, then the generalized height of $\mu$ is $\sum_{i = 1}^3 m_i - 3$. Thus $\Lambda(V)_{[0]} = \{ \nu_1, \nu_2 \}$, where we write
$$
\gamma_1 = \ffourrt1111, \quad \gamma_2 = \ffourrt0121,
$$
and for each $i$ we let $\nu_i$ be the weight such that $V_{\nu_i} = \langle e_{\gamma_i} \rangle$. Observe that if we take $s = \prod_{i = 1}^3 h_{\beta_i}(\kappa_i) \in T$, then $\nu_1(s) = \frac{\kappa_1}{\kappa_3}$ and $\nu_2(s) = \frac{\kappa_3}{\kappa_1}$; thus $\nu_1 + \nu_2 = 0$, and so $\Lambda(V)_{[0]}$ has ZLC. Set $Y = V_{[0]} = \langle e_{\gamma_1}, e_{\gamma_2} \rangle$, and
$$
\hat Y = \{ a_1 e_{\gamma_1} + a_2 e_{\gamma_2} : a_1a_2 \neq 0 \},
$$
so that $\hat Y$ is a dense open subset of $Y$. Write
$$
y_0 = e_{\gamma_1} + e_{\gamma_2} \in \hat Y.
$$

Note that $W$ acts transitively on the set $\Sigma$ of roots $\alpha$ of $H$ corresponding to the root vectors $e_\alpha$ spanning $V$. Thus if we write $W_1$ for the stabilizer in $W$ of $\gamma_1$, then $|W_1| = \frac{|W|}{|\Sigma|} = \frac{|W|}{8} = |W(A_2)|$; we then see that $W_1 = \langle w_{\beta_2}, w_{\beta_1 + \beta_2 + 2\beta_3} \rangle$. As each generator of $W_1$ stabilizes $\gamma_2$, the pointwise stabilizer in $W$ of $\{ \gamma_1, \gamma_2 \}$ is $W_1$. Now write $w^* = w_{\beta_1} w_{\beta_3}$; as $w^*$ interchanges $\gamma_1$ and $\gamma_2$, the setwise stabilizer in $W$ of $\{ \gamma_1, \gamma_2 \}$, and hence of $\Lambda(V)_{[0]}$, is $W_1\langle w^* \rangle = \langle w_{\beta_2}, w_{\beta_1} w_{\beta_3} \rangle$.

Let $A$ be the $G_2$ subgroup having simple root groups $\{ x_{\beta_1}(t) x_{\beta_3}(-t) : t \in K \}$ and $X_{\beta_2}$; set $C = A$ and $C' = Z(G) A$. Clearly we then have $C \leq C_G(y_0)$ and $C' \leq C_G(\langle y_0 \rangle)$; we shall show that in fact $C_G(y_0) = C$ and $C_G(\langle y_0 \rangle) = C'$.

By Lemma~\ref{lem: gen height zero}, if we take $g \in \Tran_G(y_0, Y)$ and set $y' = g.y_0 \in Y$, then we have $g = u_1 n u_2$ with $u_1 \in C_U(y')$, $u_2 \in C_U(y_0)$, and $n \in N_{\Lambda(V)_{[0]}}$ with $n.y_0 = y'$. In particular $G.y_0 \cap Y = N_{\Lambda(V)_{[0]}}.y_0 \cap Y$, and $C_G(y_0) = C_U(y_0) C_{N_{\Lambda(V)_{[0]}}}(y_0) C_U(y_0)$ while $C_G(\langle y_0 \rangle) = C_U(y_0) C_{N_{\Lambda(V)_{[0]}}}(\langle y_0 \rangle) C_U(y_0)$.

First, from the above the elements of $W$ which preserve $\Lambda(V)_{[0]}$ are those corresponding to elements of $C \cap N$; so we have $N_{\Lambda(V)_{[0]}}.y_0 = T.y_0$. Since any element of $T$ may be written as $h_{\beta_3}(\kappa_3) t$ where $\kappa_3 \in K^*$ and $t \in A \cap T$, by the above we have
$$
T.y_0 = \left\{ {\ts\frac{1}{\kappa_3}} e_{\gamma_1} + \kappa_3 e_{\gamma_2} : \kappa_3 \in K^* \right\}.
$$
Hence $C_{N_{\Lambda(V)_{[0]}}}(y_0) = C \cap N$; also $N_{\Lambda(V)_{[0]}}.y_0 \subseteq \hat Y$, and $N_{\Lambda(V)_{[0]}}.y_0 \cap \langle y_0 \rangle = \{ \pm y_0 \} = Z(G).y_0$, so $C_{N_{\Lambda(V)_{[0]}}}(\langle y_0 \rangle) = C' \cap N$.

Next, each of the long root subgroups in $C$ is $X_\alpha$ for some $\alpha \in \Phi_l$, and each of the short root subgroups in $C$ is diagonally embedded in $X_\alpha X_{\alpha'}$ for some $\alpha \in \Phi_l$ and $\alpha' \notin \Phi_l$. Therefore let $\Xi = \Phi^+ \setminus \Phi_l$, and set $U' = \prod_{\alpha \in \Xi} X_\alpha$; then $U = U'.(C \cap U)$ and $U' \cap (C \cap U) = \{ 1 \}$. We now observe that if $\alpha \in \Xi$ then $\nu_i + \alpha$ is a weight in $V$ for exactly one value of $i$; moreover each weight in $V$ of positive generalized height is of the form $\nu_i + \alpha$ for exactly one such root $\alpha$. Thus if we take $u = \prod x_\alpha(t_\alpha) \in U'$ satisfying $u.y_0 = y_0$, and equate coefficients of weight vectors, taking them in an order compatible with increasing generalized height, we see that for all $\alpha$ we must have $t_\alpha = 0$, so that $u = 1$; so $C_U(y_0) = C \cap U$.

Thus $C_U(y_0), C_{N_{\Lambda(V)_{[0]}}}(y_0) \leq C$ and $C_{N_{\Lambda(V)_{[0]}}}(\langle y_0 \rangle) \leq C'$, so we do indeed have $C_G(y_0) = C$ and $C_G(\langle y_0 \rangle) = C'$. Moreover $G.y_0 \cap Y = \{ b_1 e_{\gamma_1} + b_2 e_{\gamma_2} : b_1b_2 = 1 \}$.

Take $y = a_1 e_{\gamma_1} + a_2 e_{\gamma_2} \in \hat Y$. By the above, if we choose $\kappa \in K^*$ satisfying $\kappa^2 = a_1a_2$, then $\kappa^{-1}y \in T.y_0$, so there exists $h \in T$ with $h.y_0 = \kappa^{-1} y$; so $C_G(y) = C_G(\kappa^{-1}y) = C_G(h.y_0) = {}^h C$ and likewise $C_G(\langle y \rangle) = {}^h C'$. Moreover, we see that $G.y \cap Y = G.h.\kappa y_0 \cap Y = \kappa(G.y_0 \cap Y) = \{ b_1 e_{\gamma_1} + b_2 e_{\gamma_2} : b_1b_2 = a_1a_2 \}$. Since $\dim C = 14$, we have $\dim(\overline{G.y}) = \dim G - \dim C = 21 - 14 = 7$, while $\dim(\overline{G.y \cap Y}) = 1$; therefore
$$
\dim V - \dim(\overline{G.y}) = 8 - 7 = 1 \quad \hbox{and} \quad \dim Y - \dim(\overline{G.y \cap Y}) = 2 - 1 = 1.
$$
Hence $y$ is $Y$-exact. Thus the conditions of Lemma~\ref{lem: generic stabilizer from exact subset} hold; so the triple $(G, \lambda, p)$ has generic stabilizer $C/G_V \cong G_2$, while the quadruple $(G, \lambda, p, 1)$ has generic stabilizer $C'/Z(G) \cong G_2$.
\end{proof}

\begin{prop}\label{prop: C_3, omega_3, C_4, omega_4, C_5, omega_5, C_6, omega_6 modules, p = 2}
Let $G = C_3$ and $\lambda = \omega_3$, or $G = C_4$ and $\lambda = \omega_4$, or $G = C_5$ and $\lambda = \omega_5$, or $G = C_6$ and $\lambda = \omega_6$, with $p = 2$. Then the triple $(G, \lambda, p)$ and the associated first quadruple $(G, \lambda, p, 1)$ both have generic stabilizer $G_2$, or $C_3$, or $\tilde A_4.\Z_2$, or ${\tilde A_2}{}^2.{\Z_2}^2$, respectively.
\end{prop}

\begin{proof}
This is an immediate consequence of Propositions~\ref{prop: B_3, omega_3 module}, \ref{prop: D_5, omega_5, B_4, omega_4 modules}, \ref{prop: E_7, omega_7, D_6, omega_6, B_5, omega_5, A_5, omega_3, C_3, omega_3 modules} and \ref{prop: D_7, omega_7, B_6, omega_6 modules} respectively, using the exceptional isogeny $B_\ell \to C_\ell$ which exists in characteristic $2$.
\end{proof}

\begin{prop}\label{prop: E_6, omega_1, F_4, omega_4 modules}
Let $G = E_6$ and $\lambda = \omega_1$, or $G = F_4$ and $\lambda = \omega_4$. Then the triple $(G, \lambda, p)$ has generic stabilizer $F_4$ or $D_4.\Z_{(p, 3)}$ respectively; the associated first quadruple $(G, \lambda, p, 1)$ has generic stabilizer $F_4$ or $D_4.S_{(p, 3)}$ respectively.
\end{prop}

\begin{proof}
Throughout this proof we take $H$ to be the simply connected group defined over $K$ of type $E_7$, with simple roots $\beta_1, \dots, \beta_7$.

We begin with the case where $G = E_6$ and $\lambda = \omega_1$. Let $G$ have simple roots $\alpha_i = \beta_i$ for $i \leq 6$, so that $G = \langle X_\alpha : \alpha = \sum m_i \beta_i, \ m_7 = 0 \rangle < H$; then we may take $V = \langle e_\alpha : \alpha = \sum m_i \beta_i, \ m_7 = 1 \rangle < \L(H)$. Note that $Z(G) = \langle z \rangle$ where $z = h_{\beta_1}(\eta_3) h_{\beta_3}({\eta_3}^2) h_{\beta_5}(\eta_3) h_{\beta_6}({\eta_3}^2)$; since $z$ acts on $V$ as multiplication by $\eta_3$, we have $G_V = \{ 1 \}$.

We take the strictly positive generalized height function on the weight lattice of $G$ whose value at each simple root $\alpha_i$ is $1$; then the generalized height of $\lambda = \frac{1}{3}(4\alpha_1 + 3\alpha_2 + 5\alpha_3 + 6\alpha_4 + 4\alpha_5 + 2\alpha_6)$ is $8$, and as $\lambda$ and $\Phi$ generate the weight lattice it follows that the generalized height of any weight is an integer. Since $V_\lambda = \langle e_\delta \rangle$ where $\delta = \esevenrt2234321$, we see that if $\mu \in \Lambda(V)$ and $e_\alpha \in V_\mu$ where $\alpha = \sum m_i \beta_i$ with $m_7 = 1$, then the generalized height of $\mu$ is $\sum_{i = 1}^6 m_i - 8$. Thus $\Lambda(V)_{[0]} = \{ \nu_1, \nu_2, \nu_3 \}$, where we write
$$
\gamma_1 = \esevenrt1122111, \quad \gamma_2 = \esevenrt1112211, \quad \gamma_3 = \esevenrt0112221,
$$
and for each $i$ we let $\nu_i$ be the weight such that $V_{\nu_i} = \langle e_{\gamma_i} \rangle$. Observe that if we take $s = \prod_{i=1}^6 h_{\beta_i}(\kappa_i) \in T$, then $\nu_1(s) = \frac{\kappa_3}{\kappa_5}$, $\nu_2(s) = \frac{\kappa_1 \kappa_5}{\kappa_3 \kappa_6}$ and $\nu_3(s) = \frac{\kappa_6}{\kappa_1}$; thus $\nu_1 + \nu_2 + \nu_3 = 0$, and so $\Lambda(V)_{[0]}$ has ZLC. Set $Y = V_{[0]} = \langle e_{\gamma_1}, e_{\gamma_2}, e_{\gamma_3} \rangle$,
and
$$
\hat Y = \{ a_1 e_{\gamma_1} + a_2 e_{\gamma_2} + a_3 e_{\gamma_3} : a_1a_2a_3 \neq 0 \},
$$
so that $\hat Y$ is a dense open subset of $Y$. Write
$$
y_0 = e_{\gamma_1} + e_{\gamma_2} + e_{\gamma_3} \in \hat Y.
$$

Note that $W$ acts transitively on the set $\Sigma$ of roots $\alpha$ of $H$ corresponding to the root vectors $e_\alpha$ spanning $V$. Thus if we write $W_1$ for the stabilizer in $W$ of $\gamma_1$, then $|W_1| = \frac{|W|}{|\Sigma|} = \frac{|W|}{27} = |W(D_5)|$; we then see that $W_1 = \langle w_{\beta_4}, w_{\beta_2}, w_{\beta_3 + \beta_4 + \beta_5}, w_{\beta_1},$ $w_{\beta_6} \rangle$. Now the stabilizer in $W$ of any $\alpha \in \Sigma$ acts transitively on the set $\Sigma'$ of roots $\alpha' \in \Sigma$ orthogonal to $\alpha$ (this is evident if we take $\alpha = \esevenrt0000001$, as then its stabilizer in $W$ is $\langle w_{\beta_1}, w_{\beta_3}, w_{\beta_4}, w_{\beta_5}, w_{\beta_2} \rangle$, which acts transitively on the set of roots $\alpha' = \sum m_i \beta_i$ with $m_6 = 2$ and $m_7 = 1$). Thus if we write $W_2$ for the stabilizer in $W_1$ of $\gamma_2$, then $|W_2| = \frac{|W_1|}{|\Sigma'|} = \frac{|W_1|}{10} = |W(D_4)|$; we then see that $W_2 = \langle w_{\beta_4}, w_{\beta_2}, w_{\beta_3 + \beta_4 + \beta_5}, w_{\beta_1 + \beta_3 + \beta_4 + \beta_5 + \beta_6} \rangle$. As each generator of $W_2$ stabilizes $\gamma_3$, the pointwise stabilizer in $W$ of $\{ \gamma_1, \gamma_2, \gamma_3 \}$ is $W_2$. Now write $w^* = w_{\beta_3} w_{\beta_5}$ and $w^{**} = w_{\beta_1} w_{\beta_6}$; as $w^*$ interchanges $\gamma_1$ and $\gamma_2$ while fixing $\gamma_3$, and $w^{**}$ interchanges $\gamma_2$ and $\gamma_3$ while fixing $\gamma_1$, the setwise stabilizer in $W$ of $\{ \gamma_1, \gamma_2, \gamma_3 \}$, and hence of $\Lambda(V)_{[0]}$, is $W_2\langle w^*, w^{**} \rangle = \langle w_{\beta_2}, w_{\beta_4}, w_{\beta_3} w_{\beta_5}, w_{\beta_1} w_{\beta_6} \rangle$.

Let $A$ be the $F_4$ subgroup having long simple roots $\beta_2$ and $\beta_4$ and short simple root groups $\{ x_{\beta_3}(t) x_{\beta_5}(-t) : t \in K \}$ and $\{ x_{\beta_1}(t) x_{\beta_6}(-t) : t \in K \}$, and set $C = A$ and $C' = Z(G) A$. Clearly we then have $C \leq C_G(y_0)$ and $C' \leq C_G(\langle y_0 \rangle)$; we shall show that in fact $C_G(y_0) = C$ and $C_G(\langle y_0 \rangle) = C'$.

By Lemma~\ref{lem: gen height zero}, if we take $g \in \Tran_G(y_0, Y)$ and set $y' = g.y_0 \in Y$, then we have $g = u_1 n u_2$ with $u_1 \in C_U(y')$, $u_2 \in C_U(y_0)$, and $n \in N_{\Lambda(V)_{[0]}}$ with $n.y_0 = y'$. In particular $G.y_0 \cap Y = N_{\Lambda(V)_{[0]}}.y_0 \cap Y$, and $C_G(y_0) = C_U(y_0) C_{N_{\Lambda(V)_{[0]}}}(y_0) C_U(y_0)$ while $C_G(\langle y_0 \rangle) = C_U(y_0) C_{N_{\Lambda(V)_{[0]}}}(\langle y_0 \rangle) C_U(y_0)$.

First, from the above the elements of $W$ which preserve $\Lambda(V)_{[0]}$ are those corresponding to elements of $C \cap N$; so we have $N_{\Lambda(V)_{[0]}}.y_0 = T.y_0$. Since any element of $T$ may be written as $h_{\beta_1}(\kappa_1) h_{\beta_3}(\kappa_3) t$ where $\kappa_1, \kappa_3 \in K^*$ and $t \in C \cap T$, by the above we have
$$
T.y_0 = \left\{ \kappa_3 e_{\gamma_1} + {\ts\frac{\kappa_1}{\kappa_3}} e_{\gamma_2} + {\ts\frac{1}{\kappa_1}} e_{\gamma_3} : \kappa_1, \kappa_3 \in K^* \right\}.
$$
Hence $C_{N_{\Lambda(V)_{[0]}}}(y_0) = C \cap N$; also $N_{\Lambda(V)_{[0]}}.y_0 \subseteq \hat Y$, and $N_{\Lambda(V)_{[0]}}.y_0 \cap \langle y_0 \rangle = \{ {\eta_3}^i y_0 : i = 0, 1, 2 \} = Z(G).y_0$, so $C_{N_{\Lambda(V)_{[0]}}}(\langle y_0 \rangle) = C' \cap N$.

Next, take the subset $\Xi' = \{ \alpha \in \Phi : \alpha = \sum m_i \beta_i, \ m_1 + m_3 \leq m_5 + m_6 \}$ of $\Phi$; then each of the long root subgroups in $C$ is $X_\alpha$ for some $\alpha \in \Xi'$, and each of the short root subgroups in $C$ is diagonally embedded in $X_\alpha X_{\alpha'}$ for some $\alpha \in \Xi'$ and $\alpha' \notin \Xi'$. Therefore let $\Xi = \Phi^+ \setminus \Xi'$, and set $U' = \prod_{\alpha \in \Xi} X_\alpha$; then $U = U'.(C \cap U)$ and $U' \cap (C \cap U) = \{ 1 \}$. We now observe that if $\alpha \in \Xi$ then $\nu_i + \alpha$ is a weight in $V$ for exactly one value of $i$; moreover each weight in $V$ of positive generalized height is of the form $\nu_i + \alpha$ for exactly one such root $\alpha$. Thus if we take $u = \prod x_\alpha(t_\alpha) \in U'$ satisfying $u.y_0 = y_0$, and equate coefficients of weight vectors, taking them in an order compatible with increasing generalized height, we see that for all $\alpha$ we must have $t_\alpha = 0$, so that $u = 1$; so $C_U(y_0) = C \cap U$.

Thus $C_U(y_0), C_{N_{\Lambda(V)_{[0]}}}(y_0) \leq C$ and $C_{N_{\Lambda(V)_{[0]}}}(\langle y_0 \rangle) \leq C'$, so we do indeed have $C_G(y_0) = C$ and $C_G(\langle y_0 \rangle) = C'$. Moreover $G.y_0 \cap Y = \{ b_1 e_{\gamma_1} + b_2 e_{\gamma_2} + b_3 e_{\gamma_3} : b_1b_2b_3 = 1 \}$.

Take $y = a_1 e_{\gamma_1} + a_2 e_{\gamma_2} + a_3 e_{\gamma_3} \in \hat Y$. By the above, if we choose $\kappa \in K^*$ satisfying $\kappa^3 = a_1a_2a_3$, then $\kappa^{-1}y \in T.y_0$, so there exists $h \in T$ with $h.y_0 = \kappa^{-1}y$; so $C_G(y) = C_G(\kappa^{-1}y) = C_G(h.y_0) = {}^h C$ and likewise $C_G(\langle y \rangle) = {}^h C'$. Moreover, we see that $G.y \cap Y = G.h.\kappa y_0 \cap Y = \kappa(G.y_0 \cap Y) = \{ b_1 e_{\gamma_1} + b_2 e_{\gamma_2} + b_3 e_{\gamma_3} : b_1b_2b_3 = a_1a_2a_3 \}$. Since $\dim C = 52$, we have $\dim(\overline{G.y}) = \dim G - \dim C = 78 - 52 = 26$, while $\dim(\overline{G.y \cap Y}) = 2$; therefore
$$
\dim V - \dim(\overline{G.y}) = 27 - 26 = 1 \quad \hbox{and} \quad \dim Y - \dim(\overline{G.y \cap Y}) = 3 - 2 = 1.
$$
Hence $y$ is $Y$-exact. Thus the conditions of Lemma~\ref{lem: generic stabilizer from exact subset} hold; so the triple $(G, \lambda, p)$ has generic stabilizer $C/G_V \cong F_4$, while the quadruple $(G, \lambda, p, 1)$ has generic stabilizer $C'/Z(G) \cong F_4$.

To treat the case where $G = F_4$ and $\lambda = \omega_4$, we replace $G$ by $A$; we temporarily write $V'$ for the $27$-dimensional module called $V$ above. Inside $V'$ we have the submodules $X_1 = \{ \sum a_\gamma e_\gamma \in V' : \sum_i a_{\gamma_i} = 0 \}$ and $X_2 = \langle y_0 \rangle$, with the latter being trivial. If $p \neq 3$ then $V' = X_1 \oplus X_2$, and $V = X_1$; if however $p = 3$ then $X_2 < X_1$, and $V = X_1/X_2$. Thus in all cases $V = X_1/(X_1 \cap X_2)$, where $X_1 \cap X_2$ is either zero or the trivial $G$-module. As $Z(G) = \{ 1 \}$ we have $G_V = \{ 1 \}$.

We restrict the generalized height function above to the weight lattice of $G$; then $\Lambda(V)_{[0]} = \{ 0 \}$, so trivially $\Lambda(V)_{[0]}$ has ZLC, and $N_{\Lambda(V)_{[0]}} = N$. Set
$$
Y = V_{[0]} = \left\{ {\ts\sum} a_i e_{\gamma_i} + (X_1 \cap X_2) : {\ts\sum} a_i = 0 \right\},
$$
and
$$
\hat Y = \left\{ {\ts\sum} a_i e_{\gamma_i} + (X_1 \cap X_2) : {\ts\sum} a_i = 0, \ a_1a_2a_3 \neq 0, \ ({\ts\frac{a_i}{a_j}})^3 \neq 1 \hbox{ if } i \neq j \right\},
$$
so that $\hat Y$ is a dense open subset of $Y$. Take
$$
y = {\ts\sum} a_i e_{\gamma_i} + (X_1 \cap X_2) \in \hat Y.
$$

Let $A$ be the $D_4$ subgroup generated by the long root subgroups of $G$; as $Z(A) = \langle {z_1}', {z_2}' \rangle$, where ${z_1}' = h_{\beta_4}(-1) h_{\beta_3 + \beta_4 + \beta_5}(-1)$, ${z_2}' = h_{\beta_4}(-1) h_{\beta_1 + \beta_3 + \beta_4 + \beta_5 + \beta_6}(-1)$, we see that $A$ is of simply connected type. Write $n^* = n_{\alpha_3} = n_{\beta_3} {n_{\beta_5}}^{-1}$ and $n^{**} = n_{\alpha_4} = n_{\beta_1} {n_{\beta_6}}^{-1}$; if $p \neq 3$ set $C = C' = A$, while if $p = 3$ set $C = A \langle n^* n^{**} \rangle \cong D_4.\Z_3$ and $C' = A \langle n^*, n^{**} \rangle \cong D_4.S_3$. Note that if $p = 3$ then $\dim Y = 1$, and both $n^*$ and $n^{**}$ act on $Y$ as negation since for example $n^*.(\sum a_i e_{\gamma_i}) = a_1 e_{\gamma_2} + a_2 e_{\gamma_1} + a_3 e_{\gamma_3} = -(\sum a_i e_{\gamma_i}) - a_3 y_0$ as $\sum a_i = 0$. Thus we have $C \leq C_G(y)$ and $C' \leq C_G(\langle y \rangle)$; we shall show that in fact $C_G(y) = C$ and $C_G(\langle y \rangle) = C'$.

By Lemma~\ref{lem: gen height zero}, if we take $g \in \Tran_G(y, Y)$ and set $y' = g.y \in Y$, then we have $g = u_1 n u_2$ with $u_1 \in C_U(y')$, $u_2 \in C_U(y)$, and $n \in N$ with $n.y = y'$. In particular $G.y \cap Y = N.y$, and $C_G(y) = C_U(y) C_N(y) C_U(y)$ while $C_G(\langle y \rangle) = C_U(y) C_N(\langle y \rangle) C_U(y)$.

First, we note that any element of $N$ may be written as $n's$, where $n'$ is a product of elements $n_\alpha$ for various roots $\alpha$, and $s \in T$; since $s$ stabilizes $y$, and each element $n_\alpha$ can only permute the individual vectors $e_{\gamma_i}$, we see that
$$
N.y = \left\{ {\ts\sum} a_{\pi(i)} e_{\gamma_i} + (X_1 \cap X_2) : \pi \in S_3 \right\}.
$$
Moreover in the case where $p \neq 3$, suppose $\pi \in S_3$ satisfies $\sum a_{\pi(i)} e_{\gamma_i} = \kappa \sum a_i e_{\gamma_i}$ for some $\kappa \in K^*$. If $\pi$ is a transposition, say $(1 \ 2)$, then equating coefficients of $e_{\gamma_3}$ and $e_{\gamma_1}$ gives $\kappa = 1$ and then $a_1 = a_2$; if instead $\pi$ is a $3$-cycle, say $(1 \ 2 \ 3)$, then equating coefficients gives $\kappa = \frac{a_2}{a_1} = \frac{a_3}{a_2} = \frac{a_1}{a_3}$, so $(\frac{a_2}{a_1})^3 = \kappa^3 = 1$. The definition of $\hat Y$ rules out both possibilities, so we must have $\pi = 1$. Thus the only elements of $N$ which stabilize $y$ or $\langle y \rangle$ are those in $C$ or $C'$ respectively. Hence $C_N(y) = C \cap N$ and $C_N(\langle y \rangle) = C' \cap N$.

Next, let $\Xi = \Phi_s \cap \Phi^+$, and set $U' = \prod_{\alpha \in {\Phi_s}^+} X_\alpha$; then $U = U'.(C \cap U)$ and $U' \cap (C \cap U) = \{ 1 \}$. Observe that any short root element $x_{\alpha}(t) x_{\tau(\alpha)}(\e t)$ of $G$ (where $\tau$ is the graph automorphism of $E_6$, and $\e \in \{ \pm 1 \}$) sends $\sum a_i e_{\gamma_i}$ to $\sum a_i e_{\gamma_i} + t(a_i - a_j) e_{\gamma}$ for some $i \neq j$ and some root $\gamma$, so does not stabilize $y$ unless $t = 0$; moreover distinct short roots correspond to distinct roots $\gamma$. Thus if we take $u = \prod x_\alpha(t_\alpha) \in U'$ satisfying $u.y = y$, and equate coefficients of weight vectors, taking them in an order compatible with increasing generalized height, we see that for all $\alpha$ we must have $t_\alpha = 0$, so that $u = 1$; so $C_U(y) = C \cap U$.

Thus $C_U(y), C_N(y) \leq C$ and $C_N(\langle y \rangle) \leq C'$, so we do indeed have $C_G(y) = C$ and $C_G(\langle y \rangle) = C'$.

Since $\dim C = 28$, we have $\dim(\overline{G.y}) = \dim G - \dim C = 52 - 28 = 24$, while $\dim(\overline{G.y \cap Y}) = 0$ because any $N$-orbit on $Y = V_0$ is finite; therefore
$$
\dim V - \dim(\overline{G.y}) = (26 - \z_{3, p}) - 24 = 2 - \z_{3, p}
$$
and
$$
\dim Y - \dim(\overline{G.y \cap Y}) = (2 - \z_{3, p}) - 0 = 2 - \z_{3, p}.
$$
Hence $y$ is $Y$-exact. Thus the conditions of Lemma~\ref{lem: generic stabilizer from exact subset} hold; so the triple $(G, \lambda, p)$ has generic stabilizer $C/G_V \cong D_4.\Z_{(p, 3)}$, while the quadruple $(G, \lambda, p, 1)$ has generic stabilizer $C'/Z(G) \cong D_4.S_{(p, 3)}$, where the $D_4$ is of simply connected type.
\end{proof}

\begin{prop}\label{prop: F_4, omega_1 module, p = 2}
Let $G = F_4$ and $\lambda = \omega_1$ with $p = 2$. Then the triple $(G, \lambda, p)$ and the associated first quadruple $(G, \lambda, p, 1)$ both have generic stabilizer $\tilde D_4$.
\end{prop}

\begin{proof}
This is an immediate consequence of Proposition~\ref{prop: E_6, omega_1, F_4, omega_4 modules}, using the graph automorphism of $F_4$ which exists in characteristic $2$.
\end{proof}

\begin{prop}\label{prop: G_2, omega_1 module}
Let $G = G_2$ and $\lambda = \omega_1$ with $p \geq 3$ or $p = 2$. Then the triple $(G, \lambda, p)$ has generic stabilizer $A_2$ or $A_1U_5$ respectively; the associated first quadruple $(G, \lambda, p, 1)$ has generic stabilizer $A_2.\Z_2$ or $A_1 T_1 U_5$ respectively.
\end{prop}

\begin{proof}
As $Z(G) = \{ 1 \}$ we have $G_V = \{ 1 \}$. We begin with the case where $p \geq 3$; here $\dim V = 7$ and $\Lambda(V) = \Phi_s \cup \{ 0 \}$. Take an ordered basis of $V$ consisting of weight vectors $v_\mu$ for the weights $\mu = 2\alpha_1 + \alpha_2$, $\alpha_1 + \alpha_2$, $\alpha_1$, $0$, $-\alpha_1$, $-(\alpha_1 + \alpha_2)$, $-(2\alpha_1 + \alpha_2)$ respectively, such that with respect to them the simple root elements $x_{\alpha_1}(t)$ and $x_{\alpha_2}(t)$ of $G$ act by the matrices
$$
\left(
  \begin{array}{ccccccc}
    1 & -t &   &   &     &   &    \\
      &  1 &   &   &     &   &    \\
      &    & 1 & t & t^2 &   &    \\
      &    &   & 1 & 2t  &   &    \\
      &    &   &   &  1  &   &    \\
      &    &   &   &     & 1 & -t \\
      &    &   &   &     &   &  1 \\
  \end{array}
\right)
\qquad\hbox{ and }\qquad
\left(
  \begin{array}{ccccccc}
    1 &   &   &   &   &    &   \\
      & 1 & t &   &   &    &   \\
      &   & 1 &   &   &    &   \\
      &   &   & 1 &   &    &   \\
      &   &   &   & 1 & -t &   \\
      &   &   &   &   &  1 &   \\
      &   &   &   &   &    & 1 \\
  \end{array}
\right)
$$
respectively, and the corresponding negative root elements $x_{-\alpha_1}(t)$ and $x_{-\alpha_2}(t)$ act by the matrices
$$
\left(
  \begin{array}{ccccccc}
     1 &   &     &   &   &    &   \\
    -t & 1 &     &   &   &    &   \\
       &   &  1  &   &   &    &   \\
       &   & 2t  & 1 &   &    &   \\
       &   & t^2 & t & 1 &    &   \\
       &   &     &   &   &  1 &   \\
       &   &     &   &   & -t & 1 \\
  \end{array}
\right)
\qquad\hbox{ and }\qquad
\left(
  \begin{array}{ccccccc}
    1 &   &   &   &    &   &   \\
      & 1 &   &   &    &   &   \\
      & t & 1 &   &    &   &   \\
      &   &   & 1 &    &   &   \\
      &   &   &   &  1 &   &   \\
      &   &   &   & -t & 1 &   \\
      &   &   &   &    &   & 1 \\
  \end{array}
\right)
$$
respectively. Set $Y = V_0 = \langle v_0 \rangle$; let $\hat Y = \{ av_0 : a \neq 0 \}$, then $\hat Y$ is a dense open subset of $Y$. Take $y \in \hat Y$. Let $A$ be the $A_2$ subgroup having simple roots $\alpha_2$ and $3\alpha_1 + \alpha_2$, so that $A = \langle X_\alpha : \alpha \in \Phi_l \rangle$; then as $Z(A) = \langle z' \rangle$ where $z' = h_{\alpha_1}(\eta_3)$, we see that $A$ is of simply connected type. Note that $n_{\alpha_1}.y = -y$. Set $C = A$ and $C' = A \langle n_{\alpha_1} \rangle$. Clearly we then have $C \leq C_G(y)$ and $C' \leq C_G(\langle y \rangle)$; as $C'$ is a maximal subgroup of $G$ we must in fact have $C_G(y) = C$ and $C_G(\langle y \rangle) = C'$, and $G.y \cap Y = \{ \pm y \}$. Since $\dim C = 8$, we have $\dim(\overline{G.y}) = \dim G - \dim C = 14 - 8 = 6$, while $\dim(\overline{G.y \cap Y}) = 0$; therefore
$$
\dim V - \dim(\overline{G.y}) = 7 - 6 = 1 \quad \hbox{and} \quad \dim Y - \dim(\overline{G.y \cap Y}) = 1 - 0 = 1.
$$
Hence $y$ is $Y$-exact. Thus the conditions of Lemma~\ref{lem: generic stabilizer from exact subset} hold; so the triple $(G, \lambda, p)$ has generic stabilizer $C/G_V \cong A_2$, while the quadruple $(G, \lambda, p, 1)$ has generic stabilizer $C'/Z(G) \cong A_2.\Z_2$, where the $A_2$ is of simply connected type.

Now take the case where $p = 2$; here $\dim V = 6$ and $\Lambda(V) = \Phi_s$, and we may obtain matrices for the action by deleting the fourth row and column from those above. Let $v_\lambda$ be a highest weight vector in $V$, and set $Y = V_\lambda = \langle v_\lambda \rangle$; let $\hat Y = \{ a v_\lambda : a \neq 0 \}$, then $\hat Y$ is a dense open subset of $Y$. Take $y \in \hat Y$. Set $C = \langle U, X_{-\alpha_2} \rangle$ and $C' = TC$, so that $C'$ is the standard maximal parabolic subgroup corresponding to the simple root $\alpha_1$; then $C' = C_G(\langle y \rangle)$, whence $C = C_G(y)$ and $G.y \cap Y = T.y = \{ \kappa y : \kappa \in K^* \}$. Since $\dim C = 8$, we have $\dim(\overline{G.y}) = \dim G - \dim C = 14 - 8 = 6$, while $\dim(\overline{G.y \cap Y}) = 1$; therefore
$$
\dim V - \dim(\overline{G.y}) = 6 - 6 = 0 \quad \hbox{and} \quad \dim Y - \dim(\overline{G.y \cap Y}) = 1 - 1 = 0.
$$
Hence $y$ is $Y$-exact. Thus the conditions of Lemma~\ref{lem: generic stabilizer from exact subset} hold; so the triple $(G, \lambda, p)$ has generic stabilizer $C/G_V \cong A_1 U_5$, while the quadruple $(G, \lambda, p, 1)$ has generic stabilizer $C'/Z(G) \cong A_1 T_1 U_5$.
\end{proof}

\begin{prop}\label{prop: G_2, omega_2 module, p = 3}
Let $G = G_2$ and $\lambda = \omega_2$ with $p = 3$. Then the triple $(G, \lambda, p)$ has generic stabilizer $\tilde A_2$; the associated first quadruple $(G, \lambda, p, 1)$ has generic stabilizer $\tilde A_2.\Z_2$.
\end{prop}

\begin{proof}
This is an immediate consequence of Proposition~\ref{prop: G_2, omega_1 module}, using the graph automorphism of $G_2$ which exists in characteristic $3$.
\end{proof}

\begin{prop}\label{prop: A_6, omega_3 module}
Let $G = A_6$ and $\lambda = \omega_3$. Then the triple $(G, \lambda, p)$ and the associated first quadruple $(G, \lambda, p, 1)$ both have generic stabilizer $G_2$.
\end{prop}

\begin{proof}
Take $G$ to be of simply connected type, so that $G = \SL_7(K)$. We may view $V$ as the exterior cube $\bigwedge^3(V_{nat})$ of the natural module, and identify $W$ with the symmetric group $S_7$. For convenience, for $i_1, i_2, i_3 \leq 7$ write $v_{i_1i_2i_3} = v_{i_1} \wedge v_{i_2} \wedge v_{i_3}$; then $V = \{ v_{i_1i_2i_3} : 1 \leq i_1 < i_2 < i_3 \leq 7 \}$. Note that $Z(G) = \langle z \rangle$ where $z = \prod_{i = 1}^6 h_{\alpha_i}({\eta_7}^i)$; as $z$ acts on $V_{nat}$ as multiplication by $\eta_7$, it acts on $V_{nat} \otimes V_{nat} \otimes V_{nat}$ and hence on $V$ as multiplication by ${\eta_7}^3$, so $G_V = \{ 1 \}$.

We take the strictly positive generalized height function on the weight lattice of $G$ whose value at each simple root $\alpha_i$ is $1$; then the generalized height of $\lambda = \frac{1}{7}(4\alpha_1 + 8\alpha_2 + 12\alpha_3 + 9\alpha_4 + 6\alpha_5 + 3\alpha_6)$ is $6$, and as $\lambda$ and $\Phi$ generate the weight lattice it follows that the generalized height of any weight is an integer. Since $V_\lambda = \langle v_{123} \rangle$, we see that if $\mu \in \Lambda(V)$ and $v_{i_1i_2i_3} \in V_\mu$, then the generalized height of $\mu$ is $12 - (i_1 + i_2 + i_3)$. Thus $\Lambda(V)_{[0]} = \{ \nu_1, \nu_2, \nu_3, \nu_4, \nu_5 \}$, where we write
$$
x_1 = v_{147}, \quad x_2 = v_{246}, \quad x_3 = v_{345}, \quad x_4 = v_{156}, \quad x_5 = v_{237},
$$
and for each $i$ we let $\nu_i$ be the weight such that $V_{\nu_i} = \langle x_i \rangle$. Observe that if we take $s = \prod_{i = 1}^6 h_{\alpha_i}(\kappa_i) \in T$, then $\nu_1(s) = \frac{\kappa_1 \kappa_4}{\kappa_3 \kappa_6}$, $\nu_2(s) = \frac{\kappa_2 \kappa_4 \kappa_6}{\kappa_1 \kappa_3 \kappa_5}$, $\nu_3(s) = \frac{\kappa_5}{\kappa_2}$, $\nu_4(s) = \frac{\kappa_1 \kappa_6}{\kappa_4}$ and $\nu_5(s) = \frac{\kappa_3}{\kappa_1 \kappa_6}$; thus $\nu_1 + \nu_2 + \nu_3 + 2\nu_4 + 2\nu_5 = 0$, and so $\Lambda(V)_{[0]}$ has ZLC. Set $Y = V_{[0]} = \langle x_1, x_2, x_3, x_4, x_5 \rangle$,
and
$$
\hat Y = \{ a_1 x_1 + \cdots + a_5 x_5 : a_1a_2a_3a_4a_5 \neq 0 \},
$$
so that $\hat Y$ is a dense open subset of $Y$. Write
$$
y_0 = x_1 + x_2 + x_3 + x_4 + x_5 \in \hat Y.
$$

We see that $W$ acts on $\Lambda(V)$ such that if $w \in W$ and $\mu \in \Lambda(V)$ with $v_{i_1i_2i_3} \in V_\mu$ then $v_{w(i_1) w(i_2) w(i_3)} \in V_{w(\mu)}$. The pointwise stabilizer in $W$ of $\Lambda(V)_{[0]}$ is trivial, since for each of the numbers $1, \dots, 7$ there exist $i$ and $j$ distinct such that the number concerned is the intersection of the sets of three numbers appearing as subscripts in $x_i$ and $x_j$; moreover as $4$ is the only number to appear three times as a subscript in the various $x_i$, the setwise stabilizer in $W$ must in fact fix $4$, so must preserve the subsets $\{ \nu_1, \nu_2, \nu_3 \}$ and $\{ \nu_4, \nu_5 \}$, and hence is at most $S_3 \times S_2$. Since $(1 \ 2)(3 \ 5)(6 \ 7)$ interchanges $\nu_1$ and $\nu_2$, and also $\nu_4$ and $\nu_5$, while fixing $\nu_3$, and $(2 \ 3)(5 \ 6)$ interchanges $\nu_2$ and $\nu_3$ while fixing the remaining $\nu_j$, and these two permutations generate a dihedral group of order $12$, we see that the setwise stabilizer in $W$ of $\Lambda(V)_{[0]}$ is $\langle (1 \ 2)(3 \ 5)(6 \ 7), (2 \ 3)(5 \ 6) \rangle = \langle w_{\alpha_1} w_{\alpha_3 + \alpha_4} w_{\alpha_6}, w_{\alpha_2} w_{\alpha_5} \rangle$.

Let $A$ be a (simply connected) group defined over $K$ of type $G_2$, with simple roots $\beta_1$ (short) and $\beta_2$ (long). We may define a homomorphism $\psi: A \to G$, by letting $\psi(x_{\beta_1}(t))$, $\psi(x_{\beta_2}(t))$, $\psi(x_{-\beta_1}(t))$ and $\psi(x_{-\beta_2}(t))$ for $t \in K$ be the matrices given in the proof of Proposition~\ref{prop: G_2, omega_1 module} above; set $C = \psi(A)$ and $C' = Z(G) \psi(A)$. A straightforward calculation shows that the stabilizer of $y_0$ contains each element $\psi(x_{\beta_1}(t))$, $\psi(x_{\beta_2}(t))$, $\psi(x_{-\beta_1}(t))$ and $\psi(x_{-\beta_2}(t))$ for $t \in K$. Thus $C \leq C_G(y_0)$ and $C' \leq C_G(\langle y_0 \rangle)$; we shall show that in fact $C_G(y_0) = C$ and $C_G(\langle y_0 \rangle) = C'$.

By Lemma~\ref{lem: gen height zero}, if we take $g \in \Tran_G(y_0, Y)$ and set $y' = g.y_0 \in Y$, then we have $g = u_1 n u_2$ with $u_1 \in C_U(y')$, $u_2 \in C_U(y_0)$, and $n \in N_{\Lambda(V)_{[0]}}$ with $n.y_0 = y'$. In particular $G.y_0 \cap Y = N_{\Lambda(V)_{[0]}}.y_0 \cap Y$, and $C_G(y_0) = C_U(y_0) C_{N_{\Lambda(V)_{[0]}}}(y_0) C_U(y_0)$ while $C_G(\langle y_0 \rangle) = C_U(y_0) C_{N_{\Lambda(V)_{[0]}}}(\langle y_0 \rangle) C_U(y_0)$.

First, from the above the elements of $W$ which preserve $\Lambda(V)_{[0]}$ are those corresponding to elements of $C \cap N$; so we have $N_{\Lambda(V)_{[0]}}.y_0 = T.y_0$. Since any element of $T$ may be written as $h_{\alpha_1}(\kappa_1) h_{\alpha_2}(\kappa_2) h_{\alpha_3}(\kappa_3) h_{\alpha_4}(\kappa_4) t$ where $\kappa_1, \kappa_2, \kappa_3, \kappa_4 \in K^*$ and $t \in C \cap T$, by the above we have
$$
T.y_0 = \left\{ {\ts\frac{\kappa_1 \kappa_4}{\kappa_3}} x_1 + {\ts\frac{\kappa_2 \kappa_4}{\kappa_1 \kappa_3}} x_2 + {\ts\frac{1}{\kappa_2}} x_3 + {\ts\frac{\kappa_1}{\kappa_4}} x_4 + {\ts\frac{\kappa_3}{\kappa_1}} x_5  : \kappa_1, \kappa_2, \kappa_3, \kappa_4 \in K^* \right\}.
$$
Hence $C_{N_{\Lambda(V)_{[0]}}}(y_0) = C \cap N$; also $N_{\Lambda(V)_{[0]}}.y_0 \subseteq \hat Y$, and $N_{\Lambda(V)_{[0]}}.y_0 \cap \langle y_0 \rangle = \{ {\eta_7}^i y_0 : i = 0, 1, \dots, 6 \} = Z(G).y_0$, so $C_{N_{\Lambda(V)_{[0]}}}(\langle y_0 \rangle) = C' \cap N$.

Next, we calculate that the positive root subgroups of $C$ are diagonally embedded in the following products of root groups of $G$ (where we use the standard notation for the root system of $G$, and abbreviate $X_{\ve_i - \ve_j}$ to $X_{i - j}$):
\begin{eqnarray*}
\psi(X_{\beta_1})             & \subset & X_{1 - 2} X_{6 - 7} X_{3 - 4} X_{4 - 5} X_{3 - 5}, \\
\psi(X_{\beta_2})             & \subset & X_{2 - 3} X_{5 - 6}, \\
\psi(X_{\beta_1 + \beta_2})   & \subset & X_{1 - 3} X_{5 - 7} X_{2 - 4} X_{4 - 6} X_{2 - 6}, \\
\psi(X_{2\beta_1 + \beta_2})  & \subset & X_{2 - 5} X_{3 - 6} X_{1 - 4} X_{4 - 7} X_{1 - 7}, \\
\psi(X_{3\beta_1 + \beta_2})  & \subset & X_{1 - 5} X_{3 - 7}, \\
\psi(X_{3\beta_1 + 2\beta_2}) & \subset & X_{1 - 6} X_{2 - 7}.
\end{eqnarray*}
Thus if we take the subset $\Xi' = \{ \ve_1 - \ve_2, \ve_2 - \ve_3, \ve_1 - \ve_3, \ve_2 - \ve_5, \ve_1 - \ve_5, \ve_1 - \ve_6 \}$ of $\Phi^+$, then each of the positive long root subgroups of $C$ is diagonally embedded in $X_\alpha X_{\alpha'}$ for some $\alpha \in \Xi'$ and $\alpha' \notin \Xi'$, while each of the positive short root subgroups of $C$ is diagonally embedded in $X_\alpha X_{\alpha'} X_{\alpha''} X_{\alpha'''} X_{\alpha'' + \alpha'''}$ for some $\alpha \in \Xi'$ and $\alpha', \alpha'', \alpha''', \alpha'' + \alpha''' \notin \Xi'$. Therefore let $\Xi = \Phi^+ \setminus \Xi'$, and set $U' = \prod_{\alpha \in \Xi} X_\alpha$; then $U = U'.(C \cap U)$ and $U' \cap (C \cap U) = \{ 1 \}$. Now take $u \in U'$ satisfying $u.y_0 = y_0$, and equate coefficients of weight vectors in the order $v_{245}$, $v_{137}$, $v_{236}$, $v_{146}$, $v_{136}$, $v_{127}$, $v_{235}$, $v_{145}$, $v_{126}$, $v_{135}$, $v_{234}$, $v_{125}$, $v_{134}$, $v_{124}$, $v_{123}$ (which is compatible with increasing generalized height); this shows that the projection of $u$ must be trivial on each of the root groups $X_{5 - 6}$, $X_{3 - 4}$, $X_{6 - 7}$, $X_{4 - 5}$, $X_{3 - 5}$, $X_{2 - 4}$, $X_{5 - 7}$, $X_{4 - 6}$, $X_{1 - 4}$, $X_{3 - 6}$, $X_{4 - 7}$, $X_{2 - 6}$, $X_{3 - 7}$, $X_{2 - 7}$, $X_{1 - 7}$ in turn. Hence $u = 1$; so $C_U(y_0) = C \cap U$.

Thus $C_U(y_0), C_{N_{\Lambda(V)_{[0]}}}(y_0) \leq C$ and $C_{N_{\Lambda(V)_{[0]}}}(\langle y_0 \rangle) \leq C'$, so we do indeed have $C_G(y_0) = C$ and $C_G(\langle y_0 \rangle) = C'$. Moreover $G.y_0 \cap Y = \{ b_1 x_1 + \cdots + b_5 x_5 : b_1b_2b_3{b_4}^2{b_5}^2 = 1 \}$.

Take $y = a_1 x_1 + \cdots + a_5 x_5 \in \hat Y$. By the above, if we choose $\kappa \in K^*$ satisfying $\kappa^7 = a_1a_2a_3{a_4}^2{a_5}^2$, then $\kappa^{-1}v \in T.y_0$, so there exists $h \in T$ with $h.y_0 = \kappa^{-1}y$; so $C_G(y) = C_G(\kappa^{-1}y) = C_G(h.y_0) = {}^h C$ and likewise $C_G(\langle y \rangle) = {}^h C'$. Moreover, we see that $G.y \cap Y = G.h.\kappa y_0 \cap Y = \kappa(G.y_0 \cap Y) = \{ b_1 x_1 + \cdots + b_5 x_5 : b_1b_2b_3{b_4}^2{b_5}^2 = a_1a_2a_3{a_4}^2{a_5}^2 \}$. Since $\dim C = 14$, we have $\dim(\overline{G.y}) = \dim G - \dim C = 48 - 14 = 34$, while $\dim(\overline{G.y \cap Y}) = 4$; therefore
$$
\dim V - \dim(\overline{G.y}) = 35 - 34 = 1 \quad \hbox{and} \quad \dim Y - \dim(\overline{G.y \cap Y}) = 5 - 4 = 1.
$$
Hence $y$ is $Y$-exact. Thus the conditions of Lemma~\ref{lem: generic stabilizer from exact subset} hold; so the triple $(G, \lambda, p)$ has generic stabilizer $C/G_V \cong G_2$, while the quadruple $(G, \lambda, p, 1)$ has generic stabilizer $C'/Z(G) \cong G_2$.
\end{proof}


\begin{prop}\label{prop: A_7, omega_3 module}
Let $G = A_7$ and $\lambda = \omega_3$. Then the triple $(G, \lambda, p)$ has generic stabilizer $A_2.\Z_{(p, 2)}$; the associated first quadruple $(G, \lambda, p, 1)$ has generic stabilizer $A_2.\Z_2$.
\end{prop}

\begin{proof}
Take $G$ to be of simply connected type, so that $G = \SL_8(K)$. We may view $V$ as the exterior cube $\bigwedge^3(V_{nat})$ of the natural module, and identify $W$ with the symmetric group $S_8$. For convenience, for $i_1, i_2, i_3 \leq 8$ write $v_{i_1i_2i_3} = v_{i_1} \wedge v_{i_2} \wedge v_{i_3}$; then $V = \{ v_{i_1i_2i_3} : 1 \leq i_1 < i_2 < i_3 \leq 8 \}$. Note that $Z(G) = \langle z \rangle$ where $z = \prod_{i = 1}^7 h_{\alpha_i}({\eta_8}^i)$; as $z$ acts on $V_{nat}$ as multiplication by $\eta_8$, it acts on $V_{nat} \otimes V_{nat} \otimes V_{nat}$ and hence on $V$ as multiplication by ${\eta_8}^3$, so $G_V = \{ 1 \}$.

We take the generalized height function on the weight lattice of $G$ whose value at $\alpha_4$ is $0$, and at each other simple root $\alpha_i$ is $1$; then the generalized height of $\lambda = \frac{1}{8}(5\alpha_1 + 10\alpha_2 + 15\alpha_3 + 12\alpha_4 + 9\alpha_5 + 6\alpha_6 + 3\alpha_7)$ is $6$, and as $\lambda$ and $\Phi$ generate the weight lattice we see that the generalized height of any weight is an integer. Define $\sigma : \{ 1, \dots, 8 \} \to \{ 1, \dots, 7 \}$ by $\sigma(i) = i$ if $i \leq 4$ and $i - 1$ if $i \geq 5$. Since $V_\lambda = \langle v_{123} \rangle$, we see that if $\mu \in \Lambda(V)$ and $v_{i_1i_2i_3} \in V_\mu$, then the generalized height of $\mu$ is $12 - (\sigma(i_1) + \sigma(i_2) + \sigma(i_3))$. Thus $\Lambda(V)_{[0]} = \{ \nu_1, \dots, \nu_8 \}$, where we write
\begin{eqnarray*}
& x_1 = v_{148}, \quad x_2 = v_{158}, \quad x_3 = v_{247}, \quad x_4 = v_{257}, & \\
& x_5 = v_{346}, \quad x_6 = v_{356}, \quad x_7 = v_{167}, \quad x_8 = v_{238}, &
\end{eqnarray*}
and for each $i$ we let $\nu_i$ be the weight such that $x_i \in V_{\nu_i}$. Observe that if we take $s = \prod_{i = 1}^7 h_{\alpha_i}(\kappa_i) \in T$, then $\nu_1(s) = \frac{\kappa_1 \kappa_4}{\kappa_3 \kappa_7}$, $\nu_2(s) = \frac{\kappa_1 \kappa_5}{\kappa_4 \kappa_7}$, $\nu_3(s) = \frac{\kappa_2 \kappa_4 \kappa_7}{\kappa_1 \kappa_3 \kappa_6}$, $\nu_4(s) = \frac{\kappa_2 \kappa_5 \kappa_7}{\kappa_1 \kappa_4 \kappa_6}$, $\nu_5(s) = \frac{\kappa_4 \kappa_6}{\kappa_2 \kappa_5}$, $\nu_6(s) = \frac{\kappa_3 \kappa_6}{\kappa_2 \kappa_4}$, $\nu_7(s) = \frac{\kappa_1 \kappa_7}{\kappa_5}$ and $\nu_8(s) = \frac{\kappa_3}{\kappa_1 \kappa_7}$; thus given any triple $(n_1, n_2, n_3)$ of integers we have $c_1\nu_1 + \cdots + c_8\nu_8 = 0$ for $(c_1, \dots, c_8) = (n_1 - n_2 + n_3, n_1 + n_2 - n_3, n_1 + n_2, n_1 - n_2, n_1 - n_3, n_1 + n_3, n_1, n_1)$. In particular, writing \lq $(n_1, n_2, n_3) \implies (c_1, c_2, c_3, c_4, c_5, c_6, c_7, c_8)$' to indicate this relationship between triples and $8$-tuples, we have the following:
\begin{eqnarray*}
(1, 1, 1) \implies (1, 1, 2, 0, 0, 2, 1, 1), & & \phantom{-1} (2, 1, 2) \implies (3, 1, 3, 1, 0, 4, 2, 2), \\
(1, 1, 0) \implies (0, 2, 2, 0, 1, 1, 1, 1), & & \phantom{1} (2, 1, -1) \implies (0, 4, 3, 1, 3, 1, 2, 2), \\
(2, 2, 1) \implies (1, 3, 4, 0, 1, 3, 2, 2), & & \phantom{-1} (1, 0, 0) \implies (1, 1, 1, 1, 1, 1, 1, 1).
\end{eqnarray*}
It follows that any subset of $\Lambda(V)_{[0]}$ which contains $\nu_2$, $\nu_3$, $\nu_6$, $\nu_7$, $\nu_8$ and at least one of $\nu_1$ and $\nu_5$ has ZLCE. Set $Y = V_{[0]} = \langle x_1, \dots, x_8 \rangle$, and let
$$
\hat Y = \left\{ a_1 x_1 + \cdots + a_8 x_8 : (a_1a_4 - a_2a_3)(a_1a_6 - a_2a_5)(a_3a_6 - a_4a_5)a_7a_8 \neq 0 \right\},
$$
so that $\hat Y$ is a dense open subset of $Y$. Write
$$
y_0 = x_1 + x_2 + x_3 + x_6 + x_7 - x_8 \in \hat Y.
$$

Here $W$ acts on $\Lambda(V)$ such that if $w \in W$ and $\mu \in \Lambda(V)$ with $v_{i_1i_2i_3} \in V_\mu$ then $v_{w(i_1) w(i_2) w(i_3)} \in V_{w(\mu)}$. For each of the unordered pairs $\{ 1, 8 \}$, $\{ 2, 7 \}$ and $\{ 3, 6 \}$, there exist $i$ and $j$ distinct such that the pair concerned is the intersection of the sets of three numbers appearing as subscripts in $x_i$ and $x_j$; as this is true for no other unordered pairs, the setwise stabilizer in $W$ of $\Lambda(V)_{[0]}$ must permute these three pairs and hence preserve $\{ 4, 5 \}$, so it is a subgroup of $(S_2 \wr S_3) \times S_2$. Thus it must preserve $\{ \nu_7, \nu_8 \}$, so any element which preserves each of the three pairs and fixes $1$ must also fix $6$ and $7$; thus the order of the setwise stabilizer is at most $\frac{1}{4}.2^3.3!.2 = 24$. Now $(1 \ 2)(3 \ 6)(7 \ 8)$ interchanges $\nu_1$ and $\nu_3$, $\nu_2$ and $\nu_4$, and also $\nu_7$ and $\nu_8$, while fixing $\nu_5$ and $\nu_6$, and $(2 \ 3)(6 \ 7)$ interchanges $\nu_3$ and $\nu_5$, and also $\nu_4$ and $\nu_6$, while fixing the remaining $\nu_i$, and these two permutations generate a dihedral group of order $12$; moreover this group commutes with $(4 \ 5)$, which interchanges $\nu_1$ and $\nu_2$, $\nu_3$ and $\nu_4$, and also $\nu_5$ and $\nu_6$, while fixing $\nu_7$ and $\nu_8$. Thus we see that the setwise stabilizer in $W$ of $\Lambda(V)_{[0]}$ is $\langle (1 \ 2)(3 \ 6)(7 \ 8), (2 \ 3)(6 \ 7), (4 \ 5) \rangle = \langle w_{\alpha_1} w_{\alpha_3 + \alpha_4 + \alpha_5} w_{\alpha_7}, w_{\alpha_2} w_{\alpha_6}, w_{\alpha_4} \rangle$. Note that this stabilizes $\Phi_{[0]} = \langle \alpha_4 \rangle$.

Let $A$ be a simply connected group defined over $K$ of type $A_2$, with simple roots $\beta_1$ and $\beta_2$; then we may regard $V_{nat}$ as $\L(A)$, with $v_1 = e_{\beta_1 + \beta_2}$, $v_2 = e_{\beta_1}$, $v_3 = e_{\beta_2}$, $v_4 = h_{\beta_1}$, $v_5 = h_{\beta_2}$, $v_6 = f_{\beta_2}$, $v_7 = f_{\beta_1}$, $v_8 = f_{\beta_1 + \beta_2}$. Moreover the action of $A$ on its Lie algebra gives a homomorphism $\psi: A \to G$, with kernel $Z(A)$; then $\psi(A)$ is a subgroup of $G$ which is an adjoint group defined over $K$ of type $A_2$. If we take the structure constants of $A$ to be such that $[e_{\beta_1}, e_{\beta_2}] = e_{\beta_1 + \beta_2}$, then with respect to the basis $v_1, \dots, v_8$ the simple root elements $x_{\beta_1}(t)$ and $x_{\beta_2}(t)$ of $G$ act by the matrices
$$
\left(
  \begin{array}{cccccccc}
    1 &   & t &     &   &   &      &    \\
      & 1 &   & -2t & t &   & -t^2 &    \\
      &   & 1 &     &   &   &      &    \\
      &   &   &  1  &   &   &   t  &    \\
      &   &   &     & 1 &   &      &    \\
      &   &   &     &   & 1 &      & -t \\
      &   &   &     &   &   &   1  &    \\
      &   &   &     &   &   &      &  1 \\
  \end{array}
\right)
\qquad\hbox{ and }\qquad
\left(
  \begin{array}{cccccccc}
    1 & -t &   &   &     &      &   &   \\
      &  1 &   &   &     &      &   &   \\
      &    & 1 & t & -2t & -t^2 &   &   \\
      &    &   & 1 &     &      &   &   \\
      &    &   &   &  1  &   t  &   &   \\
      &    &   &   &     &   1  &   &   \\
      &    &   &   &     &      & 1 & t \\
      &    &   &   &     &      &   & 1 \\
  \end{array}
\right)
$$
respectively, and the corresponding negative root elements $x_{-\beta_1}(t)$ and $x_{-\beta_2}(t)$ act by the matrices
$$
\left(
  \begin{array}{cccccccc}
    1 &      &   &    &    &    &   &   \\
      &   1  &   &    &    &    &   &   \\
    t &      & 1 &    &    &    &   &   \\
      &  -t  &   &  1 &    &    &   &   \\
      &      &   &    &  1 &    &   &   \\
      &      &   &    &    &  1 &   &   \\
      & -t^2 &   & 2t & -t &    & 1 &   \\
      &      &   &    &    & -t &   & 1 \\
  \end{array}
\right)
\qquad\hbox{ and }\qquad
\left(
  \begin{array}{cccccccc}
     1 &   &      &    &    &   &   &   \\
    -t & 1 &      &    &    &   &   &   \\
       &   &   1  &    &    &   &   &   \\
       &   &      &  1 &    &   &   &   \\
       &   &  -t  &    &  1 &   &   &   \\
       &   & -t^2 & -t & 2t & 1 &   &   \\
       &   &      &    &    &   & 1 &   \\
       &   &      &    &    &   & t & 1 \\
  \end{array}
\right)
$$
respectively. Write $n^*$ for the element
$$
\eta_{16}
\left(
  \begin{array}{cccccccc}
   -1 &   &   &   &   &   &   &    \\
      &   & 1 &   &   &   &   &    \\
      & 1 &   &   &   &   &   &    \\
      &   &   &   & 1 &   &   &    \\
      &   &   & 1 &   &   &   &    \\
      &   &   &   &   &   & 1 &    \\
      &   &   &   &   & 1 &   &    \\
      &   &   &   &   &   &   & -1 \\
  \end{array}
\right)
$$
of $N$; then conjugation by $n^*$ interchanges the elements $\psi(x_{\beta_1}(t))$ and $\psi(x_{\beta_2}(t))$, and also the elements $\psi(x_{-\beta_1}(t))$ and $\psi(x_{-\beta_2}(t))$, so acts as a graph automorphism of $\psi(A)$, while $(n^*)^2 = z$. A straightforward calculation shows that the stabilizer of $y_0$ contains each element $\psi(x_{\beta_1}(t))$, $\psi(x_{\beta_2}(t))$, $\psi(x_{-\beta_1}(t))$ and $\psi(x_{-\beta_2}(t))$ for $t \in K$, while $n^*.y_0 = {\eta_{16}}^3 y_0$. Set $C = \psi(A)$ or $\psi(A)\langle n^* \rangle$ according as $p \geq 3$ or $p = 2$, and $C' = \psi(A)\langle n^* \rangle$. Clearly we then have $C \leq C_G(y_0)$ and $C' \leq C_G(\langle y_0 \rangle)$; we shall show that in fact $C_G(y_0) = C$ and $C_G(\langle y_0 \rangle) = C'$.

We have $U_{[0]} = X_{\alpha_4}$. Given $u = x_{\alpha_4}(t) \in U_{[0]}$ we have
$$
u.y_0 = (1 + t)x_1 + x_2 + x_3 + t x_5 + x_6 + x_7 - x_8;
$$
so the set of weights occurring in $u.y_0$ contains $\nu_2$, $\nu_3$, $\nu_6$, $\nu_7$ and $\nu_8$ and at least one of $\nu_1$ and $\nu_5$, and hence by the above has ZLCE. By Lemma~\ref{lem: gen height zero not strictly positive}, if we take $g \in \Tran_G(y_0, Y)$ and write $y' = g.y_0 \in Y$, then we may write $g = u_1 g' u_2$ with $u_1 \in C_{U_{[+]}}(y')$, $u_2 \in C_{U_{[+]}}(y_0)$, and $g' \in G_{[0]} N_{\Lambda(V)_{[0]}}$ with $g'.y_0 = y'$. In particular $G.y_0 \cap Y = G_{[0]} N_{\Lambda(V)_{[0]}}.y_0 \cap Y$, and $C_G(y_0) = C_{U_{[+]}}(y_0) C_{G_{[0]} N_{\Lambda(V)_{[0]}}}(y_0) C_{U_{[+]}}(y_0)$ while $C_G(\langle y_0 \rangle) = C_{U_{[+]}}(y_0) C_{G_{[0]} N_{\Lambda(V)_{[0]}}}(\langle y_0 \rangle) C_{U_{[+]}}(y_0)$.

First, from the above the elements of $W$ which preserve $\Lambda(V)_{[0]}$ are those corresponding to elements of $(\psi(A) \cap N)\langle n^* \rangle \langle n_{\alpha_4} \rangle$. Since $w_{\alpha_4} \in Z(W_{\Lambda(V)_{[0]}})$, we see that $G_{[0]} N_{\Lambda(V)_{[0]}}.y_0 = G_{[0]}.y_0 \cup n^*G_{[0]}.y_0$. Since any element of $G_{[0]}$ may be written as $h_{\alpha_1}(\kappa_1) h_{\alpha_2}(\kappa_2) h_{\alpha_3}(\kappa_3) h_{\alpha_5}(\kappa_5) xt$ where $\kappa_1, \kappa_2, \kappa_3, \kappa_5 \in K^*$, $t \in \psi(A) \cap T$, and $x \in \langle X_{\pm\alpha_4} \rangle$ fixes $v_i$ for $i \neq 4, 5$ and maps $v_4 \mapsto a v_4 + c v_5$ and $v_5 \mapsto b v_4 + d v_5$ for some $a, b, c, d \in K$ satisfying $ad - bc = 1$, and $n^*.x_i = {\eta_{16}}^3 x_{\pi(i)}$ where $\pi = (1 \ 2)(3 \ 6)(4 \ 5)$, by the above we have
\begin{eqnarray*}
G_{[0]}.y_0     & = & \left\{ {\ts\frac{\kappa_1}{\kappa_3}}(a + b) x_1 + \kappa_1\kappa_5 (c + d) x_2 + {\ts\frac{\kappa_2}{\kappa_1\kappa_3}} a x_3 + {\ts\frac{\kappa_2\kappa_5}{\kappa_1}} c x_4 + {\ts\frac{1}{\kappa_2\kappa_5}} b x_5 \right. \\
                &   & \left. \quad {} + {\ts\frac{\kappa_3}{\kappa_2}} d x_6 + {\ts\frac{\kappa_1}{\kappa_5}} x_7 - {\ts\frac{\kappa_3}{\kappa_1}} x_8 : \kappa_1, \kappa_2, \kappa_3, \kappa_5 \in K^*, \ ad - bc = 1 \right\}, \\
n^* G_{[0]}.y_0 & = & \left\{ {\eta_{16}}^3 ({\ts\frac{\kappa_1}{\kappa_3}}(a + b) x_2 + \kappa_1\kappa_5 (c + d) x_1 + {\ts\frac{\kappa_2}{\kappa_1\kappa_3}} a x_6 + {\ts\frac{\kappa_2\kappa_5}{\kappa_1}} c x_5 + {\ts\frac{1}{\kappa_2\kappa_5}} b x_4 \right. \\
                &   & \left. \quad {} + {\ts\frac{\kappa_3}{\kappa_2}} d x_3 + {\ts\frac{\kappa_1}{\kappa_5}} x_7 - {\ts\frac{\kappa_3}{\kappa_1}} x_8) : \kappa_1, \kappa_2, \kappa_3, \kappa_5 \in K^*, \ ad - bc = 1 \right\}.
\end{eqnarray*}
Equating the expression in the first set to $\kappa y_0$ we see that $b = c = 0$ (from $x_4$ and $x_5$), then $\kappa_1 = \kappa^{-1}\kappa_3 = \kappa\kappa_5$ (from $x_7$ and $x_8$), then $a = \kappa^2 = d^{-1}$ (from $x_1$ and $ad - bc = 1$), then $\kappa_2 = \kappa^{-2}\kappa_1$ (from $x_6$), then $\kappa_1 = \kappa^{-2}$ (from $x_3$), and finally $\kappa^8 = 1$ (from $x_2$). Likewise equating that in the second set to ${\eta_{16}}^3 \kappa y_0$ we see that $b = c = 0$ (from $x_4$ and $x_5$), after which the equations are as before, so we obtain $({\eta_{16}}^3 \kappa)^8 = 1$, and hence $\kappa^8 = -1$. Hence $C_{G_{[0]} N_{\Lambda(V)_{[0]}}}(y_0) = C \cap N$; also $G_{[0]} N_{\Lambda(V)_{[0]}}.y_0 \subseteq \hat Y$, and $G_{[0]} N_{\Lambda(V)_{[0]}}.y_0 \cap \langle y_0 \rangle = \{ {\eta_{16}}^i y_0 : i = 0, 1, \dots, 15 \} = \langle n^* \rangle.y_0$, so $C_{G_{[0]} N_{\Lambda(V)_{[0]}}}(\langle y_0 \rangle) = C' \cap N$.

Next, we calculate that the positive root subgroups of $C$ are diagonally embedded in the following products of root groups of $G$ (where we use the standard notation for the root system of $G$, and abbreviate $X_{\ve_i - \ve_j}$ to $X_{i - j}$):
\begin{eqnarray*}
\psi(X_{\beta_1})           & \subset & X_{1 - 3} X_{6 - 8} X_{2 - 4} X_{2 - 5} X_{4 - 7} X_{2 - 7}, \\
\psi(X_{\beta_2})           & \subset & X_{1 - 2} X_{7 - 8} X_{3 - 4} X_{3 - 5} X_{5 - 6} X_{3 - 6}, \\
\psi(X_{\beta_1 + \beta_2}) & \subset & X_{2 - 6} X_{3 - 7} X_{1 - 4} X_{1 - 5} X_{4 - 8} X_{5 - 8} X_{1 - 8}.
\end{eqnarray*}
Thus if we take the subset $\Xi' = \{ \ve_1 - \ve_2, \ve_1 - \ve_3, \ve_2 - \ve_6 \}$ of $\Phi^+$, then each of the positive root subgroups of $C$ is diagonally embedded in $X_\alpha X_{\alpha'} X_{\alpha''} \dots$ for some $\alpha \in \Xi'$ and $\alpha', \alpha'', \dots \notin \Xi'$. Thus if we let $\Xi = \Phi^+ \setminus (\Xi' \cup \{ \alpha_4 \})$, and set $U' = \prod_{\alpha \in \Xi} X_\alpha$, then $U_{[+]} = U'.(C \cap U_{[+]})$ and $U' \cap (C \cap U_{[+]}) = 1$. Now take $u \in U'$ satisfying $u.y_0 = y_0$, and equate coefficients of weight vectors in the order
\begin{eqnarray*}
& v_{246}, v_{256}, v_{345}, v_{147}, v_{157}, v_{237}, v_{138}, v_{137}, v_{245}, v_{156}, v_{236}, v_{128}, & \\
& v_{146}, v_{127}, v_{235}, v_{145}, v_{234}, v_{136}, v_{126}, v_{134}, v_{135}, v_{125}, v_{124}, v_{123} \phantom{,} &
\end{eqnarray*}
(which is compatible with increasing generalized height); this shows that the projection of $u$ must be trivial on each of the root groups $X_{6 - 7}$, $X_{2 - 3}$, $X_{4 - 6}$, $X_{7 - 8}$, $X_{5 - 6}$, $X_{3 - 4}$, $X_{3 - 5}$, $X_{3 - 6}$, $X_{5 - 7}$, $X_{6 - 8}$, $X_{2 - 5}$, $X_{2 - 4}$, $X_{4 - 7}$, $X_{1 - 4}$, $X_{5 - 8}$, $X_{4 - 8}$, $X_{3 - 7}$, $X_{1 - 5}$, $X_{2 - 7}$, $X_{3 - 8}$, $X_{1 - 6}$, $X_{2 - 8}$, $X_{1 - 7}$, $X_{1 - 8}$ in turn. Hence $u = 1$, so $C_{U_{[+]}}(y_0) = C \cap U_{[+]}$.

Thus $C_{U_{[+]}}(y_0), C_{G_{[0]} N_{\Lambda(V)_{[0]}}}(y_0) \leq C$ and $C_{G_{[0]} N_{\Lambda(V)_{[0]}}}(\langle y_0 \rangle) \leq C'$, so we do indeed have $C_G(y_0) = C$ and $C_G(\langle y_0 \rangle) = C'$. Moreover if we write the expression in the set $G_{[0]}.y_0$ as $b_1 x_1 + \cdots + b_8 x_8$, then $(b_1b_4 - b_2b_3)(b_1b_6 - b_2b_5)(b_3b_6 - b_4b_5)b_7b_8 = 1$; conversely given $b_1, \dots, b_8$ satisfying $(b_1b_4 - b_2b_3)(b_1b_6 - b_2b_5)(b_3b_6 - b_4b_5)b_7b_8 = 1$, if we write $\Delta_1 = b_3b_6 - b_4b_5$ and $\Delta_2 = b_1b_6 - b_2b_5$, and set $\kappa_1 = \frac{1}{\Delta_1}$, $\kappa_2 = \frac{1}{\Delta_1\Delta_2}$, $\kappa_3 = -\frac{b_8}{\Delta_1}$, $\kappa_5 = \frac{1}{b_7\Delta_1}$, $a = -\frac{b_3b_8\Delta_2}{\Delta_1}$, $b = \frac{b_5}{b_7{\Delta_1}^2\Delta_2}$, $c = b_4b_7\Delta_1\Delta_2$ and $d = -\frac{b_6}{b_8\Delta_2}$, then the expression given in the first set is equal to $b_1 x_1 + \cdots + b_8 x_8$. Treating the set $n^* G_{[0]}.y_0$ entirely similarly we see that we have $G.y_0 \cap Y = \{ b_1 x_1 + \cdots + b_8 x_8 : ((b_1b_4 - b_2b_3)(b_1b_6 - b_2b_5)(b_3b_6 - b_4b_5)b_7b_8)^2 = 1 \}$.

Take $y = a_1 x_1 + \cdots + a_8 x_8 \in \hat Y$. By the above, if we choose $\kappa \in K^*$ satisfying $\kappa^8 = (a_1a_4 - a_2a_3)(a_1a_6 - a_2a_5)(a_3a_6 - a_4a_5)a_7a_8$, then $\kappa^{-1}y \in G_{[0]}.y_0$, so there exists $h \in G_{[0]}$ with $h.y_0 = \kappa^{-1}y$; so $C_G(y) = C_G(\kappa^{-1}y) = C_G(h.y_0) = {}^h C$ and likewise $C_G(\langle y \rangle) = {}^h C'$. Moreover, we see that $G.y \cap Y = G.h.\kappa y_0 \cap Y = \kappa(G.y_0 \cap Y) = \{ b_1 x_1 + \cdots + b_8 x_8 : ((b_1b_4 - b_2b_3)(b_1b_6 - b_2b_5)(b_3b_6 - b_4b_5)b_7b_8)^2 = ((a_1a_4 - a_2a_3)(a_1a_6 - a_2a_5)(a_3a_6 - a_4a_5)a_7a_8)^2 \}$. Since $\dim C = 8$, we have $\dim(\overline{G.y}) = \dim G - \dim C = 63 - 8 = 55$, while $\dim(\overline{G.y \cap Y}) = 7$; therefore
$$
\dim V - \dim(\overline{G.y}) = 56 - 55 = 1 \quad \hbox{and} \quad \dim Y - \dim(\overline{G.y \cap Y}) = 8 - 7 = 1.
$$
Hence $y$ is $Y$-exact. Thus the conditions of Lemma~\ref{lem: generic stabilizer from exact subset} hold; so the triple $(G, \lambda, p)$ has generic stabilizer $C/G_V \cong A_2.\Z_{(p, 2)}$, while the quadruple $(G, \lambda, p, 1)$ has generic stabilizer $C'/Z(G) \cong A_2.\Z_2$, where the $A_2$ is of adjoint type.
\end{proof}

This completes the justification of the entries in Tables~\ref{table: small classical triple and first quadruple generic stab} and \ref{table: small exceptional triple and first quadruple generic stab}, and hence the proof of Theorem~\ref{thm: small triple and first quadruple generic stab}.

In concluding this section, we acknowledge that many of its results are known and available in various places in the literature. We have made no attempt to provide a list of references, both because any such list would almost certainly be incomplete, and because our work is independent of what has gone before.

\chapter{Higher quadruples not having TGS}\label{chap: non-TGS higher quadruples}

In this chapter we consider higher quadruples which do not have TGS, and establish the entries in Tables~\ref{table: large higher quadruple non-TGS}, \ref{table: small classical higher quadruple generic stab} and \ref{table: small exceptional higher quadruple generic stab}. For the most part, in Sections~\ref{sect: non-TGS large higher quadruples} and \ref{sect: small higher quadruples} we treat higher quadruples which are large and small respectively. However, it turns out to be convenient to blur the distinction between the two slightly: in a number of instances we will postpone treatment of a large higher quadruple to Section~\ref{sect: small higher quadruples}; on the other hand one of the small higher quadruples occurs in an infinite family where all the other quadruples are large, so will be treated in Section~\ref{sect: non-TGS large higher quadruples}. Throughout, given a quadruple $(G, \lambda, p, k)$ we write $V = L(\lambda)$.

As in Chapter~\ref{chap: non-TGS triples and first quadruples}, in many cases our approach will be to apply Lemma~\ref{lem: generic stabilizer from exact subset} to determine the required generic stabilizer. Again we choose $Y$ (although here it is a subset of $\Gk(V)$ rather than a subspace of $V$), and take a dense open subset $\hat Y$ of $Y$. For all $y \in \hat Y$, we show that the stabilizer $C_G(y)$ is a conjugate of a fixed subgroup $C$, and that $y$ is $Y$-exact. By Lemma~\ref{lem: generic stabilizer from exact subset} we may now conclude that the quadruple $(G, \lambda, p, k)$ has generic stabilizer $C/Z(G)$. Various methods are used to determine the stabilizer $C_G(y)$ and the transporter $\Tran_G(y, Y)$. Often we refer back to the Proposition in the corresponding section in Chapter~\ref{chap: non-TGS triples and first quadruples} which dealt with the triple to which the quadruple is associated, to allow us to use both the set-up established there and results obtained in the course of the proof.

\section{Large higher quadruples}\label{sect: non-TGS large higher quadruples}

In this section we shall treat some of the large higher quadruples listed in Table~\ref{table: large higher quadruple non-TGS}. In fact, it will be convenient at the end of this section to treat a few large higher quadruples in which the group is not simple, for use in the following section.

\begin{prop}\label{prop: A_ell, 2omega_1 module, k = 2}
Let $G = A_\ell$ for $\ell \in [2, \infty)$ and $\lambda = 2\omega_1$ with $p \geq 3$, and take $k = 2$. Then according as $\ell = 2$, or $\ell = 3$, or $\ell \geq 4$, the quadruple $(G, \lambda, p, k)$ has generic stabilizer ${\Z_2}^2.S_3$, or ${\Z_2}^3.{\Z_2}^2$, or ${\Z_2}^\ell$, respectively.
\end{prop}

\begin{proof}
We take $G = \SL_{\ell + 1}(K)$ and view $V$ as the space of $(\ell + 1) \times (\ell + 1)$ symmetric matrices over $K$, where $g \in G$ sends $A$ to $gAg^T$; write $X = \G{2}(V)$. We take $T$ to be the standard maximal torus of $G$ consisting of diagonal matrices. Set
$$
Y = \{ \langle I, \diag(\kappa_1, \dots, \kappa_{\ell + 1}) \rangle : \exists i \neq j \hbox{ with } \kappa_i \neq \kappa_j  \},
$$
and
$$
\hat Y_1 = \{ \langle I, \diag(\kappa_1, \dots, \kappa_{\ell + 1}) \rangle : \forall i \neq j, \ \kappa_i \neq \kappa_j \};
$$
then $\hat Y_1$ is a dense open subset of $Y$. As in Section~\ref{sect: localization}, we have the orbit map $\phi: G \times Y \to X$; we claim that $\phi(G \times \hat Y_1)$ contains a dense open subset of $X$.

To see this, take any subspace $\langle A, B \rangle \in X$ where $\det A = 1$ such that $A^{-1}B$ has $\ell + 1$ distinct eigenvalues; clearly the set of such subspaces is dense in $X$. Let the eigenvalues of $A^{-1}B$ be $\kappa_1, \dots, \kappa_{\ell + 1}$, with corresponding eigenvectors $v_1, \dots, v_{\ell + 1} \in K^{\ell + 1}$; write $D = \diag(\kappa_1, \dots, \kappa_{\ell + 1})$, then $\langle I, D \rangle \in \hat Y_1$. Since eigenvectors corresponding to distinct eigenvalues are linearly independent, $v_1, \dots, v_{\ell + 1}$ form a basis of $K^{\ell + 1}$; as $A$ is invertible, so do $Av_1, \dots, Av_{\ell + 1}$. For all $i$, as $A^{-1}Bv_i = \kappa_i v_i$ we have $Bv_i = \kappa_i Av_i$, so for $j \neq i$ we have
$$
\kappa_i {v_i}^T Av_j = (\kappa_i Av_i)^T v_j = (Bv_i)^T v_j = {v_i}^T Bv_j = {v_i}^T \kappa_j Av_j = \kappa_j {v_i}^T Av_j,
$$
and as $\kappa_i \neq \kappa_j$ we must have ${v_i}^T Av_j = 0$; since for any non-zero vector $v \in K^{\ell + 1}$ there exists $v' \in K^{\ell + 1}$ with $v^T v' = 0$, we must have ${v_i}^T Av_i \neq 0$, and by replacing $v_i$ by a scalar multiple we may assume ${v_i}^T Av_i = 1$. Thus for all $i$ and $j$ we have ${v_i}^T Av_j = \delta_{ij}$ and hence ${v_i}^T Bv_j = \kappa_i \delta_{ij}$. Let $R$ be the matrix whose $i$th row is ${v_i}^T$; then $RAR^T = I$ and $RBR^T = D$. Hence $(\det R)^2 = 1$, so $\det R = \pm 1$; by negating $v_1$ if necessary we may assume that $\det R = 1$. Thus $R \in G$, and $R$ sends the pair $(A, B)$ to the pair $(I, D)$, so that $\phi(R^{-1}, \langle I, D \rangle) = \langle A, B \rangle$, proving the claim.

Now take $y = \langle I, D \rangle \in \hat Y_1$, and consider $C_G(y)$. The kernel $J$ of the action of $C_G(y)$ on the subspace $y$ consists of the $g \in G$ such that $gIg^T = I$ and $gDg^T = D$; the first condition gives $g^T = g^{-1}$ and then the second gives $gD = Dg$, so $g \in C_G(D) = T$, and now the first condition again gives $g^2 = I$, so that $J = \{ s \in T : s^2 = 1 \}$. Since $J \lhd C_G(y)$ we have $C_G(y) \leq N_G(J)$. To identify $N_G(J)$, first consider $C_G(J)$. Clearly $T \leq C_G(J)$; given $g \in G \setminus T$ there exist $i$, $j$ distinct with $g_{ij} \neq 0$, and then if we take $i' \neq i, j$ then $g$ does not commute with the element of $J$ whose $i$th and $i'$th diagonal entries are $-1$ and all the other diagonal entries are $1$, so that $g \notin C_G(J)$. Hence $C_G(J) = T$; as $C_G(J) \lhd N_G(J)$ we have $N_G(J) \leq N_G(T) = N$, and as $N$ does normalise $J$ we have $N_G(J) = N$. Thus $C_G(y) \leq N$.

Write $Z_2 = \{ s \in T : s^2 \in Z(G) \}$. Let $C$ be the subgroup $Z_2\langle n_{\alpha_1}, h_{\alpha_2}(\eta_4) n_{\alpha_2} \rangle$, $Z_2\langle n_{\alpha_1} n_{\alpha_3}, h_{\alpha_2 + \alpha_3}(\eta_4) n_{\alpha_1 + \alpha_2} n_{\alpha_2 + \alpha_3} \rangle$ or $Z_2$ according as $\ell = 2$, $\ell = 3$ or $\ell \geq 4$. We shall define a dense open subset $\hat Y$ of $Y$ lying in $\hat Y_1$, and show that if $y \in \hat Y$ then $C_G(y)$ is a conjugate of $C$.

Given $y \in \hat Y_1$, we may write $y = \langle D_1, D_2 \rangle$ with $D_1 = \diag(a_1, a_2, a_3, \dots, a_{\ell + 1})$ and $D_2 = \diag(b_1, b_2, b_3, \dots, b_{\ell + 1})$. Note that by changing basis we may assume if we wish that $a_1 = b_2 = 1$, $a_2 = b_1 = 0$, in which case $a_3, \dots, a_{\ell + 1}, b_3, \dots, b_{\ell + 1} \neq 0$. Then if we take $s = \diag(\kappa_1, \dots, \kappa_{\ell + 1}) \in C_T(y)$, we must have $s.D_i \in \langle D_i \rangle$ for $i = 1, 2$, whence ${\kappa_1}^2 = {\kappa_2}^2 = \cdots = {\kappa_{\ell + 1}}^2$, and so $s^2 \in Z(G)$, whence $s \in Z_2$; conversely if $s \in Z_2$ then clearly $s.D_i \in \langle D_i \rangle$ for $i = 1, 2$, so $s \in C_T(y)$. Therefore $C_T(y) = Z_2$. It remains to consider which elements of $W$ give rise to a coset of $C_T(y)$ in $C_G(y)$.

First suppose $\ell = 2$; here we set $\hat Y = \hat Y_1$. Let $y_0 = \langle \diag(1, 0, 1), \diag(0, 1, 1) \rangle \in \hat Y$. Take $y \in \hat Y$; by the above we may assume $y = \langle \diag(1, 0, a_3), \diag(0, 1, b_3) \rangle$ with $a_3, b_3 \neq 0$. Choose $c_3, d_3, \kappa \in K^*$ satisfying ${c_3}^2 = a_3$, ${d_3}^2 = b_3$ and $\kappa^3 = (c_3d_3)^{-1}$, and write $h^{-1} = \diag(\kappa c_3, \kappa d_3, \kappa)$; then $h^{-1}.y = y_0$. As $n_{\alpha_1}$ sends the matrix $\diag(a, b, c)$ to $\diag(b, a, c)$, it interchanges $\diag(1, 0, 1)$ and $\diag(0, 1, 1)$, and thus stabilizes $y_0$; likewise as $h_{\alpha_2}(\eta_4) n_{\alpha_2}$ sends the matrix $\diag(a, b, c)$ to $\diag(a, -c, -b)$, it sends $\diag(1, 0, 1)$ to $\diag(1, 0, 1) - \diag(0, 1, 1)$ and negates $\diag(0, 1, 1)$, and thus also stabilizes $y_0$. Hence $C_G(y_0) = C$, and so $C_G(y) = C_G(h.y_0) = {}^h C$.

Now suppose $\ell \geq 3$. Take $n \in N \setminus T$, and write $n = n^* s$ where $s = \diag(\kappa_1, \dots, \kappa_{\ell + 1})$ and $n^*$ is a permutation matrix corresponding to the permutation $\pi^{-1} \in S_{\ell + 1} \setminus \{ 1 \}$. If $n.y = y$ there must exist $c_1, c_2, c_3, c_4 \in K$ with $(c_1, c_2), (c_3, c_4) \neq (0, 0)$ such that $n.D_1 = c_1 D_1 + c_2 D_2$ and $n.D_2 = c_3 D_1 + c_4 D_2$. Thus for all $i \leq \ell + 1$ we have ${\kappa_i}^2 a_i = c_1 a_{\pi(i)} + c_2 b_{\pi(i)}$ and ${\kappa_i}^2 b_i = c_3 a_{\pi(i)} + c_4 b_{\pi(i)}$, and so $c_1 a_{\pi(i)} b_i + c_2 b_{\pi(i)} b_i = c_3 a_{\pi(i)} a_i + c_4 b_{\pi(i)} a_i$, whence
$$
\left(
  \begin{array}{cccc}
    a_{\pi(1)} b_1 & b_{\pi(1)} b_1 & a_{\pi(1)} a_1 & b_{\pi(1)} a_1 \\
    a_{\pi(2)} b_2 & b_{\pi(2)} b_2 & a_{\pi(2)} a_2 & b_{\pi(2)} a_2 \\
    \vdots & \vdots & \vdots & \vdots \\
    a_{\pi(\ell + 1)} b_{\ell + 1} & b_{\pi(\ell + 1)} b_{\ell + 1} & a_{\pi(\ell + 1)} a_{\ell + 1} & b_{\pi(\ell + 1)} a_{\ell + 1} \\
  \end{array}
\right)
\left(
  \begin{array}{c}
    c_1 \\
    c_2 \\
    -c_3 \\
    -c_4 \\
  \end{array}
\right)
=
\left(
  \begin{array}{c}
    0 \\
    0 \\
    \vdots \\
    0 \\
  \end{array}
\right).
$$
Since $\pi \neq 1$, there exists $i$ with $\pi(i) \neq i$; without loss of generality we may assume $\pi(1) = 2$. We claim that, unless $\ell = 3$ and $\pi = (1\ 2)(3\ 4)$, the points $y$ for which the above matrix equation has a non-zero solution $(c_1, c_2, c_3, c_4)$ form a proper closed subvariety of $Y$; to do this it suffices to show that there is a non-zero polynomial equation in $a_1, \dots, a_{\ell + 1}, b_1, \dots, b_{\ell + 1}$ which must be satisfied.

Consider the $4 \times 4$ matrix $F$ comprising the top $4$ rows of the $(\ell + 1) \times 4$ matrix in the above equation; as $F$ has non-zero kernel we must have $\det F = 0$. If we regard $\det F$ as a polynomial of degree $8$ in the $a_i$ and $b_i$, the terms involving ${a_2}^2$ are obtained by taking either the first or third entry in row $1$ and either the third or fourth entry in row $2$; hence the coefficient of ${a_2}^2$ is
\begin{eqnarray*}
&   & b_1 a_{\pi(2)} (b_{\pi(3)} a_3 b_{\pi(4)} b_4 - b_{\pi(3)} b_3 b_{\pi(4)} a_4) \\
& + & b_1 b_{\pi(2)} (b_{\pi(3)} b_3 a_{\pi(4)} a_4 - a_{\pi(3)} a_3 b_{\pi(4)} b_4) \\
& + & a_1 b_{\pi(2)} (a_{\pi(3)} b_3 b_{\pi(4)} b_4 - b_{\pi(3)} b_3 a_{\pi(4)} b_4).
\end{eqnarray*}
If $\pi(2) > 4$, the first two of the six terms in this expression contain $a_{\pi(2)} a_3$ and $a_{\pi(2)} a_4$ and none of the other terms involves $a_{\pi(2)}$; so $\det F$ is a non-zero polynomial. Thus we may assume $\pi(2) \leq 4$; similarly we may assume $\pi(3), \pi(4) \leq 4$, so that $\pi = \pi' \pi''$ where $\pi'$ permutes $\{ 1, 2, 3, 4 \}$ and $\pi''$ permutes $\{ 5, \dots, \ell + 1 \}$. Now if $\pi(2) \neq 1$, without loss of generality we may assume $\pi(3) = 1$, in which case the coefficient of ${b_1}^2$ in the coefficient of ${a_2}^2$ above is
\begin{eqnarray*}
& a_{\pi(2)} a_3 b_{\pi(4)} b_4 - a_{\pi(2)} b_3 b_{\pi(4)} a_4 + b_{\pi(2)} b_3 a_{\pi(4)} a_4 & \\
& = \begin{cases}
{a_3}^2{b_4}^2 - a_3a_4b_3b_4 + {a_4}^2{b_3}^2 & \hbox{if } \pi' = (1\ 2\ 3), \\
2a_3a_4b_3b_4 - {a_4}^2{b_3}^2                 & \hbox{if } \pi' = (1\ 2\ 4\ 3);
\end{cases} &
\end{eqnarray*}
so $\det F$ is a non-zero polynomial. Thus we may assume $\pi(2) = 1$. If $\pi' = (1\ 2)$ the coefficient of ${a_2}^2$ is
$$
2(a_1a_3b_1b_3{b_4}^2 - a_1a_4b_1{b_3}^2b_4) + {a_4}^2{b_1}^2{b_3}^2 - {a_3}^2{b_1}^2{b_4}^2,
$$
so $\det F$ is a non-zero polynomial. Thus we may assume $\pi' = (1\ 2)(3\ 4)$, in which case we find that $\det F$ is the zero polynomial. Hence if $\ell = 3$ we have $\pi = (1\ 2)(3\ 4)$. If however $\ell \geq 4$ then arguing similarly with the first three and the fifth rows of the $(\ell + 1) \times 4$ matrix above shows that unless $\pi$ interchanges $3$ and $5$ we have a non-zero polynomial which must be satisfied. We have thus proved our claim.

Therefore if $\ell = 3$ and $\pi \notin \langle (1\ 2)(3\ 4), (1\ 3)(2\ 4) \rangle$, or if $\ell \geq 4$ and $\pi \neq 1$, the points $y \in Y$ fixed by any such $n$ form a proper closed subvariety of $Y$; we take $\hat Y_2$ to be the intersection of the complements of these proper closed subvarieties as $\pi$ runs through $S_4 \setminus \langle (1\ 2)(3\ 4), (1\ 3)(2\ 4) \rangle$ or $S_{\ell + 1} \setminus \{ 1 \}$ according as $\ell = 3$ or $\ell \geq 4$. Then $\hat Y_2$ is a dense open subset of $Y$, as therefore is $\hat Y = \hat Y_1 \cap \hat Y_2$. Thus if $\ell \geq 4$, for all $y \in \hat Y$ we have $C_G(y) = C$.

Now assume $\ell = 3$, and take $y \in \hat Y$. As above we may now assume $a_1 = b_2 = 1$, $a_2 = b_1 = 0$, in which case $a_3, a_4, b_3, b_4 \neq 0$; we must also have $a_3b_4 \neq a_4b_3$, as otherwise $b_3D_1 - a_3D_2$ would be a non-zero matrix in $y$ with two diagonal entries equal to zero, contrary to the definition of $\hat Y_1$. Take $c_1, c_2, c_3, c_4 \in K^*$ satisfying ${c_4}^8 = \frac{a_3b_3}{a_4b_4(a_3b_4 - a_4b_3)}$, ${c_3}^4 = \frac{a_4b_4}{a_3b_3}{c_4}^4$, ${c_2}^4 = \frac{b_4}{a_3}(a_3b_4 - a_4b_3){c_4}^4$ and $c_1 = \frac{1}{c_2c_3c_4}$; write $\kappa_1 = \frac{a_3{c_3}^2}{{c_1}^2}$ and $\kappa_2 = \frac{a_4{c_4}^2}{{c_1}^2}$, then we have ${\kappa_1}^2 = {\kappa_2}^2 + 1$. If we now set $h^{-1} = \diag(c_1, c_2, c_3, c_4) \in G$ and $y' = h^{-1}.y$, we have $y' = \langle {D_1}', {D_2}' \rangle$ where ${D_1}' = \diag(1, 0, \kappa_1, \kappa_2)$ and ${D_2}' = \diag(0, 1, \kappa_2, \kappa_1)$. Now with $n^* = n_{\alpha_1} n_{\alpha_3}$ we see that $n^*$ sends the matrix $\diag(a, b, c, d)$ to $\diag(b, a, d, c)$, so we have $n^*.{D_1}' = {D_2}'$ and $n^*.{D_2}' = {D_1}'$, whence $n^* \in C_G(y')$; with $n^{**} = h_{\alpha_2 + \alpha_3}(\eta_4) n_{\alpha_1 + \alpha_2} n_{\alpha_2 + \alpha_3}$ we see that $n^{**}$ sends the matrix $\diag(a, b, c, d)$ to $\diag(c, -d, a, -b)$, so we have $n^{**}.{D_1}' = \kappa_1{D_1}' - \kappa_2{D_2}'$ and $n^{**}.{D_2}' = \kappa_2 {D_1}' - \kappa_1{D_2}'$, whence $n^{**} \in C_G(y')$. Hence $C_G(y') = Z_2 \langle n^*, n^{**} \rangle = C$, so $C_G(y) = C_G(h.y') = {}^h C$.

Therefore in all cases, for all $y \in \hat Y$ there exists $h \in T$ with $C_G(y) = {}^h C$. Since $\phi(G \times \hat Y$) contains a dense open subset of $X$, the quadruple $(G, \lambda, p, k)$ has generic stabilizer $C/Z(G) \cong {\Z_2}^2.S_3$, or ${\Z_2}^3.{\Z_2}^2$, or ${\Z_2}^\ell$, according as $\ell = 2$, or $\ell = 3$, or $\ell \geq 4$.
\end{proof}

\begin{prop}\label{prop: A_1, 3omega_1 module, k = 2}
Let $G = A_1$ and $\lambda = 3\omega_1$ with $p \geq 5$, and take $k = 2$. Then the quadruple $(G, \lambda, p, k)$ has generic stabilizer ${\Z_2}^2$.
\end{prop}

\begin{proof}
We take $G = \SL_2(K)$. Recall that $V_{nat} = \langle v_1, v_2 \rangle$. As in Proposition~\ref{prop: A_1, 3omega_1 module}, we may identify $V$ with $S^3(V_{nat})$, the space of homogeneous polynomials in $v_1$ and $v_2$ of degree $3$, so that $V = \langle {v_1}^3, {v_1}^2v_2, v_1{v_2}^2, {v_2}^3 \rangle$. For convenience write $x_1 = {v_1}^3$, $x_2 = {v_1}^2v_2$, $x_3 = v_1{v_2}^2$ and $x_4 = {v_2}^3$; then with respect to the basis $x_1, x_2, x_3, x_4$ the simple root element $x_{\alpha_1}(t)$ acts as
$$
\left(
  \begin{array}{cccc}
    1 & t & t^2 &  t^3 \\
      & 1 & 2t  & 3t^2 \\
      &   &  1  &  3t  \\
      &   &     &  1   \\
  \end{array}
\right),
$$
and for $i = 1, \dots, 4$ and $\kappa \in K^*$ we have $h_{\alpha_1}(\kappa).x_i = \kappa^{5 - 2i} x_i$. Set
$$
Y = \left\{ \langle a_1x_1 + a_3x_3, a_2x_2 + a_4x_4 \rangle : (a_1, a_3), (a_2, a_4) \neq (0, 0) \right\},
$$
and
$$
\hat Y = \left\{ \langle a_1x_1 + a_3x_3, a_2x_2 + a_4x_4 \rangle :a_1a_2a_3a_4 \neq 0, \ {\ts\frac{a_2a_3}{a_1a_4}} \neq 1, -3, 9 \right\};
$$
then $\hat Y$ is a dense open subset of $Y$. Take
$$
y = \langle v^{(1)}, v^{(2)} \rangle \in \hat Y, \qquad \hbox{where} \quad v^{(1)} = a_1x_1 + a_3x_3, \ v^{(2)} = a_2x_2 + a_4x_4.
$$

Take $g \in \Tran_G(y, Y)$, and write $g = u_1nu_2$ with $u_1 \in U$, $n \in N$ and $u_2 \in U_w$ where $w = nT \in W$. Since applying $h_{\alpha_1}(\kappa)$ does not alter the value of $\frac{a_2a_3}{a_1a_4}$, we see that $T$ preserves $\hat Y$; thus we may assume $n \in \{ 1, n_{\alpha_1} \}$. Write $u_1 = x_{\alpha_1}(t)$ and $u_2 = x_{\alpha_1}(t')$, and set $t'' = tt' - 1$.

First suppose $n = 1$; then $u_2 = 1$, so $t' = 0$. Since the coefficients of $x_2$ and $x_4$ in $g.v^{(1)}$ are $2a_3t$ and $0$ respectively, and that of $x_4$ in $g.v^{(2)}$ is $a_4$, for $g.y \in Y$ we must have $t = 0$ and hence $u_1 = 1$.

Now suppose $n = n_{\alpha_1}$. We have
\begin{eqnarray*}
g.v^{(1)} & = & -(a_1t^3 + a_3t{t''}^2)x_1 - (3a_1t^2 + a_3t''(3t'' + 2))x_2 \\
          &   & {} - (3a_1t + a_3t'(3t'' + 1))x_3 - (a_1 + a_3{t'}^2) x_4, \\
g.v^{(2)} & = & -(a_2t^2t'' + a_4{t''}^3)x_1 - (a_2t(3t'' + 1) + 3a_4t'{t''}^2)x_2 \\
          &   & {} - (a_2(3t'' + 2) + 3a_4{t'}^2t'')x_3 - (a_2t' + a_4{t'}^3)x_4.
\end{eqnarray*}
For $g.y \in Y$ the projections of $g.v^{(1)}$ and $g.v^{(2)}$ on $\langle x_1, x_3 \rangle$ must be linearly dependent, as must those on $\langle x_2, x_4 \rangle$; this gives the equations
\begin{eqnarray*}
0 & = & 2a_1a_2t^3 + (3a_1a_4 - a_2a_3)tt''(2tt' - 1) + 2a_3a_4t'{t''}^3, \\
0 & = & 2a_1a_2t + (3a_1a_4 - a_2a_3)t'(2tt' - 1) + 2a_3a_4{t'}^3t''.
\end{eqnarray*}
Multiplying the second equation by $t^2$ and subtracting the first gives
$$
0 = (2tt' - 1)[(3a_1a_4 - a_2a_3)t + 2a_3a_4t't''].
$$
If the second bracket here is zero we obtain $t[3a_1a_4 - a_2a_3 + 2a_3a_4{t'}^2] = 2a_3a_4t'$; this and the second equation above are now linear in $t$, and we obtain
$$
0 = [(3a_1a_4 - a_2a_3)^2 - 4a_1a_2a_3a_4]t' = (9a_1a_4 - a_2a_3)(a_1a_4 - a_2a_3)t'.
$$
By the definition of $\hat Y$ we must have $t' = 0$, and then $t = 0$. If instead $2tt' - 1 = 0$, then $t = \frac{1}{2t'}$, and from the second equation above we have $a_1a_2 = a_3a_4{t'}^4$; so there are $4$ possibilities for $t'$, and then $t$ is determined. In this case we have
\begin{eqnarray*}
g.v^{(1)} & = & -{\ts\frac{1}{8{t'}^3}}(a_1 + a_3{t'}^2)x_1 - {\ts\frac{1}{4{t'}^2}}(3a_1 - a_3{t'}^2)x_2 \\
          &   & {} - {\ts\frac{1}{2t'}}(3a_1 - a_3{t'}^2)x_3 - (a_1 + a_3{t'}^2) x_4, \\
g.v^{(2)} & = & {\ts\frac{1}{8{t'}^2}}(a_2 + a_4{t'}^2)x_1 + {\ts\frac{1}{4t'}}(a_2 - 3a_4{t'}^2)x_2 \\
          &   & {} - {\ts\frac{1}{2}}(a_2 - 3a_4{t'}^2)x_3 - t'(a_2 + a_4{t'}^2)x_4.
\end{eqnarray*}
Suppose that there exists $h = h_{\alpha_1}(\kappa) \in T$ with $hg.y = y$: equating the projection of $hg.v^{(1)}$ on $\langle x_1, x_3 \rangle$ to a scalar multiple of $a_1x_1 + a_3x_3$, and that of $hg.v^{(2)}$ on $\langle x_2, x_4 \rangle$ to a scalar multiple of $a_2x_2 + a_4x_4$, gives
\begin{eqnarray*}
-\kappa^3 \ts{\frac{1}{8{t'}^3}} a_3(a_1 + a_3{t'}^2) & = & -\ts{\frac{1}{\kappa}} \ts{\frac{1}{2t'}} a_1(3a_1 - a_3{t'}^2), \\
\kappa \ts{\frac{1}{4t'}} a_4(a_2 - 3a_4{t'}^2) & = & -\ts{\frac{1}{\kappa^3}} a_2t'(a_2 + a_4{t'}^2),
\end{eqnarray*}
whence
$$
\ts{\frac{a_1(3a_1 - a_3{t'}^2)}{a_3(a_1 + a_3{t'}^2)}} = \ts{\frac{\kappa^4}{4{t'}^2}} = -\ts{\frac{a_2(a_2 + a_4{t'}^2)}{a_4(a_2 - 3a_4{t'}^2)}};
$$
multiplying up, substituting for ${t'}^4$ and rearranging gives
$$
(3a_1a_4 + a_2a_3)[(3a_1a_4 - a_2a_3){t'}^2 - 2a_1a_2] = 0;
$$
by the definition of $\hat Y$ the first bracket is non-zero, so we must have
$$
(3a_1a_4 - a_2a_3){t'}^2 = 2a_1a_2,
$$
and now squaring and substituting for ${t'}^4$ gives
$$
(3a_1a_4 - a_2a_3)^2 = 4a_1a_2a_3a_4,
$$
which we have seen is contrary to the definition of $\hat Y$. Therefore the elements in $\Tran_G(y, Y)$ with $t = \frac{1}{2t'}$ do not in fact stabilize $y$; so $C_G(y) \leq N$.

Thus $\Tran_G(y, Y)$ is a finite union of cosets of $T$; so
$$
\codim {\ts\Tran_G(y, Y)} = \dim G - \dim {\ts\Tran_G(y, Y)} = 3 - 1 = 2
$$
while
$$
\codim Y = \dim \G{2}(V) - \dim Y = 4 - 2 = 2.
$$
Therefore $y$ is $Y$-exact.

Now take $\kappa \in K^*$ satisfying $\kappa^8 = \frac{a_3a_4}{a_1a_2}$ and set $s = h_{\alpha_1}(\kappa)$; then
$$
s.y = \langle \kappa^4 a_1x_1 + a_3x_3, a_3x_2 + \kappa^4 a_1x_4 \rangle.
$$
Then we have $C_T(s.y) = \langle h_{\alpha_1}(\eta_4) \rangle$ and $n_{\alpha_1} \in C_N(s.y)$. Set $C = \langle h_{\alpha_1}(\eta_4), n_{\alpha_1} \rangle$; then $C_G(s.y) = C$, and hence $C_G(y) = C^s$. Thus the conditions of Lemma~\ref{lem: generic stabilizer from exact subset} hold; so the quadruple $(G, \lambda, p, k)$ has generic stabilizer $C/Z(G) \cong {\Z_2}^2$.
\end{proof}

\begin{prop}\label{prop: A_1, omega_1 + q omega_1 module, k = 2}
Let $G = A_1$ and $\lambda = \omega_1 + q \omega_1$, and take $k = 2$. Then the quadruple $(G, \lambda, p, k)$ has generic stabilizer $\Z_2$.
\end{prop}

\begin{proof}
We take $G = \SL_2(K)$. As in Proposition~\ref{prop: A_ell, omega_1 + q omega_1 and omega_1 + q omega_ell modules}, we may identify $V$ with the space of $2 \times 2$ matrices $D$ over $K$, so that $A \in G$ maps $D \mapsto AD(A^{(q)})^T$. Set
$$
Y = \left\{ \left\langle
\left(
 \begin{array}{cc}
  a_1 & a_2 \\
   0  &  0  \\
 \end{array}
\right),
\left(
 \begin{array}{cc}
   0  &  0  \\
  a_3 & a_4 \\
 \end{array}
\right)
\right\rangle: (a_1, a_2), (a_3, a_4) \neq (0, 0) \right\},
$$
and
$$
\hat Y = \left\{ \left\langle
\left(
 \begin{array}{cc}
  a_1 & a_2 \\
   0  &  0  \\
 \end{array}
\right),
\left(
 \begin{array}{cc}
   0  &  0  \\
  a_3 & a_4 \\
 \end{array}
\right)
\right\rangle: a_1a_2a_3a_4 \neq 0, \ a_1a_4 - a_2a_3 \neq 0 \right\};
$$
then $\hat Y$ is a dense open subset of $Y$. Take
$$
y = \left\langle
\left(
 \begin{array}{cc}
  a_1 & a_2 \\
   0  &  0  \\
 \end{array}
\right),
\left(
 \begin{array}{cc}
   0  &  0  \\
  a_3 & a_4 \\
 \end{array}
\right)
\right\rangle \in \hat Y.
$$

Take $g \in \Tran_G(y, Y)$, and write $g = u_1nu_2$ with $u_1 \in U$, $n \in N$ and $u_2 \in U_w$ where $w = nT \in W$. Since $T$ preserves $\hat Y$ we may assume $n \in \{ 1, n_{\alpha_1} \}$. Write $u_1 = x_{\alpha_1}(t)$ and $u_2 = x_{\alpha_1}(t')$, and set $t'' = tt' - 1$.

First suppose $n = 1$; then $u_2 = 1$, so $t' = 0$. We have
$$
g.\left(
    \begin{array}{cc}
      a_1 & a_2 \\
       0  &  0  \\
    \end{array}
  \right)
=
\left(
  \begin{array}{cc}
    a_1 + a_2t^q & a_2 \\
          0      &  0  \\
  \end{array}
\right),
\quad
g.\left(
    \begin{array}{cc}
       0  &  0  \\
      a_3 & a_4 \\
    \end{array}
  \right)
=
\left(
  \begin{array}{cc}
    a_3t + a_4t^{q + 1} & a_4t \\
     a_3 + a_4t^q       & a_4  \\
  \end{array}
\right).
$$
For $g.y \in Y$ the two top rows must be linearly dependent; if we form the matrix comprising these two top rows, and subtract $t^q$ times the second column from the first, the resulting matrix is
$$
\left(
  \begin{array}{cc}
    a_1  & a_2  \\
    a_3t & a_4t \\
  \end{array}
\right),
$$
whose determinant is $(a_1a_4 - a_2a_3)t$, so we must have $t = 0$ and hence $u_1 = 1$.

Now suppose $n = n_{\alpha_1}$. We have
\begin{eqnarray*}
g.\left(
    \begin{array}{cc}
      a_1 & a_2 \\
       0  &  0  \\
    \end{array}
  \right)
& = &
\left(
  \begin{array}{cc}
    a_1t^{q + 1} + a_2t{t''}^q & a_1t + a_2t{t'}^q \\
          a_1t^q + a_2{t''}^q  &  a_1 + a_2{t'}^q  \\
  \end{array}
\right), \\
g.\left(
    \begin{array}{cc}
       0  &  0  \\
      a_3 & a_4 \\
    \end{array}
  \right)
& = &
\left(
  \begin{array}{cc}
     a_3t^qt'' + a_4{t''}^{q + 1} &  a_3t'' + a_4{t'}^qt''    \\
      a_3t^qt' + a_4t'{t''}^q     &   a_3t' + a_4{t'}^{q + 1} \\
  \end{array}
\right).
\end{eqnarray*}
Here we need both the top rows and the bottom rows to be linearly dependent; if we form the corresponding two matrices, and subtract $t^q$ times the second column from the first, negate the first column and subtract ${t'}^q$ times the first column from the second, the resulting matrices are
$$
\left(
  \begin{array}{cc}
    a_2t   & a_1t   \\
    a_4t'' & a_3t'' \\
  \end{array}
\right),
\quad
\left(
  \begin{array}{cc}
     a_2  &  a_1  \\
    a_4t' & a_3t' \\
  \end{array}
\right),
$$
whose determinants are $-(a_1a_4 - a_2a_3)t t''$ and $-(a_1a_4 - a_2a_3)t'$, so we must have $tt'' = t' = 0$, whence $t = 0$ and hence $u_1 = u_2 = 1$.

Thus $\Tran_G(y, Y) = N$; so
$$
\codim {\ts\Tran_G(y, Y)} = \dim G - {\ts\dim \Tran_G(y, Y)} = 3 - 1 = 2
$$
while
$$
\codim Y = \dim \G{2}(V) - \dim Y = 4 - 2 = 2.
$$
Therefore $y$ is $Y$-exact.

Now take $\kappa \in K^*$ satisfying $\kappa^{4q} = -\frac{a_2a_4}{a_1a_3}$, and take $s = h_{\alpha_1}(\kappa)$; then
$$
s.y =
\left\langle
\left(
 \begin{array}{cc}
  \kappa^{2q}a_1 & a_2 \\
         0       &  0  \\
 \end{array}
\right),
\left(
 \begin{array}{cc}
   0  &         0       \\
  a_2 & -\kappa^{2q}a_1 \\
 \end{array}
\right)
\right\rangle.
$$
Then we have $C_T(s.y) = \{ \pm I \} = Z(G)$ and $n_{\alpha_1} \in C_N(s.y)$. Set $C = \langle n_{\alpha_1} \rangle$; then $C_G(s.y) = C$, and hence $C_G(y) = C^s$. Thus the conditions of Lemma~\ref{lem: generic stabilizer from exact subset} hold; so the quadruple $(G, \lambda, p, k)$ has generic stabilizer $C/Z(G) \cong \Z_2$.
\end{proof}

For the next few results we shall treat separately the cases where $p$ is coprime to $k$ and where $p$ divides $k$. Much as with the proofs of Propositions~\ref{prop: A_8, omega_3, A_7, omega_4, D_8, omega_8 modules, non-special characteristic} and \ref{prop: A_8, omega_3, A_7, omega_4, D_8, omega_8 modules, special characteristic}, for the former we shall use the approach of Section~\ref{sect: semisimple auts}, here combined with that of Section~\ref{sect: reduction}, while for the latter we shall instead use the approach of Section~\ref{sect: Lie algebra annihilators}, combined with Lemma~\ref{lem: generic stabilizer from exact subset}.

\begin{prop}\label{prop: A_2, 2omega_1 module, k = 3, A_4, omega_2 module, k = 5}
Let $G = A_2$ and $\lambda = 2\omega_1$ with $p \geq 3$, and take $k = 3$, or let $G = A_4$ and $\lambda = \omega_2$, and take $k = 5$. Then the quadruple $(G, \lambda, p, k)$ has generic stabilizer $\Z_{3/(p, 3)}.S_3$ or $\Z_{5/(p, 5)}.Dih_{10}$ respectively.
\end{prop}

\begin{proof}
Number the cases (i) and (ii) according as $G = A_2$ or $A_4$; whenever we give two choices followed by the word \lq respectively' we are taking the cases in the order (i), (ii). We shall deal separately with the cases $p \neq k$ and $p = k$.

First suppose $p \neq k$. Let $H$ be the (simply connected) group defined over $K$ of type $F_4$ or $E_8$ respectively (so that $\ell_H = 2\ell$), with simple roots $\beta_1, \dots, \beta_{2\ell}$. We have $Z(\L(H)) = \{ 0 \}$. Our strategy will be to identify the group $G^+ = G^2 = GA_{k - 1}$ as the centralizer of a semisimple automorphism of $H$, and use Lemma~\ref{lem: semisimple auts} to find the generic stabilizer in the action of $G^+$ on $\G{1}(V^+)$, where $V^+ = V \otimes V_{nat}$ with $V$ the $G$-module with high weight $2\omega_1$ or $\omega_2$ respectively and $V_{nat}$ the natural $A_{k - 1}$-module; we shall therefore employ the notation of Section~\ref{sect: semisimple auts}. Lemma~\ref{lem: reduction to projective space} will then give the result.

Define $\theta_1 \in T_H$ to be
$$
\begin{array}{ll}
h_{\beta_1}({\eta_3}^2) h_{\beta_4}({\eta_3}^2)                                                                                         & \vstrut \hbox{in case~(i),} \\
h_{\beta_1}(\eta_5) h_{\beta_2}({\eta_5}^3) h_{\beta_3}(\eta_5) h_{\beta_6}({\eta_5}^4) h_{\beta_7}({\eta_5}^2) h_{\beta_6}({\eta_5}^4) & \vstrut \hbox{in case~(ii).}
\end{array}
$$
Then ${\theta_1}^k = 1$, and $\theta_1$ sends $x_\alpha(t)$ to $x_\alpha({\eta_k}^{\height(\alpha)}t)$; so $X_\alpha < C_H(\theta_1)$ if and only if $\height(\alpha) \equiv 0$ (mod $k$). It follows that $C_H(\theta_1)$ is a connected group of type $\tilde A_2 A_2$ or ${A_4}^2$ respectively, with simple root elements $x_{\alpha_i}(t)$, where $\alpha_1, \dots, \alpha_{2\ell}$ are
$$
\begin{array}{ll}
\ffourrt0111, \ffourrt1110, \ffourrt0120, \ffourrt1122                                                                                                 & \vstrut \hbox{in case~(i),} \\
\eeightrt00111110, \eeightrt11111000, \eeightrt00011111, \eeightrt01111100, \eeightrt01121000, \eeightrt10111100, \eeightrt01011110, \eeightrt11221111 & \vstrut \hbox{in case~(ii);}
\end{array}
$$
in each case we see that $Z(C_H(\theta_1)) = \langle \theta_1 \rangle$.

Now let $\delta_1, \dots, \delta_{2\ell}$ be
$$
\begin{array}{ll}
\ffourrt0001, \ffourrt0010, -\ffourrt2342, \ffourrt1000                                                                                                 & \vstrut \hbox{in case~(i),} \\
\eeightrt10000000, \eeightrt00100000, \eeightrt00010000, \eeightrt01000000, -\eeightrt23465432, \eeightrt00000001, \eeightrt00000010, \eeightrt00000100 & \vstrut \hbox{in case~(ii);}
\end{array}
$$
and set $\theta_2 = n_{\delta_1} \dots n_{\delta_{2\ell}}$; then ${\theta_2}^k = 1$. We find that $\theta_2$ acts fixed-point-freely on both $\Phi_H$ and $\L(T_H)$. Thus $\dim C_{\L(H)}(\theta_2) = |\Phi_H|/k = 16$ or $48$ respectively; the classification of semisimple elements of $H$ (see e.g. \cite[Table~4.7.1]{GLS}) now shows that $\theta_2$ must be a conjugate of $\theta_1$.

First set $\theta = \theta_1$; then we may take $G^+ = C_H(\theta)$. We see that $e_\alpha \in \L(H)_{(i)}$ if and only if $\height(\alpha) \equiv i$ (mod $k$). Thus in $\L(H)_{(1)}$ we have a highest weight vector $e_\beta$ for $\beta = \ffourrt1342$ or $\eeightrt23465321$ respectively; the expressions above for the simple root elements of $G^+$ show that $\L(H)_{(1)}$ is the Weyl $G^+$-module with high weight $2\omega_1 \otimes \omega_1$ or $\omega_2 \otimes \omega_1$ respectively, i.e., the tensor product $V \otimes V_{nat}$. As $Z(\L(H)) = \{ 0 \}$, we may therefore take $V^+ = \L(H)_{(1)}/Z(\L(H))_{(1)}$; of course ${G^+}_{\G{1}(V^+)} = Z(G^+)$.

Now set $\theta = \theta_2$, and again take $G^+ = C_H(\theta)$ and $V^+ = \L(H)_{(1)}/Z(\L(H))_{(1)}$. We have $G^+ \cap T_H = C_{T_H}(\theta) \cong {\Z_k}^2$; indeed this group is
$$
\begin{array}{ll}
\langle h_{\beta_1}(\eta_3) h_{\beta_4}(\eta_3), h_{\beta_3}(\eta_3) h_{\beta_4}({\eta_3}^2) \rangle & \vstrut \hbox{in case~(i),} \\
\langle h_{\beta_1}({\eta_5}^4) h_{\beta_2}(\eta_5) h_{\beta_3}({\eta_5}^3) h_{\beta_4}({\eta_5}^2), & \vstrut \\
\ h_{\beta_1}({\eta_5}^2) h_{\beta_3}(\eta_5) h_{\beta_4}({\eta_5}^2) h_{\beta_6}(\eta_5) h_{\beta_7}({\eta_5}^3) h_{\beta_8}(\eta_5) \rangle & \vstrut \hbox{in case~(ii).} \\
\end{array}
$$
Moreover we find that $\L(T_H)_{(1)} = \langle h_{\delta_{2i - 1}} - \eta_3 h_{\delta_{2i}} : i = 1, 2 \rangle$ or $\langle h_{\delta_{4i - 3}} - (\eta_5 + {\eta_5}^2 + {\eta_5}^3)h_{\delta_{4i - 2}} - (\eta_5 + {\eta_5}^2)h_{\delta_{4i - 1}} - \eta_5h_{\delta_{4i}} : i = 1, 2 \rangle$ respectively. Thus $\dim \L(H)_{(1)} - \dim \L(T_H)_{(1)} = \dim G^+ - \dim (G^+ \cap T_H)$ in each case; and in each case a routine check shows that $\L(T_H)_{(1)}$ contains regular semisimple elements.

We claim that in each case we have $({W_H}^\ddagger)_{(1)} = \langle \theta T_H, w_0 \rangle$. Thus suppose $w \in W_H$ and there exists $\xi \in K^*$ such that for all $y \in \L(T_H)_{(1)}$ we have $w.y = \xi y$. For $i = 1, 2$ write $\Psi_i = \langle \delta_{2i - 1}, \delta_{2i} \rangle$ or $\langle \delta_{4i - 3}, \delta_{4i - 2}, \delta_{4i - 1}, \delta_{4i} \rangle$ respectively, so that $\Psi_i$ is of type $A_{k - 1}$. In case (i), taking $y = h_{\delta_1} - \eta_3 h_{\delta_2} \in \L(T_H)_{(1)}$ and arguing as in the paragraphs following the statement of Lemma~\ref{lem: semisimple auts} shows that $w(\beta_3)$ and $w(\beta_4)$ must be proportional outside $\{ \beta_3, \beta_4 \}$, and as $\eta_3 \neq \pm1$ that $w$ must preserve $\Psi_1$. We shall prove that the same conclusion holds in case (ii).

Take $y = h_{\delta_1} - (\eta_5 + {\eta_5}^2 + {\eta_5}^3)h_{\delta_2} - (\eta_5 + {\eta_5}^2)h_{\delta_3} - \eta_5h_{\delta_4}$ and for $j = 1, 2, 3, 4$ write $w(\beta_j) = \sum a_{ij} \beta_i$. First consider the coefficients $a_{8j}$; each lies in $\{ 0, \pm1, \pm2 \}$. Write $\rho = \eeightrt23465432$ for the high root of $\Phi_H$; then $\rho$ is the only root whose $\beta_8$-coefficient is $2$, and any root $\beta$ such that $\rho + \beta \in \Phi_H$ has $\beta_8$-coefficient equal to $-1$. As a result we see that if $\rho \in w(\Psi_1)$ then up to negation the $4$-tuple $(a_{81}, a_{82}, a_{83}, a_{84})$ must be such that either one term is $2$, or two terms are $1$ with any intermediate terms being $0$; moreover any term adjacent to the $2$ or to one of the $1$s (and not between them) is $-1$, and all other terms are $0$. If instead $\rho \notin w(\Psi_1)$ but some $a_{8j}$ is $\pm1$ then the non-zero terms in the $4$-tuple must alternate in sign. However, we know that the coefficient of $h_{\beta_8}$ in $w.y$ is $0$. In the first possibility this condition gives ${\eta_5}^{i_1} = {\eta_5}^{i_2}$ for some $0 \leq i_1 < i_2 \leq 4$; in the second it gives $\sum_{i \in S} {\eta_5}^i = 0$ where $S$ is a non-empty proper subset of $\{ 0, 1, 2, 3, 4 \}$, and using $1 + \eta_5 + \cdots + {\eta_5}^4 = 0$ we may assume $|S| = 1$ or $2$, so either ${\eta_5}^i = 0$ or ${\eta_5}^{i_1} + {\eta_5}^{i_2} = 0$ for some $0 \leq i_1 < i_2 \leq 4$. As each of these is impossible, all $a_{8j}$ must be zero, so for each $j \leq 4$ we have $w(\beta_j) \in \langle \beta_1, \dots, \beta_7 \rangle$. Now arguing similarly with the coefficients $a_{7j}$ (but ignoring the possibility that some root has $\beta_7$-coefficient equal to $2$) shows that for each $j \leq 4$ we have $w(\beta_j) \in \langle \beta_1, \dots, \beta_6 \rangle$; likewise treating the $a_{6j}$ and then the $a_{5j}$ we conclude as required that $w$ preserves $\Psi_1$.

In both cases $w$ must therefore also preserve the set of roots orthogonal to $\Psi_1$, which is $\Psi_2$. Thus $w = w_1w_2.{w_0}^j$ where each $w_i$ lies in $W(\Psi_i)$ and $j \in \{ 0, 1 \}$. For each $i$ take the group $A_{k - 1}$ with root system $\Psi_i$, and write the elements of its Lie algebra as $k \times k$ matrices in the usual way; then the corresponding basis vector of $\L(T_H)_{(1)}$ has matrix $\diag(1, {\eta_3}^2, \eta_3)$ or $\diag(1, {\eta_5}^4, {\eta_5}^3, {\eta_5}^2, \eta_5)$ respectively. As $w$ acts as a scalar, for each $i$ the element $w_i$ must be some power of $w_{\delta_{2i - 1}} w_{\delta_{2i}}$ or $w_{\delta_{4i - 3}} w_{\delta_{4i - 2}} w_{\delta_{4i - 1}} w_{\delta_{4i}}$ respectively; as the two scalars must be equal, we must have $w_1w_2 \in \langle \theta T_H \rangle$, so that $w \in \langle \theta T_H, w_0 \rangle$ as required. Note that if we write $n_0 = n_{-\ffourrt0100} n_{\ffourrt0120} n_{\ffourrt0122} n_{\ffourrt2342}$ or $n_{\beta_2} n_{\beta_3} n_{\beta_5} n_{\beta_7} n_{\rho_4} n_{\rho_6}n_{\rho_7} n_{-\rho}$ respectively (with $\rho_4 = \eeightrt01121000$, $\rho_6 = \eeightrt01122210$ and $\rho_7 = \eeightrt22343210$ in case (ii)), then $n_0$ is an involution in $N_H$ corresponding to $w_0$ which commutes with $\theta$.

Now $({N_H}^\ddagger)_{(1)} = T_H \langle \theta, n_0 \rangle$, so $C_{({N_H}^\ddagger)_{(1)}}(\theta) = C_{T_H}(\theta) \langle \theta, n_0 \rangle$. Since ${G^+}_{\G{1}(V^+)} = Z(G^+) = \langle \theta \rangle$, Lemma~\ref{lem: semisimple auts} shows that in the action of $G^+$ on $\G{1}(V^+)$ the generic stabilizer is $C_{({N_H}^\ddagger)_{(1)}}(\theta)/Z(G^+) \cong {\Z_k}^2.\Z_2 \cong \Z_k.Dih_{2k}$. Finally Lemma~\ref{lem: reduction to projective space} shows that if $p \neq k$ the quadruple $(G, \lambda, p, k)$ also has generic stabilizer $\Z_k.Dih_{2k}$.

Now suppose instead $p = k$; here $Z(G) = \{ 1 \}$. We shall write elements of both $G$ and $\L(G)$ as $k \times k$ matrices. We let $T < G$ be the subgroup of diagonal matrices, so that $N$ is the subgroup of monomial matrices; define $n_0, n_1 \in N$ by
$$
n_0 =
-\left(
  \begin{array}{ccc}
      &   & 1 \\
      & 1 &   \\
    1 &   &   \\
  \end{array}
\right),
\qquad
n_1 =
\left(
  \begin{array}{ccc}
      &   & 1 \\
    1 &   &   \\
      & 1 &   \\
  \end{array}
\right)
$$
or
$$
n_0 =
\left(
  \begin{array}{ccccc}
      &   &   &   & 1 \\
      &   &   & 1 &   \\
      &   & 1 &   &   \\
      & 1 &   &   &   \\
    1 &   &   &   &   \\
  \end{array}
\right),
\qquad
n_1 =
\left(
  \begin{array}{ccccc}
      &   &   &   & 1 \\
    1 &   &   &   &   \\
      & 1 &   &   &   \\
      &   & 1 &   &   \\
      &   &   & 1 &   \\
  \end{array}
\right)
$$
respectively, so that $n_0T$ is the long word $w_0$ of the Weyl group, and if we identify $W$ with the symmetric group $S_k$ then $n_1T$ is the $k$-cycle $(1 \ \dots \ k)$.

Recall the natural module $V_{nat}$ with basis $v_1, \dots, v_k$. In case (i) we have $V = S^2(V_{nat})$, with basis $v_{11}, v_{22}, v_{33}, v_{12}, v_{23}, v_{31}$, where we write $v_{ii} = v_i \otimes v_i$ and $v_{ij} = v_i \otimes v_j + v_j \otimes v_i$ if $i \neq j$. In case (ii) we have $V = \bigwedge^2(V_{nat})$, with basis $v_{12}, v_{23}, v_{34}, v_{45}, v_{51}, v_{14}, v_{25}, v_{31}, v_{42}, v_{53}$, where we write $v_{ij} = v_i \wedge v_j$. Note that in each case $\dim V = 2k$, so $\dim \Gk(V) = k^2$. Write
\begin{eqnarray*}
V^{(1)} & = & \langle v_{11}, v_{23} \rangle, \\
V^{(2)} & = & \langle v_{22}, v_{31} \rangle, \\
V^{(3)} & = & \langle v_{33}, v_{12} \rangle
\end{eqnarray*}
or
\begin{eqnarray*}
V^{(1)} & = & \langle v_{34}, v_{25} \rangle, \\
V^{(2)} & = & \langle v_{45}, v_{31} \rangle, \\
V^{(3)} & = & \langle v_{51}, v_{42} \rangle, \\
V^{(4)} & = & \langle v_{12}, v_{53} \rangle, \\
V^{(5)} & = & \langle v_{23}, v_{14} \rangle
\end{eqnarray*}
respectively, so that $V = V^{(1)} \oplus \cdots \oplus V^{(k)}$ and $n_1$ cycles the $V^{(i)}$. Define
$$
Y = \{ y = \langle v^{(1)}, \dots, v^{(k)} \rangle : v^{(1)} \in V^{(1)} \setminus \{ 0 \}, \ v^{(2)} = n_1.v^{(1)}, \ \dots, \ v^{(k)} = n_1.v^{(k - 1)} \};
$$
then $Y$ is a subvariety of $\Gk(V)$ of dimension $1$, whence $\codim Y = k^2 - 1 = \dim G$. Set $C = \langle n_1, n_0 \rangle$; then each $y \in Y$ is stabilized by $C$. For convenience, given $y = \langle v^{(1)}, \dots, v^{(k)} \rangle \in Y$ with $v^{(1)} = a_1v_{11} + a_2v_{23}$ or $a_1v_{34} + a_2v_{25}$ respectively, we shall write $y = y_\a$ where $\a = (a_1, a_2)$.

Define $\S = \langle I, h_0 \rangle \leq \L(T)$, where $h_0 = \diag(1, 0, -1)$ or $\diag(2, 1, 0, -1, -2)$ respectively (so that in case (i) we actually have $\S = \L(T)$). Clearly if $\alpha \in \Phi$ then $[h_0 e_\alpha] \neq 0$, so $C_{\L(G)}(\S) = \L(T)$. If $y \in Y$ we have $\S \leq \Ann_{\L(G)}(y)$. Set
$$
\hat Y = \{ y_\a \in Y : a_1a_2 \neq 0, \ a_1 \neq \pm2a_2 \};
$$
then $\hat Y$ is a dense open subset of $Y$. Take $y = y_\a \in \hat Y$.

First suppose $x \in \Ann_{\L(G)}(y)$; write $x = h + e$ where $h \in \L(T)$ and $e \in \langle e_\alpha : \alpha \in \Phi \rangle$. Clearly for each $i$ the vector $h.v^{(i)}$ lies in $V^{(i)}$; since the difference of the two weights lying in $V^{(i)}$ is not a root, the projection of $x.v^{(i)}$ on $V^{(i)}$ is equal to $h.v^{(i)}$. Thus for each $i$, the vector $h.v^{(i)}$ must be a scalar multiple of $v^{(i)}$, while for each $j \neq i$ the projection of $e.v^{(i)}$ on $V^{(j)}$ must be a scalar multiple of $v^{(j)}$. A quick calculation (needed only in case (ii)) shows that we must have $h \in \S$. Now write $e = \sum_{\alpha \in \Phi} t_\alpha e_\alpha$; then the condition on the projections of the vectors $e.v^{(i)}$ on the $V^{(j)}$ may be expressed in matrix form as $A{\bf t} = {\bf 0}$, where $A$ is an $M \times M$ matrix and ${\bf t}$ is a column vector whose entries are the various coefficients $t_\alpha$. We find that if the rows and columns of $A$ are suitably ordered then it becomes block diagonal, having $2$ or $4$ blocks respectively, with each block being a $k \times k$ matrix. In fact each block may be written in the form
$$
\left(
  \begin{array}{ccc}
    {a_1}^2 &  a_1a_2 & {a_2}^2 \\
    {a_2}^2 & {a_1}^2 &  a_1a_2 \\
     a_1a_2 & {a_2}^2 & {a_1}^2 \\
  \end{array}
\right)
\hbox{ or }
\left(
  \begin{array}{ccccc}
     {a_1}^2 &  a_1a_2  & -{a_2}^2 &          &          \\
             &  {a_1}^2 &  a_1a_2  & -{a_2}^2 &          \\
             &          &  {a_1}^2 &  a_1a_2  & -{a_2}^2 \\
    -{a_2}^2 &          &          &  {a_1}^2 &  a_1a_2  \\
     a_1a_2  & -{a_2}^2 &          &          &  {a_1}^2 \\
  \end{array}
\right)
$$
respectively, which has determinant $(2a_1 + a_2)^{2k}$. Thus in each case the definition of the set $\hat Y$ implies that each block of $A$ is non-singular, as therefore is $A$ itself; so ${\bf t}$ must be the zero vector and hence $e = 0$. Thus $x = h + e \in \S$; so $\Ann_{\L(G)}(y) = \S$.

A straightforward calculation shows that $C_T(y) = \{ 1 \}$, and $T.y \cap Y = \{ y \}$. We claim that $N.y \cap Y \subset \hat Y$, and $C_N(y) = C$. In case (i) both claims are immediate, as each element of $N$ is of the form $sc$ for $s \in T$ and $c \in C$; so assume we are in case (ii) and take $n \in \Tran_N(y, Y)$. Since the projection on $W$ of the group $C$ acts transitively on $\{ 1, \dots, 5 \}$, and $n_0T$ fixes $3$ while acting transitively on each of $\{ 1, 5 \}$ and $\{ 2, 4 \}$, there exists $c \in C$ such that $ncT$ fixes $3$ and sends $1$ to either $1$ or $2$. Then $nc$ must take $v^{(3)}$ to an element of $V^{(3)}$; according as $ncT$ sends $1$ to $1$ or $2$ it must send $5$ to $5$ or $4$, so it must be $1$, $(2\ 4)$, $(1\ 2)(4\ 5)$ or $(1\ 2\ 5\ 4)$. If it is the second or third of these elements, $nc$ does not send any other $v^{(i)}$ into any $V^{(j)}$; however, if we set
$$
n' =
-\left(
   \begin{array}{ccccc}
       &   &   & 1 &   \\
     1 &   &   &   &   \\
       &   & 1 &   &   \\
       &   &   &   & 1 \\
       & 1 &   &   &   \\
   \end{array}
 \right),
$$
then $n'.y = y_{\a'}$ where $\a' = (-a_2, a_1)$. Since $y_{\a'} \in \hat Y$ this proves the first claim; moreover the definition of $\hat Y$ ensures that $n'$ does not stabilize $y$, and so if $n \in C_N(y)$ then $n = c^{-1} \in C$, proving the second. Thus the conditions of Lemma~\ref{lem: exactness via Premet} hold, so that $\Tran_G(y, Y) \subseteq N$, and $y$ is $Y$-exact; moreover $C_G(y) = C_N(y) = C$. Therefore the conditions of Lemma~\ref{lem: generic stabilizer from exact subset} hold; so if $p = k$ the quadruple $(G, \lambda, p, k)$ has generic stabilizer $C/Z(G) \cong Dih_{2k}$.
\end{proof}

\begin{prop}\label{prop: D_5, omega_5 module, k = 4}
Let $G = D_5$ and $\lambda = \omega_5$, and take $k = 4$. Then the quadruple $(G, \lambda, p, k)$ has generic stabilizer ${\Z_{2/(p, 2)}}^2.{\Z_2}^2$.
\end{prop}

\begin{proof}
We shall deal separately with the cases $p \geq 3$ and $p = 2$.

First suppose $p \geq 3$. Let $H$ be the (simply connected) group defined over $K$ of type $E_8$, with simple roots $\beta_1, \dots, \beta_8$. We have $Z(\L(H)) = \{ 0 \}$. Our strategy will be to identify the group $G^+ = D_5A_3 = GA_{k - 1}$ as the centralizer of a semisimple automorphism of $H$, and use Lemma~\ref{lem: semisimple auts} to find the generic stabilizer in the action of $G^+$ on $\G{1}(V^+)$, where $V^+ = V \otimes V_{nat}$ with $V$ the $D_5$-module with high weight $\omega_5$ and $V_{nat}$ the natural $A_3$-module; we shall therefore employ the notation of Section~\ref{sect: semisimple auts}. Lemma~\ref{lem: reduction to projective space} will then give the result.

Define $\theta_1 \in T_H$ to be
$$
h_{\beta_1}(-1) h_{\beta_3}(-\eta_4) h_{\beta_4}(-\eta_4) h_{\beta_5}(-1) h_{\beta_7}(\eta_4) h_{\beta_8}(\eta_4).
$$
Then ${\theta_1}^4 = 1$, and $\theta_1$ sends $x_\alpha(t)$ to $x_\alpha({\eta_4}^{\height(\alpha)}t)$; so $X_\alpha < C_H(\theta_1)$ if and only if $\height(\alpha) \equiv 0$ (mod $4$). It follows that $C_H(\theta_1)$ is a connected group of type $D_5A_3$, with simple root elements $x_{\alpha_i}(t)$, where $\alpha_1, \dots, \alpha_8$ are
$$
\eeightrt01111000, \eeightrt00011110, \eeightrt11110000, \eeightrt00111100, \eeightrt00001111, \eeightrt01121111, \eeightrt10111000, \eeightrt01011100;
$$
we see that $Z(C_H(\theta_1)) = \langle \theta_1 \rangle$.

Now let $\delta_1, \dots, \delta_8$ be
$$
\eeightrt23465432, \eeightrt00000010, \eeightrt00000011, \eeightrt01122211, \eeightrt00001000, \eeightrt01000000, \eeightrt01010000, \eeightrt00110000,
$$
and set $\theta_2 = n_{\delta_1} \dots n_{\delta_8}$; then ${\theta_2}^4 = 1$, and indeed the element of $W_H$ corresponding to ${\theta_2}^2$ is the long word. We find that $\theta_2$ acts fixed-point-freely on both $\Phi_H$ and $\L(T_H)$. Thus $\dim C_{\L(H)}(\theta_2) = |\Phi_H|/4 = 60$; the classification of semisimple elements of $H$ (see e.g. \cite[Table~4.3.1]{GLS}) now shows that $\theta_2$ must be a conjugate of $\theta_1$.

First set $\theta = \theta_1$; then we may take $G^+ = C_H(\theta)$. We see that $e_\alpha \in \L(H)_{(i)}$ if and only if $\height(\alpha) \equiv i$ (mod $4$). Thus in $\L(H)_{(1)}$ we have a highest weight vector $e_\beta$ for $\beta = \eeightrt23465432$; the expressions above for the simple root elements of $G^+$ show that $\L(H)_{(1)}$ is the Weyl $G^+$-module with high weight $\omega_5 \otimes \omega_1$, i.e., the tensor product $V \otimes V_{nat}$. As $Z(\L(H)) = \{ 0 \}$, we may therefore take $V^+ = \L(H)_{(1)}/Z(\L(H))_{(1)}$; of course ${G^+}_{\G{1}(V^+)} = Z(G^+)$.

Now set $\theta = \theta_2$, and again take $G^+ = C_H(\theta)$ and $V^+ = \L(H)_{(1)}/Z(\L(H))_{(1)}$. We have $G^+ \cap T_H = C_{T_H}(\theta) \cong {\Z_2}^4$; indeed this group is
$$
\langle h_{\beta_2}(-1) h_{\beta_3}(-1), h_{\beta_2}(-1) h_{\beta_5}(-1), h_{\beta_2}(-1) h_{\beta_7}(-1), h_{\beta_4}(-1) h_{\beta_8}(-1) \rangle.
$$
Moreover we find that
\begin{eqnarray*}
\L(T_H)_{(1)} & = & \langle \eta_4 h_{\beta_2} + h_{\beta_3}, \\
              &   & \phantom{\langle} h_{\beta_3} + (1 - \eta_4) h_{\beta_4} + h_{\beta_5}, \\
              &   & \phantom{\langle} h_{\beta_2} + (1 + \eta_4) h_{\beta_4} + (1 - \eta_4) h_{\beta_6} + h_{\beta_7}, \\
              &   & \phantom{\langle} (1 - \eta_4) h_{\beta_1} + \eta_4 h_{\beta_4} + (1 + \eta_4) h_{\beta_6} + h_{\beta_8} \rangle.
\end{eqnarray*}
Thus $\dim \L(H)_{(1)} - \dim \L(T_H)_{(1)} = \dim G^+ - \dim (G^+ \cap T_H)$; and a routine check shows that $\L(T_H)_{(1)}$ contains regular semisimple elements.

We claim that we have $({W_H}^\ddagger)_{(1)} = \langle \theta T_H \rangle$. Thus suppose $w \in W_H$ and there exists $\xi \in K^*$ such that for all $y \in \L(T_H)_{(1)}$ we have $w.y = \xi y$. Taking $y = \eta_4 h_{\beta_2} + h_{\beta_3}$ and arguing as in the paragraphs following the statement of Lemma~\ref{lem: semisimple auts} shows that $w(\beta_2)$ and $w(\beta_3)$ must be proportional outside $\{ \beta_2, \beta_3 \}$, and as $\eta_4 \neq \pm1$ that $w$ must preserve the ${A_1}^2$ subsystem $\langle \beta_2, \beta_3 \rangle$. Since $\theta$ acts on $\L(T_H)_{(1)}$ as multiplication by $\eta_4$, and sends $\beta_2$ to $\beta_3$ and $\beta_3$ to $-\beta_2$, by multiplying $w$ by a power of $\theta$ we may assume $w$ fixes $\beta_2$; thus $\xi = 1$ and $w$ must also fix $\beta_3$. Next taking $y = h_{\beta_3} + (1 - \eta_4) h_{\beta_4} + h_{\beta_5}$, whose first term is now fixed by $w$, and arguing as above again shows that $w(\beta_4)$ and $w(\beta_5)$ must be proportional outside $\{ \beta_4, \beta_5 \}$. Unless $1 - \eta_4 = -1$ (which is possible if $p = 5$) we conclude as above that $w$ preserves the $A_2$ subsystem $\langle \alpha_4, \alpha_5 \rangle$; if $1 - \eta_4 = -1$ then we see that $w(\beta_4)$ and $w(\beta_5)$ must actually be equal outside $\{ \beta_4, \beta_5 \}$, and now as $w(\beta_4) + w(\beta_5)$ is a root we again draw the same conclusion. Consideration of the roots in $\langle \alpha_4, \alpha_5 \rangle$ which can be added to both $\beta_2$ and $\beta_3$ (as $\beta_4$ can), and which are orthogonal to both (as $\beta_5$ is), quickly shows that $w$ must fix both $\beta_4$ and $\beta_5$. Now taking $y = h_{\beta_2} + (1 + \eta_4) h_{\beta_4} + (1 - \eta_4) h_{\beta_6} + h_{\beta_7}$, whose first two terms are now fixed by $w$, and arguing as above once more shows that $w(\beta_6)$ and $w(\beta_7)$ must be proportional outside $\{ \beta_6, \beta_7 \}$; similarly we see that $w$ must fix both $\beta_6$ and $\beta_7$. Finally taking $y = (1 - \eta_4) h_{\beta_1} + \eta_4 h_{\beta_4} + (1 + \eta_4) h_{\beta_6} + h_{\beta_8}$, whose second and third terms are now fixed by $w$, and arguing as above yet again shows that $w(\beta_1)$ and $w(\beta_8)$ must be proportional outside $\{ \beta_1, \beta_8 \}$. Both $w(\beta_1)$ and $w(\beta_8)$ must be orthogonal to $\langle \beta_2, \beta_4, \beta_5, \beta_6 \rangle$, so must lie in the $A_4$ subsystem $\langle \eeightrt01122221, \eeightrt10000000, \eeightrt12343210, \eeightrt00000001 \rangle$; in here, the roots orthogonal to $\beta_3$ lie in the $A_3$ subsystem $\langle \eeightrt01122221, \eeightrt22343210, \eeightrt00000001 \rangle$, while those orthogonal to $\beta_7$ lie in the $A_2A_1$ subsystem $\langle \eeightrt10000000, \eeightrt12343210, \eeightrt23465432 \rangle$. As $w(\beta_1)$ is orthogonal to $\beta_7$ and can be added to $\beta_3$, while $w(\beta_8)$ is orthogonal to $\beta_3$ and can be added to $\beta_7$, we must have
$$
w(\beta_1) \in \{ \eeightrt10000000, -\eeightrt12343210 \}, \quad w(\beta_8) \in \{ \eeightrt00000001, -\eeightrt01122221, \eeightrt22343211, -\eeightrt23465431 \}.
$$
Proportionality now forces $(w(\beta_1), w(\beta_8)) = (\beta_1, \beta_8)$ or $(-\eeightrt12343210, \eeightrt22343211)$; however in the latter case we must have $1 - \eta_4 = 1$, which is impossible. Thus $w$ also fixes $\beta_1$ and $\beta_8$, so equals $1$, proving the claim.

Now $({N_H}^\ddagger)_{(1)} = T_H \langle \theta \rangle$, so $C_{({N_H}^\ddagger)_{(1)}}(\theta) = C_{T_H}(\theta) \langle \theta \rangle$. Since ${G^+}_{\G{1}(V^+)} = Z(G^+) = \langle \theta \rangle$, Lemma~\ref{lem: semisimple auts} shows that in the action of $G^+$ on $\G{1}(V^+)$ the generic stabilizer is $C_{({N_H}^\ddagger)_{(1)}}(\theta)/Z(G^+) \cong {\Z_2}^4$. Finally Lemma~\ref{lem: reduction to projective space} shows that if $p \geq 3$ the quadruple $(G, \lambda, p, k)$ also has generic stabilizer ${\Z_2}^4$.

Now suppose instead $p = 2$; here $Z(G) = \{ 1 \}$. We shall in general use the approach of Section~\ref{sect: Lie algebra annihilators}; however, this case presents some features which mean that we cannot simply apply Lemma~\ref{lem: exactness via Premet}, but instead must modify the strategy somewhat. We use the standard notation for the roots in $\Phi$, and then each weight $\nu \in \Lambda(V)$ is of the form $\frac{1}{2}\sum_{i = 1}^5 \pm\ve_i$, where the number of minus signs is even; we shall represent such a weight as a string of $5$ plus or minus signs, and write $v_\nu$ for the corresponding weight vector, so that $V = \langle v_\nu : \nu \in \Lambda(V) \rangle$ and each element $n_\alpha$ for $\alpha \in \Phi$ permutes the vectors $v_\nu$. We shall sometimes abbreviate a root $\pm\ve_i \pm \ve_j$ to $\pm i \pm j$ where it appears in a subscript.

We take the generalized height function on the weight lattice of $G$ whose value at $\alpha_4$ and $\alpha_5$ is $0$, and at $\alpha_1$, $\alpha_2$ and $\alpha_3$ is $1$; then the generalized height of $\lambda = \frac{1}{2}(\alpha_1 + 2\alpha_2 + 3\alpha_3 + \frac{3}{2}\alpha_4 + \frac{5}{2}\alpha_5)$ is $3$, and as $\lambda$, $\omega_4 = \lambda + \frac{1}{2}\alpha_4 - \frac{1}{2}\alpha_5$ and $\Phi$ generate the weight lattice it follows that the generalized height of any weight is an integer. Since $\lambda = {\ss+++++}$, we see that the generalized height of the weight $\e_1 \e_2 \e_3 \e_4 \e_5$ is $\frac{1}{2}(3\e_1 + 2\e_2 + \e_3)$ (if we regard each $\e_i$ as $\pm1$). Moreover we have $\Phi_{[0]} = \langle \alpha_4, \alpha_5 \rangle = \{ \pm\ve_4 \pm \ve_5 \}$, so that $G_{[0]} = \langle T, X_{\pm\alpha_4}, X_{\pm\alpha_5} \rangle$; thus the derived group $(G_{[0]})' = \langle X_{\pm\alpha_4}, X_{\pm\alpha_5} \rangle$ is of type $D_2$, and has centralizer $\langle X_{\pm\alpha_1}, X_{\pm\alpha_2}, X_{\pm\rho} \rangle = \langle X_{\pm\ve_i \pm \ve_j} : 1 \leq i < j \leq 3 \rangle$ of type $D_3$ (where we write $\rho = \alpha_1 + 2\alpha_2 + 2\alpha_3 + \alpha_4 + \alpha_5 = \ve_1 + \ve_2$ for the high root of $\Phi$).

For $i = 1, 2, 3, 4$ write $V^{(i)} = V^{(i), 4} \oplus V^{(i), 5}$, where
\begin{eqnarray*}
V^{(1), 4} = \langle v_{\sss{-+++-}}, v_{\sss{-++-+}} \rangle, & & V^{(1), 5} = \langle v_{\sss{+--++}}, v_{\sss{+----}} \rangle,\\
V^{(2), 4} = \langle v_{\sss{+-++-}}, v_{\sss{+-+-+}} \rangle, & & V^{(2), 5} = \langle v_{\sss{-+-++}}, v_{\sss{-+---}} \rangle,\\
V^{(3), 4} = \langle v_{\sss{++-+-}}, v_{\sss{++--+}} \rangle, & & V^{(3), 5} = \langle v_{\sss{--+++}}, v_{\sss{--+--}} \rangle,\\
V^{(4), 4} = \langle v_{\sss{---+-}}, v_{\sss{----+}} \rangle, & & V^{(4), 5} = \langle v_{\sss{+++++}}, v_{\sss{+++--}} \rangle;
\end{eqnarray*}
note that $V^{(i), 4}$ and $V^{(i), 5}$ are natural modules for $\langle X_{\pm\alpha_4} \rangle = \langle X_{\pm(4 - 5)} \rangle$ and $\langle X_{\pm\alpha_5} \rangle = \langle X_{\pm(4 + 5)} \rangle$ respectively. Then
$$
V = V^{(1)} \oplus V^{(2)} \oplus V^{(3)} \oplus V^{(4)},
$$
and indeed $V^{(1)} = V_{[0]}$ while for $i = 2, 3, 4$ we have $V^{(i)} = V_{[i - 1]} \oplus V_{[-(i - 1)]}$.

Given $\a = (a_1, a_2, a_3, a_4) \in K^4$ with $(a_1 + a_2, a_3 + a_4) \neq (0, 0)$, let $y_\a = \langle v^{(1)}, v^{(2)}, v^{(3)}, v^{(4)} \rangle$ with
\begin{eqnarray*}
v^{(1)} & \! = \! &              (a_1 + a_2) v_{\sss{-+++-}} + (a_1 + a_2) v_{\sss{-++-+}} +              (a_3 + a_4) v_{\sss{+--++}} + (a_3 + a_4) v_{\sss{+----}}, \\
v^{(2)} & \! = \! & \phantom{({} + a_2)} a_1 v_{\sss{+-++-}} + (a_1 + a_2) v_{\sss{+-+-+}} + \phantom{({} + a_4)} a_3 v_{\sss{-+-++}} + (a_3 + a_4) v_{\sss{-+---}}, \\
v^{(3)} & \! = \! & \phantom{(a_1 + {})} a_2 v_{\sss{++-+-}} + (a_1 + a_2) v_{\sss{++--+}} + \phantom{(a_3 + {})} a_4 v_{\sss{--+++}} + (a_3 + a_4) v_{\sss{--+--}}, \\
v^{(4)} & \! = \! & \phantom{(a_1 + a_2) v_{\sss{---+-}} + {}} (a_1 + a_2) v_{\sss{----+}} \phantom{{} + (a_3 + a_4) v_{\sss{+++++}}} + (a_3 + a_4) v_{\sss{+++--}}.
\end{eqnarray*}
Let $Y = \{ y_\a : (a_1 + a_2, a_3 + a_4) \neq (0, 0) \}$; then $Y$ is a subvariety of $\G{4}(V)$ of dimension $3$, whence $\codim Y = 48 - 3 = 45 = \dim G$.

Define $\S = \langle h_1, h_2, h_4 + h_5 \rangle \leq \L(T)$, where we write $h_i$ for $h_{\alpha_i}$. Here we do not have $C_{\L(G)}(\S) = \L(T)$ as in previous proofs; rather we see that if $\alpha \in \Phi \setminus \Phi_{[0]}$ then there exists $h \in \S$ with $[h e_\alpha] \neq 0$, whereas if $\alpha \in \Phi_{[0]}$ then for all $h \in \S$ we have $[h e_\alpha] = 0$, so $C_{\L(G)}(\S) = \L(G_{[0]}) = \L(T) \oplus \langle e_{\alpha_4}, e_{-\alpha_4}, e_{\alpha_5}, e_{-\alpha_5} \rangle$. If $y \in Y$ we have $\S \leq \Ann_{\L(G)}(y)$.

Given $\a = (a_1, a_2, a_3, a_4) \in K^4$ satisfying $(a_1 + a_2)(a_3 + a_4) \neq 0$, write $b_4 = \frac{a_2}{a_1 + a_2}$ and $b_5 = \frac{a_4}{a_3 + a_4}$. With this notation, set
\begin{eqnarray*}
\hat Y & = & \{ y_\a \in Y : a_1a_2a_3a_4 \neq 0, \ (a_1 + a_2)(a_3 + a_4)(a_1 + a_2 + a_3 + a_4) \neq 0, \\
       &   & \phantom{\{ y_\a \in Y : \ } a_1(a_1 + a_2) + a_3(a_3 + a_4) \neq 0, \ a_2(a_1 + a_2) + a_4(a_3 + a_4) \neq 0, \\
       &   & \phantom{\{ y_\a \in Y : \ } {b_4}^3 \neq 1, \ {b_5}^3 \neq 1, \ {\ts\frac{b_5(1 + b_5)}{b_4(1 + b_4)}} \neq 1, {\ts\frac{1}{{b_4}^3}}, {b_5}^3, {\ts\frac{1}{(1 + b_4)^3}}, (1 + b_5)^3 \};
\end{eqnarray*}
then $\hat Y$ is a dense open subset of $Y$. Take $y = y_\a \in \hat Y$.

First suppose $x \in \Ann_{\L(G)}(y)$; write $x = h + e$ where $h \in \L(G_{[0]})$ and $e \in \langle e_\alpha : \alpha \in \Phi \setminus \Phi_{[0]} \rangle$. Clearly for each $i$ the vector $h.v^{(i)}$ lies in $V^{(i)}$; since the difference of two weights lying in $V^{(i)}$ is not a root outside $\Phi_{[0]}$, the projection of $x.v^{(i)}$ on $V^{(i)}$ is equal to $h.v^{(i)}$. Thus for each $i$, the vector $h.v^{(i)}$ must be a scalar multiple of $v^{(i)}$, while for each $j \neq i$ the projection of $e.v^{(i)}$ on $V^{(j)}$ must be a scalar multiple of $v^{(j)}$. A quick calculation shows that we must have $h \in \S$. Now write $e = \sum_{\alpha \in \Phi \setminus \Phi_{[0]}} t_\alpha e_\alpha$; then the condition on the projections of the vectors $e.v^{(i)}$ on the $V^{(j)}$ may be expressed in matrix form as $A{\bf t} = {\bf 0}$, where $A$ is a $36 \times 36$ matrix and ${\bf t}$ is a column vector whose entries are the various coefficients $t_\alpha$. We find that if the rows and columns of $A$ are suitably ordered then it becomes block diagonal, having $3$ blocks, with each block being a $12 \times 12$ matrix. In fact one block may be written in the form
$$
\left(
  \begin{array}{cccccccccccc}
    a_1 s_2 & a_2 s_2 &         &         &         &         & {s_2}^2 & a_1 a_2 & a_3 s_2 &         &         & a_2 s_1 \\
    s_1 a_4 & s_1 a_3 &         &         & a_3 a_4 & {s_1}^2 &         &         &         & a_1 s_1 & a_4 s_2 &         \\
    s_1 s_2 & s_1 s_2 &         &         & a_3 s_2 &         &         & a_1 s_1 &         &         & {s_2}^2 & {s_1}^2 \\
    a_1 s_2 & a_2 s_2 &         &         &         &         & {s_2}^2 & a_1 a_2 & a_4 s_2 &         &         & a_1 s_1 \\
    s_1 a_4 & s_1 a_3 &         &         & a_3 a_4 & {s_1}^2 &         &         &         & a_2 s_1 & a_3 s_2 &         \\
    s_1 s_2 & s_1 s_2 &         &         & a_4 s_2 &         &         & a_2 s_1 &         &         & {s_2}^2 & {s_1}^2 \\
            &         &         & s_1 s_2 &         & {s_2}^2 &         &         &         &         &         & {s_2}^2 \\
            &         & s_1 s_2 &         &         &         & {s_1}^2 &         &         &         & {s_1}^2 &         \\
            &         & s_1 s_2 & s_1 s_2 & {s_1}^2 &         &         & {s_2}^2 & {s_1}^2 & {s_2}^2 &         &         \\
            &         &         & s_1 s_2 &         & {s_2}^2 &         &         & {s_1}^2 &         &         &         \\
            &         & s_1 s_2 &         &         &         & {s_1}^2 &         &         & {s_2}^2 &         &         \\
            &         & s_1 s_2 & s_1 s_2 &         &         &         &         & {s_1}^2 & {s_2}^2 &         &         \\
  \end{array}
\right)
$$
(where for reasons of space we set $s_1 = a_1 + a_2$ and $s_2 = a_3 + a_4$), and the other two are of similar form; the three determinants are $(a_1 + a_2)^4 (a_3 + a_4)^4 f(\a)^8$, where $f(\a) = (a_1 + a_2 + a_3 + a_4)^2$, $a_1(a_1 + a_2) + a_3(a_3 + a_4)$ and $a_2(a_1 + a_2) + a_4(a_3 + a_4)$. Thus the definition of the set $\hat Y$ implies that each block of $A$ is non-singular, as therefore is $A$ itself; so ${\bf t}$ must be the zero vector and hence $e = 0$. Thus $x = h + e \in \S$; so $\Ann_{\L(G)}(y) = \S$.

Now suppose $g \in \Tran_G(y, Y)$. As in the proof of Lemma~\ref{lem: exactness via Premet}, we see that $\Ad(g).\S = \Ad(g).\Ann_{\L(G)}(y) = \Ann_{\L(G)}(g.y) \geq \S$ because $g.y \in Y$, so we must have $\Ad(g).\S = \S$; therefore $\Ad(g).C_{\L(G)}(\S) = C_{\L(G)}(\S)$. However, since here we do not have $C_{\L(G)}(\S) = \L(T)$, we cannot deduce that $\Ad(g).\L(T) = \L(T)$, and so $g \in N$; instead we have $\Ad(g).\L(G_{[0]}) = \L(G_{[0]})$. However, $\Ad(g)$ must then preserve the derived subalgebra of $\L(G_{[0]})$, which is $\L((G_{[0]})') = \L(D_2)$; a quick calculation with $10 \times 10$ matrices shows that $g$ must lie in $D_3D_2$. The intersection of $\L(G_{[0]})$ with $\L(D_3)$ is then $\langle h_1, h_2 \rangle$, and it follows that $\Ad(g)$ must also preserve this; indeed an easy calculation in $D_3$ shows that $\Ad(g)$ must preserve the set $\{ h_1, h_2, h_1 + h_2 \}$, and with a little more work we find that $g \in G_{[0]} \langle n_{1 - 2}, n_{2 - 3}, n_{1 + 2}, n^* \rangle$ where $n^* = n_{1 - 2} n_{1 + 2} n_{3 - 5} n_{3 + 5}$. Since for each $i$ the element $n^*$ interchanges the ordered bases of $V^{(i), 4}$ and $V^{(i), 5}$, it sends $y_\a$ to $y_{\a'}$ where $\a' = (a_3, a_4, a_1, a_2)$, and so preserves $\hat Y$. Thus we may assume $g = g_4g_5sn$, where $g_4 \in \langle X_{\pm\alpha_4} \rangle$, $g_5 \in \langle X_{\pm\alpha_5} \rangle$, $s \in T \cap D_3$ and $n \in \langle n_{1 - 2}, n_{2 - 3}, n_{1 + 2} \rangle \cong S_4$.

For each $i$ write $v^{(i), 4}$ and $v^{(i), 5}$ for the projections of $v^{(i)}$ on $V^{(i), 4}$ and $V^{(i), 5}$ respectively. We see that $n$ permutes the $V^{(i)}$, say $n.V^{(i)} = V^{(\pi(i))}$ where $\pi \in S_4$; indeed $n$ sends the ordered bases of $V^{(i), 4}$ and $V^{(i), 5}$ to those of $V^{(\pi(i)), 4}$ and $V^{(\pi(i)), 5}$ respectively. Moreover for each $i$ the element $s$ acts on each of $V^{(i), 4}$ and $V^{(i), 5}$ as a scalar, with the two scalars being inverses of each other.

Fix $n$, and suppose $g = g_4g_5sn$ and $g' = {g_4}'{g_5}'s'n$ both lie in $\Tran_G(y, Y)$. Write $x = {g_4}'{g_4}^{-1} \in \langle X_{\pm\alpha_4} \rangle$; let the standard basis of the natural $\langle X_{\pm\alpha_4} \rangle$-module be $v_1, v_2$, and write $c_1v_1 + c_2v_2$ as $(c_1, c_2)$. Both $g_4$ and ${g_4}'$ send $n.v^{(\pi^{-1}(4)), 4}$ to a vector in $V^{(4), 4}$ in which the first basis vector has coefficient $0$; thus $x$ preserves the line $\langle (0, 1) \rangle$, so must be a lower triangular matrix. Similarly both $g_4$ and ${g_4}'$ send $n.v^{(\pi^{-1}(1)), 4}$ to a vector in $V^{(1), 4}$ in which the two basis vectors have equal coefficients; thus $x$ preserves the line $\langle (1, 1) \rangle$, so must be of the form
$\left(
  \begin{array}{cc}
    \kappa               &             \\
    \kappa + \kappa^{-1} & \kappa^{-1} \\
  \end{array}
\right)$
for some $\kappa \in K^*$. Now if $g_4$ sends $n.v^{(\pi^{-1}(2)), 4}$ to a vector in the line $\langle ({a_1}', {a_1}' + {a_2}') \rangle$, it must send $n.v^{(\pi^{-1}(3)), 4}$ to a vector in the line $\langle ({a_2}', {a_1}' + {a_2}') \rangle$, while ${g_4}'$ sends $n.v^{(\pi^{-1}(2)), 4}$ and $n.v^{(\pi^{-1}(3)), 4}$ to vectors in the lines $\langle \kappa{a_1}', \kappa{a_1}' + \kappa^{-1}{a_2}') \rangle$ and $\langle (\kappa{a_2}', \kappa^{-1}{a_1}' + \kappa{a_2}') \rangle$ respectively. Since in each case the vector in $V^{(3), 4}$ is obtained from that in $V^{(2), 4}$ by applying the transformation
$\left(
  \begin{array}{cc}
    1 & 1 \\
      & 1 \\
  \end{array}
\right)$,
we must have $\langle (\kappa{a_2}', \kappa^{-1}{a_1}' + \kappa{a_2}') \rangle = \langle (\kappa^{-1}{a_2}', \kappa{a_1}' + \kappa^{-1}{a_2}') \rangle$, whence $\kappa{a_2}'(\kappa{a_1}' + \kappa^{-1}{a_2}') = \kappa^{-1}{a_2}'(\kappa^{-1}{a_1}' + \kappa{a_2}')$ and so $\kappa^2{a_1}'{a_2}' = \kappa^{-2}{a_1}'{a_2}'$. As the four vectors $n.v^{(\pi^{-1}(i)), 4}$ lie in distinct lines in the natural $\langle X_{\pm\alpha_4} \rangle$-module, the same must be true of their images under ${g_4}'$; thus ${a_1}', {a_2}', {a_1}' + {a_2}' \neq 0$ and so $\kappa = 1$, whence $x = 1$ and ${g_4}' = g_4$. Similarly ${g_5}' = g_5$; so $g' = s'{s^{-1}}g$. Now let $g.y = y_{\a'}$, where $\a' = ({a_1}', {a_2}', {a_3}', {a_4}')$; as above we must have ${a_1}' + {a_2}' \neq 0$ and likewise ${a_3}' + {a_4}' \neq 0$. In each basis vector of $y_{\a'}$ the ratio of the coefficients of the second and fourth weight vectors is the same, namely $\frac{{a_1}' + {a_2}'}{{a_3}' + {a_4}'}$, so this must also be true in $s'{s^{-1}}.y_{\a'}$. Write $s'{s^{-1}} = h_{1 - 2}(\kappa_1) h_{2 - 3}(\kappa_2) h_{1 + 2}(\kappa_3)$; then for $i = 1, 2, 3, 4$ the element $s'{s^{-1}}$ acts on $V^{(i), 4}$ and $V^{(i), 5}$ as the scalars $\kappa$ and $\kappa^{-1}$, where $\kappa = \frac{1}{\kappa_1}$, $\frac{\kappa_1}{\kappa_2}$, $\kappa_2 \kappa_3$ and $\frac{1}{\kappa_3}$ respectively. Thus the ratio of the coefficients of the second and fourth weight vectors is multiplied by $\kappa^2$ in each case, so we must have $\frac{1}{{\kappa_1}^2} = \frac{{\kappa_1}^2}{{\kappa_2}^2} = {\kappa_2}^2 {\kappa_3}^2 = \frac{1}{{\kappa_3}^2}$, whence $\kappa_1 = \kappa_2 = \kappa_3 = 1$; so $s'{s^{-1}} = 1$ and $g' = g$. Thus for each of the $24$ elements $n$ there can be at most one element $g = g_4g_5sn \in \Tran_G(y, Y)$; so $\Tran_G(y, Y)$ is finite. Thus $\codim \Tran_G(y, Y) = \dim G$, so $y$ is $Y$-exact.

For $t_1, t_2 \in K$ define
\begin{eqnarray*}
g_1           & = & n_{2 - 3} n_{2 + 3} x_{4 - 5}(1) x_{4 + 5}(1), \\
g_2(t_1, t_2) & = & n_{1 - 2} n_{1 + 2} x_{4 - 5}(t_1) x_{4 + 5}(t_2),
\end{eqnarray*}
and set $C(t_1, t_2) = \langle g_1, g_2(t_1, t_2) \rangle \cong {\Z_2}^2$. By inspection we see that $C(b_4, b_5) \leq C_G(y)$. We shall show that in fact $C_G(y) = C(b_4, b_5)$.

Thus suppose $g \in C_G(y)$, and as above write $g = (n^*)^jg_4g_5sn$ with $j \in \{ 0, 1 \}$. Since $C(b_4, b_5)$ acts transitively on the $V^{(i)}$, by multiplying $g$ by an element of $C(b_4, b_5)$ we may assume $n$ fixes $V^{(4)}$. From the above we see that $g_4 \in \langle X_{\pm\alpha_4} \rangle$ fixes the line $\langle (0, 1) \rangle$, so must be lower triangular, say
$\left(
  \begin{array}{cc}
    \kappa &             \\
       t   & \kappa^{-1} \\
  \end{array}
\right)$
for some $\kappa \in K^*$ and $t \in K$; if $j = 0$ it permutes the lines $\langle (1, 1) \rangle$, $\langle (1 + b_4, 1) \rangle$ and $\langle (b_4, 1) \rangle$, while if $j = 1$ it sends them to $\langle (1, 1) \rangle$, $\langle (1 + b_5, 1) \rangle$ and $\langle (b_5, 1) \rangle$ in some order. Since the images of the first three lines are $\langle (\kappa, t + \kappa^{-1}) \rangle$, $\langle (\kappa(1 + b_4), t(1 + b_4) + \kappa^{-1}) \rangle$ and $\langle (\kappa b_4, tb_4 + \kappa^{-1}) \rangle$ respectively, we must have
$$
\frac{\kappa}{t + \kappa^{-1}} = c_1, \quad \frac{\kappa(1 + b_4)}{t(1 + b_4) + \kappa^{-1}} = c_2, \quad \frac{\kappa b_4}{tb_4 + \kappa^{-1}} = c_3
$$
where $\{ c_1, c_2, c_3 \} = \{ 1, 1 + b_i, b_i \}$ with $i = j + 4$. Thus $\frac{\kappa}{t + \kappa^{-1}} + \frac{\kappa(1 + b_4)}{t(1 + b_4) + \kappa^{-1}} + \frac{\kappa b_4}{tb_4 + \kappa^{-1}} = 0$, which reduces to $\kappa t^2b_4(1 + b_4) = 0$; so we must have $t = 0$, and hence $c_1 = \kappa^2$, $c_2 = \kappa^2(1 + b_4)$, $c_3 = \kappa^2b_4$ so that $b_i(1 + b_i) = c_1c_2c_3 = \kappa^6b_4(1 + b_4)$. If $j = 1$ this implies $\frac{b_5(1 + b_5)}{b_4(1 + b_4)} = \kappa^6 = {c_1}^3 \in \{ 1, {b_5}^3, (1 + b_5)^3 \}$, contrary to the definition of $\hat Y$. Thus we must have $j = 0$, and so $\kappa^6 = 1$; since ${b_4}^3 \neq 1$ we must have $c_1 = 1$, $c_2 = 1 + b_4$, $c_3 = b_4$, so that $n = 1$. Therefore by the above we must have $g = 1$ as required; so $C_G(y) = C(b_4, b_5)$.

Set $C = C(0, 1)$. Given $b_4$ and $b_5$ as above, for reasons of space write $c = (1 + \sqrt{b_4})(1 + \sqrt{b_5})$, $e = \frac{\sqrt{b_4b_5}}{b_4 + b_5}$ and $f = \frac{1}{\sqrt{b_4} + \sqrt{b_5}}$; take $h \in K$ satisfying $h^2 + h = c^2e^2$, and set $j = c(h + \frac{b_4b_5}{b_4 + b_5})$. Regarding $G$ as $\SO_{10}(K)$ and taking the standard basis of $V_{nat}$ in the order $v_1, v_2, v_3, v_4, v_5, v_{-5}, v_{-4}, v_{-3}, v_{-2}, v_{-1}$, set
$$
g =
\left(
  \begin{array}{cccccccccc}
    1 &            &         &        &         &            &            &         &            &   \\
      &      1     &         &        & ef^{-1} &            &            &         &            &   \\
      &     ce     &  h + 1  &        &         &     cf     &     cf     &    h    &     ce     &   \\
      & (1 + b_4)e &  c + j  & f^{-1} &  f^{-1} & (1 + b_5)f & (1 + b_5)f &    j    & (1 + b_5)e &   \\
      & (1 + b_5)e & b_5cf^2 &        &   b_5f  & (1 + b_5)f & (1 + b_5)f & b_5cf^2 & (1 + b_5)e &   \\
      & (1 + b_4)e & b_4cf^2 &        &   b_4f  & (1 + b_4)f &      f     & b_4cf^2 & (1 + b_4)e &   \\
      &            &         &        &         &            &      f     &         &            &   \\
      &     ce     &    h    &        &         &     cf     &            &  h + 1  &     ce     &   \\
      &            &         &        & ef^{-1} &            &   ef^{-1}  &         &      1     &   \\
      &            &         &        &         &            &            &         &            & 1 \\
  \end{array}
\right).
$$
We find that $g_1g = gg_1$ and $g_2(b_4, b_5)g = gg_2(0, 1)$, so $C(b_4, b_5)^g = C$; and $g$ preserves the relevant quadratic form, so lies in $\mathrm{O}_{10}(K)$. The matrix $g' = I + (E_{1, 1} + E_{-1, -1} + E_{1, -1} + E_{-1, 1})$ commutes with both $g_1$ and $g_2(b_4, b_5)$ and lies in $\mathrm{O}_{10}(K) \setminus \SO_{10}(K)$; therefore either $g$ or $g'g$ lies in $G$ and conjugates $C(b_4, b_5)$ to $C$. Thus the conditions of Lemma~\ref{lem: generic stabilizer from exact subset} hold; so if $p = 2$ the quadruple $(G, \lambda, p, k)$ has generic stabilizer $C/Z(G) \cong {\Z_2}^2$.
\end{proof}

For the final few results in this section we relax slightly the condition that the group acting should be simple: we allow a product of isomorphic simple groups, possibly extended by a graph automorphism. If the connected component is of the form $G_1 G_2$ or $G_1 G_2 G_3$, we shall write $\lambda = \lambda_1 \otimes \lambda_2$ or $\lambda_1 \otimes \lambda_2 \otimes \lambda_3$ to mean that $L(\lambda) = L(\lambda_1) \otimes L(\lambda_2)$ or $L(\lambda_1) \otimes L(\lambda_2) \otimes L(\lambda_3)$, where each $\lambda_i$ is a dominant weight for $G_i$.

\begin{prop}\label{prop: {A_2}^2, omega_1 otimes omega_1, k = 3}
Let $G = {A_2}^2$ and $\lambda = \omega_1 \otimes \omega_1$, and take $k = 3$; let $\tau$ be a graph automorphism of $G$ of order $2$ interchanging the simple factors of $G$. Then the quadruple $(G, \lambda, p, k)$ has generic stabilizer $\Z_{3/(p, 3)}.\Z_3$, while $(G\langle \tau \rangle, \lambda, p, k)$ has generic stabilizer $\Z_{3/(p, 3)}.\Z_3.\Z_2$.
\end{prop}

\begin{proof}
We shall deal separately with the cases $p \neq 3$ and $p = 3$.

First suppose $p \neq 3$. Let $H$ be the simply connected group defined over $K$ of type $E_6$, with simple roots $\beta_1, \dots, \beta_6$; we then have $Z(\L(H)) = \{ 0 \}$. Our strategy will be to identify the group $G^+ = {A_2}^3 = GA_{k - 1}$ as the centralizer of a semisimple automorphism of $H$, and use Lemma~\ref{lem: semisimple auts} to find the generic stabilizer in the action of $G^+$ on $\G{1}(V^+)$, where $V^+ = V \otimes V_{nat}$ with $V$ the ${A_2}^2$-module with high weight $\omega_1 \otimes \omega_1$ and $V_{nat}$ the natural $A_2$-module; we shall therefore employ the notation of Section~\ref{sect: semisimple auts}. Lemma~\ref{lem: reduction to projective space} will then give the result.

Let $\tau$ be the graph automorphism of $H$ which for all $t \in K$ interchanges $x_{\beta_1}(t)$ with $x_{\beta_6}(t)$, and $x_{\beta_3}(t)$ with $x_{\beta_5}(t)$, while fixing $x_{\beta_2}(t)$ and $x_{\beta_4}(t)$. Define $n_0 = h_{\beta_2}(-1) n_{\beta_4} n_{\beta_3 + \beta_4 + \beta_5} n_{\beta_1 + \beta_3 + \beta_4 + \beta_5 + \beta_6} n_\rho$, where $\rho = \esixrt122321$ is the high root of $\Phi_H$; then $n_0T_H$ is the long word $w_0$, and $n_0\tau$ sends each $x_\alpha(t)$ to $x_{-\alpha}(t)$.

Define $\theta_1 \in T_H$ to be
$$
h_{\beta_1}({\eta_3}^2) h_{\beta_2}({\eta_3}^2) h_{\beta_6}({\eta_3}^2).
$$
Then ${\theta_1}^3 = 1$, and $\theta_1$ sends $x_\alpha(t)$ to $x_\alpha({\eta_3}^{\height(\alpha)}t)$; so $X_\alpha < C_H(\theta_1)$ if and only if $\height(\alpha) \equiv 0$ (mod $3$). It follows that $C_H(\theta_1)$ is a connected group of type ${A_2}^3$, with simple root elements $x_{\alpha_i}(t)$, where $\alpha_1, \dots, \alpha_6$ are
$$
\esixrt000111, \esixrt011100, \esixrt101100, \esixrt010110, \esixrt001110, \esixrt111111;
$$
we see that $Z(C_H(\theta_1)) = \langle \theta_1, z \rangle$ where $z = h_{\beta_1}(\eta_3) h_{\beta_3}({\eta_3}^2) h_{\beta_5}(\eta_3) h_{\beta_6}({\eta_3}^2)$.

Now let $\delta_1, \dots, \delta_6$ be
$$
\esixrt100000, \esixrt001000, \esixrt000001, \esixrt000010, \esixrt010000, \esixrt112321,
$$
and set $\theta_2 = n_{\delta_1} \dots n_{\delta_6}$; then ${\theta_2}^3 = 1$. We find that $\theta_2$ acts fixed-point-freely on both $\Phi_H$ and $\L(T_H)$. Thus $\dim C_{\L(H)}(\theta_2) = |\Phi_H|/3 = 24$; the classification of semisimple elements of $H$ (see e.g. \cite[Table~4.7.1]{GLS}) now shows that $\theta_2$ must be a conjugate of $\theta_1$.

First set $\theta = \theta_1$; then we may take $G^+ = C_H(\theta)$, and $\tau$ acts on $G^+$ as the graph automorphism of ${A_2}^2$ while fixing pointwise the third $A_2$ factor. We see that $e_\alpha \in \L(H)_{(i)}$ if and only if $\height(\alpha) \equiv i$ (mod $3$). Thus in $\L(H)_{(1)}$ we have a highest weight vector $e_\beta$ for $\beta = \esixrt112321$; the expressions above for the simple root elements of $G^+$ show that $\L(H)_{(1)}$ is the Weyl $G^+$-module with high weight $\omega_1 \otimes \omega_1 \otimes \omega_1$, i.e., the tensor product $V \otimes V_{nat}$. As $Z(\L(H)) = \{ 0 \}$, we may therefore take $V^+ = \L(H)_{(1)}/Z(\L(H))_{(1)}$; of course ${G^+}_{\G{1}(V^+)} = Z(G^+)$.

Now set $\theta = \theta_2$, and again take $G^+ = C_H(\theta)$ and $V^+ = \L(H)_{(1)}/Z(\L(H))_{(1)}$. We have $G^+ \cap T_H = C_{T_H}(\theta) \cong {\Z_3}^3$; indeed this group is
$$
\langle h_{\beta_1}(\eta_3) h_{\beta_3}({\eta_3}^2), h_{\beta_5}({\eta_3}^2) h_{\beta_6}(\eta_3), h_{\beta_1}(\eta_3) h_{\beta_2}(\eta_3) h_{\beta_6}(\eta_3) \rangle.
$$
Moreover we find that $\L(T_H)_{(1)} = \langle h_{\delta_{2i - 1}} - \eta_3 h_{\delta_{2i}} : i = 1, 2, 3 \rangle$. Thus $\dim \L(H)_{(1)} - \dim \L(T_H)_{(1)} = \dim G^+ - \dim (G^+ \cap T_H)$; and a routine check shows that $\L(T_H)_{(1)}$ contains regular semisimple elements.

We claim that we have $({W_H}^\ddagger)_{(1)} = \langle \theta T_H \rangle$. Thus suppose $w \in W_H$ and there exists $\xi \in K^*$ such that for all $y \in \L(T_H)_{(1)}$ we have $w.y = \xi y$. For $i = 1, 2, 3$ write $\Psi_i = \langle \delta_{2i - 1}, \delta_{2i} \rangle$ and $y_i = h_{\delta_{2i - 1}} - \eta_3 h_{\delta_{2i}} \in \L(T_H)_{(1)}$. Taking $y = y_1$ and arguing as in the paragraphs following the statement of Lemma~\ref{lem: semisimple auts} shows that $w(\beta_1)$ and $w(\beta_3)$ must be proportional outside $\{ \beta_1, \beta_3 \}$, and as $\eta_3 \neq \pm1$ that $w$ must preserve $\Psi_1$. Now take $i \in \{ 2, 3 \}$. There exists $w' \in W_H$ with $w'(\delta_1) = \delta_{2i - 1}$ and $w'(\delta_2) = \delta_{2i}$, and so $w'.y_1 = y_i$, whence $w.y_i = \xi y_i$ gives $w^{w'}.y_1 = y_1$; by the above $w^{w'}$ preserves $\Psi_1$, so $w$ preserves $\Psi_i$. Thus $w = w_1w_2w_3$ where each $w_i$ lies in $W(\Psi_i)$. For each $i$, the three elements in $W(\Psi_i)$ of odd length send $y_i$ to a scalar multiple of $\eta_3 h_{\delta_{2i - 1}} - h_{\delta_{2i}}$, so we must have $w_i \in \langle w_{\delta_{2i - 1}} w_{\delta_{2i}} \rangle$; since $w$ must multiply each of the three vectors $y_i$ by the same scalar, we must have $w \in \langle \theta T_H \rangle$ as required. Note that $w_0\tau$ sends each root $\alpha$ to its negative, and therefore acts on $\L(T_H)_{(1)}$ as negation.

Now $({N_H}^\ddagger)_{(1)} = T_H \langle \theta \rangle$, so $C_{({N_H}^\ddagger)_{(1)}}(\theta) = C_{T_H}(\theta) \langle \theta \rangle$. Since ${G^+}_{\G{1}(V^+)} = Z(G^+) = \langle \theta, z \rangle$, Lemma~\ref{lem: semisimple auts} shows that in the action of $G^+$ on $\G{1}(V^+)$ the generic stabilizer is $C_{({N_H}^\ddagger)_{(1)}}(\theta)/Z(G^+) \cong {\Z_3}^2$ (and in the action of $G^+\langle \tau \rangle$ the generic stabilizer is ${\Z_3}^2.\Z_2$). Finally Lemma~\ref{lem: reduction to projective space} shows that if $p \neq 3$ the quadruple $(G, \lambda, p, k)$ also has generic stabilizer ${\Z_3}^2$ (while in the action of $G \langle \tau \rangle$ the presence of the element $n_0\tau$ means that the generic stabilizer is ${\Z_3}^2.\Z_2$).

Now suppose instead $p = 3$; here $Z(G) = \{ 1 \}$. We shall follow the strategy used in the second part of the proof of Proposition~\ref{prop: A_2, 2omega_1 module, k = 3, A_4, omega_2 module, k = 5}. Let $G$ have simple roots $\alpha_1$, $\alpha_2$, $\alpha_3$ and $\alpha_4$ (where $\alpha_1 + \alpha_2$ and $\alpha_3 + \alpha_4$ are the other two positive roots), and $V$ be the tensor product of the natural modules for the simple factors of $G$; for $i_1, i_2 \in \{ 1, 2, 3 \}$ write $v_{i_1i_2} = v_{i_1} \otimes v_{i_2}$, so that $V = \langle v_{i_1i_2} : i_1, i_2 \in \{ 1, 2, 3 \} \rangle$, and for example $x_{\alpha_1}(t)$ fixes $v_{1i_2}$ and $v_{3i_2}$ and sends $v_{2i_2}$ to $v_{2i_2} + tv_{1i_2}$, while $x_{\alpha_3}(t)$ fixes $v_{i_11}$ and $v_{i_13}$ and sends $v_{i_12}$ to $v_{i_12} + tv_{i_11}$. Since $\dim V = 9$ we have $\dim \G{3}(V) = 18$. Let $\tau$ be the automorphism of $G$ which interchanges $x_{\alpha_1}(t)$ and $x_{\alpha_2}(t)$ with $x_{\alpha_3}(t)$ and $x_{\alpha_4}(t)$ respectively; then $\tau$ acts on $V$ by sending each $v_{i_1i_2}$ to $v_{i_2i_1}$.

We shall write elements of both $G$ and $\L(G)$ as pairs of $3 \times 3$ matrices. We let $T < G$ be the subgroup of pairs of diagonal matrices, so that $N$ is the subgroup of pairs of monomial matrices; define $n_0, n_1 \in N$ by
$$
n_0 =
\left(
-\left(
  \begin{array}{ccc}
      &   & 1 \\
      & 1 &   \\
    1 &   &   \\
  \end{array}
\right),
-\left(
  \begin{array}{ccc}
      &   & 1 \\
      & 1 &   \\
    1 &   &   \\
  \end{array}
\right)
\right),
\qquad
n_1 =
\left(
\left(
  \begin{array}{ccc}
      &   & 1 \\
    1 &   &   \\
      & 1 &   \\
  \end{array}
\right),
\left(
  \begin{array}{ccc}
      &   & 1 \\
    1 &   &   \\
      & 1 &   \\
  \end{array}
\right)
\right),
$$
so that $n_0T$ is the long word $w_0$ of the Weyl group, and $n_1.v_{i_1i_2} = v_{(i_1 + 1), (i_2 + 1)}$ (with subscripts taken modulo $3$).

Write
\begin{eqnarray*}
V^{(1)} & = & \langle v_{11}, v_{23}, v_{32} \rangle, \\
V^{(2)} & = & \langle v_{22}, v_{31}, v_{13} \rangle, \\
V^{(3)} & = & \langle v_{33}, v_{12}, v_{21} \rangle,
\end{eqnarray*}
so that $V = V^{(1)} \oplus V^{(2)} \oplus V^{(3)}$ and $n_1$ cycles the $V^{(i)}$. Define
$$
Y = \{ y = \langle v^{(1)}, v^{(2)}, v^{(3)} \rangle : v^{(1)} \in V^{(1)} \setminus \{ 0 \}, \ v^{(2)} = n_1.v^{(1)}, \ v^{(3)} = n_1.v^{(2)} \};
$$
then $Y$ is a subvariety of $\G{3}(V)$ of dimension $2$, whence $\codim Y = 16 = \dim G$. Set $C = \langle n_1 \rangle$; then each $y \in Y$ is stabilized by $C$ (and also by $n_0\tau$ in the action of $G \langle \tau \rangle$). For convenience, given $y = \langle v^{(1)}, v^{(2)}, v^{(3)} \rangle \in Y$ with $v^{(1)} = a_1v_{11} + a_2v_{23} + a_3v_{32}$, we shall write $y = y_\a$ where $\a = (a_1, a_2, a_3)$.

Define $\S = \langle h_0, z_1, z_2 \rangle \leq \L(T)$, where $h_0 = (\diag(1, 0, -1), \diag(1, 0, -1))$, $z_1 = (I, 0)$ and $z_2 = (0, I)$. Clearly if $\alpha \in \Phi$ then $[h_0e_\alpha] \neq 0$, so $C_{\L(G)}(\S) = \L(T)$. Given $v^{(i)} \in V^{(i)}$ we have $h_0.v^{(i)} = (i + 1)v^{(i)}$, while $z_1.v^{(i)} = z_2.v^{(i)} = v^{(i)}$; thus if $y \in Y$ we have $\S \leq \Ann_{\L(G)}(y)$. Set
$$
\hat Y = \{ y_\a \in Y : a_1a_2a_3 \neq 0, \ a_i \neq a_{i'} \hbox{ for } i \neq i', \ {a_i}^3 \neq a_1a_2a_3, \ a_1 + a_2 + a_3 \neq 0 \};
$$
then $\hat Y$ is a dense open subset of $Y$. Take $y = y_\a \in \hat Y$.

First suppose $x \in \Ann_{\L(G)}(y)$; write $x = h + e$ where $h \in \L(T)$ and $e \in \langle e_\alpha : \alpha \in \Phi \rangle$. As before we see that the vector $h.v^{(i)}$ must be a scalar multiple of $v^{(i)}$, while for each $j \neq i$ the projection of $e.v^{(i)}$ on $V^{(j)}$ must be a scalar multiple of $v^{(j)}$. A quick calculation shows that we must have $h \in \S$. Now write $e = \sum_{\alpha \in \Phi} t_\alpha e_\alpha$; then the condition on the projections of the vectors $e.v^{(i)}$ on the $V^{(j)}$ may be expressed in matrix form as $A{\bf t} = {\bf 0}$, where $A$ is a $12 \times 12$ matrix and ${\bf t}$ is a column vector whose entries are the various coefficients $t_\alpha$. We find that if the rows and columns of $A$ are suitably ordered then it becomes block diagonal, having $2$ blocks, with each block being a $6 \times 6$ matrix. In fact each block may be written in the form
$$
\left(
  \begin{array}{cccccc}
    -{a_1}^2 &          &  a_2a_3  &  {a_3}^2 &          & -a_1a_2  \\
     a_2a_3  & -{a_1}^2 &          & -a_1a_2  &  {a_3}^2 &          \\
             &  a_2a_3  & -{a_1}^2 &          & -a_1a_2  &  {a_3}^2 \\
     {a_2}^2 &          & -a_1a_3  &          &  a_2a_3  & -{a_1}^2 \\
    -a_1a_3  &  {a_2}^2 &          & -{a_1}^2 &          &  a_2a_3  \\
             & -a_1a_3  &  {a_2}^2 &  a_2a_3  & -{a_1}^2 &
  \end{array}
\right),
$$
which has determinant ${a_1}^3(a_1 + a_2 + a_3)^9$. Thus the definition of the set $\hat Y$ implies that each block of $A$ is non-singular, as therefore is $A$ itself; so ${\bf t}$ must be the zero vector and hence $e = 0$. Thus $x = h + e \in \S$; so $\Ann_{\L(G)}(y) = \S$.

A straightforward calculation shows that $C_T(y) = \{ 1 \}$, and $T.y \cap Y = \{ y \}$. We claim that $N.y \cap Y \subset \hat Y$, and $C_N(y) = C$. Take $n = (n^{(1)}, n^{(2)}) \in \Tran_G(y, Y)$ and write $n.y = y_{\a'}$; regard the elements $n^{(i)}T$ of the Weyl group as permutations. We find that the parities of $n^{(1)}T$ and $n^{(2)}T$ must be equal, and that there exists $\pi \in S_3$ such that $\a' = (a_{\pi(1)}, a_{\pi(2)}, a_{\pi(3)})$, with $\pi = 1$ if and only if $nT \in \langle n_1T \rangle$. Since $y_{\a'} \in \hat Y$ this proves the first claim; moreover the definition of $\hat Y$ shows that if $nT \notin \langle n_1T \rangle$ then $n$ does not stabilize $y$, proving the second. Thus the conditions of Lemma~\ref{lem: exactness via Premet} hold, so that $\Tran_G(y, Y) \subseteq N$, and $y$ is $Y$-exact; moreover $C_G(y) = C_N(y) = C$. Therefore the conditions of Lemma~\ref{lem: generic stabilizer from exact subset} hold; so if $p = 3$ the quadruple $(G, \lambda, p, k)$ has generic stabilizer $C/Z(G) \cong \Z_3$ (while in the action of $G \langle \tau \rangle$ the presence of the element $n_0\tau$ means that the generic stabilizer is $\Z_3.\Z_2$).
\end{proof}

\begin{prop}\label{prop: {A_1}^3, omega_1 otimes omega_1 otimes omega_1, k = 2}
Let $G = {A_1}^3$ and $\lambda = \omega_1 \otimes \omega_1 \otimes \omega_1$, and take $k = 2$. Then the quadruple $(G, \lambda, p, k)$ has generic stabilizer $\Z_{2/(p, 2)}.\Z_2$.
\end{prop}

\begin{proof}
Let $G$ have simple roots $\alpha_1$, $\alpha_2$ and $\alpha_3$, and $V$ be the tensor product of the natural modules for the simple factors of $G$; for $i_1, i_2, i_3 \in \{ 1, 2 \}$ write $v_{i_1i_2i_3} = v_{i_1} \otimes v_{i_2} \otimes v_{i_3}$, so that $V = \langle v_{i_1i_2i_3} : i_1, i_2, i_3 \in \{ 1, 2 \} \rangle$, and for example $x_{\alpha_1}(t)$ fixes $v_{1i_2i_3}$ and sends $v_{2i_2i_3}$ to $v_{2i_2i_3} + tv_{1i_2i_3}$. Then $\Lambda(V) = \{ \pm\alpha_1 \pm \alpha_2 \pm \alpha_3 \}$. For convenience write
\begin{eqnarray*}
& x_1 = v_{111}, \quad x_2 = v_{122}, \quad x_3 = v_{212}, \quad x_4 = v_{221}, & \\
& x_5 = v_{222}, \quad x_6 = v_{211}, \quad x_7 = v_{121}, \quad x_8 = v_{112}. &
\end{eqnarray*}
Given $\a = (a_1, a_2, a_3, a_4) \in K^4 \setminus \{ (0, 0, 0, 0) \}$, we let
$$
v^{(1)} = a_1 x_1 + a_2 x_2 + a_3 x_3 + a_4 x_4, \quad v^{(2)} = a_1 x_5 + a_2 x_6 + a_3 x_7 + a_4 x_8,
$$
and set $y_\a = \langle v^{(1)}, v^{(2)} \rangle$; we let
$$
Y = \{ y_\a : \a \neq (0, 0, 0, 0) \}.
$$
Write
\begin{eqnarray*}
\hat Y & = & \{ y_\a \in Y : a_1a_2a_3a_4, a_1 \pm a_2 \pm a_3 \pm a_4, {a_1}^2 \pm {a_2}^2 \pm {a_3}^2 \pm {a_4}^2 \neq 0, \\
       &   & \phantom{ \{ y_a \in Y : \ } {\ts \frac{a_1a_2}{a_3a_4}, \frac{a_1a_3}{a_2a_4}, \frac{a_1a_4}{a_2a_3}, \frac{{a_i}^2}{{a_j}^2}, \frac{a_i(a_1 \pm a_2 \pm a_3 \pm a_4)}{a_j(a_1 \pm a_2 \pm a_3 \pm a_4)}} \neq \pm 1, \\
       &   & \phantom{ \{ y_a \in Y : \ } {\ts \sum {a_i}^8 + 6\sum {a_i}^4{a_j}^4 - 4\sum {a_i}^6{a_j}^2 + 4\sum {a_i}^4{a_j}^2{a_l}^2} \\
       &   & \phantom{ \{ y_a \in Y : \ } {} - 40{a_1}^2{a_2}^2{a_3}^2{a_4}^2 \neq 0 \},
\end{eqnarray*}
then $\hat Y$ is a dense open subset of $Y$. Take
$$
y = y_\a = \langle v^{(1)}, v^{(2)} \rangle \in \hat Y.
$$

Given $s = h_{\alpha_1}(\kappa_1) h_{\alpha_2}(\kappa_2) h_{\alpha_3}(\kappa_3) \in T$, we have
\begin{eqnarray*}
s.v^{(1)} & = & {\ts \kappa_1\kappa_2\kappa_3 a_1x_1 + \frac{\kappa_1}{\kappa_2\kappa_3}a_2x_2 + \frac{\kappa_2}{\kappa_1\kappa_3}a_3x_3 + \frac{\kappa_3}{\kappa_1\kappa_2}a_4x_4}, \\
s.v^{(2)} & = & {\ts \frac{1}{\kappa_1\kappa_2\kappa_3}a_1x_5 + \frac{\kappa_2\kappa_3}{\kappa_1}a_2x_6 + \frac{\kappa_1\kappa_3}{\kappa_2}a_3x_7 + \frac{\kappa_1\kappa_2}{\kappa_3}a_4x_8}.
\end{eqnarray*}
Thus $s \in \Tran_T(y, Y)$ if and only if ${\kappa_1}^2{\kappa_2}^2{\kappa_3}^2 = \frac{{\kappa_1}^2}{{\kappa_2}^{2}{\kappa_3}^{2}} = \frac{{\kappa_2}^2}{{\kappa_1}^{2}{\kappa_3}^{2}} = \frac{{\kappa_3}^2}{{\kappa_1}^{2}{\kappa_2}^{2}}$, which is true if and only if ${\kappa_1}^4 = {\kappa_2}^4 = {\kappa_3}^4 = \pm 1$; thus $\Tran_T(y, Y)$ is finite. Moreover $s \in C_T(y)$ if and only if $\kappa_1\kappa_2\kappa_3 = \frac{\kappa_1}{\kappa_2\kappa_3} = \frac{\kappa_2}{\kappa_1\kappa_3} = \frac{\kappa_3}{\kappa_1\kappa_2}$, which is true if and only if ${\kappa_1}^2 = {\kappa_2}^2 = {\kappa_3}^2 = \pm 1$; so if we write $h^\dagger = h_{\alpha_1}(\eta_4) h_{\alpha_2}(\eta_4) h_{\alpha_3}(\eta_4)$ then
$$
C_T(y) = \langle h_{\alpha_1}(-1), h_{\alpha_2}(-1), h_{\alpha_3}(-1), h^\dagger \rangle = Z(G) \langle h^\dagger \rangle.
$$
Also we see that $n_{\alpha_1}$, $n_{\alpha_2}$ and $n_{\alpha_3}$ send $y_\a$ to $y_{\a'}$, where $\a' = (a_2, a_1, -a_4, -a_3)$, $(a_3, -a_4, a_1, -a_2)$ and $(a_4, -a_3, -a_2, a_1)$ respectively. Thus each $n_{\alpha_i}$ preserves $\hat Y$, and it follows that $\Tran_N(y, Y)$ is finite. Moreover if we write $n^\dagger = n_{\alpha_1} n_{\alpha_2} n_{\alpha_3}$ then $n^\dagger \in C_N(y)$. Set $C = Z(G) \langle h^\dagger, n^\dagger \rangle$; we shall show that $C_G(y) = C$.

We have
\begin{eqnarray*}
s n_{\alpha_1}.v^{(1)} & = & {\ts \kappa_1\kappa_2\kappa_3 a_2x_1 + \frac{\kappa_1}{\kappa_2\kappa_3}a_1x_2 - \frac{\kappa_2}{\kappa_1\kappa_3}a_4x_3 - \frac{\kappa_3}{\kappa_1\kappa_2}a_3x_4}, \\
s n_{\alpha_1}.v^{(2)} & = & {\ts \frac{1}{\kappa_1\kappa_2\kappa_3}a_2x_5 + \frac{\kappa_2\kappa_3}{\kappa_1}a_1x_6 - \frac{\kappa_1\kappa_3}{\kappa_2}a_4x_7 - \frac{\kappa_1\kappa_2}{\kappa_3}a_3x_8};
\end{eqnarray*}
thus for $s n_{\alpha_1}$ to stabilize $y$ we require $\frac{\kappa_2}{\kappa_1\kappa_3} \frac{a_4}{a_3} = \frac{\kappa_3}{\kappa_1\kappa_2} \frac{a_3}{a_4}$ and $\frac{\kappa_1\kappa_3}{\kappa_2} \frac{a_4}{a_3} = \frac{\kappa_1\kappa_2}{\kappa_3} \frac{a_3}{a_4}$, whence $\frac{{a_3}^2}{{a_4}^2} = \frac{{\kappa_2}^2}{{\kappa_3}^2} = \frac{{a_4}^2}{{a_3}^2}$, so that $\frac{{a_3}^2}{{a_4}^2} = \pm1$, contrary to the definition of $\hat Y$. Similarly no element $s n_{\alpha_2}$ or $s n_{\alpha_3}$ can stabilize $y$; and as $n^\dagger$ does stabilize $y$ it likewise follows that no element $s n_{\alpha_i} n_{\alpha_j}$ for $i \neq j$ can stabilize $y$. Thus $C_N(y) = C$.

Now take $y \in \hat Y$ and $g \in \Tran_G(y, Y)$, and write $g = u_1nu_2$ with $u_1 \in U$, $n \in N$ and $u_2 \in U_w$ where $w = nT \in W$. Write $n = sn'$ where $s \in T$ and $n' = {n_{\alpha_1}}^{i_1} {n_{\alpha_2}}^{i_2} {n_{\alpha_3}}^{i_3}$ for some $i_1, i_2, i_3 \in \{ 0, 1 \}$, so that $n'.y \in \hat Y$. Set $y' = g.y \in Y$ and
$$
g' = s^{-1}g{n'}^{-1} = {u_1}^s.{u_2}^{{n'}^{-1}} = \prod_{i = 1}^3 x_{\alpha_i}(t_i).\prod_{i = 1}^3 x_{-\alpha_i}({t_i}');
$$
then $s^{-1}.y' = g'.(n'.y)$.

For convenience write $n'.y = y_\a = \langle v^{(1)}, v^{(2)} \rangle$ as above. For $i = 1, 2, 3$ set ${t_i}'' = t_i{t_i}' + 1$. For $i = 1, 2$ we have $g'.v^{(i)} = \sum_j \tilde a_{ij}x_j$,
where
\begin{eqnarray*}
\tilde a_{11} & = & a_1{t_1}''{t_2}''{t_3}'' + a_2{t_1}''t_2t_3 + a_3t_1{t_2}''t_3 + a_4t_1t_2{t_3}'', \\
\tilde a_{12} & = & a_1{t_1}''{t_2}'{t_3}' + a_2{t_1}'' + a_3t_1{t_2}' + a_4t_1{t_3}', \\
\tilde a_{13} & = & a_1{t_1}'{t_2}''{t_3}' + a_2{t_1}'t_2 + a_3{t_2}'' + a_4t_2{t_3}', \\
\tilde a_{14} & = & a_1{t_1}'{t_2}'{t_3}'' + a_2{t_1}'t_3 + a_3{t_2}'t_3 + a_4{t_3}'', \\
\tilde a_{15} & = & a_1{t_1}'{t_2}'{t_3}' + a_2{t_1}' + a_3{t_2}' + a_4{t_3}', \\
\tilde a_{16} & = & a_1{t_1}'{t_2}''{t_3}'' + a_2{t_1}'t_2t_3 + a_3{t_2}''t_3 + a_4t_2{t_3}'', \\
\tilde a_{17} & = & a_1{t_1}''{t_2}'{t_3}'' + a_2{t_1}''t_3 + a_3t_1{t_2}'t_3 + a_4t_1{t_3}'', \\
\tilde a_{18} & = & a_1{t_1}''{t_2}''{t_3}' + a_2{t_1}''t_2 + a_3t_1{t_2}'' + a_4t_1t_2{t_3}',
\end{eqnarray*}
and
\begin{eqnarray*}
\tilde a_{21} & = & a_1t_1t_2t_3 + a_2t_1{t_2}''{t_3}'' + a_3{t_1}''t_2{t_3}'' + a_4{t_1}''{t_2}''t_3, \\
\tilde a_{22} & = & a_1t_1 + a_2t_1{t_2}'{t_3}' + a_3{t_1}''{t_3}' + a_4{t_1}''{t_2}', \\
\tilde a_{23} & = & a_1t_2 + a_2{t_2}''{t_3}' + a_3{t_1}'t_2{t_3}' + a_4{t_1}'{t_2}'', \\
\tilde a_{24} & = & a_1t_3 + a_2{t_2}'{t_3}'' + a_3{t_1}'{t_3}'' + a_4{t_1}'{t_2}'t_3, \\
\tilde a_{25} & = & a_1 + a_2{t_2}'{t_3}' + a_3{t_1}'{t_3}' + a_4{t_1}'{t_2}', \\
\tilde a_{26} & = & a_1t_2t_3 + a_2{t_2}''{t_3}'' + a_3{t_1}'t_2{t_3}'' + a_4{t_1}'{t_2}''t_3, \\
\tilde a_{27} & = & a_1t_1t_3 + a_2t_1{t_2}'{t_3}'' + a_3{t_1}''{t_3}'' + a_4{t_1}''{t_2}'t_3, \\
\tilde a_{28} & = & a_1t_1t_2 + a_2t_1{t_2}''{t_3}' + a_3{t_1}''t_2{t_3}' + a_4{t_1}''{t_2}''.
\end{eqnarray*}
Since we require $g'.(n'.y) = s^{-1}.y'$, the right hand side of which has basis vectors lying in $\langle x_1, x_2, x_3, x_4 \rangle$ and $\langle x_5, x_6, x_7, x_8 \rangle$, the projections of $g'.v^{(1)}$ and $g'.v^{(2)}$ on $\langle x_1, x_2, x_3, x_4 \rangle$ must be linearly dependent, as must those on $\langle x_5, x_6, x_7, x_8 \rangle$; thus if for $i, j \leq 8$ we write
$$
A_{ij} = \tilde a_{1i}\tilde a_{2j} - \tilde a_{1j}\tilde a_{2i},
$$
then we must have $A_{ij} = 0$ whenever either $i, j \leq 4$, or $i, j \geq 5$.

First assume $p = 2$. Here the equations $A_{65} + A_{34} = 0$, $A_{75} + A_{24} = 0$ and $A_{85} + A_{23} = 0$ simplify to $({a_1}^2 + {a_2}^2 + {a_3}^2 + {a_4}^2){t_i}' = 0$ for $i = 1, 2, 3$ respectively; so we must have ${t_1}' = {t_2}' = {t_3}' = 0$. Now $A_{65} = 0$ and $A_{75} = 0$ reduce to $a_3t_3 = a_4t_2$ and $a_2t_3 = a_4t_1$, so $t_1 = a_2t$, $t_2 = a_3t$, $t_3 = a_4t$ for some $t \in K$; then $A_{12} - {t_1}^2A_{65} = 0$ gives ${a_1}^2t_1 + {a_2}^2t_1 + a_2a_3t_2 + a_2a_4t_3 = 0$, whence $({a_1}^2 + {a_2}^2 + {a_3}^2 + {a_4}^2)t = 0$, and so $t = 0$ and hence $t_1 = t_2 = t_3 = 0$. Thus $g' = 1$; so in this case $\Tran_G(y, Y) = \Tran_N(y, Y)$.

Now assume $p \geq 3$. To begin with, suppose $t_1 = 0$. Then ${t_1}'(A_{12} - A_{87}) = 0$ gives $2a_2a_3{t_1}'t_2 = 2a_1a_4{t_1}'{t_2}'{t_2}''$, while $t_2{t_2}''(A_{75} - A_{24}) - {t_2}'(A_{13} - A_{86}) = 0$ gives $2a_2a_3{t_1}'t_2 = -2a_1a_4{t_1}'{t_2}'{t_2}''$, so we must have ${t_1}'t_2 = {t_1}'{t_2}'{t_2}'' = 0$. If we had ${t_1}' \neq 0$ this would force $t_2 = 0 = {t_2}'$, but then $A_{65} - A_{34} = 0$ would give $({a_1}^2 - {a_2}^2 - {a_3}^2 + {a_4}^2){t_1}' = 0$, contrary to assumption; so we must have ${t_1}' = 0$. Now $A_{75} - A_{24} = 0$ gives $({a_1}^2 - {a_2}^2 - {a_3}^2 + {a_4}^2){t_2}' = 0$, so ${t_2}' = 0$, and $A_{13} - A_{86} = 0$ gives $({a_1}^2 - {a_2}^2 - {a_3}^2 + {a_4}^2)t_2{t_2}'' = 0$, so $t_2 = 0$; similarly we obtain $t_3 = {t_3}' = 0$, so $g' = 1$. Thus we may suppose $t_1 \neq 0$ (and similarly $t_2, t_3 \neq 0$).

\pagebreak

For $i = 1, 2, 3$ write $\bar t_i = \frac{{t_i}'{t_i}''}{t_i}$. Now $t_1{t_1}''(A_{65} - A_{34}) - {t_1}'(A_{12} - A_{87}) = 0$ gives
$$
2(a_1a_4(1 - \bar t_1 \bar t_2) + a_2a_3(\bar t_1 - \bar t_2)) = 0,
$$
while $t_2{t_2}''(A_{75} - A_{24}) - {t_2}'(A_{13} - A_{86}) = 0$ gives
$$
2(a_1a_4(1 - \bar t_1 \bar t_2) - a_2a_3(\bar t_1 - \bar t_2)) = 0;
$$
so $\bar t_1 \bar t_2 = 1$ and $\bar t_1 = \bar t_2$, whence $\bar t_1 = \bar t_2 = \e \in \{ \pm1 \}$ (and similarly $\bar t_3 = \e$). Next we find that
$$
t_2(2{t_1}'' - 1)[A_{12} - A_{87} - {t_1}^2(A_{65} - A_{34})] - t_1(2{t_2}'' - 1)[A_{13} - A_{86} - {t_2}^2(A_{75} - A_{24})] = 0
$$
gives
$$
2(a_2a_3({t_1}^2(2{t_2}'' - 1)^2 - {t_2}^2(2{t_1}'' - 1)^2) + a_1a_4(t_2{t_2}'{t_2}'' - t_1{t_1}'{t_1}'')) = 0,
$$
which now reduces to $(a_2a_3 - \e a_1a_4)(t_2{t_2}'{t_2}'' - t_1{t_1}'{t_1}'') = 0$; since by assumption $a_2a_3 \neq \pm a_1a_4$, we must have $t_2{t_2}'{t_2}'' = t_1{t_1}'{t_1}''$, so $(1 + 2t_2{t_2}')^2 = (1 + 2t_1{t_1}')^2$, and thus $1 + 2t_2{t_2}' = \e'(1 + 2t_1{t_1}')$ for some $\e' \in \{ \pm1 \}$.

Suppose if possible that $1 + 2t_1{t_1}' \neq 0$. Then $A_{12} - A_{87} - {t_1}^2(A_{65} - A_{34}) = 0$ gives
$$
({a_1}^2 - {a_2}^2 - {a_3}^2 + {a_4}^2)t_1 = 2\e'(a_2a_3t_2 - a_1a_4{t_2}'{t_2}'') = 2\e'(a_2a_3 - \e a_1a_4)t_2,
$$
while $A_{13} - A_{86} - {t_2}^2(A_{75} - A_{24}) = 0$ gives
$$
({a_1}^2 - {a_2}^2 - {a_3}^2 + {a_4}^2)t_2 = 2\e'(a_2a_3t_1 - a_1a_4{t_1}'{t_1}'') = 2\e'(a_2a_3 - \e a_1a_4)t_1;
$$
thus we have
$$
\frac{{a_1}^2 - {a_2}^2 - {a_3}^2 + {a_4}^2}{2(a_2a_3 - \e a_1a_4)} = \e'\frac{t_2}{t_1} = \frac{2(a_2a_3 - \e a_1a_4)}{{a_1}^2 - {a_2}^2 - {a_3}^2 + {a_4}^2},
$$
whence
$$
({a_1}^2 - {a_2}^2 - {a_3}^2 + {a_4}^2)^2 = 4(a_2a_3 - \e a_1a_4)^2,
$$
and so
$$
(({a_1}^2 - {a_2}^2 - {a_3}^2 + {a_4}^2)^2 - 4({a_1}^2{a_4}^2 + {a_2}^2{a_3}^2))^2 = 64{a_1}^2{a_2}^2{a_3}^2{a_4}^2,
$$
which upon expansion gives
$$
{\ts \sum {a_i}^8 + 6\sum {a_i}^4{a_j}^4 - 4\sum {a_i}^6{a_j}^2 + 4\sum {a_i}^4{a_j}^2{a_l}^2 - 40{a_1}^2{a_2}^2{a_3}^2{a_4}^2 = 0},
$$
contrary to the final condition in the definition of $\hat Y$. Thus we must have $1 + 2t_1{t_1}' = 0 = 1 + 2t_2{t_2}'$ (and similarly $1 + 2t_3{t_3}' = 0$); so for $i = 1, 2, 3$ we have $t_i = -\frac{1}{2{t_i}'}$, whence $\e = \bar t_i = -{{t_i}'}^2$. Thus there exist $\e_1, \e_2, \e_3 \in \{ \pm 1 \}$ and $j \in \{ 0, 1 \}$ such that for each $i$ we have ${t_i}' = \e_i{\eta_4}^j$ and $t_i = -\frac{1}{2}\e_i{\eta_4}^{-j}$; in particular there are only finitely many possibilities for each $t_i$ and ${t_i}'$, and hence for the element $g'$.

Write $\a' = ({a_1}', {a_2}', {a_3}', {a_4}') \in K^4$ where
\begin{eqnarray*}
{a_1}' & = & a_1 + \e_2 \e_3 (-1)^j a_2 + \e_1 \e_3 (-1)^j a_3 + \e_1 \e_2 (-1)^j a_4, \\
{a_2}' & = & \e_2 \e_3 (-1)^j a_1 + a_2 - \e_1 \e_2 a_3 - \e_1 \e_3 a_4, \\
{a_3}' & = & \e_1 \e_3 (-1)^j a_1 - \e_1 \e_2 a_2 + a_3 - \e_2 \e_3 a_4, \\
{a_4}' & = & \e_1 \e_2 (-1)^j a_1 - \e_1 \e_3 a_2 - \e_2 \e_3 a_3 + a_4;
\end{eqnarray*}
thus each ${a_i}'$ is of the form $\pm a_1 \pm a_2 \pm a_3 \pm a_4$. Take $c \in K^*$ with $c^2 = 2$ and let $s_1 = h_{\alpha_1}(c) h_{\alpha_2}(c) h_{\alpha_3}(c)$. Then we find that
$$
s_1g'.(n'.y) = y_{\a'} \in Y.
$$
Now suppose $g \in C_G(y)$. Since all ${t_i}'$ are non-zero we know that $n' = n_{\alpha_1} n_{\alpha_2} n_{\alpha_3} = n^\dagger$, which fixes $y$, so that in fact $y = n'.y = y_\a$; also we must have $sg' \in C_G(y)$, so that $s{s_1}^{-1}.y_{\a'} = y_\a$. Write $s{s_1}^{-1} = h_{\alpha_1}(\kappa_1) h_{\alpha_2}(\kappa_2) h_{\alpha_3}(\kappa_3)$ for $\kappa_1, \kappa_2, \kappa_3 \in K^*$; then the equations given earlier detailing the effect of an element of $T$ on points in $Y$ show that we require $\frac{\kappa_2}{\kappa_1\kappa_3} \frac{{a_3}'}{a_3} = \frac{\kappa_3}{\kappa_1\kappa_2} \frac{{a_4}'}{a_4}$ and $\frac{\kappa_1\kappa_3}{\kappa_2} \frac{{a_3}'}{a_3} = \frac{\kappa_1\kappa_2}{\kappa_3} \frac{{a_4}'}{a_4}$, whence $\frac{a_3{a_4}'}{a_4{a_3}'} = \frac{{\kappa_2}^2}{{\kappa_3}^2} = \frac{a_4{a_3}'}{a_3{a_4}'}$, so that $\frac{a_3{a_4}'}{a_4{a_3}'} = \pm 1$, contrary to the definition of $\hat Y$. Therefore the elements in $\Tran_G(y, Y)$ with $t_i = -\frac{1}{2{t_i}'}$ do not in fact stabilize $y$; so $C_G(y) \leq N$.

Thus in each case $\Tran_G(y, Y)$ is finite; so
$$
\codim {\ts\Tran_G(y, Y)} = \dim G - \dim {\ts\Tran_G(y, Y)} = 9 - 0 = 9
$$
while
$$
\codim Y = \dim \G{2}(V) - \dim Y = 12 - 3 = 9.
$$
Therefore $y$ is $Y$-exact. Moreover we have shown that $C_G(y) = C$. Thus the conditions of Lemma~\ref{lem: generic stabilizer from exact subset} hold; so the quadruple $(G, \lambda, p, k)$ has generic stabilizer $C/Z(G) \cong \Z_{2/(p, 2)}.\Z_2$.
\end{proof}

\begin{prop}\label{prop: A_2.2, omega_1 + omega_2 module, p = 3, k = 2}
Let $G = A_2$ and $\lambda = \omega_1 + \omega_2$ with $p = 3$, and take $k = 2$; let $\tau$ be a graph automorphism of $G$. Then the quadruple $(G\langle \tau \rangle, \lambda, p, k)$ has generic stabilizer $\Z_2$.
\end{prop}

\begin{proof}
We shall follow the strategy used in the second part of the proof of Proposition~\ref{prop: A_2, 2omega_1 module, k = 3, A_4, omega_2 module, k = 5}. Let $G$ have simple roots $\alpha_1$ and $\alpha_2$, and $V$ be $\L(G)/Z(\L(G))$. Since $\dim V = 7$ we have $\dim \G{2}(V) = 10$. Let $\tau$ be the automorphism of $G$ which interchanges $x_{\alpha_1}(t)$ with $x_{\alpha_2}(t)$. We have $Z(G) = \{ 1 \}$.

We shall write elements of both $G$ and $\L(G)$ as $3 \times 3$ matrices, so that $Z(\L(G)) = \langle I \rangle$. We let $T < G$ be the subgroup of diagonal matrices, so that $N$ is the subgroup of monomial matrices; define $n_0 \in N$ by
$$
n_0 =
-\left(
  \begin{array}{ccc}
      &   & 1 \\
      & 1 &   \\
    1 &   &   \\
  \end{array}
\right),
$$
so that $n_0T$ is the long word $w_0$ of the Weyl group. We see that $n_0\tau$ acts on $\L(G)$ by sending matrices to their transposes.

Given $\a = (a_1, a_2, a_3) \in K^3 \setminus \{ (0, 0, 0) \}$, define
$$
v^{(1)} =
\left(
  \begin{array}{ccc}
       &     & a_2 \\
   a_3 &     &     \\
       & a_1 &     \\
  \end{array}
\right) + Z(\L(G)),
\qquad
v^{(2)} =
\left(
  \begin{array}{ccc}
       & a_3 &     \\
       &     & a_1 \\
   a_2 &    &     \\
  \end{array}
\right) + Z(\L(G)),
$$
and set $y_\a = \langle v^{(1)}, v^{(2)} \rangle$; write $Y = \{ y_\a : \a \neq (0, 0, 0) \}$, then $Y$ is a subvariety of $\G{2}(V)$ of dimension $2$, whence $\codim Y = 8 = \dim G$. Set $C = \langle n_0\tau \rangle$; then each $y \in Y$ is stabilized by $C$.

Define $\S = \L(T)$. Clearly if $\alpha \in \Phi$ then there exists $h \in \S$ with $[he_\alpha] \neq 0$, so $C_{\L(G)}(\S) = \L(T)$. Since for $i = 1, 2$ we have $h_{\alpha_i}.v^{(1)} = v^{(1)}$ and $h_{\alpha_i}.v^{(2)} = -v^{(2)}$, if $y \in Y$ we have $\S \leq \Ann_{\L(G)}(y)$. Set
$$
\hat Y = \{ y_\a \in Y : a_1a_2a_3 \neq 0, \ {a_1}^2 + {a_2}^2 + {a_3}^2 \neq 0, \ a_i \neq \pm a_{i'} \hbox{ for } i \neq i' \};
$$
then $\hat Y$ is a dense open subset of $Y$. Take $y = y_\a \in \hat Y$.

First suppose $x \in \Ann_{\L(G)}(y)$; write $x = h + e$ where $h \in \L(T)$ and $e \in \langle e_\alpha : \alpha \in \Phi \rangle$. For $i = 1, 2$, since the difference of two weights occurring in $v^{(i)}$ is never a root, we see that the weights occurring in $e.v^{(i)}$ must be a subset of those occurring in $V^{(3 - i)}$ together with the zero weight; so $e.v^{(i)}$ must be a scalar multiple of $v^{(3 - i)}$. If we write $e = \sum_{\alpha \in \Phi} t_\alpha e_\alpha$, then this condition may be expressed in matrix form as $A{\bf t} = {\bf 0}$, where $A$ is a $6 \times 6$ matrix and ${\bf t}$ is a column vector whose entries are the various coefficients $t_\alpha$. We find that if the rows and columns of $A$ are suitably ordered then it becomes block diagonal, having $2$ blocks, with each block being a $3 \times 3$ matrix. In fact each block may be written in the form
$$
\left(
  \begin{array}{ccc}
           a_1a_3       & a_2a_3 & -{a_1}^2 - {a_2}^2 \\
     -{a_2}^2 - {a_3}^2 & a_1a_2 &       a_1a_3       \\
             a_1        &   a_2  &         a_3        \\
  \end{array}
\right),
$$
which has determinant $a_2({a_1}^2 + {a_2}^2 + {a_3}^2)^2$. Thus the definition of the set $\hat Y$ implies that each block of $A$ is non-singular, as therefore is $A$ itself; so ${\bf t}$ must be the zero vector and hence $e = 0$. Thus $x = h + e \in \S$; so $\Ann_{\L(G)}(y) = \S$.

A straightforward calculation shows that $C_T(y) = \{ 1 \}$, and $\Tran_T(y, Y) = \langle h_{\alpha_1}(-1), h_{\alpha_2}(-1) \rangle$, so that $T.y \cap Y = \{ y_{\a'}: \a' = (\pm a_1, \pm a_2, \pm a_3) \}$. We claim that $N.y \cap Y \subset \hat Y$, and $C_N(y) = C$. Take $n \in \Tran_G(y, Y)$ and let $n.y = y_{\a'}$. We find that there exists $\pi \in S_3$ such that $\a' = (\pm a_{\pi(1)}, \pm a_{\pi(2)}, \pm a_{\pi(3)})$, with $\pi = 1$ if and only if $nT = T$. Since $y_{\a'} \in \hat Y$ this proves the first claim; moreover the definition of $\hat Y$ shows that if $nT \neq T$ then $n$ does not stabilize $y$, proving the second. Thus the conditions of Lemma~\ref{lem: exactness via Premet} hold, so that $\Tran_G(y, Y) \subseteq N$, and $y$ is $Y$-exact; moreover $C_G(y) = C_N(y) = C$. Therefore the conditions of Lemma~\ref{lem: generic stabilizer from exact subset} hold; so the quadruple $(G\langle \tau \rangle, \lambda, p, k)$ has generic stabilizer $C/Z(G) \cong \Z_2$.
\end{proof}

As explained at the start of this chapter, the remaining large higher quadruples will be handled in the following section along with the small higher quadruples.

\section{Small higher quadruples}\label{sect: small higher quadruples}

In this final section we shall treat small higher quadruples along with the remaining large higher quadruples, and establish the entries in Tables~\ref{table: large higher quadruple non-TGS}, \ref{table: small classical higher quadruple generic stab} and \ref{table: small exceptional higher quadruple generic stab}, thus proving Theorems~\ref{thm: large higher quadruple generic stab} and \ref{thm: small higher quadruple generic stab}. In most cases we shall apply Lemma~\ref{lem: generic stabilizer from exact subset} to determine the generic stabilizer.

We begin with the cases where $G$ is a classical group and $V$ is the natural module. In the statement of the following result, for convenience we refer to the cases where $G = C_2$, $\lambda = \omega_1$ and $G = D_3$, $\lambda = \omega_1$; these appear in Table~\ref{table: small classical higher quadruple generic stab} as $G = B_2$, $\lambda = \omega_2$ and $G = A_3$, $\lambda = \omega_2$ respectively.

\begin{prop}\label{prop: natural modules, k arbitrary}
Let $G = A_\ell$ for $\ell \in [1, \infty)$, or $G = B_\ell$ for $\ell \in [2, \infty)$ with $p \geq 3$, or $G = C_\ell$ for $\ell \in [2, \infty)$, or $G = D_\ell$ for $\ell \in [3, \infty)$, and $\lambda = \omega_1$. Then the quadruple $(G, \lambda, p, k)$ has generic stabilizer $A_{\ell - k} A_{k - 1} T_1 U_{k(\ell + 1 - k)}$, or $B_{\frac{1}{2}(k - 1)} D_{\ell - \frac{1}{2}(k - 1)}.\Z_2$, or $C_{\frac{1}{2}(k - 1)} C_{\ell - \frac{1}{2}(k + 1)} T_1 U_{2\ell - 1}$, or $B_{\frac{1}{2}(k - 1)} B_{\ell - \frac{1}{2}(k + 1)}$, respectively if $k$ is odd, and $A_{\ell - k} A_{k - 1} T_1 U_{k(\ell + 1 - k)}$, or $D_{\frac{1}{2}k} B_{\ell - \frac{1}{2}k}.\Z_2$, or $C_{\frac{1}{2}k} C_{\ell - \frac{1}{2}k}$, or $D_{\frac{1}{2}k} D_{\ell - \frac{1}{2}k}.\Z_2$, respectively if $k$ is even.
\end{prop}

\begin{proof}
In all these cases $V$ is the natural module for $G$.

If $G = A_\ell$, then $G$ acts transitively on $k$-dimensional subspaces of $V$; if we take $y = \langle v_1, \dots, v_k \rangle$ where $v_1, \dots, v_{\ell + 1}$ is the natural basis of $V_{nat}$, then the stabilizer of $y$ is the maximal parabolic subgroup $A_{\ell - k} A_{k - 1} T_1 U_{k(\ell + 1 - k)}$ corresponding to the $k$th simple root.

Next suppose $G = C_\ell$. If $k = 2j$ is even, then $G$ acts transitively on non-singular $k$-dimensional subspaces of $V$, which form a dense open subset of $\Gk(V)$, and the stabilizer of one such is $C_j C_{\ell - j}$. If instead $k = 2j + 1$ is odd, then $G$ acts transitively on $k$-dimensional subspaces of $V$ with $1$-dimensional radical, which form a dense open subset of $\Gk(V)$, and the stabilizer of one such is $C_j C_{\ell - j - 1} T_1 U_{2\ell - 1}$ (if we take the subspace $\langle e_1, e_2, f_2, \dots, e_{j + 1}, f_{j + 1} \rangle$ where $e_1, f_1, \dots, e_\ell, f_\ell$ is the natural basis of $V_{nat}$, then the $C_j$ and $C_{\ell - j - 1}$ factors correspond to the non-singular subspaces $\langle e_2, f_2, \dots, e_{j + 1}, f_{j + 1} \rangle$ and $\langle e_{j + 2}, f_{j + 2}, \dots, e_\ell, f_\ell \rangle$ respectively, while the unipotent radical $U_{2\ell - 1}$ is the product of the root groups $X_\alpha$ for $\alpha \in \{ 2\ve_1, \ve_1 \pm \ve_i : 2 \leq i \leq \ell \}$).

Now suppose $G = B_\ell$ or $D_\ell$, and as usual write $d = \dim V$; then $G$ acts transitively on non-singular $k$-dimensional subspaces of $V$, which form a dense open subset of $\Gk(V)$, and the stabilizer in $\mathrm{O}_d(K)$ of such a subspace is $\mathrm{O}_k(K) \mathrm{O}_{d - k}(K)$, so we need to consider the intersection of this product with $G$.

First assume $p \geq 3$; then we may take $G = \SO_d(K)$, and the intersection is $\SO_k(K) \SO_{d - k}(K) \cup (\mathrm{O}_k(K) \setminus \SO_k(K))(\mathrm{O}_{d - k}(K) \setminus \SO_{d - k}(K))$. For each $a \in \N$ we have $Z(\mathrm{O}_a(K)) = \{ \pm I_a \}$, and $-I_a \in \SO_a(K)$ if and only if $a$ is even. Thus if $G = D_\ell$ and $k = 2j + 1$ is odd, the stabilizer is $\SO_{2j + 1}(K) \SO_{2\ell - 2j - 1}(K) \cup (-\SO_{2j + 1}(K))(-\SO_{2\ell - 2j - 1}(K)) = \{ \pm I_{2\ell} \} \SO_{2j + 1}(K) \SO_{2\ell - 2j - 1}(K)$; taking the quotient by $Z(G) = \{ \pm I_{2\ell} \}$ we see that the generic stabilizer is simply $B_j B_{\ell - j - 1}$. In the other cases the stabilizer does not have the form $\{ \pm I_d \} \SO_k(K) \SO_{d - k}(K)$, so the generic stabilizer is $D_j D_{\ell - j}.\Z_2$ if $G = D_\ell$ and $k = 2j$ is even, $B_j D_{\ell - j}.\Z_2$ if $G = B_\ell$ and $k = 2j + 1$ is odd, and $D_j B_{\ell - j}.\Z_2$ if $G = B_\ell$ and $k = 2j$ is even.

Finally assume $p = 2$, so that $G = D_\ell$. If $k = 2j + 1$ is odd, then the stabilizer is just $\SO_{2j + 1}(K) \SO_{2\ell - 2j - 1}(K)$, so the generic stabilizer is again simply $B_j B_{\ell - j - 1}$. If however $k = 2j$ is even, and we write the non-singular $k$-dimensional subspace and its complement as $\langle v_1, v_{-1}, \dots, v_j, v_{-j} \rangle$ and $\langle v_{j + 1}, v_{-(j + 1)}, \dots, v_\ell, v_{-\ell} \rangle$ respectively, then there is an element $n = n_{\ve_1 - \ve_\ell} n_{\ve_1 + \ve_\ell}$ of $N$ which interchanges $v_1$ with $v_{-1}$, and $v_\ell$ with $v_{-\ell}$, while fixing all other basis vectors; thus $n$ lies in the intersection required, but acts as a single transposition on the basis vectors of both the subspace and its complement, so does not lie in $D_j D_{\ell - j}$, whence the generic stabilizer is again $D_j D_{\ell - j}.\Z_2$.
\end{proof}

\begin{prop}\label{prop: natural module for B_ell, p = 2, k arbitrary}
Let $G = B_\ell$ for $\ell \in [2, \infty)$ and $\lambda = \omega_1$ with $p = 2$. Then if $k$ is odd the quadruple $(G, \lambda, p, k)$ has generic stabilizer $B_{\frac{1}{2}(k - 1)} B_{\ell - \frac{1}{2}(k + 1)} T_1 U_{2\ell - 1}$, and if $k$ is even the quadruple $(G, \lambda, p, k)$ has generic stabilizer $B_{\frac{1}{2}k} B_{\ell - \frac{1}{2}k}$.
\end{prop}

\begin{proof}
This is an immediate consequence of Proposition~\ref{prop: natural modules, k arbitrary}, using the exceptional isogeny $B_\ell \to C_\ell$ which exists in characteristic $2$.
\end{proof}

Next we consider the remaining cases which occur in infinite families.

\begin{prop}\label{prop: A_ell, ell odd, omega_2 module, k = 2}
Let $G = A_\ell$ for odd $\ell \in [5, \infty)$ and $\lambda = \omega_2$, and take $k = 2$. Then according as $\ell = 5$, or $\ell = 7$, or $\ell \geq 9$, the quadruple $(G, \lambda, p, k)$ has generic stabilizer ${A_1}^3.S_3$, or ${A_1}^4.{\Z_2}^2$, or ${A_1}^{\frac{1}{2}(\ell + 1)}$, respectively.
\end{prop}

\begin{proof}
We take $G = \SL_{\ell + 1}(K)$ and use the set-up of Proposition~\ref{prop: A_ell, 2 omega_1 and omega_2 modules}: we write $\ell = 2\ell_1 - 1$, so that $\ell_1 \geq 3$; for $1 \leq i < j \leq 2\ell_1$ we write $\bar v_{i, j} = v_i \wedge v_j$, where $v_1, \dots, v_{2\ell_1}$ is the standard basis of $V_{nat}$; we take the generalized height function on the weight lattice of $G$ whose value at each simple root $\alpha_i$ is $2$; for $i = 1, \dots, \ell_1$ we set $x_i = \bar v_{i, 2\ell_1 + 1 - i}$ and let $\nu_i$ be the weight such that $x_i \in V_{\nu_i}$; then $\Lambda(V)_{[0]} = \{ \nu_1, \dots, \nu_{\ell_1} \}$ has ZLC because $\nu_1 + \cdots + \nu_{\ell_1} = 0$, and $V_{[0]} = \langle x_1, \dots, x_{\ell_1} \rangle$; and the setwise stabilizer in $W$ of $\Lambda(V)_{[0]}$ is $\langle w_{\alpha_{\ell_1}}, w_{\alpha_{\ell_1 - 1}} w_{\alpha_{\ell_1 + 1}}, \dots, w_{\alpha_1} w_{\alpha_{2\ell_1 - 1}} \rangle$. Here however we take $Y = \G{2}(V_{[0]})$, and write
$$
\hat Y_1 = \left\{ y = \langle v^{(1)}, v^{(2)} \rangle \in Y : v^{(1)} = {\ts\sum} a_i x_i, \ v^{(2)} = {\ts\sum} b_i x_i, \ \forall i \neq j \
\left|
\begin{array}{cc}
 a_i & a_j \\
 b_i & b_j \\
\end{array}
\right| \neq 0 \right\};
$$
then $\hat Y_1$ is a dense open subset of $Y$, and the determinant condition implies that each $\nu_i$ occurs in every $y \in \hat Y_1$. We have $Z(G) = \langle z \rangle$ where $z = \eta_{2\ell_1}I$.

Let $A$ be the ${A_1}^{\ell_1}$ subgroup having simple roots $\alpha_{\ell_1}$, $\alpha_{\ell_1 - 1} + \alpha_{\ell_1} + \alpha_{\ell_1 + 1}$, \dots, $\alpha_1 + \cdots + \alpha_{2\ell_1 - 1}$; then clearly for all $y \in Y$ we have $A \leq C_G(y)$.

Take $y \in \hat Y_1$ and $g \in \Tran_G(y, Y)$, and write $y' = g.y \in Y$. By Lemma~\ref{lem: gen height zero} we have $g = u_1 n u_2$ with $u_1 \in C_U(y')$, $u_2 \in C_U(y)$, and $n \in N_{\Lambda(V)_{[0]}}$ with $n.y = y'$.

First, the identification of $W_{\Lambda(V)_{[0]}}$ above shows that we have $N_{\Lambda(V)_{[0]}} = \langle T, n_{\alpha_{\ell_1}}, n_{\alpha_{\ell_1 - 1}} n_{\alpha_{\ell_1 + 1}}, \dots, n_{\alpha_1} n_{\alpha_{2\ell_1 - 1}} \rangle = T_{\ell_1 - 1}(A \cap N).S_{\ell_1}$, where we write $T_{\ell_1 - 1} = \{ \prod_{i = 1}^{\ell_1 - 1} h_{\alpha_i}(\kappa_i) : \kappa_1, \dots, \kappa_{\ell_1 - 1} \in K^* \}$ and the symmetric group $S_{\ell_1}$ permutes the simple factors of $A$. As the elements of $N_{\Lambda(V)_{[0]}}$ permute and scale the $x_i$, we have $N_{\Lambda(V)_{[0]}}.y \subseteq \hat Y_1$.

Next, let $\Xi = \Phi^+ \setminus \Phi_A$, and set $U' = \prod_{\alpha \in \Xi} X_\alpha$; then $U = U'.(A \cap U)$ and $U' \cap (A \cap U) = \{ 1 \}$. We now observe that if $\alpha \in \Xi$ then $\nu_i + \alpha$ is a weight in $V$ for exactly one value of $i$; moreover each weight in $V$ of positive generalized height is of the form $\nu_i + \alpha$ for exactly two such roots $\alpha$. Thus if we take $u = \prod x_\alpha(t_\alpha) \in U'$ satisfying $u.y = y$, and equate coefficients of weight vectors, taking them in an order compatible with increasing generalized height, using the determinant condition in the definition of the set $\hat Y_1$ we see that for all $\alpha$ we must have $t_\alpha = 0$, so that $u = 1$; so $C_U(y) = A \cap U$. Since the previous paragraph shows that $y' = g'.y \in \hat Y_1$, likewise we have $C_U(y') = A \cap U$.

Thus $\Tran_G(y, Y) = A T_{\ell_1 - 1}.S_{\ell_1} \cong {A_1}^{\ell_1} T_{\ell_1 - 1}.S_{\ell_1}$. Hence
$$
\codim {\ts\Tran_G}(y, Y) = \dim G - \dim {\ts\Tran_G}(y, Y) = (4{\ell_1}^2 - 1) - (4\ell_1 - 1) = 4{\ell_1}^2 - 4\ell_1,
$$
while
$$
\codim Y = \dim \G{2}(V) - \dim\G{2}(V_{[0]}) = 2(2{\ell_1}^2 - \ell_1 - 2) - 2(\ell_1 - 2) = 4{\ell_1}^2 - 4\ell_1.
$$
Therefore $y$ is $Y$-exact.

We now consider stabilizers; by the above, for all $y \in \hat Y_1$ we have $Z(G) A \leq C_G(y) \leq A T_{\ell_1 - 1}.S_{\ell_1}$. Let $C$ be the subgroup $Z(G)A\langle h_{\alpha_1}(-1) n_{\alpha_1} n_{\alpha_5}, n_{\alpha_2} n_{\alpha_4} \rangle$, $Z(G)A\langle n_{\alpha_1} n_{\alpha_7} n_{\alpha_3} n_{\alpha_5}, h_{\alpha_2 + \alpha_3}(-1) n_{\alpha_1 + \alpha_2} n_{\alpha_6 + \alpha_7} n_{\alpha_2 + \alpha_3} n_{\alpha_5 + \alpha_6} \rangle$ or $Z(G)A$ according as $\ell_1 = 3$, $\ell_1 = 4$ or $\ell_1 \geq 5$. We shall define a dense open subset $\hat Y$ of $Y$ lying in $\hat Y_1$, and show that if $y \in \hat Y$ then $C_G(y)$ is a conjugate of $C$. The argument here is very similar to that in the proof of Proposition~\ref{prop: A_ell, 2omega_1 module, k = 2}.

Write $y = \langle v^{(1)}, v^{(2)} \rangle$ where $v^{(1)} = \sum a_i x_i$, $v^{(2)} = \sum b_i x_i$. Note that the determinant condition defining the set $\hat Y_1$ implies that for each $i$ we cannot have $a_i = b_i = 0$, and for each $i \neq j$ we cannot have either $a_i = a_j = 0$ or $b_i = b_j = 0$. Thus by changing basis we may assume if we wish that $a_1 = b_2 = 1$, $a_2 = b_1 = 0$, in which case $a_3, \dots, a_{\ell_1}, b_3, \dots, b_{\ell_1} \neq 0$. Then if we take $s = \diag(\kappa_1, \dots, \kappa_{2\ell_1}) \in T$, we have $s.x_i = \kappa_i \kappa_{2\ell_1 + 1 - i} x_i$, so for $s.y = y$ we require $\kappa_1 \kappa_{2\ell_1} = \kappa_2 \kappa_{2\ell_1 - 1} = \cdots = \kappa_{\ell_1} \kappa_{\ell_1 + 1}$; since $\kappa_1 \dots \kappa_{2\ell_1} = 1$ there exists $i$ such that $\kappa_1 \kappa_{2\ell_1} = {\eta_{\ell_1}}^i$, so that $z^{-i} s \in A \cap T$, whence $C_T(y) = Z(G) (A \cap T)$.

Suppose $\ell_1 = 3$; here we set $\hat Y = \hat Y_1$. Let $y_0 = \langle x_1 + x_3, x_2 + x_3 \rangle \in \hat Y$. Take $y \in \hat Y$; by the above we may assume $y = \langle x_1 + a_3x_3, x_2 + b_3 x_3 \rangle$ with $a_3, b_3 \neq 0$. Choose $\kappa \in K^*$ satisfying $\kappa^3 = (a_3b_3)^{-1}$ and set $h^{-1} = \diag(\kappa a_3, \kappa b_3, \kappa, 1, 1, 1)$; then $h^{-1}.y = y_0$. As $h_{\alpha_1}(-1) n_{\alpha_1} n_{\alpha_5}$ sends the vector $ax_1 + bx_2 + cx_3$ to $bx_1 + ax_2 + cx_3$, it interchanges $x_1 + x_3$ and $x_2 + x_3$, and thus stabilizes $y_0$; likewise as $n_{\alpha_2} n_{\alpha_4}$ sends the vector $ax_1 + bx_2 + cx_3$ to $ax_1 - cx_2 -bx_3$, it sends $x_1 + x_3$ to $(x_1 + x_3) - (x_1 + x_2)$ and negates $x_2 + x_3$, and thus also stabilizes $y_0$. Hence $C_G(y_0) = C$, and so $C_G(y) = C_G(h.y_0) = {}^h C$.

Now suppose $\ell_1 \geq 4$. Let $C_{\ell_1}$ be the subgroup with short simple root groups $\{ x_{\alpha_i}(t) x_{\alpha_{2\ell_1 + 1 - i}}(t) : t \in K \}$ for $i = 1, \dots, \ell_1 - 1$ and long simple root group $X_{\alpha_{\ell_1}}$. Take $n \in T(C_{\ell_1} \cap N) \setminus T$, and for $i = 1, \dots, \ell_1$ write $n.x_i = \kappa_i x_{\pi(i)}$ for $\kappa_i \in K^*$, where $\pi \in S_{\ell_1} \setminus \{ 1 \}$. If $n.y = y$ there must exist $c_1, c_2, c_3, c_4 \in K$ with $(c_1, c_2), (c_3, c_4) \neq (0, 0)$ such that $n.v^{(1)} = c_1 v^{(1)} + c_2 v^{(2)}$ and $n.v^{(2)} = c_3 v^{(1)} + c_4 v^{(2)}$, whence $\sum \kappa_i a_i x_{\pi(i)} = \sum (c_1 a_i + c_2 b_i)x_i$ and $\sum \kappa_i b_i x_{\pi(i)} = \sum (c_3 a_i + c_4 b_i)x_i$. Thus for all $i \leq \ell_1$ we have $\kappa_i a_i = c_1 a_{\pi(i)} + c_2 b_{\pi(i)}$ and $\kappa_i b_i = c_3 a_{\pi(i)} + c_4 b_{\pi(i)}$, whence $c_1 a_{\pi(i)} b_i + c_2 b_{\pi(i)} b_i = c_3 a_{\pi(i)} a_i + c_4 b_{\pi(i)} a_i$. These are the same equations as we had in the proof of Proposition~\ref{prop: A_ell, 2omega_1 module, k = 2} (with $\ell + 1$ there replaced by $\ell_1$ here); we therefore conclude that, unless $\ell_1 = 4$ and $\pi \in \langle (1\ 2)(3\ 4), (1\ 3)(2\ 4) \rangle$, the points $y$ for which there is a non-zero solution $(c_1, c_2, c_3, c_4)$ form a proper closed subvariety of $Y$. Again we take $\hat Y_2$ to be the intersection of the complements of these proper closed subvarieties as $\pi$ runs through $S_4 \setminus \langle (1\ 2)(3\ 4), (1\ 3)(2\ 4) \rangle$ or $S_{\ell_1} \setminus \{ 1 \}$ according as $\ell_1 = 4$ or $\ell_1 \geq 5$. Then $\hat Y_2$ is a dense open subset of $Y$, as therefore is $\hat Y = \hat Y_1 \cap \hat Y_2$. Thus if $\ell_1 \geq 5$, for all $y \in \hat Y$ we have $C_G(y) = C$.

Now assume $\ell_1 = 4$, and take $y \in \hat Y$; as above we may write $y = \langle v^{(1)}, v^{(2)} \rangle$ where $v^{(1)} = x_1 + a_3 x_3 + a_4 x_4$ and $v^{(2)} = x_2 + b_3 x_3 + b_4 x_4$, with $a_3, b_3, a_4, b_4 \neq 0$ and $a_3b_4 \neq a_4b_3$. Exactly as in the penultimate paragraph of the proof of Proposition~\ref{prop: A_ell, 2omega_1 module, k = 2}, take $c_1, c_2, c_3, c_4 \in K^*$ satisfying ${c_4}^8 = \frac{a_3b_3}{a_4b_4(a_3b_4 - a_4b_3)}$, ${c_3}^4 = \frac{a_4b_4}{a_3b_3}{c_4}^4$, ${c_2}^4 = \frac{b_4}{a_3}(a_3b_4 - a_4b_3){c_4}^4$ and $c_1 = \frac{1}{c_2c_3c_4}$; write $\kappa_1 = \frac{a_3{c_3}^2}{{c_1}^2}$ and $\kappa_2 = \frac{a_4{c_4}^2}{{c_1}^2}$, then we have ${\kappa_1}^2 = {\kappa_2}^2 + 1$. If we now set $h^{-1} = \diag(c_1, c_2, c_3, c_4, c_4, c_3, c_2, c_1) \in G$ and $y' = h^{-1}.y$, we have $y' = \langle {v^{(1)}}', {v^{(2)}}' \rangle$ where ${v^{(1)}}' = x_1 + \kappa_1 x_3 + \kappa_2 x_4$ and ${v^{(2)}}' = x_2 + \kappa_2 x_3 + \kappa_1 x_4$. Now with $n^* = n_{\alpha_1} {n_{\alpha_7}}^{-1} n_{\alpha_3} {n_{\alpha_5}}^{-1}$ we see that $n^*$ sends the vector $ax_1 + bx_2 + cx_3 + dx_4$ to $bx_1 + ax_2 + dx_3 + cx_4$, so we have $n^*.{v^{(1)}}' = {v^{(2)}}'$ and $n^*.{v^{(2)}}' = {v^{(1)}}'$, whence $n^* \in C_G(y')$; with $n^{**} = h_{\alpha_2 + \alpha_3}(-1) n_{\alpha_1 + \alpha_2} {n_{\alpha_6 + \alpha_7}}^{-1} n_{\alpha_2 + \alpha_3} {n_{\alpha_5 + \alpha_6}}^{-1}$ we see that $n^{**}$ sends the vector $ax_1 + bx_2 + cx_3 + dx_4$ to $cx_1 - dx_2 + ax_3 - bx_4$, so we have $n^{**}.{v^{(1)}}' = \kappa_1{v^{(1)}}' - \kappa_2{v^{(2)}}'$ and $n^{**}.{v^{(2)}}' = \kappa_2{v^{(1)}}' - \kappa_1{v^{(2)}}'$, whence $n^{**} \in C_G(y')$. Hence $C_G(y') = Z(G) A \langle n^*, n^{**} \rangle = C$, so $C_G(y) = C_G(h.y') = {}^h C$.

Therefore in all cases, for all $y \in \hat Y$ there exists $h \in T$ with $C_G(y) = {}^h C$. Thus the conditions of Lemma~\ref{lem: generic stabilizer from exact subset} hold; so the quadruple $(G, \lambda, p, k)$ has generic stabilizer $C/Z(G) \cong {A_1}^3.S_3$, or ${A_1}^4.{\Z_2}^2$, or ${A_1}^{\frac{1}{2}(\ell + 1)}$, according as $\ell = 5$, or $\ell = 7$, or $\ell \geq 9$.
\end{proof}

\begin{prop}\label{prop: A_ell, ell even, omega_2 module, k = 2}
Let $G = A_\ell$ for even $\ell \in [4, \infty)$ and $\lambda = \omega_2$, and take $k = 2$. Then the quadruple $(G, \lambda, p, k)$ has generic stabilizer $A_1 T_1 U_\ell$.
\end{prop}

\begin{proof}
We take $G = \SL_{\ell + 1}(K)$ and again use the set-up of Proposition~\ref{prop: A_ell, 2 omega_1 and omega_2 modules}: we write $\ell = 2\ell_1$, so that $\ell_1 \geq 2$; we identify $W$ with the symmetric group $S_{2\ell_1 + 1}$; for $1 \leq i < j \leq 2\ell_1 + 1$ we write $\bar v_{i, j} = v_i \wedge v_j$, where $v_1, \dots, v_{2\ell_1 + 1}$ is the standard basis of $V_{nat}$; we take the generalized height function on the weight lattice of $G$ whose value at each simple root $\alpha_i$ is $2$; for $1 \leq i < j \leq 2\ell_1 + 1$ we let $\nu_{i, j}$ be the weight such that $\bar v_{i, j} \in V_{\nu_{i, j}}$. We then have $\Lambda(V)_{[0]} = \{ \nu_{1, 2\ell_1 + 1}, \nu_{2, 2\ell_1}, \dots, \nu_{\ell_1, \ell_1 + 2} \}$, $\Lambda(V)_{[2]} = \{ \nu_{1, 2\ell_1}, \nu_{2, 2\ell_1 - 1}, \dots, \nu_{\ell_1, \ell_1 + 1} \}$ and $\Lambda(V)_{[+]} = \{ \nu_{i, j} \in \Lambda(V) : i + j \leq 2\ell_1 + 2 \}$. Write
$$
v_{[0]} = \bar v_{1, 2\ell_1 + 1} + \bar v_{2, 2\ell_1} + \cdots + \bar v_{\ell_1, \ell_1 + 2}, \quad
v_{[2]} = \bar v_{1, 2\ell_1} + \bar v_{2, 2\ell_1 - 1} + \cdots + \bar v_{\ell_1, \ell_1 + 1},
$$
and set
$$
y_0 = \langle v_{[0]}, v_{[2]} \rangle,
$$
so that the set of weights occurring in $y_0$ is $\Lambda(V)_{[0]} \cup \Lambda(V)_{[2]}$. We have $Z(G) = \langle z \rangle$ where $z = \eta_{2\ell_1 + 1}I$. Let $P$ be the maximal parabolic subgroup of $G$ corresponding to the $\ell_1$th simple root $\alpha_{\ell_1}$, and write $P = QL$ where $Q$ and $L$ are the unipotent radical and Levi subgroup of $P$ respectively; then $Q$ is abelian.

Take $g \in C_G(y_0)$ and write $g = u_1nu_2$ with $u_1 \in U$, $n \in N$ and $u_2 \in U_w$ where $w = nT \in W$. We have ${u_1}^{-1}.y_0 = n.(u_2.y_0)$; all weights occurring in ${u_1}^{-1}.y_0$ lie in $\Lambda(V)_{[0]} \cup \Lambda(V)_{[+]}$, and we may write $u_2.v_{[0]} = v_{[0]} + {v_{[2]}}' + v'$ and $u_2.v_{[2]} = v_{[2]} + v''$ where the weights in ${v_{[2]}}'$ lie in $\Lambda(V)_{[2]}$, and those in $v'$ or $v''$ lie in $\bigcup_{i > 2} \Lambda(V)_{[i]}$. Thus $w$ cannot send any weight in $\Lambda(V)_{[0]} \cup \Lambda(V)_{[2]}$ into $\Lambda(V)_{[-]}$.

Observe that in the $(2\ell_1 + 1)$-tuple
$$
(2\ell_1 + 1, 1, 2\ell_1, 2, \dots, \ell_1 + 2, \ell_1, \ell_1 + 1),
$$
the adjacent pairs of entries in positions $(2i - 1, 2i)$ for $i = 1, \dots, \ell_1$ sum to $2\ell_1 + 2$ and so correspond to the weights in $\Lambda(V)_{[0]}$, while the adjacent pairs of entries in positions $(2i, 2i + 1)$ for $i = 1, \dots, \ell_1$ sum to $2\ell_1 + 1$ and so correspond to the weights in $\Lambda(V)_{[2]}$. Thus if we apply $w$ to obtain the $(2\ell_1 + 1)$-tuple
$$
(w(2\ell_1 + 1), w(1), w(2\ell_1), w(2), \dots, w(\ell_1 + 2), w(\ell_1), w(\ell_1 + 1)),
$$
each adjacent pair of entries must sum to at most $2\ell_1 + 2$. Therefore in this $(2\ell_1 + 1)$-tuple, the entry $2\ell_1 + 1$ must be placed at either the extreme left or the extreme right, and the entry $1$ must be placed adjacent to it; then in the remaining $(2\ell_1 - 1)$-tuple, the entry $2\ell_1$ must be placed at either the extreme left or the extreme right, and the entry $2$ must be placed adjacent to it; then in the remaining $(2\ell_1 - 3)$-tuple, the entry $2\ell_1 - 1$ must be placed at either the extreme left or the extreme right, and the entry $3$ must be placed adjacent to it; continuing thus we see that at each of $\ell_1$ stages a choice of \lq left' or \lq right' must be made to determine the placing of $2$ entries, after which the element $w$ is determined.

Now at the $j$th stage, if the choice is \lq left' then the entries $2\ell_1 + 2 - j$ and $j$ are placed in positions $(2i - 1, 2i)$ for some $i$, while if the choice is \lq right' then the entries $j$ and $2\ell_1 + 2 - j$ are placed in positions $(2i, 2i + 1)$ for some $i$. Thus according as the choice is \lq left' or \lq right', the term $\bar v_{j, 2\ell_1 + 2 - j}$ occurs in $n.v_{[0]}$ or $n.v_{[2]}$, and hence is either absent from or present in $nu_2.v_{[2]}$. By assumption there exist $a, b \in K$ with $nu_2.v_{[2]} = a{u_1}^{-1}.v_{[0]} + b{u_1}^{-1}.v_{[2]}$; according as $a = 0$ or $a \neq 0$, the vector $nu_2.v_{[2]}$ contains either no such terms $\bar v_{j, 2\ell_1 + 2 - j}$ or all such terms, so that either all choices are \lq left' or all choices are \lq right', producing either the original $(2\ell_1 + 1)$-tuple or its reverse. Thus we must have $w \in \{ 1, w^* \}$, where $w^* = (1\ \ell_1)(2\ \ell_1 - 1) \dots (\ell_1 + 1\ 2\ell_1 + 1)(\ell_1 + 2\ 2\ell_1) \dots$. Note that $w^*$ is the long word in the Weyl group of $L$; moreover there is a corresponding element $n^* \in N$ such that either $n^*$ or $-n^*$ is a permutation matrix, and $n^*$ interchanges $v_{[0]}$ and $v_{[2]}$. Of course if $w = 1$ then $u_2 = 1$. If instead $w = w^*$ then as $u_2 \in U_w$ we see that all weights in either $n.v'$ or $n.v''$ must lie in $\Lambda_{[-]}$, which forces $v' = v'' = 0$; also if ${v_{[2]}}' \notin \langle v_{[2]} \rangle$ then some linear combination of $u_2.v_{[0]}$ and $u_2.v_{[2]}$ would contain some but not all of the weights in $\Lambda(V)_{[2]}$, and then applying $n$ would give a vector containing some but not all of the weights in $\Lambda(V)_{[0]}$, which then could not be of the form $a{u_1}^{-1}.v_{[0]} + b{u_1}^{-1}.v_{[2]}$ for some $a, b \in K$ --- so we must have $u_2.v_{[0]} = v_{[0]} + cv_{[2]}$ for some $c \in K$ and $u_2.v_{[2]} = v_{[2]}$, whence $u_2 \in C_G(y_0)$ (and as $u_2 \in U_w$ we must have $u_2 \in C_{L \cap U}(y_0)$).

We now consider the possibilities for $g$ if $w = 1$; in this case $g \in B$. Clearly $C_B(y_0) = C_T(y_0) C_U(y_0)$, so we may consider each of $C_T(y_0)$ and $C_U(y_0)$ separately.

We begin with $C_T(y_0)$. Suppose $s = \diag(\kappa_1, \dots, \kappa_{2\ell_1 + 1}) \in C_T(y_0)$; then we must have $s.v_{[0]} \in \langle v_{[0]} \rangle$ and $s.v_{[2]} \in \langle v_{[2]} \rangle$, which forces $\kappa_1 \kappa_{2\ell_1 + 1} = \kappa_2 \kappa_{2\ell_1} = \cdots = \kappa_{\ell_1} \kappa_{\ell_1 + 2}$ and $\kappa_1 \kappa_{2\ell_1} = \kappa_2 \kappa_{2\ell_1 - 1} = \cdots = \kappa_{\ell_1} \kappa_{\ell_1 + 1}$. If we take $\kappa' \in K^*$ satisfying ${\kappa'}^2 = \frac{\kappa_{2\ell_1}}{\kappa_{2\ell_1 + 1}}$ and set $\kappa = \kappa' \kappa_1 \kappa_{2\ell_1 + 1}$, then we have $\kappa_1 \kappa_{2\ell_1 + 1} = \cdots = \kappa_{\ell_1} \kappa_{\ell_1 + 2} = \kappa {\kappa'}^{-1}$ and $\kappa_1 \kappa_{2\ell_1} = \cdots = \kappa_{\ell_1} \kappa_{\ell_1 + 1} = \kappa \kappa'$; solving in terms of $\kappa_1$ gives $\kappa_{2\ell_1 + 1} = \kappa {\kappa'}^{-1} {\kappa_1}^{-1}$, $\kappa_{2\ell_1} = \kappa \kappa' {\kappa_1}^{-1}$, $\kappa_2 = {\kappa'}^{-2} \kappa_1$, $\kappa_{2\ell_1 - 1} = \kappa {\kappa'}^3 {\kappa_1}^{-1}$, $\kappa_3 = {\kappa'}^{-4} \kappa_1$, \dots, $\kappa_{\ell_1} = {\kappa'}^{2 - 2\ell_1} \kappa_1$, $\kappa_{\ell_1  + 1} = \kappa {\kappa'}^{2\ell_1 - 1} {\kappa_1}^{-1}$. Imposing the condition $\kappa_1 \kappa_2 \dots \kappa_{2\ell_1 + 1} = 1$ then gives $\kappa_1 = \kappa^{\ell_1 + 1} {\kappa'}^{\ell_1 - 1}$, whence $s = s_1 s_2$ where
\begin{eqnarray*}
s_1 & = & \diag(\kappa^{\ell_1 + 1}, \dots, \kappa^{\ell_1 + 1}, \kappa^{-\ell_1}, \dots, \kappa^{-\ell_1}), \\
s_2 & = & \diag({\kappa'}^{\ell_1 - 1}, {\kappa'}^{\ell_1 - 3}, \dots, {\kappa'}^{3 - \ell_1}, {\kappa'}^{1 - \ell_1}, {\kappa'}^{\ell_1}, {\kappa'}^{\ell_1 - 2}, \dots, {\kappa'}^{2 - \ell_1}, {\kappa'}^{-\ell_1}).
\end{eqnarray*}
Observe that the torus $T_1 = \{ \diag(\kappa^{\ell_1 + 1}, \dots, \kappa^{\ell_1 + 1}, \kappa^{-\ell_1}, \dots, \kappa^{-\ell_1}) : \kappa \in K^* \}$ is the central torus of $L$, and contains $Z(G)$.

We now turn to $C_U(y_0)$. We have $U = Q(L \cap U)$ and $Q \cap (L \cap U) = \{ 1 \}$. Each term in either $v_{[0]}$ or $v_{[2]}$ is of the form $\bar v_{i, j}$ where $1 \leq i \leq \ell_1 < j \leq 2\ell_1 + 1$; given such a term $\bar v_{i, j}$, if $q \in Q$ then $q.\bar v_{i, j} - \bar v_{i, j} \in \langle \bar v_{i', j'} : 1 \leq i' < j' \leq \ell_1 \rangle$, while if $u' \in L \cap U$ then $u'.\bar v_{i, j} - \bar v_{i, j} \in \langle \bar v_{i', j'} : 1 \leq i' \leq \ell_1 < j' \leq 2\ell_1 + 1 \rangle$. Hence $C_U(y_0) = C_Q(y_0) C_{L \cap U}(y_0) = C_{L \cap U}(y_0) C_Q(y_0)$. Thus we may consider each of $C_Q(y_0)$ and $C_{L \cap U}(y_0)$ separately.

We first consider $C_Q(y_0)$. If we take $q = I + \sum_{i = 1}^{\ell_1} \sum_{j = \ell_1 + 1}^{2\ell_1 + 1} t_{ij} E_{ij} \in C_Q(y_0)$, we must have $q.v_{[0]} = v_{[0]}$ and $q.v_{[2]} = v_{[2]}$; taking $1 \leq i < i' \leq \ell_1$ and equating coefficients of $\bar v_{i, i'}$ in these two equations gives $t_{i, 2\ell_1 + 2 - i'} = t_{i', 2\ell_1 + 2 - i}$ and $t_{i, 2\ell_1 + 1 - i'} = t_{i', 2\ell_1 + 1 - i}$ respectively. Hence $t_{i, 2\ell_1 + 2 - i} = t_{i - 1, 2\ell_1 + 1 - i} = t_{i + 1, 2\ell_1 + 3 - i} = t_{i - 2, 2\ell_1 - i} = \cdots$ and $t_{i, 2\ell_1 + 1 - i} =  t_{i + 1, 2\ell_1 + 2 - i} = t_{i - 1, 2\ell_1 - i} = t_{i + 2, 2\ell_1 + 3 - i} = \cdots$, so that $t_{i, j} = t_{i', j'}$ whenever $j - i = j' - i'$. Since $j - i \in \{1, \dots, 2\ell_1 \}$ it follows that $C_Q(y_0)$ is a connected $2\ell_1$-dimensional abelian unipotent group.

Now consider $C_{L \cap U}(y_0)$. Given $t \in K$, write
$$
x(t) = I + \sum_{i = 1}^{\ell_1 - 1} \sum_{j = i + 1}^{\ell_1} \binom{\ell_1 - i}{j - i} (-t)^{j - i} E_{ij} + \sum_{i = \ell_1 + 1}^{2\ell_1} \sum_{j = i + 1}^{2\ell_1 + 1} \binom{j - \ell_1 - 1}{j - i} t^{j - i} E_{ij},
$$
where $E_{ij}$ is the matrix unit with $(i, j)$-entry $1$ and all other entries $0$. A straightforward check shows that for $t, t' \in K$ we have $x(t) x(t') = x(t + t')$. Then
\begin{eqnarray*}
x(t).v_{[0]} & = & \sum_{j = 1}^{\ell_1 + 1} \left( \left( \sum_{i = 1}^j \binom{\ell_1 - i}{j - i} (-t)^{j - i} v_i \right) \right. \wedge \\
             &   & \phantom{\sum_{j = 1}^{\ell_1 + 1} \left( \right.} \left. \left( \sum_{i' = \ell_1 + 1}^{2\ell_1 + 2 - j} \binom{\ell_1 + 1 - j}{2\ell_1 + 2 - j - i'} t^{2\ell_1 + 2 - i' - j} v_{i'} \right) \right)
\end{eqnarray*}
(the term with $j = \ell_1 + 1$ is zero, but is included for convenience in what follows). Take $i, i'$ with $1 \leq i \leq \ell_1 < i' \leq 2\ell_1 + 1$ and $i + i' \leq 2\ell_1 + 2$. The coefficient of $\bar v_{i, i'} = v_i \wedge v_{i'}$ in $x(t).v_{[0]}$ comes from the terms with $i \leq j \leq 2\ell_1 + 2 - i'$, so is $\left( \sum_{j = i}^{2\ell_1 + 2 - i'} (-1)^{j - i} \binom{\ell_1 - i}{j - i} \binom{\ell_1 + 1 - j}{i' - \ell_1 - 1} \right) t^{2\ell_1 + 2 - i - i'}$. We see that this expression is the $t^{2\ell_1 + 2 - i - i'}$ term in $(-1 - t)^{2\ell_1 - i - i'} = (-1)^{i + i'} (1 + t)^{\ell_1 - i} (1 + t)^{\ell_1 - i'}$, so it is $1$, $t$ or $0$ according as $i + i'$ is $2\ell_1 + 2$, $2\ell_1 + 1$ or at most $2\ell_1$; thus $x(t).v_{[0]} = v_{[0]} + tv_{[2]}$. Likewise
\begin{eqnarray*}
x(t).v_{[2]} & = & \sum_{j = 1}^{\ell_1} \left( \left( \sum_{i = 1}^j \binom{\ell_1 - i}{j - i} (-t)^{j - i} v_i \right) \right. \wedge \\
             &   & \phantom{\sum_{j = 1}^{\ell_1} \left( \right.} \left. \left( \sum_{i' = \ell_1 + 1}^{2\ell_1 + 1 - j} \binom{\ell_1 - j}{2\ell_1 + 1 - j - i'} t^{2\ell_1 + 1 - i' - j} v_{i'} \right) \right).
\end{eqnarray*}
Take $i, i'$ with $1 \leq i \leq \ell_1 < i' \leq 2\ell_1 + 1$ and $i + i' \leq 2\ell_1 + 1$. The coefficient of $\bar v_{i, i'}$ in $x(t).v_{[2]}$ comes from the terms with $i \leq j \leq 2\ell_1 + 1 - i'$, so is $\left( \sum_{j = i}^{2\ell_1 + 1 - i'} (-1)^{j - i} \binom{\ell_1 - i}{j - i} \binom{\ell_1 - j}{i' - \ell_1 - 1} \right) t^{2\ell_1 + 1 - i - i'}$. This time we see that this expression is the $t^{2\ell_1 + 1 - i - i'}$ term in $(-1 - t)^{2\ell_1 - i - i'}$, so it is $1$ or $0$ according as $i + i'$ is $2\ell_1 + 1$ or at most $2\ell_1$; thus $x(t).v_{[2]} = v_{[2]}$.

Therefore $x(t) \in C_{L \cap U}(y_0)$. Moreover, as any element of $C_{L \cap U}(y_0)$ must fix $v_{[2]}$ and send $v_{[0]}$ to $v_{[0]} + t v_{[2]}$ for some $t \in K$, it now suffices to consider $u \in L \cap U$ fixing both $v_{[0]}$ and $v_{[2]}$. Write $u = I + \sum_{i = 1}^{\ell_1 - 1} \sum_{j = i + 1}^{\ell_1} t_{ij} E_{ij} + \sum_{i = \ell_1 + 1}^{2\ell_1} \sum_{j = i + 1}^{2\ell_1 + 1} t_{ij} E_{ij}$. Considering the coefficients of $\bar v_{\ell_1, \ell_1 + 1}$, $\bar v_{\ell_1 - 1, \ell_1 + 2}$, \dots, $\bar v_{1, 2\ell_1}$ in $u.v_{[0]}$, and those of $\bar v_{\ell_1 - 1, \ell_1 + 1}$, $\bar v_{\ell_1 - 2, \ell_1 + 2}$, \dots, $\bar v_{1, 2\ell_1 - 1}$ in $u.v_{[2]}$, gives $0 = t_{\ell_1 + 1, \ell_1 + 2} = t_{\ell_1 - 1, \ell_1} + t_{\ell_1 + 2, \ell_1 + 3} = t_{\ell_1 - 2, \ell_1 - 1} + t_{\ell_1 + 3, \ell_1 + 4} = \cdots = t_{1, 2} + t_{2\ell_1, 2\ell_1 + 1}$ and $0 = t_{\ell_1 - 1, \ell_1} + t_{\ell_1 + 1, \ell_1 + 2} = t_{\ell_1 - 2, \ell_1 - 1} + t_{\ell_1 + 2, \ell_1 + 3} = \cdots = t_{1, 2} + t_{2\ell_1 - 1, 2\ell_1}$, whence $t_{i, j} = 0$ whenever $j - i = 1$; then considering the coefficients of $\bar v_{\ell_1 - 1, \ell_1 + 1}$, $\bar v_{\ell_1 - 2, \ell_1 + 2}$, \dots, $\bar v_{1, 2\ell_1 - 1}$ in $u.v_{[0]}$, and those of $\bar v_{\ell_1 - 2, \ell_1 + 1}$, $\bar v_{\ell_1 - 3, \ell_1 + 2}$, \dots, $\bar v_{1, 2\ell_1 - 2}$ in $u.v_{[2]}$, gives $0 = t_{\ell_1 + 1, \ell_1 + 3} = t_{\ell_1 - 2, \ell_1} + t_{\ell_1 + 2, \ell_1 + 4} = \cdots = t_{1, 3} + t_{2\ell_1 - 1, 2\ell_1 + 1}$ and $0 = t_{\ell_1 - 2, \ell_1} + t_{\ell_1 + 1, \ell_1 + 3} = t_{\ell_1 - 3, \ell_1 - 1} + t_{\ell_1 + 2, \ell_1 + 4} = \cdots = t_{1, 3} + t_{2\ell_1 - 2, 2\ell_1}$, whence $t_{i, j} = 0$ whenever $j - i = 2$; continuing in this way we see that all $t_{i, j}$ are zero, so $u = 1$. Hence $C_{L \cap U}(y_0) = \{ x(t) : t \in K \}$.

Thus if $w = 1$ we see that $g = u_1s$ with $u_1 = x(t) u'$ with $u' \in C_Q(y_0)$, $t \in K$ and $s$ in the $2$-dimensional torus described above. If instead $w = w^*$ then $g = u_1sn^*u_2$ where $u_1$ and $s$ are as in the case $w = 1$, and $u_2 \in C_{L \cap U}(y_0)$ so $u_2 = x(t')$ for some $t' \in K$. Hence as the torus of the $A_1$ subgroup $\langle x(t), n^* : t \in K \rangle$ is $\{ \diag(\kappa^{\ell_1 - 1}, \kappa^{\ell_1 - 3}, \dots, \kappa^{3 - \ell_1}, \kappa^{1 - \ell_1}, \kappa^{\ell_1}, \kappa^{\ell_1 - 2}, \dots, \kappa^{2 - \ell_1}, \kappa^{-\ell_1}) : \kappa \in K^* \}$, we see that
$$
C_G(y_0) = \langle x(t), n^* : t \in K \rangle T_1 C_Q(y_0) \cong A_1 T_1 U_{2\ell_1}.
$$
Thus $\dim(\overline{G.y_0}) = \dim G - \dim C_G(y_0) = (4{\ell_1}^2 + 4\ell_1) - (2\ell_1 + 4) = 4{\ell_1}^2 + 2\ell_1 - 4 = \dim \G{2}(V)$, so the orbit $G.y_0$ is dense in $\G{2}(V)$. Hence the quadruple $(G, \lambda, p, k)$ has generic stabilizer $C_G(y_0)/Z(G) \cong A_1 T_1 U_\ell$.
\end{proof}

This concludes the treatment of the cases occurring in infinite families. As was the case in Section~\ref{sect: small triples and first quadruples}, although the remaining cases must be treated individually, it will be seen that there are connections between some of them which significantly reduce the amount of work involved.

\begin{prop}\label{prop: A_4, omega_2 module, k = 3}
Let $G = A_4$ and $\lambda = \omega_2$, and take $k = 3$. Then the quadruple $(G, \lambda, p, k)$ has generic stabilizer $A_1$.
\end{prop}

\begin{proof}
We take $G = \SL_5(K)$ and again use the set-up of Proposition~\ref{prop: A_ell, 2 omega_1 and omega_2 modules}: we identify $W$ with the symmetric group $S_5$; for $1 \leq i < j \leq 5$ we write $\bar v_{i, j} = v_i \wedge v_j$, where $v_1, \dots, v_5$ is the standard basis of $V_{nat}$; we take the generalized height function on the weight lattice of $G$ whose value at each simple root $\alpha_i$ is $2$; for $1 \leq i < j \leq 5$ we let $\nu_{i, j}$ be the weight such that $\bar v_{i, j} \in V_{\nu_{i, j}}$. We then have
$$
\Lambda(V)_{[-2]} = \{ \nu_{2, 5}, \nu_{3, 4} \}, \quad \Lambda(V)_{[0]} = \{ \nu_{1, 5}, \nu_{2, 4} \}, \quad \Lambda(V)_{[2]} = \{ \nu_{1, 4}, \nu_{2, 3} \}.
$$
Set
$$
v_{[-2]} = \bar v_{2, 5} + \bar v_{3, 4}, \quad v_{[0]} = \bar v_{1, 5} + \bar v_{2, 4}, \quad v_{[2]} = \bar v_{1, 4} + \bar v_{2, 3},
$$
and write
$$
y_0 = \langle v_{[-2]}, v_{[0]}, v_{[2]} \rangle,
$$
so that the set of weights occurring in $y_0$ is $\Lambda(V)_{[-2]} \cup \Lambda(V)_{[0]} \cup \Lambda(V)_{[2]}$. We have $Z(G) = \langle z \rangle$ where $z = \eta_5 I$.

Take $g \in C_G(y_0)$ and write $g = u_1nu_2$ with $u_1 \in U$, $n \in N$ and $u_2 \in U_w$ where $w = nT \in W$. We have ${u_1}^{-1}.y_0 = n.(u_2.y_0)$; all weights occurring in ${u_1}^{-1}.y_0$ lie in $\bigcup_{i \geq -2} \Lambda(V)_{[i]}$, and we may write $u_2.v_{[-2]} = v_{[-2]} + {v_{[0]}}' + {v_{[2]}}' + v'$, $u_2.v_{[0]} = v_{[0]} + {v_{[2]}}'' + v''$ and $u_2.v_{[2]} = v_{[2]} + v'''$ where ${v_{[0]}}' \in V_{[0]}$, ${v_{[2]}}', {v_{[2]}}'' \in V_{[2]}$ and $v', v'', v''' \in \bigcup_{i > 2} V_{[i]}$. Thus $w$ cannot send any weight in $\Lambda(V)_{[-2]} \cup \Lambda(V)_{[0]} \cup \Lambda(V)_{[2]}$ into $\bigcup_{i < -2} \Lambda(V)_{[i]} = \{ \nu_{4, 5}, \nu_{3, 5} \}$; therefore $w^{-1}$ must send both $\nu_{4, 5}$ and $\nu_{3, 5}$ into $\bigcup_{i \neq -2, 0, 2} \Lambda(V)_{[i]} = \{ \nu_{4, 5}, \nu_{3, 5}, \nu_{1, 3}, \nu_{1, 2} \}$. Since the only two pairs of weights in this set whose difference is a root are $\{ \nu_{4, 5}, \nu_{3, 5} \}$ and $\{ \nu_{1, 3}, \nu_{1, 2} \}$, we see that $w^{-1}$ must send $\{ \nu_{4, 5}, \nu_{3, 5} \}$ to either $\{ \nu_{4, 5}, \nu_{3, 5} \}$ or $\{ \nu_{1, 3}, \nu_{1, 2} \}$. Clearly the setwise stabilizer in $W$ of $\{ \nu_{4, 5}, \nu_{3, 5} \}$ is $\langle (1\ 2), (3\ 4) \rangle$; since $w_0 = (1\ 5)(2\ 4)$ interchanges $\{ \nu_{4, 5}, \nu_{3, 5} \}$ and $\{ \nu_{1, 3}, \nu_{1, 2} \}$, the elements of $W$ which send $\{ \nu_{4, 5}, \nu_{3, 5} \}$ to $\{ \nu_{1, 3}, \nu_{1, 2} \}$ are those in $w_0 \langle (1\ 2), (3\ 4) \rangle$. Hence $w^{-1} \in \{ 1, w_0 \}\langle (1\ 2), (3\ 4) \rangle$, so $w \in \langle (1\ 2), (3\ 4) \rangle \{ 1, w_0 \} = \{ 1, (1\ 2), (3\ 4), (1\ 2)(3\ 4), (1\ 5)(2\ 4), (1\ 5\ 2\ 4), (1\ 5)(2\ 3\ 4), (1\ 5\ 2\ 3\ 4) \}$. Arguing similarly with $g^{-1}$ we see that $w^{-1}$ must also lie in this set, so that $w \in \{ 1, (1\ 2), (3\ 4), (1\ 2)(3\ 4), (1\ 5)(2\ 4) \}$. However, if $w = (1\ 2)$ or $(3\ 4)$ then $nu_2.v_{[2]}$ contains a term $\bar v_{2, 4}$ but no term $\bar v_{1, 5}$, so cannot lie in ${u_1}^{-1}.y_0$; likewise if $w = (1\ 2)(3\ 4)$ then $nu_2.v_{[0]}$ contains a term $\bar v_{2, 5}$ but no term $\bar v_{3, 4}$, so cannot lie in ${u_1}^{-1}.y_0$. Thus $w \in \{ 1, (1\ 5)(2\ 4) \} = \langle w_0 \rangle$.

First suppose $w = 1$; then $u_2 = 1$ and $g = u_1h$ where $h \in T$, and we must have $u_1 \in C_U(y_0)$ and $h \in C_T(y_0)$. Equating to zero the coefficients of $\bar v_{1, 3}$ and $\bar v_{1, 2}$ in each of $u_1.v_{[-2]}$, $u_1.v_{[0]}$ and $u_1.v_{[2]}$, and requiring equality in the coefficients of $\bar v_{1, 5}$ and $\bar v_{2, 4}$ in $u_1.v_{[-2]}$, and in the coefficients of $\bar v_{1, 4}$ and $\bar v_{2, 3}$ in both $u_1.v_{[-2]}$ and $u_1.v_{[0]}$, shows that we must have $u_1 = x(t)$ for some $t \in K$, where we write $x(t)$ for the matrix
$$
\left(
  \begin{array}{ccccc}
    1 &   & t &   & t^2 \\
      & 1 &   & t &   \\
      &   & 1 &   &   \\
      &   &   & 1 &   \\
      &   &   &   & 1 \\
  \end{array}
\right)
\qquad \hbox{or}\qquad
\left(
  \begin{array}{ccccc}
    1 & t & \frac{3}{2}t^2 & -\frac{1}{2}t^3 & \frac{1}{4}t^4 \\
      & 1 &       3t       & -\frac{3}{2}t^2 &       t^3      \\
      &   &        1       &        -t       &       t^2      \\
      &   &                &         1       &       -2t      \\
      &   &                &                 &        1       \\
  \end{array}
\right)
$$
according as $p = 2$ or $p \geq 3$. A straightforward check shows that for $t, t' \in K$ we have $x(t) x(t') = x(t + t')$. If we write $h = \diag(\kappa_1, \dots, \kappa_5)$ then we require $\kappa_2\kappa_5 = \kappa_3\kappa_4$, $\kappa_1\kappa_5 = \kappa_2\kappa_4$ and $\kappa_1\kappa_4 = \kappa_2\kappa_3$; since in addition $\kappa_1 \kappa_2 \kappa_3 \kappa_4 \kappa_5 = 1$ there exist $\kappa \in K^*$ and $i \in \{ 0, 1, 2, 3, 4 \}$ with $h = h(\kappa) z^i$, where we write $h(\kappa) = \diag(\kappa^2, \kappa, 1, \kappa^{-1}, \kappa^{-2})$. Thus $g = x(t) h(\kappa) z^i$.

Now suppose $w = w_0$; we may write $n = h n^*$ where $h \in T$ and $n^*$ is the permutation matrix whose $(i, j)$-entry is $1$ if $i + j = 6$ and $0$ otherwise. Since $n^*$ sends $v_{[i]}$ to $-v_{[-i]}$ for $i = -2, 0, 2$ we see that $n^* \in C_G(y_0)$. Now we have $h^{-1}{u_1}^{-1}.y_0 = n^*u_2.y_0$. As $n^*$ sends $\bar v_{1, 2}$ and $\bar v_{1, 3}$ to $-\bar v_{4, 5}$ and $-\bar v_{3, 5}$ respectively, neither of which appears in any vector in $h^{-1}{u_1}^{-1}.y_0$, we see that the coefficients of $\bar v_{1, 2}$ and $\bar v_{1, 3}$ in each of $u_2.v_{[-2]}$, $u_2.v_{[0]}$ and $u_2.v_{[2]}$ must be zero. Moreover if the coefficients of $\bar v_{1, 4}$ and $\bar v_{2, 3}$ in $u_2.v_{[0]}$ were unequal, we could subtract a multiple of $u_2.v_{[2]} = v_{[2]}$ to give a vector $v_{[0]} + \kappa\bar v_{1, 4}$ for some $\kappa \in K^*$, whose image under $n^*$ would be $-v_{[0]} - \kappa\bar v_{2, 5}$, which cannot lie in $h^{-1}{u_1}^{-1}.y_0$; so we must have $u_2.v_{[0]} - v_{[0]} \in \langle v_{[2]} \rangle$. Similarly we must have $u_2.v_{[-2]} - v_{[-2]} \in \langle v_{[0]}, v_{[2]} \rangle$, so that $u_2 \in C_U(y_0)$. Since now $u_1h \in C_G(y_0)$ we see that $g = x(t) h(\kappa) z^i n^* x(t')$ for some $t, t' \in K$, some $\kappa \in K^*$ and some $i \in \{ 0, 1, 2, 3, 4 \}$.

Write $A = \langle x(t), n^* : t \in K \rangle$. It is straightforward to see that $A \cong A_1$, and that $A \cap T = \{ h(\kappa) : \kappa \in K^* \}$; thus we have $C_G(y_0) = Z(G) A$. Since $\dim(\overline{G.y_0}) = \dim G - \dim C_G(y_0) = 24 - 3 = 21 = \dim \G{3}(V)$, the orbit $G.y_0$ is dense in $\G{3}(V)$. Therefore the quadruple $(G, \lambda, p, k)$ has generic stabilizer $C_G(y_0)/Z(G) \cong A_1$.
\end{proof}

\begin{prop}\label{prop: A_4, omega_2 module, k = 4}
Let $G = A_4$ and $\lambda = \omega_2$, and take $k = 4$. Then the quadruple $(G, \lambda, p, k)$ has generic stabilizer $S_5$.
\end{prop}

\begin{proof}
We take $G = \SL_5(K)$; we have $Z(G) = \langle z \rangle$ where $z = \eta_5 I$. Recall that $V_{nat}$ has basis $v_1, \dots, v_5$. We have
$$
V = {\ts\bigwedge^2}(V_{nat}) = \langle v_{12}, v_{23}, v_{34}, v_{45}, v_{51}, v_{13}, v_{24}, v_{35}, v_{41}, v_{52} \rangle,
$$
where for $i, j \leq 5$ we write $v_{ij} = v_i \wedge v_j$. Take
$$
y_0 = \langle v_{12} - v_{23}, v_{23} - v_{34}, v_{34} - v_{45}, v_{45} - v_{51} \rangle.
$$

Define
$$
g_1 =
\left(
  \begin{array}{ccccc}
      &   &   &   & 1 \\
    1 &   &   &   &   \\
      & 1 &   &   &   \\
      &   & 1 &   &   \\
      &   &   & 1 &   \\
  \end{array}
\right),
\quad
g_2 =
\left(
  \begin{array}{ccccc}
   1 &   &    &        &    \\
     & 1 &    &        &    \\
   1 &   & -1 &   -1   &    \\
     &   &    & \pmin1 &    \\
     & 1 &    &   -1   & -1 \\
  \end{array}
\right).
$$
Since ${g_2}^2 = (g_2g_1)^4 = {g_1}^5 = ({g_2}^{-1}{g_1}^{-1}g_2g_1)^3 = 1$, we see from \cite{Atlas} that $\langle g_1, g_2 \rangle$ has presentation $(2, 4, 5; 3)$ (in the notation there); thus $\langle g_1, g_2 \rangle \cong S_5$, with $g_1$ and $g_2$ corresponding to the permutations $(1 \ 2 \ 3 \ 4 \ 5)$ and $(1 \ 2)$ respectively. We find that $g_1$ and $g_2$ both lie in $C_G(y_0)$; indeed with respect to the given basis they act on $y_0$ as
$$
\left(
  \begin{array}{cccc}
      &   &   & -1 \\
    1 &   &   & -1 \\
      & 1 &   & -1 \\
      &   & 1 & -1 \\
  \end{array}
\right)
\quad \hbox{and}\quad
\left(
  \begin{array}{cccc}
    1 &    &    &    \\
    1 & -1 &    &    \\
    1 &    & -1 &    \\
    1 &    &    & -1 \\
  \end{array}
\right)
$$
respectively. Set $C = Z(G)\langle g_1, g_2 \rangle$; then $C \leq C_G(y_0)$, and we claim that in fact $C_G(y_0) = C$.

Take $g = (a_{ij}) \in C_G(y)$. In what follows, we shall regard rows and columns as labelled with the elements of the finite field $\F_5$ of size $5$; rows or columns labelled with $i$ and $j$ will be called {\em adjacent\/} if $i - j = \pm1$ and {\em non-adjacent\/} otherwise, as will entries within a given row or column.

\pagebreak

Take the basis vector $v = v_{12} - v_{23}$ of $y_0$; for $i = 1, \dots, 5$ temporarily write $c_i = a_{i2}$ and $b_i = a_{i1} + a_{i3}$, so that the column vectors $(c_i)$ and $(b_i)$ are the second column of $g$ and the sum of the first and third columns of $g$ respectively. The coefficient of $v_{13}$ in $g.v$ is $(a_{11}a_{32} - a_{31}a_{12}) - (a_{12}a_{33} - a_{32}a_{13})$; since this must be zero, we have $a_{32}(a_{11} + a_{13}) = a_{12}(a_{31} + a_{33})$, and treating likewise the coefficients in $g.v$ of $v_{24}$, $v_{35}$, $v_{41}$ and $v_{52}$ we obtain the equations
$$
c_3b_1 = c_1b_3, \quad c_4b_2 = c_2b_4, \quad c_5b_3 = c_3b_5, \quad c_1b_4 = c_4b_1, \quad c_2b_5 = c_5b_2.
$$
Moreover the sum of the coefficients of $v_{12}$, $v_{23}$, $v_{34}$, $v_{45}$ and $v_{51}$ in $g.v$ must be zero, giving the further equation
$$
c_2b_1 + c_3b_2 + c_4b_3 + c_5b_4 + c_1b_5 = c_1b_2 + c_2b_3 + c_3b_4 + c_4b_5 + c_5b_1.
$$
In the arguments which follow, we may simultaneously cycle the entries of the vectors $(c_i)$ and $(b_i)$.

First suppose $c_1c_3c_5 \neq 0$; writing $r = \frac{b_1}{c_1}$, from the first and third equations above we would also have $r = \frac{b_3}{c_3} = \frac{b_5}{c_5}$, and from the fourth and fifth we would have $b_4 = rc_4$ and $b_2 = rc_2$, whence $(b_i) = r(c_i)$, which is impossible as the columns of $g$ are linearly independent. Thus we must have $c_1c_3c_5 = 0$; cycling we see that the vector $(c_i)$ must have at least two adjacent entries equal to zero.

Next suppose $c_1c_5 \neq 0 = c_2 = c_3 = c_4$; from the fifth, first and fourth equations above we would have $b_2 = b_3 = b_4 = 0$, and from the final equation $c_1b_5 = c_5b_1$, so writing $r = \frac{b_1}{c_1}$ we would again have $(b_i) = r(c_i)$, which is impossible. Thus if the vector $(c_i)$ has exactly two non-zero entries, they must be non-adjacent.

Next suppose $c_2c_4 \neq 0$. By the above we must have $c_1 = c_5 = 0$, and then the fourth and fifth equations above give $b_1 = b_5 = 0$; if we set $r = \frac{b_2}{c_2}$ then the second equation gives $r = \frac{b_4}{c_4}$, while the final equation becomes $c_3b_2 + c_4b_3 = c_2b_3 + c_3b_4$, which gives $b_3(c_4 - c_2) = c_3(b_4 - b_2) = rc_3(c_4 - c_2)$. If $c_4 \neq c_2$ we would have $b_3 = rc_3$, so that we would again have $(b_i) = r(c_i)$, which is impossible; so we must have $c_4 = c_2$, and hence $b_4 = b_2$.

Finally suppose $c_3 \neq 0 = c_1 = c_2 = c_4 = c_5$; then the first and third equations give $b_1 = b_5 = 0$, and the final equation gives $b_2 = b_4$.

Therefore after (simultaneous) cycling both column vectors $(c_i)$ and $(b_i)$ are of the form $(0 \ \kappa_1 \ \kappa_2 \ \kappa_1 \ 0)^T$, where $\kappa_1, \kappa_2 \in K$, and either (but not both) may be zero. Replacing $v$ by each of the other three basis vectors of $y_0$, and by the negative of the sum of all four, shows that the same is true whenever $(c_i)$ is a column of $g$ and $(b_i)$ is the sum of the two adjacent columns. We shall say that a column of the form $(0 \ \kappa_1 \ \kappa_2 \ \kappa_1 \ 0)^T$ is {\em centred\/} on the row containing the entry $\kappa_2$; this gives a map $\pi_g$ from the set $\{ 1, \dots, 5 \}$ to itself such that each column $j$ is centred on row $\pi_g(j)$, and we shall write $\pi_g$ as the $5$-tuple $(\pi_g(1), \pi_g(2), \pi_g(3), \pi_g(4), \pi_g(5))$.

Now if for some $j$ we had $\pi_g(j - 1) = \pi_g(j) = i$, it would immediately follow that we must have $\pi_g(j + 1) = i$ (by applying the above to column $j$) and $\pi_g(j - 2) = i$ (from column $j - 1$), and then that $\pi_g(j + 2) = i$ (from column $j + 1$), so that the rows non-adjacent to row $i$ would be zero, which is impossible. If $\pi_g(j - 1) = \pi_g(j + 1) = i$, considering column $j$ again would give $\pi_g(j) = i$. Thus $\pi_g$ must be injective, and hence a permutation of $\{ 1, \dots, 5\}$. We claim that $g$ must then be a scalar multiple of the element of $\langle g_1, g_2 \rangle$ corresponding to the permutation $\pi_g$. In proving this claim we shall make use of the element $g_1$ above corresponding to $(1 \ 2 \ 3 \ 4 \ 5)$, and that corresponding to $(1 \ 5)(2 \ 4)$, which we find to be
$$
\left(
  \begin{array}{ccccc}
      &   &   &   & 1 \\
      &   &   & 1 &   \\
      &   & 1 &   &   \\
      & 1 &   &   &   \\
    1 &   &   &   &   \\
  \end{array}
\right);
$$
since each of these is simply the appropriate permutation matrix, we may pre- and post-multiply $g$ by them without affecting the truth of the claim.

First suppose no two adjacent columns of $g$ are centred on adjacent rows; then using the two elements above we may assume $\pi_g = (4, 1, 3, 5, 2) = (1 \ 4 \ 5 \ 2)$. Since the second and fourth columns are centred on rows $1$ and $5$, for some $\kappa_1, \kappa_2, {\kappa_1}', {\kappa_2}' \in K$ they are $(\kappa_2 \ \kappa_1 \ 0 \ 0 \ \kappa_1)^T$ and $({\kappa_1}' \ 0 \ 0 \ {\kappa_1}' \ {\kappa_2}')^T$ respectively; since their sum is centred on row $3$ we must have ${\kappa_1}' + \kappa_2 = 0 = \kappa_1 + {\kappa_2}'$ and $\kappa_1 = {\kappa_1}'$, so that the second and fourth columns are $(\kappa_2 \ {-\kappa_2} \ 0 \ 0 \ {-\kappa_2})^T$ and $({-\kappa_2} \ 0 \ 0 \ {-\kappa_2} \ \kappa_2)^T$ respectively. Arguing exactly similarly with the other pairs of non-adjacent columns, we see that $g$ is a scalar multiple of
$$
\left(
  \begin{array}{ccccc}
          & \pmin1 &        &   -1   &   -1   \\
          &   -1   &   -1   &        & \pmin1 \\
     -1   &        & \pmin1 &        &   -1   \\
   \pmin1 &        &   -1   &   -1   &        \\
     -1   &   -1   &        & \pmin1 &        \\
  \end{array}
\right),
$$
which is the element of $\langle g_1, g_2 \rangle$ corresponding to $(1 \ 4 \ 5 \ 2)$.

Thus we may suppose $g$ has two adjacent columns centred on adjacent rows; using the two elements above we may assume $\pi_g(2) = 2$ and $\pi_g(3) = 3$, so for some $\kappa_1, \kappa_2, {\kappa_1}', {\kappa_2}' \in K$ the second column and the sum of second and fourth are $(\kappa_1 \ \kappa_2 \ \kappa_1 \ 0 \ 0)^T$ and $(0 \ {\kappa_1}' \ {\kappa_2}' \ {\kappa_1}' \ 0)^T$ respectively, so that the fourth is $(-\kappa_1 \ ({\kappa_1}' - \kappa_2) \ ({\kappa_2}' - \kappa_1) \ {\kappa_1}' \ 0)^T$. If $\pi_g(4) = 4$ we must have $\kappa_1 = {\kappa_1}' - \kappa_2 = 0$ and ${\kappa_2}' - \kappa_1 = 0$, so the second and fourth columns are $(0 \ \kappa_2 \ 0 \ 0 \ 0)^T$ and $(0 \ 0 \ 0 \ \kappa_2 \ 0)^T$ respectively; if instead $\pi_g(4) = 5$ we must have ${\kappa_1}' - \kappa_2 = {\kappa_2}' - \kappa_1 = 0$ and ${\kappa_1}' = -\kappa_1$, so the second and fourth columns are $({-\kappa_2} \ \kappa_2 \ {-\kappa_2} \ 0 \ 0)^T$ and $(\kappa_2 \ 0 \ 0 \ \kappa_2 \ 0)^T$ respectively; finally if $\pi_g(4) = 1$ we must have ${\kappa_2}' - \kappa_1 = {\kappa_1}' = 0$ and ${\kappa_1}' - \kappa_2 = 0$, so the second and fourth columns are $(\kappa_1 \ 0 \ \kappa_1 \ 0 \ 0)^T$ and $({-\kappa_1} \ 0 \ 0 \ 0 \ 0)^T$ respectively.

It is now straightforward to apply the analysis of the previous two paragraphs, using post-multiplication by the two elements above, to complete the consideration of each of these possibilities for $\pi_g$. If $\pi_g = (1, 2, 3, 4, 5) = 1$ we immediately find that $g$ is a scalar multiple of $I$. If $\pi_g = (5, 2, 3, 4, 1) = (1 \ 5)$, from the fourth and first columns we see that the first column must be $({-\kappa_2} \ 0 \ 0 \ {-\kappa_2} \ 0)^T$, then from the first and third that the third must be $(0 \ \kappa_2 \ {-\kappa_2} \ \kappa_2 \ 0)^T$, and finally from the third and fifth that the fifth must be $(0 \ {-\kappa_2} \ 0 \ 0 \ -\kappa_2)^T$; so $g$ is a scalar multiple of
$$
\left(
  \begin{array}{ccccc}
   1 &    &        &    &   \\
     & -1 &   -1   &    & 1 \\
     &    & \pmin1 &    &   \\
   1 &    &   -1   & -1 &   \\
     &    &        &    & 1 \\
  \end{array}
\right),
$$
which is the element of $\langle g_1, g_2 \rangle$ corresponding to $(1 \ 5)$. Likewise according as $\pi_g = (1, 2, 3, 5, 4) = (4 \ 5)$, $(4, 2, 3, 5, 1) = (1 \ 4 \ 5)$, $(4, 2, 3, 1, 5) = (1 \ 4)$ or $(5, 2, 3, 1, 4) = (1 \ 5 \ 4)$ we find that $g$ is a scalar multiple of
\begin{eqnarray*}
\left(
  \begin{array}{ccccc}
   -1 &   -1   &    & 1 &   \\
      & \pmin1 &    &   &   \\
      &   -1   & -1 &   & 1 \\
      &        &    & 1 &   \\
      &        &    &   & 1 \\
  \end{array}
\right),
& \quad &
\left(
  \begin{array}{ccccc}
     & \pmin1 &    & -1 &   -1   \\
     &   -1   & -1 &    & \pmin1 \\
     & \pmin1 &    &    &        \\
   1 &        & -1 & -1 &        \\
     &        &    &    & \pmin1 \\
  \end{array}
\right), \\
\left(
  \begin{array}{ccccc}
      & 1 &   & -1 &   -1   \\
      &   & 1 &    &        \\
      & 1 &   &    &        \\
   -1 &   & 1 &    &   -1   \\
      &   &   &    & \pmin1 \\
  \end{array}
\right),
& \quad &
\left(
  \begin{array}{ccccc}
   -1 & -1 &        & 1 &        \\
      &    & \pmin1 &   &        \\
      & -1 &   -1   &   & \pmin1 \\
   -1 &    & \pmin1 &   &   -1   \\
      &    &        &   & \pmin1 \\
  \end{array}
\right)
\end{eqnarray*}
respectively, which are the elements of $\langle g_1, g_2 \rangle$ corresponding to $(4 \ 5)$, $(1 \ 4 \ 5)$, $(1 \ 4)$ and $(1 \ 5 \ 4)$ respectively. We have therefore proved the claim.

Now the condition $\det g = 1$ forces the scalar involved in $g$ to be a fifth root of unity, giving $g \in Z(G)\langle g_1, g_2 \rangle = C$. Thus we do indeed have $C_G(y_0) = C$.

Since $\dim(\overline{G.y_0}) = \dim G - \dim C_G(y_0) = 24 - 0 = 24 = \dim \G{4}(V)$, the orbit $G.y_0$ is dense in $\G{4}(V)$. Thus the quadruple $(G, \lambda, p, k)$ has generic stabilizer $C_G(y_0)/Z(G) \cong S_5$.
\end{proof}

\begin{prop}\label{prop: B_3, omega_3 module, k = 2 or 3}
Let $G = B_3$ and $\lambda = \omega_3$, and take $k = 2$ or $3$. Then the quadruple $(G, \lambda, p, k)$ has generic stabilizer $A_2T_1.\Z_2$ or ${A_1}^2$ respectively.
\end{prop}

\begin{proof}
As in Proposition~\ref{prop: B_3, omega_3 module}, we take $H$ to be the (simply connected) group defined over $K$ of type $F_4$, with simple roots $\beta_1, \beta_2, \beta_3, \beta_4$; we let $G$ have simple roots $\alpha_i = \beta_i$ for $i \leq 3$, so that $G = \langle X_\alpha : \alpha = \sum m_i \beta_i, \ m_4 = 0 \rangle < H$; then we may take $V = \langle e_\alpha : \alpha = \sum m_i \beta_i, \ m_4 = 1 \rangle < \L(H)$. We have $Z(G) = \langle z \rangle$ where $z = h_{\beta_3}(-1)$.

First take $k = 2$. Here we use the set-up of Proposition~\ref{prop: B_3, omega_3 module}: we take the generalized height function on the weight lattice of $G$ whose value at each simple root $\alpha_i$ is $1$, and then $\Lambda(V)_{[0]} = \{ \nu_1, \nu_2 \}$, where we write
$$
\gamma_1 = \ffourrt1111, \quad \gamma_2 = \ffourrt0121,
$$
and for each $i$ we let $\nu_i$ be the weight such that $V_{\nu_i} = \langle e_{\gamma_i} \rangle$; we have $\nu_1 + \nu_2 = 0$, so $\Lambda(V)_{[0]}$ has ZLC; and the setwise stabilizer in $W$ of $\Lambda(V)_{[0]}$ is $\langle w_{\beta_2}, w_{\beta_1} w_{\beta_3} \rangle$. Here we take $Y = \G{2}(V_{[0]})$ and write
$$
y_0 = \langle e_{\gamma_1}, e_{\gamma_2} \rangle \in Y.
$$
By Lemma~\ref{lem: gen height zero} we have $C_G(y_0) = C_U(y_0) C_{N_{\Lambda(V)_{[0]}}}(y_0) C_U(y_0)$.

Let $A$ be the $A_2$ subgroup having simple roots $\beta_2$ and $\beta_1 + \beta_2 + 2\beta_3$ (so that $A$ is the subgroup generated by the long root subgroups of the $G_2$ subgroup seen in the proof of Proposition~\ref{prop: B_3, omega_3 module}); note that $Z(A) = \langle z \rangle$ where $z = h_{\beta_1}(\eta_3) h_{\beta_3}(\eta_3)$. Write $T_1 = C_T(A) = \{ h_{\alpha_1}(\kappa^{-2}) h_{\alpha_3} (\kappa) : \kappa \in K^* \}$, so that $Z(G) < T_1$; set $n^* = n_{\beta_1} n_{\beta_3}$, and write $C = A T_1 \langle n^* \rangle$. Clearly we then have $C \leq C_G(y_0)$; we shall show that in fact $C_G(y_0) = C$.

First, from the above the elements of $W$ which preserve $\Lambda(V)_{[0]}$ are those corresponding to elements of $C \cap N$; so $C_{N_{\Lambda(V)_{[0]}}}(y_0) = C \cap N$.

Next, let $\Xi = \Phi^+ \setminus \Phi_A$, and set $U' = \prod_{\alpha \in \Xi} X_\alpha$; then $U = U'.(C \cap U)$ and $U' \cap (C \cap U) = \{ 1 \}$. We now observe that if $\alpha \in \Xi$ then $\nu_i + \alpha$ is a weight in $V$ for exactly one value of $i$; moreover each weight in $V$ of positive generalized height is of the form $\nu_i + \alpha$ for exactly two such roots $\alpha$, one for each value of $i$. Thus if we take $u = \prod x_\alpha(t_\alpha) \in U'$ satisfying $u.y_0 = y_0$, and equate coefficients of weight vectors, taking them in an order compatible with increasing generalized height, we see that for all $\alpha$ we must have $t_\alpha = 0$, so that $u = 1$; so $C_U(y_0) = C \cap U$.

Thus $C_U(y_0), C_{N_{\Lambda(V)_{[0]}}}(y_0) \leq C$, so we do indeed have $C_G(y_0) = C$.

Since $\dim(\overline{G.y_0}) = \dim G - \dim C_G(y_0) = 21 - 9 = 12 = \dim \G{2}(V)$, the orbit $G.y_0$ is dense in $\G{2}(V)$. Thus the quadruple $(G, \lambda, p, k)$ has generic stabilizer $C_G(y_0)/Z(G) \cong A_2T_1.\Z_2$, where the $A_2$ is of simply connected type.

Now take $k = 3$. This time we take the generalized height function on the weight lattice of $G$ whose value at $\alpha_1$ and $\alpha_3$ is $0$, and at $\alpha_2$ is $1$; then the generalized height of $\lambda = \frac{1}{2}(\alpha_1 + 2\alpha_2 + 3\alpha_3)$ is $1$, and as $\lambda$ and $\Phi$ generate the weight lattice we see that the generalized height of any weight is an integer. Since $V_\lambda = \langle e_\delta \rangle$ where $\delta = \ffourrt1231$, we see that if $\mu \in \Lambda(V)$ and $e_\alpha \in V_\mu$ where $\alpha = \sum m_i \beta_i$ with $m_4 = 1$, then the generalized height of $\mu$ is $m_2 - 1$. Thus $\Lambda(V)_{[0]} = \{ \nu_1, \nu_2, \nu_3, \nu_4 \}$, where we write
$$
\gamma_1 = \ffourrt1111, \quad \gamma_2 = \ffourrt0121, \quad \gamma_3 = \ffourrt1121, \quad \gamma_4 = \ffourrt0111,
$$
and for each $i$ we let $\nu_i$ be the weight such that $V_{\nu_i} = \langle e_{\gamma_i} \rangle$. Observe that if we take $s = \prod_{i = 1}^3 h_{\beta_i}(\kappa_i) \in T$, then $\nu_1(s) = \frac{\kappa_1}{\kappa_3}$, $\nu_2(s) = \frac{\kappa_3}{\kappa_1}$, $\nu_3(s) = \frac{\kappa_1\kappa_3}{\kappa_2}$ and $\nu_4(s) = \frac{\kappa_2}{\kappa_1\kappa_3}$; thus given any pair $(n_1, n_2)$ of integers we have $c_1\nu_1 + \cdots + c_4 \nu_4 = 0$ for $(c_1, c_2, c_3, c_4) = (n_1, n_1, n_2, n_2)$, and hence $\Lambda(V)_{[0]}$ has ZLC. Take $Y = \G{3}(V_{[0]})$ and write
$$
y_0 = \langle e_{\gamma_1} + e_{\gamma_2}, e_{\gamma_3}, e_{\gamma_4} \rangle \in Y.
$$

We know that the pointwise stabilizer in $W$ of $\{ \gamma_1, \gamma_2 \}$ is $\langle w_{\beta_2}, w_{\beta_1 + \beta_2 + 2\beta_3} \rangle$; in this group the stabilizer of $\gamma_3$ is $\langle w_{\beta_1 + 2\beta_2 + 2\beta_3} \rangle$, which also stabilizes $\gamma_4$, so the pointwise stabilizer in $W$ of $\{ \gamma_1, \gamma_2, \gamma_3, \gamma_4 \}$ is $\langle w_{\beta_1 + 2\beta_2 + 2\beta_3} \rangle$. Now $\langle w_{\beta_1}, w_{\beta_3} \rangle$ acts simply transitively on $\{ \gamma_1, \gamma_2, \gamma_3, \gamma_4 \}$; as $\gamma_1$ is orthogonal to $\gamma_2$ but not to $\gamma_3$ or $\gamma_4$, and no element in the pointwise stabilizer of $\{ \gamma_1, \gamma_2 \}$ interchanges $\gamma_3$ and $\gamma_4$, it follows that the setwise stabilizer in $W$ of $\Lambda(V)_{[0]}$ is $\langle w_{\beta_1 + 2\beta_2 + 2\beta_3}, w_{\beta_1}, w_{\beta_3} \rangle$. Note that this stabilizes $\Phi_{[0]} = \langle \alpha_1, \alpha_3 \rangle = \langle \beta_1, \beta_3 \rangle$.

Let $A$ be the ${A_1}^2$ subgroup having simple root groups $\{ x_{\beta_1}(t) x_{\beta_3}(t) : t \in K \}$ and $X_{\beta_1 + 2\beta_2 + 2\beta_3}$; then $Z(A) = \langle h_{\beta_1}(-1) h_{\beta_3}(-1) \rangle$. Set $C = Z(G) A$. Clearly we have $C \leq C_G(y_0)$; we shall show that in fact $C_G(y_0) = C$.

We have $U_{[0]} = X_{\alpha_1} X_{\alpha_3}$. Given $u \in U_{[0]}$, the weights $\nu_1$ and $\nu_2$ occur in $u.(e_{\gamma_1} + e_{\gamma_2})$, while $\nu_3$ and $\nu_4$ occur in $u.e_{\gamma_3}$ and $u.e_{\gamma_4}$ respectively, so the set of weights occurring in $u.y_0$ is $\Lambda(V)_{[0]}$. By Lemma~\ref{lem: gen height zero not strictly positive}, we have $C_G(y_0) = C_{U_{[+]}}(y_0) C_{G_{[0]} N_{\Lambda(V)_{[0]}}}(y_0) C_{U_{[+]}}(y_0)$.

First, since $W_{\Lambda(V)_{[0]}} = \langle w_{\beta_1 + 2\beta_2 + 2\beta_3}, w_{\beta_1}, w_{\beta_3} \rangle$ and $\beta_1, \beta_3 \in \Phi_{[0]}$, we have $G_{[0]} N_{\Lambda(V)_{[0]}} = G_{[0]} \langle n_{\beta_1 + 2\beta_2 + 2\beta_3} \rangle$. Any element of this last group may be written as $g^*c$ where $c \in \langle x_{\beta_1}(t) x_{\beta_3}(t), x_{-\beta_1}(t) x_{-\beta_3}(t) : t \in K \rangle (\langle X_{\pm(\beta_1 + 2\beta_2 + 2\beta_3)} \rangle \cap N) < C$ and $g^* \in \langle X_{\pm\beta_3} \rangle$. Suppose then that $g^* \in C_G(y_0)$. If $g^* = x_{\beta_3}(t) h_{\beta_3}(\kappa)$ for some $t \in K$ and $\kappa \in K^*$, then we must have $t = 0$ as otherwise $g^*.e_{\gamma_4}$ has a term $e_{\gamma_2}$ but no term $e_{\gamma_1}$, so cannot lie in $y_0$; then $g^*.(e_{\gamma_1} + e_{\gamma_2}) = \kappa^{-1} e_{\gamma_1} + \kappa e_{\gamma_2}$, and for this to lie in $y_0$ we need $\kappa^2 = 1$, so that $g^* \in \langle h_{\beta_3}(-1) \rangle = Z(G)$. If instead $g^* = x_{\beta_3}(t) h_{\beta_3}(\kappa) n_{\beta_3} x_{\beta_3}(t')$ for some $t, t' \in K$ and $\kappa \in K^*$, then $g^*.e_{\gamma_3}$ has a term $e_{\gamma_1}$ but no term $e_{\gamma_2}$, so cannot lie in $y_0$. Thus $g^* \in Z(G) < C$; so $C_{G_{[0]} N_{\Lambda(V)_{[0]}}}(y_0) = C \cap G_{[0]} N_{\Lambda(V)_{[0]}}$.

Next, let $\Xi = \Phi^+ \setminus \{ \alpha_1, \alpha_3, \alpha_1 + 2\alpha_2 + 2\alpha_3 \}$, and set $U' = \prod_{\alpha \in \Xi} X_\alpha$; then $U_{[+]} = U'.(C \cap U_{[+]})$ and $U' \cap (C \cap U_{[+]}) = \{ 1 \}$. Now take $u = \prod_{\alpha \in \Xi} x_\alpha(t_\alpha) \in U'$ satisfying $u.y_0 = y_0$. The requirement that in $u.e_{\gamma_3}$ the coefficient of $e_\gamma$ for $\gamma = \ffourrt1221$ and $\ffourrt1231$ should be zero shows that $t_\alpha = 0$ for $\alpha = \alpha_2$ and $\alpha_2 + \alpha_3$ respectively; considering likewise $u.e_{\gamma_4}$ we see that the same is true for $\alpha = \alpha_1 + \alpha_2 + \alpha_3$ and $\alpha_1 + \alpha_2 + 2\alpha_3$ respectively; finally treating $u.(e_{\gamma_1} + e_{\gamma_2})$ shows that the same is true for $\alpha = \alpha_1 + \alpha_2$ and $\alpha_2 + 2\alpha_3$ respectively. Hence $u = 1$, so $C_{U_{[+]}}(y_0) = C \cap U_{[+]}$.

Thus $C_{U_{[+]}}(y_0), C_{G_{[0]} N_{\Lambda(V)_{[0]}}}(y_0) \leq C$, so we do indeed have $C_G(y_0) = C$.

Since $\dim(\overline{G.y_0}) = \dim G - \dim C_G(y_0) = 21 - 6 = 15 = \dim \G{3}(V)$, the orbit $G.y_0$ is dense in $\G{3}(V)$. Thus the quadruple $(G, \lambda, p, k)$ has generic stabilizer $C_G(y_0)/Z(G) \cong {A_1}^2$, where the ${A_1}^2$ is a central product.
\end{proof}

\begin{prop}\label{prop: C_3, omega_3 module, p = 2, k = 2 or 3}
Let $G = C_3$ and $\lambda = \omega_3$ with $p = 2$, and take $k = 2$ or $3$. Then the quadruple $(G, \lambda, p, k)$ has generic stabilizer $\tilde A_2T_1.\Z_2$ or ${\tilde A_1}{}^2$ respectively.
\end{prop}

\begin{proof}
This is an immediate consequence of Proposition~\ref{prop: B_3, omega_3 module, k = 2 or 3}, using the exceptional isogeny $B_\ell \to C_\ell$ which exists in characteristic $2$.
\end{proof}

\begin{prop}\label{prop: B_3, omega_3 module, k = 4}
Let $G = B_3$ and $\lambda = \omega_3$, and take $k = 4$. Then the quadruple $(G, \lambda, p, k)$ has generic stabilizer ${B_1}^2$ if $p \geq 3$, and semi-generic (but not generic) stabilizer ${B_1}^2$ if $p = 2$.
\end{prop}

\begin{proof}
Take $H$ to be the simply connected group defined over $K$ of type $E_6$, with simple roots $\beta_1, \dots, \beta_6$. Let $G^+$ be the $D_4$ subgroup having simple roots $\beta_3$, $\beta_4$, $\beta_2$ and $\beta_5$; for convenience we denote the positive roots of $G^+$ as
\begin{eqnarray*}
& \delta_1 = \esixrt001000, \quad \phantom{{}_{11}} \delta_2 = \esixrt000100, \quad \phantom{{}_{11}} \delta_3 = \esixrt010000, \quad \phantom{{}_{11}} \delta_4 = \esixrt000010, & \\
& \delta_5 = \esixrt001100, \quad \phantom{{}_{11}} \delta_6 = \esixrt010100, \quad \phantom{{}_{11}} \delta_7 = \esixrt000110, \quad \phantom{{}_{11}} \delta_8 = \esixrt011100, & \\
& \delta_9 = \esixrt001110, \quad \phantom{{}_{1}} \delta_{10} = \esixrt010110, \quad \phantom{{}_{1}} \delta_{11} = \esixrt011110, \quad \phantom{{}_{1}} \delta_{12} = \esixrt011210. &
\end{eqnarray*}
Let $G$ be the $B_3$ subgroup of $G^+$ having long simple roots $\delta_1$ and $\delta_2$, and short simple root group $\{ x_{\delta_3}(t) x_{\delta_4}(t) : t \in K \}$; then the other two positive short root groups are $\{ x_{\delta_6}(t) x_{\delta_7}(-t) : t \in K \}$ and $\{ x_{\delta_8}(t) x_{\delta_9}(-t) : t \in K \}$. We may take $V = \langle e_\alpha : \alpha = \sum m_i \beta_i, \ m_1 = 0, \ m_6 = 1 \rangle < \L(H)$; then $V$ is the irreducible $G^+$-module with high weight $\omega_4$, and the restriction of $V$ to $G$ is the required irreducible module with high weight $\lambda = \omega_3$. If we write
\begin{eqnarray*}
& \gamma_1 = \esixrt000001, \quad \gamma_2 = \esixrt000011, \quad \gamma_3 = \esixrt000111, \quad \gamma_4 = \esixrt001111, & \\
& \gamma_5 = \esixrt010111, \quad \gamma_6 = \esixrt011111, \quad \gamma_7 = \esixrt011211, \quad \gamma_8 = \esixrt011221, &
\end{eqnarray*}
then $V = \langle e_{\gamma_1}, \dots, e_{\gamma_8} \rangle$. We have $Z(G) = \langle z \rangle$ where $z = h_{\beta_2}(-1) h_{\beta_5}(-1)$.

Set
$$
Y = \left\{ \langle a_1 e_{\gamma_1} + a_2 e_{\gamma_2}, a_3 e_{\gamma_4} + a_4 e_{\gamma_6}, e_{\gamma_5}, e_{\gamma_7} \rangle : (a_1, a_2), (a_3, a_4) \neq (0, 0) \right\},
$$
and
$$
\hat Y = \left\{ \langle a_1 e_{\gamma_1} + a_2 e_{\gamma_2}, a_3 e_{\gamma_4} + a_4 e_{\gamma_6}, e_{\gamma_5}, e_{\gamma_7} \rangle : a_1a_2a_3a_4 \neq 0, \ a_1a_4 \neq a_2a_3 \right\},
$$
so that $\hat Y$ is a dense open subset of $Y$. Take
$$
y = \langle a_1 e_{\gamma_1} + a_2 e_{\gamma_2}, a_3 e_{\gamma_4} + a_4 e_{\gamma_6}, e_{\gamma_5}, e_{\gamma_7} \rangle \in \hat Y.
$$
Observe that if we set $s = h_{\beta_4}(\frac{a_2}{a_1}) \in T$, then
$$
s.y = \langle e_{\gamma_1} + e_{\gamma_2}, e_{\gamma_4} + a e_{\gamma_6}, e_{\gamma_5}, e_{\gamma_7} \rangle,
$$
where we write $a = \frac{a_1a_4}{a_2a_3}$.

Take $g_a = x_{\beta_2}(a) x_{-\beta_5}(1) \in G^+$; then ${g_a}^{-1}.(s.y) = \langle e_{\gamma_2}, e_{\gamma_4}, e_{\gamma_5}, e_{\gamma_7} \rangle$, whose stabilizer in $G^+$ is the ${D_2}^2.\Z_2$ subgroup with connected component having simple roots $\delta_5$, $\delta_6$, $\delta_7$ and $\delta_{11}$, and component group generated by the image of $n_{\delta_1} n_{\delta_3}$. Conjugating by $g_a$ we see that $C_{G^+}(s.y)$ has simple factors
\begin{eqnarray*}
\langle x_{\delta_5}(t) x_{\delta_8}(at), x_{-\delta_5}(t) x_{-\delta_9}(t) : t \in K \rangle, & & \\
\langle x_{\delta_6}(t), x_{-\delta_2}(-at) x_{-\delta_6}(t) x_{-\delta_7}(-at) x_{-\delta_{10}}(t) : t \in K \rangle, & & \\
\langle x_{\delta_2}(-t) x_{\delta_6}(-at) x_{\delta_7}(t) x_{\delta_{10}}(at), x_{-\delta_7}(t) : t \in K \rangle, & & \\
\langle x_{\delta_8}(t) x_{\delta_{11}}(t), x_{-\delta_9}(-at) x_{-\delta_{11}}(t) : t \in K \rangle, & &
\end{eqnarray*}
and component group generated by the image of $x_{\delta_3}(a) n_{\delta_1} n_{\delta_3} x_{\delta_3}(-a)$. Taking the intersection with $G$, we see that if we take $\kappa_1, \kappa_2 \in K^*$ satisfying ${\kappa_1}^2 = \frac{1}{a}$ and ${\kappa_2}^2 = \frac{1}{a - 1}$, and for $t \in K$ we write
\begin{eqnarray*}
   x_1^a(t) & = & x_{\delta_5}(\kappa_1 t) x_{\delta_{11}}(a\kappa_1 t), \\
x_{-1}^a(t) & = & x_{-\delta_5}(a\kappa_1 t) x_{-\delta_{11}}(\kappa_1 t), \\
   x_2^a(t) & = & x_{\delta_2}(-\kappa_2 t) x_{\delta_6}(-\kappa_2 t) x_{\delta_7}(\kappa_2 t) x_{\delta_{10}}(a\kappa_2 t), \\
x_{-2}^a(t) & = & x_{-\delta_2}(-a\kappa_2 t) x_{-\delta_6}(\kappa_2 t) x_{-\delta_7}(-\kappa_2 t) x_{-\delta_{10}}(\kappa_2 t),
\end{eqnarray*}
then we have
\begin{eqnarray*}
C_G(s.y) & = & \langle x_1^a(t), x_{-1}^a(t) : t \in K \rangle \langle x_2^a(t), x_{-2}^a(t) : t \in K \rangle \\
         & = & {}^{g_a}\left( \langle x_{\delta_5}(\kappa_1 t) x_{\delta_{11}}(a\kappa_1 t), x_{-\delta_5}(a\kappa_1 t) x_{-\delta_{11}}(\kappa_1 t) : t \in K \rangle \right. \\
         &   & \quad \left. \times \langle x_{\delta_6}((a - 1)\kappa_2 t) x_{\delta_7}(\kappa_2 t), x_{-\delta_6}(\kappa_2 t) x_{-\delta_7}((a - 1)\kappa_2 t) : t \in K \rangle \right),
\end{eqnarray*}
so that $C_G(s.y) \cong {B_1}^2$.

Now given
$$
y' = \langle b_1 e_{\gamma_1} + b_2 e_{\gamma_2}, b_3 e_{\gamma_4} + b_4 e_{\gamma_6}, e_{\gamma_5}, e_{\gamma_7} \rangle \in Y,
$$
provided $b_1b_2b_3b_4 \neq 0$ we may take $s' = h_{\beta_4}(\frac{b_2}{b_1}) \in T$ so that
$$
s'.y' = \langle e_{\gamma_1} + e_{\gamma_2}, e_{\gamma_4} + b e_{\gamma_6}, e_{\gamma_5}, e_{\gamma_7} \rangle,
$$
where we write $b = \frac{b_1b_4}{b_2b_3}$, and then the set of elements of $G^+$ sending $s.y$ to $s'.y'$ is $x_{\delta_3}(b - a) C_{G^+}(s.y)$, whose intersection with $G$ is clearly empty if $b \neq a$. Thus $\dim \overline{G.y \cap Y} = 1$, while as $\dim C_G(y) = \dim C_G(s.y) = 6$ we have $\dim \overline{G.y} = \dim G - \dim C_G(y) = 21 - 6 = 15$; therefore
$$
\dim \G{4}(V) - \dim(\overline{G.y}) = 16 - 15 = 1 \quad \hbox{and} \quad \dim Y - \dim(\overline{G.y \cap Y}) = 2 - 1 = 1.
$$
Hence $y$ is $Y$-exact.

First suppose $p \geq 3$, and set
$$
C = \langle X_{\pm(\alpha_1 + \alpha_2 + \alpha_3)} \rangle \langle x_{\alpha_2}(t) x_{\alpha_2 + 2\alpha_3}(t), x_{-\alpha_2}(t) x_{-(\alpha_2 + 2\alpha_3)}(t) : t \in K \rangle \cong {B_1}^2;
$$
then each factor of $C$ has centre $Z(G)$. Take $\kappa \in K^*$ satisfying $\kappa^4 = \frac{1}{4a} \frac{\kappa_1 - 1}{\kappa_1 + 1}$; then with $h^{-1} = h_{\alpha_3}(\kappa) x_{-\alpha_3}(-\frac{\kappa_1}{2}) x_{\alpha_3}(\frac{1}{\kappa_1})$ we have $C_G(s.y) = {}^h C$.  Thus the conditions of Lemma~\ref{lem: generic stabilizer from exact subset} hold; so the quadruple $(G, \lambda, p, k)$ has generic stabilizer $C/Z(G) = {B_1}^2$, where each $B_1$ factor is of adjoint type.

Now suppose $p = 2$; the above shows that the image of $\hat Y$ under the orbit map is dense in $\G{4}(V)$, and all points in this dense subset have stabilizer isomorphic to ${B_1}^2$, so the quadruple $(G, \lambda, p, k)$ has semi-generic stabilizer ${B_1}^2$. However, we claim that there is no generic stabilizer. If there were, then for infinitely many values $a \neq 0, 1$ the stabilizers $C_G(s.y)$ above would be conjugate, so certainly there would exist $a, b \in K \setminus \{ 0, 1 \}$ distinct and $g \in G$ such that for $i = 1, 2$ we have ${}^g \langle x_i^a(t), x_{-i}^a(t) : t \in K \rangle = \langle x_i^b(t), x_{-i}^b(t) : t \in K \rangle$. We note the following well-known facts about an $A_1$ group $A$: any two maximal tori of $A$ are conjugate; given a maximal torus $T_A$ of $A$, there are just two $1$-dimensional unipotent subgroups of $A$ normalized by $T_A$, say $U_A^+$ and $U_A^-$, which are interchanged by elements of $N_A(T_A) \setminus T_A$; the torus $T_A$ acts transitively on the non-identity elements of each of $U_A^+$ and $U_A^-$; for each non-identity element $u^+$ of $U_A^+$ there is a unique non-identity element $u^-$ of $U_A^-$ such that $u^- u^+ u^- \in N_A(T_A)$. It follows that by multiplying $g$ by an element of $C_G(s.y)$ we may assume that for $i = 1, -1, 2, -2$ we have ${}^g \langle x_i^a(t) : t \in K \rangle = \langle x_i^b(t) : t \in K \rangle$ and ${}^g x_i^a(1) = x_i^b(1)$.

We may regard $G$ as consisting of $7 \times 7$ matrices. Taking $\kappa_1, \kappa_2 \in K^*$ satisfying ${\kappa_1}^2 = \frac{1}{a}$ and ${\kappa_2}^2 = \frac{1}{a - 1}$ as above, $x_1^a(t)$ and $x_{-1}^a(t)$ are the matrices
$$
\left(
  \begin{array}{ccccccc}
          1       &   & \kappa_1 t &   & a\kappa_1 t &   &     t^2     \\
                  & 1 &            &   &             &   &             \\
                  &   &      1     &   &             &   & a\kappa_1 t \\
                  &   &            & 1 &             &   &             \\
                  &   &            &   &      1      &   &  \kappa_1 t \\
                  &   &            &   &             & 1 &             \\
    \phantom{t^2} &   &            &   &             &   &      1      \\
  \end{array}
\right), \quad
\left(
  \begin{array}{ccccccc}
        1      &   &            &   &             &   & \phantom{t^2} \\
               & 1 &            &   &             &   &               \\
   a\kappa_1 t &   &      1     &   &             &   &               \\
               &   &            & 1 &             &   &               \\
    \kappa_1 t &   &            &   &      1      &   &               \\
               &   &            &   &             & 1 &               \\
       t^2     &   & \kappa_1 t &   & a\kappa_1 t &   &       1       \\
  \end{array}
\right)
$$
respectively, while $x_2^a(t)$ and $x_{-2}^a(t)$ are the matrices
$$
\left(
  \begin{array}{ccccccc}
    1 &               &            &   &             &             &   \\
      &       1       & \kappa_2 t &   & a\kappa_2 t &     t^2     &   \\
      &               &      1     &   &             & a\kappa_2 t &   \\
      &               &            & 1 &             &  \kappa_2 t &   \\
      &               &            &   &      1      &  \kappa_2 t &   \\
      & \phantom{t^2} &            &   &             &      1      &   \\
      &               &            &   &             &             & 1 \\
  \end{array}
\right), \quad
\left(
  \begin{array}{ccccccc}
    1 &             &            &   &             &               &   \\
      &      1      &            &   &             & \phantom{t^2} &   \\
      & a\kappa_2 t &      1     &   &             &               &   \\
      &  \kappa_2 t &            & 1 &             &               &   \\
      &  \kappa_2 t &            &   &      1      &               &   \\
      &     t^2     & \kappa_2 t &   & a\kappa_2 t &       1       &   \\
      &             &            &   &             &               & 1 \\
  \end{array}
\right)
$$
respectively; the matrices $x_i^b(t)$ are obtained by replacing $a$ by $b$ throughout. For $i = 1, -1, 2, -2$ write $x_i^a(t) = I + A_i^{(1)} t + A_i^{(2)} t^2$ and $x_i^b(t) = I + B_i^{(1)} t + B_i^{(2)} t^2$, where $A_i^{(1)}$, $A_i^{(2)}$, $B_i^{(1)}$, $B_i^{(2)}$ are independent of $t$; let $D$ be the matrix representing $g$. Then for each $i$, for all $t \in K$ there exists $t' \in K$ such that $D(I + A_i^{(1)} t + A_i^{(2)} t^2)D^{-1} = I + B_i^{(1)} t' + B_i^{(2)} {t'}^2$, whence $D A_i^{(1)} D^{-1} t + D A_i^{(2)} D^{-1} t^2 = B_i^{(1)} t' + B_i^{(2)} {t'}^2$. Thus the matrices $D A_i^{(1)} D^{-1}$ and $D A_i^{(2)} D^{-1}$ must be linear combinations of the matrices $B_i^{(1)}$ and $B_i^{(2)}$, so we may write $D A_i^{(1)} D^{-1} = c_1 B_i^{(1)} + c_2 B_i^{(2)}$, $D A_i^{(2)} D^{-1} = c_3 B_i^{(1)} + c_4 B_i^{(2)}$ with $c_1, c_2, c_3, c_4 \in K$, and then for all $t \in K$ there exists $t' \in K$ such that $c_1 B_i^{(1)} t + c_2 B_i^{(2)} t + c_3 B_i^{(1)} t^2 + c_4 B_i^{(2)} t^2 = B_i^{(1)} t' + B_i^{(2)} {t'}^2$; as $B_i^{(1)}$ and $B_i^{(2)}$ are linearly independent this implies $c_1t + c_3t^2 = t'$ and $c_2t + c_4t^2 = {t'}^2$, so that $c_2t + c_4t^2 = (c_1t + c_3t^2)^2 = {c_1}^2t^2 + {c_3}^2t^4$, and as this is true for all $t$ we must have $c_2 = c_3 = 0$, ${c_1}^2 = c_4$; since by assumption $t = 1$ implies $t' = 1$, we must have $c_1 = c_4 = 1$. Therefore $D A_i^{(1)} D^{-1} = B_i^{(1)}$ and $D A_i^{(2)} D^{-1} = B_i^{(2)}$.

Now $A_1^{(2)} = B_1^{(2)} = E_{17}$, $A_{-1}^{(2)} = B_{-1}^{(2)} = E_{71}$, $A_2^{(2)} = B_2^{(2)} = E_{26}$ and $A_{-2}^{(2)} = B_{-2}^{(2)} = E_{62}$, where $E_{ij}$ is the matrix unit with $(i, j)$-entry $1$ and all other entries $0$; as $D$ must commute with each of these four matrix units, and must preserve the relevant quadratic form, we must have
$$
D =
\left(
  \begin{array}{ccccccc}
    1 &   &     &   &     &   &   \\
      & 1 &     &   &     &   &   \\
      &   & d_1 &   & d_2 &   &   \\
      &   & d_3 & 1 & d_4 &   &   \\
      &   & d_5 &   & d_6 &   &   \\
      &   &     &   &     & 1 &   \\
      &   &     &   &     &   & 1 \\
  \end{array}
\right)
$$
for some $d_1, \dots, d_6 \in K$ such that $d_1d_5 + {d_3}^2 = 0 = d_2d_6 + {d_4}^2$ and $d_1d_6 + d_2d_5 = 1$. The condition $D A_1^{(1)} D^{-1} = B_1^{(1)}$ then requires $d_1 = \frac{\sqrt{b}}{\sqrt{a}} + bd_5$, $d_6 = \frac{\sqrt{a}}{\sqrt{b}} + ad_5$, $d_2 = abd_5$ and $d_4 = ad_3$, and then the condition $D A_{-1}^{(1)} D^{-1} = B_{-1}^{(1)}$ is also satisfied; however, equating the $(4, 6)$-entries of $D A_2^{(1)} D^{-1}$ and $B_2^{(1)}$ then gives $\frac{1}{\sqrt{a - 1}} = \frac{1}{\sqrt{b - 1}}$, which is impossible. Thus there is no such matrix $D$; so the claim is proved, and the result follows.
\end{proof}

\begin{prop}\label{prop: C_3, omega_3 module, p = 2, k = 4}
Let $G = C_3$ and $\lambda = \omega_3$ with $p = 2$, and take $k = 4$. Then the quadruple $(G, \lambda, p, k)$ has semi-generic (but not generic) stabilizer ${C_1}^2$.
\end{prop}

\begin{proof}
This is an immediate consequence of Proposition~\ref{prop: B_3, omega_3 module, k = 4}, using the exceptional isogeny $B_\ell \to C_\ell$ which exists in characteristic $2$.
\end{proof}

\begin{prop}\label{prop: D_5, omega_5, B_4, omega_4 modules, k = 2}
Let $G = D_5$ and $\lambda = \omega_5$, or $G = B_4$ and $\lambda = \omega_4$, and take $k = 2$. Then the quadruple $(G, \lambda, p, k)$ has generic stabilizer $G_2 B_1$ or $A_2 T_1.\Z_2$ respectively.
\end{prop}

\begin{proof}
We begin with the case where $G = D_5$ and $\lambda = \omega_5$. We use the set-up of Proposition~\ref{prop: D_5, omega_5, B_4, omega_4 modules}: we take $H$ to be the simply connected group defined over $K$ of type $E_6$, with simple roots $\beta_1, \dots, \beta_6$; we let $G$ have simple roots $\alpha_1 = \beta_1$, $\alpha_2 = \beta_3$, $\alpha_3 = \beta_4$, $\alpha_4 = \beta_5$, $\alpha_5 = \beta_2$, so that $G = \langle X_\alpha : \alpha = \sum m_i \beta_i, \ m_6 = 0 \rangle < H$; then we may take $V = \langle e_\alpha : \alpha = \sum m_i \beta_i, \ m_6 = 1 \rangle < \L(H)$. We have $Z(G) = \langle z \rangle$ where $z = h_{\beta_1}(-1) h_{\beta_2}(\eta_4) h_{\beta_4}(-1) h_{\beta_5}(-\eta_4)$. Here we take the generalized height function on the weight lattice of $G$ whose value at $\alpha_4$ and $\alpha_5$ is $0$, and at $\alpha_1$, $\alpha_2$ and $\alpha_3$ is $1$; then the generalized height of $\lambda = \frac{1}{2}(\alpha_1 + 2\alpha_2 + 3\alpha_3 + \frac{3}{2}\alpha_4 + \frac{5}{2}\alpha_5)$ is $3$, and as $\lambda$, $\omega_4 = \lambda + \frac{1}{2}\alpha_4 - \frac{1}{2}\alpha_5$ and $\Phi$ generate the weight lattice it follows that the generalized height of any weight is an integer. Since $V_\lambda = \langle e_\delta \rangle$ where $\delta = \esixrt122321$, we see that if $\mu \in \Lambda(V)$ and $e_\alpha \in V_\mu$ where $\alpha = \sum m_i \beta_i$ with $m_6 = 1$, then the generalized height of $\mu$ is $m_1 + m_3 + m_4 - 3$. Thus $\Lambda(V)_{[0]} = \{ \nu_1, \nu_2, \nu_3, \nu_4 \}$, where we write
$$
\gamma_1 = \esixrt111111, \quad \gamma_2 = \esixrt011221, \quad \gamma_3 = \esixrt101111, \quad \gamma_4 = \esixrt011211,
$$
and for each $i$ we let $\nu_i$ be the weight such that $V_{\nu_i} = \langle e_{\gamma_i} \rangle$. Observe that if we take $s = \prod_{i = 1}^5 h_{\beta_i}(\kappa_i) \in T$, then $\nu_1(s) = \frac{\kappa_1 \kappa_2}{\kappa_4}$, $\nu_2(s) = \frac{\kappa_5}{\kappa_1}$, $\nu_3(s) = \frac{\kappa_1}{\kappa_2}$ and $\nu_4(s) = \frac{\kappa_4}{\kappa_1 \kappa_5}$; thus $\nu_1 + \nu_2 + \nu_3 + \nu_4 = 0$, and hence $\Lambda(V)_{[0]}$ has ZLC. Set
$$
Y = \left\{ \langle a_1 e_{\gamma_1} + a_2 e_{\gamma_2}, a_3 e_{\gamma_3} + a_4 e_{\gamma_4} \rangle : (a_1, a_2), (a_3, a_4) \neq (0, 0) \right\} \subset \G{2}(V_{[0]}),
$$
and
$$
\hat Y = \left\{\langle a_1 e_{\gamma_1} + a_2 e_{\gamma_2}, a_3 e_{\gamma_3} + a_4 e_{\gamma_4} \rangle : a_1a_2a_3a_4 \neq 0 \right\},
$$
so that $\hat Y$ is a dense open subset of $Y$. Write
$$
y_0 = \langle e_{\gamma_1} + e_{\gamma_2}, e_{\gamma_3} + e_{\gamma_4} \rangle \in \hat Y.
$$

In the proof of Proposition~\ref{prop: D_5, omega_5, B_4, omega_4 modules} we observed that $W$ acts transitively on the set $\Sigma$ of roots $\alpha$ of $H$ corresponding to the root vectors $e_\alpha$ spanning $V$, and that the stabilizer of any one root acts transitively on the $5$ roots orthogonal to it; so the pointwise stabilizer in $W$ of $\{ \gamma_1, \gamma_2 \}$ has size $24 = |W(A_3)|$, and we see that it is $\langle w_{\beta_3}, w_{\beta_2 + \beta_4}, w_{\beta_1 + \beta_3 + \beta_4 + \beta_5} \rangle$. In this $S_4$ subgroup the stabilizer of $\gamma_3$ contains and therefore equals the maximal subgroup $\langle w_{\beta_3}, w_{\beta_1 + \beta_2 + \beta_3 + 2\beta_4 + \beta_5} \rangle$, and as this group also stabilizes $\gamma_4$ it is the pointwise stabilizer in $W$ of $\{ \gamma_1, \gamma_2, \gamma_3, \gamma_4 \}$. Now $w_{\beta_5}$ interchanges $\gamma_2$ and $\gamma_4$ while fixing both $\gamma_1$ and $\gamma_3$, while $w_{\beta_1} w_{\beta_2 + \beta_4} w_{\beta_4 + \beta_5}$ interchanges $\gamma_1$ and $\gamma_2$, and also $\gamma_3$ and $\gamma_4$. Thus the setwise stabilizer of $\{ \gamma_1, \gamma_2, \gamma_3, \gamma_4 \}$ acts transitively on it, and as $\gamma_1$ is orthogonal to $\gamma_2$ and $\gamma_4$ but not $\gamma_3$, any element of the setwise stabilizer which fixes $\gamma_1$ must also fix $\gamma_3$, so must either fix or interchange $\gamma_2$ and $\gamma_4$. Therefore the setwise stabilizer in $W$ of $\{ \gamma_1, \gamma_2, \gamma_3, \gamma_4 \}$, and hence of $\Lambda(V)_{[0]}$, is
$$
\langle w_{\beta_3}, w_{\beta_1 + \beta_2 + \beta_3 + 2\beta_4 + \beta_5}, w_{\beta_5}, w_{\beta_1} w_{\beta_2 + \beta_4} w_{\beta_4 + \beta_5} \rangle
= \langle w_{\beta_3}, w_{\beta_5}, w_{\beta_1} w_{\beta_2 + \beta_4} w_{\beta_4 + \beta_5} \rangle.
$$
Note that this stabilizes $\Phi_{[0]} = \langle \alpha_4, \alpha_5 \rangle = \langle \beta_2, \beta_5 \rangle$.

Let $A$ be the $G_2 B_1$ subgroup with the first factor having simple root groups $\{ x_{\beta_1}(t) x_{\beta_2 + \beta_4}(t) x_{\beta_4 + \beta_5}(t) : t \in K \}$ and $X_{\beta_3}$, and the second having simple root group $\{ x_{\beta_2}(t) x_{\beta_5}(t) : t \in K \}$; then $Z(A) = \langle h_{\beta_2}(-1) h_{\beta_5}(-1) \rangle = \langle z^2 \rangle < Z(G)$. Set $C = Z(G) A$. Clearly we have $C \leq C_G(y_0)$; we shall show that in fact $C_G(y_0) = C$.

We have $U_{[0]} = X_{\alpha_4} X_{\alpha_5}$. Given $u \in U_{[0]}$, the weights $\nu_1$ and $\nu_2$ occur in $u.(e_{\gamma_1} + e_{\gamma_2})$, while $\nu_3$ and $\nu_4$ occur in $u.(e_{\gamma_3} + e_{\gamma_4})$, so the set of weights occurring in $u.y_0$ is $\Lambda(V)_{[0]}$. By Lemma~\ref{lem: gen height zero not strictly positive}, if we take $g \in \Tran_G(y_0, Y)$ and write $y' = g.y_0 \in Y$, then we have $g = u_1 g' u_2$ with $u_1 \in C_{U_{[+]}}(y')$, $u_2 \in C_{U_{[+]}}(y_0)$, and $g' \in G_{[0]} N_{\Lambda(V)_{[0]}}$ with $g'.y_0 = y'$. In particular $G.y_0 \cap Y = G_{[0]} N_{\Lambda(V)_{[0]}}.y_0 \cap Y$; moreover $C_G(y_0) = C_{U_{[+]}}(y_0) C_{G_{[0]} N_{\Lambda(V)_{[0]}}}(y_0) C_{U_{[+]}}(y_0)$.

First, since $W_{\Lambda(V)_{[0]}} = \langle w_{\beta_3}, w_{\beta_5}, w_{\beta_1} w_{\beta_2 + \beta_4} w_{\beta_4 + \beta_5} \rangle$ and $\beta_2, \beta_5 \in \Phi_{[0]}$, we have $G_{[0]} N_{\Lambda(V)_{[0]}} = G_{[0]} \langle n_{\beta_3}, n_{\beta_1} n_{\beta_2 + \beta_4} n_{\beta_4 + \beta_5} \rangle$. Any element of this last group may be written as $g^*c$ where $c \in C$ and $g^* \in \langle X_{\pm\beta_5} \rangle \{ h_{\beta_4}(\kappa_4) : \kappa_4 \in K^* \}$. Suppose first that $g^*.y_0 \in Y$. If $g^* = x_{\beta_5}(t) h_{\beta_5}(\kappa_5) h_{\beta_4}(\kappa_4)$ for some $t \in K$ and $\kappa_4, \kappa_5 \in K^*$, then we must have $t = 0$ as otherwise $g^*.(e_{\gamma_3} + e_{\gamma_4})$ has a term $e_{\gamma_2}$ but no term $e_{\gamma_1}$, whereas $g^*.(e_{\gamma_1} + e_{\gamma_2})$ has both terms $e_{\gamma_1}$ and $e_{\gamma_2}$, so $g^*.y_0$ cannot lie in $Y$. If instead $g^* = x_{\beta_5}(t) h_{\beta_5}(\kappa_5) n_{\beta_5} x_{\beta_5}(t') h_{\beta_4}(\kappa_4)$ for some $t, t' \in K$ and $\kappa_4, \kappa_5 \in K^*$, then $g^*.(e_{\gamma_1} + e_{\gamma_2})$ has a term $e_{\gamma_4}$ but no term $e_{\gamma_3}$, whereas $g^*.(e_{\gamma_3} + e_{\gamma_4})$ has a term $e_{\gamma_3}$, so $g^*.y_0$ cannot lie in $Y$. Thus we must have $g^* = h_{\beta_5}(\kappa_5) h_{\beta_4}(\kappa_4) \in T$; so $G.y_0 \cap Y = T.y_0$. Now given $y = \langle a_1 e_{\gamma_1} + a_2 e_{\gamma_2}, a_3 e_{\gamma_3} + a_4 e_{\gamma_4} \rangle \in \hat Y$, if we take $\kappa \in K^*$ satisfying $\kappa^2 = \frac{a_2a_3}{a_1a_4}$, and set $h = h_{\beta_4}(\kappa \frac{a_4}{a_3}) h_{\beta_5}(\kappa)$, then $h.y_0 = y$; so we have $G.y_0 \cap Y = \hat Y$. If we now further require $g^*.y_0 = y_0$ then as $g^*.(e_{\gamma_1} + e_{\gamma_2}) = {\kappa_4}^{-1}e_{\gamma_1} + \kappa_5 e_{\gamma_2}$ and $g^*.(e_{\gamma_3} + e_{\gamma_4}) = e_{\gamma_3} + \kappa_4{\kappa_5}^{-1} e_{\gamma_4}$, we must have $\kappa_4 = \kappa_5 = \pm 1$, whence $g^* \in \langle h_{\beta_4}(-1) h_{\beta_5}(-1) \rangle$ --- as $h_{\beta_4}(-1) h_{\beta_5}(-1) = z.h_{\beta_1}(-1) h_{\beta_2 + \beta_4}(-1) h_{\beta_4 + \beta_5}(-1). h_{\beta_2}(\eta_4) h_{\beta_5}(\eta_4) \in Z(G) (A \cap T)$, we have $g^* \in C$. Thus $C_{G_{[0]} N_{\Lambda(V)_{[0]}}}(y_0) = C \cap G_{[0]} N_{\Lambda(V)_{[0]}}$.

Next, take the $D_3D_2$ subsystem $\Psi$ of $\Phi$ consisting of roots $\sum m_i \beta_i$ with $m_4$ even; then in the $G_2$ factor of $A$, each of the long root subgroups is $X_\alpha$ for some $\alpha \in \Psi$, and each of the short root subgroups is diagonally embedded in $X_\alpha X_{\alpha'} X_{\alpha''}$ for some $\alpha \in \Psi$ and $\alpha', \alpha'' \notin \Psi$. Since the $B_1$ factor of $A$ lies in $G_{[0]}$, let $\Xi = \Phi^+ \setminus \Psi$, and set $U' = \prod_{\alpha \in \Xi} X_\alpha$; then $U_{[+]} = U'.(C \cap U_{[+]})$ and $U' \cap (C \cap U_{[+]}) = \{ 1 \}$. We now observe that if $\alpha \in \Xi$ then $\nu_i + \alpha$ is a weight in $V$ for exactly one value of $i$; moreover each weight in $V$ of positive generalized height is of the form $\nu_i + \alpha$ for exactly two such roots $\alpha$, one having $i \in \{ 1, 2 \}$ and one having $i \in \{ 3, 4 \}$. Thus if we take $u = \prod x_\alpha(t_\alpha) \in U'$ satisfying $u.y_0 = y_0$, and equate coefficients of weight vectors, taking them in an order compatible with increasing generalized height, we see that for all $\alpha$ we must have $t_\alpha = 0$, so that $u = 1$; so $C_{U_{[+]}}(y_0) = C \cap U_{[+]}$.

Thus $C_{U_{[+]}}(y_0), C_{G_{[0]} N_{\Lambda(V)_{[0]}}}(y_0) \leq C$, so we do indeed have $C_G(y_0) = C$.

Since $\dim(\overline{G.y_0}) = \dim G - \dim C_G(y_0) = 45 - 17 = 28 = \dim \G{2}(V)$, the orbit $G.y_0$ is dense in $\G{2}(V)$. Thus the quadruple $(G, \lambda, p, k)$ has generic stabilizer $C_G(y_0)/Z(G) \cong G_2 B_1$, where the $B_1$ is of adjoint type.

Before continuing we note that given $y \in \hat Y$ there exists $h \in T$ with $h.y_0 = y$; since $h$ then lies in $G_{[0]} N_{\Lambda(V)_{[0]}}$ and normalises $U_{[+]}$, by conjugating $\Tran_G(y_0, Y)$ by $h$ we see that any element of $\Tran_G(y, Y)$ is of the form $g = u_1g'u_2$ with $u_1 \in C_{U_{[+]}}(g.y)$, $u_2 \in C_{U_{[+]}}(y)$ and $g' \in G_{[0]} N_{\Lambda(V)_{[0]}}$ with $g'.y = g.y$.

To treat the case where $G = B_4$ and $\lambda = \omega_4$, we leave $H$, $V$ and $Y$ unchanged, but replace $G$ by the $B_4$ subgroup of $D_5$ which has simple root groups $X_{\beta_1}$, $X_{\beta_3}$, $X_{\beta_4}$ and $\{ x_{\beta_2}(t) x_{\beta_5}(t) : t \in K \}$; then $Z(G) = \langle z^2 \rangle$. This time we set
$$
\hat Y = \left\{\langle a_1 e_{\gamma_1} + a_2 e_{\gamma_2}, a_3 e_{\gamma_3} + a_4 e_{\gamma_4} \rangle : a_1a_2a_3a_4 \neq 0, \ (a_2a_3)^2 \neq (a_1a_4)^2 \right\},
$$
which is still a dense open subset of $Y$.

Take $y = \langle a_1 e_{\gamma_1} + a_2 e_{\gamma_2}, a_3 e_{\gamma_3} + a_4 e_{\gamma_4} \rangle \in \hat Y$. We have seen above that if we take $\kappa \in K^*$ satisfying $\kappa^2 = \frac{a_2a_3}{a_1a_4}$, and set $h = h_{\beta_4}(\kappa \frac{a_4}{a_3}) h_{\beta_5}(\kappa)$, then $h.y_0 = y$; note that then $\kappa^2 \neq 1$. The $D_5$-stabilizer of $y$ is then ${}^h C$, which has simple root groups $\{ x_{\beta_1}(t) x_{\beta_2 + \beta_4}(\frac{a_4}{a_3}t) x_{\beta_4 + \beta_5}(\frac{a_2}{a_1}t) : t \in K \}$, $X_{\beta_3}$, and $\{ x_{\beta_2}(t) x_{\beta_5}(\frac{a_2a_3}{a_1a_4}t) : t \in K \}$; thus the $B_4$-stabilizer of $y$ is the intersection of this with $G$. We therefore let $A$ be the $A_2$ subgroup having simple roots $\beta_3$ and $\beta_1 + \beta_2 + \beta_3 + 2\beta_4 + \beta_5$, let $T_1$ be the $1$-dimensional torus $\{ h_{\beta_2}(\kappa) h_{\beta_5}(\kappa) : \kappa \in K^* \}$, and write $n^* = n_{\beta_1} n_{\beta_4} n_{\beta_2 + \beta_4 + \beta_5}$; we replace $C$ by $Z(G) AT_1 \langle n^* \rangle$, and we have $C_G(y) = {}^h C = Z(G) AT_1 \langle h_{\beta_2}(\kappa) h_{\beta_5}(\kappa) h_{\beta_4}(\frac{a_2a_4}{a_1a_3}) n^* \rangle$. If we now take $\kappa' \in K^*$ satisfying ${\kappa'}^2 = \kappa$ and let $h' = h_{\beta_2}(\kappa') h_{\beta_5}(\kappa') h_{\beta_4}(\kappa \frac{a_4}{a_3})$, then $h' \in G$ and $C_G(y) = {}^{h'} C$.

Now given $g = u_1g'u_2 \in \Tran_{D_5}(y, Y)$ as above, if $g$ is to lie in $B_4$ we clearly must have each of $u_1$, $g'$ and $u_2$ in $B_4$. Thus to determine $G.y \cap Y$ it suffices to consider the elements $g'$ lying in $B_4$. Since $\langle n_{\beta_3}, n_{\beta_1} n_{\beta_2 + \beta_4} n_{\beta_4 + \beta_5} \rangle < B_4$, and the intersection of the group $G_{[0]}$ above with $B_4$ is $\langle T, x_{\beta_2}(t) x_{\beta_5}(t), x_{-\beta_2}(t) x_{-\beta_5}(t) : t \in K \rangle$, we have $g' = xn$ where $x \in \langle T, x_{\beta_2}(t) x_{\beta_5}(t), x_{-\beta_2}(t) x_{-\beta_5}(t) : t \in K \rangle$ and $n \in \langle n_{\beta_3}, n_{\beta_1} n_{\beta_2 + \beta_4} n_{\beta_4 + \beta_5} \rangle$. Now if we write $y' = \langle a_2 e_{\gamma_1} + a_1 e_{\gamma_2}, a_4 e_{\gamma_3} + a_3 e_{\gamma_4} \rangle$, then we see that $n.y \in \{ y, y' \}$. If $x = s x_{\beta_2}(t) x_{\beta_5}(t)$ for some $s \in T$ and some $t \in K$, then $s^{-1}x$ fixes $a_1 e_{\gamma_1} + a_2 e_{\gamma_2}$ and sends $a_3 e_{\gamma_3} + a_4 e_{\gamma_4}$ to $a_3 e_{\gamma_3} + a_4 e_{\gamma_4} + t(a_3 e_{\gamma_1} + a_4 e_{\gamma_2})$; as $a_1a_4 \neq a_2a_3$, for $x.y \in Y$ we must have $t = 0$. If instead $x = s x_{\beta_2}(t') x_{\beta_5}(t') n_{\beta_2} n_{\beta_5} x_{\beta_2}(t) x_{\beta_5}(t)$ for some $s \in T$ and some $t, t' \in K$, then $s^{-1}x$ sends $a_1 e_{\gamma_1} + a_2 e_{\gamma_2}$ to $-(a_1 t' e_{\gamma_1} + a_2 t' e_{\gamma_2} + a_1 e_{\gamma_3} + a_2 e_{\gamma_4})$ and $a_3 e_{\gamma_3} + a_4 e_{\gamma_4}$ to
$a_3(1 - tt') e_{\gamma_1} + a_4(1 - tt') e_{\gamma_2} - a_3 e_{\gamma_3} - a_4 e_{\gamma_4}$; again, for $x.y \in Y$ we must have $t, t' = 0$. Thus $G.y \cap Y = T \langle n_{\beta_2} n_{\beta_5} \rangle. \{ y, y' \} = \{ \langle b_1 e_{\gamma_1} + b_2 e_{\gamma_2}, b_3 e_{\gamma_3} + b_4 e_{\gamma_4} \rangle : (\frac{b_2b_3}{b_1b_4})^2 = (\frac{a_2a_3}{a_1a_4})^2 \}$.

Since $\dim C = 9$, we have $\dim(\overline{G.y}) = \dim G - \dim C = 36 - 9 = 27$, while $\dim(\overline{G.y \cap Y}) = 1$; therefore
$$
\dim \G{2}(V) - \dim(\overline{G.y}) = 28 - 27 = 1 \quad \hbox{and} \quad \dim Y - \dim(\overline{G.y \cap Y}) = 2 - 1 = 1.
$$
Hence $y$ is $Y$-exact. Thus the conditions of Lemma~\ref{lem: generic stabilizer from exact subset} hold; so the quadruple $(G, \lambda, p, k)$ has generic stabilizer $C/Z(G) = A_2 T_1.\Z_2$.
\end{proof}

\begin{prop}\label{prop: C_4, omega_4 module, p = 2, k = 2}
Let $G = C_4$ and $\lambda = \omega_4$ with $p = 2$, and take $k = 2$. Then the quadruple $(G, \lambda, p, k)$ has generic stabilizer $\tilde A_2 T_1.\Z_2$.
\end{prop}

\begin{proof}
This is an immediate consequence of Proposition~\ref{prop: D_5, omega_5, B_4, omega_4 modules, k = 2}, using the exceptional isogeny $B_\ell \to C_\ell$ which exists in characteristic $2$.
\end{proof}

\begin{prop}\label{prop: D_5, omega_5 module, k = 3}
Let $G = D_5$ and $\lambda = \omega_5$, and take $k = 3$. Then the quadruple $(G, \lambda, p, k)$ has generic stabilizer ${A_1}^2$.
\end{prop}

\begin{proof}
Again we use the set-up of Proposition~\ref{prop: D_5, omega_5, B_4, omega_4 modules}: we take $H$ to be the simply connected group defined over $K$ of type $E_6$, with simple roots $\beta_1, \dots, \beta_6$; we let $G$ have simple roots $\alpha_1 = \beta_1$, $\alpha_2 = \beta_3$, $\alpha_3 = \beta_4$, $\alpha_4 = \beta_5$, $\alpha_5 = \beta_2$, so that $G = \langle X_\alpha : \alpha = \sum m_i \beta_i, \ m_6 = 0 \rangle < H$; then we may take $V = \langle e_\alpha : \alpha = \sum m_i \beta_i, \ m_6 = 1 \rangle < \L(H)$. We have $Z(G) = \langle z \rangle$ where $z = h_{\beta_1}(-1) h_{\beta_2}(\eta_4) h_{\beta_4}(-1) h_{\beta_5}(-\eta_4)$. Here we take the generalized height function on the weight lattice of $G$ whose value at $\alpha_2$ is $1$, and at each other simple root $\alpha_i$ is $0$; then the generalized height of $\lambda = \frac{1}{2}(\alpha_1 + 2\alpha_2 + 3\alpha_3 + \frac{3}{2}\alpha_4 + \frac{5}{2}\alpha_5)$ is $1$, and as $\lambda$, $\omega_4 = \lambda + \frac{1}{2}\alpha_4 - \frac{1}{2}\alpha_5$ and $\Phi$ generate the weight lattice it follows that the generalized height of any weight is an integer. Since $V_\lambda = \langle e_\delta \rangle$ where $\delta = \esixrt122321$, we see that if $\mu \in \Lambda(V)$ and $e_\alpha \in V_\mu$ where $\alpha = \sum m_i \beta_i$ with $m_6 = 1$, then the generalized height of $\mu$ is $m_3 - 1$. Thus $\Lambda(V)_{[0]} = \{ \nu_1, \dots, \nu_8 \}$, where we write
\begin{eqnarray*}
& \gamma_1 = \esixrt001111, \quad \gamma_2 = \esixrt101111, \quad \gamma_3 = \esixrt011111, \quad \gamma_4 = \esixrt111111, & \\
& \gamma_5 = \esixrt011211, \quad \gamma_6 = \esixrt111211, \quad \gamma_7 = \esixrt011221, \quad \gamma_8 = \esixrt111221, &
\end{eqnarray*}
and for each $i$ we let $\nu_i$ be the weight such that $V_{\nu_i} = \langle e_{\gamma_i} \rangle$. Observe that if we take $s = \prod_{i = 1}^5 h_{\beta_i}(\kappa_i) \in T$, then $\nu_1(s) = \frac{\kappa_3}{\kappa_1\kappa_2}$, $\nu_2(s) = \frac{\kappa_1}{\kappa_2}$, $\nu_3(s) = \frac{\kappa_2\kappa_3}{\kappa_1\kappa_4}$, $\nu_4(s) = \frac{\kappa_1\kappa_2}{\kappa_4}$, $\nu_5(s) = \frac{\kappa_4}{\kappa_1\kappa_5}$, $\nu_6(s) = \frac{\kappa_1\kappa_4}{\kappa_3\kappa_5}$, $\nu_7(s) = \frac{\kappa_5}{\kappa_1}$ and $\nu_8(s) = \frac{\kappa_1\kappa_5}{\kappa_3}$; thus given any $4$-tuple $(n_1, n_2, n_3, n_4)$ of integers we have $c_1\nu_1 + \cdots + c_8\nu_8 = 0$ for $(c_1, \dots, c_8) = (n_1 + n_2, n_3 + n_4, n_1 + n_3 + n_4, n_2, n_4, n_1 + n_2 + n_3, n_2 + n_3, n_1 + n_4)$. In particular, writing \lq $(n_1, n_2, n_3, n_4) \implies (c_1, c_2, c_3, c_4, c_5, c_6, c_7, c_8)$' to indicate this relationship between $4$-tuples and $8$-tuples, we have the following:
\begin{eqnarray*}
(-1, 1, 1, 1) \implies (0, 2, 1, 1, 1, 1, 2, 0), & & \phantom{1} (-1, 2, 1, 1) \implies (1, 2, 1, 2, 1, 2, 3, 0), \\
(-1, 1, 1, 2) \implies (0, 3, 2, 1, 2, 1, 2, 1), & & \phantom{-1} (0, 1, 0, 1) \implies (1, 1, 1, 1, 1, 1, 1, 1).
\end{eqnarray*}
It follows that any subset of $\Lambda(V)_{[0]}$ which contains $\nu_2$, $\nu_3$, $\nu_4$, $\nu_5$, $\nu_6$ and $\nu_7$ has ZLCE. Take $Y = \G{3}(V_{[0]})$, and write
$$
y_0 = \langle e_{\gamma_2} + e_{\gamma_3}, e_{\gamma_4} + e_{\gamma_5}, e_{\gamma_6} + e_{\gamma_7} \rangle \in Y.
$$

We have seen that $W$ acts transitively on the set $\Sigma$ of roots $\alpha$ of $H$ corresponding to the root vectors $e_\alpha$ spanning $V$, and if we write $W_1$ for the stabilizer in $W$ of $\gamma_1$, then $W_1 = \langle w_{\beta_3 + \beta_4 + \beta_5}, w_{\beta_2}, w_{\beta_4}, w_{\beta_1 + \beta_3} \rangle$. Now the stabilizer in $W$ of any $\alpha \in \Sigma$ acts $5$-transitively on the set $\Sigma'$ of roots $\alpha' \in \Sigma$ orthogonal to $\alpha$ (this is evident if we take $\alpha = \esixrt000001$, as then its stabilizer in $W$ is $\langle w_{\beta_1}, w_{\beta_3}, w_{\beta_4}, w_{\beta_2} \rangle$, which acts $5$-transitively on the set of roots $\alpha' = \sum m_i \beta_i$ with $m_5 = 2$ and $m_6 = 1$). Thus if we write $W_2$ for the pointwise stabilizer in $W_1$ of $\{ \gamma_4, \gamma_6, \gamma_8 \}$, then $|W_2| = \frac{|W_1|}{5.4.3} = 2$; we then see that $W_2 = \langle w_\rho \rangle$, where we write $\rho = \beta_1 + \beta_2 + 2\beta_3 + 2\beta_4 + \beta_5$ for the highest root of $\Phi$. Since $w_\rho$ fixes all eight $\gamma_i$ we see that the pointwise stabilizer in $W$ of $\{ \gamma_1, \dots, \gamma_8 \}$ is $W_2$. Now set $W_3 = \langle w_\rho, w_{\beta_1}, w_{\beta_2}, w_{\beta_4}, w_{\beta_5} \rangle$. Then $W_3$ acts transitively on $\{ \gamma_1, \dots, \gamma_8 \}$; the stabilizer in $W_3$ of $\gamma_1$ is $\langle w_\rho, w_{\beta_4}, w_{\beta_5} \rangle$, which acts $3$-transitively on the $\gamma_i$ orthogonal to $\gamma_1$, which are $\gamma_4$, $\gamma_6$ \and $\gamma_8$; of the remaining $\gamma_i$, we see that $\gamma_2$ is orthogonal to none of $\gamma_4$, $\gamma_6$ \and $\gamma_8$, while each of $\gamma_3$, $\gamma_5$ \and $\gamma_7$ is non-orthogonal to a different one of $\gamma_4$, $\gamma_6$ \and $\gamma_8$, so that any element of the setwise stabilizer in $W$ of $\{ \gamma_1, \dots, \gamma_8 \}$ which fixes each of $\gamma_1$, $\gamma_4$, $\gamma_6$ \and $\gamma_8$ must lie in the pointwise stabilizer $W_2$. Hence the setwise stabilizer in $W$ of $\{ \gamma_1, \dots, \gamma_8 \}$, and hence of $\Lambda(V)_{[0]}$, is $\langle w_\rho, w_{\beta_1}, w_{\beta_2}, w_{\beta_4}, w_{\beta_5} \rangle$. Note that this stabilizes $\Phi_{[0]} = \langle \alpha_1, \alpha_3, \alpha_4, \alpha_5 \rangle = \langle \beta_1, \beta_2, \beta_4, \beta_5 \rangle$.

For $t \in K$ write
$$
x(t) = x_{\beta_1}(-t) x_{\beta_4}(2t) x_{\beta_5}(t) x_{\beta_2}(3t) x_{\beta_4 + \beta_5}(-t^2) x_{\beta_2 + \beta_4}(3t^2) x_{\beta_2 + \beta_4 + \beta_5}(4t^3).
$$
A straightforward calculation shows that for $t, t' \in K$ we have $x(t) x(t') = x(t + t')$. Let $A$ be the ${A_1}^2$ subgroup having simple root groups $X_\rho$ and $\{ x(t) : t \in K \}$, with the second $A_1$ factor having maximal torus $T_1 = \{ h_{\beta_1}(\kappa) h_{\beta_2}(\kappa^3) h_{\beta_4}(\kappa^4) h_{\beta_5}(\kappa^3) : \kappa \in K^* \}$ and intersection with $N$ equal to $\langle n_{\beta_1} n_{\beta_4} {n_{\beta_2 + \beta_4 + \beta_5}}^{-1} \rangle T_1$; then $Z(A) = \langle h_{\beta_1}(-1) h_{\beta_2}(-1) h_{\beta_5}(-1) \rangle$. Set $C = Z(G) A$. We find that $C \leq C_G(y_0)$; we shall show that in fact $C_G(y_0) = C$.

We have $U_{[0]} = X_{\alpha_1} X_{\alpha_3} X_{\alpha_4} X_{\alpha_5} X_{\alpha_3 + \alpha_4} X_{\alpha_3 + \alpha_5} X_{\alpha_3 + \alpha_4 + \alpha_5}$. Given $u \in U_{[0]}$, the weights $\nu_2$ and $\nu_3$ occur in $u.(e_{\gamma_2} + e_{\gamma_3})$, and $\nu_4$ and $\nu_5$ occur in $u.(e_{\gamma_4} + e_{\gamma_5})$, while $\nu_6$ and $\nu_7$ occur in $u.(e_{\gamma_6} + e_{\gamma_7})$, so the set of weights occurring in $u.y_0$ contains $\nu_2$, $\nu_3$, $\nu_4$, $\nu_5$, $\nu_6$ and $\nu_7$, and hence by the above has ZLCE. By Lemma~\ref{lem: gen height zero not strictly positive}, we have $C_G(y_0) = C_{U_{[+]}}(y_0) C_{G_{[0]} N_{\Lambda(V)_{[0]}}}(y_0) C_{U_{[+]}}(y_0)$.

First, since $W_{\Lambda(V)_{[0]}} = \langle w_\rho, w_{\beta_1}, w_{\beta_2}, w_{\beta_4}, w_{\beta_5} \rangle$ and $\beta_1, \beta_2, \beta_4, \beta_5 \in \Phi_{[0]}$, we have $G_{[0]} N_{\Lambda(V)_{[0]}} = G_{[0]} \langle n_\rho \rangle$. Since $n_\rho \in C$ it suffices to consider $C_{G_{[0]}}(y_0)$. Moreover, as $G_{[0]}$ is a subsystem subgroup of type $A_1D_3$ with the $A_1$ factor being $\langle X_{\pm\beta_1} \rangle$, to which the second $A_1$ factor of $A$ projects surjectively, we see that any element of $C_{G_{[0]}}(y_0)$ may be written as $g^*c$ where $c \in A$ and $g^* \in \langle T, X_{\pm\beta_2}, X_{\pm\beta_4}, X_{\pm\beta_5} \rangle$. Suppose then that $g^* \in C_G(y_0)$, and let $g^* = u_1nu_2$ be its Bruhat decomposition; write $w = nT$ for the corresponding Weyl group element. Since each weight $\nu_i$ for $2 \leq i \leq 7$ occurs in $u_2.y_0$, we see that $w$ cannot send any of these six weights to $\nu_1$; as $w \in \langle w_{\beta_2}, w_{\beta_4}, w_{\beta_5} \rangle$ this forces $w$ to fix $\nu_1$, so that $w \in \langle w_{\beta_4}, w_{\beta_5} \rangle$. If $w = w_{\beta_4}$ or $w_{\beta_5} w_{\beta_4}$ then $nu_2.(e_{\gamma_4} + e_{\gamma_5})$ would have a term $e_{\gamma_3}$ but no term $e_{\gamma_2}$; if $w = w_{\beta_5}$ then $nu_2.(e_{\gamma_6} + e_{\gamma_7})$ would have a term $e_{\gamma_5}$ but no term $e_{\gamma_4}$, $e_{\gamma_3}$ or $e_{\gamma_2}$; if $w = w_{\beta_4} w_{\beta_5}$ or $w_{\beta_4 + \beta_5}$ then $nu_2.(e_{\gamma_6} + e_{\gamma_7})$ would have a term $e_{\gamma_3}$ but no term $e_{\gamma_2}$ --- therefore $w = 1$. Thus $g^* = u_1h$ for some $h \in T$; so we must have $h.y_0 = y_0$ and $u_1.y_0 = y_0$. A straightforward calculation shows that $h \in C \cap T$. If we write $u_1 = \prod x_\alpha(t_\alpha)$ where the product runs over the positive roots in $\langle \beta_2, \beta_4, \beta_5 \rangle$, then we must have $t_{\beta_5} = t_{\beta_4 + \beta_5} = t_{\beta_2 + \beta_4 + \beta_5} = 0$ as otherwise $u_1.(e_{\gamma_6} + e_{\gamma_7})$, $u_1.(e_{\gamma_4} + e_{\gamma_5})$ or $u_1.(e_{\gamma_2} + e_{\gamma_3})$ would have a term $e_{\gamma_8}$; then we must have $t_{\beta_4} = t_{\beta_2 + \beta_4} = 0$ as otherwise $u_1.(e_{\gamma_4} + e_{\gamma_5})$ or $u_1.(e_{\gamma_2} + e_{\gamma_3})$ would have a term $e_{\gamma_6}$ but no term $e_{\gamma_7}$; finally we must have $t_{\beta_2} = 0$ as otherwise $u_1.(e_{\gamma_2} + e_{\gamma_3})$ would have a term $e_{\gamma_4}$ but no term $e_{\gamma_5}$ --- therefore $u_1 = 1$. Thus $g^* \in C$; so $C_{G_{[0]} N_{\Lambda(V)_{[0]}}}(y_0) = C \cap G_{[0]} N_{\Lambda(V)_{[0]}}$.

Next, let $\Xi = \{ \alpha = \sum m_i \beta_i \in \Phi^+ : m_3 = 1 \}$, and set $U' = \prod_{\alpha \in \Xi} X_\alpha$; then $U_{[+]} = U'.(C \cap U_{[+]})$ and $U' \cap (C \cap U_{[+]}) = \{ 1 \}$. We now observe that if $\alpha \in \Xi$ then $\nu_i + \alpha$ is a weight in $V$ for two values of $i$, which have the same parity and thus never correspond to terms in the same basis vector of $y_0$; moreover if we take a basis vector $e_{\gamma_{i_1}} + e_{\gamma_{i_2}}$ of $y_0$ then each weight in $V$ of positive generalized height is of the form $\nu_i + \alpha$ with $i \in \{ i_1, i_2 \}$ for either one or two such roots $\alpha$ (and for each basis vector there are two weights such that the root $\alpha$ concerned is unique, with these six roots $\alpha$ all being distinct). Thus if we take $u = \prod x_\alpha(t_\alpha) \in U'$ satisfying $u.y_0 = y_0$, and equate coefficients of weight vectors, taking them in an order compatible with increasing generalized height, we immediately see that for six roots $\alpha$ we must have $t_\alpha = 0$, after which it quickly follows that for the remaining roots $\alpha$ we must have $t_\alpha = 0$, so that $u = 1$; so $C_{U_{[+]}}(y_0) = C \cap U_{[+]}$.

Thus $C_{U_{[+]}}(y_0), C_{G_{[0]} N_{\Lambda(V)_{[0]}}}(y_0) \leq C$, so we do indeed have $C_G(y_0) = C$.

Since $\dim(\overline{G.y_0}) = \dim G - \dim C_G(y_0) = 45 - 6 = 39 = \dim \G{3}(V)$, the orbit $G.y_0$ is dense in $\G{3}(V)$. Thus the quadruple $(G, \lambda, p, k)$ has generic stabilizer $C_G(y_0)/Z(G) \cong {A_1}^2$, where the ${A_1}^2$ is a central product.
\end{proof}

\begin{prop}\label{prop: B_4, omega_4 module, k = 3}
Let $G = B_4$ and $\lambda = \omega_4$, and take $k = 3$. Then the quadruple $(G, \lambda, p, k)$ has generic stabilizer $\Z_{2/(p, 2)}.\Z_2$.
\end{prop}

\begin{proof}
We continue with the set-up of Proposition~\ref{prop: D_5, omega_5 module, k = 3}: we take $H$ to be the simply connected group defined over $K$ of type $E_6$, with simple roots $\beta_1, \dots, \beta_6$; we take the $D_5$ subgroup $\langle X_\alpha : \alpha = \sum m_i \beta_i, \ m_6 = 0 \rangle < H$; then we may take $V = \langle e_\alpha : \alpha = \sum m_i \beta_i, \ m_6 = 1 \rangle < \L(H)$. We have $Z(D_5) = \langle z \rangle$ where $z = h_{\beta_1}(-1) h_{\beta_2}(\eta_4) h_{\beta_4}(-1) h_{\beta_5}(-\eta_4)$. We write
\begin{eqnarray*}
& \gamma_1 = \esixrt001111, \quad \gamma_2 = \esixrt101111, \quad \gamma_3 = \esixrt011111, \quad \gamma_4 = \esixrt111111, & \\
& \gamma_5 = \esixrt011211, \quad \gamma_6 = \esixrt111211, \quad \gamma_7 = \esixrt011221, \quad \gamma_8 = \esixrt111221; &
\end{eqnarray*}
here in addition we write
\begin{eqnarray*}
& \delta_1 = \esixrt000011, \quad \delta_2 = \esixrt000111, \quad \delta_3 = \esixrt112221, \quad \delta_4 = \esixrt112321. &
\end{eqnarray*}
Write $\rho = \beta_1 + \beta_2 + 2\beta_3 + 2\beta_4 + \beta_5$ for the highest root in $\Phi(D_5)$. Let $A$ be the ${A_1}^2$ subgroup having simple root groups $X_\rho$ and $\{ x(t) : t \in K \}$, where as before for $t \in K$ we write
$$
x(t) = x_{\beta_1}(-t) x_{\beta_4}(2t) x_{\beta_5}(t) x_{\beta_2}(3t) x_{\beta_4 + \beta_5}(-t^2) x_{\beta_2 + \beta_4}(3t^2) x_{\beta_2 + \beta_4 + \beta_5}(4t^3);
$$
again the second $A_1$ factor has maximal torus $T_1 = \{ h_{\beta_1}(\kappa) h_{\beta_2}(\kappa^3) h_{\beta_4}(\kappa^4) h_{\beta_5}(\kappa^3) : \kappa \in K^* \}$ and intersection with $N$ equal to $\langle n_{\beta_1} n_{\beta_4} {n_{\beta_2 + \beta_4 + \beta_5}}^{-1} \rangle T_1$. Let $B$ be the Borel subgroup $\{ h_\rho(\kappa) : \kappa \in K^* \} T_1 X_\rho \{ x(t) : t \in K \}$ of $A$. Write $n_0 = {n_\rho}^{-1} n_{\beta_1} n_{\beta_4} {n_{\beta_2 + \beta_4 + \beta_5}}^{-1} \in A \cap N$.

We saw in the proof of Proposition~\ref{prop: D_5, omega_5 module, k = 3} that if we set
$$
y_0 = \langle e_{\gamma_2} + e_{\gamma_3}, e_{\gamma_4} + e_{\gamma_5}, e_{\gamma_6} + e_{\gamma_7} \rangle,
$$
then the $D_5$-orbit containing $y_0$ is dense in $\G{3}(V)$, and $C_{D_5}(y_0) = Z(D_5)A$. Given $\a = (a_1, a_2, a_3, a_4, a_5) \in K^5$, write $|\a| = {a_1}^2 + {a_3}^2 + {a_4}^2$. For $\a \in K^5$ with $|\a| = 1$ and $a_2 = 0 \neq a_5$, set
\begin{eqnarray*}
y_\a & = & \langle a_5 e_{\gamma_2} + e_{\gamma_3} - a_1 e_{\delta_3} + a_3 e_{\gamma_7} + a_4 e_{\gamma_1}, \\
     &   & \phantom{\langle} e_{\gamma_4} + e_{\gamma_5} - a_1 e_{\delta_1} - a_1 e_{\delta_4} + a_3 e_{\gamma_1} + a_3 e_{\gamma_8} + a_4 e_{\gamma_2} - a_4 e_{\gamma_7}, \\
     &   & \phantom{\langle} e_{\gamma_6} + a_5 e_{\gamma_7} - a_1 e_{\delta_2} + a_3 e_{\gamma_2} - a_4 e_{\gamma_8} \rangle.
\end{eqnarray*}
Considering coefficients of $e_{\gamma_i}$ for $i = 3, 4, 5, 6$ quickly shows that distinct such vectors $\a$ give distinct elements $y_\a$. Thus if we set
$$
Y = \{ y_\a : \a \in K^5, \ |\a| = 1, \ a_2 = 0 \neq a_5 \}
$$
then $\dim Y = 3$. Choose $\xi \in K^*$ with $\xi^2 = {a_5}^{-1}$ and define
\begin{eqnarray*}
g_\a & = & h_{\beta_5}(\xi) h_{-\beta_2}(\xi) x_{\beta_5}(a_4) x_{-\beta_2}(-a_4) x_{\beta_4 + \beta_5}(a_3) x_{-(\beta_2 + \beta_4)}(a_3) \\
     &   & {} \times x_{\beta_1 + \beta_3 + \beta_4 + \beta_5}(a_1) x_{-(\beta_1 + \beta_2 + \beta_3 + \beta_4)}(a_1)
\end{eqnarray*}
(note that the two choices for $\xi$ give elements differing by $z^2$, which fixes all points in $\G{3}(V)$); then calculation shows that
$$
g_\a.y_\a = y_0.
$$

At this point we find it convenient to switch notation. Instead of taking the root system of $D_5$ to be a subsystem of that of $E_6$, we shall use the standard notation given in Section~\ref{sect: notation}; thus we replace $\beta_1$, $\beta_3$, $\beta_4$, $\beta_5$ and $\beta_2$ by $\ve_1 - \ve_2$, $\ve_2 - \ve_3$, $\ve_3 - \ve_4$, $\ve_4 - \ve_5$ and $\ve_4 + \ve_5$ respectively, and we recall the natural module $V_{nat}$ for $D_5$. However, there is an unfortunate consequence to this change: in Section~\ref{sect: notation} we defined the action of root elements on $V_{nat}$, which implicitly determined the structure constants, and these are not the same as those given in the appendix of \cite{LSmax}, which we have been using until now. For this reason we shall largely avoid all mention of root elements from now on, but rather identify elements of $D_5$ by their action on $V_{nat}$ (the kernel of this action is $\langle z^2 \rangle$, so this is harmless). Thus with respect to the ordered basis $v_1, v_2, v_3, v_4, v_5, v_{-5}, v_{-4}, v_{-3}, v_{-2}, v_{-1}$ of $V_{nat}$, the element $g_\a$ defined above acts as
$$
\left(
  \begin{array}{ccccc|ccccc}
      1    &   &         &         &       a_1      &               &         &         &   &         \\
           & 1 &         &         &                &               &         &         &   &         \\
           &   &    1    &         &       a_3      &               &         &         &   &         \\
           &   &         &    1    &       a_4      &               &         &         &   &         \\
           &   &         &         &       a_5      &               &         &         &   &         \\
  \hline
   -\aaf15 &   & -\aaf35 & -\aaf45 & -\frac{1}{a_5} & \frac{1}{a_5} & -\aaf45 & -\aaf35 &   & -\aaf15 \\
           &   &         &         &       a_4      &               &    1    &         &   &         \\
           &   &         &         &       a_3      &               &         &    1    &   &         \\
           &   &         &         &                &               &         &         & 1 &         \\
           &   &         &         &       a_1      &               &         &         &   &    1    \\
  \end{array}
\right)
$$
The one exception to this is that we shall write a positive root element of the first $A_1$ factor of $A$ as $x_\rho(t)$; since $\rho = \ve_1 + \ve_2$, this acts on $V_{nat}$ by sending $v_{-2} \mapsto v_{-2} + tv_1$ and $v_{-1} \mapsto v_{-1} - tv_2$ and fixing all other basis vectors.

Write $V_{1, 2} = \langle v_1, v_2, v_{-2}, v_{-1} \rangle$ and $V_{3, 4, 5} = \langle v_3, v_4, v_5, v_{-5}, v_{-4}, v_{-3} \rangle$; then we have $V_{nat} = V_{1, 2} \oplus V_{3, 4, 5}$, and in the calculations which follow we will always take the basis elements of these two subspaces in the order given here. We see that $A < D_2 D_3$ where $D_2$ and $D_3$ act on $V_{1, 2}$ and $V_{3, 4, 5}$ respectively; indeed $\langle X_{\pm\rho} \rangle$ lies in $D_2$ while the second $A_1$ factor of $A$ projects non-trivially on both $D_2$ and $D_3$. In this second factor write $h(\kappa) = h_{\beta_1}(\kappa) h_{\beta_2}(\kappa^3) h_{\beta_4}(\kappa^4) h_{\beta_5}(\kappa^3)$ for $\kappa \in K^*$, and $n = n_{\beta_1} n_{\beta_4} {n_{\beta_2 + \beta_4 + \beta_5}}^{-1} = n_0 n_\rho$. We find that on $V_{1, 2}$ and $V_{3, 4, 5}$ respectively, $x(t)$ acts as
$$
\left(
  \begin{array}{cc|cc}
    1 & -t &   &   \\
      &  1 &   &   \\
  \hline
      &    & 1 & t \\
      &    &   & 1 \\
  \end{array}
\right)
\quad\hbox{and}\quad
\left(
  \begin{array}{ccc|ccc}
    1 & 2t & t^2 & 3t^2 & -2t^3 &  t^4 \\
      &  1 &  t  &  3t  & -3t^2 & 2t^3 \\
      &    &  1  &      &  -3t  & 3t^2 \\
  \hline
      &    &     &   1  &   -t  &  t^2 \\
      &    &     &      &   1   &  -2t \\
      &    &     &      &       &   1  \\
  \end{array}
\right),
$$
while $h(\kappa)$ acts as $\diag(\kappa, \kappa^{-1}, \kappa, \kappa^{-1})$ and $\diag(\kappa^4, \kappa^2, 1, 1, \kappa^{-2}, \kappa^{-4})$.

Now write
$$
v^\diamondsuit = v_5 + v_{-5},
$$
and let $G = C_{D_5}(v^\diamondsuit) = B_4$; then $Z(G) = \langle z^2 \rangle$. Since the element $n_0$ defined above fixes both $v_5$ and $v_{-5}$, while for $i \in \{ 1, 2, 3, 4 \}$ we have $n_0.v_i = v_{-i}$ and $n_0.v_{-i} = v_i$, we see that $n_0 \in G$. For $\a \in K^5$ with $|\a| = 1$ and $a_2 = 0 \neq a_5$ as above, define
$$
v_\a = g_\a.v^\diamondsuit = a_1(v_1 + v_{-1}) + a_3(v_3 + v_{-3}) + a_4(v_4 + v_{-4}) + a_5 v_5;
$$
then $v_\a$ is a vector of norm $1$ fixed by $n_0$. Write
$$
V_* = \{ v_\a : |\a| = 1, a_2 = 0 \}.
$$
Define
\begin{eqnarray*}
S & = & \{ \a \in K^5 : |\a| = 1, \ a_2 = 0, \ a_1a_3a_4a_5 \neq 0, \ {a_3}^2 + {a_4}^2 \neq 0, \\
  &   & \phantom{\{ \a \in K^5 :} 2{a_3}^2 + a_3a_5 - {a_4}^2 \neq 0, \ 4{a_3}^2 - 4a_3a_5 + {a_5}^2 + 16{a_4}^2 \neq 0, \\
  &   & \phantom{\{ \a \in K^5 :} 2{a_3}^2 - a_3a_5 + 2{a_4}^2 \neq 0 \}
\end{eqnarray*}
and set
$$
\hat Y = \{ y_\a \in Y : \a \in S \}, \quad \hat V_* = \{ v_\a \in V_* : \a \in S \};
$$
then $\hat Y$ and $\hat V_*$ are dense open subsets of $Y$ and $V_*$ respectively.

Take $y_\a \in \hat Y$ and suppose $g \in \Tran_G(y_\a, Y)$; write $g.y_\a = y_{\a'}$ and set $g' = g_{\a'} g {g_\a}^{-1} \in D_5$. Then $g'.y_0 = g_{\a'} g {g_\a}^{-1}.y_0 = g_{\a'} g.y_\a = g_{\a'}.y_{\a'} = y_0$, and $g'.v_\a = g_{\a'} g {g_\a}^{-1}.v_\a = g_{\a'} g.v^\diamondsuit = g_{\a'}.v^\diamondsuit = v_{\a'}$ since $g \in G$; so any element of $\Tran_G(y_\a, Y)$ is of the form ${g_{\a'}}^{-1} g' g_\a$, where $g' \in C_{D_5}(y_0) = Z(D_5)A$ and $g'.v_\a = v_{\a'} \in V_*$. In particular, taking $\a' = \a$ we see that $C_G(y_\a) = C_{Z(D_5)A}(v_\a)^{g_\a}$. We shall show that there is a dense open subset $S'$ of $S$ such that if $\a \in S'$ then $\Tran_{Z(D_5)A}(v_\a, V_*)$ is finite, as is then $\Tran_G(y_\a, Y)$, and we shall identify $C_{Z(D_5)A}(v_\a)$.

Take $g' \in \Tran_{Z(D_5)A}(v_\a, V_*)$, so that $g'.v_\a = v_{\a'}$ for some $\a'$ with $|\a'| = 1$ and ${a_2}' = 0$; thus the coefficients in $g'.v_\a$ of $v_2$, $v_{-2}$ and $v_{-5}$ must all be zero, while for $i \in \{ 1, 3, 4 \}$ those of $v_i$ and $v_{-i}$ must be equal to each other. We have $g' = z^i {g_1}' {g_2}'$, where $i \in \{ 0, 1, 2, 3 \}$, ${g_1}' \in \langle X_{\pm\rho} \rangle$ and ${g_2}'$ lies in the second $A_1$ factor. Write $\e = (-1)^i$, so that $z^i$ acts on $V_{nat}$ as multiplication by $\e$. According as ${g_1}'$ lies in the Borel subgroup $B$ or not we have ${g_1}' = h_\rho(\kappa_1) x_\rho(t_1)$ or $x_\rho({t_1}') n_\rho h_\rho(\kappa_1) x_\rho(t_1)$, where $t_1, {t_1}' \in K$ and $\kappa_1 \in K^*$; likewise according as ${g_2}'$ lies in $B$ or not we have ${g_2}' = h(\kappa_2) x(t_2)$ or $x({t_2}') n h(\kappa_2) x(t_2)$, where $t_2, {t_2}' \in K$ and $\kappa_2 \in K^*$.

First suppose ${g_1}', {g_2}' \in B$. From the coefficients of $v_2$ and $v_{-2}$ we immediately see that $t_1 = t_2 = 0$; those of $v_1$ and $v_{-1}$ give ${\kappa_1}^2 {\kappa_2}^2 = 1$, and those of $v_4$ and $v_{-4}$ give ${\kappa_2}^4 = 1$. Thus there are finitely many such elements $g'$. Moreover if $\a' = \a$, the coefficient of $v_5$ shows that $\e = 1$, and then those of $v_1$ and $v_4$ that $\kappa_1\kappa_2 = {\kappa_2}^2 = 1$; so $\kappa_1 = \kappa_2 = \pm1$, and $g' = 1$ or $z^2$.

Next suppose ${g_2}' \in B$ but ${g_1}' \notin B$. Here the coefficients of $v_2$, $v_{-2}$ and $v_{-5}$ give ${\kappa_1}^2 t_1 {t_1}' = 1$, $t_1t_2 = -1$ and $a_3{t_2}^2 - a_4 t_2 = 0$ respectively; since the second of these implies that $t_2 \neq 0$, the third gives $t_2 = \frac{a_4}{a_3}$, and then the second again gives $t_1 = -\frac{a_3}{a_4}$. Now the coefficients of $v_4$ and $v_{-4}$ give ${\kappa_2}^2[2a_3 {t_2}^3 - a_4(3{t_2}^2 - 1)+ a_5 t_2] = {\kappa_2}^{-2}[-2a_3 t_2 + a_4]$, which reduces to ${\kappa_2}^4({a_3}^2 + a_3 a_5 - {a_4}^2) = -{a_3}^2$, so ${\kappa_2}^4$ is determined, as is thus $\kappa_2$ up to a power of $\eta_4$; the coefficients of $v_1$ and $v_{-1}$ give ${\kappa_1}^2 {t_1}^2 = {\kappa_2}^2$, so $\kappa_1$ is determined up to a sign; finally the first equation given then determines ${t_1}'$. Thus there are finitely many such elements $g'$. Moreover if $\a' = \a$, the coefficient of $v_5$ gives $\e(3a_3{t_2}^2 -3a_4 t_2 + a_5) = a_5$, which reduces to $\e = 1$; now the coefficient of $v_{-3}$ gives ${\kappa_2}^{-4} a_3 = a_3$, so ${\kappa_2}^4 = 1$; but then the equation above determining ${\kappa_2}^4$ gives $2{a_3}^2 + a_3 a_5 - {a_4}^2 = 0$, contrary to the definition of the set $S$. Therefore no elements of this type fix $v_\a$.

Now suppose ${g_1}' \in B$ but ${g_2}' \notin B$. Here the coefficients of $v_2$ and $v_{-2}$ give $t_1 t_2 = -1$ and ${\kappa_2}^2 t_2 {t_2}' = 1$ respectively; using the second of these, the coefficient of $v_{-5}$ gives ${\kappa_2}^4 a_3 {t_2}'{}^2 + {\kappa_2}^2 a_4 {t_2}' = 0$ and hence $a_3 + a_4 t_2 = 0$, so that $t_2 = -\frac{a_3}{a_4}$, whence $t_1 = \frac{a_4}{a_3}$ and ${t_2}' = -\frac{a_4}{{\kappa_2}^2 a_3}$. Now the coefficients of $v_4$ and $v_{-4}$ give $2{\kappa_2}^2 a_3 {t_2}' ({t_2}'{}^2 + 1) + a_4 ({t_2}'{}^2 - {t_2}^2 + 3) + a_5 t_2 = 0$; substituting for $t_2$ and ${t_2}'$ gives $\frac{1}{{\kappa_2}^4} = \frac{{a_3}^2}{{a_4}^4}({a_4}^2 - {a_3}^2 - a_3 a_5)$, so ${\kappa_2}^4$ is determined, as is thus $\kappa_2$ up to a power of $\eta_4$; now ${t_2}'$ is determined, and finally the coefficients of $v_1$ and $v_{-1}$ give ${\kappa_1}^2 = {\kappa_2}^2 {t_2}^2$, so $\kappa_1$ is determined up to a sign. Thus there are finitely many such elements $g'$. Moreover if $\a' = \a$, the coefficient of $v_5$ gives $\e(3{\kappa_2}^4 a_3 {t_2}'{}^2 + 3{\kappa_2}^2 a_4 {t_2}' + a_5) = a_5$, which reduces to $\e = 1$; now the coefficient of $v_3$ gives ${\kappa_2}^4 a_3 {t_2}'{}^4 = a_3$, so $\frac{{a_4}^4}{{\kappa_2}^4 {a_3}^4} = 1$ and thus ${\kappa_2}^4 = \frac{{a_4}^4}{{a_3}^4}$; but then the equation above determining ${\kappa_2}^4$ gives $2{a_3}^2 + a_3 a_5 - {a_4}^2 = 0$, contrary to the definition of the set $S$. Therefore no elements of this type fix $v_\a$.

Finally suppose ${g_1}', {g_2}' \notin B$; this is the most complicated case. Here the coefficients of $v_2$ and $v_{-2}$ give ${t_1}'(t_1 t_2 + 1) = \frac{1}{{\kappa_1}^2} t_2$ and ${t_2}'(t_1 t_2 + 1) = \frac{1}{{\kappa_2}^2} t_1$; the fact that elements of $D_2$ preserve norms of vectors in $V_{1, 2}$ gives $\kappa_1 \kappa_2 (t_1 t_2 + 1) = \e' \in \{ \pm1 \}$, and thus ${t_1}' = \e' \frac{\kappa_2}{\kappa_1} t_2$ and ${t_2}' = \e' \frac{\kappa_1}{\kappa_2} t_1$; it follows that $1 - {\kappa_2}^2 t_2 {t_2}' = \e' \kappa_1 \kappa_2$. Using these equations, the coefficients of $v_3$ and $v_{-3}$ give
$$
{\kappa_1}^4[a_3({t_1}^4 + 1) - 2a_4t_1({t_1}^2 - 1) + a_5{t_1}^2] = {\kappa_2}^4[a_3({t_2}^4 + 1) - 2a_4t_2({t_2}^2 - 1) + a_5{t_2}^2],
$$
those of $v_4$ and $v_{-4}$ give
\begin{eqnarray*}
&   & {\kappa_1}^2[2a_3({t_1}^3 - t_2) + a_4({t_1}^3 t_2 - 3{t_1}^2 - 3t_1 t_2 + 1) - a_5 t_1(t_1t_2 - 1)] \\
& = & {\kappa_2}^2[2a_3({t_2}^3 - t_1) + a_4(t_1 {t_2}^3 - 3{t_2}^2 - 3t_1 t_2 + 1) - a_5 t_2(t_1t_2 - 1)],
\end{eqnarray*}
and that of $v_{-5}$ gives $f_1(t_1, t_2) = 0$, where
$$
f_1(t_1, t_2) = a_3({t_1}^2 + {t_2}^2) + a_4(t_1 + t_2)(t_1 t_2 - 1) - a_5 t_1 t_2.
$$
Squaring the second of these and using the first to eliminate $\frac{{\kappa_1}^4}{{\kappa_2}^4}$ produces an equation which eventually simplifies to $(2{a_3}^2 + a_3 a_5 - {a_4}^2)f_2(t_1, t_2) = 0$, where writing $t_{j_1, j_2} = {t_1}^{j_1} {t_2}^{j_2} - {t_1}^{j_2} {t_2}^{j_1}$ for convenience we have
\begin{eqnarray*}
f_2(t_1, t_2) & = & (2a_3 - a_5)[t_{6, 4} + 2 t_{5, 3} - 2 t_{3, 1} - t_{2, 0}] \\
              &   & {} + 2a_4[t_{6, 5} - t_{6, 3} - 3 t_{5, 2} - 5 t_{4, 3} - 3 t_{4, 1} - 5 t_{3, 2} - t_{3, 0} + t_{1, 0}].
\end{eqnarray*}
As $2{a_3}^2 + a_3 a_5 - {a_4}^2 \neq 0$ by definition of the set $S$, we must have $f_2(t_1, t_2) = 0$. We may write $f_1(t_1, t_2) = \sum_{j = 0}^2 P_j(t_1) {t_2}^j$ and $f_2(t_1, t_2) = \sum_{j = 0}^6 Q_j(t_1) {t_2}^j$, where the various $P_j$ and $Q_j$ are polynomials; multiplying $f_2(t_1, t_2)$ by $P_2(t_1)^5$ and repeatedly replacing $P_2(t_1) {t_2}^2$ by $-P_1(t_1) t_2 - P_0(t_1)$ gives the equation $R_1(t_1) t_2 + R_0(t_1) = 0$, where
\begin{eqnarray*}
R_1 & = & (-3 {P_2}^2 P_1 {P_0}^2 + 4 P_2 {P_1}^3 P_0 -{P_1}^5)Q_6 + ({P_2}^3 {P_0}^2 - 3 {P_2}^2 {P_1}^2 P_0 + P_2 {P_1}^4)Q_5 \\
    &   & {} + (2 {P_2}^3 P_1 P_0 - {P_2}^2 {P_1}^3)Q_4 + (-{P_2}^4 P_0 + {P_2}^3 {P_1}^2) Q_3 - {P_2}^4 P_1 Q_2 + {P_2}^5 Q_0, \\
R_0 & = & (-{P_2}^2 {P_0}^3 + 3 P_2 {P_1}^2 {P_0}^2 - {P_1}^4 P_0)Q_6 + (-2 {P_2}^2 P_1 {P_0}^2 + P_2 {P_1}^3 P_0)Q_5 \\
    &   & {} + ({P_2}^3 {P_0}^2 - {P_2}^2 {P_1}^2 P_0)Q_4 + {P_2}^3 P_1 P_0 Q_3 - {P_2}^4 P_0 Q_2 + {P_2}^5 Q_0;
\end{eqnarray*}
now multiplying $f_1(t_1, t_2)$ by $R_1(t_1)^2$ and replacing $R_1(t_1) t_2$ by $-R_0(t_1)$ gives the polynomial equation $f_3(t_1) = 0$, where
$$
f_3 = P_2 {R_0}^2 - P_1 R_0 R_1 + P_0 {R_1}^2.
$$
The coefficient in each term of $f_3$ is a polynomial in $a_3$, $a_4$ and $a_5$ (indeed, a homogeneous polynomial of degree $13$, since each $P_j$ and $Q_j$ is homogeneous of degree $1$). If we view each such coefficient as a polynomial in $a_4$, we find that the coefficient of $t^4$ has constant term ${a_3}^9(2a_3 - a_5)^2(4{a_3}^2 - 3{a_5}^2)$. Since this is not identically zero, there is a dense open subset $S'$ of $S$ where the coefficient of $t^4$ in $f_3$ is non-zero, and thus $f_3$ is not the zero polynomial; so if $\a \in S'$ then $t_1$ is a root of a non-zero polynomial, and hence there are only finitely many possibilities for $t_1$. Since interchanging $t_1$ and $t_2$ fixes $f_1$ and negates $f_2$, there are also only finitely many possibilities for $t_2$. As elements of $D_3$ preserve the norms of vectors in $V_{3, 4, 5}$, and ${a_3}^2 + {a_4}^2 \neq 0$ by the definition of the set $S$, we see that in the equations above obtained from coefficients of $v_j$ and $v_{-j}$ for $j \in \{ 3, 4 \}$ we cannot have both sides of both equations being zero; thus the value of $\frac{{\kappa_1}^4}{{\kappa_2}^4}$ is determined, as is thus $\frac{\kappa_1}{\kappa_2}$ up to a power of $\eta_4$; using $\kappa_1 \kappa_2 (t_1 t_2 + 1) = \e'$ we see that there are finitely many possibilities for each of $\kappa_1$ and $\kappa_2$, and then ${t_1}'$ and ${t_2}'$ are both determined. Thus once more there are finitely many elements $g'$. At this point we have indeed proved that if $\a \in S'$ then $\Tran_{Z(D_5)A}(v_\a, V_*)$ is finite, as is then $\Tran_G(y_\a, Y)$. Thus
$$
\codim {\ts\Tran_G(y_\a, Y)} = \dim G - \dim {\ts\Tran_G(y_\a, Y)} = 36 - 0 = 36
$$
while
$$
\codim Y = \dim \G{2}(V) - \dim Y = 39 - 3 = 36.
$$
Therefore $y_\a$ is $Y$-exact.

Now suppose in this final case that $\a' = \a$. Here it is convenient to note that $g'$ acts as $z^i{g_2}' = z^i x({t_2}') n h(\kappa_2) x(t_2)$ on $V_{3, 4, 5}$. If temporarily we write $v$ for the projection of $v_\a$ on $V_{3, 4, 5}$, we have
$$
g'.v = v \iff h(\kappa_2) x(t_2).v = n^{-1} x({t_2}')^{-1}z^{-i}.v = nx(-{t_2}').\e v;
$$
subtracting $\e$ times the matrix representing $nx(-{t_2}')$ on $V_{3, 4, 5}$ from that representing $h(\kappa_2) x(t_2)$ gives
$$
\left(
  \begin{array}{ccc|ccc}
    {\kappa_2}^4 & 2{\kappa_2}^2 t_2 & {\kappa_2}^4{t_2}^2 & 3{\kappa_2}^4{t_2}^2 &       -2{\kappa_2}^4{t_2}^3     &        {\kappa_2}^4{t_2}^4 - \e        \\
                 &    {\kappa_2}^2   &   {\kappa_2}^2 t_2  &   3{\kappa_2}^2 t_2  &   -3{\kappa_2}^2{t_2}^2 - \e    &    2{\kappa_2}^2{t_2}^3 - 2\e{t_2}'    \\
                 &                   &       1 - \e        &                      &         -3t_2 - 3\e{t_2}'       &         3{t_2}^2 - 3\e{t_2}'{}^2       \\
  \hline
                 &                   &                     &        1 - \e        &          -t_2 - \e{t_2}'        &          {t_2}^2 - \e{t_2}'{}^2        \\
                 &        -\e        &      \e{t_2}'       &       3\e{t_2}'      & {\kappa_2}^{-2} + 3\e{t_2}'{}^2 &  -2{\kappa_2}^{-2} t_2 + 2\e{t_2}'{}^3 \\
        -\e      &     2\e{t_2}'     &    -\e{t_2}'{}^2    &    -3\e{t_2}'{}^2    &          -2\e{t_2}'{}^3         &     {\kappa_2}^{-4} - \e{t_2}'{}^4     \\
  \end{array}
\right),
$$
which then must send $v$ to the zero vector in $V_{3, 4, 5}$. From the coefficients of $v_5$ and $v_{-5}$ we immediately see that we must have $\e = 1$ (and considering the action on $V_{1, 2}$ now gives $\e' = -1$, so that ${t_1}' = -\frac{\kappa_2}{\kappa_1}t_2$ and $t_1 = -\frac{\kappa_2}{\kappa_1}{t_2}'$, and $1 + \kappa_1 \kappa_2 = {\kappa_2}^2 t_2{t_2}'$); thus that of $v_{-5}$ gives $0 = a_3({t_2}^2 - {t_2}'{}^2) - a_4(t_2 + {t_2}') = (t_2 + {t_2}')[a_3(t_2 - {t_2}') - a_4]$. Suppose if possible that $t_2 + {t_2}' \neq 0$; then $a_3(t_2 - {t_2}') - a_4 = 0$, so $t_2 - {t_2}' = \frac{a_4}{a_3}$. From the coefficients of $v_3$ and $v_{-3}$ we have $a_3({t_2}^4 + 1 - {\kappa_2}^{-4}) - 2a_4t_2({t_2}^2 - 1) + a_5{t_2}^2 = 0$ and $a_3({t_2}'{}^4 + 1 - {\kappa_2}^{-4}) + 2a_4{t_2}'({t_2}'{}^2 - 1) + a_5{t_2}'{}^2 = 0$; subtracting and dividing by $t_2 + {t_2}'$ gives $a_3(t_2 - {t_2}')({t_2}^2 + {t_2}'{}^2) - 2a_4({t_2}^2 - t_2{t_2}' + {t_2}'{}^2 - 1) + a_5(t_2 - {t_2}') = 0$, which on substituting for $t_2 - {t_2}'$ reduces to $2{a_3}^2 + a_3a_5 - {a_4}^2 = 0$, contrary to the definition of the set $S$. Thus we must have ${t_2}' = -t_2$. We therefore have
\begin{eqnarray*}
a_3({t_2}^4 + 1) - 2a_4t_2({t_2}^2 - 1) + a_5{t_2}^2 & = & {\kappa_2}^{-4}a_3, \\
2a_3{t_2}^3 - a_4(3{t_2}^2 - 1) + a_5 t_2            & = & {\kappa_2}^{-2}(-2a_3t_2 + a_4),
\end{eqnarray*}
the second of these equations coming from the coefficient of $v_4$. Squaring the second and using the first to eliminate ${\kappa_2}^{-4}$ gives an equation which reduces to $(2{a_3}^2 + a_3 a_5 - {a_4}^2)t_2(2a_4{t_2}^2 + (2a_3 - a_5)t_2 - 2a_4) = 0$; so
$$
t_2(2a_4{t_2}^2 + (2a_3 - a_5)t_2 - 2a_4) = 0.
$$
Adding this to the first of the two displayed equations above produces $a_3({t_2}^2 + 1)^2 = {\kappa_2}^{-4}a_3$; so ${\kappa_2}^{-2} = \e''({t_2}^2 + 1)$ where $\e'' \in \{ \pm1 \}$. Substituting for ${\kappa_2}^{-2}$ in the second of the two displayed equations gives
$$
2(1 + \e'')a_3{t_2}^3 - (3 + \e'')a_4{t_2}^2 + (2\e'' a_3 + a_5)t_2 + (1 - \e'')a_4 = 0.
$$

If $p = 2$ this last equation is simply $a_5 t_2 = 0$, so we have $t_2 = 0$, from which it immediately follows that ${t_2}' = 0 = t_1 = {t_1}'$ and $\kappa_2 = 1 = \kappa_1$; thus $g' = n_0$ (and of course $z = 1$ in this case), so $C_{Z(D_5)A}(v_\a) = \langle n_0 \rangle$. Now assume $p \neq 2$.

First suppose $\e'' = 1$; then we have $4a_3{t_2}^3 - 4a_4{t_2}^2 + (2a_3 + a_5)t_2 = 0$. If $t_2 \neq 0$ then we have $4a_3{t_2}^2 - 4a_4 t_2 + (2a_3 + a_5) = 0$ and $2a_4{t_2}^2 + (2a_3 - a_5)t_2 - 2a_4 = 0$; multiplying the first of these by $a_4$ and the second by $2a_3$, and subtracting, gives $(-4{a_3}^2 + 2a_3 a_5 - 4{a_4}^2)t_2 + (6a_3 a_4 + a_4 a_5) = 0$, so that $t_2 = \frac{a_4(6a_3 + a_5)}{2(2{a_3}^2 - a_3 a_5 + 2{a_4}^2)}$; substituting in the second of the two equations above and clearing denominators produces $0 = 2(2{a_3}^2 + a_3 a_5 - {a_4}^2)(4{a_3}^2 - 4a_3 a_5 + {a_5}^2 + 16{a_4}^2)$, contrary to the definition of the set $S$. Thus we must have $t_2 = 0$, whence ${t_2}' = 0 = t_1 = {t_1}'$ and $\kappa_2 = \pm1 = -\kappa_1$; so $g' = n_0$ or $z^2 n_0$.

Now suppose $\e'' = -1$; then we have $-2a_4{t_2}^2 - (2a_3 - a_5)t_2 + 2a_4 = 0$, and $-1 - \frac{1}{{\kappa_2}^2} = {t_2}^2 = -\frac{\kappa_1}{\kappa_2} - \frac{1}{{\kappa_2}^2}$, so that $\kappa_1 = \kappa_2$ and hence ${t_1}' = -t_2$ and $t_1 = -{t_2}'$. Writing simply $t$ for $t_2$, we have the element which acts on $V_{1, 2}$ as
$$
\frac{1}{t^2 + 1}
\left(
  \begin{array}{cc|cc}
    t^2 &   t  &  -t  &  1  \\
     t  & -t^2 &  -1  & -t  \\
  \hline
    -t  &  -1  & -t^2 &  t  \\
     1  &  -t  &   t  & t^2 \\
  \end{array}
\right)
$$
and on $V_{3, 4, 5}$ as
$$
\frac{1}{(t^2 + 1)^2}
\left(
  \begin{array}{ccc|ccc}
      t^4 &     -2t^3   &       t^2     &      3t^2     &      2t     &   1   \\
    -2t^3 & -t^4 + 3t^2 &    t^3 - t    &   3t^3 - 3t   &   3t^2 - 1  &  2t   \\
     3t^2 &  3t^3 - 3t  & t^4 - t^2 + 1 &     -9t^2     &  3t^3 - 3t  &   t^2 \\
  \hline
      t^2 &   t^3 - t   &      -t^2     & t^4 - t^2 + 1 &   t^3 - t   &  3t^2 \\
     2t   &   3t^2 - 1  &    t^3 - t    &   3t^3 - 3t   & -t^4 + 3t^2 & -2t^3 \\
      1   &      2t     &       t^2     &      3t^2     &     -2t^3   &   t^4 \\
  \end{array}
\right).
$$
The quadratic satisfied by $t$ has discriminant $4{a_3}^2 - 4a_3 a_5 + {a_5}^2 + 16{a_4}^2 \neq 0$, so it has distinct roots; the two choices for the root give two such elements, both being involutions commuting with, and conjugate to, $n_0$. Call one of them $x_\a$; then the other is $n_0 x_\a$. Thus here $g' = x_\a$, $z^2 x_\a$, $n_0 x_\a$ or $z^2 n_0 x_\a$. Therefore in the case $p \neq 2$ we have $C_{Z(D_5)A}(v_\a) = \langle z^2, n_0, x_\a \rangle$.

Thus according as $p = 2$ or $p \neq 2$ we have $C_G(y_\a) = C_{Z(D_5)A}(v_\a)^{g_\a} = \langle n_0 \rangle$ or $\langle z^2, n_0, {x_\a}^{g_\a} \rangle$ (since $g_\a$ commutes with $n_0$, and of course with $z$). If $p = 2$ we need say no more, so assume $p \neq 2$. Here a simple check shows that there is a single conjugacy class in $G$ of involutions lying in the same $D_5$-class as $n_0$, and it contains $h_0 = h_\rho(-1)$; the centralizer of this in $G$ is of type $D_2 B_2$, in which there is a single conjugacy class of involutions $x$ such that both $x$ and $x h_0$ lie in the same $G$-class as $h_0$. Since $n_0$ is such an element $x$, we see that if we set $C = \langle z^2, h_0, n_0 \rangle$ then $C_G(y_\a)$ is $G$-conjugate to $C$ (and this last statement is also true for $p = 2$). Thus the conditions of Lemma~\ref{lem: generic stabilizer from exact subset} hold; so the quadruple $(G, \lambda, p, k)$ has generic stabilizer $C/Z(G) \cong \Z_{2/(p, 2)}.\Z_2$.
\end{proof}

\begin{prop}\label{prop: C_4, omega_4 module, p = 2, k = 3}
Let $G = C_4$ and $\lambda = \omega_4$ with $p = 2$, and take $k = 3$. Then the quadruple $(G, \lambda, p, k)$ has generic stabilizer $\Z_2$.
\end{prop}

\begin{proof}
This is an immediate consequence of Proposition~\ref{prop: B_4, omega_4 module, k = 3}, using the exceptional isogeny $B_\ell \to C_\ell$ which exists in characteristic $2$.
\end{proof}

\begin{prop}\label{prop: E_6, omega_1 module, k = 2}
Let $G = E_6$ and $\lambda = \omega_1$, and take $k = 2$. Then the quadruple $(G, \lambda, p, k)$ has generic stabilizer $D_4.S_3$.
\end{prop}

\begin{proof}
We use the set-up of Proposition~\ref{prop: E_6, omega_1, F_4, omega_4 modules}: we take $H$ to be the simply connected group defined over $K$ of type $E_7$, with simple roots $\beta_1, \dots, \beta_7$; we let $G$ have simple roots $\alpha_i = \beta_i$ for $i \leq 6$, so that $G = \langle X_\alpha : \alpha = \sum m_i \beta_i, \ m_7 = 0 \rangle < H$; then we may take $V = \langle e_\alpha : \alpha = \sum m_i \beta_i, \ m_7 = 1 \rangle < \L(H)$. We have $Z(G) = \langle z \rangle$ where $z = h_{\beta_1}(\eta_3) h_{\beta_3}({\eta_3}^2) h_{\beta_5}(\eta_3) h_{\beta_6}({\eta_3}^2)$. We take the strictly positive generalized height function on the weight lattice of $G$ whose value at each simple root $\alpha_i$ is $1$, and then $\Lambda(V)_{[0]} = \{ \nu_1, \nu_2, \nu_3 \}$, where we write
$$
\gamma_1 = \esevenrt1122111, \quad \gamma_2 = \esevenrt1112211, \quad \gamma_3 = \esevenrt0112221,
$$
and for each $i$ we let $\nu_i$ be the weight such that $V_{\nu_i} = \langle e_{\gamma_i} \rangle$; we know that $\Lambda(V)_{[0]}$ has ZLC; and the setwise stabilizer in $W$ of $\Lambda(V)_{[0]}$ is $\langle w_{\beta_2}, w_{\beta_4}, w_{\beta_3} w_{\beta_5}, w_{\beta_1} w_{\beta_6} \rangle$. Here however we take $Y = \G{2}(V_{[0]})$, and write
$$
\hat Y = \left\{ y = \langle v^{(1)}, v^{(2)} \rangle \in Y : v^{(1)} = {\ts\sum} a_i e_{\gamma_i}, \ v^{(2)} = {\ts\sum} b_i e_{\gamma_i}, \ \forall i \neq j \
\left|
\begin{array}{cc}
 a_i & a_j \\
 b_i & b_j \\
\end{array}
\right| \neq 0 \right\};
$$
then $\hat Y$ is a dense open subset of $Y$, and the determinant condition implies that each $\nu_i$ occurs in every $y \in \hat Y$.

Let $A$ be the $D_4$ subgroup having simple roots $\beta_4$, $\beta_2$, $\beta_3 + \beta_4 + \beta_5$ and $\beta_1 + \beta_3 + \beta_4 + \beta_5 + \beta_6$; then $Z(A) = \langle z_1, z_2 \rangle$ where $z_1 = h_{\beta_3}(-1) h_{\beta_5}(-1)$ and $z_2 = h_{\beta_1}(-1) h_{\beta_6}(-1)$. We see that $V_{[0]}$ is the fixed point space of $A$ in its action on $V$, so clearly for all $y \in Y$ we have $A \leq C_G(y)$; let $C = Z(G) A \langle n_{\beta_3} n_{\beta_5}, n_{\beta_1} n_{\beta_6} \rangle$. Write $T_2 = \{ h_{\beta_1}(\kappa) h_{\beta_3}(\kappa') h_{\beta_5}({\kappa'}^{-1}) h_{\beta_6}(\kappa^{-1}) : \kappa, \kappa' \in K^* \}$, then for all $y \in Y$ we have $T_2 C \subseteq \Tran_G(y, Y)$.

Take $y \in \hat Y$ and $g \in \Tran_G(y, Y)$, and write $y' = g.y \in Y$. By Lemma~\ref{lem: gen height zero} we have $g = u_1 n u_2$ with $u_1 \in C_U(y')$, $u_2 \in C_U(y)$, and $n \in N_{\Lambda(V)_{[0]}}$ with $n.y = y'$.

First, from the above the elements of $W$ which preserve $\Lambda(V)_{[0]}$ are those corresponding to elements of $T_2 C \cap N$; so we have $N_{\Lambda(V)_{[0]}}.y = T_2.y$. As the elements of $N_{\Lambda(V)_{[0]}}$ permute and scale the $e_{\gamma_i}$, we have $N_{\Lambda(V)_{[0]}}.y \subseteq \hat Y$.

Next, let $\Xi = \Phi^+ \setminus \Phi_A$, and set $U' = \prod_{\alpha \in \Xi} X_\alpha$; then $U = U'.(A \cap U)$ and $U' \cap (A \cap U) = \{ 1 \}$. We now observe that if $\alpha \in \Xi$ then $\nu_i + \alpha$ is a weight in $V$ for exactly one value of $i$; moreover each weight in $V$ of positive generalized height is of the form $\nu_i + \alpha$ for exactly two such roots $\alpha$. Thus if we take $u = \prod x_\alpha(t_\alpha) \in U'$ satisfying $u.y = y$, and equate coefficients of weight vectors, taking them in an order compatible with increasing generalized height, using the determinant condition in the definition of the set $\hat Y_1$ we see that for all $\alpha$ we must have $t_\alpha = 0$, so that $u = 1$; so $C_U(y) = A \cap U$. Since the previous paragraph shows that $y' = g'.y \in \hat Y$, likewise we have $C_U(y') = A \cap U$.

Thus $\Tran_G(y, Y) = T_2 C$; so
$$
\codim {\ts\Tran_G}(y, Y) = \dim G - \dim {\ts\Tran_G}(y, Y) = 78 - 30 = 48,
$$
while
$$
\codim Y = \dim \G{2}(V) - \dim\G{2}(V_{[0]}) = 50 - 2 = 48.
$$
Therefore $y$ is $Y$-exact.

Now we may write $y = \langle e_{\gamma_1} + a_3 e_{\gamma_3}, e_{\gamma_2} + b_3 e_{\gamma_3} \rangle$ with $a_3, b_3 \neq 0$. Take $\kappa, \kappa' \in K^*$ satisfying $\kappa^6 = -\frac{1}{a_3 b_3}$ and ${\kappa'}^2 = b_3$, and take $h = h_{\beta_1}(\kappa) h_{\beta_3}(\kappa^2\kappa') h_{\beta_5}({\kappa^{-2}\kappa'}^{-1})$ $h_{\beta_6}(\kappa^{-1}) \in T_2$; then we find that ${}^h (n_{\beta_3} n_{\beta_5}) = n_{\beta_3} n_{\beta_5} h_{\beta_3}(\kappa^{-3}{b_3}^{-1}) h_{\beta_5}(\kappa^3 b_3) \in C_G(y)$ and ${}^h (n_{\beta_1} n_{\beta_6}) = n_{\beta_1} n_{\beta_6} h_{\beta_1}(\kappa') h_{\beta_6}({\kappa'}^{-1}) \in C_G(y)$, whence ${}^h C \leq C_G(y)$. Conversely $C_G(y) \leq T_2 C = {}^h (T_2 C)$. Given $s = h_{\beta_1}(\kappa) h_{\beta_3}(\kappa') h_{\beta_5}({\kappa'}^{-1}) h_{\beta_6}(\kappa^{-1}) \in T_2$ we have $s.(e_{\gamma_1} + a_3 e_{\gamma_3}) = {\kappa'}^2 e_{\gamma_1} + \kappa^{-2} a_3 e_{\gamma_3}$ and $s.(e_{\gamma_2} + b_3 e_{\gamma_3}) = \kappa^2 {\kappa'}^{-2} e_{\gamma_2} + \kappa^{-2} b_3 e_{\gamma_3}$, so $s \in C_G(y)$ requires $\kappa^2 {\kappa'}^2 = 1 = \kappa^{-4} {\kappa'}^2$, whence $\kappa^6 = 1$ and ${\kappa'}^2 = \kappa^4$ and so $s \in \langle z, z_1, z_2 \rangle < C$; hence $C_G(y) = {}^h C$.

Thus the conditions of Lemma~\ref{lem: generic stabilizer from exact subset} hold; so the quadruple $(G, \lambda, p, k)$ has generic stabilizer $C/Z(G) \cong D_4.S_3$, where the $D_4$ is of simply connected type.
\end{proof}

\begin{prop}\label{prop: E_6, omega_1, A_5, omega_2 modules, k = 3}
Let $G = E_6$ and $\lambda = \omega_1$, or $G = A_5$ and $\lambda = \omega_2$, and take $k = 3$. Then the quadruple $(G, \lambda, p, k)$ has generic stabilizer $A_2.\Z_{3/(p, 3)}.S_3$ or $T_1.\Z_{3/(p, 3)}.S_3$ respectively.
\end{prop}

\begin{proof}
We begin with the case where $G = E_6$ and $\lambda = \omega_1$. We use the set-up of Proposition~\ref{prop: E_6, omega_1, F_4, omega_4 modules}: we take $H$ to be the simply connected group defined over $K$ of type $E_7$, with simple roots $\beta_1, \dots, \beta_7$; we let $G$ have simple roots $\alpha_i = \beta_i$ for $i \leq 6$, so that $G = \langle X_\alpha : \alpha = \sum m_i \beta_i, \ m_7 = 0 \rangle < H$; then we may take $V = \langle e_\alpha : \alpha = \sum m_i \beta_i, \ m_7 = 1 \rangle < \L(H)$. We have $Z(G) = \langle z \rangle$ where $z = h_{\beta_1}(\eta_3) h_{\beta_3}({\eta_3}^2) h_{\beta_5}(\eta_3) h_{\beta_6}({\eta_3}^2)$. Here we take the generalized height function on the weight lattice of $G$ whose value at $\alpha_2$ and $\alpha_4$ is $1$, and at each other simple root $\alpha_i$ is $0$; then the generalized height of $\lambda = \frac{1}{3}(4\alpha_1 + 3\alpha_2 + 5\alpha_3 + 6\alpha_4 + 4\alpha_5 + 2\alpha_6)$ is $3$, and as $\lambda$ and $\Phi$ generate the weight lattice it follows that the generalized height of any weight is an integer. Since $V_\lambda = \langle e_\delta \rangle$ where $\delta = \esevenrt2234321$, we see that if $\mu \in \Lambda(V)$ and $e_\alpha \in V_\mu$ where $\alpha = \sum m_i \beta_i$ with $m_7 = 1$, then the generalized height of $\mu$ is $m_2 + m_4 - 3$. Thus $\Lambda(V)_{[0]} = \{ \nu_{11}, \dots, \nu_{33} \}$, where we write
\begin{eqnarray*}
& \gamma_{11} = \esevenrt0112111, \quad \gamma_{12} = \esevenrt0112211, \quad \gamma_{13} = \esevenrt0112221, & \\
& \gamma_{21} = \esevenrt1112111, \quad \gamma_{22} = \esevenrt1112211, \quad \gamma_{23} = \esevenrt1112221, & \\
& \gamma_{31} = \esevenrt1122111, \quad \gamma_{32} = \esevenrt1122211, \quad \gamma_{33} = \esevenrt1122221, &
\end{eqnarray*}
and for each $(i, j)$ we let $\nu_{ij}$ be the weight such that $V_{\nu_{ij}} = \langle e_{\gamma_{ij}} \rangle$. Observe that if we take $s = \prod_{i = 1}^6 h_{\beta_i}(\kappa_i) \in T$ then $\nu_{11}(s) = \frac{\kappa_4}{\kappa_1 \kappa_5}$, $\nu_{12}(s) = \frac{\kappa_5}{\kappa_1 \kappa_6}$, $\nu_{13}(s) = \frac{\kappa_6}{\kappa_1}$, $\nu_{21}(s) = \frac{\kappa_1 \kappa_4}{\kappa_3 \kappa_5}$, $\nu_{22}(s) = \frac{\kappa_1 \kappa_5}{\kappa_3 \kappa_6}$, $\nu_{23}(s) = \frac{\kappa_1 \kappa_6}{\kappa_3}$, $\nu_{31}(s) = \frac{\kappa_3}{\kappa_5}$, $\nu_{32}(s) = \frac{\kappa_3 \kappa_5}{\kappa_4 \kappa_6}$ and $\nu_{33}(s) = \frac{\kappa_3 \kappa_6}{\kappa_4}$; thus given any $5$-tuple $(n_1, n_2, n_3, n_4, n_5)$ of integers we have $c_{11} \nu_{11} + \cdots + c_{33} \nu_{33} = 0$ for $(c_{11}, c_{12}, c_{13}, c_{21}, c_{22}, c_{23}, c_{31}, c_{32}, c_{33}) = (n_1, n_3 + n_4, n_2 + n_5, n_3 + n_5, n_2, n_1 + n_4, n_2 + n_4, n_1 + n_5, n_3)$. In particular, writing \lq $(n_1, n_2, n_3, n_4, n_5) \implies (c_{11}, c_{12}, c_{13}, c_{21}, c_{22}, c_{23}, c_{31}, c_{32}, c_{33})$' to indicate this relationship between $5$-tuples and $9$-tuples, and for convenience writing $\bar 1$ for $-1$, we have the following:
\begin{eqnarray*}
     (0, 0, 0, 1, 1) \implies (0, 1, 1, 1, 0, 1, 1, 1, 0), & &           (0, 1, 1, 0, 0) \implies (0, 1, 1, 1, 1, 0, 1, 0, 1), \\
(1, 1, 1, \bar 1, 0) \implies (1, 0, 1, 1, 1, 0, 0, 1, 1), & &           (1, 1, 0, 0, 0) \implies (1, 0, 1, 0, 1, 1, 1, 1, 0), \\
(1, 1, 1, 0, \bar 1) \implies (1, 1, 0, 0, 1, 1, 1, 0, 1), & &           (1, 0, 1, 0, 0) \implies (1, 1, 0, 1, 0, 1, 0, 1, 1), \\
     (1, 0, 0, 0, 0) \implies (1, 0, 0, 0, 0, 1, 0, 1, 0), & & (1, 1, 1, \bar 1, \bar 1) \implies (1, 0, 0, 0, 1, 0, 0, 0, 1), \\
     (0, 1, 0, 0, 0) \implies (0, 0, 1, 0, 1, 0, 1, 0, 0), & &           (0, 0, 0, 1, 0) \implies (0, 1, 0, 0, 0, 1, 1, 0, 0), \\
     (0, 0, 1, 0, 0) \implies (0, 1, 0, 1, 0, 0, 0, 0, 1). & &           (0, 0, 0, 0, 1) \implies (0, 0, 1, 1, 0, 0, 0, 1, 0).
\end{eqnarray*}
By taking sums of these it follows that any subset of $\Lambda(V)_{[0]}$ whose complement is a subset of $\{ \nu_{11}, \nu_{22}, \nu_{33} \}$, $\{ \nu_{21}, \nu_{32}, \nu_{13} \}$, $\{ \nu_{31}, \nu_{12}, \nu_{23} \}$, $\{ \nu_{11}, \nu_{32}, \nu_{23} \}$, $\{ \nu_{21}, \nu_{12}, \nu_{33} \}$ or $\{ \nu_{31}, \nu_{22}, \nu_{13} \}$ has ZLCE.

Take $Y = \G{3}(V_{[0]})$. Given vectors $v^{(1)} = \sum a_{ij} e_{\gamma_{ij}}$, $v^{(2)} = \sum b_{ij} e_{\gamma_{ij}}$ and $v^{(3)} = \sum c_{ij} e_{\gamma_{ij}}$ in $V_{[0]}$, define the following $3 \times 3$ matrices $J_{i, j} = J_{i, j}(v^{(1)}, v^{(2)}, v^{(3)})$: for $j = 1, 2, 3$ set
$$
J_{1, j} =
\left(
  \begin{array}{ccc}
    a_{1j} & a_{2j} & a_{3j} \\
    b_{1j} & b_{2j} & b_{3j} \\
    c_{1j} & c_{2j} & c_{3j} \\
  \end{array}
\right), \quad
J_{2, j} =
\left(
  \begin{array}{ccc}
    a_{j1} & a_{j2} & a_{j3} \\
    b_{j1} & b_{j2} & b_{j3} \\
    c_{j1} & c_{j2} & c_{j3} \\
  \end{array}
\right);
$$
now for $i = 1, 2$ define the $9 \times 9$ matrices $J_i = J_i(v^{(1)}, v^{(2)}, v^{(3)})$ by
$$
J_i =
\left(
  \begin{array}{ccc}
        0     &  J_{i, 1} & J_{i, 2} \\
    -J_{i, 1} &     0     & J_{i, 3} \\
    -J_{i, 2} & -J_{i, 3} &     0    \\
  \end{array}
\right).
$$
We find that in the case where
\begin{eqnarray*}
v^{(1)} & = & a_{33} e_{\gamma_{33}} + a_{12} e_{\gamma_{12}} + a_{21} e_{\gamma_{21}}, \\
v^{(2)} & = & b_{11} e_{\gamma_{11}} + b_{23} e_{\gamma_{23}} + b_{32} e_{\gamma_{32}}, \\
v^{(3)} & = & e_{\gamma_{22}} + e_{\gamma_{31}} + e_{\gamma_{13}}
\end{eqnarray*}
then
$$
\det J_1 = -\det J_2 = (a_{12} b_{23} - a_{33} b_{11})(a_{21} b_{32} - a_{12} b_{23})(a_{33} b_{11} - a_{21} b_{32}),
$$
so that $\det J_1$ and $\det J_2$ are not identically zero. Observe that if we take $D = (d_{ij}) \in \GL_3(K)$, and for $i = 1, 2, 3$ we set ${v^{(i)}}' = d_{i1} v^{(1)} + d_{i2} v^{(2)} + d_{i3} v^{(3)}$, then for each $i$ and $j$ we have $J_{i, j}({v^{(1)}}', {v^{(2)}}', {v^{(3)}}') = D J_{i, j}(v^{(1)}, v^{(2)}, v^{(3)})$, whence for $i = 1, 2$ we have
$$
J_i({v^{(1)}}', {v^{(2)}}', {v^{(3)}}') =
\left(
  \begin{array}{ccc}
    D & 0 & 0 \\
    0 & D & 0 \\
    0 & 0 & D \\
  \end{array}
\right) J_i(v^{(1)}, v^{(2)}, v^{(3)}),
$$
so that $\det J_i({v^{(1)}}', {v^{(2)}}', {v^{(3)}}') = (\det D)^3 \det J_i(v^{(1)}, v^{(2)}, v^{(3)})$.  Therefore if we take $y \in Y$ and write $y = \langle v^{(1)}, v^{(2)}, v^{(3)} \rangle$, then although the individual determinants of the matrices $J_i(v^{(1)}, v^{(2)}, v^{(3)})$ depend on the choice of basis, whether or not they are zero does not. Thus if for each $i$ we set $\Delta_i = \det J_i(v^{(1)}, v^{(2)}, v^{(3)})$, we may define
$$
\hat Y_1 = \left\{ \langle v^{(1)}, v^{(2)}, v^{(3)} \rangle \in Y : \Delta_1 \Delta_2 \neq 0 \right\};
$$
then $\hat Y_1$ is a dense open subset of $Y$. (In fact a lengthy calculation shows that
$$
\Delta_1 + \Delta_2 = 0,
$$
so we could replace the condition in the definition of $\hat Y_1$ by simply \lq $\Delta_1 \neq 0$'.) Note that if $v^{(1)}$, $v^{(2)}$, $v^{(3)}$ are such that two weights differing by a root in $\Phi_{[0]}$ both fail to occur in any $v^{(i)}$, then one of the columns of either $J_1(v^{(1)}, v^{(2)}, v^{(3)})$ or $J_2(v^{(1)}, v^{(2)}, v^{(3)})$ is zero. Hence if $y \in \hat Y_1$ then the set of weights occurring in $y$ must meet any pair of weights differing by a root in $\Phi_{[0]}$; it follows that the complement of this set is a subset of $\{ \nu_{11}, \nu_{22}, \nu_{33} \}$, $\{ \nu_{21}, \nu_{32}, \nu_{13} \}$, $\{ \nu_{31}, \nu_{12}, \nu_{23} \}$, $\{ \nu_{11}, \nu_{32}, \nu_{23} \}$, $\{ \nu_{21}, \nu_{12}, \nu_{33} \}$ or $\{ \nu_{31}, \nu_{22}, \nu_{13} \}$.

In the proof of Proposition~\ref{prop: E_6, omega_1, F_4, omega_4 modules} we observed that the pointwise stabilizer in $W$ of $\{ \gamma_{13}, \gamma_{22}, \gamma_{31} \}$ is $W_1 = \langle w_{\beta_4}, w_{\beta_2}, w_{\beta_3 + \beta_4 + \beta_5}, w_{\beta_1 + \beta_3 + \beta_4 + \beta_5 + \beta_6} \rangle \cong W(D_4)$. Now if we write $\delta = \esixrt112221$, then the stabilizer in $W_1$ of $\beta_1$ contains $W_2 = \langle w_{\beta_4}, w_{\beta_2}, w_{\delta} \rangle \cong W(A_3)$, of index $8$, while the $W_1$-orbit of $\beta_1$ contains the eight roots $\esixrt100000$, $\esixrt101110$, $\esixrt111110$, $\esixrt111210$, $-\esixrt000001$, $-\esixrt001111$, $-\esixrt011111$ and $-\esixrt011211$, so has size at least $8$; thus the stabilizer in $W_1$ of $\beta_1$ is $W_2$. Since $\gamma_{12} = \gamma_{22} - \beta_1$ and $\gamma_{23} = \gamma_{13} + \beta_1$, we see that the pointwise stabilizer in $W$ of $\{ \gamma_{12}, \gamma_{13}, \gamma_{22}, \gamma_{23}, \gamma_{31} \}$ is $W_2$. Similarly if we write $\rho = \esixrt122321$, then the stabilizer in $W_2$ of $\beta_3$ contains $W_3 = \langle w_{\beta_2}, w_{\rho} \rangle \cong W(A_2)$, of index $4$, while the $W_2$-orbit of $\beta_3$ contains the four roots $\esixrt001000$, $\esixrt001100$, $\esixrt011100$ and $-\esixrt111221$, so has size at least $4$; thus the stabilizer in $W_2$ of $\beta_3$ is $W_3$. Since $\gamma_{21} = \gamma_{31} - \beta_3$, $\gamma_{32} = \gamma_{22} + \beta_3$, $\gamma_{33} = \gamma_{23} + \beta_3$ and $\gamma_{11} = \gamma_{21} - \beta_1$, we see that the pointwise stabilizer in $W$ of $\{ \gamma_{11}, \dots, \gamma_{33} \}$ is $W_3$. Now set $W_4 = \langle w_{\beta_1}, w_{\beta_3}, w_{\beta_5}, w_{\beta_6} \rangle \cong W({A_2}^2)$, and write $w^* = w_{\beta_1 + \beta_3 + \beta_4} w_{\beta_3 + \beta_4 + \beta_5} w_{\beta_4 + \beta_5 + \beta_6}$; then $W_4$ commutes with $W_3$, and $\langle w^* \rangle$ normalizes each of $W_3$ and $W_4$. Moreover $W_4$ acts transitively on $\{ \gamma_{11}, \dots, \gamma_{33} \}$, so given $w$ in the setwise stabilizer in $W$ of $\{ \gamma_{11}, \dots, \gamma_{33} \}$, there exists $w' \in W_4$ such that $w' w$ stabilizes $\gamma_{11}$; then $\gamma_{ij}$ is orthogonal to $\gamma_{11}$ only for $i, j \neq 1$, and the stabilizer in $W_4$ of $\gamma_{11}$ is $\langle w_{\beta_3}, w_{\beta_6} \rangle$, which acts transitively on $\{ \gamma_{22}, \gamma_{23}, \gamma_{32}, \gamma_{33} \}$, so there exists $w'' \in \langle w_{\beta_3}, w_{\beta_6} \rangle$ such that $w'' w' w$ stabilizes both $\gamma_{11}$ and $\gamma_{22}$; as $\gamma_{ij}$ is orthogonal to both $\gamma_{11}$ and $\gamma_{22}$ only for $(i, j) = (3, 3)$, we see that $w'' w' w$ also stabilizes $\gamma_{33}$; as $w^*$ interchanges $\gamma_{23}$ and $\gamma_{32}$ while fixing $\gamma_{11}$, $\gamma_{22}$ and $\gamma_{33}$, either $w'' w' w$ or $w^* w'' w' w$ stabilizes each of $\gamma_{11}$, $\gamma_{22}$, $\gamma_{23}$, $\gamma_{32}$ and $\gamma_{33}$; as each of the remaining $\gamma_{ij}$ is uniquely determined by which of $\gamma_{22}$ and $\gamma_{32}$ are orthogonal to it, we see that either $w'' w' w$ or $w^* w'' w' w$ lies in $W_3$. Thus the setwise stabilizer in $W$ of $\{ \gamma_{11}, \dots, \gamma_{33} \}$, and hence of $\Lambda(V)_{[0]}$, is $W_3 W_4 \langle w^* \rangle = \langle w_{\beta_2}, w_\rho, w_{\beta_1}, w_{\beta_3}, w_{\beta_1 + \beta_3 + \beta_4} w_{\beta_3 + \beta_4 + \beta_5} w_{\beta_4 + \beta_5 + \beta_6} \rangle \cong W({A_2}^3).\Z_2$. Note that this stabilizes $\Phi_{[0]} = \langle \alpha_1, \alpha_3, \alpha_5, \alpha_6 \rangle = \langle \beta_1, \beta_3, \beta_5, \beta_6 \rangle$.

Let $A$ be the $A_2$ subgroup having simple roots $\beta_2$ and $\rho - \beta_2$; then $Z(A) = \langle z' \rangle$ where $z' = h_{\beta_1}(\eta_3) h_{\beta_3}({\eta_3}^2) h_{\beta_5}({\eta_3}^2) h_{\beta_6}(\eta_3)$. We see that $V_{[0]}$ is the fixed point space of $A$ in its action on $V$, so clearly for all $y \in Y$ we have $A \leq C_G(y)$. Write $n^* = n_{\beta_1 + \beta_3 + \beta_4} n_{\beta_3 + \beta_4 + \beta_5} n_{\beta_4 + \beta_5 + \beta_6}$, and let $G_1$ be the derived group $(G_{[0]})' = \langle X_{\pm\alpha_1}, X_{\pm\alpha_3}, X_{\pm\alpha_5}, X_{\pm\alpha_6} \rangle \cong {A_2}^2$; then for all $y \in Y$ we have $A G_1 \langle n^* \rangle \subseteq \Tran_G(y, Y)$. Write $h^\dagger = h_{\beta_1}(\eta_3) h_{\beta_3}(\eta_3) h_{\beta_5}(\eta_3) h_{\beta_6}(\eta_3)$ and $n^\dagger = n_{\beta_1} n_{\beta_3} n_{\beta_5} n_{\beta_6}$, and set $C' = Z(G_1) \langle h^\dagger, n^\dagger, n^* \rangle$; let $C = C' A$, and then as $Z(G_1) = \langle h_{\beta_1}(\eta_3) h_{\beta_3}({\eta_3}^2), h_{\beta_5}(\eta_3) h_{\beta_6}({\eta_3}^2) \rangle = Z(G) Z(A)$ we have $C = Z(G) A \langle h^\dagger, n^\dagger, n^* \rangle$. Take $y =  \langle v^{(1)}, v^{(2)}, v^{(3)} \rangle \in \hat Y_1$; we shall show that $\Tran_G(y, Y) = A G_1 \langle n^* \rangle$, and that there is a dense open subset $\hat Y$ of $Y$ contained in $\hat Y_1$ such that if in fact $y \in \hat Y$ then $C_G(y) = {}^x C$ for some $x \in G$.

We have $U_{[0]} = X_{\alpha_1} X_{\alpha_3} X_{\alpha_1 + \alpha_3} X_{\alpha_5} X_{\alpha_6} X_{\alpha_5 + \alpha_6}$. If we take the root element $u = x_{\alpha_1}(t)$ for some $t \in K$, and write
$$
M =
\left(
  \begin{array}{ccc}
     1 & t &   \\
       & 1 &   \\
       &   & 1 \\
  \end{array}
\right),
$$
then for each $j$ we have $J_{1, j}(u.v^{(1)}, u.v^{(2)}, u.v^{(3)}) = J_{1, j}(v^{(1)}, v^{(2)}, v^{(3)}) M$, so that
$$
J_1(u.v^{(1)}, u.v^{(2)}, u.v^{(3)}) = J_1(v^{(1)}, v^{(2)}, v^{(3)})
\left(
  \begin{array}{ccc}
     M &   &   \\
       & M &   \\
       &   & M \\
  \end{array}
\right);
$$
however for $j = 1, 3$ we have $J_{2, j}(u.v^{(1)}, u.v^{(2)}, u.v^{(3)}) = J_{2, j}(v^{(1)}, v^{(2)}, v^{(3)})$, while $J_{2, 2}(u.v^{(1)}, u.v^{(2)}, u.v^{(3)}) = J_{2, 2}(v^{(1)}, v^{(2)}, v^{(3)}) + t J_{2, 1}(v^{(1)}, v^{(2)}, v^{(3)})$, so that
$$
J_2(u.v^{(1)}, u.v^{(2)}, u.v^{(3)}) =
\left(
  \begin{array}{ccc}
     I &    &   \\
       &  I &   \\
       & tI & I \\
  \end{array}
\right)
J_2(v^{(1)}, v^{(2)}, v^{(3)})
\left(
  \begin{array}{ccc}
     I &   &    \\
       & I & tI \\
       &   &  I \\
  \end{array}
\right).
$$
Similar equations hold for any root element $u = x_\alpha(t)$ where $\alpha \in \Phi_{[0]}$. Therefore $U_{[0]}$ preserves $\hat Y_1$; so given $u \in U_{[0]}$, by the above the set of weights occurring in $u.y$ has ZLCE. By Lemma~\ref{lem: gen height zero not strictly positive}, if we take $g \in \Tran_G(y, Y)$ and write $y' = g.y \in Y$, then we have $g = u_1 g' u_2$ with $u_1 \in C_{U_{[+]}}(y')$, $u_2 \in C_{U_{[+]}}(y)$, and $g' \in G_{[0]} N_{\Lambda(V)_{[0]}}$ with $g'.y = y'$. In particular $G.y \cap Y = G_{[0]} N_{\Lambda(V)_{[0]}}.y \cap Y$; moreover $C_G(y) = C_{U_{[+]}}(y) C_{G_{[0]} N_{\Lambda(V)_{[0]}}}(y) C_{U_{[+]}}(y)$.

First, since $W_{\Lambda(V)_{[0]}} = \langle w_{\beta_2}, w_\rho, w_{\beta_1}, w_{\beta_3}, w_{\beta_1 + \beta_3 + \beta_4} w_{\beta_3 + \beta_4 + \beta_5} w_{\beta_4 + \beta_5 + \beta_6} \rangle$ and $\beta_1, \beta_3 \in \Phi_{[0]}$, we have $G_{[0]} N_{\Lambda(V)_{[0]}} = G_{[0]} \langle n_{\beta_2}, n_\rho, n^* \rangle = G_1 (A \cap N) \langle n^* \rangle$. Any element of this last group may be written as $n' g^* c$ where $c \in A$, $g^* \in G_1$ and $n' \in \{ 1, n^* \}$; as $c.y = y$ it suffices to consider $n' g^*.y$. The above shows that applying any root element in $G_1$ has no effect on the determinants $\Delta_i$, so the same is true of $g^*$. We find that $n^*$ sends $e_{\gamma_{ij}}$ to $e_{\gamma_{ji}}$. Thus for $i = 1, 2$ and $j = 1, 2, 3$ we have $J_{i, j}(n^*.v^{(1)}, n^*.v^{(2)}, n^*.v^{(3)}) = J_{3 - i, j}(v^{(1)}, v^{(2)}, v^{(2)})$, and so $J_i(n^*.v^{(1)}, n^*.v^{(2)}, n^*.v^{(3)}) = J_{3 - i}(v^{(1)}, v^{(2)}, v^{(2)})$; so applying $n'$ permutes the determinants $\Delta_i$. Thus $G_{[0]} N_{\Lambda(V)_{[0]}}.y \subset \hat Y_1$. If we now further require the element $n' g^* c$ to stabilize $y$, we must have $n' g^*.y = y$. Since $V_{[0]}$ is the $G_1$-module with high weight $\omega_1 \otimes \omega_1$, using Proposition~\ref{prop: {A_2}^2, omega_1 otimes omega_1, k = 3} we see that there is a dense open subset $\hat Y_2$ of $Y$ each point of which has $G_1 \langle n^* \rangle$-stabilizer a conjugate of $C'$. Set $\hat Y = \hat Y_1 \cap \hat Y_2$; then if $y \in \hat Y$ we see that $C_{G_{[0]} N_{\Lambda(V)_{[0]}}}(y) = {}^x C' (A \cap N)$ for some $x \in G_1$.

Next, let $\Xi = \Phi^+ \setminus (\Phi_{[0]} \cup \Phi_A)$, and set $U' = \prod_{\alpha \in \Xi} X_\alpha$; then $U_{[+]} = U'.(A \cap U_{[+]})$ and $U' \cap (A \cap U_{[+]}) = \{ 1 \}$. We now observe that if $\alpha \in \Xi$ then $\nu_{ij} + \alpha$ is a weight in $V$ for exactly two pairs $(i, j)$; moreover each weight in $V$ of positive generalized height is of the form $\nu_i + \alpha$ for exactly six such roots $\alpha$. Indeed $\Xi$ is the union of three $W(G_1)$-orbits of size nine, distinguished by the coefficients of $\beta_2$ and $\beta_4$; likewise the nine weights in $V$ of positive generalized height form three sets of three, distinguished by the coefficients of $\beta_2$ and $\beta_4$ in the corresponding roots. If we now take a product of root elements corresponding to the nine roots in the orbit, and require it to stabilize $y$, equating coefficients of the corresponding three weight vectors in all three basis vectors of $y$ gives $9$ linear equations which may be expressed in matrix form using one of the matrices $J_i(v^{(1)}, v^{(2)}, v^{(3)})$ above. For example, one such orbit consists of the roots $\sum m_i \alpha_i$ with $m_2 = 0$ and $m_4 = 1$; here the three weights $\nu_{ij} + \alpha$ are those corresponding to the roots $\delta_1 = \esevenrt1123211$, $\delta_2 = \esevenrt1123221$ and $\delta_3 = \esevenrt1123321$. If we set $u = x_{\alpha_4}(t_1) x_{\alpha_3 + \alpha_4}(t_2) x_{\alpha_1 + \alpha_3 + \alpha_4}(t_3) x_{\alpha_4 + \alpha_5}(t_4) x_{\alpha_3 + \alpha_4 + \alpha_5}(t_5) x_{\alpha_1 + \alpha_3 + \alpha_4 + \alpha_5}(t_6)$ $x_{\alpha_4 + \alpha_5 + \alpha_6}(t_7) x_{\alpha_3 + \alpha_4 + \alpha_5 + \alpha_6}(t_8) x_{\alpha_1 + \alpha_3 + \alpha_4 + \alpha_5 + \alpha_6}(t_9)$, then we find that
\begin{eqnarray*}
u.{\ts\sum a_{ij} e_{\gamma_{ij}}} & = & {\ts\sum a_{ij} e_{\gamma_{ij}}} + (a_{32}t_1 - a_{22}t_2 + a_{12}t_3 + a_{31}t_4 - a_{21}t_5 + a_{11}t_6) e_{\delta_1} \\
                                   &   & {\phantom{{\ts\sum a_{ij} e_{\gamma_{ij}}}}} + (a_{33}t_1 - a_{23}t_2 + a_{13}t_3 - a_{31}t_7 + a_{21}t_8 - a_{11}t_9) e_{\delta_2} \\
                                   &   & {\phantom{{\ts\sum a_{ij} e_{\gamma_{ij}}}}} - (a_{33}t_4 - a_{23}t_5 + a_{13}t_6 + a_{32}t_7 - a_{22}t_8 + a_{12}t_9) e_{\delta_3}.
\end{eqnarray*}
Equating to zero the coefficients of $e_{\delta_1}$, $e_{\delta_2}$ and $e_{\delta_3}$ in each $u.v^{(i)}$ now gives the equation $J_1(v^{(1)}, v^{(2)}, v^{(3)}) {\bf t} = {\bf 0}$, where ${\bf t} = ( t_9 \ {-t_8} \ t_7 \ t_6 \ {-t_5} \ t_4 \ t_3 \ {-t_2} \ t_1 )^T$; since the matrix concerned has non-zero determinant we see that $t_i = 0$ for $i = 1, \dots, 9$. Thus if we take $u = \prod x_\alpha(t_\alpha) \in U'$ satisfying $u.y = y$, and equate coefficients of weight vectors, taking them in an order compatible with increasing generalized height, we see that for all $\alpha$ we must have $t_\alpha = 0$, so that $u = 1$; so $C_{U_{[+]}}(y) = A \cap U_{[+]}$. Since the previous paragraph shows that $y' = g'.y \in \hat Y_1$, likewise we have $C_{U_{[+]}}(y') = A \cap U_{[+]}$.

As $G_{[0]} N_{\Lambda(V)_{[0]}} \subseteq \Tran_G(y, Y)$ and $C_{U_{[+]}}(y), C_{U_{[+]}}(y') \leq A$ we do indeed have $\Tran_G(y, Y) = A G_1 \langle n^* \rangle \cong {A_2}^3.\Z_2$; so
$$
\codim {\ts\Tran_G(y, Y)} = \dim G - \dim {\ts\Tran_G(y, Y)} = 78 - 24 = 54
$$
while
$$
\codim Y = \dim \G{3}(V) - \dim Y = 72 - 18 = 54.
$$
Therefore $y$ is $Y$-exact. Now suppose further $y \in \hat Y$. As $C_{G_{[0]} N_{\Lambda(V)_{[0]}}}(y), C_{U_{[+]}}(y) \leq {}^x C' A$ for some $x \in G_1$ we have $C_G(y) = {}^x C' A = {}^x C$. Thus the conditions of Lemma~\ref{lem: generic stabilizer from exact subset} hold; so the quadruple $(G, \lambda, p, k)$ has generic stabilizer $C/Z(G) \cong A_2.\Z_{3/(p, 3)}.S_3$, where the $A_2$ is of simply connected type.

Now we may replace $G$ by the $A_5$ subgroup $\langle X_\alpha : \alpha = \sum m_i \beta_i, \ m_2 = m_7 = 0 \rangle$ of $E_6$, and $V$ by $\langle e_\alpha : \alpha = \sum m_i \beta_i, \ m_2 = m_7 = 1 \rangle$, leaving $Y$ and its subsets unchanged. We replace $A$ by the intersection of that above with $G$, which is the $T_1$ subgroup $\{ h_{\beta_1}(\kappa) h_{\beta_3}(\kappa^2) h_{\beta_4}(\kappa^3) h_{\beta_5}(\kappa^2) h_{\beta_6}(\kappa) : \kappa \in K^* \}$. We again let $C = C' A$; then $Z(C) = \langle h_{\beta_1}(\eta_3) h_{\beta_3}({\eta_3}^2) h_{\beta_5}({\eta_3}^2) h_{\beta_6}(\eta_3), h_{\beta_1}(-1) h_{\beta_4}(-1) h_{\beta_6}(-1) \rangle = Z(G)$. As $G_1 \langle n^* \rangle < G$, for all $y \in \hat Y_1$ we have $\Tran_G(y, Y) = A G_1 \langle n^* \rangle \cong {A_2}^2 T_1.\Z_2$; so
$$
\codim {\ts\Tran_G(y, Y)} = \dim G - \dim {\ts\Tran_G(y, Y)} = 35 - 17 = 18
$$
while
$$
\codim Y = \dim \G{3}(V) - \dim Y = 36 - 18 = 18.
$$
Therefore $y$ is $Y$-exact. Also if $y \in \hat Y$ then $C_G(y) = {}^x C' A = {}^x C$ for some $x \in G_1$. Thus the conditions of Lemma~\ref{lem: generic stabilizer from exact subset} hold; so the quadruple $(G, \lambda, p, k)$ has generic stabilizer $C/Z(G) \cong T_1.\Z_{3/(p, 3)}.S_3$.
\end{proof}

\begin{prop}\label{prop: F_4, omega_4, C_3, omega_2 modules, k = 2}
Let $G = F_4$ and $\lambda = \omega_4$, or $G = C_3$ and $\lambda = \omega_2$, and take $k = 2$. Then the quadruple $(G, \lambda, p, k)$ has generic stabilizer $A_2$ or $T_1$ respectively if $p \neq 3$, and $A_2.\Z_2$ or $T_1.\Z_2$ respectively if $p = 3$.
\end{prop}

\begin{proof}
We use the set-up of Proposition~\ref{prop: E_6, omega_1, A_5, omega_2 modules, k = 3}, but modify the notation slightly: we take $H$ to be the simply connected group defined over $K$ of type $E_7$, with simple roots $\beta_1, \dots, \beta_7$; we let $G^+ = E_6$ have simple roots $\alpha_i = \beta_i$ for $i \leq 6$, so that $G^+ = \langle X_\alpha : \alpha = \sum m_i \beta_i, \ m_7 = 0 \rangle < H$; then we may take $V^+ = \langle e_\alpha : \alpha = \sum m_i \beta_i, \ m_7 = 1 \rangle < \L(H)$; we write
\begin{eqnarray*}
& \gamma_{11} = \esevenrt0112111, \quad \gamma_{12} = \esevenrt0112211, \quad \gamma_{13} = \esevenrt0112221, & \\
& \gamma_{21} = \esevenrt1112111, \quad \gamma_{22} = \esevenrt1112211, \quad \gamma_{23} = \esevenrt1112221, & \\
& \gamma_{31} = \esevenrt1122111, \quad \gamma_{32} = \esevenrt1122211, \quad \gamma_{33} = \esevenrt1122221, &
\end{eqnarray*}
and let $Y^+ = \G{3}(\langle e_{\gamma_{ij}} : 1 \leq i, j \leq 3 \rangle)$. We let $A$ be the $A_2$ subgroup having simple roots $\beta_2$ and $\rho - \beta_2$, where $\rho = \esixrt122321$ is the high root of $G^+$; we write $G_1 = \langle X_{\pm\alpha_1}, X_{\pm\alpha_3}, X_{\pm\alpha_5}, X_{\pm\alpha_6} \rangle \cong {A_2}^2$, and set $n^* = n_{\beta_1 + \beta_3 + \beta_4} n_{\beta_3 + \beta_4 + \beta_5} n_{\beta_4 + \beta_5 + \beta_6}$. The proof of Proposition~\ref{prop: E_6, omega_1, A_5, omega_2 modules, k = 3} showed that there is a dense open subset $\hat Y_1^+$ of $Y^+$, defined by the non-vanishing of a determinant, such that if $y^+ \in \hat Y_1^+$ then $\Tran_{G^+}(y^+, Y^+) = A G_1 \langle n^* \rangle$.

We saw in the proof of Proposition~\ref{prop: E_6, omega_1, F_4, omega_4 modules} that if we set $v_0 = e_{\gamma_{22}} + e_{\gamma_{31}} + e_{\gamma_{13}}$ then $C_{G^+}(v_0)$ is the $F_4$ subgroup having long simple roots $\beta_2$ and $\beta_4$ and short simple root groups $\{ x_{\beta_3}(t) x_{\beta_5}(-t) : t \in K \}$ and $\{ x_{\beta_1}(t) x_{\beta_6}(-t) : t \in K \}$. We now take $G$ to be this $F_4$ subgroup; then $Z(G) = \{ 1 \}$. We have $A \langle n^* \rangle \leq G$, and if we write $\tilde A_2$ for the subgroup having simple root groups $\{ x_{\beta_3}(t) x_{\beta_5}(-t) : t \in K \}$ and $\{ x_{\beta_1}(t) x_{\beta_6}(-t) : t \in K \}$ then $G \cap G_1 = \tilde A_2$. Inside $V^+$ we have the submodules $X_1 = \{ \sum a_\gamma e_\gamma \in V^+ : a_{\gamma_{22}} + a_{\gamma_{31}}  + a_{\gamma_{13}} = 0 \}$ and $X_2 = \langle v_0 \rangle$, with the latter being trivial. If $p \neq 3$ then $V^+ = X_1 \oplus X_2$, and $V = X_1$; if however $p = 3$ then $X_2 < X_1$, and $V = X_1/X_2$. Thus in all cases $V = X_1/(X_1 \cap X_2)$, where $X_1 \cap X_2$ is either zero or the trivial $G$-module.

Write $\bar V = \{ \sum a_{ij} e_{\gamma_{ij}} + (X_1 \cap X_2) : a_{\gamma_{22}} + a_{\gamma_{31}}  + a_{\gamma_{13}} = 0 \}$, and let $Y = \G{2}(\bar V)$; then $\dim Y = 12 - 2\z_{3, p}$. Given $y \in Y$, write $y = \langle v^{(1)} + (X_1 \cap X_2), v^{(2)} + (X_1 \cap X_2) \rangle$ and set $y^+ = \langle v^{(1)}, v^{(2)}, v_0 \rangle \in Y^+$. We observed in the proof of Proposition~\ref{prop: E_6, omega_1, A_5, omega_2 modules, k = 3} that the determinant defining the set $\hat Y_1^+$ is not identically zero for such elements of $Y^+$; thus there is a dense open subset $\hat Y_1$ of $Y$ such that if $y \in \hat Y_1$ then $y^+ \in \hat Y_1^+$.

Take $y \in \hat Y_1$; then $A \tilde A_2 \langle n^* \rangle \subseteq \Tran_G(y, Y)$. Conversely, given $g \in \Tran_G(y, Y)$ we have $g.v_0 = v_0$ and $g.y \in Y$, so $g.y^+ \in Y^+$, whence $g \in \Tran_{G^+}(y^+, Y^+) = A G_1 \langle n^* \rangle$; thus $g \in G \cap A G_1 \langle n^* \rangle = A \tilde A_2 \langle n^* \rangle$. Therefore $\Tran_G(y, Y) = A \tilde A_2 \langle n^* \rangle$; so
$$
\codim {\ts\Tran_G(y, Y)} = \dim G - \dim {\ts\Tran_G(y, Y)} = 52 - 16 = 36
$$
while
$$
\codim Y = \dim \G{2}(V) - \dim Y = (48 - 2\z_{3, p}) - (12 - 2\z_{3, p}) = 36.
$$
Therefore $y$ is $Y$-exact.

Now we certainly have $A \leq C_G(y)$, so it remains to consider $C_{\tilde A_2 \langle n^* \rangle}(y)$. Observe that $\bar V$ is an irreducible module for $\tilde A_2$ with high weight $\omega_1 + \omega_2$, of dimension $8 - \z_{3, p}$. By Proposition~\ref{prop: classical Lie algebras, k = 2} we know that the action of $\tilde A_2$ on $\G{2}(L(\omega_1 + \omega_2))$ has trivial generic stabilizer; thus if we extend the action to $\tilde A_2 \langle n^* \rangle$ the generic stabilizer must be either $1$ or $\Z_2$. By \cite[Table~4.3.1]{GLS} there is a single conjugacy class of outer involutions in $\tilde A_2 \langle n^* \rangle$, represented by $n^*$ and with centralizer $B_1$, so $\dim (n^*)^{\tilde A_2} = 5$. We saw in the proof of Proposition~\ref{prop: E_6, omega_1, A_5, omega_2 modules, k = 3} that $n^*$ sends each $e_{\gamma_{ij}}$ to $e_{\gamma_{ji}}$; thus its fixed point space in $\bar V$ is spanned by the images of $e_{\gamma_{33}}$, $e_{\gamma_{23}} + e_{\gamma_{32}}$, $e_{\gamma_{11}}$, $e_{\gamma_{12}} + e_{\gamma_{21}}$ and $-2e_{\gamma_{22}} + e_{\gamma_{31}} + e_{\gamma_{13}}$ (the last of these being $v_0$ if $p = 3$), and so has dimension $5 - \z_{3, p}$. Thus in the notation of Chapter~\ref{chap: TGS quadruples} we have $\d = (5 - \z_{3, p}, 3)$, whence using Proposition~\ref{prop: B value when t = 2} we have $\codim C_{\G{2}(L(\omega_1 + \omega_2))}(n^*) = B_{\d, 2} = 6 - \z_{3, p}$. Thus if $p \neq 3$ we have $\codim C_{\G{2}(L(\omega_1 + \omega_2))}(n^*) > \dim (n^*)^{\tilde A_2}$; so the generic stabilizer for the action of $\tilde A_2 \langle n^* \rangle$ does not meet $(n^*)^{\tilde A_2}$, and so must be trivial. If however $p = 3$ we see from Proposition~\ref{prop: A_2.2, omega_1 + omega_2 module, p = 3, k = 2} that the generic stabilizer is $\langle n^* \rangle \cong \Z_2$. Therefore there is a dense open subset $\hat Y$ of $Y$, which we may assume lies in $\hat Y_1$, such that if $y \in \hat Y$ then according as $p \neq 3$ or $p = 3$ we have $C_G(y) = A$ or $A \langle n^* \rangle^x$ for some $x \in \tilde A_2$. Thus the conditions of Lemma~\ref{lem: generic stabilizer from exact subset} hold; so according as $p \neq 3$ or $p = 3$ the quadruple $(G, \lambda, p, k)$ has generic stabilizer $C/Z(G) \cong A_2$ or $A_2.\Z_2$.

Now as in the proof of Proposition~\ref{prop: E_6, omega_1, A_5, omega_2 modules, k = 3} we may replace $G^+$ by the $A_5$ subgroup $\langle X_\alpha : \alpha = \sum m_i \beta_i, \ m_2 = m_7 = 0 \rangle$ of $E_6$, and then $G$ by the intersection of that above with $G^+$, which is the $C_3$ subgroup with long simple root $\beta_4$ and short simple root groups $\{ x_{\beta_3}(t) x_{\beta_5}(-t) : t \in K \}$ and $\{ x_{\beta_1}(t) x_{\beta_6}(-t) : t \in K \}$, and $A$ by the intersection of that above with $G^+$, which is the $T_1$ subgroup $\{ h_{\beta_1}(\kappa) h_{\beta_3}(\kappa^2) h_{\beta_4}(\kappa^3) h_{\beta_5}(\kappa^2) h_{\beta_6}(\kappa) : \kappa \in K^* \}$. We also replace $V^+$ by $\langle e_\alpha : \alpha = \sum m_i \beta_i, \ m_2 = m_7 = 1 \rangle$, and then $X_1$ by the intersection of that above with $V^+$, leaving $X_2$ unchanged; again we have $V = X_1/(X_1 \cap X_2)$, and we leave $Y$ and its subsets unchanged. As $\tilde A_2 \langle n^* \rangle < G$, for all $y \in \hat Y_1$ we have $\Tran_G(y, Y) = A \tilde A_2 \langle n^* \rangle \cong \tilde A_2T_1.\Z_2$; so
$$
\codim {\ts\Tran_G(y, Y)} = \dim G - \dim {\ts\Tran_G(y, Y)} = 21 - 9 = 12
$$
while
$$
\codim Y = \dim \G{2}(V) - \dim Y = (24 - 2\z_{3, p}) - (12 - 2\z_{3, p}) = 12.
$$
Therefore $y$ is $Y$-exact. Also if $y \in \hat Y$ then according as $p \neq 3$ or $p = 3$ we have $C_G(y) = A$ or $A \langle n^* \rangle^x$ for some $x \in \tilde A_2$. Thus the conditions of Lemma~\ref{lem: generic stabilizer from exact subset} hold; so according as $p \neq 3$ or $p = 3$ the quadruple $(G, \lambda, p, k)$ has generic stabilizer $C/Z(G) \cong T_1$ or $T_1.\Z_2$.
\end{proof}

\begin{prop}\label{prop: F_4, omega_1, B_3, omega_2 modules, p = 2, k = 2}
Let $G = F_4$ and $\lambda = \omega_1$ with $p = 2$, or $G = B_3$ and $\lambda = \omega_2$ with $p = 2$, and take $k = 2$. Then the quadruple $(G, \lambda, p, k)$ has generic stabilizer $\tilde A_2$ or $T_1$ respectively.
\end{prop}

\begin{proof}
This is an immediate consequence of Proposition~\ref{prop: F_4, omega_4, C_3, omega_2 modules, k = 2}, using the graph automorphism of $F_4$ and the exceptional isogeny $B_\ell \to C_\ell$ which exist in characteristic $2$.
\end{proof}

\begin{prop}\label{prop: E_7, omega_7, D_6, omega_6, A_5, omega_3, C_3, omega_3 modules, k = 2}
Let $G = E_7$ and $\lambda = \omega_7$, or $G = D_6$ and $\lambda = \omega_6$, or $G = A_5$ and $\lambda = \omega_3$, or $G = C_3$ and $\lambda = \omega_3$ with $p \geq 3$, and take $k = 2$. Then the quadruple $(G, \lambda, p, k)$ has generic stabilizer $D_4.\Z_{2/(p, 2)}.\Z_2$, or ${A_1}^3.\Z_{2/(p, 2)}.\Z_2$, or $T_2.\Z_{2/(p, 2)}.\Z_2$, or ${\Z_2}^4$, respectively.
\end{prop}

\begin{proof}
We begin with the case where $G = E_7$ and $\lambda = \omega_7$. We use the set-up of Proposition~\ref{prop: E_7, omega_7, D_6, omega_6, B_5, omega_5, A_5, omega_3, C_3, omega_3 modules}: we take $H$ to be the (simply connected) group defined over $K$ of type $E_8$, with simple roots $\beta_1, \dots, \beta_8$; we let $G$ have simple roots $\alpha_i = \beta_i$ for $i \leq 7$, so that $G = \langle X_\alpha : \alpha = \sum m_i \beta_i, \ m_8 = 0 \rangle < H$; then we may take $V = \langle e_\alpha : \alpha = \sum m_i \beta_i, \ m_8 = 1 \rangle < \L(H)$. We have $Z(G) = \langle z \rangle$ where $z = h_{\beta_2}(-1) h_{\beta_5}(-1) h_{\beta_7}(-1)$. Here we take the generalized height function on the weight lattice of $G$ whose value at $\alpha_2$, $\alpha_5$ and $\alpha_7$ is $0$, and at each other simple root $\alpha_i$ is $1$; then the generalized height of $\lambda = \frac{1}{2}(2\alpha_1 + 3\alpha_2 + 4\alpha_3 + 6\alpha_4 + 5\alpha_5 + 4\alpha_6 + 3\alpha_7)$ is $8$, and as $\lambda$ and $\Phi$ generate the weight lattice it follows that the generalized height of any weight is an integer. Since $V_\lambda = \langle e_\delta \rangle$ where $\delta = \eeightrt23465431$, we see that if $\mu \in \Lambda(V)$ and $e_\alpha \in V_\mu$ where $\alpha = \sum m_i \beta_i$ with $m_8 = 1$, then the generalized height of $\mu$ is $m_1 + m_3 + m_4 + m_6 - 8$. Thus $\Lambda(V)_{[0]} = \{ \nu_1, \dots, \nu_8 \}$, where we write
\begin{eqnarray*}
& \gamma_1 = \eeightrt11232211, \quad \gamma_2 = \eeightrt11232221, \quad \gamma_3 = \eeightrt11233211, \quad \gamma_4 = \eeightrt11233221, & \\
& \gamma_5 = \eeightrt12232211, \quad \gamma_6 = \eeightrt12232221, \quad \gamma_7 = \eeightrt12233211, \quad \gamma_8 = \eeightrt12233221, &
\end{eqnarray*}
and for each $i$ we let $\nu_i$ be the weight such that $V_{\nu_i} = \langle e_{\gamma_i} \rangle$. Observe that if we take $s = \prod_{i = 1}^7 h_{\beta_i}(\kappa_i) \in T$ then $\nu_1(s) = \frac{\kappa_4 \kappa_6}{\kappa_2 \kappa_5 \kappa_7}$, $\nu_2(s) = \frac{\kappa_4 \kappa_7}{\kappa_2 \kappa_5}$, $\nu_3(s) = \frac{\kappa_5}{\kappa_2 \kappa_7}$, $\nu_4(s) = \frac{\kappa_5 \kappa_7}{\kappa_2 \kappa_6}$, $\nu_5(s) = \frac{\kappa_2 \kappa_6}{\kappa_5 \kappa_7}$, $\nu_6(s) = \frac{\kappa_2 \kappa_7}{\kappa_5}$, $\nu_7(s) = \frac{\kappa_2 \kappa_5}{\kappa_4 \kappa_7}$ and $\nu_8(s) = \frac{\kappa_2 \kappa_5 \kappa_7}{\kappa_4 \kappa_6}$; thus given any $5$-tuple $(n_1, n_2, n_3, n_4, n_5)$ of integers we have $c_1 \nu_1 + \cdots + c_8 \nu_8 = 0$ for $(c_1, \dots, c_8) = (n_1 + n_5, n_2, n_3, n_4 + n_5, n_4, n_3 + n_5, n_2 + n_5, n_1)$. In particular, writing \lq $(n_1, n_2, n_3, n_4, n_5) \implies (c_1, c_2, c_3, c_4, c_5, c_6, c_7, c_8)$' to indicate this relationship between $5$-tuples and $8$-tuples, we have the following:
\begin{eqnarray*}
(0, 0, 0, 0, 1) \implies (1, 0, 0, 1, 0, 1, 1, 0), & & \phantom{1} (1, 1, 1, 1, -1) \implies (0, 1, 1, 0, 1, 0, 0, 1), \\
(1, 0, 0, 0, 0) \implies (1, 0, 0, 0, 0, 0, 0, 1), & & \phantom{-1} (0, 1, 0, 0, 0) \implies (0, 1, 0, 0, 0, 0, 1, 0), \\
(0, 0, 1, 0, 0) \implies (0, 0, 1, 0, 0, 1, 0, 0), & & \phantom{-1} (0, 0, 0, 1, 0) \implies (0, 0, 0, 1, 1, 0, 0, 0).
\end{eqnarray*}
By taking sums of these it follows that any subset of $\Lambda(V)_{[0]}$ whose complement either is a subset of $\{ \nu_2, \nu_3, \nu_5, \nu_8 \}$ or $\{ \nu_1, \nu_4, \nu_6, \nu_7 \}$, or is of the form $\{ \nu_j, \nu_{9 - j} \}$ for some $j \leq 4$, has ZLCE.

Take $Y = \G{2}(V_{[0]})$. Given vectors $v^{(1)} = \sum a_i e_{\gamma_i}$ and $v^{(2)} = \sum b_i e_{\gamma_i}$ in $V_{[0]}$, define the following $4 \times 4$ matrices $J_i = J_i(v^{(1)}, v^{(2)})$:
$$
J_1 =
\left(
  \begin{array}{cccc}
    a_1 & a_2 & a_3 & a_4 \\
    b_1 & b_2 & b_3 & b_4 \\
    a_5 & a_6 & a_7 & a_8 \\
    b_5 & b_6 & b_7 & b_8 \\
  \end{array}
\right), \quad
J_2 =
\left(
  \begin{array}{cccc}
    a_1 & a_5 & a_2 & a_6 \\
    b_1 & b_5 & b_2 & b_6 \\
    a_3 & a_7 & a_4 & a_8 \\
    b_3 & b_7 & b_4 & b_8 \\
  \end{array}
\right), \quad
J_3 =
\left(
  \begin{array}{cccc}
    a_1 & a_3 & a_5 & a_7 \\
    b_1 & b_3 & b_5 & b_7 \\
    a_2 & a_4 & a_6 & a_8 \\
    b_2 & b_4 & b_6 & b_8 \\
  \end{array}
\right).
$$
Observe that if we take $D = (d_{ij}) \in \GL_2(K)$, and for $i = 1, 2$ we set ${v^{(i)}}' = d_{i1} v^{(1)} + d_{i2} v^{(2)}$, then for each $i$ we have
$$
J_i({v^{(1)}}', {v^{(2)}}') =
\left(
  \begin{array}{cc}
    D & 0 \\
    0 & D \\
  \end{array}
\right) J_i(v^{(1)}, v^{(2)}),
$$
so that $\det J_i({v^{(1)}}', {v^{(2)}}') = (\det D)^2 \det J_i(v^{(1)}, v^{(2)})$. Therefore if we take $y \in Y$ and write $y = \langle v^{(1)}, v^{(2)} \rangle$, then although the individual determinants of the matrices $J_i(v^{(1)}, v^{(2)})$ depend on the choice of basis, the ratio of any two of these determinants does not. Thus if for each $i$ we set $\Delta_i = \det J_i(v^{(1)}, v^{(2)})$, we may define
$$
\hat Y_1 = \left\{ \langle v^{(1)}, v^{(2)} \rangle \in Y : \Delta_1 \Delta_2 \Delta_3 \neq 0, \ \forall i \neq j \ ({\ts\frac{\Delta_i}{\Delta_j}})^3 \neq 1 \right\};
$$
then $\hat Y_1$ is a dense open subset of $Y$. (In fact a straightforward calculation shows that
$$
\Delta_1 + \Delta_2 + \Delta_3 = 0.)
$$
Note that if $v^{(1)}$, $v^{(2)}$ are such that two weights differing by a root in $\Phi_{[0]}$ both fail to occur in either $v^{(i)}$, then one of the columns of $J_1(v^{(1)}, v^{(2)})$, $J_2(v^{(1)}, v^{(2)})$ or $J_3(v^{(1)}, v^{(2)})$ is zero. Hence if $y \in \hat Y_1$ then the set of weights occurring in $y$ must meet any pair of weights differing by a root in $\Phi_{[0]}$; it follows that the complement of this set either is a subset of $\{ \nu_2, \nu_3, \nu_4, \nu_8 \}$ or $\{ \nu_1, \nu_5, \nu_6, \nu_7 \}$, or is of the form $\{ \nu_j, \nu_{9 - j} \}$ for some $j \leq 4$.

In the proof of Proposition~\ref{prop: E_7, omega_7, D_6, omega_6, B_5, omega_5, A_5, omega_3, C_3, omega_3 modules} we observed that the pointwise stabilizer in $W$ of $\{ \gamma_3, \gamma_6 \}$ is $W_1 = \langle w_{\beta_6}, w_{\beta_1}, w_{\beta_2 + \beta_4 + \beta_5}, w_{\beta_3}, w_{\beta_4}, w_{\beta_5 + \beta_6 + \beta_7} \rangle \cong W(E_6)$; in fact if we write $\delta = \esevenrt1123321$ then ${W_1}^{w_\delta} = \langle w_{\beta_1}, \dots, w_{\beta_6} \rangle$. Now ${W_1}^{w_\delta}$ acts transitively on the set $\Sigma$ of roots $\alpha$ of the form $\sum m_i \beta_i$ with $m_7 = 1$ and $m_8 = 0$, so the stabilizer in ${W_1}^{w_\delta}$ of $\beta_7$ has order $\frac{|{W_1}^{w_\delta}|}{|\Sigma|} = \frac{|W_1|}{27} = |W(D_5)|$; we see then that this stabilizer is $\langle w_{\beta_1}, \dots, w_{\beta_5} \rangle$, and hence if we write $\rho_6 = \esevenrt0112221$ then the stabilizer in $W_1$ of $w_\delta(\beta_7) = \beta_7$ is $W_2 = \langle w_{\beta_1}, \dots, w_{\beta_5} \rangle^{w_\delta} = \langle w_{\rho_6}, w_{\beta_1}, w_{\beta_3}, w_{\beta_4}, w_{\beta_2 + \beta_4 + \beta_5} \rangle$. Since $\gamma_5 = \gamma_6 - \beta_7$ and $\gamma_4 = \gamma_3 + \beta_7$, the pointwise stabilizer in $W$ of $\{ \gamma_3, \gamma_4, \gamma_5, \gamma_6 \}$ is $W_2$. Next, if we write $\rho_4 = \esevenrt0112100$, the stabilizer in $W_2$ of $\beta_2$ contains $W_3 = \langle w_{\beta_3}, w_{\beta_1}, w_{\rho_4}, w_{\rho_6} \rangle \cong W(D_4)$, of index $10$, while the $W_2$-orbit of $\beta_2$ contains the eight roots $\sum m_i \beta_i$ with $(m_2, m_5, m_6) = (1, 0, 0)$ or $(0, -1, 0)$, along with $\esevenrt1223221$ and $-\esevenrt1123321$, so has size at least $10$; thus the stabilizer in $W_2$ of $\beta_2$ is $W_3$. Since $\gamma_1 = \gamma_5 - \beta_2$, $\gamma_2 = \gamma_6 - \beta_2$, $\gamma_7 = \gamma_3 + \beta_2$ and $\gamma_8 = \gamma_4 + \beta_2$, we see that the pointwise stabilizer in $W$ of $\{ \gamma_1, \dots, \gamma_8 \}$ is $W_3$. Now set $W_4 = \langle w_{\beta_2}, w_{\beta_5}, w_{\beta_7} \rangle \cong W({A_1}^3)$, and write $w^* = w_{\beta_2 + \beta_4} w_{\beta_4 + \beta_5}$ and $w^{**} = w_{\beta_5 + \beta_6} w_{\beta_6 + \beta_7}$; then $W_4$ commutes with $W_3$, and $\langle w^*, w^{**} \rangle \cong S_3$ normalizes each of $W_3$ and $W_4$. Moreover $W_4$ acts simply transitively on $\{ \gamma_1, \dots, \gamma_8 \}$, so given $w$ in the setwise stabilizer in $W$ of $\{ \gamma_1, \dots, \gamma_8 \}$, there exists $w' \in W_4$ such that $w' w$ stabilizes $\gamma_1$; then $\gamma_i - \gamma_1$ is a root only for $i \in \{ 2, 3, 5 \}$, and $\langle w^*, w^{**} \rangle$ stabilizes $\gamma_1$ while acting as $S_3$ on $\{ \gamma_2, \gamma_3, \gamma_5 \}$, so there exists $w'' \in \langle w^*, w^{**} \rangle$ such that $w'' w' w$ stabilizes each of $\gamma_1$, $\gamma_2$, $\gamma_3$ and $\gamma_5$; as $\gamma_4 = \gamma_2 + \gamma_3 - \gamma_1$, $\gamma_6 = \gamma_2 + \gamma_5 - \gamma_1$, $\gamma_7 = \gamma_3 + \gamma_5 - \gamma_1$ and $\gamma_8 = \gamma_2 + \gamma_3 + \gamma_5 - 2\gamma_1$, we see that $w'' w' w \in W_3$. Thus the setwise stabilizer in $W$ of $\{ \gamma_1, \dots, \gamma_8 \}$, and hence of $\Lambda(V)_{[0]}$, is $W_3 W_4 \langle w^*, w^{**} \rangle = \langle w_{\beta_3}, w_{\beta_1}, w_{\beta_2}, w_{\beta_2 + \beta_4} w_{\beta_4 + \beta_5}, w_{\beta_5 + \beta_6} w_{\beta_6 + \beta_7} \rangle \cong W(D_4{A_1}^3).S_3$. Note that this stabilizes $\Phi_{[0]} = \langle \alpha_2, \alpha_5, \alpha_7 \rangle = \langle \beta_2, \beta_5, \beta_7 \rangle$.

Let $A$ be the $D_4$ subgroup having simple roots $\beta_3$, $\beta_1$, $\rho_4$ and $\rho_6$; then $Z(A) = \langle z_1, z_2 \rangle$ where $z_1 = h_{\beta_2}(-1) h_{\beta_5}(-1)$ and $z_2 = h_{\beta_5}(-1) h_{\beta_7}(-1)$. We see that $V_{[0]}$ is the fixed point space of $A$ in its action on $V$, so clearly for all $y \in Y$ we have $A \leq C_G(y)$. Write $n^* = n_{\beta_2 + \beta_4} n_{\beta_4 + \beta_5}$ and $n^{**} = n_{\beta_5 + \beta_6} n_{\beta_6 + \beta_7}$, and let $G_1$ be the derived group $(G_{[0]})' = \langle X_{\pm\alpha_2}, X_{\pm\alpha_5}, X_{\pm\alpha_7} \rangle \cong {A_1}^3$; then for all $y \in Y$ we have $A G_1 \langle n^*, n^{**} \rangle \subseteq \Tran_G(y, Y)$. Write $h^\dagger = h_{\beta_2}(\eta_4) h_{\beta_5}(\eta_4) h_{\beta_7}(\eta_4)$ and $n^\dagger = n_{\beta_2} n_{\beta_5} n_{\beta_7}$, and set $C' = Z(G_1) \langle h^\dagger, n^\dagger \rangle$; let $C = C' A$, and then as $Z(G_1) = \langle h_{\beta_2}(-1), h_{\beta_5}(-1), h_{\beta_7}(-1) \rangle = Z(G) Z(A)$ we have $C = Z(G) A \langle h^\dagger, n^\dagger \rangle$. Take $y =  \langle v^{(1)}, v^{(2)} \rangle \in \hat Y_1$; we shall show that $\Tran_G(y, Y) = A G_1 \langle n^*, n^{**} \rangle$, and that there is a dense open subset $\hat Y$ of $Y$ contained in $\hat Y_1$ such that if in fact $y \in \hat Y$ then $C_G(y) = {}^x C$ for some $x \in G_1$.

We have $U_{[0]} = X_{\alpha_2} X_{\alpha_5} X_{\alpha_7}$. If we take the root element $u = x_{\alpha_2}(t)$ for some $t \in K$, and write
$$
A_1 =
\left(
  \begin{array}{cccc}
     1 &   &   &   \\
       & 1 &   &   \\
     t &   & 1 &   \\
       & t &   & 1 \\
  \end{array}
\right), \quad
A_2 =
\left(
  \begin{array}{cccc}
     1 & t &   &   \\
       & 1 &   &   \\
       &   & 1 & t \\
       &   &   & 1 \\
  \end{array}
\right), \quad
A_3 =
\left(
  \begin{array}{cccc}
     1 &   & t &   \\
       & 1 &   & t \\
       &   & 1 &   \\
       &   &   & 1 \\
  \end{array}
\right),
$$
then $J_1(u.v^{(1)}, u.v^{(2)}) = A_1 J_1(v^{(1)}, v^{(2)})$, $J_2(u.v^{(1)}, u.v^{(2)}) = J_2(v^{(1)}, v^{(2)}) A_2$ and $J_3(u.v^{(1)}, u.v^{(2)}) = J_3(v^{(1)}, v^{(2)}) A_3$. Similar equations hold for any root element $u = x_\alpha(t)$ where $\alpha \in \Phi_{[0]}$. Therefore $U_{[0]}$ preserves $\hat Y_1$; so given $u \in U_{[0]}$, by the above the set of weights occurring in $u.y$ has ZLCE. By Lemma~\ref{lem: gen height zero not strictly positive}, if we take $g \in \Tran_G(y, Y)$ and write $y' = g.y \in Y$, then we have $g = u_1 g' u_2$ with $u_1 \in C_{U_{[+]}}(y')$, $u_2 \in C_{U_{[+]}}(y)$, and $g' \in G_{[0]} N_{\Lambda(V)_{[0]}}$ with $g'.y = y'$. In particular $G.y \cap Y = G_{[0]} N_{\Lambda(V)_{[0]}}.y \cap Y$; moreover $C_G(y) = C_{U_{[+]}}(y) C_{G_{[0]} N_{\Lambda(V)_{[0]}}}(y) C_{U_{[+]}}(y)$.

First, since $W_{\Lambda(V)_{[0]}} = \langle w_{\beta_3}, w_{\beta_1}, w_{\beta_2}, w_{\beta_2 + \beta_4} w_{\beta_4 + \beta_5}, w_{\beta_5 + \beta_6} w_{\beta_6 + \beta_7} \rangle$ and $\beta_2 \in \Phi_{[0]}$, we have $G_{[0]} N_{\Lambda(V)_{[0]}} = G_{[0]} \langle n_{\beta_3}, n_{\beta_1}, n^*, n^{**} \rangle = G_1 (A \cap N) \langle n^*, n^{**} \rangle$. Any element of this last group may be written as $n' g^* c$ where $c \in A$, $g^* \in G_1$ and $n' \in \{ 1, n^*, n^{**}, n^* n^{**}, n^{**} n^*, n^* n^{**} n^* \}$; as $c.y = y$ it suffices to consider $n' g^*.y$. The above shows that applying any root element in $G_1$ has no effect on the determinants $\Delta_i$, so the same is true of $g^*$. We find that $n^*$ interchanges $e_{\gamma_3}$ and $e_{\gamma_5}$, and also $e_{\gamma_4}$ and $e_{\gamma_6}$, while fixing the other $e_{\gamma_i}$; likewise $n^{**}$ interchanges $e_{\gamma_2}$ and $e_{\gamma_3}$, and also $e_{\gamma_6}$ and $e_{\gamma_7}$, while fixing the other $e_{\gamma_i}$. Thus if we set
$$
M =
\left(
  \begin{array}{cccc}
     1 &   &   &   \\
       &   & 1 &   \\
       & 1 &   &   \\
       &   &   & 1 \\
  \end{array}
\right),
$$
and write $\pi^* = (1\ 2)$ and $\pi^{**} = (2\ 3)$, then for each $i$ we have $J_i(n^*.v^{(1)}, n^*.v^{(2)}) = J_{\pi^*(i)}(v^{(1)}, v^{(2)}) M$ and $J_i(n^{**}.v^{(1)}, n^{**}.v^{(2)}) = J_{\pi^{**}(i)}(v^{(1)}, v^{(2)}) M$; so applying $n'$ permutes the determinants $\Delta_i$. Thus $G_{[0]} N_{\Lambda(V)_{[0]}}.y \subset \hat Y_1$. If we now further require the element $n' g^* c$ to stabilize $y$, it must preserve the triple ratio $\Delta_1 : \Delta_2 : \Delta_3$ of determinants; the last condition in the definition of $\hat Y_1$ implies that we must have $n' = 1$, and so $g^*.y = y$. Since $V_{[0]}$ is the $G_1$-module with high weight $\omega_1 \otimes \omega_1 \otimes \omega_1$, using Proposition~\ref{prop: {A_1}^3, omega_1 otimes omega_1 otimes omega_1, k = 2} we see that there is a dense open subset $\hat Y_2$ of $Y$ each point of which has $G_1$-stabilizer a conjugate of $C'$. Set $\hat Y = \hat Y_1 \cap \hat Y_2$; then if $y \in \hat Y$ we see that $C_{G_{[0]} N_{\Lambda(V)_{[0]}}}(y) = {}^x C' (A \cap N)$ for some $x \in G_1$.

Next, let $\Xi = \Phi^+ \setminus (\Phi_{[0]} \cup \Phi_A)$, and set $U' = \prod_{\alpha \in \Xi} X_\alpha$; then $U_{[+]} = U'.(A \cap U_{[+]})$ and $U' \cap (A \cap U_{[+]}) = \{ 1 \}$. We now observe that if $\alpha \in \Xi$ then $\nu_i + \alpha$ is a weight in $V$ for exactly two values of $i$; moreover each weight in $V$ of positive generalized height is of the form $\nu_i + \alpha$ for exactly four such roots $\alpha$. Indeed $\Xi$ is the union of $12$ $W(G_1)$-orbits of size $4$; each such orbit is orthogonal to precisely one of $\alpha_2$, $\alpha_5$ and $\alpha_7$, and all $4$ roots $\alpha$ in the orbit give the same two weights $\nu_i + \alpha$. If we now take a product of root elements corresponding to the four roots in the orbit, and require it to stabilize $y$, equating coefficients of the corresponding two weight vectors in both basis vectors of $y$ gives $4$ linear equations which may be expressed in matrix form using one of the matrices $J_i(v^{(1)}, v^{(2)})$ above. For example, one such orbit is $\{ \alpha_6, \alpha_5 + \alpha_6, \alpha_6 + \alpha_7, \alpha_5 + \alpha_6 + \alpha_7 \}$, which is orthogonal to $\alpha_2$; here the two weights $\nu_i + \alpha$ are those corresponding to the roots $\delta_1 = \eeightrt11233321$ and $\delta_2 = \eeightrt12233321$. If we set $u = x_{\alpha_6}(t_1) x_{\alpha_5 + \alpha_6}(t_2) x_{\alpha_6 + \alpha_7}(t_3) x_{\alpha_5 + \alpha_6 + \alpha_7}(t_4)$, then we find that
$$
u.{\ts\sum a_i e_{\gamma_i}} = {\ts\sum a_i e_{\gamma_i}} + (a_4t_1 - a_2t_2 + a_3t_3 - a_1t_4) e_{\delta_1} + (a_8t_1 - a_6t_2 + a_7t_3 - a_5t_4) e_{\delta_2}.
$$
Equating to zero the coefficients of $e_{\delta_1}$ and $e_{\delta_2}$ in each $u.v^{(i)}$ now gives the equation $J_1(v^{(1)}, v^{(2)}) {\bf t} = {\bf 0}$, where ${\bf t} = ( {-t_4} \ {-t_2} \ t_3 \ t_1 )^T$; since the matrix concerned has non-zero determinant we see that $t_i = 0$ for $i = 1, 2, 3, 4$. Thus if we take $u = \prod x_\alpha(t_\alpha) \in U'$ satisfying $u.y = y$, and equate coefficients of weight vectors, taking them in an order compatible with increasing generalized height, we see that for all $\alpha$ we must have $t_\alpha = 0$, so that $u = 1$; so $C_{U_{[+]}}(y) = A \cap U_{[+]}$. Since the previous paragraph shows that $y' = g'.y \in \hat Y_1$, likewise we have $C_{U_{[+]}}(y') = A \cap U_{[+]}$.

As $G_{[0]} N_{\Lambda(V)_{[0]}} \subseteq \Tran_G(y, Y)$ and $C_{U_{[+]}}(y), C_{U_{[+]}}(y') \leq A$ we do indeed have $\Tran_G(y, Y) = A G_1 \langle n^*, n^{**} \rangle \cong D_4{A_1}^3.S_3$; so
$$
\codim {\ts\Tran_G(y, Y)} = \dim G - \dim {\ts\Tran_G(y, Y)} = 133 - 37 = 96
$$
while
$$
\codim Y = \dim \G{2}(V) - \dim Y = 108 - 12 = 96.
$$
Therefore $y$ is $Y$-exact. Now suppose further $y \in \hat Y$. As $C_{G_{[0]} N_{\Lambda(V)_{[0]}}}(y), C_{U_{[+]}}(y) \leq {}^x C' A$ for some $x \in G_1$ we have $C_G(y) = {}^x C' A = {}^x C$. Thus the conditions of Lemma~\ref{lem: generic stabilizer from exact subset} hold; so the quadruple $(G, \lambda, p, k)$ has generic stabilizer $C/Z(G) \cong D_4.\Z_{2/(p, 2)}.\Z_2$, where the $D_4$ is of simply connected type.

Next as in the proof of Proposition~\ref{prop: E_7, omega_7, D_6, omega_6, B_5, omega_5, A_5, omega_3, C_3, omega_3 modules} we may replace $G$ by the $D_6$ subgroup $\langle X_\alpha : \alpha = \sum m_i \beta_i, \ m_1 = m_8 = 0 \rangle$ of $E_7$, and $V$ by $\langle e_\alpha : \alpha = \sum m_i \beta_i, \ m_1 = m_8 = 1 \rangle$, leaving $Y$ and its subsets unchanged. We replace $A$ by the intersection of that above with $G$, which is the ${A_1}^3$ subgroup with simple roots $\beta_3$, $\rho_4$ and $\rho_6$; since $Z(A) = \langle z_1, z_2, z_3 \rangle$ where $z_1 = h_{\beta_3}(-1)$, $z_2 = h_{\beta_2}(-1) h_{\beta_3}(-1) h_{\beta_5}(-1)$ and $z_3 = h_{\beta_2}(-1) h_{\beta_3}(-1) h_{\beta_7}(-1)$, we see that $A$ is of simply connected type. We again let $C = C' A$; then $Z(C) = \langle h_{\beta_2}(-1), h_{\beta_3}(-1), h_{\beta_5}(-1), h_{\beta_7}(-1) \rangle = Z(G) Z(A)$, where $Z(G) = \langle h_{\beta_2}(-1) h_{\beta_3}(-1), h_{\beta_3}(-1) h_{\beta_5}(-1) h_{\beta_7}(-1) \rangle$. As $G_1 \langle n^*, n^{**} \rangle < G$, for all $y \in \hat Y_1$ we have $\Tran_G(y, Y) = A G_1 \langle n^*, n^{**} \rangle \cong {A_1}^3{A_1}^3.S_3$; so
$$
\codim {\ts\Tran_G(y, Y)} = \dim G - \dim {\ts\Tran_G(y, Y)} = 66 - 18 = 48
$$
while
$$
\codim Y = \dim \G{2}(V) - \dim Y = 60 - 12 = 48.
$$
Therefore $y$ is $Y$-exact. Also if $y \in \hat Y$ then $C_G(y) = {}^x C' A = {}^x C$ for some $x \in G_1$. Thus the conditions of Lemma~\ref{lem: generic stabilizer from exact subset} hold; so the quadruple $(G, \lambda, p, k)$ has generic stabilizer $C/Z(G) \cong {A_1}^3.\Z_{2/(p, 2)}.\Z_2$, where the ${A_1}^3$ is a central product.

Now we replace $G$ by the $A_5$ subgroup $\langle X_\alpha : \alpha = \sum m_i \beta_i, \ m_1 = m_3 = m_8 = 0 \rangle$, and $V$ by $\langle e_\alpha : \alpha = \sum m_i \beta_i, \ m_1 = m_8 = 1, \ m_3 = 2 \rangle$, again leaving $Y$ and its subsets unchanged. We replace $A$ by the intersection of that above with $G$, which is the $T_2$ subgroup $\{ h_{\beta_2}(\kappa_1) h_{\beta_4}({\kappa_1}^2) h_{\beta_5}(\kappa_1 \kappa_2) h_{\beta_6}({\kappa_2}^2) h_{\beta_7}(\kappa_2) : \kappa_1, \kappa_2 \in K^* \}$. We let $C = C' A$. As $G_1 \langle n^*, n^{**} \rangle < G$, for all $y \in \hat Y_1$ we have $\Tran_G(y, Y) = A G_1 \langle n^*, n^{**} \rangle \cong T_2{A_1}^3.S_3$; so
$$
\codim {\ts\Tran_G(y, Y)} = \dim G - \dim {\ts\Tran_G(y, Y)} = 35 - 11 = 24
$$
while
$$
\codim Y = \dim \G{2}(V) - \dim Y = 36 - 12 = 24.
$$
Therefore $y$ is $Y$-exact. Also if $y \in \hat Y$ then $C_G(y) = {}^x C' A = {}^x C$ for some $x \in G_1$. Thus the conditions of Lemma~\ref{lem: generic stabilizer from exact subset} hold; so the quadruple $(G, \lambda, p, k)$ has generic stabilizer $C/Z(G) \cong T_2.\Z_{2/(p, 2)}.\Z_2$.

Finally for $p \geq 3$ we replace $G$ by the $C_3$ subgroup with simple root groups $\{ x_{\beta_2 + \beta_4}(t) x_{-(\beta_4 + \beta_5)}(t) : t \in K \}$, $\{ x_{\beta_5 + \beta_6}(t) x_{-(\beta_6 + \beta_7)}(t) : t \in K \}$ and $X_{\beta_7}$, and $V$ by its submodule which has highest weight $\omega_3$, again leaving $Y$ and its subsets unchanged. We replace $A$ by the intersection of that above with $G$, which is the ${\Z_2}^2$ subgroup $\langle h_{\beta_2}(-1) h_{\beta_5}(-1), h_{\beta_5}(-1) h_{\beta_7}(-1) \rangle$. We let $C = C' A = Z(G) A \langle h^\dagger, n^\dagger \rangle$. As $G_1 \langle n^*, n^{**} \rangle < G$, for all $y \in \hat Y_1$ we have $\Tran_G(y, Y) = A G_1 \langle n^*, n^{**} \rangle \cong {\Z_2}^2{A_1}^3.S_3$; so
$$
\codim {\ts\Tran_G(y, Y)} = \dim G - \dim {\ts\Tran_G(y, Y)} = 21 - 9 = 12
$$
while
$$
\codim Y = \dim \G{2}(V) - \dim Y = 24 - 12 = 12.
$$
Therefore $y$ is $Y$-exact. Also if $y \in \hat Y$ then $C_G(y) = {}^x C' A = {}^x C$ for some $x \in G_1$. Thus the conditions of Lemma~\ref{lem: generic stabilizer from exact subset} hold; so the quadruple $(G, \lambda, p, k)$ has generic stabilizer $C/Z(G) \cong {\Z_2}^4$.
\end{proof}

\begin{prop}\label{prop: B_5, omega_5 module, k = 2}
Let $G = B_5$ and $\lambda = \omega_5$, and take $k = 2$. Then the quadruple $(G, \lambda, p, k)$ has generic stabilizer $\Z_{2/(p, 2)}.\Z_2$.
\end{prop}

\begin{proof}
We use the set-up of (the relevant part of) Proposition~\ref{prop: E_7, omega_7, D_6, omega_6, A_5, omega_3, C_3, omega_3 modules, k = 2}: again we take $H$ to be the (simply connected) group defined over $K$ of type $E_8$, with simple roots $\beta_1, \dots, \beta_8$; we take the $D_6$ subgroup $\langle X_\alpha : \alpha = \sum m_i \beta_i, \ m_1 = m_8 = 0 \rangle < H$; then we may take $V = \langle e_\alpha : \alpha = \sum m_i \beta_i, \ m_1 = m_8 = 1 \rangle < \L(H)$; we have $Z(D_6) = \langle z_1, z_2 \rangle$ where $z_1 = h_{\beta_2}(-1) h_{\beta_3}(-1)$, $z_2 = h_{\beta_3}(-1) h_{\beta_5}(-1) h_{\beta_7}(-1)$. We write
\begin{eqnarray*}
& \gamma_1 = \eeightrt11232211, \quad \gamma_2 = \eeightrt11232221, \quad \gamma_3 = \eeightrt11233211, \quad \gamma_4 = \eeightrt11233221, & \\
& \gamma_5 = \eeightrt12232211, \quad \gamma_6 = \eeightrt12232221, \quad \gamma_7 = \eeightrt12233211, \quad \gamma_8 = \eeightrt12233221; &
\end{eqnarray*}
here in addition we write
\begin{eqnarray*}
& \delta_1 = \eeightrt11121111, \quad \delta_2 = \eeightrt11122111, \quad \delta_3 = \eeightrt11122211, \quad \delta_4 = \eeightrt11122221, & \\
& \delta_5 = \eeightrt12243211, \quad \delta_6 = \eeightrt12243221, \quad \delta_7 = \eeightrt12243321, \quad \delta_8 = \eeightrt12244321. &
\end{eqnarray*}
Let $A$ be the ${A_1}^3$ subgroup with simple roots $\beta_3$, $\rho_4 = \eeightrt01121000$ and $\rho_6 = \eeightrt01122210$, and $G_1$ be the ${A_1}^3$ subgroup with simple roots $\beta_2$, $\beta_5$ and $\beta_7$; then $z_2 \in A$. Write $h^\dagger = h_{\beta_2}(\eta_4) h_{\beta_5}(\eta_4) h_{\beta_7}(\eta_4)$ and $n^\dagger = n_{\beta_2} n_{\beta_5} n_{\beta_7}$, and $h_0 = h_{\beta_3}(\eta_4) h_{\rho_4}(\eta_4) h_{\rho_6}(\eta_4) h^\dagger \in A h^\dagger$ and $n_0 = {n_{\beta_3}}^{-1} {n_{\rho_4}}^{-1}{n_{\rho_6}}^{-1} n^\dagger \in A n^\dagger$; then ${h_0}^2 = {n_0}^2 = [h_0, n_0] = z_1$. Set $n^* = n_{\beta_2 + \beta_4} n_{\beta_4 + \beta_5}$ and $n^{**} = n_{\beta_5 + \beta_6} n_{\beta_6 + \beta_7}$.

We saw in the proof of Proposition~\ref{prop: E_7, omega_7, D_6, omega_6, A_5, omega_3, C_3, omega_3 modules, k = 2} that if we let $Y_0 = \G{2}(\langle e_{\gamma_1}, \dots, e_{\gamma_8} \rangle)$ then there is a dense open subset $\hat Y_0$ of $Y_0$ such that if $y \in \hat Y_0$ then $\Tran_{D_6}(y, Y_0) = A G_1 \langle n^*, n^{**} \rangle$ and $C_{D_6}(y)$ is a $G_1$-conjugate of $Z(D_6) A \langle h^\dagger, n^\dagger \rangle = A \langle h^\dagger, n^\dagger \rangle$. Given $\c = (c_1, c_2, c_3) \in K^3$, define
$$
y_\c = \langle e_{\gamma_1} + c_1 e_{\gamma_4} + c_2 e_{\gamma_6} + c_3 e_{\gamma_7}, e_{\gamma_8} + c_1 e_{\gamma_5} + c_2 e_{\gamma_3} + c_3 e_{\gamma_2} \rangle \in Y_0;
$$
write ${Y_0}' = \{ y_\c : \c \in K^3 \}$ and $\hat Y_0{}' = {Y_0}' \cap \hat Y_0$. From the proofs of Propositions~\ref{prop: {A_1}^3, omega_1 otimes omega_1 otimes omega_1, k = 2} and \ref{prop: E_7, omega_7, D_6, omega_6, A_5, omega_3, C_3, omega_3 modules, k = 2} we see that $\hat Y_0{}' \neq \emptyset$, and if $y_\c \in \hat Y_0{}'$ then $C_{G_1}(y_\c) = Z(G_1) \langle h^\dagger, n^\dagger \rangle$, so that $C_{D_6}(y_\c) = A Z(G_1) \langle h^\dagger, n^\dagger \rangle = A \langle h_0, n_0 \rangle$; moreover $\Tran_{G_1}(y_\c, {Y_0}')$ is finite, so that $\Tran_{D_6}(y_\c, {Y_0}')$ is a finite union of left cosets of $A$, each of which lies in $A G_1 \langle n^*, n^{**} \rangle$. Note that as $A \lhd A G_1 \langle n^*, n^{**} \rangle = \Tran_{D_6}(y_\c, Y_0)$, each left coset of $A$ in $\Tran_{D_6}(y_\c, {Y_0}')$ is also a right coset.

Given $\a = (a_1, a_2, a_3, a_4, a_5, a_6) \in K^6$, write $|\a| = \sum_{i = 1}^6 {a_i}^2$. For $\a \in K^6$ with $|\a| = 1$ and $a_6 \neq 0$, and $\c \in K^3$ as above, set
\begin{eqnarray*}
y_{\a, \c} & = & \langle e_{\gamma_1} + c_1 e_{\gamma_4} + a_6 c_2 e_{\gamma_6} + a_6 c_3 e_{\gamma_7} - a_5 e_{\gamma_5} - a_5 c_1 e_{\gamma_8} - a_2 e_{\delta_1} - a_1 c_1 e_{\delta_2} \\
           &   & \quad {} + a_4 e_{\delta_3} + a_3 c_1 e_{\delta_4} + a_3 e_{\delta_5} - a_4 c_1 e_{\delta_6} - a_1 e_{\delta_7} + a_2 c_1 e_{\delta_8}, \\
           &   & \phantom{\langle} a_6 e_{\gamma_8} + a_6 c_1 e_{\gamma_5} + c_2 e_{\gamma_3} + c_3 e_{\gamma_2} - a_5 c_2 e_{\gamma_7} - a_5 c_3 e_{\gamma_6} - a_1 c_3 e_{\delta_1} - a_2 c_2 e_{\delta_2} \\
           &   & \quad {} + a_3 c_2 e_{\delta_3} + a_4 c_3 e_{\delta_4} - a_4 c_2 e_{\delta_5} + a_3 c_3 e_{\delta_6} + a_2 c_3 e_{\delta_7} - a_1 c_2 e_{\delta_8} \rangle.
\end{eqnarray*}
An easy check shows that, if $\a' \in K^6$ with $|\a'| = 1$ and ${a_6}' \neq 0$, and $\c' \in K^3$, then $y_{\a', \c'} = y_{\a, \c}$ if and only if either $\a' = \a$, $\c' = \c$ or ${a_i}' = a_i$ for $i \leq 5$, ${a_6}' = -a_6$, ${c_1}' = c_1$, ${c_2}' = -c_2$, ${c_3}' = -c_3$. Thus if we set
$$
Y = \{ y_{\a, \c} : \a \in K^6, \c \in K^3, |\a| = 1, a_6 \neq 0 \}
$$
then $\dim Y = 8$. Choose $\xi \in K^*$ with $\xi^2 = {a_6}^{-1}$, and define
\begin{eqnarray*}
g_\a & = & h_{\beta_2}(\xi) h_{-\beta_3}(\xi) x_{\beta_2}(a_5) x_{-\beta_3}(-a_5) x_{\beta_2 + \beta_4}(-a_4) x_{-(\beta_3 + \beta_4)}(a_4) \\
     &   & {} \times x_{\beta_2 + \beta_4 + \beta_5}(a_3) x_{-(\beta_3 + \beta_4 + \beta_5)}(-a_3) \\
     &   & {} \times x_{\beta_2 + \beta_4 + \beta_5 + \beta_6}(-a_2) x_{-(\beta_3 + \beta_4 + \beta_5 + \beta_6)}(a_2) \\
     &   & {} \times x_{\beta_2 + \beta_4 + \beta_5 + \beta_6 + \beta_7}(a_1) x_{-(\beta_3 + \beta_4 + \beta_5 + \beta_6 + \beta_7)}(-a_1)
\end{eqnarray*}
(note that the two choices for $\xi$ give elements differing by $z_1$, which fixes all points in $\G{2}(V)$); then calculation shows that
$$
g_\a.y_{\a, \c} = y_\c.
$$

At this point we find it convenient to switch notation. Instead of taking the root system of $D_6$ to be a subsystem of that of $E_8$, we shall use the standard notation given in Section~\ref{sect: notation}; thus we replace $\beta_7$, $\beta_6$, $\beta_5$, $\beta_4$, $\beta_2$ and $\beta_3$ by $\ve_1 - \ve_2$, $\ve_2 - \ve_3$, $\ve_3 - \ve_4$, $\ve_4 - \ve_5$, $\ve_5 - \ve_6$ and $\ve_5 + \ve_6$ respectively, and we recall the natural module $V_{nat}$ for $D_6$. However, just as in the proof of Proposition~\ref{prop: B_4, omega_4 module, k = 3} there is an unfortunate consequence to this change: in Section~\ref{sect: notation} we defined the action of root elements on $V_{nat}$, which implicitly determined the structure constants, and these are not the same as those given in the appendix of \cite{LSmax}, which we have been using until now. For this reason we shall avoid all mention of root elements from now on, but rather identify elements of $D_6$ by their action on $V_{nat}$ (the kernel of this action is $\langle z_1 \rangle$, so this is harmless). Thus with respect to the ordered basis $v_1, v_2, v_3, v_4, v_5, v_6, v_{-6}, v_{-5}, v_{-4}, v_{-3}, v_{-2}, v_{-1}$ of $V_{nat}$, the element $g_\a$ defined above acts as
$$
\left(
  \begin{array}{cccccc|cccccc}
      1    &         &         &         &         &         a_1         &               &         &         &         &         &         \\
           &    1    &         &         &         &         a_2         &               &         &         &         &         &         \\
           &         &    1    &         &         &         a_3         &               &         &         &         &         &         \\
           &         &         &    1    &         &         a_4         &               &         &         &         &         &         \\
           &         &         &         &    1    &         a_5         &               &         &         &         &         &         \\
           &         &         &         &         &         a_6         &               &         &         &         &         &         \\
  \hline
   -\aaf16 & -\aaf26 & -\aaf36 & -\aaf46 & -\aaf56 & a_6 - \frac{1}{a_6} & \frac{1}{a_6} & -\aaf56 & -\aaf46 & -\aaf36 & -\aaf26 & -\aaf16 \\
           &         &         &         &         &         a_5         &               &    1    &         &         &         &         \\
           &         &         &         &         &         a_4         &               &         &    1    &         &         &         \\
           &         &         &         &         &         a_3         &               &         &         &    1    &         &         \\
           &         &         &         &         &         a_2         &               &         &         &         &    1    &         \\
           &         &         &         &         &         a_1         &               &         &         &         &         &    1    \\
  \end{array}
\right)
$$

For $i = 1, 2, 3$ write $V_{2i - 1, 2i} = \langle v_{2i - 1}, v_{-(2i - 1)}, v_{2i}, v_{-2i} \rangle$; then we have $V_{nat} = V_{1, 2} \oplus V_{3, 4} \oplus V_{5, 6}$. We see that $A G_1 \langle n^*, n^{**} \rangle = {D_2}^3.S_3$, where the three $D_2$ factors have roots $\pm\ve_{2i - 1} \pm \ve_{2i}$ for $i = 1, 2, 3$ and thus act on $V_{2i - 1, 2i}$, and the $S_3$ permutes these three summands.

Now write
$$
v^\diamondsuit = v_6 + v_{-6},
$$
and let $G = C_{D_6}(v^\diamondsuit) = B_5$; then $Z(G) = \langle z_1 \rangle$. Since the elements $h_0$ and $n_0$ defined above have the property that for each $i \in \{ 1, \dots, 6 \}$ we have $h_0.v_i = (-1)^i v_i$ and $h_0.v_{-i} = (-1)^i v_{-i}$, and $n_0.v_i = v_{-i}$ and $n_0.v_{-i} = v_i$, we see that $h_0, n_0 \in G$. For $\a \in K^6$ with $|\a| = 1$, define
$$
v_\a = g_\a.v^\diamondsuit = {\ts\sum_{i = 1}^6} a_i(v_i + v_{-i});
$$
then $v_\a$ is a vector of norm $1$ fixed by $n_0$ (and $v^\diamondsuit = v_\a$ where $\a = (0, 0, 0, 0, 0, 1)$). Write
$$
V_* = \{ v_\a : |\a| = 1 \}.
$$
Define
\begin{eqnarray*}
S_A    & = & \{ \a \in K^6 : |\a| = 1, \ {a_{2i - 1}}^2 + {a_{2i}}^2 \neq 0, 1 \hbox { for } i = 1, 2, 3 \}, \\
{S_A}' & = & \{ \a \in S_A : a_6 \neq 0 \}, \\
S_C    & = & \{ \c \in K^3 : y_\c \in \hat Y_0{}', \ c_1c_2c_3 \neq 0 \},
\end{eqnarray*}
and set
$$
\hat Y = \{ y_{\a, \c} \in Y : \a \in {S_A}', \c \in S_C \}, \quad \hat V_* = \{ v_\a \in V_* : \a \in {S_A}' \};
$$
then $\hat Y$ and $\hat V_*$ are open dense subsets of $Y$ and $V_*$ respectively.

Take $y_{\a, \c} \in \hat Y$ and suppose $g \in \Tran_G(y_{\a, \c}, Y)$; write $g.y_{\a, \c} = y_{\a'', \c''}$ and set $g' = g_{\a''} g {g_\a}^{-1} \in D_6$. Then $g'.y_\c = g_{\a''} g {g_\a}^{-1}.y_\c = g_{\a''} g.y_{\a, \c} = g_{\a''}.y_{\a'', \c''} = y_{\c''}$, and $g'.v_\a = g_{\a''} g {g_\a}^{-1}.v_\a = g_{\a''} g.v^\diamondsuit = g_{\a''}.v^\diamondsuit = v_{\a''}$ since $g \in G$; so any element of $\Tran_G(y_{\a, \c}, Y)$ is of the form ${g_{\a''}}^{-1} g' g_\a$, where $g' \in \Tran_{D_6}(y_\c, {Y_0}')$ and $g'.v_\a = v_{\a''} \in V_*$. By the above $\Tran_{D_6}(y_\c, {Y_0}') = \bigcup_{j = 1}^n A x_j$ for some finite set $\{ x_1, \dots, x_n \}$ of elements of ${D_2}^3.S_3$. Take $j \in \{ 1, \dots, n \}$ and write $x_j.y_\c = y_{\c'} \in {Y_0}'$. Since $\a \in S_A$, the projection of $v_\a$ on each of $V_{1, 2}$, $V_{3, 4}$ and $V_{5, 6}$ is a non-singular vector, so as $x_j \in {D_2}^3.S_3$ the same is true of $x_j.v_\a$. Since for $i = 1, 2, 3$ the $A_1$ subgroup with roots $\pm(\ve_{2i - 1} + \ve_{2i})$ acts simply transitively on the set of non-singular vectors of a given norm lying in $V_{2i - 1, 2i}$, we see that the coset $Ax_j = x_jA$ has intersection with $\Tran_{D_6}(v_\a, V_*)$ of dimension $3$. For each element $g'$ lying in this intersection, we have $g' g_\a.y_{\a, \c} = y_{\c'}$ and $g' g_\a.v^\diamondsuit = v_{\a'}$ for some $\a' \in S_A$; for those with $\a' \in {S_A}'$ we have ${g_{\a'}}^{-1} g' g_\a \in G$ with ${g_{\a'}}^{-1} g' g_\a.y_{\a, \c} = y_{\a', \c'}$. Since distinct elements $g'$ give distinct vectors $v_{\a'}$ and thus distinct subspaces $y_{\a', \c'}$ (note that the definition of the set $S_C$ implies that no component of $\c'$ can be zero, so we cannot have ${c_2}' = -{c_2}'$, ${c_3}' = -{c_3}'$), we see that the elements of $\Tran_G(y_{\a, \c}, Y)$ arising from the coset $Ax_j$ form a $3$-dimensional variety. Since this is true for each $j \in \{ 1, \dots, n \}$, we have $\dim \Tran_G(y_{\a, \c}, Y) = 3$. Thus
$$
\codim {\ts\Tran_G(y_{\a, \c}, Y)} = \dim G - \dim {\ts\Tran_G(y_{\a, \c}, Y)} = 55 - 3 = 52
$$
while
$$
\codim Y = \dim \G{2}(V) - \dim Y = 60 - 8 = 52.
$$
Therefore $y_{\a, \c}$ is $Y$-exact.

Now suppose $g \in C_G(y_{\a, \c})$; as above if we set $g' = {}^{g_\a} g \in D_6$ then $g'$ fixes both $y_\c$ and $v_\a$. The first of these conditions implies $g' \in A \langle h_0, n_0 \rangle$, in which $n_0$ fixes $v_\a$; the simple transitivity of $A$ on triples of non-singular vectors of given norms in $V_{1, 2}$, $V_{3, 4}$ and $V_{5, 6}$ gives $C_A(v_\a) = \{ 1 \}$, while if $p \neq 2$ the coset $Ah_0$ contains a unique element fixing $v_\a$. Indeed, in this case a straightforward calculation shows that if for $i = 1, 2, 3$ we write the basis elements of $V_{2i - 1, 2i}$ in the order $v_{2i - 1}, v_{-2i}, v_{2i}, v_{-(2i - 1)}$, and take $\kappa_i \in K^*$ with ${\kappa_i}^2 = {a_{2i - 1}}^2 + {a_{2i}}^2$, then the element $x_\a$ of $A$ acting on $V_{2i - 1, 2i}$ as
$$
{\ts \frac{1}{\kappa_i}}
\left(
  \begin{array}{cc|cc}
      a_{2i}   & -a_{2i - 1} &             &            \\
    a_{2i - 1} &    a_{2i}   &             &            \\
    \hline
               &             &    a_{2i}   & a_{2i - 1} \\
               &             & -a_{2i - 1} &   a_{2i}   \\
  \end{array}
\right)
$$
commutes with $n_0$, and the element ${h_0}^{x_\a}$ of $Ah_0$ fixes $v_\a$. Thus $C_G(y_{\a, \c}) = \langle h_0, n_0 \rangle^{x_\a g_\a}$ (and this also holds if $p = 2$ since then $h_0 = 1$). Observe that $x_\a g_\a.v^\diamondsuit = x_\a.v_\a = v_{\a'}$ where $\a' = (0, \kappa_1, 0, \kappa_2, 0, \kappa_3)$; as $\a \in S_A$ we have ${\kappa_1}^2 + {\kappa_2}^2 \neq 0$. Now the $D_3$ with roots $\pm\ve_2 \pm \ve_4, \pm\ve_4 \pm \ve_6, \pm\ve_2 \pm \ve_6$ commutes with $h_0$. Take $\kappa \in K^*$ with $\kappa^2 = {\kappa_1}^2 + {\kappa_2}^2$, and let ${x_\a}'$ and ${x_\a}''$ be the elements of $\langle X_{\pm(\ve_2 - \ve_4)} \rangle$ and $\langle X_{\pm(\ve_4 - \ve_6)} \rangle$ which act on $\langle v_2, v_4, v_{-4}, v_{-2} \rangle$ and $\langle v_4, v_6, v_{-6}, v_{-4} \rangle$ respectively as
$$
{\ts \frac{1}{\kappa}}
\left(
  \begin{array}{cc|cc}
   \kappa_2 & -\kappa_1 &           &          \\
   \kappa_1 &  \kappa_2 &           &          \\
    \hline
            &           &  \kappa_2 & \kappa_1 \\
            &           & -\kappa_1 & \kappa_2 \\
  \end{array}
\right)
\quad \hbox{and} \quad
\left(
  \begin{array}{cc|cc}
   \kappa_3 & -\kappa  &          &          \\
    \kappa  & \kappa_3 &          &          \\
    \hline
            &          & \kappa_3 &  \kappa  \\
            &          & -\kappa  & \kappa_3 \\
  \end{array}
\right);
$$
then ${x_\a}'$ and ${x_\a}''$ both commute with both $h_0$ and $n_0$, and ${x_\a}'.v_{\a'} = v_{\a''}$ where $\a'' = (0, 0, 0, \kappa, 0, \kappa_3)$, while ${x_\a}''.v_{\a''} = v^\diamondsuit$. Therefore $C_G(y_{\a, \c}) = \langle h_0, n_0 \rangle^x$ where $x = {x_\a}'' {x_\a}' x_\a g_\a$, and as $x.v^\diamondsuit = {x_\a}'' {x_\a}' x_\a g_\a.v^\diamondsuit = {x_\a}'' {x_\a}' x_\a.v_\a = {x_\a}'' {x_\a}'.v_{\a'} = {x_\a}''.v_{\a''} = v^\diamondsuit$ we have $x \in G$. Thus if we let $C = \langle h_0, n_0 \rangle$ then the conditions of Lemma~\ref{lem: generic stabilizer from exact subset} hold; so the quadruple $(G, \lambda, p, k)$ has generic stabilizer $C/Z(G) \cong \Z_{2/(p, 2)}.\Z_2$.
\end{proof}

\begin{prop}\label{prop: C_5, omega_5 module, p = 2, k = 2}
Let $G = C_5$ and $\lambda = \omega_5$ with $p = 2$, and take $k = 2$. Then the quadruple $(G, \lambda, p, k)$ has generic stabilizer $\Z_2$.
\end{prop}

\begin{proof}
This is an immediate consequence of Proposition~\ref{prop: B_5, omega_5 module, k = 2}, using the exceptional isogeny $B_\ell \to C_\ell$ which exists in characteristic $2$.
\end{proof}

\begin{prop}\label{prop: G_2, omega_1 module, k = 2}
Let $G = G_2$ and $\lambda = \omega_1$ with $p \geq 3$ or $p = 2$, and take $k = 2$. Then the quadruple $(G, \lambda, p, k)$ has generic stabilizer $A_1T_1.\Z_2$ or $A_1 \tilde A_1$ respectively.
\end{prop}

\begin{proof}
We begin with the case where $p \geq 3$; here $\dim V = 7$ and $\Lambda(V) = \Phi_s \cup \{ 0 \}$. We take an ordered basis of $V$ consisting of weight vectors $v_\mu$ for the weights $\mu = 2\alpha_1 + \alpha_2$, $\alpha_1 + \alpha_2$, $\alpha_1$, $0$, $-\alpha_1$, $-(\alpha_1 + \alpha_2)$, $-(2\alpha_1 + \alpha_2)$ respectively, such that with respect to them the simple root elements $x_{\alpha_1}(t)$ and $x_{\alpha_2}(t)$ and the corresponding negative root elements $x_{-\alpha_1}(t)$ and $x_{-\alpha_2}(t)$ act by the matrices given in the proof of Proposition~\ref{prop: G_2, omega_1 module}.

We take the generalized height function on the weight lattice of $G$ whose value at $\alpha_1$ is $0$, and at $\alpha_2$ is $1$; then the generalized height of $\lambda = 2\alpha_1 + \alpha_2$ is $1$, and as $\Phi$ generates the weight lattice we see that the generalized height of any weight is an integer. We have $\Lambda(V)_{[0]} = \{ \pm\alpha_1, 0 \}$; evidently both $\{ \pm\alpha_1 \}$ and $\{ \pm\alpha_1, 0 \}$ have ZLCE. Take $Y = \G{2}(V_{[0]})$ and write
$$
y_0 = \langle v_{\alpha_1}, v_{-\alpha_1} \rangle \in Y.
$$
Clearly the setwise stabilizer in $W$ of $\Lambda(V)_{[0]}$ is $\langle w_{\alpha_1}, w_{3\alpha_1 + 2\alpha_2} \rangle$. Note that this stabilizes $\Phi_{[0]} = \langle \alpha_1 \rangle$.

Let $A$ be the $A_1$ subgroup having simple root $3\alpha_1 + 2\alpha_2$; then we have $Z(A) = \langle h_{3\alpha_1 + 2\alpha_2}(-1) \rangle$. Write $T_1 = \{ h_{\alpha_1}(\kappa) : \kappa \in K^* \}$ and set $C = AT_1 \langle n_{\alpha_1} \rangle$. Clearly we have $C \leq C_G(y_0)$; we shall show that in fact $C_G(y_0) = C$.

We have $U_{[0]} = X_{\alpha_1}$. Given $u \in U_{[0]}$, the weights $\alpha_1$ and $-\alpha_1$ occur in $u.v_{\alpha_1}$ and $u.v_{-\alpha_1}$ respectively, so the set of weights occurring in $u.y_0$ contains $\pm\alpha_1$, and hence has ZLCE. By Lemma~\ref{lem: gen height zero not strictly positive}, we have $C_G(y_0) = C_{U_{[+]}}(y_0) C_{G_{[0]} N_{\Lambda(V)_{[0]}}}(y_0) C_{U_{[+]}}(y_0)$.

First, since $W_{\Lambda(V)_{[0]}} = \langle w_{\alpha_1}, w_{3\alpha_1 + 2\alpha_2} \rangle$ and $\alpha_1 \in \Phi_{[0]}$, we have $G_{[0]} N_{\Lambda(V)_{[0]}} = G_{[0]} \langle n_{3\alpha_1 + 2\alpha_2} \rangle$. Any element of this last group may be written as $g^*c$ where $c \in \langle X_{\pm(3\alpha_1 + 2\alpha_2)} \rangle \cap N < C$ and $g^* \in \langle X_{\pm\alpha_1} \rangle$. Suppose then that $g^* \in C_G(y_0)$. If $g^* = x_{\alpha_1}(t) h_{\alpha_1}(\kappa)$ for some $t \in K$ and $\kappa \in K^*$, then we must have $t = 0$ as otherwise $g^*.v_{-\alpha_1}$ has a term $v_0$; thus $g^* \in T_1$. If instead $g^* = x_{\alpha_1}(t) h_{\alpha_1}(\kappa) n_{\alpha_1} x_{\alpha_1}(t')$ for some $t, t' \in K$ and $\kappa \in K^*$, then we must have $t = 0$ as otherwise $g^*.v_{\alpha_1}$ has a term $v_0$, and then we must have $t' = 0$ as otherwise $g^*.v_{-\alpha_1}$ has a term $v_0$; thus $g^* \in T_1 n_{\alpha_1}$. Therefore we have $g^* \in T_1 \langle n_{\alpha_1} \rangle < C$; so $C_{G_{[0]} N_{\Lambda(V)_{[0]}}}(y_0) = C \cap G_{[0]} N_{\Lambda(V)_{[0]}}$.

Next, let $\Xi = \Phi^+ \setminus \{ \alpha_1, 3\alpha_1 + 2\alpha_2 \}$, and set $U' = \prod_{\alpha \in \Xi} X_\alpha$; then $U_{[+]} = U'.(C \cap U_{[+]})$ and $U' \cap (C \cap U_{[+]}) = \{ 1 \}$. Now take $u = \prod_{\alpha \in \Xi} x_\alpha(t_\alpha) \in U'$ satisfying $u.y_0 = y_0$. The requirement that in $u.v_{\alpha_1}$ the coefficients of $v_{2\alpha_1 + \alpha_2}$ and $v_{\alpha_1 + \alpha_2}$ should be zero shows that $t_\alpha = 0$ for $\alpha = \alpha_1 + \alpha_2$ and $\alpha_2$ respectively; considering likewise $u.v_{-\alpha_1}$ we see that the same is true for $\alpha = 3\alpha_1 + \alpha_2$ and $2\alpha_1 + \alpha_2$ respectively. Hence $u = 1$, so $C_{U_{[+]}}(y_0) = C \cap U_{[+]}$.

Therefore $C_G(y_0) = (C \cap U_{[+]})(C \cap G_{[0]} N_{\Lambda(V)_{[0]}})(C \cap U_{[+]}) \leq C$, so that we do indeed have $C_G(y_0) = C$.

Since $\dim(\overline{G.y_0}) = \dim G - \dim C_G(y_0) = 14 - 4 = 10 = \dim \G{2}(V)$, the orbit $G.y_0$ is dense in $\G{2}(V)$. Thus the quadruple $(G, \lambda, p, k)$ has generic stabilizer $C_G(y_0)/Z(G) \cong A_1 T_1.\Z_2$, where the $A_1$ is of simply connected type.

Now take the case where $p = 2$; here $\Lambda(V) = \Phi_s$. Again write
$$
y_0 = \langle v_{\alpha_1}, v_{-\alpha_1} \rangle.
$$
Let $C$ be the $A_1 \tilde A_1$ subgroup having simple roots $3\alpha_1 + 2\alpha_2$ and $\alpha_1$; clearly we have $C \leq C_G(y_0)$, and as $C$ is a maximal subgroup we must have $C_G(y_0) = C$. Since $\dim(\overline{G.y_0}) = \dim G - \dim C_G(y_0) = 14 - 6 = 8 = \dim \G{2}(V)$, the orbit $G.y_0$ is dense in $\G{2}(V)$. Thus the quadruple $(G, \lambda, p, k)$ has generic stabilizer $C_G(y_0)/Z(G) \cong A_1 \tilde A_1$.
\end{proof}

\begin{prop}\label{prop: G_2, omega_1 module, k = 3}
Let $G = G_2$ and $\lambda = \omega_1$ with $p \geq 3$ or $p = 2$, and take $k = 3$. Then the quadruple $(G, \lambda, p, k)$ has generic stabilizer $A_1$ or $A_1 U_2$ respectively.
\end{prop}

\begin{proof}
As in the proof of Proposition~\ref{prop: G_2, omega_1 module, k = 2}, we shall assume throughout that the basis of weight vectors $v_\mu$ of $V$ is chosen such that the elements of $G$ act by the matrices given in the proof of Proposition~\ref{prop: G_2, omega_1 module} (where if $p = 2$ the fourth row and column are deleted).

Again we begin with the case where $p \geq 3$; here $\Lambda(V) = \Phi_s \cup \{ 0 \}$. We take the same generalized height function on the weight lattice of $G$ as in Proposition~\ref{prop: G_2, omega_1 module, k = 2}, so that $\Lambda(V) = \Lambda(V)_{[-1]} \cup \Lambda(V)_{[0]} \cup \Lambda(V)_{[1]}$ where
\begin{eqnarray*}
\Lambda(V)_{[-1]} & = & \{ -(2\alpha_1 + \alpha_2), -(\alpha_1 + \alpha_2) \}, \\
\Lambda(V)_{[0]}  & = & \{ -\alpha_1, 0, \alpha_1 \}, \\
\Lambda(V)_{[1]}  & = & \{ \alpha_1 + \alpha_2, 2\alpha_1 + \alpha_2 \}.
\end{eqnarray*}
Here we set
$$
Y = \left\{ \langle v_{-(\alpha_1 + \alpha_2)}, a_1 v_{-\alpha_1} + a_2 v_0 + a_3 v_{\alpha_1}, v_{\alpha_1 + \alpha_2} \rangle : (a_1, a_2, a_3) \in K^3 \setminus \{ (0, 0, 0) \} \right\},
$$
so that $\dim Y = 2$; we let
$$
\hat Y = \left\{ \langle v_{-(\alpha_1 + \alpha_2)}, a_1 v_{-\alpha_1} + a_2 v_0 + a_3 v_{\alpha_1}, v_{\alpha_1 + \alpha_2} \rangle : a_1a_2a_3 \neq 0, \ {\ts\frac{a_1a_3}{{a_2}^2} \neq \frac{1}{4}} \right\},
$$
and then $\hat Y$ is a dense open subset of $Y$.

Take $y = \langle v^{(1)}, v^{(2)}, v^{(3)} \rangle \in \hat Y$, where
$$
v^{(1)} = v_{-(\alpha_1 + \alpha_2)}, \quad v^{(2)} = a_1 v_{-\alpha_1} + a_2 v_0 + a_3 v_{\alpha_1}, \quad v^{(3)} = v_{\alpha_1 + \alpha_2};
$$
note that the condition $\frac{a_1a_3}{{a_2}^2} \neq \frac{1}{4}$ implies that $X_{\alpha_1}.v^{(2)}$ does not contain any vector in $\langle v_{-\alpha_1} \rangle$. Take $g \in \Tran_G(y, Y)$, and write $y' = g.y$ and $g = u_1nu_2$ with $u_1 \in U$, $n \in N$ and $u_2 \in U_w$ where $w = nT \in W$; write $u_2 = \prod x_\alpha(t_\alpha)$ where the product takes the relevant roots $\alpha$ in order of increasing height. We have ${u_1}^{-1}.y' = n.(u_2.y)$; the weights $-(\alpha_1 + \alpha_2)$, $-\alpha_1$ and $\alpha_1 + \alpha_2$ occur in $u_2.v^{(1)}$, $u_2.v^{(2)}$ and $u_2.v^{(3)}$ respectively, so $w$ cannot send any of these three weights to $-(2\alpha_1 + \alpha_2)$ as this does not occur in ${u_1}^{-1}.y'$, whence $w \in \langle w_{\alpha_2} \rangle \{ 1, w_{3\alpha_1 + \alpha_2}, w_{2\alpha_1 + \alpha_2} \}$. Thus $n = hn^*$ where $h \in T$ and
$$
n^* \in \left\{ 1, n_{\alpha_2}, n_{\alpha_1} n_{\alpha_2} n_{\alpha_1}, n_{\alpha_2} n_{\alpha_1} n_{\alpha_2} n_{\alpha_1}, n_{\alpha_1} n_{\alpha_2} n_{\alpha_1} n_{\alpha_2} n_{\alpha_1}, n_{\alpha_2} n_{\alpha_1} n_{\alpha_2} n_{\alpha_1} n_{\alpha_2} n_{\alpha_1} \right\}.
$$
Note that, in addition to each of the vectors $n^*u_2.v^{(1)}$, $n^*u_2.v^{(2)}$ and $n^*u_2.v^{(3)}$ having no $v_{-(2\alpha_1 + \alpha_2)}$ term, some non-zero linear combination of the three vectors must equal $h^{-1}{u_1}^{-1}.v_{\alpha_1 + \alpha_2}$ and therefore lie in $V_{[+]}$, so that the projections on $V_{[0]}$ of the three vectors must be linearly dependent.

If $n^* = n_{\alpha_2}$ then the projections on $V_{[0]}$ of the vectors $n^*u_2.v^{(i)}$ are $-v_{-\alpha_1}$, $a_2 v_0 - a_3 t_{\alpha_2} v_{\alpha_1}$ and $-v_{\alpha_1}$, which are linearly independent. If $n^* = n_{\alpha_1} n_{\alpha_2} n_{\alpha_1}$ then $n^*u_2.v^{(1)}$ has zero projection on both $V_{[0]}$ and $V_{[+]}$; the coefficient of $v_{-(2\alpha_1 + \alpha_2)}$ in $n^*u_2.v^{(2)}$ is $a_1 {t_{\alpha_1}}^2 + a_2 t_{\alpha_1} + a_3$, so this expression must be zero, whence $t_{\alpha_1}, 2a_1 t_{\alpha_1} + a_2 \neq 0$; now the projections on $V_{[0]}$ of $n^*u_2.v^{(2)}$ and $n^*u_2.v^{(3)}$ are $(a_1 t_{\alpha_1} t_{2\alpha_1 + \alpha_2} + a_1 t_{3\alpha_1 + \alpha_2} - a_2 t_{2\alpha_1 + \alpha_2}) v_{-\alpha_1} + (2a_1 t_{\alpha_1} + a_2) v_0$ and $-t_{\alpha_1} v_{-\alpha_1}$, which are linearly independent. If $n^* = n_{\alpha_2} n_{\alpha_1} n_{\alpha_2} n_{\alpha_1}$ then the coefficient of $v_{-(2\alpha_1 + \alpha_2)}$ in $n^*u_2.v^{(2)}$ is again $a_1 {t_{\alpha_1}}^2 + a_2 t_{\alpha_1} + a_3$, so once more $2a_1 t_{\alpha_1} + a_2 \neq 0$; now the projections on $V_{[0]}$ of the vectors $n^*u_2.v^{(i)}$ are $v_{-\alpha_1}$, $(2a_1t_{\alpha_1} + a_2) v_0 + a_1 t_{2\alpha_1 + \alpha_2} v_{\alpha_1}$ and $v_{\alpha_1}$, which are linearly independent. If $n^* = n_{\alpha_1} n_{\alpha_2} n_{\alpha_1} n_{\alpha_2} n_{\alpha_1}$ then the coefficient of $v_{-(2\alpha_1 + \alpha_2)}$ in $n^*u_2.v^{(3)}$ is $t_{\alpha_1}$, which thus must be zero; now the projections on $V_{[0]}$ of the vectors $n^*u_2.v^{(i)}$ are ${t_{\alpha_1 + \alpha_2}}^2 v_{-\alpha_1} - 2t_{\alpha_1 + \alpha_2} v_0 + v_{\alpha_1}$, $(a_1 t_{2\alpha_1 + \alpha_2} + a_2 t_{\alpha_1 + \alpha_2}) v_{-\alpha_1} - a_2 v_0$ and $v_{-\alpha_1}$, which are linearly independent. Thus we must have $n^* \in \{ 1, n_{\alpha_2} n_{\alpha_1} n_{\alpha_2} n_{\alpha_1} n_{\alpha_2} n_{\alpha_1} \}$.

First suppose $n^* = 1$; then $u_2 = 1$ and so $g = u_1 h = h {u_1}'$ where ${u_1}' = {u_1}^h$. From the coefficient of $v_{2\alpha_1 + \alpha_2}$ in ${u_1}'.v^{(3)}$ we see that the projection of ${u_1}'$ on the root group $X_{\alpha_1}$ must be trivial; hence ${u_1}'.v^{(2)} - v^{(2)} \in V_{[+]}$, and it follows that we must have ${u_1}' \in C_U(y)$. Equating to zero the coefficient of $v_{2\alpha_1 + \alpha_2}$ in ${u_1}'.v^{(1)}$ and ${u_1}'.v^{(2)}$, and requiring the projection of ${u_1}'.v^{(1)}$ on $V_{[0]}$ to be a scalar multiple of $v^{(2)}$, shows that ${u_1}' = x_{(a_1, a_2, a_3)}(t)$ for some $t \in K$, where we write
$$
x_{(a_1, a_2, a_3)}(t) = x_{\alpha_2}(a_3 t) x_{\alpha_1 + \alpha_2}(-{\ts\frac{a_2}{2}} t) x_{2\alpha_1 + \alpha_2}(a_1 t) x_{3\alpha_1 + \alpha_2}({\ts\frac{a_1a_2}{2a_3}} t) x_{3\alpha_1 + 2\alpha_2}({\ts\frac{a_1a_2}{2}} t^2).
$$

Now suppose $n^* = n_{\alpha_2} n_{\alpha_1} n_{\alpha_2} n_{\alpha_1} n_{\alpha_2} n_{\alpha_1}$. Since $n^*.v_{2\alpha_1 + \alpha_2} = -v_{-(2\alpha_1 + \alpha_2)}$, the coefficient of $v_{2\alpha_1 + \alpha_2}$ in each vector $u_2.v^{(i)}$ must be zero; again taking $i = 3$ shows that $t_{\alpha_1} = 0$, and now it follows that we must have $u_2.v^{(1)} = v^{(1)} + v + av^{(3)}$, $u_2.v^{(2)} = v^{(2)} + bv^{(3)}$ and $u_2.v^{(3)} = v^{(3)}$ for some $v \in V_{[0]}$ and some $a, b \in K$. As $n^*$ interchanges $v^{(1)}$ and $v^{(3)}$ and preserves $V_{[0]}$, we have $n^*u_2.v^{(1)} = av^{(1)} + n^*.v + v^{(3)}$, $n^*u_2.v^{(2)} = bv^{(1)} + n^*.v^{(2)}$ and $n^*u_2.v^{(3)} = v^{(1)}$; since some linear combination of these must be in $V_{[+]}$, we must have $n^*.v \in \langle n^*.v^{(2)} \rangle$, and so $v \in \langle v^{(2)} \rangle$. Therefore $u_2 \in C_U(y)$; as $n^*.y \in \hat Y$, the previous paragraph shows that we must have ${u_1}' \in C_U(n^*.y)$, so $y' = n^*.y$ and ${u_1}' \in C_U(y')$, $u_2 \in C_U(y)$.

Write $n_0 = h_{\alpha_2}(-\frac{a_1}{a_3}) n_{\alpha_2} n_{\alpha_1} n_{\alpha_2} n_{\alpha_1} n_{\alpha_2} n_{\alpha_1}$; then $n_0.v^{(2)} = -v^{(2)}$, and so $n_0 \in C_N(y)$. Therefore we have
$$
{\ts\Tran_G(y, Y) = T \{ x_{(a_1, a_2, a_3)}(t) : t \in K \}(\{ 1 \} \cup \{ n_0 x_{(a_1, a_2, a_3)}(t') : t' \in K \})}.
$$
Thus $G.y \cap Y = T.y$; we see that $C_T(y) = \{ h_{3\alpha_1 + 2\alpha_2}(\kappa) : \kappa \in K^* \}$, and
$$
T.y = \left\{ \langle v_{-(\alpha_1 + \alpha_2)}, b_1 v_{-\alpha_1} + b_2 v_0 + b_3 v_{\alpha_1}, v_{\alpha_1 + \alpha_2} \rangle \in \hat Y : {\ts\frac{b_1b_3}{{b_2}^2}} = {\ts\frac{a_1a_3}{{a_2}^2}} \right\}.
$$
Moreover we have $C_G(y) = C_U(y) C_N(y) C_U(y) = \langle x_{(a_1, a_2, a_3)}(t), n_0 : t \in K \rangle$. Since $\dim C_G(y) = 3$, we have $\dim(\overline{G.y}) = \dim G - \dim C_G(y) = 14 - 3 = 11$, while $\dim(\overline{G.y \cap Y}) = 1$; therefore
$$
\dim \G{3}(V) - \dim(\overline{G.y}) = 12 - 11 = 1 \quad \hbox{and} \quad \dim Y - \dim(\overline{G.y \cap Y}) = 2 - 1 = 1.
$$
Therefore $y$ is $Y$-exact.

To prove the conjugacy of stabilizers requires a little more work. Write
$$
C = \langle x_{\alpha_2}(t) x_{2\alpha_1 + \alpha_2}(t), x_{-\alpha_2}(t) x_{-(2\alpha_1 + \alpha_2)}(t) : t \in K \rangle;
$$
then $Z(C) = \langle h_{\alpha_2}(-1) h_{2\alpha_1 + \alpha_2}(-1) \rangle = \{ 1 \}$. Set $c = \frac{4a_1a_3}{{a_2}^2}$, so that $c \neq 0, 1$, and let $\xi$ be a root of the cubic $x(x + 3)^2 = c(3x + 1)^2$, so that $\xi \neq 0, 1, -3, -\frac{1}{3}$. Let ${a_1}^* = \frac{(\xi - 1)(3\xi + 1)}{\xi(\xi + 3)}.\frac{a_1}{a_2}$ and ${a_2}^* = \frac{\xi}{4{a_1}^*}$; then with $x = x_{-\alpha_1}(\frac{2{a_1}^*}{4{a_1}^*{a_2}^* - 1}) x_{\alpha_1}(\frac{1}{2{a_1}^*})$ we have
$$
x^{-1}.y = \langle v_{-(2\alpha_1 + \alpha_2)} + 2{a_2}^* v_{-(\alpha_1 + \alpha_2)}, {a_1}^* v_{-\alpha_1} + v_0 + {a_2}^* v_{\alpha_1}, 2{a_1}^* v_{\alpha_1 + \alpha_2} + v_{2\alpha_1 + \alpha_2} \rangle,
$$
and
\begin{eqnarray*}
x_{(a_1, a_2, a_3)}(t)^x & = & x_{\alpha_2}({a_1}^* t^*) x_{2\alpha_1 + \alpha_2}({a_2}^* t^*) \qquad \hbox{where } t^* = {\ts\frac{(\xi - 1)^2}{\xi(3\xi + 1)}} t, \\
(h_{3\alpha_1 + 2\alpha_2}({\ts\frac{\xi - 1}{4{{a_1}^*}^2}}) n_0)^x & = & h_{3\alpha_1 + \alpha_2}({\ts\frac{{a_2}^*}{{a_1}^*}}) h_{\alpha_2}(-1) n_{\alpha_2} n_{\alpha_1} n_{\alpha_2} n_{\alpha_1} n_{\alpha_2} n_{\alpha_1},
\end{eqnarray*}
so that
$$
C_G(x^{-1}.y) = \langle x_{\alpha_2}({a_1}^* t^*) x_{2\alpha_1 + \alpha_2}({a_2}^* t^*), x_{-\alpha_2}({a_2}^* t^*) x_{-(2\alpha_1 + \alpha_2)}({a_1}^* t^*) : t^* \in K \rangle.
$$
Now if we take $\kappa \in K^*$ satisfying $\kappa^2 = \frac{{a_2}^*}{{a_1}^*}$ and set $h = h_{3\alpha_1 + \alpha_2}(\kappa)$, then we have $C_G(x^{-1}.y) = {}^h C$, so that $C_G(y) = {}^{xh} C$. Thus the conditions of Lemma~\ref{lem: generic stabilizer from exact subset} hold; so the quadruple $(G, \lambda, p, k)$ has generic stabilizer $C/Z(G) \cong A_1$, where the $A_1$ is of adjoint type.

Now take the case where $p = 2$; here $\Lambda(V) = \Phi_s$. Write $y_0 = \langle v^{(1)}, v^{(2)}, v^{(3)} \rangle$ where
$$
v^{(1)} = v_{-(\alpha_1 + \alpha_2)}, \quad v^{(2)} = v_{-\alpha_1} + v_{\alpha_1 + \alpha_2}, \quad v^{(3)} = v_{2\alpha_1 + \alpha_2}.
$$
Take $g \in C_G(y_0)$, and write $g = u_1nu_2$ with $u_1 \in U$, $n \in N$ and $u_2 \in U_w$ where $w = nT \in W$; write $u_2 = \prod x_\alpha(t_\alpha)$ where the product takes the relevant roots $\alpha$ in order of increasing height. We have ${u_1}^{-1}.y_0 = n.(u_2.y_0)$; the weights $-(\alpha_1 + \alpha_2)$, $-\alpha_1$ and $2\alpha_1 + \alpha_2$ occur in $u_2.v^{(1)}$, $u_2.v^{(2)}$ and $u_2.v^{(3)}$ respectively, so $w$ cannot send any of these three weights to $-(2\alpha_1 + \alpha_2)$ as this does not occur in ${u_1}^{-1}.y_0$, whence $w \in \langle w_{\alpha_2}, w_{3\alpha_1 + \alpha_2} \rangle$. Thus $n = hn^*$ where $h \in T$ and
$$
n^* \in \left\{ 1, n_{\alpha_2}, n_{\alpha_1} n_{\alpha_2} n_{\alpha_1}, n_{\alpha_2} n_{\alpha_1} n_{\alpha_2} n_{\alpha_1}, n_{\alpha_1} n_{\alpha_2} n_{\alpha_1} n_{\alpha_2}, n_{\alpha_2} n_{\alpha_1} n_{\alpha_2} n_{\alpha_1} n_{\alpha_2} \right\}.
$$
Note that, in addition to each of the vectors $n^*u_2.v^{(1)}$, $n^*u_2.v^{(2)}$ and $n^*u_2.v^{(3)}$ having no $v_{-(2\alpha_1 + \alpha_2)}$ term, some non-zero linear combination of the three vectors must equal $h^{-1}{u_1}^{-1}.v_{2\alpha_1 + \alpha_2}$ and therefore lie in $\langle v_{2\alpha_1 + \alpha_2} \rangle$.

If $n^* = n_{\alpha_1} n_{\alpha_2} n_{\alpha_1}$ or $n_{\alpha_2} n_{\alpha_1} n_{\alpha_2} n_{\alpha_1}$ then the coefficient of $v_{-(2\alpha_1 + \alpha_2)}$ in $n^*u_2.v^{(1)}$ is $t_{2\alpha_1 + \alpha_2}$, which thus must be zero; now $n^*u_2.v^{(i)}$ has non-zero projection on $\langle v_{\alpha_1}, v_{\alpha_1 + \alpha_2}, v_{2\alpha_1 + \alpha_2} \rangle$ only for $i = 2$, and the projection is $v_{\alpha_1 + \alpha_2} + v_{2\alpha_1 + \alpha_2}$ or $v_{\alpha_1} + v_{2\alpha_1 + \alpha_2}$ respectively, so no non-zero linear combination of the three vectors lies in $\langle v_{2\alpha_1 + \alpha_2} \rangle$. If $n^* = n_{\alpha_1} n_{\alpha_2} n_{\alpha_1} n_{\alpha_2}$ or $n_{\alpha_2} n_{\alpha_1} n_{\alpha_2} n_{\alpha_1} n_{\alpha_2}$ then the coefficient of $v_{-(2\alpha_1 + \alpha_2)}$ in $n^*u_2.v^{(2)}$ is $1 + t_{2\alpha_1 + \alpha_2}$, which thus must be zero; now $n^*u_2.v^{(i)}$ has non-zero projection on $\langle v_{\alpha_1}, v_{\alpha_1 + \alpha_2}, v_{2\alpha_1 + \alpha_2} \rangle$ only for $i = 1$, and the projection is $v_{\alpha_1 + \alpha_2} + v_{2\alpha_1 + \alpha_2}$ or $v_{\alpha_1} + v_{2\alpha_1 + \alpha_2}$ respectively, so no non-zero linear combination of the three vectors lies in $\langle v_{2\alpha_1 + \alpha_2} \rangle$. Thus we must have $n^* \in \{ 1, n_{\alpha_2} \}$.

First suppose $n^* = 1$; then $u_2 = 1$ and $g = u_1 h$, and we must have $h \in C_T(y_0)$ and $u_1 \in C_U(y_0)$. The requirement that $h.v^{(2)} \in \langle v^{(2)} \rangle$ shows that $h \in \{ h_{\alpha_2}(\kappa) : \kappa \in K^* \}$, so this is $C_T(y_0)$. From the coefficient of $v_{\alpha_1}$ in $u_1.v^{(i)}$ for $i = 1, 2$ we see that the projection of $u_1$ on the root groups $X_{2\alpha_1 + \alpha_2}$ and $X_{\alpha_1}$ must be trivial; for the coefficients of $v_{-\alpha_1}$ and $v_{\alpha_1 + \alpha_2}$ in $u_1.v^{(1)}$ to be equal we require $u_1 \in \{ x_{\alpha_2}(t^2) x_{\alpha_1 + \alpha_2}(t) : t \in K \} X_{3\alpha_1 + \alpha_2} X_{3\alpha_1 + 2\alpha_2}$, so this is $C_U(y_0)$.

Now suppose $n^* = n_{\alpha_2}$; by the previous paragraph we may write $g = uhg'$, where $u \in C_U(y_0)$, $h \in T$ and $g' = x_{\alpha_1}(t_1) x_{\alpha_1 + \alpha_2}(t_2) x_{2\alpha_1 + \alpha_2}(t_3) n_{\alpha_2} x_{\alpha_2}(t)$ for some $t_1, t_2, t_3, t \in K$, and we must have $hg' \in C_G(y_0)$. The coefficients of $v_{\alpha_1}$ in $g'.v^{(1)}$ and $g'.v^{(2)}$ are $t_3t + {t_1}^2$ and $1 + t_3$, so we must have $t_3 = 1$ and $t = {t_1}^2$; the coefficients of $v_{-\alpha_1}$ and $v_{\alpha_1 + \alpha_2}$ in $g'.v^{(2)}$ are $0$ and ${t_2}^2$, so we must have $t_2 = 0$; thus $g' = x_{\alpha_1}(t_1) x_{2\alpha_1 + \alpha_2}(1) n_{\alpha_2} x_{\alpha_2}({t_1}^2)$ and we see that now $g' \in C_G(y_0)$, whence also $h \in C_T(y_0)$. Observe that $g' = x_{3\alpha_1 + \alpha_2}(t_1) x_{2\alpha_1 + \alpha_2}(1) n_{\alpha_2} x_{\alpha_2}({t_1}^2) x_{\alpha_1 + \alpha_2}(t_1)$.

Thus if we write $n^\dagger = x_{2\alpha_1 + \alpha_2}(1) n_{\alpha_2}$, and $x(t) = x_{\alpha_2}(t^2) x_{\alpha_1 + \alpha_2}(t)$ for $t \in K$, then we have
$$
C_G(y_0) = X_{3\alpha_1 + \alpha_2} X_{3\alpha_1 + 2\alpha_2} \left\{ x(t)h_{\alpha_2}(\kappa), x(t) h_{\alpha_2}(\kappa) n^\dagger x(t') : t, t' \in K, \ \kappa \in K^* \right\}.
$$
Calculation now shows that if $t \in K^*$ we have $x(t).x(\frac{1}{t})^{n^\dagger}.x(t) = h_{\alpha_2}(t^2) n^\dagger$; so $\{ x(t)h_{\alpha_2}(\kappa), x(t) h_{\alpha_2}(\kappa) n^\dagger x(t') : t, t' \in K, \ \kappa \in K^* \} = \langle x(t), n^\dagger : t \in K \rangle$ is an $A_1$ subgroup.

Since $\dim(\overline{G.y_0}) = \dim G - \dim C_G(y_0) = 14 - 5 = 9 = \dim \G{3}(V)$, the orbit $G.y_0$ is dense in $\G{3}(V)$. Thus the quadruple $(G, \lambda, p, k)$ has generic stabilizer $C_G(y_0)/Z(G) \cong A_1 U_2$, where the $A_1$ is of adjoint type.
\end{proof}

\begin{prop}\label{prop: G_2, omega_2 module, p = 3, k = 2 or 3}
Let $G = G_2$ and $\lambda = \omega_2$ with $p = 3$, and take $k = 2$ or $3$. Then the quadruple $(G, \lambda, p, k)$ has generic stabilizer $\tilde A_1 T_1.\Z_2$ or $A_1$ respectively.
\end{prop}

\begin{proof}
This is an immediate consequence of Propositions~\ref{prop: G_2, omega_1 module, k = 2} and \ref{prop: G_2, omega_1 module, k = 3} respectively, using the graph automorphism of $G_2$ which exists in characteristic $3$.
\end{proof}

This completes the justification of the entries in Tables~\ref{table: large higher quadruple non-TGS}, \ref{table: small classical higher quadruple generic stab} and \ref{table: small exceptional higher quadruple generic stab}, and hence the proof of Theorems~\ref{thm: large higher quadruple generic stab} and \ref{thm: small higher quadruple generic stab}.

\appendix

\backmatter

\bibliographystyle{amsplain}


\printindex

\end{document}